\documentclass[12pt,a4paper,twoside]{article}
\usepackage{amssymb, amsmath, amsfonts, stmaryrd}
\usepackage[mathscr]{euscript}
\usepackage{mathrsfs}
\usepackage{theorem}
\usepackage[english]{babel}
\usepackage{bm}
\usepackage[all]{xy}
\usepackage[a4paper,bookmarks=false]{hyperref}

\newif\ifthesis\thesisfalse

\long\def\nodo#1{\relax}
\long\def\comment#1\endcomment{\relax}
\def\makepoint#1{\medbreak\noindent{\bf #1.\ }}
\def\zeropoint{\setcounter{subsection}{-1}}
\def\zerosubpoint{\setcounter{subsubsection}{-1}}
\def\nxpoint{\refstepcounter{subsection}%
  \makepoint{\thesubsection}}
\def\nxsubpoint{\refstepcounter{subsubsection}%
  \makepoint{\thesubsubsection}}

\let\ptref=\refpoint

\def\markbothsame#1{\markboth{#1}{#1}}
\def\pttocinfo#1{\addcontentsline{toc}{subsection}{\thesubsection. #1}%
  \markright{\textbf{\thesubsection.} #1}}
\def\mysection#1{\section{#1}\markbothsame{Chapter \textbf{\thesection.} #1}}

\def\nxpointtocb#1{\subsection{#1}\pttocinfo{#1}}
\let\nxpointtoc=\nxpointtocm
\let\nxpointtocb=\nxpointtocm

\def\smallmatrix#1#2#3#4{\genfrac{(}{.}{0pt}{}{#1}{#3}\,\genfrac{.}{)}{0pt}{}{#2}{#4}}
\def\legendre#1#2{\genfrac{(}{)}{0.4pt}{}{#1}{#2}}
\def\setmyindent{\hangindent=3.9em\parindent=\hangindent}
\def\myindent#1{\setmyindent\noindent
\hbox to\hangindent{\hfill #1\hskip0.4emplus0.4emminus0.2em}\ignorespaces}

\let\myunderline=\mathbf
\let\stdfsets=\mathbf
\let\cst=\underline
\def\st#1{{\stdfsets{#1}}}
\def\stn{\st{n}}
\def\stm{\st{m}}
\def\stp{\st{p}}
\def\stk{\st{k}}
\def\str{\st{r}}
\def\typ{\negthinspace:\negthinspace}
\def\card{\operatorname{card}}
\def\charact{\operatorname{char}}
\def\id{\operatorname{id}}
\def\ev{\operatorname{ev}}
\def\Id{\operatorname{Id}}
\def\gr{\operatorname{gr}}
\def\pt{\operatorname{pt}}
\def\sk{\operatorname{sk}}
\def\cosk{\operatorname{cosk}}
\def\pr{\operatorname{pr}}
\def\ch{\operatorname{ch}}
\def\Ob{\operatorname{Ob}}
\def\Ar{\operatorname{Ar}}
\def\Hom{\operatorname{Hom}}
\def\End{\operatorname{End}}
\def\END{\operatorname{END}}
\def\Aut{\operatorname{Aut}}

\def\Ind{\operatorname{Ind}}
\def\Ker{\operatorname{Ker}}
\def\Coker{\operatorname{Coker}}
\def\Im{\operatorname{Im}}
\def\Coim{\operatorname{Coim}}
\def\Spec{\operatorname{Spec}}
\def\Stab{\operatorname{Stab}}
\def\Proj{\operatorname{Proj}}
\def\Aff{\operatorname{Aff}}
\def\Ch{\operatorname{Ch}}
\def\Pic{\operatorname{Pic}}
\def\Oct{\operatorname{Oct}}
\def\Tot{\operatorname{Tot}}
\def\Ho{\operatorname{Ho}}
\def\HO{\operatorname{HO}}
\def\Ext{\operatorname{Ext}}

\def\Perf{\operatorname{Perf}}
\def\Coh{\operatorname{Coh}}
\def\Cl{\operatorname{Cl}}
\def\SCl{\operatorname{SCl}}
\def\GCl{\operatorname{GCl}}
\def\GSCl{\operatorname{GSCl}}
\def\FinCl{\operatorname{FinCl}}
\def\FinSCl{\operatorname{FinSCl}}
\def\GFinCl{\operatorname{GFinCl}}
\def\GFinSCl{\operatorname{GFinSCl}}
\def\Lim{\operatorname{Lim}}
\def\genLim#1{\operatorname{\vtop{\ialign{##\cr$\Lim$\cr\noalign{\nointerlineskip\kern0.5ex}\hfil$#1$\hfil\cr\noalign{\nointerlineskip\kern-0.5ex}\cr}}}}
\def\projLim{\genLim\longleftarrow}
\def\injLim{\genLim\longrightarrow}
\def\quotprojlim{\operatorname{\mbox{``$\projlim$''}}}

\def\diag{\operatorname{diag}}
\def\div{\operatorname{div}}
\def\sgn{\operatorname{sgn}}
\def\cl{\operatorname{cl}}
\def\supp{\operatorname{supp}}
\def\tr{\operatorname{tr}}
\def\rank{\operatorname{rank}}
\def\conv{\operatorname{conv}}
\def\Bilin{\operatorname{Bilin}}
\def\Polylin{\operatorname{Polylin}}
\def\Inv{\operatorname{Inv}}
\def\iHom{\operatorname{\myunderline{Hom}}}
\def\iEnd{\operatorname{\myunderline{End}}}
\def\iEND{\operatorname{\myunderline{END}}}
\def\iBilin{\operatorname{\myunderline{Bilin}}}
\def\Sing{\operatorname{Sing}}
\def\cA{\mathcal{A}}
\def\cB{\mathcal{B}}
\def\cC{\mathcal{C}}
\def\cD{\mathcal{D}}
\def\cE{\mathcal{E}}
\def\cS{\mathcal{S}}
\def\cU{\mathcal{U}}
\def\cF{\mathcal{F}}
\def\cG{\mathcal{G}}
\def\cL{\mathcal{L}}
\def\cI{\mathcal{I}}
\def\cO{\mathcal{O}}
\def\cT{\mathcal{T}}
\def\cP{\mathcal{P}}
\def\cX{\mathcal{X}}
\def\sA{\mathscr{A}}
\def\sB{\mathscr{B}}
\def\sE{\mathscr{E}}
\def\sF{\mathscr{F}}
\def\sG{\mathscr{G}}
\def\sL{\mathscr{L}}
\def\sR{\mathscr{R}}
\def\sX{\mathscr{X}}
\def\sO{\mathscr{O}}
\def\sP{\mathscr{P}}
\def\ga{\mathfrak{a}}
\def\gb{\mathfrak{b}}
\def\gm{\mathfrak{m}}
\def\gp{\mathfrak{p}}
\def\gq{\mathfrak{q}}
\def\gS{\mathfrak{S}}
\def\gP{\mathfrak{P}}
\def\rad{\mathfrak{r}}
\def\fnorm{|\negthinspace\cdot\negthinspace|}
\def\vnorm{\|\negthinspace\cdot\negthinspace\|}
\def\nin{\not\in}
\def\Cop{C^\circ}
\def\bbQ{\mathbb{Q}}
\def\bbZ{\mathbb{Z}}
\def\bbN{\mathbb{N}}
\def\bbF{\mathbb{F}}
\def\bbR{\mathbb{R}}
\def\bbC{\mathbb{C}}
\def\bbH{\mathbb{H}}
\def\bbA{\mathbb{A}}
\def\bbP{\mathbb{P}}
\def\bbG{\mathbb{G}}
\def\bbT{\mathbb{T}}
\def\Synt{\mathbb{S}}
\def\dL{\mathbb{L}}
\def\dR{\mathbb{R}}
\def\sitV{\mathfrak{V}}
\def\sitW{\mathfrak{W}}
\def\Hr{\stackrel{r}}
\def\Hl{\stackrel\ell}

\def\simr{\stackrel{r}\sim}
\def\siml{\stackrel\ell\sim}
\def\sims{\stackrel{s}\sim}
\def\simss{\stackrel{ss}\sim}
\def\Loplus{\underline{\underline\oplus}}
\def\Lotimes{\underline{\underline\otimes}}
\def\Loslash{\underline{\underline\oslash}}
\def\af{{}^a\negthinspace f}
\def\bbDelta{{\bm\Delta}}
\def\bpi{{\bm\pi}}
\def\bphi{{\bm\phi}}
\def\Zinfty{\bbZ_\infty}
\def\Zninfty{\bbZ_{(\infty)}}
\def\barZinfty{\bar\bbZ_\infty}

\def\Sigmainf{\Sigma_\infty}
\def\Fone{\bbF_1}
\def\Fpm{\bbF_{\pm1}}
\def\Fempty{\bbF_\emptyset}
\def\Finfty{\bbF_\infty}
\def\SpecZ{{\Spec\bbZ}}
\def\CompZ{\widehat\SpecZ}
\def\CompuZ{\widehat{\Spec^u\bbZ}\strut}
\def\quotedinjlim{\operatorname{\mbox{``$\varinjlim$''}}}
\def\quotedprojlim{\operatorname{\mbox{``$\varprojlim$''}}}
\def\simto{\stackrel\sim\to}
\def\Unit{{\bm 1}}
\def\bu{{\bm e}}
\def\wedgenl{\bigwedge\nolimits}
\let\boxe=\boxempty
\let\phi=\varphi
\let\epsilon=\varepsilon
\let\zeroword=\emptyset
\let\emptyset=\varnothing
\let\injlim=\varinjlim
\let\projlim=\varprojlim
\let\longto=\longrightarrow
\let\textcat=\textsl
\def\textstk#1{{\rm\textbf{#1}}}
\def\catSets{{\textcat{Sets\/}}}
\def\catGrps{{\textcat{Grps\/}}}
\def\catTop{{\textcat{Top\/}}}
\def\catCat{{\textcat{Cat\/}}}
\def\catMon{{\textcat{Mon}}}
\def\catAb{{\textcat{Ab}}}
\def\catPt{{\langle{\bm *}\rangle}}
\def\catN{{\underline{\mathbb N}}}
\def\catAlg{\operatorname{\textcat{Alg\/}}}
\def\catCommAlg{\operatorname{\textcat{CommAlg\/}}}
\def\catMonads{\operatorname{\textcat{Monads\/}}}
\def\catGenR{\operatorname{\textcat{GenR\/}}}
\def\catFunct{\operatorname{\textcat{Funct\/}}}
\def\catCart{\operatorname{\textcat{Cart\/}}}
\def\catDesc{\operatorname{\textcat{Desc\/}}}
\def\catEndof{\operatorname{\textcat{Endof\/}}}
\def\catEndofalg{\catEndof_{alg}}
\def\catDiagFunct{\operatorname{\textcat{DiagFunct\/}}}
\def\catDiagEndof{\operatorname{\textcat{DiagEndof\/}}}
\def\catInnFunct{\operatorname{\textcat{InnFunct\/}}}
\def\catInnEndof{\operatorname{\textcat{InnEndof\/}}}
\let\catDelta=\bbDelta
\def\monMon{\mathbb{M}}
\def\ZpLat{{\textcat{$\bbZ_p$-Lat}}}
\def\catMod#1{{\textcat{$#1$-Mod}}}
\def\catGradMod#1{{\textcat{$#1$-GrMod}}}
\def\catFlMod#1{{\textcat{$#1$-Fl.Mod}}}
\def\catVect#1{{\textcat{$#1$-Vect}}}
\def\catAlgebr#1{{\textcat{$#1$-Alg}}}
\def\catCommAlgebr#1{{\textcat{$#1$-CommAlg}}}
\def\catQCoh#1{{\textcat{$#1$-QCoh}}}
\def\ZinfLat{{\textcat{$\Zinfty$-Lat\/}}}
\def\ZinfFlat{\catFlMod{\Zinfty}}
\def\ZpFlat{\catFlMod{\bbZ_p}}
\def\ZinfMod{\catMod{\Zinfty}}
\def\RVect{\catVect{\bbR}}
\def\stSETS{{\textstk{SETS}}}
\def\stSHEAVES{{\textstk{SHEAVES}}}
\def\stLCSETS{{\textstk{LCSETS}}}
\def\stMODo{{\textstk{MOD}}}
\def\stMOD#1{{\textstk{$#1$-MOD}}}
\def\stCART{{\textstk{CART}}}
\def\stQCOH{{\textstk{QCOH}}}
\def\stks{\mathfrak{s}}
\def\stkc{\mathfrak{c}}
\def\stkAr{\operatorname{\mathfrak{Ar}}}
\def\univU{\mathcal{U}}
\def\univV{\mathcal{V}}
\def\propQ{\mathscr{Q}}
\def\propP{\mathscr{P}}
\def\propW{\mathscr{W}}
\theorembodyfont{\sl}

\newtheorem{ThD}[subsubsection]{Theorem}
\newtheorem{Thz}{Theorem.}

\newtheorem{PropD}[subsubsection]{Proposition}
\newtheorem{Propz}{Proposition.}

\newtheorem{LemmaD}[subsubsection]{Lemma}
\newtheorem{CorD}[subsubsection]{Corollary}

\newtheorem{DefD}[subsubsection]{Definition}

\theorembodyfont{\rm}

\newtheorem{RemD}[subsubsection]{Remark}
\newtheorem{NotatD}[subsubsection]{Notation}

\newtheorem{Proof}{Proof.}

\numberwithin{equation}{subsubsection}

\title{New Approach to Arakelov Geometry}
\author{Nikolai Durov}
\begin{document}

\ifthesis
\else
\maketitle
\fi

\pagestyle{myheadings}


\section*{Introduction}
\setcounter{section}{0}
\setcounter{subsection}{0}
\markbothsame{Introduction}

The principal aim of this work is to provide an alternative algebraic
framework for Arakelov geometry, and to demonstrate its usefulness by
presenting several simple applications. This framework,
called {\em theory of generalized rings and schemes}, appears to be
useful beyond the scope of Arakelov geometry, providing a uniform
description of classical scheme-theoretical algebraic geometry
(``schemes over $\Spec\bbZ$''), Arakelov geometry (``schemes over
$\Spec\Zinfty$ and $\CompZ$''), tropical geometry (``schemes over
$\Spec\bbT$ and $\Spec\bbN$'') and the geometry over the so-called
field with one element (``schemes over $\Spec\Fone$'').
Therefore, we develop this theory a bit further than it is strictly
necessary for Arakelov geometry.

The approach to Arakelov geometry developed in this work is completely 
{\em algebraic}, in the sense that it doesn't require the combination
of scheme-theoretical algebraic geometry and complex differential geometry,
traditionally used in Arakelov geometry since the works of Arakelov himself.

However, we show that our models $\bar\sX/\CompZ$ of algebraic
varieties $X/\bbQ$ define both a model $\sX/\Spec\bbZ$ in the usual
sense and a (possibly singular) Banach (co)metric on (the smooth locus
of) the complex analytic variety $X(\bbC)$. This metric cannot be
chosen arbitrarily; however, some classical metrics like the
Fubini--Study metric on $\bbP^n$ do arise in this way. It is
interesting to note that ``good'' models from the algebraic point of
view (e.g.\ finitely presented) usually give rise to not very nice
metrics on $X(\bbC)$, and conversely, nice smooth metrics like
Fubini--Study correspond to models with ``bad'' algebraic properties
(e.g.\ not finitely presented).

Our algebraic approach has some obvious advantages over the classical one.
For example, we never need to require $X$ to be smooth or proper, 
and we can deal with singular metrics.

In order to achieve this goal we construct a theory of {\em
generalized rings}, commutative or not, which include classical rings
(always supposed to be associative with unity), then define {\em
spectra\/} of such (commutative) generalized rings, and construct {\em
generalized schemes\/} by patching together spectra of generalized
rings. Of course, these generalized schemes are {\em generalized
ringed spaces}, i.e.\ topological spaces, endowed with a sheaf of
generalized rings. Then the ``compactified'' $\Spec\bbZ$, denoted by
$\CompZ$, is constructed as a (pro-)generalized scheme, and our
models $\bar\sX/\CompZ$ are (pro-)generalized schemes as well.

All the ``generalized'' notions we discuss are indeed generalizations of
corresponding ``classical'' notions. More precisely, ``classical'' objects
(e.g.\ commutative rings) always constitute a full subcategory of the
category of corresponding ``generalized'' objects (e.g.\ commutative 
generalized rings). In this way we can always treat for example a classical
scheme as a generalized scheme, since no new morphisms between classical
schemes arise in the larger category of generalized schemes.

In particular, the category of (commutative) generalized rings
contains all classical (commutative) rings like
$\bbZ$, $\bbQ$, $\bbR$, $\bbC$, $\bbZ_p$, \dots, as well as some new
objects, such as $\Zinfty$ (the ``archimedian valuation ring'' of~$\bbR$,
similar to $p$-adic integers $\bbZ_p\subset\bbQ_p$), 
$\Zninfty$ (the ``non-completed localization at $\infty$'', or
the ``archimedian valuation ring'' of~$\bbQ$), $\barZinfty$
(``the integral closure of $\Zinfty$ in~$\bbC$''). Furthermore,
once these ``archimedian valuation rings'' are constructed, we can
define some other generalized rings, such as $\Fpm:=\Zinfty\cap\bbZ$,
or the ``field with one element''~$\Fone$. Tropical numbers $\bbT$
and other semirings are also generalized rings, thus 
fitting nicely into this picture as well.

In this way we obtain not only a ``compact model'' $\CompZ$ of $\bbQ$
(called also ``compactification of~$\Spec\bbZ$''), and models
$\sX\to\CompZ$ of algebraic varieties $X/\bbQ$, but a
geometry over ``the field with one element'' as well. For example,
$\CompZ$ itself is a pro-generalized scheme over~$\Fone$ and~$\Fpm$.

In other words, we obtain rigorous definitions both of the
``archimedian local ring'' $\Zinfty$ and of the ``field with one element''
$\Fone$. They have been discussed in mathematical folklore for quite a
long time, but usually only in a very informal fashion.

We would like to say a few words here about some applications of
the theory of generalized rings and schemes presented in this work.
Apart from defining generalized rings, their spectra, and
generalized schemes, we discuss some basic properties of generalized
schemes, essentially transferring some results of EGA~I and~II to our
case. For example, we discuss projective (generalized) schemes and morphisms,
study line and vector bundles, define Picard groups and so on.

Afterwards we do some homological (actually homotopic) algebra over
generalized rings and schemes, define perfect simplicial objects and
cofibrations (which replace perfect complexes in this theory), 
define $K_0$ of perfect simplicial objects and vector bundles,
briefly discuss higher algebraic $K$-theory (Waldhausen's construction
seems to be well-adapted to our situation), and construct Chow rings
and Chern classes using the $\gamma$-filtration on~$K_0$,
in the way essentially known since Grothendieck's proof of
Riemann--Roch theorem.

We apply the above notions to Arakelov geometry as well. For example,
we compute Picard group, Chow ring and Chern classes of vector bundles
over $\CompZ$, and construct the moduli space of such vector bundles.
In particular, we obtain the notion of {\em (arithmetic) degree\/}
$\deg\sE\in\log\bbQ^\times_+$ of a vector bundle $\sE$ over $\CompZ$;
it induces an isomorphism $\deg:\Pic(\CompZ)\to\log\bbQ^\times_+$.
We also prove that any affine or projective algebraic variety $X$
over~$\bbQ$ admits a finitely presented model $\sX$ over $\CompZ$, and
show (under some natural conditions) that rational points $P\in
X(\bbQ)$ extend to uniquely determined sections
$\sigma_P:\CompZ\to\sX$. We show that when $X$ is a closed subvariety
of the projective space $\bbP^n_\bbQ$, and its model $\sX$ is chosen
accordingly (e.g.\ $\sX$ is the ``scheme-theoretical closure'' of $X$
in $\bbP^n_{\CompZ}$), then the (arithmetic) degree of the pullback
$\sigma_P^*(\sO_\sX(1))$ of the ample line bundle of $\sX$
equals the logarithmic height of point $P\in X(\bbQ)\subset\bbP^n(\bbQ)$.

There are several reasons to believe that our ``algebraic'' Arakelov
geometry can be related to its more classical variants, based on
K\"ahler metrics, differential forms and Green currents, as developed
first by Arakelov himself, and then in the series of works
of H.~Gillet, C.~Soul\'e and their colloborators. The simplest reason
is that our ``algebraic'' models give rise to some (co)metrics,
and classical metrics like the Fubini--Study do appear in this way.
More sophisticated arguments involve comparison with the non-archimedian
variant of classical Arakelov geometry, developed in \cite{BGS}
and~\cite{GS}. This non-archimedian Arakelov geometry is quite similar
to (classical) archimedian Arakelov geometry, and at the same time
admits a natural interpretation in terms of models of algebraic varieties
over discrete valuation rings. Analytic torsion corresponds in this
picture to torsion in the special fiber (i.e.\ lack of flatness).

Therefore, one might hope to transfer eventually
the results of these two works to archimedian context, using our
theory of generalized schemes, thus establishing a direct connection
between our ``algebraic'' and classical ``analytic'' variant of Arakelov
geometry. 

{\bf Acknowledgements.} 
First of all, I would like to thank my scientific advisor,
Gerd Faltings, for his constant attention to this work, as well as
for teaching me a lot of arithmetic geometry during my graduate studies
in Bonn. I would also like to thank Alexandr Smirnov for some very
interesting discussions and remarks, some of which have considerably
influenced this work. I'd like to thank Christophe Soul\'e for bringing
very interesting works \cite{GS} and~\cite{BGS} on non-archimedian
Arakelov geometry to my attention.

This work has been written during my graduate studies in Bonn,
as a part of the IMPRS (International Max Planck Research School for
Moduli Spaces and their applications), a joint program of the
Bonn University and the Max-Planck-Institut f\"ur Mathematik.
I'd like to thank all people involved from these organizations for
their help and support. Special thanks go to Christian Kaiser,
the coordinator of the IMPRS, whose
help was indispensable throughout all my stay in Bonn.


\cleardoublepage

\tableofcontents

\cleardoublepage


\section*{Overview}
\setcounter{section}{0}
\setcounter{subsection}{0}
\markbothsame{Overview}
\addcontentsline{toc}{section}{Overview}

We would like to start with a brief overview of the rest of this work,
discussing chapter after chapter. Purely technical definitions
and statements will be omitted or just briefly mentioned, while those
notions and results, which we consider crucial for the
understanding of the remainder of this work, 
will be explained at some length.

\nxpointtoc{Motivation}
Chapter~{\bf 1} is purely motivational. Here we discuss proper smooth
models both of functional and number fields, and indicate why
non-proper models (e.g.\ the affine line $\bbA^1_k$ as a model of
$k(t)$, or $\Spec\bbZ$ as a model of $\bbQ$) are not sufficient for
some interesting applications. We also introduce some notations.  For
example, $\CompZ$ denotes the ``compactification'' of $\Spec\bbZ$. Its
closed points must correspond to all valuations of $\bbQ$, archimedian
or not, i.e.\ we expect $\CompZ=\Spec\bbZ\cup\{\infty\}$ as a set,
where $\infty$ denotes a new ``archimedian point''.  We denote by
$\Zinfty\subset\bbQ_\infty:=\bbR$ and $\Zninfty\subset\bbQ$ the
completed and non-completed local rings of $\CompZ$ at~$\infty$,
analogous to classical notations $\bbZ_p\subset\bbQ_p$ and $\bbZ_{(p)}
\subset\bbQ$.

\nxsubpoint
It is important to notice here that $\CompZ$, $\Zinfty$
and $\Zninfty$ are not defined in this chapter. Instead, they are used 
in an informal way to describe the properties we would expect these
objects to have. In this way we are even able to explain 
the classical approach to Arakelov geometry, which insists on
defining an Arakelov model $\bar\sX$ of a smooth 
projective algebraic variety $X/\bbQ$
as a flat proper model $\sX\to\Spec\bbZ$ together with a metric on $X(\bbC)$
subject to certain restrictions (e.g.\ being a K\"ahler metric).

\nxsubpoint
Another thing discussed in this chapter is that the problem of constructing
models over $\CompZ$ can be essentially reduced to the problem of
constructing $\Zninfty$-models of algebraic varieties $X/\bbQ$,
or $\Zinfty$-models of algebraic varieties $X/\bbR$. In other words,
we need a notion of a $\Zinfty$-structure on an algebraic variety~$X$ 
over~$\bbR$; if $X=\Spec A$ is affine, this is the same thing as a
$\Zinfty$-structure on an $\bbR$-algebra~$A$. So we see that
a proper understanding of ``compactified'' models of algebraic varieties
over~$\bbQ$ must include an understanding of $\Zinfty$-structures on
$\bbR$-algebras and vector spaces.

\nxpointtoc{$\Zinfty$-structures}
Chapter~{\bf 2} is dedicated to a detailed study of $\Zinfty$-structures
on real vector spaces and algebras. We start from the simplest cases
and extend our definitions step by step, arriving at the end
to the ``correct'' definition of $\catMod\Zinfty$, the category
of $\Zinfty$-modules. In this way we learn what the $\Zinfty$-modules
are, without still having a definition of $\Zinfty$ itself.

The main method employed here to obtain ``correct'' definitions
is the comparison with the $p$-adic case.

\nxsubpoint ($\Zinfty$-lattices: classical description.)
The first step is to describe $\Zinfty$-structures on a finite-dimensional
real vector space~$E$, i.e.\ {\em $\Zinfty$-lattices\/} $A\subset E$.

The classical solution is this. In the $p$-adic case a $\bbZ_p$-lattice
$A$ in a finite-dimensional $\bbQ_p$-vector space~$E$ defines a maximal
compact subgroup $K_A:=\Aut_{\bbZ_p}(A)\cong GL(n,\bbZ_p)$ in
locally compact group $G:=\Aut_{\bbQ_p}(E)\cong GL(n,\bbQ_p)$, all
maximal compact subgroups of~$G$ arise in this way, and $K_A=K_{A'}$
iff $A'$ and $A$ are similar, i.e.\ $A'=\lambda A$ for some $\lambda\in
\bbQ_p^\times$.

Therefore, it is reasonable to expect similarity classes of
$\Zinfty$-lattices inside real vector space~$E$ to be in one-to-one
correspondence with maximal compact subgroups $K$ of locally compact
group $G:=\Aut_\bbR(E)\cong GL(n,\bbR)$. Such maximal compact subgroups
are exactly the orthogonal subgroups $K_Q\cong O(n,\bbR)$, defined
by positive definite quadratic forms $Q$ on~$E$, and $K_Q=K_{Q'}$
iff $Q$ and $Q'$ are proportional, i.e.\ $Q'=\lambda Q$ for some $\lambda>0$.

In this way the classical answer is that a $\Zinfty$-structure on
a finite dimensional real space~$E$ is just a positive definite quadratic
form on~$E$, and similarly, a $\barZinfty$-structure on a finite dimensional
complex vector space is a positive definite hermitian form. This point of
view, if developed further, explains why classical Arakelov geometry
insists on equipping (complex points of) all varieties and vector bundles
involved with hermitian metrics.

\nxsubpoint ($\Zinfty$-structures on finite $\bbR$-algebras.)
Now suppose that $E$ is a finite $\bbR$-algebra. We would like to describe
$\Zinfty$-structures on this algebra, i.e.\ $\Zinfty$-lattices $A\subset E$,
compatible with the multiplication and unit of~$E$. In the $p$-adic case
this would actually mean $1\in A$ and $A\cdot A\subset A$, but if we want
to obtain ``correct'' definitions in the archimedian case, we must
re-write these conditions for $A$ in terms of corresponding maximal
compact subgroup $K_A\subset G$.

And here a certain problem arises. These compatibility conditions
cannot be easily expressed in terms of maximal compact subgroups of
automorphism groups even in the $p$-adic case. However, if we consider
maximal compact submonoids of endomorphism monoids instead, this problem
disappears.

\nxsubpoint (Maximal compact submonoids: $p$-adic case.)
Thus we are induced to describe $\Zinfty$-lattices $A$ 
in a finite-dimensional real space~$E$ 
in terms of maximal compact submonoids $M_A$ of locally compact monoid
$M:=\End_{\bbR}(E)\cong M(n,\bbR)$. When we study the corresponding 
$p$-adic problem, we see that all maximal compact submonoids of
$\End_{\bbQ_p}(E)$ are of form $M_A:=\{\phi:\phi(A)\subset A\}\cong
\End_{\bbZ_p}(A)\cong M(n,\bbZ_p)$ for a $\bbZ_p$-lattice $A\subset E$,
and that $M_A=M_{A'}$ iff $A$ and $A'$ are similar, i.e.\ 
in the $p$-adic case maximal compact submonoids of $\End(E)$
are classified by similarity classes of lattices $A\subset E$,
exactly in the same way as maximal compact subgroups of $\Aut(E)$.

\nxsubpoint (Maximal compact submonoids: archimedian case.)
However, in the archimedian case there are much more maximal compact submonoids
$M_A$ inside $M:=\End_\bbR(E)$ than positive definite quadratic forms. Namely,
we can take any symmetric compact convex body $A\subset E$
(a convex {\em body\/} is always required to be absorbent, 
i.e.\ $E=\bbR\cdot A$), and put
$M_A:=\{\phi\in\End_\bbR(E):\phi(A)\subset A\}$. Any such $M_A$ is a maximal
compact submonoid in $M=\End_\bbR(E)$, all maximal compact submonoids
in~$M$ arise in this way, and $M_A=M_{A'}$ iff $A'=\lambda A$ for some
$\lambda\in\bbR^\times$.

\nxsubpoint\label{sp:intro.zinflat} ($\Zinfty$-lattices.)
This leads us to define a {\em $\Zinfty$-lattice~$A$}
in a finite-dimensional real space~$E$ as a symmetric compact convex body~$A$.
It is well-known that any such~$A$ is the unit ball with respect to some
Banach norm $\vnorm$ on~$E$; in other words, {\em a $\Zinfty$-structure
on~$E$ is essentially a Banach norm on~$E$.}

Next, we can define the {\em category $\ZinfLat$ of $\Zinfty$-lattices}
as follows. Objects of $\ZinfLat$ are couples $A=(A_{\Zinfty},A_\bbR)$,
where $A_\bbR$ is a finite-dimensional real vector space, and
$A_{\Zinfty}\subset A_\bbR$ is a symmetric compact convex body in~$A_\bbR$.
Morphisms $f:A\to B$ are couples $(f_{\Zinfty},f_\bbR)$, where $f_\bbR:A_\bbR
\to B_\bbR$ is an $\bbR$-linear map, and 
$f_{\Zinfty}:A_{\Zinfty}\to B_{\Zinfty}$ is its restriction.

Of course, this definition is again motivated by the $p$-adic case.
Notice, however, that in the $p$-adic case we might define the category
of $\bbZ_p$-lattices without any reference to ambient $\bbQ_p$-vector
spaces, while we are still not able to describe $\Zinfty$-lattices
without reference to a real vector space.

Another interesting observation is that $\ZinfLat$ is essentially
the category of finite-dimensional Banach vector spaces, with
$\bbR$-linear maps of norm $\leq1$ as morphisms. While this description
establishes a connection to Banach norms, we don't insist on using it
too much, since the $p$-adic case suggests that we should
concentrate our attention on $A_{\Zinfty}$, not on ambient space~$A_\bbR$.

\nxsubpoint\label{sp:intro.zinf.alg} 
($\Zinfty$-structures on finite $\bbR$-algebras.)
Now we are able to define a $\Zinfty$-structure~$A$ on a finite 
$\bbR$-algebra~$E$, i.e.\ a $\Zinfty$-lattice $A\subset E$, compatible
with the multiplication and unit of~$E$. The key idea here is to
express this compatibility in terms of corresponding maximal compact
submonoid $M_A\subset\End(E)$ first in the $p$-adic case, and then
to transfer the conditions on $M_A$ thus obtained to the archimedian case.

The final result is that we must consider symmetric compact convex
bodies $A\subset E$, such that $1\in\partial A$ and $A\cdot A\subset A$.
In the language of Banach norms this means $\|1\|=1$ and $\|xy\|\leq
\|x\|\cdot\|y\|$, i.e.\ we recover the notion of a finite-dimensional
real Banach algebra.

\nxsubpoint (From lattices to torsion-free modules.)
Next step is to embed $\ZinfLat$ into a larger category of
{\em torsion-free\/ $\Zinfty$-modules\/ $\ZinfFlat$.}
Comparing to the $p$-adic case, we see that $\ZinfFlat$ might
be constructed as the category $\Ind(\ZinfLat)$ of ind-objects
over $\ZinfLat$ (cf.\ SGA~4~I). This category consists of
``formal'' inductive limits $\quotedinjlim M_\alpha$ 
of $\Zinfty$-lattices,
taken along filtered ordered sets or small categories, with morphisms
given by $\Hom(\quotedinjlim M_\alpha, \quotedinjlim N_\beta)=
\projlim_\alpha\injlim_\beta\Hom(M_\alpha,N_\beta)$.

\nxsubpoint (Direct description of $\ZinfFlat$.)
However, the category $\ZinfFlat$ of torsion-free $\Zinfty$ admits
a more direct description, similar to~\ptref{sp:intro.zinflat}. 
Namely, one can define it as the category of couples $(A_{\Zinfty},A_\bbR)$,
where $A_\bbR$ is a real vector space (not required more to be 
finite-dimensional), and $A_{\Zinfty}\subset A_\bbR$ is a symmetric
convex body (not required to be compact, but still required to be absorbent:
$A_\bbR=\bbR\cdot A_{\Zinfty}$). Morphisms are defined in the same way
as for $\Zinfty$-lattices, and it is immediate from this construction
that $\ZinfLat$ is a full subcategory of $\ZinfFlat$.

\nxsubpoint (Inductive and projective limits in $\ZinfFlat$.)
We show that arbitrary inductive and projective limits exist in
$\ZinfFlat$. For example, product $A\times B$ equals
$(A_{\Zinfty}\times B_{\Zinfty},A_\bbR\times B_\bbR)$, and the
direct sum (i.e.\ coproduct) $A\oplus B$ can be computed as
$(\conv(A_{\Zinfty}\cup B_{\Zinfty}),A_\bbR\oplus B_\bbR)$,
where $\conv(S)$ denotes the convex hull of a subset~$S$ in a
real vector space.

\nxsubpoint (Tensor structure on $\ZinfFlat$.)
We define an ACU $\otimes$-structure on $\ZinfFlat$ and $\ZinfLat$,
having the following property. An algebra $(A,\mu,\epsilon)$ 
in $\ZinfLat$ is the same thing as a couple $(A_{\Zinfty},A_\bbR)$,
consisting of a finite dimensional $\bbR$-algebra $A_\bbR$,
and a $\Zinfty$-lattice $A_{\Zinfty}\subset A_\bbR$, compatible with
the algebra structure of~$A_{\bbR}$ in the sense of~\ptref{sp:intro.zinf.alg}.

This tensor structure is extremely natural in several other respects.
For example, when translated into the language of (semi)norms on real
vector spaces, it corresponds to Grothendieck's projective tensor product
of (semi)norms.

The unit object for this tensor structure is $\Zinfty:=([-1,1],\bbR)$.
One should think of this $\Zinfty$ as ``$\Zinfty$, considered as a left
module over itself'', not as ``the ring~$\Zinfty$''.

\nxsubpoint (Underlying set of a torsion-free $\Zinfty$-module.)
Once we have a unit object $\Zinfty$, we can define the
``forgetful functor'' $\Gamma:\ZinfFlat\to\catSets$,
$A\mapsto\Hom(\Zinfty,A)$. It is natural to say that $\Gamma(A)$
is the {\em underlying set of~$A$}.
Direct computation shows that $\Gamma(A)=A_{\Zinfty}$ for $A=(A_{\Zinfty},
A_\bbR)$, thus suggesting that we ought to concentrate our attention 
on~$A_{\Zinfty}$, not on auxiliary vector space~$A_\bbR$.

\nxsubpoint (Free $\Zinfty$-modules. Octahedral combinations.)
Since arbitrary direct sums exist in $\ZinfFlat$, we can construct
{\em free $\Zinfty$-modules\/} $\Zinfty^{(S)}$, by taking the direct sum
of $S$ copies of $\Zinfty$. It is immediate that $L_{\Zinfty}:S\mapsto
\Zinfty^{(S)}$ is a left adjoint to $\Gamma=\Gamma_{\Zinfty}:
\ZinfFlat\to\catSets$. One can describe this $\Zinfty^{(S)}=(\Sigma_\infty(S),
\bbR^{(S)})$ explicitly. Its vector space component is simply $\bbR$-vector
space $\bbR^{(S)}$ freely generated by~$S$. Its standard 
basis elements will be denoted by $\{s\}$, $s\in S$. Then
the symmetric convex subset $\Sigma_\infty(S)\subset\bbR^{(S)}$ 
consists of all {\em octahedral (linear) combinations\/} of these basis
vectors:
\begin{align}
\Sigma_\infty(S)&=\conv\bigl(\pm\{s\}\,:\,s\in S\bigr)\\
&=\Bigl\{\sum_s\lambda_s\{s\}\,:\,
\sum_s|\lambda_s|\leq1,\text{ almost all $\lambda_s=0$}\Bigr\}
\end{align}
In other words, $\Sigma_\infty(S)$ is the {\em standard octahedron\/}
in $\bbR^{(S)}$.

\nxsubpoint (Notation: finite sets and basis elements.)
We would like to mention here some notation, used throughout this work.
If $M$ is a free ``object'', generated by a set~$S$
(e.g.\ $M=R^{(S)}$ for some ring~$R$), we denote by $\{s\}$ the 
``basis element'' of $M$ corresponding to $s\in S$. Thus $s\mapsto\{s\}$
is the natural embedding $S\to M$.

Another notation: we denote by $\stn$ the {\em standard finite set\/}
$\{1,2,\ldots,n\}$, where $n\geq0$ is any integer. For example,
$\st0=\emptyset$, and $\st2=\{1,2\}$. Furthermore, whenever we have
a functor $\Sigma$ defined on the category of sets, we write $\Sigma(n)$
instead of $\Sigma(\stn)$, for any $n\geq0$. In this way $\bbR^{(n)}=\bbR^n$
is the standard $n$-dimensional real vector space, with standard basis $\{k\}$,
$1\leq k\leq n$, and $\Sigma_\infty(n)\subset\bbR^n$ is the standard
$n$-dimensional octahedron.

\nxsubpoint (From $\ZinfFlat$ to $\catMod\Zinfty$.)
Our next step is to recover the category $\catMod\Zinfty$ of
all $\Zinfty$-modules, starting from the category $\ZinfFlat$ of
torsion-free $\Zinfty$-modules. We use adjoint functors
$L_{\Zinfty}:\catSets\leftrightarrows\ZinfFlat:\Gamma$ for this.
Namely, we observe that this pair of adjoint functors defines
a {\em monad\/} structure $(\Sigma_\infty,\mu,\epsilon)$ 
on endofunctor $\Sigma_\infty:=\Gamma L_{\Zinfty}:\catSets\to\catSets$.
However, functor $\Gamma$ happens not to be monadic, i.e.\
the arising ``comparison functor'' $I:\ZinfFlat\to\catSets^{\Sigma_\infty}$
from $\ZinfFlat$ into the category of $\Sigma_\infty$-algebras
in $\catSets$ (which we prefer to call $\Sigma_\infty$-{\em modules})
is not an equivalence of categories, but just a fully faithful functor.

The $p$-adic analogy suggests to define $\catMod\Zinfty:=
\catSets^{\Sigma_\infty}$. Then category $\ZinfFlat$ can be identified
with a full subcategory of $\catMod\Zinfty$ with the aid of
functor~$I$, and the forgetful functor
$\Gamma_{\Zinfty}:\catMod\Zinfty\to\catSets$ is now monadic by construction,
similarly to the forgetful functor on any category defined by an algebraic
system (e.g.\ category $\catMod R$ of modules over an associative
ring~$R$).

\nxsubpoint (Explicit description of $\Zinfty$-modules.)
We can obtain a more explicit description of $\Zinfty$-modules,
i.e.\ of objects $M=(M,\alpha)\in\Ob\catMod\Zinfty=
\Ob\catSets^{\Sigma_\infty}$. Indeed, by definition $M=(M,\alpha)$
consists of a set~$M$, together with a ``$\Zinfty$-structure'',
i.e.\ a map $\alpha:\Sigma_\infty(M)\to M$, subject to certain conditions.
Since $\Sigma_\infty(M)$ consists of formal octahedral combinations
$\sum_x\lambda_x\{x\}$, $\sum_x|\lambda_x|\leq1$, of elements of~$M$,
such a map~$\alpha$ should be thought of as a way of evaluating such
formal octahedral combinations. Thus we write $\sum_x\lambda_xx$
instead of $\alpha(\sum_x\{x\})$. We use this notation for finite sums as
well, e.g.\ 
$\lambda x+\mu y$ actually means $\alpha(\lambda\{x\}+\mu\{y\})$,
for any $x$, $y\in M$, $\lambda$, $\mu\in\bbR$, $|\lambda|+|\mu|\leq1$.
Notice, however, that the $+$ sign in expression $\lambda x+\mu y$
is completely formal: in general we don't get any addition operation on~$M$.

As to the conditions for~$\alpha$, they are exactly all the usual relations
satisfied by octahedral combinations of elements of real vector spaces. 
For example, $\nu(\lambda x+\mu y)=\nu\lambda\cdot x+\nu\mu\cdot y$.

In this way we see that a $\Zinfty$-module is nothing else than a set~$M$,
together with a way of evaluating octahedral combinations of its elements,
in such a way that all usual relations hold in~$M$.

\nxsubpoint\label{sp:intro.tors.zinfmod} (Torsion modules.)
Category $\catMod\Zinfty$ is strictly larger than the category
$\ZinfFlat$ of torsion-free $\Zinfty$-modules. In fact, it contains
objects like $\Finfty:=
\Coker(\gm_\infty\rightrightarrows\bbZ_\infty)$, where $\gm_\infty=
((-1,1),\bbR)$ in $\ZinfFlat$, and the two morphisms
$\gm_\infty\to\Zinfty$ are the natural inclusion and the zero
morphism. Notice that this cokernel is zero in $\ZinfFlat$, but not in
$\catMod\Zinfty$: in fact, this $\Finfty$ is a three-element set with a
certain $\Zinfty$-module structure.  This gives an example of a
non-trivial torsion $\Zinfty$-module.

\nxsubpoint (Arakelov affine line.)
Among other things discussed in Chapter~{\bf 2}, we would like to mention
the study of ``Arakelov affine line'' $\Spec\Zinfty[T]$, carried
both from the point of view of (co)metrics and from that of
(generalized) schemes.

When we compute the naturally arising cometric on $\bbA^1_\bbR$ or 
$\bbA^1_\bbC$, coming from this $\Zinfty$-structure $\Zinfty[T]$ on $\bbR[T]$,
it turns out to be identically zero outside unit disk 
$\{\lambda:|\lambda|<1\}$, i.e.\ all points outside this disk are at infinite
distance from each other. Inside the disk we obtain a continuous piecewise
smooth cometric, very similar to Poincar\'e model
of hyperbolic plane in the unit disk.

On the other hand, we can use the $\otimes$-structure defined on $\ZinfFlat$
(and actually on all of $\catMod\Zinfty$) to define $\Zinfty$-algebras~$A$,
modules over them, ideals and prime ideals inside them and so on,
thus obtaining a definition of prime spectrum $\Spec A$. For example,
topological space $\Spec\Zinfty$ looks like the spectrum of a DVR.
When we study the ``Arakelov affine line'' $\Spec\Zinfty[T]$ from this
point of view, we observe some unexpected phenomena, e.g.\ Krull dimension
of this topological space turns out to be infinite.

\nxpointtoc{$\otimes$-categories, algebras and monads}
Chapter {\bf 3} collects some general definitions and constructions,
related to $\otimes$-categories, algebras and monads. It is quite technical,
but nevertheless quite important for the next two chapters.
Most results collected here are well-known and can be found
in~\cite{MacLane}; however, we want to fix the terminology, and
discuss the generalization to the topos case.

\nxsubpoint (AU $\otimes$-categories and $\otimes$-actions.)
We discuss AU $\otimes$-categories, i.e.\ categories $\cA$,
equipped with a tensor product functor $\otimes:\cA\times\cA\to\cA$
and a unit object $\Unit\in\Ob\cA$, together with some {\em associativity
constraint\/} $\alpha_{X,Y,Z}:(X\otimes Y)\otimes Z\cong X\otimes(Y\otimes Z)$
and {\em unit constraints\/} $\Unit\otimes X\cong X\cong X\otimes\Unit$,
subject to certain axioms (e.g.\ the pentagon axiom). Roughly speaking,
these axioms ensure that multiple tensor products $X_1\otimes X_2\otimes\cdots
\otimes X_n$ are well-defined and have all the usual properties
of multiple tensor products, apart from commutativity. 

After that we discuss {\em external (left) $\otimes$-actions\/}
$\obslash:\cA\times\cB\to\cB$ of an AU $\otimes$-category~$\cA$ on a category
$\cB$. Here we have to impose some external associativity and unit
constraints $(X\otimes Y)\obslash M\cong X\obslash(Y\obslash M)$ and
$\Unit\obslash M\cong M$, subject to similar relations.

Of course, any AU $\otimes$-category $\cA=(\cA,\otimes)$ admits a
natural $\otimes$-action on itself, given by $\obslash:=\otimes$.

\nxsubpoint (Algebras and modules.)
Whenever we have an AU $\otimes$-category~$\cA$, we can consider
{\em algebras\/} $A=(A,\mu,\epsilon)$, $A\in\Ob\cA$, $\mu:A\otimes A\to A$,
$\epsilon:\Unit\to A$, always supposed to be associative with unity
(but not commutative -- in fact, commutativity doesn't make sense
without a commutativity constraint on~$\cA$). Thus we obtain
a {\em category of algebras in~$\cA$}, denoted by $\catAlg(\cA)$.

Next, if $\obslash:\cA\times\cB\to\cB$ is an external $\otimes$-action,
and $A$ is an algebra in~$\cA$, we denote by $\catMod A$ or $\cB^A$
the category of {\em $A$-modules in~$\cB$}, consisting of couples
$(M,\alpha)$, $M\in\Ob\cB$, $\alpha:A\obslash M\to M$, subject to classical
relations $\alpha\circ(\epsilon\obslash 1_M)=1_M$ and 
$\alpha\circ(1_A\obslash\alpha)=\alpha\circ(\mu\obslash 1_M)$.

\nxsubpoint (Monads over a category~$\cC$.)
Now if $\cC$ is an arbitrary category, the category of endofunctors
$\cA:=\catEndof(\cC)=\catFunct(\cC,\cC)$ admits a natural AU 
$\otimes$-structure, given by composition of functors: $F\otimes G:=F\circ G$.
Furthermore, there is a natural $\otimes$-action of~$\cA$ on~$\cC$, defined
by $F\obslash X:=F(X)$.

Then a {\em monad\/} $\Sigma=(\Sigma,\mu,\epsilon)$ over a category~$\cC$
is simply an algebra in this category of endofunctors~$\cA$, i.e.\ 
$\Sigma:\cC\to\cC$ is an endofunctor and $\mu:\Sigma^2\to\Sigma$,
$\epsilon:\Id_\cC\to\Sigma$ are natural transformations, subject to
associativity and unit relations: $\mu\circ(\Sigma\star\mu)=\mu\circ(\mu
\star\Sigma)$ and $\mu\circ(\epsilon\star\Sigma)=\id_\Sigma=\mu\circ(\Sigma
\star\epsilon)$, or equivalently, $\mu_X\circ\Sigma(\mu_X)=\mu_X\circ
\mu_{\Sigma(X)}:\Sigma^3(X)\to\Sigma(X)$, and $\mu_X\circ\epsilon_{\Sigma(X)}=
\id_{\Sigma(X)}=\mu_X\circ\Sigma(\epsilon_X)$, for any $X\in\Ob\cC$.

Similarly, the category $\cC^\Sigma$ of $\Sigma$-modules (in~$\cC$)
is defined as the category of $\Sigma$-modules in~$\cC$ with respect to
the external $\otimes$-action of~$\cA$ on~$\cC$ just discussed.
In this way $\cC^\Sigma$ consists of couples $M=(M,\alpha)$, with
$M\in\Ob\cC$, $\alpha:\Sigma(M)\to M$, such that $\alpha\circ\Sigma(\alpha)=
\alpha\circ\mu_M:\Sigma^2(M)\to M$, and $\alpha\circ\epsilon_M=\id_M$.

\nxsubpoint (Monads and adjoint functors.)
Whenever we have a monad $\Sigma$ over a category~$\cC$, we get a 
forgetful functor $\Gamma_\Sigma:\cC^\Sigma\to\cC$, $(M,\alpha)\mapsto M$,
which admits a left adjoint $L_\Sigma:\cC\to\cC^\Sigma$, $S\mapsto(\Sigma(S),
\mu_S)$, such that $\Gamma_\Sigma\circ L_\Sigma=\Sigma$. Conversely,
given any two adjoint functors $L:\cC\leftrightarrows\cD:\Gamma$,
we obtain a canonical monad structure on $\Sigma:=\Gamma\circ L:\cC\to\cC$,
together with a ``comparison functor'' $I:\cD\to\cC^\Sigma$, such that
$\Gamma=\Gamma_\Sigma\circ I$. Functor $\Gamma$ is said to be {\em monadic\/}
if $I$ is an equivalence of categories.

\nxsubpoint (Examples: $R$-modules, $\Zinfty$-modules\dots)
For example, if $\cC=\catSets$, $\cD=\catMod R$ is the category of
left modules over a ring~$R$ (always supposed to be associative with unity),
then the forgetful functor $\Gamma:\catMod R\to\catSets$ turns out to 
be monadic. The corresponding monad $\Sigma_R:\catSets\to\catSets$
transforms a set $S$ into the underlying set of free $R$-module $R^{(S)}$
generated by~$S$, i.e.\ into the set of all formal $R$-linear combinations
of basis elements $\{s\}$, $s\in S$. In this way an $R$-module~$M$
is just a set~$M$ together with a method $\alpha:\Sigma_R(M)\to M$ of
evaluating formal $R$-linear combinations of its elements, subject to some
natural conditions.

Similarly, $\catMod\Zinfty$ was defined to be $\catSets^{\Sigma_\infty}$
for a certain monad $\Sigma_\infty$ over $\catSets$, so the set
$\Sigma_\infty(S)$ of formal octahedral combinations of elements of~$S$
should be thought of as the set of all formal $\Zinfty$-linear combinations
of elements of~$S$, or as the underlying set of free $\Zinfty$-module
$\Zinfty^{(S)}$. This was the way we've defined $\Sigma_\infty$ in the
first place.

This analogy between $\catMod R=\catSets^{\Sigma_R}$ and
$\catMod\Zinfty=\catSets^{\Sigma_\infty}$ suggests that our category of
generalized rings, which is expected to contain all classical
rings~$R$ as well as exotic objects like~$\Zinfty$, might be constructed
as a certain full subcategory of the category of monads over~$\catSets$.
Then we might put $\catMod\Sigma:=\catSets^\Sigma$ for any monad
$\Sigma$ from this full subcategory, thus obtaining a reasonable
definition of $\Sigma$-modules without any special considerations.

\nxsubpoint (Topos case: inner endofunctors and inner monads.)
Apart from the things just discussed, Chapter~{\bf 3} contains
some technical definitions and statements, related to
{\em inner endofunctors\/} and {\em inner monads\/} over a topos~$\cE$.
These notions are used later to transfer the definitions of generalized
rings and modules over them from the category of sets into arbitrary topoi,
e.g.\ categories of sheaves of sets over a topological space~$X$.
This is necessary to obtain a monadic interpretation of sheaves
of generalized rings and of modules over them. Without such notions we
wouldn't be able to discuss generalized ringed spaces and in particular 
generalized schemes. However, we suggest to the reader to skip
these technical topos-related pages during the first reading.

\nxpointtoc{Algebraic monads as non-commutative generalized rings}
Chapter~{\bf 4} is dedicated to the study of {\em algebraic\/}
endofunctors and monads over $\catSets$ (and algebraic {\em inner\/}
endofunctors and monads over a topos $\cE$ as well). This
notion of {\em algebraicity\/} is actually the first condition
we need to impose on monads over $\catSets$ in order to define
the category of generalized rings. In fact, {\em algebraic monads
are non-commutative generalized rings.} For example, any monad
$\Sigma_R$ defined by a classical ring~$R$ 
(as usual, associative with unity) is algebraic, as well as monad 
$\Sigma_\infty$, used to define $\catMod\Zinfty$.

\nxsubpoint (Algebraic monads vs.\ algebraic systems.)
Another important remark is that {\em an algebraic monad over $\catSets$
is essentially the same thing as an algebraic system.} 
We discuss this correspondence in more detail below. In some sense
algebraic systems are something like presentations (by a system 
of generators and relations) of algebraic monads. Thus different 
(but equivalent) algebraic systems may correspond to isomorphic
algebraic monads, and algebraic monads provide an invariant way of
describing algebraic systems. In this way the study of algebraic monads
might be thought of as the study of algebraic systems from a
categorical point of view, i.e.\ a categorical approach to {\em
universal algebra}.

\nxsubpoint (Algebraic endofunctors over $\catSets$.)
We say that an endofunctor $\Sigma:\catSets\to\catSets$ is
{\em algebraic\/} if it commutes with filtered inductive limits:
$\Sigma(\injlim_\alpha S_\alpha)\cong\injlim_\alpha\Sigma(S_\alpha)$.
Since any set is a filtered inductive limit of its finite subsets,
we obtain $\Sigma(S)=\injlim_{(\stn\stackrel\phi\to S)\in\catN_{/S}}
\Sigma(n)$, where $\catN=\{\st0,\st1,\ldots,\stn,\ldots\}_{n\geq0}$ 
denotes the category of {\em standard finite sets}, considered as a full
subcategory of $\catSets$. Therefore, {\em any algebraic endofunctor
$\Sigma:\catSets\to\catSets$ is completely determined by its 
restriction $\Sigma|_{\catN}:\catN\to\catSets$.} In fact, this
restriction functor $\Sigma\mapsto\Sigma|_{\catN}$ induces an equivalence
between the category of algebraic endofunctors $\cA_{alg}\subset\cA=
\catEndof(\catSets)$ and $\catFunct(\catN,\catSets)=\catSets^{\catN}$.
This means that {\em an algebraic endofunctor $\Sigma$ is essentially
the same thing as a countable collection of sets $\{\Sigma(n)\}_{n\geq0}$,
together with maps $\Sigma(\phi):\Sigma(n)\to\Sigma(m)$, defined for
any $\phi:\stn\to\stm$, such that $\Sigma(\psi\circ\phi)=\Sigma(\psi)\circ
\Sigma(\phi)$ and $\Sigma(\id_\stn)=\id_{\Sigma(n)}$.}

\nxsubpoint (Algebraic monads.)
On the other hand, if two endofunctors $\Sigma$ and $\Sigma'$
commute with filtered inductive limits, the same is true for
their composition $\Sigma\otimes\Sigma'=\Sigma\circ\Sigma'$,
i.e.\ {\em $\cA_{alg}\cong\catSets^{\catN}$ is a full $\otimes$-subcategory
of $\cA=\catEndof(\catSets)$.} Therefore, we can define an
{\em algebraic monad\/} $\Sigma=(\Sigma,\mu,\epsilon)$ as an algebra
in $\cA_{alg}$. Of course, an algebraic monad is just a monad,
such that its underlying endofunctor commutes with filtered inductive 
limits.

\nxsubpoint (Why algebraic?)
Notice that monads $\Sigma_R$ defined by an associative ring~$R$
are algebraic, i.e.\ $\Sigma_R(S)=\bigcup_{\text{finite $I\subset S$}}
\Sigma_R(I)$, just because any element of $\Sigma_R(S)=R^{(S)}$,
i.e.\ any formal $R$-linear combination of elements of~$S$,
involves only finitely many elements of~$S$, hence comes from $\Sigma_R(I)$
for some finite subset $I\subset S$. Similarly, $\Sigma_\infty$ is
algebraic, because any octahedral combination $\sum_{s\in S}\lambda_s\{s\}
\in\Sigma_\infty(S)$ has only finitely many $\lambda_s\neq0$.

We can express this by saying that we consider only ``operations''
depending on finitely many arguments. For example, we might remove
the requirement ``$\lambda_s\neq0$ only for finitely many $s\in S$'' and
consider ``infinite octahedral combinations'' $\sum\lambda_s\{s\}$,
with the only requirement $\sum_s|\lambda_s|\leq1$. In this way we
obtain a larger monad $\hat\Sigma_\infty\supset\Sigma_\infty$,
which coincides with $\Sigma_\infty$ on finite sets, but is different on
larger sets. A $\hat\Sigma_\infty$-structure on a set~$M$ is essentially
a way of computing such ``infinite octahedral combinations'' of elements
of~$M$. This is definitely not an algebraic operation, and $\hat\Sigma_\infty$
is not an algebraic monad.

Therefore, word {\em ``algebraic''\/} means here something like
``expressible in terms of operations involving only finitely many 
arguments''.

\nxsubpoint (Algebraic monads and operations.)
Now let $\Sigma$ be an algebraic endofunctor (over $\catSets$),
$M$ be a set, $\alpha:\Sigma(M)\to M$ be any map of sets. For example,
we might take an algebraic monad $\Sigma$ and a $\Sigma$-module~$M$.

Let $t\in\Sigma(n)$ for some integer~$n$. Take any $n$-tuple
$x=(x_1,\ldots,x_n)\in M^n$ of elements of~$M$. It can be considered
as a map $\tilde x:\stn\to M$, $k\mapsto x_k$, hence we get a map
$\Sigma(\tilde x):\Sigma(n)\to\Sigma(M)$. Now we can apply 
$\alpha\circ\Sigma(\tilde x)$ to~$t$, thus obtaining an element of~$M$:
\begin{equation}
t(x_1,\ldots,x_n)={[t]}_M(x_1,\ldots,x_n):=(\alpha\circ\Sigma(\tilde x))(t)
\end{equation}

In this way any $t\in\Sigma(n)$ defines an $n$-ary operation
${[t]}_M:M^n\to M$ on~$M$. That's why we say that $\Sigma(n)$ is
{\em the set of $n$-ary operations of~$\Sigma$}, and call its elements
{\em $n$-ary operations of~$\Sigma$.} Furthermore, we say that 
${[t]}_M:M^n\to M$ is the {\em value\/} of operation~$t$ on~$M$.
Of course, when $n=0$, we speak about {\em constants} and their values
(any constant $c\in\Sigma(0)$ has a value ${[c]}_M\in M$),
and when $n=1,2,3,\ldots$ we obtain {\em unary, binary, ternary, \dots\ 
operations}.

One can show that giving a map $\alpha:\Sigma(M)\to M$ is actually 
{\em equivalent\/} to giving a family of 
``evaluation maps'' $\alpha^{(n)}:\Sigma(n)\times M^n\to
M$, $(t,x_1,\ldots,x_n)\mapsto{[t]}_M(x_1,\ldots,x_n)$, $n\geq0$,
satisfying natural compatibility relations, which can be written as
\begin{equation}
{[\phi_*(t)]}_M(x_1,\ldots,x_n)={[t]}_M(x_{\phi(1)},\ldots,x_{\phi(m)}),
\quad\forall t\in\Sigma(m), \phi:\stm\to\stn
\end{equation}
Here $\phi_*$ is a shorthand for $\Sigma(\phi)$.

\nxsubpoint\label{sp:intro.elem.algmon} 
(Elementary description of algebraic monads.)
The above description is applicable to maps 
$\mu_n:\Sigma(\Sigma(n))\to\Sigma(n)$. We see that such a map is completely
determined by a sequence of ``evaluation'' or ``substitution'' 
maps $\mu_n^{(k)}:\Sigma(k)\times\Sigma(n)^k\to\Sigma(n)$,
$(t,t_1,\ldots,t_k)\mapsto t(t_1,\ldots,t_k)=
{[t]}_{\Sigma(n)}(t_1,\ldots,t_k)$. In this way we obtain an
``elementary'' description of an algebraic monad~$\Sigma$,
consisting of 
a sequence of sets $\{\Sigma(n)\}_{n\geq0}$, transition maps
$\Sigma(\phi):\Sigma(m)\to\Sigma(n)$, defined for all $\phi:\stm\to\stn$,
evaluation maps $\mu_n^{(k)}:\Sigma(k)\times\Sigma(n)^k\to\Sigma(n)$,
and a ``unit element'' $\bu:=\epsilon_{\st1}(1)\in\Sigma(1)$.

In fact, this collection completely determines the algebraic monad~$\Sigma$,
and conversely, if we start from such a collection, satisfying some
natural compatibility conditions 
(certain ``associativity'' and ``unit'' axioms), we obtain a uniquely
determined algebraic monad~$\Sigma$.

Of course, $\Sigma$-modules $M=(M,\alpha)$ also admit such an
``elementary'' description, consisting of a set~$M$, and a collection
of ``evaluation maps'' $\alpha^{(n)}:\Sigma(n)\times M^n\to M$,
subject to certain compatibility, associativity and unit conditions
(e.g.\ $[\bu]_M=\id_M$).

Therefore, we might have defined algebraic monads and modules over them
in such an ``elementary'' fashion. We didn't do this just because
arising definitions and especially relations seem to be quite complicated
and not too enlightening, unless one knows that they come from 
the definition of algebraic monads, i.e.\ they are the algebra relations 
in~$\cA_{alg}$, written in explicit form.

\nxsubpoint (Unit and basis elements.)
Notice that the unit element $\bu=\epsilon_1(1)\in\Sigma(1)$ actually
completely determines the unit $\epsilon:\Id_{\catSets}\to\Sigma$,
since $\epsilon_X(x)=(\Sigma(i_x))(\bu)$ for any $x\in X$, where
$i_x:\st1\to X$ is the map $1\mapsto x$. 

Recall that we denote by $\{k\}_\stn$ or simply by $\{k\}$ the
``basis element'' $\epsilon_n(k)\in\Sigma(n)$, $1\leq k\leq n$.
For example, $\bu=\{1\}_{\st1}$.
These basis elements have some natural properties, e.g.\ 
$\phi_*\{k\}_\stm=\{\phi(k)\}_\stn$, for any $\phi:\stm\to\stn$, 
$1\leq k\leq m$. More interesting properties come from
the unit axioms for $\Sigma$ and $\Sigma$-modules. For example,
$\mu\circ(\Sigma\star\epsilon)=\id$ implies (and actually is equivalent to)
\begin{equation}
t(\{1\},\{2\},\ldots,\{n\})={[t]}_{\Sigma(n)}(\{1\}_\stn,\ldots,\{n\}_\stn)=t,
\quad\forall t\in\Sigma(n).
\end{equation}
The other unit axiom $\mu\circ(\epsilon\star\Sigma)=\id$ is actually equivalent
to $\bu(t)=t$ for all $n\geq0$, $t\in\Sigma(n)$.

Unit axiom for a $\Sigma$-module~$M$ is equivalent to
${[\{k\}_\stn]}_M=\pr_k:M^n\to M$, or just to ${[\bu]}_M=\id_M$.

\nxsubpoint (Special notation for unary operations.)
If $u\in\Sigma(1)$ is a unary operation, and $t\in\Sigma(n)$ is an $n$-ary
operation, we usually write $ut$ or $u\cdot t$ instead of ${[u]}_{\Sigma(n)}(t)
\in\Sigma(n)$, and similarly $ux:={[u]}_M(x)$ for any $\Sigma$-module $M$
and any $x\in M$. This notation is unambiguous because of associativity
relations
$(ut)(x_1,\ldots,x_n)=u\cdot t(x_1,\ldots,x_n)$ for any $u\in\Sigma(1)$,
$t\in\Sigma(n)$, $x_i\in M$ or $\Sigma(m)$.

In this way we obtain on set $|\Sigma|:=\Sigma(1)$ a monoid structure with
identity~$\bu$, and a monoid action of $|\Sigma|$ on the
underlying set of any $\Sigma$-module~$M$.

\nxsubpoint (Special notation for binary operations.)
Of course, when we have a binary operation, denoted by a sign like
$+$, $*$, \dots, usually written in infix form, we write $x+y$,
$x*y$ etc.\ instead of $[+](x,y)$, $[*](x,y)$ etc.

Since $t=t(\{1\},\ldots,\{n\})$ for any $t\in\Sigma(n)$, we can write
$[+]=\{1\}+\{2\}\in\Sigma(2)$ when we want to point out the corresponding
element of~$\Sigma(2)$.

\nxsubpoint (Free modules.)
Notice that $L_\Sigma(n)=(\Sigma(n),\mu_n)$ is the free $\Sigma$-module
of rank~$n$, i.e.\ $\Hom_\Sigma(L_\Sigma(n),M)\cong M^n$ for any
$\Sigma$-module~$M$. Of course, the map $\Hom_\Sigma(L_\Sigma(n),M)\to M^n$
is given by evaluating $f:L_\Sigma(n)\to M$ on basis elements $\{k\}_\stn$.
We also denote $L_\Sigma(n)$ simply by $\Sigma(n)$, when no confusion 
can arise.

\nxsubpoint (Set of unary operations.)
For any algebraic monad~$\Sigma$ the set $|\Sigma|:=\Sigma(1)$ 
has two natural structures: that of a monoid, and that of a $\Sigma$-module
(free of rank one). If $\Sigma=\Sigma_R$ for a classical associative
ring~$R$, then $|\Sigma|$ is the underlying set of~$R$, its monoid structure
is given by the multiplication of~$R$, and its module structure is the
natural left $R$-module structure on~$R$.

In this way $|\Sigma|$ plays the role of the underlying set of 
algebraic monad~$\Sigma$. Notice that it is always a monoid 
under multiplication, but in general it doesn't have any addition, i.e.\
{\em multiplication is in some sense more fundamental than addition.}

Another interesting observation is that, while the addition of a
classical ring $R$ indeed corresponds to a binary operation
$[+]=(1,1)\in\Sigma_R(2)=R^2$, the multiplication of~$R$ doesn't
correspond to any element of $\Sigma_R(n)$. Instead, it is built in
the ``composition maps'' $\mu_n^{(k)}$, i.e.\ it is part of the
more fundamental structure of algebraic monad.

\nxsubpoint (Matrix description. Comparison to Shai Haran's approach.)
Since $\Hom_\Sigma(\Sigma(n),\Sigma(m))\cong\Sigma(m)^n$,
we {\em define\/} the set of $m\times n$-matrices over $\Sigma$
by $M(m,n;\Sigma):=\Sigma(m)^n$. Composition of morphisms defines
maps $M(n,k;\Sigma)\times M(k,m;\Sigma)\to M(n,m;\Sigma)$, i.e.\
we obtain a well-defined {\em product of matrices}. Putting here $m=1$,
we get maps $\Sigma(n)^k\times\Sigma(k)\to\Sigma(n)$, which are
nothing else than the ``structural maps'' $\mu_n^{(k)}$ of~$\Sigma$,
up to a permutation of the two arguments.
This matrix language is often quite convenient. For example,
the ``associativity relations'' for maps $\mu_n^{(k)}$ are
essentially equivalent to associativity of matrix products,
and the ``unit relations'' actually mean that the identity matrix
$I_n:=(\{1\}_\stn,\ldots,\{n\}_\stn)\in M(n,n;\Sigma)$ is indeed
left and right identity with respect to matrix multiplication.

One might actually try to describe generalized rings in terms
of collections of ``matrix sets'' $\{M(m,n;\Sigma)\}_{m,n\geq0}$,
together with composition (i.e.\ matrix multiplication) maps as above,
and direct sum maps $M(m,n;\Sigma)\times M(m',n';\Sigma)\to 
M(m+m',n+n';\Sigma)$, $(u,v)\mapsto u\oplus v$. Pursuing this path one
essentially recovers Shai Haran's notion of an {\em $\bbF$-algebra}
(or rather its non-commutative counterpart, since we don't say
anything about tensor products at this point), defined in \cite{ShaiHaran}.
Furthermore, this notion of $\bbF$-algebra is more general than
our notion of generalized ring, since Shai Haran never requires
$M(m,n;\Sigma)=M(1,n;\Sigma)^m$.

In fact, once we impose this condition (which is quite natural if we want
to have $\Sigma(m)\oplus\Sigma(n)=\Sigma(m+n)$ for directs sums, i.e.\ 
coproducts of free modules) on (non-commutative) $\bbF$-algebras, 
we recover our notion of generalized ring. However, Shai Haran
doesn't impose such restrictions. 
Actually, he doesn't define a module over an $\bbF$-algebra~$A$
as a set~$Q$ with some additional structure (in fact, any algebraic structure
on a set is described by some algebraic monad, so our theory of generalized
rings is the largest algebraic theory of ring-like
objects, which admit a notion of module over them), but considers infinite
collections $\{M(m,n;Q)\}_{m,n\geq0}$ instead, 
thought of as ``$m\times n$-matrices
with entries in~$Q$''. This more complicated notion of module
roughly corresponds to our ``$\Sigma$-bimodules in~$\cA_{alg}$'',
i.e.\ to algebraic endofunctors, equipped with compatible left and 
right actions of~$\Sigma$.

It is not clear whether one might transfer more
sophisticated constructions of our work to Shai Haran's more general
case, since these constructions involve forgetful functors
$\Gamma:\catMod\Sigma\to\catSets$, i.e.\ expect $\Sigma$-modules
to be sets with additional structure, and heavily rely on formulas like
$\Sigma(m)\oplus\Sigma(n)=\Sigma(m+n)$ for categorial coproducts.

Unfortunately, we cannot say much more about this right now, since we haven't
found any publications of Shai Haran on this topic apart from his
original preprint \cite{ShaiHaran}, which deals only with the basic
definitions of his theory.

\nxsubpoint (Initial and final algebraic monads.)
Notice that the category of algebraic monads has an initial object~$\Fempty$,
given by the only monad structure on $\Id_{\catSets}$, as well as a 
final object~$\st1$, given by the constant functor with value~$\st1$,
equipped with its only monad structure.

This initial algebraic monad~$\Fempty$ has the property $\Fempty(n)=\stn$,
i.e.\ a free $\Fempty$-module of rank~$n$ consists only of basis
elements $\{k\}_\stn$, $1\leq k\leq n$. Furthermore, any set admits a
unique $\Fempty$-module structure, i.e.\ $\catMod\Fempty=\catSets$.

As to $\catMod{\st1}$, the only sets which admit a $\st1$-module structure
are the one-element sets, hence $\catMod{\st1}$ is equivalent to the
``final category'' $\st1$.

\nxsubpoint (Projective limits of algebraic monads. Submonads.)
Projective limits of algebraic monads can be computed componentwise:
$(\projlim_\alpha\Sigma_\alpha)(n)=\projlim_\alpha\Sigma_\alpha(n)$,
for any $n\geq0$. 
For example, $(\Sigma\times\Sigma')(n)=\Sigma(n)\times\Sigma'(n)$. 
An immediate consequence is that an algebraic monad
homomorphism $f:\Sigma'\to\Sigma$ is a {\em monomorphism\/} iff
all maps $f_n:\Sigma'(n)\to\Sigma(n)$ are injective. (Notice that
injectivity of $f_1=|f|:|\Sigma'|\to|\Sigma|$ does not suffice.)
Therefore, $\Sigma'$ is an {\em (algebraic) submonad\/} of $\Sigma$
iff $\Sigma'(n)\subset\Sigma(n)$ for all $n\geq0$, $\bu_{\Sigma'}=\bu_\Sigma$,
and the ``composition maps'' $\mu^{(k)}_n$ of $\Sigma'$ are restrictions
of those of~$\Sigma$.

We can easily compute intersections of algebraic submonads inside an
algebraic monad. For example, we can intersect $\Zinfty\subset\bbR$
and $\bbQ\subset\bbR$ (here we denote by $\Zinfty$ the algebraic monad
previously denoted by~$\Sigma_\infty$, and identify classical associative
rings~$R$ with corresponding algebraic monads~$\Sigma_R$), thus obtaining
the ``non-completed local ring at~$\infty$'':
\begin{equation}
\Zninfty:=\Zinfty\cap\bbQ
\end{equation}

\nxsubpoint (Classical associative rings as algebraic monads.)
In fact, the functor $R\mapsto\Sigma_R$, transforming a classical ring
into corresponding algebraic monad, is fully faithful, and we have
$\catMod R=\catMod{\Sigma_R}$. Therefore, we can safely identify
$R$ with $\Sigma_R$ and treat classical associative rings as
algebraic monads, i.e.\ (non-commutative) generalized rings.

\nxsubpoint (Underlying set of an algebraic monad.)
We denote by $\|\Sigma\|$ the
``total'' or ``graded'' underlying set of an algebraic monad~$\Sigma$, given by
\begin{equation}
\|\Sigma\|:=\bigsqcup_{n\geq0}\Sigma(n)
\end{equation}
This is a $\bbN_0$-graded set, i.e.\ a set together with a fixed decomposition
into a disjoint union indexed by~$\bbN_0$, or equivalently, a set $\|\Sigma\|$
together with a ``degree map'' $r:\|\Sigma\|\to\bbN_0$. In our case
we say that $r$ is the {\em arity map\/} of~$\Sigma$.

\nxsubpoint (Free algebraic monads.)
The underlying set functor $\Sigma\mapsto\|\Sigma\|$ from the category
of algebraic monads into the category $\catSets_{/\bbN_0}$ of $\bbN_0$-graded
sets admits a left adjoint, called the {\em free algebraic monad functor\/}
and denoted by $S\mapsto\langle S\rangle$ or $S\mapsto\Fempty\langle S\rangle$.
If $S$ is a finite set $\{f_1,\ldots,f_n\}$, consisting of elements $f_i$
of ``degrees'' or ``arities'' $r_i=r(f_i)$, we also write
$\Fempty\langle f_1^{[r_1]},\ldots,f_n^{[r_n]}\rangle$ or just
$\langle f_1^{[r_1]},\ldots, f_n^{[r_n]}\rangle$.

This $\Sigma=\Fempty\langle S\rangle$ is something like an algebra
of polynomials over~$\Fempty$ in non-commuting indeterminates from~$S$.
Notice, however, that the structure of $\Sigma$ depends on the
choice of arities of indeterminates from~$S$: if all of them are unary,
the result is indeed very much like an algebra of polynomials in
non-commuting variables, while non-unary free algebraic monads often
exhibit more complicated behaviour.

Existence of free algebraic monads $\Sigma=\Fempty\langle S\rangle$
is shown as follows. We explicitly define $\Sigma(X)$ as the set of
all terms, constructed from free variables $\{x\}$, $x\in X$,
with the aid of ``formal operations'' $f\in S$. In other words,
$\Sigma(X)$ is the set of all {\em terms\/} over~$X$, defined by
structural induction as follows:
\begin{itemize}
\item Any $\{x\}$, $x\in X$, is a term.
\item If $f\in S$ is a formal generator of arity~$r$, and
$t_1$, \dots, $t_r$ are terms, then $f\,t_1\,\ldots\, t_r$ is also a term.
\end{itemize}
Of course, this is just the formal construction of free algebraic systems,
borrowed from mathematical logic. When we don't want to be too formal,
we write $f(t_1,\ldots,t_n)$ instead of $f\,t_1\,\ldots\,t_n$,
and $t_1*t_2$ instead of $*\,t_1\,t_2$, if $*$ is a binary operation,
traditionally written in infix form.

\nxsubpoint (Generators of an algebraic monad.)
Given any (graded) subset $S\subset\|\Sigma\|$, we can always find the
smallest algebraic submonad $\Sigma'\subset\Sigma$, containing~$S$
(i.e.\ such that $\|\Sigma'\|\supset S$), for example by taking 
the intersection of all such algebraic submonads. Another description:
$\Sigma'$ is the image of the natural homomorphism $\Fempty\langle S\rangle
\to\Sigma$ from the free algebraic monad generated by~$S$ into~$\Sigma$,
induced by the embedding $S\to\|\Sigma\|$. Therefore, $\Sigma'(n)$ consists
of all operations which can be obtained by applying operations
from~$S$ to the basis elements $\{k\}_\stn$ finitely many times.

If $\Sigma'=\Sigma$, i.e.\ if $\Fempty\langle S\rangle\to\Sigma$ is surjective,
we say that {\em $S$ generates $\Sigma$.}

\nxsubpoint (Relations and strict quotients of algebraic monads.)
Notice that epimorphisms of algebraic monads needn't be surjective,
as illustrated by epimorphism $\bbZ\to\bbQ$. However, an algebraic monad 
homomorphism $f:\Sigma\to\Sigma'$ is a {\em strict\/} epimorphism
(i.e.\ coincides with the cokernel of its kernel pair) iff
all components $f_n:\Sigma(n)\to\Sigma'(n)$ are surjective. Therefore,
strict quotients $\Sigma\twoheadrightarrow\Sigma'$ of $\Sigma$ are
in one-to-one correspondence with compatible (algebraic) equivalence
relations $R\subset\Sigma\times\Sigma$. Any such equivalence relation
is completely determined by $\|R\|\subset\|\Sigma\times\Sigma\|$,
i.e.\ by the collection of equivalence relations $R(n)\subset\Sigma(n)\times
\Sigma(n)$. {\em Compatibility\/} of such a family of equivalence relations
with the algebraic monad structure of $\Sigma$ essentially means compatibility
with all maps $\Sigma(\phi):\Sigma(n)\to\Sigma(m)$ and $\mu^{(k)}_n:
\Sigma(k)\times\Sigma(n)^k\to\Sigma(n)$.

Given any ``system of equations'' $E\subset\|\Sigma\times\Sigma\|$,
we can construct the smallest compatible algebraic 
equivalence relation $R=\langle E\rangle\subset\Sigma\times\Sigma$ 
containing $E$, e.g.\ by taking the intersection of all such equivalence
relations. The corresponding strict quotient $\Sigma\twoheadrightarrow
\Sigma/\langle E\rangle$ is universal among all homomorphisms 
$\Sigma\stackrel f\to\Sigma'$, such that all relations from $f(E)$ are
satisfied in~$\Sigma'$.

\nxsubpoint (Presentations of algebraic monads.)
This is applicable in particular to subsets $E\subset\|\Fempty\langle S\rangle
\times\Fempty\langle S\rangle\|$, where $S\stackrel r\to\bbN_0$ is any
$\bbN_0$-graded set. Corresponding strict quotient $\Fempty\langle S\rangle
/\langle E\rangle$ will be denoted by $\Fempty\langle S\,|\,E\rangle$
or $\langle S\,|\,E\rangle$, and called the {\em free algebraic monad
generated by~$S$ modulo relations from~$E$.} When $S$ or $E$ are finite sets,
we can replace $S$ or $E$ in $\langle S\,|\,E\rangle$ by an explicit
list of generators (with arities) or relations. Furthermore,
when we write relations from~$E$, we often replace standard ``free variables''
$\{1\}$, $\{2\}$, \dots, by lowercase letters like $x$, $y$, \ldots,
thus writing for example $\Fpm=\langle 0^{[0]},$ $-^{[1]}\,|\,-0=0,$
$-(-x)=x\rangle$.

Whenever $\Sigma\cong\langle S\,|\,E\rangle$, we say that $(S,E)$
is a {\em presentation\/} of algebraic monad~$\Sigma$. Clearly, any
algebraic monad~$\Sigma$ admits a presentation: we just have to choose
any system of generators $S\subset\|\Sigma\|$ (e.g.\ $\|\Sigma\|$ itself),
and take any system of equations $E$ generating the kernel
$R\subset\Fempty\langle S\rangle\times\Fempty\langle S\rangle$ of the
canonical surjection $\Fempty\langle S\rangle\twoheadrightarrow\Sigma$
(e.g.\ $E=R$). If both $S$ and $E$ can be chosen to be finite,
we say that {\em $\Sigma$ is finitely presented\/} (absolutely, 
i.e.\ over $\Fempty$).

\nxsubpoint (Inductive limits of algebraic monads.)
Presentations of algebraic monads are very handy when we need to compute
inductive limits of algebraic monads. For example, the coproduct
$\Sigma\boxtimes\Sigma'$ of two algebraic monads 
$\Sigma=\langle S\,|\,E\rangle$ and $\Sigma'=\langle S'\,|\,E'\rangle$
can be computed as $\langle S, S\,|\,E, E'\rangle$. Since the cokernel
of a couple of morphisms can be computed as a strict quotient
modulo a suitable compatible algebraic equivalence relation, and
filtered inductive limits of algebraic monads can be computed componentwise,
similarly to the projective limits, we can conclude that
{\em arbitrary inductive and projective limits exist in the category
of algebraic monads.}

\nxsubpoint (Endomorphism monad of a set.)
Let $X$ be any set. We denote by $\END(X)$ its {\em endomorphism monad},
constructed as follows. We put $(\END(X))(n):=\Hom_\catSets(X^n,X)$,
$\bu:=\id_X\in(\END(X))(1)$, and define $\mu_n^{(k)}:\Hom(X^k,X)\times
\Hom(X^n,X)^k=\Hom(X^k,X)\times\Hom(X^n,X^k)\to\Hom(X^n,X)$ to be
the composition of maps: $(f,g)\mapsto g\circ f$. It is easy to check
that this indeed defines an algebraic monad $\END(X)$, which acts
on~$X$ in a canonical way. Furthermore, giving an action of an
algebraic monad~$\Sigma$ on set~$X$ turns out to be the same thing
as giving algebraic monad homomorphism $\Sigma\to\END(X)$,
since a $\Sigma$-action on~$X$ is essentially a rule that
transforms ``formal operations'' $t\in\Sigma(n)$ into their ``values''
${[t]}_X:X^n\to X$.

We transfer some classical terminology to our case. For example,
a $\Sigma$-module~$X$ is said to be {\em exact\/} or {\em faithful\/} 
if corresponding homomorphism $\Sigma\to\END(X)$ is injective.

\nxsubpoint (Presentations of algebraic monads and algebraic systems.)
Let $\Sigma=\langle S\,|\,E\rangle$ be a presentation of an algebraic
monad~$\Sigma$. Clearly, algebraic monad homomorphisms $\Sigma\to\Sigma'$
are in one-to-one correspondence to graded maps $\phi:S\to\|\Sigma'\|$, such
that the images under $\phi$ of the relations from~$E$ hold in $\Sigma'$.
Applying this to $\Sigma'=\END(X)$, we see that {\em a $\Sigma$-module
structure on a set~$X$ is the same thing as assignment of a map
${[f]}_X:X^r\to X$ to each generator $f\in S$ of arity $r=r(f)$,
in such a way that all relations from~$E$ hold.} Since $S$ is a set of
operations with some arities, and $E$ is a set of relations between expressions
involving operations from~$S$ and free variables, we see that $(S,E)$
is an algebraic system, and a $\Sigma$-module is just an ``algebra''
for this algebraic system $(S,E)$. In other words, ``algebras''
under an algebraic system $(S,E)$ are exactly the 
$\langle S\,|\,E\rangle$-modules. Since any algebraic monad admits a
presentation, the converse is also true: the category of modules
over an algebraic monad can be described as the category of ``algebras''
for some algebraic system.

In this way we see that algebraic systems are nothing else than presentations
of algebraic monads, and algebraic monads are just a convenient
categorical way of describing algebraic systems, i.e.\ {\em
theory of algebraic monads is a category-theoretic exposition of
universal algebra.}

The reader might think that algebraic monads are not so useful, because
they are just algebraic systems in another guise. If this be the case,
and the reader is still not convinced of the convenience of using algebraic
monads, we suggest to do following two exercises: 
\begin{itemize}
\item Describe the intersection of two algebraic
submonads of an algebraic monads in terms of their presentations.
\item Find a presentation of $\END(X)$, and prove the correspondence
between homomorphisms $\Sigma\to\END(X)$ and $\Sigma$-module structures
on~$X$, in terms of presentations of these algebraic monads.
\end{itemize}

\nxsubpoint (Examples: $\bbZ$ and $\Fone$.)
For example, we can write
\begin{multline}
\bbZ=\langle 0^{[0]}, -^{[1]}, +^{[2]}\,|\,
0+x=x,\,(x+y)+z=x+(y+z),\\x+y=y+x,\,x+(-x)=0\rangle
\end{multline}
Strictly speaking, we should write $++\{1\}\{2\}\{3\}=+\{1\}+\{2\}\{3\}$
instead of $(x+y)+z=x+(y+z)$,\dots

In any case the above equality actually means that the category
$\catMod\bbZ$ consists of sets~$X$, endowed by a constant $0_X$,
a unary operation $-_X:X\to X$, and a binary operation $+_X:X^2\to X$,
satisfying the above relations, i.e.\ abelian groups.

On the other hand, this algebraic monad $\bbZ=\Sigma_\bbZ$ can be
described explicitly, since $\bbZ(n)=\Sigma_\bbZ(n)=\bbZ^n$ is the
free $\bbZ$-module of rank~$n$.

Of course, we can remove the commutativity relation from the above presentation
of~$\bbZ$, and add relations $x+0=x$ and $(-x)+x=0$ instead, thus obtaining
an algebraic monad $\bbG$, such that $\catMod\bbG$ is the category of groups.
Sets $\bbG(n)$ are free groups in $n$ generators $\{1\}$, \dots, $\{n\}$,
and we have a natural homomorphism $\bbG\to\bbZ$. We'll see later that
algebraic monad $\bbZ$ is {\em commutative}, while $\bbG$ is not.

Another example: the ``field with one element'' $\Fone$ can be defined by
\begin{equation}
\Fone=\Fempty\langle 0^{[0]}\rangle\quad,
\end{equation}
i.e.\ it is the free algebraic monad generated by one constant. We see that
$\catMod\Fone$ is the category of sets~$X$ with one pointed element
$0_X\in X$, and $\Fone$-homomorphisms are just maps of pointed sets,
i.e.\ maps $f:X\to Y$, such that $f(0_X)=0_Y$.

On the other hand, we can describe sets $\Fone(n)$ explicitly:
any such set consists of $n+1$ elements, namely, $n$ basis elements
$\{k\}_\stn$, and a constant~$0$.

We hope that these examples illustrate how algebraic systems and
algebraic monads are related to each other.

\nxsubpoint (Modules over algebraic monads.)
After discussing the properties of algebraic monads themselves, and
their relation to algebraic systems, we consider some properties
of the category $\catMod\Sigma$ of modules over an algebraic monad~$\Sigma$.
These are actually properties of ``algebras'' under a fixed algebraic system;
we prefer to prove these properties again, using the theory of algebraic 
monads, instead of referring to well-known properties of such ``algebras'',
partly because we want to illustrate how the theory of modules over
algebraic monads works. Anyway, our proofs are quite short, compared
to those classical proofs given in terms of algebraic systems.

For example, we show that arbitrary projective and filtered inductive
limits of $\Sigma$-modules exist and can be computed on the level of
underlying sets. (Notice that the statement about filtered inductive
limits uses algebraicity of monad $\Sigma$ in an essential way.) 
Furthermore, arbitrary inductive limits of
$\Sigma$-modules exist. Monomorphisms of $\Sigma$-modules are just
injective $\Sigma$-homomorphisms, and {\em strict\/} epimorphisms
are the surjective $\Sigma$-homomorphisms. Since we have a notion
of free $\Sigma$-modules and free $\Sigma$-modules, we have reasonable
notions of a system of generators and of a presentation of
$\Sigma$-module~$M$, and we can define
{\em finitely generated\/} and {\em finitely presented\/} 
$\Sigma$-modules as the $\Sigma$-modules isomorphic to a strict quotient
of some free module $\Sigma(n)$, resp.\ to cokernel of a couple of
homomorphisms $\Sigma(m)\rightrightarrows\Sigma(n)$.

These notions seem to have all the usual properties we expect from them,
e.g.\ most of the properties of left modules over a classical associative
ring, with some notable exceptions:
\begin{itemize}
\item Not all epimorphisms are strict, i.e.\ surjective;
\item Direct sums (i.e.\ coproducts) $M\oplus N$ are not isomorphic to
direct products $M\times N$;
\item Direct sum of two monomorphisms needn't be a monomorphism.
\end{itemize}

\nxsubpoint (Scalar restriction and scalar extension.)
Given any algebraic monad homomorphism $\rho:\Sigma\to\Xi$, we obtain
a {\em scalar restriction functor\/} $\rho_*:\catMod\Xi\to\catMod\Sigma$,
such that $\Gamma_\Sigma\circ\rho_*=\Gamma_\Xi$, where $\Gamma_\Sigma$
and $\Gamma_\Xi$ are the forgetful functors into the category of sets.
Conversely (cf.~\ptref{sp:monhom.catfunct}), any functor $H:\catMod\Xi\to
\catMod\Sigma$, such that $\Gamma_\Sigma\circ H=\Gamma_\Xi$, equals
$\rho_*$ for some algebraic monad homomorphism $\rho:\Sigma\to\Xi$.
For example, the forgetful functor from the category $\catMod\bbZ$
of abelian groups into the category $\catMod\bbG$ of groups is
the scalar restriction functor with respect to canonical homomorphism
$\bbG\to\bbZ$.

Furthermore, scalar restriction functor $\rho_*$ admits a left adjoint
$\rho^*$, called the {\em scalar extension\/} or {\em base change\/}
functor. When no confusion can arise, we denote it by $\Xi\otimes_\Sigma-$,
$-\otimes_\Sigma\Xi$ or $-_{(\Xi)}$, even if it is not immediately related
to any tensor product. For example, if $G$ is a group, then
$\bbZ\otimes_\bbG G=G^{ab}=G/[G,G]$. Another example: $\bbR\otimes_{\Zinfty}-$
transforms a torsion-free $\Zinfty$-module $A=(A_{\Zinfty},A_\bbR)$
into real vector space $A_{(\bbR)}=A_\bbR$. Finally, scalar extension
functor with respect to $\Fempty\to\Sigma$ is just the
free $\Sigma$-module functor $L_\Sigma:\catSets\to\catMod\Sigma$.

\nxsubpoint (Flatness and unarity.)
Of course, these scalar restriction and extension functors retain almost all
their classical properties, e.g.\ scalar restriction commutes with 
arbitrary projective limits (hence is left exact), and scalar extension
commutes with arbitrary inductive limits (hence is right exact).
Notice, however, that scalar restriction $\rho_*:\catMod\Xi\to\catMod\Sigma$
needn't preserve direct sums (i.e.\ coproducts), as illustrated by
the case $\rho:\Fempty\to\bbZ$. In fact, {\em $\rho_*$ commutes with
finite direct sums iff $\rho_*$ is right exact iff $\rho$ is {\bf unary}}.
This notion of {\em unarity\/} is fundamental for the theory of generalized
rings, but completely absent from the classical theory of rings,
since {\em all homomorphisms of classical rings are unary.}

Similarly, we say that $\rho$ is flat, or that $\Xi$ is flat over~$\Sigma$,
if $\rho^*$ is (left) exact.

\nxsubpoint (Algebraic monads over topoi.)
Chapter~{\bf 4} treats the topos case as well. There is a natural approach
to algebraic monads over a topos~$\cE$, based on {\em algebraic inner
endofunctors.} However, the category of algebraic inner endofunctors
turns out to be canonically equivalent to $\cE^{\bbN}=\catFunct(\bbN,\cE)$,
i.e.\ an algebraic inner endofunctor~$\Sigma$ over~$\cE$ can be described as
a sequence $\{\Sigma(n)\}_{n\geq0}$ of objects of~$\cE$ together
with some transition morphisms $\Sigma(\phi):\Sigma(n)\to\Sigma(m)$,
defined for each $\phi:\stn\to\stm$. Furthermore, the ``elementary''
description of algebraic monads recalled in~\ptref{sp:intro.elem.algmon}
transfers verbatim to topos case, if we replace all sets by objects of 
topos~$\cE$.

Applying this to the topos $\cE$ of sheaves of sets over a site~$\cS$
or a topological space~$X$, we obtain a notion of a ``sheaf of algebraic 
monads'' over $\cS$ or~$X$: it is a collection of sheaves of sets 
$\{\Sigma(n)\}_{n\geq0}$, together with maps $\Sigma(\phi):
\Sigma(n)\to\Sigma(n)$, a unit section $\bu\in\Gamma(X,\Sigma(1))$
and ``evaluation maps'' $\mu_n^{(k)}:\Sigma(k)\times\Sigma(n)^k\to
\Sigma(n)$, such that for any open subset $U\subset X$ 
(or any object $U\in\Ob\cS$) $\Gamma(U,\Sigma)$ becomes an algebraic monad
(over $\catSets$). Sheaves of modules over such $\Sigma$ are defined 
similarly. We can also define generalized ringed spaces (ringed spaces
with a sheaf of generalized rings), their morphisms and so on.

\nxpointtoc{Commutativity}
Chapter~{\bf 5} is dedicated to another fundamental property we
require from our generalized rings, namely, {\em commutativity.} In
fact, we define a (commutative) generalized ring as a commutative
algebraic monad.  A classical ring is commutative as an algebraic
monad iff it is commutative in the classical sense, so our terminology
extends the classical one.

\nxsubpoint\label{sp:intro.def.comm} (Definition of commutativity.)
Given an algebraic monad~$\Sigma$, two operations $t\in\Sigma(n)$,
$t'\in\Sigma(m)$, a $\Sigma$-module~$X$, and a matrix $(x_{ij})_{1\leq i\leq m,
1\leq j\leq n}$, $x_{ij}\in X$, we can consider two elements $x'$, $x''\in X$,
defined as follows. Let $x_{i.}:=t(x_{i1},\ldots,x_{in})$ be the element 
of~$X$, obtained by applying $t$ to $i$-th row of $(x_{ij})$, and
$x_{.j}:=t'(x_{1j},\ldots,x_{mj})$ be the element, obtained by applying
$t'$ to $j$-th column of~$(x_{ij})$. Then $x':=t'(x_{1.},\ldots,x_{m.})$,
and $x'':=t(x_{.1},\ldots,x_{.n})$. Now we say that {\em $t$ and $t'$
commute on~$X$} if $x'=x''$ for any choice of $x_{ij}\in X$.
We say that {\em $t$ and $t'$ commute (in $\Sigma$)} if they commute
on any $\Sigma$-module~$X$.

Notice that $t$ and $t'$ commute iff the commutativity relation
is fulfilled for the ``universal matrix'' $x_{ij}:=\{(i-1)n+j\}$
in $X=\Sigma(nm)$. In this way commutativity of $t$ and $t'$ is equivalent
to one equation $x'=x''$ in $\Sigma(nm)$, which will be denoted by $[t,t']$.

Given any subset $S\subset\|\Sigma\|$, we can consider its 
{\em commutant\/} $S'\subset\|\Sigma\|$, i.e.\ the set of all operations 
of~$\Sigma$ commuting with all operations from~$S$. This subset $S'\subset
\|\Sigma\|$ turns out to be the underlying set of an algebraic submonad
of $\Sigma$, which will be also denoted by~$S'$. We have some classical
formulas like $S'''=S'$, $S''\supset S$ and so on.

Finally, the {\em center\/} of an algebraic monad $\Sigma$ is the
commutant of~$\|\Sigma\|$ itself, i.e.\ the set of all operations
commuting with each operation of~$\Sigma$. We say that $\Sigma$
is {\em commutative}, or a {\em generalized ring}, if
it coincides with its center, i.e.\ if any two operations of
$\Sigma$ commute. If $S\subset\|\Sigma\|$ is any system of generators
of~$\Sigma$, it suffices to require that any $t$, $t'\in S$ commute.

\nxsubpoint (Why this definition of commutativity?)
We would like to present some argument in favour of the above definition
of commutativity. Let $\Sigma$ be an algebraic monad,
$M$ and $N$ be two $\Sigma$-modules. Then the set 
$\Hom(M,N)=\Hom_{\catSets}(M,N)=N^M$ of {\em all\/}
maps $\phi:M\to N$ admits a natural $\Sigma$-module structure,
namely, the product $\Sigma$-structure on $N^M$. We can describe this
structure by saying that {\em operations $t\in\Sigma(n)$ act on maps
$\phi:M\to N$ pointwise:}
\begin{equation}\label{eq:intro.ptwise.act}
t(\phi_1,\ldots,\phi_n):x\mapsto t\bigl(\phi_1(x),\ldots,\phi_n(x)\bigr)
\end{equation}
So far we have used only $\Sigma$-module structure on~$N$, but
$M$ might be just an arbitrary set. Now if $M$ is a $\Sigma$-module,
we get a subset $\Hom_\Sigma(M,N)\subset\Hom(M,N)=N^M$, and we
might ask whether this subset is a $\Sigma$-submodule.

Looking at the case of modules over classical associative ring~$R$,
we see that $\Hom_R(M,N)\subset\Hom(M,N)=N^M$ is an $R$-submodule
for all $R$-modules $M$ and~$N$ if and only if $R$ is commutative,
and in this case the induced $R$-module structure on $\Hom_R(M,N)$
is its usual $R$-module structure, known from commutative algebra.

Therefore, it makes sense to say that an algebraic monad~$\Sigma$ is
commutative iff $\Hom_\Sigma(M,N)\subset\Hom(M,N)=N^M$ is a
$\Sigma$-submodule for all choices of $\Sigma$-modules $M$ and~$N$.
When we express this condition in terms of operations of~$\Sigma$,
taking~\eqref{eq:intro.ptwise.act} into account, we obtain exactly the
definition of commutativity given above in~\ptref{sp:intro.def.comm}.

\nxsubpoint (Consequences: bilinear maps, tensor products and inner Homs.)
Whenever $M$ and~$N$ are two modules over a generalized ring
(i.e.\ commutative algebraic monad) $\Sigma$, we denote by
$\iHom_\Sigma(M,N)$, or simply by $\Hom_\Sigma(M,N)$, the set
$\Hom_\Sigma(M,N)$, considered as a $\Sigma$-module with respect to
the structure induced from $\Hom(M,N)=N^M$.

Looking at the case of classical commutative rings, we see that we might
expect formulas like
\begin{equation}
\Hom_\Sigma(M\otimes_\Sigma N,P)\cong\Bilin_\Sigma(M,N;P)\cong
\Hom_\Sigma(M,\iHom_\Sigma(N,P))
\end{equation}
This expectation turns out to be correct, i.e.\ there is
a tensor product functor $\otimes_\Sigma$ on $\catMod\Sigma$ with the
above property, as well as a notion of {\em $\Sigma$-bilinear maps
$\Phi:M\times N\to P$}, fitting into the above formula.

This notion of $\Sigma$-bilinear map is very much like the classical one:
a map $\Phi:M\times N\to P$ is said to be {\em bilinear\/}
if the map $s_\Phi(x):N\to P$, $y\mapsto\Phi(x,y)$, is $\Sigma$-linear
(i.e.\ a $\Sigma$-homomorphism) for each $x\in M$, and
$d_\Phi(y):M\to P$, $x\mapsto\Phi(x,y)$, is $\Sigma$-linear for
each $y\in N$. Multilinear maps are defined similarly. 

In this way we obtain an ACU $\otimes$-structure on $\catMod\Sigma$,
which admits inner Homs $\iHom_\Sigma$. The unit object with respect to
this tensor product is $|\Sigma|=L_\Sigma(1)$, the free $\Sigma$-module
of rank one.

When we apply these definitions to a classical commutative ring~$R$,
we recover classical notions of $R$-bilinear maps, tensor product of
$R$-modules and so on.  On the other hand, when we apply them to
$\catMod\Zinfty$, we recover the definitions of Chapter~{\bf 2},
originally based on properties of maximal compact submonoids and
symmetric compact convex sets. This shows that the tensor product we
consider is indeed a very natural one.

\nxsubpoint (Algebras over a generalized ring.)
Let $\Lambda$ be a generalized ring. We can define a (non-commutative)
$\Lambda$-algebra 
as a generalized ring~$\Sigma$ together with a {\em central\/}
homomorphism $f:\Lambda\to\Sigma$ (i.e.\ all operations from the
image of~$f$ must lie in the center of~$\Sigma$). Of course, if
$\Sigma$ is commutative, all homomorphisms $\Lambda\to\Sigma$ are central,
i.e.\ {\em a (commutative) $\Lambda$-algebra is a generalized ring
$\Sigma$ with a homomorphism $\Lambda\to\Sigma$.}

\nxsubpoint (Variants of free algebras.)
We have several notions of a free algebra, depending on whether we impose
some commutativity relations or not:
\begin{itemize}
\item
Let $\Lambda$ be an algebraic monad, and consider the category of
$\Lambda$-monads, i.e.\ algebraic monad homomorphisms $\Lambda\to\Sigma$
with source~$\Lambda$. We have a (graded) underlying set functor
$\Sigma\to\|\Sigma\|$, which admits a left adjoint $S\mapsto\Lambda\langle
S\rangle$. Indeed, if $\Lambda=\langle T\,|\,E\rangle$ is any presentation
of $\Lambda$, then $\Lambda\langle S\rangle$ can be constructed as
$\langle S,T\,|\,E\rangle$.  Notice that this $\Lambda\langle S\rangle$
is something like a ``very non-commutative'' polynomial algebra over~$S$:
not only the indeterminates from~$S$ are not required to commute
between themselves, they are even not required to commute with
operations from~$\Lambda$!
\item Let $\Lambda$ be a generalized ring, i.e.\ a commutative algebraic
monad, and consider the category of all $\Lambda$-algebras, commutative
or not, i.e.\ central homomorphisms $\Lambda\to\Sigma$. The graded
underlying set functor $\Sigma\mapsto\|\Sigma\|$ still admits a left adjoint,
which will be denoted by $S\mapsto\Lambda\{S\}$. If 
$\Lambda=\langle T\,|\,E\rangle$, then $\Lambda\{S\}=\langle T,S\,|\,
E,[S,T]\rangle$, where $[S,T]$ denotes the set of all commutativity relations
$[s,t]$, $s\in S$, $t\in T$. In other words, now the indeterminates from~$S$
still don't commute among themselves, but at least commute with operations
from~$\Lambda$.
\item Finally, we can consider the category of commutative algebras
$\Lambda\to\Sigma$ over a generalized ring~$\Lambda$, and repeat
the above reasoning. We obtain free commutative algebras
$\Lambda[S]=\langle T,S\,|\,E,[S,T],[S,S]\rangle$, which are quite similar
to classical polynomial algebras, especially when all indeterminates
from~$S$ are unary.
\end{itemize}

\nxsubpoint (Presentations of $\Lambda$-algebras.)
In any of the three cases listed above we can impose some relations on
free algebras, thus obtaining the notion of a presentation of a
$\Lambda$-algebra. Consider for example the case of commutative 
$\Lambda$-algebras, with $\Lambda=\langle T\,|\,E\rangle$ any
generalized ring. Then, given any $\bbN_0$-graded set~$S$ and any
set of relations $R\subset\|\Lambda[S]\times\Lambda[S]\|$, we can
construct the strict quotient
\begin{equation}
\Lambda[S\,|\,R]=\Lambda[S]/\langle R\rangle=
\langle T,S\,|\,E,\,[S,T],\,[S,S],\,R\rangle
\end{equation}
This strict quotient is easily seen to be a commutative $\Lambda$-algebra,
and whenever we are given a $\Lambda$-algebra isomorphism
$\Sigma\cong\Lambda[S\,|\,R]$, we say that {\em $(S,R)$ is a 
presentation of (commutative) $\Lambda$-algebra~$\Sigma$.}

\nxsubpoint (Finitely generated/presented algebras.)
If $S$ can be chosen to be finite, we say that {\em $\Sigma$ is
a finitely generated $\Lambda$-algebra}, or {\em a $\Lambda$-algebra
of finite type}, and if both $S$ and~$R$ can be chosen to be finite,
we speak about {\em finitely presented\/} $\Lambda$-algebras.
These notions are quite similar to their classical counterparts.

\nxsubpoint (Pre-unary and unary algebras.)
Apart from notions of finite generation and presentation, which have well-known
counterparts in classical commutative algebra, we obtain notions of
{\em pre-unary\/} and {\em unary\/} $\Lambda$-algebras. Namely, we
say that a $\Lambda$-algebra $\Sigma$ is {\em pre-unary\/} if it admits
a system of unary generators, and {\em unary\/} if it admits a presentation
$\Sigma=\Lambda[S\,|\,R]$, 
consisting only of unary generators (i.e.\ $S\subset|\Sigma|=\Sigma(1)$) 
and unary relations (i.e.\ $R\subset|\Lambda[S]\times\Lambda[S]|$).

These notions are very important for the theory of generalized rings,
but they don't have non-trivial counterparts in classical algebra,
since any algebra over a classical ring is automatically unary.

\nxsubpoint (Tensor products, i.e.\ coproducts of $\Lambda$-algebras.)
Since any commutative $\Lambda$-algebra~$\Sigma_i$ admits some presentation 
$\Lambda[S_i\,|\,E_i]$ over~$\Lambda$, we can easily construct
{\em coproducts\/} in the category of commutative $\Lambda$-algebras,
called also {\em tensor products\/} of $\Lambda$-algebras, by
\begin{equation}
\Sigma_1\otimes_\Lambda\Sigma_2=\Lambda[S_1,S_2\,|\,E_1,E_2]
\end{equation}
Using existence of filtered inductive limits and arbitrary projective
limits of generalized rings, which actually coincide with those
computed in the category of algebraic monads, we conclude that
{\em arbitrary inductive and projective limits exist in the
category of generalized rings.}

\nxsubpoint (Unary algebras as algebras 
in $\otimes$-category $\catMod\Lambda$.)
Notice that we might try to generalize another classical definition
of a $\Lambda$-algebra and say that a $\Lambda$-algebra~$A$ is just
an algebra in ACU $\otimes$-category $\catMod\Lambda$, with respect
to our tensor product $\otimes=\otimes_\Lambda$, and define $A$-modules
accordingly.

However, the category of such algebras in $\catMod\Lambda$ turns out
to be {\em equivalent\/} to the category of {\em unary\/} $\Lambda$-algebras,
and we obtain same categories of modules under this correspondence.

This allows us to construct for example tensor, symmetric or
exterior algebras of a $\Lambda$-module~$M$ first by applying
the classical constructions inside $\otimes$-category $\catMod\Lambda$,
and then taking corresponding unary $\Lambda$-algebras. Monoid and
group algebras can be also constructed in this way; in particular,
they are unary.

Unary algebras have some other equivalent descriptions, which correspond
to some properties always fulfilled in the classical case.
For example, $\Sigma$ is a (commutative) unary $\Lambda$-algebra
iff the scalar restriction functor $\rho_*$ with respect to
$\Lambda\stackrel\rho\to\Sigma$ commutes with finite direct sums
iff it is exact iff it admits a right adjoint $\rho^!$ iff the
``projection formula'' $\rho_*(M\otimes_\Sigma\rho^*N)\cong\rho_*M
\otimes_\Lambda N$ holds.

\nxsubpoint (``Affine base change'' theorem.)
Another important result of Chapter~{\bf 5}, called by us the
{\em affine base change theorem} (cf.~\ptref{th:aff.base.change}),
essentially asserts the following. Suppose we are given two commutative
$\Lambda$-algebras $\Lambda'$ and $\Sigma$, and put 
$\Sigma':=\Lambda'\otimes_\Lambda\Sigma$. Then we can start from a
$\Sigma$-module~$M$, compute its scalar restriction to $\Lambda$,
and afterwards extend scalars to~$\Lambda'$. On the other hand,
we can first extend scalars to $\Sigma'$, and then restrict them
to $\Lambda'$. The affine base change theorem claims that
{\em the two $\Lambda'$-modules thus obtained are canonically isomorphic,
provided either $\Sigma$ is unary or $\Lambda'$ is flat over~$\Lambda$.}

We also consider an example, which shows that the statement is not true
without additional assumptions on $\Sigma$ or~$\Lambda'$.

\nxsubpoint
We have discussed at some length some consequences of commutativity for
categories of modules and algebras over a generalized ring~$\Sigma$.
Now we would like to discuss the commutativity relations themselves.
At first glance commutativity of two operations of arity $\geq2$
seems to be not too useful. We are going to demonstrate that
in fact such commutativity relations are very powerful, and that
they actually imply classical associativity, commutativity and
distributivity laws, thus providing a common generalization of all such
laws.

\nxsubpoint (Commutativity for operations of low arity.)
Let's start with some simple cases. Here $\Sigma$ is some fixed algebraic 
monad.
\begin{itemize}
\item Two constants $c$, $c'\in\Sigma(0)$ commute iff $c=c'$. In particular,
{\em a generalized ring contains at most one constant.}
\item Two unary operations $u$, $u'\in|\Sigma|=\Sigma(1)$ commute iff
$uu'\{1\}=u'u\{1\}$, i.e.\ $uu'=u'u$ in monoid $|\Sigma|$. This is classical
commutativity of multiplication.
\item An $n$-ary operation $t$ commutes with a constant $c$ iff
$t(c,c,\ldots,c)=c$. For example, a constant $c$ and a unary operation
$u$ commute iff $uc=c$ in $\Sigma(0)$. A binary operation $*$ and
constant $c$ commute iff $c*c=c$.
\item An $n$-ary operation $t$ commutes with unary $u$ iff
$u\cdot t(\{1\},\ldots,\{n\})=t(u\{1\},\ldots,u\{2\})$. For example,
if $t$ is a binary operation~$*$, this means $u(\{1\}*\{2\})=u\{1\}*u\{2\}$,
or $u(x*y)=ux*uy$, i.e.\ classical distributivity relation.
\item Two binary operations $+$ and $*$ commute iff $(x+y)*(z+w)=(x*z)+(y*w)$.
Notice that, while any operation of arity $\leq1$ automatically commutes
with itself, this is not true even for a binary operation~$*$, since
$(x*y)*(z*w)=(x*z)*(y*w)$ is already a non-trivial condition.
\end{itemize}

\nxsubpoint (Zero.)
We say that a {\em zero\/} of an algebraic monad~$\Sigma$ is just a
central constant $0\in\Sigma(0)$. Since it commutes with any other
constant~$c$, we must have $c=0$, i.e.\ a zero is automatically 
the only constant of~$\Sigma$. Furthermore, we have $u0=0$ for any
unary~$u$, $0*0=0$ for any binary~$*$, and $t(0,\ldots,0)=0$ for
an $n$-ary~$t$, i.e.\ all operations fix the zero.

Of course, a generalized ring either doesn't have any constants
(then we say that it is a {\em generalized ring without zero}),
or it has a zero (and then it is a {\em generalized ring with zero}).
We have a universal generalized ring with zero $\Fone=\Fempty[0^{[0]}]=
\Fempty\langle 0^{[0]}\rangle=\langle 0^{[0]}\rangle$, called
{\em the field with one element.} Furthermore, $\Fone$-algebras
are exactly the algebraic monads with a central constant, i.e.\ with zero.

\nxsubpoint (Addition.)
Let $\Sigma$ be an algebraic monad with zero~$0$. We say that a binary
operation $+^{[2]}\in\Sigma(2)$ is a {\em pseudoaddition\/}
iff $\{1\}+0=\{1\}=0+\{1\}$. This can be also written as
$\bu+0=\bu=0+\bu$ or $x+0=x=0+x$. An {\em addition\/} is simply a
central pseudoaddition. If $\Sigma$ admits an addition $+$, any
other pseudoaddition $*$ must coincide with $+$: indeed, we can put $y=z=0$
in commutativity relation $(x+y)*(z+w)=(x*z)+(y*z)$, and obtain
$x*w=(x+0)*(0+w)=(x*0)+(0*w)=x+w$.

So let $\Sigma$ be an algebraic monad with zero~$0$
and addition~$+$. Commutativity of $+$ and unary operations $u\in|\Sigma|$
means $u(x+y)=ux+uy$, i.e.\ distributivity. Now consider the commutativity
condition for $+$ and itself: $(x+y)+(z+w)=(x+z)+(y+w)$. Putting here $x=w=0$,
we obtain $y+z=z+y$, i.e.\ classical commutativity of addition.
Next, if we put $z=0$, we get $(x+y)+w=x+(y+w)$, i.e.\ classical associativity
of addition. Therefore, classical 
commutativity, associativity and distributivity laws are just consequences
of our ``generalized'' commutativity law.

\nxsubpoint (Semirings.)
Given any (classical) semiring~$R$, we can consider corresponding algebraic
monad~$\Sigma_R$, which will be an algebraic monad with addition
(hence also with zero). Conversely, whenever $\Sigma$ is an algebraic
monad with addition~$+$, we get an addition $+$ on monoid $|\Sigma|=\Sigma(1)$,
commuting with the action of all unary operations, i.e.\ $R:=|\Sigma|$ is 
a semiring. One can check that $\Sigma=\Sigma_R$, i.e.\ classical semirings
correspond exactly to algebraic monads with addition. For example,
semirings $\bbN_0$ and tropical numbers~$\bbT$ can be treated as
generalized rings. This reasoning also shows that any (generalized)
algebra~$\Sigma$ over a semiring~$\Lambda$ 
is also an algebraic monad with addition (since $\Lambda\to\Sigma$
is central, it maps the zero and addition of $\Lambda$ into zero and
addition of $\Sigma$), i.e.\ is also given by a classical semiring.
Furthermore, any such $\Lambda$-algebra $\Sigma$ is automatically unary.

Of course, we have a universal algebraic monad with addition, namely,
\begin{align}
\bbN_0=&\Fone[+^{[2]}\,|\,\bu+0=\bu=0+\bu]\\
=&\Fempty[0^{[0]},+^{[2]}\,|\,\bu+0=\bu=0+\bu]\\
=&\langle 0^{[0]}, +^{[2]}\,|\,[+,+],\,\bu+0=\bu=0+\bu\rangle
\end{align}

\nxsubpoint (Symmetry.)
A {\em symmetry\/} in an algebraic monad~$\Sigma$ is simply a central
unary operation $-$, such that $-^2=\bu$ in monoid $|\Sigma|$,
i.e.\ $-(-\bu)=\bu$, or $-(-x)=x$. Centrality of $-$ actually means
$-t(\{1\},\ldots,\{n\})=t(-\{1\},\ldots,-\{n\})$ for any $n$-ary operation~$t$.
In particular, $-c=c$ for any constant~$c$.

We have a universal algebraic monad with zero and symmetry, namely,
\begin{equation}
\Fpm=\Fone[-^{[1]}\,|\,-(-\bu)=\bu]=\Zinfty\cap\bbZ
\end{equation}

Using the notion of symmetry, we obtain a finite presentation of~$\bbZ$:
\begin{align}
\bbZ=&\Fpm[+^{[2]}\,|\,x+0=x=0+x,\,x+(-x)=0]\\
=&\Fempty[0^{[0]},-^{[1]},+^{[2]}\,|\,x+0=x=0+x,\,x+(-x)=0]
\end{align}

Reasoning as in the case of semirings, we see that classical rings
can be described as algebraic monads with addition~$+$ and symmetry~$-$,
compatible in the sense that $x+(-x)=0$. For example, if $\Lambda$
is a classical commutative ring, and $\Sigma$ is a generalized 
$\Lambda$-algebra, then algebraic monad $\Sigma$ admits a zero, a symmetry
and an addition (equal to the images of those of~$\Lambda$), compatible
in the above sense, hence $\Sigma$ is a classical ring, i.e.\ 
a classical $\Lambda$-algebra. In other words, we cannot obtain any ``new''
generalized algebras over a classical commutative ring.

\nxsubpoint (Examples of generalized rings.)
Let's list some generalized rings.
\begin{itemize}
\item Classical commutative rings, e.g.\ $\bbZ$, $\bbQ$, $\bbR$, $\bbC$,\dots
\item Classical commutative semirings: $\bbN_0$, $\bbR_{\geq0}$, $\bbT$,\dots
\item Algebraic submonads of generalized rings, e.g.\ $\Zinfty\subset\bbR$,
$\Zninfty\subset\bbQ$, or $A_N:=\Zninfty\cap\bbZ[1/N]\subset\bbQ$.
\item {\em Archimedian valuation rings\/} $V\subset K$, defined as follows.
Let $K$ be a classical field, equipped with an archimedian valuation $\fnorm$.
Then
\begin{equation}
V(n)=\{(\lambda_1,\ldots,\lambda_n)\in K^n\,:\,\sum_i|\lambda_i|\leq 1\}
\end{equation}
For example, $\Zninfty\subset\bbQ$ and $\Zinfty\subset\bbR$ are special
cases of this construction.
\item Strict quotients of generalized rings. Consider for example
the three-point $\Zinfty$-module~$Q$, obtained by identifying all inner
points of segment $|\Zinfty|=[-1,1]$ (we've denoted this~$Q$
by~$\Finfty$ in~\ptref{sp:intro.tors.zinfmod}), and denote by
$\Finfty$ the image of induced homomorphism $\Zinfty\to\END(Q)$.
Then $\Finfty$ is commutative, i.e.\ a generalized ring, being a strict
quotient of~$\Zinfty$. Furthermore, it is generated over~$\Fpm$ by
one binary operation $*$, the image of $(1/2)\{1\}+(1/2)\{2\}\in\Zinfty(2)$:
\begin{equation}
\Finfty=\Fpm[*^{[2]}\,|\,x*(-x)=0,\,x*x=x,\,x*y=y*x,\,(x*y)*z=x*(y*z)]
\end{equation}
\item Cyclotomic extensions of $\Fone$:
\begin{equation}
\bbF_{1^n}=\Fone[\zeta^{[1]}\,|\,\zeta^n=\bu]
\end{equation}
Modules over $\bbF_{1^n}$ are sets~$X$ with a marked element $0_X$
and a permutation $\zeta_X:X\to X$, such that $\zeta_X^n=\id_X$ and 
$\zeta_X(0_X)=0_X$. This is applicable to $\Fpm=\bbF_{1^2}$.
\item Affine generalized rings $\Aff_R\subset R$: take any classical 
commutative (semi)ring~$R$,
and put $\Aff_R(n):=\{\lambda_1\{1\}+\cdots+\lambda_n\{n\}\in R^n\,|\,
\sum_i\lambda_i=1\}$. Then $\catMod{\Aff_R}$ is just the category 
of affine $R$-modules.
\item Generalized ring $\Delta:=\Aff_\bbR\cap\bbR_{\geq0}\subset\bbR$.
Set $\Delta(n)$ is the standard simplex in $\bbR^n$, and $\catMod\Delta$
consists of abstract convex sets.
\end{itemize}

\nxsubpoint (Some interesting tensor products.)
We compute some interesting tensor products of generalized rings: 
$\bbZ\otimes\bbZ=\bbZ$, 
and similarly $\bbZ\otimes\Zinfty=\bbR$,
$\bbZ\otimes\Zninfty=\bbQ$, if the tensor products are computed over
$\Fpm$, $\Fone$ or $\Fempty$. Furthermore, $\Finfty\otimes_{\Fpm}\Finfty=
\Finfty$ and $\Zninfty\otimes_{\Fpm}\Zninfty=\Zninfty$.

\nxsubpoint (Tensor and symmetric powers.)
Given any module $M$ over any generalized ring $\Lambda$, we can consider
its tensor and symmetric powers $T^n(M)=M^{\otimes n}$ and 
$S^n(M)=T^n(M)/\gS_n$. They seem to retain all their usual properties,
e.g.\ they commute with scalar extension, and a symmetric power of a free
$\Lambda$-module~$M$ is again free, 
with appropriate products of base elements of~$M$ as a basis.

\nxsubpoint (Exterior powers, alternating generalized rings and determinants.)
On the other hand, exterior powers $\wedge^n(M)$ are more tricky to define.
First of all, $M$ has to be a module over 
a commutative $\Fpm$-algebra~$\Sigma$, i.e.\ generalized ring $\Sigma$
must admit a zero and a symmetry. Even if this condition is fulfilled,
these exterior powers $\wedge^n(M)$ don't retain some of their classical
properties, e.g.\ exterior powers of a free module are not necessarily free,
and in particular $\wedge^n L_\Sigma(n)\not\cong|\Sigma|$. 
In order to deal with these problems we introduce and study {\em
alternating\/} generalized rings. Over them exterior powers of a free module
are still free with the ``correct'' basis, e.g.\ if $M$ is free with
basis $e_1$, \dots, $e_n$, then $\wedge^n(M)$ is freely generated by
$e_1\wedge\ldots\wedge e_n$, so we write $\det M:=\wedge^n(M)$, and
use it to define {\em determinants\/} $\det u\in|\Sigma|$ of endomorphisms
$u\in\End_\Sigma(M)$. Then $\det(uv)=\det(u)\det(v)$, $\det(\id_M)=\bu$,
so any automorphism has invertible determinant. However, 
invertibility of $\det(u)$ doesn't imply bijectivity of~$u$ in general case.
We need to introduce some additional conditions $(DET)$, which assure that
invertibility $\det(u)$ indeed implies invertibility of~$u$.

Our general impression of this theory of exterior powers and determinants
is that one should avoid using exterior powers in the generalized ring 
context, and use symmetric powers instead whenever possible. Then
we don't have to impose any complicated conditions on the base generalized
ring~$\Lambda$.

\nxsubpoint (Topos case.)
Of course, definitions of commutativity and alternativity of algebraic monads
can be transferred to the topos case, so we have a reasonable notion
of a sheaf of generalized rings over a topological space~$X$ or a site~$\cS$.
Furthermore, usually we consider only generalized commutatively ringed
spaces or sites, i.e.\ require the structural sheaf~$\sO$ of algebraic monads
to be commutative. Then we get a tensor product $\otimes_\sO$ of
(sheaves of) $\sO$-modules, inner Homs $\iHom_\sO(\sF,\sG)$,
functors $i_!$, $i^*$, $i_*$, related to an open immersion $i:U\to X$,
and even the classical description of the category of $\sO$-modules
over~$X$ in terms of categories of modules over an open subset $U\subset X$
and its closed complement $Y=X-U$.

\nxpointtoc{Localization and generalized schemes}
Chapter~{\bf 6} is dedicated to {\em unary localization},
{\em localization theories}, {\em spectra of generalized rings}
and {\em generalized schemes}. We transfer some of results 
of EGA~I and~II to our situation, e.g.\ notions of morphisms of finite
presentation or of finite type, affine and projective morphisms, and
show existence of
fibered products of generalized schemes. We discuss some examples,
e.g.\ projective spaces and projective bundles.

While this ``generalized algebraic geometry'', i.e.\ theory of 
generalized schemes, seems to be quite similar to its classical counterpart,
there is a considerable distinction related to closed subschemes.
Namely, a strict epimorphism (i.e.\ surjective homomorphism) of generalized
rings $f:A\to B$ doesn't define a {\em closed\/} map
${}^af:\Spec B\to\Spec A$. In fact, even the diagonal of $\bbA^1_{\Fone}=
\Spec\Fone[T^{[1]}]$ over $\Spec\Fone$ is not closed. 

We partially avoid this problem by defining a ``closed'' immersion as an
affine morphism $i:Y\to X$, such that $\sO_X\to i_*\sO_Y$ is a strict
epimorphism, a ``closed'' (generalized) subscheme as a subobject in the
category of generalized schemes, given by a ``closed'' immersion,
and a ``separated'' morphism $X\to S$ by requiring 
$\Delta_{X/S}:X\to X\times_SX$ to be a ``closed'' immersion. In this way
affine and projective spaces over $\Fone$ are ``separated'',
even if their diagonals are not topologically closed.

\nxsubpoint (Multiplicative systems and unary localizations.)
Let $A$ be a generalized ring, $M$ be an $A$-module.
Recall that $|A|$ has a natural (commutative) monoid structure,
and the underlying set of $M$ admits a canonical monoid action of $|A|$.
Now a {\em multiplicative system\/} $S$ in $A$ is simply a submonoid
of $|A|$. Given any subset $S\subset|A|$, we can consider the
multiplicative system (i.e.\ submonoid) $\langle S\rangle\subset|A|$,
generated by~$S$. Next, we can define the {\em localization $A[S^{-1}]$ of $A$}
as the universal commutative $A$-algebra $i_A^S:A\to A[S^{-1}]$, such that
all elements from~$S$ become invertible in $A[S^{-1}]$. Similarly,
we define the localization $M[S^{-1}]$ of $A$-module $M$ by considering the
universal $A$-module homomorphism $i_M^S:M\to M[S^{-1}]$, such that
all $s\in S$ act on $M[S^{-1}]$ by bijections.

When $S$ is a multiplicative system, we write $S^{-1}A$ and $S^{-1}M$ instead
of $A[S^{-1}]$ and $M[S^{-1}]$. We can construct $S^{-1}M$ in the classical
fashion, as the set of couples $(x,s)\in M\times S$, modulo 
equivalence relation $(x,s)\sim(y,t)$ iff $utx=usy$ for some $u\in S$.
Sets $(S^{-1}A)(n)=S^{-1}(A(n))$ admit a similar description. Next,
if $S$ is not a multiplicative system, we can always replace $S$ by 
$\langle S\rangle$, and apply the above constructions.

In this way (unary) localizations always exist and have almost all of
the classical properties. For example, $M[S^{-1}]$ admits a natural
$A[S^{-1}]$-module structure, and coincides with the scalar extension
$M\otimes_AA[S^{-1}]$. Furthermore, $A[S^{-1}]$ is always a flat $A$-algebra,
i.e.\ $M\mapsto M[S^{-1}]$ is exact,
and $\catMod{A[S^{-1}]}$ is identified by the scalar restriction functor
with the full subcategory of $\catMod A$, consisting of all $A$-modules~$N$,
such that all $[s]_N:N\to N$, $s\in S$, are bijections.

The reason why (unary) 
localization theory extends so nicely to the generalized ring
context is that almost all constructions and proofs in this theory
are ``multiplicative'', i.e.\ involve only the multiplication of $|A|$
and the action of $|A|$ on~$M$.

\nxsubpoint (Ideals.)
Now consider the {\em ideals\/} $\ga$ of~$A$, i.e.\ the $A$-submodules
$\ga\subset|A|$. They constitute a lattice with respect to inclusion,
with $\inf(\ga,\gb)$ given by the intersection $\ga\cap\gb$, and
$\sup(\ga,\gb)$ equal to $\ga+\gb=(\ga,\gb)=\Im(\ga\oplus\gb\to|A|)$.
The largest element of this lattice is the unit ideal $(1)=(\bu)=|A|$,
and the smallest element is the {\em initial ideal\/} $\emptyset_A=L_A(0)
\subset L_A(1)=|A|$. If $A$ admits a zero~$0$, then $\emptyset_A=(0)$
is called the {\em zero ideal}; if not, we say that $\emptyset_A$
is the {\em empty ideal}.

We have also a notion of {\em principal ideal\/} $(a)=a|A|\subset|A|$,
for any $a\in|A|$, and of {\em product of ideals\/} $\ga\gb$,
defined for example as $\Im(\ga\otimes_A\gb\to A)$; this is the
ideal generated by all products $ab$, $a\in\ga$, $b\in\gb$.

So far ideals in $A$ have all their properties known from commutative algebra.
However, when we try to compute strict quotients $A/\ga$ (for example first
as algebras in $\catMod A$, and then as unary algebras over~$A$), the
picture is not so nice any longer. For example, if $A$ admits a zero
(this is actually a prerequisite for defining $A/\ga$),
then the preimage of zero under $A\to A/\ga$ might be strictly larger that
$\ga$, i.e.\ {\em quotient algebras and modules lose some of their
classical properties.}

\nxsubpoint (Prime spectra.)
Next, we define a {\em prime ideal\/} $\gp\subset|A|$ as an ideal,
such that $|A|-\gp$ is a multiplicative system. We denote by
$\Spec A$ or $\Spec^pA$ the {\em prime spectrum\/} of~$A$, defined
as the set of all prime ideals in~$|A|$. For any subset $S\subset|A|$
we denote by $V(S)$ the set of prime ideals $\gp\supset S$. Clearly,
$V(S)=V(\ga)$, where $\ga=(S)$ is the ideal generated by~$S$,
and one can show that $V((1))=\emptyset$, $V(\emptyset)=\Spec A$,
$V(\bigcup_\alpha S_\alpha)=\bigcap_\alpha V(S_\alpha)$ and
$V(\ga\gb)=V(\ga\cap\gb)=V(\ga)\cup V(\gb)$. This means that these
$V(S)$ are the closed sets for some topology on $\Spec A$, called the
{\em Zariski topology}. Principal open subsets $D(f):=\Spec A-V(\{f\})=
\{\gp\,:\,\gp\not\ni f\}$ constitute a base of Zariski topology
on~$\Spec A$, and any $A$-module~$M$ defines a {\em sheaf\/} $\tilde M$
on $\Spec A$, characterized by $\Gamma(D(f),M)=M_f=S_f^{-1}M$, where
$S_f$ is the multiplicative system $\{1,f,f^2,\ldots\}$. 

Furthermore, $\tilde A(n)=\widetilde{A(n)}$ also turn out to be sheaves,
and their collection defines a sheaf of generalized rings~$\tilde A=
\sO_{\Spec A}$, and all $\tilde M$ are $\sO_{\Spec A}$-modules.
In this way $\Spec A$ becomes a generalized ringed space, and
any $A$-module $M$ defines a sheaf of modules~$\tilde M$ over~$\Spec A$;
such sheaves of modules are called {\em quasicoherent.}

This construction is contravariant in~$A$ in the usual manner,
i.e.\ any homomorphism $\phi:A\to B$ defines a continuous map
${}^a\phi:\Spec B\to\Spec A$, given by $\gp\mapsto\phi^{-1}(\gp)$.
This map naturally extends to a morphism of generalized ringed spaces,
having all the elementary properties of EGA~I, e.g.\ direct images
of quasicoherent sheaves correspond to scalar restriction with respect
to~$\phi$, and pullbacks of quasicoherent sheaves are given by scalar
extension.

In this way we obtain a theory of prime spectra of generalized rings,
i.e.\ {\em affine generalized schemes\/}, having almost all elementary
properties of EGA~I. Of course, when we apply our theory to a
classical ring~$A$, we recover classical prime spectra $\Spec A$.

\nxsubpoint (Drawbacks of prime spectra.)
Despite all their advantages and similarity to classical case,
prime spectra don't behave very nicely when applied to generalized
rings without sufficiently many unary operations. For example,
$|\Aff_\bbZ|=\{\bu\}$, and the only prime ideal of $\Aff_\bbZ$ is the
initial ideal $\emptyset_{\Aff_\bbZ}$, i.e.\ $\Spec\Aff_\bbZ$ is
a one-element set. On the other hand, one would rather expect
$\Spec\Aff_\bbZ$ to be similar to $\Spec\bbZ$, but with a different
(smaller) structural sheaf.

\nxsubpoint (Solution: open pseudolocalizations and localization theories.)
We try to deal with this problem as follows. We say that a generalized
ring homomorphism $\rho:A\to A'$ is a {\em pseudolocalization} if
$A'$ is flat over~$A$, and $\rho_*:\catMod{A'}\to\catMod A$ is fully faithful.
For example, $A\to S^{-1}A$ is a pseudolocalization for any multiplicative
system~$S$. Next, we say that $\rho:A\to A'$ is an 
{\em open pseudolocalization} if it is a finitely presented pseudolocalization.
For example, $A\to A_f=A[f^{-1}]$ is an open pseudolocalization for any 
$f\in|A|$. There might be more open pseudolocalizations even in the
classical case; notice that $\Aff_{\bbZ}\to\Aff_{\bbZ[f^{-1}]}$ is an
open pseudolocalization for any $f\in\bbZ$.

After that, we define a {\em localization theory $\cT^?$} as follows.
For any generalized ring~$A$ we must have a collection $\cT^?_A$ of open 
pseudolocalizations of~$A$, considered as a full subcategory of the
category of $A$-algebras. We require all unary open localizations
$A\to A_f$ to be contained in $\cT^?_A$, and we impose some natural
restrictions, e.g.\ stability of open pseudolocalizations from $\cT^?$
under composition and base change.

We have a {\em minimal\/} localization theory, namely, the
{\em unary localization theory\/} $\cT^u$, consisting only of 
unary localizations $A\to A_f$. We have a maximal localization theory,
called the {\em total localization theory\/} $\cT^t$, which contains
all open pseudolocalizations. Any other localization theory is
somewhere in between: $\cT^u\subset\cT^?\subset\cT^t$.

\nxsubpoint (Spectra with respect to a localization theory.)
Once we fix a localization theory $\cT^?$, we can construct the
{\em $\cT^?$-spectrum\/} $\Spec^?A$ of a generalized ring~$A$ as follows.
We denote by $\cS^?_A$ the category, opposite to the category $\cT^?_A$
of open pseudolocalizations $A\to A'$ from $\cT^?$ with source~$A$.
Since all pseudolocalizations are epimorphisms of generalized rings,
all morphisms in $\cS^?_A$ are monomorphisms, and there is at most
one morphism between any two objects of~$\cS^?_A$. In other words,
$\cS^?_A$ is equivalent to category defined by an ordered set.

For example, if we choose the unary localization theory $\cT^u$,
then $\cS^?_A$ can be thought of as the ordered set of principal open
subsets $D(f)\subset\Spec A$, since $A\to A_g$ factorizes through
$A\to A_f$ iff $D(g)\subset D(f)$. In this way objects of $\cS^?_A$
are something like ``abstract principal open subsets'' of $\Spec^?A$.

For any $A$-module~$M$ we define a presheaf of sets $\tilde M:(\cS^?_A)^0\to
\catSets$, given by $\tilde M:(A\to A')\mapsto M_{(A')}$. Furthermore,
we have a ``structural presheaf'' of generalized rings 
$\sO:(A\to A')\mapsto A'$, and $\tilde M$ admits a natural $\sO$-module
structure.

Since $\cT^?$ is closed under pushouts and composition,
$\cS^?_A$ has fibered products, which can be thought of as
``intersections of abstract principal open subsets''. 
We introduce on~$\cS^?_A$ 
the finest Grothendieck topology, such that all $\tilde M$
become sheaves. In other words, a family of morphisms
$(U_\alpha\to U)$ is a cover in~$\cS^?_A$ iff the following diagram
\begin{equation}
\tilde M(U)\longrightarrow\prod_\alpha\tilde M(U_\alpha)\rightrightarrows
\prod_{\alpha,\beta}\tilde M(U_\alpha\cap U_\beta)
\end{equation}
is left exact for any $A$-module~$M$.

For example, if we work with $\cT^u$, then $\{D(f_\alpha)\to D(f)\}_\alpha$
is a cover in $\cS^u_A$ (where $D(f)$ denotes the object of $\cS^u_A$
corresponding to $A\to A_f$) iff $D(f)=\bigcup D(f_\alpha)$ in $\Spec^p A$.

In general we obtain a Grothendieck topology on $\cS^?_A$, and define
the {\em strong spectrum $\Spec^?_s A$} as the corresponding topos,
i.e.\ the category of sheaves of sets on $\cS^?_A$. It is a generalized
ringed topos with structural sheaf $\sO$, and all $\tilde M$ are
(sheaves of) $\sO$-modules.

\nxsubpoint (Finite spectra.)
In fact, it turns out to be more convenient to use {\em finite spectra
$\Spec^?A=\Spec^?_fA$}, defined by the Grothendieck topology on~$\cS^?_A$,
generated by all {\em finite\/} covers $\{U_\alpha\to U\}$ with
the same property as above. For example, these finite spectra
$\Spec^?A$ are always quasicompact, as well as their ``abstract principal
open subsets'' (i.e.\ objects of $\cS^?_A$), hence they are coherent,
and in particular admit enough points by a general result of SGA~4.
Since they are obviously generated by open objects (i.e.\ subobjects
of the final object), namely, by the objects $U\in\Ob\cS^?_A$, identified
with corresponding representable sheaves $h_U:U'\mapsto\Hom(U',U)$,
we can conclude that {\em topos $\Spec^?A$ is given by a uniquely 
determined topological space}, which will be denoted by $\Spec^?A$ as
well.

In this way we obtain a generalized ringed space $\Spec^?A$, and can
forget about sites and topoi, used during its construction.
Unary spectrum $\Spec^uA$ always coincides with the prime spectrum $\Spec^pA$
constructed before, so this notion of $\cT^?$-spectra indeed
generalizes that of prime spectra.

Notice that, while classical rings~$A$ 
sometimes admit open pseudolocalizations not of the form $A\to A_f$
(in fact, $A\to A'$ is an open pseudolocalization of classical rings
iff $\Spec A'\to\Spec A$ is an open embedding of classical schemes),
i.e.\ categories $\cS^?_A$ might be distinct from $\cS^u_A$, arising
topoi $\Spec^?A$ are always equivalent, and corresponding topological
spaces are canonically homeomorphic, i.e.\ we can write $\Spec^?A=\Spec^uA=
\Spec A$ for a classical ring~$A$, regardless of the localization theory
chosen.

This property is due to the fact that $\Spec^uA'\to\Spec^uA$ is
an open map for any open pseudolocalization of classical rings
$A\to A'$, since open pseudolocalizations are flat and finitely presented.
Furthermore, if we have some localization theory~$\cT^?$,
such that $\Spec^?A'\to\Spec^?A$ is open 
for any open pseudolocalization $A\to A'$,
then generalized ringed space $\Spec^?A$ is automatically isomorphic
to~$\Spec^t A$. In this way, any ``reasonable'' localization theory~$\cT^?$,
having the property that $\Spec^? B\to\Spec^? A$ is an open map for any
finitely presented flat $A$-algebra~$B$, yields the same spectra as
the total localization theory, i.e.\ the {\em total localization theory\/}
$\cT^t$ is a very natural candidate for constructing ``correct'' spectra.

\nxsubpoint (Unary spectra vs.\ total spectra.)
We see that unary spectra $\Spec^u\negthinspace A=\Spec^p\negthinspace A$ 
have the advantage
of admitting a direct description in terms of prime ideals, similar
to the classical case, while total spectra $\Spec^tA$ can be expected
to have better properties for non-classical generalized rings, yielding
the same results for classical rings. For example, $\Spec^t\bbZ=\Spec^u\bbZ=
\Spec\bbZ$, $\Spec^u\Aff_\bbZ$ consists only of one point, but
$\Spec^t\Aff_\bbZ$ must be much larger, since it contains ``abstract
principal open subsets'' corresponding to open pseudolocalizations
$\Aff_{\bbZ}\to\Aff_{\bbZ[f^{-1}]}$.

Nevertheless, we lack a direct description of points of the total
spectrum $\Spec^tA$, so the reader might find more convenient
to consider only unary spectra $\Spec^u$, and skip remarks concerning
other localization theories while reading Chapter~{\bf 6} for the first time.

\nxsubpoint (Generalized schemes.)
Once we have a notion of spectra $\Spec A=\Spec^?A$ (for some
fixed localization theory $\cT^?$), we can define an {\em affine
generalized scheme} as a generalized ringed space isomorphic to some
$\Spec A$, and a {\em generalized scheme\/} as a generalized ringed space
$(X,\sO_X)$, which admits an open cover by affine generalized schemes.

Morphisms of generalized schemes can be defined either as local
morphisms of generalized locally ringed spaces (notice that the
notions of local generalized rings, and local homomorphisms of such,
depend on the choice of $\cT^?$), or as generalized ringed space
morphisms $f=(f,\theta):(Y,\sO_Y)\to(X,\sO_X)$, such that for any
open affine $\Spec B=V\subset Y$ and $\Spec A=U\subset X$,
such that $f(V)\subset U$, the restriction $f|_V:V\to U$ is induced by
a generalized ring homomorphism $A\to B$.

\nxsubpoint (Basic properties of generalized schemes.)
We show that most basic properties of schemes and sheaves of modules over
schemes, known from EGA~I, transfer to our case without much modifications.
For example, fibered products of generalized schemes exist and can
be constructed in the classical fashion, by considering appropriate 
affine covers and putting $\Spec A\times_{\Spec C}\Spec B:=
\Spec(A\otimes_CB)$ for affine generalized schemes. Quasicoherent sheaves
(of $\sO_X$-modules) constitute a full subcategory of $\catMod{\sO_X}$,
closed under finite projective and arbitrary inductive limits.
Quasicoherence is a local property (this is not so trivial as it might seem).
Sheaves of $\sO_X$-modules of finite type and of finite presentation
can be defined in a natural way, and have some reasonable properties,
e.g.\ any finitely presented $\sF$ is quasicoherent, and a quasicoherent
$\sF$ is finitely presented iff for any open affine $U\subset X$
the $\Gamma(U,\sO)$-module $\Gamma(U,\sF)$ is finitely presented iff
this condition holds for all $U$ from an affine open cover of~$X$.

We also have some classical filtered inductive limit properties,
for example, $\Hom_{\sO_X}(\sF,-)$ commutes with filtered inductive limits
of quasicoherent $\sO_X$-modules, provided $\sF$ is a finitely presented
$\sO_X$-module, and $X$ is quasicompact and quasiseparated.

Now let us mention some notable distinctions from the classical case. 

\nxsubpoint (Quasicoherent $\sO_X$-algebras, affine morphisms, unarity.)
We have a notion of a {\em quasicoherent $\sO_X$-algebra~$\sA$.} It is
more tricky than one might expect: $\sA$ has to be a sheaf 
of generalized rings on~$X$, equipped with a homomorphism $\sO_X\to\sA$,
such that all components $\sA(n)$ are quasicoherent $\sO_X$-modules.
Then we can construct {\em relative spectra\/} $\Spec\sA\to X$,
and define {\em affine morphisms\/} accordingly. We obtain a notion
of {\em unary affine morphisms}, corresponding to unary quasicoherent
algebras. This notion has no counterpart in classical algebraic geometry,
where all affine morphisms are unary.

\nxsubpoint (``Closed'' immersions are not closed.)
Another distinction is that generalized scheme morphisms
$\Spec B\to\Spec A$, induced by strict epimorphisms (i.e.\ surjective
generalized ring homomorphisms) $A\twoheadrightarrow B$ are not necessarily
closed. Roughly speaking, if $A$ admits a zero~$0$, and $B$ is obtained
from $A$ by imposing several relations of the form $f_i=0$, $f_i\in|A|$,
then $\Spec B$ can be identified with the closed subset $V(f_1,\ldots,f_n)$
of $\Spec A$, the complement of the union of principal open subsets $D(f_i)$.
However, in general $B$ is obtained by imposing
relations $f_i=g_i$, and we cannot replace such relations with $f_i-g_i=0$,
in contrast with the classical case.

We cope with this problem by defining ``closed'' immersions
$i:Y\to X$ as affine morphisms, defined by quasicoherent strict quotients
of $\sO_X$, and ``closed'' subschemes $Y\subset X$ as subobjects of~$X$
in the category of generalized schemes, defined by a ``closed'' immersion.
We say that $X\to S$ is ``separated'' if $\Delta_{X/S}:X\to X\times_SX$
is ``closed''. Surprisingly, these notions retain some of 
their classical properties,
e.g.\ any affine or projective morphism is automatically ``separated'',
even if its diagonal is not closed in the topological sense. 

\nxsubpoint (Absence of residue fields, and non-injective monomorphisms.)
Another notable distinction is the absence of a reasonable theory of
``generalized residue fields'' $\kappa(x)$ of points $x\in X$ of a generalized
scheme~$X$. Even when such ``residue fields'' can be constructed
(e.g.\ $\Finfty$ is a natural candidate for the residue field of
$\Spec\Zninfty$ at~$\infty$), they lack their usual properties,
and cannot be used to transfer classical proofs to our situation.

For example, one of the classical application of residue fields is the
statement about surjectivity of $|X\times_SY|\to|X|\times_{|S|}|Y|$,
where $X$ and $Y$ are $S$-schemes, and $|Z|$ denotes the underlying set of~$Z$.
In particular, a monomorphism of classical schemes is always injective on
underlying sets.

These properties fail for generalized schemes. For example, $\Spec\bbZ\to
\Spec\Fone$ is a monomorphism of generalized schemes
(since $\bbZ\otimes_{\Fone}\bbZ=\bbZ$, i.e.\ $\Spec\bbZ\times_{\Spec\Fone}
\Spec\bbZ=\Spec\bbZ$), but it is obviously not injective, since
$\Spec^p\Fone$ consists of one point.

\nxsubpoint (Classical schemes as generalized schemes.)
This statement about monomorphicity of $\Spec\bbZ\to\Spec\Fone$ 
and $\Spec\bbZ\to\Spec\Fempty$, while discarding all hopes to
obtain a non-trivial ``square of $\Spec\bbZ$'', has some benign
consequences as well. Namely, the category of classical schemes,
which can be also described as generalized schemes over $\Spec\bbZ$,
turns out to constitute a full subcategory of the category of all
generalized schemes. In this way we can safely treat classical schemes
as generalized schemes, since no new morphisms arise. Furthermore,
product of two classical schemes $X$ and~$Y$, computed over $\Spec\Fone$
or $\Spec\Fempty$, coincides with their product over $\Spec\bbZ$,
so we can write $X\times Y$ without any risk of confusion.

\nxsubpoint\label{sp:intro.graded.genrg} 
(Projective geometry and graded generalized rings.)
A considerable part of Chapter~{\bf 6} is dedicated to the study
of {\em graded generalized rings\/} and {\em modules},
{\em (relative) projective spectra}, and related notions, i.e.\ 
to ``generalized projective geometry''. We don't want to explain this
theory in any detail right now, since it is quite parallel to its classical
counterpart, with the usual exception of notions of unarity and
pre-unarity. We would like, however, to present a motivation for our
definition of a graded generalized ring.

Suppose we are given a generalized scheme~$X$, a commutative monoid
$\Delta$, and a homomorphism $\phi:\Delta\to\Pic(X)$. Let's denote
by $\sO[\lambda]$ the line bundle corresponding to $\lambda\in\Delta$,
and suppose we have some canonical isomorphisms $\sO[\lambda+\delta]=
\sO[\lambda]\otimes_{\sO_X}\sO[\delta]$. For example, we might start
from one ``ample'' line bundle $\sL$ on~$X$, and put
$\Delta=\bbZ$, $\sO[n]:=\sL^{\otimes n}$.

In the classical situation we would define a $\Delta$-graded ring
$R_*=\Gamma_*(X,\phi)$ simply by putting $R_\lambda:=\Gamma(X,\sO[\lambda])$,
using isomorphisms $\sO[\lambda]\otimes\sO[\delta]\to\sO[\lambda+\delta]$
to define multiplication on $R_*:=\bigoplus R_\lambda$. In the generalized
situation we should care about operations of higher arities as well.

Thus we put $R_\lambda(n):=\Gamma(X,L_{\sO_X}(n)\otimes_{\sO_X}\sO[\lambda])$,
where $L_{\sO_X}(n)=\sO_X^{(n)}=\sO_X(n)$ 
is the free $\sO_X$-module of rank~$n$. Twisting the ``composition maps''
$\sO_X(k)\times\sO_X(n)^k\to\sO_X(n)$ with the aid of $\sO[\lambda]$ and
$\sO[\mu]$, and taking global sections, we arrive to some
``graded composition maps'' $R_\lambda(k)\times R_\mu(n)^k\to 
R_{\lambda+\mu}(n)$.

In this way a $\Delta$-graded algebraic monad $R$ may be defined as 
a collection of sets $R_\lambda(n)$, $n\geq0$, $\lambda\in\Delta$,
together with some transition maps $R_\lambda(\phi):R_\lambda(n)\to 
R_\lambda(m)$, defining an algebraic endofunctor $R_\lambda$,
a ``unit element'' $\bu\in R_0(1)$, and some ``graded composition maps''
$\mu^{(k,\lambda)}_{n,\mu}:
R_\lambda(k)\times R_\mu(n)^k\to R_{\lambda+\mu}(n)$,
satisfying some natural ``graded'' analogues of the axioms of an algebraic 
monad. Of course, commutativity also has its graded counterpart,
so we can speak about commutative graded algebraic monads, i.e.\
graded generalized rings. 

Informally, $R_\lambda(n)$ should be thought
of as ``the degree $\lambda$ part of the free $R$-module of rank~$n$''.

\nxsubpoint (Sections of projective bundles.)
Once this ``generalized projective geometry'' is constructed,
we apply it to symmetric algebras $S_{\sO_S}(\sF)$ of 
quasicoherent sheaves~$\sF$ (e.g.\ locally free $\sO_S$-modules,
i.e.\ vector bundles) on a generalized scheme~$S$, thus defining
{\em projective bundles\/} $\bbP_S(\sF)=\Proj S_{\sO_S}(\sF)$ over~$S$.

These projective bundles seem to retain most of their classical properties.
For example, sections of $\bbP_S(\sF)$ over~$S$ are in one-to-one
correspondence to those strict quotients of $\sF$, which are line bundles
over~$S$. We can apply this to compute $A$-valued points
of $\bbP^n=\bbP^n_{\Fempty}=\Proj\Fempty[T_0^{[1]},\ldots,T_n^{[1]}]$
over any generalized ring~$A$ with $\Pic(A)=0$: elements of $\bbP^n(A)$
correspond to surjective homomorphisms from free $A$-module $A(n+1)$
into $|A|=A(1)$, considered modulo multiplication by elements 
from~$|A|^\times$. Since $\Hom_A(A(n+1),|A|)\cong|A|^{n+1}$, we see that
$\bbP^n(A)$ consists of $(n+1)$-tuples $(X_0:X_1:\ldots:X_n)$ of elements
of~$|A|$, generating together the unit ideal of~$|A|$ (i.e.\ ``coprime''),
modulo multiplication by invertible elements from~$|A|^\times$,
exactly as in the classical case.

For example, $\bbP^n(\Fone)$ consists of $2^{n+1}-1$ points
$(X_0:\ldots:X_n)$, with $X_i\in\{0,1\}$, and not all $X_i=0$.
In particular, $\bbP^1(\Fone)=\{(0:1),(1:1),(1:0)\}=\{0,1,\infty\}$.
On the other hand, the underlying topological space
of $\bbP^1_{\Fone}=\Proj\Fone[X,Y]$ consists also of three points:
two closed points, corresponding to graded ideals $(X)$ and $(Y)$,
and one generic point $\xi$, corresponding to $(0)$.

It is interesting to note that $\bbF_1$-points $0$ and $\infty$
are supported at the closed points of $\bbP^1_{\Fone}$, while
{\em $1=(1:1)$ is supported at the generic point $\xi\in\bbP^1_{\Fone}$}.
This is possible only because $\bbP^1_{\Fone}$ is ``separated''
over $\Spec\Fone$, but not separated in the classical sense:
otherwise any section would have been a closed map.

\nxpointtoc{Applications to Arakelov geometry}
Chapter {\bf 7} deals with the applications of the theory of generalized
rings and schemes to Arakelov geometry. In particular, we construct
the ``compact model'' $\CompZ$ of~$\bbQ$, show that any algebraic
variety $X/\bbQ$ admits a finitely presented model $\sX$ over $\CompZ$,
prove an ``archimedian valuative criterion of properness'',
use it to extend rational points $P\in X(\bbQ)$ to uniquely
determined sections $\sigma_P$ of a model $\sX/\CompZ$, provided $X$
and $\sX$ are projective, and finally relate the arithmetic degree
of the pullback $\sigma_P^*\sO_\sX(1)$ of the ample line bundle on~$\sX$
to the logarithmic height of point~$P$.

Notice, however, that the development of an intersection theory and
a theory of Chern classes is postponed until Chapter~{\bf 10}.

\nxsubpoint (Construction of $\CompZ\strut^{(N)}$.)
Fix any integer $N>1$, and consider generalized rings
$B_N:=\bbZ[N^{-1}]\subset\bbQ$, $A_N:=B_N\cap\Zninfty$. One can show that
$1/N\in|A_N|$, and that $A_N[(1/N)^{-1}]=B_N$, i.e.\ $\Spec B_N$
is isomorphic to principal open subsets both of $\Spec\bbZ$ and of
$\Spec A_N$. Therefore, we can patch together $\Spec\bbZ$ and $\Spec A_N$
along their principal open subsets isomorphic to $\Spec B_N$,
and obtain a generalized scheme, denoted by $\CompZ\strut^{(N)}$.

\nxsubpoint (Structure of $\CompZ\strut^{(N)}$.)
One can use any localization theory $\cT^?$ to construct $\CompZ\strut^{(N)}$.
Let's describe the structure of $\CompZ\strut^{(N)}$, constructed with
the aid of unary localization theory~$\cT^u$. Then $\CompZ\strut^{(N)}=
\Spec\bbZ\cup\{\infty\}$ as a set. The topology on $\CompZ\strut^{(N)}$
is such that the generic point $\xi$ of $\Spec\bbZ$ is still generic,
$\infty$ is a closed point, as well as all points $p$ with $p\mid N$,
while points $p$ with $p\nmid N$ are not closed. Their closure consists
of $p$ and~$\infty$. Furthermore, while 
the local rings of $\CompZ\strut^{(N)}$ at points $p$ are equal
to $\bbZ_{(p)}$ as expected, the local ring at~$\infty$ equals
$A_N$, not $\Zninfty$ as one expects from the ``true'' compactification
of $\Spec\bbZ$.

In this way $\CompZ\strut^{(N)}$ is just a crude approximation
of $\CompZ$, which depends on the choice of $N>1$ to a considerable extent.

\nxsubpoint (Morphisms $f^{NM}_N:\CompZ\strut^{(NM)}\to\CompZ\strut^{(N)}$.)
Given any integers $N>1$ and $M\geq1$, we construct a natural
generalized scheme morphism 
$f^{NM}_N:\CompZ\strut^{(NM)}\to\CompZ\strut^{(N)}$, which is an 
isomorphism on the open subschemes of $\CompZ^{(N)}$ and $\CompZ^{(NM)}$,
isomorphic to $\Spec\bbZ$. In other words, $f^{NM}_N$ changes something
only over the archimedian point~$\infty$, and we see later that
$f^{NM}_N$ is projective, i.e.\ it behaves like a sort of blow-up
over the archimedian point $\infty$ of $\CompZ\strut^{(N)}$. Notice that
$f^{NM}_N$ is bijective, i.e.\ this is something like a blow-up
of a cusp-like singularity on a projective curve.

\nxsubpoint (Compactification $\CompZ$.)
Since these $f^{NM}_N$ are ``transitive'' in an obvious sense,
we get a projective system of generalized schemes $\CompZ\strut^{(N)}$
over the filtered set of integers $N>1$, ordered by division.
It is very natural to define the ``true'' compactification $\CompZ$
as $\quotedprojlim\CompZ\strut^{(N)}$.

Unfortunately, this projective limit doesn't exist in the category
of generalized schemes. However, we can compute it either in the
category of pro-generalized schemes, or in the category of
generalized ringed spaces. When we are concerned only with finitely
presented objects (sheaves, schemes,\dots) over $\CompZ$,
these two approaches tend to yield equivalent answers, so we can
freely choose the most convenient approach in each particular situation.

\nxsubpoint ($\CompZ$ as a generalized ringed space.)
When we consider $\CompZ$ as a generalized ringed space,
we obtain something very similar to one's idea of the ``smooth''
compactification of $\Spec\bbZ$. For example, $\CompZ=\Spec\bbZ\cup\{\infty\}$,
but now all points except the generic one $\xi$ are closed, and
the local ring at $\infty$ is indeed equal to $\Zninfty$, as one would
expect. Therefore, we can say that our $\CompZ$ is the ``true''
compactification of $\Spec\bbZ$.

\nxsubpoint (Finitely presented sheaves on $\CompZ$.)
Notice that we have an open embedding $i:\Spec\bbZ\to\CompZ$,
and natural morphisms $\hat\xi:\Spec\bbQ\to\CompZ$ (``generic point'')
and $\hat\eta:\Spec\Zninfty\to\CompZ$) (``archimedian point''),
both in the pro-generalized scheme and in the generalized ringed space
descriptions of $\CompZ$. Therefore, if we start with a finitely presented
$\sO_{\CompZ}$-module~$\sF$, compute its pullbacks with respect to
these three maps, and take their global sections, we obtain
a finitely presented $\bbZ$-module~$F_\bbZ$, a finite-dimensional
$\bbQ$-vector space~$F_\bbQ$, and a finitely presented $\Zninfty$-module
$F_\infty$. Furthermore, since $\hat\xi$ factorizes through both
$i$ and~$\hat\eta$, we get $F_\bbZ\otimes_\bbZ\bbQ\cong F_\bbQ\cong
F_\infty\otimes_{\Zninfty}\bbQ$. 

In this way we obtain a functor
from the category of finitely presented $\sO_{\CompZ}$-modules
into the category $\cC$ of triples $(M_\bbZ,M_\infty,\theta)$,
consisting of a finitely generated $\bbZ$-module~$M_\bbZ$,
a finitely presented $\Zninfty$-module~$M_\infty$, and an isomorphism
$\theta:M_\bbZ\otimes_\bbZ\bbZ\simto M_\infty\otimes_{\Zninfty}\bbQ$
between their scalar extensions to~$\bbQ$. In fact, {\em this functor
is an equivalence of categories.} Furthermore, {\em other categories
of finitely presented objects over~$\CompZ$, e.g.\ finitely presented
schemes or finitely presented sheaves on finitely presented schemes,
admit a similar description in terms of corresponding objects over
$\Spec\bbZ$, $\Spec\Zninfty$, and isomorphisms between their base change to
$\Spec\bbQ$.}

\nxsubpoint (Implication for finitely presented models.)
An immediate implication is that constructing a finitely presented model
$\bar\sX/\CompZ$ of an algebraic variety $X/\bbQ$ is equivalent to
constructing finitely presented models $\sX/\bbZ$ and $\sX_\infty/\Zninfty$
of this~$X$. Since existence and properties of models over $\bbZ$ 
have been intensively studied for quite a long time, we have to
concentrate our efforts on models over~$\Zninfty$.

Before studying models of algebraic varieties, we discuss some further
properties of $\CompZ$ and $\CompZ\strut^{(N)}$:

\nxsubpoint (Finite presentation of $\CompZ\strut^{(N)}$.)
Generalized scheme $\CompZ\strut^{(N)}$ is finitely presented over
$\Fpm=\Fone[-^{[1]}\,|\,-(-\bu)=\bu]$, hence also over $\Fone$ and~$\Fempty$.
Indeed, since $\CompZ\strut^{(N)}=\Spec\bbZ\cup\Spec A_N$ and
$\bbZ=\Fpm[+^{[2]}\,|\,x+0=x=0+x$, $x+(-x)=0]$, all we have to check is
that $A_N$ is a finitely presented $\Fpm$-algebra. We show this
by presenting an explicit presentation of $A_N$ over~$\Fpm$.
Namely, $A_N$ is generated over $\Fpm$ by $(s_p^{[p]})_{p\mid N}$,
where $p$ runs through all prime divisors of~$N$, and
$s_n$ denotes the {\em averaging operation\/}
$s_n:=(1/n)\{1\}+\cdots+(1/n)\{n\}$, subject to following
{\em idempotency}, {\em symmetry\/} and {\em cancellation\/} relations:
\begin{align}
&s_n(\{1\},\{1\},\ldots,\{1\})=\{1\}\\
&s_n(\{1\},\{2\},\ldots,\{n\})=
s_n\bigl(\{\sigma(1)\},\ldots,\{\sigma(n)\}\bigr),
\quad\forall\sigma\in\gS_n\\
&s_n(\{1\},\ldots,\{n-1\},-\{n-1\})=s_n(\{1\},\ldots,\{n-2\},0,0)
\end{align}

For example, $A_2=\Fpm[*^{[2]}\,|\,x*x=x$, $x*y=y*x$, $x*(-x)=0]$
is generated by one binary operation $*:=s_2$, similarly to $\bbZ$ and
$\Finfty$. There are some interesting matrices over~$A_2$, which satisfy
the braid relations; cf.~\ptref{sp:genrg.braid.repr} for more details.

\nxsubpoint (Line bundles over $\CompZ\strut^{(N)}$.)
Another interesting result is that $\Pic(\Spec A_N)=0$; since
$\Pic(\Spec\bbZ)=0$, we see that $\Pic(\CompZ\strut^{(N)})=
\bbZ^{\times}\backslash B_N^{\times}/|A_N|^\times=B_N^\times/\{\pm1\}=
B_{N,+}^\times$. We prefer to write this group of positive invertible
elements of $B_N=\bbZ[N^{-1}]$ in additive form:
\begin{equation}
\Pic(\CompZ\strut^{(N)})=\log B_{N,+}^\times=\bigoplus_{p\mid N}\bbZ\cdot\log p
\end{equation}
For any $\lambda\in B_{N,+}^\times$ we denote by $\sO(\log\lambda)$
the line bundle on $\CompZ\strut^{(N)}$ 
corresponding to $\log\lambda$ under this
isomorphism. In this way $\sO(0)=\sO$ and $\sO(\log\lambda+\log\mu)=
\sO(\log\lambda)\otimes_\sO\sO(\log\mu)$.

We see that $\Pic(\CompZ\strut^{(N)})$ is a free abelian group of rank~$r$,
generated by $\sO(\log p_i)$, where $p_1$, \dots, $p_r$ are the distinct
prime divisors of~$N$. In this respect $\CompZ\strut^{(N)}$ is somewhat 
similar to $P^r:=(\bbP^1)^r$, which also has $\Pic\cong\bbZ^r$. 
Notice, however, that $\Aut(P^r)$ acts transitively on the canonical
basis of $\Pic(P^r)$, i.e.\ all standard generators of $\Pic(P^r)$
have the same properties, while $\Aut(\CompZ\strut^{(N)})=1$ acts trivially
on generators $\sO(\log p_i)$ of $\Pic(\CompZ\strut^{(N)})$.

Another interesting consequence is that we can compute $\Pic(\CompZ)=
\Pic(\quotedprojlim\CompZ\strut^{(N)})=\injlim_{N>1}\log B_{N,+}^\times=
\log\bbQ^\times_+$. We see that $\Pic(\CompZ)$ is a free abelian group
with a countable system of generators $\sO(\log p)$, $p\in\bbP$.
If $\sL\cong\sO(\log\lambda)$, we say that $\log\lambda$ is the
{\em (arithmetic) degree\/} of line bundle $\sL$ over~$\CompZ$.

\nxsubpoint (Ampleness of $\sO(\log N)$ on $\CompZ\strut^{(N)}$.)
Comparison with $(\bbP^1)^r$ leads us to believe that
$\sO(\log N)$, being the tensor product of the standard generators of
the Picard group of $\CompZ\strut^{(N)}$, must be an ample line bundle.
In order to check this we have to compute graded generalized ring
$R_*:=\Gamma(\CompZ\strut^{(N)},\sO(\log N))$, as explained
in~\ptref{sp:intro.graded.genrg}, and check whether $\Proj R_*$ is
isomorphic to $\CompZ\strut^{(N)}$. This turns out to be the case,
i.e.\ {\em $\CompZ\strut^{(N)}$ is a projective generalized scheme}
(absolutely, i.e.\ over $\Spec\Fempty$, as well as over $\Spec\Fone$ and
$\Spec\Fpm$). 

\nxsubpoint ($\CompZ$ as an infinite resolution of singularities.)
In this way we arrive to the following picture. Each $\CompZ\strut^{(N)}$
is a projective finitely presented scheme over~$\Fpm$, containing
$\Spec\bbZ$ as an open subscheme. All $f^{NM}_N:\CompZ\strut^{(NM)}\to
\CompZ\strut^{(N)}$ are projective finitely presented morphisms of
generalized schemes as well, identical on open subsets isomorphic to 
$\Spec\bbZ$. We have already discussed that each individual $f^{NM}_N$ 
looks like a blow-up of a complicated cusp-like singularity at~$\infty$.
Therefore, one might consider the ``true'' or ``smooth''
compactification $\CompZ$ as the result of an infinite resolution of
singularities of any ``non-smooth'' compactification $\CompZ\strut^{(N)}$.
Some further discussion may be found in~\ptref{sp:smirnov.descr.compz}.

\nxsubpoint (Existence of models.)
After discussing properties of $\CompZ$ and $\CompZ\strut^{(N)}$,
we show existence of finitely presented models over $\CompZ$ or $\Zninfty$  
of affine or projective algebraic varieties $X/\bbQ$. Since a finitely
presented model $\bar\sX/\CompZ$ of $X/\bbQ$ is essentially the same thing
as finitely presented models $\sX/\bbZ$ and $\sX_\infty/\Zninfty$, we
need to consider only models over~$\Zninfty$. Then we have two basic
ways of constructing models:
\begin{itemize}
\item If $X=\Spec\bbQ[T_1,\ldots,T_n]/(f_1,\ldots,f_m)$ for some
polynomicals $f_i\neq0$, we may multiply $f_i$ by some non-zero rational
numbers so as to have $f_i\in|\Zninfty[T_1^{[1]},\ldots,T_n^{[1]}]|$, and
put $\sX_\infty:=\Spec\Zninfty[T_1^{[1]},\ldots,T_n^{[1]}\,|\,f_1=0$, 
\dots, $f_m=0]$. This is obviously a finitely presented $\Zninfty$-model
of~$X$, and the case of a projective~$X$ is dealt with similarly,
writing $X=\Proj\bbQ[T_1,\ldots,T_n]/(F_1,\ldots,F_m)$ for some
homogeneous polynomials $F_i$.
\item If $X$ is a closed subvariety of some algebraic variety $P$, known
to admit a $\Zninfty$-model~$\sP_\infty$ (e.g.\ $\bbA^n$ or $\bbP^n$), then
we can take for $\sX_\infty$ the ``scheme-theoretical closure''
of~$X\subset P\subset\sP_\infty$ in~$\sP_\infty$.
\end{itemize}

The first of these two approaches always produces finitely presented models,
which, however, usually have some $\Zninfty$-torsion, i.e.\
``torsion in the archimedian fiber'', which may be thought of some sort
of ``embedded analytic torsion''. On the other hand, the second approach
always produces models without $\Zninfty$-torsion, but they are usually
not finitely presented, but only of finite type. In particular,
such models cannot be immediately extended to models over $\CompZ$. 

All we can do
is represent such a model as a projective limit of finitely presented
models $\sX_\infty^{(\alpha)}$, extend each of them to a 
finitely presented model $\bar\sX^{(\alpha)}$ over~$\CompZ$, and compute the
projective limit $\bar\sX:=\quotedprojlim\bar\sX^{(\alpha)}$ in the
category of pro-generalized schemes. In this way we can always obtain
a torsion-free pro-finitely presented model $\bar\sX/\CompZ$ in the
category of pro-generalized schemes over~$\CompZ$.

This construction is quite similar to ``non-archimedian Arakelov theory''
constructions of \cite{BGS} and~\cite{GS}, where the projective system
of all models of a fixed variety~$X/K$ is considered, instead
of choosing fixed models of algebraic varieties $X$, $Y$, \dots, 
which might be poorly adapted to interesting morphisms $f:X\to Y$.

\nxsubpoint (Application: heights of points on projective varieties.)
We give an interesting application of the above constructions.
Namely, let $X\subset\bbP^n_\bbQ$ be a projective variety over~$\bbQ$,
embedded into $\bbP^n_\bbQ$ as a closed subvariety. Denote by $\sX$
the ``scheme-theoretical closure'' of $X$ in $\bbP^n_{\CompZ}$;
it is a torsion-free model of $X$ of finite type, but in general not
of finite presentation. Next, let $\sO_\sX(1)$ be the natural
ample line bundle on $\sX$, i.e.\ the pullback of the Serre line bundle
on~$\sP:=\bbP^n_{\CompZ}$ with respect to ``closed'' embedding
$\sX\to\sP$.

We show that any rational point $P\in X(\bbQ)$ extends to a unique
section $\sigma_P:\CompZ\to\sX$, using an ``archimedian valuative
criterion'' for this. Furthermore, we compute the arithmetic degree
of line bundle $\sigma_P^*\sO_\sX(1)$ over~$\CompZ$: it turns out to
be equal to the {\em logarithmic height\/} of point $P\in X(\bbQ)\subset
\bbP^n(\bbQ)$, defined by the classical formula
\begin{equation}
h(P)=\log\prod_{v\in\bbP\cup\{\infty\}}\max\bigl({|X_0|}_v,{|X_1|}_v,\ldots,
{|X_n|}_v\bigr)
\end{equation}
If we choose homogeneous coordinates $P=(X_0:X_1:\cdots:X_n)$ of~$P$
to be coprime integers, then $h(P)=\log\max(|X_0|,\ldots,|X_n|)$.

This theorem seems to be a nice counterpart of the classical theorem
of Arakelov geometry, relating heights of rational points to
arithmetic degrees of pullbacks of Serre bundles.

\nxpointtoc{Homological and homotopic algebra}
Chapter {\bf 8} opens the ``homological'', or rather ``homotopic''
part of this work, consisting of the three last chapters.
It begins with a detailed introduction or ``plan'' 
(cf.~\ptref{p:great.plan}) for 
this chapter, and for the remaining ``homotopic'' chapters as well.
Therefore, we won't repeat that introduction here; instead, we would
like to explain how homotopic algebra replaces homological algebra,
for non-additive categories like $\catMod\Sigma$,
where $\Sigma$ is some generalized ring. This is especially useful because
homotopic algebra is not widely known among specialists in arithmetic
geometry.

\nxsubpoint (Model categories.)
One of the most fundamental notions of homotopic algebra is that
of a {\em model category}, due to Quillen (cf.~\cite{Quillen};
cf.\ also \cite{DS} for a short introduction to model categories).
A {\em model category~$\cC$} is simply a category~$\cC$, with arbitrary
inductive and projective limits (actually Quillen originally required
only existence of finite limits), together with three distinguished 
classes of morphisms, called {\em fibrations}, {\em cofibrations\/}
and {\em weak equivalences}. Furthermore, if a morphism is both a 
weak equivalence and a fibration (resp.\ cofibration), it is said to be
an {\em acyclic fibration\/} (resp.\ {\em acyclic cofibration}).

These three distinguished classes of morphisms are subject to certain axioms.
For example, each of them is stable under composition and retracts,
and weak equivalences satisfy the {\em 2-out-of-3 axiom:}
``if two of $g\circ f$, $f$ and $g$ are weak equivalences, so is the third''.
Furthermore, there is an important {\em factorization axiom},
which claims that any morphism~$f$ can be factorized into a cofibration,
followed by an acyclic fibration, as well as into an acycloc cofibation,
followed by a fibration (usually one can even choose such factorizations
functorially in~$f$). This axiom is especially useful for constructions
in model categories, e.g.\ for constructing derived functors.

Finally, the last important model category axiom says that acyclic cofibrations
have the {\em LLP (left lifting property)} with respect to fibrations, or 
equivalently, fibrations have the {\em RLP (right lifting property)}
with respect to acyclic cofibrations, and that cofibrations have the
LLP with respect to acyclic fibrations. Recall that a morphism $i:A\to B$
is said to have the LLP with respect to $p:X\to Y$, or equivalently,
$p$ is said to have the RLP with respect to~$i$, if for any
two morphisms $u:A\to X$ and $v:B\to Y$, such that $pu=vi$, one
can find a ``lifting'' $h:B\to X$, such that $hi=u$ and $ph=v$:
\begin{equation}
\xymatrix{
A\ar[d]^{i}\ar@{-->}[r]^{u}&X\ar[d]^{p}\\
B\ar@{-->}[r]^{v}\ar@{.>}[ur]|{\exists h}&Y}
\end{equation}

\nxsubpoint (Example: non-negative chain complexes.)
Let us make an example to illustrate the above definition. Let $\sA$
be any abelian category with sufficiently many projective objects, 
e.g.\ $\cA=\catMod R$ for a classical ring~$R$.
Then the category $\cC:=\Ch(\cA)$ of {\em non-negative\/} chain
complexes $K=(\cdots\stackrel\partial\to K_1\stackrel\partial\to K_0\to0)$
admits a natural model category structure. Weak equivalences for this
model category structure are just the quasi-isomorphisms, i.e.\
chain maps $f:K\to L$, such that all $H_n(f):H_n(K)\to H_n(L)$ are
isomorphisms in~$\cA$. A chain map $f:K\to L$ is a fibration iff
all its components $f_n:K_n\to L_n$ with $n>0$ are epimorphisms
(surjective $R$-linear maps, if $\cA=\catMod R$). Finally,
$f:K\to L$ is a cofibration iff all its components $f_n:K_n\to L_n$ 
are monomorphisms with projective cokernels.

Then the lifting axiom becomes a reformulation of a well-known
statement from homological algebra about chain maps from
complexes of projective modules into acyclic complexes.

\nxsubpoint (Fibrant and cofibrant objects and replacements.)
Let $\cC$ be a model category, and denote by $\emptyset_\cC$ and
$e_\cC$ its initial and final objects. An object $X\in\Ob\cC$
is said to be {\em cofibrant\/} if $\emptyset_\cC\to X$ is a cofibration,
and {\em fibrant\/} if $X\to e_\cC$ is a fibration. We denote by
$\cC_c$ (resp.\ $\cC_f$, $\cC_{cf}$) full subcategories of~$\cC$,
consisting of cofibrant (resp.\ fibrant, resp.\ fibrant and cofibrant)
objects. Finally, a {\em cofibrant replacement\/} of an object~$X\in\Ob\cC$
is a weak equivalence $P\to X$ with $P\in\Ob\cC_c$ (sometimes we consider
{\em strong\/} cofibrant replacements, i.e.\ require
$P\to X$ to be an acyclic fibration), and similarly a 
{\em fibrant replacement\/} is a weak equivalence $X\to Q$ with fibrant
target. Factorization axiom of model categories, when applied to
$\emptyset_\cC\to X$ and $X\to e_\cC$, shows existence of (co)fibrant
replacements.

For example, if $\cC=\Ch(\cA)$ is the category of non-negatively graded
chain complexes as above, then $\cC_f=\cC$, but $\cC_c$ consists of
complexes of projective objects, and a cofibrant replacement $P\to X$
is a quasi-isomorphism from a complex~$P$ consisting of projective objects
to given complex~$X$, i.e.\ a {\em projective resolution.}

This explains how cofibrant replacements are used in general:
under some additional restrictions, we can construct left derived functor
$\dL F$ of a functor $F:\cC\to\cD$ between model categories by putting
$\dL F(X):=F(P)$, where $P\to X$ is any cofibrant replacement of~$X$.

One can also say that {\em fibrations and cofibrations, and the
factorization axiom are used for constructions in model categories.}

\nxsubpoint (Homotopic categories and derived functors.)
Whenever we have a model category~$\cC$, or just any category
$\cC$ with a distinguished class of weak equivalences, we denote by
$\Ho\cC$ the {\em homotopic category\/} of~$\cC$, equal by definition
to the localization of $\cC$ with respect to the set of weak equivalences.
Canonical functor $\cC\to\Ho\cC$ is denoted by~$\gamma_\cC$ or simply
by~$\gamma$. One of its
fundamental properties is that {\em $f$ is a weak equivalence in model
category~$\cC$ iff $\gamma(f)$ is an isomorphism in~$\Ho\cC$.}

Given a functor $F:\cC\to\cD$ between two model categories, we define
its {\em left derived functor\/} $\dL F:\Ho\cC\to\Ho\cD$ as the functor,
which is as close from the left to making the obvious diagram commutative
as possible. More precisely, we must have a natural transformation
$\xi:\dL F\circ\gamma_C\to\gamma_D\circ F$, such that for any functor
$H:\Ho\cC\to\Ho\cD$ and any natural transformation $\eta:H\circ\gamma_C\to
\gamma_D\circ F$ there is a unique natural transformation
$\zeta:H\to\dL F$, such that $\eta=\xi\circ(\zeta\star\gamma_\cC)$.

For example, if we have a functor $F:\cA\to\cB$ between abelian categories
with enough projectives, and denote by the same letter $F$
its extension $F:\Ch(\cA)\to\Ch(\cB)$ to non-negative chain complexes,
then $\Ho\Ch(\cA)$ is equivalent to full subcategory $\cD^{\leq0}(\cA)$
of the derived category of~$\cA$, and similarly for $\Ho\Ch(\cB)$,
and then $\dL F$ corresponds to a left derived functor in the sense of
Verdier.

We see that the homotopic category~$\Ho\cC$ plays 
the role of the derived category (or at least of its non-positive part). 
It depends only on weak equivalences of~$\cC$, hence the same is true
for derived functors as well. In fact, one can have different model
structures on the same category~$\cC$, having the same weak equivalences;
in such situations we get same $\Ho\cC$ and $\dL F$.

However, the correct choice of (co)fibrations is crucial for showing
existence of derived functors, because of the following statement,
due to Quillen: {\em if a functor $F:\cC\to\cD$ between two model categories
transforms acyclic cofibrations between cofibrant objects into 
weak equivalences, then $\dL F:\Ho\cC\to\Ho\cD$ exists and can be computed
with the aid of cofibrant replacements, i.e.\ $(\dL F)(\gamma X)\cong
\gamma F(P)$ for any cofibrant replacement $P\to X$.}

This corresponds to the classical computation of derived functors via
projective resolutions.

\nxsubpoint (Simplicial objects and Dold--Kan correspondence.)
We see that in order to obtain ``homological algebra''
over a generalized ring~$\Sigma$ we have to do two things:
(a) construct a counterpart of $\Ch(\catMod\Sigma)$; this is necessary
since the ``non-negative chain complex'' description makes sense only
over abelian categories, and (b) introduce a reasonable model category 
structure on this replacement for $\Ch(\catMod\Sigma)$; we need this
to be able to define and construct homotopic categories and derived functors.

Of course, the principal requirement for these constructions is that
we should recover $\Ch(\catMod\Sigma)$ with its model category structure
explained above whenever $\Sigma$ is a classical ring.

It turns out that the natural replacement for $\Ch(\cA)$ for a non-abelian
category $\cA$ (e.g.\ $\cA=\catMod\Sigma$) is the category
$s\cA=\catFunct(\catDelta^0,\cA)$ of {\em simplicial objects\/} over~$\cA$.
This is due to the following fact: {\em whenever $\cA$ is an abelian
category, there are natural equivalences $K:\Ch(\cA)\leftrightarrows s\cA:N$
between the category of non-negative chain complexes over~$\cA$ and
the category of simplicial objects over~$\cA$.} This correspondence
between chain complexes and simplicial objects is called {\em Dold--Kan
correspondence.}

\nxsubpoint (Simplicial objects.)
Since the simplicial objects play a fundamental role in the homotopic algebra
parts of our work, we would like to fix some notations,
mostly consistent with those of~\cite{GZ}.

We denote by $[n]$ the {\em standard finite ordered set\/}
$\{0,1,\ldots,n\}$, endowed with its natural linear order. Notice that
$[n]$ is an $n+1$-element ordered set, while $\stn=\{1,2,\ldots,n\}$
is an $n$-element unordered set. One shouldn't confuse these two notations.

Next, we denote by $\catDelta$ the category of all standard finite ordered sets
$[n]$, $n\geq0$, considered as a full subcategory of the category of
ordered sets (i.e.\ morphisms $\phi:[n]\to[m]$ are the non-decreasing maps).
We denote by $\partial^i_n:[n-1]\to[n]$, $0\leq i\leq n>0$, the 
{\em $i$-th face map},
defined as increasing injection not taking value~$i$, and by 
$\sigma^i_n:[n+1]\to[n]$, $0\leq i\leq n$, the {\em $i$-th degeneracy map},
i.e.\ the non-decreasing surjection, which takes value~$i$ twice.

Now a {\em simplicial object\/} over a category~$\cA$ is simply a contravariant
functor $X:\catDelta^0\to\cA$. We usually write $X_n$ instead of $X([n])$,
$n\geq0$, and put $d_i^{n,X}:=X(\partial^i_n):X_n\to X_{n-1}$,
$s_i^{n,X}:=X(\sigma^i_n):X_n\to X_{n+1}$. Since any morphism in $\catDelta$
can be written as a product of face and degeneracy maps, a simplicial
object can be described in terms of a sequence of objects $(X_n)_{n\geq0}$ and
a collection of face and degeneracy morphisms $d_i^{n,X}:X_n\to X_{n-1}$,
$s_i^{n,X}:X_n\to X_{n+1}$, $0\leq i\leq n$, subject to certain relations.

\nxsubpoint (Simplicial sets.)
Our next step is to introduce a model category structure on $s\catMod\Sigma$,
compatible via Dold--Kan correspondence with one we had on $\Ch(\catMod\Sigma)$
whenever $\Sigma$ is a classical ring. We consider first the ``basic''
case $\Sigma=\Fempty$, i.e.\ we need a model structure on
the category $s\catMod\Fempty=s\catSets$ of simplicial sets.

This category has indeed a classical model category structure, which 
can be described as follows. Denote by $\Delta(n)\in\Ob s\catSets$ the
``standard $n$-dimensional simplex'', i.e.\ the contravariant functor
$\catDelta^0\to\catSets$, $[k]\mapsto\Hom_\catDelta([k],[n])$,
represented by~$[n]$. We denote by $\dot\Delta(n)\subset\Delta(n)$
the ``boundary'' of $\Delta(n)$; actually $\dot\Delta(n)_k$ can be described
as the set of all non-decreasing non-surjective maps $\phi:[k]\to[n]$.

We denote by $\Lambda_k(n)\subset\dot\Delta(n)$, $0\leq k\leq n$, 
the simplicial set,
obtained from $\dot\Delta(n)$ by removing its $k$-th $(n-1)$-dimensional face.
Set $\Lambda_k(n)_m$ consists of all order-preserving maps $\phi:[m]\to[n]$,
such that $[n]-\phi([m])$ is distinct from both $\emptyset$ and $\{k\}$.

Now consider the set of {\em standard cofibrations\/} $I:=\{\dot\Delta(n)\to
\Delta(n)\}_{n\geq0}$ and the set of {\em standard acyclic cofibrations\/}
$J:=\{\Lambda_k(n)\to\Delta(n)\}_{0\leq k\leq n>0}$. We {\em define\/}
acyclic fibrations in $s\catSets$ as the set of morphisms, which have
the RLP with respect to all standard cofibrations from~$I$, and
fibrations as the set of morphisms having the RLP with respect to~$J$.
After that cofibrations (resp.\ acyclic cofibrations) 
are defined as morphisms having the LLP with respect to all acyclic
fibrations (resp.\ fibrations). In particular, all morphisms from~$I$
(resp.\ $J$) are indeed cofibrations (resp.\ acyclic cofibrations).
Finally, we define {\em weak equivalences\/} in $s\catSets$ as 
morphisms, which can be decomposed into an acyclic cofibration,
followed by an acyclic fibration.

It is a theorem of Quillen that the cofibrations, fibrations and
weak equivalences thus defined do satisfy the model category axioms,
so we obtain a natural model category structure on $s\catSets$.
This model category structure admits a ``topological'' description:
if $|X|$ denotes the geometric realization of simplicial set~$X$,
then $f:X\to Y$ is a weak equivalence iff $|f|:|X|\to|Y|$ is a weak
equivalence in topological sense, i.e.\ iff all $\pi_n(|f|,x):
\pi_n(|X|;x)\to\pi_n(|Y|;f(x))$ are bijective, for all $n\geq0$
and all choices of base point $x\in|X|$.

Furthermore, cofibrations in $s\catSets$ admit a very simple description:
{\em $f:X\to Y$ is a cofibration of simplicial sets iff
$f$ is a monomorphism, i.e.\ iff all $f_n:X_n\to Y_n$ are injective.}

\nxsubpoint (Model category structure on $s\catMod\Sigma$.)
Now we are ready to describe our model category structure on $s\catMod\Sigma$,
for any algebraic monad~$\Sigma$. First of all, notice that we have
adjoint {\em free $\Sigma$-module\/} and {\em forgetful\/} functors
$L_\Sigma:s\catSets\leftrightarrows s\catMod\Sigma:\Gamma_\Sigma$,
equal to simplicial extensions of corresponding functors 
$\catSets\leftrightarrows\catMod\Sigma$. It is extremely natural to expect
that $L_\Sigma$ preserves cofibrations and acyclic cofibrations,
especially if we hope to construct a left derived $\dL L_\Sigma$ afterwards.
In particular, we can apply $L_\Sigma$ to ``standard generators''
$I$ and $J$ of the model structure on $s\catSets$: we see that
all morphisms from $L_\Sigma(I)$ (resp.\ $L_\Sigma(J)$) must be
cofibrations (resp.\ acyclic cofibrations) in $s\catMod\Sigma$.

A natural next step is to consider the model structure on $s\catMod\Sigma$,
generated by these sets $L_\Sigma(I)$ and $L_\Sigma(J)$ in the same way
as the model category structure on $s\catSets$ was generated by $I$ and~$J$.
This means that we {\em define\/} fibrations in $s\catMod\Sigma$ as
morphisms with the RLP with respect to all morphisms from $L_\Sigma(J)$,
and so on.

We apply a theorem of Quillen (\cite[2.4]{Quillen}, th.~4) to show
that we indeed obtain a model category structure on $s\catMod\Sigma$
in this way. Adjointness of $L_\Sigma$ and~$\Gamma_\Sigma$ shows
that {\em $f:X\to Y$ is a fibration (resp.\ acyclic fibration)
in $s\catMod\Sigma$ iff $\Gamma_\Sigma(f)$ is one in $s\catSets$.}
Furthermore, it turns our that {\em $f:X\to Y$ is a weak equivalence in
$s\catMod\Sigma$ iff $\Gamma_\Sigma(f)$ is one in $s\catSets$.}
Since cofibrations are characterized by their LLP with respect to
acyclic fibrations, this completely determines the model category 
structure of $s\catMod\Sigma$.

\nxsubpoint (Derived category of $\Lambda$-modules.)
One can check that the model category structure just defined on
$s\catMod\Lambda$ indeed corresponds to the natural model category structure
on $\Ch(\catMod\Lambda)$ under Dold--Kan correspondence 
whenever $\Lambda$ is a classical ring. Therefore,
it is extremely natural to put $\cD^{\leq0}(\Lambda)=\cD^{\leq0}
(\catMod\Lambda):=\Ho s\catMod\Lambda$. If $\Lambda$ admits a zero
(i.e.\ is an $\Fone$-algebra), then we have a natural {\em suspension functor}
$\Sigma:s\catMod\Lambda\to s\catMod\Lambda$, which has a left derived
functor $\Sigma=L\Sigma:\cD^{\leq0}(\Lambda)\to\cD^{\leq0}(\Lambda)$.
This suspension functor is actually a generalization of the degree translation
functor $T$ on $\Ch(\Lambda)$. Therefore, it is natural to define
$\cD^-(\Lambda)$ as the {\em stable homotopic category\/} of $s\catMod\Lambda$,
obtained by formally inverting functor $L\Sigma$ on $\Ho s\catMod\Lambda$.

\nxsubpoint (Derived category of $\Fone$-modules.)
It is interesting to note that $s\catMod\Fone$ is exactly the category
of pointed simplicial sets, with its model category structure extensively
used in algebraic topology. Therefore, $\cD^{\leq0}(\Fone)=\Ho
s\catMod\Fone$ is the topologists' homotopic category,
and $\cD^-(\Fone)$ is the {\em stable homotopic category}. In this way
a lot of statements in algebraic topology become statements
about homological algebra over~$\Fone$.  

\nxsubpoint (Derived scalar extension.)
Given any algebraic monad homomorphism $\rho:\Sigma\to\Xi$, we have
a scalar extension functor $\rho^*=\Xi\otimes_\Sigma-:\catMod\Sigma\to
\catMod\Xi$, which can be extended componentwise to a functor
$\rho^*=s\rho^*:s\catMod\Sigma\to s\catMod\Xi$. This functor happens
to transform generating cofibrations $L_\Sigma(I)$ of $s\catMod\Sigma$
into cofibrations $L_\Xi(I)$ in $s\catMod\Xi$, and similarly for $L_\Sigma(J)$.
Starting from this fact, we show that $\rho^*$ preserves cofibrations and
acyclic cofibrations, hence by Quillen's theorem 
it admits a left derived $\dL\rho^*:\cD^{\leq0}
(\Sigma)\to\cD^{\leq0}(\Xi)$, which can be computed with the aid of
cofibrant replacements: $\dL\rho^*X=\rho^*P$ for any cofibrant replacement
$P\to X$. If $\Sigma$ and $\Xi$ are classical rings, this $\dL\rho^*$
corresponds to the usual left derived functor, since $P\to X$ corresponds
via Dold--Kan to a projective resolution of a chain complex 
of $\Sigma$-modules. On the other hand, if $\Sigma$ (hence also $\Xi$) is 
an algebraic monad with zero, $\dL\rho^*$ commutes with suspension functors,
hence it induces a ``stable left derived functor'' $\dL\rho^*:\cD^-(\Sigma)
\to\cD^-(\Xi)$, which also corresponds via Dold--Kan to the classical
left derived functor whenever $\Sigma$ and $\Xi$ are classical rings.

\nxsubpoint (Derived tensor product.)
Now let $\Lambda$ be a generalized ring, i.e.\ a commutative algebraic monad.
Then we have a tensor product $\otimes=\otimes_\Lambda:\catMod\Lambda\times
\catMod\Lambda\to\catMod\Lambda$, and we want to construct
a {\em derived tensor product\/} $\Lotimes:\cD^{\leq0}(\Lambda)\times
\cD^{\leq0}(\Lambda)\to\cD^{\leq0}(\Lambda)$. We show that
such a derived tensor product exists, and can be computed as follows.
We extend $\otimes$ to simplicial $\Lambda$-modules by putting
$(X\otimes_\Lambda Y)_n:=X_n\otimes_\Lambda Y_n$ for all $n\geq0$.
After that we define $X\Lotimes Y:=P\otimes_\Lambda Q$, where
$P\to X$ and $Q\to Y$ are arbitrary cofibrant replacements.

If $\Lambda$ is a classical commutative ring, this construction yields
the usual derived product, because of {\em Eilenberg--Zilber theorem},
which says that taking the diagonal of a bisimplicial object over an
abelian category corresponds via Dold--Kan (up to a homotopy equivalence) to
computing the total complex of a bicomplex. Notice, however, that
in the classical situation it is enough to replace {\em one\/} of
the arguments $X$ and~$Y$ by any its projective resolution, while
in the generalized context this usually does not suffice: we {\em must\/}
take cofibrant replacements of both $X$ and~$Y$.

If $\Lambda$ admits a zero, this derived tensor product $\Lotimes$
commutes with suspension in each variable, hence it extends to the
stable homotopic category, yielding a functor $\Lotimes:\cD^-(\Lambda)
\times\cD^-(\Lambda)\to\cD^-(\Lambda)$.

An important remark is that $\Lotimes$ defines an ACU $\otimes$-structure
on $\cD^{\leq0}(\Lambda)$ (or on $\cD^-(\Lambda)$, when $\Lambda$ 
admits a zero), which has inner Homs $\dR\Hom$; these inner Homs
are indeed constructed as some right derived functors.

\nxpointtoc{Homotopic algebra over topoi}
The aim of Chapter {\bf 9} is to transfer the previous results
to the case of modules over a generalized ringed topos $\cE=(\cE,\sO)$,
e.g.\ topos of sheaves of sets over a site $\cS$ or a topological space~$X$.
In particular, results of this chapter are applicable to generalized ringed
spaces, e.g.\ generalized schemes.

\nxsubpoint (Stacks.)
We start by recalling some general properties of stacks, which can be found
for example in \cite{Giraud}. Let us explain here why we need to use stacks
at all. The main idea is that a stack (over a site or a topos) is
the ``correct'' counterpart of a sheaf of categories.

Consider the following example. Let $\cS$ be a site with fibered products
(e.g.\ the category of open subsets $U\subset X$ of a topological space~$U$.)
Then a presheaf of sets $\sF$ over~$\cS$ is simply a contravariant functor
$\cS^0\to\catSets$, i.e.\ for any $U\in\Ob\cS$ we have a set $\sF(U)$,
and for any $V\to U$ we have a ``restriction map''
$\rho_U^V:\sF(U)\to\sF(V)$, often denoted by $x\mapsto x|_V$. Of course,
these restriction maps must be ``transitive''.

Now the sheaf condition for $\sF$ can be expressed as follows.
Suppose we have a property $P$ of sections of $\sF$ (i.e.\
for any $U$ and any $x\in\sF(U)$ we can determine whether $x$ has this
property). Suppose that $P$ is compatible with restrictions
(i.e.\ if $x$ has $P$, then so has any $x|_V$), that a section $x\in\sF(U)$
with property $P$ is uniquely determined (i.e.\ if both $x$ and $x'\in\sF(U)$
have property $P$ for some~$U$, then $x=x'$), and finally that
a section $x$ with property $P$ exists ``locally'' over some~$U$,
i.e.\ there is some cover $(U_\alpha\to U)$ in $\cS$, such that
in each $\sF(U_\alpha)$ there is a (necessarily unique) element
$x_\alpha$ with property $P$. The sheaf condition then says that
there is a unique element $x\in\sF(U)$, restricting to $x_\alpha$
on each $U_\alpha$. In other words, {\em the sheaf condition means
that from local existence and uniqueness of a section with some property
we can deduce global existence and uniqueness.}

Let's try to generalize this to ``sheaves of categories''. First of all,
a ``presheaf of categories'' $\cC$ over~$\cS$ must be a collection of
categories $\cC_U=\cC(U)$, defined for each $U\in\Ob\cS$, together with
some ``pullback'' or ``restriction'' functors 
$\phi^*:\cC(U)\to\cC(V)$, defined for
each $V\stackrel\phi\to U$ (we usually denote $\phi^*X$ by $X|_V$).
However, these pullback functors are usually defined only up to a canonical
isomorphism, hence $(\phi\circ\psi)^*$ might be not equal to 
$\psi^*\circ\phi^*$, but must be canonically isomorphic to it. Of course,
these isomorphisms $(\phi\circ\psi)^*\cong\psi^*\circ\phi^*$ must be
compatible on triple compositions $\chi\circ\phi\circ\psi$.

There is a better equivalent way of describing ``presheaves of categories'',
which doesn't require to fix choices of all pullback functors~$\phi^*$. 
Namely, one considers the category $\cC$, with objects given by
$\Ob\cC:=\bigsqcup_{U\in\Ob\cS}(\cC(U))$, and morphisms given by 
$\Hom_\cC((V,Y),(U,X)):=
\bigsqcup_{\phi\in\Hom_\cS(V,U)}\Hom_{\cC(V)}(Y,\phi^*X)$, together with
the natural projection functor $p:\cC\to\cS$. 
This $\cS$-category $\cC\stackrel p\to\cS$ turns out to contain 
all data we need,
e.g.\ $\cC(U)=\cC_U$ is the pre-image of $U\subset\Ob\cS$, considered as 
a point subcategory of~$\cS$, but determines $\phi^*$ only up to isomorphism,
i.e.\ we don't have to fix any choice of~$\phi^*$. Transitivity
conditions for $\phi^*$ translate into some conditions for $\cC\stackrel p\to
\cS$, and when these conditions hold, we say that
{\em $\cC$ is a fibered category over~$\cS$.} (cf.\ SGA~1 or \cite{Giraud} 
for more details). However, we often use the ``na{\"\i}ve'' description
in terms of categories $\cC(U)$ and pullback functors $\phi^*$,
just to convey the main ideas of the constructions we perform.

Next step is to define ``sheaves of categories''. Let $\cC$ be a ``presheaf
of categories'', i.e.\ a fibered category $\cC\to\cS$. Let $P$ be a 
property of objects of $\cC$, compatible with restrictions, and
determining an object of fiber $\cC(U)=\cC_U$ uniquely up to a unique
isomorphism. Notice that we cannot expect objects of categories to
be singled out by some category-theoretic property, since any isomorphic
object must have the same property; anyway, the restriction functors
are always defined up to isomorphism, so compatibility with restrictions
makes sense only for such properties. Now suppose that objects with property
$P$ exist ``locally'' over some $U\in\Ob\cS$, i.e.\ one can find a
cover~$(U_\alpha\to U)$ in~$\cS$, and objects $X_\alpha\in\Ob\cC(U_\alpha)$,
having property $P$. However, uniqueness up to isomorphism doesn't
imply that the restrictions of $X_\alpha$ and $X_\beta$ to
$U_{\alpha\beta}:=U_\alpha\times_UU_\beta$ 
(i.e.\ to $U_\alpha\cap U_\beta$, if we work
with open subsets of a topological space) coincide; they rather imply existence
of isomorphisms $\theta_{\alpha\beta}$ between these two restrictions,
and uniqueness of these isomorphisms translates into ``transitivity''
$\theta_{\alpha\gamma}=\theta_{\beta\gamma}\circ\theta_{\alpha\beta}$
for the restrictions of these isomorphisms to $U_{\alpha\beta\gamma}$.

In this way we obtain a {\em descent datum\/} for $\cC$ over $U_\alpha\to U$.
We expect it to be {\em efficient}, i.e.\ to determine 
an object $X\in\Ob\cC(U)$, unique up to a unique isomorphism, such that
$X|_{U_\alpha}\cong X_\alpha$. Replacing ``presheaves of categories''
with technically more convenient ``fibered categories over~$\cS$'',
we arrive to the correct definition: {a stack over a site~$\cS$ is 
a fibered category $\cC\stackrel p\to\cS$, such that any descent datum
for~$\cC$ over any cover of~$\cS$ is efficient} (cf.~\cite{Giraud}).

\nxsubpoint (Kripke--Joyal semantics.)
After discussing some basic properties of stacks, we outline a variant
of the so-called {\em Kripke--Joyal semantics}, adapted for proving statements
about stacks. The key idea is that one can transfer classical proofs
of statements about sets, maps and categories to the topos case
(i.e.\ to the case of sheaves over a site~$\cS$), by replacing in a 
regular fashion sets by sheaves of sets, i.e.\ objects of corresponding
topos~$\cE=\tilde\cS$, maps of sets by morphisms in~$\cE$,
categories by stacks over~$\cS$, and functors between categories by
{\em cartesian\/} functors between stacks. However, there is one
important restriction: the proof being transferred must be 
{\em intuitionistic}, i.e.\ it cannot involve the logical law of excluded
middle, and any methods of proofs such as {\em reductio ad absurdum},
based on this logical law.

For example, most of elementary algebra (ring and module theory) is
intuitionistic, and that's the reason why most statements about
rings and modules can be generalized to statements about sheaves of
rings and sheaves of modules over them.

\nxsubpoint (Model stacks.)
Then we use Kripke--Joyal semantics to transfer model category axioms
to the topos case, thus obtaining an interesting notion of
a {\em model stack~$\cC$} over a site $\cS$ or a topos~$\cE$. It is interesting
to note that individual fibers $\cC(U)$ of a model stack~$\cC$ are
not necessarily model categories, since the lifting properties hold
in~$\cC$ only locally (i.e.\ after passing to a suitable cover).

We study this new notion of model stack, and prove some counterparts of
Quillen's results about model categories.

\nxsubpoint (Pseudomodel stacks.)
However, when we try to construct a model structure on stack
$\stks\stSETS$, given by $(\stks\stSETS)(U):=s\widetilde{\cS_{/U}}$
(simplicial sheaves of sets over~$U$), we fail to transfer Quillen's
proof. The problem is that in the case of model categories we knew that
arbitrary direct sums (coproducts) of (acyclic) cofibrations 
are (acyclic) cofibrations again, since they could be characterized
by the LLP with respect to acyclic or all fibrations. However,
in model stacks cofibrations are characterized by {\em local\/} LLP
with respect to acyclic fibrations. This is sufficient to conclude
only that {\em finite\/} direct sums of cofibrations are cofibrations.
On the other hand, Quillen's ``small object argument'', used to
prove the factorization axiom in $s\catSets$, expects the class of
cofibrations to be closed under arbitrary direct sums, 
(possibly infinite) sequential compositions and retracts. Therefore,
we are unable to prove the factorization axiom for $\stks\stSETS$.

The only way to deal with this problem is to weaken either the 
factorization or the lifting axiom of model stacks. We decided to weaken
the lifting axiom, which is used to obtain ``homotopic'' description
of derived (i.e.\ homotopic) categories, but preserve the factorization
axiom, because it is crucial for proving existence of derived functors.

In this way we arrive to the axioms of a {\em pseudomodel stack.}
We construct a very natural structure of pseudomodel stack on $\stks\stSETS$.
Fibrations and acyclic fibrations are still characterized by their
{\em local\/} RLP with respect to ``standard generators''
$\underline J$ and $\underline I$ (considered here as maps of {\em constant\/}
simplicial sheaves over~$\cS$), while cofibrations and acyclic cofibrations
are defined as the smallest local classes of morphisms, containing
$\underline I$ (resp.\ $\underline J$) and stable under arbitrary direct
sums, pushouts, sequential composition and retracts.

This notion of pseudomodel stack seems to be sufficient for our purpose.
At least, it leads to a reasonable definition of derived categories and
derived functors, and we are able to transfer and apply Quillen's
results about existence of derived functors.

\nxsubpoint (Pseudomodel category structure on $\stks\stMOD\sO$.)
Now let $\cE=(\cE,\sO)$ be a generalized ringed topos.
Consider the stack $\cC:=\stks\stMOD\sO$, characterized by
$\cC(U):=s\catMod{\sO|_U}$. We manage to construct a pseudomodel category
structure on this stack, generated by $L_\sO(\underline I)$
and $L_\sO(\underline J)$ in a similar way to what we did for $\stks\stSETS$.

If topos $\cE$ has enough points (e.g.\ it is given by some topological
space~$X$), then {\em $f:X\to Y$ is a fibration (resp.\ acyclic fibration,
weak equivalence) in $\cC(U)=s\catMod{\sO|_U}$ iff $f_p:X_p\to Y_p$
has this property in $s\catMod{\sO_p}$ for any point $p$ of $\cE_{/U}$.}
We can define {\em pointwise cofibrations} and {\em pointwise acyclic 
cofibrations} in a similar fashion. Any cofibration is a pointwise cofibration,
but not the other way around. However, we can use these pointwise
(acyclic) cofibrations, together with the above (pointwise) fibrations,
acyclic fibrations and weak equivalences to define another
pseudomodel structure on stack $\cC$, called the {\em pointwise
pseudomodel structure}.

Since the weak equivalences are the same, these two pseudomodel structures
lead to the same derived categories and derived functors.

\nxsubpoint (Derived categories and functors.)
If $\cC$ is a pseudomodel stack, we have a class of weak equivalences 
in its ``fiber category'' $\cC(e_\cE)$ over the final object of~$\cE$
(if $\cE$ is given by a topological space~$X$, this corresponds to
taking $\cC(X)=\cC_X$). Therefore, we can define {\em homotopic categories\/}
$\Ho\cC(e_\cE)$ as the localization of $\cC(e_\cE)$ with respect to 
weak equivalences, and define left or right derived functors accordingly.

For example, if $\cC=\stks\stMOD\sO$, then $\cC(e_\cE)=s\catMod{\sO}$
is the category of simplicial sheaves of $\sO$-modules. Using the
pointwise description of weak equivalences, we see that for a classical
sheaf of rings $\sO$ our weak equivalences in $s\catMod\sO$ still
correspond to quasi-isomorphisms of chain complexes, hence we can still
write $\cD^{\leq0}(\sO):=\Ho\cC(e_\cE)$. If $\sO$ admits a zero,
we obtain a ``stable'' homotopic category $\cD^-(\sO)$ as well.

\nxsubpoint (Derived pullbacks and tensor products.)
We manage to generalize our results on existence of derived scalar extensions
and tensor products to the topos case, thus proving existence of
left derived pullbacks with respect to morphisms of generalized ringed
spaces, and of derived tensor products. When sheaves of generalized 
rings involved are classical, we recover Verdier's derived pullbacks
and derived tensor products.

\nxsubpoint (Derived symmetric powers.)
The last very important thing done in Chapter {\bf 9} is the
proof of existence of {\em derived symmetric powers $\dL S^n_\sO$}, 
which turn out to admit a description in terms of cofibrant replacements,
as usual. The proof is quite long and technical.
(It would be complicated even in the classical situation, since
$S^n_\sO$ are non-additive functors even then.) 
In order to avoid some topos-related technicalities (such as defining
a {\em constructible\/} pseudomodel structure on $\stks\stMOD\sO$),
we give a complete proof only for the {\em pointwise\/} pseudomodel 
structure, thus assuming the base topos to have enough points.

Nevertheless, this result
is very important for the next chapter, and it seems to be new
even for classical rings. Notice that we obtain only functors
$\dL S^n_\sO:\cD^{\leq0}(\sO)\to\cD^{\leq0}(\sO)$, but
no ``stable'' versions $\cD^-(\sO)\to\cD^-(\sO)$. This is sufficient
for our applications.

\nxpointtoc{Perfect cofibrations and intersection theory}
In the final Chapter~{\bf 10} we combine the previous results to
obtain a reasonable notion of perfect simplicial sheaf of
modules and perfect cofibrations, then construct $K_0$ of
perfect cofibrations (and briefly discuss higher $K$-theory),
and finally transfer Grothendieck's construction of Chow rings
and Chern classes to our case.

This chapter (and also this work) ends with an application of this
intersection theory to the ``compactified $\CompZ$'', constructed
in Chapter~{\bf 7}. We manage to construct the moduli spaces of 
vector bundles on~$\CompZ$, compute $K_0$ of vector bundles and
perfect cofibrations, and use these results to compute the Chow ring
of~$\CompZ$.

\nxsubpoint (Cofibrations as presentations of their cokernels.)
Recall that in the classical case of modules over a classical ring~$\Lambda$
cofibrations in $s\catMod\Lambda\cong\Ch(\catMod\Lambda)$ correspond
to monomorphisms of non-negative chain complexes $u:A\to B$, with
$P=\Coker u$ a complex of projective modules. Therefore, one might
think of $u:A\to B$ as a ``presentation'' of this complex~$P$.

This language is quite convenient. For example, if $u'$ is a pushout of~$u$,
then $\Coker u'\cong\Coker u$, i.e.\ ``$u'$ and $u$ represent the
same complex~$P$''. This suggests that we might use $u$ as a representative
of $P$ in some $K_0$-group, and require $[u']=[u]$ for any pushout $u'$ of~$u$.
Any composable sequence of cofibrations
$A\stackrel u\to B\stackrel v\to C$ defines a short exact sequence
$0\to\Coker v\to\Coker vu\to\Coker u\to 0$, i.e.\ {\em composition of
cofibrations corresponds to short exact sequences}. Therefore, we might
expect to have $[vu]=[u]+[v]$ in some $K_0$-group. Similarly, longer
finite composable sequences of cofibrations correspond to finite filtrations
on the quotient complex, and sequential compositions of cofibrations correspond
to (infinite) discrete exhausitive filtrations.

It is convenient to extend this way of thinking about cofibrations to
the case of cofibrations $u:A\to B$ of simplicial modules over any
generalized ring (or even algebraic monad) $\Lambda$. Actually
this is already helpful in the two previous chapters. For example,
for any two morphisms $i:K\to L$ and $u:A\to B$ in $s\catMod\Lambda$
we construct a new morphism $u\boxe i\colon
?\,\to B\otimes L$, such that
$\Coker u\boxe i=\Coker u\otimes\Coker i$ in classical situation. In this
way this ``box product'' of morphisms may be thought of as corresponding 
to the tensor product of their ``virtual cokernels''. 

The proof of existence
of derived tensor products, given in Chapter~{\bf 8}, actually involves
the following statement: {\em $u\boxe i$ is a cofibration whenever
$u$ and $i$ are cofibrations; if at least one of $u$ or $i$ is acyclic,
then $u\boxe i$ is also acyclic.} Similarly, our proof of existence
of derived symmetric powers in Chapter~{\bf 9} involves some
``symmetric powers'' $\rho_n(u):?\,\to S^n(B)$, having property
$\Coker\rho_n(u)\cong S^n(\Coker(u))$ in classical situation.

\nxsubpoint (Perfect cofibrations.)
Let $\cE=(\cE,\sO)$ be any generalized ringed topos. Recall that
we've defined in Chapter~{\bf 9} the cofibrations in $s\catMod\sO$
as the smallest {\em local\/} class of morphisms, 
containing ``standard generators''
from $L_\sO(\underline I)$, and closed under some operations,
such as pushouts, retracts and finite or infinite sequential compositions.

This leads us to define {\em perfect\/} cofibrations as the set
of morphisms $u:A\to B$ in $s\catMod\sO$, which can be {\em locally\/}
obtained from morphisms of $L_\sO(\underline I)$ by means of pushouts,
finite compositions and retracts. In classical situation perfect
cofibrations $u$ correspond to injective maps of complexes with perfect
cokernel.

Next, we say that $X$ is a perfect simplicial object if $\emptyset\to X$
is a perfect cofibration. This is again compatible with classical terminology.

\nxsubpoint ($K_0$ of perfect cofibrations.)
After that we define $K_0$ of perfect cofibrations between perfect objects
as the free abelian group, generated by such cofibrations $[u]$,
modulo certain relations. Namely, we impose relations $[\id_A]=0$,
$[u']=[u]$ for any pushout $u'$ of~$u$, and $[vu]=[u]+[v]$ for any
composable $u$ and~$v$. Another important relation is $[u']=[u]$
if $u'$ is isomorphic to $u$ in the derived category $\cD^{\leq0}(\cE,\sO)$.

We denote by $[A]$ the class of $[\emptyset\to A]$ in $K_0$. Since
$[u]=[B]-[A]$ for any perfect cofibration $u:A\to B$ between perfect
simplicial objects, our $K_0$ is generated by classes~$[A]$ 
of perfect simplicial objects.

Under some restrictions (such as quasicompactness of~$\cE$) we
are able to show that $K_0$ is already generated by {\em constant\/}
perfect simplicial objects, i.e.\ by {\em vector bundles},
similarly to the classical case.

The last important remark is that operations $\boxe$ and $\rho_n$
preserve perfect cofibrations and perfect cofibrations between perfect objects,
and are compatible with the relations of~$K_0$ (notice that we have to
apply here our results on existence of derived tensor products and
symmetric powers). In this way we obtain a {\em multiplication\/} on $K_0$,
given by $[u]\cdot[v]=[u\boxe v]$, or equivalently, $[A]\cdot[B]=[A\otimes B]$,
and {\em symmetric power operations\/} $s^n:[u]\mapsto[\rho^n(u)]$,
$[A]\mapsto[S^n(A)]$. This multiplication and symmetric power operations
determine a {\em pre-$\lambda$-ring structure\/} on $K_0$,
since the generating function $\lambda_t(x):=\sum_{n\geq0}\lambda^n(x)t^n$ 
for the ``external power operations'' $\lambda^n$ can be expressed
in terms of symmetric powers by $\lambda_t(x):=s_{-t}(x)^{-1}$.

Notice that we never need to require $\sO$ to be an alternating sheaf of 
$\Fpm$-algebras, since we use symmetric powers to define $\lambda$-structure
on~$K_0$, which always have good properties, even when $\sO$ doesn't
have symmetry or zero.

\nxsubpoint (Waldhausen's construction.)
Our construction of $K_0$ turns out to be a modification of Waldhausen's
construction of~$K_0$, adapted to the case of categories with cofibrations
and weak equivalences, but without a zero object. Since Waldhausen's
construction yields a definition of higher $K$-groups as well,
we can expect higher algebraic $K$-theory of generalized rings
to admit a description in terms of a suitable modification of Waldhausen's
construction as well.

\nxsubpoint (Chow rings and Chern classes.)
Once we obtain a pre-$\lambda$-ring $K^0=K_0(\cE,\sO)$, we can apply
Grothendieck's reasoning: construct $\gamma$-operations and
$\gamma$-filtration, define the Chow ring $CH(\cE)_\bbQ$ as the associated
graded of $K^0_\bbQ$ with respect to the $\gamma$-filtration, and define
Chern classes of elements of $K^0$ (and in particular, of vector bundles)
with the aid of $\gamma$-operations again.

The only complication here is that we don't prove that pre-$\lambda$-ring
$K^0=K_0(\cE,\sO)$ is a $\lambda$-ring. Instead, we replace $K^0$
by its largest quotient $K^0_\lambda$, which is a $\lambda$-ring.

\nxsubpoint (Intersection theory of $\CompZ$.)
Chapter~{\bf 10}, hence also this work, ends with the application of
the above constructions to $\CompZ$, already constructed in Chapter~{\bf 7}.
We show that any finitely generated projective $\Zninfty$-module is free,
and use this result to classify vector bundles and constant perfect
cofibrations between vector bundles over $\CompZ$. We use this classification
to compute $K^0=K^0(\CompZ)$, and the Chow ring $CH(\CompZ,\bbQ)$.

The answer is that $K^0(\CompZ)=\bbZ\oplus\log\bbQ_+^*$. The isomorphism
is given here by $\sE\mapsto(\rank\sE,\deg\sE)$, where $\rank\sE$ 
denotes the rank of a vector bundle~$\sE$, and $\deg\sE=\deg\det(\sE)$
is its ``arithmetic degree'' in $\Pic(\CompZ)\cong\log\bbQ_+^*$.
The Chow ring $CH(\CompZ)$, i.e.\ the associated graded of $K^0(\CompZ)$
with respect to $\gamma$-filtration, also turns out to be 
$\bbZ\oplus\log\bbQ_+^*$, with graded components 
$CH^0(\CompZ)=\bbZ$ and $CH^1(\CompZ)=
\log\bbQ_+^*=\Pic(\CompZ)$. The only non-trivial Chern class
$c_1:K^0\to CH^1(\CompZ)$ turns out to be the arithmetic degree map.

This result illustrates that the intersection theory we constructed
yields very natural results when applied to the simplest ``Arakelov
variety'' $\CompZ$, leading us to believe that it will work equally
well in other situations, related to Arakelov geometry or not.


\cleardoublepage

\mysection{Motivation: Looking for a compactification of $\Spec\bbZ$}
\label{sect:motivation}

\nxpointtoc{Original motivation} Arakelov geometry appeared as an 
attempt to transfer some proofs valid for algebraic varieties (especially 
curves and abelian varieties) over a 
functional base field $K$ to the number field case. For example, people 
wanted to prove statements like the Shafarevich conjecture or the 
Mordell conjecture.

\nxsubpoint Suppose we are given some functional field $K$. This means that 
$K$ is finitely generated of transcendence degree one over some base field~$k$,
usually assumed to be either finite or algebraically closed. Then there 
is a unique smooth projective curve $C$ over $k$ such that $K=k(C)$; in 
other words, $\Spec K$ is the generic point $\xi$ of $C$.

\nxsubpoint Given an algebraic variety $X/K$ (usually supposed to be 
smooth and projective), one might look for {\em models\/} $\sX\to C$ of $X$: 
by definition, these are flat projective schemes with generic fiber $\sX_\xi$ 
isomorphic to $X$. Now, $\sX$ is a (projective) algebraic variety over~$k$. 
For example, if $X/K$ is a curve, 
$\sX/k$ will be a projective surface. Now we 
can use intersection theory on~$\sX$ to obtain some interesting bounds 
and estimates. For example, any rational point $P\in X(K)$ lifts to a 
section $\sigma_P:C\to\sX$ by the valuative criterion of properness. 
Hence for any two rational points $P$, $Q\in X(K)$ of a projective curve 
$X/K$ one can compute the {\em intersection numbers} 
$(P\cdot Q) = (\sigma_P(C)\cdot\sigma_Q(C))$ and systematically use 
them, e.g.\ to obtain some bounds for the total number of rational points 
of $X/K$.

\nxsubpoint Observe that for such reasoning we need $\sX$ to be proper 
over~$k$, so the properness of our model $C$ of $K$ over $k$ is 
very important. If we replace $C$ by an affine curve $\Cop\subset C$, 
no argument of this sort will be possible.

\nxsubpoint Note that any functional field $K$ contains a subfield 
$K'=k(t)$ isomorphic to the field of rational functions over~$k$, and 
$K/K'$ is a finite extension. Hence, once a smooth proper model~$C'$ 
of $K'/k$ is constructed, one can take the normalization of $C'$ in $K$ 
as a smooth proper model of $K$. This means that it is sufficient to 
consider functional fields of form $K=k(T)$ to learn how to construct 
smooth proper models of {\em all\/} functional fields. 

Similarly, any number field $K$ is by definition a finite extension of 
$K'=\bbQ$, so it should be sufficient to construct (in some sense) 
a ``smooth proper model of $\bbQ$''. Such a model is usually thought of 
as a ``compactification of $\Spec\bbZ$''.

\nxsubpoint Consider first the functional case, i.e.\ $K=k(T)$. 
The first candidate for a model $C$ of $K$ over $k$ is the affine line 
$\bbA_k^1=\Spec k[T]$: it is a smooth (hence regular) curve 
(i.e.\ one-dimensional scheme), and has 
the field of rational functions equal to $K=k(T)$. However, it is 
not proper, so one has to consider the projective line $C:=\bbP_k^1$. 
The affine line $\bbA_k^1$ is an open subscheme of $\bbP_k^1$; its 
complement is exactly one point -- the {\em point at infinity~$\infty$}.

One sees that each closed point $p$ of $C=\bbP_k^1$ defines a 
discrete valuation $v_p : K\to\bbZ\cup\{+\infty\}$ with the valuation 
ring equal to the local ring $\sO_{C,p}$. This also defines a norm 
$\fnorm_p$ on the field~$K$ by the rule $|x|_p:=\rho^{-v_p(x)}$ where 
$\rho$ is any real number greater than one (one can take $\rho=e$ or 
$\rho=e^{[\kappa(p):k]}$). Conversely, all norms on $K$ which are 
trivial on $k$ (condition automatically fulfilled when $k$ is algebraic 
over a finite field) come from points of $\bbP_k^1$, and this correspondence 
is one-to-one provided we identify equivalent norms (i.e.~norms such that 
$\fnorm_1=\fnorm_2^\alpha$ for some positive~$\alpha$).

If we consider just the affine line $\bbA_k^1=\Cop\subset C$, we 
do not have any point corresponding to the valuation 
$\fnorm_\infty$ given by $|f(T)/g(T)|_\infty=\rho^{\deg f - \deg g}$, that 
corresponds to the infinite point on the projective line. Note that 
it is the only valuation $\fnorm$ of~$K$ trivial on~$k$, for which 
$|T|_\infty>1$. For all other valuations we have $|f|_p\leq1$ for any 
polynomial $f\in k[T]$.

Another way to see that $\Cop$ is not complete (= proper over~$k$) is 
this. On a complete curve~$C$ the degree of the divisor~$\div(f)$ of a 
non-zero rational function $f\in k(C)=K$ is zero. In other words, we have
\begin{equation}
\sum_{P\in C, P\neq\xi} v_p(f)\cdot[k(P):k] = 0
\end{equation}
After elevating some $\rho>1$ into this power we obtain the ``product formula'':
\begin{equation}
\prod_{P\in C, P\neq\xi} |f|_P = 1
\end{equation}

Of course, once we omit the point $\infty$, both these formulas cease to be true.

\nxsubpoint The first obvious candidate for a smooth model of $\bbQ$ is 
$\Spec\bbZ$. Indeed, this is a one-dimensional regular scheme with 
the field of rational functions equal to~$\bbQ$. However, $\Cop=\Spec\bbZ$ 
is the counterpart of the affine line $\bbA_k^1=\Spec k[T]$, 
not of the projective line $\bbP_k^1$. For example, both $\Spec\bbZ$ and 
$\bbA_k^1=\Spec k[T]$ are affine, spectra of principal ideal domains. 
Closed points of $\bbA_k^1$ correspond to irreducible unitary polynomials 
$\pi\in k[T]$, and those of $\Spec\bbZ$ -- to prime numbers $p\in\bbZ$. 
In both cases these closed points correspond to all valuations of 
the field of rational functions but one. For $\bbA_k^1$ the omitted 
valuation is the 
valuation $\fnorm_\infty$ given by the degrees of polynomials. For 
$\Spec\bbZ$ it is the only archimedian valuation $\fnorm_\infty$ defined 
by the usual absolute value of rational (or real) numbers.

If we normalize the $p$-adic norms on $\bbQ$ by requiring $|p|_p=p^{-1}$, 
we have a product formula for any $f\in\bbQ^*$, {\em provided we take 
into account the remaining archimedian norm $\fnorm_\infty$}:
\begin{equation}\label{eq:product.formula}
|f|_\infty\cdot\prod_{(p)\in\Spec\bbZ,(p)\neq(0)} |f|_p = 1
\end{equation}

All this means that $\Cop=\Spec\bbZ$ is a counterpart of $\bbA_k^1$, not 
of $\bbP_k^1$, and to construct the ``compactification'' $C$ of $\Spec\bbZ$ 
one has to take into account the archimedian valuation $\fnorm_\infty$, 
which should be thought of as corresponding to some ``point at infinity'' 
$\infty$, which is the complement of $\Cop$ in $C$. 

\nxsubpoint
Up to now this picture has been just a fancy way of describing the essential 
ideas of Arakelov geometry. In one of subsequent sections
we will indeed construct the compactification $\CompZ$ of $\Spec\bbZ$ as a 
``generalized scheme''. It happens to be a topological space with a local 
sheaf of ``generalized rings''. {\bf We shall assign precise mathematical 
meaning to these notions in the remaining part of this work.}

\nxpoint
Let's describe now how one could describe a smooth projective curve $C$ 
(say, $C=\bbP_k^1$) in terms of the complement $\Cop$ of a point $\infty$ 
in $C$ and the valuation $\fnorm_\infty$ on $K=k(C)$ corresponding to 
the omitted point. 
\nxsubpoint
Let's treat the functional field case first. Then $\fnorm_\infty$ is 
non-archimedian, so we can reconstruct the local ring 
$\sO_{C,\infty}=\{f\in K : |f|_\infty\leq 1\}$ as well as its maximal ideal 
$\gm_{C,\infty}=\{f\in K : |f|_\infty < 1\}$. Then one observes that 
$\Cop\to C$ and $C_\infty:=\Spec\sO_{C,\infty}\to C$ form a covering for 
the flat topology (fpqc). Indeed, they are flat and cover the whole of~$C$. 

So we can apply faithfully flat descent here: We assign to an 
object $\sX$ (say, a quasi-projective scheme, or a vector bundle, 
or a quasi-coherent sheaf on such a scheme) over~$C$ a triple 
$(\sX^\circ, \sX_\infty, \sigma)$, where $\sX^\circ$ and $\sX_\infty$ 
are objects of the same sort over $\Cop$ and $C_\infty$, 
respectively (the pullbacks of $\sX$), and  
$\sigma:\sX^\circ_{(K)}\simto(\sX_\infty)_{(K)}$ is an isomorphism  
of the pullbacks of these 
objects onto $\Cop\times_C C_\infty=\Spec K$ (this is the canonical 
isomorphism for a triple constructed from some $\sX$ over $C$).

Then faithfully flat descent assures us that for some kinds of objects 
(e.g.\ quasi-projective schemes or quasi-coherent sheaves on such schemes) 
the functor thus constructed from the category of objects of this kind on $C$ 
into the category of triples described above is an equivalence of categories.

\nxsubpoint 
So, for example, to describe a quasi-projective scheme $\sX$ over $C$ we can 
use a quasi-projective scheme $\sX^\circ$ over $\Cop$ together with a 
``$\sO_{C,\infty}$-structure'' on its generic fiber $X:=\sX^\circ_\xi=
\sX^\circ_{(K)}$, i.e.\ a quasi-projective scheme $\sX_\infty$ over 
$\sO_{C,\infty}$ equipped by an isomorphism of $K$-schemes 
$\sigma:X\simto (\sX_\infty)_{(K)}$.

Similarly, a vector bundle $\sE$ on $C$ is essentially the same thing 
as a vector bundle $\sE^\circ$ on $C^\circ$ equipped by a 
$\sO_{C,\infty}$-structure on its generic fiber $\sE_\xi$, which is 
a finite-dimensional $K$-vector space.

\nxsubpoint
Note that we have used here the fact that
$\Cop\times_C\Cop\cong\Cop$ and $C_\infty\times_C C_\infty\cong C_\infty$: 
otherwise we would also need some descent data on these fibered products 
as well.

\nxsubpoint
Geometrically, $\Cop$ is $C$ without one closed point $\infty$, 
$C_\infty=\Spec\sO_{C,\infty}$ consists of two points $\{\infty,\xi\}$, 
and their fibered product over $C$ is $\Spec K=\{\xi\}$. So we just 
glue $\Cop$ and $C_\infty$ by their generic points to reconstruct $C$:
\begin{equation}
\xymatrix{
\Spec K=\{\xi\}\ar[rr]\ar[d]&&C_\infty=\{\xi,\infty\}\ar[d]\\
\Cop=C-\{\infty\}\ar[rr]&&C
}
\end{equation}

\nxpoint
Sometimes it is more convenient to work with the completions 
$\hat{C}_\infty=\Spec\hat\sO_{C,\infty}$ and 
$\Cop\times_C\hat{C}_\infty = \Spec \hat{K}$, where 
$\hat K=\hat\sO_{C,\infty}\otimes_{\sO_{C,\infty}}K$ is the 
completion of $K$ with respect to $\fnorm_\infty$. In this case 
two morphisms $\Cop\to C$ and $\hat{C}_\infty\to C$ form a faithfully 
flat family, so we might have used faithfully flat descent here. 
However, $\hat{C}\times_C\hat{C}=
\Spec\hat\sO_{C,\infty}\otimes_{\sO_{C,\infty}}\hat\sO_{C,\infty}$ is 
{\em not\/} isomorphic to $\hat{C}$, so in this case our descent data 
would include an isomorphism of the two pullbacks to this scheme 
of the object $\hat{\sX}_\infty$ defined over $\hat{C}_\infty$. 

\nxsubpoint 
Since this scheme usually is not noetherian and not too much can be said 
about its structure, people usually proceed in the na{\"\i}ve way. They 
map an object $\sX$ over $C$ into a triple 
$(\sX^\circ, \hat\sX_\infty, \sigma)$, where $\sX^\circ$ is an object of 
the same sort over $\Cop$, $\hat\sX_\infty$ is over $\hat{C}_\infty$, 
and $\sigma$ is an isomorphism of the pullbacks of these two objects 
on $\Cop\times_C\hat{C}_\infty=\Spec\hat{K}$. In other words, some 
of descent data is omitted.

\nxsubpoint
So the functor that maps an object $\sX$ over $C$ into a triple as above 
is {\em not\/} necessarily an equivalence of categories now. 
However, it is still faithful, and it often turns out to be fully faithful 
or even an equivalence of categories, when restricted to objects~$\sX$ 
with some additional finiteness or flatness conditions over the 
infinite point $\infty$.

\nxsubpoint
Let's show how this works for vector bundles $\sE$ over $C$. We map 
such a vector bundle into a vector bundle $\sE^\circ$ over $\Cop$ 
together a $\hat\sO_{C,\infty}$-structure $\hat{E}_\infty$ 
on the finite dimensional $\hat K$-vector space 
$\hat E:=\sE^\circ_{(\hat{K})}=E\otimes_K\hat{K}$, where
$E:=\sE^\circ_\xi$ is the generic fiber of $\sE^\circ$. Now observe that 
the corresponding $\sO_{C,\infty}$-structure $E_\infty$ on $E$ is completely 
determined by this data since $E_\infty=\hat{E}_\infty\cap E\subset \hat{E}$, 
so we get some descent datum for $\Cop,C_\infty\to C$, that uniquely 
determines the original vector bundle~$\sE$.

Note that any choice of a $\hat\sO_{C,\infty}$-structure 
$\hat E_\infty$ on $\hat E$ corresponds to a unique 
$\sO_{C,\infty}$-structure 
on $E$ since the intersection $\hat E_\infty\cap E$ necessarily has a rank 
equal to $\dim_K E$. The real reason behind this is that, apart from 
general faithfully flat descent theory, in this case we also have 
$\sO_{C,\infty}=\hat\sO_{C,\infty}\cap K\subset\hat{K}$. In other words, 
the following square is cartesian, hence bicartesian in the category of 
rings:
\begin{equation}
\xymatrix{
\sO_{C,\infty}\ar[rr]\ar[d]&&K\ar[d]\\
\hat\sO_{C,\infty}\ar[rr]&&
  \hat\sO_{C,\infty} \otimes_{\sO_{C,\infty}} K = \hat{K}
}
\end{equation}

This means that our functor is still an equivalence 
of categories, at least for the categories of vector bundles. In other words, 
{\em any choice of a $\hat\sO_{C,\infty}$-%
structure on the completed generic fiber of $\sE$ is algebraic}. 

\nxpoint\label{p:firstcompz}
Let's try to transport these arguments to the case $K=\bbQ$, 
$\Cop=\Spec Z$. So, a flat quasi-projective scheme $\sX$ over the 
compactification $C=\CompZ$ should correspond to a flat quasi-projective 
scheme $\sX^\circ$ over $\Cop=\Spec Z$ together with a ``$\Zinfty$-structure'' 
on its completed generic fiber $X_\bbR:=\sX_\xi\otimes_\bbQ\bbR$, 
which is an algebraic scheme over the completed field of rational 
functions $\hat{K}=\bbQ_\infty=\bbR$. In other 
words, to construct a model over $\CompZ$ of some (flat) variety $X/\bbQ$ 
we need to construct its model $\sX^\circ$ over $\Spec\bbZ$, and 
to construct a $\Zinfty$-model $\hat\sX_\infty$ of $X_\bbR:=X_{(\bbR)}$.

Since there is not much mystery about models over $\Spec\bbZ$ (they are 
just some schemes in the usual sense), we will be mostly concerned with 
describing models over $\Zinfty$.

\nxsubpoint
Here we denote by $\Zinfty\subset\bbQ_\infty=\bbR$ 
the conjectural counterpart of 
the completed localization $\hat\sO_{C,\infty}$, and by 
$\Zninfty=\Zinfty\cap\bbQ\subset\bbQ$ the counterpart of the (non-completed) 
localization $\sO_{C,\infty}$. Maybe ${\hat\bbZ}_\infty$ and $\Zinfty$ 
would be a better choice of notation; our choice is motivated by the 
$p$-adic analogy: $\bbZ_p\subset\bbQ_p$ is the completed localization of $\bbZ$
at a prime~$p$, while the usual localization $\bbZ_p\cap\bbQ$ is denoted by~%
$\bbZ_{(p)}$.

Of course, for now $\Zinfty$ and $\Zninfty$ are just some formal symbols 
that we use in our considerations. On the contrary, our 
analogue of $\hat{K}$ is $\Zinfty\otimes_\bbZ\bbQ=\bbQ_\infty=\bbR$, so it 
has a well-known mathematical meaning. So we will try to describe different 
objects over $\Zinfty$ as corresponding objects over its ``fraction field'' 
$\bbR$ together with some additional structure.

\nxsubpoint
Actually it would be better to study the non-completed localization 
$\Zninfty$ instead, since it should give a theory capable of handling 
objects with ``torsion over $\infty$'', and without any additional 
``algebraicity'' restrictions (i.e.\ any descent datum would be efficient). 

We will return to this question later.

\nxsubpoint\label{ss:firstbarzinfty}
We introduce some more fancy notation: we denote by $\barZinfty$ the 
``algebraic closure of $\Zinfty$ in $\bbC$''. Then, by Galois descent 
a $\Zinfty$-structure on, say, some algebraic variety $X_\bbR$ over $\bbR$ 
should be essentially the same thing as a $\barZinfty$-structure on 
its complexification $X_\bbC$, invariant under complex conjugation.

\nxpointtoc{Vector bundles over $\CompZ$}
Let's discuss first how a vector bundle $\sE$ over $\CompZ$ looks like. 
Firstly, we need some vector bundle $\sE^\circ$ over $\Spec\bbZ$. 
Vector bundles over $\bbZ$ correspond simply to free $\bbZ$-modules of finite 
rank; so $\sE^\circ$ corresponds to some $E_\bbZ$, that can be 
identified with a lattice in $E:=E_\bbZ\otimes\bbQ\cong\sE_\xi$. 
Secondly, we need to define a ``$\Zinfty$-structure'' or a 
``$\Zinfty$-lattice'' in the real vector space $E_\bbR=E_{(\bbR)}$. 
So we have to discuss what this might mean.

\nxsubpoint\label{p:zilat}
This question is usually answered as follows. A $\Zinfty$-lattice 
(resp.\ $\barZinfty$-lattice) in some finite-dimensional real 
(resp.\ complex) vector space~$E$ is just a positive-definite quadratic 
(resp.\ hermitian) form on this space. 
{\em Our approach will be slightly different, 
but let's adopt this point of view for now.}

\nxsubpoint\label{p:zilatexpl}
Here are some reasons for such a definition. Given a finite-dimensional 
$\bbQ_p$-vector space, any $\bbZ_p$-lattice $\Lambda\subset E$ defines 
a maximal compact (for the $p$-adic topology) subgroup 
$G_\Lambda:=GL_{\bbZ_p}(\Lambda)=\{g\in GL(E):g(\Lambda)=\Lambda\}$ 
in the locally compact group $GL(E)=GL_{\bbQ_p}(E)$. Conversely, all 
maximal compact subgroups in $GL(E)$ are of this form, 
and $G_\Lambda=G_{\Lambda'}$ iff $\Lambda$ and $\Lambda'$ are similar, 
i.e.\ $\Lambda'=c\Lambda$ for some $c\in\bbQ_p^*$. This means that 
$\bbZ_p$-lattices considered up to similitude are in one-to-one correspondence 
with maximal compact subgroups of~$GL(E)$.

Now, if $E$ is a real vector space, $GL(E)$ is again a locally compact group, 
and its maximal compact subgroups are the orthogonal groups $G_Q:=O(Q)$, where 
$Q$ is any positive-definite quadratic form on~$E$. Again, $G_Q=G_{Q'}$ iff 
$Q'=cQ$ for some $c>0$, so it seems quite natural to assume that 
$\Zinfty$-lattices in $E$ correspond to positive-definite quadratic forms.

Similarly, if $E$ is a (finite-dimensional) complex vector space, 
the maximal compact subgroups of $GL(E)=GL_\bbC(E)$ are the unitary groups 
of positive-definite hermitian forms, so $\barZinfty$-lattices in~$E$ 
should correspond to positive-definite hermitian forms.

\nxsubpoint\label{ss:bunoncompz}
Hence a vector bundle $\sE$ over $\CompZ$ corresponds to a free $\bbZ$-module 
$E_\bbZ$ of finite rank, together with a positive definite quadratic form 
on the real vector space $E_\bbR:=E_\bbZ\otimes\bbR$. Vector bundles over 
$\hat{C}_\infty$, i.e.\ free $\bbZ_\infty$-modules of finite rank, correspond 
to quadratic real vector spaces. Note that these categories are {\em not\/} 
additive (they don't have neither direct products nor coproducts), 
for any reasonable choice of morphisms.

\nxpointtoc{Usual description of Arakelov varieties}\label{p:cls.arak.vbun}
Now let's try to describe more complicated objects. Suppose we start 
with a {\em complete} (i.e.\ proper over~$\bbC$) algebraic variety $X/\bbC$, and we are given some proper model $\sX/\barZinfty$, whatever this might mean. 
Now we would like to describe a vector bundle $\sE$ on $\sX$.

\nxsubpoint
First of all, the generic fiber $E=\sE_\xi=\sE_{(\bbC)}$ of~$\sE$ 
is just a vector bundle on~$X$, and we know what this means. 
Next, if we believe that 
$\barZinfty$ is something like a valuation ring, then ``by the valuative 
criterion of properness'' any $\bbC$-point $P\in X(\bbC)$ should lift 
to a unique section $\sigma_P:\Spec\barZinfty\to\sX$, so we can consider the 
pullback $\sE_P:=\sigma_P^*(\sE)$. This is a vector bundle on $\Spec\barZinfty$
with generic fiber $\sE_P\otimes\bbC\cong E_P$, i.e.\ we get a 
$\barZinfty$-lattice in each fiber $E_P$. By our convention this means that 
we have a hermitian form on each $E_P$, and it is quite natural to assume 
that these metrics depend continuously or even smoothly on the point~$P$.

In this way we see that a vector bundle $\sE$ on our model $\sX$ gives rise 
to an (algebraic) vector bundle~$E$ on~$X$ equipped by a smooth hermitian 
metric, when considered as a holomorphic vector bundle over the complex 
analytic variety $X(\bbC)$.

\nxsubpoint
If we start with a complete algebraic variety $X/\bbR$ and a 
vector bundle $\sE$ on its proper model $\sX/\Zinfty$, we can 
extend these data to an algebraic variety $X_{(\bbC)}$ over~$\bbC$ together 
with a vector bundle $\sE_{(\barZinfty)}$ on its proper model 
$\sX_{(\barZinfty)}$. Then, after reasoning as above, we end up with 
an (algebraic) vector bundle~$E$ on~$X$ and a hermitian metric~$h_E$
on the corresponding holomorphic vector bundle on complex analytic variety~%
$X(\bbC)$, compatible with complex conjugation on $X(\bbC)$.

\nxsubpoint
If we are lucky, it may turn out that this collection of data completely 
determines the original vector bundle~$\sE$. Of course, we cannot 
expect any such pair $(E,h_E)$ to be ``algebraic'', i.e.\ to correspond 
to a vector bundle~$\sE$ on~$\sX$. The metrics expected to be 
algebraic are usually called {\em admissible}.

\nxsubpoint
Suppose that $X/\bbR$ is complete and {\em smooth}, and that it has a
proper and smooth model $\sX/\Zinfty$. Then the sheaf of 
relative K\"ahler differentials $\Omega^1_{\sX/\Zinfty}$ is locally free, 
i.e.\ it is a vector bundle. Hence the above arguments can be applied to 
$\sE=\Omega^1_{\sX/\Zinfty}$. The generic fiber~$E$ of $\sE$ in this 
case is $\Omega^1_{X/\bbR}$, i.e.\ it is the cotangent bundle $T^*X$ of~$X$.

In this way a smooth proper model $\sX$ of a complete smooth algebraic 
variety $X/\bbR$ gives us a hermitian metric on $T^*X(\bbC)$, invariant under 
complex conjugation, i.e.\ a cometric on $X(\bbC)$. Since this cometric 
is supposed to be non-degenerate, it corresponds by duality to a hermitian 
metric on~$X(\bbC)$ (i.e.\ on~$TX(\bbC)$), compatible with complex conjugation.

\nxsubpoint
Note that in this reasoning it is important that not only $X$ is supposed 
to be smooth, but $\sX$ as well. One might expect an $\sX$ with 
some singularities in the special fiber to give a singular (co)metric 
on $X(\bbC)$.

\nxsubpoint
Again, we can hope that this data ($X$ plus a hermitian metric $h_X$ on 
$X(\bbC)$, compatible with complex conjugation) completely determines~$\sX$. 
At least, this is what people usually assume while doing Arakelov geometry. 

So an arithmetic (or Arakelov) variety over $\Spec\bbZ$ is usually 
{\em defined\/} to be 
a flat proper scheme $\sX/\bbZ$ with smooth generic fiber $X=\sX_\xi/\bbQ$, 
equipped by a conjugation-invariant hermitian metric~$h_X$ on~$X(\bbC)$. Some 
additional ``admissibility'' requirements are usually imposed on $h_X$, 
the most common of them being the requirement to be a K\"ahler metric.

\nxsubpoint
Of course, such an approach cannot work for arithmetic varieties 
not supposed to be generically smooth and proper. We are going to 
provide a more direct approach and define an arithmetic variety~$\sX$ 
as some ``generalized scheme''. In particular, it will be a topological 
space, so the ``fiber over infinity'' will receive a very concrete meaning, 
and the informal considerations of this chapter will get a 
rigorous interpretation.

\clearpage

\mysection{$\Zinfty$-Lattices and flat $\Zinfty$-modules}
\label{sect:zinfty.lat.mod}

In this chapter we are going to study in more detail $\Zinfty$-lattices 
and $\Zinfty$-structures on real vector spaces. At first we restrict 
ourselves to $\Zinfty$-lattices in finite-dimensional vector spaces. 
We resume the considerations of~\ptref{p:zilat} and \ptref{p:zilatexpl}, and 
try to extend them so as to obtain a notion of a $\Zinfty$-structure 
on a finite $\bbR$-algebra $A$. This leads to a considerable change of 
point of view with respect to one recalled in~\ptref{p:zilat}.

\nxpointtoc{Lattices stable under multiplication}
Recall that our currently adopted approach consists in describing  
sublattices~$A$ of a finite-dimensional (real or $p$-adic) vector space~$E$ 
in terms of maximal compact subgroups $G_A=\{g:g(A)=A\}$ of~$GL(E)$; 
cf.~\ptref{p:zilatexpl}.

Suppose we are given an algebra structure on~$E$; we want to describe 
$\bbZ_p$-structures on~$E$, i.e.\ those sublattices~$A$ of $E$ which are 
stable under multiplication of~$E$. Of course, we would like to obtain 
a description in terms of corresponding subgroups $G_A$.

\nxsubpoint
Let's try the $p$-adic case first. So $E$ is a finite $\bbQ_p$-algebra, 
and we are looking for finite flat $\bbZ_p$-algebras $A$, such that 
$A_{(\bbQ_p)}\cong E$; in other words, $A$ must be a sublattice of $E$ 
stable under multiplication: $A\cdot A=A$.

\nxsubpoint\label{ss:propsconstr}
From the abstract point of view, a $\bbZ_p$-algebra $A$ is a $\bbZ_p$-module 
$A$ equipped by two morphisms $\eta_A:\bbZ_p\to A$ (the unit) and 
$\mu_A:A\otimes A\to A$ (the multiplication), subject to some conditions. 
If $A_{(\bbQ_p)}\cong E$ as a $\bbQ_p$-algebra, 
then we must have $(\eta_A)_{(\bbQ_p)}=\eta_E$ and similarly for the 
multiplication. 
Hence, if we start from several automorphisms 
$\gamma_1$, \dots, $\gamma_n\in GL(E)$ that preserve~$A$, and construct from 
them some new element $\gamma$ of $GL(E)$ using only tensor products and 
maps $\mu_E$ and $\eta_E$, the resulting element must also preserve $A$, since 
the same construction can be done first in the category of $\bbZ_p$-modules, 
yielding some $\gamma':A\to A$, and then $\gamma=\gamma'_{(\bbQ_p)}$ will 
preserve~$A$.

\nxsubpoint\label{ss:defgamma12}
Proceeding in this way, for any two $\gamma_1,\gamma_2:E\to E$ we define 
two new maps $\gamma_1\llcorner\gamma_2$ and $\gamma_1\lrcorner\gamma_2$ 
from~$E$ to~$E$:
\begin{equation}
\xymatrix{
\gamma_1\llcorner\gamma_2:
E \ar[rr]^{1_E\otimes\eta_E} && 
E\otimes E \ar[r]^{\gamma_1\otimes\gamma_2} & 
E\otimes E \ar[r]^<>(.5){\mu_E} & E}
\end{equation}
\begin{equation}\xymatrix{
\gamma_1\lrcorner\gamma_2:
E \ar[rr]^{\eta_E\otimes 1_E} && 
E\otimes E \ar[r]^{\gamma_1\otimes\gamma_2} & 
E\otimes E \ar[r]^<>(.5){\mu_E} & E}
\end{equation}
In other words, $(\gamma_1\llcorner\gamma_2)(x)=\gamma_1(x)\cdot\gamma_2(1)$ 
and $(\gamma_1\lrcorner\gamma_2)(x)=\gamma_1(1)\cdot\gamma_2(x)$. 

\nxsubpoint\label{ss:gamma12stab} 
Clearly, these are examples of constructions considered 
in~\ptref{ss:propsconstr}, i.e.\ if both $\gamma_1$ and $\gamma_2$ preserve 
some $\bbZ_p$-subalgebra $A\subset E$, the same is true for 
$\gamma_1\llcorner\gamma_2$ and $\gamma_1\lrcorner\gamma_2$.

So we are tempted to describe sublattices $A\subset E$ that are 
$\bbZ_p$-subalgebras by requiring $\gamma_1\llcorner\gamma_2$, 
$\gamma_1\lrcorner\gamma_2\in G_A$ for all $\gamma_1$, $\gamma_2\in G_A$.

However, these two new elements need not be invertible, i.e.\ in some 
cases they are not in $G=GL(E)$. Even when they lie in~$G$, all we can 
say about them is that $(\gamma_1\llcorner\gamma_2)(A)\subset A$ and 
similarly for $\gamma_1\lrcorner\gamma_2$, 
condition insufficient to conclude $\gamma_1\llcorner\gamma_2\in G_A$.

Observe that we can avoid both these problems by considering from the 
very beginning the {\em compact submonoids\/} 
$M_A:=\{g\in\End(E):g(A)\subset A\}\cong\End_{\bbZ_p}(A)$ of $\End(E)$.

Then conditions $\gamma_1, \gamma_2\in M_A$ imply $\gamma_1\llcorner\gamma_2$, 
$\gamma_1\lrcorner\gamma_2\in M_A$ whenever the sublattice $A\subset E$ is 
a $\bbZ_p$-algebra.

\nxpoint
Before proceeding further we would like to check whether these $M_A$ 
still determine sublattices $A\subset E$ up to similitude, and whether 
they can be described as maximal compact submonoids of $\End(E)$. The 
answer to both these questions is {\em positive,} so the situation 
is completely similar to one we had before when we considered 
maximal compact subgroups of~$GL(E)$:

\label{th:maxsubmpadic}
\begin{Thz}
Let $E$ be a finite-dimensional $\bbQ_p$-vector space. For any 
$\bbZ_p$-sublattice $A\subset E$ denote by $M_A$ 
the submonoid of $\End(E)$ consisting of endomorphisms of~$E$ that preserve 
$A$, i.e.\ $M_A=\{g\in\End(E):g(A)\subset A\}$. Then:
\begin{itemize}
\item[a)] The $M_A$ are maximal compact submonoids in $\End(E)$ (with respect to 
the $p$-adic topology);
\item[b)] All maximal compact submonoids of $\End(E)$ are of this form;
\item[c)] $M_A=M_{A'}$ iff $A$ and $A'$ are similar, i.e.\ iff 
$A'=c\cdot A$ for some $c\in\bbQ_p^*$.
\end{itemize}
\end{Thz}

\nxsubpoint
Let us check first that any compact submonoid $M\subset\End(E)$ stabilizes 
some lattice~$A$ in~$E$, i.e.\ that $M\subset M_A$. For this take any 
sublattice $A_0\subset E$ and consider the $\bbZ_p$-submodule 
$A$ of $E$ generated by $S:=M\cdot A_0=\{g(x):g\in M, x\in A_0\}$. 
Note that both $M$ and $A_0$ are compact, hence $S$ is compact, 
hence $S$ and $A$ are contained in a compact sublattice 
$p^{-n}A_0$ for some $n>0$, hence the closure $\bar{A}$ is a compact 
$\bbZ_p$-submodule of~$E$. On the other hand, it contains $A_0$, hence 
it is open, and all open compact $\bbZ_p$-submodules in $E$ are 
lattices. So we have constructed a sublattice $\bar{A}$ in $E$, stable 
under all elements $g\in M$, since $gS=gM\cdot A_0\subset M\cdot A_0=S$. 
In other words, we have $M\subset M_{\bar{A}}$.

\nxsubpoint
Clearly, if $A$ and $A'$ are similar, then $M_A=M_{A'}$. Let us prove that 
$M_A\subset M_{A'}$ implies that $A$ and $A'$ are similar, hence  
$M_A=M_{A'}$; this will complete the proof of the theorem.

So suppose that $M_A\subset M_{A'}$. Case $E=0$ is trivial, so assume $E\neq0$.
After rescaling $A'$, we can assume that 
$A'$ is contained in $A$, but not in $pA$. Choose any element 
$u\in A'-pA\subset A-pA$. According to Lemma~\ptref{l:transonlat} below, 
for any element $v\in A$ we can find $g\in M_A$, for which $g(u)=v$.
Then, since $u\in A'$ and $g\in M_A\subset M_{A'}$, we have $v=g(u)\in A'$. 
So we have proved $A\subset A'$, hence $A=A'$.

\begin{LemmaD}\label{l:transonlat}
If $A$ is a sublattice in $E$, then for any $u\in A-pA$, $v\in A$ 
there is an element $g\in M_A\cong\End_{\bbZ_p}(A)$, such that $g(u)=v$. 
In other words, for any $u\in A-pA$ we have $A=M_A\cdot u=\{g(u):g\in M_A\}$.
\end{LemmaD}

\begin{Proof}
Indeed, $A$ is a free $\bbZ_p$-module of finite rank, so $u\in A-pA$ implies 
that we can complete $u$ to some base $e_1=u$, $e_2$, \dots, $e_n$ of $A$. 
Then we can define an element $g\in\End_{\bbZ_p}(A)$ with required properties 
by putting $g(e_1):=v$, $g(e_k):=0$ for $k>1$.
\end{Proof}

\nxsubpoint\label{th:padicsubalg}
Now we would like to check whether the conditions 
$\gamma_1\llcorner\gamma_2$, $\gamma_1\lrcorner\gamma_2\in M_A$ for all
$\gamma_1$, $\gamma_2\in M_A$, necessary for a sublattice $A$ in a 
finite-dimensional $\bbQ_p$-algebra~$E$ to be a 
$\bbZ_p$-subalgebra, are in fact sufficient. 

The answer is again positive, modulo some rescaling:

\begin{Thz}
Let $A$ be a $\bbZ_p$-sublattice in a finite-dimensional 
$\bbQ_p$-algebra~$E$. Suppose that $1\in A-pA$ (any lattice is 
similar to exactly one lattice of this sort). 
The following conditions are equivalent:
\begin{itemize}
\item[(i)] $A$ is a $\bbZ_p$-subalgebra in~$E$;
\item[(ii)] For any $\gamma_1$ and $\gamma_2$ in $M_A$ the elements 
$\gamma_1\llcorner\gamma_2$ and $\gamma_1\lrcorner\gamma_2$ defined in 
\ptref{ss:defgamma12} also belong to $M_A$;
\item[(iii)] For any $\gamma\in M_A$ we have $1_E\llcorner\gamma\in M_A$,
where $1_E$ denotes the identity of $\End(E)$.
\end{itemize}
\end{Thz}

\begin{Proof}
(i)$\Rightarrow$(ii) has been explained in~\ptref{ss:gamma12stab}, and 
(ii)$\Rightarrow$(iii) is evident since $1_E\in M_A$. Let's prove 
(iii)$\Rightarrow$(i). Let $u$ and $v$ be two elements of $A$; we want to 
prove that $uv\in A$. According to Lemma~\ptref{l:transonlat}, we can 
find an element $\gamma\in M_A$, such that $\gamma(1)=v$. By assumption 
$1_E\llcorner\gamma\in M_A$, hence $uv=1_E(u)\cdot\gamma(1)=
(1_E\llcorner\gamma)(u)$ belongs to~$A$, q.e.d.
\end{Proof}

\nxpointtoc{Maximal compact submonoids of $\End(E)$}
Now it is natural to discuss the archimedian case. Let $E$ be a 
finite-dimensional real vector space. Any positive-definite quadratic 
form~$Q$ on~$E$ defines a maximal compact submonoid 
$M_Q:=\{g\in\End(E):\forall x, Q(g(x))\leq Q(x)\}$. 
Is it true that all maximal compact submonoids of $\End(E)$ are of this form?

Contrary to what we might have expected after looking to the $p$-adic case, the answer is {\em negative}. There are many other maximal compact submonoids of 
$\End(E)$. Now we are going to describe them.

\nxsubpoint\label{sp:convexdefs}
Recall that a subset $A$ of a (not necessarily finite-dimensional) 
real vector space~$E$ is said to be {\em convex\/}, if 
$\lambda x_1+(1-\lambda)x_2\in A$ whenever 
$x_1$, $x_2\in A$ and $0\leq\lambda\leq 1$. This condition is 
equivalent to $A$ being closed under baricentric combinations of its 
points, i.e.\ $\lambda_1 x_1+\cdots+\lambda_n x_n\in A$ whenever 
$x_1$, \dots, $x_n\in A$, all $\lambda_k\geq 0$ and $\sum_k\lambda_k=1$.

For any subset $S\subset E$ there is a smallest convex subset 
$\conv(S)\subset E$, containing $S$. This subset is called the 
{\em convex hull\/} of~$S$; note that it can be described as the set of 
all baricentric combinations of elements of~$S$.

We say that $A$ is {\em (central) symmetric}, if $A=-A$, 
i.e.\ if it is stable under multiplication by $-1$.

Finally, we say that $A$ is a {\em convex body}, if it is convex, and if its 
affine span (i.e.\ the smallest affine subspace of~$E$ containing~$A$) 
is equal to~$E$. For finite-dimensional~$E$ this is equivalent to the 
non-emptiness of the interior of~$A$.

A symmetric convex set $A$ is a convex body iff its linear span 
$\bbR\cdot A = \bigcup_{\lambda\geq0}\lambda A$ is equal to the whole of~$E$. 
Subsets $A\subset E$ with the latter property are also called {\em absorbent}.
In the finite-dimensional case this is equivalent to $0$ being an interior 
point of~$A$, i.e.\ to $A$ being a neighborhood of the origin.

In the infinite-dimensional case a similar description of convex bodies 
can be obtained, provided we equip~$E$ with its ``algebraic'' topology, i.e.\ 
with its finest locally convex topology (cf.~\cite{EVT}, ch.~II, \S4, n.~2).

\nxsubpoint\label{p:norms}
Note that compact symmetric convex bodies~$A$ in a finite-dimensional 
real space~$E$ are in one-to-one correspondence with {\em norms\/} on~$E$, 
i.e.\ maps $\vnorm:E\to\bbR_{\geq0}$ satisfying following conditions:
\begin{enumerate}
\item $\|\lambda x\|=|\lambda|\cdot\|x\|$ for any $\lambda\in\bbR$, $x\in E$;
\item $\|x+y\|\leq\|x\|+\|y\|$ for any $x,y\in E$;
\item $\|x\|=0$ iff $x=0$.
\end{enumerate}

Indeed, any such norm defines a compact symmetric convex body, namely 
$A_{\vnorm}:=\{x\in E:\|x\|\leq 1\}$. Conversely, any such body~$A$ defines a 
norm~$\vnorm_A$, namely $\|x\|_A:=\inf\{\lambda>0:\lambda^{-1}x\in A\}$.

Recall that this is a one-to-one correspondence, and that all these norms 
define the same topology on~$E$, hence all of them are continuous.

Now we are in position to state an analogue of Theorem~\ptref{th:maxsubmpadic}:

\begin{ThD}\label{th:maxsubmrealc}
Let $E$ be a finite-dimensional real vector space. For any symmetric compact 
convex body $A\subset E$ we denote by~$M_A$ the compact submonoid 
of\/ $\End(E)$ defined by $M_A:=\{g\in\End(E) : g(A)\subset A\}$. Then:
\begin{itemize}
\item[a)] The $M_A$ are maximal compact submonoids of\/~$\End(E)$ 
(with respect to the real topology);
\item[b)] All maximal compact submonoids of\/~$\End(E)$ are of this form;
\item[c)] $M_A=M_{A'}$ iff $A$ and $A'$ are similar, i.e.\ iff 
$A'=c\cdot A$ for some $c\in\bbR^*$.
\end{itemize}
\end{ThD}

\nxsubpoint
Let us check that any compact submonoid $M\subset\End(E)$ stabilizes 
some symmetric compact convex body $A$, i.e.\ that $M\subset M_A$. 
For this take any such body $A_0\subset E$ (e.g.\ the closed unit ball with 
respect to some metric on~$E$) and consider the convex hull~$A$ of 
the set $S:=M\cdot A_0=\{g(x):g\in M, x\in A_0\}$. Note that 
$S$ is compact since both $M$ and $A_0$ are compact, hence $S$ is bounded, 
hence the same is true for its convex hull~$A$. Now $\bar{A}$ is a compact 
convex subset; it is symmetric since $S=-S$, and it is a neighborhood of 
the origin since $\bar{A}\supset A_0$. So we have constructed a 
symmetric compact convex body $\bar{A}$, clearly stable under $M$ 
since $gS=gMA_0\subset MA_0=S$ for any $g\in M$.

\nxsubpoint
Clearly, if $A$ and $A'$ are similar, then $M_A=M_{A'}$. Let us prove 
that $M_A\subset M_{A'}$ implies that $A$ and $A'$ are similar, hence 
$M_A=M_{A'}$; this will complete the proof of the theorem.

So suppose that $M_A\subset M_{A'}$. Case $E=0$ is trivial, so we 
assume $E\neq0$. Put $\lambda:=\sup_{x\in A'}\|x\|_A$. Since $A'$ is 
compact and $\vnorm_A$ continuous, this supremum is actually achieved 
at some point $u\in A'$, hence $\lambda=\|u\|_A$ is finite and $>0$. 
After rescaling $A'$, we can assume $\lambda=1$, i.e.\  
$A'\subset A$, and $\|u\|_A=\lambda=1$. 
Now, according to Lemma~\ptref{l:transonconv} below, 
for any element $v\in A$ we can find some $g\in M_A$, for which $g(u)=v$.
Then, since $u\in A'$ and $g\in M_A\subset M_{A'}$, we have $v=g(u)\in A'$. 
So we obtain $A\subset A'$, hence $A=A'$.

\begin{LemmaD}\label{l:transonconv}
If $A$ is a symmetric convex body in~$E$, then for any two 
$u$, $v\in A$, such that $\|u\|_A=1$, there is an endomorphism 
$g\in M_A\subset\End(E)$, such that $g(u)=v$.
\end{LemmaD}

\begin{Proof}
Indeed, our condition for $u$ means that $u$ belongs to the boundary 
$\partial A$ of~$A$, i.e.\ that it doesn't belong to the interior~$A^0$ 
of~$A$. 
Hence by Hahn--Banach theorem (cf.\ \cite{EVT}, ch.~II, \S3, cor.~2 of th.~1) 
there is a linear form $\phi:E\to\bbR$, such that $\phi(x)\leq\|x\|_A$ for 
all $x\in E$ and $\phi(u)=\|u\|_A=1$, hence $-1\leq\phi(A)\leq1=\phi(u)$. 
We can define our endomorphism~$g$ 
by $g(x) := \phi(x)\cdot v$. The image of~$A$ under~$g$ is 
$[-1,1]\cdot v\subset A$, hence $g\in M_A$, and by construction $g(u)=v$.
\end{Proof}

\nxpointtoc{Category of $\Zinfty$-lattices}
Thus we are led to define a {\em $\Zinfty$-lattice} in a finite-dimensional 
real vector space~$E$ as a symmetric compact convex body $A$ in $E$: 

\begin{DefD}\label{def:zinflat}
We define the category of $\Zinfty$-lattices $\ZinfLat$ as follows. 
Its objects are pairs $A=(A_{\Zinfty},A_\bbR)$, consisting of a 
finite-dimensional real vector space $A_\bbR$ and a compact 
symmetric convex body $A_{\Zinfty}\subset A_\bbR$; 
when no confusion can arise, we denote 
$A_{\Zinfty}$ by the same letter~$A$ as the whole pair.

The morphisms from $A=(A_{\Zinfty},A_\bbR)$ to $B=(B_{\Zinfty},B_\bbR)$ 
are pairs $f=(f_{\Zinfty}, f_\bbR)$, where $f_\bbR: A_\bbR\to B_\bbR$ is 
an $\bbR$-linear map, and $f_{\Zinfty}: A_{\Zinfty}\to B_{\Zinfty}$ is 
required to coincide with the restriction of~$f_\bbR$ to $A_{\Zinfty}$. 
When no confusion can arise, we denote $f_{\Zinfty}$ simply by~$f$.

Composition of morphisms is defined in the natural way.
\end{DefD}

Note that a morphism $f=(f_{\Zinfty},f_\bbR)$ is determined by any of 
its components, since $f_{\Zinfty}$ is the restriction of $f_\bbR$ to 
$A_{\Zinfty}$; on the other hand, $A_{\Zinfty}$ is absorbent, i.e.\ 
any element of $A_\bbR$ is of form $\lambda x$ for some $\lambda>0$, 
$x\in A_{\Zinfty}$, hence $f_\bbR(\lambda x)=\lambda f_{\Zinfty}(x)$.

We could also describe $\Hom_\ZinfLat(A,B)$ as the set of $\bbR$-linear maps 
$f:A_\bbR\to B_\bbR$, such that $f(A)\subset B$. Another possible 
description is this: $\ZinfLat$ is the category of finite-dimensional 
normed real spaces, and the morphisms $f: A \to B$ are the 
$\bbR$-linear maps $A_\bbR\to B_\bbR$ of norm $\leq1$, i.e.\ we require 
$\|f(x)\|_B \leq \|x\|_A$ for all $x\in A_\bbR$.

This choice of morphisms is motivated by the $p$-adic case, in which 
we had $M_A\cong\End_{\bbZ_p}(A)$ for any $\bbZ_p$-lattice $A$ in a 
$\bbQ_p$-vector space. In our case we also have 
$\End_\ZinfLat(A)\cong M_A$, where $M_A$ is the corresponding maximal 
compact submonoid in $\End(A_\bbR)$ (cf.~\ptref{th:maxsubmrealc}).

\nxsubpoint
We are going to embed $\ZinfLat$ into the category of ``flat $\Zinfty$-modules''
$\ZinfFlat$ as a full subcategory, and later we will embed this latter 
into the category of (all) $\Zinfty$-modules $\catMod\Zinfty$, again as 
a full subcategory. That's why we will say that morphisms $f:A\to B$ of 
$\ZinfLat$ are {\em $\Zinfty$-linear maps} or {\em $\Zinfty$-homomorphisms}, 
and will denote $\Hom_{\ZinfLat}(A,B)$ simply by $\Hom_{\Zinfty}(A,B)$.

\nxsubpoint\label{p:limszinflat}
Note that $\ZinfLat$ is already a nicer category than the category of 
quadratic vector spaces. At least, it has finite inductive and projective 
limits, and in particular finite direct sums (i.e.\ coproducts) and products. 
For example, the product $A\times B$ is equal to 
$(A_{\Zinfty}\times B_{\Zinfty}, A_\bbR\times B_\bbR)$, and the 
direct sum $A\oplus B$ is given by $(\conv(A_{\Zinfty}\cup B_{\Zinfty}), 
A_\bbR\oplus B_\bbR)$. Note that both $\bbR$-vector spaces $(A\times B)_\bbR$ 
and $(A\oplus B)_\bbR$ can be identified with $A_\bbR\oplus B_\bbR$, 
while the corresponding convex subsets are different: all we can say is 
$(A\oplus B)_{\Zinfty}\subset (A\times B)_{\Zinfty}$, i.e.\ we have a 
canonical {\em monomorphism} $A\oplus B\to A\times B$. The fact that 
this is a monomorphism will be explained later as a manifestation 
of {\em hypoadditivity} of $\Zinfty$. Note that in an additive category it 
would have been an isomorphism, hence $\ZinfLat$ is not additive.

In the language of normed vector spaces the direct sum corresponds to norm 
$\|x+y\|_{A\oplus B}=\|x\|_A+\|y\|_B$, and the direct product to norm
$\|x+y\|_{A\times B}=\sup(\|x\|_A,\|y\|_B)$ on the same vector space 
$A_\bbR\oplus B_\bbR$, where we take $x\in A_\bbR$, $y\in B_\bbR$.

\nxpoint
Before proceeding further let us study the $\Zinfty$-lattices $A$ 
which are stable under multiplication in a given finite-dimensional 
real algebra~$E$. This should give us an understanding of finite flat 
$\Zinfty$-algebras. Here is an analogue of Theorem~\ptref{th:padicsubalg}:

\label{th:realsubalg}
\begin{Thz}
Let $A$ be a $\Zinfty$-sublattice in a finite-dimensional 
$\bbR$-algebra~$E$, i.e.\ $A\subset E$ is a symmetric compact convex body. 
Suppose that $1\in \partial A$, i.e.\ $\|1\|_A=1$ (any $\Zinfty$-lattice in~$E$
is similar to exactly one $\Zinfty$-lattice of this sort). 
The following conditions are equivalent:
\begin{itemize}
\item[(i)] $A$ is closed under multiplication of~$E$, i.e.\ 
$\mu_E(A\times A)\subset A$, or, in other words, $A$ is a multiplicative 
submonoid of~$E$;
\item[(ii)] For any $\gamma_1$ and $\gamma_2$ in $M_A$ the elements 
$\gamma_1\llcorner\gamma_2$ and $\gamma_1\lrcorner\gamma_2$ defined in 
\ptref{ss:defgamma12} also belong to $M_A$;
\item[(iii)] For any $\gamma\in M_A$ we have $1_E\llcorner\gamma\in M_A$,
where $1_E$ denotes the identity of $\End(E)$.
\item[(iv)] $E$ is a Banach algebra with respect to $\vnorm_A$, i.e.\ 
$\|1\|_A=1$ and ${\|xy\|}_A\leq{\|x\|}_A\cdot{\|y\|}_A$ for all $x,y\in E$.
\end{itemize}
\end{Thz}

\begin{Proof}
(i)$\Rightarrow$(ii): Clear from definitions (cf.~\ptref{ss:defgamma12}); note that we use $1\in A$ here.
(ii)$\Rightarrow$(iii): Evident since $1_E\in M_A$.
(iii)$\Rightarrow$(i): Same as in the proof of~\ptref{th:padicsubalg}, 
with Lemma~\ptref{l:transonlat} replaced by \ptref{l:transonconv}.
(iv)$\Rightarrow$(i): Follows from $A=\{x:\|x\|_A\leq1\}$.
(i)$\Rightarrow$(iv): If one of $x$ and $y$ is zero, the statement is trivial; 
assume $x\neq 0$, $y\neq 0$. Put $x':=x/\|x\|_A$, $y':=y/\|y\|_A$. Then 
$\|x'\|_A=\|y'\|_A=1$, hence $x',y'\in A$, so by (i) we get $x'y'\in A$, i.e.\ 
$\|x'y'\|_A\leq1$, hence 
$\|xy\|_A={\|x\|}_A\cdot{\|y\|}_A\cdot{\|x'y'\|}_A\leq{\|x\|}_A\cdot{\|y\|}_A$.
\end{Proof}

In other words, the $\Zinfty$-lattices~$A$ of~$E$, which are at the same time 
``finite $\Zinfty$-subalgebras'', are exactly the symmetric compact convex 
absorbent (multiplicative) submonoids of~$E$.

\nxsubpoint
In particular, $\bbR$ contains exactly one such ``$\Zinfty$-subalgebra'', 
namely, $\Zinfty:=([-1,1], \bbR)$. This is our first definition 
of~$\Zinfty$ itself; later it will be improved.

\nxsubpoint
Note that multiplicative submonoids appeared again in the statement of 
\ptref{th:realsubalg}. This suggests some connection to maximal 
compact submonoids $M_A$ considered in~\ptref{th:maxsubmrealc}. Indeed, 
for any $\Zinfty$-lattice $A$ in~$E$, the corresponding maximal compact 
submonoid $M_A$ is easily seen to be symmetric and convex, hence 
it satisfies the conditions of~\ptref{th:realsubalg} for $\End(E)$, 
i.e.\ $M_A\cong\End_{\Zinfty}(A)$ is a ``finite $\Zinfty$-subalgebra'' 
of $\End(E)$, as one might have expected.

\nxsubpoint\label{ss:intcloszinf}
If the finite-dimensional $\bbR$-algebra $E$ is commutative and semisimple, 
all the ``finite $\Zinfty$-subalgebras'' in $E$ are contained in a 
maximal one $A_{max}$; it is natural to call this $A_{max}$ the 
{\em integral closure of $\Zinfty$ in $E$}. To construct this 
``integral closure'' consider all embeddings $\sigma_i:E\to\bbC$; then 
$A_{max}=\{x\in E:\forall i, |\sigma_i(x)|\leq 1\}$. Indeed, 
$A_{max}$ is a closed symmetric convex body, stable under multiplication; 
it is compact since $\bigcap_i\Ker\sigma_i=0$, so it is a ``finite 
$\Zinfty$-subalgebra'' of~$E$. If $A$ is another one, then 
by compactness of~$A$, there is some $C>0$, such that $|\sigma_i(x)|\leq C$ 
for all $x\in A$ and $i$. Since $A$ is closed under multiplication, 
we get $|\sigma_i(x^n)|=|\sigma_i(x)|^n\leq C$ for any $n\geq1$, 
hence $|\sigma_i(x)|\leq1$ for each~$i$, i.e.\ $A\subset A_{max}$.

\nxsubpoint
In particular, we can apply this to $E=\bbC$. The integral 
closure of $\Zinfty$ in $\bbC$ will be denoted by 
$\barZinfty$ (cf.~\ptref{ss:firstbarzinfty}). So we have 
$\barZinfty=(\{z\in\bbC:|z|\leq1\}, \bbC)$, with the multiplication 
induced by that of~$\bbC$.

\nxsubpoint
Any commutative semisimple $\bbR$-algebra $E$ is isomorphic to 
$\bbR^r\times\bbC^s$ for some $r,s\geq0$. 
The considerations of \ptref{ss:intcloszinf} show that the integral closure 
of $\Zinfty$ in $E$ coincides with $\Zinfty^r\times\barZinfty^s$, as 
one might have expected.

\nxsubpoint
Resuming the considerations of \ptref{ss:bunoncompz} in our new context, 
we see that a vector bundle~$\sE$ on $\CompZ$ corresponds to 
a free $\bbZ$-module $E_\bbZ$ of finite rank, together with 
a $\Zinfty$-lattice, i.e.\ a symmetric compact convex body~$A$, in 
$E_\bbR:=E_\bbZ\otimes\bbR$. In other words, we are given both a 
$\bbZ$-lattice $E_\bbZ$ and a symmetric convex body~$A$ in the same 
real vector space~$E_\bbR$. The ``global sections'' of $\sE$ correspond 
to the intersection $E_\bbZ\cap A$; this is a finite set 
(``an $\Fone$-module''), and there are a lot of interesting theorems 
offering bounds for the number of lattice points inside a convex set, 
starting with Minkowski theorem. So one might think of these theorems 
as some ``Riemann--Roch theorems'' for $\CompZ$.

\nxsubpoint
In particular, we have the trivial line bundle $\sO_{\CompZ}$, defined 
by $\bbZ$ and $\Zinfty$ in $\bbR$. If $K/\bbQ$ is a number field, 
we can consider the normalization of $\CompZ$ in~$K$; as a 
vector bundle over $\CompZ$, it will be given by the integral closure 
$\sO_K$ of $\bbZ$ in $K$, and by the integral closure $\sO_{K,\infty}$ 
of $\Zinfty$ in $K\otimes_\bbQ\bbR$. We have a description 
of $\sO_{K,\infty}$ in~\ptref{ss:intcloszinf}; when we compute the 
``global sections'' of our vector bundle over $\CompZ$, we get 
$\sO_K\cap\sO_{K,\infty}=\{x\in\sO_K:|\sigma(x)|\leq1$ for any 
$\sigma:K\to\bbC\}$. By classical algebraic number theory this set 
is known to be the union of zero and the subgroup of roots of unity 
in~$\sO_K$, which is a finite cyclic group.

\nxsubpoint\label{ss:maxsubmcomplc}
What about a complex analogue of~\ptref{th:maxsubmrealc}? In other 
words, what are the maximal compact submonoids of $\End_\bbC(E)$, when 
$E$ is a finite-dimensional $\bbC$-vector space? One can check that 
an analogue of~\ptref{th:maxsubmrealc} still holds, provided we consider 
{\em balanced\/} compact convex bodies $A\subset E$, i.e.\ 
we require $\lambda A\subset A$ for any $\lambda\in\bbC$ with $|\lambda|\leq1$.
In other words, $\barZinfty\cdot A\subset A$, i.e.\ $A$ is a 
``$\barZinfty$-module''. 

Note that balanced compact convex bodies $A\subset E$ are 
in one-to-one correspondence with the {\em complex\/} norms $\vnorm$ on~$E$, 
i.e.\ norms on~$E$, such that $\|\lambda x\|=|\lambda|\cdot\|x\|$ for any 
$\lambda\in\bbC$, $x\in E$.

Then we define $M_A:=\{g\in\End_\bbC(E):g(A)\subset A\}$ as before, 
and these are precisely the maximal compact submonoids of $\End_\bbC(E)$. 
The proof is essentially that of~\ptref{th:maxsubmrealc}, but we 
have to apply complex Hahn--Banach in the lemma (cf.~\cite{EVT}, 
ch.~II, \S8, n.~3, cor.~1 of th.~1).

We do not develop this point of view further since these 
``$\barZinfty$-lattices'' will be later described as $\Zinfty$-lattices 
equipped with a $\barZinfty$-module structure, hence their theory will 
be a particular case of the theory of modules over a $\Zinfty$-algebra.

\nxsubpoint
We have already observed (cf.~\ptref{ss:intcloszinf}), that  
any {\em commutative\/} semisimple $\bbR$-algebra~$E$ contains a unique 
maximal compact submonoid (note that it automatically will be 
symmetric and convex, otherwise we would replace it by the closure of 
its symmetric convex hull). What about non-commutative 
semisimple $\bbR$-algebras? Clearly, it is enough to consider the case 
of a simple $\bbR$-algebra~$E$, and there are three kinds of these.
If $E=M(n,\bbR)=\End_\bbR(\bbR^n)$, its maximal compact submonoids 
correspond to $\Zinfty$-lattices in~$\bbR^n$ (cf.~\ptref{th:maxsubmrealc}), 
and the maximal compact submonoids of $E=M(n,\bbC)=\End_\bbC(\bbC^n)$ 
correspond to $\barZinfty$-lattices in~$\bbC^n$ (cf.~\ptref{ss:maxsubmcomplc}).

So only the quaternionic case $E=M(n,\bbH)$ remains. Note that  
there is the largest compact submonoid (hence a $\Zinfty$-subalgebra) 
$\Zinfty^\bbH:=\{q\in\bbH:|q|\leq1\}$ in~$\bbH$, 
where $|q|:=\sqrt{q\bar{q}}$ is 
the quaternionic absolute value. Then it seems quite plausible  
that maximal compact submonoids of $M(n,\bbH)\cong\End_\bbH(\bbH^n)$ 
are exactly the stabilizers $M_A$ of those $\Zinfty$-lattices $A\subset\bbH^n$
which are $\Zinfty^\bbH$-submodules, i.e.\ $\Zinfty^\bbH\cdot A\subset A$. 
In the language of norms they correspond to quaternionic norms $\vnorm$ 
on $\bbH^n$, i.e.\ we require $\|\lambda x\|=|\lambda|\cdot\|x\|$ for any 
$\lambda\in\bbH$, $x\in\bbH^n$. Actually the proof of~\ptref{th:maxsubmrealc} 
still works, provided we have a quaternionic Hahn--Banach theorem.

\nxpoint
When $A$ is a $\bbZ_p$-lattice in a finite-dimensional $\bbQ_p$-algebra~$E$, 
stable under multiplication $\mu_E:E\times E\to E$ of~$E$, then 
$A$ is a $\bbZ_p$-algebra, with the multiplication $\mu_A:A\times A\to A$ 
induced by that of~$E$. Note that $\mu_A$, being the restriction 
of $\mu_E$, is automatically a $\bbZ_p$-bilinear map. 
Conversely, if we are given a $\bbZ_p$-bilinear map $\mu_A:A\times A\to A$ 
on some $\bbZ_p$-lattice~$A$, subject to some additional conditions, 
we obtain a $\bbZ_p$-algebra structure on~$A$.

In other words, a finite flat $\bbZ_p$-algebra can be roughly described as 
a $\bbZ_p$-lattice~$A$ equipped by a $\bbZ_p$-bilinear map~$\mu_A$. We would like 
to obtain a similar description of finite flat $\Zinfty$-algebras. 
For this we have to define $\Zinfty$-bilinear maps as follows:

\begin{DefD}\label{def:zinflatbilin}
Let $M$, $N$, $P$ be three $\Zinfty$-lattices (i.e.\ three objects of 
$\ZinfLat$). A {\em $\Zinfty$-bilinear map $M\times N\to P$} is by 
definition a pair $\Phi=(\Phi_{\Zinfty},\Phi_\bbR)$, consisting of an 
$\bbR$-bilinear map $\Phi_\bbR:M_\bbR\times N_\bbR\to P_\bbR$, and 
a map of sets $\Phi_{\Zinfty}:M_{\Zinfty}\times N_{\Zinfty}\to P_{\Zinfty}$, 
required to coincide with the restriction of~$\Phi_\bbR$. When no 
confusion can arise, we denote $\Phi_{\Zinfty}$ by the same letter~$\Phi$.
We denote by $\Bilin_{\Zinfty}(M,N;P)$ the set of all $\Zinfty$-bilinear maps 
$M\times N\to P$. 
\end{DefD}

\nxsubpoint
One can make for $\Zinfty$-bilinear map remarks similar to those 
made for $\Zinfty$-linear maps (i.e.\ morphisms of $\ZinfLat$) 
after~\ptref{def:zinflat}. In particular, $\Phi=(\Phi_{\Zinfty},\Phi_\bbR)$ 
is completely determined by either of its components. It can be described 
as an $\bbR$-linear map~$\Phi_\bbR:M_\bbR\times N_\bbR\to P_\bbR$, 
such that $\Phi_\bbR(M_{\Zinfty}\times N_{\Zinfty})\subset P_{\Zinfty}$. 
In the language of normed vector spaces this means that we consider 
$\bbR$-bilinear maps of norm $\leq 1$, i.e.\ ${\|\Phi_\bbR(x,y)\|}_P\leq 
{\|x\|}_M\cdot{\|y\|}_N$ for all $x\in M_\bbR$, $y\in N_\bbR$.

Note that, similar to what we have for modules over a commutative ring, 
a $\Zinfty$-bilinear form can be described as a map 
$\Phi:M\times N\to P$ (i.e.\ $M_{\Zinfty}\times N_{\Zinfty}\to P_{\Zinfty}$), 
such that for any $x\in M$ the map $s_\Phi(x):y\mapsto\Phi(x,y)$ is 
$\Zinfty$-linear (i.e.\ extends to a morphism $N\to P$ in $\ZinfLat$), 
as well as the map $d_\Phi(y):x\mapsto\Phi(x,y)$ for any $y\in N$.

\begin{PropD}\label{p:extensprod} (Existence of tensor products.) 
For any two $\Zinfty$-lattices $M$ and $N$ the functor 
$\Bilin_{\Zinfty}(M,N;-)$ is representable by some $\Zinfty$-lattice 
$M\otimes_{\Zinfty}N$, i.e.\ we have $\Bilin_{\Zinfty}(M,N;P)\cong
\Hom_{\Zinfty}(M\otimes_{\Zinfty}N,P)$, functorially in~$P$. 
In other words, we have a universal $\Zinfty$-bilinear map 
$\Phi_0:M\times N\to M\otimes_{\Zinfty}N$, such that for any other 
$\Zinfty$-bilinear map $\Psi:M\times N\to P$ there is a unique 
$\Zinfty$-linear map $\psi:M\otimes_{\Zinfty}N\to P$ satisfying 
$\Psi=\psi\circ\Phi_0$.
\end{PropD}

\begin{Proof}
We define our tensor product $M\otimes N$ in $\ZinfLat$ by 
$(M\otimes M)_\bbR := M_\bbR\otimes_\bbR N_\bbR$, 
$(M\otimes N)_{\Zinfty}:=\conv(S)$, where 
$S:=\{x\otimes y : x\in M_{\Zinfty}, y\in N_{\Zinfty}\}$. Clearly, 
$(M\otimes N)_{\Zinfty}$ is a symmetric convex body in $(M\otimes N)_\bbR$; 
it is compact, since it is the convex hull of a compact subset~$S$ of 
a finite-dimensional real vector space (this is a simple corollary of the  
Carath\'eodory theorem, which asserts that each point of $\conv(S)$ can be 
written as a baricentric combination of at most $\dim (M\otimes N)_\bbR+1$ 
points of~$S$).
Now the required universal property follows immediately 
from that of $M_\bbR\times N_\bbR\to (M\otimes N)_\bbR$, and from the 
definition of $(M\otimes N)_{\Zinfty}$.
\end{Proof}

\nxsubpoint\label{sp:defpolylin}
Sets of $\Zinfty$-polylinear maps $\Polylin_{\Zinfty}(M_1,M_2,\ldots, M_n;P)$ 
and multiple tensor products $M_1\otimes M_2\otimes\cdots\otimes M_n$ 
are defined essentially in the same way. The universal property of 
multiple tensor products gives us immediately a chain of isomorphisms 
$M_1\otimes(M_2\otimes M_3)\cong M_1\otimes M_2\otimes M_3\cong
(M_1\otimes M_2)\otimes M_3$, i.e.\ we have an associativity constraint 
for our tensor product. It clearly satisfies the pentagon axiom, since 
all vertices of the pentagon diagram are canonically isomorphic to 
a quadruple tensor product. Similarly, the natural commutativity 
isomorphisms $M_1\otimes M_2\cong M_2\otimes M_1$ obviously satisfy 
the hexagon axiom. Finally, $\Zinfty$ is a unit for this tensor product, 
i.e.\ $\Zinfty\otimes M\cong M$: this can be seen either from the 
explicit construction of $\Zinfty\otimes M$ given in~\ptref{p:extensprod}, 
or from its universal property: one has then to check 
$\Bilin_{\Zinfty}(\Zinfty,M;N)\cong\Hom_{\Zinfty}(M,N)$.

Therefore, $\ZinfLat$ is an ACU (associative, commutative, with unity) 
$\otimes$-category with respect to the tensor product introduced above. 
We could proceed and define $\Zinfty$-algebras and modules over these 
algebras in terms of this $\otimes$-structure; however, it will be 
more convenient to do this after extending $\ZinfLat$ to a larger category. 
It will be sufficient for now to remark that this definition of 
algebras in $\ZinfLat$ is equivalent to one we had in~\ptref{th:realsubalg}.

\nxsubpoint
Note that {\em inner Homs exist in $\ZinfLat$}. In other words, for any 
two $\Zinfty$-lattices $N$ and $P$ the functor 
$\Hom_{\Zinfty}(-\otimes N,P)\cong\Bilin_{\Zinfty}(-,N;P)$ is representable 
by some $\Zinfty$-lattice $\iHom_{\Zinfty}(N,P)$. This ``inner Hom'' 
$H:=\iHom_{\Zinfty}(N,P)$ can be constructed as follows. We put 
$H_\bbR:=\Hom_\bbR(N_\bbR, P_\bbR)$ and $H_{\Zinfty}:=
\{\phi:N_\bbR\to P_\bbR: \phi(N_{\Zinfty})\subset P_{\Zinfty}\}$. 
The corresponding 
norm $\vnorm_H$ on~$H_\bbR$ is defined by the usual rule for norms of 
linear maps between normed vector spaces: 
${\|\phi\|}_H=\sup_{{\|x\|}_N\leq 1}{\|\phi(x)\|}_P$. Since 
$H_{\Zinfty}=\{\phi:{\|\phi\|}_H\leq1\}$, this also shows that 
$H_{\Zinfty}$ is a compact symmetric convex body in $H_\bbR$, i.e.\ 
$H=(H_{\Zinfty},H_\bbR)$ is indeed an object of~$\ZinfLat$. Now the 
verification of the required property for $\iHom(N,P)$ is straightforward.

\nxsubpoint
This gives us a lot of canonical morphisms and isomorphisms, valid in 
any ACU $\otimes$-category with inner Homs. For example, 
\begin{align}
&\iHom(M,\iHom(N,P))\cong\iHom(M\otimes N,P)\cong\iHom(N,\iHom(M,P))\\
&\Gamma(\iHom(M,N))\cong\Hom_{\Zinfty}(M,N),
\quad \iHom(\Zinfty,M)\cong M
\end{align}
Here we denote by $\Gamma=\Gamma_{\Zinfty}$ the ``functor of global sections'' 
in $\ZinfLat$, given by $\Gamma(-):=\Hom(\Zinfty,-)$. Note that 
in our situation $\Gamma(M)=M_{\Zinfty}$ for any $M=(M_{\Zinfty},M_\bbR)$ in 
$\ZinfLat$.

\nxsubpoint
We can define ``inner Bilins'' as well: $\iBilin(M,N;P):=\iHom(M\otimes N,P)$.
Clearly, $\Gamma(\iBilin(M,N;P))=\Bilin(M,N;P)$. In the language of 
normed vector spaces $B:=\iBilin(M,N;P)$ corresponds to the space~$B_\bbR$ 
of $\bbR$-bilinear maps $\Phi:M_\bbR\times N_\bbR\to P_\bbR$, equipped by the 
usual norm: ${\|\Phi\|}_B=
\sup_{{\|x\|}_M\leq1,\,{\|y\|}_N\leq1}{\|\Phi(x,y)\|}_P$.

\nxsubpoint
Once we have the unit object and inner Homs, we can define {\em dual objects\/}
by $M^*:=\iHom_{\Zinfty}(M,\Zinfty)$. In our situation we get the dual 
convex body $M^*_{\Zinfty}$ to~$M_{\Zinfty}$ 
inside the dual space $(M^*)_\bbR=(M_\bbR)^*$. This corresponds to the 
usual norm ${\|\phi\|}_{M^*}=\sup_{{\|x\|}_M\leq 1}|\phi(x)|$ on the space 
of linear functionals $M^*_\bbR$.

We also have canonical $\Zinfty$-homomorphisms $M\to M^{**}$, which 
in our situation turn out to be isomorphisms, i.e.\ all $\Zinfty$-lattices are 
{\em reflexive}. In fact, this is an immediate consequence of the 
Hahn--Banach theorem, which says that the dual to the dual norm is the 
original one.

\nxsubpoint
Is the $\otimes$-category $\ZinfLat$ {\em rigid?} It already has inner Homs, 
and all its objects are reflexive. However, the remaining 
property fails: the map $M^*\otimes N\to\iHom(M,N)$ is 
a monomorphism, but usually not an isomorphism. To see this take 
$M=N=\Zinfty\oplus\Zinfty$; then the element $1_M\in\Hom(M,M)
=\Gamma(\iHom(M,M))$ doesn't come from any element of $\Gamma(M^*\otimes M)$, 
since any~$f$ from the image of this set has $|\tr f|\leq1$.

\nxsubpoint
We are going to embed $\ZinfLat$ into larger categories. However, 
$\ZinfLat$ can be always retrieved as the full subcategory consisting 
of the reflexive objects of any of these larger categories, subject 
to some additional finiteness conditions.

\nxpointtoc{Torsion-free $\Zinfty$-modules}
Now we want to construct the category of ``flat'' or ``torsion-free'' 
$\Zinfty$-modules $\ZinfFlat$, which will contain $\ZinfLat$ as 
a full subcategory.

\nxsubpoint
As usual, we consider the $p$-adic case first. How can one reconstruct 
the category of flat $\bbZ_p$-modules $\ZpFlat$, starting from 
its full subcategory of $\bbZ_p$-lattices (i.e.\ free modules of finite rank) 
$\ZpLat$? Observe that $\ZpFlat$ has arbitrary inductive limits, hence 
the inclusion $\ZpLat\to\ZpFlat$ induces a functor $\Ind(\ZpLat)\to\ZpFlat$ 
from the category of ind-objects over $\ZpLat$ (i.e.\ the category 
of filtered inductive systems of objects of $\ZpLat$; we refer to 
SGA~4~I for a more detailed discussion of ind-objects) to~$\ZpFlat$. 
This functor is essentially surjective, since any torsion-free $\bbZ_p$-module 
can be represented as a filtered inductive limit of its finitely generated 
submodules, which are automatically $\bbZ_p$-lattices. It is quite easy to 
check that it is in fact an equivalence of categories, i.e.\ 
$\ZpFlat\cong\Ind(\ZpLat)$.

\nxsubpoint
So we are tempted to construct $\ZinfFlat$ as $\Ind(\ZinfLat)$. However, 
we would like to have a more explicit description of this category, 
so we give a direct definition of a category $\ZinfFlat$ that contains 
$\ZinfLat$ as a full subcategory and is closed under arbitrary inductive 
limits. This will give a functor $\Ind(\ZinfLat)\to\ZinfFlat$. The 
reader is invited to check later, if so inclined, that this functor 
is indeed an equivalence of categories, to convince himself that 
the definition of $\ZinfFlat$ given below is the ``correct'' one.

\begin{DefD}\label{def:zinfflat}
We define the category of flat (or torsion-free) $\Zinfty$-modules $\ZinfFlat$ 
as follows. 
Its objects are pairs $A=(A_{\Zinfty},A_\bbR)$, consisting of a 
real vector space $A_\bbR$ and a symmetric convex body 
$A_{\Zinfty}\subset A_\bbR$ (in particular, $A_{\Zinfty}$ is required to be 
absorbent; cf.~\ptref{sp:convexdefs}). When no confusion can arise, we denote 
$A_{\Zinfty}$ by the same letter~$A$ as the whole pair.

The morphisms $f=(f_{\Zinfty}, f_\bbR)$ of $\ZinfFlat$ and their 
composition are defined in the same way as in~\ptref{def:zinflat}.
\end{DefD}

\nxsubpoint
In contrast with \ptref{def:zinflat}, now we do not require $A_\bbR$ to 
be finite-dimen\-sional, and we do not require $A_{\Zinfty}$ to be compact 
or even closed (with respect to the ``algebraic topology'' on $A_\bbR$ 
in the infinite-dimensional case).

\nxsubpoint
Any flat $\Zinfty$-module $A$ defines a {\em seminorm\/} 
$\vnorm_A:A_\bbR\to\bbR_{\geq0}$ 
by the same formula as before: 
${\|x\|}_A=\inf\{\lambda> 0:\lambda^{-1}x\in A_{\Zinfty}\}$. 
However, in general this will be just a seminorm, i.e.\ it has only 
properties 1.\ and 2.\ of~\ptref{p:norms}. This seminorm is still continuous 
with respect to the algebraic topology on~$A_\bbR$. However, it doesn't 
determine $A_{\Zinfty}$ uniquely, unless we know $A_{\Zinfty}$ to be closed: 
all we can say is 
$\{x:{\|x\|}_A<1\}\subset A_{\Zinfty} \subset \{x:{\|x\|}_A\leq 1\}$
(cf.~\cite{EVT}, ch.~II, \S2, prop.~22).

\nxpoint\label{p:projlimzflat}
{\em Arbitrary projective limits exist in $\ZinfFlat$}. Indeed, 
any such limit $L=\projlim{}_{\alpha\in I} F_\alpha$ 
can be computed as follows. 
(Here we consider generalized projective limits, i.e.\ limits in the sense 
of~\cite{MacLane}; so $I$ is actually a small category, and 
$F:\alpha\mapsto F_\alpha$ is a functor $I\to\ZinfFlat$; however, 
for simplicity we write 
generalized projective limits in the same way as the ``classical'' projective 
limits along an ordered set.)
First of all, compute $V:=\projlim{}_{\alpha} F_{\alpha,\bbR}$ in the 
category $\RVect$ 
of $\bbR$-vector spaces, and denote by $\pi_\alpha:V\to F_{\alpha,\bbR}$
the canonical projections. 
Then put $L_{\Zinfty}:=\bigcap_\alpha \pi_\alpha^{-1}(F_{\alpha,\Zinfty})$. 
This is clearly a symmetric convex subset in $V$. It need not be absorbent, 
so we take $L_\bbR:=\bbR\cdot L_{\Zinfty}\subset V$ and put 
$L:=(L_{\Zinfty},L_\bbR)$. The projections $\pi'_\alpha:L\to F_\alpha$ are 
given by the restrictions of $\pi_\alpha$ onto $L_\bbR$; they are 
well-defined (i.e.\ they define morphisms in $\ZinfFlat$) since 
by construction $\pi_\alpha(L_{\Zinfty})\subset F_{\alpha,\Zinfty}$.

We have to check that $L$ with $\{\pi'_\alpha:L\to F_\alpha\}$ have the 
universal property required for projective limits. So suppose
$\{\rho_\alpha:M\to F_\alpha\}$ is any compatible system of morphisms from 
some object~$M$ into our projective system. By definition of~$V$ we have 
a unique $\bbR$-linear map $\phi:M_\bbR\to V$, such that $\pi_\alpha\circ\phi
=\rho_{\alpha,\bbR}$ for all $\alpha$. Since $\rho_{\alpha,\bbR}(M_{\Zinfty})
\subset F_{\alpha,\Zinfty}$, we see that $\phi(M_{\Zinfty})$ is contained 
in each $\pi_\alpha^{-1}(F_{\alpha,\Zinfty})$, hence in their intersection 
$L_{\Zinfty}$. Now $M_\bbR=\bbR\cdot M_{\Zinfty}$, hence 
$\phi(M_\bbR)\subset\bbR\cdot L_{\Zinfty}=L_\bbR\subset V$. This gives 
us a $\Zinfty$-homomorphism $\phi':M\to L$, such that $\pi'_\alpha\circ
\phi'=\rho_\alpha$ for all $\alpha$; its uniqueness is immediate.

This construction can be expressed in the language of norms, at least 
if all $F_{\alpha,\Zinfty}$ are closed, as follows. We introduce on the 
projective limit $V$ of given normed vector spaces the $\sup$-norm by 
${\|x\|}_{\sup}:=\sup_\alpha{\|\pi_\alpha(x)\|}_{F_\alpha}$. Then we 
consider the subspace $L_\bbR:=\{x\in V:{\|x\|}_{\sup}<+\infty\}$ of elements 
with finite $\sup$-norm. Then $L_\bbR$ is a normed vector space with 
respect to the restriction of the $\sup$-norm, and 
$L_{\Zinfty}=\{x:{\|x\|}_{\sup}\leq 1\}$.

\nxsubpoint\label{ss:tensRlex}
Note that, if $I$ is finite, then $L_\bbR=V$ in the above notations. 
In other words, functor $\rho^*:X\mapsto X_\bbR$ 
commutes with finite projective limits, i.e.\ 
it is left exact. To show this take any $x\in V$ and observe 
that its $\sup$-norm is finite (being the supremum of a finite set of 
real numbers), hence we can find some $\lambda>{\|x\|}_{\sup}$. Now 
it follows from definitions that $\pi_\alpha(\lambda^{-1}x)$ belongs 
to $F_{\alpha,\Zinfty}$ for each $\alpha$, hence $\lambda^{-1}x\in L_{\Zinfty}$
and $x\in\bbR\cdot L_{\Zinfty}=L_\bbR$.

\nxsubpoint
In particular, finite products in $\ZinfFlat$ are computed componentwise: 
$A\times B=(A_{\Zinfty}\times B_{\Zinfty}, A_\bbR\times B_\bbR)$; cf.~%
\ptref{p:limszinflat}.

\nxsubpoint\label{sp:descr.mono}
Another consequence is this: {\em a morphism $f=(f_{\Zinfty},f_\bbR)$ is 
a monomorphism (in $\ZinfFlat$) iff $f_\bbR$ is injective iff 
$f_{\Zinfty}$ is injective}. 
Indeed, we have just seen that $X\mapsto X_\bbR$ 
is left exact, hence transforms monomorphisms into monomorphisms. So if $f$ is monic, then $f_\bbR$ is monic in $\RVect$, i.e.\ injective. Since $f_{\Zinfty}$ 
is a restriction of $f_\bbR$, injectivity of $f_\bbR$ implies injectivity of $f_{\Zinfty}$. Finally, this last injectivity clearly implies that $f$ is monic.

\nxsubpoint
This means that {\em subobjects} $M=(M_{\Zinfty},M_\bbR)$ of a fixed 
object $L$ can be characterized by requiring $M_\bbR$ to be an 
$\bbR$-linear subspace of~$L_\bbR$, and $M_{\Zinfty}$ to 
be contained in $L_{\Zinfty}$.

For example, the subobjects (i.e.\ ``ideals'') of~$\Zinfty=([-1,1],\bbR)$ are 
$0=(0,0)$, $\lambda\Zinfty=([-\lambda,\lambda],\bbR)$ and 
$\lambda\gm_\infty=((-\lambda,\lambda),\bbR)$ for $0<\lambda\leq 1$. Of 
these $0$ and $\lambda\Zinfty$ lie in $\ZinfLat$, and $\lambda\gm_\infty$ 
does not, since the open interval $(-\lambda,\lambda)$ is not compact 
in~$\bbR$. Note that among all ideals $\neq\Zinfty$ there is a maximal one, 
namely, $\gm_\infty$, so $\Zinfty$ is something like a local ring.

\nxpoint\label{p:indlimzflat}
{\em Arbitrary inductive limits exist in $\ZinfFlat$}. 
Let us compute such a limit $L=\injlim{}_\alpha F_\alpha$. First of all, 
put $L_\bbR:=\injlim{}_\alpha F_{\alpha,\bbR}$ in $\RVect$. Denote 
by $\sigma_{\alpha,\bbR}: F_{\alpha,\bbR} \to L_\bbR$ the canonical 
maps into the inductive limit, and put 
$L_{\Zinfty}:=\conv\bigl(\bigcup_\alpha\sigma_{\alpha,\bbR}(F_{\alpha,\Zinfty})
\bigr)$. Clearly, $L_{\Zinfty}$ is a symmetric and convex subset 
of~$L_\bbR$. To check that it is absorbent observe that 
$L_{\bbR}$ is generated as an $\bbR$-vector space 
by the union of~$\sigma_{\alpha,\bbR}(F_{\alpha,\bbR})$, and each 
$F_{\alpha,\bbR}$ is generated by $F_{\alpha,\Zinfty}$. Hence the union 
of~$\sigma_{\alpha,\bbR}(F_{\alpha,\Zinfty})$ generates $L_\bbR$; 
{\em a fortiori\/} this is true for $L_{\Zinfty}$.

In this way we get an object $L=(L_{\Zinfty},L_\bbR)$ of $\ZinfFlat$. Clearly, 
$\sigma_{\alpha,\bbR}(F_{\alpha,\Zinfty})\subset L_{\Zinfty}$, so we 
get a compatible system of morphisms $\sigma_\alpha:F_\alpha\to L$. 
Let's check the universal property of this system. If 
$\tau_\alpha:F_\alpha\to M$ is another compatible system of morphisms, 
we get a unique $\bbR$-linear map $\phi_\bbR:L_\bbR\to M_\bbR$, such 
that $\tau_{\alpha,\bbR}=\phi_\bbR\circ\sigma_{\alpha,\bbR}$, since 
$L_\bbR$ is an inductive limit of $F_{\alpha,\bbR}$. Now 
$\tau_{\alpha,\bbR}(F_{\alpha,\Zinfty})$ is contained in $M_{\Zinfty}$, 
hence each $\sigma_{\alpha,\bbR}(F_{\alpha,\Zinfty})$ is contained 
in the convex set $\phi_\bbR^{-1}(M_{\Zinfty})$, therefore, the 
convex hull $L_{\Zinfty}$ of their union 
is contained in this convex set as well. This 
means that $\phi_\bbR$ does indeed induce 
a morphism $\phi:L\to M$ in $\ZinfFlat$, such that $\tau_\alpha=\phi\circ
\sigma_\alpha$, and this $\phi$ is clearly unique, q.e.d.

\nxsubpoint\label{ss:zinffiltind} Note that {\em filtered\/} inductive 
limits of~$\ZinfFlat$ can be computed componentwise in the category of sets, 
i.e.\ in this case we have $L_{\Zinfty}=\injlim F_{\alpha,\Zinfty}$ and 
$L_\bbR=\injlim F_{\alpha,\bbR}$ in the category of sets. This can be 
either checked directly (define $L:=(L_{\Zinfty},L_\bbR)$ and check that it 
lies in $\ZinfFlat$), or deduced from the description given 
in~\ptref{p:indlimzflat}.

\nxsubpoint\label{ss:infdirsum} (Infinite direct sums and norms.) 
In particular, 
arbitrary (not necessarily finite) direct sums (i.e.\ coproducts) 
$L=\bigoplus_\alpha F_\alpha$ exist in $\ZinfFlat$. We see that 
$L_\bbR=\bigoplus_\alpha F_{\alpha,\bbR}$, and 
$L_{\Zinfty}=\conv\bigl(\bigcup_\alpha F_{\alpha,\Zinfty}\bigr)$, 
where the union and the convex hull are computed inside~$L_\bbR$. 
This convex hull can be clearly described as the set of ``octahedral'' 
combinations $\sum_\alpha \lambda_\alpha x_\alpha$, where all $x_\alpha\in
F_{\alpha,\Zinfty}$, almost all $\lambda_\alpha\in\bbR$ are equal to zero and 
$\sum_\alpha|\lambda_\alpha|\leq 1$. (Actually, baricentric combinations 
would suffice here.) The corresponding norm $\vnorm_L$ is the ``$L_1$-norm'' 
$\|\sum_\alpha x_\alpha\|\strut_L=\sum_\alpha {\|x_\alpha\|}_{F_\alpha}$. 
One checks that if all $F_\alpha$ are closed (i.e.\ each $F_{\alpha,\Zinfty}$ 
is closed in $F_{\alpha,\bbR}$ with respect to the algebraic topology) 
and do not contain any non-trivial real vector spaces (i.e.\ each 
$F_{\alpha,\Zinfty}$ is given by a {\em norm\/} on $F_{\alpha,\bbR}$, and not 
just a seminorm), 
then the same is true for $L$, hence $L_{\Zinfty}=\{x:{\|x\|}_L\leq1\}$.

\nxsubpoint\label{sp:rvectinzflat} (Right exactness of $X\mapsto X_\bbR$.) 
Our construction 
of inductive limits shows that the functor $\rho^*:X\mapsto X_\bbR$ commutes 
with arbitrary inductive limits, and, in particular, it is right exact, 
hence exact (cf.~\ptref{ss:tensRlex}).
This leads us to believe that it might have a right adjoint 
$\rho_*:\RVect\to\ZinfFlat$. This adjoint indeed exists and can be 
defined by $\rho_*(V):=(V,V)$ for any $\bbR$-vector space~$V$, i.e.\ 
any $\bbR$-vector space can be considered as a flat $\Zinfty$-module. 
Note that $\rho_*$ is fully faithful, so $\RVect$ can be identified with a 
full subcategory of $\ZinfFlat$. It is also interesting to remark that 
in the language of (semi)norms $\rho^*V$ corresponds to $V$ together with 
the seminorm which is identically zero; i.e.\ we equip $V$ with its 
coarsest (non-separated) topology. 

\nxsubpoint\label{sp:zinfstrepi} (Strict epimorphisms and quotients.)
Another immediate consequence is this: {\em a morphism $f=(f_{\Zinfty},f_\bbR):
M\to P$ is a strict epimorphism iff both $f_{\Zinfty}$ and $f_\bbR$ 
are subjective iff $f_\bbR:M_\bbR\to P_\bbR$ is surjective and 
$P_{\Zinfty}=f_\bbR(M_{\Zinfty})$}. By definition, $f$ is a strict 
epimorphism iff it is a cokernel of its kernel pair 
$M\times_P M\rightrightarrows M$. The 
construction of inductive limits given above implies that if $f:M\to P$ 
is cokernel of a pair of morphisms, then $P_\bbR$ is the corresponding 
quotient in $\RVect$, and $P_{\Zinfty}=\conv\bigl(f(M_{\Zinfty})\bigr)=
f(M_{\Zinfty})$. Conversely, if this is the case, then clearly 
$f$ is the cokernel of its kernel pair.

Thus we see that strict quotients $P$ of a fixed object $M$ are in 
one-to-one correspondence with quotients $\pi_\bbR:M_\bbR\to P_\bbR$ of the 
corresponding $\bbR$-vector space $M_\bbR$. For each of these strict 
quotients we have $P_{\Zinfty}=\pi_\bbR(M_{\Zinfty})$.

In the language of (semi)norms strict quotients correspond to 
quotient norms: ${\|y\|}_P=\inf_{x\in\pi_\bbR^{-1}(y)}{\|x\|}_M$.

\nxsubpoint\label{sp:zinfstrmono} (Strict monomorphisms and subobjects.)
Dually, $f:N\to M$ is a strict monomorphism, if it is the kernel of some 
pair of morphisms, and then it is the kernel of its cokernel pair. 
The description of projective limits given in~\ptref{p:projlimzflat} 
shows immediately that {\em $f:N\to M$ is a strict monomorphism iff 
$f_\bbR$ is injective and $N_{\Zinfty}=f_\bbR^{-1}(M_{\Zinfty})$}. Hence 
strict subobjects of $M$ are given by $N=(N_{\Zinfty},N_\bbR)$, where 
$N_\bbR\subset M_\bbR$ is an $\bbR$-linear subspace, and 
$N_{\Zinfty}=M_{\Zinfty}\cap N_\bbR$. We see that strict subobjects of 
$M$ are in one-to-one correspondence with $\bbR$-subspaces of $M_\bbR$. 
In the language of (semi)norms this means that we consider 
the restriction of $\vnorm_M$ to $N_\bbR\subset M_\bbR$.

\nxsubpoint\label{sp:zinf.coim} (Coimages.)
We see that any morphism $f:M\to N$ decomposes uniquely into 
a strict epimorphism $\pi:M\to Q$ (namely, the cokernel of the 
kernel pair $M\times_N M\rightrightarrows M$ of~$f$), followed by 
a monomorphism $i:Q\to N$. Uniqueness of such a decomposition is a general 
fact. To show existence one has to check that the morphism 
$i:Q\to N$ induced by~$f$ is a monomorphism. We know that injectivity of 
$i_\bbR$ suffices for this; now we apply exact functor~$\rho^*$ and 
see that $i_\bbR:\Coim f_\bbR\to N_\bbR$ is indeed a monomorphism in $\RVect$. 
We will say that $Q$ is the {\em image} of $M$ in $N$ under $f$ (even 
if it would be more accurate to say that it is the {\em coimage} of~$f$).

\nxsubpoint (Inductive limits in $\ZinfMod$ are different.)
We are going to embed $\ZinfFlat$ into a larger category $\ZinfMod$ later. 
This inclusion functor will have a left adjoint; hence it commutes 
with arbitrary projective limits, i.e.\ the projective limits 
computed in~\ptref{p:projlimzflat} will remain such when considered 
in $\ZinfMod$. This is in general {\em not\/} true for arbitrary 
inductive limits, even if we consider only finite inductive limits. 
For example, if we compute the cokernel of the inclusion and of the zero 
morphisms $\gm_\infty\rightrightarrows\Zinfty$, we obtain $0$ in 
$\ZinfFlat$, but in $\ZinfMod$ this cokernel ``$\Zinfty/\gm_\infty$'' 
will be a non-trivial torsion $\Zinfty$-module (cf.~\ptref{sp:def.finf}).

However, {\em filtered\/} inductive limits will be the same in 
$\ZinfMod$ and in $\ZinfFlat$, because in both categories they 
can be described as in~\ptref{ss:zinffiltind}.

\nxsubpoint (Zero object.)
Our category $\ZinfFlat$ has a {\em zero object\/} $0:=(0,0)$, i.e.\ 
an object that is both initial and final at the same time. Hence 
we have a pointed element in each set of morphisms $\Hom_{\Zinfty}(M,N)$, 
namely, the {\em zero morphism} $0_{MN}:M\to0\to N$. Note that 
$0$ lies actually in $\ZinfLat$, so the same is true in this smaller 
category.

\nxpoint\label{p:zinfflat.polylin} (Multilinear algebra.)
Bilinear and polylinear $\Zinfty$-maps between objects of $\ZinfFlat$ 
are defined exactly in the same way as it has been done 
in~\ptref{def:zinflatbilin} and~\ptref{sp:defpolylin}, once we replace 
``$\Zinfty$-lattices'' with ``flat $\Zinfty$-modules'' in these definitions. 
Moreover, double and multiple tensor products can be still constructed 
in the same way as in~\ptref{p:extensprod}. For example, 
$M\otimes N=M\otimes_{\Zinfty}N$ is defined to be $L=(L_{\Zinfty},L_\bbR)$, 
where $L_\bbR=M_\bbR\otimes_\bbR N_\bbR$, and $L_{\Zinfty}$ is the convex hull 
of the set of all products 
$\{x\otimes y: x\in M_{\Zinfty}, y\in N_{\Zinfty}\}$. Clearly, $L_{\Zinfty}$ 
generates $L_\bbR$ as an $\bbR$-vector space, so this is indeed an object 
of~$\ZinfFlat$.

\zerosubpoint\nxsubpoint (Relation to Grothendieck's tensor product
of seminorms.)
Notice that the tensor product $L=M\otimes N$ just discussed
(cf.\ also \ptref{p:extensprod}) admits a direct description in terms
of seminorms on $L_\bbR=M_\bbR\otimes N_\bbR$, $M_\bbR$ and $N_\bbR$.
Namely, if $\vnorm$, $\vnorm_1$ and $\vnorm_2$ denote the seminorms
defined by $L$, $M$ and $N$, respectively, then
\begin{equation}
\|z\|=\inf_{z=\sum_{i=1}^nx_i\otimes y_i}
\sum_{i=1}^n{\|x_i\|}_1\cdot{\|y_i\|}_2,\quad\text{for any $z\in L_\bbR$.}
\end{equation}
This is nothing else than the projective tensor product of seminorms defined 
by Grothendieck in~\cite{Gr0}. However, this doesn't describe our
tensor product completely since $L_{\Zinfty}$ isn't completely
determined by~$\vnorm$.

\nxsubpoint\label{sp:acu.zinfflat} (ACU $\otimes$-structure on $\ZinfFlat$.) 
In this way we get an ACU $\otimes$-structure on $\ZinfFlat$, compatible 
with one introduced on $\ZinfLat$ before. In particular, $\Zinfty$ 
is the unit object with respect to this tensor product, so we still 
have a ``global sections functor'' 
$\Gamma=\Hom_{\Zinfty}(\Zinfty,-):\ZinfFlat\to\catSets$, that maps 
a flat $\Zinfty$-module $M=(M_{\Zinfty},M_\bbR)$ into its 
``underlying set'' $M_{\Zinfty}\cong\Hom_{\Zinfty}(\Zinfty,M)$. This 
terminology is again motivated by the $p$-adic case, where 
$\Hom_{\bbZ_p}(\bbZ_p,-)$ maps a flat $\bbZ_p$-module~$M$ into its 
underlying set.

\nxsubpoint
{\em Inner Homs\/ $\iHom_{\Zinfty}(N,P)$ exist in $\ZinfFlat$.}
We construct such an inner Hom $H=\iHom(N,P)$ as follows. Put 
$H':=\Hom_\bbR(N_\bbR,P_\bbR)$, then 
$H_{\Zinfty}:=\{f\in H' : f(N_{\Zinfty})\subset P_{\Zinfty}\}$, 
$H_\bbR := \bbR\cdot H_{\Zinfty}$ and $H:=(H_{\Zinfty},H_\bbR)$. 
Clearly, $H_{\Zinfty}$ is a symmetric convex subset of~$H'$, hence 
it is a symmetric convex body in the real vector subspace $H_\bbR\subset H'$, 
so $H$ is indeed an object of~$\ZinfFlat$. Required universal property 
$\Hom(M,H)\cong\Hom(M\otimes N, P)\cong\Bilin(M,N;P)$ can be deduced now from 
the similar property of $H'$ in $\RVect$ in the same 
manner as it has been done for projective limits in~\ptref{p:projlimzflat}.

\nxsubpoint\label{sp:r.ex.tens}
Existence of inner Homs implies that {\em tensor products in $\ZinfFlat$ 
commute with arbitrary inductive limits in each variable.} In particular, 
$(M'\oplus M'')\otimes N\cong(M'\otimes N)\oplus(M''\otimes N)$.

\nxsubpoint
Note that $\Gamma(\iHom_{\Zinfty}(N,P))\cong\Hom_{\Zinfty}(N,P)$, as 
one expects for general reasons. Another interesting observation is this: 
For any $M=(M_{\Zinfty},M_\bbR)$ we have $M\otimes_{\Zinfty}\bbR=
(M_\bbR, M_\bbR)=M_\bbR$, if we identify $\RVect$ with a full subcategory 
of~$\ZinfFlat$ as in~\ptref{sp:rvectinzflat}. So we can write 
$M\otimes_{\Zinfty}\bbR$ or $M_{(\bbR)}$ instead of $\rho^*M$ or $M_\bbR$. 
This leads us to believe that $\rho_*:\RVect\to\ZinfFlat$ must have a 
right adjoint $\rho^!:\ZinfFlat\to\RVect$, given by 
$\rho^!(M):=\iHom_{\Zinfty}(\bbR,M)$, and this is indeed the case. 
This $\bbR$-vector space $\rho^!M$ can be described more explicitly 
as the largest $\bbR$-vector subspace of $M_\bbR$ contained in~$M_{\Zinfty}$.

\nxsubpoint (Duals.)
Of course, we define the {\em dual $M^*$} of some flat $\Zinfty$-module~$M$ 
by $M^*:=\iHom_{\Zinfty}(M,\Zinfty)$. We have a lot of canonical maps 
like $M^*\otimes N\to\iHom(M,N)$ and $M\to M^{**}$, coming from the 
general theory of $\otimes$-categories. One can check that 
$M$ is {\em reflexive}, i.e.\ $M\to M^{**}$ is an isomorphism iff 
$M_\bbR$ is a Banach space with respect to $\vnorm_M$, and 
$M_{\Zinfty}$ is closed, i.e.\ $M_{\Zinfty}=\{x:{\|x\|}_M\leq 1\}$.
For any $M$ its double dual $M^{**}$ is something like the {\em completion} 
of~$M$: $(M^{**})_\bbR$ is the completion of $M_\bbR$ with respect to 
$\vnorm_M$, and $(M^{**})_{\Zinfty}$ is the closure of the image of 
$M_{\Zinfty}$ in this completed space.

One should be careful with these definitions. For example, 
$\Zinfty^{**}=\Zinfty$, but $\bbR^{**}=0$.

\nxsubpoint (Hilbert spaces.)
We have seen that Banach spaces have an interior characterization in 
$\ZinfFlat$. Hilbert spaces have such a characterization as well: they 
correspond to flat $\Zinfty$-modules $H$ equipped by a {\em perfect\/}
symmetric pairing $\phi:H\times H\to\Zinfty$, i.e.\ $\phi$ is 
a symmetric $\Zinfty$-bilinear form, such that the induced map 
$\tilde\phi: H\to H^*$ is an isomorphism. Actually, to obtain Hilbert 
spaces we must also require $\phi$ to be positive definite, otherwise 
$-\phi$ would also do.

\nxsubpoint\label{sp:freezinfmod} (Free $\Zinfty$-modules.)
For any flat $\Zinfty$-module $M$ and any set $S$ we can consider 
the {\em product $M^S$} and the {\em direct sum $M^{(S)}$} 
of the constant family with value $M$ and index set~$S$. If $S$ 
is a standard finite set $\stn=\{1,2,\ldots,n\}$, $n\geq0$, we write 
$M^{(n)}$ and $M^n$ instead of $M^{(\stn)}$ and $M^\stn$. Note that 
we have a canonical {\em monomorphism\/} $M^{(S)}\to M^S$, which 
usually is {\em not\/} an isomorphism, even for a finite~$S$.
By definition, we have
\begin{align}
\Hom_{\Zinfty}(N,M^S)\cong\Hom_{\Zinfty}(N,M)^S\cong
\Hom_{\catSets}\bigl(S,\Hom_{\Zinfty}(N,M)\bigr)\\
\Hom_{\Zinfty}(M^{(S)},N)\cong\Hom_{\Zinfty}(M,N)^S\cong
\Hom_{\catSets}\bigl(S,\Hom_{\Zinfty}(M,N)\bigr)
\end{align}
From these formulas we deduce a lot of other canonical isomorphisms, 
for example, $P\otimes_{\Zinfty} M^{(S)}\cong (P\otimes_{\Zinfty} M)^{(S)}$, 
$\iHom_{\Zinfty}(M^{(S)},N)\cong\iHom_{\Zinfty}(M,N)^S$, 
$(M^{(S)})^*\cong (M^*)^S$ and so on. We can specialize these constructions 
to the case $M=\Zinfty$. We obtain $M^{(S)}\cong\Zinfty^{(S)}\otimes M$, 
$\Zinfty^{(S)}\otimes\Zinfty^{(T)}\cong\Zinfty^{(S\times T)}$, 
$M^S\cong\iHom(\Zinfty^{(S)},M)$, $(\Zinfty^{(S)})^*\cong\Zinfty^S$, and 
in particular
\begin{equation}
\Hom_{\Zinfty}(\Zinfty^{(S)},M)\cong\Gamma(M)^S\cong
\Hom_{\catSets}(S,\Gamma(M))
\end{equation}
This last equality means that the ``forgetful functor'' 
$\Gamma:\ZinfFlat\to\catSets$ has a left adjoint 
$L_{\Zinfty}:\catSets\to\ZinfFlat$, given by $S\mapsto \Zinfty^{(S)}$. 
That's the reason why we say that any (flat) $\Zinfty$-module~$M$, 
isomorphic to some $\Zinfty^{(S)}$, is {\em a free $\Zinfty$-module
(of rank $\card S$)}. The existence of this left adjoint to $\Gamma$ 
will be crucial for our definition of $\ZinfMod$.

\nxsubpoint\label{sp:explfreezinfmod}
Since free $\Zinfty$-modules $\Zinfty^{(S)}$ will be very important to us 
later, let's describe them more explicitly. According to~\ptref{ss:infdirsum},
$\Zinfty^{(S)}=(\Sigmainf(S),\bbR^{(S)})$, where $\Sigmainf(S)$ is the 
convex hull of all $\pm\{s\}$, $s\in S$. Here we denote by $\{s\}$ 
the basis element of $\bbR^{(S)}$ corresponding to some $s\in S$. 
Hence $\bbR^{(S)}$ is the set of formal linear combinations 
$\sum_s\lambda_s\{s\}$ of elements of~$S$ (i.e.\ all but finitely 
many $\lambda_s\in\bbR$ are equal to zero), and $\Sigmainf(S)$ can 
be described as the set of all formal {\em octahedral combinations\/} 
$\sum_s\lambda_s\{s\}$ of elements of~$S$, i.e.\ we require 
in addition that $\sum_s|\lambda_s|\leq 1$. This means that 
$\Sigmainf(S)$ is closed in $\bbR^{(S)}$, and the corresponding 
norm is the $L_1$-norm: $\bigl\|\sum_s\lambda_s\{s\}\bigr\|=\sum_s|\lambda_s|$.

If $S$ is infinite, $\bbR^{(S)}$ is not complete with respect to this 
norm, hence $\Zinfty^{(S)}$ is not reflexive. On the other hand, 
if $S$ is finite, e.g.\ $S=\stn=\{1,2,\ldots,n\}$, then 
$\Zinfty^{(n)}=(\Sigmainf(n),\bbR^n)$ lies in $\ZinfLat$, hence is reflexive. 
In this case $\Sigmainf(n)$ is the convex hull of $\pm e_i$, where 
$\{e_i\}_{i=1}^n$ is the standard base of $\bbR^n$. In other words, 
$\Sigmainf(n)$ is the {\em standard $n$-dimensional octahedron} 
with vertices $\pm e_i$.

\nxsubpoint (Finitely generated modules.)
We say that a flat $\Zinfty$-module $M$ is {\em finitely generated}, or 
is {\em of finite type}, if it is a strict quotient of some 
free $\Zinfty$-module $\Zinfty^{(n)}$ of finite rank.
The description of strict quotients given in~\ptref{sp:zinfstrepi} 
combined with the description of~$\Zinfty^{(n)}$ given 
in~\ptref{sp:explfreezinfmod}, yields the following characterization 
of finitely generated flat $\Zinfty$-modules: $M$ is finitely generated 
iff $M_\bbR$ is finite-dimensional, and $M_{\Zinfty}$ is the convex hull 
of a finite number of vectors $\pm x_i$ in $M_\bbR$. In other words, 
$M_{\Zinfty}$ has to be a symmetric convex polyhedron in a 
finite-dimensional real vector space~$M_\bbR$. This implies that all 
flat $\Zinfty$-modules of finite type are $\Zinfty$-lattices, but not 
conversely. In particular, $\barZinfty$ is a $\Zinfty$-lattice, but 
it is not finitely generated.
However, any $\Zinfty$-lattice, as well as any 
flat $\Zinfty$-module, is a filtered inductive limit of 
its finitely generated submodules.

\nxpointtoc{Category of torsion-free algebras and modules}\label{p:tf.alg.mod}
Our ACU $\otimes$-structure on category $\ZinfFlat$ allows us to define flat 
$\Zinfty$-algebras and ($\Zinfty$-flat) modules over these algebras 
as algebras and modules for this $\otimes$-structure.

\nxsubpoint
For example, an (associative) algebra $A$ in $\ZinfFlat$ 
is a flat $\Zinfty$-module 
$A$, together with multiplication $\mu=\mu_A:A\otimes_{\Zinfty} A\to A$ 
and unit $\eta=\eta_A:\Zinfty\to A$, subject to associativity and 
unit relations: $\mu\circ(\mu\otimes 1_A)=\mu\circ(1_A\otimes\mu)$, 
$\mu\circ(\eta\otimes 1_A)=1_A=\mu\circ(1_A\otimes\eta)$.

We can give an alternative description of this structure. Observe for this
that $\mu_A$ corresponds to a $\Zinfty$-bilinear map 
$\mu'=\mu'_A:A\times A\to A$,
and $\eta_A\in\Hom_{\Zinfty}(\Zinfty,A)=\Gamma(A)\cong A_{\Zinfty}$ corresponds
to some element $\eta'=\eta'_A\in A_{\Zinfty}\subset A_\bbR$. Then 
$\mu'_\bbR$ defines on $A_\bbR$ the structure of an $\bbR$-algebra with 
unity $\eta'$, and $A_{\Zinfty}\subset A_\bbR$ is a symmetric 
convex absorbent multiplicative {\em submonoid\/} of $A_\bbR$, containing 
the unity $\eta'$. Clearly, this is an equivalent description of 
flat $\Zinfty$-algebras; for $\Zinfty$-lattices it coincides with one 
we had in~\ptref{th:realsubalg}.

\nxsubpoint
For example, $\Zinfty$, $\barZinfty$, $\bbR$, $\bbC$ and $\Zinfty^n$ 
are $\Zinfty$-algebras; 
this is not the case for $\Zinfty^{(n)}\subset\Zinfty^n$, since it does not 
contain the unity. Note that the tensor product of two flat $\Zinfty$-algebras 
is an $\Zinfty$-algebra again. In particular, for any flat $\Zinfty$-algebra~%
$A$ we obtain that $A\otimes_{\Zinfty}\bbR=A_{(\bbR)}\cong A_\bbR$ is a 
$\Zinfty$-algebra, and even an $\bbR$-algebra; we have already
seen this before.

\nxsubpoint\label{sp:inn.end.alg}
Another example of a (non-commutative) flat $\Zinfty$-algebra is given 
by the ``inner Ends'' $A:=\iEnd_{\Zinfty}(M):=\iHom_{\Zinfty}(M,M)$, with 
the multiplication given by composition of endomorphisms; more formally, 
we consider canonical morphism 
$\iHom(M',M'')\otimes\iHom(M,M')\to\iHom(M,M'')$, existing in any 
ACU $\otimes$-category with inner Homs, for $M''=M'=M$.

\nxsubpoint
A {\em (left) $A$-module structure} on a flat $\Zinfty$-module $M$ 
is by definition 
a morphism $\alpha=\alpha_M:A\otimes M\to M$, such that 
$\alpha\circ(1_A\otimes\alpha)=\alpha\circ(\mu\otimes 1_M)$ and 
$\alpha\circ(\eta\otimes 1_M)=1_M$. Of course, $\alpha$ corresponds to 
a $\Zinfty$-bilinear map $\alpha':A\times M\to M$. It induces an 
$A_\bbR$-module structure $\alpha'_\bbR$ 
on the real vector space $M_\bbR$, such that 
$A_{\Zinfty}\cdot M_{\Zinfty}=M_{\Zinfty}$. In other words, 
$M_{\Zinfty}$ is a symmetric convex body in the $A_\bbR$-module $M_\bbR$, 
stable under action of $A_{\Zinfty}\subset A_\bbR$. Note that 
in this way we obtain an action of monoid $A_{\Zinfty}$ on the 
set $M_{\Zinfty}$.

\nxsubpoint
Both these descriptions of $A$-modules are useful. For example, 
the first of them, combined with~\ptref{sp:r.ex.tens}, shows us that 
the direct sum of two $A$-modules is again an $A$-module, while 
the second one shows that the product of two $A$-modules is an $A$-module.

\nxsubpoint
Another possible description of an $A$-module is this. Note that 
$\alpha_M:A\otimes M\to M$ corresponds by adjointness to some 
$\gamma:A\to\iEnd(M)$, and $\alpha_M$ is an $A$-module structure on~$M$ 
iff $\gamma$ is a $\Zinfty$-algebra homomorphism with respect to 
the $\Zinfty$-structure on~$\iEnd(M)$ considered in~\ptref{sp:inn.end.alg}.

\nxsubpoint\label{def:catAmod}
Of course, (left) $A$-module homomorphisms, or $A$-linear maps 
$f:M=(M,\alpha_M)\to N=(N,\alpha_N)$ are defined to be 
$\Zinfty$-linear maps $f:M\to N$, such that $\alpha_N\circ(1_A\otimes f)=
f\circ\alpha_M$. In our situation this means that $f_\bbR:M_\bbR\to N_\bbR$ 
is an $A_\bbR$-linear map, and $f_\bbR(M_{\Zinfty})\subset N_{\Zinfty}$.

Now we can define the category $\catFlMod{A}$ of {\em (left) $A$-modules 
which are flat over~$\Zinfty$}. Any $\Zinfty$-algebra homomorphism 
$f:A\to B$ gives us a ``scalar restriction'' functor 
$f_*:\catFlMod{B}\to\catFlMod{A}$, which turns out to have 
both left and right adjoints $f^*,f^!:\catFlMod{A}\to\catFlMod{B}$. 
If $A$ is commutative, we obtain a natural notion of $A$-bilinear maps, 
hence of tensor products $\otimes_A$ and inner Homs $\iHom_A$.
We do not develop these notions here, because we would like to remove the 
$\Zinfty$-flatness restriction first. We will return to these questions 
later, in the context of generalized rings.

\nxsubpoint
Given any flat $\Zinfty$-module $M$, we can construct its 
{\em tensor algebra $T(M)$} in the usual way. Namely, we consider the tensor 
powers $T^n(M)=T_{\Zinfty}^n(M):=M\otimes M\otimes\cdots\otimes M$ 
(tensor product of $n$ copies of~$M$; $T^1(M)=M$ and $T^0(M)=\Zinfty$), 
and put $T(M):=\bigoplus_{n\geq0} T^n(M)$. Canonical maps 
$T^n(M)\otimes T^m(M)\to T^{n+m}(M)$ induce the multiplication on 
$T(M)$; thus $T(M)$ becomes a graded associative flat $\Zinfty$-algebra. 
We have a canonical embedding $M\to T(M)$ that maps $M$ into $T^1(M)\cong M$. 
The tensor algebra together with this embedding have the usual 
universal property with respect to $\Zinfty$-linear maps 
$\phi:M\to A$ from $M$ into associative flat $\Zinfty$-algebras $A$: 
any such map lifts to a unique homomorphism of $\Zinfty$-algebras 
$\phi^\sharp:T(M)\to A$.

Note that tensor products of free modules are free, so, if 
$M$ is a free $\Zinfty$-module, then all $T^n(M)$ and $T(M)$ are also free.

\nxsubpoint
The {\em symmetric algebra $S(M)$} can be defined by its universal property 
with respect to $\Zinfty$-linear maps of $M$ into {\em commutative} 
algebras: any such map $\phi:M\to A$ should lift to a unique 
$\Zinfty$-algebra homomorphism $\phi^\sharp:S(M)\to A$. This symmetric 
algebra does exist, and it is a commutative graded $\Zinfty$-algebra: 
$S(M)=\bigoplus_{n\geq0}S^n(M)$. We have a canonical $\Zinfty$-algebra 
map $\pi:T(M)\to S(M)$. It is a strict epimorphism in $\ZinfFlat$, 
and it respects the grading; the induced maps $\pi_n:T^n(M)\to S^n(M)$ 
are also strict epimorphisms. The {\em symmetric powers $S^n(M)$} have 
the usual universal property with respect to symmetric polylinear 
maps from $M^n$ into a variable $\Zinfty$-module $N$. They can be constructed 
as the image of $T_{\Zinfty}^n(M)\subset T_\bbR^n(M_{\bbR})$ under 
the canonical surjection $\pi_{n,\bbR}:T_\bbR^n(M_\bbR)\to S_\bbR^n(M_\bbR)$, 
i.e.\ $S^n(M)=(S^n(M)_{\Zinfty}, S^n(M)_\bbR)$, where 
$S^n(M)_\bbR=S_\bbR^n(M_\bbR)$, and $S^n(M)_{\Zinfty}=\pi_{n,\bbR}
(T^n(M)_{\Zinfty})$ (recall the description of strict quotients given 
in~\ptref{sp:zinfstrepi}).

\nxsubpoint
The {\em exterior algebra $\Lambda(M)$} is defined similarly, but we 
consider $\Zinfty$-linear maps $\phi:M\to A$ into associative algebras~$A$, 
subject to additional condition $\phi(x)^2=0$ for any $x$ in $M_{\Zinfty}$ 
(or in~$M_\bbR$). Again, $\Lambda(M)=\bigoplus_{n\geq0}\Lambda^n(M)$ 
is graded supercommutative, and it is a strict quotient of $T(M)$, 
and any individual graded component $\Lambda^n(M)$ is strict quotient of 
$T^n(M)$; together with property $(\Lambda^n(M))_{\bbR}=
\Lambda_\bbR^n(M_\bbR)$ 
these conditions determine $\Lambda^n(M)$ and $\Lambda(M)$ uniquely. 
Finally, 
all the $\Lambda^n(M)$ have the usual universal property with respect to 
alternating polylinear maps from $M^n$ to a variable $\Zinfty$-module $N$.

\nxsubpoint
{\em Monoid algebras\/} give us another example of $\Zinfty$-algebras. 
Given any monoid $M$ with unity $e$, we define the {\em monoid algebra}
$\Zinfty[M]$ to be the free $\Zinfty$-module $\Zinfty^{(M)}$, 
with unity $\{e\}$ and multiplication induced by that of $M$ by means of 
the canonical isomorphism $\Zinfty^{(M)}\otimes\Zinfty^{(M)}\cong
\Zinfty^{(M\times M)}$. In other words, $\Zinfty[M]$ consists of 
formal octahedral combinations of its basis elements $\{m\}$, and 
the multiplication is defined by the usual requirement $\{m\}\cdot\{m'\}=
\{mm'\}$. According to~\ptref{sp:explfreezinfmod}, we have 
$\Zinfty[M]_\bbR=\bbR[M]$, and the corresponding symmetric convex body in 
$\bbR[M]$ is given by the $L_1$-norm 
$\bigl\|\sum_m \lambda_m\{m\}\bigr\|=\sum_m|\lambda_m|$.

Note that $\Zinfty[M]$ has the usual universal property of monoid algebras, 
namely, $\Zinfty$-algebra homomorphisms $f:\Zinfty[M]\to A$ are 
in one-to-one correspondence with monoid homomorphisms $f^\flat:M\to A^\times$,
where $A^\times$ denotes $A_{\Zinfty}$, considered as a  
monoid under multiplication.

\nxsubpoint
Now we can combine together several of the above constructions and 
construct the {\em polynomial algebras}. Given any set $X$, 
we define the polynomial algebra $\Zinfty[X]$ in variables from $X$ 
to be a commutative algebra, such that 
$\Hom_\catSets(X,\Gamma(A))\cong\Hom_{\Zinfty-alg}(\Zinfty[X],A)$ 
for any {\em commutative} (flat) $\Zinfty$-algebra $A$. If $X$ is 
finite, say $X=\stn=\{1,2,\ldots,n\}$, we denote the corresponding 
elements of $\Zinfty[X]$ by $T_i$, and write $\Zinfty[T_1,T_2,\ldots, T_n]$ 
instead of $\Zinfty[\stn]$.

We have two ways of constructing polynomial algebras:

\noindent a) We can construct $\Zinfty[X]$ as the symmetric algebra of the 
free module $\Zinfty^{(X)}$. Indeed, for any commutative $\Zinfty$-algebra 
we have $\Hom_{\Zinfty-alg}(S(\Zinfty^{(X)}),A)\cong
\Hom_{\Zinfty}(\Zinfty^{(X)},A)\cong\Hom_\catSets(X,\Gamma(A))=\Gamma(A)^X$. 
In particular, we obtain a grading on $\Zinfty[X]$.

\noindent b) We can construct first the free commutative monoid 
$N(X)$ generated by $X$, i.e.\ the set of monomials in $X$, and then define 
$\Zinfty[X]$ to be the monoid algebra $\Zinfty[N(X)]$. Indeed, 
for any commutative algebra~$A$ we have
$\Hom_{\Zinfty-alg}(\Zinfty[N(X)],A)\cong\Hom_{\catMon}(N(X),A^\times)
\cong\Hom_\catSets(X,\Gamma(A))=\Gamma(A)^X$.
In this way we see that $\Zinfty[X]$ is free as a $\Zinfty$-module, 
i.e.\ the symmetric algebra of a free module is free. This is 
actually true for any individual symmetric power $S^n(\Zinfty^{(X)})$, 
since it admits the set $N_n(X)\subset N(X)$ of all monomials of degree~$n$ 
as a basis.

The corresponding norm on $\Zinfty[X]_\bbR=\bbR[X]$ is again the 
$L_1$-norm with respect to the basis given by the monomials $N(X)$.

\nxpointtoc{Arakelov affine line}
We would like to consider an interesting example -- namely, 
the {\em Arakelov affine line} $\bbA_{\Zinfty}^1$. Of course, 
$\bbA_{\Zinfty}^1=\Spec\Zinfty[T]$. However, we haven't yet developed 
the theory of ``generalized rings'' like $\Zinfty$ and $\Zinfty[T]$ 
and their spectra, so we will try to study this affine line 
in terms of $\Zinfty$-sections and norms, as described 
in~\ptref{p:cls.arak.vbun}.

\nxsubpoint
First of all, let's take some $\bbR$-rational point $P_\lambda$ in the 
generic fiber $\bbA_\bbR^1=\Spec\bbR[T]=\Spec \Zinfty[T]_{(\bbR)}$. 
It is given by some real number $\lambda$, and the corresponding 
section $\sigma'_\lambda:\Spec\bbR\to\Spec\bbR[T]$ is given by 
the evaluation at~$\lambda$ map 
$\pi_\lambda:\bbR[T]\to\bbR$, $F(T)\mapsto F(\lambda)$. Now we want 
to compute the scheme-theoretic closure of $\{P_\lambda\}$ in 
$\bbA_{\Zinfty}^1$. The $p$-adic case tells us that we have to consider 
$\Spec C_\lambda$ for this, where $C_\lambda$ is the image of $\Zinfty[T]$ in 
$\bbR$ under $\pi_\lambda$. Since $\Zinfty[T]$ is the convex hull of 
$\pm T^k$ in $\bbR[T]$, we see that $C_\lambda=(C_\lambda,\bbR)$ is given 
by $C_\lambda=\conv\{\pm\lambda^k:k\geq0\}$. If $|\lambda|>1$, we get 
$C_\lambda=\bbR$, i.e.\ $P_\lambda$ is closed in $\bbA_{\Zinfty}^1$. 
If $|\lambda|\leq1$, then $C_\lambda=\Zinfty$, so $P_\lambda$ lifts to 
a $\Zinfty$-section $\sigma_\lambda$ of $\bbA_{\Zinfty}^1$, 
in perfect accordance with the $p$-adic case.

\nxsubpoint\label{sp:a1zinfty.infdata}
Now we would like to extract some infinitesimal data at $P_\lambda$ 
when $P_\lambda$ lifts to a $\Zinfty$-section, i.e.\ when $|\lambda|\leq1$. 
More precisely, we want to obtain some cometric at these points.

In the usual situation we would have 
$\sigma_\lambda^*(\Omega^1)\cong\gp_\lambda/\gp_\lambda^2$, 
where $\gp_\lambda$ is the ideal defining 
the image of this section, i.e.\ the kernel of 
$\Zinfty[T]\to C_\lambda=\Zinfty$. Clearly, we must have 
$(\gp_\lambda)_{(\bbR)}=\gp'_\lambda=
\Ker\pi_\lambda=(T-\lambda)\cdot\bbR[T]$, so let's 
take $\gp_\lambda:=\Zinfty[T]\cap\gp'_\lambda$, and let's 
assume that $\gp_\lambda\to E=\gp_\lambda/\gp_\lambda^2$ is a strict 
epimorphism, and that $E_\bbR=\gp'_\lambda/{\gp'_\lambda}^2\cong\bbR$, 
where the isomorphism between $E_\bbR$ and $\bbR$ maps the class of
$T-\lambda$ into $1$. Note that $T-\lambda$ usually corresponds to $dT$, 
evaluated at $P_\lambda$, so we have chosen the value of $dT$ at $P_\lambda$ 
as our basis for $E_\bbR=\sigma_\lambda^*\Omega^1_{\bbA^1/\bbR}$.

\nxsubpoint
Let's compute the induced $\Zinfty$-structure. Roughly speaking, 
we have to compute the quotient norm $\vnorm_\lambda$ 
of the restriction of the $L_1$-norm 
of $\bbR[T]$ onto hyperplane $\gp'_\lambda=\{F:F(\lambda)=0\}$, and the 
quotient is computed with respect to the map $\gp'_\lambda\to E_\bbR=\bbR$ 
given by $F\mapsto F'(\lambda)$. We are actually interested in the value 
of our cometric on $(dT)_{\lambda}\in E_\bbR$, which has been identified with 
$1\in\bbR$.

\nxsubpoint\label{sp:comp.ar.aff}
The result of this computation is the following. For each $n\geq1$ denote by 
$\alpha_n$ the only positive root of polynomial 
$\lambda^{n+1}+(n+1)\lambda-n$, i.e.\ the only positive 
$\lambda$ for which $\lambda^n=n\lambda^{-1}-(n+1)$, 
and put $\alpha_0:=0$. Then we have 
$0=\alpha_0<\alpha_1<\alpha_2<\cdots<\alpha_n<\cdots<1$. The claim is that 
for $\alpha_{n-1}\leq\lambda\leq\alpha_n$ we have 
$\|dT\|_\lambda=(1+\lambda^n)/(n\lambda^{n-1})=\inf_{k\geq1}
(1+\lambda^k)/(k\lambda^{k-1})$ and 
$E=(n\lambda^{n-1})/(1+\lambda^n)\cdot\Zinfty$. For $-1<\lambda\leq0$ 
we have $\vnorm_{-\lambda}=\vnorm_\lambda$, and for $|\lambda|\geq1$ we have 
$\vnorm_\lambda=0$.

\nxsubpoint\label{sp:comp.ar.aff.prop}
Let's check first the properties of $\alpha_n$ that have been implicitly 
stated above. For any integer $n\geq1$ put 
$G_n(\lambda):=\lambda^{n+1}+(n+1)\lambda-n$. Since $G_n(0)=-n<0$ and 
$G_n(1)=2>0$, this polynomial has a root $\alpha_n$ 
between 0 and 1. On the other hand, 
$G'_n(\lambda)>0$ for $\lambda>0$, so $\alpha_n$ is the only positive root 
of $G_n$. We have just seen that $0<\alpha_n<1$; let's check 
$\alpha_{n-1}<\alpha_n$ for $n\geq2$. For this note that 
$\lambda G_{n-1}(\lambda)-G_n(\lambda)=n\lambda^2-(n-1)\lambda-(n+1)\lambda+n=
n(\lambda-1)^2>0$ for $\lambda<1$. Apply this for $\lambda=\alpha_{n-1}<1$: 
then $G_{n-1}(\lambda)=0$, so we obtain $G_n(\alpha_{n-1})<0$. Since 
$G_n$ is strictly increasing, we get $\alpha_{n-1}<\alpha_n$.

We have shown that $0=\alpha_0<\alpha_1<\alpha_2<\cdots<\alpha_n<\cdots<1$. 
Put $\phi_k(\lambda):=(1+\lambda^k)/(k\lambda^{k-1})$ for any $\lambda>0$ and 
any integer $k\geq1$, and put $\phi(\lambda):=\inf_{k\geq1}\phi_k(\lambda)$. 
It has been implicitly stated that for any $0<\lambda<1$ we can find 
an integer $n\geq1$, for which $\alpha_{n-1}\leq\lambda\leq\alpha_n$, and then 
$\phi(\lambda)=\phi_n(\lambda)$, i.e.\ 
$\phi_k(\lambda)\geq\phi_n(\lambda)$ for all $k\geq1$.

To show existence of~$n$ we have to check that $\sup_n\alpha_n=1$. Let's prove 
that $n/(n+2)<\alpha_n<1$. We know that $G_n(1)=2>0$, and for $\lambda=n/(n+2)$
we have $G_n(\lambda)=(n+1+\lambda^n)\lambda-n<(n+2)\lambda-n=0$, hence 
$\alpha_n$ is indeed between $n/(n+2)$ and $1$.

Finally, let's show that for $0<\lambda\leq\alpha_n$ 
we have $\phi_n(\lambda)\leq
\phi_{n+1}(\lambda)$, and for $\lambda\geq\alpha_n$ we have 
$\phi_n(\lambda)\geq\phi_{n+1}(\lambda)$. Then for $\alpha_{n-1}\leq\lambda
\leq\alpha_n$ we deduce by induction $\phi_n(\lambda)\leq\phi_{n+1}(\lambda)
\leq\phi_{n+2}(\lambda)\leq\cdots$ and $\phi_n(\lambda)\leq\phi_{n-1}(\lambda)
\leq\cdots\leq\phi_1(\lambda)$, hence $\phi_k(\lambda)\geq\phi_n(\lambda)$ for 
any integer $k\geq1$ and $\phi(\lambda)=\phi_n(\lambda)$.

Let's compute $\phi_n(\lambda)-\phi_{n+1}(\lambda)=
\frac{1+\lambda^n}{n\lambda^{n-1}}-\frac{1+\lambda^{n+1}}{(n+1)\lambda^n}=
\bigl(n(n+1)\lambda^n\bigr)^{-1}\cdot 
\bigl((n+1)\lambda(1+\lambda^n)-n(1+\lambda^{n+1})\bigr)
=\bigl(n(n+1)\lambda^n\bigr)^{-1}\cdot G_n(\lambda)$, so the sign 
of this expression coincides with that of $G_n(\lambda)$, and $G_n$ is 
strictly increasing for $\lambda>0$. This proves our inequality.

\nxsubpoint
Recall that we have to restrict the $L_1$-norm of $\bbR[T]$ onto 
$\gp'_\lambda$, and compute the quotient norm of this restriction with 
respect to $\gp'_\lambda\to\gp'_\lambda/{\gp'_\lambda}^2=E_\bbR\cong\bbR$. 
Instead of doing this, we can first compute the quotient norm on 
$\bbR[T]/{\gp'_\lambda}^2\cong\bbR^2$, and then restrict it to $E_\bbR$. 
This is clear in the language of norms; since they do not define completely
the corresponding $\Zinfty$-structure, we would like to prove this 
statement directly for $\Zinfty$-modules.

\begin{LemmaD}\label{l:lem.monoepi}
Let $i:N\to M$ be a strict monomorphism and $\pi:M\to M'$ be a strict 
epimorphism in $\ZinfFlat$. Suppose we are given a decomposition 
$\pi\circ i:N\stackrel\sigma\to N'\stackrel{j}\to M'$ of $\pi\circ i$ into an 
epimorphism $\sigma$, followed by a monomorphism $j$ 
(cf.~\ptref{sp:zinf.coim}, where existence of such decompositions with 
a strict $\sigma$ is shown). 
Suppose that the square, obtained by applying $\rho^*:\ZinfFlat\to\RVect$ to 
the commutative square~\eqref{leq:lem.monoepi} below, 
is bicartesian in $\RVect$ (since the 
vertical arrows are epic and horizontal arrows are monic, this condition 
is equivalent to $\Im i_\bbR\supset\Ker\pi_\bbR$):
\begin{equation}\label{leq:lem.monoepi}
\xymatrix{
N\ar[r]^{i}\ar[d]_{\sigma}&M\ar[d]_{\pi}\\
N'\ar[r]^{j}&M' }
\end{equation}
Then $j:N'\to M'$ is a strict monomorphism, $\sigma:N\to N'$ is 
a strict epimorphism, and this square is bicartesian in $\ZinfFlat$. 
The decomposition $j\circ\sigma$ of $\pi\circ i$ with the above properties 
is unique.
\end{LemmaD}

\begin{Proof}
a) Let's check that the square is cartesian. The explicit construction of 
projective limits given in~\ptref{p:projlimzflat} shows that this is 
equivalent to proving that it becomes cartesian after applying $\rho^*$, 
something that we know already, and that $N_{\Zinfty}=i_\bbR^{-1}(M_{\Zinfty})
\cap\sigma_\bbR^{-1}(N'_{\Zinfty})$. This is clear, since 
$\sigma_\bbR^{-1}(N'_{\Zinfty})$ contains $N_{\Zinfty}$ (this is true 
for any morphism), and $i_\bbR^{-1}(M_{\Zinfty})=N_{\Zinfty}$ since 
$i$ is a strict monomorphism (cf.~\ptref{sp:zinfstrmono}).

b) Let's check that the square is cocartesian. Since we know already that 
it becomes cocartesian after applying $\rho^*$, this amounts to check  
$M'_{\Zinfty}=\conv\bigl(j_\bbR(N'_{\Zinfty})\cup\pi_\bbR(M_{\Zinfty})\bigr)$.
Again, $j_\bbR(N'_{\Zinfty})\subset M'_{\Zinfty}$ for general reasons, and 
$\pi_\bbR(M_{\Zinfty})=M'_{\Zinfty}$ because $\pi$ is a strict epimorphism 
(cf.~\ptref{sp:zinfstrepi}).

c) By assumption, $j$ is a monomorphism, i.e.\ $j_\bbR$ is injective. 
We have to show that $j$ 
is a strict monomorphism, i.e.\ that $j_\bbR^{-1}(M'_{\Zinfty})\subset
N'_{\Zinfty}$, the opposite inclusion being evident. So let's take some 
$y'\in N'_\bbR$, such that $x':=j_\bbR(y')\in M'_{\Zinfty}$. Since 
$\pi$ is a strict epimorphism, $\pi_\bbR(M_{\Zinfty})=M'_{\Zinfty}$, so 
we can find some $x\in M_{\Zinfty}$, such that $\pi_\bbR(x)=x'$. We know 
that after application of $\rho^*$ our square becomes cartesian in 
$\RVect$, hence in $\catSets$, so $N_\bbR\cong M_\bbR\times_{M'_\bbR}N'_\bbR$,
and, since $x$ and $y'$ have same image in $M'_\bbR$, we can find some 
element $y\in N_\bbR$, such that $i_\bbR(y)=x$ and $\sigma_\bbR(y)=y'$. 
Since $i$ is a strict monomorphism and $i_\bbR(y)=x\in M_{\Zinfty}$, we see 
that $y\in N_{\Zinfty}$. Since $\sigma_\bbR(N_{\Zinfty})\subset N'_{\Zinfty}$, 
we see that $y'=\sigma_\bbR(y)\in\sigma_\bbR(N_{\Zinfty})\subset N'_{\Zinfty}$.
This reasoning is valid for any $y'\in j_\bbR^{-1}(M'_{\Zinfty})\supset 
N'_{\Zinfty}$, so we see that both $j_\bbR^{-1}(M'_{\Zinfty})$ and 
$\sigma_\bbR(N_{\Zinfty})$ coincide with $N'_{\Zinfty}$, i.e.\ 
$j$ is a strict monomorphism, and $\sigma$ is a strict epimorphism, q.e.d.
\end{Proof}

\nxsubpoint
An immediate consequence of this lemma is that a $\Zinfty$-structure 
$V_{\Zinfty}$
on an $\bbR$-vector space $V$ induces a canonical $\Zinfty$-structure on 
any its subquotient $V'/V''$, $V''\subset V'\subset V$, which can be defined 
either by first restricting the original $\Zinfty$-structure to $V'$ 
($V'_{\Zinfty}:=V_{\Zinfty}\cap V'$), and then taking its quotient 
($(V'/V'')_{\Zinfty}:=\Im (V'_{\Zinfty}\to V'/V'')$), or by taking 
first the quotient $V/V''$, and then restricting the quotient 
$\Zinfty$-structure to $V'/V''\subset V/V''$. To see this one has 
just to apply the lemma for $N=V'$, $M=V$, $N'=V'/V''$, $M'=V/V''$ 
with the induced $\Zinfty$-structures described above; for $V'/V''$ we 
take the quotient of that of $V'$, and obtain that $V'/V''\to V/V''$ 
is strict, i.e.\ this $\Zinfty$-structure of $V'/V''$ coincides 
with that induced from $V/V''$.

\nxsubpoint
In our situation we apply this lemma for $M=\Zinfty[T]$, 
$N=\gp_\lambda=\gp'_\lambda\cap\Zinfty[T]$, $N'=\gp_\lambda/\gp_\lambda^2=E$ 
and $M'=D:=\Zinfty[T]/\gp_\lambda^2$. The latter notation is understood as 
the strict quotient $D$ of $\Zinfty[T]$, such that 
$D_\bbR=\bbR[T]/{\gp'_\lambda}^2$. Note that $D_\bbR$ can be identified 
with $\bbR^2$, and the projection $\pi_\bbR:\bbR[T]\to D_\bbR$ with the map 
$F(T)\mapsto(F(\lambda),F'(\lambda))$. Then $E_\bbR=\bbR$ is identified 
with the coordinate axis $0\times\bbR$ in $D_\bbR$.

Now we have to compute the strict quotient $\Zinfty$-structure on $D_\bbR$, 
given by $D_{\Zinfty}=\pi_\bbR(\Zinfty[T])$. Since 
$\Zinfty[T]=\conv(\{\pm T^n\}_{n\geq0})$,
we see that
\begin{equation}\label{eq:infcvhull}
D_{\Zinfty}=\conv\bigl(\{\pm(1,0)\}\cup
\{\pm(\lambda^n, n\lambda^{n-1})\}_{n\geq1}\bigr).
\end{equation}
After computing this 
convex hull we will have $E_{\Zinfty}=D_{\Zinfty}\cap E_\bbR$, since 
$D\to E$ is a strict monomorphism according to the above lemma.

The following theorem is quite helpful for computing infinite convex hulls 
like~\eqref{eq:infcvhull}:

\begin{ThD} (``octahedral Carath\'eodory theorem'') 
Let $S$ be a subset of an $n$-dimensional real vector space~$V$. Denote 
by $M$ the $\Zinfty$-submodule of~$V$, generated by~$S$, 
i.e.\ the smallest $\Zinfty$-submodule containing $S$ 
(i.e.\ such that $S\subset M_{\Zinfty}$); 
in other words, $M=(M_{\Zinfty},M_\bbR)$, where 
$M_{\Zinfty}=\conv(S\cup-S)$ is the symmetric convex hull of~$S$
and $M_\bbR=\bbR\cdot M_{\Zinfty}$ is the $\bbR$-span of~$S$.

Then any point $x\in M_{\Zinfty}$ can be written as an octahedral combination 
of at most $\dim M_\bbR\leq\dim V=n$ linearly independent elements of~$S$:
\begin{equation}\label{eq:octcar}
x=\sum_{i=1}^k \lambda_i x_i,\quad x_i\in S,\quad 
\sum_{i=1}^k |\lambda_i|\leq1.
\end{equation}
\end{ThD}

\begin{Proof}
Clearly, any $x\in M_{\Zinfty}$ can be written as an octahedral combination 
of a finite subset of~$S$, i.e.\ it has a representation \eqref{eq:octcar} 
for some $k\geq0$. Fix some $x$ and choose such a representation with 
minimal $k$. We are going to show that the $x_i$ in such a minimal 
representation must be linearly independent; since they lie in $M_\bbR$, 
this will imply that $k\leq\dim M_\bbR\leq n$. Indeed, suppose 
we have a non-trivial linear relation $\sum_{i=1}^k \mu_i x_i=0$. 
After renumbering the $x_i$, $\lambda_i$ and $\mu_i$, we can assume 
that $\lambda_i>0$ for $1\leq i\leq r$ and $\lambda_i<0$ for 
$r+1\leq i\leq k$ for some $0\leq r\leq k$ (note that $\lambda_i=0$ would 
contradict minimality of $k$). Changing the signs of all $\mu_i$ if 
necessary, we can assume that $\sum_{i=1}^r\mu_i\geq\sum_{i=r+1}^k\mu_i$. 

Observe that this choice of signs guarantees that the set of indices 
$i$ for which $\lambda_i\mu_i>0$ is not empty: otherwise we would have 
$\mu_i\leq 0$ for $i\leq r$ and $\mu_i\geq 0$ for $i>r$, hence 
$\sum_{i=1}^r\mu_i\leq 0\leq\sum_{i=r+1}^k\mu_i$; together with the opposite 
inequality this implies $\mu_i=0$ for all $i$. This contradicts the 
non-triviality of chosen linear relation between the $x_i$.

Now put $\rho:=\min_{\lambda_i\mu_i>0}\lambda_i/\mu_i$. We have just 
checked that this minimum is computed along a non-empty set, 
hence $\rho$ is well-defined; clearly, $\rho>0$. Put 
$\lambda'_i:=\lambda_i-\rho\mu_i$. Then $x=\sum_{i=1}^k\lambda'_i x_i$, 
and at least one of $\lambda'_i$ is zero by the choice of $\rho$. 
Let's check that this combination is still octahedral, i.e.\ 
$\sum_{i=1}^k|\lambda'_i|\leq1$; this will give a contradiction 
with the choice of $k$, and thus it will prove linear independence of $x_i$.

From the choice of $\rho$ we see that all $\lambda'_i$, if non-zero, have 
the same signs as $\lambda_i$: $\lambda'_i\geq0$ for $i\leq r$ and 
$\lambda'_i\leq0$ for $i>r$. Hence $\sum_i|\lambda'_i|=
\sum_{i\leq r}\lambda'_i-\sum_{i>r}\lambda'_i=\sum_i|\lambda_i|-\rho\delta$, 
where $\delta:=\sum_{i\leq r}\mu_i-\sum_{i>r}\mu_i\geq 0$ according to 
our choice of signs of $\mu_i$. We also have $\rho>0$; hence 
$\sum_i|\lambda'_i|\leq\sum_i|\lambda_i|\leq 1$, i.e.\ our new expression 
$x=\sum_i\lambda'_i x_i$ is an octahedral combination with 
smaller number of non-zero terms, q.e.d.
\end{Proof}

\nxsubpoint
This theorem has an obvious $p$-adic counterpart. Namely, if $M$ 
is the $\bbZ_p$-submodule generated by an arbitrary subset $S$ 
of a finite-dimensional $\bbQ_p$-vector space $V$, then any element of $M$ 
can be written as a $\bbZ_p$-linear combination of several 
$\bbQ_p$-linearly independent elements of $S$; hence the total number of 
these elements doesn't exceed $\dim_{\bbQ_p}M_{(\bbQ_p)}\leq\dim_{\bbQ_p}V$.

\nxsubpoint
Let's apply this theorem to compute $D_{\Zinfty}=\conv(\pm(\lambda^n,
n\lambda^{n-1}))\subset D_\bbR=\bbR^2$. We see that $D_{\Zinfty}=
\bigcup_{k,l\geq 0} A_{kl}$, where $A_{kl}:=\conv\bigl(\pm(\lambda^k,
k\lambda^{k-1})$, $\pm(\lambda^l,l\lambda^{l-1})\bigr)$ is the 
$\Zinfty$-sublattice generated by these two generators of $D_{\Zinfty}$. 
Therefore, $E_{\Zinfty}=D_{\Zinfty}\cap (0\times\bbR)=\bigcup_{k,l\geq0}
B_{kl}$, where $B_{kl}:=A_{kl}\cap (0\times\bbR)$ 
is a compact symmetric convex subset of $0\times\bbR\cong\bbR$, 
hence $B_{kl}=0\times[-b_{kl},+b_{kl}]$ for some real number $b_{kl}\geq0$.

\nxsubpoint
If $\lambda=0$, then $D_{\Zinfty}=\conv((\pm1,0),(0,\pm1))=\Zinfty^{(2)}$ 
is the standard two-dimensional octahedron, and $E_{\Zinfty}=\Zinfty$ 
as stated in~\ptref{sp:comp.ar.aff}. Case $\lambda<0$ is reduced to 
$\lambda>0$ by means of the automorphism $T\mapsto-T$ of 
$\bbR[T]$ that preserves $\Zinfty[T]$, maps $\lambda$ into $-\lambda$ 
and $dT$ into $-dT$, hence $\|dT\|_{\lambda}=\|-dT\|_{-\lambda}=
\|dT\|_{-\lambda}$. So we can assume $\lambda>0$.

\nxsubpoint
Let's compute the numbers $b_{kl}\geq0$, always assuming $\lambda>0$. 
Clearly, $b_{kk}=0$ for any $k\geq0$ and $b_{kl}=b_{lk}$, so we 
can assume $k>l\geq0$. Now, $y\in B_{kl}$ is equivalent to 
$(0,y)\in A_{kl}$, i.e.\ to the existence of $u,v\in\bbR$ with $|u|+|v|\leq1$, 
such that $(0,y)=u(\lambda^k,k\lambda^{k-1})+v(\lambda^l,l\lambda^{l-1})$. 
This is equivalent to $\lambda^ku+\lambda^lv=0$ and $y=k\lambda^{k-1}u+
l\lambda^{l-1}v$. Since both $\lambda^k$ and $\lambda^l$ are non-zero, 
the first of these conditions is equivalent to $u=w\lambda^l$, 
$v=-w\lambda^k$ for some real $w$. Then $|u|+|v|\leq1$ can be rewritten 
as $|w|\leq(\lambda^k+\lambda^l)^{-1}$, and our equation for $y$ becomes 
$y=(k-l)\lambda^{k+l-1}w$, so the maximal possible value of~$y$ is 
$b_{kl}=(k-l)\lambda^{k+l-1}/(\lambda^k+\lambda^l)=
(k-l)\lambda^{k-1}/(\lambda^{k-l}+1)$.

\nxsubpoint
If $\lambda\geq1$, then for $l=1$ and $k\gg1$ we have
$b_{k1}>k/3$, hence $E_{\Zinfty}$, being the union of all 
$[-b_{kl},b_{kl}]$, coincides with the whole $\bbR$; the corresponding 
seminorm on $E_\bbR=\bbR$ is identically zero, as stated 
in~\ptref{sp:comp.ar.aff}. So we restrict ourselves to the case 
$0<\lambda<1$.

In this case $b_{kl}=(k-l)\lambda^{k-1}/(1+\lambda^{k-l})\leq 
b_{k0}=k\lambda^{k-1}/(1+\lambda^k)=\phi_k(\lambda)^{-1}$ for any 
$0\leq l<k$, hence $E_{\Zinfty}$ is equal to the union of 
$[-b_{k0},b_{k0}]=[-\phi_k(\lambda)^{-1},\phi_k(\lambda)^{-1}]$ for 
all integer $k\geq1$. We know that for any $0<\lambda<1$ the sequence 
$\phi_k(\lambda)$ has a minimal element $\phi(\lambda)=\phi_n(\lambda)$, 
where $n$ is determined by $\alpha_{n-1}\leq\lambda\leq\alpha_n$ 
(cf.~\ptref{sp:comp.ar.aff.prop}). Therefore, 
$E_{\Zinfty}=[-\phi(\lambda)^{-1},\phi(\lambda)^{-1}]=\phi(\lambda)^{-1}\cdot
\Zinfty$, and $\|dT\|_\lambda=\|1\|_E=\phi(\lambda)$. This finishes the 
proof of~\ptref{sp:comp.ar.aff}.

\nxsubpoint
Thus we have obtained a singular cometric $\vnorm_\lambda$ on the 
real affine line, given by $\|dT\|_\lambda=\phi(|\lambda|)$ for 
$|\lambda|<1$ and $\vnorm_\lambda=0$ for $|\lambda|\geq 1$. We see 
that this cometric is continuous and piecewise smooth for $|\lambda|<1$; 
we want to study its asymptotic behavior as $\lambda\to1-$. We will 
show that $\|dT\|_\lambda=\phi(\lambda)=e^{\kappa}(1-\lambda)+
O((1-\lambda)^2)$ as $\lambda\to1-$, for some positive constant 
$\kappa>0$. This will imply the continuity of $\vnorm_\lambda$ at all 
$\lambda\in\bbR$.

Denote by $\kappa$ the only real solution of $\kappa=1+e^{-\kappa}$. 
Numerically, $\kappa=1,2784645\ldots$. We will need the asymptotic 
behavior of $\alpha_n$ as $n\to+\infty$. Notice for this that 
for any $\lambda>0$ we have $G'_{n-1}(\lambda)=n\lambda^{n-1}+n>n$, 
hence by the mean value theorem we get $|\lambda-\alpha_{n-1}|\leq n^{-1}
|G_{n-1}(\lambda)|$ for any $\lambda>0$. Now put $\lambda:=1-\kappa/n$; 
then $n\log\lambda=-\kappa+O(n^{-1})$, hence $\lambda^n=e^{-\kappa}+O(n^{-1})$
and $G_{n-1}(\lambda)=\lambda^n+n\lambda-(n-1)=e^{-\kappa}+O(n^{-1})+
(n-\kappa)-(n-1)=O(n^{-1})$ according to the choice of $\kappa$, hence 
$\bigl|(1-\kappa/n)-\alpha_{n-1}\bigr|=O(n^{-2})$ and 
$\alpha_{n-1}=1-\kappa/n+O(n^{-2})$. Replacing $n$ by $n+1$, we get 
$\alpha_n=1-\kappa/n+O(n^{-2})$ as well, so $\alpha_{n-1}\leq\lambda\leq
\alpha_n$ implies $\lambda=1-\kappa/n+O(n^{-2})$, hence 
$n=\kappa(1-\lambda)^{-1}+O(1)$.

Now let's estimate $\phi(\lambda)=\phi_n(\lambda)=
(1+\lambda^n)/(n\lambda^{n-1})=(\lambda/n)(1+\lambda^{-n})$. 
We know that $\lambda=1-\kappa/n+O(n^{-2})$, 
hence $\lambda^n=e^{-\kappa}+O(n^{-1})$ as before, hence 
$\phi(\lambda)=n^{-1}(1+e^{\kappa}+O(n^{-1}))=(1+e^\kappa)n^{-1}+O(n^{-2})$. 
Now we combine this with $1+e^\kappa=\kappa e^\kappa$ and with the 
asymptotic expression for $n$. We get $\phi(\lambda)=e^{\kappa}(1-\lambda)+
O((1-\lambda)^2)$ for $\lambda\to1-$, as claimed above.

\nxsubpoint (Complex Arakelov affine line.)
Of course, we can consider the {\em complex Arakelov affine line}
$\bbA^1_{\barZinfty}=\bbA^1_{\Zinfty}\otimes_{\Zinfty}\barZinfty=
\Spec \barZinfty[T]$, where $\barZinfty[T]=\barZinfty\otimes_{\Zinfty}
\Zinfty[T]$ is the polynomial algebra over $\barZinfty$ in one variable 
$T$, easily seen to satisfy the usual universal property in the category 
of $\barZinfty$-algebras with one pointed element.

Then we might ask about the induced cometrics $\vnorm_\lambda$ for 
all $\lambda\in\bbA^1(\bbC)=\bbC$. However, the answer remains essentially 
the same as before: $\|dT\|_\lambda=\phi(|\lambda|)$ for 
$|\lambda|<1$, and $\|dT\|_\lambda=0$ for $|\lambda|\geq1$. To prove this 
we reduce first to the case $\lambda\geq0$, by observing that for 
any complex $\epsilon$ with $|\epsilon|=1$ (e.g.\ $|\lambda|/\lambda$) 
the map $T\mapsto\epsilon T$ induces an automorphism of $\barZinfty[T]$, 
hence of $\bbA^1_{\barZinfty}$, so we must have 
$\|dT\|_\lambda=\|\epsilon\,dT\|_{\epsilon\lambda}=\|dT\|_{\epsilon\lambda}=
\|dT\|_{|\lambda|}$ if we take $\epsilon:=|\lambda|/\lambda$.

Now we want to say that our previous computations for a real $\lambda$ and 
the corresponding point $P_\lambda\in\bbA^1(\bbR)$ remain valid after 
extending everything to $\bbC$ (or, more precisely, to $\barZinfty$). 
To prove this we consider again the strict epimorphism 
$\Zinfty[T]\to F=\Zinfty[T]/\gp_\lambda^2$ and the strict monomorphism 
$E=\gp_\lambda/\gp_\lambda^2\to F$, and apply our complexification functor 
$\barZinfty\otimes_{\Zinfty}-$ to them. Clearly, this functor preserves 
strict epimorphisms (i.e.\ cokernels of pairs of morphisms), being right exact.
If we check that $E_{(\barZinfty)}\to F_{(\barZinfty)}$ remains a strict 
monomorphism, we'll obtain that the subquotient $\barZinfty$-structure on 
$E_\bbC$ coincides with $E_{(\barZinfty)}$; since 
$E_{\Zinfty}=\phi(\lambda)^{-1}\cdot\Zinfty$, we get $E_{\barZinfty}=
\phi(\lambda)^{-1}\cdot\barZinfty$ and $\|z\,dT\|_\lambda=
|z|\cdot\phi(|\lambda|)$ as claimed above.

So we want to check that $j:E\to F$ remains a strict monomorphism after 
tensor multiplication with $\barZinfty$. For this we show that $j$ 
admits a left inverse $p:F\to E$. Recall that $E_\bbR\cong\bbR$ and that 
$E_{\Zinfty}$ is closed, hence $E\cong\bbR$ or $E\cong\Zinfty$. 
In the first case existence of a left inverse is evident; in the second 
we apply Hahn--Banach theorem to linear form 
$1_E: E_\bbR\to E_\bbR=\bbR$ and extend it to a linear form 
$p:F_\bbR\to\bbR=E_\bbR$ 
with norm $\leq1$, i.e.\ to a morphism $p:F\to E$. Once existence of 
$p$ is proved, we observe that $j\otimes 1_P$ also admits a left inverse 
$p\otimes 1_P$ for any $P$ (e.g.\ $P=\barZinfty$), hence it is a 
strict monomorphism, being the kernel of the pair $1_{E\otimes P}$ and 
$(jp)\otimes 1_P:E\otimes P\rightrightarrows E\otimes P$.

\nxsubpoint (Projective line.) Of course, now we can take two copies of 
$\bbA^1_{\Zinfty}$ and glue them along open subsets $\bbG_{m,\Zinfty}=
\Spec\Zinfty[T,T^{-1}]$ by means of the map given by $T\mapsto T^{-1}$ 
into an {\em Arakelov projective line~$\bbP^1_{\Zinfty}$}. 
The corresponding cometric $\vnorm_\lambda$ for $|\lambda|<1$ 
will be the same, since these points 
lift to $\Zinfty$-sections already in $\bbA^1_{\Zinfty}$. For 
$|\lambda|>1$ we observe that $T\mapsto T^{-1}$ defines an automorphism 
of $\bbP^1$, hence $\|dT\|_{\lambda}=\|d(T^{-1})\|_{\lambda^{-1}}
=|\lambda|^2\cdot\|dT\|_{\lambda^{-1}}=|\lambda|^2\phi(|\lambda|^{-1})$. 
In particular, for $\lambda\to1+$ we have 
$\|dT\|_\lambda=e^\kappa(\lambda-1)+O((\lambda-1)^2)$. We cannot 
make a detailed analysis for $|\lambda|=1$ now, but it seems natural to assume 
``by continuity'' that $\|dT\|_\lambda=0$ for $|\lambda|=1$.

So we have obtained a continuous singular cometric on $\bbP^1(\bbC)$ 
with singularities along the unit circle. The corresponding metric 
will have a ``simple pole'' along this circle. Note that the singular 
locus is not Zariski closed in this situation; it is actually Zariski dense.

\nxsubpoint (Hyperbolic plane?) 
Once we have a metric on $\bbC$, we can define the length $l(\gamma)$ of 
a piecewise $C^1$-curve $\gamma:[a,b]\to\bbC$ by the usual formula 
$l(\gamma)=\int_a^b \|\gamma'(t)\|_{\gamma(t)}dt$. Then we can define the 
distance between two points of $\bbC$ as the infimum of lengths of 
curves connecting these two points, thus obtaining a metric space structure 
on $\bbC$ and on $\bbP^1(\bbC)$. Note that the unit circle is at infinite 
distance from any point either inside or outside this circle. 
When $\lambda\to1-$, the distance between $0$ and $\lambda$ is 
$\rho(0,\lambda)=e^{-\kappa}\log(1-\lambda)+O(1)$. 
In this respect the interior of the unit circle, 
equipped with this metric, behaves like the 
usual Poincar\'e model of the hyperbolic plane inside the unit circle. 
Of course, our metric is only piecewise smooth.

\nxsubpoint (Origin of singularities.)
Let's give some indications why Arakelov affine and projective lines 
constructed here turned out to have singular metrics. Later 
we will construct projective spaces $\bbP_{\Zinfty}(V)=
\Proj S_{\Zinfty}(V)$ for any $\Zinfty$-lattice~$V$, and in particular we 
will get some (co)metrics on the complex points of $\bbP(V_\bbR)$.

The $\Zinfty$-sections of $\bbP(V)$ will correspond to strict quotients of 
$V$, isomorphic to $\Zinfty$, i.e.\ to strict epimorphisms 
$V\to\Zinfty$ modulo $\pm1$. By duality, they correspond to strict 
monomorphisms $\Zinfty\to V^*$ modulo sign, i.e.\ to the boundary
points $\partial V^*=V^*-\gm_\infty V^*$, again modulo sign.

Now, the projective line constructed in this section corresponds to 
$V=\Zinfty^{(2)}$, hence $V^*=\Zinfty^2$ is the square with vertices 
$(\pm1,\pm1)$. Its boundary is smooth outside these vertices, which correspond
exactly to $\lambda=\pm1$. So we can think of singularity of our 
cometric at these points as a consequence of non-smoothness of 
$\partial V^*$ at corresponding points. Our other critical points 
$\alpha_n$, where the (co)metric is not smooth, condense at $|\lambda|=1$, 
so they could be thought of as ``echoes'' of this singularity.

So we are inclined to think that if $V^*$ has a smooth boundary, then 
the corresponding (co)metric on the projective space will be 
everywhere smooth and non-singular. Some evidence to this is given by the 
quadratic and hermitian $V$'s (i.e.\ $V_{\Zinfty}=\{x:Q(x)\leq1\}$, 
where $Q$ is a positive definite quadratic or hermitian form 
on $V_\bbR$ or $V_\bbC$), which will be shown later 
to give rise to Fubiny--Study metrics on corresponding projective spaces.

\nxpointtoc{Spectra of flat $\Zinfty$-algebras}\label{p:spec.flat.zalg}
We will obtain later the theory of spectra of (arbitrary) $\Zinfty$-algebras 
as a part of theory of spectra of generalized rings. For several reasons 
this seems to be a more natural approach than to develop the theory 
only for one special case. However, we would like to show the 
impatient reader a glimpse of this special case of the theory.
The exposition will be somewhat sketchy here, since we are going to give 
the rigorous proofs later in a more general context.

\nxsubpoint
It will be more convenient to think of some object $M=(M_{\Zinfty},M_\bbR)$ 
of $\ZinfFlat$ as a set $M=M_{\Zinfty}$, equipped by a flat $\Zinfty$-module 
structure. By definition, a (flat) $\Zinfty$-module structure on some 
set $M$ is simply an embedding $j_M:M\to M_\bbR$ of $M$ into some 
real vector space $M_\bbR$, such that $j_M(M)$ becomes a symmetric convex 
body in $M_\bbR$. We say that a map of sets $f:M\to N$ is compatible with 
given $\Zinfty$-structures on them if $f$ extends to a $\bbR$-linear map 
$f_\bbR:M_\bbR\to N_\bbR$, necessarily unique.

Of course, this is just another description of $\ZinfFlat$; we have only 
switched our focus of interest to $M_{\Zinfty}$. In particular, a 
(commutative flat) $\Zinfty$-algebra $A$ can be now described as a 
(commutative) monoid, equipped by a flat $\Zinfty$-module structure, 
with respect to which the multiplication becomes $\Zinfty$-bilinear.

\nxsubpoint
The notions of a multiplicative system $S\subset A$ and an invertible 
element $s\in A$ involve only the multiplicative (i.e.\ monoid) 
structure of $A$, hence they are well-defined in our situation.

\nxsubpoint (Ideals.)
An {\em ideal\/} of $A$ is by definition an $A$-submodule $\ga\subset A$, 
i.e.\ a $\Zinfty$-submodule, stable under the multiplication of $A$: 
$A\cdot\ga\subset A$. The {\em intersection\/} of two ideals $\ga\cap\gb$ is 
understood in the usual way. The {\em sum\/} $\ga+\gb=(\ga,\gb)$ is understood 
as the ideal generated by $\ga\cup\gb$, i.e.\ the smallest ideal 
containing both $\ga$ and $\gb$; it consists of all octahedral combinations 
of elements of~$\ga$ and~$\gb$, i.e.\ 
$\ga+\gb=\{\lambda x+\mu y:$ $|\lambda|+|\mu|\leq1$, $x\in\ga$, $y\in\gb\}$. 
Finally, the {\em product\/} $\ga\gb$ is by definition the ideal generated 
by all products $xy$, $x\in\ga$, $y\in\gb$. It can be also described as 
the image of $\ga\otimes_{\Zinfty}\gb\to A$.

An ideal $\gp$ in $A$ is {\em prime\/} if its complement is a multiplicative 
system, i.e.\ if $1\nin\gp$, and $s,t\nin\gp$ implies $st\nin\gp$. 
We say that $\gm$ is a {\em maximal ideal\/} if it $\neq A$ and if it is 
not contained in any ideal $\neq\gm$ and $\neq A$. Finally, we say that 
$A$ is {\em local\/} if it has exactly one maximal ideal.

Note that {\em any maximal ideal\/ $\gm$ is prime.} Indeed, if $s,t\nin\gm$, 
then $(\gm,s)=(\gm,t)=A$ by maximality, hence 
$(\gm,st)\supset(\gm^2,s\gm,t\gm,st)=(\gm,s)(\gm,t)=A$ and $st\nin\gm$.

Also note that the Krull theorem still holds: {\em Any ideal $\ga\neq A$ 
is contained in some maximal ideal\/ $\gm$}. For this we apply the Zorn lemma 
to the set of ideals of $A$, containing $\ga$ and distinct from $A$, 
and observe that filtered unions of ideals are ideals.

Some corollaries of these theorems also hold. For example, the union of all 
maximal ideals of~$A$ coincides with the complement of the set of invertible 
elements of~$A$.

Another interesting statement is that {\em $\ga\gb\subset\gp$ iff
$\ga\subset\gp$ or $\gb\subset\gp$, for a prime ideal\/ $\gp$}.

\nxsubpoint\label{sp:loczflat} (Localization.)
Given a multiplicative subset $S\subset A$ and an $A$-module $M$, we can 
construct their {\em localizations\/} $S^{-1}A$ and $S^{-1}M$ as follows. 
As a set, $S^{-1}M=M\times S/\sim$, where $(x,s)\sim(y,t)$ iff 
$utx=usy$ for some $u\in S$, and similarly for $S^{-1}A$, so this is the 
usual definition of localizations. The multiplication of $S^{-1}A$ and 
the action of $S^{-1}A$ on $S^{-1}M$ are defined in the usual way. To obtain 
$\Zinfty$-structures on $S^{-1}M$ and $S^{-1}A$, we write them down 
as filtered inductive limits: $S^{-1}M=\varinjlim_{s\in\cS} M_{[s]}$. 
Here the index category $\cS$ is defined as follows: $\Ob\cS = S$, and 
morphisms $[s]\to [s']$ are given by elements $t\in S$, such that $s'=ts$. 
Our inductive system is given by $M_{[s]}=M$ for all $s$, and the transition 
maps $M_{[s]}\to M_{[st]}$ are given by multiplication by $t$. They are 
clearly $\Zinfty$-linear, and filtered inductive limits in $\ZinfFlat$ 
coincide with those computed in $\catSets$ (cf.~\ptref{ss:zinffiltind}), 
so we indeed get a flat $\Zinfty$-structure on $S^{-1}M$ and $S^{-1}A$.

Then $S^{-1}A$ becomes a $\Zinfty$-algebra and even an $A$-algebra, 
and $S^{-1}M$ is an $S^{-1}A$-module; from their inductive limit description 
one can obtain their usual universal properties: $S^{-1}A$ is the universal 
$A$-algebra in which all elements of $S$ become invertible, and 
$M\to S^{-1}M$ is the universal $A$-linear map of $M$ into a module, 
on which all elements of $S$ act bijectively.

Of course, we introduce the usual notations $M_\gp$ and $M_f$ if 
$S$ is the complement of a prime ideal $\gp$, resp.\ if 
$S=S_f=\{1,f,f^2,\ldots\}$. If $g$ divides some power of~$f$, we get 
canonical homomorphisms $M_g\to M_f$.

In short, this notion of localization seems to possess all usual properties 
of localizations of commutative rings and modules over them.

\nxsubpoint (Definition of $\Spec A$.)
Now we define the {\em (prime) spectrum\/} $\Spec A$ 
to be the set of all prime ideals of $A$. 
For any $M\subset A$ we define $V(M)\subset\Spec A$ by 
$V(M):=\{\gp:\gp\supset M\}$. We have the usual properties: 
$V(0)=\Spec A$, $V(1)=\emptyset$, $V(\bigcup M_\alpha)=\bigcap V(M_\alpha)$, 
$V(M)=V(\langle M\rangle)$, where $\langle M\rangle$ is the ideal generated 
by $M$, and finally 
$V(\ga\gb)=V(\ga\cap\gb)=V(\ga)\cup V(\gb)$. These properties show that 
the $\{V(M)\}$ are the closed subsets of $X:=\Spec A$ for some topology on 
$X$, called the {\em Zariski topology}. Clearly, the principal open 
subsets $D(f):=X-V(\{f\})=\{\gp:f\nin\gp\}$ form a base of this topology. 
Note that $\bigcup D(f_\alpha)=X$ iff the $f_\alpha$ generate the unit 
ideal $A$, hence $X$ is quasicompact.

Any $\Zinfty$-algebra homomorphism $\phi: A\to B$ induces a continuous 
map ${}^a\phi:\Spec B\to\Spec A$ by the usual rule $\gp\mapsto\phi^{-1}(\gp)$.
If $B=S^{-1}A$, $\Spec B$ can be identified by means of this map with 
the subset $\{\gp: S\cap\gp=\emptyset\}\subset X$, equipped with the induced 
topology. In particular, $\Spec A_f$ is homeomorphic to $D(f)\subset X$.

\nxsubpoint (Quasicoherent sheaves.)
Any $A$-module $M$ defines a presheaf of $\Zinfty$-modules $\widetilde{M}$ 
on $X=\Spec A$ with respect to the base $\{D(f)\}$, given by 
$\Gamma(D(f),\widetilde{M}):=M_f$. These presheaves are actually sheaves, 
something that will be checked later in a more general context. 
In particular, we get the {\em structure sheaf} $\sO_X:=\widetilde{A}$. It is
a sheaf of $\Zinfty$-algebras, and each $\widetilde{M}$ is a sheaf of 
$\sO_X$-modules.
Of course, we have the usual functoriality 
with respect to $\Zinfty$-algebra homomorphisms $\phi:A\to B$; for 
quasicoherent sheaves (i.e.\ those isomorphic to some $\widetilde{M}$) 
this functoriality involves tensor products over~$A$ which haven't 
been explained in detail yet (cf.~\ptref{def:catAmod}).

\nxsubpoint (Flat $\Zinfty$-schemes.)
Now we define {\em flat $\Zinfty$-schemes\/} to be topological spaces equipped 
with a sheaf of $\Zinfty$-algebras, each point of which has a neighborhood 
isomorphic to $\Spec A$ for some flat $\Zinfty$-algebra~$A$. Flat 
$\Zinfty$-schemes, isomorphic to some $\Spec A$, are called {\em affine}. 
{\em Morphisms\/} of flat $\Zinfty$-schemes are just {\em local\/} morphisms 
of topological spaces equipped by local sheaves of flat $\Zinfty$-algebras. 
As usual, morphisms $X\to\Spec A$ are in one-to-one correspondence with 
$\Zinfty$-algebra homomorphisms $A\to\Gamma(X,\sO_X)$.
Note that the category of flat $\Zinfty$-schemes has fibered products, 
the existence of which is shown in the same way as in EGA~I, once 
we have the tensor products~$\otimes_A$.

\nxsubpoint (Schemes over $\Spec\bbR$.)
Any $\bbR$-algebra $A$ has a natural flat $\Zinfty$-algebra structure; 
$\Spec A$ in this case coincides with the usual spectrum of an affine 
ring. In other words, the category of $\bbR$-schemes can be identified 
with a full subcategory of the category of flat $\Zinfty$-schemes. 

\nxsubpoint
For example, $S:=\Spec\Zinfty=\{0,\gm_\infty\}$, and the only non-trivial 
open subset is $\{0\}$, so $\Spec\Zinfty$ looks exactly like $\Spec\bbZ_p$. 
The structural sheaf is given by $\Gamma(S,\sO_S)=\Zinfty$, 
$\Gamma(\{0\},\sO_S)=\bbR$ and $\Gamma(\emptyset,\sO_S)=0$.
Spectrum of $\barZinfty$ gives the same topological space, but of course 
with a different structure sheaf.

For any flat $\Zinfty$-algebra $A$ we have a canonical homomorphism 
$\Zinfty\to A$, hence there is a unique morphism from 
any flat $\Zinfty$-scheme $X$ into $\Spec\Zinfty$, i.e.\ 
$\Spec\Zinfty$ is the final object of this category. 
The generic fiber $X_\xi=X_\bbR$, i.e.\ the pullback of the open subset 
$\{0\}\cong\Spec\bbR$ with respect to this morphism $X\to\Spec\Zinfty$, 
is an $\bbR$-scheme in the usual sense. 
For example, the generic (or open) fiber of $\Spec A$ is $\Spec A_\bbR$. 
Not so much can be said about the closed fiber.

\nxsubpoint
In this way we see that any prime ideal~$\gp$ of~$A$ falls into 
one of two subcategories:

a) Ideals over $0$, i.e.\ $\gp\cap\Zinfty=0$. These are in one-to-one 
correspondence with (usual) prime ideals $\gp'$ in the $\bbR$-algebra $A$. 
This correspondence is given by $\gp\mapsto\gp_\bbR$, $\gp'\mapsto 
\gp'\cap A$; in other words, we consider on $\gp'$ the $\Zinfty$-structure 
induced from $A_\bbR$, i.e.\ $\gp$ is the strict subobject of $A$, 
corresponding to $\gp'\subset A_\bbR$.

b) Ideals over $\gm_\infty$, i.e.\ $\gp\cap\Zinfty=\gm_\infty$, hence 
$\gp\supset\gm_\infty$ and $\gp\supset\gm_\infty A$. In general we cannot 
say too much about these ideals. One of the reasons is that we cannot describe 
them as prime ideals in $A/\gm_\infty A$, since this quotient is zero 
in $\ZinfFlat$, and is $\Zinfty$-torsion in the category $\ZinfMod$ 
that will be constructed later.

\nxsubpoint\label{sp:examp.spec.zinfty2}
Let's consider $X=\Spec A$ for $A=\Zinfty^2=\Zinfty\times\Zinfty$. Contrary 
to what one might expect, $X\neq\Spec\Zinfty\sqcup\Spec\Zinfty$, because 
we can have prime ideals that contain both idempotents $e_1=(1,0)$ and 
$e_2=(0,1)$. In our case the open fiber consists of two points 
(closed in this fiber), given by ideals $0\times\Zinfty$ and $\Zinfty\times0$. 
However, the closed fiber consists of {\em three\/} points 
$0\times\gm_\infty$, $\gm_\infty\times0$ and $\gm:=A-\{(\pm1,\pm1)\}$. 
The latter is the only maximal ideal of~$A$, i.e.\ $A$ is local, and 
$\gm$ defines the only closed point of~$X$. The complement of this closed 
point is an open subscheme, isomorphic to $\Spec\Zinfty\sqcup\Spec\Zinfty$. 
We'll give an explanation for this phenomenon in~\ptref{sp:prop.unary.env}.

\nxsubpoint
Now we would like to consider the affine line $\bbA^1_{\Zinfty}=\Spec
\Zinfty[T]$ from this point of view. The prime ideals $\gp$ in the 
open fiber are in one-to-one correspondence with $\gp'\in\Spec\bbR[T]$, 
so we obtain immediately their complete description: there is the 
zero ideal $(0)$, ideals $\gp_\lambda:=\Zinfty[T]\cap(T-\lambda)\bbR[T]=
\{F(T)\in\Zinfty[T]:F(\lambda)=0\}$ for each real $\lambda\in\bbR$, 
and $\gp_z:=\Zinfty[T]\cap(T-z)(T-\bar{z})\bbR[T]=
\{F(T)\in\Zinfty[T]:F(z)=0\}$ for each complex~$z$ with $\Im z>0$. All these 
ideals are closed in this fiber, with the only exception of the zero ideal.

The closed fiber, i.e.\ prime ideals $\gp\subset\Zinfty[T]$ containing 
$\gm_\infty\cdot\Zinfty[T]$, seem to be more complicated to describe. 
Here are some examples of them:
\begin{itemize}
\item The only maximal ideal $\gm:=\Zinfty[T]-\{\pm1\}$; i.e.\ 
$\Zinfty[T]$ is {\em local}.
\item Prime ideal $\gq:=\Zinfty[T]-\{\pm T^k\}_{k\geq0}$.
\item Prime ideals $\gq_n:=\{F(T)\in\Zinfty[T]:|F(\zeta_n)|<1\}$, where 
$\zeta_n$ is a primitive root of unity of degree~$n$. In particular, we have 
$\gq_+:=\gq_1=\{F(T):|F(1)|<1\}$ and $\gq_-:=\gq_2$. Note that $\gm$ 
can be defined in the same way with $\zeta_0:=0$. For odd $n$ the 
complement of $\gq_n$ can be described as the set of polynomials 
$F(T)=\sum_k a_kT^k$ with $a_{qn+r}\neq0$ only for at most one value of 
remainder $0\leq r<n$, and 
with all $a_{qn+r}$ of same sign and of sum equal to $\pm1$. When $n$ 
is even, the description is slightly more complicated: we require 
$a_{qn/2+r}\neq0$ for at most one value of $0\leq r<n/2$, and all 
$(-1)^q a_{qn/2+r}$ of same sign and of sum equal to $\pm1$ for this value 
of~$r$.
\end{itemize}

Note that $\gq_n\subset\gq_{nn'}$ for all integer $n$ and odd $n'$.
This means that the minimal between these ideals are the $\gq_{2^k}$ for 
$k\geq0$; so they should correspond to irreducible components of the 
special fiber, if we admit that 
this list of prime ideals exhausts the whole of 
$\Spec\Zinfty[T]$. Note that the closed fiber has infinite Krull dimension 
in any case.

\nxsubpoint (Projective $\Zinfty$-schemes.) 
We can also define {\em graded\/} flat $\Zinfty$-algebras 
$A=\bigoplus_{n\geq0} A_n$, either in terms of such direct sum decompositions, 
understood as in~\ptref{ss:infdirsum}, subject to usual conditions 
$A_n\cdot A_m\subset A_{n+m}$ and $1\in A_0$, or in terms of sequences 
$(A_n)_{n\geq0}$ of flat $\Zinfty$-modules and $\Zinfty$-bilinear 
multiplication maps $\mu_{nm}:A_n\times A_m\to A_{n+m}$, subject to usual 
associativity and commutativity relations. Both descriptions are easily seen 
to be equivalent, and graded $\Zinfty$-flat $A$-modules $M$ can be described 
in a similar way. If $S\subset A$ is a multiplicative system, consisting 
only of homogeneous elements, we can define a natural grading on 
$S^{-1}A$ and $S^{-1}M$, and consider their degree zero part. Another 
possible construction of $(S^{-1}M)_0$ is by means of the 
inductive limit $\varinjlim_{s\in\cS}M_{[s]}$, where this time we take 
$M_{[s]}:=M_{\deg s}$ (cf.~\ptref{sp:loczflat}). If $f\in A$ is homogeneous 
of positive degree, we denote the degree zero part of $M_f$ by $M_{(f)}$, 
and similarly for $A_f$.

Now we can define the {\em projective spectrum\/} $\Proj A$ to be 
the set of homogeneous prime ideals $\gp\subset A$, not containing the 
ideal $A^+:=\bigoplus_{n>0}A_n$, with the topology induced by that of 
$\Spec A\supset\Proj A$. Next, the reasoning of EGA~II~2 shows us 
that the open subsets $D_+(f):=D(f)\cap\Proj A$ for homogeneous $f\in A^+$ 
form a base of topology on $\Proj A$, and that each $D_+(f)$ is 
homeomorphic to $\Spec A_{(f)}$. Modulo some compatibility issues, 
checked in the same way as in {\em loc.~cit.}, we obtain in this way a 
flat $\Zinfty$-scheme structure on $\Proj A$.

\nxsubpoint\label{sp:proj.vb.zinflat} (Projective spaces.) 
In particular, for any $\Zinfty$-lattice $V$ we consider the projective 
space $\bbP_{\Zinfty}(V):=\Proj S_{\Zinfty}(V)$, and similarly over 
$\barZinfty$. Note that the $\Zinfty$ or $\barZinfty$-structure 
of these projective spaces give rise to metrics on their complex-valued 
points in the same way as in~\ptref{sp:a1zinfty.infdata}, at least if 
we believe in the ``valuative criterion of properness'' in our situation. 
Locally this means considering the subquotient $\Zinfty$-structure 
on $\gm_P/\gm_P^2$, where $\gm_P$ is the maximal ideal corresponding 
to some closed point of the open fiber.

If we start with a hermitian $\barZinfty$-module $V$, we end up with 
a metric on the complex projective space $\bbP(V_{(\bbC)})\cong\bbP^{n-1}(\bbC)$, 
clearly equivariant under the action of $\Aut_{\barZinfty}(V)=U(n)$. There is 
only one such metric, up to multiplication by a positive constant, namely, 
the Fubini--Study metric on the projective space. 

We will use similar arguments later to show that any 
smooth projective variety 
with K\"ahler metric induced by the Fubini--Study metric on an ambient 
projective space admits a description in terms of flat $\Zinfty$-schemes.

\nxsubpoint
However, for several reasons we do not want to develop further the theory 
of (flat) $\Zinfty$-schemes right now. We are going to develop the 
theory of generalized rings first, and then we shall construct 
a theory of their spectra and generalized schemes, which will contain 
the theory sketched above as a special case. We would like to list 
some reasons for this.

\begin{itemize}
\item Up to now we have considered only flat $\Zinfty$-modules and 
flat $\Zinfty$-schemes. However, even if we are interested only in 
($\Zinfty$-flat) models of algebraic varieties over~$\bbR$, their closed 
subschemes, required for example to construct a reasonable intersection 
theory, in general need not be $\Zinfty$-flat. On the other hand, the 
theory of (abstract) $\Zinfty$-modules and $\Zinfty$-schemes doesn't seem 
to be more simple than that of modules and schemes over an arbitrary 
(hypoadditive) generalized ring, so we don't lose much if we do all 
constructions in the more general setting.
\item The reader may have observed that the localization theory of 
flat $\Zinfty$-algebras and modules seems to enjoy almost all 
the usual properties known from commutative algebra. However, this 
is not the case with the other important construction of commutative algebra, 
namely, the theory of quotient modules and rings. Even if we consider the 
category of all (abstract) $\Zinfty$-modules, we get $Q:=\Zinfty/\frac12\Zinfty
\cong\Zinfty/\gm_\infty$, so the kernel of $\Zinfty\to Q$ is equal to 
$\gm_\infty$, and not to $\frac12\Zinfty$ (cf.~\ptref{sp:def.finf}); 
another similar phenomenon is 
given by $\bbR/\Zinfty$, which is equal to zero even in $\ZinfMod$. 
However, if we consider strict quotients of a $\Zinfty$-algebra in the 
larger category of generalized rings, the situation becomes more 
natural. In this way we are led to consider closed (generalized) subschemes 
of a flat $\Zinfty$-scheme that are not flat themselves.
\item If we develop separately the theory of $\Zinfty$-schemes and 
the usual theory of (Grothendieck) schemes, we will still be forced to 
describe objects (say, schemes) over $\CompZ$ as corresponding objects 
over $\Spec\bbZ$ and $\Spec\Zinfty$, the pullbacks of which to $\Spec\bbR$ 
agree (cf.~\ptref{p:firstcompz}). On the other hand, the theory of 
generalized rings and schemes enables us to construct $\CompZ$ and, say, 
models $\sX\to\CompZ$ as generalized schemes, since 
both the theory of commutative rings and the theory of $\Zinfty$-algebras 
are special cases of this general theory.
\item Furthermore, we will be able to consider $\sX$ as a generalized 
scheme over $\Fpm:=\Gamma(\CompZ,\sO_{\CompZ})=\bbZ\cap\Zinfty:=
\bbZ\times_\bbR\Zinfty$, where this intersection (i.e.\ fibered product) 
is computed in the category of generalized rings. Among other things, 
this will lead to a reasonable construction of intersection theory on $\sX$, 
with $\Fpm$ for the base ring.
\end{itemize}

\nxpointtoc{Abstract $\Zinfty$-modules}
Now we are going to use our knowledge of the category $\ZinfFlat$ of flat 
(or torsion-free) $\Zinfty$-modules to construct the category 
$\ZinfMod$ of all (abstract) $\Zinfty$-modules. We shall also need 
for this the definition and some basic properties of monads; they can 
be found either in the next chapter or in \cite{MacLane}.

\nxsubpoint
Recall that we have constructed a pair of adjoint functors: 
the {\em forgetful functor} $\Gamma_{\Zinfty}:\ZinfFlat\to\catSets$, 
$A\mapsto A_{\Zinfty}\cong\Hom_{\Zinfty}(\Zinfty,A)$, 
and its left adjoint $L_{\Zinfty}:\catSets\to\ZinfFlat$, 
$S\mapsto\Zinfty^{(S)}$ (cf.~\ptref{sp:freezinfmod}). Denote by 
$\xi:\Id_\catSets\to\Gamma_{\Zinfty}L_{\Zinfty}$ and 
$\eta:L_{\Zinfty}\Gamma_{\Zinfty}\to\Id_{\ZinfFlat}$ the  
natural transformations defining adjointness between $L_{\Zinfty}$ 
and~$\Gamma_{\Zinfty}$.

By general theory of monads this pair of adjoint functors gives us a 
monad $\Sigmainf=\Sigma_{\Zinfty}=(\Sigmainf,\mu,\epsilon)$ 
on $\catSets$, where $\Sigmainf:=\Gamma_{\Zinfty}L_{\Zinfty}$ 
is an endofunctor on $\catSets$, $\epsilon:=\xi:\Id_\catSets\to\Sigmainf$ 
is the unit, and $\mu:=\Gamma_{\Zinfty}*\eta*L_{\Zinfty}:\Sigmainf^2\to
\Sigmainf$ is the multiplication of this monad.

\nxsubpoint\label{sp:prel.def.cat.alg}
Once we have a monad $\Sigmainf$ on the category of sets, we can 
consider the category of $\Sigmainf$-algebras or $\Sigmainf$-sets, 
usually denoted by $\catSets^{\Sigmainf}$. By definition, its objects 
are sets $X$, equipped with a $\Sigmainf$-structure $\alpha:
\Sigmainf(X)\to X$, i.e.\ a map of sets subject to conditions 
$\alpha\circ\epsilon_X=\id_X$ and $\alpha\circ\mu_X=\alpha\circ\Sigmainf(\alpha):\Sigmainf^2(X)\to X$. Morphisms $f:(X,\alpha_X)\to(Y,\alpha_Y)$ 
are maps of sets $f:X\to Y$, compatible with given $\Sigmainf$-structures, 
i.e.\ $f\circ\alpha_X=\alpha_Y\circ\Sigmainf(f)$.

Clearly, we have a forgetful functor $\Gamma':\catSets^{\Sigmainf}\to
\catSets$, $(X,\alpha)\mapsto X$. It has a left adjoint $L'$, given by 
$L'(S):=(\Sigmainf(S),\mu_S)$, such that $\Gamma'L'$ is still equal to 
$\Sigmainf$.

Finally, $\Gamma_{\Zinfty}:\ZinfFlat\to\catSets$ canonically factorizes 
through $\Gamma'$, 
yielding a functor $I:\ZinfFlat\to\catSets^{\Sigmainf}$, 
such that $\Gamma'I=\Gamma_{\Zinfty}$; recall that 
$I:A\mapsto\bigl(\Gamma_{\Zinfty}(A),\Gamma_{\Zinfty}(\eta_A)\bigr)$. Now, 
by definition $\Gamma$ would be a {\em monadic functor\/} if 
$I$ had been an equivalence of categories; in this case 
we would be able to identify $\ZinfFlat$ with $\catSets^{\Sigmainf}$, 
$\Gamma_{\Zinfty}$ with $\Gamma'$, and $L_{\Zinfty}$ with $L'$.

\nxsubpoint\label{sp:notallzinfmodflat}
However, this is not the case: {\em $I:\ZinfFlat\to\catSets^{\Sigmainf}$ 
is fully faithful, but not an equivalence of categories.} 
We will prove this statement later; let's discuss its consequences first.

\begin{DefD}\label{def:zinfmod}
We define the category $\ZinfMod$ of (abstract) $\Zinfty$-modu\-les to be 
the category $\catSets^{\Sigmainf}$ of sets equipped with a 
$\Sigmainf$-structure. We identify $\ZinfFlat$ with a full subcategory 
of $\ZinfMod$ by means of functor~$I$ constructed above. A 
{\em $\Zinfty$-structure $\alpha$} on some set $X$ is by definition 
the same thing as a $\Sigmainf$-structure, i.e.\ a map 
$\alpha:\Sigmainf(X)\to X$, subject to two conditions recalled 
in~\ptref{sp:prel.def.cat.alg}. A map $f:X\to Y$ between two 
$\Zinfty$-modules is {\em $\Zinfty$-linear}, or a 
{\em $\Zinfty$-homomorphism}, if it respects the $\Zinfty$-structures, 
i.e.\ if it defines a morphism in $\ZinfMod$.
\end{DefD}

\nxsubpoint 
Note that this definition is more ``algebraic'' than  
\ptref{def:zinflat} and \ptref{def:zinfflat}: instead of considering 
pairs $(A_{\Zinfty},A_\bbR)$, consisting of a symmetric convex body in 
an ambient real space, we consider sets $A=A_{\Zinfty}$, equipped with 
some $\Zinfty$-structure, i.e.\ $\Sigma_{\infty}$-structure.

\nxsubpoint
The reason for this definition is the following. If we start from any 
commutative ring~$K$ (say, $K=\bbZ_p$), and consider the category 
$\catMod{K}$ of $K$-modules together with the forgetful functor 
$\Gamma_K:\catMod{K}\to\catSets$, then this functor turns out to be 
monadic, i.e.\ it has a left adjoint $L_K:S\mapsto K^{(S)}$, hence 
it defines a monad $\Sigma_K:=\Gamma_KL_K$ on $\catSets$, and induced 
functor $I_K:\catMod{K}\to\catSets^{\Sigma_K}$ is an equivalence.

However, if we start from the category $\ZpFlat$ of flat $\bbZ_p$-modules, 
we will end up with the same monad $\Sigma_{\bbZ_p}$, but induced functor 
$\ZpFlat\to\catSets^{\Sigma_{\bbZ_p}}\cong\catMod{\bbZ_p}$ won't be 
an equivalence, but just a fully faithful functor. In this way we can 
reconstruct first $\Sigma_{\bbZ_p}$, and then $\catMod{\bbZ_p}\cong
\catSets^{\Sigma_{\bbZ_p}}$, starting from the category $\ZpFlat$ and 
the forgetful functor $\ZpFlat\to\catSets$. This is exactly what we have 
done in~\ptref{def:zinfmod} for~$\Zinfty$.

\nxsubpoint
Recall the explicit description of 
$\Sigmainf(S)=\Gamma_{\Zinfty}(\Zinfty^{(S)})$ given 
in~\ptref{sp:explfreezinfmod}: $\Sigmainf(S)\subset\bbR^{(S)}$ consists 
of all formal octahedral combinations $\sum_{s\in S}\lambda_s\{s\}$ 
of elements of $S$, where almost all $\lambda_s\in\bbR$ are equal to zero 
and $\sum_s|\lambda_s|\leq1$. Hence a $\Sigmainf$-structure $\alpha$ 
on some set~$S$ is a map $\alpha:\Sigmainf(S)\to S$, i.e.\ a way to 
evaluate formal octahedral combinations of elements of~$S$. We denote 
$\alpha\bigl(\sum_i\lambda_i\{s_i\}\bigr)$ by 
$\bigl<\sum_i\lambda_i s_i\bigr>_\alpha$,
or even by $\sum_i\lambda_is_i$, when no confusion can arise.
Here $s_i\in S$, almost all $\lambda_i\in\bbR$ are zero, 
and $\sum_i|\lambda_i|\leq1$, as usual.

If we start from an object $A=(A_{\Zinfty},A_\bbR)$ of $\ZinfFlat$ and 
construct a $\Sigmainf$-structure $\alpha_A$ on $A_{\Zinfty}$ as explained 
in~\ptref{sp:prel.def.cat.alg} (this corresponds to considering $I(A)$), 
then $\alpha_A=\Gamma(\eta_A):\Sigmainf(A_{\Zinfty})\to A_{\Zinfty}$, 
and we see immediately that the image in $A_\bbR$ 
of $\langle\sum_i\lambda_i s_i\rangle_{\alpha_A}\in A_{\Zinfty}$ 
coincides with the corresponding linear combination computed in $A_\bbR$; 
this explains our notation.

\nxsubpoint\label{sp:mu.zinfty}
In particular, we can apply this observation to $A:=\Zinfty^{(S)}$. 
In this case it follows from definitions that $\alpha_A=\mu_S$, so 
we obtain an explicit description of $\mu_S:\Sigmainf^2(S)\to\Sigmainf(S)$. 
Namely, $\mu_S\bigl(\sum_i\lambda_i\{\sum_j\mu_{ij}\{s_j\}\}\bigr)$ 
equals $\sum_{i,j}\lambda_i\mu_{ij}\{s_j\}$. An explicit description of 
$\epsilon_S:S\to\Sigmainf(S)$ is even easier to obtain: 
$\epsilon_S$ maps any $s\in S$ into corresponding basis element 
$\{s\}\in\Sigmainf(S)\subset\bbR^{(S)}$.

Now, our conditions for a map $\alpha:\Sigma(S)\to S$ to be a 
$\Sigmainf$-structure are $\alpha\circ\epsilon_S=1_S$ and 
$\alpha\circ\Sigmainf(\alpha)=\alpha\circ\mu_S$. The first condition 
translates into $\alpha(\{s\})=s$, i.e.\ $\langle s\rangle_\alpha=s$, for 
any $s\in S$; this requirement seems to be quite natural. 

The second condition translates into 
\begin{equation}
\Bigl<\sum_{i=1}^n \lambda_i\bigl<\sum_{j=1}^m \mu_{ij} s_j\bigr>_\alpha
\Bigr>_\alpha =
\Bigl<\sum_{j=1}^m \bigl(\sum_{i=1}^n \lambda_i\mu_{ij}\bigr)s_j\Bigr>_\alpha
\quad,
\end{equation}
where $s_j\in S$, $\sum_i|\lambda_i|\leq 1$ and
$\sum_j|\mu_{ij}|\leq 1$ for any $i$. If we remove the angular brackets, 
the equation we obtain looks like the usual distributivity relation, 
written for arbitrary {\em octahedral\/} combination of octahedral 
combinations of some $s_j\in S$.

\nxsubpoint
Since any formal octahedral combination involves only finitely many elements 
of a set~$X$, we see that a $\Sigmainf$-structure $\alpha$ on $X$ is 
completely determined by the family of maps 
$\alpha_n:\Sigma(n)\times X^n\to X$, $n\geq0$, 
given by $\alpha_n(\lambda_1,\ldots,\lambda_n; x_1,\ldots,x_n):=
\langle\lambda_1x_1+\cdots+\lambda_nx_n\rangle_\alpha=
\alpha(\lambda_1\{x_1\}+\cdots+\lambda_n\{x_n\})$. Here 
$\Sigma(n)=\{(\lambda_1,\ldots,\lambda_n):\sum_i|\lambda_i|\leq1\}$ is 
the standard octahedron in~$\bbR^n$; cf.~\ptref{sp:explfreezinfmod}.

Of course, these maps $\alpha_n$, $n\geq0$, have to satisfy some compatibility 
relations (e.g.\ invariance under all permutations of arguments~$\sigma$) 
in order to define together some $\alpha:\Sigmainf(X)\to X$:
\begin{gather}
\alpha_n(\lambda_1,\ldots,\lambda_n; x_1,\ldots,x_n)=
 \alpha_n(\lambda_{\sigma(1)},\ldots,\lambda_{\sigma(n)}; 
  x_{\sigma(1)},\ldots,x_{\sigma(n)})\text{, $\forall\sigma$};\\
\label{eq:extby0}
\alpha_{n+1}(\lambda_1,\ldots,\lambda_n,0; x_1, \ldots, x_n, x_{n+1})=
 \alpha_n(\lambda_1,\ldots,\lambda_n; x_1,\ldots,x_n);\\
\alpha_{n+1}(\lambda_1,\ldots,\lambda_n,\lambda_{n+1}; x_1,\ldots, x_n, x_n)=
 \alpha_n(\lambda_1,\ldots,\lambda_n+\lambda_{n+1};x_1,\ldots,x_n)
\end{gather}

Resulting $\alpha$ is a $\Sigmainf$-structure iff the unit and associativity 
relations are fulfilled:
\begin{gather}
\alpha_1(1;x)=x\\
\begin{split}
&\alpha_n\bigl(\lambda_1,\ldots,\lambda_n;\,
 \alpha_m(\mu_{11},\ldots,\mu_{1m}; x_1,\ldots, x_m),\ldots,
 \alpha_m(\mu_{n1},\ldots,\mu_{nm}; x_1,\ldots)\bigr)\\ 
& =
\alpha_m\bigl(\sum_{i=1}^n\lambda_i\mu_{i1},\ldots, 
 \sum_{i=1}^n\lambda_i\mu_{im};\, x_1,\ldots, x_m\bigr)
\end{split}
\end{gather}

We have thus obtained a description of $\Zinfty$-modules in terms of 
a countable family of operations $\alpha_n$.

\nxsubpoint
Note that {\em all $\alpha_n$, hence also $\alpha$, are completely determined 
by $\alpha_2:\Sigmainf(2)\times X^2\to X$}. Indeed, for $\alpha_0$ and 
$\alpha_1$ we have $\alpha_1(\lambda; x)=\alpha_2(\lambda,0;x,x)$ and 
$\alpha_0=\alpha_2(0,0;x,x)$ for any $x\in X$ (note that $X$ has to be 
non-empty, otherwise we won't be able to retrieve its zero element $\alpha_0$).
For $\alpha_n$, $n\geq3$, we prove this statement by induction, using 
identity
\begin{multline}
\alpha_n(\lambda_1,\ldots,\lambda_{n-2},\lambda_{n-1},\lambda_n;\,
x_1,\ldots,x_{n-2},$ $x_{n-1},x_n)=\\
\alpha_{n-1}(\lambda_1,\ldots,\lambda_{n-2},\mu;\,
x_1,\ldots, x_{n-2},\alpha_2(\lambda_{n-1}/\mu, 
\lambda_n/\mu;\, x_{n-1}, x_n))
\end{multline}
when $\mu:=|\lambda_{n-1}|+|\lambda_n|\neq0$; 
if it is zero, we use~\eqref{eq:extby0} instead.

We see that we might describe a $\Zinfty$-module $X$ as a triple 
$(X,\alpha_0,\alpha_2)$, where $X$ is a set, $\alpha_0\in X$ is its 
marked element (``zero''), and $\alpha_2:\Sigmainf(2)\times X^2\to X$ is 
the binary octahedral combination evaluation map discussed above. These 
data have to satisfy several axioms, which do not seem to be especially 
enlightening, so we don't list them here. The most interesting among them 
is the associativity/distributivity axiom for~$\alpha_2$:
\begin{equation}
\alpha_2(\lambda'\nu',\mu;\, x, \alpha_2(\nu,\rho;\, y,z))=
\alpha_2(\lambda',\mu\rho;\, \alpha_2(\nu',\rho';\, x,y), z),
\end{equation}
for any $x,y,z\in X$, whenever 
$\mu\nu=\lambda'\rho'$, $|\lambda'\nu'|+|\mu|\leq1$, 
$|\lambda'|+|\mu\rho|\leq1$, $|\nu|+|\rho|\leq1$ and $|\nu'|+|\rho'|\leq1$.

In this form the definition of (abstract) $\Zinfty$-modules is very similar 
to that of abstract convex sets, which in fact can be obtained in the 
same way starting from the monad $\bbDelta$ that maps 
any set $S$ into the simplex with vertices $\{s\}$, i.e.\ the convex hull 
of all basis elements $\{s\}$ of $\bbR^{(S)}$.

\nxsubpoint\label{sp:zinfmod.specel}
We see that any $\Zinfty$-module $X$ has a marked element, namely, 
$\alpha_0:=\alpha_X(0)$; we denote it by $0_X$ or $0$, and call it 
the {\em zero element of~$X$.} We have also a map $\alpha_1:\Zinfty\times X\to 
X$, which defines an action of $\Zinfty$, considered as a multiplicative 
monoid, on~$X$, i.e.\ 
$\alpha_1(\lambda,\alpha_1(\mu,x))=\alpha_1(\lambda\mu,x)$ for any 
$x\in X$ and $|\lambda|$, $|\mu|\leq1$. Of course, $\alpha_1(\lambda,x)$ 
is usually written as $\lambda\cdot x$ or $\lambda x$, so the above formula 
can be written simply as $\lambda(\mu x)=(\lambda\mu)x$. Note that 
$0_X$ is a fixed point for this action, i.e.\ $\lambda\cdot 0_X=0_X$ for any 
$|\lambda|\leq1$, and $0\cdot x=0_X$ for any $x\in X$. However, this 
structure of a set with marked point, equipped by an action of monoid 
$\Zinfty=[-1,1]$, doesn't determine uniquely $\alpha_2$, i.e.\ the structure 
of a $\Zinfty$-module on~$X$.

\nxsubpoint\label{sp:zinfmod.listprop}
We see that a $\Zinfty$-module $X$ can be described as an algebraic system 
$X=(X,\alpha_0,\alpha_2)$, consisting of a set $X$, its element $\alpha_0$ 
and a map $\alpha_2:\Sigmainf(2)\times X^2\to X$, subject to a finite set 
of axioms. Therefore, the category of $\Zinfty$-modules has all the 
usual properties of categories defined by algebraic systems. We list 
some of these properties below without proof; most of them can be easily 
checked directly, and all of them will be proved later in~\ptref{p:algmon.mod} 
for the categories of modules over arbitrary generalized rings.
\begin{itemize}
\item Arbitrary projective limits exist in $\ZinfMod$; they can be essentially 
computed in $\catSets$ (i.e.\ $\Gamma_{\Zinfty}:\ZinfMod\to\catSets$ commutes 
with arbitrary projective limits).
\item Arbitrary inductive limits exist in $\ZinfMod$ (this is the only 
complicated statement in this list).
\item Filtered inductive limits exist in $\ZinfMod$; they can be computed in 
$\catSets$, i.e.\ $\Gamma_{\Zinfty}$ commutes with filtered inductive limits.
\item A $\Zinfty$-linear map $f:M\to N$ is a monomorphism in $\ZinfMod$ 
iff it is injective as a map of sets (i.e.\ iff $\Gamma_{\Zinfty}(f)$ is 
injective).
\item A $\Zinfty$-linear map $f:M\to N$ is a strict epimorphism iff it is 
surjective as a map of sets.
\item Given an injective map of sets $f:N\to M$ and a $\Zinfty$-module 
structure on $M$, there is at most one $\Zinfty$-module structure 
(called the {\em induced structure}) on~$N$, 
compatible with $f$. Therefore, subobjects (i.e.\ $\Zinfty$-submodules) 
of $M$ are in one-to-one correspondence with those subsets $N$ of $M$, which 
admit an induced $\Zinfty$-module structure.
\item Similarly, given a surjective map of sets $f:M\to Q$ and a 
$\Zinfty$-module structure on $M$, there is at most one $\Zinfty$-module 
structure on $Q$, compatible with $f$. Strict quotients of $M$ are 
in one-to-one correspondence with certain quotient sets $Q$ of $M$.
\item Any $\Zinfty$-homomorphism $f:M\to N$ decomposes into a 
strict epimorphism $p:M\to I$, followed by a monomorphism $i:I\to N$. 
As a set, $I$ coincides with $f(M)$; its $\Zinfty$-module structure can 
be described either as one induced from $N$ or as one induced from $M$. 
We denote $I$ by $f(M)$ or $\Im f$ and say that it is the {\em image} of~$f$.
\end{itemize}

\nxsubpoint\label{sp:def.finf}
Consider the $\Zinfty$-module $Q=\Finfty$, constructed as 
follows. Its underlying set consists of three elements, denoted 
$0$, $+1$ and $-1$. The $\Sigmainf$-structure $\alpha=\alpha_Q$ on
$Q=\{0,+1,-1\}$ is given by
\begin{equation}
\alpha_Q\bigl(\lambda_{-1}\{-1\}+\lambda_0\{0\}+\lambda_1\{+1\}\bigr) =
 \begin{cases}
 +1&\text{if $\lambda_1-\lambda_{-1}=1$}\\
 -1&\text{if $\lambda_1-\lambda_{-1}=-1$}\\
  0&\text{if $|\lambda_1-\lambda_{-1}|<1$}
\end{cases}\end{equation}
Here $|\lambda_{-1}|+|\lambda_0|+|\lambda_1|\leq1$, so in any case 
$|\lambda_1-\lambda_{-1}|\leq1$, hence the cases listed above exhaust all 
possibilities. We have to check that this $\alpha_Q$ satisfies 
the axioms for a $\Sigmainf$-structure; for this we consider the 
surjective map $\phi:\Zinfty\to Q$, which maps $\pm1\in\Zinfty=[-1,1]$ into 
$\pm1\in Q$, and all other elements of $\Zinfty$ into $0$. It is easy to check
that $\phi\circ\alpha_{\Zinfty}=\alpha_Q\circ\Sigmainf(\phi)$; since 
$\phi$ is surjective, this implies that $\alpha_Q$ is indeed a 
$\Zinfty$-structure on $Q$, and that $\phi:\Zinfty\to Q$ is a surjective 
$\Zinfty$-linear map, i.e.\ a strict epimorphism in $\ZinfMod$.

We see that $Q$ is obtained from $\Zinfty$ by identifying all elements of 
$\gm_\infty$ with zero; this explains our alternative notation 
$\Finfty$ and $\Zinfty/\gm_\infty$ for~$Q$.

Note that for any non-zero object $A=(A_{\Zinfty},A_\bbR)$ of $\ZinfFlat$ 
the underlying set of corresponding object $I(A)$ of $\ZinfMod$ 
is equal to $A_{\Zinfty}$, hence it is infinite. Since $Q\neq0$ and 
it is finite, we see that it cannot be isomorphic to any $I(A)$, i.e.\ 
it doesn't lie in the essential image of $I:\ZinfFlat\to\ZinfMod$, 
hence $I$ cannot be an equivalence of categories.

\nxsubpoint\label{sp:ff.zinfflat.mod}
Clearly, $I:\ZinfFlat\to\ZinfMod$ is faithful, since $\Gamma'\circ I=\Gamma
:=\Gamma_{\Zinfty}$ is faithful; 
in order to prove~\ptref{sp:notallzinfmodflat} we have to check 
that $I$ is fully faithful, i.e.\ that {\em for any two flat $\Zinfty$-modules 
$A$ and $B$, any map of underlying sets $f:\Gamma(A)\to\Gamma(B)$ 
compatible with induced $\Sigmainf$-structures 
actually defines a morphism $f'=(f'_{\Zinfty},f'_\bbR):A\to B$ in 
$\ZinfFlat$ with $f'_{\Zinfty}=f$ (cf.~\ptref{def:zinfflat}).}

Recall that the $\Sigmainf$-structure on $\Gamma(A)$ is given by 
$\alpha_A:=\Gamma(\eta_A):\Gamma(\Zinfty^{(A)})=\Sigmainf\Gamma(A)\to
\Gamma(A)$, and similarly for~$B$; so the compatibility of $f$ with 
these $\Sigmainf$-structures means commutativity of the following square:
\begin{equation}\label{eq:ffproof1}
\xymatrix{
\Gamma(\Zinfty^{(A)})\ar[rr]^{\Sigmainf(f)}
\ar[d]^{\Gamma(\eta_A)}&&\Gamma(\Zinfty^{(B)})\ar[d]^{\Gamma(\eta_B)}\\
\Gamma(A)\ar[rr]^{f}&&\Gamma(B)}
\end{equation}
Now consider the following diagram with right exact rows:
\begin{equation}
\xymatrix@+1pc{
R_A\ar[r]<2pt>\ar[r]<-3pt>\ar[d]_{R_f}&
\Zinfty^{(A)}\ar[d]_{\Zinfty^{(f)}}\ar[r]^{\eta_A}&
A\ar@{-->}[d]_{f'}\\
R_B\ar[r]<2pt>\ar[r]<-3pt>&\Zinfty^{(B)}\ar[r]^{\eta_B}&B}
\end{equation}
Note that $\eta_A$ is a strict epimorphism in $\ZinfFlat$, 
hence it is cokernel of its kernel pair $R_A\rightrightarrows 
\Zinfty^{(A)}$, and similarly for $\eta_B$; therefore, once we 
construct the two left vertical arrows in $\ZinfFlat$, we'll obtain 
immediately an $f':A\to B$ in $\ZinfFlat$ completing this diagram; 
since $\Gamma(\eta_A)$ is surjective and~\eqref{eq:ffproof1} is commutative, 
necessarily $\Gamma(f')=f$ as required. Now, the middle vertical 
arrow $L_{\Zinfty}(f)=\Zinfty^{(f)}$ exists for any map of sets~$f$; to show 
existence of $R_f$ we observe that $\Zinfty^{(f)}\times\Zinfty^{(f)}$ 
maps the underlying set 
of strict subobject $R_A$ of $\Zinfty^{(A)}\times\Zinfty^{(A)}$ into 
the underlying set of $R_B\subset\Zinfty^{(B)}\times\Zinfty^{(B)}$
because of commutativity of~\eqref{eq:ffproof1}, hence 
it induces some $R_f:R_A\to R_B$ with required properties.

This finishes the proof of~\ptref{sp:notallzinfmodflat}.

\nxsubpoint
Henceforth we shall identify $\ZinfFlat$ with a full subcategory of $\ZinfMod$ 
by means of the functor~$I$ constructed in~\ptref{sp:prel.def.cat.alg}. 
Recall that the category of $\bbR$-vector spaces $\RVect$ has been identified 
in~\ptref{sp:rvectinzflat} with a full subcategory of $\ZinfFlat$, hence 
it is now identified with a full subcategory of $\ZinfMod$. In other words, 
we have a canonical $\Zinfty$-structure $\alpha_V$ on any 
$\bbR$-vector space~$V$; it is easy to see that $\alpha_V:\Sigmainf(V)\to V$ 
simply evaluates any given formal octahedral combination of elements of 
$V$ in the natural way by means of the $\bbR$-vector space structure on~$V$.

Let us denote the functor just constructed by $\rho_*:\RVect\to\ZinfMod$. 
It is not very hard to prove that it admits both a left adjoint $\rho^*$ 
and a right adjoint $\rho^!$, similarly to what we had before for flat 
$\Zinfty$-modules. We will see later that formulas 
$\rho^*M=\bbR\otimes_{\Zinfty}M$ and $\rho^!M=\iHom_{\Zinfty}(\bbR,M)$ 
are still valid in this context, so we can denote $\rho^*M$ also 
by~$M_{(\bbR)}$ or even~$M_\bbR$.

\nxsubpoint\label{sp:constr.tensR}
Existence of $M_{(\bbR)}$ and of a $\Zinfty$-linear map $i_M:M\to M_{(\bbR)}$, 
universal with respect to $\Zinfty$-linear maps from $M$ into 
$\bbR$-vector spaces, can be shown directly in several different ways. 
For example, we can construct $M_{(\bbR)}$ as the quotient of free 
$\bbR$-vector space $\bbR^{(M)}$ by the $\bbR$-vector subspace generated 
by all elements of the form $\{\lambda_1x_1+\cdots+\lambda_nx_n\}-
\lambda_1\{x_1\}-\cdots-\lambda_n\{x_n\}$, for all $n\geq0$ and all octahedral 
combinations of $n$ elements of~$M$. 

Second construction: we can characterize 
$\bbR$-vector spaces as those $\Zinfty$-modules, on which all elements 
of $S:=\Zinfty-\{0\}$ act bijectively; then the required universal property 
of $M_{(\bbR)}$ is that of $S^{-1}M$, so we have to compute $S^{-1}M$, 
and this can be done by means of filtered inductive limits 
in the same way as in~\ptref{sp:loczflat}. In particular, this means that 
all elements of $M_{(\bbR)}$ can be written in form $x/s$ for some 
$x\in M$ and $s\in\Zinfty$, $s\neq0$, and $x/s=y/t$ iff there is some 
$u\neq0$, such that $utx=usy$.

Finally, we can develop first the theory of tensor products in $\ZinfMod$, 
and construct $M_{(\bbR)}$ as $\bbR\otimes_{\Zinfty}M$.

\begin{DefD}
We say that a $\Zinfty$-module $M$ is {\em a torsion-free $\Zinfty$-module} 
if the canonical map $i_M:M\to M_{(\bbR)}$ is injective, i.e.\ a 
monomorphism in $\ZinfMod$. We say that $M$ is {\em a (pure) torsion 
$\Zinfty$-module} if $M_{(\bbR)}=0$. Finally, we denote by $M_{tf}$ the 
image of $M$ in $M_{(\bbR)}$; it is the largest torsion-free quotient of~$M$.
\end{DefD}

Note that by the universal property of $M_{(\bbR)}$, we have that 
$M$ is torsion-free iff there exists an embedding (=monomorphism) of 
$M$ into an $\bbR$-vector space. This shows that $M_{tf}\subset M_{(\bbR)}$ 
is indeed torsion-free. Any map $f:M\to N$ from $M$ to a torsion-free $N$ 
induces an $\bbR$-linear map $f_{(\bbR)}:M_{(\bbR)}\to N_{(\bbR)}$; since 
$i_N:N\to N_{(\bbR)}$ is injective, and $f_{(\bbR)}\circ i_M=i_N\circ f$,
this implies immediately that $f$ factorizes through $i_M(M)=M_{tf}$, as 
stated above. In other words, $M\mapsto M_{tf}$ is a left adjoint to 
the embedding of the full subcategory of torsion-free $\Zinfty$-modules 
into~$\ZinfMod$.

\nxsubpoint
Let us show that the torsion-free $\Zinfty$-modules are exactly those 
which lie in the essential image of $I:\ZinfFlat\to\ZinfMod$, i.e.\ 
which are isomorphic to some $I(A)$. This will characterize internally 
$\ZinfFlat$ as a subcategory of $\ZinfMod$, and it will show that $
M\mapsto M_{tf}$ is a left adjoint to~$I$.

Clearly, any $\Zinfty$-module of form $I(A)=A_{\Zinfty}$ for some 
object $A=(A_{\Zinfty},A_{\bbR})$ embeds into some $\bbR$-vector space, namely,
$A_\bbR$, so it must be torsion-free. Conversely, if $M$ is torsion-free, 
we put $A_\bbR:=M_{(\bbR)}$ and $A_{\Zinfty}:=M_{tf}=i_M(M)\cong M$. 
Then $A_{\Zinfty}\subset A_\bbR$ is stable under all octahedral combinations 
of its elements, hence it is symmetric and convex; it also generates the 
whole of $A_\bbR=M_{(\bbR)}$, as it easily follows either from the universal 
property or from any of the three constructions of $M_{(\bbR)}$ given
in~\ptref{sp:constr.tensR}, hence $A\in\Ob\ZinfFlat$ and $M\cong I(A)$.

\nxsubpoint (Torsion $\Zinfty$-modules.)
By definition, $M$ is a torsion $\Zinfty$-module iff 
$M_{(\bbR)}=S^{-1}M=0$, where $S=\Zinfty-\{0\}$. Our description of 
localizations shows immediately that this condition is equivalent to 
the following: For any $x\in M$ there is some $\lambda\in S$
(i.e.\ a non-zero $\lambda\in\Zinfty$), such that $\lambda x=0$ in~$M$. 
For example, the $\Zinfty$-module $\Finfty$ considered 
in~\ptref{sp:def.finf} is torsion, since $(1/2)x=0$ for any $x\in\Finfty$.

\nxsubpoint (Limits in $\ZinfFlat$ and $\ZinfMod$.)
Since $I:\ZinfFlat\to\ZinfMod$ has a left adjoint $M\mapsto M_{tf}$, 
it commutes with arbitrary projective limits, i.e.\ arbitrary projective 
limits of torsion-free $\Zinfty$-modules can be still computed 
as in~\ptref{p:projlimzflat}. This is in general not true for 
arbitrary inductive limits: we have only the formula 
$I(\injlim_\alpha M_\alpha)=(\injlim_\alpha I(M_\alpha))_{tf}$, which 
follows immediately from the fact that $I$ is fully faithful and admits a 
left adjoint $M\mapsto M_{tf}$.

However, this is still true for filtered inductive limits, since in
both categories they can be essentially computed in $\catSets$, and
for arbitrary direct sums as well.  In the latter case we are to check
that arbitrary direct sums of torsion-free modules are still
torsion-free. For this we first write a direct sum over an arbitrary
index set~$I$ as a filtered inductive limit of direct sums over all
finite subsets $J\subset I$, thus reducing to the case of a finite
direct sum, and then by induction it suffices to check the statement
for a direct sum $M\oplus N$.

So let $M$ and $N$ be torsion-free $\Zinfty$-modules, $P:=M\oplus N$
be their direct sum in $\catMod\Zinfty$. Notice that $i_1:M\to M\oplus
N$ is injective, since it admits a left inverse $\pi_1:M\oplus N\to
M\oplus 0\cong M$, and similarly for $i_2:N\to M\oplus N$, so we can
identify $M$ and $N$ with subsets of $M\oplus N$. Let us check that
$P$ is torsion-free, i.e.\ that $P\to S^{-1}P$ is injective. Let $z$,
$z'\in P$ be two elements of~$P$, such that $z/1=z'/1$ in $S^{-1}P$,
i.e.\ $\nu z=\nu z'$ for some $0<\nu\leq 1$. We'll check
in~\ptref{sp:elem.dirsum} that elements $z$, $z'\in P=M\oplus N$ can
be always represented as octahedral combinations $z=\lambda x+\mu y$,
$z'=\lambda'x'+\mu'y'$ with $x$, $x'\in M$, $y$, $y'\in N$. Applying
$\pi_1:P\to M$ to $\nu z=\nu z'$, we obtain $\nu\lambda
x=\nu\lambda'x'$, hence $\lambda x= \lambda'x'$, $M$ being
torsion-free. Furthermore, if $|\lambda|\geq|\lambda'|$, we get
$\lambda\cdot(\lambda'/\lambda)x'=\lambda x$, whence
$x=(\lambda'/\lambda)x'$, unless both $\lambda=\lambda'=0$.
Substituting this into $z=\lambda x+\mu y$, we get $z=\lambda'x'+\mu
y$, i.e.\ we can assume $\lambda=\lambda'$, $x=x'$. Case
$|\lambda'|\geq|\lambda|$ is dealt with similarly, and if both
$\lambda=\lambda'=0$, then $z$ and $z'\in N$, hence $z=z'$, $N$ being
torsion-free, so this case can be excluded. Next, we can consider
cases $|\mu|\geq|\mu'|$, $|\mu'|\geq|\mu|$ and $\mu=\mu'=0$, and
deduce from $N$ being torsion-free that we can assume $\mu=\mu'$,
$y=y'$ as well. We see that $z=\lambda x+\mu y=z'$, so $P\to S^{-1}P$
is indeed injective and $P=M\oplus N$ is torsion-free, q.e.d.

\nxsubpoint (Linear algebra, $\otimes$-structure etc.)
Given any three $\Zinfty$-modules $M$, $N$ and $P$, we say that a map 
$\Phi: M\times N\to P$ is {\em $\Zinfty$-bilinear} if 
for any $x\in M$ the map $s_\Phi(x):y\mapsto\Phi(x,y)$ is a $\Zinfty$-linear 
map $N\to P$, and, similarly, for any $y\in N$ the map 
$d_\Phi(y):x\mapsto\Phi(x,y)$ is $\Zinfty$-linear.

Proceeding as in~\ptref{p:zinfflat.polylin}, we construct tensor products 
and inner Homs in $\ZinfMod$, obtaining an ACU $\otimes$-structure on this 
category. This leads to a natural definition of a $\Zinfty$-algebra 
and of a module over such an algebra; we can construct their localization 
theory and their spectra exactly in the same way we did it 
in~\ptref{p:spec.flat.zalg}.

Again, we shall return to all these questions later in a more general 
context. We would like to finish this quite long chapter by remarking 
that $I:\ZinfFlat\to\ZinfMod$ commutes with 
tensor products and inner Homs (the reader is invited to check this easy 
statement once tensor products and inner Homs over arbitrary
generalized rings are defined in~\ptref{p:mod.over.genrg}), 
so almost all our previous considerations and computations inside $\ZinfFlat$ 
will retain their validity in~$\ZinfMod$.


\cleardoublepage

\let\nxpointtoc=\nxpointtocb

\mysection{Generalities on monads}

In this chapter we want to collect some general facts on monads, which will be 
of use to us in the remaining part of this work. Almost all of these facts are 
very well known and can be found, for example, in~\cite{MacLane}. We recall 
also some proofs, especially in the cases where we need to generalize them 
later to the case of inner monads.

After that we study the concept of inner endofunctors and inner monads on a 
Cartesian closed category (e.g.\ a topos). This theory doesn't seem to be 
widely known; however, most definitions and proofs appear just as 
natural generalizations of their counterparts in the theory of monads.

Later we shall usually consider monads on the category of sets, and 
inner monads on topoi; we shall see that the latter intuitively correspond 
to sheaves of monads over the category of sets, 
hence their significance for our further developments.

\nxpointtoc{AU $\otimes$-categories}
Let $\cA$ be an AU (associative with unity) $\otimes$-category. This means 
that we are given a bifunctor $\otimes:\cA\times\cA\to\cA$, an associativity 
constraint $\alpha$, i.e.\ a family of functorial isomorphisms 
$\alpha_{X,Y,Z}:(X\otimes Y)\otimes Z\simto X\otimes(Y\otimes Z)$ satisfying 
the pentagon axiom for quadruple tensor products, and that we are given 
an unit object $\Unit_\cA$ together with functorial isomorphisms 
$\Unit_\cA\otimes X\cong X\cong X\otimes\Unit_\cA$, compatible with the 
associativity constraint in a natural way which need not be explained 
explicitly here. By a well-known theorem of MacLane, this data allows 
us to define multiple tensor products of any (possibly empty) finite sequence 
of objects of $\cA$, and to establish canonical isomorphisms between 
tensor products of such multiple products. In short, we get a notion 
of tensor product on $\cA$, enjoying all usual properties of tensor product 
of modules over a commutative ring, {\em apart from commutativity}, since 
we don't impose any commutativity constraint. We'll profit by this remark 
by working with this tensor product in the customary way, without bothering 
to write down explicitly all arising canonical isomorphisms.

\nxsubpoint\label{sp:ext.tens.act} (External $\otimes$-action.) 
Once the AU $\otimes$-category $\cA$ is fixed, we can consider its 
{\em external (left) $\otimes$-action,} or {\em $\obslash$-action,} 
on some category $\cB$. This means that 
we are given some ``external tensor product'' $\obslash:\cA\times\cB\to\cB$ 
(which will be also denoted by the same symbol $\otimes$ as the tensor 
product of $\cA$, when no confusion can arise), together with some external  
associativity constraint, i.e.\ a family of functorial isomorphisms 
$(X\otimes Y)\obslash M\simto X\obslash(Y\obslash M)$, satisfying a variant of 
the pentagon axiom with respect to quadruple products of sort 
$X\otimes Y\otimes Z\obslash M$, where $X$, $Y$, $Z$ are in $\cA$ and 
$M$ is in $\cB$. We also require some functorial isomorphisms 
$\Unit_\cA\obslash M\simto M$ to be given and to be compatible with given 
external associativity constraint.

As before, we don't use these canonical isomorphisms explicitly, but simply 
write down multiple tensor products $X_1\otimes\cdots\otimes X_n\obslash M$, 
for any $X_i$ in $\cA$ and $M$ in $\cB$, and use different canonical 
isomorphisms between them without explaining or even naming them.

We have also the notion of a {\em right\/} external $\otimes$-action 
(also called $\oslash$-action) 
of $\cA$ on $\cB$, given by some bifunctor $\oslash:\cB\times\cA\to\cB$ 
together with some external associativity and unit constraints. This 
case can be reduced to one considered before by permuting the 
order of arguments to $\otimes$ and~$\oslash$.

\nxsubpoint 
Of course, we have both a left and a right external $\otimes$-action 
of $\cA$ on itself. Another example, even more trivial, is given by taking 
for $\cB$ the empty or the point category.

\nxsubpoint\label{sp:tens.functor}
We have the notion of a {\em $\otimes$-functor $F:\cA\to\cA'$} between 
two AU $\otimes$-categories $\cA$ and $\cA'$. By definition, this is 
a functor $F:\cA\to\cA'$, together with a family of functorial isomorphisms 
$F(X\otimes_\cA Y)\simto F(X)\otimes_{\cA'}F(Y)$ and an isomorphism 
$F(\Unit_\cA)\simto\Unit_{\cA'}$, compatible with associativity and unit 
constraints on $\cA$ and $\cA'$.

If we are given some external $\otimes$-action of $\cA$ on some $\cB$, and 
of $\cA'$ on $\cB'$, we can speak about {\em external $\otimes$-functors 
(or $\obslash$-functors) $G:\cB\to\cB'$, compatible with~$F$}. 
By definition, this means that we 
are given some functor $G:\cB\to\cB'$, together with functorial isomorphisms 
$G(X\obslash_{\cB} M)\simto F(X)\obslash_{\cB'} G(M)$, compatible with 
external associativity and unit constraints on $\cB$ and $\cB'$ in a 
natural way.

In particular, we can apply this definition in case $\cA'=\cA$, $F=\Id_{\cA}$. 
Then we are given external $\otimes$-actions of $\cA$ on two categories 
$\cB$ and $\cB'$, and study $\obslash$-functors $G:\cB\to\cB'$ compatible with 
these $\otimes$-actions.

Another special case is given by $\cB=\cA$, $\cB'=\cA'$. Then we see that 
$F$ is compatible with itself.

\nxsubpoint
We have also a notion of natural transformation $\zeta:G\to G'$ of two 
$\obslash$-functors $G,G':\cB\to\cB'$, compatible with same $\otimes$-functor 
$F:\cA\to\cA'$. This is simply a natural transformation $\zeta:G\to G'$ of 
corresponding underlying functors, such that the following square is 
commutative for all $X\in\Ob\cA$ and $Y\in\Ob\cB$, where the horizontal 
arrows come from the $\obslash$-structure on $G$ and $G'$:
\begin{equation}
\xymatrix@+1pc{
G(X\obslash M)\ar[d]^{\zeta_{X\obslash M}}\ar[r]^{\sim}&
F(X)\obslash G(M)\ar[d]^{1_{F(X)}\obslash\zeta_M}\\
G'(X\obslash M)\ar[r]^{\sim}&F(X)\obslash G'(M)}
\end{equation}
Natural transformations of $\oslash$-functors are defined similarly.

\nxsubpoint\label{sp:examp.tens.cat}
Let's mention some important examples of $\otimes$-categories and 
$\otimes$-actions.

a) Take $\cA:=\catSets$, and $\otimes:=\times$, i.e.\ we consider the 
the cartesian products on the category of sets. Of course, this is 
an AU (even ACU) $\otimes$-category, and it acts on itself.

b) Take $\cA:=\catMod{K}$ for any commutative ring~$K$, and 
let $\otimes:=\otimes_K$ be the usual tensor product of $K$-modules. 
Again, this is clearly an AU (even ACU) $\otimes$-category, acting on 
itself.

c) Let $K'$ be any (associative) $K$-algebra. Let $\cA:=\catMod{K}$ as before, 
and $\cB:=\catMod{K'}$. Then any tensor product $X\otimes_K M$, where 
$X$ is a $K$-module and $M$ is a (left) $K'$-module, has a canonical 
$K'$-module structure, so we get an external tensor product 
$\obslash:\cA\times\cB\to\cB$. It satisfies indeed all requirements 
for an external $\otimes$-action, so we get an external $\otimes$-action 
of $\catMod{K}$ on $\catMod{K'}$.

d) If $F:\cA\to\cA'$ is an $\otimes$-functor, then it induces both a 
left and a right $\otimes$-action of $\cA$ on $\cA'$, given by 
$X\obslash M:=F(X)\otimes_{\cA'} M$ (resp.\ $M\oslash X:=M\otimes_{\cA'} F(X)$).
If $K'$ was commutative in the previous example, we see that that 
example can be recovered in this manner starting from 
the $\otimes$-functor $F:\catMod{K}\to\catMod{K'}$ 
given by $X\mapsto X_{(K')}=K'\otimes_K M$.

\nxsubpoint\label{sp:def.alg} (Algebras.) 
By definition, an {\em algebra\/} in an AU $\otimes$-category $\cA$ 
is a triple $A=(A,\mu,\epsilon)$, consisting of an object $A$ of $\cA$, 
a ``multiplication'' morphism $\mu: A\otimes A\to A$ and a ``unit'' 
morphism $\epsilon:\Unit\to A$, subject to usual associativity and unit 
axioms: $\mu\circ(1_A\otimes\mu)=\mu\circ(\mu\otimes 1_A):A\otimes A\otimes A
\to A$ and $\mu\circ(1_A\otimes\epsilon)=1_A=\mu\circ(\epsilon\otimes1_A)$. 
Note that in the first of these axioms we have to 
identify $A\otimes(A\otimes A)$ with $(A\otimes A)\otimes A$ 
by means of the associativity constraint, 
and, similarly, we have implicitly used the unit constraint in the 
second axiom. So AU $\otimes$-categories are actually the minimal 
possible context sufficient to consider algebras. Note that we cannot 
impose any commutativity conditions on $A$ without a commutativity constraint 
on $\cA$.

An {\em algebra homomorphism $f: (A,\mu_A,\epsilon_A) \to 
(B,\mu_B,\epsilon_B)$ in $\cA$} is by definition 
a morphism $f:A\to B$, such that $\mu_B\circ(f\otimes f)=f\circ\mu_A$ 
and $f\circ\epsilon_A=\epsilon_B$. In this way we construct the 
{\em category of algebras in $\cA$,} denoted by $\catAlg(\cA)$.

This category has an initial object, called {\em unit\/} or {\em initial 
algebra}, given by the natural algebra structure on $\Unit_\cA$ coming from 
the canonical isomorphism $\Unit_\cA\otimes\Unit_\cA\cong\Unit_\cA$. Indeed, 
the only algebra homomorphism from $\Unit_\cA$ to an arbitrary algebra 
$A=(A,\mu,\epsilon)$ is defined by $\epsilon:\Unit_\cA\to A$. If $\cA$ has 
a final object $e_\cA$, it admits a unique algebra structure, thus becoming 
the final object in $\catAlg(\cA)$.

\nxsubpoint
Of course, this definition is motivated by the situation of 
\ptref{sp:examp.tens.cat}b), where we recover the notion of an 
associative $K$-algebra with unity. In other situations we obtain 
other interesting notions, e.g.\ in \ptref{sp:examp.tens.cat}a) 
we recover the category of monoids.

\nxsubpoint\label{sp:proj.lim.alg} (Projective limits of algebras.)
Note that all the projective limits which exist in $\cA$, exist also 
in $\catAlg(\cA)$, and they can be computed in~$\cA$. Indeed, if we have 
to compute $\projlim(A_\alpha,\mu_\alpha,\epsilon_\alpha)$, we compute 
first $A:=\projlim A_\alpha$ in $\cA$, denote the canonical projections by 
$\pi_\alpha:A\to A_\alpha$, and define $\mu:A\otimes A\to A$ and $\epsilon:
\Unit_\cA\to A$ starting from compatible families of morphisms 
$\mu_\alpha(\pi_\alpha\otimes\pi_\alpha):A\otimes A\to A_\alpha$ and 
$\epsilon_\alpha:\Unit_\cA\to A$. We see that the forgetful functor 
$\catAlg(\cA)\to\cA$ commutes with arbitrary projective limits, 
and in particular is left exact, if finite projective limits exist in $\cA$. 
In this situation it preserves monomorphisms, and since in addition it is 
faithful, we see that $\phi:A'\to A$ is a monomorphism in $\catAlg(\cA)$ 
iff $\phi$ is a monomorphism in $\cA$.

\nxsubpoint\label{sp:subalg} (Subalgebras.)
Given an algebra $A=(A,\mu,\epsilon)$ in $\catAlg(\cA)$ and a monomorphism
$i: A'\to A$, we see that there is at most one algebra structure 
$(\mu',\epsilon')$ on $A'$ 
compatible with $i$, since we must have $i\mu'=\mu(i\otimes i)$ and 
$i\epsilon'=\epsilon$. Note that whenever $\mu'$ and $\epsilon'$ satisfying 
these relations exist, they satisfy the associativity and unit axioms 
automatically, thus giving required algebra structure on $A'$. In particular, 
we can apply this to embeddings $A'\to A$ of subobjects $A'\subset A$. 
If a compatible algebra structure on $A'$ exists, we call it the 
{\em induced structure\/}, and say that $A'$ is a {\em subalgebra} of $A$. 
If finite projective limits exist in $\cA$, then the subalgebras of a fixed 
algebra $A$ are the same thing as the subobjects of $A$ in $\catAlg(\cA)$.

\nxsubpoint\label{sp:def.mod} (Modules over an algebra.) 
Suppose we are given a (left) external $\otimes$-action of $\cA$ on $\cB$, 
and an algebra $A=(A,\mu,\epsilon)$ in $\cA$. Then a {\em (left) $A$-module 
in $\cB$} is by definition a pair $M=(M,\alpha)$, consisting of an 
object $M$ of $\cB$, and a left action $\alpha$ of $A$ on $M$, i.e.\ a
morphism $\alpha:A\obslash M\to M$, such that $\alpha\circ(\mu\obslash 1_M)=
\alpha\circ(1_A\obslash\alpha)$ and $\alpha\circ(\epsilon\obslash 1_M)=1_M$. 
Note that we have again implicitly used here the external associativity and 
unity constraint for the $\otimes$-action of $\cA$ on $\cB$, so our 
context seems to be the most general context in which the notion of 
an $A$-module makes sense.

Of course, a {\em morphism of $A$-modules $f:(M,\alpha_M)\to(N,\alpha_N)$ 
in $\cB$} is just a morphism $f:M\to N$ in $\cB$, compatible with given 
$A$-actions, i.e.\ $f\circ\alpha_M=\alpha_N\circ(1_A\obslash f)$. In this way 
we obtain the category of $A$-modules in $\cB$, denoted by $\cB^A$.

If we have a {\em right\/} $\otimes$-action $\oslash:\cB\times\cA\to\cB$ of 
$\cA$ on $\cB$, we can define {\em right\/} $A$-modules $M=(M,\alpha)$ in 
$\cB$ in a similar way. In this case $\alpha:M\oslash A\to M$ is required 
to satisfy $\alpha\circ(1_M\oslash\mu)=\alpha\circ(\alpha\oslash 1_A)$ 
and $\alpha\circ(1_M\oslash\epsilon)=1_M$. This yields the category $\cB_A$ 
of right $A$-modules in $\cB$.
Usually we will state and prove the statements only for left modules, 
leaving their counterparts for right modules to the reader.

In particular, we can always take the canonical left and right 
$\otimes$-action of $\cA$ on itself. This leads to a definition of 
left and right $A$-modules in $\cA$, which are usually called simply 
left and right $A$-modules. For example, the multiplication $\mu$ gives 
both a left and a right $A$-module structure on $A$ itself.

\nxsubpoint
Of course, this definition is motivated by the situation of 
\ptref{sp:examp.tens.cat}b) again, where we recover the usual notions of 
a left or right $A$-module over an associative $K$-algebra $A$. In the 
situation of \ptref{sp:examp.tens.cat}a) we recover the notion of 
a (left or right) action of a monoid on a set, and 
in~\ptref{sp:examp.tens.cat}c) 
we obtain the notion of an $A\otimes_K K'$-module.

\nxsubpoint\label{sp:submod.proj.lim} (Submodules and projective limits.)
Similarly to what we had in \ptref{sp:proj.lim.alg}, all projective limits 
existing in $\cB$ exist also in $\cB^A$ and can be essentially computed in 
$\cB$; if finite projective limits exist in $\cB$, then $\phi:N\to M$ is 
a monomorphism in $\cB^A$ iff it is a monomorphism in $\cB$. In this 
situation we have a notion of {\em submodules\/} $M'$ of 
a module $M=(M,\alpha)$, 
i.e.\ subobjects $M'\subset M$ in $\cB$, such that $\alpha(1_A\obslash i)$ 
factorizes into $i\alpha'$, thus yielding an $A$-module structure on $M'$; 
here $i:M'\to M$ is the canonical embedding. Again, when finite projective 
limits exist in $\cB$, the set of submodules of $M$ is actually the set of 
subobjects of $M$ in $\cB^A$.

\nxsubpoint\label{sp:gen.scal.restr} (Scalar restriction.)
We have a very natural notion of scalar restriction in this setup. Namely, 
if $f:A\to A'$ is an algebra homomorphism in $\cA$, and we are given 
a (left) $A'$-module $N=(N,\alpha')$ in some category $\cB$ on which 
$\cA$ $\otimes$-acts, then we can define its {\em scalar restriction 
with respect to $f$}, denoted by $f^*N$, to be the same object $N$ of $\cB$ 
equipped with (left) $A$-action $\alpha:=\alpha'\circ(f\obslash 1_N):
A\obslash N\to N$. As usual, this gives us a functor 
$f^*:\cB^{A'}\to\cB^{A}$.

\nxsubpoint\label{prop:triv.sc.ext} (Scalar extension.)
In some cases the scalar restriction functor has a left adjoint, called the 
{\em scalar extension\/} functor. In general it can be not so easy to 
construct it, but we'll consider now a very simple special case. 
Fix an algebra $A=(A,\mu,\epsilon)$ and consider the scalar restriction 
$\epsilon^*:\cB^A\to\cB^{\Unit_\cA}$ with respect to the only homomorphism 
$\epsilon:\Unit_\cA\to A$ from the initial algebra $\Unit_\cA$ into $A$. 
Clearly, $\cB^{\Unit_\cA}$ is canonically equivalent (even isomorphic) 
to $\cB$, since any object of $\cB$ admits exactly one $\Unit_\cA$-module 
structure, and $\epsilon^*$ can be identified with the forgetful 
functor $\Gamma=\Gamma_A:\cB^A\to\cB$, $(M,\alpha)\mapsto M$.

\begin{Propz} 
For any algebra $A=(A,\mu,\epsilon)$ in $\cA$ and any 
$\obslash$-action of $\cA$ on $\cB$ the forgetful functor 
$\Gamma=\Gamma_A:\cB^A\to\cB$, $(M,\alpha)\mapsto M$, has a left adjoint 
$L=L_A:\cB\to\cB^A$, given by $X\mapsto (A\obslash X, \mu\obslash 1_X)$. 
The functorial isomorphism $\Hom_{\cB^A}(L(X),M)\simto\Hom_\cB(X,\Gamma(M))$ 
maps an $A$-module homomorphism $\phi:A\obslash X\to M$ into 
$\phi^\flat: X\to M$ given by $\phi^\flat:=\phi\circ(\epsilon\obslash 1_X)$, 
and its inverse maps a morphism $\psi:X\to M$ into 
$\psi^\sharp:=\alpha\circ(1_A\obslash\psi)$, where $\alpha:A\obslash M\to M$ 
is the multiplication of $M$.
Moreover, the adjointness natural transformations $\xi:\Id_{\cB}\to\Gamma L$ 
and $\eta:L\Gamma\to\Id_{\cB^A}$ can be described explicitly as follows: 
$\xi_X=\epsilon\obslash 1_X: X\to A\obslash X$ for any $X\in\Ob\cB$, and
$\eta_{(M,\alpha)}=\alpha: (A\obslash M,\mu\obslash 1_M)\to (M, \alpha)$.
\end{Propz}
\begin{Proof}
We have to check that $\psi^\sharp=\alpha\circ(1_A\obslash\psi)$ 
is indeed an $A$-module homomorphism $(A\obslash X,\mu\obslash1_X)\to (M,\alpha)$ for any $\psi:X\to M$. This means $\psi^\sharp\circ(\mu\obslash 1_X)=
\alpha\circ(1_A\obslash\psi^\sharp)$, i.e.\ 
$\alpha\circ(1_A\obslash\psi)\circ(\mu\obslash1_X)=
\alpha\circ(1_A\obslash\alpha)\circ(1_{A\otimes A}\obslash\psi)$. This is clear 
since $\alpha\circ(1_A\obslash\alpha)=\alpha\circ(\mu\obslash1_M)$ and 
$(1_A\obslash\psi)\circ(\mu\obslash1_X)=\mu\obslash\psi=(\mu\obslash1_M)\circ 
(1_{A\otimes A}\obslash\psi)$. Next, we have to check that the maps
$\phi\mapsto\phi^\flat$ and $\psi\mapsto\psi^\sharp$ are inverse to 
each other. This is also quite straightforward: 
$(\phi^\flat)^\sharp=\alpha\circ(1_A\obslash\phi^\flat)=
\alpha\circ(1_A\obslash\phi)\circ(1_A\otimes\epsilon\obslash 1_X)=
\phi\circ(\mu\obslash1_X)\circ(1_A\otimes\epsilon\obslash 1_X)=
\phi\circ\bigl((\mu\circ(1_A\otimes\epsilon))\obslash1_X\bigr)=\phi$
and $(\psi^\sharp)^\flat=\psi^\sharp\circ(\epsilon\obslash 1_X)=
\alpha\circ(1_A\obslash\psi)\circ(\epsilon\obslash1_X)=
\alpha\circ(\epsilon\obslash\psi)=\alpha\circ(\epsilon\obslash1_M)\circ\psi=
\psi$. This completes the proof of adjointness of $L$ and $\Gamma$. 
The formulas for $\xi$ and $\eta$ follow immediately from 
$\xi_X=(1_{L(X)})^\flat$ and $\eta_{(M,\alpha)}=(1_M)^\sharp$.
\end{Proof}

\nxsubpoint\label{sp:mod.fapplic} (Application of functors.)
Given an $\otimes$-functor $F:\cA\to\cA'$, we can apply it to any 
algebra $A=(A,\mu,\epsilon)$ in $\cA$, obtaining an algebra 
$F(A)=\bigl(F(A),F(\mu),F(\epsilon)\bigr)$, where we identify for simplicity 
$F(A\otimes A)$ with $F(A)\otimes F(A)$, and $F(\Unit_{\cA})$ with 
$\Unit_{\cA'}$. Clearly, $F$ transforms algebra homomorphisms into 
algebra homomorphisms, hence it induces a functor $\catAlg(F)$ 
from the category $\catAlg(\cA)$ of algebras in $\cA$ 
into the category $\catAlg(\cA')$ of algebras in $\cA'$; 
this functor will be usually denoted by the same letter~$F$.

Similarly, if we have an $\otimes$-action of $\cA$ on $\cB$ and of 
$\cA'$ on $\cB'$, and an external $\otimes$-functor 
$G:\cB\to\cB'$ compatible with $F$, then we can apply $G$ to any 
$A$-module $M=(M,\alpha)$ in $\cB$, thus obtaining an $F(A)$-module 
$G(M)=\bigl(G(M),G(\alpha)\bigr)$ in $\cB'$. It is clear again that 
this gives us actually a functor $\tilde{G}$ (or simply $G$) 
from the category $\cB^A$ of $A$-modules 
in~$\cB$ into the category ${\cB'}^{F(A)}$ of $F(A)$-modules in $\cB'$.

Furthermore, any natural transformation of $\obslash$-functors 
$\zeta:G\to G'$ induces an $F(A)$-module homomorphism 
$\zeta_M:G(M)\to G'(M)$ in $\cB'$ for any $A$-module $M$ in $\cB$, 
thus effectively defining a natural transformation $\tilde\zeta:\tilde{G}\to 
\tilde G'$ of functors $\cB^A\to{\cB'}^{F(A)}$.

\nxpointtoc{Categories of functors}\label{p:cat.funct}
Recall that the category of all categories $\catCat$ is actually a 
(strictly associative) 2-category. 
In other words, the set of functors $\catFunct(\cC,\cD)=
\Hom_\catCat(\cC,\cD)$ between two arbitrary categories $\cC$ and $\cD$ 
is not just a set, but a category, since we have the notion of 
a natural transformation between two functors from $\cC$ to $\cD$, 
and the composition map $\circ:\catFunct(\cD,\cE)\times\catFunct(\cC,\cD)
\to\catFunct(\cC,\cE)$ is not just a map of sets, but a functor.

This means that, apart from the usual composition of natural transformations 
$\xi'\circ\xi=\xi'\xi:F\to F''$, defined for any natural transformations 
$\xi':F'\to F''$, $\xi:F\to F'$ between functors $F,F',F'':\cC\to\cD$, 
we have also the notion of {\em $\star$-composition} 
$\eta\star\xi: GF\to G'F'$ for any natural transformations 
$\xi:F\to F':\,\cC\to\cD$ and $\eta:G\to G':\,\cD\to\cE$. The requirement 
for $\circ$ to be a functor means that $\id_G\star\id_F=\id_{GF}$ 
and $(\eta'\star\xi')(\eta\star\xi)=(\eta'\eta)\star(\xi'\xi)$.

We adopt the usual convention and denote $\eta\star\id_F$ simply by 
$\eta\star F$, and $\id_G\star\xi$ by $G\star\xi$. Hence we have 
$\eta\star\xi=(\eta\star F')(G\star\xi)=(G'\star\xi)(\eta\star F)$.

\nxsubpoint 
Let's make these operations with natural transformations of functors more 
explicit. First of all, a natural transformation $\xi:F\to G$ between 
two functors from $\cC$ to $\cD$ is by definition a collection of 
morphisms in~$\cD$ $(\xi_X)_{X\in\Ob\cC}$, indexed by objects of $\cC$, 
where $\xi_X:F(X)\to G(X)$ are required to satisfy 
$\xi_Y\circ F(\phi)=G(\phi)\circ\xi_X$ for any morphism $\phi:X\to Y$ in $\cC$.
Therefore, we would like to have explicit descriptions of our operations 
with natural transformations in terms of such families.

a) First of all, the usual composition $\xi'\xi=\xi'\circ\xi:F\to F''$ 
of two natural transformations $\xi:F\to F'$ and $\xi':F'\to F''$ of functors 
from $\cC$ to $\cD$ is computed in the obvious way: 
$(\xi'\xi)_X:=\xi'_X\circ\xi_X:F(X)\to F''(X)$.

b) Now, given a functor $G:\cD\to\cE$, we would like to obtain a natural 
transformation $G\star\xi=\id_G\star\xi: GF\to GF'$. 
This is done in a very natural way: $(G\star\xi)_X:=G(\xi_X):GF(X)\to GF'(X)$.

c) Similarly, given a natural transformation $\eta:G\to G'$, we have to 
define $\eta\star F=\eta\star\id_F:GF\to G'F$. Again, the most natural way 
of doing this is correct: $(\eta\star F)_X:=\eta_{F(X)}:GF(X)\to G'F(X)$.

d) Now we combine these two particular cases into a general formula for 
$\eta\star\xi:GF\to G'F'$ for any two $\xi:F\to F'$ and $\eta:G\to G'$. 
We must have $\eta\star\xi=(\eta\star F')(G\star\xi)=
(G'\star\xi)(\eta\star F)$, hence 
$(\eta\star\xi)_X=\eta_{F'(X)}\circ G(\xi_X)=G'(\xi_X)\circ\eta_{F(X)}$:
\begin{equation}
\xymatrix@+1pc{
GF(X)\ar[r]^{G(\xi_X)}\ar[d]^{\eta_{F(X)}}\ar[rd]^{(\eta\star\xi)_X}
&GF'(X)\ar[d]^{\eta_{F'(X)}}\\
G'F(X)\ar[r]^{G'(\xi_X)}&G'F'(X)}
\end{equation}
These two descriptions of $(\eta\star\xi)_X$ coincide, i.e.\ the above square 
is commutative, simply because of compatibility of $\eta_{F(X)}$ and 
$\eta_{F'(X)}$ with the morphism $\xi_X:F(X)\to F'(X)$ in $\cD$.

e) Once the equivalence of these two descriptions of $\eta\star\xi$ is checked,
we deduce immediately the formula $(\eta'\star\xi')(\eta\star\xi)=
\eta'\eta\star\xi'\xi$, which actually means that $\catCat$ is a 
(strictly associative) 2-category.

\nxsubpoint\label{sp:tens.str.endof}
Now we can exploit this 2-category structure on $\catCat$ to obtain 
some important examples of AU $\otimes$-categories and their external 
$\otimes$-actions. First of all, consider the category of endofunctors 
$\catEndof(\cC):=\catFunct(\cC,\cC)$ on any category $\cC$. Let's 
take the composition functor $\circ:\catEndof(\cC)\times\catEndof(\cC)\to
\catEndof(\cC)$ as the tensor product on $\catEndof(\cC)$. In other words, 
we put $G\otimes F:=GF$ for any two endofunctors $G$ and $F$ on $\cC$, 
and for any two natural transformations $\eta:G\to G'$ and $\xi:F\to F'$ 
of such endofunctors we put $\eta\otimes\xi:=\eta\star\xi:GF\to G'F'$. 
Of course, the identity functor $\Id_\cC$ of $\cC$ is the unity for this 
$\otimes$-structure.

In this way we obtain an AU $\otimes$-structure on $\catEndof(\cC)$. Clearly, 
it is even strictly associative (i.e.\ we have $(X\otimes Y)\otimes Z=
X\otimes(Y\otimes Z)$, and the associativity constraint is identity), since 
$\catCat$ is strictly associative (i.e.\ $(XY)Z=X(YZ)$).

\nxsubpoint\label{sp:tens.act.funct}
Moreover, for any other category $\cD$ we have a left $\otimes$-action of 
$\cA:=\catEndof(\cC)$ on $\catFunct(\cD,\cC)$, and a right $\otimes$-action 
of $\cA$ on $\catFunct(\cC,\cD)$, given in both cases by the composition 
functor $\circ:\catEndof(\cC)\times\catFunct(\cD,\cC)\to\catFunct(\cD,\cC)$ 
(resp.\ $\circ:\catFunct(\cC,\cD)\times\catEndof(\cC)\to\catFunct(\cC,\cD)$). 
Therefore, in both cases we have $F\obslash H=FH$ (resp.\ $H\oslash F=HF$) 
for any $F\in\Ob\catEndof(\cC)$ and $H:\cD\to\cC$ (resp.\ $H:\cC\to\cD$), 
and the action on natural transformations is given by the $\star$-product.

\nxsubpoint\label{sp:cat.fapplic}
Any functor $H:\cD\to\cD'$ defines functors $H_*:\catFunct(\cC,\cD)\to
\catFunct(\cC,\cD')$ and $H^*:\catFunct(\cD',\cC)\to\catFunct(\cD,\cC)$, 
given again by the composition functors in $\catCat$, i.e.\ 
$H_*:G\mapsto HG, \eta\mapsto H\star\eta$ and 
$H^*:G\mapsto GH, \eta\mapsto \eta\star H$. It is immediate that 
$H_*$ (resp.\ $H^*$) is compatible with right (resp.\ left) 
external $\otimes$-action of $\cA=\catEndof(\cC)$ on corresponding categories, 
i.e.\ $H_*$ is an $\oslash$-functor, and $H^*$ is an $\obslash$-functor.

Moreover, any natural transformation $\zeta:H\to H'$ induces a 
natural transformation of $\obslash$-functors $\zeta^*:H^*\to {H'}^*$, given 
by $(\zeta^*)_G:=G\star\zeta:GH\to GH'$, and a natural transformation 
of $\oslash$-functors $\zeta_*:H_*\to H'_*$, given by
$(\zeta_*)_G:=\zeta\star G:HG\to H'G$.

\nxsubpoint\label{sp:tens.act.cat}
So far this reasoning was valid in any 2-category. Let's do something more 
specific now. Put $\cD:=\catPt$ --- the point (or final) category, consisting 
of exactly one object $*$ and only one morphism --- the identity of $*$. 
Then we get a right external $\otimes$-action of $\cA=\catEndof(\cC)$ 
on $\catFunct(\cC,\catPt)$, which is not really interesting since 
this category is equivalent (even isomorphic) to $\catPt$, and a 
left $\otimes$-action of $\cA$ on $\catFunct(\catPt,\cC)$. This is more 
interesting, since $\catFunct(\catPt,\cC)$ is canonically isomorphic 
to $\cC$ itself.

Therefore, we get an external $\otimes$-action of $\cA=\catEndof(\cC)$ on 
$\cC$, $\obslash:\catEndof(\cC)\times\cC\to\cC$, which is easily seen to 
be the evaluation map: $F\obslash X=F(X)$ for any endofunctor $F$ and any 
object $X$ of $\cC$, and $\xi\obslash f=\xi_Y\circ F(f)=F'(f)\circ\xi_X:
F(X)\to F'(Y)$ for 
any natural transformation $\xi:F\to F'$ and any morphism $f:X\to Y$ in $\cC$.
In particular, $\xi\obslash X=\xi_X$ and $F\obslash f=F(f)$.

The situation of \ptref{sp:cat.fapplic} also admits a similar special case, 
if we put $\cD:=\catPt$, and rename $\cD'$ to $\cD$. We see that for any 
object $Z\in\Ob\cD$ the corresponding ``evaluation functor'' 
$Z^*:\catFunct(\cD,\cC)\to\cC$, $Z^*:G\mapsto G(Z), \eta\mapsto\eta_Z$ 
is an $\obslash$-functor with respect to the external $\otimes$-action of 
$\cA$ on $\catFunct(\cD,\cC)$ and on~$\cC$.
In this case any morphism $\zeta:Z\to Z'$ in $\cD$ gives rise to 
a natural $\obslash$-transformation $\zeta^*:Z^*\to {Z'}^*$
of corresponding evaluation functors, given by $\zeta^*_G:=G(\zeta)$.

\nxsubpoint\label{sp:set.th.issues} (Set-theoretical issues.)
Strictly speaking, the category of all categories is no more legal 
mathematical object than the set of all sets. Therefore, we encounter some 
set-theoretical complications, which really do not affect anything 
important for us. 
However, we would like to explain here how we are going to avoid them, 
without paying too much attention to similar issues in future.

To do this we work in Tarski--Grothendieck set theory, i.e.\ we accept the 
usual Zermelo--Frenkel axiomatics for set theory and mathematics 
(recall that all mathematical objects are sets in this picture), and 
we accept Grothendieck's universe axiom as well, 
which tells us that any set can be embedded into an {\em universe}, i.e.\ 
a set~$\univU$ closed under all usual set-theoretical operations with its 
elements (e.g.\ the union of a family of sets $\in\univU$ indexed by a 
set $\in\univU$ also belongs to $\univU$). Usually we fix such an universe 
from the very beginning, choosing it as to contain all sets of interest to 
us (e.g.\ the set of integers $\bbZ$). Then all usual constructions 
produce sets $\in\univU$, and we restrict ourselves to working with such sets 
(cf.\ appendix to SGA~4~I).

A set is said to be {\em $\univU$-small\/} (or just {\em small}, when 
$\univU$ is fixed) if it belongs to $\univU$, or in some cases if we can 
establish a bijection between this set and an element of $\univU$. 
A category is said to be $\univU$-small, or just small when $\univU$ is fixed, 
if it belongs to $\univU$, i.e.\ if both its set of objects and its set 
of morphisms are $\univU$-small. There is also a notion of a 
$\univU$-category~$\cC$, where the only requirement is for all individual 
Hom-sets $\Hom_\cC(x,y)$ to be $\univU$-small for the less restrictive 
understanding of $\univU$-smallness (cf.\ SGA~4~I for more details).   

When we consider the categories $\catSets$, $\catAb$, {\dots} of,
 say, sets, abelian groups and so on, 
we actually consider the categories $\univU$-$\catSets$, $\univU$-$\catAb$ 
\dots of $\univU$-small sets, $\univU$-small abelian groups 
and so on (i.e.\ we require these sets, abelian groups, {\dots} to belong to 
$\univU$). In this way we obtain $\univU$-categories, which, however, are 
not $\univU$-small.

Now, if we consider the category $\univU$-$\catCat$ of $\univU$-small 
categories, this perfectly makes sense, and all previous constructions can 
be applied. For example, for any two $\univU$-small categories $\cC$ and $\cD$ 
the category of functors $\catFunct(\cC,\cD)$ is also $\univU$-small.

However, this approach is too restrictive for our purposes, since we cannot 
even consider the category of endofunctors on the category of sets. 
Therefore, we choose a larger universe $\univV\ni\univU$, and consider the 
category $\catCat:=\univV$-$\univU$-$\catCat$ of $\univV$-small 
$\univU$-categories. Then all previous constructions begin to make sense; 
however, observe that even in this case 
$\catFunct(\cC,\cD)$ usually will not be an $\univU$-category, but only
a $\univV$-category. For example, $\catEndof(\catSets)$ is not an 
$\univU$-category.

We'll deal with this minor problem later by restricting our attention to 
certain full subcategories of $\catFunct(\cC,\cD)$, which will be actually
$\univU$-categories. Then the output of our constructions won't involve the 
larger universe $\univV$, and won't depend on its choice, 
even if it is required at some intermediate stages. 
One example of such full $\univU$-subcategory of a category of functors 
is given by the subcategory of {\em algebraic\/} endofunctors in 
$\catEndof(\catSets)$, which will be studied in the next chapter.

\nxpointtoc{Monads}
Now we combine the AU $\otimes$-structure on the categories of endofunctors 
$\catEndof(\cC)$ studied in \ptref{p:cat.funct} with the definition~%
\ptref{sp:def.alg} of algebras in $\otimes$-categories and 
of modules over such 
algebras. Of course, we thus obtain just the definition of monads. 
Our choice of this approach is motivated by the theory of inner endofunctors 
and inner monads to be developed later. The statements for this second 
case will be almost the same, but direct proofs tend to be more 
lengthy and slightly more complicated.

\begin{DefD}\label{def:monad}
A {\em monad $\Sigma$ on} (or {\em over}) {\em some category $\cC$}) 
is an algebra (cf.~\ptref{sp:def.alg}) 
in the category $\catEndof(\cC)$ with respect to 
its AU $\otimes$-structure defined in~\ptref{sp:tens.str.endof}. In other 
words, a monad $\Sigma$ over $\cC$ is a triple $\Sigma=(\Sigma,\mu,\epsilon)$ 
consisting of an endofunctor $\Sigma:\cC\to\cC$ and two natural 
transformations: the {\em multiplication} $\mu:\Sigma^2=\Sigma\Sigma\to\Sigma$
and the {\em unit} $\epsilon:\Id_\cC\to\Sigma$, required to satisfy the 
associativity axiom
$\mu\circ(\Sigma\star\mu)=\mu\circ(\mu\star\Sigma):\Sigma^3\to\Sigma$ 
and the unit axiom 
$\mu\circ(\Sigma\star\epsilon)=\id_\Sigma=\mu\circ(\epsilon\star\Sigma)$. 

A {\em morphism of monads} $\phi:\Sigma\to\Xi$ is simply a morphism of algebras
in $\catEndof(\cC)$, i.e.\ a natural transformation of underlying endofunctors 
$\phi:\Sigma\to\Xi$, such that $\phi\circ\epsilon_\Sigma=\epsilon_\Xi$ and 
$\phi\circ\mu_\Sigma=\mu_\Xi\circ(\phi\star\phi)$. The category of all monads 
over $\cC$ will be denoted by $\catMonads(\cC)$. Hence 
$\catMonads(\cC)=\catAlg(\catEndof(\cC))$.
\end{DefD}

\nxsubpoint
The definition of a monad $\Sigma=(\Sigma,\mu,\epsilon)$ over some 
category $\cC$ can be made even more explicit in terms of individual 
components of natural transformations $\mu:\Sigma^2\to\Sigma$ and 
$\epsilon:\Id_\cC\to\Sigma$, i.e.\ morphisms $\mu_X:\Sigma^2(X)\to\Sigma(X)$ 
and $\epsilon_X:X\to\Sigma(X)$ parametrized by $X\in\Ob\cC$. The requirement 
for $\mu$ and $\epsilon$ to be natural transformations translates into 
$\Sigma(f)\circ\mu_X=\mu_Y\circ\Sigma^2(f)$ and $\Sigma(f)\circ\epsilon_X=
\epsilon_Y\circ f$ for any morphism $f:X\to Y$ in $\cC$, and our axioms 
translate into $\mu_X\circ\Sigma(\mu_X)=\mu_X\circ\mu_{\Sigma(X)}:
\Sigma^3(X)\to\Sigma(X)$ and $\mu_X\circ\Sigma(\epsilon_X)=
\id_{\Sigma(X)}=\mu_X\circ\epsilon_{\Sigma(X)}$ for all $X\in\Ob\cC$.

Similarly, a morphism of monads $\phi:\Sigma\to\Xi$ is simply a collection 
of morphisms $\phi_X:\Sigma(X)\to\Xi(X)$, parametrized by $X\in\Ob\cC$, 
such that $\phi_Y\circ\Sigma(f)=\Xi(f)\circ\phi_X$ for any morphism 
$f:X\to Y$ in $\cC$, and $\phi_X\circ\epsilon_{\Sigma,X}=\epsilon_{\Xi,X}$ 
and $\phi_X\circ\mu_{\Sigma,X}=\mu_{\Xi,X}\circ\Xi(\phi_X)\circ\phi_{\Sigma(X)}
=\mu_{\Xi,X}\circ\phi_{\Xi(X)}\circ\Sigma(\phi_X)$ for all $X\in\Ob\cC$. 
Here, of course, $\Xi(\phi_X)\circ\phi_{\Sigma(X)}=(\phi\star\phi)_X=
\phi_{\Xi(X)}\circ\Sigma(\phi_X)$, so only one equality has to be checked in 
the last axiom.

\nxsubpoint\label{sp:submonads} (Submonads and projective limits of monads.)
Note that all projective limits that exist in $\cC$, exist also in 
$\cA=\catEndof(\cC)$, since they can be computed componentwise: 
$(\projlim F_\alpha)(X):=\projlim F_\alpha(X)$. 
We say that $F'$ is a {\em subfunctor\/} 
of $F$, if $F'(X)\subset F(X)$ for all $X\in\Ob\cC$. When $\cC$ has 
finite projective limits, the subfunctors of some endofunctor $F$ are actually 
its subobjects in $\cA$, and we are in position to 
apply~\ptref{sp:proj.lim.alg} and~\ptref{sp:subalg} to $\catMonads(\cC)=
\catAlg(\cA)$.

We see that all projective limits existing in $\cC$ exist also in 
$\catMonads(\cC)$, and they can be computed componentwise. Moreover, 
we obtain the notion of a {\em submonad\/} $\Sigma'$ of a monad 
$\Sigma$: by definition it is a subfunctor $\Sigma'\subset\Sigma$ stable 
under multiplication and unit of $\Sigma$, or equivalently a monad $\Sigma'$, 
such that $\Sigma'(X)\subset\Sigma(X)$ for all $X\in\Ob\cC$, and 
$\Sigma'\to\Sigma$ is a morphism of monads.

\nxsubpoint
Once we have the definition of a monad $\Sigma$ over $\cC$ as an algebra 
in AU $\otimes$-category $\cA:=\catEndof(\cC)$, we can use the 
left (resp.\ right) $\otimes$-action of this category on 
$\catFunct(\cD,\cC)$ and $\cC$ (resp.\ on $\catFunct(\cC,\cD)$) explained 
in~\ptref{sp:tens.act.funct} and~\ptref{sp:tens.act.cat} to define the 
categories of left (resp.\ right) $\Sigma$-modules in 
$\catFunct(\cD,\cC)$ and $\cC$ (resp.\ in $\catFunct(\cC,\cD)$) according 
to~\ptref{sp:def.mod}; the corresponding categories of modules are denoted 
by $\catFunct(\cD,\cC)^\Sigma$ and $\cC^\Sigma$ (resp.\ 
$\catFunct(\cC,\cD)_\Sigma$), according to conventions of {\em loc.cit.} 
In all of these situations we have a faithful forgetful functor $\Gamma$, 
acting from corresponding category of modules into the underlying category.

\nxsubpoint\label{sp:cat.left.sigmamod} (Category of $\Sigma$-modules.)
Let's make these definitions more explicit, starting from 
the category~$\cC^\Sigma$ of (left) $\Sigma$-modules in $\cC$. Its objects 
are pairs $M=(M,\alpha)$, consisting of an object $M$ of $\cC$, and a 
{\em $\Sigma$-structure} or {\em $\Sigma$-action} $\alpha$ on~$M$, 
i.e.\ a morphism $\alpha:\Sigma\obslash M=\Sigma(M)\to M$, such that 
$\alpha\circ\mu_M=\alpha\circ\Sigma(\alpha):\Sigma^2(M)\to M$ and 
$\alpha\circ\epsilon_M=\id_M$. Objects of $\cC^\Sigma$, i.e.\ these pairs 
$M=(M,\alpha)=(M,\alpha_M)$ are called {\em $\Sigma$-modules} or 
{\em $\Sigma$-objects}; if $\cC=\catSets$, we also call them 
{\em $\Sigma$-sets.} Another classical terminology is ``$\Sigma$-algebras''; 
we avoid it since in our setup these are modules over algebra $\Sigma$, 
not algebras.

Morphisms $f:(M,\alpha_M)\to(N,\alpha_N)$ in $\cC^\Sigma$ are simply those 
morphisms $f:M\to N$ in $\cC$, which agree with given $\Sigma$-structures, 
i.e.\ such that $f\circ\alpha_M=\alpha_N\circ\Sigma(f)$. 
When no confusion can arise, the subset $\Hom_{\cC^\Sigma}(M,N)\subset
\Hom_\cC(M,N)$ is denoted simply by $\Hom_\Sigma(M,N)$.

\nxsubpoint
Let's consider the category $\catFunct(\cD,\cC)^\Sigma$ now. 
First of all, any functor $H:\cD'\to\cD$ induces an $\obslash$-functor 
$H^*:\catFunct(\cD,\cC)\to\catFunct(\cD',\cC)$, hence a functor 
$H^*:\catFunct(\cD,\cC)^\Sigma\to\catFunct(\cD',\cC)^\Sigma$, and 
natural transformations $\zeta:H\to H'$ induce natural transformations 
$\zeta^*:H^*\to {H'}^*$
(cf.\ \ptref{sp:cat.fapplic} and \ptref{sp:mod.fapplic}).

Similarly, the evaluation functor $Z^*:\catFunct(\cD,\cC)\to\cC$, 
$M\mapsto M(Z)$, defined by any $Z\in\Ob\cD$, is an $\obslash$-functor 
(cf.\ \ptref{sp:tens.act.cat}), hence it induces a functor 
$Z^*:\catFunct(\cD,\cC)^\Sigma\to\cC^\Sigma$, and morphisms 
$\zeta:Z\to Z'$ induce natural transformations $\zeta^*:Z^*\to {Z'}^*$ 
of these functors.

In other words, for any $M=(M,\alpha)\in\Ob\catFunct(\cD,\cC)^\Sigma$, 
where $M:\cD\to\cC$ is a functor and $\alpha:M\to\Sigma M$ is a natural 
transformation, such that $\alpha\circ(\mu\star M)=
\alpha\circ(\Sigma\star\alpha):\Sigma^2M\to M$, we get
an object $Z^*M$ of $\cC^\Sigma$, easily seen to be 
$M(Z)=\bigl(M(Z),\alpha_Z\bigr)$, and any $\zeta:Z\to Z'$ induces a 
$\Sigma$-morphism $\zeta^*M=M(\zeta): M(Z)\to M(Z')$.

In this way we obtain a functor $\tilde M:\cD\to\cC^\Sigma$, which 
is usually denoted with the same letter $M$, since it coincides with~$M$ 
on the level of the underlying objects; applying 
this for different $\Sigma$-modules $M$ and their homomorphisms, we 
construct a functor $\catFunct(\cD,\cC)^\Sigma\to\catFunct(\cD,\cC^\Sigma)$, 
which is now easily seen to be an equivalence, and even an isomorphism of 
categories. One can say that $\cC^\Sigma$ ``represents'' the functor 
$\catFunct(-,\cC)^\Sigma$ on the category of all categories. Let's state 
this result separately:

\begin{PropD}\label{prop:left.sigma.funct}
For any monad $\Sigma$ over a category $\cC$ and any category $\cD$ 
we have a canonical\/ {\em isomorphism} of categories 
$\catFunct(\cD,\cC)^\Sigma\to\catFunct(\cD,\cC^\Sigma)$ which maps 
$(M,\alpha:M\to\Sigma M)$ into the functor $Z\mapsto(M(Z),\alpha_Z)$. 
This isomorphism of categories is compatible with corresponding forgetful 
functors, i.e.\ it identifies $\catFunct(\cD,\cC)^\Sigma\to
\catFunct(\cD,\cC)$ with $(\Gamma_\Sigma)_*:\catFunct(\cD,\cC^\Sigma)\to
\catFunct(\cD,\cC)$, where $\Gamma_\Sigma:\cC^\Sigma\to\cC$ is the forgetful 
functor for $\cC^\Sigma$. Moreover, the functor $\catFunct(\cD,\cC)^\Sigma
\to\catFunct(\cD',\cC)^\Sigma$ induced by any functor $H:\cD'\to\cD$ 
is identified with $H_*:\catFunct(\cD,\cC^\Sigma)\to\catFunct(\cD',\cC^\Sigma)$,
and this statement extends to natural transformations arising from 
any $\zeta:H\to H'$. A similar statement is also true for evaluation functors 
$Z^*:\catFunct(\cD,\cC)^\Sigma\to\cC^\Sigma$ arising from objects $Z$ 
of~$\cD$.
\end{PropD}
\begin{Proof} The functor mentioned in the statement has been constructed 
before. It remains to check that it is an isomorphism of categories, but 
this is quite straightforward, once we understand that this statement 
actually means that to give a left $\Sigma$-module structure on a 
functor $F:\cD\to\cC$ is the same thing as to give a left $\Sigma$-structure 
on each individual $F(X)$, depending functorially on $X\in\Ob\cD$. 
Another possibility --- construct the inverse functor directly, by 
first introducing a left $\Sigma$-functor structure on 
$\Gamma:\cC^\Sigma\to\cC$, and then applying $\tilde{F}^*$ to it for any 
$\tilde{F}:\cD\to\cC^{\Sigma}$, thus obtaining a left $\Sigma$-functor 
structure on $\Gamma\tilde{F}:\cD\to\cC$.
\end{Proof}

\nxsubpoint\label{sp:sigma.free.obj}
In particular, we see that the category $\cA^\Sigma$ 
of left $\Sigma$-modules in $\cA:=\catEndof(\cC)$ 
is isomorphic to $\catFunct(\cC,\cC^\Sigma)$. Since the multiplication $\mu$ 
defines a left $\Sigma$-action on $\Sigma$ itself, we obtain a functor 
$L=L_\Sigma:\cC\to\cC^{\Sigma}$, such that 
$L(X)=\bigl(\Sigma(X),\mu_X\bigr)$ for any object $X$ in $\cC$, and 
$L(f)$ is given by $\Sigma(f)$ for any morphism $f:X\to Y$ in $\cC$. 
Moreover, $L(X)$ can be described as $(\Sigma\obslash X, \mu\obslash\id_X)$ 
in terms of the $\obslash$-action of $\cA$ on $\cC$, so we are in position 
to apply~\ptref{prop:triv.sc.ext}:

\begin{Propz}
The forgetful functor $\Gamma=\Gamma_\Sigma:\cC^\Sigma\to\cC$, 
$(M,\alpha)\mapsto M$, has a left adjoint $L=L_\Sigma:\cC\to\cC^{\Sigma}$, 
given by $X\mapsto\bigl(\Sigma(X),\mu_X\bigr)$. The adjointness isomorphism 
$\Hom_{\cC^\Sigma}\bigl(L(X),(M,\alpha)\bigr)\simto\Hom_\cC(X,M)$ transforms 
a $\Sigma$-morphism $\phi:L(X)=\Sigma(X)\to M$ into 
$\phi^\flat:=\phi\circ\epsilon_X$, and its inverse transforms 
$\psi:X\to\Gamma(M)=M$ into $\psi^\sharp:=\alpha\circ\Sigma(\psi)$.

Moreover, $\Gamma L=\Sigma$, and the natural transformations 
$\xi:\Id_\cC\to\Gamma L=\Sigma$ and $\eta:L\Gamma\to\Id_{\cC^\Sigma}$ 
defined by adjointness of $L$ and $\Gamma$  
are given by $\xi=\epsilon$ and $\eta_{(M,\alpha)}=\alpha$.

Finally, the pair of adjoint functors $L$ and $\Gamma$ 
completely determines the monad $\Sigma=(\Sigma,\mu,\epsilon)$ since 
$\Sigma=\Gamma L$, $\epsilon=\xi$ and 
$\mu=\Gamma\star\eta\star L:\Gamma L\Gamma L=\Sigma^2\to \Gamma L=\Sigma$.
\end{Propz}
\begin{Proof} All statements but the last one follow immediately from 
\ptref{prop:triv.sc.ext}. Let's check the formula for $\mu$: 
$(\Gamma\star\eta\star L)_X=\Gamma(\eta_{L(X)})=
\Gamma(\eta_{(\Sigma(X),\mu_X)})=\Gamma(\mu_X)=\mu_X$.
\end{Proof}

\nxsubpoint\label{sp:monad.freeobj}
Since $L_\Sigma$ is a left adjoint to the forgetful functor $\Gamma_\Sigma$, 
the $\Sigma$-object $L_\Sigma(X)=(\Sigma(X),\mu_X)$ will be called 
{\em the free $\Sigma$-object (or $\Sigma$-module) generated by~$X$}, 
and $\Sigma$-objects isomorphic to some $L_\Sigma(X)$ will be called 
{\em free.} Since the underlying object of $L_\Sigma(X)$ is equal to 
$\Sigma(X)$, this $\Sigma$-object will be often denoted simply by 
$\Sigma(X)$, when no confusion can arise. In other words, when we consider 
$\Sigma(X)$ as a $\Sigma$-object, we endow it with the $\Sigma$-structure 
given by $\mu_X$, unless otherwise specified.

\nxsubpoint\label{sp:monad.from.adjoints}
We have seen that the monad $\Sigma$ over $\cC$ is completely determined 
by the pair of adjoint functors $L:\cC\to\cC^\Sigma$ and 
$\Gamma:\cC^\Sigma\to\cC$. 

This construction can be generalized as follows. 
Suppose we are given two adjoint functors $F:\cC\to\cD$ and $G:\cD\to\cC$ 
with adjointness transformations $\xi:\Id_\cC\to GF$ and 
$\eta:FG\to\Id_\cD$, i.e.\ we have $(\xi\star G)(G\star\eta)=\id_G$ and 
$(\eta\star F)(F\star\xi)=\id_F$. Put $\Sigma:=GF\in\Ob\catEndof(\cC)$, 
$\mu:=G\star\eta\star F:\Sigma^2=GFGF\to\Sigma=GF$ and $\epsilon:=\xi:
\Id_\cC\to\Sigma$. Then it is easy to see that $\Sigma:=(\Sigma,\mu,\epsilon)$ 
is a monad over $\cC$. Indeed, $\mu(\Sigma\star\mu)=(G\star\eta\star F)
(GFG\star\eta\star F)=G\star\eta\star\eta\star F=\mu(\mu\star\Sigma)$, 
$\mu(\Sigma\star\epsilon)=(G\star\eta\star F)(GF\star\xi)=G\star\id_F=
\id_\Sigma$, and $\mu(\epsilon\star\Sigma)=\id_\Sigma$ is checked similarly.

Moreover, $\alpha:=G\star\eta:\Sigma G\to G$ determines a left $\Sigma$-module 
structure on $G:\cD\to\cC$: indeed, $\alpha(\Sigma\star\alpha)=
(G\star\eta)(GFG\star\eta)=G\star\eta\star\eta=
(G\star\eta)(G\star\eta\star FG)=\alpha(\mu\star G)$ and 
$\alpha(\epsilon\star G)=(G\star\eta)(\xi\star G)=\id_G$. Therefore, 
according to~\ptref{prop:left.sigma.funct} our left $\Sigma$-functor 
$(G,G\star\eta)\in\Ob\catFunct(\cD,\cC)^\Sigma$ determines in a unique way 
a functor $\tilde{G}:\cD\to\cC^\Sigma$, such that 
$\tilde{G}(Z)=\bigl(G(Z),G(\eta_Z)\bigr)$ for any $Z\in\Ob\cD$. Clearly, 
$G=\Gamma_\Sigma\tilde{G}$, i.e.\ we have shown that $G:\cD\to\cC$ factorizes 
through $\cC^\Sigma$ in a natural way. It is also easy to see that 
$\tilde{G}F=L_\Sigma$.

\begin{DefD}\label{def:monadic.funct}
A functor $G:\cD\to\cC$ is called\/ {\em monadic}, if it admits a left 
adjoint $F:\cC\to\cD$, and if the induced functor $\tilde{G}:\cD\to\cC^\Sigma$ 
is an equivalence of categories, where $\Sigma$ is the monad over $\cC$ 
defined by $F$ and $G$ as described above.
\end{DefD}

Clearly, in this case we can replace $\cD$ with $\cC^\Sigma$, 
$G$ with $\Gamma_\Sigma$, and $F$ with the left adjoint $L_\Sigma$ of 
$\Gamma_\Sigma$.

Note that the monadicity of $G$ doesn't depend on the choice of its left 
adjoint $F$, since such an adjoint is unique up to a unique isomorphism.

\nxsubpoint
The concept of monadicity is very important since it allows us to replace 
some categories $\cD$ with the categories of $\Sigma$-objects inside a 
simpler category $\cC$ with respect to a certain monad $\Sigma$ over $\cC$, 
thus providing a uniform description of these categories $\cD$. 
There are several interesting criteria of monadicity, e.g.\ the
{\em Beck's monadicity theorem} (cf.~\cite{MacLane}). However, what is 
really important for us is that if $\cC$ is the category of sets and
$\cD$ is a category defined by some algebraic structure 
(e.g.\ the category of groups, monoids, rings, left modules over a fixed ring, 
\dots), then the forgetful functor $F:\cD\to\catSets$ is monadic, 
i.e.\ any of these ``algebraic'' categories is equivalent (usually even 
isomorphic) to $\catSets^\Sigma$ for a suitable monad $\Sigma$ over $\catSets$.
In fact, such categories are defined by a certain special class of monads, 
which we call {\em algebraic}; they will be studied in some detail in 
the next chapter.

\nxsubpoint\label{prop:univ.lower.sigma}
We have seen in \ptref{prop:left.sigma.funct} 
that $\cC^\Sigma$ in a certain sense represents the category of 
left $\Sigma$-functors
$\catFunct(\cD,\cC)^\Sigma$ for variable category $\cD$. 
Is a similar result true for the categories of right $\Sigma$-functors 
$\catFunct(\cC,\cD)_\Sigma$? The answer turns out to be positive:
\begin{Propz}
Given a monad $\Sigma=(\Sigma,\mu,\epsilon)$ over a category $\cC$, 
denote by $\cC_\Sigma$ 
the category with the same objects as $\cC$, and with morphisms given by
$\Hom_{\cC_\Sigma}(X,Y):=\Hom_\cC(X,\Sigma(Y))\cong
\Hom_{\cC^\Sigma}(L_\Sigma(X),L_\Sigma(Y))$. The composition of morphisms 
is given by the second description of\/ $\Hom_{\cC_\Sigma}$; in terms 
of the first description it is given by 
$(\psi,\phi)\mapsto\mu_Z\circ\Sigma(\psi)\circ\phi$,
where $\phi\in\Hom_{\cC_\Sigma}(X,Y)$ and $\psi\in\Hom_{\cC_\Sigma}(Y,Z)$. 
In other words, $\cC_\Sigma$ is equivalent to the full subcategory of 
free objects in $\cC^\Sigma$.

Consider the functor $I=I_\Sigma:\cC\to\cC_\Sigma$, which is identical on 
objects and maps $\phi:X\to Y$ into $\epsilon_Y\circ\phi\in
\Hom_{\cC_\Sigma}(X,Y)$. Then the category $\catFunct(\cC,\cD)_\Sigma$ 
of right $\Sigma$-functors $\cC\to\cD$ is canonically isomorphic to 
the category $\catFunct(\cC_\Sigma,\cD)$, 
and under this isomorphism the forgetful functor is identified with 
$I^*:\catFunct(\cC_\Sigma,\cD)\to\catFunct(\cC,\cD)$.

More explicitly, $F=(F,\alpha)\in\Ob\catFunct(\cC,\cD)_\Sigma$, where 
$F:\cC\to\cD$ and $\alpha:F\Sigma\to F$, corresponds under this 
isomorphism to the functor 
$\tilde{F}:\cC_\Sigma\to\cD$ which coincides with $F$ on objects, and 
transforms $\phi:X\to \Sigma(Y)$ into $\tilde{F}(\phi):=\alpha_Y\circ F(\phi):
F(X)\to F(Y)$.
\end{Propz}
\begin{Proof}
Let's check first the equivalence of the two descriptions of composition in 
$\cC_\Sigma$. In the notations of~\ptref{sp:sigma.free.obj} this amounts 
to check 
$\psi^\sharp\circ\phi^\sharp=(\mu_Z\circ\Sigma(\psi)\circ\phi)^\sharp$, 
i.e.\ $\mu_Z\circ\Sigma(\psi)\circ\mu_Y\circ\Sigma(\phi)=
\mu_Z\circ\Sigma(\mu_Z)\circ\Sigma^2(\psi)\circ\Sigma(\phi)$; this equality 
follows immediately from $\Sigma(\psi)\circ\mu_Y=\mu_{\Sigma(Z)}\circ
\Sigma^2(\psi)$ and from the associativity condition $\mu_Z\circ\Sigma(\mu_Z)=
\mu_Z\circ\mu_{\Sigma(Z)}$.

Now we want to construct a functor $\catFunct(\cC_\Sigma,\cD)\to
\catFunct(\cC,\cD)_\Sigma$, inverse to one constructed above. To do this 
we construct first a right $\Sigma$-structure $\beta:I\Sigma\to I$ on 
$I:\cC\to\cC_\Sigma$, and then apply to it the $\oslash$-functor
$\tilde{F}_*$ for any $\tilde{F}:\cC_\Sigma\to\cD$, thus obtaining 
a right $\Sigma$-functor $\tilde{F}_*(I,\beta)=
(\tilde{F}I,\tilde{F}\star\beta)$; this construction, being functorial 
in $\tilde{F}$, yields the required functor $\catFunct(\cC_\Sigma,\cD)\to
\catFunct(\cC,\cD)_\Sigma$.

So we have to construct $\beta:I\Sigma\to I$. For any $X\in\Ob\cC$
we must produce some $\beta_X\in\Hom_{\cC_\Sigma}(\Sigma(X),X)=
\Hom_\cC(\Sigma(X),\Sigma(X))$; of course, we take $\id_{\Sigma(X)}$. 
It remains to check that $\beta$ is indeed a right $\Sigma$-action on $I$, 
and that our two functors between $\catFunct(\cC,\cD)_\Sigma$ and 
$\catFunct(\cC_\Sigma,\cD)$ are inverse to each other; 
we leave these verifications to the reader.
\end{Proof}

\nxsubpoint
Note that $I=I_\Sigma:\cC\to\cC_\Sigma$ admits a right adjoint 
$K=K_\Sigma:\cC_\Sigma\to\cC$, given by $K(X):=\Sigma(X)$, 
$K(\phi):=\mu_Y\circ\Sigma(\phi)$ for any $\phi:X\to\Sigma(Y)$. Clearly, 
$KI=\Sigma$, and $IK:\cC_\Sigma\to\cC_\Sigma$ coincides with $\Sigma$ 
on objects and transforms $\phi:X\to\Sigma(Y)$ into 
$\epsilon_{\Sigma(Y)}\circ\mu_Y\circ\Sigma(\phi)$. The adjointness 
natural transformations $\xi:\Id_\cC\to KI=\Sigma$ and 
$\eta:IK\to\Id_{\cC_\Sigma}$ are computed as follows: $\xi=\epsilon$ 
and $\eta_X\in\Hom_{\cC_\Sigma}(\Sigma(X),X)=\Hom_\cC(\Sigma(X),\Sigma(X))$ 
is given again by $\id_{\Sigma(X)}$.

It is easy to see that this pair of adjoint functors $I:\cC\to\cC_\Sigma$ 
and $K:\cC_\Sigma\to\cC$ determines by the recipe 
of~\ptref{sp:monad.from.adjoints} the original monad $\Sigma=(\Sigma,\mu,
\epsilon)$: indeed, we have already seen that $KI=\Sigma$ and $\xi=\epsilon$, 
and $K\star\eta\star I=\mu$ is quite easy to check: 
$(K\star\eta\star I)_X=K(\eta_{I(X)})=K(\id_{\Sigma(X)})=\mu_X$.

\nxsubpoint
Moreover, if we are given two adjoint functors $F:\cC\to\cD$ and 
$G:\cD\to\cC$ with adjointness natural transformations 
$\xi:\Id_\cC\to GF$ and $\eta:FG\to\Id_\cD$, defining a monad 
$\Sigma=(GF,G\star\eta\star F,\xi)$ as described 
in~\ptref{sp:monad.from.adjoints}, we obtain a right $\Sigma$-action 
$\eta\star F:F\Sigma\to F$ on $F$, hence by~\ptref{prop:univ.lower.sigma} 
there is a unique functor $\tilde{F}:\cC_\Sigma\to\cD$, such that 
$\tilde{F}I=F$ and $\tilde{F}\star\beta=\eta\star F$. 

In this respect the situation is quite similar (or rather dual) to 
that of~\ptref{sp:monad.from.adjoints}, where we obtained a functor 
$\tilde{G}:\cD\to\cC^\Sigma$. We see that among all pairs of adjoint 
functors $(F,G)$, defining the same monad $\Sigma$ on $\cC$, 
we have an initial object $(I_\Sigma,K_\Sigma)$ and a final object 
$(L_\Sigma,\Gamma_\Sigma)$.

\nxsubpoint\label{sp:Q.sigma}
Proposition \ptref{prop:univ.lower.sigma} is important for us since it 
allows us to identify the category $\cA_\Sigma$ of right $\Sigma$-modules 
in $\cA:=\catEndof(\cC)$ with the category $\catFunct(\cC_\Sigma,\cC)$ 
of functors from $\cC_\Sigma$ to $\cC$.

Combining this with \ptref{prop:left.sigma.funct}, we see that, given 
a monad $\Sigma$ over $\cC$ and a monad $\Xi$ over $\cD$, the category 
$\catFunct(\cD,\cC)^\Sigma_\Xi$ of $\Sigma$-$\Xi$-bimodules in 
$\catFunct(\cD,\cC)$ (i.e.\ we consider the category of 
triples $(M,\alpha,\beta)$, where $M:\cD\to\cC$, $\alpha:\Sigma M\to M$ 
is a left $\Sigma$-action, $\beta:M\Xi\to M$ is a right $\Xi$-action, 
and these actions commute: $\alpha(\Sigma\star\beta)=\beta(\alpha\star\Xi)$) 
is isomorphic to the category $\catFunct(\cD_\Xi,\cC^\Sigma)$.

In particular, the category $\cA^\Sigma_\Sigma$ 
of $\Sigma$-bimodules in $\cA$ is isomorphic to 
$\catFunct(\cC_\Sigma, \cC^\Sigma)$. In other words, a $\Sigma$-bimodule in 
$\cA$ is essentially the same thing as a functor from the subcategory of 
free objects in $\cC^\Sigma$ into the whole of $\cC^\Sigma$. Since 
$\Sigma$ is canonically a bimodule over itself, we obtain a functor 
$Q=Q_\Sigma:\cC_\Sigma\to\cC^\Sigma$, such that $Q_\Sigma I_\Sigma=L_\Sigma$ 
and $\Gamma_\Sigma Q_\Sigma=K_\Sigma$; explicitly, 
$Q(X)=L(X)=(\Sigma(X),\mu_X)$ on objects, and $Q(\phi)=\mu_Y\circ\Sigma(\phi)$ 
for any $\phi\in\Hom_{\cC_\Sigma}(X,Y)=\Hom_\cC(X,\Sigma(Y))
\cong\Hom_{\cC^\Sigma}(L(X),L(Y))$. We see again 
that this functor $Q:\cC_\Sigma\to\cC^\Sigma$ is fully faithful 
(cf.\ \ptref{prop:univ.lower.sigma}). Different functors between 
categories $\cC$, $\cC^\Sigma$ and $\cC_\Sigma$ fit into the 
following commutative diagram:
\begin{equation}
\xymatrix@+1pc{
&\cC_\Sigma\ar[r]^{Q_\Sigma}\ar[rrd]|{K_\Sigma}&
\cC^\Sigma\ar[rd]^{\Gamma_\Sigma}\\
\cC\ar[ru]^{I_\Sigma}\ar[rru]|{L_\Sigma}\ar[rrr]^{\Sigma}&&&\cC}
\end{equation}

\nxsubpoint\label{sp:gen.scalrestr.mod} (Scalar restriction.)
Given a morphism $\rho:\Sigma\to\Xi$ of monads over $\cC$, we obtain 
the {\em scalar restriction functors\/} $\rho^*:\cC^\Xi\to\cC^\Sigma$, 
$\rho^*:\catFunct(\cD,\cC)^\Xi\to\catFunct(\cD,\cC)^\Sigma$ and 
$\rho^*:\catFunct(\cC,\cD)_\Xi\to\catFunct(\cC,\cD)_\Sigma$ according to 
the general recipe of~\ptref{sp:gen.scal.restr}. It is easy to see that 
the second of these functors can be identified with $(\rho^*)_*:
\catFunct(\cD,\cC^\Xi)\to\catFunct(\cD,\cC^\Sigma)$, $H\mapsto\rho^*\circ H$ 
(cf.~\ptref{prop:left.sigma.funct}). Similarly, the third of these functors 
can be identified according to \ptref{prop:univ.lower.sigma} 
with a functor $\catFunct(\cC_\Xi,\cD)\to\catFunct(\cC_\Sigma,\cD)$, 
and this construction is functorial in $\cD$, hence ``by Yoneda'' 
this functor has to be induced by some functor $\rho_*:\cC_\Sigma\to\cC_\Xi$, 
acting in the opposite direction.

Let's describe these two functors $\rho^*:\cC^\Xi\to\cC^\Sigma$ and 
$\rho_*:\cC_\Sigma\to\cC_\Xi$ more explicitly. Clearly, 
$\rho^*:\cC^\Xi\to\cC^\Sigma$ transforms a $\Xi$-module 
$(N,\alpha:\Xi(N)\to N)$ into $\Sigma$-module $\rho^*(N):=
(N,\alpha\circ\rho_N)$, and acts identically on morphisms. 
It is also easy to check that $\rho_*:\cC_\Sigma\to\cC_\Xi$ acts identically 
on $\Ob\cC_\Sigma=\Ob\cC=\Ob\cC_\Xi$, and transforms a morphism 
$f\in\Hom_{\cC_\Sigma}(X,Y)=\Hom_\cC(X,\Sigma(Y))$ into 
$\rho_Y\circ f\in\Hom_{\cC_\Xi}(X,Y)$. Note that 
$\Gamma_\Sigma\rho^*=\Gamma_\Xi$ and $\rho_*I_\Sigma=I_\Xi$.

\nxsubpoint (Base change.)
Since $Q_\Sigma$ establishes an equivalence between $\cC_\Sigma$ and the 
full subcategory of free $\Sigma$-modules in $\cC^\Sigma$, and 
similarly for $Q_\Xi$, we see that 
$\rho_*:\cC_\Sigma\to\cC_\Xi$ can be considered as a partially defined 
functor $\tilde{\rho}_*:\cC^\Sigma\dashrightarrow\cC^\Xi$, defined only 
for free $\Sigma$-modules. Clearly, $\tilde{\rho}_*$ transforms 
$L_\Sigma(X)=\Sigma(X)$ into $\Xi(X)$, and $f:\Sigma(X)\to\Sigma(Y)$ into 
$(\rho_Y\circ f^\flat)^\sharp=
\mu_{\Xi,Y}\circ\Xi(\rho_Y\circ f\circ\epsilon_{\Sigma,X})$ 
(cf.~\ptref{sp:sigma.free.obj}). We claim that 
{\em $\tilde{\rho}_*:\cC^\Sigma\dashrightarrow\cC^\Xi$ is a partially 
defined left adjoint to $\rho^*:\cC^\Xi\to\cC^\Sigma$.} 
Indeed, $\Hom_{\cC^\Xi}(\tilde{\rho}_*L_\Sigma(X),N)=
\Hom_{\cC^\Xi}(L_\Xi(X),N)\cong\Hom_\cC(X,\Gamma_\Xi(N))=
\Hom_\cC(X,\Gamma_\Sigma\rho^*(N))\cong\Hom_{\cC^\Sigma}(L_\Sigma(X),\rho^*N)$. 
In some situations $\tilde{\rho}_*$ can be extended to a well-defined 
{\em base change} (or {\em scalar extension}) {\em functor\/} 
$\rho_*:\cC^\Sigma\to\cC^\Xi$, 
left adjoint to $\rho^*:\cC^\Xi\to\cC^\Sigma$.

\begin{PropD}\label{prop:ex.base.change}
If cokernels of pairs of morphisms exist in $\cC^\Xi$, 
then $\tilde{\rho}_*$ extends to a well-defined {\em scalar extension} 
(or {\em base change}) functor $\rho_*:\cC^\Sigma\to\cC^\Xi$, left adjoint 
to $\rho^*:\cC^\Xi\to\cC^\Sigma$. Moreover, in this situation 
$\rho^*$ is monadic, i.e.\ it induces an equivalence and even an isomorphism 
of categories $\cC^\Xi\to(\cC^\Sigma)^{\Xi/\Sigma}$, where 
$\Xi/\Sigma:=\rho^*\rho_*$ is the monad over $\cC^\Sigma$ 
defined by $\rho_*$ and $\rho^*$ (cf.~\ptref{sp:monad.from.adjoints}).
\end{PropD}

\begin{Proof}
According to Lemma~\ptref{l:mod.coker.free} below, any 
$\Sigma$-module $M=(M,\alpha)$ in $\cC$ can be represented (even in a 
functorial way) as a cokernel of two $\Sigma$-morphisms between two 
free $\Sigma$-modules: $M=\Coker(p,q:L_\Sigma(R)\rightrightarrows 
L_\Sigma(X))$.
If $\rho_*$ exists, it must preserve arbitrary inductive limits and 
transform any free $\Sigma$-object $L_\Sigma(Z)$ into 
free $\Xi$-object $L_\Xi(Z)$, hence $\rho_*(M)$ must be isomorphic 
to $\Coker(\tilde{\rho}_*(p),\tilde{\rho}_*(q))$. We put 
$\rho_*(M):=\Coker(\tilde{\rho}_*(p),\tilde{\rho}_*(q):
L_\Xi(R)\rightrightarrows L_\Xi(X))$. It remains to check 
$\Hom_{\cC^\Xi}(\rho_*M,N)\cong\Hom_{\cC^\Sigma}(M,\rho^*N)$, but this 
is quite clear: $\Hom_{\cC^\Xi}(\rho_*M,N)=
\Hom_{\cC^\Xi}\bigl(\Coker(\tilde{\rho}_*(p),\tilde{\rho}_*(q)), N\bigr)\cong$
$\Ker\bigl(\Hom_{\cC^\Xi}(L_\Xi(X),N)\rightrightarrows
\Hom_{\cC^\Xi}(L_\Xi(R),N)\bigr)$, and since $\tilde{\rho}_*$ is 
a partial left adjoint to $\rho^*$ and $L_\Xi(X)=\tilde{\rho}_*L_\Sigma(X)$, 
this kernel turns out to be isomorphic to
$\Ker\bigl(\Hom_{\cC^\Sigma}(L_\Sigma(X),\rho^*N)\rightrightarrows
\Hom_{\cC^\Sigma}(L_\Sigma(R),\rho^*N)\bigr)\cong
\Hom_{\cC^\Sigma}(M,\rho^*N)$, hence $\rho^*M$ represents indeed whatever 
it has to represent. The monadicity of $\rho^*$ will not be used in the sequel,
so we leave it as an exercise for the reader; an important point here is 
to construct natural transformations $\zeta=(\rho_*)\star\eta_\Sigma
:L_\Xi\Gamma_\Sigma\cong\rho_*L_\Sigma\Gamma_\Sigma\to\rho_*$ and
$\Gamma_\Xi\star\zeta:\Xi\Gamma_\Sigma\to
\Gamma_\Sigma(\Xi/\Sigma):\cC^\Sigma\to\cC$, and to use them 
for constructing the inverse functor $(\cC^\Sigma)^{\Xi/\Sigma}\to\cC^\Xi$. 
Another possibility is to use Beck's monadicity theorem to check that 
the monadicity of $\rho^*$ follows from the existence of a left adjoint to 
$\rho^*$ and from the monadicity of $\Gamma_\Sigma$ and that of 
$\Gamma_\Sigma\rho^*=\Gamma_\Xi$.
\end{Proof}

\begin{LemmaD}\label{l:mod.coker.free}
Any $\Sigma$-module 
$M=(M,\alpha)\in\Ob\cC^\Sigma$ can be represented (even in a functorial way) 
as a cokernel of a pair of morphisms between two free $\Sigma$-objects: 
$M\cong\Coker(p,q:L_\Sigma(R)\rightrightarrows L_\Sigma(X))$. Moreover, 
one can take $X:=M$, $R:=\Sigma(M)$, $p:=\mu_M$, $q:=\Sigma(\alpha):
\Sigma^2(M)\to\Sigma(M)$, and $\alpha$ for the 
strict epimorphism $L_\Sigma(X)=\Sigma(M)\to M$.
\end{LemmaD}
\begin{Proof}
We want to check that $\xymatrix{\Sigma^2(M)\ar@/_2pt/[r]\ar@/^2pt/[r]^{p,q}&
\Sigma(M)\ar[r]^{\alpha}&M}$ is right exact in $\cC^\Sigma$. First of all, 
notice that $\sigma:=\epsilon_M:M\to\Sigma(M)$ and $\tau:=\epsilon_{\Sigma(M)}:
\Sigma(M)\to\Sigma^2(M)$ provide a splitting of this diagram in~$\cC$, i.e.\ 
we have $\alpha\circ\sigma=\id_M$, $p\circ\tau=\id_{\Sigma(M)}$ and 
$q\circ\tau=\sigma\circ\alpha$. 
This together with $\alpha\circ p=\alpha\circ q$ 
implies the right exactness of this diagram 
in~$\cC$, and, since this splitting is preserved by any functor, 
the right exactness of this diagram after an application of~$\Sigma$.

Now, given any $\Sigma$-morphism $\phi:(\Sigma(M),\mu_M)\to (N,\beta)$, 
such that $\phi\circ p=\phi\circ q$, 
we want to show the existence and uniqueness of 
a $\Sigma$-morphism $\psi:M\to N$, such that $\phi=\psi\circ\alpha$. 
Right exactness in~$\cC$ implies existence and uniqueness of such a 
morphism~$\psi$ in~$\cC$; it remains to check that this is a $\Sigma$-morphism,
i.e.\ that $\psi\circ\alpha=\beta\circ\Sigma(\psi)$. Since $\Sigma(\alpha)$ 
is a (split) epimorphism in $\cC$, this follows from 
$\psi\circ\alpha\circ\Sigma(\alpha)=\psi\circ\alpha\circ\mu_M=
\phi\circ\mu_M=\beta\circ\Sigma(\phi)=
\beta\circ\Sigma(\psi)\circ\Sigma(\alpha)$.
\end{Proof}

\nxsubpoint\label{sp:monhom.catfunct}
We have seen that any monad homomorphism $\rho:\Sigma\to\Xi$ induces a 
scalar restriction functor $\rho^*:\cC^\Xi\to\cC^\Sigma$, such that 
$\Gamma_\Sigma\rho^*=\Gamma_\Xi$. We claim that 
{\em any functor $H:\cC^\Xi\to\cC^\Sigma$, such that 
$\Gamma_\Sigma H=\Gamma_\Xi$, is equal to $\rho^*$ for a uniquely 
determined monad homomorphism $\rho:\Sigma\to\Xi$.}
To show this we first observe that by~\ptref{prop:left.sigma.funct} such 
functors~$H:\cC^\Xi\to\cC^\Sigma$ are in one-to-one correspondence with 
left $\Sigma$-structures $\beta:\Sigma\Gamma_\Xi\to\Gamma_\Xi$ 
on~$\Gamma_\Xi$. For example, $H=\rho^*$ corresponds to 
$\beta=(\Gamma_\Xi\star\eta_\Xi)\circ(\rho\star\Gamma_\Xi)$, where 
$\eta_\Xi:L_\Xi\Gamma_\Xi\to\Id_{\cC^\Xi}$ and 
$\xi_\Xi=\epsilon_\Xi:\Id_\cC\to\Xi=\Gamma_\Xi L_\Xi$ are the usual 
adjointness natural transformations. Notice that $\rho$ is uniquely 
determined by $\beta$ since $\rho=\rho\circ\mu_\Sigma\circ
(\Sigma\star\epsilon_\Sigma)=\mu_\Xi\circ(\rho\star\rho\epsilon_\Sigma)=
\mu_\Xi\circ(\rho\star\epsilon_\Xi)=\mu_\Xi\circ(\rho\star\Xi)\circ
(\Sigma\star\epsilon_\Xi)=(\beta\star L_\Xi)\circ(\Sigma\star\epsilon_\Xi)$. 
In general we put $\rho:=(\beta\star L_\Xi)\circ(\Sigma\star\epsilon_\Xi)$ 
for an arbitrary left $\Sigma$-action $\beta$ on $\Gamma_\Xi$. 
We check that $\rho:\Sigma\to\Xi$ is a monad homomorphism using the 
fact that $\beta$ is a left $\Sigma$-action. For example, we 
use $\beta\circ(\epsilon_\Sigma\star\Gamma_\Xi)=\id_{\Gamma_\Xi}$ to 
prove $\rho\circ\epsilon_\Sigma=(\beta\star L_\Xi)\circ (\Sigma\star
\epsilon_\Xi)\circ\epsilon_\Sigma=(\beta\star L_\Xi)\circ(\epsilon_\Sigma\star
\epsilon_\Xi)=(\beta\star L_\Xi)\circ(\epsilon_\Sigma\star\Gamma_\Xi L_\Xi)
\circ\epsilon_\Xi=\epsilon_\Xi$. The compatibility with multiplication 
$\mu_\Xi\circ(\rho\star\rho)=\rho\circ\mu_\Sigma$ is checked similarly 
but longer, using $\mu_\Xi=\Gamma_\Xi\star\eta_\Xi\star L_\Xi$, 
$(\Gamma_\Xi\star\eta_\Xi)\circ(\epsilon_\Xi\star\Gamma_\Xi)=\id_{\Gamma_\Xi}$ 
and $\beta\circ(\Sigma\star\beta)=\beta\circ(\mu_\Sigma\star\Gamma_\Xi)$. 
Finally, we have to check that 
$(\Gamma_\Xi\star\eta_\Xi)\circ(\rho\star\Gamma_\Xi)$ will give back 
the original $\beta$. Indeed, this expression equals 
$(\Gamma_\Xi\star\eta_\Xi)\circ(\beta\star L_\Xi\Gamma_\Xi)\circ
(\Sigma\star\epsilon_\Xi\star\Gamma_\Xi)=(\beta\star\eta_\Xi)\circ
(\Sigma\star\epsilon_\Xi\star\Gamma_\Xi)=\beta\circ
(\Sigma\Gamma_\Xi\star\eta_\Xi)\circ(\Sigma\star\epsilon_\Xi\star\Gamma_\Xi)=
\beta$.

\nxpointtoc{Examples of monads}
Now we would like to give some examples of monads which will be important for 
us later. The most important of them are those defined by pairs of 
adjoint functors $F:\cC\to\cD$ and $G:\cD\to\cC$ 
(cf.~\ptref{sp:monad.from.adjoints}), where usually 
$\cC:=\catSets$, $\cD$ is some ``algebraic category'' (category of sets 
with some algebraic structure, e.g.\ the category of groups, rings and so on), 
$G=\Gamma:\cD\to\catSets$ is the forgetful functor, and $F$ is a left 
adjoint of~$G$. This gives us a monad $\Sigma$ on the endofunctor 
$GF$ over the category of sets, and the forgetful functor $G$ 
usually turns out to be monadic (cf.~\ptref{def:monadic.funct}), 
i.e.\ the induced functor 
$\tilde{G}:\cD\to\catSets^\Sigma$ turns out to be an equivalence 
and even an isomorphism of categories, thus identifying $\cD$ with 
$\catSets^\Sigma$. Another important source of monads is obtained by 
considering submonads of previously constructed monads 
(cf.~\ptref{sp:submonads}).

\nxsubpoint\label{sp:mon.words} (The monad of words.)
Let's take $\cC:=\catSets$, $\cD:=\catMon$ be the category of monoids, 
$G:\catMon\to\catSets$ be the forgetful functor, and $F:\catSets\to\catMon$ 
be its left adjoint, which maps any set $X$ into the free monoid $F(X)$ 
generated by $X$. Recall that $F(X)$ consists of all finite words 
$x_1\ldots x_n$ in alphabet $X$, i.e.\ all $x_i\in X$, $n\geq0$, the 
multiplication on $F(X)$ is given by the concatenation of words, and 
the unit of $F(X)$ is the empty word~$\zeroword$, 
i.e.\ the only word of length zero.
According to \ptref{sp:monad.from.adjoints}, this defines a monad 
$W=(W,\mu,\epsilon)$ over the category of sets, which will be called 
{\em the monad of words.}

We see that $W(X)=\bigsqcup_{n\geq0}X^n$ for any set~$X$. Sequences 
$(x_1,x_2,\ldots,x_n)$ lying in 
$X^n\subset W(X)$ are usually called {\em words of length~$n$ over~$X$}, 
and are written simply as $x_1x_2\ldots x_n$; the individual components $x_i$ 
of a word are usually called {\em letters.} When some of $x_i$ are 
replaced by more complicated expressions, or when there is some other source 
of potential confusion, we enclose some or all of them into braces, thus 
writing $\{x_1\}\{x_2\}\ldots\{x_n\}$ instead of $x_1x_2\ldots x_n$ or 
$(x_1,x_2,\ldots,x_n)$. For example, $\{x\}$ is the one-letter word 
corresponding to an element $x\in X$. The {\em empty word\/} will be usually 
denoted by $\zeroword$, or on some occasions by $\{\}$.

Now the unit $\epsilon:\Id_\catSets\to W$ and the multiplication 
$\mu:W^2(X)\to W(X)$ can be described explicitly as follows. 
$\epsilon_X:X\to W(X)$ maps any $x\in X$ into the corresponding one-letter 
word $\{x\}\in W(X)$. As to the multiplication $\mu_X:W^2(X)\to W(X)$, 
it ``removes one layer of braces'', i.e.\ maps 
$\{x_1x_2\ldots x_n\}\{y_1y_2\ldots\}\ldots\{z_1\ldots z_m\}$ into 
$x_1x_2\ldots x_ny_1y_2\ldots z_1\ldots z_m$. The axioms of a monad can 
be checked now directly. For example, the associativity condition 
$\mu_X\circ\mu_{W(X)}=\mu_X\circ W(\mu_X):W^3(X)\to W(X)$ means simply that 
if we have an expression with two layers of braces, and remove 
first the inner layer and then the outer layer, or conversely, then 
the resulting word in $W(X)$ is the same in both cases. So we might 
have defined the monad $W=(W,\mu,\epsilon)$ directly, without any reference 
to the category of monoids.

\nxsubpoint\label{sp:expl.mona.mono}
Of course, the forgetful functor $G:\catMon\to\catSets$ is monadic, i.e.\ 
the induced functor $\tilde{G}:\catMon\to\catSets^W$ is an equivalence, 
and even an isomorphism of categories. In other words, the structure of
a monoid on some set $X$ is essentially the same thing as a $W$-structure 
$\alpha:W(X)\to X$, and $W$-morphisms are exactly monoid morphisms.

Let's show how this works in this situation. If we have a monoid structure 
$(X,*,e)$ on $X$, where $*:X\times X\to X$ is the multiplication and 
$e\in X$ is the unit, then the corresponding $W$-structure 
$\alpha:W(X)\to X$ simply maps a word $x_1x_2\ldots x_n\in W(X)$ into 
the corresponding product $x_1*x_2*\cdots*x_n$ (this actually follows from 
the description of $\tilde{G}$ given in~\ptref{sp:monad.from.adjoints}). 
Conversely, if we are given a $W$-structure $\alpha:W(X)\to X$ on a set $X$,
we can recover a monoid structure on $X$ by setting $e:=\alpha(\{\})$ 
and $x*y:=\alpha(\{x\}\{y\})=\alpha(xy)$. The associativity and unit 
axioms for this choice of operations on $X$ follow from the axioms 
$\alpha\circ W(\alpha)=\alpha\circ\mu_X$ and $\alpha\circ\epsilon_X=\id_X$ 
for $\alpha$. For example, 
$(x*y)*z=\alpha(\{\alpha(xy)\}\{z\})=\alpha(\{\alpha(xy)\}\{\alpha(z)\})=
(\alpha\circ W(\alpha))(\{xy\}\{z\})=\alpha(\mu_X(\{xy\}\{z\}))=
\alpha(xyz)$, and $x*(y*z)=\alpha(xyz)$ is shown similarly, thus proving 
the associativity of $*$.

The fact that $W$-morphisms $f:(X,\alpha)\to (Y,\beta)$ are exactly monoid 
homomorphisms is also checked similarly. For example, if $f$ is a $W$-morphism,
i.e.\ $\beta\circ W(f)=f\circ\alpha$, then $f(x *_X y)=
f(\alpha(xy))=(\beta\circ W(f))(\{x\}\{y\})=
\beta(\{f(x)\}\{f(y)\})=f(x)*_Yf(y)$.

\nxsubpoint
This construction generalizes to any (Grothendieck) topos $\cE$, where 
we can also put $W_\cE(X):=\bigsqcup_{n\geq0}X^n$ for any $X\in\Ob\cE$. 
Then $\cE^{W_\cE}$ is of course just the category $\catMon(\cE)$ of monoids 
in $\cE$. If $\cE$ is the category $\tilde{\cS}$ of sheaves (of sets) 
on some site $\cS$, then $\cE^{W_\cE}$ is the category of sheaves of monoids 
on~$\cS$. We'll see later that both $W$ and $W_\cE$ are examples of 
algebraic monads, and that $W_\cE$ is the pullback of $W$ with respect 
to the canonical morphism of topoi $\cE\to\catSets$.

\nxsubpoint
Given a monad $\Sigma=(\Sigma,\mu,\epsilon)$ over a category $\cC$ 
with finite coproducts, and an object $U\in\Ob\cC$, we can construct 
a new monad $\Sigma_U=(\Sigma_U,\mu',\epsilon')$ with the underlying functor 
given by $\Sigma_U(X):=\Sigma(U\sqcup X)$ and 
$\Sigma_U(f):=\Sigma(\id_U\sqcup f)$.
Informally, we want to add some 
``constants'' from $U$ to the ``variables'' from $X$. Let's define 
the multiplication and the unit of this monad. Denote for this by 
$i_X:U\to U\sqcup X$ and $j_X:X\to U\sqcup X$ the natural embeddings, and 
by $k_X: U\sqcup\Sigma(U\sqcup X)\to\Sigma(U\sqcup X)$ the morphism 
$k_X:=\langle\Sigma(i_X)\circ\epsilon_U,\id\rangle$
with components $\Sigma(i_X)\circ \epsilon_U=\epsilon_{U\sqcup X}\circ i_X:
U\to\Sigma(U\sqcup X)$ and $\id_{\Sigma(U\sqcup X)}$, i.e.\ 
$k_X\circ i_{\Sigma_U(X)}=\Sigma(i_X)\circ\epsilon_U=\epsilon_{U\sqcup X}\circ
i_X$ and $k_X\circ j_{\Sigma_U(X)}=\id$. 
Now put 
$\epsilon'_X:=\Sigma(j_X)\circ\epsilon_X=\epsilon_{U\sqcup X}\circ j_X:
X\to\Sigma_U(X)$, and define $\mu'_X:\Sigma(U\sqcup\Sigma(U\sqcup X))\to
\Sigma(U\sqcup X)$ by $\mu'_X:=\mu_{U\sqcup X}\circ\Sigma(k_X)$.

Let's check first that $\epsilon'$ is a unit for $\mu'$. Indeed, 
$\mu'_X\circ\Sigma_U(\epsilon'_X)=\mu_{U\sqcup X}\circ\Sigma(k_X)\circ
\Sigma(\id_U\sqcup\epsilon'_X)=\mu_{U\sqcup X}\circ\Sigma\bigl(
\langle\epsilon_{U\sqcup X}\circ i_X,
\epsilon_{U\sqcup X}\circ j_X\rangle\bigr)=
\mu_{U\sqcup X}\circ\Sigma(\epsilon_{U\sqcup X})=\id_{\Sigma(U\sqcup X)}$ and 
$\mu'_X\circ\epsilon'_{\Sigma_U(X)}=\mu_{U\sqcup X}\circ\Sigma(k_X)\circ
\Sigma(j_{\Sigma_U(X)})\circ\epsilon_{\Sigma_U(X)}=\mu_{U\sqcup X}\circ
\Sigma(\id)\circ\epsilon_{\Sigma(U\sqcup X)}=\id_{\Sigma(U\sqcup X)}$.

Next we have to check the associativity axiom for $\mu'$. Before doing this 
let's introduce $\lambda_X:U\sqcup\Sigma^2(U\sqcup X)\to\Sigma^2(U\sqcup X)$ 
by $\lambda_X:=\langle \epsilon_{\Sigma(U\sqcup X)}\circ\epsilon_{U\sqcup X}
\circ i_X, \id\rangle$. Now the associativity condition 
$\mu'_X\circ\Sigma_U(\mu'_X)=\mu'_X\circ\mu'_{\Sigma_U(X)}$ 
will be a consequence of the 
commutativity of the outer circuit of the following diagram:
$$\xymatrix@+12pt{
\Sigma(U\sqcup\Sigma_U^2(X))
\ar[r]^{\Sigma(k_{\Sigma_U(X)})}
\ar[d]^{\Sigma(\id_U\sqcup\Sigma(k_X))}
&\Sigma^2(U\sqcup\Sigma_U(X))
\ar[r]^{\mu_{U\sqcup\Sigma_U(X)}}
\ar[d]^{\Sigma^2(k_X)}
&\Sigma(U\sqcup\Sigma_U(X))
\ar[d]_{\Sigma(k_X)}
\\\Sigma(U\sqcup\Sigma^2(U\sqcup X))
\ar[r]^{\Sigma(\lambda_X)}
\ar[d]^{\Sigma(\id_U\sqcup\mu_{U\sqcup X})}
&\Sigma^3(U\sqcup X)
\ar[r]^{\mu_{\Sigma(U\sqcup X)}}
\ar[d]^{\Sigma(\mu_{U\sqcup X})}
&\Sigma^2(U\sqcup X)
\ar[d]_{\mu_{U\sqcup X}}
\\
\Sigma(U\sqcup\Sigma(U\sqcup X))
\ar[r]^{\Sigma(k_X)}
&\Sigma^2(U\sqcup X)
\ar[r]^{\mu_{U\sqcup X}}
&\Sigma(U\sqcup X)
}$$
The upper right square is commutative simply because the natural transformation
$\mu$ is compatible with $k_X$, and the commutativity of the lower right 
square follows from the associativity of~$\mu$. The two remaining squares 
are obtained by applying $\Sigma$ to some simpler diagrams, and it is 
sufficient to check the commutativity of these simpler diagrams. Both 
of them assert equality of some morphisms from certain coproducts 
$U\sqcup\ldots$, hence it is sufficient to check the equality of the 
restrictions of these morphisms to individual components of the coproduct.

Now restricting to $U$ in the (counterpart of) the upper left square 
yields $\lambda_X\circ i_{\Sigma^2(U\sqcup X)}=
\epsilon_{\Sigma(U\sqcup X)}\circ\epsilon_{U\sqcup X}\circ i_X=
\Sigma(\epsilon_{U\sqcup X}\circ i_X)\circ\epsilon_U=
\Sigma(k_X\circ i_{\Sigma_U(X)})\circ\epsilon_U=
\Sigma(k_X)\circ k_{\Sigma_U(X)}\circ i_{\Sigma^2_U(X)}$, and restricting 
to $\Sigma^2_U(X)$ yields $\Sigma(k_X)$ for both paths; this proves the 
commutativity of this square.

Restricting to $U$ in the lower left square we get 
$\mu_{U\sqcup X}\circ\lambda_X\circ i_{\Sigma^2(U\sqcup X)}=
\mu_{U\sqcup X}\circ\epsilon_{\Sigma(U\sqcup X)}\circ\epsilon_{U\sqcup X}\circ
i_X=\epsilon_{U\sqcup X}\circ i_X=k_X\circ i_{\Sigma(U\sqcup X)}$, and 
restricting to $\Sigma^2(U\sqcup X)$ we get $\mu_{U\sqcup X}$ for both 
paths. This finishes the proof of the commutativity of the whole diagram, 
hence also that of the associativity of~$\mu'$.

\nxsubpoint
Note that the above construction is functorial in~$U$, so we get 
canonical monad homomorphisms $\zeta_f:\Sigma_U\to\Sigma_V$ 
for any $f:U\to V$, given by $\zeta_{f,X}:=\Sigma(f\sqcup\id_X):
\Sigma(U\sqcup X)\to\Sigma(V\sqcup X)$, 
and homomorphisms $\zeta_U:\Sigma\to\Sigma_U$, 
$\zeta_{U,X}:=\Sigma(j_X):\Sigma(X)\to\Sigma(U\sqcup X)$ as well. 

\nxsubpoint\label{sp:words.with.const}
In particular, we can apply this to the monad of words $W$ on the category 
of sets, or to its topos counterpart $W_\cE$, thus obtaining a family of new 
monads $W\langle U\rangle=W_U$ parametrized by sets (resp.\ objects of~$\cE$)
$U$. The alternative notation $W\langle U\rangle$ reflects the fact that 
this is something like ``the algebra of polynomials over $W$ in 
non-commuting variables from~$U$''. More precisely, $\catSets^{W_U}$ 
is the category of pairs $(M,\phi)$, consisting of a monoid $M$ and some 
``choice of constants'', i.e.\ an arbitrary map of sets $\phi:U\to M$.

Clearly, the elements of $W_U(X)=W(U\sqcup X)$ are simply words in the 
alphabet $U\sqcup X$. Usually we don't enclose 
the letters from $U$ in braces, even if we do it with some other letters 
of the same word, thus obtaining expressions like 
$u_1\{x_1\}u_2\{x_2\}\{x_3\}$. The unit morphism $\epsilon'_X:X\to W_U(X)$ 
still maps any $x\in X$ into the corresponding ``variable'' 
(one-letter word) $\{x\}$, and the multiplication 
$\mu'_X:W_U^2(X)\to W_U(X)$ is given by the concatenation of words and 
subsequent identification of ``constants'', i.e.\ letters from~$U$, 
coming from different layers. 
For example, $\mu'_X(u_1\{u_2\{x_1\}\{x_2\}\}\{u_3\})=u_1u_2\{x_1\}\{x_2\}u_3$.
This property is actually the reason which allows us to omit braces around 
letters from $U$ in almost all situations. 

\nxsubpoint
The importance of these monads $W_U$ obtained from the monad of words $W$ 
is due to the fact that 
{\em any algebraic monad appears as a subquotient of a monad of this form.}
We'll return to this issue in the next chapter.

\nxsubpoint\label{sp:mon.from.ring} (Monads defined by rings.) 
Now fix any ring~$R$ (required to be associative, with unity, but not 
necessarily commutative) and consider the forgetful functor 
$\Gamma_R:\catMod{R}\to\catSets$ and its left adjoint $L_R:\catSets\to
\catMod{R}$, which transforms any set $X$ into the corresponding 
free module $R^{(X)}$. As usual, this defines a monad 
$\Sigma_R=(\Sigma_R,\mu,\epsilon)$ over $\catSets$ 
with the underlying endofunctor equal to 
$\Gamma_R L_R$. In other words, $\Sigma_R(X)$ is simply the underlying 
set of $R^{(X)}$, i.e.\ the set of all formal $R$-linear combinations 
$\lambda_1\{x_1\}+\cdots+\lambda_n\{x_n\}$, where $n\geq0$, all 
$\lambda_i\in R$, $x_i\in X$, and where we have denoted by $\{x\}$ 
the basis element corresponding to $x\in X$. Another description: 
$\Sigma_R(X)=\{\lambda:X\to R\,|\, \lambda(x)=0$ for almost all 
$x\in X\}$.

Clearly, the unit $\epsilon_X:X\to\Sigma_R(X)$ maps any $x\in X$ into 
the corresponding basis element $\{x\}$, and 
$\mu_X:\Sigma_R^2(X)\to\Sigma_R(X)$ computes linear combinations of 
formal linear combinations (cf.~\ptref{sp:mu.zinfty}): 
$\mu_X\bigl(\sum_i\lambda_i\{\sum_j\mu_{ij}\{x_j\}\}\bigr)$ 
equals $\sum_{i,j}\lambda_i\mu_{ij}\{x_j\}$. 

Now a $\Sigma_R$-structure $\alpha:\Sigma_R(X)\to X$ on some set $X$ can 
be considered as a way of evaluating formal linear combinations of elements 
of~$X$. In particular, we can define an addition and an action of $R$ on 
$X$ by setting $x+y:=\alpha(\{x\}+\{y\})$ and $\lambda\cdot x:=
\alpha(\lambda\{x\})$. Similarly to what we had in~\ptref{sp:expl.mona.mono}, 
one can check directly that this addition and $R$-action satisfy the axioms for
 a left $R$-module, thus obtaining a one-to-one correspondence between 
$\Sigma_R$-structures and left $R$-module structures. This shows that 
the categories $\catMod{R}$ and $\catSets^{\Sigma_R}$ are equivalent and 
actually isomorphic, hence $\Gamma_R$ is monadic. 

\nxsubpoint\label{sp:underl.set.monad} (The underlying set of a monad.)
Can we recover the ring~$R$, or at least its underlying set, starting from 
the corresponding monad~$\Sigma_R$? The answer is {\em positive}, since 
$R=\Sigma_R(\st1)$, where $\st1=\{1\}$ is the standard one-element set 
(i.e.\ a final object of $\catSets$). Moreover, we can recover the 
multiplication on $R$, as well as the action of $R$ on the underlying set 
of any $\Sigma_R$-module~$X$.

To do this observe that $R_s:=L_{\Sigma_R}(\st1)$ is actually the ring~$R$, 
considered as a left module over itself. Now $R\cong \End_R(R_s)^{op}=
\End_{\Sigma_R}(L_{\Sigma_R}(\st1))^{op}$, and this is an isomorphism of 
monoids, if we consider the multiplication on $\End$'s defined by 
composition of endomorphisms. Moreover, the action of~$R$ on any 
$R$-module (i.e.\ $\Sigma_R$-module) $M$ can be recovered from the canonical 
isomorphism $\Hom_\catSets(\st1,M)\cong
\Hom_{\Sigma_R}(L_{\Sigma_R}(\st1),M)=\Hom_R(R_s,M)$, and from the 
right action of $\End_R(R_s)$ on this set, given again by the composition of 
homomorphisms.

This construction can be generalized to any monad $\Sigma$ over $\catSets$. 
For any such monad we define its {\em underlying set\/~$|\Sigma|$} 
(sometimes denoted simply by $\Sigma$) by $|\Sigma|:=\Sigma(\st1)$. 
Clearly, $|\Sigma|\cong\Hom_\catSets(\st1,\Sigma(\st1))\cong 
\End_\Sigma(L_\Sigma(\st1))\cong\End_{\catSets_\Sigma}(\st1)$, 
and we introduce a monoid structure on 
$|\Sigma|$ by transporting the natural monoid structure of 
$\End_\Sigma(L_\Sigma(\st1))$ given by composition of endomorphisms and 
taking the opposite. This defines the {\em underlying monoid\/} of~$\Sigma$, 
which will be also denoted by $|\Sigma|$. Clearly, $\{1\}:=\epsilon_\st1(1)
\in|\Sigma|$ is the unit for this multiplication.

Moreover, for any $\Sigma$-module $M$ we get a canonical right action of 
monoid $\End_\Sigma(L_\Sigma(\st1))$ on $\Hom_\Sigma(L_\Sigma(\st1),M)\cong
\Hom_\catSets(\st1,M)\cong M$, hence a left action of the opposite monoid 
$|\Sigma|\cong\End_\Sigma(L_\Sigma(\st1))^{op}$. All these constructions 
are functorial: for example, any $\Sigma$-homomorphism $f:M\to N$ is 
at the same time a map of $|\Sigma|$-sets, any homomorphism of monads 
$\rho:\Sigma\to\Xi$ induces a monoid homomorphism 
$|\rho|:=\rho_\st1: |\Sigma|\to|\Xi|$, and the scalar restriction functors 
$\rho^*:\catSets^\Xi\to\catSets^\Sigma$ are compatible with these constructions
as well, i.e.\ the $|\Sigma|$-action on $\rho^*N$ coincides with the 
scalar restriction of the $|\Xi|$-action on $N$ with respect to $|\rho|$.

Of course, in any particular case there is some extra information which 
cannot be recovered by this general construction. For example, the addition 
of $R$ or of an $R$-module is not recovered in this way from $\Sigma_R$.

\nxsubpoint (Functoriality of $\Sigma_R$.)
Any ring homomorphism $\rho:R\to S$ induces a homomorphism 
$\Sigma_\rho:\Sigma_R\to\Sigma_S$ of corresponding monads, explicitly given 
by the collection of maps $\Sigma_{\rho,X}:\Sigma_R(X)\to\Sigma_S(X)$, 
$\sum_{i=1}^n\lambda_i\{x_i\}\mapsto\sum_{i=1}^n\rho(\lambda_i)\{x_i\}$. 
It is easy to see that this natural transformation is indeed compatible 
with the multiplication and the unit of these monads, and that 
$\Sigma_\rho^*:\catSets^{\Sigma_S}=\catMod{S}\to\catSets^{\Sigma_R}=
\catMod{R}$ is just the usual scalar restriction functor 
$\rho^*:\catMod{S}\to\catMod{R}$; by uniqueness of adjoints we see that 
the scalar extension functor $(\Sigma_\rho)_*$ constructed 
in~\ptref{prop:ex.base.change} is isomorphic to the usual scalar 
extension functor $\rho_*:\catMod{R}\to\catMod{S}$, $M\mapsto S\otimes_R M$.

Note that $|\Sigma_\rho|=\rho$, hence $\rho$ can be recovered from 
$|\Sigma_\rho|$. Hence the functor $R\mapsto\Sigma_R$ from the category of 
rings into the category of monads over $\catSets$ is at least faithful. 
We claim that {\em the functor $R\mapsto\Sigma_R$ is fully faithful, 
hence it identifies the category of rings with a full subcategory of 
the category of monads.} In other words, any monad homomorphism 
$\zeta:\Sigma_R\to\Sigma_S$ is of form $\Sigma_\rho$ for a (uniquely 
determined) ring homomorphism $\rho:R\to S$. Clearly, we must have 
$\rho=|\zeta|:R\to S$. We postpone some details of the proof 
until~\ptref{p:hypadd}, where we develop a suitable machinery for such tasks. 
Up to these details the proof goes as follows.

We know already that $\rho:=|\zeta|$ is a monoid 
homomorphism, and we'll check later that it respects addition as well. 
Thus we are reduced to check $\Sigma_{\rho,X}=\zeta_X:
\Sigma_R(X)\to\Sigma_S(X)$ for any set~$X$. 
Since both $\Sigma_R$ and $\Sigma_S$ commute with 
filtered inductive limits and any set is a filtered inductive limit of 
its finite subsets, we can assume $X$ to be finite, e.g.\ 
$X=\stn=\{1,2,\ldots,n\}$. Then we construct canonical projections 
$\pr_i: \Sigma_R(\stn)\to\Sigma_R(\st1)=|\Sigma_R|=R$, 
using only the monad structure 
of $\Sigma_R$, and combine them to obtain a canonical isomorphism 
$\Sigma_R(\stn)\simto R^n$. Since this construction uses only the monad 
structure of $\Sigma_R$, it has to be compatible with $\zeta$, i.e.\ 
it fits into a commutative diagram relating $\zeta_\stn$ with 
$(\zeta_\st1)^n=\rho^n$. This proves $\zeta_\stn=\Sigma_{\rho,\stn}$. 
This also proves that $\rho$ respects addition as well, since the only map 
$\st2=\{1,2\}\to\st1$ induces a map 
$R^2\cong\Sigma_R(\st2)\to R=\Sigma_R(\st1)$, which is nothing else than 
the addition of~$R$. 
Now the functoriality of this construction with respect to 
$\zeta:\Sigma_R\to\Sigma_S$ and the identification of $\zeta_\st2$ with 
$\rho^2$ proves the compatibility of~$\rho$ with the addition on~$R$ and~$S$.

\nxsubpoint\label{sp:mon.as.genr} (Monads as generalized rings.)
We have embedded the category of rings into 
the category of monads over $\catSets$ as a full subcategory. 
We will often identify a ring $R$ with the 
corresponding monad $\Sigma_R$, thus writing $L_R$ instead of $L_{\Sigma_R}$, 
and on some occasion even $R(X)$ instead of $\Sigma_R(X)$. 
Moreover, we can treat arbitrary monads over $\catSets$ as some sort of 
non-commutative generalized rings, and transport some notations usual 
for rings into this more general context. For example, we write 
$\catMod{\Sigma}$ instead of $\catSets^\Sigma$, $\Hom_\Sigma$ instead of 
$\Hom_{\catSets^\Sigma}$, and even $\Sigma^{(X)}$ instead of 
$\Sigma(X)$ or $L_\Sigma(X)$.

However, one should be careful with these identifications, since to 
obtain the correct category of non-commutative generalized rings we have 
to impose the algebraicity condition, satisfied by all monads we considered 
so far. For example, infinite projective limits  
computed in the category of all monads (cf.~\ptref{sp:submonads}) 
and in the category of algebraic monads are different.

Notice that {\em $R\mapsto\Sigma_R$ is left exact, i.e.\ it commutes with 
finite projective limits.} Indeed, given such a limit $R=\projlim R_\alpha$, 
we first observe that $\Sigma_R(X)\cong\projlim\Sigma_{R_\alpha}(X)$ is 
true for a finite $X$, since in this case $\Sigma_R(\stn)\cong R^n$ and 
similarly for $\Sigma_{R_\alpha}$, and then extend this result to the 
case of an arbitrary~$X$ by representing $X$ as a filtered inductive limit 
of finite sets and observing that filtered inductive limits commute with 
finite projective limits.

Nevertheless, if we compute for example $\hat\bbZ_p:=\projlim\bbZ/p^n\bbZ$ in 
the category of monads, we obtain a monad $\hat\bbZ_p$ which is not 
algebraic and in particular is different from $\bbZ_p$, regardless of 
the fact that $|\hat\bbZ_p|=\bbZ_p$. Actually the category of 
$\hat\bbZ_p$-modules turns out to be quite similar to that of 
{\em topological\/} $\bbZ_p$-modules. 
In particular, it doesn't have the usual properties 
of a category of modules over a ring. 

\nxsubpoint\label{sp:examp.submon} (Subrings and submonads.)
Clearly, if $R'$ is a subring of $R$, then $\Sigma_{R'}(X)\subset\Sigma_R(X)$ 
for any set~$X$, i.e.\ $\Sigma_{R'}$ is a {\em submonad\/} of~$\Sigma_R$ 
(cf.~\ptref{sp:submonads}). However, not all submonads of $\Sigma_R$ come 
from subrings of $R$. Let's give some examples of such submonads.

a) Consider the submonad $\bbZ_{\geq0}$ of $\bbZ$, defined by 
$\bbZ_{\geq0}^{(X)}:=\{$formal linear combinations 
of elements of $X$ with {\em non-negative\/} 
integer coefficients$\}$. Clearly, $\bbZ_{\geq0}$ is a subfunctor of 
$\bbZ=\Sigma_\bbZ$, and it is compatible with the multiplication and the 
unit of $\Sigma_{\bbZ}$, since all basis vectors $\{x\}$ belong to 
$\bbZ_{\geq0}^{(X)}$, and a linear combination with non-negative integer 
coefficients of several such formal linear combinations is again a 
formal linear combination with such coefficients. Hence $\bbZ_{\geq0}$ is 
indeed a submonad of $\bbZ$. Notice that $|\bbZ_{\geq0}|$ is indeed the 
multiplicative monoid of all non-negative integers, as suggested by 
our notation. It is easy to see that $\catMod{\bbZ_{\geq0}}$ 
is actually the category of commutative monoids, and that the 
scalar restriction functor for the inclusion $\bbZ_{\geq0}\to\bbZ$ is 
simply the forgetful functor from the category of abelian groups 
to the category of commutative monoids.
We have also a 
unique (surjective) homomorphism of monads 
$W\to\bbZ_{\geq0}$, and the corresponding 
scalar restriction functor is simply the embedding of the category of 
commutative monoids into the category of all monoids.

b) Now consider the submonad $\Zinfty\subset\bbR$, defined by the following 
property: $\Zinfty^{(X)}\subset\bbR^{(X)}$ consists of all {\em octahedral\/}
formal combinations $\sum_i\lambda_i\{x_i\}$ of elements of~$X$, i.e.\ 
we require $\sum_i|\lambda_i|\leq1$, $\lambda_i\in\bbR$. Clearly, any 
octahedral combination of formal octahedral combinations is an octahedral 
combination again, and all basis vectors lie in $\Zinfty^{(X)}$, hence 
$\Zinfty$ is indeed a submonad of~$\bbR$. Of course, $\catMod{\Zinfty}$ 
coincides with the category denoted in the same way in~\ptref{def:zinfmod}, 
where it has been defined as the category of $\Sigma_\infty$-modules, 
for a monad $\Sigma_\infty$ isomorphic to our $\Zinfty$.

Clearly, $|\Zinfty|$ is the segment $[-1,1]\subset\bbR$ with the multiplication
induced by that of~$\bbR$.

c) Similarly, we can consider $\barZinfty\subset\bbC$, such that 
$\barZinfty^{(X)}\subset\bbC^{(X)}$ consists of all formal $\bbC$-linear 
combinations $\sum_i\lambda_i\{x_i\}$, subject to condition 
$\sum_i|\lambda_i|\leq1$, $\lambda_i\in\bbC$. Again, $|\barZinfty|$ is just 
the unit disk in~$\bbC$ with the multiplication induced by that of~$\bbC$.

d) Given two submonads $\Sigma'$, $\Sigma''$ of the same monad $\Sigma$, 
we can construct a new submonad $\Sigma'\cap\Sigma'':=
\Sigma'\times_\Sigma\Sigma''\subset\Sigma$. Since projective limits of monads 
are computed componentwise (cf.~\ptref{sp:submonads}), 
$(\Sigma'\cap\Sigma'')(X)=\Sigma'(X)\cap\Sigma''(X)\subset\Sigma(X)$.

In particular, we can define the ``non-completed localization of $\bbZ$ 
at infinity'' $\Zninfty:=\Zinfty\cap\bbQ\subset\bbR$, similarly to the 
$p$-adic case $\bbZ_{(p)}=\bbZ_p\cap\bbQ\subset\bbQ_p$. Clearly, 
$\Zninfty^{(X)}$ consists of all formal octahedral combinations of elements
of~$X$ with rational coefficients. This generalized ring will be important 
in the construction of $\CompZ$, where it will appear as the local ring 
at infinity.

We can intersect further and consider $\bbZ[1/p]\cap\Zninfty$ for 
any prime~$p$. This corresponds to considering formal octahedral combinations 
with $p$-rational numbers (i.e.\ rational numbers of form $x/p^n$) as 
coefficients. Later we'll construct $\CompZ$ by gluing together 
$\Spec\bbZ$ and $\Spec(\bbZ[1/p]\cap\Zninfty)$ along their principal 
open subsets isomorphic to $\Spec\bbZ[1/p]$.

e) Put $\Fpm:=\bbZ\cap\Zninfty\subset\bbQ$, or equivalently 
$\Fpm:=\bbZ\cap\Zinfty\subset\bbR$. Clearly, if we construct $\CompZ$ 
as described above, $\Fpm$ must be the (generalized) ring of 
global sections of the structural sheaf of~$\CompZ$. We see that 
$\Fpm^{(X)}$ consists of all formal octahedral combinations of elements 
of~$X$ with {\em integer\/} coefficients. This means that there can be 
at most one non-zero coefficient, and it must be equal to $\pm1$. 
In other words, $\Fpm^{(X)}$ consists of $0$, basis elements 
$\{x\}$, and their opposites~$-\{x\}$, i.e.\ 
$\Fpm^{(X)}=0\sqcup X\sqcup-X$. Hence an $\Fpm$-structure $\alpha$ 
on a set $X$ is a map $0\sqcup X\sqcup-X\to X$, identical on $X$. 
Its restriction to $0$ yields a marked (or zero) element $\alpha(0)$ of $X$, 
also denoted by $0$ or~$0_X$. The restriction of~$\alpha$ to $-X$ gives 
a map $-_X:X\to X$, which can be shown to be an involution preserving $0_X$. 
Therefore, $\catMod\Fpm$ is the category of sets with a marked point 
and an involution preserving this marked point.

We introduce an alternative notation $\bbF_{1^2}$ for $\Fpm$. This 
is due to the fact that we'll construct later generalized rings 
$\bbF_{1^n}$, such that $\catMod{\bbF_{1^n}}$ is the category of sets 
with a marked point and a permutation of order~$n$ which fixes the 
marked point.

f) We can also consider the ``field with one element'' 
$\Fone:=\bbZ_{\geq0}\cap\Fpm=\bbZ_{\geq0}\cap\Zinfty$. The positivity 
condition implies that $\Fone^{(X)}=0\sqcup X$, and $\catMod\Fone$ 
is simply the category of sets with a marked point.

g) Here are some more examples of submonads: ${\bm\Delta}\subset\bbR$ 
corresponds to baricentric combinations $\sum_i\lambda_i\{x_i\}$, 
$\lambda_i\in\bbR_{\geq0}$, $\sum_i\lambda_i=1$; then 
$\catMod{{\bm\Delta}}$ is the category of abstract convex sets. 
We also define $\bbR_{\geq0}\subset\bbR$ by considering linear combinations 
with non-negative real coefficients. For any ring $R$ we define 
$\Aff_R\subset R$, obtained by considering formal $R$-linear combinations 
with the sum of coefficients equal to one; then $\catMod{\Aff_R}$ is 
the category of (abstract) affine spaces over~$R$.

h) Finally, the identity functor $\Id_\catSets$ admits a natural monad 
structure, thus becoming the initial object in the category of monads over 
$\catSets$ (cf.~\ptref{sp:def.alg}). 
This monad will be denoted by $\Fempty$, and sometimes called
{\em the field without elements}. 
Note that the category of monads has a final object $\st1$ as well, such that 
$\st1^{(X)}=\{0\}$ for any set~$X$. 
Clearly, this is just the monad corresponding to the classical 
trivial ring~$0$.

\nxsubpoint\label{sp:ex.nonalg.mon}
All monads of the above examples a)--h) turn out to be algebraic, 
i.e.\ they commute with filtered inductive limits of sets.
Let's construct a monad $\hat\bbZ_\infty\supset\Zinfty$ which is {\em not\/} 
algebraic. For this we take 
$\hat\bbZ_\infty^{(X)}$ to be the set of all expressions 
$\sum_x\lambda_x\{x\}$, where all $\lambda_x$ are real, and 
$\sum_x|\lambda_x|\leq1$. (In particular, only countably many of 
$\lambda_x$ are $\neq0$, otherwise this sum would be $+\infty$.) 
In several respects the relationship between $\Zinfty$ and 
$\hat\bbZ_\infty$ is similar to the relationship between 
$\bbZ_p$ and $\hat\bbZ_p$ (cf.~\ptref{sp:mon.as.genr}). 

\nxsubpoint
All of the above monads, maybe with the only exception of $W$ and $W_U$, 
behave themselves very much like ordinary rings. Let's consider an example 
where this is definitely not so. For this consider the monad $\Sigma$ 
defined by the forgetful functor from the category of commutative rings 
to the category of sets. This functor turns out to be monadic, so 
$\catMod\Sigma$ is actually the category of commutative rings. 
Furthermore, $\Sigma(X)=\bbZ[X]$ is the set of polynomials in variables 
from~$X$ with integer coefficients. In particular, $|\Sigma|$ is the 
set $\Sigma[T]$ of polynomials in one variable $T$, and the monoid 
structure on $|\Sigma|$ is given by the application of polynomials:
$f*g=f(g)$. This is definitely not what somebody expects from the 
multiplicative monoid of a ring.

\nxpointtoc{Inner functors}
Given a monad $\Sigma$ over a topos~$\cE$, we would like to be able to 
localize it, i.e.\ to obtain a monad $\Sigma_S$ over $\cE/S$ for any 
$S\in\Ob\cE$. However, this requires some additional structure on functor 
$\Sigma$. We suggest one possible choice of such a structure, sufficient 
for our purposes.

\begin{DefD}\label{def:diag.funct}
Let $P:\cC\to\cD$ be a functor between two cartesian categories 
(i.e.\ categories with finite products). A {\em diagonal structure} 
$(\rho,\theta)$ on $P$ consists of a morphism 
$\theta:e_\cD\to P(e_\cC)$, where $e_\cC$ and $e_\cD$ are the final 
objects of $\cC$ and $\cD$, and a family of morphisms 
$\{\rho_{X,Y}:P(X)\times P(Y)\to P(X\times Y)\}$ 
in $\cD$, 
parametrized by couples of objects $X$, $Y$ of $\cC$, required 
to satisfy the following axioms:
\begin{itemize}
\item[$\rho0)$] 
$\rho_{X,Y}:P(X)\times P(Y)\to P(X\times Y)$ is functorial 
in~$X$ and~$Y$, i.e. for any two morphisms $f:X\to X'$ and $g:Y\to Y'$ 
we have $P(f\times g)\circ\rho_{X,Y}=\rho_{X',Y'}\circ(P(f)\times P(g))$;
\item[$\rho1)$] 
$\xymatrix@C+5pt{P(X)\ar[r]^<>(.5){\theta\times\id_{P(X)}}&P(e_\cC)\times P(X)
\ar[r]^<>(.5){\rho_{e_\cC,X}}&P(e_\cC\times X)\cong P(X)}$ 
is the identity of $P(X)$ 
for any $X\in\Ob\cC$;
\item[$\rho2)$]
For any $X$, $Y$, $Z\in\Ob\cC$ we have 
$\rho_{X\times Y,Z}\circ(\rho_{X,Y}\times\id_{P(Z)})=
\rho_{X,Y\times Z}\circ(\id_{P(X)}\times\rho_{Y,Z})$:
$$\xymatrix{
P(X)\times P(Y)\times P(Z)\ar[rr]^<>(.5){\rho_{X,Y}\times\id_{P(Z)}}
\ar[d]^{\id_{P(X)}\times\rho_{Y,Z}}
&&P(X\times Y)\times P(Z)\ar[d]^{\rho_{X\times Y,Z}}\\
P(X)\times P(Y\times Z)\ar[rr]^{\rho_{X,Y\times Z}}
&&P(X\times Y\times Z)}$$
\end{itemize}
A {\em diagonal (inner) functor $P=(P,\rho,\theta):\cC\to\cD$} is 
a functor $P:\cC\to\cD$ endowed with a diagonal structure. 
A {\em natural transformation} $\xi:(P,\rho,\theta)\to(P',\rho',\theta')$ 
of diagonal functors is simply a natural transformation $\xi:P\to P'$, 
compatible with diagonal structures: $\xi_{e_\cC}\circ\theta=\theta'$ 
and $\xi_{X\times Y}\circ\rho_{X,Y}=\rho'_{X,Y}\circ(\xi_X\times\xi_Y)$.
The category of all diagonal functors between two cartesian categories 
will be denoted by $\catDiagFunct(\cC,\cD)$; if $\cD=\cC$, we also write 
$\catDiagEndof(\cC)$.
\end{DefD}

\nxsubpoint
Any left exact functor $P:\cC\to\cD$ admits a unique diagonal structure, 
given by canonical isomorphisms $P(X\times Y)\cong P(X)\times P(Y)$ 
and $P(e_\cC)\cong e_\cD$. In particular, $\Id_\cC$ can be considered 
as a diagonal endofunctor on~$\cC$.

\nxsubpoint
If $(P,\rho,\theta):\cC\to\cD$ and $(Q,\rho',\theta'):\cD\to\cE$ are 
two diagonal functors, we obtain a diagonal structure $(\rho'',\theta'')$ on 
$QP:\cC\to\cE$ as follows. We put $\theta'':=Q(\theta)\circ\theta':
e_\cE\to QP(e_\cC)$, and $\rho''_{X,Y}:=Q(\rho_{X,Y})\circ
\rho'_{P(X),P(Y)}$. This actually defines a {\em functor\/} 
$\circ:\catDiagFunct(\cD,\cE)\times\catDiagFunct(\cC,\cD)\to 
\catDiagFunct(\cC,\cE)$. Combining together these categories of 
diagonal functors, we obtain a strictly associative 2-category 
with cartesian categories for objects and diagonal functors for 
morphisms. We can replace the 2-category of all categories with this 
new 2-category in our previous considerations, thus obtaining an  
AU $\otimes$-structure on $\catDiagEndof(\cC)$, its $\obslash$-action 
on $\catDiagFunct(\cD,\cC)$ and $\cC$, and its $\oslash$-action on 
$\catDiagFunct(\cC,\cD)$.

\nxsubpoint
Of course, a {\em diagonal (inner) monad $\Sigma$} over a cartesian 
category $\cC$ is simply an algebra $\Sigma=(\Sigma,\mu,\epsilon)$ 
in $\catDiagEndof(\cC)$. Natural transformations 
$\mu:\Sigma^2\to\Sigma$ and $\epsilon:\Id_\cC\to\Sigma$ define a monad 
structure on the underlying endofunctor $\Sigma$, and the category 
$\cC^\Sigma$ is the same for the diagonal monad $\Sigma$ and for its 
underlying usual monad, also denoted by~$\Sigma$.

We see that a diagonal monad $\Sigma$ over $\cC$ is just a usual monad 
$\Sigma=(\Sigma,\mu,\epsilon)$ over $\cC$, together with a 
diagonal structure $(\rho,\theta)$ on $\Sigma$, compatible in a certain sense. 

For example, the compatibility of $\epsilon$ with $\theta$ actually means 
$\theta=\epsilon_{e_\cC}:e_\cC\to\Sigma(e_\cC)$, 
so we don't have to specify $\theta$ separately. Compatibility of 
$\epsilon$ with $\rho$ means $\rho_{X,Y}\circ(\epsilon_X\times\epsilon_Y)=
\epsilon_{X\times Y}$. Next, the compatibility of $\mu$ with 
$\theta=\epsilon_{e_\cC}$ means $\rho_{e_\cC,X}\circ(\epsilon_{e_\cC}\times 
\id_{\Sigma(X)})=\id_{\Sigma(X)}$, and the remaining compatibility of 
$\mu$ with $\rho$ means $\mu_{X\times Y}\circ\Sigma(\rho_{X,Y})\circ
\rho_{\Sigma(X),\Sigma(Y)}=\rho_{X,Y}\circ(\mu_X\times\mu_Y)$:
\begin{equation}\label{eq:comp.mu.rho}
\xymatrix@C+20pt{
\Sigma^2(X)\times\Sigma^2(Y)\ar[d]^{\mu_X\times\mu_Y}
\ar[r]^<>(.5){\rho_{\Sigma(X),\Sigma(Y)}}
&\Sigma(\Sigma(X)\times\Sigma(Y))\ar[r]^<>(.5){\Sigma(\rho_{X,Y})}
&\Sigma^2(X\times Y)\ar[d]^{\mu_{X\times Y}}\\
\Sigma(X)\times\Sigma(Y)\ar[rr]^{\rho_{X,Y}}&&\Sigma(X\times Y)}
\end{equation}

\begin{DefD}\label{def:p-inner.funct}
Given a diagonal functor $(P,\rho,\theta):\cC\to\cD$, we define a 
{\em $P$-inner functor} $F=(F,\alpha):\cC\to\cD$ as follows: 
$F$ is a functor $\cC\to\cD$, and $\alpha$ is a {\em $P$-inner structure 
on~$F$}, i.e.\ a family of morphisms 
$\{\alpha_{X,Y}:P(X)\times F(Y)\to F(X\times Y)\}_{X,Y\in\Ob\cC}$,
subject to the following conditions, similar to $\rho0)$--$\rho2)$ 
of~\ptref{def:diag.funct}:
\begin{itemize}
\item[$\alpha0)$] 
$\alpha_{X,Y}:P(X)\times F(Y)\to F(X\times Y)$ is functorial 
in~$X$ and~$Y$, i.e. for any two morphisms $f:X\to X'$ and $g:Y\to Y'$ 
we have $F(f\times g)\circ\alpha_{X,Y}=\alpha_{X',Y'}\circ(P(f)\times F(g))$;
\item[$\alpha1)$] 
$\xymatrix@C+5pt{F(X)\ar[r]^<>(.5){\theta\times\id_{F(X)}}&P(e_\cC)\times F(X)
\ar[r]^<>(.5){\alpha_{e_\cC,X}}&F(e_\cC\times X)\cong F(X)}$ 
is the identity of $F(X)$ for any $X\in\Ob\cC$;
\item[$\alpha2)$]
For any $X$, $Y$, $Z\in\Ob\cC$ we have 
$\alpha_{X\times Y,Z}\circ(\rho_{X,Y}\times\id_{F(Z)})=
\alpha_{X,Y\times Z}\circ(\id_{P(X)}\times\alpha_{Y,Z})$:
$$\xymatrix{
P(X)\times P(Y)\times F(Z)\ar[rr]^<>(.5){\rho_{X,Y}\times\id_{F(Z)}}
\ar[d]^{\id_{P(X)}\times\alpha_{Y,Z}}
&&P(X\times Y)\times F(Z)\ar[d]^{\alpha_{X\times Y,Z}}\\
P(X)\times F(Y\times Z)\ar[rr]^{\alpha_{X,Y\times Z}}
&&F(X\times Y\times Z)}$$
\end{itemize}
A {\em natural transformation} of $P$-inner functors 
$\zeta:(F,\alpha)\to(F',\alpha')$ is a natural transformation 
$\zeta:F\to F'$, such that $\zeta_{X\times Y}\circ\alpha_{X,Y}=
\alpha'_{X,Y}\circ(\id_{P(X)}\times\zeta_Y):P(X)\times F(Y)\to F'(X\times Y)$.
The category of $P$-inner functors will be denoted by 
$\catInnFunct_{/P}(\cC,\cD)$.
\end{DefD}

\nxsubpoint 
Given a $P$-inner functor $(F,\alpha):\cC\to\cD$ and a 
$Q$-inner functor $(G,\alpha'):\cD\to\cE$, we have a canonical 
$QP$-inner structure $\alpha''$ on $GF:\cC\to\cE$, given by 
$\alpha''_{X,Y}=G(\alpha_{X,Y})\circ\alpha'_{P(X),F(Y)}:
QP(X)\times GF(Y)\to GF(X\times Y)$. In this way we construct the composition 
functors $\circ:\catInnFunct_{/Q}(\cD,\cE)\times\catInnFunct_{/P}(\cC,\cD)
\to\catInnFunct_{/QP}(\cC,\cE)$.

\nxsubpoint
Notice that $P=(P,\rho,\theta)$ 
is itself a $P$-inner functor, since the axioms 
$\alpha0)$--$\alpha2)$ specialize to $\rho0)$--$\rho2)$ when we put 
$F:=P$, $\alpha:=\rho$. Hence the statements proved for inner functors 
over a diagonal functor can be transferred to diagonal inner functors 
themselves. 

\nxsubpoint\label{sp:def.plain.monad} 
If $\cC=\cD$, we can take $P=\Id_\cC$, thus obtaining the 
{\em category of plain inner endofunctors\/ $\catInnEndof(\cC)$}, which 
consists of endofunctors $F:\cC\to\cC$, endowed with a plain inner 
structure $\alpha$, $\alpha_{X,Y}:X\times F(Y)\to F(X\times Y)$. 
We have an AU $\otimes$-structure on $\catInnEndof(\cC)$; the algebras with 
respect to this $\otimes$-structure will be called 
{\em (plain) inner monads.} 
Clearly, a plain inner monad $\Sigma=(\Sigma,\mu,\epsilon,\alpha)$ 
is simply a usual monad $(\Sigma,\mu,\epsilon)$ over~$\cC$, together with 
some plain inner (i.e.\ $\Id_\cC$-inner) structure $\alpha$ on $\Sigma$, 
subject to some compatibility conditions similar to those considered before. 
Namely, we must have 
\begin{gather}
\label{eq:plain.mon1}
\epsilon_{X\times Y}=\alpha_{X,Y}\circ(\id_X\times\epsilon_Y) \\
\label{eq:plain.mon2}
\mu_{X\times Y}\circ\Sigma(\alpha_{X,Y})\circ\alpha_{X,\Sigma(Y)}=
\alpha_{X,Y}\circ(\id_X\times\mu_Y):X\times\Sigma^2(Y)\to\Sigma(X\times Y)
\end{gather}

\begin{PropD}\label{prop:uniq.plain.str}
Any endofunctor $F$ over $\catSets$ admits a unique plain inner structure, 
i.e.\ the $\otimes$-categories $\catInnEndof(\catSets)$ and 
$\catEndof(\catSets)$ are isomorphic. In particular, any monad over 
$\catSets$ uniquely extends to a plain inner monad.
\end{PropD}
\begin{Proof}
For any set $X$ and any element $x\in X$ denote by 
$\tilde{x}:e_{\catSets}=\st1=\{1\}\to X$ the only map with image~$x$. 
Then $\tilde{x}\times\id_Y:Y\to X\times Y$ is the map $y\mapsto (x,y)$ for any 
set~$Y$. A plain inner structure $\{\alpha_{X,Y}:
X\times F(Y)\to F(X\times Y)\}$ 
is functorial with respect to $X$ and $Y$; applying this to $\tilde{x}$ 
and taking $\alpha1)$ into account we see that 
$\alpha_{X,Y}\circ(\tilde{x}\times\id_{F(Y)})=F(\tilde{x}\times\id_Y):
F(Y)\to F(X\times Y)$, 
i.e.\ $\alpha_{X,Y}(x,w)=(F(\tilde{x}\times\id_Y))(w)$ for any $x\in X$, 
$w\in F(Y)$. This shows the uniqueness of $\alpha$, and also gives 
a way of defining it. It is immediate that the family $\{\alpha_{X,Y}\}$ 
defined by this formula indeed satisfies $\alpha0)$--$\alpha2)$.
\end{Proof}

\nxsubpoint\label{sp:plain.str.words}
For example, the plain inner structure on the monad of words $W$ is given 
by $\alpha_{X,Y}(x,y_1y_2\ldots y_n)=\{x,y_1\}\{x,y_2\}\ldots\{x,y_n\}$, 
where we write $\{x,y\}$ instead of $\{(x,y)\}$. For the monad $W_U$ 
of words with constants from~$U$ the situation is similar, but the 
constants are not affected, e.g.\ 
$\alpha_{X,Y}(x,y_1u_1y_2y_3u_2)=\{x,y_1\}u_1\{x,y_2\}\{x,y_3\}u_2$.

\nxsubpoint\label{sp:gamma.descr.inn} 
(Alternative description of inner structures.) 
A $P$-inner structure $\{\alpha_{X,Y}:P(X)\times F(Y)\to F(X\times Y)\}$ 
on a functor $F:\cC\to\cD$ admits an alternative description in terms 
of a natural transformation $\gamma$ of 
functors $\cC^0\times\cC^0\times\cC\to\catSets$, 
$\gamma_{X,Y,Z}:\Hom_\cC(X\times Y, Z)\to\Hom_\cD(P(X)\times F(Y), F(Z))$. 
Namely, for any $\phi:X\times Y\to Z$ we put 
$\gamma_{X,Y,Z}(\phi):=F(\phi)\circ\alpha_{X,Y}$, and conversely, 
any such $\gamma$ defines (by Yoneda) a family $\alpha_{X,Y}$, functorial 
in~$X$ and~$Y$; actually, 
$\alpha_{X,Y}=\gamma_{X,Y,X\times Y}(\id_{X\times Y})$. Moreover, 
the axioms $\alpha0)$--$\alpha2)$ can be rewritten directly in terms 
of~$\gamma$:
\begin{itemize}
\item[$\gamma0)$] 
$\gamma_{X,Y,Z}:\Hom_\cC(X\times Y, Z)\to\Hom_\cD(P(X)\times F(Y), F(Z))$ 
is functorial in $X$, $Y$ and $Z\in\Ob\cC$;
\item[$\gamma1)$]
The following diagram is commutative for any~$X$ and $Y$ in $\cC$:
\begin{equation}
\xymatrix@C+5pt{
\Hom_\cC(e_\cC\times X,Y)\ar[d]^{\sim}\ar[r]^<>(.5){\gamma_{e_\cC,X,Y}}
&\Hom_\cD(P(e_\cC)\times F(X),F(Y))\ar[d]^{(\theta\times\id_{F(X)})^*}\\
\Hom_\cC(X,Y)\ar[r]^{F(\cdot)}&\Hom_\cD(F(X),F(Y))}
\end{equation}
\item[$\gamma2)$]
For any morphisms $\phi:U\times X\to Y$ and $\psi:V\times Y\to Z$ in~$\cC$ 
define $\chi:V\times U\times X\to Z$ by $\chi:=\psi\circ(\id_V\times\phi)$. 
Then $\gamma_{V\times U,X,Z}(\chi)\circ(\rho_{V,U}\times\id_{F(X)})=
\gamma_{V,Y,Z}(\psi)\circ(\id_{P(V)}\times\gamma_{U,X,Y}(\phi))$:
\begin{equation}
\xymatrix@C+12pt{
P(V)\times P(U)\times F(X)\ar[rr]^<>(.5){\id_{P(V)}\times\gamma_{U,X,Y}(\phi)}
\ar[d]^{\rho_{V,U}\times\id_{F(X)}}
&&P(V)\times F(Y)\ar[d]^{\gamma_{V,Y,Z}(\psi)}\\
P(V\times U)\times F(X)\ar[rr]^<>(.5){\gamma_{V\times U,X,Z}(\chi)}&&F(Z)}
\end{equation}
\end{itemize}
The verification of the equivalence of $\alpha0)$--$\alpha2)$ and 
$\gamma0)$--$\gamma2)$ doesn't present any difficulties. For example, 
$\alpha2)$ is obtained from $\gamma2)$ by setting $\phi=\id$, $\psi=\id$, 
and the converse follows from the explicit formula for $\gamma$ and 
the functoriality of~$\alpha$.

\nxsubpoint (Inner structures over cartesian closed categories.)
Now suppose that both $\cC$ and $\cD$ are {\em cartesian closed,} 
i.e.\ they are cartesian, and the functor $\Hom(-\times Y,Z)$ is representable 
by some object $\iHom(Y,Z)$ ({\em local\/} or {\em inner\/} Hom), for 
any couple of objects $Y$ and $Z$ from $\cC$ or from~$\cD$. The canonical 
bijection $\Hom(X\times Y,Z)\simto\Hom(X,\iHom(Y,Z))$ will be denoted 
$\phi\mapsto\phi^\flat$, and its inverse $\psi\mapsto\psi^\sharp$. 
We also have the {\em evaluation morphisms\/} 
$\ev_{Y,Z}:=(\id_{\iHom(Y,Z)})^\sharp:\iHom(Y,Z)\times Y\to Z$, 
and the {\em functor of global sections\/} 
$\Gamma_\cC(-):=\Hom_\cC(e_\cC,-):\cC\to\catSets$, 
and similarly for~$\cD$. Clearly, 
$\Gamma_\cC(\iHom_\cC(Y,Z))\cong\Hom_\cC(Y,Z)$.

\nxsubpoint 
Consider now a $P$-inner structure $\alpha$ on some functor 
$F:\cC\to\cD$. We know already that it admits a description in terms of 
morphisms $\gamma_{X,Y,Z}:\Hom_\cC(X\times Y,Z)
\to\Hom_\cD(P(X)\times F(Y),F(Z))$, functorial in $X$, $Y$ and~$Z$ 
(cf.~\ptref{sp:gamma.descr.inn}). Since we have assumed $\cC$ and $\cD$ to be 
cartesian closed, we can interpet $\gamma_{X,Y,Z}$ as a morphism 
$\Hom_\cC(X,\iHom_\cC(Y,Z))\to\Hom_\cD(P(X),\iHom_\cD(F(Y),F(Z))$. 
By Yoneda such a functorial family of morphisms is determined by its 
value on $\id_{\iHom_\cC(Y,Z)}$, so we put $\beta_{Y,Z}:=
\gamma_{\iHom(Y,Z),Y,Z}(\ev_{Y,Z})^\flat:P(\iHom_\cC(Y,Z))\to 
\iHom_\cD(F(Y),F(Z))$. 
This gives another equivalent description 
of a $P$-inner structure on a functor between two cartesian closed categories,
and in fact it is possible to restate the axioms in terms of these morphisms 
$\beta$. For example, $\beta0)$ states that 
$\beta_{Y,Z}:P(\iHom_\cC(Y,Z))\to\iHom_\cD(F(Y),F(Z))$ is functorial in
$Y$ and~$Z$, axiom $\beta1)$ relates $\beta$ with the global section 
functors, and $\beta2)$ relates $\beta$ with the composition morphisms
$o_{X,Y,Z}:\iHom(Y,Z)\times\iHom(X,Y)\to\iHom(X,Z)$. Actually, it 
states $\beta_{X,Z}\circ P(o_{X,Y,Z})\circ\rho_{\iHom(Y,Z),\iHom(X,Y)}=
o_{F(X),F(Y),F(Z)}\circ(\beta_{Y,Z}\times\beta_{X,Y})$. We omit the 
verification of equivalence of this set of axioms $\beta0)$--$\beta2)$ with 
our other sets of axioms. Let's just remark that this description 
shows that an inner structure can be thought of as an extension of 
$F(-):\Hom_\cC(Y,Z)\to\Hom_\cD(F(Y),F(Z))$ to the local Homs. 
This also shows again that in case $\cC=\cD=\catSets$ there is always a 
unique plain inner structure, since in this case local Homs coincide with 
the usual ones.

\nxsubpoint
Here is one application of (plain) inner endofuctors: 
{\em Any plain inner endofunctor $(F,\alpha)$ over a category $\cC$ with 
finite projective limits defines a family of endofunctors 
$F_{/S}:\cC/S\to\cC/S$, parametrized by objects $S$ of~$\cC$.}
For any object $p:X\to S$ of $\cC/S$ we construct $F_{/S}(p)\in\Ob\cC/S$ as 
follows. We consider the graph $\Gamma'_p=(p,\id_X):X\to S\times X$, and 
define $F_{/S}(p)$
from the first line of the following diagram, involving a cartesian square:
\begin{equation}
\xymatrix{F_{/S}(p)\ar@{-->}[r]\ar@{-->}[d]
&S\times F(X)\ar[r]^<>(.5){\pr_1}\ar[d]^{\alpha_{S,X}}&S\\
F(X)\ar[r]^<>(.5){F(\Gamma'_p)}&F(S\times X)
}
\end{equation}
This construction is functorial in~$F=(F,\alpha)$ and $S$, but in general 
it doesn't respect composition of functors, hence it cannot be used for 
the localization of monads. Notice that we don't even obtain a natural 
plain inner structure on functors~$F_{/S}$. However, we'll see in the 
next chapter that this construction behaves very nicely when $\cC=\cE$ is a 
topos and $F$ is an algebraic inner endofunctor or monad over~$\cE$; 
this is exactly the case we need.

\nxsubpoint\label{sp:sigma.str.hom}
Here is another useful application of plain inner monads. Given a monad 
$\Sigma$ over $\catSets$, any $\Sigma$-object $M=(M,\sigma)$ and 
any set $I$, the local Hom $\iHom(I,M)$ equals $M^I$, hence it has a 
canonical $\Sigma$-structure, namely, the product $\Sigma$-structure 
(cf.~\ptref{sp:submod.proj.lim}). We want to generalize this to the case 
of a monad $\Sigma$ over a cartesian closed category $\cC$, an object 
$M=(M,\sigma)\in\Ob\cC^\Sigma$, and any $I\in\Ob\cC$. However, to obtain 
a $\Sigma$-structure on $H:=\iHom_\cC(I,M)$ we need $\Sigma$ to be 
a plain inner monad, i.e.\ we need a compatible plain inner structure 
$\alpha$ on~$\Sigma$. In this case we define a $\Sigma$-structure 
$\tau:\Sigma(H)\to H$ by the following diagram:
\begin{equation}
\xymatrix{
I\times\Sigma(H)\ar[r]^{\alpha_{I,H}}&\Sigma(I\times H)
\ar[r]^<>(.5){\Sigma(\ev')}&\Sigma(M)\ar[d]^{\sigma}\\
\Sigma(H)\times I\ar[u]^{\sim}\ar[rr]^{\tau^\sharp}&&M}
\end{equation}
Here $\ev':I\times H=I\times\iHom(I,M)\to M$ is essentially the evaluation 
map (up to a permutation of arguments), and the left vertical arrow is 
the natural symmetry.

It is easy to check that in case $\cC=\catSets$, when any monad $\Sigma$ 
admits a unique compatible plain inner structure 
(cf.\ \ptref{prop:uniq.plain.str}), we obtain again the product 
$\Sigma$-structure on $\iHom(I,M)=M^I$. In the general case we have to 
check that $\tau$ is indeed a $\Sigma$-structure on $H$, i.e.\ that 
$\tau\circ\epsilon_H=\id_H$ and $\tau\circ\Sigma(\tau)=
\tau\circ\mu_H$. Both these statements follow from the definition of $\tau$ 
and from the compatibility conditions for $\alpha$ recalled 
in~\ptref{sp:def.plain.monad}; we suggest the reader to skip 
their verification presented below.

a) First of all, denote by $\alpha'_X:\Sigma(X)\times I\to\Sigma(X\times I)$ 
the morphism deduced from $\alpha_{I,X}:I\times\Sigma(X)\to\Sigma(I\times X)$ 
by permuting the factors in direct products. Then 
$\tau^\sharp=\sigma\circ\Sigma(\ev)\circ\alpha'_H:\Sigma(H)\times I\to M$, 
where $\ev=\ev_{H,I}=\id_H^\sharp:H\times I\to M$ is the evaluation map.
Now we compute $(\tau\circ\epsilon_H)^\sharp=
\tau^\sharp\circ(\epsilon_H\times\id_I)=
\sigma\circ\Sigma(\ev)\circ\alpha'_H\circ(\epsilon_H\times\id_I)$; 
according to~\eqref{eq:plain.mon1}, this equals 
$\sigma\circ\Sigma(\ev)\circ\epsilon_{H\times I}=\sigma\circ\epsilon_M\circ\ev=
\ev=(\id_H)^\sharp$. This proves $\tau\circ\epsilon_H=\id_H$.

b) Notice that $\ev\circ(\tau\times\id_I)=\tau^\sharp$. Now compute 
$(\tau\circ\Sigma(\tau))^\sharp=\tau^\sharp\circ(\Sigma(\tau)\times\id_I)=
\sigma\circ\Sigma(\ev)\circ\alpha'_H\circ(\Sigma(\tau)\times\id_I)=
\sigma\circ\Sigma(\ev)\circ\Sigma(\tau\times\id_I)\circ\alpha'_{\Sigma(H)}=
\sigma\circ\Sigma(\tau^\sharp)\circ\alpha'_{\Sigma(H)}=
\sigma\circ\Sigma(\sigma)\circ\Sigma^2(\ev)\circ\Sigma(\alpha'_H)\circ
\alpha'_{\Sigma(H)}$. Now we use~\eqref{eq:plain.mon2} together with 
$\sigma\circ\Sigma(\sigma)=\sigma\circ\mu_M$, thus obtaining 
$\sigma\circ\mu_M\circ\Sigma^2(\ev)\circ\Sigma(\alpha'_H)\circ
\alpha'_{\Sigma(H)}=\sigma\circ\Sigma(\ev)\circ\mu_{H\times I}\circ
\Sigma(\alpha'_H)\circ\alpha'_{\Sigma(H)}=\sigma\circ\Sigma(\ev)\circ
\alpha'_H\circ(\mu_H\times\id_I)=\tau^\sharp\circ(\mu_H\times\id_I)=
(\tau\circ\mu_H)^\sharp$. This proves $\tau\circ\Sigma(\tau)=
\tau\circ\mu_H$, hence $(H,\mu_H)\in\Ob\cC^\Sigma$.

\nxsubpoint\label{sp:diag.from.plain}
Given any diagonal inner structure $\rho$ on a monad $\Sigma$ over $\cC$, 
we can retrieve a plain inner structure $\alpha$ on the same monad by putting 
$\alpha_{X,Y}:=\rho_{X,Y}\circ(\epsilon_X\times\id_{\Sigma(Y)})$. 
Conversely, we might try to extend 
a plain inner structure $\alpha$ to a diagonal inner structure $\rho$ by 
defining $\rho_{X,Y}:\Sigma(X)\times\Sigma(Y)\to\Sigma(X\times Y)$ from 
the following commutative diagram:
\begin{equation}
\xymatrix@C+4pt{
\Sigma(X)\times\Sigma(Y)\ar[d]^{\sim}
\ar[rrr]^{\rho_{X,Y}}&&&\Sigma(X\times Y)\\
\Sigma(Y)\times\Sigma(X)\ar[r]^{\alpha_{\Sigma(Y),X}}
&\Sigma(\Sigma(Y)\times X)\ar[r]^{\sim}
&\Sigma(X\times\Sigma(Y))\ar[r]^<>(.5){\Sigma(\alpha_{X,Y})}
&\Sigma^2(X\times Y)\ar[u]_{\mu_{X\times Y}}
}
\end{equation}

Strictly speaking, we have to check that these formulas indeed define some 
plain inner structure $\alpha$ starting from 
a diagonal inner structure $\rho$,
and study the exact conditions under which a plain inner monad extends to 
a diagonal inner monad in this manner; it is true in general that we 
always obtain a diagonal inner structure on the underlying endofunctor 
of~$\Sigma$, which need not be compatible with the monad structure. 
The verifications are straightforward but lengthy; we 
omit them and consider an example instead.

Let $\cC=\catSets$, $\Sigma=W$ be the monad of words (cf.\ 
\ptref{sp:mon.words}), and consider the only compatible plain inner 
structure $\alpha$ on~$W$ (cf.~\ptref{prop:uniq.plain.str}); 
$\alpha_{X,Y}:X\times W(Y)\to W(X\times Y)$ is given by 
$(x, y_1y_2\ldots y_n)\mapsto\{x,y_1\}\{x,y_2\}\ldots\{x,y_n\}$ 
(cf.~\ptref{sp:plain.str.words}). Then the corresponding diagonal inner 
structure $\rho_{X,Y}:W(X)\times W(Y)\to W(X\times Y)$ maps 
a pair of words $(x_1x_2\ldots x_m, y_1y_2\ldots y_n)$ into the word 
$\{x_1,y_1\}\{x_1,y_2\}\ldots\{x_1,y_n\}\{x_2,y_1\}\ldots\{x_m,y_n\}$, 
i.e.\ we list all pairs of letters $(x_i,y_j)$ of the two original words 
in the lexicographical order with respect to the pair of indices $(i,j)$. 
Notice that this is indeed a diagonal inner structure on the underlying 
endofunctor~$W$, but it is {\em not\/} compatible with the monad 
structure. For example, put $X:=\{a,b\}$, $Y:=\{x,y,z\}$ and 
consider the two maps $W^2(X)\times W^2(Y)\to W(X\times Y)$ which 
have to be equal by~\eqref{eq:comp.mu.rho}. The images of 
$(\{\{a\}\{b\}\},\{\{x\}\{y\}\}\{\{z\}\})\in W^2(X)\times W^2(Y)$ 
under these two maps are distinct: 
$\{a,x\}\{a,y\}\{a,z\}\{b,x\}\{b,y\}\{b,z\}\neq
\{a,x\}\{a,y\}\{b,x\}\{b,y\}\{a,z\}\{b,z\}$.

\nxsubpoint\label{sp:comm.diag.str}
Notice that the diagonal inner structure $\rho$ on endofunctor~$\Sigma$ 
constructed above in 
general need not be {\em commutative}, i.e.\ the following diagram in general 
does {\em not\/} commute:
\begin{equation}
\xymatrix{
\Sigma(X)\times\Sigma(Y)\ar[r]^<>(.5){\rho_{X,Y}}\ar[d]^{\sim}&
\Sigma(X\times Y)\ar[d]^{\sim}\\
\Sigma(Y)\times\Sigma(X)\ar[r]^<>(.5){\rho_{Y,X}}&\Sigma(Y\times X)}
\end{equation}

It is easy to see that the family $\rho'_{X,Y}:\Sigma(X)\times\Sigma(Y)\to 
\Sigma(X\times Y)$, obtained from $\rho$ by permuting factors in direct 
products, is another diagonal structure on~$\Sigma$, also restricting to 
the same plain inner structure $\alpha$.

We say that a diagonal (resp.\ plain) inner monad is {\em commutative\/} if $\rho'=\rho$ (resp.\ if this is true for diagonal inner structure constructed 
above), i.e.\ if this diagram commutes. Since any monad over $\catSets$ 
extends uniquely to a plain inner monad, the notion of commutativity makes 
sense in this context as well. For example, the monad of words is 
not commutative, while the monad $\Sigma_R$ defined by an associative 
ring~$R$ is commutative iff $R$ is commutative, since 
in this case $\rho_{X,Y}$ maps $\bigl(\sum_i\lambda_i\{x_i\}$, 
$\sum_j\mu_j\{y_j\}\bigr)$ into $\sum_{i,j}\lambda_i\mu_j\{x_i,y_j\}$.

An interesting statement is that 
{\em if the diagonal inner structure~$\rho$ defined by some plain inner 
structure~$\alpha$ on a monad $\Sigma$ is commutative, then it is 
compatible with the monad structure of~$\Sigma$}. It seems to be known 
in similar contexts (cf.~\ptref{sp:inn.mon.terms} below), so 
we omit the proof. We won't use it anyway.

\nxsubpoint\label{sp:gen.locinnhom}
Let $\Sigma=(\Sigma,\mu,\epsilon,\alpha)$ be a plain inner monad over 
a cartesian closed category~$\cC$, and let $M=(M,\sigma)$ and 
$N=(N,\tau)$ be two objects of $\cC^\Sigma$. Then we can construct 
``local inner Homs'' $\iHom_{\cC,\Sigma}(M,N)=\iHom_\Sigma(M,N)\subset
\iHom(M,N)$, such that $\Gamma_\cC(\iHom_\Sigma(M,N))=\Hom_\Sigma(M,N)$. 
To do this consider an arbitrary morphism $i:H'\to H:=\iHom(M,N)$, 
which corresponds by adjointness to some $i^\sharp: H'\times M\to N$. 
We say that {\em $H'$ acts by $\Sigma$-morphisms from $M$ to $N$} if 
the following diagram commutes:
\begin{equation}
\xymatrix{
H'\times\Sigma(M)\ar[d]^{\alpha_{H',M}}\ar[rr]^{\id_{H'}\times\sigma}
&&H'\times M\ar[d]^{i^\sharp}\\
\Sigma(H'\times M)\ar[r]^<>(.5){\Sigma(i^\sharp)}&\Sigma(N)\ar[r]^{\tau}&N}
\end{equation}
In particular, we can take $H'=H$, $i=\id_H$ and obtain two morphisms 
$H\times\Sigma(M)\rightrightarrows N$, which correspond by adjointness 
to some $p,q:H\rightrightarrows\iHom(\Sigma(M),N)$. 
In general the above diagram commutes 
for some $i:H'\to H$ iff $p\circ i=q\circ i$, i.e.\ iff $i$ 
factorizes through the kernel $H_0=\iHom_\Sigma(M,N):=\Ker(p,q)\subset H$. 
In other words, $H_0$ is the largest subobject of $\iHom(M,N)$ which 
acts by $\Sigma$-homomorphisms, so it is a natural candidate for 
$\iHom_\Sigma(M,N)$. Notice that $\Gamma_\cC(H_0)$ corresponds to 
morphisms $i:e_\cC\to\iHom(M,N)$ which factorize through $H_0$, i.e.\ 
which make the above diagram commutative; using $\alpha1)$ we see that 
this is equivalent to $i^\sharp:M\to N$ being a $\Sigma$-homomorphism, i.e.\ 
$\Gamma_\cC(\iHom_\Sigma(M,N))\cong\Hom_\Sigma(M,N)$.

Of course, one can check directly that the composition morphisms 
$o_{M,N,P}:\iHom(N,P)\times\iHom(M,N)\to\iHom(M,P)$ respect these subobjects 
$\iHom_\Sigma$ once we fix $\Sigma$-structures on $M$, $N$ and $P$; 
we'll prove this later for algebraic monads over a topos in another way.

Notice that in general $\iHom_\Sigma(M,N)$ doesn't inherit a 
$\Sigma$-structure from $\iHom(M,N)$ (cf.~\ptref{sp:sigma.str.hom}), i.e.\ 
it is a subobject of $\iHom(M,N)$, but not necessarily a $\Sigma$-stable 
subobject. Actually, we'll prove later that $\iHom_\Sigma(M,N)$ is a 
$\Sigma$-subobject of $\iHom(M,N)$ for all choices of $M$ and $N$ from 
$\cC^\Sigma$ iff $\Sigma$ is commutative. 

\nxsubpoint\label{sp:inn.mon.terms}
A terminological remark: if we consider on a cartesian category $\cC$ 
the ACU $\otimes$-structure defined by the direct product 
($A\otimes B:=A\times B$), then our notion of plain inner structure 
corresponds (up to some minor details) 
to what people call {\em tensorial strength}, and 
plain inner monads -- to {\em strong monads}. Furthermore, our 
diagonal inner monads correspond to {\em monoidal\/} (or {\em tensorial\/}) 
{\em monads}, and our commutativity condition seems to be known either under 
the name of {\em commutativity\/} or {\em symmetry\/} condition in this 
setup. 


\cleardoublepage

\mysection{Algebraic monads and algebraic systems}

This chapter is dedicated to the study of a very important class of 
endofunctors and monads --- the {\em algebraic\/} endofunctors and 
monads over the category of sets, and their topos counterparts --- 
algebraic (plain inner) endofunctors and monads over a (Grothendieck) 
topos~$\cE$. These notions admit both category-theoretic descriptions 
(in terms of some functors and monads) and algebraic descriptions 
(as sets or sheaves with some operations on them, in the spirit 
of universal algebra). Both these approaches are important: for example, 
the first of them enables us to construct (arbitrary) projective limits 
of algebraic monads, while the second is equally indispensable for the 
construction of some inductive limits (e.g.\ tensor products). That's 
why we carefully develop some language and some notation to simplify 
transition between these two alternative descriptions. 

Let us emphasize once more that the notion of an algebraic monad 
(resp.\ algebraic inner monad over a topos) is extremely important 
in the sequel since in our setup it corresponds to the notion of a
``(non-commutative) generalized ring'' (or a sheaf of such generalized rings); 
in this sense the content of this chapter should be viewed as some 
generalization of the basic theory of (sheaves of) 
associative rings and algebras.

{\bf Notations.} All endofunctors and monads considered in this chapter, 
unless otherwise specified,  
will be defined over a base category $\cC$, 
which will be supposed to be equal either to $\catSets$, or to a topos~$\cE$. 
Sometimes $\cE$ is the topos of sheaves (of sets) on some site $\cS$; then 
we write $\cE=\tilde\cS$. We also denote by $\cA$ the AU 
$\otimes$-category $\catEndof(\cC)$.

\nxpointtoc{Algebraic endofunctors on $\catSets$}
The definition itself is very simple:
\begin{DefD}\label{def:alg.endof}
We say that an endofunctor $\Sigma:\catSets\to\catSets$ is 
{\em algebraic} if it commutes with filtered inductive limits. The 
full $\otimes$-subcategory of $\cA=\catEndof(\catSets)$, 
consisting of all algebraic endofunctors, 
will be denoted $\catEndofalg(\catSets)$ or $\cA_{alg}$. 
Thus an {\em algebraic monad} is simply an algebra in $\cA_{alg}$, i.e.\ 
a monad over $\catSets$, the underlying functor of which commutes with 
arbitrary inductive limits.
\end{DefD}
Clearly, if $F$ and $G$ commute with filtered inductive limits, the same 
is true for $F\otimes G=F\circ G$, so $\cA_{alg}$ is indeed a 
$\otimes$-subcategory of $\cA$.

\nxsubpoint
Notice that any set~$X$ is a filtered inductive limit of the ordered set 
of all its finite subsets $I\subset X$, hence 
an algebraic endofunctor~$\Sigma$ is completely determined by its values 
on finite sets: $\Sigma(X)=\injlim_{\text{finite $I\subset X$}}\Sigma(I)$. 
On the other hand, any finite set is isomorphic to some 
{\em standard finite set} $\stn=\{1,2,\ldots,n\}$, $n\geq0$. Let us denote 
by $\catN$ the full subcategory of $\catSets$, consisting of all 
standard finite sets. Sometimes we identify the objects of $\catN$ with 
the set of non-negative integers, and write $n$ instead of $\stn$. Let 
us denote by $J=J_\catSets:\catN\to\catSets$ the inclusion functor.

We see that $J$ is fully faithful, and that its essential image is the 
category of all finite sets. Since any algebraic endofunctor $\Sigma$ is 
completely determined by its restriction to the category of finite sets, 
which is equivalent to $\catN$, we see that the canonical functor 
$J^*:\cA\to\catFunct(\catN,\catSets)=\catSets^\catN$, given by 
$F\mapsto F\circ J$, induces an equivalence between $\cA_{alg}$ and a 
full subcategory of $\catSets^\catN$. Let's make a more precise statement:

\begin{PropD}\label{prop:j.descr.algfun}
Consider the inclusion functor $J:\catN\to\catSets$, the corresponding 
restriction functor $J^*:\cA=\catEndof(\catSets)\to\catSets^\catN$, 
$F\mapsto F\circ J$, and its left adjoint (i.e.\ left Kan extension of~$J$) 
$J_!:\catSets^\catN\to\cA$, defined by the usual formula 
$(J_!G)(X)=\injlim_{\catN/X}G(\stn)$, where $\catN/X$ is the category of all 
maps from standard finite sets into a set~$X$.
Then $J_!$ is fully faithful, and its essential 
image equals $\cA_{alg}$, hence it induces an equivalence between 
$\catSets^\catN$ and $\cA_{alg}$, and $J^*|_{\cA_{alg}}$ is the 
quasi-inverse equivalence. Moreover, $J^*$ commutes with arbitrary limits 
of functors, and $J_!$ commutes with arbitrary inductive and finite 
projective limits, hence $\cA_{alg}\subset\cA$ is stable under these 
types of limits in~$\cA$.
\end{PropD}
\begin{Proof}
The existence of $J_!$ follows from the existence of inductive limits of sets, 
and the explicit formula for $(J_!G)(X)$ is well-known and can be either 
checked directly or found in SGA~4~I, for example. Full faithfulness 
of~$J_!$ is a consequence of that of~$J$, since $(J_!G)(\stm)$ is computed 
as the limit of $G(\stn)$ along the category $\catN/\stm$ of all maps 
$\stn\to\stm$, which has a final object $\id_\stm$, hence 
$(J_!G)(\stm)\cong G(\stm)$, i.e.\ $G\to J^*J_!G$ is an isomorphism, hence 
$J_!$ is fully faithful. Notice that $\catN/X$ is a filtered category, 
and $\injlim_{\catN/X}\stn\cong X$, 
hence for any algebraic endofunctor~$F$ we have 
$F(X)\cong \injlim_{\catN/X} F(\stn)=(J_!J^*F)(X)$, i.e.\ any algebraic 
functor lies indeed in the essential image of~$J_!$. To prove the converse 
we use the alternative description of $(J_!G)(X)$ in terms of inductive 
limits and finite products, given in Lemma~\ptref{l:alt.descr.algfun} below, 
and the fact that filtered inductive limits in $\catSets$ commute with 
finite projective limits; this shows immediately that $J_!G$ is indeed 
algebraic. This fact also shows that $J_!$ commutes with finite projective 
limits, and it has to commute with arbitrary inductive limits since 
it has a right adjoint; this proves the remaining statements.
\end{Proof}

\begin{LemmaD}\label{l:alt.descr.algfun}
Given any functor $G:\catN\to\catSets$, its left Kan extension 
$J_!G:\catSets\to\catSets$ can be computed as follows: 
$(J_!G)(X)$ is the cokernel of the pair of maps $p$, $q$ from 
$H_1(X):=\bigsqcup_{\phi:\stm\to\stn} G(\stm)\times X^n$ to 
$H_0(X):=\bigsqcup_{n\geq0} G(\stn)\times X^n$, defined by the following 
two requirements:
\begin{itemize}
\item[a)] The restriction of $p$ to the component 
of $H_1(X)$ indexed by $\phi:\stm\to\stn$ is the map 
$\id_{G(\stm)}\times X^\phi: G(\stm)\times X^n\to G(\stm)\times X^m$, 
where $X^\phi=\iHom(\phi,\id_X):\iHom(\stn,X)=X^n\to\iHom(\stm,X)=X^m$ is 
the canonical map $(x_1,\ldots, x_n)\mapsto (x_{\phi(1)},\ldots,x_{\phi(m)})$.
\item[b)] The restriction of $q$ to the same component is the map 
$G(\phi)\times\id_{X^n}:G(\stm)\times X^n\to G(\stn)\times X^n$.
\end{itemize}
\end{LemmaD}
\begin{Proof}
Recall that {\em any\/} inductive limit $\injlim_\cI G$ of a functor 
$G:\cI\to\cC$ can be written as the cokernel of two maps $p$ and $q$ 
from $\bigsqcup_{(\phi:s\to t)\in\Ar\cI} G(s)$ to 
$\bigsqcup_{s\in\Ob\cI} G(s)$, 
the restrictions of which to the component indexed by some $\phi:s\to t$
is the identity morphism of $G(s)$ (resp.\ the morphism 
$G(\phi):G(s)\to G(t)$).
Now we apply this for the computation of $(J_!G)(X)=\injlim_{\catN/X}G(\stn)$. 
In this case $\cI=\catN/X$, $\Ob\cI=\{$maps $\psi:\stn\to X\}=
\sqcup_{n\geq0}X^n=W(X)$, and $\Ar\cI=\{\stm\stackrel{\phi}\to\stn
\stackrel{\psi}\to X\}=\sqcup_{\phi:\stm\to\stn}X^n$. Taking the required 
sums along these index sets yields the objects $H_0(X)$ and $H_1(X)$ 
described in the statement of the lemma, and 
$p,q:H_1(X)\rightrightarrows H_0(X)$ 
also turn out to coincide with those mentioned in the statement.
\end{Proof}

\nxsubpoint\label{sp:max.alg.subf}
Notice that the injective maps $\stn\to X$ form a cofinal set in $\catN/X$, 
hence it is sufficient to take the inductive limit along only such maps 
in order to compute $(J_!G)(X)$. 
In particular, for any endofunctor $F$ we have 
$(J_!J^*F)(X)=\injlim_{\text{finite }I\subset X}F(I)\to F(X)$. If $X$ is 
finite, this map is an isomorphism; in general it is always {\em injective}, 
i.e.\ $J_!J^*F\to F$ is a {\em monomorphism}. This follows from the fact that 
any monomorphism with non-empty source splits in $\catSets$, hence all 
maps $F(I)\to F(X)$ for non-empty finite subsets $I\subset X$ are injective, 
hence this is true for the filtered inductive limit of the $F(I)$. 
Let us denote $J_!J^*F$ by $F_{alg}$. Clearly, $F_{alg}\subset F$ 
is the largest algebraic subfunctor of $F$, and any morphism 
$F'\to F$ from an algebraic $F'$ factorizes through~$F_{alg}$. We also 
know that $F\mapsto F_{alg}$ commutes with arbitrary inductive and 
finite projective limits (cf.~\ptref{prop:j.descr.algfun}). Moreover, 
arbitrary (infinite) projective limits can be computed in $\cA_{alg}$ by 
first computing it (componentwise) in $\cA$, and then taking the largest 
algebraic subfunctor; alternatively, we might first restrict everything 
to $\catN$ by applying $J^*$, compute (componentwise)
the projective limit in $\catSets^\catN$, 
and then extend resulting functor $\catN\to\catSets$ 
to an algebraic endofunctor by means of~$J_!$.

\nxsubpoint
Observe that if $\Sigma$ is a monad over~$\catSets$, then 
$\Sigma_{alg}\subset\Sigma$ is a {\em submonad\/} of $\Sigma$ 
(cf.~\ptref{sp:submonads}). Indeed, $\Sigma_{alg}^2$ is algebraic, hence 
$\Sigma_{alg}^2\to\Sigma^2\stackrel\mu\to\Sigma$ factorizes through 
$\Sigma_{alg}\subset\Sigma$, and similarly $\epsilon:\Id_\catSets\to\Sigma$ 
factorizes through $\Sigma_{alg}$ since $\Id_\catSets$ is algebraic. 
Clearly, $\Sigma_{alg}\to\Sigma$ is universal among all homomorphisms 
from algebraic monads into $\Sigma$. In other words, the inclusion 
functor from the category of algebraic monads into the category of all 
monads admits a right adjoint $\Sigma\mapsto\Sigma_{alg}$, hence 
it commutes with arbitrary inductive limits of algebraic monads.

\nxsubpoint
We see that the algebraicity of a monad $\Sigma$ means that 
$\Sigma(X)=\bigcup_{\text{finite }I\subset X}\Sigma(I)$. If we think of 
$\Sigma(X)$ as ``the set of all formal $\Sigma$-linear combinations of 
elements of~$X$'', this condition means that any such $\Sigma$-linear 
combination actually involves only finitely many elements of~$X$. This 
explains why all monads defined by algebraic systems (in terms of some 
operations and relations between them), e.g.\ the monads 
$W$ and $\Sigma_R=R$ for an associative ring~$R$, or any of the 
monads listed in~\ptref{sp:examp.submon}, are actually algebraic in the 
sense of~\ptref{def:alg.endof}. However, monads $\hat\bbZ_\infty$   
(cf.~\ptref{sp:ex.nonalg.mon}) and $\hat\bbZ_p$ (cf.~\ptref{sp:mon.as.genr}) 
are clearly {\em not\/} algebraic; their largest algebraic subfunctors 
(which turn out to be submonads) are $\bbZ_\infty$ and 
$\bbZ_p$, respectively.

\nxsubpoint
Notice that these monads like $\hat\bbZ_\infty$ still have a similar property:
any formal $\hat\bbZ_\infty$-linear combination of elements of a set~$X$ 
involves only at most countably many elements of~$X$. This means that if 
we replace $\catN$ by the category of all ordinals $\leq\omega$, these 
monads will lie in the essential image of the corresponding left Kan 
extension $J_{\leq\omega,!}$. This leads to the idea of 
{\em cardinal filtration} on the category of all endofunctors, when we 
consider those endofunctors which lie in the essential image of 
$J_{\leq\alpha,!}$ or $J_{<\alpha,!}$ for some ``small'' cardinal 
$\alpha$ (e.g.\ belonging to the chosen universe~$\univU$). The set of 
all endofunctors which belong to some step of the cardinal filtration 
forms a full subcategory of $\cA$, stable under $\otimes$ (i.e.\ 
composition of functors) and arbitrary projective and inductive limits 
(along small index categories). We call such endofunctors {\em accessible}
(cf.\ SGA~4~I); 
clearly, their category $\cA_{acc}$ is a very good (rigorous) replacement 
for $\cA=\catEndof(\catSets)$, which contains all interesting endofunctors 
and at the same time avoids all set-theoretical complications 
(e.g.\ it is a $\univU$-category).

\nxsubpoint
Since $J^*$ induces an equivalence between categories $\cA_{alg}$ and 
$\catSets^\catN$, and since it is nothing else than the restriction functor, 
we usually denote an algebraic endofunctor $F$ and its restriction $J^*F$ by 
the same letter~$F$. We also often write $F(n)$ instead of $F(\stn)$ or 
$(J^*F)(\stn)$. In this sense an algebraic endofunctor $F$ is something 
like a (co)simplicial object: we get a collection of sets $\{F(n)\}_{n\geq0}$, 
together with some transition maps $F(\phi):F(m)\to F(n)$, defined for any 
map of finite sets $\phi:\stm\to\stn$. In this context we think of 
$J_!F$ (denoted again by~$F$) as some sort of extension of $F$ to arbitrary 
sets.

\nxsubpoint
In particular, we see that the category $\cA_{alg}\cong\catSets^\catN$ 
really makes sense, i.e.\ it is not ``too large'' (it is a 
$\univU$-category in the language of universes), so in this way we 
avoid the set-theoretical complications of~\ptref{sp:set.th.issues}.

\nxpoint\label{p:op.algendof} (Operations with algebraic endofunctors.)
Let's express some category-theoretic operations 
with algebraic endofunctors in terms of corresponding ``cosimplicial objects''.

\nxsubpoint\label{sp:simpl.descr.algendf}
We have already seen that an algebraic endofunctor~$F$ corresponds to 
a ``cosimplicial set'', i.e.\ to a collection of sets $\{F(n)\}_{n\geq0}$, 
together with some maps $F(\phi):F(m)\to F(n)$, defined for any map 
of finite sets $\phi:\stm\to\stn$, subject to usual conditions 
$F(\psi\circ\phi)=F(\psi)\circ F(\phi)$ and $F(\id_\stn)=\id_{F(n)}$. 
Moreover, a natural transformation $\beta:F\to G$ from an algebraic 
endofunctor~$F$ into another endofunctor~$G$ is given by a collection of 
maps $\beta_n:F(n)\to G(n)$, such that 
$G(\phi)\circ\beta_m=\beta_n\circ F(\phi)$ for any $\phi:\stm\to\stn$.

\nxsubpoint\label{sp:algmap.SXtoY}
Consider an algebraic endofunctor $\Sigma$, two arbitrary sets $X$ and~$Y$, 
and some map $\alpha:\Sigma(X)\to Y$ (if $X=Y$, we say that any such map is 
a {\em pre-action of $\Sigma$ on~$X$}). Now recall the expression for 
$\Sigma(X)$ as $\Coker(p,q:H_1(X)\rightrightarrows H_0(X))$, given 
in~\ptref{l:alt.descr.algfun}. It shows immediately that to give such an 
$\alpha:\Sigma(X)\to Y$ is equivalent to giving a family of maps
$\{\alpha^{(n)}:\Sigma(n)\times X^n\to Y\}_{n\geq0}$, subject to condition 
\begin{equation}\label{eq:preact.cond0}
\alpha^{(m)}\circ(\id_{\Sigma(m)}\times X^\phi)=
\alpha^{(n)}\circ(\Sigma(\phi)\times\id_{X^n}):\Sigma(m)\times X^n\to Y
\text{ for all }\phi:\stm\to\stn
\end{equation}
We see that $\Sigma(n)$ parametrizes something like ``$n$-ary operations 
from~$X$ to~$Y$''; if we agree to write $[t]_\alpha(x_1,x_2,\ldots,x_n)$ 
instead of $\alpha^{(n)}(t;x_1,x_2,\ldots,x_n)$, we see that giving an 
$\alpha:\Sigma(X)\to Y$ is equivalent to giving a collection of maps 
(``operations'') $[t]_\alpha:X^n\to Y$ for each ``formal $n$-ary operation'' 
$t\in\Sigma(n)$, $n\geq0$, subject to the only condition
\begin{multline}\label{eq:preact.cond0a}
\bigl[(F(\phi))(t)\bigr]_\alpha(x_1,x_2,\ldots,x_n)=
[t]_\alpha(x_{\phi(1)},x_{\phi(2)},\ldots,x_{\phi(m)})\\
\quad\text{for any $t\in\Sigma(m)$, $\phi:\stm\to\stn$, and 
$x_1,\ldots, x_n\in X$.}
\end{multline}
When we are given a pre-action $\alpha:\Sigma(X)\to X$ on some set $X$, 
we sometimes write $[t]_X$ instead of $[t]_\alpha$. Note that a map 
$f:X\to Y$ is compatible with some pre-actions $\alpha:\Sigma(X)\to X$ 
and $\beta:\Sigma(Y)\to Y$ (i.e.\ $f\circ\alpha=\beta\circ\Sigma(f)$) 
iff $f\bigl([t]_\alpha(x_1,x_2,\ldots,x_n)\bigr)=
[t]_\beta\bigl(f(x_1),f(x_2),\ldots,f(x_n)\bigr)$ for all choices of~$t$ 
and~$x_i$.

\nxsubpoint\label{sp:algnat.SFtoG}
Now suppose $F$, $G$ and $\Sigma$ are three algebraic endofunctors, and we 
want to describe all natural transformations $\alpha:\Sigma F\to G$. 
We know that this is equivalent to describing families of morphisms 
$\alpha_n:\Sigma F(n)\to G(n)$, functorial in $n$. On the other hand, 
we have just seen that any such $\alpha_n$ is given by a family of maps  
$\alpha_n^{(k)}:\Sigma(k)\times F(n)^k\to G(n)$, indexed by $k,n\geq0$, 
subject to compatibility conditions with respect to maps 
$\phi:\stm\to\stn$ and $\psi:\st{k}\to\st{l}$, similar to those considered 
above.

\nxsubpoint\label{sp:algmap.TSXtoY}
Let's describe maps $\alpha:\Xi\Sigma(X)\to Y$, where $\Xi$ and 
$\Sigma$ are algebraic endofunctors, and $X$ and $Y$ are two arbitrary sets. 
We know that $\alpha$ is determined by a collection of maps 
$\alpha^{(n)}:\Xi(n)\times\Sigma(X)^n\to Y$. On the other hand, functors 
$X\mapsto\Xi(n)\times\Sigma(X)^n$ are clearly algebraic again (for any fixed 
$n\geq0$), hence $\alpha^{(n)}$ is in its turn given by a collection of 
maps $\alpha^{(n,m)}:\Xi(n)\times\Sigma(m)^n\times X^m\to Y$, subject to 
some compatibility conditions with respect to maps $\phi:\stm\to\stm'$ and 
$\psi:\stn\to\stn'$.

\nxsubpoint\label{sp:algnat.TSFtoG}
This result generalizes to natural transformations $\alpha:\Xi\Sigma F\to G$, 
where $\Xi$, $\Sigma$, $F$ and $G$ are algebraic endofunctors. We see that 
such natural transformations are given by compatible families of maps 
$\{\alpha_n^{(k,m)}:\Xi(k)\times\Sigma(m)^k\times F(n)^m\to 
G(n)\}_{k,m,n\geq0}$. Of course, similar descriptions can be obtained for 
longer compositions of algebraic endofunctors: 
a map $\alpha:\Sigma_1\cdots\Sigma_s(X)\to Y$ is given by 
$\alpha^{(n_1,\ldots,n_s)}:\Sigma_1(n_1)\times\Sigma_2(n_2)^{n_1}\times\cdots
\times\Sigma_s(n_s)^{n_{s-1}}\times X^{n_s}\to Y$, and 
a natural transformation $\alpha:\Sigma_1\cdots\Sigma_s F\to G$ is given by 
$\alpha_n^{(n_1,\ldots,n_s)}:\Sigma_1(n_1)\times\Sigma_2(n_2)^{n_1}\times\cdots
\times\Sigma_s(n_s)^{n_{s-1}}\times F(n)^{n_s}\to G(n)$. If we put 
$M(m,n;\Sigma):=\Sigma(m)^n$ and consider this as ``the set of 
$m\times n$-matrices with entries in~$\Sigma$'', then the above maps can 
be seen as some ``matrix multiplication rules'' 
(notice the inverse order of arguments!) 
$M(n_1,n_0;\Sigma_1)\times M(n_2,n_1;\Sigma_2)\times\cdots\times
M(n_s,n_{s-1};\Sigma_s)\times M(n,n_s; F)\to M(n,n_0; G)$; this observation  
will be used later to compare our approach with that of M.~J.~Shai Haran.

\nxsubpoint\label{sp:lim.alg.funct}
Of course, inductive and projective limits of algebraic functors are 
computed in terms of corresponding ``cosimplicial sets'' componentwise: 
we have $(\injlim\Sigma_\alpha)(n)=\injlim\Sigma_\alpha(n)$, and 
$(\projlim\Sigma_\alpha)(n)=\projlim\Sigma_\alpha(n)$. For example, 
$(\Sigma_1\times_\Sigma\Sigma_2)(n)=\Sigma_1(n)\times_{\Sigma(n)}\Sigma_2(n)$. 
Notice that the formula $(\injlim\Sigma_\alpha)(X)=\injlim\Sigma_\alpha(X)$ is 
still valid for an arbitrary set~$X$, while the similar formula for 
projective limits in general doesn't hold unless the limit is finite 
(cf.~\ptref{sp:max.alg.subf}).

\nxsubpoint
In particular, the algebraic subfunctors $F'$ of a given algebraic functor 
$F$ are given by collections of subsets $\{F'(n)\subset F(n)\}_{n\geq0}$, 
stable under all maps $F(\phi):F(m)\to F(n)$. Of course, these are 
exactly the subobjects of~$F$ in~$\cA_{alg}$. Notice that after we extend 
$F'$ and $F$ to endofunctors on~$\catSets$, we still have 
$F'(X)\subset F(X)$ for any set $X$. 
Another interesting consequence: 
an algebraic equivalence relation~$R$ on some algebraic functor 
$F$ is given by a family of equivalence relations 
$R(n)\subset F(n)\times F(n)$, compatible with all maps~$F(\phi)$. In this 
case we can construct the quotient $F/R$, given by $(F/R)(n)=F(n)/R(n)$. 
It is again an algebraic functor, and, if we extend everything to arbitrary 
sets~$X$, we still have $(F/R)(X)=F(X)/R(X)$, and $F\to F/R$ is a 
strict epimorphism with kernel $R$ both in $\cA$ and~$\cA_{alg}$.

\nxsubpoint\label{sp:ker.im.morph.algfun}
Notice that any morphism $\rho:F\to G$ of algebraic endofunctors 
has an algebraic kernel pair $R:=F\times_G F$ as well as an algebraic image
$I=\rho(F)\cong F/R$. Of course, these notions can be computed componentwise: 
$R(n)\subset F(n)\times F(n)$ is the equivalence relation defined by 
$\rho_n:F(n)\to G(n)$, and $F(n)/R(n)\cong\rho_n(F(n))=I(n)\subset G(n)$. 
Moreover, these formulas extend to arbitrary sets $X$, i.e.\ 
$R(X)\subset F(X)\times F(X)$ is still the kernel of $\rho_X: F(X)\to G(X)$ 
and so on.

\nxsubpoint\label{sp:nat.idtoalg}
What about natural transformations $\epsilon:\Id_{\catSets}\to F$ from 
the identity functor into an arbitrary (algebraic) endofunctor~$F$? 
Since $\Id_{\catSets}=\Hom(\st1,-)$, we see by Yoneda that such an 
$\epsilon$ is uniquely determined by an element $\bu:=\epsilon_{\st1}(1)\in 
F(1)$, which can be chosen arbitrarily. In this case we get 
$\epsilon_X(x)=(F(\tilde{x}))(\bu)$ for any set $X$ and any $x\in X$, 
where $\tilde{x}:\st1\to X$ is the only map with image~$x$. 
In particular, for an algebraic~$F$ the maps $\epsilon_n:\stn\to F(n)$ 
can be reconstructed from $\bu\in F(1)$ in this way.

\nxsubpoint\label{sp:expl.SX}
According to~\ptref{l:alt.descr.algfun}, we see that for any algebraic 
endofuctor~$\Sigma$ and any set~$X$ we have a surjection 
$H_0(X)=\bigsqcup_{n\geq0}\Sigma(n)\times X^n\twoheadrightarrow \Sigma(X)$. 
Consider its individual component $\Sigma(n)\times X^n\to\Sigma(X)$; clearly, 
it maps $(t;x_1,\ldots,x_n)$ into $(\Sigma(x))(t)\in\Sigma(X)$, where 
$x:\stn\to X$ is the map $k\mapsto x_k$. Let's denote this element 
by $t(\{x_1\},\ldots,\{x_n\})$. For example, if we take $X:=\stn$, 
$x:=\id_\stn$, $x_k=k$, we get $t(\{1\},\ldots,\{k\})=t\in\Sigma(\stn)$. 
Notice that {\em any\/} element of $\Sigma(X)$ can be written in form 
$t(\{x_1\},\ldots,\{x_n\})$ for some $n\geq0$; in fact, $\Sigma(X)$ 
consists of all such expressions modulo the equivalence relation $\equiv$ 
generated by $(\phi_*t)(\{x_1\},\ldots,\{x_n\})\equiv
t(\{x_{\phi(1)}\},\ldots,\{x_{\phi(n)}\})$ for any $\phi:\stm\to\stn$, 
$t\in\Sigma(m)$ and $x_1,\ldots, x_n\in X$, where we put 
$\phi_*t:=(\Sigma(\phi))(t)$. Moreover, since the category $\catN/X$ 
is filtering, we see that two elements $t(\{x_1\},\ldots,\{x_n\})$ and 
$t'(\{y_1\},\ldots,\{y_m\})$ are equal in $\Sigma(X)$ iff 
we can find some $\phi:\stn\to\stp$, $\psi:\stm\to\stp$ and some 
elements $z_1,\ldots,z_p\in X$, such that $\phi_*t=\psi_*t'$, 
$x_i=z_{\phi(i)}$, and $y_j=z_{\psi(j)}$. We can even assume these elements 
$z_k$ to be pairwise distinct.

Notice that $\Sigma(f):\Sigma(X)\to\Sigma(Y)$ 
maps $t(\{x_1\},\ldots,\{x_n\})\in\Sigma(X)$ into 
$t\bigl(\{f(x_1)\},\ldots,\{f(x_n)\}\bigr)$, for any map $f:X\to Y$.

\nxsubpoint\label{sp:expl.SXtoY}
In particular, suppose we have a map $\alpha:\Sigma(X)\to Y$, 
described by the collection $\{\alpha^{(n)}:\Sigma(n)\times X^n\to Y\}$ as 
in~\ptref{sp:algmap.SXtoY}. By definition, $\alpha^{(n)}$ is just the 
composite map $\Sigma(n)\times X^n\to H_0(X)\twoheadrightarrow\Sigma(X)
\stackrel\alpha\to Y$, hence $\alpha$ maps an element 
$t(\{x_1\},\ldots,\{x_n\})\in\Sigma(X)$ into $\alpha^{(n)}(t;x_1,\ldots,x_n)=
[t]_\alpha(x_1,\ldots,x_n)$.

\nxpointtoc{Algebraic monads}\label{p:algmon}
Now we want to combine together some of previous considerations to 
obtain a ``na{\"\i}ve'' description of an algebraic monad 
$\Sigma=(\Sigma,\mu,\epsilon)$ (i.e.\ an algebra in~$\cA_{alg}$, 
or equivalently, a monad over $\catSets$ with an algebraic underlying 
endofunctor) in terms of a sequence of sets $\{\Sigma(n)\}_{n\geq0}$ 
and some maps between these sets and their products, subject to some 
conditions. This na{\"\i}ve approach has its advantages: for example, 
it immediately generalizes to the case of an arbitrary category~$\cC$, 
and it is completely clear how to apply left exact functors to such 
objects (e.g.\ pullback, direct image or global sections functors between 
topoi), or how to define sheaves of them over any site $\cS$. However, 
we'll develop a (plain inner) monad interpretation in the 
topos case as well, since it also has its advantages.

\nxsubpoint\label{sp:algmon.predescr}
So let's fix an algebraic monad $\Sigma=(\Sigma,\mu,\epsilon)$ over 
$\catSets$, and a set $X$ with some $\Sigma$-action $\alpha:\Sigma(X)\to X$. 
We are going to obtain a ``cosimplicial'' description of these things.

First of all, an algebraic endofunctor $\Sigma$ and natural transformations 
$\mu:\Sigma^2\to\Sigma$ and $\epsilon:\Id_\catSets\to\Sigma$ are uniquely 
determined by the following data:

a) A collection of sets $\{\Sigma(n)\}_{n\geq0}$ and maps 
$\Sigma(\phi):\Sigma(m)\to\Sigma(n)$, defined for any map $\phi:\stm\to\stn$, 
subject to conditions $\Sigma(\id_\stn)=\id_{\Sigma(n)}$ and 
$\Sigma(\psi\circ\phi)=\Sigma(\psi)\circ\Sigma(\phi)$ 
(cf.~\ptref{sp:simpl.descr.algendf}).

b) A collection of ``multiplication'' or ``evaluation'' maps 
$\mu_n^{(k)}:\Sigma(k)\times\Sigma(n)^k\to\Sigma(n)$, subject to conditions 
$\mu_n^{(k)}\circ(\id_{\Sigma(k)}\times\Sigma(n)^\phi)=
\mu_n^{(k')}\circ(\Sigma(\phi)\times\id_{\Sigma(n)^{k'}})$ for any 
$\phi:\stk\to\stk'$, and 
$\mu_n^{(k)}\circ(\id_{\Sigma(k)}\times\Sigma(\psi)^k)=\mu_m^{(k)}$ 
for any $\psi:\stm\to\stn$ (cf.~\ptref{sp:algnat.SFtoG}). 
Usually we write $[t]_{\Sigma(n)}(x_1,x_2,\ldots,x_k)$ or even 
$t(x_1,x_2,\ldots,x_k)$ instead of $\mu_n^{(k)}(t;x_1,x_2,\ldots,x_k)$ 
for any $t\in\Sigma(k)$ and $x_1$, $\ldots$, $x_k\in\Sigma(n)$. In this 
case the above requirements can be rewritten as 
$[\phi_*(t)]_{\Sigma(n)}(x_1,\ldots,x_{k'})=
[t]_{\Sigma(n)}(x_{\phi(1)},\ldots,x_{\phi(k)})$ and 
$\psi_*\bigl([t]_{\Sigma(m)}(x_1,\ldots,x_k)\bigr)
=[t]_{\Sigma(n)}(\psi_*x_1,\ldots,\psi_*x_k)$, where we write 
$\phi_*(t)$ or $\phi_*t$ instead of $(\Sigma(\phi))(t)$. 
If we put $M(m,n;\Sigma):=\Sigma(m)^n$ (``the set of $m\times n$-matrices 
over~$\Sigma$''), then ${(\mu_n^{(k)})}^m:\Sigma(k)^m\times\Sigma(n)^k\to
\Sigma(n)^m$ can be interpreted as a ``matrix multiplication map'' 
$M(k,m;\Sigma)\times M(n,k;\Sigma)\to M(n,m;\Sigma)$
(cf.~\ptref{sp:algnat.TSFtoG}).

c) An element $\bu\in\Sigma(1)$, called the {\em identity\/} of $\Sigma$. 
Then $\epsilon_X:X\to\Sigma(X)$ maps any $x\in X$ into 
$(\Sigma(\tilde{x}))(\bu)$, where $\tilde{x}:\st1\to X$ is the map with 
image~$x$ (cf.~\ptref{sp:nat.idtoalg}). We usually denote the element 
$\epsilon_X(x)\in\Sigma(X)$ by $\{x\}_X$ or $\{x\}$. In particular, 
we have the ``basis elements'' $\{k\}=\{k\}_\stn\in\Sigma(n)$, 
$1\leq k\leq n$. For example, $\{1\}_{\st1}=\bu$. 
Clearly, for any map $f:X\to Y$ and any $x\in X$ we have 
$f_*\{x\}=(\Sigma(f))(\epsilon_X(x))=\{f(x)\}$.

d) Similarly, a pre-action $\alpha:\Sigma(X)\to X$ is given by a collection 
of maps $\alpha^{(n)}:\Sigma(n)\times X^n\to X$, 
subject to relation \eqref{eq:preact.cond0}, and we usually 
write $[t]_\alpha(x_1,\ldots,x_n)$ or $[t]_X(x_1,\ldots,x_n)$ or even 
$t(x_1,\ldots,x_n)$ instead of $\alpha^{(n)}(t;x_1,\ldots,x_n)$
(cf~\ptref{sp:algmap.SXtoY}). Notice that $\mu_X:\Sigma^2(X)\to\Sigma(X)$ 
defines a pre-action of $\Sigma$ on $\Sigma(X)$, so the notation 
$[t]_{\Sigma(X)}(f_1,\ldots,f_n)$ makes sense for any $t\in\Sigma(n)$ 
and $f_1$, $\ldots$, $f_n\in\Sigma(X)$. If we take $X=\stn$, 
we recover again the maps $\mu_n^{(k)}$ of~b), so our notation 
$[t]_{\Sigma(n)}$ is at least consistent.

\nxsubpoint\label{sp:expl.action}
Now we want to express in terms of these data the monad axioms 
$\mu\circ(\epsilon\star\Sigma)=\id_\Sigma=\mu\circ(\Sigma\star\epsilon)$ 
and $\mu\circ(\mu\star\Sigma)=\mu\circ(\Sigma\star\mu)$, as well as the 
axioms $\alpha\circ\epsilon_X=\id_X$ and $\alpha\circ\Sigma(\alpha)=
\alpha\circ\mu_X$ for a pre-action $\alpha:\Sigma(X)\to X$ to be an action.
Let's start with the action axioms:

a) Condition $\alpha\circ\epsilon_X=\id_X$ is equivalent to 
$\alpha^{(1)}_1(\bu;x)=x$ for any $x\in X$, i.e.\ $[\bu]_X(x)=x$; in other 
words, $[\bu]_X:X\to X$ has to be the identity map. This condition together 
with functoriality in~$\stn$ implies that ${[\{k\}_\stn]}_X:X^n\to X$ is the 
projection onto the $k$-th component.

b) Condition $\alpha\circ\Sigma(\alpha)=\alpha\circ\mu_X$ translates into 
the commutativity of the following diagram, which can be thought of as 
a sort of ``associativity of matrix multiplication'' 
$M(k,1;\Sigma)\times M(n,k;\Sigma)\times M(1,n;X)\to M(1,1;X)$ 
(cf.~\ptref{sp:algnat.TSFtoG}):
\begin{equation}
\xymatrix@C+4pt{
\Sigma(k)\times\Sigma(n)^k\times X^n
\ar[rr]^<>(.5){\id_{\Sigma(k)}\times(\alpha^{(n)})^k}
\ar[d]^{\mu_n^{(k)}\times\id_{X^n}}
&&\Sigma(k)\times X^k\ar[d]^{\alpha^{(k)}}\\
\Sigma(n)\times X^n\ar[rr]^<>(.5){\alpha^{(n)}}&&X}
\end{equation}
In our operational notation this can be written as follows:
\begin{multline}\label{eq:act.assoc.0b}
\bigl[{[t]}_{\Sigma(n)}(t_1,\ldots,t_k)\bigr]_X(x_1,\ldots,x_n)=
{[t]}_X\bigl({[t_1]}_X(x_1,\ldots,x_n),\ldots,{[t_k]}_X(x_1,\ldots)\bigr)\\
\text{for $t\in\Sigma(k)$, $t_1,\ldots, t_k\in\Sigma(n)$, 
$x_1,\ldots,x_n\in X$}
\end{multline}
It is useful to think about $t\in\Sigma(n)$ as some sort of formal polynomial 
or formal operation with $n$ arguments, and to interpret the above 
condition as some sort of substitution rule.

c) Recall that $f:X\to Y$ is compatible with given pre-actions 
$\alpha:\Sigma(X)\to X$ and $\beta:\Sigma(Y)\to Y$ iff 
$f\circ\alpha^{(n)}=\beta^{(n)}\circ(\id_{\Sigma(n)}\times f^n)$, i.e.\ 
$f\bigl({[t]}_X(x_1,\ldots,x_n)\bigr)={[t]}_Y\bigl(f(x_1),\ldots,f(x_n)\bigr)$ 
(cf.~\ptref{sp:algmap.SXtoY}).

\nxsubpoint\label{sp:expl.algmon.cond}
We continue by expressing the conditions for $\Sigma=(\Sigma,\mu,\epsilon)$ 
to be a monad. Notice that these conditions are some equalities between 
natural transformations of algebraic endofunctors, hence it is enough 
to require these pairs of natural transformations to coincide
on standard finite sets.

First of all, $\mu_n\circ\Sigma(\mu_n)=\mu_n\circ\mu_{\Sigma(n)}$ and 
$\mu_n\circ\epsilon_{\Sigma(n)}=\id_{\Sigma(n)}$ are equivalent to 
saying that the pre-action $\mu_n:\Sigma^2(n)\to\Sigma(n)$ of $\Sigma$ 
on $\Sigma(n)$ is actually an action, so we obtain conditions 
\ptref{sp:expl.action}a) and b) for $(X,\alpha)=(\Sigma(n),\mu_n)$, 
i.e.\ we have to replace $\alpha^{(k)}$ with $\mu_n^{(k)}$. 
We obtain the following two conditions:

a) $\mu_n^{(1)}(\bu;t)=t$ for any $t\in\Sigma(n)$ and any $n\geq0$. 
In other words, we require $[\bu]_{\Sigma(n)}(t)=t$, i.e.\ 
$[\bu]_{\Sigma(n)}=\id_{\Sigma(n)}$.

b) The associativity $\mu\circ(\Sigma\star\mu)=\mu\circ(\mu\star\Sigma)$ 
translates into some sort of ``associativity of matrix multiplication'' 
$M(k,1;\Sigma)\times M(n,k;\Sigma)\times M(m,n;\Sigma)\to M(m,1;\Sigma)$, 
or, if the reader prefers,  
$M(k,s;\Sigma)\times M(n,k;\Sigma)\times M(m,n;\Sigma)\to M(m,s;\Sigma)$, 
where we put $M(m,n;\Sigma):=\Sigma(m)^n$
(cf.~\ptref{sp:algnat.TSFtoG}):
\begin{equation}
\xymatrix@C+4pt{
\Sigma(k)\times\Sigma(n)^k\times\Sigma(m)^n
\ar[rr]^<>(.5){\id_{\Sigma(k)}\times(\mu_m^{(n)})^k}
\ar[d]^{\mu_n^{(k)}\times\id_{\Sigma(m)^n}}
&&\Sigma(k)\times\Sigma(m)^k\ar[d]^{\mu_m^{(k)}}\\
\Sigma(n)\times\Sigma(m)^n\ar[rr]^<>(.5){\mu_m^{(n)}}&&\Sigma(m)}
\end{equation}
Of course, this can be rewritten in form~\eqref{eq:act.assoc.0b}, 
where we put $X:=\Sigma(m)$.

c) Now only the axiom $\mu_n\circ\Sigma(\epsilon_n)=\id_{\Sigma(n)}$ 
remains. According to~\ptref{sp:expl.SX}, $\Sigma(\epsilon_n)$ 
maps $t=t(\{1\},\ldots,\{n\})\in\Sigma(n)$ into 
$t(\{\{1\}\},\ldots,\{\{n\}\})\in\Sigma^2(n)$, and~\ptref{sp:expl.SXtoY} 
shows that $\mu_n$ maps this element to 
${[t]}_{\Sigma(n)}(\{1\},\ldots,\{n\})$,
hence the last axiom boils down to 
\begin{equation}\label{eq:mon.unit.1c}
{[t]}_{\Sigma(n)}\bigl(\{1\}_\stn,\ldots,\{n\}_\stn\bigr)=t
\text{ for any $t\in\Sigma(n)$}
\end{equation}
This is something like the identity $F(T_1,T_2,\ldots,T_n)=F$ 
for any polynomial~$F$ from $\bbZ[T_1,T_2,\ldots,T_n]$ 
(actually, if we consider the monad defined by the category of commutative 
rings, we see that this identity is indeed a special case 
of~\eqref{eq:mon.unit.1c}). 
Another possible interpretation: 
$I_n:=(\{1\},\ldots,\{n\})\in\Sigma(n)^n=M(n,n;\Sigma)=:M(n,\Sigma)$ 
is the ``identity matrix'', i.e.\ it induces the identity map 
$\Sigma(n)\to\Sigma(n)$.

d) Finally, let us express the conditions for a natural transformation 
$\rho:\Sigma\to\Xi$ to be a homomorphism of algebraic monads. 
Of course, $\rho$ is given by a collection of maps 
$\rho_n:\Sigma(n)\to\Xi(n)$, such that 
$\Xi(\phi)\circ\rho_m=\rho_n\circ\Sigma(\phi)$ for any $\phi:\stm\to\stn$. 
Compatibility with the unit $\rho\circ\epsilon_\Sigma=\epsilon_\Xi$ 
translates into $\rho_1(\bu_\Sigma)=\bu_\Xi$, and compatibility with 
the multiplication $\rho\circ\mu_\Sigma=\mu_\Xi\circ(\rho\star\rho)$
translates into $\rho_n\circ\mu_{\Sigma,n}^{(k)}=\mu_{\Xi,n}^{(k)}\circ
(\rho_k\times\rho_n^k)$, or, equivalently, into 
\begin{multline}\label{eq:cond.algmon.hom}
\rho_n\bigl({[t]}_{\Sigma(n)}(t_1,\ldots,t_k)\bigr)=
{[\rho_k(t)]}_{\Xi(n)}\bigl(\rho_n(t_1),\ldots,\rho_n(t_k)\bigr)\\
\text{for any $t\in\Sigma(k)$, $t_1,\ldots,t_k\in\Sigma(n)$}
\end{multline}
Sometimes we identify $\rho$ with the ``graded'' map
$\|\rho\|:=\sqcup\rho_n:\|\Sigma\|=\bigsqcup_{n\geq0}\Sigma(n)\to
\|\Xi\|=\bigsqcup_{n\geq0}\Xi(n)$, 
and write $\rho(t)$ or even $\rho t$ instead of $\rho_n(t)$ for some 
$t\in\Sigma(n)$.

\nxsubpoint
Consider the case when we require all $\Sigma(n)$ to be modules over a 
commutative ring~$K$, all $\Sigma(\phi)$ to be $K$-linear, and all 
$\mu_n^{(k)}:\Sigma(k)\times\Sigma(n)^k\to\Sigma(n)$ to be 
$K$-polylinear. Then $\mu_n^{(k)}$ induces some 
$\tilde\mu_n^{(k)}:\Sigma(k)\otimes\Sigma(n)^{\otimes k}\to\Sigma(n)$, 
and we can rewrite the above conditions for $\Sigma$ to be an algebraic 
monad in terms of these maps $\tilde\mu_n^{(k)}$, thus obtaining nothing 
else than the definition of an {\em operad over~$K$} (up to some minor
points: in fact, we obtain an equivalent description of operads with
compatible cosimplicial structure). Of course, 
the definition of an operad makes sense in any ACU $\otimes$-category 
$\cC=(\cC,\otimes)$; if we put $(\cC,\otimes):=\catMod\Fempty=
(\catSets,\times)$, we obtain back the definition of an algebraic monad, 
hence algebraic monads can be thought of as 
{\em (cosimplicial) operads over~$\Fempty$.}

\nxsubpoint\label{sp:underl.monoid}
Notice that $\mu_1^{(1)}:\Sigma(1)\times\Sigma(1)\to\Sigma(1)$ and 
$\bu\in\Sigma(1)$ define a monoid structure on the underlying set 
$|\Sigma|:=\Sigma(1)$ of an algebraic monad~$\Sigma$, if we put 
$\lambda\cdot\mu=\lambda\mu:=\mu_1^{(1)}(\lambda;\mu)=[\lambda]_{\Sigma(1)}
(\mu)=[\lambda]_{|\Sigma|}(\mu)$; in particular, 
$[\lambda]_{|\Sigma|}:|\Sigma|\to|\Sigma|$ is just the left multiplication 
by~$\lambda$ map. Indeed, the associativity $(\lambda\mu)\nu=
\lambda(\mu\nu)$ is a special case of~\ptref{sp:expl.algmon.cond}b), 
and $\bu\lambda=\lambda$ and $\lambda\bu=\lambda$ are special cases 
of~\ptref{sp:expl.algmon.cond}a) and~c). Of course, this is exactly the 
monoid structure on $|\Sigma|$ considered in~\ptref{sp:underl.set.monad}. 
If $M=(M,\alpha)$ is a $\Sigma$-module, we recover again the monoid action 
of $|\Sigma|$ on $M$ by taking $\alpha^{(1)}:\Sigma(1)\times M\to M$. 
Then $\lambda\cdot x=\lambda x:=\alpha^{(1)}(\lambda;x)={[\lambda]}_M(x)$ 
for any $\lambda\in|\Sigma|$, $x\in M$. That's why we usually denote 
the action of unary operations $\lambda\in|\Sigma|$ on an element 
$x\in M$ simply by $\lambda x$.

\nxsubpoint
We have already seen in~\ptref{sp:expl.SX} that $\Sigma(X)$ consists of all 
formal expressions of form $t(\{x_1\},\ldots,\{x_n\})$, $n\geq0$, 
$t\in\Sigma(n)$, $x_1,\ldots, x_n\in X$, modulo a certain equivalence 
relation. On the other hand, since we have agreed to write 
$\{x\}$ for $\epsilon_X(x)\in\Sigma(X)$, and since 
$\mu_X$ is a canonical $\Sigma$-structure on $\Sigma(X)$, 
we see that ${[t]}_{\Sigma(X)}(\{x_1\},\ldots,\{x_n\})$ actually means 
$\mu_X^{(n)}(t;\epsilon_X(x_1),\ldots,\epsilon_X(x_n))$. Since 
$\mu_X\circ\Sigma(\epsilon_X)=\id_{\Sigma(X)}$, this expression is 
indeed equal to $t(\{x_1\},\ldots,\{x_n\})$ of~\ptref{sp:expl.SX}, 
so our notation is consistent in this respect.

\nxsubpoint\label{sp:algmon.mod.alg}
We see that a $\Sigma$-module structure $\alpha:\Sigma(X)\to X$ is 
given by a collection of ``operation maps'' 
${[t]}_X:X^n\to X$, parametrized by $t\in\bigsqcup_{n\geq0}\Sigma(n)$. 
These maps have to satisfy some relations like $[\bu]_X=\id_X$, 
${[\phi_*t]}_X={[t]}_X\circ X^\phi$ for any $\phi:\stm\to\stn$, 
i.e.\ ${[\phi_*t]}_X(x_1,\ldots,x_n)={[t]}_X(x_{\phi(1)},\ldots,x_{\phi(m)})$, 
and the ``compound operation'' relations~\eqref{eq:act.assoc.0b}. 
Conversely, \ptref{sp:expl.action} shows that these relations are 
sufficient for a collection of maps ${[t]}_X$ to define a $\Sigma$-action 
on~$X$. Since all these relations are algebraic, we see that a 
$\Sigma$-structure on a set~$X$ is an example of an algebraic structure,  
hence $\catMod{\Sigma}=\catSets^\Sigma$ has all the properties 
of a category defined by an algebraic structure.

\nxsubpoint\label{sp:def.absend.ring}
Given any set~$X$, we can construct its ``absolute endomorphism ring'' 
$\END(X)$. It is an algebraic monad defined by 
$(\END(X))(n):=\iHom_\catSets(X^n,$ $X)$; we take $\bu_\Xi$ to be equal to 
$\id_X$, and $\mu_n^{(k)}:\iHom(X^k,X)\times\iHom(X^n,X)^k\cong
\iHom(X^k,X)\times\iHom(X^n,X^k)\to\iHom(X^n,X)$ are the usual 
composition maps. The conditions of~\ptref{sp:expl.algmon.cond} 
obviously hold, and we have a canonical action of~$\END(X)$ on $X$. 
Furthermore, the considerations of~\ptref{sp:algmon.mod.alg} show that 
giving an action $\alpha:\Sigma(X)\to X$ of an algebraic monad $\Sigma$ 
on~$X$ is equivalent to giving a monad homomorphism (``representation'') 
$\rho:\Sigma\to\END(X)$, and then the $\Sigma$-action on~$X$ can be 
recovered from the canonical $\END(X)$-action on~$X$ by means of 
scalar restriction along~$\rho$.

We say that $X$ is a {\em faithful\/ $\Sigma$-module\/} if the corresponding 
homomorphism $\rho:\Sigma\to\END(X)$ is injective, i.e.\ 
if ${[t]}_X={[t']}_X$ implies $t=t'$ for any $t,t'\in\Sigma(n)$ and any 
$n\geq0$. In general case we denote the image $\rho(\Sigma)$ of $\rho$ 
by $\Sigma_X$; clearly, this is the only strict quotient of~$\Sigma$ 
which acts faithfully on~$X$.

\nxsubpoint (Terminology and notations.)
Given any algebraic monad~$\Sigma$, we call the elements of $\Sigma(0)$ 
{\em constants (of~$\Sigma$)}, elements of $|\Sigma|=\Sigma(1)$ --- 
{\em unary operations}, elements of $\Sigma(2)$ --- {\em binary operations}, 
\dots, elements of $\Sigma(n)$ --- {\em $n$-ary operations (of~$\Sigma$)}. 
We say that $\|\Sigma\|:=\bigsqcup_{n\geq0}\Sigma(n)$ is the set of all 
operations of~$\Sigma$; if $t\in\Sigma(n)$, we say that {\em the arity of~$t$ 
is equal to~$n$}, and sometimes write $t^{[n]}$ or $r(t)=n$ 
to emphasize this fact. 
Cleary, an operation $t$ of arity~$n$ induces a map 
${[t]}_X:X^n\to X$ on any $\Sigma$-module~$X$. We often write 
$t_X(x_1,\ldots,x_n)$, $t(x_1,\ldots,x_n)$ or even $t\,x_1\,\ldots\,x_n$ 
instead of ${[t]}_X(x_1,\ldots,x_n)$. In particular, if $\lambda\in|\Sigma|$ 
is a unary operation, we write $\lambda_X(x)$ or $\lambda x$ instead of 
${[\lambda]}_X(x)$; this is especially useful when $\lambda$ is denoted 
by some sign like $-$; then we write $-x$ instead of ${[-]}_X(x)$, 
and $-_X:X\to X$ instead of~${[-]}_X$. 
Notice that $(\lambda t)(x_1,\ldots,x_n)=\lambda(t(x_1,\ldots,x_n))$ 
for any $\lambda\in|\Sigma|$, $t\in\Sigma(n)$ 
(this is a special case of~\eqref{eq:act.assoc.0b}), so the notation 
$\lambda t(x_1,\ldots,x_n)$ is unambiguous.

Similarly, if a binary operation 
is denoted by some sign like $+$, $*$, $\times$ \dots, which is 
usually written in infix form, we write $x_1+x_2$ instead of 
${[+]}_X(x_1,x_2)$; notice that $[+]={[+]}_{\Sigma(2)}(\{1\},\{2\})=
\{1\}+\{2\}$ in $\Sigma(2)$. We obey the usual precedence rules while 
using such notations, inserting parentheses when necessary; 
if an operation $+$ is associative, we write $x+y+z$ instead of 
$(x+y)+z$ or $x+(y+z)$. For example, any $R$-linear combination 
$t=(\lambda_1,\ldots,\lambda_n)\in\Sigma_R(n)=R^n$, where $\lambda_i\in R$, 
$R$ is any associative ring (cf.~\ptref{sp:mon.from.ring}), 
can be safely written now in form 
$t=\lambda_1\{1\}+\lambda_2\{2\}+\cdots+\lambda_n\{n\}$, where 
$\lambda_i\in R=|\Sigma_R|$ are interpreted as unary operations of 
$\Sigma_R$, and $[+]=(1,1)\in\Sigma_R(2)$ is of course a binary operation 
of $\Sigma_R$. In this situation we also have a unary operation 
$[-]:=-1\in|\Sigma_R|$, and a constant $0\in\Sigma_R(0)$.

\nxsubpoint\label{sp:sigma.injtoinj}
In particular, for any constant $c\in\Sigma(0)$ and any $\Sigma$-module 
$X$ we get the {\em value\/} $c_X={[c]}_X\in X$ of this constant in~$X$. 
Of course, any $\Sigma$-homomorphism $f:X\to Y$ has to preserve all 
constants: $f(c_X)=c_Y$. For example, if $i:\st0=\emptyset\to\st1$ is the 
canonical embedding, then the unary operation $i_*c=c_{|\Sigma|}\in\Sigma(1)$ 
acts on any $\Sigma$-module~$X$ by mapping all elements of~$X$ into 
$c_X$, i.e.\ ${[i_*c]}_X(x)=c_X$. 
This shows that $\Sigma(i):\Sigma(0)\to\Sigma(1)$ is {\em injective}: 
indeed, $i_*c=i_*c'$ implies ${[i_*c]}_{\Sigma(0)}={[i_*c']}_{\Sigma(0)}$, 
hence $c={[i_*c]}_{\Sigma(0)}(x)={[i_*c']}_{\Sigma(0)}(x)=c'$, where we 
can take any auxiliary element $x\in\Sigma(0)$ (e.g.\ $x=c$). 

This fact, combined with the observation that any monomorphism (i.e.\ 
injective map) $f:X\to Y$ with non-empty source splits in $\catSets$ 
(i.e.\ admits a left inverse), hence $\Sigma(f):\Sigma(X)\to\Sigma(Y)$ 
has to be injective when $X$ is non-empty, shows that $\Sigma(f)$ 
is in fact injective even if~$X$ is empty. In particular, all maps 
$\Sigma(\phi):\Sigma(m)\to\Sigma(n)$ for injective $\phi:\stm\to\stn$ 
are injective; if we put $\Sigma(\infty):=\Sigma(\bbZ_{\geq0})=
\injlim_{n\geq0}\Sigma(n)=\bigcup_n\Sigma(n)$, we see that 
$\Sigma(n)$ is a subset of $\Sigma(\infty)$, and these subsets form a 
filtration on~$\Sigma(\infty)$, so we can combine all $\Sigma(n)$ 
into one large filtered set $\Sigma(\infty)$.

\nxsubpoint\label{sp:mon.wo.const}
Notice that $\Sigma(0)=L_\Sigma(\emptyset)$ is an initial object of 
$\catMod\Sigma=\catSets^\Sigma$. We say that $\Sigma$ is a 
{\em monad without constants\/} if $\Sigma(0)=\emptyset$; in this case 
the empty set is (the only) initial object of $\catMod\Sigma$, 
in particular, the empty set admits a (unique) $\Sigma$-action.
We have already seen some examples of monads without constants like 
$\Fempty$, ${\bm\Delta}$ or $\Aff_R$ (cf.~\ptref{sp:examp.submon}). 
We can also consider the submonad $W_+\subset W$ given by $W_+(X)=W(X)$ if 
$X\neq\emptyset$, $W_+(\emptyset)=\emptyset$; this is also a monad 
without constants, and $\catMod{W_+}$ is the category of semigroups.

\nxsubpoint\label{sp:mon.with.zero}
If $\Sigma(0)$ consists of exactly one constant, usually denoted by $0$, 
we say that $\Sigma$ is a {\em monad with zero.} In this case 
$\Sigma(0)=\{0\}$ is both a initial and a final object, i.e.\ 
a {\em zero object\/} $0:=\Sigma(0)$ of $\catMod\Sigma$, and conversely, 
since $\Sigma(0)$ is always initial and $\st1$ is always final in 
$\catMod\Sigma$, if $\catMod\Sigma$ has a zero object, then 
$\Sigma$ is a monad with zero. In this case any $\Sigma$-module~$X$ 
has a pointed element $0_X=0$, which has to be respected by any 
$\Sigma$-homomorphism; in particular, $X$ cannot be empty. 
Moreover, any set $\Hom_\Sigma(M,N)$ has a pointed element --- the 
{\em zero homomorphism $0_{MN}:M\to0\to N$}. Almost all monads 
considered by us up to now have been monads with zero, 
with the exception of $W_U$ and of those listed in~\ptref{sp:mon.wo.const}.

\nxsubpoint\label{sp:prop.triv.mon} (Triviality of a monad.)
Recall that the category of all monads over $\catSets$ admits a final object 
--- the {\em trivial\/} or {\em final monad\/} $\st1$, 
such that $\st1(X)$ is the one-point set $\st1$ for any set~$X$. 
Clearly, this monad is algebraic, so it is 
a final object in the category of algebraic monads. We say that 
$\Sigma$ is {\em trivial\/} if it is isomorphic to~$\st1$; in this case 
we write $\Sigma\cong\st1$ or even $\Sigma=\st1$. In some cases we write 
$0$ or $\st0$ instead of~$\st1$, since $\st1$ is the monad defined by the 
trivial ring~$0$. Let's list 
some criteria for triviality of a monad:
\begin{Propz}
The following conditions for an algebraic monad~$\Sigma$ are equivalent:
(i) $\Sigma$ is trivial;
\quad(ii) $\catMod\Sigma$ is a trivial category, i.e.\ all its objects are 
final;
\quad(ii') The only sets which admit a $\Sigma$-structure are the one-element
sets;
\quad(iii) $\epsilon_X:X\to\Sigma(X)$ is {\em not\/} injective for 
at least one set~$X$, and $\Sigma(0)\neq\emptyset$;
\quad(iii') $\epsilon_{\st2}:\st2\to\Sigma(2)$ is not injective, i.e.\ 
$\{1\}_{\st2}=\{2\}_{\st2}$, and $\Sigma(0)\neq\emptyset$;
\quad(iv) $\bu\in\Sigma(1)$ comes from a constant $c\in\Sigma(0)$, i.e.\ 
$\bu={[c]}_{\Sigma(1)}$.
\end{Propz}
\begin{Proof}
(i)$\Rightarrow$(ii): Indeed, for any $\Sigma$-module~$X$ the action 
$\alpha:\Sigma(X)=\st1\to X$ has to be surjective, and this is possible only 
if $X\cong\st1$; (ii)$\Leftrightarrow$(ii') is evident.
(ii')$\Rightarrow$(i): All sets $\Sigma(n)$ admit a $\Sigma$-structure, so 
we must have $\Sigma(n)\cong\st1$, i.e.\ $\Sigma=\st1$. 
The implications (i)$\Rightarrow$(iii')$\Rightarrow$(iii) are also evident. 
(iii)$\Rightarrow$(iii'): If we have $\epsilon_X(x)=\epsilon_X(y)$ for some 
$x\neq y\in X$, consider the map $f:X\to\st2$, such that $f(x)=1$ and 
$f(z)=2$ otherwise; applying $\Sigma(f)$ to our equality we get (iii'). 
(iii')$\Rightarrow$(ii'): 
Since $\{1\}$ and $\{2\}$ act on any $\Sigma$-module~$X$ 
by means of corresponding projections $X^2\to X$, we see that 
$\{1\}=\{2\}$ implies the equality of these projections, hence 
$X$ consists of at most one element; it cannot be empty since $\Sigma(0)\neq
\emptyset$. Now (i)$\Rightarrow$(iv) is again trivial, and 
(iv)$\Rightarrow$(ii') follows from the fact that if $\bu=[c]_{|\Sigma|}$ 
for some constant~$c$, then for any $\Sigma$-module~$X$ 
the identity map $\id_X={[\bu]}_X$ 
has to coincide with the constant map which maps all elements of~$X$ 
into the value $c_X\in X$ of the constant, hence~$X=\{c_X\}$ is an 
one-element set.
\end{Proof}

\nxpointtoc{Algebraic submonads and strict quotients}\label{sp:projlim.algmon}
Projective limits of algebraic monads are computed in the category of 
algebraic endofunctors (cf.\ \ptref{sp:proj.lim.alg} 
and \ptref{sp:max.alg.subf}), i.e.\ 
$(\projlim\Sigma_\alpha)(n)=\projlim(\Sigma_\alpha(n))$. In particular, 
$\rho:\Sigma'\to\Sigma$ is a monomorphism of algebraic monads iff 
all maps $\rho_n:\Sigma'(n)\to\Sigma(n)$ are injective. We see that 
the subobjects of $\Sigma$ in the category of algebraic monads 
are exactly the {\em algebraic submonads} $\Sigma'\subset\Sigma$, 
i.e.\ those algebraic subfunctors which admit a (necessarily unique) 
monad structure compatible with that of~$\Sigma$. According 
to~\ptref{sp:lim.alg.funct}, such a subfunctor is determined by 
a collection of subsets $\Sigma'(n)\subset\Sigma(n)$, stable under  
maps $\Sigma(\phi):\Sigma(m)\to\Sigma(n)$ for all $\phi:\stm\to\stn$. 
Clearly, an algebraic subfunctor $\Sigma'\subset\Sigma$ is a submonad 
iff $\bu\in\Sigma'(1)\subset\Sigma(1)$ and 
$\mu_n^{(k)}\bigl(\Sigma'(k)\times\Sigma'(n)^k\bigr)\subset\Sigma'(n)$ for all 
$n$, $k\geq0$. For example, if we put $\Sigma_+(0):=\emptyset$, 
$\Sigma_+(n):=\Sigma(n)$ for $n\geq1$, we obtain an algebraic submonad 
$\Sigma_+\subset\Sigma$, which is the largest algebraic submonad of~$\Sigma$ 
without constants. Notice that since an algebraic submonad 
$\Sigma'\subset\Sigma$ is completely 
determined by the subset $\|\Sigma'\|:=\bigsqcup_{n\geq0}\Sigma'(n)
\subset\|\Sigma\|$, all algebraic submonads of 
a fixed monad~$\Sigma$ constitute a set.

\nxsubpoint\label{sp:img.algmon}
Notice that the image $\rho(\Sigma)\subset\Xi$ of a homomorphism of 
algebraic monads $\rho:\Sigma\to\Xi$, computed in the category of 
algebraic functors (cf.~\ptref{sp:ker.im.morph.algfun} 
and~\ptref{sp:expl.algmon.cond}d)), is actually an algebraic submonad of 
$\Xi$, not just an algebraic subfunctor, and we obtain a decomposition 
$\rho:\Sigma\twoheadrightarrow\rho(\Sigma)\hookrightarrow\Xi$ of 
$\rho$ into a surjective monad homomorphism followed by an inclusion 
of a submonad. In this case the kernel $R:=\Sigma\times_\Xi\Sigma$ 
in the category of algebraic monads coincides with the same kernel 
computed in the category of (algebraic or all) endofunctors, hence 
$\Sigma\to\rho(\Sigma)$ is a strict epimorphism in all these categories, 
and $\rho(\Sigma)\cong \Sigma/R$. Recall that 
$(\rho(\Sigma))(n)=\rho_n(\Sigma(n))\subset\Xi(n)$ 
(cf.~\ptref{sp:ker.im.morph.algfun}). We say that $\rho(\Sigma)$ 
is the {\em image\/} of $\rho:\Sigma\to\Xi$; clearly, 
$\Sigma\to\rho(\Sigma)\to\Xi$ is the (necessarily unique) decomposition 
of~$\rho$ into a strict epimorphism followed by a monomorphism.

\nxsubpoint\label{sp:gen.submon} (Submonad generated by a set of operations.)
Given a monad $\Sigma$ and any subset $U=\bigsqcup_{n\geq0} U_n$ 
of the set $\|\Sigma\|=\bigsqcup_{n\geq0}\Sigma(n)$ of all operations 
of~$\Sigma$, we can consider the {\em submonad\/~$\Sigma'$ of\/~$\Sigma$, 
generated by~$U$}, i.e.\ the smallest submonad~$\Sigma'$ of~$\Sigma$, 
which contains $U$ (more precisely, such that $U\subset\|\Sigma'\|$). 
Such a submonad always exists: indeed, we can simply take the intersection 
$\bigcap_\alpha\Sigma_\alpha$ 
of all submonads $\Sigma_\alpha\subset\Sigma$ which contain~$U$. 
Recall that projective limits of algebraic monads are computed 
componentwise, hence $\Sigma'(n)=(\bigcap_\alpha\Sigma_\alpha)(n)=
\bigcap_\alpha\Sigma_\alpha(n)\subset\Sigma(n)$. We usually denote 
this submonad $\Sigma'$ by $\langle U\rangle$ or $\Fempty\langle U\rangle$; 
when $U$ is a finite set $\{u_1,\ldots,u_n\}$, we also write 
$\Fempty\langle u_1,u_2,\ldots,u_n\rangle$; if we want to emphasize the 
arities $r_i:=r(u_i)$ of the $u_i$, we write 
$\Fempty\langle u_1^{[r_1]},\ldots, u_n^{[r_n]}\rangle$.

If $\langle U\rangle=\Sigma$, we say that {\em $U$ generates $\Sigma$}, 
or that {\em $U$ is a system of generators of\/~$\Sigma$}; by definition, 
this means that any algebraic submonad of~$\Sigma$ containing~$U$ 
is necessarily equal to~$\Sigma$. This implies the usual consequences: 
for example, if two homomorphisms $\rho_1$, $\rho_2:\Sigma\to\Xi$ 
coincide on~$U$, then they are equal, since $\Ker(\rho_1,\rho_2)$ 
is an algebraic submonad of~$\Sigma$ containing~$U$. Applying this 
to $\Xi=\END(X)$ we see that a $\Sigma$-structure on a set~$X$ is 
completely determined by the collection of maps 
${[u]}_X:X^{r(u)}\to X$ for $u\in U$.

\nxsubpoint\label{sp:expl.descr.gen.subm} 
(Explicit description of~$\Fempty\langle U\rangle$.) 
Let $\Sigma'=\Fempty\langle U\rangle\subset\Sigma$ be the 
submonad of~$\Sigma$ generated by some set $U\subset\|\Sigma\|$. 
Notice that the subsets $\Sigma'(n)\subset\Sigma(n)$ have to 
satisfy the following two requirements:
\begin{itemize}
\item[1)] All ``variables'' 
$\{k\}_\stn\in\Sigma(n)$, $1\leq k\leq n$, are contained in $\Sigma'(n)$.
\item[2)] If $u\in U_k$ is any generator of arity~$k$, and 
$t_1$, \dots, $t_k\in\Sigma'(n)$ for some $k\geq0$, then 
${[u]}_{\Sigma(n)}(t_1,\ldots,t_k)$ also belongs to $\Sigma'(n)$.
\end{itemize}
Conversely, consider the subsets $\Sigma''(n)\subset\Sigma(n)$ of operations 
which can be obtained by applying the above rules a finite number of times. 
Clearly, $U_n\subset\Sigma''(n)\subset\Sigma'(n)$. 
In fact, this collection actually defines a submonad of~$\Sigma$ 
containing~$U$, hence $\Sigma'=\Sigma''$ and we obtain a description 
of $\Sigma'=\Fempty\langle U\rangle$. Indeed, we have $\bu=\{1\}_\st1\in
\Sigma''(1)$, and the substitution property 
``${[t]}_{\Sigma(n)}(t_1,\ldots,t_k)\in\Sigma''(n)$ whenever 
$t\in\Sigma''(k)$, $t_i\in\Sigma''(n)$'' is shown starting from requirement~2) 
by induction on the number of applications of rules~1) and 2) used for 
constructing~$t$. The remaining property ``$\phi_*t=(\Sigma(\phi))(t)\in 
\Sigma''(n)$ for any $t\in\Sigma''(m)$ and $\phi:\stm\to\stn$'' follows 
from the substitution property since $\phi_*t={[t]}_{\Sigma(n)}
(\{\phi(1)\}_\stn,\ldots,\{\phi(m)\}_\stn)$.

\nxsubpoint\label{sp:rel.gen.algmon}
The above definitions admit a useful generalization. Suppose that 
we are given some algebraic monad homomorphism $\rho:\Sigma_0\to\Sigma$, 
and some set $U\subset\|\Sigma\|$. Then we denote by 
$\Sigma_0\langle U\rangle$ the smallest algebraic submonad of $\Sigma$ 
that contains both $\rho(\Sigma_0)$ and $U$; we say that this is 
{\em the submonad\/ of~$\Sigma$ generated over~$\Sigma_0$ by~$U$}.
Clearly, 
$\Sigma_0\langle U\rangle=\bigl<U\cup\|\rho(\Sigma_0)\|\bigr>$, and 
$\rho:\Sigma_0\to\Sigma$ factorizes through $\Sigma_0\langle U\rangle\subset 
\Sigma$. Of course, we write $\Sigma_0\langle u_1,\ldots, u_n\rangle$ 
when $U$ is a finite set. This notion has similar properties to those 
considered above. For example, if $\Sigma$ is generated by~$U$
over~$\Sigma_0$, a $\Sigma$-structure on any set~$X$ 
is completely determined by its scalar restriction 
along~$\rho:\Sigma_0\to\Sigma$ 
together with the action of the generators ${[u]}_X:X^{r(u)}\to X$, 
$u\in U$.

\nxsubpoint\label{sp:fingen.un.algmon} (Finite generation and pre-unarity.)
Given an algebraic monad homomorphism $\rho:\Sigma\to\Xi$, we say that 
{\em $\Xi$ is finitely generated over~$\Sigma$} if there is a finite set 
$U\subset\|\Xi\|$, such that $\Xi=\Sigma\langle U\rangle$. In this case 
we can choose a finite subsystem of generators~$V_0\subset V$
from any other system~$V$ of generators of~$\Xi$ over~$\Sigma$: 
indeed, according 
to~\ptref{sp:expl.descr.gen.subm}, any operation $t\in\|\Xi\|$ 
can be expressed in terms of a finite number of operations from $\rho(\Sigma)$
and $V$; let's denote by $V_0(t)$ the finite set of involved operations from 
$V$, and put $V_0:=\bigcup_{u\in U}V_0(u)$.

We say that {\em $\Xi$ is generated over~$\Sigma$ by operations of 
arity~$\leq r$}, or that {\em $\rho:\Sigma\to\Xi$ is of arity $\leq r$} 
if $\Xi=\Sigma\langle V\rangle$ for some subset 
$V\subset\|\Xi\|_{\leq r}:=\bigsqcup_{n\leq r}\Xi(n)\subset\|\Xi\|$; clearly, 
in this case we can always take $V=\|\Xi\|_{\leq r}$. On the other hand, 
if at the same time $\Xi$ is finitely generated over~$\Sigma$, we can 
find a finite system of generators of arity $\leq r$.

Finally, we say that {\em $\Xi$ is pre-unary over~$\Sigma$}, or that 
{\em $\Xi$ is a pre-unary extension of~$\Sigma$}, or that 
{\em $\rho:\Sigma\to\Xi$ is a pre-unary homomorphism} if 
$\Xi=\Sigma\langle V\rangle$ for some set of unary operations 
$V\subset\Xi(1)$. In this case we can again take $V=|\Xi|=\Xi(1)$; 
if $\Xi$ is finitely generated over~$\Sigma$, we can choose 
a finite set of unary generators $V\subset|\Xi|$.

\nxsubpoint (Examples.)
For example, $\bbZ=\Sigma_\bbZ$ 
is finitely generated (absolutely, i.e.\ over $\Fempty$) 
by operations of arity $\leq 2$:
we can choose $0\in\Sigma_\bbZ(0)$, $[-]:=-1\in\bbZ=\Sigma_\bbZ(1)$ and 
$[+]:=(1,1)\in\bbZ^2=\Sigma_\bbZ(2)$ as a system of generators of~$\bbZ$. 
Any ring~$R$ is pre-unary over~$\bbZ$; we'll show later they are actually 
{\em unary\/} and that all unary 
extensions of~$\bbZ$ come from associative rings. An associative ring~$R$ 
is (absolutely) finitely generated iff it is finitely generated over~$\bbZ$ 
in the usual sense. Notice that $\bbZ$ and $\Zinfty$ are {\em not\/} 
pre-unary extensions of $\Fone$ or $\Fpm$, but are generated in arity $\leq2$. 
Another interesting example: while $R\times R$ ($=\Sigma_R\times\Sigma_R$) 
is a unary extension of~$R=\Sigma_R$ for any associative ring~$R$, 
algebraic monad $\Zinfty\times\Zinfty$ is {\em not\/} a pre-unary extension 
of~$\Zinfty$.

\nxsubpoint\label{sp:strquot.algmon} 
(Strict quotients and equivalence relations.) 
Recall that the cokernels and the fibered products are computed 
in~$\cA_{alg}$ componentwise, hence a morphism $\rho:\Sigma\to\Xi$ in 
$\cA_{alg}\cong\catSets^\catN$ is a strict epimorphism iff all its 
components $\rho_n:\Sigma(n)\to\Xi(n)$ are strict epimorphisms in 
$\catSets$, i.e.\ surjective maps.
We see that the strict quotients of an algebraic endofunctor~$\Sigma$ 
are given by ``surjective'' maps $\rho:\Sigma\to\Xi$, i.e.\ 
natural transformations~$\rho$, such that all $\rho_n:\Sigma(n)\to\Xi(n)$ 
are surjective, or equivalently, such that $\|\rho\|:\|\Sigma\|\to\|\Xi\|$ 
is surjective. Such (algebraic) strict quotients are in one-to-one 
correspondence with algebraic equivalence relations 
$R\subset\Sigma\times\Sigma$, given by some systems of equivalence relations 
$R(n)\subset\Sigma(n)\times\Sigma(n)$, compatible with all maps~$\Sigma(\phi)$
as usual (cf.~\ptref{sp:ker.im.morph.algfun}); 
in this case $\Xi=\Sigma/R$ and $\Xi(n)=\Sigma(n)/R(n)$.

Now suppose~$\Sigma$ has an algebraic monad structure. We say that its 
strict quotient $\Sigma/R$ (in $\cA_{alg}$) is {\em compatible} with 
the monad structure of~$\Sigma$, or that equivalence relation~$R$ 
is {\em compatible} with the monad structure of~$\Sigma$, if 
$\Xi=\Sigma/R$ admits a (necessarily unique) monad structure, for which 
$\Sigma\to\Sigma/R$ becomes a monad homomorphism.

Clearly, in this case $\Sigma/R$ is a quotient of $\Sigma$ modulo~$R$ 
in the category of algebraic monads as well, so $\Sigma\to\Sigma/R$ is 
still a strict epimorphism in this category. Conversely, an 
algebraic monad homomorphism $\rho:\Sigma\to\Xi$ is a strict epimorphism 
iff $\rho$ is surjective, i.e.\ is of the form just discussed, and then 
$\Xi\cong\Sigma/R$ in both categories for some compatible algebraic 
equivalence relation~$R$ on~$\Sigma$. To see this we just consider 
the decomposition of any algebraic monad homomorphism 
$\rho:\Sigma\to\Xi$ into a surjection $\tilde\rho:\Sigma\to\rho(\Sigma)$ 
and an embedding of a submonad $i:\rho(\Sigma)\to\Xi$ 
(cf.~\ptref{sp:img.algmon}), and $R:=\Sigma\times_\Xi\Sigma$ is the 
kernel of both $\rho$ and $\tilde\rho$, hence it is necessarily 
a compatible equivalence relation, and $\rho(\Sigma)\cong\Sigma/R$ both 
as a monad and as an algebraic endofunctor. This also shows that 
$\rho$ is a strict epimorphism iff $\rho(\Sigma)=\Xi$, i.e.\ iff 
$\rho$ is surjective.

\nxsubpoint
When an algebraic equivalence relation $R\subset\Sigma\times\Sigma$ on 
an algebraic monad is compatible with the monad structure?
First of all, it is given by a collection of equivalence relations 
$R(n)\subset\Sigma(n)\times\Sigma(n)$, compatible with maps 
$\phi_*=\Sigma(\phi):\Sigma(m)\to\Sigma(n)$ for all $\phi:\stm\to\stn$; 
in other words, 
\begin{multline}\label{eq:cond.comp.eqrel0}
t\equiv_{R(m)} t'\text{ implies }\phi_*(t)\equiv_{R(n)}\phi_*(t')\\
\text{for any $\phi:\stm\to\stn$, and any $t$, $t'\in\Sigma(m)$.}
\end{multline}
Clearly, $\Xi:=\Sigma/R$ admits a compatible monad structure iff 
all the maps $\Sigma(k)\times\Sigma(n)^k\stackrel{\mu_n^{(k)}}\to
\Sigma(n)\to\Xi(n)=\Sigma(n)/R(n)$ factorize through the corresponding 
canonical projections $\Sigma(k)\times\Sigma(n)^k\to\Xi(k)\times\Xi(n)^k$. 
This means the following:
\begin{multline}\label{eq:cond.comp.eqrel1}
{[t]}_{\Sigma(n)}(t_1,\ldots,t_k)\equiv_{R(n)}
{[t']}_{\Sigma(n)}(t'_1,\ldots,t'_k)\\
\text{whenever $t\equiv_{R(k)}t'$ and $t_i\equiv_{R(n)}t'_i$ 
for $1\leq i\leq k$}
\end{multline}
Notice that this condition actually implies \eqref{eq:cond.comp.eqrel0} 
in view of the formula $\phi_*t={[t]}_{\Sigma(n)}(\{\phi(1)\},\ldots,
\{\phi(m)\})$, hence it is sufficient for the collection 
of equivalence relations $\{R(n)\subset\Sigma(n)\times\Sigma(n)\}_{n\geq0}$ to 
define a compatible algebraic equivalence relation on~$\Sigma$. 

\nxsubpoint\label{sp:eqrel.gen.by.eqs} 
(Compatible equivalence relations generated by a set of equations.)
First of all, notice that the intersection $R=\bigcap_\alpha R_\alpha$ 
of any family of compatible algebraic equivalence relations $R_\alpha$ 
on an algebraic monad~$\Sigma$ is again compatible. This can be either checked 
directly from~\eqref{eq:cond.comp.eqrel1} and 
$R(n)=\bigcap_\alpha R_\alpha(n)$, or deduced from the fact that $R$ 
is the kernel of $\Sigma\to\prod_\alpha \Sigma/R_\alpha$.

This means that, given any set of ``equations'' $E=\bigsqcup_{n\geq0}E_n$, 
$E_n\subset\Sigma(n)\times\Sigma(n)$, we can find the smallest compatible 
algebraic equivalence relation $R=\langle E\rangle$ on~$\Sigma$, which
contains~$E$, by simply taking the intersection of all such equivalence 
relations. Alternatively, we might consider all statements of form 
``$t_1\equiv t_2$'', which can be obtained by a finite number of 
applications of the following rules:
\begin{itemize}
\item[0)] $t\equiv t'$ if $(t,t')\in E$;
\item[1)] $t\equiv t$ for any $t\in\Sigma(n)$;
\item[2)] $t\equiv t'$ implies $t'\equiv t$;
\item[3)] $t\equiv t'$ and $t'\equiv t''$ imply $t\equiv t''$;
\item[4)] The substitution rule~\eqref{eq:cond.comp.eqrel1}.
\end{itemize}
We say that $\langle E\rangle$ is {\em the (compatible algebraic) 
equivalence relation on~$\Sigma$ generated by (the equations, or relations 
from)~$E$}.

Clearly, an algebraic monad homomorphism $\rho:\Sigma\to\Xi$ factorizes 
through $\Sigma\to\Sigma/\langle E\rangle$ iff the kernel 
$R:=\Sigma\times_\Xi\Sigma$ contains~$E$, i.e.\ iff all the equations from 
$E$ are fulfilled in~$\Xi$: $\rho_n(t)=\rho_n(t')$ whenever $(t,t')\in E_n$. 
Putting here $\Xi:=\END(X)$ and taking into account the surjectivity of 
$\Sigma\to\Sigma/\langle E\rangle$, we see that 
$\catMod{\Sigma/\langle E\rangle}$ can be identified with the full 
subcategory of $\catMod\Sigma$, consisting of all $\Sigma$-modules~$X$, 
such that ${[t]}_X={[t']}_X$ whenever $(t,t')\in E$.

When $E$ is a finite set of pairs: $E=\{(t_1,u_1),\ldots,(t_s,u_s)\}$, 
we also write $\langle t_1=u_1,\ldots, t_s=u_s\rangle$ instead of 
$\langle E\rangle$.

\nxsubpoint\label{sp:fingen.unary.eqrel} 
(Finitely generated and unary equivalence relations.)
We introduce for compatible equivalence relations $R$ on algebraic monads 
$\Sigma$ and corresponding strict quotients $\Xi:=\Sigma/R$ a terminology 
similar to that of~\ptref{sp:fingen.un.algmon}. For example, we say that 
$R$ is {\em finitely generated\/}, or that {\em $\Xi$ is obtained from 
$\Sigma$ by imposing finitely many relations (or conditions, or equations)}, 
if $R$ is generated by some finite set of equations~$E$. We say that 
$R$ is {\em generated in arity $\leq k$}, or {\em generated by equations or 
relations of arity $\leq k$}, or that {\em $\Xi$ is obtained from $\Sigma$ 
by imposing relations in arity $\leq k$}, if $R$ is generated by some 
subset $E\subset\bigsqcup_{0\leq n\leq k}\Sigma(n)\times\Sigma(n)$; 
if $R$ is also finitely generated, then such a set of equations can be chosen 
to be finite. Finally, we say that $R$ is {\em unary\/} if it is generated 
by some subset of $\Sigma(1)\times\Sigma(1)$.

\nxpointtoc{Free algebraic monads}\label{p:free.algmon}
Now we would like to construct a very important class of algebraic monads, 
namely, the {\em free algebraic monads.} For this purpose we 
fix some ``graded set'' $U=\bigsqcup_{n\geq0}U_n$, consider the 
category of couples $(\Sigma,\phi)$, where $\Sigma$ is an algebraic monad 
and $\phi=\bigsqcup_{n\geq0}\phi_n:U\to\|\Sigma\|$ is a graded map of sets 
(i.e.\ essentially a collection of maps $\phi_n:U_n\to\Sigma(n)$), 
and define the {\em free algebraic monad\/~$\Fempty\langle U\rangle$} 
to be the initial object of this category. 

\nxsubpoint
By definition, 
we must have a graded map $j:U\to\|\Fempty\langle U\rangle\|$, such that 
any graded map $\phi:U\to\Sigma$ induces a uniquely determined monad 
homomorphism $\rho:\Fempty\langle U\rangle\to\Sigma$, for which 
$\phi=\|\rho\|\circ j$. In other words, monad homomorphisms 
$\rho:\Fempty\langle U\rangle\to\Sigma$ are in one-to-one correspondence 
with graded maps $\phi:U\to\|\Sigma\|$, i.e.\ with families of maps 
$\{\phi_n:U_n\to\Sigma(n)\}_{n\geq0}$. Applying this to 
$\Sigma=\END(X)$, we see that $\catMod{\Fempty\langle U\rangle}$ 
is equivalent to the category of sets $X$, equipped with an 
arbitrary family of maps ${[u]}_X:X^{r(u)}\to X$, defined for each 
$u\in U$, where $r:U\to\bbZ_{\geq0}$ is the arity map (i.e.\ 
$r(u)=n$ iff $u\in U_n$).

\nxsubpoint
Of course, we have to show the existence of free algebraic monads 
$\Fempty\langle U\rangle$. For this purpose we construct such a monad 
as a certain submonad $W_{U,r}$ of the monad $W_U$ of words with constants 
from~$U$ (cf.~\ptref{sp:words.with.const}). Recall that 
$W_U(X)=W(U\sqcup X)$ consists of all words in alphabet $U\sqcup X$, 
usually written in form $u_1\{x_1\}\{x_2\}u_2\{x_3\}$ (we omit 
braces around constants, i.e.\ letters from~$U$). Informally speaking, 
$\Fempty\langle U\rangle=W_{U,r}\subset W_U$ is constructed as the 
set of all valid expressions, constructed with the aid of operations 
from~$U$ from the variables from~$X$, and written in the prefix notation. 
For example, if $r(u_1)=1$ and $r(u_2)=2$, then 
$\{x\}$, $u_1\{z\}$ and $u_2u_1\{x\}u_2\{y\}\{z\}$ are elements of 
$W_{U,r}(\{x,y,z\})$. We present a more formal construction below, 
which has the advantage of allowing a straightforward generalization to 
the topos case.

\nxsubpoint
Let us denote by $\Synt$ the set of all maps $\phi$ from the set 
$\{-\infty\}\cup\bbZ_{\geq0}$ into itself, which satisfy the following 
conditions: 1) $\phi(-\infty)=-\infty$; 2) $\phi(x)\neq-\infty$ for at least  
one $x$; 3) $\phi(x+1)=\phi(x)+1$ whenever $\phi(x)\neq-\infty$. Clearly, 
this set is closed under composition and contains the identity map, 
hence $\Synt$ is a monoid. We call $\Synt$ the {\em syntax monoid}. 
Its elements $\phi$ are classified by pairs of integers $r$, $s\geq0$, 
where $r$ is the smallest $x\geq0$ for which $\phi(x)\geq0$, and 
$s=\phi(r)$. We denote the corresponding element of $\Synt$ by~$G_{r,s}$. 
Notice that $G_{0,s}$ is the only element $\phi\in\Synt$, such that 
$\phi(0)=s$. In particular, $G_{0,0}$ is the identity of~$\Synt$. 
We can compute $G_{r,s}G_{r',s'}$ explicitly: if $r\leq s'$, it is equal to 
$G_{r',s+s'-r}$, otherwise it equals~$G_{r+r'-s',s}$.

\nxsubpoint
Now we define the {\em verification maps} $d_X:W_U(X)=W(U\sqcup X)\to\Synt$ 
as follows. Recall that $W_U(X)$ is a free monoid generated by $U\sqcup X$, 
so $d_X$ is completely determined by its values on generators $U\sqcup X$, 
if we also require it to be a monoid homomorphism. We put 
$d_X(u):=G_{r(u),1}$ for any $u\in U$, i.e.\ $d_X(u)=G_{n,1}$ for $u\in U_n$, 
and $d_X(x):=G_{0,1}$ for any $x\in X$; then $d_X(z_1z_2\ldots z_s)=
d_X(z_1)d_X(z_2)\ldots d_X(z_s)$ for any $z_i\in U\sqcup X$.

Put $W_{U,r}(X):=d_X^{-1}(G_{0,1})\subset W_U(X)$, i.e.\ this is the set of 
all words $t\in W_U$, such that $d_X(t)=G_{0,1}$, or equivalently 
such that $(d_X(t))(0)=1$; we call such words~$t$ {\em valid expressions} 
or {\em terms} with respect to the set of operations $U$.
Clearly, $d_Y\circ W_U(f)=d_X$ for 
any map of sets $f:X\to Y$, hence $(W_U(f))(W_{U,r}(X))\subset W_{U,r}(Y)$, 
so $W_{U,r}$ is indeed a subfunctor of~$W_U$, and even an algebraic subfunctor,
since it commutes with filtered inductive limits.

\nxsubpoint
We have to check that $W_{U,r}\subset W_U$ is indeed a submonad, i.e.\ 
that if we replace in some term~$t$ some variables $\{x_i\}$ 
(or even $\{t_i\}$) with some other terms~$t_i$, then 
the resulting word~$\tilde t$ is necessarily a term. This is clear, since 
the expression for $d(\tilde t)\in\Synt$ is obtained from that of $d(t)$ 
by replacing $d(\{x_i\})$ with $d(t_i)$; but $d(t_i)=G_{0,1}=d(\{x_i\})$ 
since all $t_i$ are terms, hence $d(\tilde t)=d(t)=G_{0,1}$, hence 
$\tilde t$ is indeed a term.

\nxsubpoint\label{sp:struct.ind} (Structural induction.)
Notice that we have some general rules for constructing new terms:
\begin{itemize}
\item[1)] Any variable $x\in X$ defines a one-letter term 
$\{x\}\in W_{U,r}(X)$.
\item[2)] If $u\in U_n$ is an operation of arity~$n$ and $t_1$, \dots, $t_n$ 
are terms, then $u\,t_1\,\ldots\,t_n$ is also a term.
\end{itemize}
Conversely, one checks in the usual way that any term $t\in W_{U,r}(X)$ is 
necessarily non-empty, and that it is either of form 1), when the first 
letter of~$t$ belongs to~$X$, or of form 2) with uniquely determined 
$u\in U$ and terms of smaller length $t_1$, \dots, $t_n\in W_{U,r}(X)$, 
when the first letter of~$t$ is an operation from~$U$; 
moreover, if $t=u\,t_1\ldots\,t_n$ for some operation $u\in U$ and some 
terms $t_i$, then necessarily $n=r(u)$, and such a list of $t_i$ is unique.

The fact that any term of~$W_{U,r}(X)$ can be uniquely written 
either in form 1) or 2) allows us to prove statements by induction on length 
of a term~$t$, proving them first for case $t=\{x\}$, and then proving them 
for $t=u\,t_1\,\ldots\,t_n$, assuming them to be already proved for all $t_i$; 
this sort of induction is called {\em structural induction}.

\nxsubpoint
In particular, for any $u\in U_n$ we get an element 
$j(u)=j_n(u):=u\{1\}\{2\}\ldots\{n\}\in W_{U,r}(n)$, so we get a graded map 
$j:U\to\|W_{U,r}\|$. Since $j$ is injective, we often identify $U$ 
with the subset $j(U)\subset\|W_{U,r}\|$ and write $\langle u\rangle$ or 
even just~$u$ instead of~$j(u)$. 
Notice that
\begin{equation}\label{eq:action.free.gen}
{\bigl[j(u)\bigr]}_{W_{U,r}(X)}(t_1,t_2,\ldots,t_n)=u\,t_1\,t_2\,\ldots\,t_n
\text{ for any $u\in U_n$, $t_i\in W_{U,r}(X)$.}
\end{equation}
This shows that $j(U)\subset\|W_{U,r}\|$ generates $W_{U,r}$ 
(cf.\ \ptref{sp:expl.descr.gen.subm}).

\nxsubpoint
Now let's check that $(W_{U,r},j)$ satisfies the universal property required 
from~$\Fempty\langle U\rangle$, thus proving the existence of free algebraic 
monads. We have to show that for any graded map $\phi:U\to\|\Sigma\|$ there 
is a unique monad homomorphism $\rho:W_{U,r}\to\Sigma$, such that 
$\phi=\|\rho\|\circ j$, i.e.\ we require $\rho_n(j_n(u))=\phi_n(u)$. 
Taking \eqref{eq:action.free.gen} into account, we obtain
\begin{multline}\label{eq:ext.to.freemon}
\rho_X(u\, t_1\,\ldots\, t_n)={[\phi(u)]}_{\Sigma(X)}
\bigl(\rho_X(t_1),\ldots,\rho_X(t_n)\bigr)\\
\text{ for any $u\in U_n$, $t_i\in W_{U,r}(X)$.}
\end{multline}
Since we must also have $\rho_X(\{x\})=\{x\}$ for any $x\in X$, we 
prove by structural induction in $t\in W_{U,r}(X)$ the existence and 
uniqueness of such maps $\rho_X:W_{U,r}(X)\to\Sigma(X)$. They are clearly 
functorial in~$X$, so we get a natural transformation of algebraic 
functors $\rho:W_{U,r}\to\Sigma$. It remains to prove that $\rho$ is indeed a 
monad homomorphism; we have to check~\eqref{eq:cond.algmon.hom} for this, 
and this is easily shown by structural induction in~$t$, 
using~\eqref{eq:ext.to.freemon} for the induction step.

\nxsubpoint
We have already seen that $\Fempty\langle U\rangle=W_{U,r}$ is generated 
by $U\cong j(U)\subset\|W_{U,r}\|$. This justifies the notation 
$\Fempty\langle U\rangle$. Moreover, if we have any subset 
$U\subset\|\Sigma\|$, it induces a canonical homomorphism 
$\rho:\Fempty\langle U\rangle\to\Sigma$, and its image is exactly 
the submonad $\Sigma'$ of $\Sigma$ generated by~$U$, which has been also 
denoted before by $\Fempty\langle U\rangle$. In particular, $U$ generates 
$\Sigma$ iff $\rho:\Fempty\langle U\rangle\to\Sigma$ is surjective, 
i.e.\ a strict epimorphism. 
If it is an isomorphism, we say that {\em $U$ freely generates~$\Sigma$}, 
or that {\em $U$ is a system of free generators of~$\Sigma$}. 
For example, the free algebraic monad $\Fempty\langle U\rangle$ is 
freely generated by $j(U)\cong U$.

\nxsubpoint
In general we tend to use capital letters for the elements of $U$ when 
$\Sigma$ is freely generated by~$U$, to distinguish this situation from 
the case when $U$ is just any system of generators of~$\Sigma$.
Of course, we follow the conventions of~\ptref{sp:gen.submon} when 
$U$ is a finite set, thus writing for example 
$\Fempty\langle N^{[0]},T^{[1]},U^{[2]}\rangle$ to denote the free monad 
generated by one constant, one unary, and one binary generator.

\nxsubpoint\label{sp:pres.algmon} (Algebraic monads from algebraic systems.)
Now suppose given a graded set of operations $U=\sqcup_{n\geq0}U_n$, 
or equivalently, a set~$U$ together with an arity map $r:U\to\bbZ_{\geq0}$, 
and any set of ``equations'' or ``relations'' 
$E\subset\|\Fempty\langle U\rangle\times\Fempty\langle U\rangle\|$. 
Then we can construct the compatible algebraic equivalence relation 
$\langle E\rangle$ on $\Fempty\langle U\rangle$ generated by~$E$, 
and consider the quotient $\Sigma:=\Fempty\langle U\rangle/\langle E\rangle$. 
This quotient will be denoted $\Fempty\langle U|E\rangle$ 
or even $\langle U|E\rangle$; when~$U$ 
and~$E$ are finite we adopt the conventions of~\ptref{sp:gen.submon} 
and~\ptref{sp:eqrel.gen.by.eqs}, thus writing expressions like 
$\Fempty\langle 0^{[0]}, \zeta^{[1]}\,|\,\zeta 0=0, 
\zeta^n\{1\}=\{1\}\rangle$ (this is actually $\bbF_{1^n}$) or 
$\Fempty\langle U^{[2]}\,|\, U\{1\}\{2\}=U\{2\}\{1\}, 
U\{1\}U\{2\}\{3\}=UU\{1\}\{2\}\{3\}\rangle$. We usually prefer to replace 
the free variables like $\{1\}$, $\{2\}$, \dots, with some (arbitrarily 
chosen) letters like $x$, $y$, \dots, and of course we use the infix notation 
and parentheses when appropriate, thus writing the second of the above 
examples in form $\Fempty\langle+^{[2]}\,|\,x+y=y+x,\,x+(y+z)=(x+y)+z\rangle$. 
Another useful convention: we often write unary equations without 
naming the free variable $\bu=\{1\}_\st1$ explicitly, using the monoid 
structure of $|\Sigma|=\Sigma(1)$ instead. Thus we write 
$\zeta^n=\bu$ instead of $\zeta\zeta\cdots\zeta\{1\}=\{1\}$, and 
$\phi^2=\phi+\bu$ instead of $\phi\phi\{1\}=\phi\{1\}+\{1\}$; 
we can also write $-^2=\bu$ instead of $-(-x)=x$ (here $\zeta$, $\phi$ and 
$-$ are unary, and $+$ is binary). This is possible because of the 
associativity relation $t(u_1,\ldots,u_k)\cdot x=t(u_1(x),\ldots,u_k(x))$, 
true for any $t\in\Sigma(k)$, $u_i\in|\Sigma|$, and for any $x$ from a 
$\Sigma$-module~$X$, e.g.\ for~$\bu\in|\Sigma|$.

We call such couples $(U,E)$ as above {\em algebraic systems}, and 
$\Fempty\langle U|E\rangle$ is the {\em algebraic monad defined by algebraic 
system~$(U,E)$}. If $\Sigma\cong\Fempty\langle U|E\rangle$, we say that 
$(U,E)$ is a {\em presentation of algebraic monad~$\Sigma$}. Of course, 
in this case the image of $U$ in~$\Sigma$ generates~$\Sigma$, and all 
relations from~$E$ are fulfulled in~$\Sigma$, and all other relations between 
operations of~$\Sigma$ can be deduced from~$E$, so we have indeed something 
very similar to the usual description of an algebra in terms of a list of 
generators and relations.

Clearly, algebraic monad homomorphisms $\rho:\Fempty\langle U|E\rangle\to\Xi$ 
are in one-to-one correspondence with graded maps $\phi:U\to\|\Xi\|$, such 
that the image under~$\rho$ of any equation from~$E$ is fulfilled in~$\Xi$. 
Taking here $\Xi=\END(X)$ for some set~$X$, we see that 
$\catMod{\Fempty\langle U|E\rangle}$ is exactly the category of sets~$X$ with 
an algebraic structure of species defined by the algebraic system $(U,E)$. 
In this way the study of any algebraic structure is reduced to the study of 
modules over some algebraic monad. Of course, different algebraic systems can 
correspond to isomorphic algebraic monads.

\nxsubpoint\label{sp:algmon.algsyst}
Conversely, any algebraic monad~$\Sigma$ admits some presentation 
$(U,E)$, since one can always take for $U$ any system of generators 
of~$\Sigma$ (e.g.\ $U=\|\Sigma\|$), and then take for~$E$ any system of 
generating equations for the equivalence relation~$R$ on
the free monad~$\Fempty\langle U\rangle$, defined by the canonical surjection
$\Fempty\langle U\rangle\to\Sigma$ (we can always take $E=R$). In this sense 
the study of categories of sets equipped with some algebraic structures, 
i.e.\ {\em universal algebra}, is nothing else than the study of 
categories of modules over algebraic monads. For example, we 
see immediately that the forgetful functors from any of these categories 
into $\catSets$ is monadic.
However, the category of algebraic monads themselves 
is more convenient for different 
category-theoretic operations (e.g.\ computation of projective limits) than 
the category of algebraic systems (since even the definition of a morphism 
of algebraic systems is quite difficult to handle when written without monads).

\nxsubpoint
The above constructions generalize to the case when we replace $\Fempty$ with 
an arbitrary algebraic monad~$\Sigma_0$. For example, for any graded set 
we can construct an algebraic monad~$\Sigma_0\langle U\rangle$ over $\Sigma_0$ 
(i.e.\ a homomorphism of algebraic monads 
$\Sigma_0\to\Sigma_0\langle U\rangle$), 
which is an initial object in the category of triples $(\Sigma,\rho,\phi)$, 
consisting of an (algebraic) monad~$\Sigma$, a homomorphism 
$\rho:\Sigma_0\to\Sigma$, and a graded map $\phi:U\to\|\Sigma\|$. Then 
the category of $\Sigma_0\langle U\rangle$-modules consists exactly of 
$\Sigma_0$-modules~$X$ together with arbitrarily chosen maps 
${[u]}_X:X^{r(u)}\to X$ for all $u\in U$. This actually shows that the 
canonical map $j:U\to\|\Sigma_0\langle U\rangle\|$ is injective unless 
$\Sigma_0=\st1$ or $\Sigma_0\cong\st1_+\subset\st1$.

We have to show the existence of such free monads $\Sigma_0\langle U\rangle$, 
but this is simple: take any presentation $(U',E')$ of $\Sigma_0$ and put 
$\Sigma_0\langle U\rangle:=\Fempty\langle U'\sqcup U|E'\rangle$. Then 
everything follows from the universal properties of all constructions involved,
or from the observation that in this way we get a correct category of 
$\Sigma_0\langle U\rangle$-modules, using~\ptref{sp:monhom.catfunct} 
to obtain monad homomorphisms from functors between categories of modules.

\nxsubpoint\label{sp:relpres.algmon}
Of course, we can construct a monad $\Sigma_0\langle U|E\rangle=
\Sigma_0\langle U\rangle/\langle E\rangle$ over
an algebraic monad $\Sigma_0$, starting from an arbitrary graded set~$U$ 
and an arbitrary set of relations $E\subset\|\Sigma_0\langle U\rangle^2\|$. 
It has properties similar to those considered before for $\Sigma_0=\Fempty$, 
and we introduce similar terminology and notations. For example, given 
a monad $\Sigma$ over $\Sigma_0$, i.e.\ a homomorphism of algebraic monads 
$\rho:\Sigma_0\to\Sigma$, we say that $(U,E)$ is a {\em presentation of 
$\Sigma$ over~$\Sigma_0$} if $\Sigma\cong\Sigma_0\langle U|E\rangle$. Of 
course, such a presentation always exists; it can found by the same 
reasoning as in~\ptref{sp:algmon.algsyst}.

\nxsubpoint\label{sp:homalgmon.finpres.un} (Finite presentation and unarity.)
Given any algebraic monad $\Sigma$ over $\Sigma_0$, or equivalently a 
homomorphism $\rho:\Sigma_0\to\Sigma$, we say that {\em $\Sigma$ is 
finitely presented over~$\Sigma_0$}, or that {\em $\rho$ is of finite 
presentation}, if there exists a presentation 
$\Sigma\cong\Sigma_0\langle U|E\rangle$ with both $U$ and $E$ finite. 
We say that {\em $\Sigma$ is generated in arity $\leq r$ with relations 
in arity $\leq s$ over~$\Sigma_0$} if we can find a presentation with 
all operations in~$U$ of arity $\leq r$ and all relations in~$E$ of 
arity~$\leq s$. If $r\leq s$, we can replace $U$ by any other system~$U'$ of 
generators of arity $\leq r$, since all the relations expressing  
operations from~$U$ in terms of new generators from~$U'$ will be in arity 
$\leq r\leq s$.

Finally, we say that $\Sigma$ is {\em unary\/} over $\Sigma_0$, if 
$U$ can be chosen inside $|\Sigma|=\Sigma(1)$, and $E$ inside 
$|\Sigma_0\langle U\rangle^2|$. If we have another set $U'$ of unary 
generators of such a $\Sigma$ over $\Sigma_0$, then the induced 
equivalence relation on $\Sigma_0\langle U'\rangle$ will be necessarily 
unary, since we can replace in all unary equations from~$E$ the operations 
from~$U$ with their expressions in terms of unary operations from~$U'$ 
(and maybe some operations of other arities from $\Sigma_0$, but this 
doesn't affect anything).

\nxsubpoint (Examples.)
We know that $\Fone=\Fempty\langle0^{[0]}\rangle$, so it is finitely presented 
over~$\Fempty$, but of course not unary. Similarly, 
$\Fpm=\Fone\langle-^{[1]}\,|\,-(-x)=x$, $-0=0\rangle=
\Fempty\langle0^{[0]},-^{[1]}\,|\,-(-x)=x$, $-0=0\rangle$ is finitely presented
over both $\Fempty$ and $\Fone$, and unary over $\Fone$ (since 
we might replace the equation for constants $-0=0$ with unary equation 
$-0_{\Fone(1)}\{1\}=0_{\Fone(1)}\{1\}$), but not over~$\Fempty$.
A more interesting example: $\bbZ$ 
is finitely presented over~$\Fempty$, since 
$\bbZ=\Fempty\bigl<0^{[0]},-^{[1]},+^{[2]}\;|\;
x+y=y+x$, $x+(-x)=0$, $x+0=x$, $(x+y)+z=x+(y+z)\bigr>$. It is also 
finitely presented over~$\Fone$ and~$\Fpm$, but not unary over any of these 
monads.

\nxsubpoint\label{sp:nc.tensprod} (Non-commutative tensor products.)
Presentations of algebraic monads can be used to show existence of 
``non-commutative tensor products'' $\Sigma_1\boxtimes_{\Sigma}\Sigma_2$, 
i.e.\ pushouts of pairs of morphisms $\rho_i:\Sigma\to\Sigma_i$ in the 
category of algebraic monads. Indeed, we have just to take any 
presentations $\Sigma_i=\Sigma\langle U_i|E_i\rangle$ of~$\Sigma_i$, and 
put $\Sigma_1\boxtimes_\Sigma\Sigma_2:=\Sigma\langle U_1,U_2|E_1,E_2\rangle=
\Sigma\langle U_1\sqcup U_2|E_1\sqcup E_2\rangle$. Alternatively, we might 
take $\Sigma_1\langle U_2|\rho_{1,*}(E_2)\rangle$, where 
$\rho_{1,*}:\|\Sigma_0\langle U_2\rangle^2\|\to
\|\Sigma_1\langle U_2\rangle^2\|$ is the canonical map induced by~$\rho_1$. 
Again, the required universal property of~$\Sigma_1\boxtimes_\Sigma\Sigma_2$ 
follows from the universal properties of monads 
of form $\Sigma\langle U|E\rangle$. We see that 
$\catMod{\bigl(\Sigma_1\boxtimes_\Sigma\Sigma_2\bigr)}$ 
is isomorphic to the category of 
sets~$X$, equipped with both a $\Sigma_1$-structure and a $\Sigma_2$-structure,
restricting to the same $\Sigma$-structures. Actually, since this is a 
category of algebraic systems, we could use it to define 
$\Sigma_1\boxtimes_\Sigma\Sigma_2$ as the monad defined by the forgetful 
functor from this category, applying~\ptref{sp:monhom.catfunct} to construct 
monad homomorphisms from functors between categories of modules.

\nxsubpoint\label{sp:nc.finpres.un.prop} 
(Finite presentation, unarity and pushouts.)
Of course, the notions of finite generation and presentation have the usual 
properties with respect to pushouts, i.e.\ NC-tensor products. They are 
proved essentially in the usual way, well-known for usual associative 
and especially commutative algebras, so we just list them in the order they 
can be proved, omitting the proofs themselves.
\begin{itemize}
\item Finitely generated and finitely presented homomorphisms of algebraic 
monads are stable under pushouts.
\item A strict epimorphism $\Sigma\twoheadrightarrow\Xi$ is of finite 
presentation iff its kernel $R=\Sigma\times_\Xi\Sigma$ is finitely generated 
as an equivalence relation on~$\Sigma$.
\item For any homomorphism $\rho:\Sigma\to\Xi$ the codiagonal 
(or multiplication) map $\nabla:\Xi\boxtimes_\Sigma\Xi\to\Xi$ is a 
strict epimorphism (i.e.\ a surjection); if $\rho$ is of finite type 
(i.e.\ if $\Xi$ is finitely generated over~$\Sigma$), then $\nabla$ is 
of finite presentation.
\item In the situation $\Sigma\to\Sigma'\to\Xi$, if $\Xi$ is finitely 
presented over~$\Sigma$, and $\Sigma'$ is finitely generated over~$\Sigma$, 
then $\Xi$ is finitely presented over~$\Sigma'$.
\item In the same situation if $\Xi$ is finitely generated over~$\Sigma$, 
then it is finitely generated over~$\Sigma'$.
\end{itemize}
Moreover, we have similar statements about unary and pre-unary 
homomorphisms:
\begin{itemize}
\item Pre-unary and unary homomorphisms are stable under pushouts.
\item A strict epimorphism is always pre-unary; it is unary iff its kernel 
is a unary equivalence relation (cf.~\ptref{sp:fingen.unary.eqrel}).
\item If $\Xi$ is pre-unary over $\Sigma$, then the codiagonal map 
$\nabla:\Xi\boxtimes_\Sigma\Xi\to\Xi$ is unary.
\item In the situation $\Sigma\to\Sigma'\to\Xi$, if $\Xi$ is unary 
over~$\Sigma$, and $\Sigma'$ is pre-unary over~$\Sigma$, then 
$\Xi$ is unary over $\Sigma'$.
\item In the same situation, if $\Xi$ is pre-unary over~$\Sigma$, the 
same is true over~$\Sigma'$.
\end{itemize}

\nxsubpoint\label{sp:indlim.algmon} 
(Arbitrary inductive limits of algebraic monads.)
Since $\Fempty$ is an initial object of the category of algebraic monads, 
we see that we have coproducts $\Sigma_1\boxtimes\Sigma_2:=
\Sigma_1\boxtimes_{\Fempty}\Sigma_2$ and cokernels of pairs of morphisms 
in this category as well. Since 
filtered inductive limits of algebraic monads also exist (they can be 
computed componentwise), we deduce first the existence of infinite coproducts, 
and then the existence of arbitrary (small) inductive limits of algebraic 
monads. So presentations and algebraic systems enable us to compute 
inductive limits of algebraic monads, while projective limits and 
submonads are better described in terms of algebraic monads themselves 
(cf.~\ptref{sp:projlim.algmon}). 

\nxsubpoint\label{sp:catalgmon.mon} (Algebraic monads as algebraic structures.)
We have seen that the forgetful functor $\Sigma\mapsto\|\Sigma\|$ from 
the category of algebraic monads into the category of $\bbZ_{\geq0}$-graded 
sets $\catSets_{/\bbZ_{\geq0}}$ 
admits a left adjoint $U\mapsto\Fempty\langle U\rangle$, hence 
we get a monad $\monMon$ over $\catSets_{/\bbZ_{\geq0}}$. Moreover, 
$\monMon:U\mapsto\|\Fempty\langle U\rangle\|$ clearly commutes with filtered 
inductive limits of graded sets, so it is something like an algebraic monad 
over~$\catSets_{/\bbZ_{\geq0}}$. On the other hand, we have seen 
in~\ptref{sp:expl.algmon.cond} that the structure of an algebraic monad 
on a graded set $\|\Sigma\|$ is itself algebraic, so we might expect the 
forgetful functor $\Sigma\to\|\Sigma\|$ to be monadic, i.e.\ the category 
of algebraic monads to be equivalent to $(\catSets_{/\bbZ_{\geq0}})^\monMon$, 
and this is indeed the case.

In fact, we might fix some set~$S$ from the very beginning, 
and consider monads $\Sigma$ over the category $\catSets_{/S}$ of 
$S$-graded sets $X=\bigsqcup_{s\in S}X_s$, or equivalently, 
of maps $X\stackrel{r}\to S$. Algebraic monads and algebraic endofunctors 
$\Sigma:\catSets_{/S}\to\catSets_{/S}$
are defined by the same requirement to commute with filtered inductive limits;
they are completely determined by their restriction to the category 
$\catN_{/S}$ of finite $S$-graded sets. Since $\Ob(\catN_{/S})\cong
\bigsqcup_{n\geq0}S^n=W(S)$, we can index objects of $\catN_{/S}$ by words 
or sequences $(s_1,\ldots,s_n)$ in alphabet~$S$. We see that an 
algebraic endofunctor $\Sigma$ over $\catSets_{/S}$ is given by a collection 
of $S$-graded sets $\Sigma(s_1,\ldots,s_n)=\bigsqcup_{t\in S}
\Sigma(s_1,\ldots,s_n)_t$, parametrized by $(s_1,\ldots,s_n)\in W(S)$, 
and some maps between these sets. A pre-action $\alpha$ of $\Sigma$ on 
some $S$-graded set is given then by a collection of maps 
$\Sigma(s_1,\ldots,s_n)_t\times X_{s_1}\times\cdots\times X_{s_n}\to X_t$, 
and a monad structure on~$\Sigma$ is given by a collection of graded identity 
elements $\bu_s\in\Sigma(s)_s$, $s\in S$, and pre-actions of~$\Sigma$ on 
all $S$-graded sets $\Sigma(s_1,\ldots,s_n)$. Of course, all these data 
are subject to some compatibility conditions similar to those considered 
before, but more cumbersome to write down, so we decided not to adopt 
this approach from the very beginning and treat the simplest case, even if 
almost all statements and constructions generalize to the $S$-graded case.

Notice that the underlying set $\|\Sigma\|=\bigsqcup\Sigma(s_1,\ldots,s_n)_t$ 
of an algebraic monad $\Sigma$ over $\catSets_{/S}$ is itself 
$W(S)\times S$-graded, so the category of such monads can be described itself 
as a category of modules over some monad $\monMon_S$ on the category of 
$W(S)\times S$-graded sets. This observation allows one in principle to deduce 
properties of categories of algebraic monads from the properties of 
categories of modules over algebraic monads, at the cost of making 
everything less explicit.

\nxsubpoint (Graded algebraic monads.)
One might expect that in the case when $S$ is a commutative group or at least 
a commutative monoid, algebraic monads over $\catSets_{/S}$ are something like 
$S$-graded algebraic monads (over $\catSets$). However, this is usually 
not the case. For example, the degree translation functors 
$T_t:X\mapsto X(t)$, $X(t)_s=X_{s+t}$, are usually expected to transform 
graded modules into graded modules, but in our situation we don't obtain 
a $\Sigma$-structure on $X(t)$ from a $\Sigma$-structure on~$X$ unless 
we are given a compatible family of morphisms 
$\Sigma(s_1,\ldots,s_n)_s\to\Sigma(s_1+t,\ldots,s_n+t)_{s+t}$ 
for each $t\in S$;
if $S$ is a group, then all this morphisms clearly have to be isomorphisms.

For example, we can construct some canonical maps 
$\theta:\monMon(r_1,\ldots,r_m)_n=
\Fempty\bigl<\langle1\rangle^{[r_1]},\ldots,\langle m\rangle^{[r_m]}\bigr>(n)
\to\monMon(r_1+1,\ldots,r_m+1)_{n+1}$, $t\mapsto\tilde t$ 
by structural induction in~$t$: 
we map $\{k\}_\stn$ into $\{k\}_\st{n+1}$, and 
$t=\langle k\rangle\,t_1\,\ldots\,t_{r_k}$ into $\tilde t:=
\langle k\rangle\,\tilde t_1,\ldots,\tilde t_{r_k}\,\{n+1\}_\st{n+1}$, i.e.\ 
we add an extra argument $\{n+1\}$ to each operation $\langle k\rangle$. 
Of course, these maps $\theta$ are injective, but not bijective, so 
$\monMon$ is not $\bbZ$-graded; however, we can take the inductive limit 
along all these maps, thus obtaining a $\bbZ$-graded monad $\monMon_+$ 
and a monomorphism $\rho:\monMon\to\monMon_+$. Then the $\monMon_+$-modules 
are something like algebraic monads, but they admit degree translation 
in both directions, and we obtain functors $\rho^*$ and $\rho_*$ between 
such things and usual algebraic monads.

\nxpointtoc{Modules over an algebraic monad}\label{p:algmon.mod}
Let's fix an algebraic monad~$\Sigma$ over $\cC=\catSets$. We want to 
study some basic properties of the category of $\Sigma$-modules 
$\catMod\Sigma=\cC^\Sigma$, as well as the categories of 
(left, right or two-sided) $\Sigma$-modules in $\cA_{alg}$. These latter 
categories of ``algebraic modules'' over $\Sigma$ actually behave themselves 
more like complexes of modules over an ordinary ring, while 
$\catMod\Sigma$ is a close counterpart of category of modules over an 
ordinary (associative) ring. In particular, we are going to prove 
all properties listed and used before 
for $\Sigma=\Zinfty$ in~\ptref{sp:zinfmod.listprop}.

\nxsubpoint\label{sp:projlim.mod} (Projective limits.)
Of course, arbitrary projective limits $\projlim M_\alpha$ exist in the 
category $\catMod\Sigma=\catSets^\Sigma$ of modules over an algebraic 
monad~$\Sigma$, and they are essentially computed in the category of sets 
(cf.~\ptref{sp:submod.proj.lim}). If we identify $(\projlim M_\alpha)^n$ 
with $\projlim M_\alpha^n$, then ${[t]_M}:M^n\to M$ is identified with 
$\projlim {[t]}_{M_\alpha}$, for any $t\in\Sigma(n)$. 

\nxsubpoint\label{sp:submod} (Submodules.)
In particular, $f:N\to M$ is a monomorphism in $\catMod\Sigma$ iff it is 
a monomorphism in $\cC=\catSets$, i.e.\ an injective map. Therefore, 
the subobjects in $\catMod\Sigma$ of a $\Sigma$-module~$M=(M,\eta_M)$ 
are given by the {\em submodules} $N$ of $M$, i.e.\ those subsets 
$N\subset M$, which admit a (necessarily unique) $\Sigma$-structure, 
compatible with that of~$M$, i.e.\ such that the inclusion $i:N\to M$ is a 
$\Sigma$-homomorphism (cf.~\ptref{sp:submod.proj.lim}). Clearly, 
$N\subset M$ is a submodule of~$M$ iff the image of 
$\Sigma(N)\stackrel{\Sigma(i)}\to\Sigma(M)\stackrel{\eta_M}\to M$ is contained 
in~$N$; in this case it has to be equal to~$N$ since 
$\eta_M\circ\epsilon_M=\id_M$.

Since $\Sigma(N)$ consists of expressions $t(\{x_1\},\ldots,\{x_n\})$ with 
$t\in\Sigma(n)$ and $x_i\in N$ (cf.~\ptref{sp:expl.SX}), and 
$\eta_M:\Sigma(M)\to M$ maps such an expression into ${[t]}_M(x_1,\ldots,x_n)$,
we see that $N\subset M$ is a submodule of~$M$ iff ${[t]}_M(N^n)\subset N$ 
for any $n\geq0$, $t\in\Sigma(n)$, i.e.\ $N$ has to be stable under all 
operations of~$\Sigma$. Actually, the set of all operations of~$\Sigma$, 
under which $N$ is stable, forms an algebraic submonad of~$\Sigma$ 
for any subset $N\subset M$, hence it is sufficient to check the stability 
of~$N$ under a set of generators of~$\Sigma$.

\nxsubpoint\label{sp:module.of.maps} (Module structure on sets of maps.)
Given any set~$S$ and any $\Sigma$-module~$M$, we obtain a canonical 
$\Sigma$-structure on the set $H:=\Hom(S,M)=\Hom_\catSets(S,M)$ of all maps 
from~$S$ to~$M$ (cf.~\ptref{sp:sigma.str.hom}). 
This structure is nothing else than the product $\Sigma$-structure 
on $\Hom_\catSets(S,M)\cong M^S$; the description of projective limits 
of $\Sigma$-modules given in~\ptref{sp:projlim.mod} shows that 
all operations $t\in\Sigma(n)$ act on maps $S\to M$ pointwise, i.e.
\begin{multline}
\bigl({[t]}_{\Hom(S,M)}(f_1,\ldots,f_n)\bigr)(s)=
{[t]}_M\bigl(f_1(s),\ldots,f_n(s)\bigr)\\
\text{for any $t\in\Sigma(n)$, any maps $f_i:S\to M$ and any $s\in S$.}
\end{multline}
When $N$ is another $\Sigma$-module, $\Hom_\Sigma(N,M)$ is a subset of 
$\Hom(N,M)$, but in general this is not a $\Sigma$-submodule of 
$\Hom(N,M)$ with respect to the $\Sigma$-structure just considered, 
at least if we don't suppose $\Sigma$ to be commutative 
(cf.~\ptref{sp:modstr.homsigma}). 

\nxsubpoint\label{sp:modhom.img} (Image of a $\Sigma$-homomorphism.)
Let $f:M=(M,\eta_M)\to N=(N,\eta_N)$ be an arbitrary $\Sigma$-homomorphism. 
Let's denote by 
$I:=f(M)\subset N$ its image as a map of sets, so we get the canonical 
decomposition $f:M\stackrel\pi\twoheadrightarrow I\stackrel i\hookrightarrow N$
of~$f$ into a surjection followed by an embedding. Notice that any 
epimorphism in~$\catSets$ admits a section, hence it is respected by~$\Sigma$,
and we have seen in~\ptref{sp:sigma.injtoinj} that $\Sigma$ preserves 
injectivity of maps as well, so $\Sigma(f)=\Sigma(i)\circ\Sigma(\pi)$ 
is the canonical decomposition of~$\Sigma(f)$. Now 
$\eta_N\circ\Sigma(i)\circ\Sigma(\pi)=i\circ\pi\circ\eta_M$ since 
$f$ is a $\Sigma$-homomorphism, i.e.\ the outer circuit of the following 
diagram is commutative:
\begin{equation}
\xymatrix@C+8pt{
\Sigma(M)\ar[r]_{\Sigma(\pi)}\ar@/^1pc/[rr]^{\Sigma(f)}\ar[d]^{\eta_M}&
\Sigma(I)\ar[r]_{\Sigma(i)}\ar@{-->}[d]^{\eta_I}&
\Sigma(N)\ar[d]^{\eta_N}\\
M\ar[r]^{\pi}\ar@/_1pc/[rr]_{f}&I\ar[r]^{i}&N}
\end{equation}
The existence of the middle arrow and the commutativity of this diagram 
now follow from the fact that the rows are the canonical decompositions 
of $\Sigma(f)$ and of~$f$. Existence of~$\eta_I$ shows that $I$ is a submodule
of~$N$, with $\Sigma$-structure given by~$\eta_I$, and the commutativity of 
the diagram shows that both $\pi:M\to I$ and $i:I\to N$ are 
$\Sigma$-homomorphisms. In this way we see that any $\Sigma$-homomorphism 
can be uniquely decomposed into a surjective homomorphism followed by 
an embedding of a submodule.

\nxsubpoint (Free $\Sigma$-modules.)
Recall that $\mu_S:\Sigma^2(S)\to\Sigma(S)$ is a $\Sigma$-structure 
on~$\Sigma(S)$, thus defining a $\Sigma$-module 
$L_\Sigma(S)=(\Sigma(S),\mu_S)$, which will be usually denoted by~$\Sigma(S)$.
Moreover, we have seen in~\ptref{sp:sigma.free.obj} that $L_\Sigma$ is a 
left adjoint to the forgetful functor $\Gamma_\Sigma:\catMod\Sigma\to\catSets$,
i.e.\ there is a canonical bijection $\Hom_\Sigma(\Sigma(S),M)\cong
\Hom_\catSets(S,M)$ between $\Sigma$-homomorphisms $f:\Sigma(S)\to M$ and 
arbitrary maps $f^\flat:S\to M$. Of course, we have a canonical embedding 
$\epsilon_S:S\to\Sigma(S)$, and the above correspondence is given by 
$f^\flat=f\circ\epsilon_S$. We say that a $\Sigma$-module $M$ is {\em free},
if it is isomorphic to some $\Sigma(S)$ (cf.~\ptref{sp:monad.freeobj}); 
the image of~$S$ in~$M$ is called a {\em system of free generators of~$M$}. 
Finally, we say that $M$ is {\em free of (finite) rank~$n$} 
if it is isomorphic to~$\Sigma(n)$.

\nxsubpoint (Matrices.)
In~\ptref{prop:univ.lower.sigma} and \ptref{sp:Q.sigma} 
we have constructed a certain category $\cC_\Sigma=\catSets_\Sigma$ and 
a fully faithful functor $Q_\Sigma:\catSets_\Sigma\to\catMod\Sigma$, 
which transforms a set~$S$ into the corresponding free module 
$(\Sigma(S),\mu_S)$. Clearly, the essential image of~$Q_\Sigma$ consists 
of all free $\Sigma$-modules, and $Q_\Sigma$ induces an equivalence of 
categories $\catSets_\Sigma$ and the category of free $\Sigma$-modules. 
Recall that $\Hom_{\cC_\Sigma}(S,T)=\Hom_\cC(S,\Sigma(T))\cong
\Hom_\Sigma(\Sigma(S),\Sigma(T))$.

We also have the full subcategory $\catN_\Sigma\subset\catSets_\Sigma$, 
given by the standard finite sets and corresponding free modules; 
of course, the essential image under~$Q_\Sigma$ of this category consists 
of all free $\Sigma$-modules of finite rank. We see that 
$\Hom_{\catN_\Sigma}(\stn,\stm)=\Hom_\catSets(\stn,\Sigma(m))=\Sigma(m)^n\cong
\Hom_\Sigma(\Sigma(n),\Sigma(m))$. This is the reason why we put 
$M(m,n;\Sigma):=\Sigma(m)^n$ and call this ``the set of $m\times n$-matrices 
with entries in~$\Sigma$'' (cf.~\ptref{sp:algmon.predescr},b)). 
The composition of homomorphisms 
$\Hom_\Sigma(\Sigma(n),\Sigma(m))\times\Hom_\Sigma(\Sigma(k),\Sigma(n))\to
\Hom_\Sigma(\Sigma(k),\Sigma(m))$ defines the ``matrix multiplication'' 
$M(m,n;\Sigma)\times M(n,k;\Sigma)\to M(m,k;\Sigma)$, which is essentially 
given by the maps $\mu_m^{(n)}:\Sigma(n)\times\Sigma(m)^n\to\Sigma(m)$, 
up to a permutation of arguments (cf.~\ptref{sp:expl.algmon.cond},b) 
and~\ptref{prop:univ.lower.sigma}). We have the ``identity matrices'' 
$I_n=(\{1\}_\stn,\ldots,\{n\}_\stn)\in M(n,n;\Sigma)$, which correspond to 
$\id_{\Sigma(n)}$ (cf.~\ptref{sp:expl.algmon.cond},c)). 
Moreover, if we put $M(1,n;X):=X^n$ for any set~$X$ 
(``the set of rows over~$X$''), then a $\Sigma$-action $\alpha$ on~$X$ 
corresponds to a family of maps~$M(n,k;\Sigma)\times M(1,n;X)\to M(1,k;X)$, 
essentially given by the $\alpha^{(n)}:\Sigma(n)\times X^n\to X$ 
(cf.~\ptref{sp:expl.action},b)), which can be understood as some sort of 
multiplication of a row by a matrix.

\nxsubpoint (Invertible matrices.)
Of course, we have the groups of {\em invertible matrices} 
$GL_n(\Sigma)=GL(n,\Sigma)\subset M(n,n;\Sigma)$, which correspond to 
$\Aut_\Sigma(\Sigma(n))\subset\End_\Sigma(\Sigma(n))$. Clearly, a 
square matrix $f=(f_1,f_2,\ldots,f_n)$, $f_i\in\Sigma(n)$ 
is invertible iff there is another square matrix $g=(g_1,\ldots,g_n)$, 
such that $f_i(g_1,\ldots,g_n)=\{i\}_\stn$ and
$g_i(f_1,\ldots,f_n)=\{i\}_\stn$ for all~$1\leq i\leq n$.

\nxsubpoint (Initial object and ideals.)
Notice that $\st0:=\Sigma(0)$ is the initial object of $\catMod\Sigma$. If 
$\Sigma$ is a monad without constants, $\Sigma(0)$ is actually an empty set, 
so we prefer to denote it by $\emptyset$ and call it the 
{\em empty $\Sigma$-module}; when $\Sigma$ is a monad with zero, then 
$\st0=\Sigma(0)$ is a one-element set, hence a zero (both initial and final) 
object of $\catMod\Sigma$. Any module~$M$ contains a smallest submodule, 
equal to the image of $\Sigma(0)\to M$; it is called the 
{\em initial submodule} of~$M$, or the {\em empty} or {\em zero submodule} 
of~$M$, and it is denoted by $\emptyset_M$, $\emptyset$ or $0$ 
(depending on the situation). In any case it is a quotient of~$\Sigma(0)$; 
if $\Sigma$ has at most one constant, then it is always isomorphic 
to~$\Sigma(0)$. 

On the other hand, $\Sigma(0)\to\Sigma(1)$ is always injective 
(cf.\ \ptref{sp:sigma.injtoinj}), so the set of all {\em ideals} 
(i.e.\ $\Sigma$-submodules of $|\Sigma|=\Sigma(1)$) has a smallest element 
--- the {\em initial, empty} or {\em zero ideal}, which is always isomorphic 
to the initial object $\Sigma(0)$.

\nxsubpoint\label{sp:subm.gen.subset} (Submodule generated by a subset.)
If $M=(M,\eta_M)$ is a $\Sigma$-module, any subset $S\subset M$ induces 
a $\Sigma$-homomorphism $f:\Sigma(S)\to M$; according to~\ptref{sp:modhom.img},
the image $\langle S\rangle:=f(\Sigma(S))\subset M$ is a submodule of~$M$, 
clearly containing~$S$. We claim that {\em $\langle S\rangle$ is the smallest 
submodule of~$M$ containing~$S$, i.e.\ the submodule of~$M$, generated by~$S$.}
Indeed, if a submodule $N\subset M$ contains $S$, then 
$\langle S\rangle=f(\Sigma(S))=\eta_M(\Sigma(S))\subset\eta_M(\Sigma(N))=N$ 
(here we have identified $\Sigma(S)$ and $\Sigma(N)$ with their images in 
$\Sigma(M)$). Since $\Sigma(S)$ consists of all expressions 
$t(\{z_1\},\ldots,\{z_n\})$ with $t\in\Sigma(n)$ and $z_i\in S$ 
(cf.~\ptref{sp:expl.SX}), we see that $\langle S\rangle$ consists of all 
elements of the form $t_M(z_1,\ldots,z_n)\in M$ with $t\in\Sigma(n)$ and 
$z_i\in S$. In particular, any element of $\langle S\rangle$ can be expressed 
in terms of a finite subset $S_0\subset S$, i.e.\ belongs to some 
$\langle S_0\rangle$ with finite $S_0\subset S$.

\nxsubpoint\label{sp:fingen.mod} (Finitely generated modules.)
Of course, we say that {\em $S\subset M$ generates~$M$} if 
$\langle S\rangle=M$, i.e.\ if $f:\Sigma(S)\to N$ is surjective. 
In this case we also say that {\em $S$ is a system of generators of~$M$}. 
If $M$ admits a finite system of generators, i.e.\ if there exists a 
surjective homomorphism $\Sigma(n)\to M$, then we say that 
{\em $M$ is finitely generated\/} or {\em of finite type.} 
In this case any system $T$ of generators of~$M$ contains a finite subsystem 
$T_0\subset M$ of generators of~$M$, since the finite set $S$ is contained 
in $\langle T_0\rangle$ for a finite $T_0\subset T$.

\nxsubpoint (Compatible equivalence relations and strict quotients.)
We say that an equivalence relation~$R\subset M\times M$ 
on the underlying set of a $\Sigma$-module~$M=(M,\eta_M)$ 
is {\em compatible\/} with its 
$\Sigma$-structure, or that the quotient set $M/R$ is {\em compatible\/} 
with the $\Sigma$-structure of~$M$, if $M/R$ admits a (necessarily unique) 
$\Sigma$-structure, for which $\pi:M\to M/R$ becomes a homomorphism.
Clearly, in this case $R=M\times_{M/R}M$ is the kernel of $\pi$, and 
$\pi:M\to M/R$ is the cokernel of $R\rightrightarrows M$ both in 
$\catSets$ and in $\catMod\Sigma$, hence $M/R$ is a strict quotient of~$M$. 
Conversely, the kernel $R=M\times_N M$ of any $\Sigma$-homomorphism 
$f:M\to N$ is always compatible with the~$\Sigma$-structure on~$M$, since 
the image $f(M)\cong M/R$ admits a $\Sigma$-structure 
(cf.~\ptref{sp:modhom.img}), so if $f$ is a strict epimorphism, we must 
have $N\cong M/R$ in $\catMod\Sigma$, hence in $\catSets$ as well because 
of the compatibility of~$R$. We see that strict epimorphisms of 
$\catMod\Sigma$ are precisely the surjective $\Sigma$-homomorphisms, 
and that the strict quotients of~$M$ in $\catMod\Sigma$ are exactly the 
quotients $M/R$ of~$M$ modulo compatible equivalence relations~$R$; in 
particular, strict quotients of~$M$ constitute a set. 
Observe that in~\ptref{sp:strquot.algmon} we have obtained similar results 
for the category of algebraic monads; one might expect this 
because of~\ptref{sp:catalgmon.mon}.

Notice that $R$ is compatible with the $\Sigma$-structure on~$M$ iff all maps
$M^n\stackrel{{[t]}_M}\to M\stackrel\pi\to M/R$ factorize through 
$M^n\to(M/R)^n$. This means the following:
\begin{multline}\label{eq:comprel.mod}
{[t]}_M(x_1,\ldots,x_n)\equiv_R{[t]}_M(y_1,\ldots,y_n)
\text{ whenever }x_i\equiv_Ry_i\text{ for all }1\leq i\leq n,\\
\text{for any $t\in\Sigma(n)$, $x_1$, \dots, $x_n$, $y_1$, \dots, $y_n\in M$}
\end{multline}
Of course, it is enough to require this for $t$ from a system of generators 
of~$\Sigma$.

\nxsubpoint (Equivalence relation generated by a set of equations.)
Notice that the intersection $R=\bigcap R_\alpha$ of any family of 
compatible equivalence relations on a $\Sigma$-module~$M$ is again a 
compatible equivalence relation, since $R$ is the kernel of 
$M\to\prod_\alpha M/R_\alpha$. This means that for any set of equations 
(or relations) $E$ on~$M$, i.e.\ for any subset $E\subset M\times M$, 
we can find the smallest compatible equivalence relation $\langle E\rangle$ 
containing~$E$, simply by taking the intersection of all such relations. 
We have an alternative description of $R=\langle E\rangle$: namely, 
$x\equiv_R y$ iff this relation can be obtained after a finite number of  
applications of rules 0) $x\equiv y$ for any $(x,y)\in E$; 1) $x\equiv x$; 
2) $x\equiv y$ implies $y\equiv x$; 3) $x\equiv y$ and $y\equiv z$ 
imply $x\equiv z$; 4) the substitution rule~\eqref{eq:comprel.mod}.

Since the kernel~$R'$ of any $\Sigma$-homomorphism $f:M\to N$ is a 
compatible equivalence relation, it contains~$E$ iff it contains 
$\langle E\rangle$, i.e.\ iff $f$ factorizes through 
$\pi:M\to M/\langle E\rangle$. In other words, $\pi$ is universal among 
all $\Sigma$-homomorphisms $f:M\to N$, such that $f(x)=f(y)$ for any 
$(x,y)\in E$.

\nxsubpoint\label{sp:coker.mod} (Cokernels of pairs of morphisms.)
Suppose $M$ is a $\Sigma$-module, and $p$, $q:S\rightrightarrows M$ 
are two maps of sets. 
Consider the category of all $\Sigma$-homomorphisms $f:M\to N$, such that 
$f\circ p=f\circ q$. We claim that it has an initial object: indeed, 
$f\circ p=f\circ q$ iff the kernel of $f$ contains the set $E=(p,q)(S)$ 
of all pairs $(p(s),q(s))$, hence $M\to M/\langle E\rangle$ has 
the required universal property. This applies in particular when 
$S$ is a $\Sigma$-module, and $p$ and $q$ are $\Sigma$-homomorphisms; 
we see that {\em cokernels of pairs of morphisms exist in~$\catMod\Sigma$.}

\nxsubpoint\label{sp:findirsum.mod} (Finite direct sums.)
Let's show that {\em finite direct sums}, i.e.\ {\em coproducts}
exist in $\catMod\Sigma$. Since this category has an initial object 
$\Sigma(0)$, we have to show the existence of the direct sum $M_1\oplus M_2$ 
of two $\Sigma$-modules. First of all, if both $M_1$ and $M_2$ are free, 
then their direct sum exists since $\Sigma(S_1\sqcup S_2)$ satisfies the 
universal property required from $\Sigma(S_1)\oplus\Sigma(S_2)$. 
Next, we know that any $M_i$ can be written as a cokernel of a pair 
of morphisms $p_i$, $q_i:\Sigma(E_i)\rightrightarrows\Sigma(S_i)$ 
between two free modules (cf.~\ptref{l:mod.coker.free}). 
Since inductive limits commute with other inductive limits, we see that 
$M_1\oplus M_2\cong\Coker(p_1,q_1)\oplus\Coker(p_2,q_2)\cong
\Coker\bigl(p_1\oplus p_2, q_1\oplus q_2:
\Sigma(E_1\sqcup E_2)\rightrightarrows\Sigma(S_1\sqcup S_2)\bigr)$. 
Now the last cokernel is representable (cf.~\ptref{sp:coker.mod}), hence 
$M_1\oplus M_2$ is representable as well.

\nxsubpoint\label{sp:elem.dirsum} (Elements of direct sums.)
Notice that if $S$ generates $M$, and $T$ generates $N$, then 
$S\sqcup T$ generates $M\oplus N$; in particular, if both 
$M$ and $N$ are finitely generated, $M\oplus N$ is also finitely generated. 
This also means that $M\sqcup N$ generates $M\oplus N$; 
taking~\ptref{sp:subm.gen.subset} into account, we see that any element 
of $M\oplus N$ can be written in form $t(x_1,\ldots, x_n,y_1,\ldots, y_m)$ 
for some $m$, $n\geq0$, $t\in\Sigma(n)$, $x_i\in M$ and 
$y_j\in N$; we can even assume $x_i\in S$ and $y_j\in T$.

Contrary to the case of modules over an associative ring, we cannot 
express in general any element of $M\oplus N$ in form 
$t'(x,y)$ for some $t'\in\Sigma(2)$, $x\in M$ and $y\in N$. However, 
if such a statement is true for all direct sum decompositions of form 
$\Sigma(n+m)=\Sigma(n)\oplus\Sigma(m)$, i.e.\ if any $t\in\Sigma(n+m)$ 
can be represented in form $t=t'\bigr(t_1(\{1\},\ldots,\{n\}),
t_2(\{n+1\},\ldots,\{n+m\})\bigr)$ for some $t'\in\Sigma(2)$, 
$t_1\in\Sigma(n)$ and $t_2\in\Sigma(m)$, then our description of 
arbitrary elements of $M\oplus N$ together with the associativity relations 
shows that any element of $M\oplus N$ can be written in form 
$t'(x,y)$ with $x\in M$, $y\in N$ in this case, i.e.\ the validity 
of the statement for $\Sigma(n+m)=\Sigma(n)\oplus\Sigma(m)$, for all 
$n$, $m\geq0$, implies its validity in general. 

This remark is applicable in particular to $\Zinfty$, $\Zninfty$ 
and $\barZinfty$, since any formal octahedral combination
$t=\lambda_1x_1+\cdots+\lambda_nx_n+\mu_1y_1+\cdots+\mu_my_m$, 
$\sum_i|\lambda_i|+\sum_j|\mu_j|\leq 1$, can be re-written as
$\lambda x+\mu y$ with $x$ an octahedral combination of $x_i$ 
and $y$ an octahedral combination of $y_j$. Indeed, put $\lambda:=\sum_i
|\lambda_i|$, $\mu:=\sum_j|\mu_j|$. If $\lambda=0$ or $\mu=0$, the
statement is trivial; otherwise we put $x:=\sum_i(\lambda^{-1}\lambda_i)x_i$
and $y:=\sum_j(\mu^{-1}\mu_j)y_j$.

\nxsubpoint\label{sp:filtindlim.mod} (Filtered inductive limits.)
{\em Filtered inductive limits $\injlim M_\alpha$ of 
$\Sigma$-modules exist, and they can be computed in $\catSets$.}
Indeed, since $\Sigma$ is algebraic, it commutes 
with filtered inductive limits, hence 
$\Sigma(\injlim M_\alpha)\cong\injlim\Sigma(M_\alpha)$. Now put 
$M:=\injlim M_\alpha$ (in $\catSets$), and 
$\eta_M:=\injlim\eta_{M_\alpha}:\Sigma(M)\to M$, where 
$\eta_{M_\alpha}:\Sigma(M_\alpha)\to M_\alpha$ is the $\Sigma$-structure 
of~$M_\alpha$.

\nxsubpoint\label{sp:indlim.mod} (Arbitrary inductive limits.)
Notice that arbitrary sums $\bigoplus_{i\in I}M_i$ exist in $\catMod\Sigma$, 
since such a sum can be written as the filtered inductive limit of 
$\bigoplus_{i\in J}M_i$ for all finite $J\subset I$; one could also use 
directly the same reasoning as in~\ptref{sp:findirsum.mod}. Since 
cokernels of pairs of morphisms also exist in $\catMod\Sigma$ 
(cf.~\ptref{sp:coker.mod}), we see that {\em arbitrary inductive limits 
exist in $\catMod\Sigma$}.

\nxsubpoint (Lattices of submodules and of strict quotients.)
Existence of arbitrary sums and products, together with existence of 
images of homomorphisms (cf.\ \ptref{sp:modhom.img}), shows that 
any family of subobjects $N_\alpha$ of $M$ has both an infimum and a 
supremum in the ordered set of subobjects of~$M$: one simply computes 
$\inf N_\alpha$ as $\bigcap N_\alpha$, and $\sup N_\alpha$, denoted also 
by $\sum_\alpha N_\alpha$, is the image of the canonical homomorphism 
$\bigoplus N_\alpha\to M$. Similar statements are true for strict quotients 
of~$M$ as well. 

\nxsubpoint\label{sp:ex.scalext} (Scalar extension.)
Combining~\ptref{sp:indlim.mod} and~\ptref{prop:ex.base.change}, we see that 
for any homomorphism of algebraic monads~$\rho:\Sigma\to\Xi$ 
the scalar restriction functor $\rho^*:\catMod\Xi\to\catMod\Sigma$ admits 
a left adjoint --- the {\em scalar extension} or {\em base change} functor 
$\rho_*:\catMod\Sigma\to\catMod\Xi$. We denote $\rho_*M$ also by 
$\Xi\otimes_\Sigma M$; however, this notation still has to be justified. 
Of course, in situation $\Sigma\stackrel\rho\to\Xi\stackrel\sigma\to\Lambda$ 
we have $(\sigma\rho)^*=\rho^*\sigma^*$ and 
$(\sigma\rho)_*\cong\sigma_*\rho_*$, so $\Lambda\otimes_\Xi
(\Xi\otimes_\Sigma M)$ can be identified with $\Lambda\otimes_\Sigma M$, i.e.\ 
we have some sort of ``associativity''. 
Also note that $\Hom_\Xi(\Xi\otimes_\Sigma M,N)\cong\Hom_\Sigma(M,\rho^*N)$,
and $\Xi\otimes_\Sigma\Sigma(S)\cong\Xi(S)$. This is actually a special case 
of the ``associativity'', since $\Sigma\otimes_{\Fempty}S=\Sigma(S)$ for any 
set~$S$. We see that $\otimes_{\Fempty}$ corresponds to the 
left $\otimes$-action of $\cA_{alg}\subset\cA$ on $\cC$, denoted 
by $\otimes$ or $\obslash$ in~\ptref{sp:tens.act.cat}.
Of course, 
$\rho_*$ commutes with arbitrary inductive limits, and in particular it is 
right exact; for example, $\Xi\otimes_\Sigma(M\oplus M')\cong 
(\Xi\otimes_\Sigma M)\oplus(\Xi\otimes_\Sigma M')$. Similarly, 
the scalar restriction functor $\rho^*$ commutes with arbitrary 
projective limits and in particular it is left exact, but in general 
not (right) exact, contrary to what one might expect. In fact, we'll 
show in the next chapter that {\em when $\Sigma$ is commutative, 
$\rho^*$ is (right) exact iff $\rho:\Sigma\to\Xi$ is unary}, and in this 
case $\rho^*$ commutes with arbitrary inductive limits, and even admits 
a right adjoint~$\rho^!$.

\nxsubpoint (Some counterexamples.)
Up to now we have seen that $\catMod\Sigma$ behaves in most respects like 
the category of (left) modules over an associative ring. However, some 
properties do not extend to this case. Consider for example the 
algebraic monad~$\Sigma$, such that $\catMod\Sigma$ is the category of 
commutative rings. Then $\bbZ\to\bbQ$ is both a monomorphism and an 
epimorphism in $\catMod\Sigma$, but not a strict epimorphism and not an 
isomorphism. Also note that the direct sums lose most of their properties: 
for example, direct sum of two monomorphisms need not be a monomorphism, 
since both $\bbZ\to\bbQ$ and $\id_{\bbZ/2\bbZ}$ are injective, but their 
coproduct in the category of commutative rings, i.e.\ the tensor product 
over~$\bbZ$, is not injective. Direct sums don't commute with finite 
direct products as well, since this is not true even in 
$\catMod\Fempty=\catSets$. On the other hand, filtered inductive limits 
are still left exact, since both filtered inductive limits and 
finite projective limits of $\Sigma$-modules are computed in~$\catSets$.

\nxsubpoint (Finite presentation.)
We say that a compatible equivalence relation~$R$ on a $\Sigma$-module~$M$ 
is {\em finitely generated}, if it is generated by a finite set of equations 
$E\subset M\times M$. Clearly, $R$ is finitely generated iff 
$\pi:M\to M/R$ is the universal coequalizer of a pair of maps  
$p',q':\stn\rightrightarrows M$ (cf.~\ptref{sp:coker.mod}), or equivalently, 
if $\pi:M\to M/R$ is the cokernel of a pair of morphisms 
$p,q:\Sigma(n)\rightrightarrows M$.

Given a module~$M$, we say that a pair of sets $(S,E)$, 
where $S\subset M$ and $E\subset\Sigma(S)\times\Sigma(S)$ is a 
{\em presentation\/} of~$M$ if $f:\Sigma(S)\to M$ is surjective and the kernel
of~$f$ is generated by~$E$; then $M\cong\Sigma(S)/\langle E\rangle$. 
Equivalently, we can say that a presentation of~$M$ is a pair of maps of sets 
$p',q':E\rightrightarrows\Sigma(S)$, or a pair of homomorphisms 
$p,q:\Sigma(E)\rightrightarrows\Sigma(S)$, such that~$M$ is the cokernel 
in $\catMod\Sigma$ of this pair of maps.

We say that $M$ is {\em of finite presentation}, or {\em finitely presented}, 
if it admits a finite presentation $(S,E)$, i.e.\ if it is a coequalizer 
of a pair of maps of sets $p',q':\stm\rightrightarrows\Sigma(n)$, 
or equivalently, if it is a cokernel of a pair of homomorphisms 
$p,q:\Sigma(m)\rightrightarrows\Sigma(n)$ between free modules of finite rank. 
Of course, such a pair of morphisms can be given by two matrices 
$P$, $Q\in\Sigma(n)^m=M(n,m;\Sigma)$.

Notice that the isomorphism classes of all finitely generated 
(resp.\ finitely presented) modules over a fixed algebraic monad~$\Sigma$ 
constitute a set, since any such module is isomorphic to a strict quotient 
of some~$\Sigma(n)$. In other words, the category of finitely generated 
(resp.\ finitely presented) modules is equivalent to a small category.

These notions of finite presentation and of finite type have all the usual 
properties, which can be shown essentially in the same way as over 
associative rings. For example, the cokernel of a pair of homomorphisms from 
a finitely generated module into a finitely presented module is itself 
finitely presented; conversely, the kernel of a surjective homomorphism from 
a finitely generated module to a finitely presented module is finitely 
generated as a compatible equivalence relation. Another simple fact: 
any finite inductive limit (e.g.\ finite sum) of finitely generated 
(resp.\ finitely presented) modules is itself finitely generated 
(resp.\ finitely presented).

\nxsubpoint\label{sp:mod.filtlim.fpres}
Any $\Sigma$-module is a filtered union of its finitely generated submodules, 
hence it is a filtered inductive limit of a system of finitely generated 
modules with injective transition morphisms. Moreover, {\em any 
$\Sigma$-module~$M$ is a filtered inductive limit of finitely presented
$\Sigma$-modules.} To show this we consider the category $\cD$ of finitely 
presented $\Sigma$-modules (or any its small subcategory, equivalent to 
the whole of~$\cD$), and observe that the category $\cD/M$ of homomorphisms 
$N\to M$ from finitely presented modules into~$M$ is filtering, being 
stable under finite inductive limits, and $\injlim_{\cD/M} N\cong M$.

\nxsubpoint (Direct factors and projective modules.)
Recall that $N$ is a {\em direct factor} of an object $M$ of some 
category (e.g.\ $\catMod\Sigma$) if 
we have two morphisms $i:N\to M$ and $\sigma:M\to N$, such that 
$\sigma i=\id_N$. Then $i$ is a strict monomorphism, and 
$\sigma$ is a strict epimorphism, so~$N$ can be considered both as a submodule 
and as a strict quotient of~$M$. Moreover, in this situation 
$p:=i\sigma\in\End_\Sigma(M)$ is an idempotent, i.e.\ $p^2=p$, 
$i:N\to M$ is the kernel of $p,\id_M:M\rightrightarrows M$, and 
$\sigma:M\to N$ is the cokernel of the same pair. Conversely, if 
$p^2=p$ is an idempotent in $\End(M)$, then we might reconstruct 
$N$ either as the kernel or as the cokernel of $p$ and $\id_M$ 
(if any of these exists in the category under consideration; of 
course, they both exist in~$\catMod\Sigma$), thus 
obtaining a direct factor of~$M$. In this way we obtain a bijection between 
direct factors of~$M$ and idempotents in~$\End(M)$. Notice that the property 
of being a direct factor is preserved by all functors, and since direct 
factors can be described both as kernels and as cokernels, we see that 
a direct factor of a left (resp.\ right) exact functor is again left 
(resp.\ right) exact. Another consequence: since 
$N=\Coker(p,\id_M:M\rightrightarrows M)$, 
if $M$ is finitely generated (resp.\ finitely presented), the same is true 
for any its direct factor~$N$.

Let's apply this to {\em projective} modules~$P$, i.e.\ such $\Sigma$-modules 
$P$, for which $\Hom_\Sigma(P,-)$ preserves strict epimorphisms,
i.e.\ transforms surjective homomorphisms into surjective maps. 
Clearly, any free $\Sigma$-module $\Sigma(S)$ is projective, since 
$\Hom_\Sigma(\Sigma(S),-)=\Hom_\catSets(S,-)$ (if~$S$ is infinite, we need 
the axiom of choice here), and any direct factor of a projective module is 
projective, so we see that {\em direct factors of free modules are 
projective.} Conversely, if $P$ is projective, and $S$ is any system of 
generators of~$P$, then the surjection $\sigma:\Sigma(S)\to P$ admits a 
section, hence $P$ is a direct factor of a free module~$\Sigma(S)$.

\nxsubpoint 
Given any $\Sigma$-module $M$ and any {\em filtered\/} inductive limit 
$N=\injlim N_\alpha$, we have a canonical map 
$i=i_M:\injlim\Hom_\Sigma(M,N_\alpha)\to\Hom_\Sigma(M,\injlim N_\alpha)$. 
We claim that {\em if $M$ is finitely presented (resp.\ finitely generated), 
then $i$ is always bijective (resp.\ injective).} Indeed, in this 
case we can find a finite presentation of~$M$, i.e.\ a right exact sequence 
$\Sigma(m)\rightrightarrows\Sigma(n)\to M$ (resp.\ a surjective homomorphism 
$\Sigma(n)\to M$). Consider the following diagram (resp.\ its left square):
\begin{equation}
\xymatrix{
\injlim\Hom_\Sigma(M,N_\alpha)\ar[r]\ar[d]^{i}&
\injlim\Hom_\Sigma(\Sigma(n),N_\alpha)\ar@/_1pt/[r]\ar@/^3pt/[r]
\ar[d]^{i_{\Sigma(n)}}&
\injlim\Hom_\Sigma(\Sigma(m),N_\alpha)\ar[d]^{i_{\Sigma(m)}}\\
\Hom_\Sigma(M,\injlim N_\alpha)\ar[r]&
\Hom_\Sigma(\Sigma(n),\injlim N_\alpha)\ar@/_1pt/[r]\ar@/^3pt/[r]&
\Hom_\Sigma(\Sigma(m),\injlim N_\alpha)}
\end{equation}
Its rows are left exact (resp.\ the horizontal arrows of the left square are 
injective), and the vertical arrows $i_{\Sigma(n)}$ and $i_{\Sigma(m)}$ 
are isomorphisms, since $\Hom_\Sigma(\Sigma(n),N)\cong N^n$, and 
filtered inductive limits commute with finite products in $\catSets$. 
This implies the bijectivity (resp.\ injectivity) of $i=i_M$.

Conversely, {\em if for some $\Sigma$-module~$M$ the functor 
$\Hom_\Sigma(M,-)$ commutes with filtered inductive limits, then 
$M$ is finitely presented.} Indeed, write $M$ as a filtered inductive 
limit $\injlim N_\alpha$ of some finitely presented modules 
(cf.~\ptref{sp:mod.filtlim.fpres}), and denote by $j_\alpha:N_\alpha\to M$ 
the canonical homomorphisms $N_\alpha\to\injlim N_\alpha\cong M$. We see that 
$\id_M\in\Hom_\Sigma(M,M)\cong\Hom_\Sigma(M,\injlim N_\alpha)$ has to come 
from some $\sigma\in\Hom_\Sigma(M,N_\alpha)$, i.e.\ 
$j_\alpha\circ\sigma=\id_M$. This means that $M$ is a direct factor of a 
finitely presented module~$N_\alpha$, hence $M$ is finitely presented itself.

\nxpointtoc{Categories of algebraic modules}\label{p:cat.algmod}
Now we would like to prove some basic properties of {\em algebraic 
modules over an algebraic monad~$\Sigma$}, i.e.\ the (left or right) 
$\Sigma$-modules in $\cA_{alg}\subset\cA$. 

\nxsubpoint (Left $\Sigma$-modules.)
Let's begin with the case of a left algebraic $\Sigma$-module~$F$, i.e.\ 
an algebraic endofunctor $F$, together with a $\Sigma$-action 
$\alpha:\Sigma F\to F$, such that $\alpha\circ(\epsilon\star F)=\id_F$ and 
$\alpha\circ(\Sigma\star\alpha)=\alpha\circ(\mu\star F)$ 
(cf.~\ptref{sp:cat.left.sigmamod}). We know that $F$ is given by 
a collection of sets $\{F(n)\}_{n\geq0}$ and of maps 
$F(\phi):F(m)\to F(n)$ for each $\phi:\stm\to\stn$ 
(cf.~\ptref{sp:simpl.descr.algendf}), and $\alpha:\Sigma F\to F$ is 
given by maps $\alpha^{(k)}_n:\Sigma(k)\times F(n)^k\to F(n)$, subject 
to certain compatibility conditions (cf.~\ptref{sp:algnat.SFtoG}). Recall that 
we put ${[t]}_{F(n)}(x_1,\ldots,x_k)=\alpha^{(k)}_n(t;x_1,\ldots,x_n)$ 
for any $t\in\Sigma(k)$ and any $x_i\in F(n)$.

Now we can write down the conditions for such a collection of data to 
define a left $\Sigma$-module~$F$. Actually, equality of natural 
transformations of algebraic functors can be checked on standard finite sets, 
so we end up with requiring the $\{\alpha^{(k)}_n:\Sigma(k)\times F(n)^k\to
F(n)\}_{k\geq0}$ to define a $\Sigma$-module structure on each $F(n)$, 
and requiring $F(\phi):F(m)\to F(n)$ to be a $\Sigma$-homomorphism for 
each $\phi:\stm\to\stn$.

We see that a left algebraic $\Sigma$-module~$F$ is essentially the same thing 
as a functor $F:\catN\to\catMod\Sigma=\cC^\Sigma$, $\stn\mapsto F(n)$. 
We could deduce this result directly from~\ptref{prop:j.descr.algfun} 
and~\ptref{prop:left.sigma.funct}: indeed, we know that the restriction 
to $\cA_{alg}\subset\cA=\catFunct(\cC,\cC)$ of the restriction functor 
$J^*:\catFunct(\cC,\cC)\to\catFunct(\catN,\cC)$ is an equivalence of 
categories (cf.~\ptref{prop:j.descr.algfun}), clearly compatible with the 
left $\obslash$-action of $\cA_{alg}$ on $\cA_{alg}\subset\cA$ and on 
$\catFunct(\catN,\cC)$, hence $\cA_{alg}^\Sigma\cong\catFunct(\catN,\cC)^\Sigma
\cong\catFunct(\catN,\cC^\Sigma)=\catMod\Sigma^\catN$ 
by~\ptref{prop:left.sigma.funct}.

\nxsubpoint (Matrix interpretation.)
Let $F$ be an algebraic left $\Sigma$-module. Then we put 
$M(n,m;F):=F(n)^m$, and interpret this as ``the set of $n\times m$-matrices 
with entries in~$F$'' as before. The action maps $\alpha^{(k)}_n:\Sigma(k)
\times F(n)^k\to F(n)$ can be interpreted as some matrix multiplication maps 
$M(k,m;\Sigma)\times M(n,k;F)\to M(n,m;F)$ (with opposite order of arguments), 
and then the associativity and unit relations for $\alpha$ can be understood 
as associativity of matrix multiplication 
$M(m,s;\Sigma)\times M(k,m;\Sigma)\times M(n,k;F)\to M(n,s;F)$, and the 
requirement for $I_k=(\{1\},\ldots,\{k\})\in M(k,k;\Sigma)$ to act 
identically on $M(n,k;F)=F(n)^k$.

\nxsubpoint\label{sp:def.absend.ring.endof} 
(Monad of endomorphisms of an algebraic endofunctor.)
If $F$ is an algebraic left $\Sigma$-module, any $t\in\Sigma(k)$ induces 
a family of compatible maps ${[t]}_{F(n)}:F(n)^k\to F(n)$, i.e.\ 
a natural transformation ${[t]}_F:F^k\to F$. These natural transformations 
have to satisfy a substitution property similar to that discussed 
in~\ptref{sp:algmon.mod.alg}, and conversely, once we have a family of 
such natural transformations, satisfying this property as well as the 
requirement for ${[\{s\}_\stn]}_F:F^k\to F$ to coincide with the projection 
onto the $k$-th component, we obtain a left $\Sigma$-structure on~$F$.

This enables us to repeat the reasoning of \ptref{sp:def.absend.ring}, 
and define for any algebraic endofunctor~$F$ its 
``absolute endomorphism ring'' $\END(F)$ by putting 
$(\END(F))(n):=\Hom_\cA(F^n,F)$, with the multiplication of~$\END(F)$ 
given by composition of natural transformations. Then $\END(F)$ 
canonically acts on $F$ from the left, i.e.\ we have a natural 
transformation $\END(F)\otimes F\to F$, and giving a left $\Sigma$-module 
structure on~$F$ is equivalent to giving a homomorphism of algebraic monads 
$\rho:\Sigma\to\END(F)$; then the left $\Sigma$-structure on~$F$ is recovered 
from the canonical $\END(F)$-structure by means of scalar restriction 
along~$\rho$. 

Similarly to~\ptref{sp:def.absend.ring}, we denote by $\Sigma_F$ the 
image of~$\Sigma$ in~$\END(F)$, and say that $F$ is a {\em faithful\/} 
algebraic left $\Sigma$-module if $\rho:\Sigma\to\END(F)$ is injective. 
For example, $\Sigma$ is always a faithful module over itself.

\nxsubpoint\label{sp:limits.leftmod}(Limits, subobjects, strict quotients\dots)
Once the equivalence of categories $\cA_{alg}^\Sigma\cong
\catFunct(\catN,\catMod\Sigma)$ is established, we can deduce properties 
of this category from those of~$\catMod\Sigma$. For example, 
arbitrary projective and inductive limits exist in $\cA_{alg}^\Sigma$, 
and they are computed componentwise, e.g.\ $(\injlim F_\alpha)(n)=
\injlim F_\alpha(n)$, the forgetful functor $\cA_{alg}^\Sigma\to\cA_{alg}$ 
is left exact and commutes with filtered inductive limits, and so on. 
In particular, subobjects~$F'$ of~$F$ are given by families of 
submodules $F'(n)\subset F(n)$, compatible with all maps $F(\phi)$, 
strict quotients of~$F$ are parametrized by compatible equivalence 
relations~$R\subset F\times F$, which can be interpeted as families 
of compatible equivalence relations $R(n)\subset F(n)\times F(n)$, 
and the strict quotient~$F/R$ is computed componentwise: $(F/R)(n)=F(n)/R(n)$.

\nxsubpoint\label{sp:scal.change.lmod} (Scalar restriction and base change.)
Recall that any homomorphism of algebraic monads $\rho:\Sigma\to\Xi$ 
induces scalar restriction functors $\rho^*:\catMod\Xi=\cC^\Xi\to
\catMod\Sigma$ and $\rho^*_\cA:\cA^\Xi\to\cA^\Sigma$; 
identifying $\cA^\Sigma$ with $\catFunct(\cC,\cC^\Sigma)$ 
(cf.~\ptref{prop:left.sigma.funct}), 
we see that this functor $\cA^\Xi\to\cA^\Sigma$ is essentially given by 
$F\mapsto\rho^*\circ F$, for any $F:\cC\to\cC^\Xi$. When we restrict 
our attention to full subcategories of algebraic left modules $\cA_{alg}^\Sigma
\subset\cA^\Sigma$, $\cA_{alg}^\Sigma\cong\catFunct(\catN,\catMod\Sigma)$, 
and similarly for $\cA_{alg}^\Xi$, we see that the scalar restriction 
functor $\cA_{alg}^\Xi\cong\catFunct(\catN,\catMod\Xi)\to\cA_{alg}^\Sigma$ 
is given again by $F\mapsto\rho^*\circ F$, i.e.\ the scalar restriction 
of left algebraic modules is computed componentwise. 
Since $\rho^*:\catMod\Xi\to\catMod\Sigma$ admits a left adjoint 
$\rho_*$ (cf.~\ptref{sp:ex.scalext}), we see that the scalar restriction 
functor between categories of algebraic left modules admits a left 
adjoint as well, given by $F\mapsto\rho_*\circ F$. Of course, 
we say that $\rho_*F$ is the {\em scalar extension} or {\em base change} 
of~$F$ with respect to $\rho:\Sigma\to\Xi$, and denote $\rho_*F$ also 
by $\Xi\otimes_\Sigma F$; we have just seen that 
$(\rho_*F)(n)=\rho_*(F(n))$, or $(\Xi\otimes_\Sigma F)(n)=
\Xi\otimes_\Sigma F(n)$.

\nxsubpoint (Kan extensions, skeleta and coskeleta.)
Since $\cA_{alg}^\Sigma\cong\catFunct(\catN,\cC^\Sigma)=(\catMod\Sigma)^\catN$,
and arbitrary inductive and projective limits exist in $\catMod\Sigma$, 
we see that any functor $K:\cI\to\catN$ admits both left and right 
Kan extensions $K_!$ and $K_*:(\catMod\Sigma)^\cI\to\cA_{alg}^\Sigma$, 
i.e.\ both left and right adjoints to $K^*:F\mapsto F\circ K$. 
Usually we apply this for the embeddings $\cI\to\catN$ of some subcategories 
$\cI$ of $\catN$. For example, if we take for $\cI$ the subcategory of 
$\catN$, which contains all objects of~$\catN$, but only the identity 
morphisms, we see that the functor $K^*$ which maps any 
left algebraic $\Sigma$-module~$F$ into the collection of $\Sigma$-modules 
$\{F(n)\}_{n\geq0}$ (without any transition morphisms), admits 
both a left and a right adjoint, and we can compute these adjoints by 
usual formulas for Kan extensions. For example, $K_!G$ is given by 
$(K_!G)(n)=\injlim_{\cI/\stn} G(i)$, where $\cI/\stn$ is the category of 
couples $(i,\phi)$, $i\in\Ob\cI$, $\phi:K(i)\to\stn$; in our case 
$K_!$ transforms a collection of $\Sigma$-modules $M=\{M_n\}_{n\geq0}$, 
into $K_!M$, given by $(K_!M)(n):=\bigoplus_{\phi:\stm\to\stn}M_m$, 
and similarly, $(K_*M)(n)=\prod_{\phi:\stn\to\stm}M_m$. 
Notice that left Kan extensions commute with base change, while 
right Kan extensions commute with scalar restriction.

Another important example: when $\cI=\{\st1\}\subset\catN$, we see that
$F\mapsto F(1)$, $\cA_{alg}^\Sigma\to\catMod\Sigma$, admits both a 
left and a right adjoint; they are given by $(K_!M)(n)=\bigoplus_\stn M
=M^{(n)}$ (the direct sum of $n$ copies of~$M$), and 
$(K_*M)(n)=M$ for any~$n\geq 0$.

Finally, let's take for $\cI$ the full subcategory $\catN_{\leq n}$, 
consisting of all objects $\{\stm\}_{m\leq n}$; in this manner we obtain 
the {\em skeleta} $\sk_n:=K_{\leq n,!}K_{\leq n}^*$ and 
the {\em coskeleta} functors $\cosk_n:=K_{\leq n,*}K_{\leq n}^*$, 
computed almost in the same way as for cosimplicial objects 
(remember, $\cA_{alg}^\Sigma=(\catMod\Sigma)^\catN$ is something like 
the category of cosimplicial $\Sigma$-modules, but not exactly, 
since we allow all maps between finite sets in~$\catN$, not just 
the non-decreasing ones, and consider the empty set $\st0\in\Ob\catN$ 
as well.)

\nxsubpoint (Free algebraic left $\Sigma$-modules.)
An important consequence is that the forgetful functor 
$F\mapsto\|F\|=\bigsqcup_{n\geq0}F(n)$ from the category of algebraic 
left $\Sigma$-modules into the category of graded sets admits a left adjoint. 
Indeed, we can decompose this functor into two functors: the first 
maps $F$ into the collection of $\Sigma$-modules $\{F(n)\}_{n\geq0}$, 
and it admits both a left and a right adjoint by the theory of Kan 
extensions, and the second maps such a collection $\{M_n\}$ 
into the collection of underlying sets of $M_n$; it also admits a left 
adjoint, given by $S=\bigsqcup_{n\geq0} S_n\mapsto\{\Sigma(S_n)\}_{n\geq0}$.

\nxsubpoint (Right $\Sigma$-modules.)
Now we would like to obtain similar descriptions for the category 
$\cA_{alg,\Sigma}$ of right algebraic $\Sigma$-modules, i.e.\ of 
algebraic endofunctors $G$, equipped with a right $\Sigma$-action 
$\beta:G\Sigma\to G$. Of course, $G$ itself is given by a collection of 
sets $\{G(n)\}_{n\geq0}$ and of maps $G(\phi):G(m)\to G(n)$, defined 
for each $\phi:\stm\to\stn$ (cf.~\ptref{sp:simpl.descr.algendf}), 
and $\beta$ is given by a collection of maps 
$\beta^{(k)}_n:G(k)\times\Sigma(n)^k\to G(n)$. We can put 
$M(n,m;G):=G(n)^m$ as usual, and interpet the above maps as 
some matrix multiplication rules $M(k,s;G)\times M(n,k;\Sigma)\to M(n,s;G)$. 
Then the requirements for $\beta$ to be a right $\Sigma$-action 
translate into some ``associativity of matrix multiplication'' 
$M(k,s;G)\times M(n,k;\Sigma)\times M(m,n;\Sigma)\to M(m,s;G)$ and 
into the requirement for the identity matrix $I_n\in\Sigma(n)^n$ to 
act identically on~$G(n)$; the verification is similar to that 
of~\ptref{sp:expl.algmon.cond},b) and c). 

Notice that $\beta^{(k)}_n:G(k)\times\Sigma(n)^k\to G(n)$ 
can be interpreted as a map from $\Sigma(n)^k\cong
\Hom_\Sigma(\Sigma(k),\Sigma(n))$ into $\Hom_\catSets(G(k),G(n))$, so we 
essentially obtain a functor $\tilde G$ from the category $\catN_\Sigma$ 
of standard free $\Sigma$-modules of finite rank into the category of sets. 
We'll check in a moment that this is indeed a functor, but this can be 
also seen directly in terms of the ``associativity of matrix multiplication'' 
just discussed.

\nxsubpoint\label{sp:altdescr.rmod} (Alternative description of right modules.)
Recall that in~\ptref{prop:univ.lower.sigma} we have shown that 
the category $\cA_\Sigma$ of right $\Sigma$-modules in $\cA=\catFunct(\cC,\cC)$
is equivalent to $\catFunct(\cC_\Sigma,\cC)$, where $\cC_\Sigma$ has the 
same objects as $\cC$, but morphisms given by $\Hom_{\cC_\Sigma}(X,Y)=
\Hom_\cC(X,\Sigma(Y))\cong\Hom_\Sigma(\Sigma(X),\Sigma(Y))$, so 
$\cC_\Sigma$ is essentially the category of free $\Sigma$-modules. 
We have also constructed a fully faithful functor $Q_\Sigma:\cC_\Sigma\to
\cC^\Sigma=\catMod\Sigma$, that transforms $X$ into the free module 
$L_\Sigma(X)=(\Sigma(X),\mu_X)$, and a functor 
$I_\Sigma:\cC\to\cC_\Sigma$ with a canonical right $\Sigma$-action; 
then the equivalence $\catFunct(\cC_\Sigma,\cC)\cong\catFunct(\cC,\cC)_\Sigma
=\cA_\Sigma$ is given by $I_\Sigma^*: G\mapsto G\circ I_\Sigma$, 
with the right $\Sigma$-structure on $G\circ I_\Sigma$ 
induced by that of~$I_\Sigma$.

In this way the full subcategory $\cA_{alg,\Sigma}\subset\cA_\Sigma$ 
of algebraic right $\Sigma$-modules corresponds to a certain 
full subcategory of $\catFunct(\cC_\Sigma,\cC)$. We claim that
{\em the full subcategory of $\catFunct(\cC_\Sigma,\cC)$, equivalent to 
$\cA_{alg,\Sigma}$, coincides with the essential image of the 
left Kan extension $J_{\Sigma,!}:\catFunct(\catN_\Sigma,\cC)\to
\catFunct(\cC_\Sigma,\cC)$, where $\catN_\Sigma$ is the full 
subcategory of $\cC_\Sigma$ with standard finite sets for objects 
(i.e. $\catN_\Sigma$ is essentially the category of ``standard free 
$\Sigma$-modules of finite rank''), and 
$J_\Sigma:\catN_\Sigma\to\cC_\Sigma$ is the natural embedding. 
The quasi-inverse equivalence is given by the restriction to $\cA_{alg,\Sigma}$
of the restriction functor $J^*_\Sigma$.}

Since arbitrary inductive limits exist in $\catFunct(\cC_\Sigma,\cC)$, 
the left Kan extension $J_{\Sigma,!}$, i.e.\ the left adjoint to 
$J^*_\Sigma$, exists for general reasons; it is fully faithful since 
$J_\Sigma$ is fully faithful itself, so we have $J^*_\Sigma J_{!,\Sigma}(G')
\cong G'$ for any functor $G':\catN_\Sigma\to\cC$. Now we have to check 
that some functor $G:\cC_\Sigma\to\cC$ is algebraic, i.e.\ the 
underlying functor $I_\Sigma^*(G)=G\circ I_\Sigma:\cC\to\cC$ 
of the corresponding right 
$\Sigma$-module commutes with arbitrary inductive limits, 
iff $G$ is isomorphic to some $J_{\Sigma,!}(G')$. This follows now 
from~\ptref{prop:j.descr.algfun} and lemma~\ptref{l:kanext.rightmod} below, 
once we take into account $J_\Sigma\circ\bar I_\Sigma=I_\Sigma\circ J$, and 
the conservativity of $I_\Sigma^*$. Indeed, if 
$G=J_{\Sigma,!}(G')$, we get $I_\Sigma^*(G)\cong J_!\bar I_\Sigma^*(G')$ 
from \ptref{l:kanext.rightmod}, so $I_\Sigma^*(G)$ is algebraic 
by~\ptref{prop:j.descr.algfun}. Conversely, if $I_\Sigma^*(G)$ is algebraic 
for some $G:\cC_\Sigma\to\cC$, we obtain 
$I_\Sigma^*(G)\cong J_!J^*I_\Sigma^*(G)=J_!\bar I_\Sigma^*J_\Sigma^*(G)
\cong I_\Sigma^* J_{\Sigma,!}J_\Sigma^*(G)$, hence 
$J_{\Sigma,!}J_\Sigma^*(G)\to G$ is an isomorphism because of the 
conservativity of~$I_\Sigma^*$.

\begin{LemmaD}\label{l:kanext.rightmod}
In the above notations $I_\Sigma^* J_{\Sigma,!}(G')\cong J_!\bar 
I_\Sigma^*(G')$, 
where $\bar I_\Sigma:\catN\to\catN_\Sigma$ is the restriction of 
$I_\Sigma:\cC\to\cC_\Sigma$ to $\catN$, and $G':\catN_\Sigma\to\cC$ 
is an arbitrary functor. Moreover, both inductive limits used 
in computation of left Kan extensions $J_{\Sigma,!}$ and $J_!$ 
are in fact filtering.  
\end{LemmaD}
\begin{Proof}
First of all, a canonical morphism $\kappa_{G'}:
J_!\bar I_\Sigma^*(G')\to I_\Sigma^* J_{\Sigma,!}(G')$ is constructed by 
adjointness from $\bar I_\Sigma^*(G')\simto J^*I_\Sigma^*J_{\Sigma,!}(G')
=\bar I_\Sigma^*J_\Sigma^*J_{\Sigma,!}(G')$, obtained in its turn by 
applying $\bar I_\Sigma$ to the adjointness morphism 
$G'\to J_\Sigma^*J_{\Sigma,!}(G')$, which is in fact an isomorphism, 
$J_{\Sigma,!}$ being fully faithful. We have to show that $\kappa_{G'}$ is 
also an isomorphism; for this we write down the expressions for left 
Kan extensions in terms of inductive limits:
$(J_!\bar I_\Sigma^*(G'))(S)=\injlim_{\catN/S} G'(I_\Sigma(\stn))$, and 
$(I_\Sigma^* J_{\Sigma,!}(G'))(S)=\injlim_{\catN_\Sigma/I_\Sigma(S)} 
G'$. We see that in the first case we compute the inductive limit of 
$G'(I_\Sigma(\stn))$ along the category $\catN/S$ of all maps $i:\stn\to S$, 
i.e.\ essentially the limit of $G'$ along the category $I_\Sigma(\catN/S)$ 
of all morphisms $I_\Sigma(i):I_\Sigma(\stn)\to I_\Sigma(S)$; if we identify 
$\cC_\Sigma$ with a full subcategory of $\catMod\Sigma$ by means of $Q_\Sigma$,
we see that the first index category consists of all morphisms 
$\Sigma(n)\to\Sigma(S)$ of form $\Sigma(i)$ for some $i:\stn\to S$. 
Now observe that the second inductive limit is computed for the same functor 
$G'$, but along the category $\catN_\Sigma/I_\Sigma(S)$, which essentially 
consists of {\em all\/} $\Sigma$-homomorphisms $\Sigma(n)\to\Sigma(S)$.

So we have to show that the natural embedding $I_\Sigma(\catN/S)\to 
\catN_\Sigma/I_\Sigma(S)$ induces an isomorphism of inductive limits along 
these index categories, i.e.\ that the subcategory of homomorphisms of form 
$\Sigma(i):\Sigma(n)\to\Sigma(S)$ is {\em cofinal\/} in the category of all
$\Sigma$-homomorphisms $\Sigma(n)\to\Sigma(S)$.

First of all, notice that
$\Sigma(S)=\injlim_{S_0\subset S}\Sigma(S_0)$, 
where the limit is taken along all finite $S_0\subset S$, 
since $\Sigma$ is algebraic, hence 
$\Hom_\Sigma(\Sigma(n),\Sigma(S))\cong\Sigma(S)^n\cong\injlim_{S_0\subset S}
\Sigma(S_0)^n\cong\injlim_{S_0\subset S}\Hom_\Sigma(\Sigma(n),\Sigma(S_0))$. 
This shows that any $\Sigma$-homomorphism $f:\Sigma(n)\to\Sigma(S)$ 
factorizes through $\Sigma(S_0)\subset\Sigma(S)$ for some finite 
$S_0\subset S$; replacing $S_0$ by an isomorphic standard finite set $\stm$, 
we see that any $f:\Sigma(n)\to\Sigma(S)$ factorizes through some 
$\Sigma(i):\Sigma(m)\to\Sigma(S)$ for some {\em injective} $i:\stm\to S$.

Now this remark implies immediately that such homomorphisms form a cofinal 
set in $\catN_\Sigma/I_\Sigma(S)$, and that both these categories are 
in fact filtering.
\end{Proof}

\nxsubpoint\label{sp:alg.bimod} (Algebraic bimodules.)
Now suppose we are given two algebraic monads $\Sigma$ and $\Lambda$, 
and we want to describe the category $\cA_{alg,\Sigma}^\Lambda\subset
\cA_\Sigma^\Lambda$ of {\em algebraic $(\Lambda,\Sigma)$-bimodules}, 
i.e.\ of algebraic functors $G:\cC\to\cC$, equipped both with a 
left $\Lambda$-action $\alpha:\Lambda G\to G$ and a right $\Sigma$-action 
$\beta:G\Sigma\to G$, these two actions being required to commute: 
$\alpha\circ(\Lambda\star\beta)=\beta\circ(\alpha\star\Sigma)$.

Recall that the category of all $(\Sigma,\Lambda)$-bimodules is equivalent 
to the category of functors $\catFunct(\cC_\Sigma,\cC^\Lambda)$, 
since $\cA_\Sigma\cong\catFunct(\cC_\Sigma,\cC)$ 
by~\ptref{prop:univ.lower.sigma}, and this 
equivalence is compatible with the left $\obslash$-actions of $\cA$, 
so we get $\cA_\Sigma^\Lambda=(\cA_\Sigma)^\Lambda\cong
\catFunct(\cC_\Sigma,\cC)^\Lambda\cong\catFunct(\cC_\Sigma,\cC^\Lambda)$ 
by~\ptref{prop:left.sigma.funct}.

Similarly, the category $\cA_{alg,\Sigma}\subset\cA_\Sigma$ is equivalent 
to $\catFunct(\catN_\Sigma,\cC)$ by \ptref{sp:altdescr.rmod}, and this 
equivalence is again compatible with the left $\obslash$-action of $\cA$, 
hence $\cA_{alg,\Sigma}^\Lambda\cong\catFunct(\catN_\Sigma,\cC)^\Lambda\cong 
\catFunct(\catN_\Sigma,\cC^\Lambda)=\catFunct(\catN_\Sigma,\catMod\Lambda)$. 
In the ``cosimplicial set'' setting this means that $\Lambda$ controls the 
category $\catMod\Lambda$ where our ``cosimplicial objects'' take values, 
and $\Sigma$ affects the ``index category'' $\catN$, replacing the category 
of standard finite sets $\catN$ with a richer category $\catN_\Sigma$. 
It might be interesting, for example, to consider the case $\Lambda=\Sigma=
{\bm\Delta}$.

Of course, if we take $\Lambda=\Fempty$, then $\cA_{alg,\Sigma}^{\Fempty}=
\cA_{alg,\Sigma}$, so all the properties proved below for algebraic bimodules 
can be specialized to the case of algebraic right modules.

\nxsubpoint (Matrix description of algebraic bimodules.)
Of course, an algebraic bimodule $(G,\alpha,\beta)$ with 
$\alpha:\Lambda G\to G$ and $\beta:G\Sigma\to G$ can be described in terms 
of a collection of sets $\{G(n)\}_{n\geq0}$, maps $G(\phi):G(m)\to G(n)$ 
for each $\phi:\stm\to\stn$, and some maps 
$\alpha^{(k)}_n:\Lambda(k)\times G(n)^k\to G(n)$ and 
$\beta^{(k)}_n:G(k)\times\Sigma(n)^k\to G(n)$, subject to all conditions 
listed before for left and right algebraic modules, as well as a certain 
commutativity relation for $\alpha$ and $\beta$. We have a matrix 
interpretation as before; then this commutativity relation can be 
interpreted again as some associativity of matrix multiplication 
$M(k,s;\Lambda)\times M(n,k;G)\times M(m,n;\Sigma)\to M(m,s;G)$. 
For example, if $G$ is an algebraic $\Sigma$-bimodule, i.e.\ a 
$(\Sigma,\Sigma)$-bimodule, then we obtain some multiplication maps 
of matrices over $G$ with matrices over $\Sigma$ from both sides, 
with products being matrices over~$G$.

\nxsubpoint (Operations with algebraic bimodules.)
Equivalence of categories $\cA_{alg,\Sigma}^\Lambda\cong
\catFunct(\catN_\Sigma,\catMod\Lambda)$ shows 
that arbitrary inductive and projective limits, subobjects, 
strict quotients, and so on, can be defined and computed componentwise 
in $\catMod\Lambda$, similarly to \ptref{sp:limits.leftmod}. Moreover, 
any algebraic monad homomorphism $\sigma:\Lambda\to\Lambda'$ induces 
a scalar restriction functor $\sigma^*:\cA_{alg,\Sigma}^{\Lambda'}\to
\cA_{alg,\Sigma}^\Lambda$, and its left adjoint, the scalar extension or 
base change functor $\sigma_*$, both of which can be computed componentwise 
by composing functors from $\catN_\Sigma$ into $\catMod\Lambda$ or 
$\catMod{\Lambda'}$ with scalar restriction or extension functors between 
$\catMod\Lambda$ and $\catMod{\Lambda'}$. Therefore, we can fix $\Lambda$ 
and study what happens when we change $\Sigma$.

\nxsubpoint (Scalar restriction and extension.)
Now let's take some algebraic monad homomorphism $\rho:\Sigma\to\Xi$, 
and study the scalar restriction functor $\rho^*:\cA_{alg,\Xi}^\Lambda\to
\cA_{alg,\Sigma}^\Lambda$. We know that $\rho^*:\cA_\Xi^\Lambda\cong 
\catFunct(\cC_\Xi,\cC^\Lambda)\to\cA_\Sigma^\Lambda\cong
\catFunct(\cC_\Sigma,\cC^\Lambda)$ is given by pre-composing functors 
$\cC_\Xi\to\cC^\Lambda$ with a certain ``base change'' functor 
$\rho_*:\cC_\Sigma\to\cC_\Xi$ (cf.~\ptref{sp:gen.scalrestr.mod}). 
This implies that $\rho^*:\cA_{alg,\Xi}^\Lambda\cong\catFunct(\catN_\Xi,
\cC^\Lambda)\to\cA_{alg,\Sigma}^\Lambda\cong\catFunct(\catN_\Sigma,
\cC^\Lambda)$ is given by pre-composing functors 
$\catN_\Xi\to\cC^\Lambda$ with the restriction $\bar\rho_*:\catN_\Sigma\to
\catN_\Xi$ of $\rho_*$ to $\catN_\Sigma$, i.e.\ 
$\rho^*$ is essentially equal to $(\bar\rho_*)^*$.

Now we see that $\rho_*=(\bar\rho_*)^*$ admits both a left and a right 
adjoint, namely, the Kan extensions $\rho_!:=(\bar\rho_*)_!$ and 
$\rho_?:=(\bar\rho_*)_*:\catFunct(\catN_\Sigma,\catMod\Lambda)\to
\catFunct(\catN_\Xi,\catMod\Lambda)$, computed by means of appropriate 
inductive and projective limits in $\catMod\Lambda$.

\nxsubpoint\label{sp:tensprod.bimod} (Tensor products of bimodules.)
Given an algebraic $(\Xi,\Sigma)$-bimodule $M$ and an algebraic 
$(\Sigma,\Lambda)$-bimodule~$N$, we would like to construct a new 
algebraic $(\Xi,\Lambda)$-bimodule $M\otimes_\Sigma N$, and similarly, 
for any $\Sigma$-module $X$ we would like to obtain a $\Xi$-module 
$M\otimes_\Sigma X$. Interpreting all bimodules involved as functors between 
certain categories, we see that all we have to do is to extend in some 
way $M:\catN_\Sigma\to\catMod\Xi$ to a functor $\tilde M:\catMod\Sigma\to
\catMod\Xi$; then we'll put $M\otimes_\Sigma N:=\tilde M\circ N$ and 
$M\otimes_\Sigma X:=\tilde M(X)$. Since the canonical functor 
$Q'_\Sigma:\catN_\Sigma\to\catMod\Sigma$ is fully faithful (with the essential 
image equal to the subcategory of free $\Sigma$-modules of finite rank), 
we have a natural choice of $\tilde M$, given by the 
left Kan extension: $\tilde M:=Q'_{\Sigma,!}(M)$. This defines our 
tensor product $M\otimes_\Sigma N$. 

Clearly, $M\mapsto\tilde M=Q'_{\Sigma,!}M$ commutes with arbitrary inductive 
limits in~$M$, hence the same can be said about $M\otimes_\Sigma N$. 
Another useful observation: since $Q'_{\Sigma,!}=Q_{\Sigma,!}\circ
J_{\Sigma,!}$, we see that we might first extend $M:\catN_\Sigma\to\cC^\Lambda$
to a functor $M':\cC_\Sigma\to\cC^\Lambda$ by means of $J_{\Sigma,!}$, 
and then extend this $M'$ to $\tilde M:\cC^\Sigma\to\cC^\Lambda$. Since 
$M=M'$ is the algebraic functor corresponding to~$M$ 
(cf.~\ptref{sp:altdescr.rmod}), we see that for any set $S$ we have 
$\tilde M(\Sigma(S))=M(S)$, i.e.\ $M\otimes_\Sigma\Sigma(S)=M(S)$.

From this and similar observations we deduce that all scalar restrictions 
and extensions with respect to right and left algebraic module structures 
can be written in terms of this tensor product in the usual way. For example, 
if $\rho:\Sigma\to\Xi$ is a homomorphism, $N$ is any $\Sigma$-module, 
and $P$ is any $\Xi$-module, then the scalar extension $\rho_*N$ 
is isomorphic to $\Xi\otimes_\Sigma N$, where $\Xi$ is considered as a 
$(\Xi,\Sigma)$-bimodule, and the scalar restriction $\rho^*P$ can 
be identified with $\Xi\otimes_\Xi P$, where this time $\Xi$ is 
considered as a $(\Sigma,\Xi)$-bimodule.

\nxpointtoc{Addition. Hypoadditivity and hyperadditivity}\label{p:hypadd}
We have seen that the multiplication on the underlying set $|\Sigma|$ of 
an algebraic monad $\Sigma$ is always defined, so in some sense multiplication 
is more fundamental than addition. Now we would like to study different 
versions of addition and additivity of algebraic monads.

We fix some algebraic monad~$\Sigma$ and a constant $0\in\Sigma(0)$. Since 
$\Fone=\Fempty\langle 0^{[0]}\rangle$, we see that this is equivalent to 
fixing some $\rho:\Fone\to\Sigma$, i.e.\ $\Sigma$ is an algebraic monad 
over~$\Fone$. In most cases $\Sigma$ is supposed to be an {\em algebraic 
monad with zero}, i.e.\ with exactly one constant 
(cf.~\ptref{sp:mon.with.zero}), so we have only one way of 
choosing~$0\in\Sigma(0)$.

\nxsubpoint\label{sp:compar.maps} (Comparison maps.)
Once the zero constant $0\in\Sigma(0)$ is fixed, we can define the 
{\em comparison maps} $\pi_n=\pi_{\Sigma,n}:\Sigma(n)\to\Sigma(1)^n=
|\Sigma|^n$ as follows. 
The $k$-th component $\pi_n^k$ of $\pi_n$ is given by 
$\pi_n^k:t\mapsto t(0_{\Sigma(1)},\ldots,\{1\}_\st1,\ldots,0_{\Sigma(1)})$ 
for any $t\in\Sigma(n)$, i.e.\ we substitute $0=0_{\Sigma(1)}\in\Sigma(1)$ 
for all formal variables $\{i\}_\stn$ in~$t=t(\{1\}_\stn,\ldots,\{n\}_\stn)$, 
with the only exception of~$\{k\}_\stn$.

If we are given a monad homomorphism $\rho:\Sigma\to\Xi$, and if we choose 
the constant $\rho_0(0)\in\Xi(0)$ to construct $\pi_{\Xi,n}$, then we 
get $\pi_{\Xi,n}\circ\rho_n=|\rho|^n\circ\pi_{\Sigma,n}$, i.e.\ the 
comparison maps are functorial in~$\Sigma$.

Notice that the choice of $0\in\Sigma(0)$ gives us a homomorphism 
$\st1\to\Sigma(0)$ from the final to the initial $\Sigma$-module 
(when $\Sigma$ is a monad with zero, this is an isomorphism), so for 
any direct sum $M\oplus N$ in $\catMod\Sigma$ we get canonical 
homomorphisms $M\oplus N\to M\oplus\st1\to M\oplus\Sigma(0)\cong M$ and 
$M\oplus N\to N$, which define together a {\em comparison homomorphism}
$\pi_{M,N}:M\oplus N\to M\times N$. Of course, this construction can 
be generalized to arbitrary finite and even infinite sums and products, 
thus defining comparison homomorphisms $\bigoplus_i M_i\to\prod_i M_i$.

For example, the comparison homomorphism for $n$ copies of $\Sigma(1)$ 
is a homomorphism $\Sigma(n)=\Sigma(1)\oplus\cdots\oplus\Sigma(1)\to
\Sigma(1)^n$, 
which is exactly the comparison map $\pi_n$ considered before.

\begin{DefD}\label{def:ps.add} (Addition and pseudoaddition.)
We say that a binary operation $[+]\in\Sigma(2)$ is a {\em pseudoaddition} 
(with respect to constant $0\in\Sigma(0)$, when $\Sigma$ has more than 
one constant), if $\pi_2([+])=(\bu,\bu)\in\Sigma(1)^2$, in other words, 
if $\{1\}+0=\{1\}=0+\{1\}$ in $|\Sigma|$, i.e.\ if $x+0=x=0+x$ for any 
$\Sigma$-module~$X$ and any $x\in X$. We say that 
$[+]$ is {\em the addition of~$\Sigma$} if $[+]$ is the only pseudoaddition 
of an algebraic monad with zero~$\Sigma$. In this case we say that
{\em $\Sigma$ is a monad with addition.}
\end{DefD}
Of course, any monad homomorphism $\rho:\Sigma\to\Xi$ transforms a 
pseudoaddition into a pseudoaddition, and if both $\Sigma$ and $\Xi$ 
are monads with addition, then any $\rho:\Sigma\to\Xi$ respects this addition.

We can define pseudoadditions of higher arities $t\in\Sigma(n)$ by requiring 
$\pi_n(t)=(\bu,\bu,\ldots,\bu)\in\Sigma(1)^n$; if we have a binary 
pseudoaddition $[+]$, we can construct pseudoadditions of all higher arities 
by taking $(\cdots((\{1\}+\{2\})+\{3\})+\cdots)+\{n\}$.

\nxsubpoint (Alternative definition of addition.)
Arguably better definitions of addition and zero are these:
a {\em zero\/} is simply any {\em central\/} constant in an algebraic 
monad~$\Sigma$ in the sense of~\ptref{def:comm}
(such central constant is necessarily unique, cf.~\ptref{sp:comm.lower.arity}),
and an {\em addition\/} in an algebraic monad with zero $\Sigma$ is any
central pseudoaddition (it is automatically the only pseudoaddition
of~$\Sigma$ by~\ptref{sp:comm.and.add}, i.e.\ is an addition in the
sense of~\ptref{def:ps.add}). The only reason why we don't adopt and
study these stronger definitions here is that 
they require the commutativity notions of the next Chapter, 
so we postpone a more detailed study of the
relationship between (pseudo)addition and commutativity 
until~\ptref{sp:comm.and.add}.

\begin{DefD}\label{def:hypohyperadd}
We say that an algebraic monad~$\Sigma$ is {\em hypoadditive} 
(resp.\ {\em additive, hyperadditive}) if all comparison maps 
$\pi_n:\Sigma(n)\to\Sigma(1)^n$ are injective (resp.\ bijective, surjective).
\end{DefD}
Clearly, $\Sigma$ is additive iff it is both hypoadditive and hyperadditive. 
Notice that additivity and hypoadditivity imply that $\pi_0:\Sigma(0)\to\st1$ 
is injective, so $\Sigma$ has at most one constant, so it has to be a monad 
with zero, and we don't need to specify our choice of $0\in\Sigma(0)$ while 
discussing hypoadditivity or additivity.

Of course, there can be monads which are neither hypoadditive nor 
hyperadditive; some of them, like $\Finfty$, will play important role in 
the sequel. Notice that any submonad $\Sigma'$ of a hypoadditive 
monad~$\Sigma$, 
containing the zero of $\Sigma$, is automatically hypoadditive, and 
similarly, any strict quotient of a hyperadditive monad is hyperadditive; 
this follows immediately from the functoriality of comparison maps.

\nxsubpoint\label{sp:hyperadd.psadd}
If $\Sigma$ is additive, then it has an addition $[+]=\pi_2^{-1}(\bu,\bu)$, 
i.e.\ $\Sigma$ is then a monad with addition. Conversely, if 
$\Sigma$ is hypoadditive and admits a pseudoaddition $[+]$, it is 
necessarily unique (since $\pi_2$ is injective), so $\Sigma$ becomes 
a monad with addition; moreover, this addition is automatically 
commutative, since $\pi_2(\{2\}+\{1\})=(\bu,\bu)=\pi_2(\{1\}+\{2\})$, 
and associative, since both $(\{1\}+\{2\})+\{3\}$ and $\{1\}+(\{2\}+\{3\})$
are mapped into $(\bu,\bu,\bu)$ by $\pi_3$. We see that 
$\{1\}+\{2\}+\cdots+\{n\}$ is the unique pseudoaddition of arity~$n$, and 
for any choice of unary operations $\lambda_1$, \dots, $\lambda_n\in|\Sigma|$ 
we get an $n$-ary operation $\lambda_1\{1\}+\cdots+\lambda_n\{n\}\in\Sigma(n)$,
mapped into $(\lambda_1,\ldots,\lambda_n)$ by~$\pi_n$, so {\em a hypoadditive 
monad~$\Sigma$ with pseudoaddition $[+]$ is in fact additive}, and in this case
any $n$-ary operation $t\in\Sigma(n)$ can be uniquely written 
in form $\lambda_1\{1\}+\cdots+\lambda_n\{n\}$ with $\lambda_i\in|\Sigma|$.

One shows, essentially in the same way, that {\em an algebraic monad~$\Sigma$ 
is hyperadditive iff it admits a pseudoaddition}. Indeed, if $\Sigma$ 
is hyperadditive, it admits a pseudoaddition, $\pi_2:\Sigma(2)\to|\Sigma|^2$ 
being surjective; conversely, if $[+]$ is a pseudoaddition, then for 
any $\lambda=(\lambda_1,\ldots,\lambda_n)\in|\Sigma|^n$ we can construct an 
$n$-ary operation $(\cdots((\lambda_1\{1\}+\lambda_2\{2\})+\lambda_3\{3\})+
\cdots)+\lambda_n\{n\}$, mapped into $\lambda$ by $\pi_n$.

\nxsubpoint\label{sp:addmon.semirings} (Additive monads and semirings.)
This implies that an additive monad $\Sigma$ is generated (over~$\Fempty$) 
by its only constant 
$0\in\Sigma(0)$, the addition $[+]\in\Sigma(2)$, and the set of its unary 
operations~$\Sigma(1)$. We have also some relations: $0+0=0$, the 
commutativity $\{1\}+\{2\}=\{2\}+\{1\}$ and the associativity 
$(\{1\}+\{2\})+\{3\}=\{1\}+(\{2\}+\{3\})$ of the addition,
$\lambda\cdot 0=0$ for any $\lambda\in|\Sigma|$; the addition 
$[+]_{|\Sigma|}:|\Sigma|^2\to|\Sigma|$ gives a family of relations of form 
$\lambda=\lambda'+\lambda''$ for all $\lambda'$, $\lambda''\in|\Sigma|$, 
we have also a family of relations $\lambda=\lambda'\lambda''$ coming 
from the monoid structure of~$|\Sigma|$,
and finally, we have relations $\lambda(\{1\}+\{2\})=\lambda\{1\}+\lambda\{2\}$
for any $\lambda\in|\Sigma|$, since both sides are mapped by~$\pi_2$
into~$(\lambda,\lambda)$. Notice that this list of generators and relations 
is a presentation of~$\Sigma$, since these relations are sufficient 
to rewrite any $t\in\Sigma(n)$ in form $\lambda_1\{1\}+\cdots+\lambda_n\{n\}$.

Conversely, if we have any algebraic monad~$\Sigma$ generated by such a list 
of generators and relations, i.e.\ if we are given some set $|\Sigma|$ 
together with a constant $0\in|\Sigma|$ and two binary operations 
$+,\times:|\Sigma|^2\to|\Sigma|$, subject to a family of relations which 
turn out to be exactly the axioms for a semiring structure 
$(0,+,\times)$ on $|\Sigma|$, then $\Sigma$ is automatically additive, 
since the list of relations given above allows us to write any $t\in\Sigma(n)$ 
in form $\lambda_1\{1\}+\cdots+\lambda_n\{n\}$ with $\lambda_i\in|\Sigma|$,
and $\catMod\Sigma$ is the category of $|\Sigma|$-modules, i.e.\ 
additively written commutative monoids with a biadditive action of~$|\Sigma|$.

We see that additive algebraic monads~$\Sigma$ are in one-to-one 
correspondence with semirings~$|\Sigma|$.

\nxsubpoint\label{sp:addmon.rings} (Additive monads and rings.)
In particular, any monad~$R=\Sigma_R$, defined by a ring~$R$, is additive. 
Conversely, if $\Sigma$ is additive, and if we have a unary operation 
$[-]\in\Sigma(1)$, such that $[-]^2=\bu$ (giving such a $[-]$ is equivalent 
to giving a homomorphism $\Fpm=\Fone\langle-^{[1]}|-^2=\bu$, 
$-0=0\rangle\to\Sigma$), and such that $\{1\}+(-\{1\})=0$ (this condition 
doesn't follow from the others), then $R:=|\Sigma|$ is a ring, and 
$\Sigma\cong\Sigma_R$.

\begin{PropD}\label{prop:crit.add} (Criteria of additivity.)
The following conditions for an algebraic monad~$\Sigma$ are equivalent:
(i) $\Sigma$ is additive; (ii) All comparison maps $\pi_n:\Sigma(n)\to
\Sigma(1)^n$ are isomorphisms; (iii) All comparison homomorphisms 
$\pi_{M,N}:M\oplus N\to M\times N$ are isomorphisms;
(iii') All comparison homomorphisms $\bigoplus_{i=1}^kM_i\to\prod_{i=1}^kM_i$ 
are isomorphisms;
(iv) Comparison homomorphisms $\pi_{\Sigma(n),\Sigma(1)}:\Sigma(n+1)
\cong\Sigma(n)\oplus\Sigma(1)\to\Sigma(n)\times\Sigma(1)$ are isomorphisms 
for all $n\geq0$; (v) It is possible to introduce on each 
$\Hom_\Sigma(M,N)$ a commutative monoid structure (usually written in additive 
form), such that all composition maps $\Hom_\Sigma(N,P)\times\Hom_\Sigma(M,N)
\to\Hom_\Sigma(M,P)$ become biadditive; (vi) $\catMod\Sigma$ is equivalent 
to the category of modules over a semiring~$R$.
\end{PropD}
\begin{Proof}
(i)$\Leftrightarrow$(ii) is true by definition, and (ii)$\Leftrightarrow$(iv) 
is shown by induction in $n\geq0$, using $\pi_{n+1}=
(\pi_n\times\id_{\Sigma(1)})\circ\pi_{\Sigma(n),\Sigma(1)}$. 
Now (iii)$\Leftrightarrow$(iii') by an easy induction in~$k$, and 
(iii)$\Rightarrow$(iv) is evident. We know (i)$\Rightarrow$(vi) already. 
Let's show (vi)$\Rightarrow$(v). When we have (vi), we can 
define an addition on each $\Hom_\Sigma(M,N)\cong\Hom_R(M',N')$ pointwise, 
where $M'$ and $N'$ are the $R$-modules corresponding to $M$ and $N$: 
$(\phi_1+\phi_2)(x):=\phi_1(x)+\phi_2(x)$; it is immediate that 
sum of two $R$-homomorphisms is an $R$-homomorphism again, that 
this addition on $\Hom_R(M',N')\cong\Hom_\Sigma(M,N)$ 
is commutative and associative, with 
the zero map $M'\to0\to N'$ for zero, and that the composition 
maps are biadditive with respect to this addition, so we obtain 
(vi)$\Rightarrow$(v). 

It remains to show (v)$\Rightarrow$(iii).
First of all, (v) implies that the final object $\st1$ of $\catMod\Sigma$ 
is an initial object as well, i.e.\ $\catMod\Sigma$ admits a zero object. 
Indeed, for any~$M$ the set $\Hom_\Sigma(\st1,M)$ is a monoid under addition, 
so it cannot be empty; on the other hand, $\End_\Sigma(\st1)$ consists 
of exactly one element $\id_\st1$, so it has to be the zero of this 
additive monoid, and biadditivity of $\Hom_\Sigma(\st1,M)\times
\End_\Sigma(\st1)\to\Hom_\Sigma(\st1,M)$ implies that for any $\phi:\st1\to M$
we have $\phi=\phi\circ\id_\st1=\phi\circ0_\st1=0$, so $\Hom_\Sigma(\st1,M)=0$,
$\st1$ is the zero object of $\catMod\Sigma$, and for any $M$ and $N$ the 
composite homomorphism $M\to\st1\to N$ is the zero of $\Hom_\Sigma(M,N)$. 
In particular, we have shown that $\Sigma(0)\cong\st1$, i.e.\ $\Sigma$ is 
a monad with zero.

Now let's take any two $\Sigma$-modules $M_1$ and $M_2$ and prove that 
$\pi_{M_1,M_2}:M_1\oplus M_2\to M_1\times M_2$ is an isomorphism. 
Let's denote the projections $M_1\times M_2\to M_i$ by $p_i$, and the 
composite maps $M_i\to M_1\oplus M_2\to M_1\times M_2$ by $\lambda_i$.
These maps are defined by $p_i\lambda_i=\id_{M_i}$, $p_i\lambda_{3-i}=0$, 
and $\pi_{M_1,M_2}$ is an isomorphism iff $M_1\times M_2$ together with these 
maps $\lambda_i$ satisfies the universal property for $M_1\oplus M_2$. 
Notice that $\lambda_1p_1+\lambda_2p_2=\id_{M_1\times M_2}$, since 
$p_i(\lambda_1p_1+\lambda_2p_2)=(p_i\lambda_i)p_i+0\cdot p_{3-i}=p_i$.
Now for any $f:M_1\times M_2\to N$ we have $f=f(\lambda_1p_1+\lambda_2p_2)=
f_1p_1+f_2p_2$, so $f$ is completely determined by its ``components'' 
$f_i:=f\lambda_i:M_i\to N$. Conversely, given any $f_i:M_i\to N$, we can 
define $f:M_1\times M_2\to N$ by $f:=f_1p_1+f_2p_2$, and then clearly 
$f\lambda_i=f_i$, so $M_1\times M_2$ is indeed the direct sum of $M_1$ and 
$M_2$, and $\pi_{M_1,M_2}$ is an isomorphism.
\end{Proof}

This proposition, together with some more similar reasoning, implies the 
following corollary:
\begin{CorD}
The following conditions for an algebraic monad~$\Sigma$ are equivalent:
(i) $\Sigma\cong\Sigma_R$ for some associative ring~$R$;
(ii) $\Sigma$ is additive, and there is a unary operation $[-]\in\Sigma(1)$, 
such that $\{1\}+(-\{1\})=0$;
(iii) $\catMod\Sigma$ is equivalent to the category of modules over an 
associative ring~$R$;
(iv) $\catMod\Sigma$ is additive;
(v) $\catMod\Sigma$ is abelian.
\end{CorD}

\nxsubpoint (Hypoadditivity.)
Some of the above conditions generalize to the hypoadditive case. For example, 
{\em $\Sigma$ is hypoadditive iff all comparison maps~$\Sigma(n+1)\cong
\Sigma(n)\oplus\Sigma(1)\to\Sigma(n)\times\Sigma(1)$ are injective.} 
One might expect all comparison maps $M\oplus N\to M\times N$ to be 
injective for a hypoadditive~$\Sigma$; 
however, this is not true, as shown by the following example. 
Let's take $\Sigma:=\Zinfty$; since $\Zinfty$ is a submonad of an 
additive monad~$\bbR$, containing the zero constant of~$\bbR$, it is 
hypoadditive. Consider the equivalence relation~$R$ on $\Zinfty(2)=
\{\lambda\{1\}+\mu\{2\}\;|\;|\lambda|+|\mu|\leq1\}$, which identifies 
all interior points of~$\Zinfty(2)$, i.e.\ all elements with 
$|\lambda|+|\mu|<1$, with zero. It is easy to see that 
$R$ is compatible with the $\Zinfty$-structure of~$\Zinfty(2)$, so we 
can compute the strict quotient $M:=\Zinfty(2)/R$, and let's denote by~$\pi$ 
the canonical projection $\Zinfty(2)\to M$. Clearly, $M\oplus M$ is 
a strict quotient of $\Zinfty(4)=\Zinfty(2)\oplus\Zinfty(2)$; let's 
denote the projection by $\pi'$ and its kernel by $R'$. It is easy to check 
that~$R'$ identifies all interior points of $\Zinfty(4)$ together, and doesn't
identify any boundary points. In particular, 
$x:=\pi'(\frac12\{1\}+\frac12\{3\})$ and $y:=\pi'(\frac12\{2\}+\frac12\{4\})$ 
are distinct in~$M\oplus M$, but the comparison map $M\oplus M\to M\times M$ 
maps both of them into $(0,0)$, so it cannot be injective.

\nxsubpoint (Importance of hypoadditivity.)
Notice that {\em any homomorphism $\rho:\Sigma\to\Xi$ of hypoadditive monads 
is completely determined by the map $|\rho|=\rho_1:|\Sigma|\to|\Xi|$.} 
Indeed, the functoriality of comparison maps  
$\pi_{\Xi,n}\circ\rho_n=|\rho|^n\circ\pi_{\Sigma,n}$ shows that 
$\rho_n$ is completely determined by $|\rho|$ when $\pi_{\Xi,n}$ is injective. 
This can be seen as a manifestation of the following fact: {\em most properties
of a hypoadditive monad~$\Sigma$ can be seen at the level of $|\Sigma|$.}
Conversely, there are some non-hypoadditive monads, like $\Aff_R\subset R$ 
(cf.~\ptref{sp:examp.submon},g)), which are highly non-trivial 
($\catMod{\Aff_R}$ is the category of ``affine spaces over~$R$''), but 
still have $|\Aff_R|=\st1$, so almost all structure of~$\Aff_R$ 
is outside of $|\Aff_R|$. This remark will have important consequences in 
our theory of spectra of generalized rings, i.e.\ commutative algebraic monads.

\nxsubpoint (Examples.)
a) Of course, all monads~$\Sigma_R$, defined by some associative ring~$R$, 
are additive, hence all their submonads $\Sigma'\subset R$, containing the 
zero constant of $\Sigma_R$, are hypoadditive. This applies to almost all 
examples from~\ptref{sp:examp.submon}. In particular, $\Zinfty$, $\barZinfty$, 
$\Zninfty$, $\Fone$ and $\Fpm$ are hypoadditive. However, 
$\Aff_R\subset R$, ${\bm\Delta}\subset\bbR$ and $\Fempty\subset\Fone$ are 
not hypoadditive, since all of them are monads without constants. 

b) Consider the algebraic monad~$\Sigma$, defined by the category of 
commutative rings. Then $\Sigma(n)\cong\bbZ[T_1,\ldots,T_n]$. If we fix 
some constant $c\in\Sigma(0)=\bbZ$, say, $c=0$, then the comparison maps 
look like $\pi_2:F(X,Y)\mapsto(F(T,0),F(0,T))$. Clearly, all these maps 
are surjective, so $\Sigma$ is hyperadditive. It admits a lot of different 
pseudoadditions like $X+Y$, $X+Y+XY$ or even $X+Y+X^2Y$; this shows that 
a pseudoaddition need not be commutative or associative.

c) We'll see later that most non-commutative tensor products like 
$\bbZ\boxtimes_{\Fone}\bbZ$ are hyperadditive. 
This means that we should not expect to read all properties of such monads 
from~$\Sigma(1)$.

d) All additive monads are generated in arity $\leq2$, and the 
same is true for all hypoadditive monads listed in a). Actually, all 
hypoadditive monads we encountered so far are generated in arity $\leq2$, 
and can be embedded into some additive monads. However, 
it is not clear why this should be true in general.
Consider the hypoadditive submonad 
$\bbZ_{\sqrt\infty}\subset\bbR$, such that $\bbZ_{\sqrt\infty}(n)$ 
consists of all formal linear combinations $\lambda_1\{1\}+\cdots+
\lambda_n\{n\}$ with $\lambda_i\in\bbR$ and $\sum_i\sqrt{|\lambda_i|}\leq1$. 
It is easy to see that such linear combinations are stable under substitution, 
so $\bbZ_{\sqrt\infty}$ is indeed a hypoadditive submonad of $\bbR$. 
However, it still turns out to be generated in arity $\leq2$. Another
example is given by $A_3:=\Zinfty\cap\bbZ[1/3]\subset\bbQ$, which is 
hypoadditive and generated over $\Fpm$ by the ternary ``averaging operation''
$s_3:=(1/3)\{1\}+(1/3)\{2\}+(1/3)\{3\}$, but can be shown not to be generated
in arity $\leq 2$.

\nxsubpoint\label{sp:examp.finfty}
Let's construct now a monad with zero, which is neither hypoadditive 
nor hyperadditive. Consider for this the equivalence relation $R$ on 
$\Zinfty(1)=[-1,1]$, which identifies all interior points of $\Zinfty(1)$ 
with zero. We know that~$R$ is compatible with the $\Zinfty$-structure on 
$\Zinfty(1)$, so we obtain a strict quotient $Q=\Zinfty/R\cong\{-1,0,1\}$ of 
$\Zinfty(1)$, which has been denoted by $\Finfty$ in~\ptref{sp:def.finf}. 
The $\Zinfty$-structure on~$Q$ defines a homomorphism of algebraic monads 
$\rho:\Zinfty\to\END(Q)$, cf.~\ptref{sp:def.absend.ring}; 
let's denote its image by~$\Finfty$. 
Clearly, $\Finfty$ is a strict quotient of $\Zinfty$, and two elements 
$t$, $t'\in\Zinfty(n)$ have the same image in $\Finfty(n)$ iff the 
induced maps ${[t]}_Q$ and ${[t']}_Q:Q^n\to Q$ coincide. In particular, 
all $\Finfty(n)\subset\END(Q)(n)=\Hom_\catSets(Q^n,Q)$ are finite.
It is easy to see that all elements $t=\lambda_1\{1\}+\cdots+\lambda_n\{n\}\in
\Zinfty(n)$ with $\sum_i|\lambda_i|<1$ are identified with~$0$ in 
$\Finfty(n)$, and all remaining elements are identified iff their 
sequences of signs $(\sgn(\lambda_1),\ldots,\sgn(\lambda_n))$ coincide, 
so we can write corresponding elements of $\Finfty(n)$ in form 
$+?\{1\}+?\{2\}-?\{3\}+?\{5\}$, where the question marks replace arbitrary 
positive real numbers with sum one; clearly, elements of $\Finfty(n)$ 
correspond to faces of the standard octahedron~$\Zinfty(n)$. In particular, 
$\card\Finfty(n)=3^n=\card|\Finfty|^n$; however, the comparison maps 
$\pi_n:\Finfty(n)\to|\Finfty|^n$ are not bijective, since they map 
all elements different from $\pm\{i\}$ for $1\leq i\leq n$ into zero.

\nxsubpoint\label{sp:finfty.gens0}
Notice that $\Finfty$ is generated by one constant $0^{[0]}$, one 
unary operation $-^{[1]}$, and one binary operation 
$*^{[2]}:=+?\{1\}+?\{2\}$. We have some relations between them, like 
$-0=0$, $0*0=0$, $-(-x)=x$, $x*0=0$, $x*x=x$, $x*(-x)=0$, 
$(-x)*(-y)=-(x*y)$, $x*y=y*x$ and 
$(x*y)*z=x*(y*z)$. This system of generators and relations constitutes 
a presentation of $\Finfty$, since it already enables us  
to represent any non-zero element of $\Sigma(n)$, 
where $\Sigma$ is the generalized 
ring generated by these operations and relations, 
in form $\pm\{i_1\}*\cdots*\pm\{i_k\}$ 
with some $1\leq i_1<\cdots<i_k\leq n$, hence the surjection
$\Sigma\twoheadrightarrow\Finfty$ is an isomorphism, and 
$\Finfty=\Sigma$ is finitely presented (over~$\Fempty$).

\nxpointtoc{Algebraic monads over a topos}\label{p:algmon.over.topos}
Most of the results of this chapter generalize, when properly understood, to 
the topos case. We start with a na{\"\i}ve approach:

\nxsubpoint\label{sp:algmon.anycat} 
(Algebraic monads over an arbitrary category.)
Recall that in \ptref{sp:expl.algmon.cond} we have described an algebraic 
monad $\Sigma$ over $\catSets$ as a sequence of sets $\{\Sigma(n)\}_{n\geq0}$, 
together with some maps $\Sigma(\phi):\Sigma(m)\to\Sigma(n)$ for each 
$\phi:\stm\to\stn$, an element $\bu\in\Sigma(1)$, and some 
multiplication maps $\mu_n^{(k)}:\Sigma(k)\times\Sigma(n)^k\to\Sigma(n)$, 
subject to some identities between certain maps, constructed from these maps 
by means of finite products and composition. Of course, this definition 
generalizes immediately to an arbitrary cartesian category~$\cC$, 
if we require all $\Sigma(n)$ to be objects of~$\cC$ instead of being sets, 
all $\Sigma(\phi)$ and $\mu_n^{(k)}$ to be morphisms in~$\cC$, and 
choose some $\bu\in\Gamma_\cC(\Sigma(1))=\Hom_\cC(e_\cC,\Sigma(1))$, where 
$e_\cC$ is the final object of~$\cC$. Given any such ``algebraic monad'' 
over~$\cC$, we can define its action on a object~$X$ of~$\cC$ as a collection 
of morphisms $\alpha^{(k)}:\Sigma(k)\times X^k\to X$, satisfying the 
conditions of~\ptref{sp:expl.action}; this defines 
{\em the category $\cC^\Sigma=\catMod\Sigma$ of $\Sigma$-modules in~$\cC$.} 
Furthermore, we can define left and right algebraic modules, algebraic monad 
homomorphisms and so on simply by copying their descriptions in terms of 
sets or sequences of sets and maps between products of these sets.
Notice that ``algebraic endofunctors over~$\cC$'' are nothing else than 
functors~$\catN\to\cC$.

In fact, we can generalize these definitions even to the case when involved 
finite products are not representable in~$\cC$, by identifying $\cC$ 
with a full subcategory of 
the category $\hat\cC=\catFunct(\cC^0,\catSets)$ of presheaves of sets. 
Notice that giving an algebraic monad in $\hat\cC$ is equivalent to giving 
a presheaf of algebraic monads over~$\cC$, i.e.\ a functor from $\cC^0$ 
to the category of algebraic monads over~$\catSets$.

\nxsubpoint (Application of left exact functors.)
The advantage of this approach is that we can apply any left exact functor 
$h:\cC\to\cC'$ to any of these collections of data. For example, if 
$\Sigma$ is an algebraic monad over~$\cC$, then $h\Sigma$, defined by 
$(h\Sigma)(n):=h\Sigma(n)$, is an algebraic monad over~$\cC'$, and 
for any $\Sigma$-module~$X$ we obtain a canonical $h\Sigma$-structure 
on~$hX$. In this way any such left exact functor $h:\cC\to\cC'$ induces a 
functor $h=h^\Sigma:\cC^\Sigma\to{\cC'}^{h\Sigma}$ between corresponding 
categories of modules. 

\nxsubpoint\label{sp:appl.top.funct} (Application to topoi.)
This is especially useful for the left exact functors that arise from 
morphisms of topoi $f:\cE'\to\cE$, i.e.\ the pullback functor 
$f^*:\cE\to\cE'$ and the direct image functor $f_*:\cE'\to\cE$. In this 
way $f^*$ induces a functor $f^{*,\Sigma}:\cE^\Sigma\to{\cE'}^{f^*\Sigma}$ 
for any algebraic monad~$\Sigma$ over~$\cE$, and $f_*$ induces a functor 
$f_*^{\Sigma'}:{\cE'}^{\Sigma'}\to\cE^{f_*\Sigma'}$ for any algebraic 
monad~$\Sigma'$ over~$\cE'$. Moreover, canonical adjointness morphisms 
$\Sigma\to f_*f^*\Sigma$ and $f^*f_*\Sigma'\to\Sigma'$ are easily seen to 
be monad homomorphisms, so we can combine $f_*^{f^*\Sigma}$ with the 
scalar restriction along $\Sigma\to f_*f^*\Sigma$, thus obtaining a 
functor ${\cE'}^{f^*\Sigma}\to\cE^\Sigma$.

We can apply this in partucular to the canonical morphism $q:\cE\to\catSets$ 
from an arbitrary topos~$\cE$ into the point topos~$\catSets$. Then 
$q_*X=\Gamma_\cE(X)=\Hom_\cE(e_\cE,X)$ for any $X\in\Ob\cE$, i.e.\ 
$q_*$ is the functor of global sections, and $q^*:S\mapsto S_\cE=\bigsqcup_{s\in S}e_\cE$ is the ``functor of constant objects''. We see that for any 
algebraic monad~$\Sigma$ over~$\cE$ we have its global sections monad 
$\Gamma_\cE(\Sigma)$ over~$\catSets$, that any algebraic monad~$\Xi$ over 
$\catSets$ defines a constant algebraic monad $\Xi_\cE:=q^*\Xi$ over~$\cE$, 
and so on. Notice that we could also apply to~$\Sigma$ the functor of sections 
$\Gamma_S:X\mapsto\Hom_\cE(S,X)$ over any object~$S\in\Ob\cE$.

Another application: if $p:\catSets\to\cE$ is any point of~$\cE$, and 
$p^*:X\mapsto X_p$ is the corresponding fiber functor, then we can compute 
the fiber at~$p$ of any algebraic monad $\Sigma$ over~$\cE$ by putting 
$\Sigma_p:=p^*\Sigma$.

Finally, if $\cS$ is a (small) site, we have the topos of presheaves~$\hat\cS$ 
and the topos of sheaves~$\tilde\cS$, the inclusion functor $i:\tilde\cS\to
\hat\cS$, and its left adjoint, the ``sheafification functor'' 
$a:\hat\cS\to\tilde\cS$. Both of them are left exact, and 
$ai\cong\Id_{\tilde\cS}$ since $i$ is fully faithful; this allows us to treat 
any algebraic monad~$\Sigma$ in~$\tilde\cS$, i.e.\ a sheaf of algebraic monads 
over~$\cS$, as a presheaf of monads $i\Sigma$, and conversely, to 
``sheafify'' any presheaf of algebraic monads, obtaining a sheaf of algebraic 
monads, i.e.\ an algebraic monad in~$\tilde\cS$.

\nxsubpoint\label{sp:monadic.gammas} (Monadicity of forgetful functors.)
Given any algebraic monad~$\Sigma$ over a topos~$\cE$, we have a 
forgetful functor $\Gamma_\Sigma:\cE^\Sigma\to\cE$ (not to be confused with 
the global sections functor $\Gamma_\cE:\cE\to\catSets$). We claim that 
{\em $\Gamma_\Sigma:\cE^\Sigma\to\cE$ is monadic, i.e.\ it admits a left 
adjoint $L_\Sigma:\cE\to\cE^\Sigma$, and induced functor 
$\cE^\Sigma\to\cE^{\tilde\Sigma}$, where $\tilde\Sigma:=\Gamma_\Sigma 
L_\Sigma$ is a monad over~$\cE$, is an equivalence of categories.} 
We are going to show this in several steps.

\nxsubpoint\label{sp:top.ex.lsigma}
First of all, let's show the existence of a left adjoint~$L_{\tilde\Sigma}$ 
to~$\Gamma_\Sigma$. Any topos~$\cE$ is equivalent to the category~$\tilde\cS$ 
of presheaves on some site~$\cS$, which can be chosen to be small. 
Consider first the category of presheaves $\hat\cS\supset\tilde\cS$. 
The collection of sheaves $\Sigma=\{\Sigma(n)\}$ is an ``algebraic monad'' 
over~$\hat\cS$ as well, so we have a forgetful functor 
$\Gamma_{\hat\Sigma}:\hat\cS^\Sigma\to\hat\cS$. On the other hand, 
$\Sigma$ defines an algebraic monad $\Sigma_{[S]}$ over~$\catSets$ 
for any object $S\in\Ob\cS$, once we apply to $\Sigma$ the left exact 
functor $F\mapsto F(S)$. We see that a $\Sigma$-structure $\alpha$ on a 
presheaf $F\in\Ob\hat\cS$, given by a collection of morphisms 
of presheaves $\alpha^{(k)}:\Sigma(k)\times F^k\to F$, is nothing else 
than a collection of $\Sigma_{[S]}$-structures on each set~$F(S)$, 
compatible with the ``restriction maps'' $F(S)\to F(S')$ and 
$\Sigma_{[S]}\to\Sigma_{[S']}$ for any morphism $\phi:S'\to S$ in $\cS$. 
Now we can define a left adjoint~$L_{\hat\Sigma}$ to $\Gamma_{\hat\Sigma}$ 
by putting $(L_{\hat\Sigma}F)(S):=L_{\Sigma_{[S]}}(F(S))$ for any 
presheaf of sets~$F$. Clearly, $\hat\Sigma:=\Gamma_{\hat\Sigma}L_{\hat\Sigma}$ 
is given by $(\hat\Sigma F)(S):=\Sigma_{[S]}(F(S))$. The monadicity of 
$\Gamma_{\hat\Sigma}:\hat\cS^\Sigma\to\hat\cS$ follows now pointwise 
from the monadicity of all~$\Gamma_{\Sigma_{[S]}}$.

We extend these results from $\hat\cS$ to $\tilde\cS$ in the usual way 
by putting $L_{\tilde\Sigma}:=aL_{\hat\Sigma}i$ and $\tilde\Sigma:=
\Gamma_{\tilde\Sigma}L_{\tilde\Sigma}=
a\hat\Sigma i$, where $a:\hat\cS\to\tilde\cS$ and $i:\tilde\cS\to\hat\cS$ 
denote the ``sheafification'' and inclusion functors as usual; in the first 
expression $a$ actually denotes the extension $a^\Sigma:\hat\cS^{i\Sigma}\to
\tilde\cS^{ai\Sigma}=\tilde\cS^\Sigma$; it exists because of the left exactness
of~$a$ and full faithfullness of~$i$. Now the adjointness of 
$L_{\tilde\Sigma}$ and $\Gamma_{\tilde\Sigma}$ is immediate: 
$\Hom_{\tilde\cS^\Sigma}(L_{\tilde\Sigma}F,G)=
\Hom_{\tilde\cS^\Sigma}(aL_{\hat\Sigma}iF,G)=
\Hom_{\hat\cS^\Sigma}(L_{\hat\Sigma}iF,iG)=
\Hom_{\hat\cS}(iF,\Gamma_{\hat\Sigma}iG)=\Hom_{\hat\cS}(iF,
i\Gamma_{\tilde\Sigma}G)=\Hom_{\tilde\cS}(F,\Gamma_{\tilde\Sigma}G)$.

It remains to show the monadicity of~$\Gamma_{\tilde\Sigma}$; however, 
we postpone this proof until~\ptref{sp:topos.same.algmon}.

\nxsubpoint (Algebraic endofunctors on a topos.)
Recall that an ``algebraic endofunctor'' $\Sigma$ on a topos $\cE$ is given 
by a collection of objects $\{\Sigma(n)\}_{n\geq0}$ and morphisms 
$\Sigma(\phi):\Sigma(m)\to\Sigma(n)$ for each $\phi:\stm\to\stn$, i.e.\ 
by a functor $\Sigma\in\Ob\cE^\catN=\Ob\catFunct(\catN,\cE)$. We would like 
to extend $\Sigma$ to a ``true'' endofunctor $\tilde\Sigma:\cE\to\cE$, 
such that any morphism $\alpha:\tilde\Sigma(X)\to Y$ would be given by 
a collection of morphisms $\alpha^{(k)}:\Sigma(k)\times X^k\to Y$ 
as in~\ptref{sp:algmap.SXtoY}. Clearly, this requirement determines 
$\tilde\Sigma(X)$ uniquely up to a unique isomorphism, and we see 
that $\tilde\Sigma(X)$ can be constructed in the same way as 
in~\ptref{l:alt.descr.algfun}, i.e.\ $\tilde\Sigma(X)=\Coker(p,q:H_1(X)
\rightrightarrows H_0(X))$, where $H_0(X)=\bigsqcup_{n\geq0}\Sigma(n)\times 
X^n$ and $H_1(X)=\bigsqcup_{\phi:\stm\to\stn}\Sigma(m)\times X^n$.

When $\cE=\tilde\cS$, we might first construct an extension of $\Sigma$ 
to $\hat\Sigma:\hat\cS\to\hat\cS$ with the correct universal property 
in~$\hat\cS$; clearly, $\hat\Sigma F$ is still computed by the same formulas 
as above in the category of presheaves, and all limits in this category are 
computed componentwise, so we have $(\hat\Sigma F)(S)=\Sigma_{[S]}(F(S))$, 
where $\Sigma_{[S]}$ is the algebraic endofunctor over~$\catSets$, 
constructed from the collection of sets $(\Sigma(n))(S)$ as before. Then 
we can put $\tilde\Sigma:=a\hat\Sigma i:\tilde\cS\to\tilde\cS$; clearly, 
$\Hom_{\tilde\cS}(\tilde\Sigma F, G)\cong\Hom_{\hat\cS}(\hat\Sigma iF,iG)$ 
can be still described as in~\ptref{sp:algmap.SXtoY}, so we get an 
alternative construction of the extension~$\tilde\Sigma$. In particular, 
our notations are compatible with those of~\ptref{sp:top.ex.lsigma}. 
Moreover, we see that the canonical morphism $a\hat\Sigma\to a\hat\Sigma ia=
\tilde\Sigma a$ is an isomorphism. This follows from the fact that 
$\hat\Sigma F=\Coker(H_1(F)\rightrightarrows H_0(F))$ is computed 
for any presheaf~$F$ by means of inductive limits and finite products, and 
$a:\hat\cS\to\tilde\cS$ commutes with these types of limits.

This observation can be used to compute the values of $\tilde\Sigma$ on 
finite constant sheaves $\stn_\cE$: since $\stn_\cE=a\stn_{\hat\cS}$, and 
$\stn_{\hat\cS}$ is the constant presheaf $S\mapsto\stn$, we see that 
$\tilde\Sigma(\stn_{\tilde\cS})=a\hat\Sigma(\stn_{\hat\cS})=\Sigma(n)$, 
since $\hat\Sigma(\stn_{\hat\cS})$ is the presheaf $S\mapsto\Sigma_{[S]}(\stn)=
\Sigma_{[S]}(n)=(\Sigma(n))(S)$, equal to~$\Sigma(n)$. One checks similarly 
that $\Sigma(\phi)=\tilde\Sigma(\phi_\cE)$ for any map of 
finite sets~$\phi:\stm\to\stn$.  

\nxsubpoint (Universal description of extensions.)
We have just constructed a functor $\tilde J_\cE:
\cE^\catN=\catFunct(\catN,\cE)\to\catEndof(\cE)$, $\Sigma\mapsto\tilde\Sigma$, 
and a restriction functor $J_\cE^*:\catEndof(\cE)\to\cE^\catN$ in the opposite 
direction, where $J_\cE:\catN\to\cE$ is given by $\stn\mapsto\stn_\cE$, 
such that $J_\cE^*\tilde J_\cE\cong\Id_{\cE^\catN}$. This implies that 
$\tilde J_\cE$ is faithful, so it induces an equivalence between 
$\cE^\catN$ and some subcategory $\catEndof_{alg}(\cE)$ 
of~$\catEndof(\cE)$, which merits 
to be called {\em the category of algebraic endofunctors on~$\cE$}, 
and the quasi-inverse equivalence is given by the restriction of~$J_\cE^*$. 
In this respect everything is quite similar to~\ptref{l:alt.descr.algfun}; 
however, $\tilde J_\cE$ cannot be the left Kan extension $J_{\cE,!}$ of 
$J_\cE$, since $J_\cE$ is in general not fully faithful, so we cannot expect 
$J_\cE^*J_{\cE,!}\cong\Id$. However, a similar description can be obtained 
if we replace $\catEndof(\cE)$ with the category of {\em (plain) inner 
endofunctors\/}~$\catInnEndof(\cE)$ (cf.~\ptref{sp:def.plain.monad}):

\begin{PropD}
Denote by $J_\cE:\catN\to\cE$ the canonical functor, and consider the 
restriction functor $J_\cE^*:\catInnEndof(\cE)\to\cE^\catN$. This restriction 
functor admits a left adjoint $\tilde J_\cE:\cE^\catN\to\catInnEndof(\cE)$, 
such that the composite functor $\cE^\catN\to\catInnEndof(\cE)\to\catEndof(\cE)$
coincides with the extension functor $\tilde J_\cE:\Sigma\mapsto\tilde\Sigma$ 
constructed before. Moreover, $\tilde J_\cE:\cE^\catN\to\catInnEndof(\cE)$ 
is fully faithful, i.e.\ $J_\cE^*\tilde J_\cE\cong\Id$, so the category 
$\cE^\catN$ of ``algebraic endofunctors'' over~$\cE$ is equivalent to 
the essential image of $\tilde J_\cE$,
a full subcategory $\catInnEndof_{alg}(\cE)$ of $\catInnEndof(\cE)$, 
which will be called {\em the category of algebraic (plain) inner endofunctors
on~$\cE$}. In particular, $\catEndof_{alg}(\cE)\subset\catEndof(\cE)$ 
is equivalent to this category, hence any algebraic endofunctor on~$\cE$ 
admits a canonical (algebraic) plain inner structure.
\end{PropD}
\begin{Proof}
First of all, let's construct a canonical plain inner structure on 
each endofunctor $\tilde\Sigma:\cE\to\cE$, obtained by extending some 
$\Sigma:\catN\to\cE$; this would define a functor $\tilde J_\cE:
\cE^\catN\to\catInnEndof(\cE)$, such that $J_\cE^*\tilde J_\cE\cong\Id$. 
So we need to construct a functorial family of morphisms 
$\alpha_{X,Y}:X\times\tilde\Sigma(Y)\to\tilde\Sigma(X\times Y)$,
defined for any $X$, $Y\in\Ob\cE$, satisfying the conditions 
of~\ptref{def:p-inner.funct}. Since $\tilde\Sigma(Y)=\Coker(H_1(Y)
\rightrightarrows H_0(Y))$, and similarly for $\tilde\Sigma(X\times Y)$, 
we get $X\times\tilde\Sigma(Y)=\Coker(X\times H_1(Y)\rightrightarrows 
X\times H_0(Y))$, since $j_X^*:S\mapsto X\times S$ commutes with arbitrary 
inductive limits, having a right adjoint $\iHom_\cE(X,-)$, so
we see that it is sufficient to define some morphisms 
$X\times H_0(Y)=X\times\bigsqcup_{n\geq0}(\Sigma(n)\times Y^n)=
\bigsqcup_{n\geq0}(X\times\Sigma(n)\times Y^n)\to H_0(X\times Y)=
\bigsqcup_{n\geq0}(\Sigma(n)\times(X\times Y)^n)$. Of course, we can 
construct such a morphism by defining its $n$-th component
$X\times\Sigma(n)\times Y^n\to\Sigma(n)\times (X\times Y)^n$ with the aid 
of the diagonal morphism~$X\to X^n$. We define 
$X\times H_1(Y)\to H_1(X\times Y)$ in a similar way, thus obtaining a 
morphism $\alpha_{X,Y}:X\times\tilde\Sigma(Y)\to\tilde\Sigma(X\times Y)$ 
after taking cokernels. This construction is clearly functorial in~$X$ and~$Y$;
in order to check the remaining conditions for~$\alpha$ we simply observe 
that in case $\cE=\tilde\cS$ we might construct $\alpha$ first on 
the category of presheaves $\hat\cS$, by considering the unique plain 
inner structure on each~$\Sigma_{[S]}$ over~$\catSets$, and then 
``sheafify'' this plain inner structure on~$\hat\Sigma$, thus obtaining 
a plain inner structure on~$\tilde\Sigma$, easily seen to coincide with 
one constructed before in an invariant fashion.

We already know that $J_\cE^*\tilde\Sigma\cong\Sigma$, i.e.\ 
$J_\cE^*\tilde J_\cE\cong\Id$, so it remains to show that $\tilde J_\cE$ 
is indeed a left adjoint to $J_\cE^*$, i.e.\ that for any 
$\Sigma:\catN\to\cE$ and any $\Sigma'=(\Sigma',\alpha')\in\Ob
\catInnEndof(\cE)$ we have
$\Hom_{\catInnEndof(\cE)}(\tilde\Sigma,\Sigma')\cong
\Hom_{\cE^\catN}(\Sigma,J_\cE^*\Sigma')$, where $\tilde\Sigma=(\tilde\Sigma,
\alpha):=\tilde J_\cE\Sigma$ as before.

So we have to show that natural transformations of inner endofunctors 
$\tilde\gamma:\tilde\Sigma\to\Sigma'$, $\tilde\gamma_X:\tilde\Sigma(X)\to
\Sigma'(X)$, are in one-to-one correspondence with natural transformations 
$\gamma:\Sigma\to J_\cE^*\Sigma'$, $\gamma_n:\Sigma(n)\to(J_\cE^*\Sigma')(n)=
\Sigma'(\stn_\cE)$. Since $\tilde\Sigma(X)=\Coker(H_1(X)\rightrightarrows
H_0(X))$, we see that giving a $\tilde\gamma_X:\tilde\Sigma(X)\to\Sigma'(X)$ 
is equivalent to giving a family of morphisms 
$\tilde\gamma_X^{(n)}:\Sigma(n)\times X^n\to\Sigma'(X)$, such that 
$\tilde\gamma_X^{(n)}\circ(\Sigma(\phi)\times\id_{X^n})=
\tilde\gamma_X^{(m)}\circ(\id_{\Sigma(m)}\times X^\phi)$ for any 
$\phi:\stm\to\stn$ (cf.~\ptref{sp:algmap.SXtoY}). Let's denote by $\theta_X^n:
X^n\times\Sigma(n)\to\Sigma'(X)$ the same morphism as $\tilde\gamma_X^{(n)}$ 
with interchanged arguments. Since $X^n=\iHom_\cE(\stn_\cE,X)$, we have 
some evaluation morphisms $\ev_{\stn,X}:X^n\times\stn_\cE\to X$. We know 
that $\tilde\Sigma(\stn_\cE)\cong\Sigma(n)$, and $\tilde\gamma$ is supposed 
to be compatible with plain inner structures; this proves the commutativity 
of the following diagram:
\begin{equation}
\xymatrix@C+10pt{
X^n\times\Sigma(n)\ar[r]^{\alpha_{X^n,\stn_\cE}}
\ar[d]^{\id_{X^n}\times\gamma_n}\ar@{-->}[rrd]^<>(.75){\theta_X^n}
&\tilde\Sigma(X^n\times\stn_\cE)\ar[d]\ar[r]^<>(.5){\tilde\Sigma(\ev_{\stn,X})}
&\tilde\Sigma(X)\ar[d]^{\tilde\gamma_X}\\
X^n\times\Sigma'(\stn_\cE)\ar[r]^{\alpha'_{X^n,\stn_\cE}}
&\Sigma'(X^n\times\stn_\cE)\ar[r]_<>(.5){\Sigma'(\ev_{\stn,X})}
&\Sigma'(X)
}
\end{equation}
Notice that the upper row defines the canonical morphism 
$X^n\times\Sigma(n)\to\tilde\Sigma(X)$, given by $X^n\times\Sigma(n)\to 
H_0(X)\to\tilde\Sigma(X)$, hence the diagonal arrow has to be equal to 
$\theta_X^n$. This shows that the collection of the $\theta_X^n$, hence 
also $\tilde\gamma_X$, is completely determined by the collection of 
$\gamma_n$, i.e.\ by~$\gamma$. Conversely, given any $\gamma:\Sigma\to
J_\cE^*\Sigma'$, we can construct the morphisms $\theta_X^n:X^n\times\Sigma(n)
\to\Sigma'(X)$ from the above diagram; 
they are easily seen to satisfy necessary compatibility 
conditions, so they define some $\tilde\gamma_X:\tilde\Sigma(X)\to\Sigma'(X)$,
clearly functorial in~$X$.
\end{Proof}

\nxsubpoint (Algebraic inner endofunctors and filtered inner inductive limits.)
In the case $\cE=\catSets$ we had an equivalent description of 
the category of algebraic endofunctors $\catEndof_{alg}(\catSets)\subset
\catEndof(\catSets)$: it was the full subcategory, consisting of those 
endofunctors, which commute with filtered inductive limits. In the topos 
case all algebraic endofunctors are easily seen to commute with 
arbitrary filtered inductive limits; however, this property is insufficient 
to be an algebraic endofunctor over a topos~$\cE$, since in general 
an arbitrary object~$X$ of~$\cE$ cannot be represented as a filtered inductive 
limit of finite constant objects $\stn_\cE$.

It is possible to save the situation by considering the category of 
algebraic inner endofunctors $\catInnEndof_{alg}(\cE)\subset\catInnEndof(\cE)$:
these inner endofunctors can be characterized by their property to commute 
with filtered {\em inner} (or {\em local\/}) inductive limits, 
i.e.\ inductive limits 
along filtered inner index categories $\cI$ in~$\cE$ 
(recall that this means that $\cI$ has an object of objects $\cI_0$ and 
an object of morphisms $\cI_1$ instead of the sets of objects and morphisms,
together with some morphisms $i:\cI_0\to\cI_1$, $s$, $t:\cI_1\to\cI_0$, 
$m:\cI_1\times_{t,s}\cI_1\to\cI_1$, subject to certain conditions which 
need not be written here). For example, any $X\in\Ob\cE$ can be written 
as an inner filtered inductive limit along a certain inner category 
$\cI_X$, with $\cI_{X,0}=\bigsqcup_{n\geq0} X^n$ and 
$\cI_{X,1}=\bigsqcup_{\phi:\stm\to\stn}X^n$. If we write down explicitly 
the requirement for an inner endofunctor $\Sigma'=(\Sigma',\alpha')$ to 
commute with this particular inner inductive limit, we obtain 
$\Sigma'\cong\tilde J_\cE J_\cE^*\Sigma'$, and one can check that 
any functor of form $\tilde J_\cE\Sigma$ commutes with filtered inner 
inductive limits, essentially in the same way as 
in~\ptref{l:alt.descr.algfun}. 

Such an approach has its advantages. For example, it shows immediately 
that {\em the composite of two algebraic inner endofunctors is again 
algebraic}.
However, we don't want to follow this path, since it would require to explain 
in more detail the theory of inner categories and limits; instead, 
we prove the above statement about composites of algebraic endofunctors 
directly:

\begin{PropD}
The composite of two algebraic inner endofunctors is again algebraic, i.e.\ 
$\catInnEndof_{alg}(\cE)$ is a full $\otimes$-subcategory of 
$\catInnEndof(\cE)$, and $\catEndof_{alg}(\cE)$ is a $\otimes$-subcategory 
of $\catEndof(\cE)$. Therefore, we have a canonical AU $\otimes$-structure 
on $\cE^\catN\cong\catInnEndof_{alg}(\cE)$, and a canonical 
$\obslash$-action of this category on~$\cE$. Moreover, 
$\catInnEndof_{alg}(\cE)$ is stable under arbitrary inductive and finite 
projective limits of~$\catInnEndof(\cE)$, and $\tilde J_\cE$ commutes 
with these types of limits.
\end{PropD}
\begin{Proof}
First of all, notice that $\catInnEndof(\cE)$ admits arbitrary 
projective and inductive limits, which can be computed componentwise. 
Indeed, to define a plain inner structure $\alpha_{X,Y}:X\times
(\projlim\Sigma_\iota)(Y)\to(\projlim\Sigma_\iota)(X\times Y)$ we 
simply consider the morphism with components $X\times\projlim\Sigma_\iota(Y)
\to X\times\Sigma_\iota(Y)\to\Sigma_\iota(X\times Y)$, coming from the 
plain inner structures on~$\Sigma_\iota$. The case of inductive limits is 
similar, once we observe that the functor 
$S\mapsto X\times S$ commutes with arbitrary inductive limits,
$\cE$ being cartesian closed.

Now notice that $\cE^\catN$ also admits arbitrary projective and inductive 
limits, which can be also computed componentwise, hence 
$J_\cE^*:\catInnEndof(\cE)\to\cE^\catN$ commutes with all limits, 
and $\tilde J_\cE$ commutes with arbitrary inductive limits, having a 
right adjoint~$J_\cE^*$. We claim that, similarly to what we had 
in~\ptref{prop:j.descr.algfun}, {\em $\tilde J_\cE$ commutes with 
finite projective limits as well}, hence {\em $\catInnEndof_{alg}(\cE)$ 
is stable under arbitrary inductive and finite projective limits 
of~$\catInnEndof(\cE)$}. We have to show that 
$\widetilde{\projlim\Sigma_\iota}\to\projlim\tilde\Sigma_\iota$ is an 
isomorphism for any finite projective limit $\projlim\Sigma_\iota$ 
in~$\cE^\catN$. We can assume $\cE=\tilde\cS$ for some site~$\cS$; in this case
our statement can be first shown for the extension $\hat\Sigma$ 
of~$\Sigma$ to the category of presheaves $\hat\cS$, 
using~\ptref{prop:j.descr.algfun} componentwise; then we extend this 
result to $\tilde\Sigma=a\hat\Sigma i$, using the left exactness of~$a$.

Notice that {\em $\tilde J_\cE$ transforms a constant functor 
$\catN\to\cE$ with value $Z$ into a constant endofunctor~$\cE\to\cE$ 
with the same value}. Indeed, this is clear for $\cE=\catSets$, 
constant functors being algebraic, then the statement can be 
shown for $\hat\cS$, applying the result for~$\catSets$ componentwise, 
and then we can extend it to the general case $\cE=\tilde\cS$ by 
sheafification. 

Now let's show that the composite $\tilde\Xi\tilde\Sigma$ of two 
algebraic inner endofunctors is again algebraic. Indeed, for any 
$X\in\Ob\cE$ we have $\tilde\Xi(X)=\Coker(H_1(X)\rightrightarrows H_0(X))$, 
hence $\tilde\Xi\tilde\Sigma=\Coker(H_1(\tilde\Sigma)\rightrightarrows 
H_0(\tilde\Sigma))$, where $H_0(\tilde\Sigma)=\bigsqcup_{n\geq0}
\Xi(n)\times\tilde\Sigma^n$ and $H_1(\tilde\Sigma)=\bigsqcup_{\phi:\stm\to\stn}
\Xi(m)\times\tilde\Sigma^n$ are computed in~$\catInnEndof(\cE)$. 
This allows us to compute~$\tilde\Xi\tilde\Sigma$ in
$\catInnEndof(\cE)$ by means of finite products and some inductive limits, 
starting from algebraic inner endofunctor $\tilde\Sigma$ and constant 
inner endofunctors with values~$\Xi(n)$, which are also algebraic; 
hence $\tilde\Xi\tilde\Sigma$ is algebraic as well.
\end{Proof}

\nxsubpoint\label{sp:topos.same.algmon}
Once we know that $\catInnEndof_{alg}(\cE)\subset\catInnEndof(\cE)$ 
is stable under composition, arbitrary inductive, and finite projective limits,
and that morphisms $\tilde\Sigma(X)\to Y$ admit a description similar 
to that given in~\ptref{sp:algmap.SXtoY} for any algebraic endofunctor~$\Sigma$
on~$\cE$ and any two objects~$X$ and~$Y$ of~$\cE$, 
we can extend all statements from~\ptref{p:op.algendof} 
and~\ptref{p:algmon} to the topos case. In particular, we see that 
an algebraic inner monad~$\Sigma$ over~$\cE$, i.e.\ an algebra 
in $\catInnEndof_{alg}(\cE)$, can be described in the same way 
as in~\ptref{sp:algmon.predescr} and~\ptref{sp:expl.algmon.cond}, 
so we obtain the same notion of an algebraic monad over~$\cE$ as 
in~\ptref{sp:algmon.anycat}, and the same notions of modules over 
such an algebraic monad. This proves the monadicity of 
$\Gamma_\Sigma:\cE^\Sigma\to\cE$, claimed in~\ptref{sp:monadic.gammas}.

\nxsubpoint
Most statements and definitions given before in this chapter for algebraic 
monads over~$\catSets$ generalize directly to the topos case, at least 
those of them which are intuitionistic, i.e.\ don't involve the logical law 
of excluded middle and the axiom of choice. However, we have to interpret 
these statements in the so-called {\em Kripke--Joyal semantics}, where 
existence is always understood as local existence, i.e.\ existence after 
restricting everything to a suitable cover. 
Consider for example~\ptref{sp:expl.SX}, 
which tells us that any element $\xi$ of $\Sigma(X)$ can be written in form 
$t(\{x_1\},\ldots,\{x_n\})$ for some $n\geq0$, $t\in\Sigma(n)$ and 
$x_i\in X$. In the topos case we can denote by 
$t(\{x_1\},\ldots,\{x_n\})$ the image of $(t,x_1,\ldots,x_n)\in
\Sigma(n)(S)\times X(S)^n$ in $\Sigma(X)(S)$, for any $S\in\Ob\cE$ 
(or in $\Ob\cS$, when $\cE=\tilde\cS$) and any $t\in\Sigma(n)(S)$ and 
$x_i\in X(S)$. However, the above statement has to be interpreted as follows:
{\em for any $\xi\in (\Sigma(X))(S)$ we can find a cover $\{\phi_\alpha:
S_\alpha\to S\}$, some integers $n_\alpha\geq0$, elements 
$t_\alpha\in(\Sigma(n_\alpha))(S_\alpha)$ and $x_{\alpha,i}\in X(S_\alpha)$, 
$1\leq i\leq n_\alpha$, such that $\phi_\alpha^*\xi=t_\alpha(
\{x_{\alpha,1}\},\ldots,\{x_{\alpha,n_\alpha}\})$ in $(\Sigma(X))(S_\alpha)$ 
for any~$\alpha$.}

\nxsubpoint
However, some statements over $\catSets$ have been proved before 
with the aid of the axiom of choice. We need to find other proofs of such 
statements, or to extend these results over $\catSets$ first to the categories 
of presheaves $\hat\cS$, applying already known results componentwise, and 
then extend these statements to arbitrary $\cE=\tilde\cS$ by sheafification.

For example, we have seen that algebraic monads (even endofunctors) 
over~$\catSets$ preserve epimorphisms and monomorphisms 
(cf.~\ptref{sp:sigma.injtoinj} and~\ptref{sp:modhom.img}); 
however, our proofs used the axiom of choice, so they cannot be generalized 
to the topos case directly. Nevertheless, we claim that {\em 
any algebraic endofunctor~$\Sigma$ over a topos~$\cE$ preserves monomorphisms 
and epimorphisms.}

The statement about monomorphisms is easily deduced from the known case 
by sheafification, since both $i:\tilde\cS\to\hat\cS$ and $a:\hat\cS\to
\tilde\cS$ are left exact, and $\tilde\Sigma=a\hat\Sigma i$. 
The statement about epimorphisms cannot be easily shown in this way, since 
$i$ doesn't respect epimorphisms; however, it can be deduced from the fact 
that any finite product of epimorphisms in a topos is an epimorphism again, 
hence for any epimorphism $f:X\to Y$ the induced morphism 
$H_0(f):H_0(X)=\bigsqcup_{n\geq0}\Sigma(n)\times X^n\to H_0(Y)$ is again 
an epimorphism, and $H_0(Y)\to\Sigma(Y)$ is a (strict) epimorphism as well, 
hence $H_0(X)\to\Sigma(X)\to\Sigma(Y)$ is an epimorphism, hence 
$\Sigma(f):\Sigma(X)\to\Sigma(Y)$ has to be an epimorphism as well.

Since all epimorphisms in~$\cE$ are strict, and any morphism in~$\cE$ 
decomposes uniquely into an epimorphism followed by a monomorphism, we 
are in position to generalize~\ptref{sp:modhom.img} to the topos case, 
proving that the image of any morphism of $\Sigma$-modules in~$\cE$ admits 
a $\Sigma$-structure. Most statements after~\ptref{sp:modhom.img} can 
be deduced from this fact in the same way as before. For example, 
arbitrary projective and inductive limits exist in~$\cE^\Sigma$.
 
\nxsubpoint
Some statements that involve the law of excluded middle also need to be 
reconsidered. For example, an algebraic monad~$\Sigma$ over a topos~$\cE$ 
is said to be {\em a monad without constants\/} (resp.\ {\em with zero,
with at most one constant}) iff $\Sigma(0)$ is isomorphic to the initial object
$\emptyset_\cE$ of~$\cE$ (resp.\ iff the canonical morphism 
$\Sigma(0)\to e_\cE$ is an isomorphism, or a monomorphism for the last case). 
We see that $\Sigma(0)$ can be a non-trivial open object of~$\cE$ (i.e.\ 
a subobject of the final object~$e_\cE$, different from both $\emptyset_\cE$ 
and $e_\cE$); in this case $\Sigma$ is a monad with at most one constant, 
but neither a monad without constants nor a monad with zero.

\nxsubpoint (Topos-specific properties.)
Finally, there are some topos-specific properties, very similar to those 
discussed in SGA~4~IV. For example, any morphism of topoi $f:\cE'\to\cE$ 
induces functors $f^{*,\Sigma}:\cE^\Sigma\to{\cE'}^{f^*\Sigma}$ and 
$f_*^\Sigma:{\cE'}^{f^*\Sigma}\to\cE^\Sigma$ for any algebraic monad 
$\Sigma$ over~$\cE$ (cf.~\ptref{sp:appl.top.funct}), which are easily seen 
to be adjoint. Since they coincide with $f^*$ and $f_*$ on the level of 
underlying objects, we usually omit the monad~$\Sigma$ in these notations.

In particular, we can apply this to localization morphisms 
$j_X:\cE_{/X}\to\cE$, for any object~$X$ of~$\cE$, and any algebraic 
monad~$\Sigma$ on~$\cE$. We usually write $\Sigma_{|X}$ or $\Sigma|X$ 
instead of $j_X^*\Sigma$; so we get a pair of adjoint functors 
$j_X^*=j_X^{*,\Sigma}:\cE^\Sigma\to(\cE_{/X})^{\Sigma_{|X}}$ and 
$j_{X,*}=(j_X)_*^\Sigma$ in the opposite direction, enjoying all the usual 
properties. Of course, we usually write $M_{|X}$ or $M|X$ instead of 
$j_X^{*,\Sigma}M$. Moreover, $j_X^{*,\Sigma}$ admits a left adjoint 
$j_{X,!}^\Sigma$, which has all the properties listed 
in SGA~4 IV~11.3.1; these properties can be shown in the same way as in 
SGA~4, and the existence of $j_{X,!}^\Sigma$ can be shown either 
by a reference to SGA~4 III~1.7, where the statement has been shown for 
arbitrary algebraic structures defined by finite projective limits, 
or directly: clearly, $j_{X,!}^\Sigma L_{\Sigma|X}=L_\Sigma j_{X,!}$, 
so we see how to define $j_{X,!}^\Sigma$ on free $\Sigma_{|X}$-modules; 
now by~\ptref{l:mod.coker.free} an arbitrary $\Sigma_{|X}$-module~$N$ 
can be represented as a cokernel of a pair of morphisms 
$p$, $q$ between some free $\Sigma_{|X}$-modules $L_{\Sigma|X}(R)$ and 
$L_{\Sigma|X}(S)$; 
since $j_{X,!}$ has to commute with arbitrary inductive limits, 
$j_{X,!}^\Sigma N$ has to be equal to the cokernel 
of $j_{X,!}^\Sigma(p)$ and $j_{X,!}^\Sigma(q)$, and this cokernel is 
easily seen to have the universal property required from $j_{X,!}^\Sigma N$, 
hence $j_{X,!}^\Sigma$ is well-defined.

\nxsubpoint (Open and closed subtopoi.)
Let $\cE$ be a topos, $\Sigma$ an algebraic monad over~$\cE$, 
$U\subset e_\cE$ an open object of~$\cE$, $\cU$ the open subtopos of~$\cE$ 
defined by~$U$, and $\cF$ the complementary closed subtopos of~$\cE$ 
(cf.~SGA~4 IV~9.3). Let's denote the corresponding embeddings by 
$j:\cU\to\cE$ and $i:\cF\to\cE$; recall that $\cU\cong\cE_{/U}$, and 
$i_*:\cF\to\cE$ is a fully faithful functor, the essential image of 
which consists of all $X\in\Ob\cE$, such that $j^*X=X_{|U}$ is isomorphic 
to the final object of $\cU=\cE_{/U}$. We obtain a left exact 
{\em gluing functor\/} $\rho:=i^*j_*:\cU\to\cF$, 
and we know that $\cE$ is canonically equivalent to the ``comma category'' 
$(\cF,\cU,\rho)$, the objects of which are triples $(A,B,f)$, with 
$A\in\Ob\cF$, $B\in\Ob\cU$ and $f:A\to\rho B$ a morphism in~$\cF$. 
Recall that this equivalence $\Phi:\cE\to(\cF,\cU,\rho)$ transforms 
$X\in\Ob\cE$ into $(i^*X,j^*X,h(X):i^*X\to \rho j^*X)$, where $h(X)$ 
is obtained by applying~$i^*$ to adjointness morphism $X\to j_*j^*X$ 
(cf.~SGA~4 IV~9.5). All functors $j_!$, $j^*$, $j_*$, $i^*$ and $i_*$ 
admit descriptions in terms of this equivalence: 
$j_!:B\mapsto(\emptyset_\cF,B,\emptyset_\cF\to\rho B)$, 
$j^*:(A,B,f)\mapsto B$, $j_*:B\mapsto (\rho B, B,\id_{\rho B})$, 
$i^*:(A,B,f)\mapsto A$, and $i_*:A\mapsto (A,e_\cU,A\to\rho e_\cU=e_\cF)$.

Since all functors $i^*$, $j^*$ and $\rho$ are left exact, we can extend 
the above results to the category of modules over an algebraic monad~$\Sigma$; 
we see that $\cE^\Sigma$ is equivalent to the comma category 
$(\cF^{i^*\Sigma},\cU^{j^*\Sigma},\rho)$, similarly to SGA~4 IV~14. 
Moreover, {\em when $\Sigma$ is an algebraic monad with zero, 
$i_*=i_*^\Sigma$ admits a right adjoint $i^!=i^{!,\Sigma}:\cE^\Sigma\to
\cF^{i^*\Sigma}$, given by $i^!:(A,B,f)\mapsto\Ker(f:A\to\rho B)$.} 
Here, of course, $\Ker(f)$ means the kernel of the pair, consisting of 
$f$ and of the zero morphism~$0:A\to\rho B$. We check that $i^!$ 
is a right adjoint to $i_*$, that $i_*i^!X\to X$ is a monomorphism, and 
that $i^!$ commutes with scalar restriction, essentially in the 
same way as in SGA~4 IV~14.5.

\nxsubpoint (Sheaves of homomorphisms.)
Let $\Sigma$ be an algebraic monad over a topos~$\cE$, and $M$, 
$N$ be two $\Sigma$-modules. Then we can construct the 
``sheaf of local homomorphisms'' $\iHom_\Sigma(M,N)\subset\iHom_\cE(M,N)$ 
in two different ways. On one hand, we can repeat the reasoning 
of SGA~4 IV~12 and define $\iHom_\Sigma(M,N)$ by requiring  
$\Hom_\cE(S,\iHom_\Sigma(M,N))\cong\Hom_{\Sigma|S}(M_{|S},N_{|S})$; 
this notion turns out to have almost all properties of SGA~4 IV~12, 
which can be proved essentially in the same way. On the other hand, 
$\Sigma$ admits a canonical plain inner structure, so we 
might use~\ptref{sp:gen.locinnhom} instead. {\em Both these definitions yield 
the same subobject of~$\iHom_\cE(M,N)$}. We don't give a detailed proof 
since we are not going to use this fact anyway; it can be shown, for example, 
by reduction to the case~$\cE=\hat\cS$.

\nxsubpoint (Small families of generators.)
Notice that {\em $\cE^\Sigma$ admits a small family of generators.} 
Indeed, any topos admits a small family of generators, so let 
${\mathcal G}\subset\Ob\cE$ be such a family for~$\cE$. Then the corresponding 
free $\Sigma$-modules $\{\Sigma(S)\}_{S\in\mathcal G}$ are easily seen 
to constitute a small family of generators of~$\cE^\Sigma$.

\nxsubpoint (Algebraic bimodules.)
Of course, we can consider algebraic bimodules (and left or right modules) 
over~$\cE$, i.e.\ 
bimodules in $\catInnEndof_{alg}(\cE)$. Most statements made 
in~\ptref{p:cat.algmod} generalize directly to the topos case. For example, 
we have matrix interpretations for all these notions, if we put 
$M(m,n;\Sigma):=\Sigma(m)^n$ as before; of course, now these are objects 
of~$\cE$, not just sets. However, $(\Sigma,\Lambda)$-bimodules cannot 
be interpreted as functors $\catN_\Lambda\to\cE^\Sigma$ when 
$\Lambda$ is non-constant, at least if we don't want to replace 
$\catN_\Lambda$ with an inner category in~$\cE$.

\nxsubpoint (Local $\END$'s.)
For any $X\in\Ob\cE$ we can construct an algebraic monad~$\iEND(X)$ over~$\cE$,
such that giving an action of some algebraic monad~$\Sigma$ on~$X$ 
is equivalent to giving a monad homomorphism $\rho:\Sigma\to\iEND(X)$. 
Indeed, we can put $(\iEND(X))(n):=\iHom_\cE(X^n,X)$, similarly 
to~\ptref{sp:def.absend.ring}. We can also consider the algebraic monad 
$\END(X):=\Gamma_\cE(\iEND(X))$ over~$\catSets$; then 
$(\END(X))(n)=\Hom_\cE(X^n,X)$, and monad homomorphisms 
$\Sigma\to\END(X)$ are in one-to-one correspondence with actions 
of {\em constant\/} algebraic monads~$\Sigma$ on~$X$. These notions 
generalize to the case of algebraic endofunctor~$X$ over~$\cE$, 
similarly to~\ptref{sp:def.absend.ring.endof}.


\cleardoublepage

\mysection{Commutative monads}

This chapter is dedicated to the notion of {\em commutativity\/} of 
algebraic monads, and to the properties of {\em generalized rings}, 
i.e.\ commutative algebraic monads, and categories of modules over them. 
We'll see that they have indeed properties very similar to classical 
commutative rings, so our terminology is justified. Moreover, sometimes 
we even consider arbitrary algebraic monads as {\em (non-commutative) 
generalized rings}; however, they can have properties very different from 
those of classical (non-commutative) rings.

We do some basic linear algebra over generalized rings; we construct a 
theory of exterior powers and determinants, at least for {\em alternating\/}
generalized rings, and define and compute some $K$-groups. 
We are also able to construct a reasonable theory of traces, notwithstanding 
the absence of addition. The properties of these theories are illustrated 
by a series of important examples.

\zeropoint\nxpoint
{\bf Notation.} Given an algebraic monad homomorphism $\rho:\Sigma\to\Xi$, 
we denote by ${}^a\rho^*:\catMod\Sigma\to\catMod\Xi$ or even by 
$\rho^*$  the base change functor 
(already denoted by $\rho_*$ in~\ptref{sp:ex.scalext}), and by 
${}^a\rho_*:\catMod\Xi\to\catMod\Sigma$ or even $\rho_*$ the scalar 
restriction functor (denoted before by $\rho^*$). This change of notation 
may seem slightly confusing, since ${}^a\rho^*$ is covariant in~$\rho$, 
and ${}^a\rho_*$ is contravariant; the idea, of course, is that we think 
of pullback and direct image functors for the corresponding morphism 
of (generalized) schemes ${}^a\rho:\Spec\Xi\to\Spec\Sigma$. 
This new notation has its advantages: for example, ${}^a\rho^*$ is now
a left adjoint to ${}^a\rho_*$.

Unless otherwise specified, all monads are supposed to be algebraic. 
We denote them by letters like $\Sigma$, $\Xi$, \dots, and also by 
$A$, $B$, $C$, \dots, especially when they are commutative; thus our notation 
reminds of classical rings. We also denote free $A$-module $L_A(S)$ by 
$A(S)$ or $A^{(S)}$; when $S=\stn=\{1,2,\ldots,n\}$, we write 
$A(n)$ or $A^{(n)}$, but never $A^n$; the latter notation actually means 
the product of $n$ copies of $A=A^{(1)}$, which is usually not isomorphic to 
$A^{(n)}$, at least for a non-additive~$A$.

We preserve the notation $|A|$ for $A(1)$, considered both as a monoid 
(with identity $\bu\in|A|$), and an $A$-module. We denote by 
$\|A\|$ the (graded) set of operations $\bigsqcup_{n\geq0}A(n)$.

\nxpointtoc{Definition of commutativity}

\begin{DefD}\label{def:comm} Let $\Sigma$ be an algebraic monad.
\begin{itemize}
\item[a)] We say that two operations $t\in\Sigma(n)$ and $t'\in\Sigma(m)$
{\em commute on a $\Sigma$-module~$X$} if for any collection of elements 
$\{x_{ij}\}$, $1\leq i\leq n$, $1\leq j\leq m$ of~$X$ we have 
$t(x_{1.},\ldots,x_{n.})=t'(x_{.1},\ldots,x_{.m})$, where we put 
$x_{i.}=t'(x_{i1},\ldots,x_{im})$ and $x_{.j}=t(x_{1j},\ldots,x_{nj})$.
\item[b)] We say that two operations $t$, $t'\in\|\Sigma\|$ 
{\em commute\/} if they commute on any $\Sigma$-module~$X$; clearly, 
it suffices to require them to commute on the set of free generators of 
$\Sigma(\stn\times\stm)\cong\Sigma(mn)$.
\item[c)] For any subset $S\subset\|\Sigma\|$ we denote by 
$S'\subset\|\Sigma\|$ the {\em commutant of~$\Sigma$}, i.e.\ the set of 
all operations of~$\Sigma$, commuting with all operations from~$S$. 
We'll show in a moment that $S'$ is the underlying set of a certain submonad 
of $\Sigma$, which will be also denoted by $S'$. Finally, if 
$\Lambda\subset\Sigma$ is a submonad, we denote by $\Lambda'$ the 
commutant of its underlying set $\|\Lambda\|\subset\|\Sigma\|$.
\item[d)] We say that $t\in\|\Sigma\|$ is a {\em central operation\/} 
of~$\Sigma$ if it commutes with all operations of~$\Sigma$. 
All such operations form a submonad $\Sigma'\subset\Sigma$, called the 
{\em center\/} of~$\Sigma$.
\item[e)] We say that $\Sigma$ is {\em commutative}, if any two its operations 
commute, i.e.\ if $\Sigma$ coincides with its center $\Sigma'$.
\end{itemize}
\end{DefD}

\nxsubpoint
Let's show that the commutant $S'$ of any subset $S\subset\|\Sigma\|$ is 
indeed a submonad of~$\Sigma$. Since $S'=\bigcap_{t\in S}\{t\}'$, we see 
that it is sufficient to show this for a one-element set $S=\{t\}$, with 
$t\in\Sigma(n)$.
Given any two $\Sigma$-modules $M$ and $N$, and any map of sets 
$f:M\to N$, we say that an operation $t'\in\Sigma(m)$ {\em commutes with~$f$}, 
if ${[t']}_N\circ f^m=f\circ{[t']}_M:M^m\to N$. It follows immediately from 
the definitions that $t'$ commutes with $t$ iff it commutes with all maps 
${[t]}_{\Sigma(k)}:\Sigma(k)^n\to\Sigma(k)$, $k\geq0$; actually, $k=mn$ 
would suffice for $t'\in\Sigma(m)$. Thus we are reduced to proving the 
following statement:

\nxsubpoint\label{sp:stab.of.subset}
{\em The operations of~$\Sigma$, commuting with a given map $f:M\to N$ 
between two $\Sigma$-modules $M$ and $N$, constitute a submonad 
$\Sigma_0\subset\Sigma$. The map $f$ is a $\Sigma_0$-homomorphism, and 
$\Sigma_0$ is the largest submonad of $\Sigma$ with this property.}
First of all, given any set $X$ and any subset $Y\subset X$, we 
can consider the submonad $(Y:Y)$ of $\END(X)$, given by 
$(Y:Y)(n)=\{f:X^n\to X\,|\,f(Y^n)\subset Y\}$ (cf.~\ptref{sp:def.absend.ring}),
i.e.\ $(Y:Y)$ consists of those maps $f:X^n\to X$, which stabilize~$Y$. 
It is immediate that such maps are stable under composition, hence  
$(Y:Y)$ is indeed a submonad of~$\END(X)$.

Now suppose that $X$ is a $\Sigma$-module, i.e.\ we have a monad homomorphism 
$\rho:\Sigma\to\END(X)$. Consider the pullback $\Sigma_0:=\rho^{-1}(Y:Y)$; 
clearly, this is the submonad $\Sigma_0$ of~$\Sigma$, consisting of 
all operations $t\in\Sigma(n)$ that stabilize $Y$ (i.e.\ such that 
${[t]}_X(Y^n)\subset Y$), hence $\Sigma_0$ is the largest submonad 
of~$\Sigma$, such that $Y$ is a $\Sigma_0$-submodule of~$X$.

Finally, let's consider any map $f:M\to N$ of $\Sigma$-modules as above. 
Its graph $\Gamma_f\subset M\times N$ is a subset of the $\Sigma$-module 
$M\times N$, and we see that an operation $t\in\Sigma(n)$ commutes with~$f$, 
i.e.\ $f\circ{[t]}_M={[t]}_N\circ f^n$, iff $t$ stabilizes the graph 
$\Gamma_f\subset M\times N$, i.e.\ ${[t]}_{M\times N}(\Gamma_f^n)\subset
\Gamma_f$, and we have just seen that all such operations constitute a 
submonad $\Sigma_0\subset\Sigma$.

\nxsubpoint\label{sp:prop.commut}
Notice that our notion of commutant of sets $S\subset\|\Sigma\|$ 
enjoys the usual properties, e.g.\ $S\subset T$ implies $T'\subset S'$, 
and we have $S\subset S''$ and $S'=S'''$ for any~$S$. We have seen that 
$S'$ is always a submonad of~$\Sigma$. This implies that for any subset 
$S\subset\|\Sigma\|$ its commutant $S'$ coincides with the commutant 
$\langle S\rangle'$ of the algebraic submonad $\langle S\rangle\subset\Sigma$ 
generated by~$S$ (cf.~\ptref{sp:gen.submon}). Indeed, 
$S\subset\langle S\rangle$ (we deliberately confuse submonads of~$\Sigma$ 
with their underlying sets), hence $\langle S\rangle'\subset S'$. Conversely, 
$S\subset S''$ implies $\langle S\rangle\subset S''$, $S''$ being a submonad, 
hence $S'=S'''\subset\langle S\rangle'$. In particular, if $S$ is a set of 
generators of~$\Sigma$, then $S'$ coincides with the center of~$\Sigma$. 
Therefore, $\Sigma$ is commutative iff $S'=\Sigma=\langle S\rangle$, i.e.\ 
iff $S\subset S'$, i.e.\ iff any two operations from $S$ commute.

\nxsubpoint
Given any two operations $t\in\Sigma(n)$, $t'\in\Sigma(m)$, we denote by 
$[t,t']\in\Sigma(mn)^2$ the {\em commutator\/} of $t$ and $t'$, i.e.\ 
the relation between the free generators of
$\Sigma(\stm\times\stn)\cong\Sigma(mn)$, described in~\ptref{def:comm}. 
In this way $[t,t']$ holds in~$\Sigma$ iff $t$ and $t'$ commute. Notice that 
the actual element $[t,t']\in\Sigma(mn)^2$ depends on the choice of bijection 
$\stm\times\stn\cong\stm\stn$, and on the order of $t$ and $t'$, but 
the compatible equivalence relation $\langle[t,t']\rangle$ generated by 
$[t,t']$ (cf.~\ptref{sp:eqrel.gen.by.eqs}) doesn't depend on these choices, 
since $\Sigma/\langle[t,t']\rangle$ is the largest strict quotient of~$\Sigma$,
such that the images of $t$ and $t'$ commute in this quotient. Similarly, 
for any two subsets $S$, $T\subset\|\Sigma\|$ we denote by $[S,T]$ the 
set of all relations $[s,t]$ for $s\in S$, $t\in T$; sometimes we  
denote the corresponding compatible equivalence relation $\langle[S,T]\rangle$ 
on~$\Sigma$ simply by $[S,T]$, especially when we form quotients like 
$\Sigma/[S,T]$. Clearly, this is the largest strict quotient of~$\Sigma$, 
on which all operations from the image of~$S$ commute with those from the 
image of~$T$. In particular, $\Sigma_{ab}:=\Sigma/[\Sigma,\Sigma]$ is the 
largest commutative quotient of~$\Sigma$, i.e.\ $\Sigma\mapsto\Sigma_{ab}$ 
is the left adjoint to the inclusion functor from the category of 
commutative algebraic monads into the category of all algebraic monads. 
Notice that for any two subsets $S$, $T\subset\|\Sigma\|$ we have 
$\bigl<[\langle S\rangle,\langle T\rangle]\bigr>=\langle[S,T]\rangle$; this 
is shown by the same reasoning as in~\ptref{sp:prop.commut}.

\nxsubpoint\label{sp:catgenr} (Category of generalized rings.)
We denote by $\catGenR$ the {\em category of generalized rings}, i.e.\ 
the full subcategory of the category $\catMonads_{alg}(\catSets)$ of algebraic 
monads, consisting of {\em commutative\/} algebraic monads. 
Clearly, $\catGenR$ 
is stable under subobjects, strict quotients, and arbitrary projective limits 
of~$\catMonads_{alg}(\catSets)$; in particular, arbitrary projective limits 
exist in~$\catGenR$. As to inductive limits, their existence in~$\catGenR$ 
follows from their existence in~$\catMonads_{alg}(\catSets)$ 
(cf.~\ptref{sp:indlim.algmon}), and from the existence of a left adjoint 
$\Sigma\to\Sigma_{ab}$ to the inclusion functor $\catGenR\to\catMonads_{alg}
(\catSets)$. Indeed, we have just to compute $\injlim_\alpha\Sigma_\alpha$ 
in the category of algebraic monads, and then take 
$(\injlim_\alpha\Sigma_\alpha)_{ab}$. For example, pushouts exist 
in~$\catGenR$: given two homomorphisms of generalized rings 
$\rho_i:\Sigma\to\Sigma_i$, we construct their pushout, i.e.\ the 
``tensor product of $\Sigma$-algebras'' $\Sigma_1\otimes_\Sigma\Sigma_2$, 
by computing first the non-commutative tensor product 
$\Sigma_1\boxtimes_\Sigma\Sigma_2$ (cf.~\ptref{sp:nc.tensprod}), 
and then taking quotient by the compatible equivalence relation, 
generated by all commutators between pairs of operations
of this non-commutative tensor product. Since $\Sigma_1\boxtimes_\Sigma
\Sigma_2$ is generated by $\lambda_1(\Sigma_1)\cup\lambda_2(\Sigma_2)$, 
where $\lambda_i:\Sigma_i\to\Sigma_1\boxtimes_\Sigma\Sigma_2$ is the canonical 
``inclusion'' homomorphism, we see that this equivalence relation is generated 
by commutators between elements of this sort. Now the elements of 
$\lambda_1(\Sigma_1)$ will automatically commute between themselves in 
$\Sigma_1\boxtimes_\Sigma\Sigma_2$, and similarly for $\lambda_2(\Sigma_2)$; 
hence $\Sigma_1\otimes_\Sigma\Sigma_2=(\Sigma_1\boxtimes_\Sigma\Sigma_2)/
[\lambda_1(\Sigma_1),\lambda_2(\Sigma_2)]$. 

Notice that $\catGenR$ is stable under filtered inductive limits 
of the category of algebraic monads $\catMonads_{alg}(\catSets)$, 
and such inductive limits can be computed componentwise: 
$(\injlim\Sigma_\alpha)(n)=\injlim\Sigma_\alpha(n)$.

\nxsubpoint\label{sp:comm.lower.arity} 
(Commutativity for operations of lower arities.) 
Now we want to study the shape of commutativity relations between operations 
of lower arities, and their immediate consequences.

a) Commutativity between two constants $c$, $c'\in\Sigma(0)$ means simply that 
$c=c'$, i.e.\ {\em a commutative monad has at most one constant}, hence 
is either a monad without constants or a monad with zero. So we can denote 
all constants of commutative monads by the same symbol, say $0\in\Sigma(0)$.

b) Commutativity between a constant $0\in\Sigma(0)$ and a unary operation 
$u\in\Sigma(1)$ means that $u(0)=0$, i.e.\ the zero constant of a commutative 
monad (if it exists) is necessarily fixed by all unary operations of~$\Sigma$.

c) Similarly, commutativity between a constant $0\in\Sigma(0)$ and 
an $n$-ary operation $t\in\Sigma(n)$ means $t(0,0,\ldots,0)=0$, 
so the zero constant of a commutative monad has to be respected by 
operations of all arities.

d) Commutativity between two unary operations $u$, $u'\in|\Sigma|=\Sigma(1)$ 
means that $uu'\{1\}=u'u\{1\}$, i.e.\ $uu'=u'u$ inside the monoid 
$|\Sigma|$ (cf.~\ptref{sp:underl.monoid} and~\ptref{sp:underl.set.monad}), 
so $|\Sigma|$ has to be a commutative monoid for any commutative algebraic 
monad~$\Sigma$.

e) Commutativity between a unary operation $u$ and an $n$-ary operation 
$t$ means $t(u\{1\},\ldots,u\{n\})=ut(\{1\},\ldots,\{n\})$. In particular, 
commutativity between unary $u$ and binary $+$ actually means some sort of 
distributivity: $u(\{1\}+\{2\})=u\{1\}+u\{2\}$.

f) Finally, the commutativity between two binary operations $+$ and $*$ means 
$(\{1\}+\{2\})*(\{3\}+\{4\})=(\{1\}*\{3\})+(\{2\}*\{4\})$, i.e.\ 
$(x+y)*(z+t)=(x*z)+(y*t)$ for free variables $x$, $y$, $z$, $t$ 
(or for any four elements of any $\Sigma$-module~$X$). Notice that any 
operation of arity $\leq1$ automatically commutes with itself, while 
this is not true for operations of higher arities: a binary operation $+$ 
commutes with itself iff $(x+y)+(z+t)=(x+z)+(y+t)$.

\nxsubpoint\label{sp:comm.and.add} (Commutativity and addition.)
Now we want to study the relationship between commutativity and 
additivity notions of~\ptref{p:hypadd}. Notice that these notions make sense 
only if we fix some constant $0\in\Sigma(0)$; however, commutative monads 
have at most one constant, so we don't have to specify it explicitly.

a) Any two pseudoadditions $+$ and $*$ of a commutative monad~$\Sigma$ 
coincide. Indeed, we have $0+x=x=x+0$ and $0*x=x=x*0$ by definition of a 
pseudoaddition, and $(x+y)*(z+t)=(x*z)+(y*t)$ because of the commutativity 
between $+$ and $*$. We obtain $x+y=(x*0)+(0*y)=(x+0)*(0+y)=x*y$. 
A similar proof is valid for pseudoadditions $t$ and $t'$ of same arity 
$n\geq3$: we just have to apply the commutativity relation between $t$ 
and~$t'$ to the set of elements $x_{ij}=0$ for $i\neq j$, 
$x_{ii}=\{i\}$ in $\Sigma(n)$. This implies that any pseudoaddition in 
a commutative monad is unique, i.e.\ an addition, and it is automatically 
commutative and associative: we just apply uniqueness to $+$ and $\{2\}+\{1\}$ 
for the first property, and to $(\{1\}+\{2\})+\{3\}$ and $\{1\}+(\{2\}+\{3\})$ 
for the second; alternatively, we might put $z=0$ in $(x+y)+(z+t)=(x+z)+(y+t)$.
Actually, all these statements are true for a central pseudoaddition $+$ 
in an arbitrary algebraic monad~$\Sigma$.

b) Now we show that {\em any hyperadditive commutative monad~$\Sigma$ 
is in fact additive}; in other words, {\em if a commutative monad 
admits a pseudoaddition, it is necessarily additive.} 
Indeed, any hyperadditive monad~$\Sigma$ admits 
a pseudoaddition $+$ (cf.~\ptref{sp:hyperadd.psadd}), and commutativity 
implies that $+$ is an associative and commutative addition. 
Now let $t\in\Sigma(n)$ be an arbitrary operation of~$\Sigma$, and 
let $\lambda=(\lambda_1,\ldots,\lambda_n)$ be the image of~$t$ under the 
comparison map $\pi_n:\Sigma(n)\to|\Sigma|^n$. Let us denote by $\tilde+$ 
the $n$-ary pseudoaddition (necessarily unique), constructed from $+$. 
We apply the commutativity relation between $t$ and $\tilde+$ to the 
collection $x_{ij}=0$ if $i\neq j$, $x_{ii}=\{i\}$ of elements of~$\Sigma(n)$.
We have $x_{i.}=[\tilde+](0,\ldots,\{i\},\ldots,0)=0+\cdots+\{i\}+\cdots+0=
\{i\}$, $x_{.j}=t(0,\ldots,\{j\},\ldots,0)=\lambda_j\{j\}$ by the definition 
of comparison map~$\pi_n$, hence 
$t=t(\{1\},\ldots,\{n\})=t(x_{1.},\ldots,x_{n.})=[\tilde+](x_{.1},\ldots,
x_{.n})=\lambda_1\{1\}+\cdots+\lambda_n\{n\}$, i.e.\ $t$ is uniquely 
determined by its image under~$\pi_n$, hence $\pi_n$ is bijective and 
$\Sigma$ is additive. Notice that the same proof actually shows that 
{\em if an algebraic monad admits a central pseudoaddition, 
it is automatically additive.}

c) We know that any additive monad~$\Sigma$ is of form $\Sigma_R$ for 
a uniquely determined semiring~$R$; in fact, $R=|\Sigma|$ 
(cf.~\ptref{sp:addmon.semirings}). We claim that {\em an additive 
monad~$\Sigma$ is commutative iff the semiring $|\Sigma|$ is commutative}. 
Indeed, the necessity of this condition has been shown 
in~\ptref{sp:comm.lower.arity},d); its sufficience follows from the 
fact that $\Sigma$ is generated by $0^{[0]}$, $+^{[2]}$ and all its 
unary operations~$|\Sigma|$ (cf.~\ptref{sp:addmon.semirings}), and the 
relations between these generators imply all necessary commutativity 
relations between them, with the only exception of commutativity among 
unary operations (cf.~\ptref{sp:comm.lower.arity} 
and~\ptref{sp:addmon.semirings}). Therefore, the commutativity of $|\Sigma|$ 
implies that of~$\Sigma$.

d) Similarly, classical commutative rings~$R$ correspond to commutative 
additive monads~$\Sigma$ that admit a {\em symmetry~$[-]$}, i.e.\ 
a unary operation $-\in\Sigma(1)$, such that $\bu+(-\bu)=0$ 
(this automatically implies $-^2=\bu$; cf.~\ptref{sp:addmon.rings}). 
Notice that if these properties hold in~$\Sigma$ (i.e.\ $\Sigma$ admits 
a pseudoaddition and a symmetry), and if we have any homomorphism 
of commutative monads $\rho:\Sigma\to\Xi$, they hold in~$\Xi$ as well, 
$\rho(+)$ and $\rho(-)$ being a pseudoaddition and a symmetry in~$\Xi$, 
hence if~$\Sigma$ is given by a classical ring, the same is true for~$\Xi$. 
Of course, we obtain similar results for semirings as well.

e) In general a commutative algebraic monad need not be additive. 
For example, $\Zinfty\subset\bbR$, as well as $\Fempty\subset\Fone\subset
\Fpm\subset\Zninfty\subset\bbQ$ are examples of hypoadditive commutative 
monads that are not additive (they are commutative, being submonads of 
commutative monad~$\bbR$). In fact, a commutative algebraic monad needn't 
be even hypoadditive, as shown by~$\Finfty$ (cf.~\ptref{sp:examp.finfty}; 
notice that $\Finfty$ is commutative, being a strict quotient of~$\Zinfty$).

\nxsubpoint\label{sp:alg.over.genr} (Algebras over a generalized ring.)
Let $\Lambda$ be a generalized ring, i.e.\ a commutative algebraic monad. 
A {\em (commutative) $\Lambda$-algebra~$\Sigma$} 
is by definition a homomorphism of generalized rings $\rho:\Lambda\to\Sigma$. 
We denote the category of (commutative) $\Lambda$-algebras by 
$\catCommAlgebr\Lambda$. Clearly, arbitrary projective and inductive limits 
exist in this category; for example, coproducts of $\Lambda$-algebras 
are given by pushouts of~$\catGenR$, i.e.\ the ``tensor products'' 
$\Sigma_1\otimes_\Lambda\Sigma_2$ (cf.~\ptref{sp:catgenr}).

Similarly, we define a {\em (not necessarily commutative) $\Lambda$-algebra 
$\Sigma$} to be a {\em central\/} homomorphism of algebraic monads 
$\rho:\Lambda\to\Sigma$, i.e.\ we require $\rho(\Lambda)$ to lie in the 
center of~$\Sigma$. Clearly, if~$\Lambda$ is not commutative, any such 
$\rho$ factorizes through $\Lambda_{ab}$, so it doesn't really make much 
sense to allow~$\Lambda$ to be non-commutative. Therefore, we shall consider 
algebras only over commutative monads, i.e.\ generalized rings~$\Lambda$. 
We denote the category of (all) $\Lambda$-algebras by~$\catAlgebr\Lambda$. 
Notice that this category has arbitrary projective and inductive limits: 
all projective and filtered inductive limits are computed in the usual 
way (i.e.\ componentwise), $\id:\Lambda\to\Lambda$ is the initial object, 
and the pushouts are given by the non-commutative tensor products 
$\Sigma_1\boxtimes_\Sigma\Sigma_2$ of~\ptref{sp:nc.tensprod}. We can 
also consider the commutative tensor products $\Sigma_1\otimes_\Sigma\Sigma_2
:=(\Sigma_1\boxtimes_\Sigma\Sigma_2)/[\Sigma_1,\Sigma_2]$; they also lie in
$\catAlgebr\Lambda$, being strict quotients of corresponding non-commutative 
tensor products, and they have a universal property with respect to all 
$\Lambda$-algebras~$\Xi$ that fit into a commutative square with given
$\Sigma\to\Sigma_1$ and $\Sigma\to\Sigma_2$, and such that the images of 
$\Sigma_1$ and $\Sigma_2$ in~$\Xi$ commute.

\nxsubpoint\label{sp:alg.over.addring} 
(Algebras over an additive generalized ring.)
If $\Lambda$ is additive, i.e.\ given by a commutative semiring $|\Lambda|$, 
then any $\Lambda$-algebra $\Sigma$ is additive as well: indeed, 
the structural homomorphism $\rho:\Lambda\to\Sigma$ has to be central, 
hence $\rho(+)$ is a central pseudoaddition of~$\Sigma$, and we have seen 
in~\ptref{sp:comm.and.add},b) that this implies additivity of~$\Sigma$, 
so $\Sigma$ is given by a semiring $|\Sigma|$ and a central homomorphism 
$|\rho|:|\Lambda|\to|\Sigma|$. Similarly, if $\Lambda$ is additive and admits 
a symmetry $-$, i.e.\ if $\Lambda$ corresponds to a classical commutative 
ring $|\Lambda|$, then any $\Lambda$-algebra $\Sigma$ admits a central 
pseudoaddition and a symmetry as well, hence it is given by a classical ring 
$|\Sigma|$ and a central homomorphism $|\rho|:|\Lambda|\to|\Sigma|$, i.e.\ 
a classical $|\Lambda|$-algebra. This explains our terminology. 

Here is another result of the same kind: if $\Lambda$ is a generalized ring 
with zero, then any $\Lambda$-algebra $\Sigma$ is also a monad with zero. 
Indeed, $\Sigma$ has at least one constant, namely, $\rho(0)$; and this 
constant is central, so any other constant must be equal to it 
(cf.~\ptref{sp:comm.lower.arity},a)).

\nxsubpoint (Scalar restriction and base change.)
If $\sigma:\Lambda\to\Lambda'$ is a homomorphism of generalized rings, 
we can define the {\em scalar restriction} of a $\Lambda'$-algebra 
$\rho':\Lambda'\to\Sigma'$ by considering $\rho'\circ\sigma:\Lambda\to\Sigma'$.
We denote this $\Lambda$-algebra by ${}^a\sigma_*(\Sigma')$, 
$\sigma_*\Sigma'$ or ${}_{\sigma}\Sigma'$. This scalar restriction functor 
has an obvious left adjoint --- the {\em base change functor\/~$\sigma^*$}, 
given by the commutative tensor product: $\sigma^*\Sigma:=
\Lambda'\otimes_\Lambda\Sigma$. These considerations are valid both for the 
categories of all (non-commutative) algebras and for the categories of 
commutative algebras.

\nxsubpoint\label{sp:descr.centr.el} (Central elements.)
We see directly from~\ptref{def:comm} that $t\in\Sigma(n)$ is a central element
iff the maps ${[t]}_X:X^n\to X$ are $\Sigma$-homomorphisms for all 
$\Sigma$-modules~$X$ (or just for all $X=\Sigma(m)$). Another description:
the center of $\Sigma$ is the largest submonad $\Sigma_0\subset\Sigma$, 
such that for all $t\in\Sigma(n)$ and all $\Sigma$-modules~$X$ 
(or just for all $X=\Sigma(m)$) the maps ${[t]}_X:X^n\to X$ are 
$\Sigma_0$-homomorphisms. In particular, $\Sigma$ is commutative iff 
all maps ${[t]}_X:X^n\to X$ are $\Sigma$-homomorphisms, and 
$\rho:\Lambda\to\Sigma$ is central iff ${[t]}_X:X^n\to X$ is a 
$\Lambda$-homomorphism for any $t\in\Sigma(n)$ and any $\Sigma$-module~$X$ 
(again, $X=\Sigma(m)$ suffice).

We can generalize this a bit: the actions of two submonads 
$\Sigma_1$, $\Sigma_2\subset\Xi$ commute on a $\Xi$-module~$X$ 
(i.e.\ any two operations $t\in\|\Sigma_1\|$, $t'\in\|\Sigma_2\|$ commute 
on~$X$) iff for any $t\in\Sigma_1(n)$ the map ${[t]}_X:X^n\to X$ is a
$\Sigma_2$-homomorphism. If this is true for any $X$ (or for any $X=\Xi(m)$),
then $\Sigma_1$ and $\Sigma_2$ commute inside~$\Xi$.

We can apply this to $\Xi=\Sigma_1\boxtimes_\Lambda\Sigma_2$, where 
$\Sigma_1$ and $\Sigma_2$ are $\Lambda$-algebras. We see that a 
$\Xi$-module~$X$ is just a set~$X$ with a $\Sigma_1$ and a 
$\Sigma_2$-structure, restricting to same $\Lambda$-structure 
(this actually follows from the universal property of $\Xi$, applied to 
$\Xi\to\END(X)$), and a $\Sigma_1\otimes_\Lambda\Sigma_2$-module~$X$ 
is a set~$X$ with two {\em commuting} $\Sigma_1$ and $\Sigma_2$-structures, 
restricting to the same $\Lambda$-structure, i.e.\ we require 
$\Sigma_1$ to act by $\Sigma_2$-homomorphisms on~$X$, or equivalently, 
$\Sigma_2$ to act by $\Sigma_1$-homomorphisms.

\nxsubpoint\label{sp:gen.pres.alg} (Generators and relations.)
Let's fix a generalized ring $\Lambda$. Then for any graded set 
$S=\bigsqcup_{n\geq0}S_n$ we can construct the free (non-commutative) 
$\Lambda$-algebra $\Lambda\{S\}$ (resp.\ the free commutative 
$\Lambda$-algebra $\Lambda[S]$), enjoying the usual universal property with 
respect to non-commutative (resp.\ commutative) $\Lambda$-algebras~$\Sigma$ 
and graded maps $S\to\|\Sigma\|$. Indeed, we just have to put 
$\Lambda\{S\}:=\Lambda\langle S\,|\,[S,\Lambda]\rangle=\Lambda\langle S\rangle
/\langle[S,\Lambda]\rangle$ (cf.~\ptref{sp:pres.algmon}), resp.\ 
$\Lambda[S]:=\Lambda\langle S\,|\,[S,\Lambda],[S,S]\rangle$. If all operations 
from~$S$ are unary, then $\Lambda\{S\}$ is ``the ring of polynomials in 
non-commuting variables from~$S$'' (that are nevertheless assumed to commute 
with constants from~$\Lambda$), and $\Lambda[S]$ is the ring of polynomials 
in commuting variables from~$S$.

Suppose that we are also given a system of ``equations'' or ``relations'' 
$E\subset\|\Lambda\{S\}\|^2$ (resp.\ $E\subset\|\Lambda[S]\|^2$); since 
$\Lambda\langle S\rangle\to\Lambda\{S\}$ (resp.\ $\cdots\to\Lambda[S]$) 
is surjective, we can assume that $E$ comes from a system of relations 
$\tilde E\subset\|\Lambda\langle S\rangle\|^2$, which will be usually denoted 
by the same letter~$E$. Then we can construct a (non-commutative) algebra 
$\Lambda\{S|E\}=\Lambda\{S\}/\langle E\rangle=\Lambda\langle S\,|\,
\tilde E$, $[\Lambda,S]\rangle$, resp.\ a commutative algebra 
$\Lambda[S|E]=\Lambda[S]/\langle E\rangle=\Lambda\langle S\,|\,\tilde E$, 
$[\Lambda,S]$, $[S,S]\rangle=\Lambda\{S\,|\,E,[S,S]\}$, having the usual 
universal property. When a non-commutative (resp.\ commutative) 
$\Lambda$-algebra~$\Sigma$ is represented in form $\Lambda\{S|E\}$ 
(resp.\ $\Lambda[S|E]$), we say that {\em $(S,E)$ is a presentation 
of the non-commutative (resp.\ commutative) $\Lambda$-algebra~$\Sigma$}.

\nxsubpoint\label{sp:fingen.preun.alg} 
(Finitely generated and pre-unary algebras.)
We say that a subset $S\subset\|\Sigma\|$ {\em generates the 
non-commutative (resp.\ commutative) $\Lambda$-algebra~$\Sigma$} 
if the induced homomorphism of $\Lambda$-algebras $\Lambda\{S\}\to\Sigma$ 
(resp.\ $\Lambda[S]\to\Sigma$) is surjective, i.e.\ a strict epimorphism. 
Since $\Lambda\langle S\rangle\to\Lambda\{S\}\to\Lambda[S]$ are also strict 
epimorphisms, we see that this condition is equivalent to the surjectivity 
of~$\Lambda\langle S\rangle\to\Sigma$, i.e.\ to~$S$ being a set of 
generators of~$\Sigma$ over~$\Lambda$ in the sense 
of~\ptref{sp:rel.gen.algmon}. In particular, $S$ is a set of generators of 
a commutative $\Lambda$-algebra~$\Sigma$ iff it generates $\Sigma$ as 
a non-commutative algebra. This means that we can safely transfer all 
terminology from~\ptref{sp:rel.gen.algmon} and~\ptref{sp:fingen.un.algmon}; 
we obtain the notions of finitely generated $\Lambda$-algebras 
(or $\Lambda$-algebras of finite type), pre-unary $\Lambda$-algebras, 
and $\Lambda$-algebras, generated in arity $\leq n$, both in the commutative 
and non-commutative case.

\nxsubpoint\label{sp:finpres.un.alg} (Finitely presented and unary algebras.)
We can go further on and define the finitely presented and unary 
$\Lambda$-algebras to be those (commutative or non-commutative) 
$\Lambda$-algebras that admit a finite presentation $(S,E)$ 
(resp.\ admit a unary presentation $(S,E)$, i.e.\ a presentation with 
all operations from~$S$ and all relations from~$E$ unary). 
Relation $\Lambda[S|E]=\Lambda\{S\,|\,E$, $[S,S]\rangle$ 
shows that a commutative $\Lambda$-algebra is finitely 
presented iff it is finitely presented as a non-commutative algebra, and 
similarly for unarity (note that $[u,u']$ is a unary relation when both 
$u$ and $u'$ are unary). However, these notions of finite presentation and 
unarity do not coincide with those of~\ptref{sp:homalgmon.finpres.un}; 
henceforth {\em finite presentation and unarity of algebras over generalized 
rings will be understood as just defined.} Nevertheless, they have essentially 
the same properties as listed before in~\ptref{sp:nc.finpres.un.prop}, 
when they are properly understood. For example, in the commutative case 
we have to replace the NC-tensor product $\boxtimes$ with its commutative 
counterpart $\otimes$, and in the non-commutative case we compute base 
change of a $\Lambda$-algebra $\Sigma$ only with respect to homomorphisms 
of generalized rings $\Lambda\to\Lambda'$; such base change functors preserve 
finite generation and presentation, unarity and pre-unarity, due to the 
usual property $\Lambda'\otimes_\Lambda\Lambda\{S|E\}=\Lambda'\{S|E\}$ 
and its commutative version. Another consequence is that if a 
$\Lambda'$-algebra~$\Sigma'$ 
is finitely presented (resp.\ unary) over $\Lambda$, and 
$\Lambda'$ is finitely generated (resp.\ pre-unary) over $\Lambda$, then 
$\Sigma'$ is finitely presented (resp.\ unary) over~$\Lambda'$.

\nxsubpoint\label{sp:examp.pres.comalg} (Examples.)
Presentations of $\Lambda$-algebras (commutative or not) defined 
in~\ptref{sp:gen.pres.alg} tend to be shorter than those 
of~\ptref{sp:relpres.algmon}, since some commutativity relations need not 
be written explicitly. Let's consider some examples.

a) We have $\Fone=\Fempty[0^{[0]}]=\Fempty\{0^{[0]}\}$ and 
$\Fpm=\Fone[-^{[1]}\,|\,-^2=\bu]=\Fone\{-^{[1]}\,|\,-^2=\bu\}$, and in 
general $\bbF_{1^n}=\Fone[\zeta_n^{[1]}\,|\,\zeta_n^n=\bu]$ 
(notice that we don't have to require $\zeta_n0=0$ explicitly, since 
it is a consequence of~$[\zeta_n,0]$). 

b) For any $n$, $m\geq1$ we have a 
canonical embedding $\bbF_{1^n}\to\bbF_{1^{mn}}$, given by 
$\zeta_n\mapsto\zeta_{mn}^m$; taking the (filtered) inductive limit along 
all such maps we obtain a new generalized ring $\bbF_{1^\infty}=\injlim
\bbF_{1^n}$. Clearly, $\bbF_{1^\infty}=\Fone[\zeta_1,\zeta_2,\ldots\,|$
$\zeta_1=\bu,$ $\zeta_n=\zeta_{nm}^m]$.

c) Similarly, for any $n$, $m\geq1$ we have a canonical surjection 
$\bbF_{1^{nm}}\to\bbF_{1^n}$, given by $\zeta_{nm}\mapsto\zeta_n$. We 
can compute the corresponding projective limit 
$\bbF_{1^{-\infty}}=\projlim\bbF_{1^n}$, and we have a special unary operation 
$\zeta\in\bbF_{1^{-\infty}}(1)$, defined by the projective system of the 
$\zeta_n$. However, $\zeta$ doesn't generate this generalized ring 
over~$\bbF_1$, and in fact $\bbF_1[\zeta]\subset\bbF_{1^{-\infty}}$ is 
isomorphic to the algebra of polynomials $\bbF_1[T^{[1]}]$.

d) Now we can write down shorter presentations for $\bbZ_{\geq0}$ and 
$\bbZ$: we have $\bbZ_{\geq0}=\Fone[+^{[2]}\,|\,0+\bu=\bu=\bu+0]=
\Fempty[0^{[0]},+^{[2]}\,|\,0+\bu=\bu=\bu+0]$, and 
$\bbZ=\bbZ_{\geq0}[-^{[1]}\,|\,(-\bu)+\bu=0]$.

e) We also have $\Finfty=\Fpm\bigl[*^{[2]}\,|\,\bu*(-\bu)=0$, 
$\bu*\bu=\bu$, $x*y=y*x$, $(x*y)*z=x*(y*z)\bigr]$; the remaining 
relations of~\ptref{sp:finfty.gens0} follow from commutativity.

\nxsubpoint (Some tensor products.)
Since $\Lambda[S|E]\otimes_\Lambda\Lambda[S'|E']=\Lambda[S,S'|E,E']$, we 
can use our results on commutative algebras and their generators to compute 
some tensor products, i.e.\ some pushouts in~$\catGenR$. For example, 
$\Fone=\Fempty[0^{[0]}]$, hence $\Fone\otimes_{\Fempty}\Fone\cong
\Fempty[0^{[0]},{0'}^{[0]}]\cong\Fempty[0^{[0]}]=\Fone$, since 
the commutativity relation $[0,0']$ implies $0=0'$ 
by~\ptref{sp:comm.lower.arity},a). Since the coproduct $\Fone\otimes\Fone$ 
is isomorphic to $\Fone$ itself, we see that $\Fempty\to\Fone$ is an 
epimorphism in~$\catGenR$, but of course not an epimorphism in 
the category of all algebraic monads, since $\Fone\boxtimes_{\Fempty}\Fone=
\Fempty\langle 0^{[0]},{0'}^{[0]}\rangle$, and this is definitely not 
isomorphic to $\Fone$ (in fact, the category of modules over this 
NC-tensor product is the category of sets with two marked points).

\nxsubpoint\label{sp:nc.epi} (NC-epimorphisms.)
We see that a homomorphism $\rho:\Lambda\to\Sigma$ in $\catGenR$ 
can be an epimorphism in~$\catGenR$ without being an epimorphism in the 
category of all algebraic monads. To distinguish between these two situations 
we say that a homomorphism of generalized rings $\rho:\Lambda\to\Sigma$ 
is an {\em NC-epimorphism\/} if it is 
an epimorphism in the category of all algebraic monads. Of course, any 
NC-epimorphism is automatically an epimorphism (in~$\catGenR$). Applying the 
definition to homomorphisms $\Sigma\to\END(X)$ we see that $\rho$ is 
an NC-epimorphism iff any $\Lambda$-module~$X$ admits at most one compatible 
$\Sigma$-structure (this is the behavior one expects from localization 
morphisms like $\Lambda\to\Lambda[S^{-1}]$, and in the next chapter 
we'll see that this expectation turns out to be well-founded). 
To show the sufficience of this condition we observe 
that the category of $\Sigma\boxtimes_\Lambda\Sigma$-modules is the 
category of sets, equipped with two different $\Sigma$-structures, restricting
to the same $\Lambda$-structure; if any two such $\Sigma$-structures 
necessarily coincide, we have 
$\catMod{\bigl(\Sigma\boxtimes_\Lambda\Sigma\bigr)}\cong\catMod\Sigma$, 
hence $\Sigma\boxtimes_\Lambda\Sigma\cong\Sigma$ 
(any monad is completely determined by its category of modules, together 
with the forgetful functor).

Notice that while epimorphisms in $\catGenR$ and in the category of 
algebraic monads are different, monomorphisms and strict epimorphisms coincide,
so we can freely use these notions in both contexts, as well as notions 
of submonads and strict quotients, without any risk of confusion.

\nxsubpoint (Semirings and $\bbZ_{\geq0}$-algebras.)
We have seen in~\ptref{sp:alg.over.addring} that any 
$\bbZ_{\geq0}$-algebra~$\Sigma$ admits a central pseudoaddition and a central 
constant, hence it is additive. Conversely, if $\Sigma$ is an additive 
algebraic monad, it has exactly one constant~$0$, so this constant is 
automatically central, and its addition~$+$ is also automatically central, 
a statement easily checked by applying the comparison map $\pi_{2n}$ to 
both sides of the commutativity relation $[+,t]$ for an arbitrary 
$t\in\Sigma(n)$. We know that 
$\bbZ_{\geq0}=\Fempty[0^{[0]},+^{[2]}\,|\,0+\bu=\bu=\bu+0]$, so we obtain 
a unique central homomorphism $\bbZ_{\geq0}\to\Sigma$, i.e.\ any additive 
algebraic monad~$\Sigma$ admits a unique $\bbZ_{\geq0}$-algebra structure. 
On the other hand, the list of generators and relations given 
in~\ptref{sp:addmon.semirings} shows that any such $\Sigma$ is generated 
as a (non-commutative) $\bbZ_{\geq0}$-algebra 
by the set of its unary operations $|\Sigma|$, subject to unary relations
of form $\lambda'+\lambda''=\lambda$ and $\lambda'\lambda''=\lambda$; 
all other relations of {\em loc.cit.}\ (associativity, commutativity and 
distributivity) are either consequences of similar relations in~$\bbZ_{\geq0}$,
or consequences of implied commutativity relations between $\bbZ_{\geq0}$ and 
the generators from~$|\Sigma|$ (cf.~\ptref{sp:comm.lower.arity}). 
Therefore, any such $\Sigma$ is a unary $\bbZ_{\geq0}$-algebra. 
We have just proved the following statement:

\begin{Propz} The following categories are canonically isomorphic:
a) the category of additive algebraic monads; b) the category of 
(non-commutative) $\bbZ_{\geq0}$-algebras; c) the category of 
(non-commutative) unary $\bbZ_{\geq0}$-algebras.
\end{Propz}

\nxsubpoint\label{sp:rings.z-alg} (Rings and $\bbZ$-algebras.)
Of course, the above equivalence of categories establishes an equivalence 
between the category of (unary) commutative $\bbZ_{\geq0}$-algebras and 
the category of additive generalized rings. Another consequence: 
{\em If $\Lambda$ is an additive generalized ring, any 
$\Lambda$-algebra~$\Sigma$ is automatically additive and unary 
(over~$\Lambda$).} Indeed, $\bbZ_{\geq0}\to\Lambda\to\Sigma$ defines a 
$\bbZ_{\geq0}$-algebra structure on~$\Sigma$, hence $\Sigma$ is additive and 
unary over~$\bbZ_{\geq0}$, hence it is unary over~$\Lambda$ as well, 
$\Lambda$ being pre-unary over~$\bbZ_{\geq0}$ (cf.~\ptref{sp:finpres.un.alg}). 
This is applicable to $\bbZ=\bbZ_{\geq0}[-^{[1]}\,|\,(-\bu)+\bu=0]$: we 
see that all $\bbZ$-algebras are additive and unary (over~$\bbZ$), and that 
the category of $\bbZ$-algebras consists of those $\bbZ_{\geq0}$-algebras,
i.e.\ additive monads, that admit a symmetry~$[-]$; clearly, this is just 
the category of classical semirings with symmetry, i.e.\ of classical rings. 
Similarly, any algebra~$\Sigma$ over a commutative $\bbZ$-algebra~$\Lambda$, 
corresponding to a commutative classical ring~$|\Lambda|$, is automatically 
additive and unary over~$\Lambda$, hence given by classical ring 
$|\Sigma|$ and a central homomorphism $|\Lambda|\to|\Sigma|$, i.e.\ 
a $|\Lambda|$-algebra in the classical sense. This means that when we 
work over a classical base ring~$\Lambda$, all arising algebras turn out to 
be unary, additive and classical, so we don't obtain anything new over 
such~$\Lambda$. This also explains why we cannot observe the importance 
of unarity in the classical case, since in this case it turns out to be 
always fulfilled.

\nxsubpoint\label{sp:ex.aff.grp.sch} (Some examples: affine group schemes.)
Consider, for example, the functor 
$M_n:\Sigma\mapsto M(n,n;\Sigma)=\Sigma(n)^n$ on the category of
commutative $\Fempty$-algebras, i.e.\ generalized rings. It is easily seen 
to be representable by a certain commutative $\Fempty$-algebra $A_n$, 
namely, $A_n:=\Fempty[T_1^{[n]},\ldots,T_n^{[n]}]$. Of course, 
$A_n$ is not unary over $\Fempty$ for $n\geq2$ (indeed, if $A_2$ would be 
unary, its universal property would enable us to express any central binary 
operation in terms of some unary operations, hence any generalized ring 
generated in arity $\leq2$ would be also generated in arity $\leq1$; 
however, this is not true for~$\bbZ$). But if we extend the base to~$\bbZ$, 
we obtain a commutative $\bbZ$-algebra $A_{n,\bbZ}=\bbZ\otimes A_n=
\bbZ\otimes_{\Fempty}A_n$, representing the restriction of the same functor 
to the category of commutative $\bbZ$-algebras; of course, this algebra 
is additive and unary, and in fact it equals 
$\bbZ[T_{11},T_{12},\ldots,T_{nn}]$; so a non-unary algebra becomes unary 
over~$\bbZ$.

Another example is given by $GL_n:\Sigma\mapsto GL_n(\Sigma)\subset 
M_n(\Sigma)$. It is representable by $A'_n=\Fempty\bigl[T_1^{[n]},\ldots,
T_n^{[n]},U_1^{[n]},\ldots,U_n^{[n]}\,\big|\,U_i(T_1,\ldots,T_n)=\{i\}_n=
T_i(U_1,\ldots,U_n)\bigr]$. Again, $A'_n$ is not unary over~$\Fempty$ for 
$n\geq2$, but it becomes unary over $\bbZ$; of course, $A'_{n,\bbZ}$ is 
the usual Hopf algebra of $GL_{n,\bbZ}$.

Proceeding in this way one can show that the only unary affine group schemes 
defined over~$\Fempty$ are the (split) diagonal groups $D_{\Fempty}(A)$, 
parametrized by abelian groups~$A$, e.g.\ the split tori 
$\bbG_m^n\cong D(\bbZ^n)$.

\nxsubpoint\label{sp:tens.sqr.z} (Tensor square of~$\bbZ$.)
We have seen in~\ptref{sp:rings.z-alg} that for any generalized ring~$\Sigma$ 
there is at most one homomorphism $\bbZ\to\Sigma$. By definition this 
is equivalent to saying that $\Fempty\to\bbZ$ is an epimorphism in~$\catGenR$, 
or to $\bbZ\otimes_{\Fempty}\bbZ\cong\bbZ$. This statement can be also 
shown directly starting from the presentation of $\bbZ\otimes\bbZ$ obtained 
from the presentation $\bbZ=\Fempty[0^{[0]},-^{[1]},+^{[2]}\,|\,
\bu+0=\bu=0+\bu$, $(-\bu)+\bu=0]$: we see that $\bbZ\otimes\bbZ$ is generated 
by two constants $0$ and $0'$, two binary operations $+$ and $+'$, and 
two unary operations $-$ and $-'$, and we use~\ptref{sp:comm.lower.arity} 
to show first $0=0'$, then $[+]=[+']$, and finally deduce $[-]=[-']$, hence 
$\bbZ\otimes\bbZ\cong\bbZ$. Yet another proof: use the well-known lemma 
which states that any two commuting group structures on the same set~$H$ 
coincide and are abelian (cf.\ e.g.\ SGA~3 II~3.10).

Of course, once we know that $\Fempty\to\bbZ$ 
is an epimorphism, we deduce that $\Fone\to\bbZ$ and $\Fpm\to\bbZ$ are 
epimorphisms as well, and we can compose them with strict epimorphisms 
$\bbZ\to\bbZ/n\bbZ$, thus obtaining statements like $\bbZ/n\bbZ\otimes_{\Fone}
\bbZ/n\bbZ\cong\bbZ/n\bbZ$. This is a phenomenon of commutativity: one easily 
checks that a countable set admits different $\bbF_p$-vector space structures, 
even if we fix in advance the zero and the symmetry, hence 
$\bbF_p\boxtimes_{\Fpm}\bbF_p\not\cong\bbF_p$.

Another important consequence is that for any two $\bbZ$-algebras 
$\Sigma$ and $\Sigma'$ (commutative or not) we have $\Sigma\otimes_\bbZ\Sigma'
\cong\Sigma\otimes_{\Fempty}\Sigma'$, so we can freely use the usual 
notation $\Sigma\otimes\Sigma'$ for such tensor products.

\nxpointtoc{Topos case}
Now we want to indicate a way to transfer all our previous results to 
the topos case. Everything here turns out to be quite similar to the 
classical case.

\nxsubpoint\label{sp:comm.diagrams} (Commutativity diagrams.)
Notice that the commutativily of an algebraic monad~$\Sigma$ is equivalent 
to the commutativity of the following diagrams for all $k$, $n$, $m\geq0$:
\begin{equation}\label{eq:comm.comm.diagr}
\xymatrix@C+1pc{
\Sigma(n)\times\Sigma(m)\times\Sigma(k)^{\stn\times\stm}\ar[r]\ar[d]&
\Sigma(n)\times\Sigma(k)^n\ar[d]^{\mu_k^{(n)}}\\
\Sigma(m)\times\Sigma(k)^m\ar[r]^{\mu_k^{(m)}}&\Sigma(k)
}
\end{equation}
Of course, the two remaining arrows are also given by products of 
$\mu_k^{(n)}$, $\mu_k^{(m)}$, so the commutativity of this diagram makes 
sense also if we are given an ``algebraic monad'' $\Sigma$ over an 
arbitrary cartesian category~$\cC$ in the sense of~\ptref{sp:algmon.anycat}. 
We thus obtain the notion of a ``commutative algebraic monad'' or 
a ``generalized ring'' over such~$\cC$, and this commutativity property 
is stable under any left exact functors $h:\cC\to\cC'$, i.e.\ if~$\Sigma$  
is commutative, the same is true for $h\Sigma$. We define the notion 
of a central homomorphism $\rho:\Lambda\to\Sigma$ by means of similar diagrams 
(comparing two morphisms from $\Lambda(n)\times\Sigma(m)\times\Sigma(k)^
{\stm\times\stn}$ into $\Sigma(k)$), so we obtain the notion of a  
$\Lambda$-algebra (commutative or not) as well. 
These notions have essentially the same properties as before: 
we can either repeat the proofs or embed $\cC$ into the category 
of presheaves $\hat\cC$ as a full subcategory, and prove all statements for 
$\hat\cC$ componentwise.

\nxsubpoint (Topos case.) 
Of course, we immediately obtain the notion of a commutative algebraic monad 
or a generalized ring in a topos~$\cE$, i.e.\ the notion of a sheaf of 
generalized rings on a site~$\cS$, if we consider $\cE=\tilde\cS$. 
Since the pullback, direct image and global section functors are left exact, 
they preserve commutativity of algebraic monads. Conversely, if 
$\Sigma$ is an algebraic monad over~$\cE$, such that $\Sigma_{[U]}=
\Gamma(U,\Sigma)$ is commutative for all $U$ from a system of generating 
objects of~$\cE$, then~$\Sigma$ is commutative. If we are given a 
homomorphism $\rho:\Lambda\to\Lambda'$ of generalized rings in~$\cE$, 
we obtain a scalar restriction functor $\rho_*:\cE^{\Lambda'}\to\cE^\Lambda$; 
as usual, it admits a left adjoint, the base change functor $\rho^*$, 
denoted also by $\Lambda'\otimes_\Lambda-$; one can deduce its existence 
from the presheaf case by applying the ``sheafification'' functor 
$a:\hat\cS\to\tilde\cS$.

\nxsubpoint\label{sp:gen.ringsp.top} (Generalized ringed spaces and topoi.)
This means that we can define a {\em generalized ringed topos 
$\cE=(\cE,\Lambda)$} as a topos~$\cE$, together with a generalized ring 
$\Lambda$ in~$\cE$ (the ``structural sheaf'' of~$\cE$). A {\em morphism\/} 
$f:(\cE',\Lambda')\to(\cE,\Lambda)$ is defined in the usual way: 
$f=(\phi,\theta)$, with $\phi:\cE'\to\cE$ a morphism of topoi, and 
$\theta:\Lambda\to\phi_*\Lambda'$ a homomorphism of generalized rings. 
We obtain a pair of adjoint functors $f^*:\cE^\Lambda\to{\cE'}^{\Lambda'}$ and 
$f_*:{\cE'}^{\Lambda'}\to\cE^\Lambda$ in the usual way: $f_*$ by composing 
$\phi_*^{\Lambda'}:{\cE'}^{\Lambda'}\to\cE^{\phi_*\Lambda'}$ with 
scalar restriction with respect to $\theta$, and $f^*$ by composing 
$\phi^{*,\Lambda}:\cE^\Lambda\to{\cE'}^{\phi^*\Lambda}$ with the 
base change functor $\Lambda'\otimes_{\phi^*\Lambda}-$ with respect to 
$\theta^\sharp:\phi^*\Lambda\to\Lambda'$. Sometimes we denote the 
underlying morphism of topoi $\phi$ by the same letter~$f$; then we 
write $f^{-1}$ for the pullback functor $\phi^*$ of sheaves of sets, and 
$f^*$ for the functor just defined.

In this way we obtain the 2-category of generalized ringed topoi. Of course, 
we can consider topoi defined by some sites, thus obtaining the notions 
of generalized ringed sites and their morphisms. In particular, we can 
consider topoi defined by topological spaces; then we obtain the 
category of {\em generalized ringed spaces}.

We postpone until the next chapter the definitions of generalized locally 
ringed spaces and topoi, since they depend considerably on the choice
of {\em theory of spectra} (or {\em localization theory}): 
$\cE=(\cE,A)$ is locally ringed iff for any $X\in\Ob\cE$ 
the morphism of ringed topoi $(\cE_{/X},A_{/X})\to(\catSets,A(X))$ factorizes 
through $\Spec A(X)$. 

\nxsubpoint (Commutativity in terms of plain inner structures.)
Recall that an algebraic monad~$\Sigma$ over a topos~$\cE$ admits also a 
description in terms of certain plain inner monads 
(cf.~\ptref{sp:topos.same.algmon}), so we get a monad 
$\Sigma:\cE\to\cE$ over~$\cE$, and a canonical plain inner structure~$\alpha$,
$\alpha_{X,Y}:X\times\Sigma(Y)\to\Sigma(X\times Y)$. We can construct 
a diagonal inner structure~$\rho$ on endofunctor~$\Sigma$, 
$\rho_{X,Y}:\Sigma(X)\times\Sigma(Y)\to\Sigma(X\times Y)$, and then ask 
whether this $\rho_{X,Y}$ is commutative (cf.~\ptref{sp:comm.diag.str}); 
one knows that in this case~$\rho$ is even compatible with the monad 
structure of~$\Sigma$ {\em (loc.~cit.)} It turns out that 
{\em $\rho$ is commutative iff~$\Sigma$ is commutative in the sense 
of~\ptref{sp:comm.diagrams}}. One checks this directly for $\cE=\catSets$, 
using the explicit description of the plain inner structure maps
$\alpha_{X,Y}:X\times\Sigma(Y)\to\Sigma(X\times Y)$ and the construction 
of~$\rho$ from~$\alpha$, and then generalizes this result first to the 
categories of presheaves $\hat\cS$, and then by sheafification to 
arbitrary $\cE\cong\tilde\cS$.

\nxpointtoc{Modules over generalized rings}\label{p:mod.over.genrg}
Now we want to study some properties of the category~$\catMod\Lambda$ 
of modules over a generalized ring~$\Lambda$; sometimes we also need 
another generalized ring~$\Lambda'$ and a homomorphism 
$\rho:\Lambda\to\Lambda'$. Of course, all general properties obtained 
in~\ptref{p:algmon.mod} hold in the commutative case as well, e.g.\ 
$\catMod\Lambda$ has arbitrary inductive and projective limits and so on. 
We want to study properties specific to the commutative case. The most 
important among them is the existence of a canonical ACU 
$\otimes$-structure on $\catMod\Lambda$, and of inner Homs with respect 
to this tensor product. This implies a theory of bilinear and multilinear 
forms, which admit a more direct description as well, and in particular of 
symmetric forms; we postpone the study of alternating forms, since 
it requires some additional structure on~$\Lambda$. We can also use 
this $\otimes$-structure to consider the category of algebras 
(commutative or all) in $\catMod\Lambda$; this category turns out to 
be equivalent to the category of unary $\Lambda$-algebras.

\nxsubpoint\label{sp:modstr.homsigma} 
($\Sigma$-structure on $\Hom_\Sigma(N,P)$.) 
Recall that for any algebraic monad~$\Sigma$, any $\Sigma$-module~$P$ and 
any set~$N$ we have a canonical $\Sigma$-structure on the set
$\Hom(N,P)=\Hom_{\catSets}(N,P)\cong P^N$, namely, the product structure; 
in other words, operations of~$\Sigma$ are applied to functions 
$N\to P$ pointwise (cf.~\ptref{sp:module.of.maps}). Now suppose that 
$N$ is also a $\Sigma$-module; then we have a subset $\Hom_\Sigma(N,P)
\subset\Hom(N,P)\cong P^N$, in general not stable under all operations 
of~$\Sigma$. We claim the following:

\begin{Propz} An algebraic monad~$\Sigma$ is commutative iff 
$\Hom_\Sigma(N,P)$ is a $\Sigma$-submodule of $\Hom_\catSets(N,P)$ for 
any two $\Sigma$-modules~$N$ and~$P$ (or just for $N$ and $P$ free of finite 
rank, i.e.\ $\cong\Sigma(k)$). If this be the case, we denote 
by~$\iHom_\Sigma(N,P)$ the set $\Hom_\Sigma(N,P)$, endowed with the  
$\Sigma$-structure induced from $\Hom(N,P)\supset\Hom_\Sigma(N,P)$.
\end{Propz}
We deduce this statement from its more precise form:

\begin{PropD} 
An operation $t\in\Sigma(n)$ is central iff 
$\Hom_\Sigma(N,P)\subset\Hom(N,P)$ is stable under $t$ for all 
$\Sigma$-modules~$N$ and~$P$ (or just for $N$ and $P$ free of finite rank). 
This condition is true for fixed~$t$ and~$P$ and variable~$N$ iff 
$t_P:P^n\to P$ is a $\Sigma$-homomorphism.
\end{PropD}
\begin{Proof} Put $H:=\Hom(N,P)\cong P^N$, $H':=\Hom_\Sigma(N,P)\subset H$. 
Then by definition $t_H:\Hom(N,P)^n\cong\Hom(N,P^n)\to\Hom(N,P)$ is given 
by composing maps $N\to P^n$ with $t_P:P^n\to P$. Now ${H'}^n\subset H^n$ 
is identified with $\Hom_\Sigma(N,P^n)\subset\Hom(N,P^n)$. Therefore, if 
$t_P:P^n\to P$ is a $\Sigma$-homomorphism, then $t_H$ maps ${H'}^n$ into 
$H'=\Hom_\Sigma(N,P)$, the composite of two $\Sigma$-homomorphisms being 
a $\Sigma$-homomorphism. Conversely, taking $N:=P^n$ and considering
the image under~$t_H$ of $\id_N\in\Hom_\Sigma(N,P^n)\cong{H'}^n$, equal to 
$t_P$, we see that the stability of $\Hom_\Sigma(P^n,P)\subset\Hom(P^n,P)$ 
under~$t$ implies that $t_P$ is a $\Sigma$-homomorphism.

This proves the second statement of our proposition; notice that it would 
suffice to require the stability under~$t$ of $\Hom_\Sigma(N,P)\subset
\Hom(N,P)$ not for all~$N$, but only for all $N=\Sigma(k)$, $k\geq0$: 
indeed, stability for such $N$'s implies stability for arbitrary filtered 
inductive limits of these, in particular for all $\Sigma(S)$; now 
we can represent any~$N$ as a cokernel of a pair of morphisms 
$\Sigma(R)\rightrightarrows\Sigma(S)$ (cf.~\ptref{l:mod.coker.free}), and 
a simple diagram chasing shows that the stability for~$\Sigma(R)$ and 
$\Sigma(S)$ implies the stability for~$N$.

Now the first statement follows from the second, in its stronger form 
given above, and from the following observation, already mentioned before 
on several occasions (cf.~\ptref{sp:descr.centr.el}): 
{\em $t$ is central iff all $t_P:P^n\to P$ are $\Sigma$-homomorphisms 
iff all $t_{\Sigma(k)}:\Sigma(k)^n\to\Sigma(k)$ are $\Sigma$-homomorphisms}.
\end{Proof}

\begin{DefD} (Bilinear forms.) Given three $\Lambda$-modules $M$, $N$ and~$P$ 
over a generalized ring~$\Lambda$ and a map $\Phi:M\times N\to P$, 
we say that {\em $\Phi$ is $\Lambda$-bilinear\/} if for any 
$x\in M$ the map $s_\Phi(x):N\to P$, $y\mapsto\Phi(x,y)$, is a 
$\Lambda$-homomorphism, and for any $y\in N$ the map $d_\Phi(y):M\to P$, 
$x\mapsto\Phi(x,y)$ is a $\Lambda$-homomorphism as well. We denote 
by~$\Bilin_\Lambda(M,N;P)$ the set of all $\Lambda$-bilinear maps 
$M\times N\to P$.
\end{DefD}

\nxsubpoint
Notice that $\Bilin_\Lambda(M,N;P)\subset\Hom(M\times N,P)\cong P^{M\times N}$,
and the latter set is canonically isomorphic to $\Hom(M,\Hom(N,P))$, and 
this isomorphism is clearly $\Lambda$-linear (i.e.\ a $\Lambda$-homomorphism). 
We see that 
this bijection transforms a $\Phi\in\Hom(M\times N,P)$ into 
$s_\Phi\in\Hom(M,\Hom(N,P))$. Now the first condition in the 
definition of bilinearity actually requires $s_\Phi\in
\Hom(M,\Hom_\Lambda(N,P))\subset\Hom(M,\Hom(N,P))$; on the other hand, we 
have a canonical $\Lambda$-structure on $\Hom(N,P)\cong P^N$, and the 
second condition for the bilinearity of~$\Phi$ is easily seen to 
be equivalent to $s_\Phi\in\Hom_\Lambda(M,\Hom(N,P))$. Now we know that 
$\Hom_\Lambda(N,P)$ is a $\Lambda$-submodule of $\Hom(N,P)$, so we get 
$\Bilin_\Lambda(M,N;P)\cong\Hom(M,\Hom_\Lambda(N,P))\cap
\Hom_\Lambda(M,\Hom(N,P))=\Hom_\Lambda(M,\iHom_\Lambda(N,P))$, as expected. 

The latter set admits a canonical $\Lambda$-structure, so we obtain 
a $\Lambda$-module $\iBilin_\Lambda(M,N;P)$, canonically isomorphic 
to~$\iHom_\Lambda(M,\iHom_\Lambda(N,P))$. This $\Lambda$-structure 
is actually induced by that of $P^{M\times N}\cong\Hom(M\times N,P)\supset
\Bilin_\Lambda(M,N;P)$.

\nxsubpoint\label{sp:tensprod.genrg.mod} (Tensor products.)
We have just proved that $\iHom_\Lambda(N,P)$ represents 
$\Bilin_\Lambda(-,N;P)$; what about the representability of 
$\Bilin_\Lambda(M,N;-)$? 

\begin{Propz} For any $\Lambda$-modules $M$ and $N$ the functor\/ 
$\Bilin_\Lambda(M,N;-)$ is representable by a $\Lambda$-module 
$M\otimes_\Lambda N$, called the {\em tensor product of~$M$ and~$N$}, 
and a $\Lambda$-bilinear map $\otimes:M\times N\to M\otimes N$.
\end{Propz}
\begin{Proof} We construct $M\otimes_\Lambda N$ in the usual way. 
Consider first the free $\Lambda$-module $L:=\Lambda(M\times N)$; let's 
denote its basis elements $\{(x,y)\}$, $x\in M$, $y\in N$, simply 
by~$\{x,y\}$. Now consider the set of relations $E\subset L\times L$, 
given by the union of the following two families of relations: 
\begin{multline}
t\bigl(\{x_1,y\},\ldots,\{x_n,y\}\bigr)\equiv
\{t(x_1,\ldots,x_n),y\}\\
\text{for any $t\in\Lambda(n)$, $x_1$, \dots, $x_n\in M$ and $y\in N$;}
\end{multline}
\begin{multline}
t\bigl(\{x,y_1\},\ldots,\{x,y_n\}\bigr)\equiv
\{x,t(y_1,\ldots,y_n)\}\\
\text{for any $t\in\Lambda(n)$, $x\in M$, and $y_1$, \dots, $y_n\in N$.}
\end{multline}
We consider the compatible equivalence relation $\equiv$ on~$L$, generated 
by this set of relations, and put $M\otimes_\Lambda N:=L/\text{$\equiv$}$. We have 
a canonical map $M\times N\to\Lambda(M\times N)=L\to M\otimes_\Lambda N$, 
and its bilinearity and universal property follow immediately from the 
construction.
\end{Proof}

\nxsubpoint (Multilinear maps and closed ACU $\otimes$-structure.)
Of course, we can define $\Lambda$-multilinear maps 
$M_1\times\cdots\times M_n\to N$ in a similar fashion, by requiring them 
to become $\Lambda$-linear (i.e.\ $\Lambda$-homomorphisms) after we fix  
arbitrary values for all arguments but one. The 
set of multilinear maps from the product of fixed $\Lambda$-modules 
$M_i$ into a variable $\Lambda$-module~$N$ is easily seen to be representable 
by a {\em multiple tensor product\/} 
$M_1\otimes_\Lambda M_2\otimes_\Lambda\cdots\otimes_\Lambda M_n$.

Clearly, $(M_1\otimes_\Lambda M_2)\otimes_\Lambda M_3\cong
M_1\otimes_\Lambda M_2\otimes_\Lambda M_3\cong 
M_1\otimes_\Lambda(M_2\otimes_\Lambda M_3)$, so we get 
associativity isomorphisms 
for this tensor product on~$\catMod\Lambda$, easily seen to fit into 
the pentagon diagram, since all vertices of such a diagram are canonically 
isomorphic to a quadruple tensor product. We have a canonical bijection 
$\Bilin_\Lambda(M_1,M_2;N)\simto\Bilin_\Lambda(M_2,M_1;N)$, hence 
canonical isomorphisms $M_1\otimes_\Lambda M_2\cong M_2\otimes_\Lambda M_1$, 
fitting into the hexagon diagram.

We have inner Hom's as well: $\Hom_\Lambda(M\otimes_\Lambda N,P)\cong 
\Bilin_\Lambda(M,N;P)\cong\Hom_\Lambda(M,\iHom_\Lambda(N,P))$. The 
free $\Lambda$-module $|\Lambda|=\Lambda(1)$ is easily seen to be the 
unit for this tensor product: indeed, $\Hom_\Lambda(\Lambda(1),P)\cong
\Hom_\catSets(\st1,P)\cong P$, and this bijection is easily seen to preserve 
the $\Lambda$-structures, so we have $\iHom_\Lambda(\Lambda(1),P)\cong P$ 
in $\catMod\Lambda$, and $\Hom_\Lambda(M\otimes_\Lambda\Lambda(1),P)\cong 
\Hom_\Lambda(M,$ $\iHom_\Lambda(\Lambda(1),P))\cong\Hom_\Lambda(M,P)$ 
functorially in~$P$, hence $M\otimes_\Lambda\Lambda(1)\cong M$ by Yoneda.

We have shown that {\em $\catMod\Lambda$ admits a closed ACU 
$\otimes$-structure.} When $\Lambda$ is given by a classical ring $|\Lambda|$, 
this structure coincides with that given by classical multilinear algebra.
Notice for this that it is sufficient to check bilinearity for 
a set of generators of~$\Lambda$, e.g.\ $0\in\Lambda(0)$, $[+]\in\Lambda(2)$ 
and $|\Lambda|$; this shows equivalence of our definition of bilinearity 
with the classical one, and isomorphicity of our tensor products and 
inner Homs to their classical counterparts follows. 
On the other hand, for $\Lambda=\Zinfty$ our new definitions 
are compatible with those given in~\ptref{def:zinflatbilin}  
and~\ptref{p:extensprod} for $\ZinfLat$, and in~\ptref{p:zinfflat.polylin} 
and~\ptref{sp:acu.zinfflat} for $\ZinfFlat$, both full subcategories 
of~$\catMod\Zinfty$. This gives a motivation for defining tensor product 
in~$\catMod\Lambda$ the way we have just done it.

\nxsubpoint (Basic properties of tensor products.)
We have some simple properties of our tensor product and inner Hom on 
$\catMod\Lambda$; most of them follow formally from the fact that 
$\catMod\Lambda$ is a closed ACU $\otimes$-category, with arbitrary inductive 
and projective limits. For example, $M\otimes_\Lambda N$ commutes with 
arbitrary inductive limits in each variable (since it is commutative 
and has a right adjoint~$\iHom_\Lambda$), and $\iHom_\Lambda(N,P)$ commutes 
with inductive limits in~$N$ and projective limits in~$P$. We also 
have $\Lambda(S)\otimes_\Lambda M\cong M^{(S)}$ for any set~$S$ and any 
$\Lambda$-module~$M$, where $M^{(S)}$ denotes the sum of~$S$ copies of~$M$ 
(apply $\Hom_\Lambda(-,N)$ to both sides), hence also 
$M^{(S)}\otimes_\Lambda N\cong (M\otimes_\Lambda N)^{(S)}$ (associativity), 
$M^{(S)}\otimes_\Lambda N^{(T)}\cong (M\otimes_\Lambda N)^{(S\times T)}$, 
$\Lambda(S)\otimes_\Lambda\Lambda(T)\cong\Lambda(S\times T)$, and in 
particular $\Lambda(m)\otimes_\Lambda\Lambda(n)\cong\Lambda(mn)$. 
Using adjointness between $\otimes_\Lambda$ and $\iHom_\Lambda$, we get 
$\iHom_\Lambda(N^{(S)},P^T)\cong\iHom_\Lambda(N,P)^{S\times T}$; 
in particular, $\iHom_\Lambda(\Lambda(n),P)\cong P^n$.

Another easy consequence is that $\Bilin_\Lambda(\Lambda(m),\Lambda(n);P)
\cong P^{mn}$, and in particular $\Bilin_\Lambda(\Lambda(m),\Lambda(n);
|\Lambda|)\cong|\Lambda|^{mn}$, i.e.\ bilinear forms are given by 
$m\times n$-matrices (understood in the classical sense) with entries 
in~$|\Lambda|$.

\nxsubpoint\label{sp:alg.in.Lmod} (Algebras in $\catMod\Lambda$.)
Since $\catMod\Lambda$ is an ACU $\otimes$-category, we can consider algebras 
inside this category, thus obtaining the category $\catAlg(\catMod\Lambda)$ 
of algebras in $\catMod\Lambda$, and its full subcategory 
$\catCommAlg(\catMod\Lambda)$, consisting of commutative algebras. Moreover, 
once we have an algebra $A$ in $\catMod\Lambda$, we have the category of 
(left) $A$-modules $(\catMod\Lambda)^A$ (cf.~\ptref{sp:def.alg} 
and~\ptref{sp:def.mod}). For example, when $\Lambda=\Zinfty$, we recover 
the $\Zinfty$-algebras and modules over them, defined in~\ptref{p:tf.alg.mod}. 
On the other hand, if $\Lambda$ corresponds to a classical ring~$|\Lambda|$, 
we recover the classical notion of a $|\Lambda|$-algebra and of a module 
over such an algebra.

We are going to show the following: {\em The category of algebras 
in~$\catMod\Lambda$ is equivalent to the category of unary $\Lambda$-algebras 
(cf.~\ptref{sp:alg.over.genr} and~\ptref{sp:finpres.un.alg}), and similarly 
for commutative algebras. Moreover, if a unary $\Lambda$-algebra $A$ 
corresponds to an algebra $|A|$ inside $\catMod\Lambda$ under this 
equivalence, then $(\catMod\Lambda)^{|A|}$ is canonically equivalent 
(even isomorphic) to $\catMod A$, and this equivalence is compatible with 
the underlying set functors.}

Therefore, all $\Zinfty$-algebras constructed before can be re-interpreted 
in terms of unary $\Zinfty$-algebras, i.e.\ in terms of some algebraic monads.

\nxsubpoint
Let's fix an algebra $A$ in $\catMod\Lambda$. By definition, this means that 
we are given $\mu:A\otimes_\Lambda A\to A$ and $\epsilon:|\Lambda|\to A$, 
satisfying the usual axioms. On the other hand, such $\Lambda$-homomorphisms 
$\mu$ are in one-to-one correspondence with $\Lambda$-bilinear maps 
$\tilde\mu:A\times A\to A$, and $\epsilon:|\Lambda|\to A$ corresponds to 
an element $\tilde\epsilon\in A$. We see that the algebra axioms can be 
now reformulated by requiring $A$ to be both a $\Lambda$-module and 
a monoid (commutative, if we want $A$ to be a commutative algebra), 
such that the monoid multiplication $A\times A\to A$ is $\Lambda$-bilinear. 
Similarly, an $A$-module~$M$ is just a $\Lambda$-module~$M$ with a 
$\Lambda$-bilinear monoid action~$A\times M\to M$.

For example, if $\Lambda=\bbZ$ we recover classical rings and modules over 
them, and if $\Lambda=\Fempty$ we recover monoids and sets with monoid 
actions. 

\nxsubpoint\label{sp:rho!.for.unary} (Some functors.)
Now we construct some functors. First of all, we have a forgetful functor 
$\rho_*:(\catMod\Lambda)^A\to\catMod\Lambda$, $(M,\alpha)\mapsto M$. 
It admits a left adjoint $\rho^*$, given by $N\mapsto (A\otimes_\Lambda N, 
\mu_A\otimes\id_N)$ (this is true in any AU $\otimes$-category, 
cf.~\ptref{prop:triv.sc.ext}). Moreover, existence of inner Homs allows 
us to construct a right adjoint $\rho^!$ to~$\rho_*$: we put 
$\rho^!N:=\iHom_\Lambda(A,N)$ for any $N\in\Ob\catMod\Lambda$, and 
define $\alpha_{\rho^!N}:A\otimes_\Lambda\iHom_\Lambda(A,N)\to
\iHom_\Lambda(A,N)$ by requiring its adjoint $\alpha_{\rho^!N}^\sharp:
A\otimes_\Lambda A\otimes_\Lambda\iHom_\Lambda(A,N)\to N$ 
to be the composite map of 
$\mu_A\otimes\id_{\iHom(A,N)}$ and the evaluation map 
$\ev:A\otimes_\Lambda\iHom_\Lambda(A,N)\to N$. In other words, an 
element $a\in A$ acts on $\iHom_\Lambda(A,N)$ by pre-composing maps 
$A\to N$ with the right multiplication $R_a:A\to A$. This shows that 
we've got a $\Lambda$-bilinear monoid action of~$A$ on $\rho^!N$, i.e.\ 
an $A$-module structure. Verification of  
$\Hom_A(M,\rho^!N)\cong\Hom_\Lambda(\rho_*M,N)$ can be done now essentially 
in the usual way; let us remark that $\rho_*\rho^!N=\iHom_\Lambda(A,N)\to N$ 
is given by the evaluation at the identity of~$A$, and 
$M\to\rho^!\rho_*M=\iHom_\Lambda(A,M)$ is equal to $\bar\alpha^\flat$, where 
$\bar\alpha:M\otimes_\Lambda A\simto A\otimes_\Lambda M\stackrel\alpha\to M$ 
is the $A$-action on~$M$ with permuted arguments.

Notice that arbitrary inductive and projective limits exist in 
$(\catMod\Lambda)^A$ and are essentially computed in $\catMod\Lambda$, 
scalar restriction $\rho_*$ commutes with all these limits (and in particular 
it is exact), $\rho^*$ commutes with inductive limits, and 
$\rho^!$ with projective limits.

\nxsubpoint (Forgetful functor and construction of an algebraic monad.)
Now we can compose $\rho_*:(\catMod\Lambda)^A\to\catMod\Lambda$ with the 
forgetful functor $\Gamma_\Lambda:\catMod\Lambda\to\catSets$; thus 
$\Gamma_A:=\Gamma_\Lambda\circ\rho_*:(\catMod\Lambda)^A\to\catSets$ is the 
``underlying set'' functor, and it admits a left adjoint $L_A=\rho^*\circ 
L_\Gamma$, hence it defines a monad $\Sigma=\Lambda_A:=\Gamma_A L_A$ 
over~$\catSets$, and a comparison functor $(\catMod\Lambda)^A\to\catMod\Sigma$.

We want to show that this functor is an equivalence (even an isomorphism) of 
categories, and that~$\Lambda_A$ is an algebraic monad 
(cf.~\ptref{def:alg.endof}). The last statement 
is evident, since $\Lambda_A=\Gamma_\Lambda\rho_*\rho^*L_\Gamma$, and all 
functors involved commute with filtered inductive limits (a priori this 
is not clear only for $\rho_*$, but we have just seen that it admits a right 
adjoint~$\rho^!$). Notice that this expression actually implies 
$\Lambda_A(S)=\Gamma_\Lambda(A\otimes_\Lambda\Lambda(S))\cong\Gamma_\Lambda
(A^{(S)})$, where $A^{(S)}$ denotes the direct sum in~$\catMod\Lambda$ 
of $S$ copies of~$A$. In particular, $|\Lambda_A|=A$ as a set.

\nxsubpoint\label{sp:gammaA.monadic} (Monadicity of~$\Gamma_A$.)
So we are reduced to showing the monadicity of~$\Gamma_A$. Consider for 
this the following unary $\Lambda$-algebra $\Xi$: $\Xi$ is generated 
as a non-commutative $\Lambda$-algebra by unary generators $[a]$, $a\in A$, 
subject to unary relations $[e]=\bu$, where $e\in A$ is the unity of~$A$,
$[a][a']=[aa']$ for any $a$, $a'\in A$, and 
$\bigl[t(a_1,\ldots,a_n)\bigr]=t\bigl([a_1],\ldots,[a_n]\bigr)$ 
for any $t\in\Lambda(n)$, $a_i\in A$.

We see that a $\Xi$-module structure $\Xi\to\END(X)$ on a set~$X$ 
is the same thing as a $\Lambda$-module structure on~$X$, plus a family 
of unary operations $a_X:X\to X$, one for each $a\in A$, commuting 
with all operations from~$\Lambda$, i.e.\ all $a_X$ have to be 
$\Lambda$-linear, and such that $(aa')_X=a_X\circ a'_X$, i.e.\ we get 
a monoid action $A\times X\to X$, linear in the second argument, 
and the remaining family of relations expresses its linearity with respect 
to the first argument. 

Therefore, a $\Xi$-module structure is exactly 
a $\Lambda$-module structure plus an $A$-module structure, i.e.\ 
$\catMod\Xi$ is canonically equivalent (and even isomorphic) to 
$(\catMod\Lambda)^A$, and this equivalence is compatible with the forgetful 
functors to~$\catSets$, hence $\Xi$ has to coincide with the monad~$\Lambda_A$ 
constructed before, and we have thus shown $(\catMod\Lambda)^A\to
\catMod{\Lambda_A}$ to be an isomorphism of categories.

Moreover, by construction we have a $\Lambda$-algebra structure on 
$\Lambda_A=\Xi$, i.e.\ a central homomorphism $\rho:\Lambda\to\Lambda_A$, 
and $\rho_*$ obviously coincides with the functor 
$\rho_*:(\catMod\Lambda)^A\to\catMod\Lambda$ considered before, so 
the left adjoints $\rho^*$ also coincide. We have also proved that 
$\Lambda_A$ is unary over~$\Lambda$, and that $\rho_*$ admits a 
right adjoint~$\rho^!$. Another consequence: $|\Lambda_A|=A$.

\nxsubpoint (Unary envelopes of $\Lambda$-algebras.)
Let $\Sigma$ be an arbitrary $\Lambda$-algebra. Then $|\Sigma|$ is 
a $\Sigma$-module, hence it has a natural $\Lambda$-module structure as well 
by scalar restriction. On the other hand, $|\Sigma|$ has a monoid structure, 
given by the composition of unary operations of~$\Sigma$, and we see that 
the multiplication $|\Sigma|\times|\Sigma|\to|\Sigma|$ is $\Lambda$-bilinear:
$\Lambda$-linearity in the first argument is evident, while 
$\Lambda$-linearity in the second argument uses commutativity of 
all operations from the image of $\Lambda$ with the unary operations 
of~$\Sigma$, a consequence of the centrality of $\Lambda\to\Sigma$. 
Therefore, any such $\Sigma$ defines an algebra $|\Sigma|$ in $\catMod\Lambda$,
commutative for a commutative~$\Sigma$.

Conversely, we have just seen that any algebra $A$ in $\catMod\Lambda$ defines 
a $\Lambda$-algebra $\Lambda_A$, such that $\catMod{\Lambda_A}\cong
(\catMod\Lambda)^A$, so we have a functor in the opposite direction. 
We claim that {\em $A\mapsto\Lambda_A$ is a left adjoint 
to~$\Sigma\mapsto|\Sigma|$.} 
Indeed, this follows immediately from the description of~$\Lambda_A$ 
as the $\Lambda$-algebra, generated by unary operations $[a]$, $a\in A$, 
subject to certain unary relations, cf.~\ptref{sp:gammaA.monadic}.

We also know that $|\Lambda_A|\cong A$, i.e.\ {\em the functor $A\to\Lambda_A$
is fully faithful.} We know that all $\Lambda$-algebras in its essential 
image are unary. Now we claim that {\em all unary algebras lie in the 
essential image of $A\mapsto\Lambda_A$}, hence {\em the category of 
unary $\Lambda$-algebras is equivalent to the category of algebras 
in~$\catMod\Lambda$}, as announced in~\ptref{sp:alg.in.Lmod}. Indeed, we have 
to show that the homomorphism of $\Lambda$-algebras 
$\rho:\Lambda_{|\Sigma|}\to\Sigma$ is an isomorphism for a unary~$\Sigma$; 
this follows from the fact that {\em 
if a homomorphism of unary $\Lambda$-algebras $\rho:\Sigma'\to\Sigma$ 
induces an isomorphism between their sets of unary operations, then~$\rho$ 
is itself an isomorphism}.

For a general $\Lambda$-algebra~$\Sigma$ the homomorphism 
$\Lambda_{|\Sigma|}\to\Sigma$ is not an isomorphism; we usually denote 
$\Lambda_{|\Sigma|}$ simply by~$|\Sigma|$ and say that it is the 
{\em unary envelope\/} of~$\Sigma$, since it has a universal property with 
respect to homomorphisms from unary $\Lambda$-algebras into~$\Sigma$.

\nxsubpoint\label{sp:prop.unary.env} (Properties.)
Notice that $A\mapsto\Lambda_A$ commutes with arbitrary inductive limits, 
and $\Sigma\mapsto|\Sigma|$ commutes with arbitrary projective and 
filtered inductive limits. However, the first functor does not commute 
in general with projective limits, e.g.\ products: all we can say is 
that $\Lambda_{A\times B}$ is the unary envelope of 
$\Lambda_A\times\Lambda_B$. For example, for $\Lambda=\Zinfty$, and 
$A=B=|\Zinfty|$, we see that $\Zinfty\times\Zinfty$ computed in the 
category of unary $\Zinfty$-algebras (cf.~\ptref{sp:examp.spec.zinfty2}) 
differs from 
the same product computed in the category of all $\Zinfty$-algebras. 
This explains why we've got $\Spec\Zinfty\times\Zinfty\neq
\Spec\Zinfty\sqcup\Spec\Zinfty$ in~\ptref{sp:examp.spec.zinfty2}: 
this equality would be true 
if we had computed $\Zinfty\times\Zinfty$ in the larger category. 
This also explains why unary $\Zinfty$-algebras are insufficient for our 
purpose, so we need to enlarge our category. 

\begin{PropD}\label{prop:cond.unarity}
Let $\Sigma$ be a $\Lambda$-algebra, given by a central homomorphism 
$\rho:\Lambda\to\Sigma$. The following conditions are equivalent:
\begin{itemize}
\item[$a)$] $\Sigma$ is a unary $\Lambda$-algebra;
\item[$b)$] $\Sigma$ is isomorphic to $\Lambda_A$ for some algebra~$A$ 
in~$\catMod\Lambda$;
\item[$c)$] Scalar restriction $\rho_*:\catMod\Sigma\to\catMod\Lambda$ 
admits a right adjoint~$\rho^!$;
\item[$d)$] $\rho_*$ commutes with finite direct sums;
\item[$d')$] $\rho_*$ commutes with finite direct sums of 
free $\Sigma$-modules of finite rank.
\end{itemize}
\end{PropD}
\begin{Proof}
We have already proved $a)\Leftrightarrow b)$, $b)\Rightarrow c)$ has 
been shown in~\ptref{sp:rho!.for.unary}, and 
$c)\Rightarrow d)\Rightarrow d')$ is trivial. Let's show $d')\Rightarrow b)$. 
Put $A:=|\Sigma|$, and consider the canonical homomorphism 
$\Lambda_A\to\Sigma$. Now $\Lambda_A(n)$ can be identified with $A^{(n)}$, 
i.e.\ $(\rho_*\Sigma(1))^{(n)}$, and $\Sigma(n)$ with 
$\rho_*(\Sigma(1)^{(n)})$. Then $\Lambda_A(n)\to\Sigma(n)$ is identified with 
the canonical morphism $(\rho_*\Sigma(1))^{(n)}\to\rho_*(\Sigma(1)^{(n)})$, 
and $d')$ implies that all these maps are bijective, hence 
$\Lambda_A\cong\Sigma$, q.e.d.
\end{Proof}

\nxsubpoint\label{sp:pullback.tensfunct} (Base change and tensor product.)
Let $\rho:\Lambda\to\Lambda'$ be a homomorphism of generalized rings. 
We claim that {\em the base change functor $\rho^*:\catMod\Lambda\to
\catMod{\Lambda'}$ is a $\otimes$-functor}, i.e.\ {\em we have canonical 
isomorphisms $\rho^*(M\otimes_\Lambda N)\simto\rho^*M\otimes_{\Lambda'}
\rho^*N$ for all $\Lambda$-modules~$M$ and~$N$.}

First of all, for any three $\Lambda'$-modules $M'$, $N'$, $P'$ we have 
a natural map $\Bilin_{\Lambda'}(M',N';P')\to\Bilin_\Lambda(\rho_*M',\rho_*N'; 
\rho_*P')$, given by considering any $\Lambda'$-bilinear map $M'\times N'\to 
P'$ as a $\Lambda$-bilinear map with respect to $\Lambda$-module structures 
defined by scalar restriction. Put $M':=\rho^*M$, $N':=\rho^*N$; then 
$\Bilin_{\Lambda'}(\rho^*M,\rho^*N;P')\to\Bilin_\Lambda(\rho_*\rho^*M, 
\rho_*\rho^*N;\rho_*P')\to\Bilin_\Lambda(M,N;\rho_*P')$, where the second 
arrow is defined with the aid of canonical homomorphisms $M\to\rho_*\rho^*M$ 
and $N\to\rho_*\rho^*N$.

Using the universal property of tensor products we get natural maps 
$\Hom_{\Lambda'}(\rho^*M\otimes_{\Lambda'}\rho^*N,P')\cong
\Bilin_{\Lambda'}(\rho^*M,\rho^*N;P')\to\Bilin_\Lambda(M,N;\rho_*P')\cong 
\Hom_\Lambda(M\otimes_\Lambda N,\rho_*P')\cong\Hom_{\Lambda'}
(\rho^*(M\otimes_\Lambda N),P')$; by Yoneda the composite arrow is induced 
by a canonical $\Lambda'$-homomorphism $\rho^*(M\otimes_\Lambda N)\to
\rho^*M\otimes_{\Lambda'}\rho^*N$. We want to show that it is an 
isomorphism, i.e.\ that for all $M$, $N$ and $P'$ the map 
$\Bilin_{\Lambda'}(\rho^*M,\rho^*N;P')\to 
\Bilin_{\Lambda}(M,N;\rho_*P')$ is a bijection.

Notice that the source of this map $\Bilin_{\Lambda'}(\rho^*M,\rho^*N;P')$ is 
isomorphic to $\Hom_{\Lambda'}(\rho^*M,\iHom_{\Lambda'}(\rho^*N,P'))
\cong\Hom_\Lambda(M,\rho_*\iHom_{\Lambda'}(\rho^*N,P'))$, 
and the target is isomorphic to $\Hom_\Lambda(M,\iHom_\Lambda(N,\rho_*P'))$; 
by Yoneda our map is induced by some $\Lambda$-homomorphism 
$\rho_*\iHom_{\Lambda'}(\rho^*N,P')\to\iHom_\Lambda(N,\rho_*P')$, and all 
we have to do is to check that it is an isomorphism. Now this is trivial: 
taking global sections, i.e.\ looking at underlying sets, we get 
the adjointness map $\Hom_{\Lambda'}(\rho^*N,P')\to\Hom_\Lambda(N,\rho_*P')$, 
which is of course bijective.

One can construct in a similar way canonical homomorphisms, relating 
$\rho^*$ with multiple tensor products; expressing multiple tensor products 
in terms of repeated binary tensor products we see that these homomorphisms 
are in fact isomorphisms, and their existence implies the compatibility 
of our isomorphisms $\rho^*(M\otimes_\Lambda N)\simto\rho^*M\otimes_{\Lambda'}
\rho^*N$ with the associativity and commutativity constraints, hence 
$\rho^*$ is a $\otimes$-functor.

\nxsubpoint (Base change for unary algebras.)
An immediate consequence is that $\rho^*A$ is an algebra 
in~$\catMod{\Lambda'}$ for any algebra~$A$ in~$\catMod\Lambda$. 
Conversely, if~$A'$ is an algebra in $\catMod{\Lambda'}$, 
then $\rho_*A'$ admits a canonical algebra structure: this can 
be either seen directly from the description of algebras as monoids 
with bilinear multiplication, or deduced from the existence of canonical 
homomorphisms $\rho_*M'\otimes_\Lambda\rho_*N'\to
\rho_*(M'\otimes_{\Lambda'}N')$.

We want to show that this base change is compatible with the general 
base change for $\Lambda$-algebras, i.e.\ that $\Lambda'_{\rho^*A}$ 
is canonically isomorphic to~$\Lambda'\otimes_\Lambda\Lambda_A$; 
this would allow us to denote $\Lambda_A$ by~$A$ and still be able to 
use notations like $\Lambda'\otimes_\Lambda A$ or $A_{(\Lambda')}$ without 
ambiguity. To prove this statement we observe that $\Sigma:=\Lambda_A$ is 
a unary $\Lambda$-algebra, hence $\Sigma':=\Lambda'\otimes_\Lambda\Sigma$ 
is a unary $\Lambda'$-algebra, hence it is given by some algebra $A'$ 
in~$\catMod{\Lambda'}$, namely, $A':=|\Sigma'|$; all we have to show is
$A'\cong\rho^*A$. This follows from the ``base change theorem'' 
(cf.~\ptref{th:aff.base.change} below), 
applied to $\Sigma$-module $|\Sigma|$: we see that 
doing base change to $\Sigma'$ and then restricting scalars to $\Lambda'$ 
yields the same result as first restricting scalars to $\Lambda$, and then 
applying base change~$\rho^*$; now the first path yields~$A'=|\Sigma'|$, 
and the second --- $\rho^*|\Sigma|=\rho^*A$.

\nxsubpoint (Tensor algebras.) 
Given any $\Lambda$-module~$M$, we can construct its {\em tensor algebra\/}
$T_\Lambda(M)$. We construct $T_\Lambda(M)=T(M)$ first as an algebra 
in~$\catMod\Lambda$, by putting $T^n(M):=M^{\otimes n}$ 
({\em tensor powers\/} of~$M$) and 
$T(M):=\bigoplus_{n\geq0} T^n(M)$, and defining the multiplication 
$T(M)\otimes_\Lambda T(M)\to T(M)$ with the aid of canonical maps 
$T^n(M)\otimes_\Lambda T^m(M)\to T^{n+m}(M)$, using the fact that tensor 
products commute with infinite direct sums. One shows in the usual fashion 
that $M\to T_\Lambda(M)$ is universal among all $\Lambda$-homomorphisms 
from~$M$ into algebras in~$\catMod\Lambda$.

We denote the unary $\Lambda$-algebra, defined by $T_\Lambda(M)$, by the 
same expression $T_\Lambda(M)$. Using the existence of unary envelopes 
one shows that $T_\Lambda(M)$ is universal among $\Lambda$-homomorphisms 
$M\to|\Sigma|$, with~$\Sigma$ a variable $\Lambda$-algebra.

Finally, if $\rho:\Lambda\to\Lambda'$ is a homomorphism of generalized 
rings, we have a canonical isomorphism $T_{\Lambda'}(\rho^*M)\cong
\Lambda'\otimes_\Lambda T_\Lambda(M)$, obtained for example by 
comparing the universal properties of the two sides.

\nxsubpoint (Symmetric algebras.)
Notice that the symmetric group $\gS_n$ acts on $T^n(M)$ by permuting the 
components. We can define the {\em symmetric powers\/} of~$M$ by putting 
$S_\Lambda^n(M):=T_\Lambda^n(M)/\gS_n$, i.e.\ we take the largest strict 
quotient of~$T^n(M)$, on which $\gS_n$ acts trivially (in other words, 
the kernel of $T^n(M)\to S^n(M)$ is the compatible equivalence relation, 
generated by equations $\sigma(t)\equiv t$, $t\in T^n(M)$, $\sigma\in\gS_n$).
These symmetric powers parametrize the spaces of symmetric multilinear maps 
$M^n\to P$, i.e.\ $\Polylin_\Lambda(M^n;P)^{\gS_n}\cong\Hom_\Lambda
(S_\Lambda(M),P)$.

Now we can define canonical maps $S^n(M)\otimes S^m(M)\to S^{n+m}(M)$, induced 
by corresponding maps for tensor powers by passage to strict quotients, 
and construct the {\em symmetric algebra\/} $S_\Lambda^n(M):=\bigoplus_{n\geq0}
S_\Lambda^n(M)$, considered both as a commutative algebra inside 
$\catMod\Lambda$, and as a commutative unary $\Lambda$-algebra; clearly, 
$S_\Lambda(M)$ is a strict quotient of the tensor algebra~$T_\Lambda(M)$.

We see that $S_\Lambda(M)$ has a universal property with respect to 
$\Lambda$-homo\-mor\-phisms $M\to|\Sigma|$, where $\Sigma$ runs through 
commutative $\Lambda$-algebras. This universal property implies stability 
of the construction under any base change~$\rho:\Lambda\to\Lambda'$. 
In particular, $S_{\Lambda'}^n(\rho^*M)\cong \rho^*(S_\Lambda^n(M))$.

\nxsubpoint (Monoid and group algebras.)
For any monoid~$G$ we can construct the corresponding {\em monoid algebra 
$\Lambda[G]$} (called a {\em group algebra\/} if~$G$ is a group), as follows.
We construct first the corresponding algebra in~$\catMod\Lambda$ by defining 
a multiplication on the free module $\Lambda(G)=\bigoplus_{g\in G}\Lambda\{g\}$
in the usual fashion, and then consider the corresponding unary 
$\Lambda$-algebra. It clearly has a universal property with respect to 
monoid homomorphisms $G\to|\Sigma|$, for variable $\Lambda$-algebra~$\Sigma$
(again, for a unitary $\Sigma$ this is shown in the classical way, and 
in general case we replace $\Sigma$ by its unary envelope). 
This gives another description of $\Lambda[G]$, namely, 
$\Lambda[G]=\Lambda\{[g], g\in G\,|\,[e]=\bu$, $[gh]=[g][h]\}$, and shows the 
compatibility with base change: $\Lambda'\otimes_\Lambda\Lambda[G]\cong
\Lambda'[G]$. Clearly, $\Lambda[G]$ is commutative for a commutative~$G$; 
considering homomorphisms $\Lambda[G]\to\END(X)$, we see that 
$\Lambda[G]$-modules are exactly $\Lambda$-modules, equipped with a 
$\Lambda$-linear action of~$G$.

\nxsubpoint (Unary algebras over~$\Fempty$.) 
Recall that $\catMod\Fempty$ is just the category of sets, with direct sums 
(i.e.\ coproducts) given by disjoint unions, and tensor products equal 
to direct products (since any map $M\times N\to P$ is bilinear). 
This means that unary $\Fempty$-algebras, i.e.\ algebras in~$\catMod\Fempty
\cong\catSets$, are nothing else than {\em monoids\/}~$G$ (commutative algebras
corresponding to commutative monoids). Clearly, in this case the unary 
$\Fempty$-algebra defined by~$G$ coincides with the monoid algebra 
$\Fempty[G]$, i.e.\ all unary algebras over~$\Fempty$ are monoid algebras. 
Of course, for any $\Lambda$ we have $\Lambda[G]=\Lambda\otimes_{\Fempty}
\Fempty[G]=\Lambda\otimes G$. For example, we have 
$(\Fempty[G])(n)=\bigsqcup_{i\in\stn} G=G\times\stn$, with $(g,i)\in G\times 
\stn$ corresponding to $g\{i\}$. Computing base change to~$\Lambda$, we 
obtain $(\Lambda[G])(n)\cong\Lambda(G\times\stn)$.

\nxsubpoint (Unary polynomial algebras.)
We can use the above results to describe the unary polynomial algebras 
$\Sigma=\Lambda[(T_i^{[1]})_{i\in I}]$. Indeed, an inspection of 
universal properties shows that it is isomorphic both to the symmetric algebra 
of the free $\Lambda$-module $\Lambda(I)$, hence it admits a natural 
grading, and to the monoid algebra $\Lambda[\bbZ_{\geq0}^{(I)}]$ of the 
free commutative monoid $\bbZ_{\geq0}^{(I)}$, generated by~$I$. For example, 
the set of unary operations of $\Fempty[T_1^{[1]},\ldots,T_n^{[1]}]$ 
is canonically isomorphic to $\bbZ_{\geq0}^n$, i.e.\ it is the set of 
monomials in the~$T_i$; the set of $m$-ary operations is given by the 
product $\bbZ_{\geq0}^n\times\stm$, i.e.\ it consists of all expressions 
$T_1^{k_1}\cdots T_n^{k_n}\{i\}$, $k_j\geq0$, $1\leq i\leq m$. 

Notice that the grading on the polynomial algebra, coming from its 
symmetric algebra description, defines an increasing filtration, so it 
makes sense to speak about the degree of a polynomial with respect to all 
or to some group of variables. However, in general we cannot extract from 
a polynomial its leading term, or any other coefficient, except the free term.

This isomorphism between polynomial algebras and symmetric algebras of 
free modules shows that symmetric powers of free modules are free, and 
(in the finite-dimensional case) of correct rank: $S_\Lambda^k(\Lambda(n))$ 
is free of rank~$\binom{n+k-1}{k}$. 

\nxsubpoint\label{sp:un.poly.alg.nc} 
(Unary polynomial algebras in non-commuting variables.)
Similarly, we see by comparing the universal properties that 
the free algebra $\Lambda\{S\}$, when all free generators from~$S$ are 
supposed to be unary, is canonically isomorphic to the 
monoid algebra of~$W(S)$ (the free monoid over~$S$, i.e.\ the set of 
monomials in non-commuting variables from~$S$) over~$\Lambda$. We deduce 
$(\Lambda\{S\})(n)\cong\Lambda(W(S)\times\stn)$, and in particular 
$|\Lambda\{S\}|\cong\Lambda(W(S))$.

\nxsubpoint (Commutativity and algebraic bimodules.) 
Finally, let us mention an application of commutativity to the category of 
algebraic $(\Lambda,\Lambda)$-bimodules. Namely, we have seen 
in~\ptref{sp:alg.bimod} that any such bimodule~$M$ is given by a 
functor $\tilde M:\catN_\Lambda\to\catMod\Lambda$, where $\catN_\Lambda$ 
is the category of standard free $\Lambda$-modules~$\Lambda(n)$ of finite rank.
Therefore, for any $n$, $m\geq0$ we have a map of sets 
$\Hom_\Lambda(\Lambda(n),\Lambda(m))\to\Hom_\Lambda(\tilde M(n),\tilde M(n))$. 
When $\Lambda$ is commutative, we obtain $\Lambda$-structures on these two 
sets, so it makes sense to consider the case when all these maps of sets 
are $\Lambda$-linear. When $\tilde M$ has this property, we say that 
$M$ or $\tilde M$ is a {\em central\/} algebraic bimodule. This definition 
generalizes to the categories of algebraic $(\Sigma,\Xi)$-bimodules, with 
$\Sigma$ and $\Xi$ two $\Lambda$-algebras; then we speak about 
{\em $\Lambda$-central\/} algebraic bimodules. Of course, these notions 
admit some matrix descriptions as well.

\nxsubpoint (Comparison to~\cite{ShaiHaran}.)
Since $\Lambda(n)\otimes_\Lambda\Lambda(m)\cong\Lambda(nm)$ and 
$\Lambda(n)\oplus\Lambda(m)\cong\Lambda(n+m)$ for any generalized 
ring~$\Lambda$, and we have defined the matrix sets $M(n,m;\Lambda)$ 
by $M(n,m;\Lambda):=\Lambda(n)^m\cong\Hom_\Lambda(\Lambda(m),\Lambda(n))$, 
functors $\oplus$ and $\otimes$ induce maps 
$\oplus:M(n,m;\Lambda)\times M(n',m';\Lambda)\to M(n+n',m+m';\Lambda)$ and 
$\otimes:M(n,m;\Lambda)\times M(n',m';\Lambda)\to M(nn',mm';\Lambda)$. 
We have the composition maps $\circ:M(n,m;\Lambda)\times M(m,k;\Lambda)\to 
M(n,k;\Lambda)$ as well, and when $\Lambda$ is a (commutative) 
$\Fone$-algebra, we have also some maps $M(n,m;\Fone)\to M(n,m;\Lambda)$. 

This collection of data, consisting of the sets $M(n,m;\Lambda)$, the 
$\circ$, $\otimes$ and $\oplus$-operations between them, and the maps  
$M(n,m;\Fone)\to M(n,m;\Lambda)$, is easily checked to satisfy all the 
conditions for an ``$\bbF$-algebra'' of~\cite{ShaiHaran} 
(up to a minor point --- Shai Haran requires all $M(n,m;\Fone)\to
M(n,m;\Lambda)$ to be injective; this excludes the trivial monad $\st1$ 
and $\st1_+\subset\st1$), and our (commutative) $\Fpm$-algebras define 
``$\bbF_\pm$-algebras'' in the sense of {\em loc.~cit.} 
Notice that the approach of Shai Haran is more general, since it 
never requires $M(n,m;\Lambda)\cong M(n,1;\Lambda)^m$; actually, our 
category of commutative $\Fone$-algebras (almost) corresponds to the 
category of ``$\bbF$-algebras'', satisfying this additional condition.

Shai Haran's approach has the obvious advantage of being more general and 
more symmetric (one can ``transpose'' all matrices and define $\Lambda^t$ 
by $M(n,m;\Lambda^t):=M(m,n;\Lambda)$ for any $\bbF$-algebra~$\Lambda$). 
However, there are considerable disadvantages of such a generality as well: 
his approach is more categorial, 
while our approach is more algebraic and has a direct 
connection to algebraic systems. This enables us to construct  
free objects, and compute inductive limits 
(both in the categories of algebras and modules) quite explicitly, and transfer
more statements from the classical case. At the same time our 
approach seems sufficiently general, at least for the present moment, 
to deal with objects arising from arithmetics and Arakelov geometry 
(e.g.\ the compactification of $\Spec\bbZ$). That's why we prefer not  
to replace our approach by the more general one.

\nxpointtoc{Flatness and unarity}
The notions of unarity and of flatness (defined below) seem to be 
in some sense ``orthogonal'', like properness and smoothness of morphisms 
with respect to \'etale cohomology.

\begin{DefD}
A homomorphism of generalized rings $\rho:\Sigma\to\Xi$ is {\em flat} 
(resp.\ {\em faithfully flat}), if the corresponding base change functor 
$\rho^*:\catMod\Sigma\to\catMod\Xi$ is (left) exact (resp.\ left exact 
and faithful). A $\Lambda$-algebra $\Sigma$ is said to be {\em flat}
(resp.\ {\em faithfully flat}) if the corresponding central homomorphism 
$\Lambda\to\Sigma$ has this property.
\end{DefD}

Clearly, $\rho$ is faithfully flat iff it is flat, and $\rho^*$ is 
conservative. If $\rho$ is already known to be flat, the latter condition
is equivalent to $\rho^*N\neq\rho^*M$ for any submodule $N\neq M$
of any $\Sigma$-module~$M$.

We have the usual properties. For example, in the situation 
$\Sigma\to\Sigma'\to\Sigma''$, we have transitivity of flatness and of 
faithful flatness, and if $\Sigma''$ is flat (resp.\ faithfully flat) 
over~$\Sigma$ and faithfully flat over~$\Sigma'$, then $\Sigma'$ is also flat 
(resp.\ faithfully flat) over~$\Sigma$. Stability of flatness under base 
change, at least for algebras over a generalized ring, follows from the 
following statement:

\begin{ThD}\label{th:aff.base.change} (``Base change theorem'')
Let $\sigma:\Lambda\to\Lambda'$ be a homomorphism of generalized rings, 
$\Sigma$ be a $\Lambda$-algebra, given by a central homomorphism 
$\rho:\Lambda\to\Sigma$. Put $\Sigma':=\Lambda'\otimes_\Lambda\Sigma$, 
and denote by $\rho':\Lambda'\to\Sigma'$ and $\tau:\Sigma\to\Sigma'$ the 
canonical homomorphisms, so as to obtain a cocartesian square:
\begin{equation}
\xymatrix@C+5pt{
\Sigma'&\ar[l]_{\rho'}\Lambda'\\
\Sigma\ar[u]_{\tau}&\Lambda\ar[l]_{\rho}\ar[u]_{\sigma}
}
\end{equation}
Then we have a canonical natural transformation 
$\gamma:\sigma^*\rho_*\to\rho'_*\tau^*$ of functors $\catMod\Sigma\to
\catMod{\Lambda'}$. 
When either a) $\rho$ is unary, or b) $\sigma$ is flat, this 
natural transformation 
$\gamma:\sigma^*\rho_*\to\rho'_*\tau^*$ is an isomorphism.
\end{ThD}
\begin{Proof}
First of all, let us construct $\gamma:\sigma^*\rho_*\to\rho'_*\tau^*$. 
By adjointness it is sufficient to construct $\gamma^\flat:\rho_*\to
\sigma_*\rho'_*\tau^*=\rho_*\tau_*\tau^*$, and to do this we just 
apply~$\rho_*$ to the adjointness natural transformation~$\Id\to\tau_*\tau^*$.

We want to prove that $\gamma$ is an isomorphism under some conditions. 
Let us fix some presentation $\Sigma=\Lambda\{S|E\}$; in case a) we choose 
a unary presentation, i.e.\ $S\subset|\Sigma|$ and $E\subset|\Lambda\{S\}|^2$. 
Now let's take some $\Sigma$-module~$M$ and put $\bar M:=\rho_*M$; we 
have some additional structure on this $\Lambda$-module~$\bar M$, namely, 
for each~$s\in S$ we have a $\Lambda$-linear map $s_M:\bar M^{r(s)}\to\bar M$ 
($\Lambda$-linearity follows from the centrality of $\Lambda\to\Sigma$), 
and appropriate composites of these maps satisfy the relations from~$E$. 
Moreover, the category of $\Sigma$-modules is in fact equivalent to the 
category of $\Lambda$-modules, equipped with these additional $\Lambda$-linear 
maps, and since $\Sigma'=\Lambda'\otimes_\Lambda\Sigma=\Lambda'\{S|E\}$, we 
have a similar description of $\Sigma'$-modules in terms of $\Lambda'$-modules 
with some extra $\Lambda'$-linear maps.

We construct this additional structure on the $\Lambda'$-module 
$\sigma^*\bar M$ by putting $s_{\sigma^*M}:=\sigma^*(s_M)\in
\Hom_{\Lambda'}(\sigma^*(\bar M^{r(s)}),\sigma^*\bar M)$. Notice that 
in case a) we have $r(s)=1$, and in case b) $\sigma^*$ commutes with finite 
products, so in both cases we have $\sigma^*(\bar M^{r(s)})\cong
(\sigma^*\bar M)^{r(s)}$, hence $s_{\sigma^*M}$ acts from 
$\sigma^*\bar M^{r(s)}$ to~$\sigma^*\bar M$ as required.

We claim that these $\Lambda'$-homomorphisms $s_{\sigma^*M}$ satisfy all the 
relations from~$E$, thus defining a $\Sigma'$-module structure on 
$\sigma^*\bar M$; resulting $\Sigma'$-module will be provisionally denoted 
by~$\tau^?M$. Indeed, this is clear in case~b), since all relations from~$E$ 
equate some morphisms between products of several copies of a $\Lambda$- or 
a $\Lambda'$-module, and $\sigma^*$ commutes with such products, $\sigma$ 
being flat, hence the validity of~$E$ for $s_{\sigma^*M}=\sigma^*(s_M)$ 
follows from the validity of~$E$ for~$s_M$; one might also say that we have 
a homomorphism of $\Lambda$-algebras $\END_\Lambda(\bar M)\to
\END_{\Lambda'}(\sigma^*\bar M)$. In case a) we have only unary relations, 
i.e.\ $E\subset|\Lambda\{S\}|^2$, and according to~\ptref{sp:un.poly.alg.nc} 
any element of~$|\Lambda\{S\}|$ can be written in form 
$t(u_1,\ldots, u_m)$, for some $t\in\Lambda(m)$, $m\geq0$, 
and $u_i\in W(S)$. Let us denote by~$s$ the image in~$|\Sigma|$ of such an 
element. All we have to show is that $s=t(u_1,\ldots,u_n)$ in~$\Sigma$ 
implies similar relation between $s_{\sigma^*M}$ and $u_{j,\sigma^*M}$. 
Now $\sigma^*:\End_\Lambda(\bar M)\to\End_{\Lambda'}(\sigma^*\bar M)$ is 
clearly a monoid homomorphism, so we can assume that $S$ is closed under 
multiplication (in~$|\Sigma|$), and that all $u_j\in S$. We know that 
$s_M=t(u_{1,M},\ldots,u_{m,M})$ in $\iEnd_\Lambda(\bar M)$; now our statement 
follows from the $\Lambda$-linearity of~$\iEnd_\Lambda(\bar M)\to
\iHom_\Lambda(\bar M,\sigma_*\sigma^*\bar M)\cong
\sigma_*\iEnd_{\Lambda'}(\sigma^*\bar M)$, 
already shown in~\ptref{sp:pullback.tensfunct}.

Now let's show that $\tau^?M$ has the universal property required 
from~$\tau^*M$. Indeed, giving a $\Sigma$-homomorphism 
$f:M\to\tau_*N$ is equivalent to giving a $\Lambda$-homomorphism 
$\bar f:\bar M\to\sigma_*\bar N$, where $\bar N:=\rho'_*N$, such that 
for all generators~$s\in S$ the following diagram commutes:
\begin{equation}
\xymatrix@C+1pc{
\bar M^{r(s)}\ar[r]^{\bar f^{r(s)}}\ar[d]^{s_M}&
\sigma_*\bar N^{r(s)}\ar[d]^{\sigma_*(s_N)}\\
\bar M\ar[r]^{\bar f}&\sigma_*\bar N}
\end{equation}
Using adjointness between $\sigma^*$ and $\sigma_*$, and the fact that either 
$r(s)=1$ or $\sigma$ is flat, we see that the commutativity of the above 
diagrams for~$\bar f$ is equivalent to the commutativity of similar diagrams 
for~$\bar f^\sharp:\sigma^*\bar M\to\bar N$:
\begin{equation}
\xymatrix@C+1pc{
\sigma^*\bar M^{r(s)}\ar[r]^{(\bar f^\sharp)^{r(s)}}\ar[d]^{\sigma^*(s_M)}&
\bar N^{r(s)}\ar[d]^{s_N}\\
\sigma^*\bar M\ar[r]^{\bar f^\sharp}&\bar N}
\end{equation}
On the other hand, giving such a $\Lambda'$-homomorphism $\bar f^\sharp:
\sigma^*\bar M\to\bar N$ is equivalent to giving a~$\Sigma'$-homomorphism 
$f^\sharp:\tau^?M\to N$, and we get the required universal property 
for~$\tau^?M$.

We see that $\tau^*M\cong\tau^?M$, hence $\rho'_*\tau^*M\cong\rho'_*\tau^?M
=\sigma^*\bar M=\sigma^*\rho_*M$, i.e.\ 
$\gamma_M:\sigma^*\rho_*M\to\rho'_*\tau^*M$ is an isomorphism, q.e.d.
\end{Proof}

\nxsubpoint 
Notice that the conclusion of the theorem is not true without any additional 
assumptions about~$\sigma$ or~$\rho$. Indeed, let's take $\Lambda=\Fempty$, 
$\Sigma=\Lambda'=\bbZ$; then $\Sigma'=\bbZ\otimes_{\Fempty}\bbZ\cong\bbZ$ 
by~\ptref{sp:tens.sqr.z}, and consider $M:=\Sigma(1)=\bbZ$. Then 
$\sigma^*\rho_*\bbZ$ equals $\bbZ[\bbZ]$, the monoid algebra over~$\bbZ$ of 
the multiplicative monoid of~$\bbZ$, i.e.\ the ring of Dirichlet polynomials 
$\sum_{n\in\bbZ}c_n\cdot n^\omega$. On the other hand, both~$\tau$ 
and~$\rho'$ are isomorphisms, hence $\rho'_*\tau^*\bbZ=\bbZ$, and 
$\gamma_\bbZ:\bbZ[\bbZ]\to\bbZ$ is the map $\sum_nc_nn^\omega\mapsto\sum_n
nc_n$. Clearly, it is not an isomorphism. This example also shows that 
$\bbZ$ is neither unary nor flat over~$\Fempty$, $\Fone$ and~$\Fpm$. 

\nxsubpoint\label{sp:bc.stab.flat}
An important consequence is that {\em if $\Sigma$ is flat over~$\Lambda$, 
then $\Sigma'=\Lambda'\otimes_\Lambda\Sigma$ is flat over~$\Lambda'$}, i.e.\ 
flatness of algebras is stable under base change. Indeed, in the 
notations of~\ptref{th:aff.base.change} we have $\rho'_*\tau^*\cong 
\sigma^*\rho_*$; since $\rho'_*$, $\rho_*$ and $\sigma^*$ are left exact, 
and $\rho'_*$ is in addition conservative, we see that $\tau^*$ is 
also left exact, hence $\tau$ is flat. One shows similarly that if in addition 
$\sigma^*$ is conservative, then the same is true for $\tau^*$, i.e.\ 
{\em faithful flatness is also stable under base change.}

\begin{DefD}
We say that a $\Lambda$-module~$M$ is flat (resp.\ faithfully flat) if 
the functor $M\otimes_\Lambda-:\catMod\Lambda\to\catMod\Lambda$ is 
(left) exact (resp.\ left exact and faithful).
\end{DefD}
This definition is motivated by the following easy observation: a unary 
$\Lambda$-algebra~$A$ is flat (resp.\ faithfully flat) iff $|A|$ is 
flat (resp.\ faithfully flat) as a~$\Lambda$-module.

Clearly, $|\Lambda|=\Lambda(1)$ is always faithfully flat, and if 
$\Lambda$ is a monad with zero, then the zero module $\Lambda(0)$ is flat. 
However, in general we cannot say that all free $\Lambda$-modules are flat, 
since the direct sum of two monomorphisms is not necessarily a monomorphism. 
Actually, flatness of free modules seems to be deeply related to 
the existence of injectives in~$\catMod\Lambda$.


\nxsubpoint (Projection formula.)
Notice that the projection formula $\rho_*(M'\otimes_{\Lambda'}\rho^*N)\cong
\rho_*M\otimes_\Lambda N$, for $M'$ a $\Lambda'$-module, 
$N$ a $\Lambda$-module, and $\rho:\Lambda\to\Lambda'$ a homomorphism 
of generalized rings, in general does {\em not\/} hold unless $\rho$ 
is unary. Indeed, the validity of projection formula for 
$M':=\Lambda'(1)$, $N:=\Lambda(n)$ means $\rho_*\Lambda'(n)\cong
\Lambda'(1)\otimes_\Lambda\Lambda(n)\cong\Lambda'(1)^{\oplus n}$, 
hence the unarity of~$\rho$ by~\ptref{prop:cond.unarity},$d')$. Conversely, if 
$\rho$ is unary, we can prove the projection formula first for~$N=\Lambda(1)$ 
(trivial), then for $N=\Lambda(S)=\Lambda(1)^{(S)}$ (using the fact that 
$\rho_*$ commutes in this case with direct sums), and then for arbitrary~$N$, 
writing it as a cokernel of two morphisms between free modules.

\nxsubpoint (Products of flat algebras.)
We have seen that in general we don't have a lot of flat modules and 
algebras in $\catMod\Lambda$: all we can say is that $\Lambda(1)$ is flat, 
and that $\Lambda(0)$ is flat when $\Lambda$ admits a zero. 
On the other hand, we can construct more (non-unary) $\Lambda$-algebras: 
we know that $\Lambda$ is flat over itself, and {\em when $\Lambda$ is a monad 
with zero, the product $\Sigma_1\times\Sigma_2$ of two flat 
$\Lambda$-algebras is itself a flat $\Lambda$-algebra.} This follows 
immediately from the following fact:

\begin{ThD}\label{th:mod.under.prod} 
If $\Sigma_1$ and $\Sigma_2$ are two algebraic monads with zero, 
then $\catMod{(\Sigma_1\times\Sigma_2)}$ is canonically equivalent to 
$\catMod{\Sigma_1}\times\catMod{\Sigma_2}$.
\end{ThD}
\begin{Proof}
Put $\Sigma:=\Sigma_1\times\Sigma_2$; by definition $\Sigma(n)=\Sigma_1(n)
\times\Sigma_2(n)$, so any operation $t\in\Sigma(n)$ can be written in form 
$(t^{(1)},t^{(2)})$, with $t^{(i)}\in\Sigma_i(n)$. Now we can construct a 
functor $F:\catMod{\Sigma_1}\times\catMod{\Sigma_2}\to\catMod\Sigma$: 
given any $\Sigma_1$-module~$M_1$ and $\Sigma_2$-module~$M_2$, we put 
$F(M_1,M_2):=M_1\times M_2$, with the~$\Sigma$-structure 
$\alpha^{(k)}:\Sigma(k)\times(M_1\times M_2)^k\to M_1\times M_2$ defined 
componentwise by means of the maps $\Sigma_i(k)\times M_i^k\to M_i$.

Let's construct a functor~$G$ in the opposite direction. Denote for this 
the only constant of~$\Sigma$ by~$0$, so we have $0=(0^{(1)},0^{(2)})$, 
where $0^{(i)}$ is the zero of~$\Sigma_i$. We also have 
$\bu=(\bu^{(1)},\bu^{(2)})$ for the identity of~$\Sigma$, and we consider 
the elements $\bu_1$ and $\bu_2$ in $|\Sigma|$, given by 
$\bu_1=(\bu^{(1)},0^{(2)})$ and $\bu_2=(0^{(1)},\bu^{(2)})$. Clearly, 
$|\Sigma|=|\Sigma_1|\times|\Sigma_2|$ as a monoid, hence $\bu_1$ and 
$\bu_2$ are central idempotents in~$|\Sigma|$, such that $\bu_1\bu_2=\bu_2\bu_1
=0$.

Now for an arbitrary $\Sigma$-module~$M$ we put $M_i:=\bu_iM\subset M$. 
Clearly, $M_i$ is a $\Sigma$-submodule of~$M$, $\bu_i$ being central 
in~$\Sigma$; $M_i$ consists of all elements~$x$ of~$M$, such that $\bu_ix=x$. 
We define a $\Sigma_i$-action on~$M_i$ in the natural way 
($t\in\Sigma_1(n)$ acts on~$M_1$ by means of ${[(t,0^{(2)})]}_{M_1}:M_1^n\to
M_1$), and put $G(M):=(M_1,M_2)$.

We have to construct some isomorphisms $GF(M_1,M_2)\cong(M_1,M_2)$ and 
$M\simto FG(M)$, defining an adjunction between $G$ and~$F$. First of them 
is obvious: we observe that $\bu_1(M_1\times M_2)=M_1\times 0\cong M_1$, 
and similarly for the second component of~$GF(M_1,M_2)$. The 
$\Sigma$-homomorphism $\theta_M:M\to\bu_1M\times\bu_2M$ is also easy to 
construct: put $\theta_M(z):=(\bu_1z,\bu_2z)$. One checks immediately that 
these natural transformations do define an adjunction between~$G$ and~$F$, 
and we have only to check that $\theta_M$ is an isomorphism.

Consider for this the element $\phi:=(\{1\}_\st2^{(1)},\{2\}_\st2^{(2)})\in
\Sigma(2)$. It is clearly central in~$\Sigma$, both its components being 
central, hence it defines a $\Sigma$-homo\-mor\-phism $\phi_M:M\times M\to M$ 
for any $\Sigma$-module~$M$. Put $M_1:=\bu_1M$, $M_2:=\bu_2M$, and 
denote by $\phi'_M$ the composition of the embedding $M_1\times M_2\to
M\times M$ with $\phi_M$. We claim that $\phi'_M$ is an inverse for~$\theta_M$.

a) $\theta_M\phi'_M=\id_{M_1\times M_2}$. Let's take $x\in M_1$, $y\in M_2$; 
then $x=\bu_1x$, $y=\bu_2y$, and we want to show $\bu_1\phi(x,y)=x$, 
$\bu_2\phi(x,y)=y$. Clearly, it suffices to check $\bu_1\phi(\bu_1\{1\},\{2\})=
\bu_1\{1\}$; we check this componentwise: $(\bu_1\phi(\bu_1\{1\},\{2\}))^{(1)}=
\bu_1^{(1)}\phi^{(1)}(\bu_1^{(1)}\{1\},\{2\})=\bu^{(1)}\bu^{(1)}\{1\}=
\{1\}=\bu_1^{(1)}\{1\}$, and $(\bu_1\phi(\bu_1\{1\},\{2\}))^{(2)}=
\bu_1^{(2)}\phi^{(2)}(\bu_1^{(2)}\{1\},\{2\})=0^{(2)}\{2\}=0^{(2)}=
\bu_1^{(2)}\{1\}$.

b) $\phi'_M\theta_M=\id_M$. We have to check that $\phi(\bu_1z,\bu_2z)=z$ 
for any $z\in M$, i.e.\ that $\phi(\bu_1,\bu_2)=\bu$ in $|\Sigma|$. 
We check this componentwise: $(\phi(\bu_1,\bu_2))^{(1)}=\phi^{(1)}
(\bu_1^{(1)},\bu_2^{(1)})=\bu_1^{(1)}=\bu^{(1)}$, and similarly for the 
other component. So $\phi'_M$ is indeed an inverse to~$\theta_M$, hence 
$F$ and $G$ are adjoint equivalences, q.e.d.
\end{Proof}

\nxsubpoint 
One might try to generalize the above statement as follows. Suppose 
$\sigma_i:\Sigma_i\to\Sigma_0$, $i=1$, $2$, are two strict epimorphisms. 
Put $\Sigma:=\Sigma_1\times_{\Sigma_0}\Sigma_2$, and denote by~$\rho_i$ 
the canonical projections $\Sigma\to\Sigma_i$. Then $\rho_1^*$ and 
$\rho_2^*$ induce a functor~$G$ from $\catMod\Sigma$ into the category 
of triples $(M_1,M_2,\phi)$, where $M_i\in\catMod{\Sigma_i}$, and 
$\phi:\sigma_1^*M_1\simto\sigma_2^*M_2$ is a $\Sigma_0$-isomorphism. 
One can construct a functor~$F$ in the opposite direction by mapping 
a triple as above into the subset $M\subset M_1\times M_2$, consisting 
of all pairs $(x,y)\in M_1\times M_2$, such that $\phi\gamma_1(x)=\gamma_2(y)$,
where $\gamma_i:M_i\to\sigma_{i,*}\sigma_i^*M_i$ is the canonical map. 
It is easy to check that~$M$ is indeed a $\Sigma$-module (with the action 
of $\Sigma\subset\Sigma_1\times\Sigma_2$ defined componentwise), and that 
$G$ and~$F$ are adjoint. It seems plausible that they are always equivalences; 
in this case we might be able to prove that $\Lambda$-flatness of 
$\Sigma_0$, $\Sigma_1$ and $\Sigma_2$ implies that of $\Sigma$.

In any case similar considerations explain why in~\ptref{th:mod.under.prod} 
we need 
to suppose that $\Sigma_1$ and $\Sigma_2$ admit some constants, i.e.\ that 
both $\Sigma_i\to\st1$ are surjective. Without this assumption the theorem 
is false, as illustrated by the case $\Sigma_1=\Sigma_2=\st1_+\subset\st1$:
then $(\emptyset,\st1)$ and $(\st1,\emptyset)$ are not isomorphic in 
$\catMod{\Sigma_1}\times\catMod{\Sigma_2}$, but their images under~$G$ 
are both equal to~$\emptyset$, hence isomorphic.

\nxsubpoint (Flatness of algebraic bimodules.)
Recall that any algebraic $(\Lambda,\Lambda)$-bimodule~$M$ 
induces a functor $\tilde M=M\otimes_\Lambda-:\catMod\Lambda\to
\catMod\Lambda$ (cf.~\ptref{sp:tensprod.bimod}), such that 
$M\otimes_\Lambda\Lambda(n)\cong M(n)$. It is natural to say that 
{\em $M$ is flat\/} if {\em this functor $M\otimes_\Lambda-$ is exact.} 
One defines similarly flatness for algebraic $(\Sigma,\Lambda)$-bimodules: 
in this case we require exactness of a functor $\catMod\Lambda\to
\catMod\Sigma$. Clearly, a $\Lambda$-algebra~$\Sigma$ is flat iff it is 
flat as a $(\Sigma,\Lambda)$-bimodule.

\nxpointtoc{Alternating monads and exterior powers}
Now we would like to study some simple property (called {\em alternativity}) 
of generalized rings that enables us to construct a reasonable theory of 
exterior powers and determinants.

To do this we fix some commutative $\Fpm$-algebra~$K$; if we need another 
such algebra, we denote it by~$K'$. Since $\Fpm=\Fempty[0^{[0]},-^{[1]}\,
|\,-^2=\bu]$, we see that this is equivalent to fixing a generalized 
ring~$K$ with zero~$0$ and a unary operation~$-$, such that $-(-\bu)=\bu$.

\nxsubpoint (Bilinear forms and matrices.)
Let $\Lambda$ be a generalized ring. We know that $\Lambda$-bilinear maps 
$\Phi:\Lambda(m)\times\Lambda(n)\to X$ are in one-to-one correspondence with 
collections $(\Phi_{ij})$ of elements of~$X$, indexed 
by~$(i,j)\in\stm\times\stn$, since 
$\Bilin_\Lambda(\Lambda(m),\Lambda(n);X)\cong\Hom_\Lambda(
\Lambda(m)\otimes_\Lambda\Lambda(n),X)\cong\Hom_\Lambda(\Lambda(mn),X)\cong
X^{mn}=X^{\stm\times\stn}$. Clearly, this correspondence is given by 
evaluating~$\Phi$ on base elements: $\Phi_{ij}=\Phi(\{i\},\{j\})$. 
Therefore, such bilinear maps are parametrized in the usual manner by 
(classical) $m\times n$-matrices with entries in~$X$.

Given any elements $t\in\Lambda(m)$ and $t'\in\Lambda(n)$, we can construct 
a new element $t\otimes t'\in\Lambda(m)\otimes_\Lambda\Lambda(n)\cong
\Lambda(\stm\times\stn)$. For any $\Lambda$-module~$X$ the map 
${[t\otimes t']}_X:X^{\stm\times\stn}\to X$ is computed in the same way
as in~\ptref{def:comm}, first applying $t'$ to the rows of a matrix 
$(x_{ij})$, and then $t$ to the elements of~$X$ thus obtained, or the 
other way around: $[t\otimes t']((x_{ij}))=[t'\otimes t]((x_{ji}))$.
Clearly, if a bilinear map~$\Phi$ is given by a matrix 
$(\Phi_{ij})$ as above, ${[t\otimes t']}_X((\Phi_{ij}))=\Phi(t,t')$

\begin{DefD}\label{def:altform}
Let $K$ be a commutative $\Fpm$-algebra.
\begin{itemize}
\item[a)] A $K$-bilinear map $\Phi:M\times M\to X$ is {\em skew-symmetric} 
if $\Phi(x,y)=-\Phi(y,x)$ for any $x$, $y\in M$.
\item[b)] A $K$-bilinear map $\Phi:M\times M\to X$ is {\em alternating} 
if it is skew-symmetric, and $\Phi(x,x)=0$ for any $x\in M$.
\item[c)] A $K$-bilinear map $\Phi:K(n)\times K(n)\to X$ is said to be 
{\em pseudo-alternating} if it is given by an {\em alternating matrix} 
$(\Phi_{ij})_{1\leq i,j\leq n}$, $\Phi_{ij}:=\Phi(\{i\},\{j\})$, i.e.\ 
if $\Phi_{ii}=0$ and $\Phi_{ij}=-\Phi_{ji}$ for all $1\leq i,j\leq n$.
\end{itemize}
\end{DefD}

\nxsubpoint\label{sp:rem.altforms}
Clearly, any alternating $K$-bilinear map $\Phi:K(n)\times K(n)\to X$ 
is pseudo-alternating; the converse is in general not true. However, 
it is easy to see that any pseudo-alternating map is skew-symmetric: indeed, 
we have $\Phi(t,t')=[t\otimes t']((\Phi_{ij}))=[t'\otimes t]((\Phi_{ji}))=
[t'\otimes t]((-\Phi_{ij}))=-[t'\otimes t]((\Phi_{ij}))=-\Phi(t',t)$ 
(notice the way we used here the commutativity of $t$, $t'$ and $-$ between 
themselves). We see that a skew-symmetric map $\Phi:K(n)\times K(n)\to X$ 
is pseudo-alternating iff $\Phi(\{i\},\{i\})=0$ for all~$1\leq i\leq n$.

Observe that there is a universal pseudo-alternating map 
$\Phi:K(n)\times K(n)\to X$ for a fixed~$n$ and variable~$X$. To obtain it we 
put $X:=K(n(n-1)/2)$, order the pairs $(i,j)$ with $1\leq i<j\leq n$ in some 
way (e.g.\ lexicographically), and put $\Phi_{ij}:=\{i,j\}$  
(the basis element corresponding to~$(i,j)$) if $i<j$, $\Phi_{ii}:=0$, and 
$\Phi_{ij}:=-\{j,i\}$ if $i>j$. We usually denote this free $K$-module with 
this base by $\bigwedge_K^2K(n)$, and the universal pseudo-alternating map 
by $\wedge:K(n)\times K(n)\to\bigwedge_K^2K(n)$; we write 
$x\wedge y$ instead of~$\wedge(x,y)$.

\begin{DefD}\label{def:altop}
Let $K$ be a commutative $\Fpm$-algebra.
\begin{itemize}
\item[a)] A $K$-linear map $f:K(m)\to K(n)$ is said to be 
{\em alternating} if for any pseudo-alternating bilinear map 
$\Phi:K(n)\times K(n)\to X$ (or just for the universal 
pseudo-alternating bilinear map $\wedge$) the bilinear map 
$\Phi\circ(f\times f):K(m)\times K(m)\to X$ is also 
pseudo-alternating.
\item[b)] An operation $t\in K(n)$ is said to be {\em alternating} 
if the corresponding map $t:K(1)\to K(n)$ is alternating, i.e.\ 
if for any pseudo-alternating map $\Phi:K(n)\times K(n)\to X$ 
(or just for the universal one) we have $\Phi(t,t)=0$, or equivalently,
if for any alternating $n\times n$-matrix $(\Phi_{ij})$ (or just for the 
universal one, given by $\Phi_{ij}=-\Phi_{ji}=\{i,j\}$ for $i<j$, 
$\Phi_{ii}=0$) we have $[t\otimes t]((\Phi_{ij}))=0$.
\item[c)] For any operation $t\in K(n)$ we denote by $\wedge t$ the 
corresponding {\em alternativity relation} for~$t$, i.e.\ the 
relation of arity $n(n-1)/2$, obtained by equating to zero the value 
of $[t\otimes t]$ on the universal alternating $n\times n$-matrix 
with entries in~$K(n(n-1)/2)$. Clearly, $\wedge t$ is fulfilled in~$K$ iff
$t$ is alternating.
\item[d)] We say that a commutative $\Fpm$-algebra is {\em alternating} 
iff all its operations are alternating. 
\item[e)] For any $n\geq0$ we denote by $K^{alt}(n)$ the subset 
of~$K(n)$, consisting of all alternating $n$-arity operations. We'll see in a 
moment that this collection of sets defines a submonad $K^{alt}\subset K$, 
clearly the largest alternating submonad of~$K$.
\item[f)] We denote by $K_{alt}$ or $K/_{alt}$ the quotient of~$K$ 
with respect to all alternativity relations~$\wedge t$, for all $t\in\|K\|$. 
This is the largest alternating strict quotient of~$K$.
\end{itemize}
\end{DefD}

\nxsubpoint\label{sp:rem.altmon}
Notice that if $K'$ is a submonad of~$K$, and $t\in K'(n)$, then alternativity 
of~$t$ with respect to~$K$ is equivalent to its alternativity with respect 
to~$K'$, since both conditions require the same element 
$\wedge t\in K'(n(n-1)/2)\subset K(n(n-1)/2)$ to be equal to zero.

Let $f:K(m)\to K(n)$ be a $K$-linear map. 
Since $\Hom_K(K(m),K(n))\cong K(n)^m$, $f$ is given by a collection 
$(f_1,\ldots,f_m)\in K(n)^m$ of $n$-ary operations. Now let 
$\Phi:K(n)\times K(n)\to X$ be pseudo-alternating, and put 
$\Psi:=\Phi\circ(f\times f)$. Clearly, $\Psi$ is skew-symmetric, so it is 
pseudo-alternating iff all $\Psi_{ii}=\Phi(f(\{i\}),f(\{i\}))=\Phi(f_i,f_i)$ 
are equal to zero, hence {\em $f$ is alternating iff all its components 
$f_i$ are alternating.} It is evident that the composite of two alternating 
maps, say, $t:K(1)\to K(m)$ and $(t_1,\ldots, t_m):K(m)\to K(n)$, is 
alternating again; we deduce that ${[t]}_{K(n)}(t_1,\ldots,t_m)$ is alternating
whenever $t\in K(m)$ and all $t_i\in K(n)$ are alternating. 

Hence $K^{alt}\subset K$ is closed under composition, 
and obviously contains all projections $\{i\}_\stn$; 
therefore, $K^{alt}$ is indeed an algebraic submonad of~$K$, and 
$K^{alt}(n)=K(n)$ for $n\leq1$, all operations of arity $\leq1$ being 
automatically alternating. Therefore, our operations $0\in K(0)$ 
and $-\in K(1)$, such that $-^2=\bu$, lie in~$K^{alt}$, hence 
$K^{alt}$ is a commutative $\Fpm$-algebra, clearly alternating.

An immediate consequence is this: {\em 
if $S\subset\|K\|$ generates $K$, then $K$ is alternating iff all 
operations from~$S$ are alternating}. Since all unary operations are 
automatically alternating, we see that {\em any pre-unary algebra over 
an alternating $\Fpm$-algebra is itself alternating}, and in particular 
{\em any commutative unary $\Fpm$-algebra is alternating.} 

\nxsubpoint\label{sp:alt.binop} (Alternativity and binary operations.)
Clearly, all constants and unary operations are alternating. Let's consider 
the case of a binary operation~$*$. We see that $\wedge*$ is the unary 
relation $(0*\bu)*(-\bu*0)=0$. For an addition~$+$ this is equivalent to
$\bu+(-\bu)=0$, i.e.\ to $-$ being a symmetry for~$+$. Since any classical 
commutative ring is generated over~$\Fpm$ by its unary operations and its 
addition, we see that any classical commutative ring is alternating.

\nxsubpoint (Exterior square of a module.)
For any commutative $\Fpm$-algebra~$K$ and any $K$-module~$M$ we can 
construct a universal alternating map $\wedge:M\times M\to\bigwedge_K^2M$, 
simply by taking the strict quotient of $M\times M$ modulo relations 
$x\otimes y\equiv-y\otimes x$ and $x\otimes x\equiv 0$. This 
$K$-module $\bigwedge_K^2M$ is called the {\em exterior square of\/~$M$}.

When $M=K(n)$ is free and $K$ is alternating, all pseudo-alternating 
bilinear maps $\Phi:K(n)\times K(n)\to X$ are automatically alternating, 
condition $\Phi(t,t)=0$ for any pseudo-alternating~$\Phi$ 
being actually equivalent to 
the alternativity of $t\in K(n)$. Therefore, in this case the universal 
alternating map on $K(n)\times K(n)$ coincides with the universal 
pseudo-alternating map of~\ptref{sp:rem.altforms}, hence 
$\bigwedge_K^2K(n)\cong K(n(n-1)/2)$, and $\{i\}\wedge\{j\}$, 
$1\leq i<j\leq n$ form a base of $\bigwedge_K^2K(n)$, as one would expect. 
Notice that in general this is not true when $K$ is not alternating.

\nxsubpoint\label{sp:ext.sqr.sum} (Exterior square of a direct sum.)
We want to show that over an alternating~$K$ we have 
\begin{equation}\label{eq:ext.sqr.sum}
\wedgenl_K^2(M_1\oplus M_2)\cong\wedgenl_K^2M_1\oplus
\left(M_1\otimes_KM_2\right)\oplus\wedgenl_K^2M_2
\end{equation}
First of all, consider bilinear maps $\Phi:M\times M'\to X$, where 
$M=M_1\oplus M_2$ and $M'=M'_1\oplus M'_2$. Since 
$\Bilin_K(M,M';X)\cong\Hom_K(M\otimes_K M',X)\cong\Hom_K(\bigoplus_{i,j}
M_i\otimes_K M'_j,X)\cong\prod_{i,j}\Bilin_K(M_i,M'_j;X)$, we see that 
giving such a $\Phi$ is equivalent to giving four bilinear maps
$\Phi_{ij}:M_i\times M'_j\to X$; of course, these $\Phi_{ij}$ are 
just restrictions of~$\Phi$ to corresponding components of direct sums.

Now consider the case $M=M'$, $M_i=M'_i$, and $\Phi:M\times M\to X$ 
skew-symmetric. Then $\Phi_1:=\Phi_{11}:M_1\times M_1\to X$ and 
$\Phi_2:=\Phi_{22}$ are clearly also skew-symmetric, and 
$\Phi_{21}$ is completely determined by $\Phi_{12}$ since 
$\Phi_{21}(y,x)=-\Phi_{12}(x,y)$. Therefore, a skew-symmetric map~$\Phi$ 
on~$M$ gives rise to two skew-symmetric maps $\Phi_i$ on~$M_i$ and 
a bilinear map $\Phi_{12}:M_1\times M_2\to X$. Conversely, any such collection 
$(\Phi_1,\Phi_2,\Phi_{12})$ allows us to reconstruct a bilinear map~$\Phi$ 
(we recover $\Phi_{21}$ from $\Phi_{12}$ as above), clearly skew-symmetric.

We would like to obtain a similar description of alternating maps on~$M$; 
this would prove~\eqref{eq:ext.sqr.sum}, since both sides would represent 
the same functor. Clearly, if $\Phi$ is alternating, then both $\Phi_1$ 
and $\Phi_2$ are alternating. Let's prove the converse, i.e.\ let's 
prove $\Phi(z,z)=0$ for a $z\in M=M_1\oplus M_2$, assuming both $\Phi_1$ 
and $\Phi_2$ to be alternating.

We know that $z$ can be written in form $z=t(x_1,\ldots, x_n,y_1,\ldots, y_m)$,
for some $n$, $m\geq0$, $t\in K(n+m)$, $x_i\in M_1$, and $y_j\in M_2$ 
(cf.~\ptref{sp:elem.dirsum}). Consider the map $f:K(n+m)\to M$, defined by 
these elements $x_i$ and $y_j$. Clearly, $\Psi:=\Phi\circ(f\times f)$ 
is pseudo-alternating: indeed, it is obviously skew-symmetric, 
and the diagonal elements of its matrix are equal either to some 
$\Phi(x_i,x_i)=\Phi_1(x_i,x_i)$ or to some $\Phi(y_j,y_j)=\Phi_2(y_j,y_j)$;
in both cases we get zero, since both $\Phi_1$ and $\Phi_2$ have been 
assumed to be alternating. Now the alternativity of~$K$ means that 
any pseudo-alternating form is alternating; in particular, $\Psi$ 
is alternating, hence $\Phi(z,z)=\Phi(f(t),f(t))=\Psi(t,t)=0$ as required.

\begin{DefD}
We say that a multilinear map $\Phi:M^n\to X$ is {\em skew-symmetric}, 
if it becomes a skew-symmetric bilinear map after 
an arbitrary choice of all arguments but two, or equivalently, if
\begin{equation}\label{eq:multilin.skewsym}
\Phi(x_1,\ldots,x_i,\ldots,x_j,\ldots,x_n)=
-\Phi(x_1,\ldots,x_j,\ldots,x_i,\ldots,x_n)\text{ for any $i<j$}
\end{equation}
Similarly, we say that $\Phi$ is {\em alternating\/} iff it is skew-symmetric, 
and it vanishes whenever any two of its arguments coincide:
\begin{equation}\label{eq:multilin.alt}
\Phi(x_1,\ldots,x_i,\ldots,x_j,\ldots,x_n)=0\text{ if $x_i=x_j$, $i<j$}
\end{equation}
Finally, when $M=K(m)$, we say that $\Phi$ is {\em pseudo-alternating\/} iff 
it is skew-symmetric, and vanishes on collections $(x_j)$ of base elements 
of~$K(m)$, such that $x_i=x_j$ for some $i\neq j$.
\end{DefD}

\nxsubpoint
Clearly, \eqref{eq:multilin.skewsym} for all $j=i+1$ suffices for $\Phi$ to 
be skew-symmetric, since it implies 
\begin{equation}
\Phi(x_{\sigma(1)},\ldots,x_{\sigma(n)})=\sgn(\sigma)\cdot\Phi(x_1,\ldots,x_n)
\text{ for any $\sigma\in\gS_n$}
\end{equation}
In order to show this one has just to decompose a permutation 
$\sigma\in\gS_n$ into a product of elementary transpositions.
Similarly, it would suffice to consider $j=i+1$ in 
\eqref{eq:multilin.skewsym} and $i=1$, $j=2$ in~\eqref{eq:multilin.alt} 
to establish that $\Phi$ is alternating.

Finally, one deduces from the bilinear map case that {\em 
if $K$ is alternating, any pseudo-alternating multilinear map
$\Phi:K(m)^n\to X$ is automatically alternating.}

\nxsubpoint (Exterior powers of a module.)
For any $K$-module~$M$ and any $n\geq0$ we can construct a universal 
alternating map $\wedge:M^n\to\bigwedge_K^nM$ by computing the strict quotient 
of~$M^{\otimes n}$ modulo relations coming from~\eqref{eq:multilin.skewsym}
and~\eqref{eq:multilin.alt}. This $K$-module $\bigwedge_K^nM$ is called 
the {\em $n$-th exterior power of~$M$.}

Notice that $\bigwedge_K^mM\otimes_K\bigwedge_K^nN$ has a universal property 
among all multilinear maps $M^m\times N^n\to X$, alternating with 
respect to each group of variables. Putting $N:=M$, $X:=\bigwedge_K^{n+m}M$, 
we obtain a canonical homomorphism $\bigwedge_K^mM\otimes_K\bigwedge_K^nM\to
\bigwedge_K^{n+m}M$. This yields a graded algebra structure on 
$\bigwedge_KM:=\bigoplus_{n\geq0}\bigwedge_K^nM$. This algebra is called 
{\em the exterior algebra of~$M$}. Clearly, it is a strict quotient of the 
tensor algebra of~$M$; it has a universal property with respect to 
$K$-linear homomorphisms $f:M\to A$ from~$M$ into algebras~$A$ in 
$\catMod K$, such that $f(x)f(y)=-f(y)f(x)$ and $f(x)^2=0$ for any 
$x$, $y\in M$. In other words, we require $M\times M\stackrel{f\times f}\to
A\times A\stackrel\mu\to A$ to be alternating.

\nxsubpoint (Exterior algebra of a direct sum.)
Notice that $A:=\bigwedge_K^nM=A^+\oplus A^-$, 
$A^+:=\bigoplus_{\text{$n$ even}}\bigwedge_K^nM$,
$A^-:=\bigoplus_{\text{$n$ odd}}\bigwedge_K^nM$, is a supercommutative 
algebra in $\catMod K$, i.e.\ $xy=(-1)^{\deg(x)\deg(y)}yx$ for any 
$x$, $y\in A^\pm$. We define the tensor product $A\otimes_K B$ of 
two supercommutative algebras by the usual rule 
$(x\otimes y)(x'\otimes y'):=(-1)^{\deg(x')\deg(y)}(xx'\otimes yy')$, 
$\deg(x\otimes y):=\deg(x)+\deg(y)$. One checks in the usual way that 
the resulting superalgebra is supercommutative, and that $A\otimes_K B$ 
is the coproduct in the category of supercommutative superalgebras 
in $\catMod K$.

It is easy to see that the exterior algebra $\bigwedge_KM$ together with the 
canonical map $M\to\bigwedge_KM$ is universal among all pairs $(A,\phi)$, 
where $A$ is a supercommutative algebra, and $\phi:M\to A^-$ is a $K$-linear 
map from $M$ into the odd part of~$A$, such that $\phi(x)^2=0$ for any 
$x\in M$; this is equivalent to requiring $M\times M\to A^-\times A^-\to A^+$ 
to be alternating.

Now suppose that $K$ is alternating. Let $M_1$ and $M_2$ be two $K$-modules, 
$A_i:=\bigwedge_KM_i$, $f_i:M_i\to A_i$ the canonical embedding. Then we 
have a commutative superalgebra $A:=A_1\otimes_KA_2$, and a canonical 
$K$-linear map $f:M:=M_1\oplus M_2\to A$, such that the restriction of~$f$ 
to $M_i$ equals $M_i\stackrel{f_i}\to A_i\to A$. One checks immediately 
that $f$ maps $M$ into the odd part of~$A$, hence 
$M\times M\stackrel{f\times f}\to A\times A\to A$ is skew-symmetric; moreover, 
the restriction of this bilinear map to $M_i\times M_i$ is alternating, 
hence it is alternating as well (cf.~\ptref{sp:ext.sqr.sum}; we use 
alternativity of~$K$ here). Now we see immediately that $(A,f)$ has the 
universal property required from~$\bigwedge_KM$, so we have proved 
\begin{equation}\label{eq:extalg.sum}
\wedgenl_K(M_1\oplus M_2)\cong
\bigl(\wedgenl_KM_1\bigr)\otimes_K\bigl(\wedgenl_KM_2\bigr)
\end{equation}
over an alternating monad~$K$. Taking individual graded components 
we obtain (under the same assumption)
\begin{equation}\label{eq:extpow.sum}
\wedgenl_K^n(M_1\oplus M_2)\cong\bigoplus_{p+q=n}
\bigl(\wedgenl_K^pM_1\bigr)\otimes_K\bigl(\wedgenl_K^qM_2\bigr)
\end{equation}

\nxsubpoint\label{sp:extalg.freemod} (Exterior algebra of a free module.)
Notice that $\bigwedge_KK(1)= K(1)\oplus K(1)$. When $K$ is alternating, 
one deduces from~\eqref{eq:extalg.sum} and~\eqref{eq:extpow.sum} that
\begin{align}
\wedgenl^r_K(M\oplus K(1))&\cong\wedgenl_K^rM\oplus\wedgenl_K^{r-1}M\\
\wedgenl_K(M\oplus K(1))&\cong\wedgenl_KM\otimes_K(K(1)\oplus K(1))
\end{align}
From this one shows by induction that $\bigwedge_KK(n)$ is a free $K$-module 
of rank~$2^n$. Its basis elements $e_I$ are parametrized by subsets 
$I\subset\stn$ in the usual way: $I=\{i_1,i_2,\ldots,i_r\}$, 
$1\leq i_1<i_2<\cdots<i_r\leq n$ corresponds to $e_I:=\{i_1\}\wedge\{i_2\}
\wedge\cdots\wedge\{i_r\}$. Considering individual graded pieces of the 
exterior algebra we see that the exterior power $\bigwedge_K^rK(n)$ is 
a free $K$-module of rank $\binom nr$, and its basis elements are parametrized
by $r$-element subsets $I\subset\stn$. One can extend these results to 
the exterior powers of $K(S)$, $S$ any linearly ordered set, by observing that 
exterior powers and algebras commute with filtered inductive limits.

\nxsubpoint (Exterior algebra and pullbacks.)
Notice that the exterior algebra $\bigwedge_KM$, considered as a unary 
$K$-algebra, is universal among all pairs $(\Sigma,f)$, where 
$\Sigma$ is a $K$-algebra, and $f:M\to|\Sigma|$ is a $K$-linear map, 
such that $f(x)^2=0$ for all $x\in M$. 
Indeed, the universal property among all unary~$\Sigma$ has been already 
discussed, and in the general case we can replace $\Sigma$ by its 
unary envelope.

Another easy observation: when $K$ is alternating, a skew-symmetric map 
$\Phi:M\times M\to X$ is alternating iff $\Phi(x,x)=0$ for all $x$ from 
a system of generators of~$M$. Combining this with our previous remark we see 
that (over an alternating~$K$) the exterior algebra~$\bigwedge_KM$ is 
universal among all $K$-linear maps $f:M\to|\Sigma|$, $\Sigma$ a $K$-algebra,
such that $f(x)^2=0$ for all $x\in S$, where $S$ is any fixed system of 
generators of~$M$.

Now let $\rho:K\to K'$ be a homomorphism of alternating $\Fpm$-algebras, 
$M$ be a $K$-module, and $S$ any system of generators of~$M$ 
(e.g.\ $M$ itself). Then $\bigwedge_{K'}\rho^*M$ is universal among all 
$K'$-linear maps 
$f:\rho^*M\to|\Sigma'|$, with $\Sigma'$ a $K'$-algebra, such that 
$f(\bar s)^2=0$ for all $s\in S$, where $\bar s$ denotes the image of 
$s\in S\subset M$ under the canonical map $M\to\rho^*M$ 
(we use here the fact that the image of~$S$ in~$\rho^*M$ is a system of 
generators). Clearly, such maps are in one-to-one correspondence with 
$K$-linear maps $f^\flat:M\to\rho_*|\Sigma|=|\rho_*\Sigma|$, such that 
$(f^\flat(s))^2=0$ for any $s\in S$. Such $f^\flat$ are in one-to-one 
correspondence with $K$-algebra homomorphisms $\bigwedge_KM\to\rho_*\Sigma$. 
Comparing the universal properties involved, we see that, for any
$\rho$ and $M$ as above we have
\begin{equation}
\wedgenl_{K'}\rho^*M\cong\rho^*\bigl(\wedgenl_KM\bigr)=K'\otimes_K
\wedgenl_KM
\end{equation}
Considering individual graded pieces we obtain
\begin{equation}
\wedgenl_{K'}^r\rho^*M\cong\rho^*\bigl(\wedgenl_K^rM\bigr)
\end{equation}
Notice that the alternativity of~$K'$ is important in this reasoning.

\nxsubpoint (Determinants.)
Let $K$ be an alternating monad (i.e.\ an alternating commutative 
$\Fpm$-algebra), $M$ be a free $K$-module of rank~$n$, and 
$u\in\End_K(M)$. We have seen in~\ptref{sp:extalg.freemod} that 
$\bigwedge_K^nM$ is a free $K$-module of rank~1, hence the 
endomorphism $\wedge^nu$ of this module defines an element 
$\det(u)\in|K|\cong\End_K(\bigwedge_K^nM)$, such that 
$(\wedge^nu)(\alpha)=\det(u)\cdot\alpha$ for any $\alpha\in\bigwedge_K^nM$. 
Of course, we say that $\det(u)$ is the {\em determinant\/} of~$u$; this 
definition can be extended to the case of any $K$-module $M$ and any integer
$n>0$, such that $\bigwedge_K^nM$ is free of rank~one. 
When $M=K(n)$, $\End_K(K(n))\cong K(n)^n=M(n,n;K)$ is the set of 
$n\times n$-matrices over~$K$, so we have a notion of determinant for such 
matrices as well.

Functoriality of exterior powers implies $\det(\id_M)=\bu$ and 
$\det(v\circ u)=\det(v)\cdot\det(u)$, for any $u$, $v\in\End_K(M)$. 
Hence {\em if $u$ is invertible, $\det(u)$ is also invertible, and\/ 
$\det(u^{-1})=\det(u)^{-1}$}. Unfortunately, the converse doesn't seem to 
be true in general; we'll discuss this point later in more detail.

It is also easy to check that $\det(u\oplus v)=\det(u)\det(v)$ for 
any $u\in M(n,n;K)$ and $v\in M(m,m;K)$, so we obtain indeed a reasonable 
theory of determinants, up to the point already mentioned above. Notice that 
these determinants are not so easy to compute as one might think. 
For example, if $A=(T,U)\in K(2)^2$ is a $2\times2$-matrix, we have 
\begin{equation}\label{eq:det.2x2}
\det(A)=T(U(0,\bu),U(-\bu,0))=U(T(0,\bu),T(-\bu,0)) 
\end{equation}

\nxsubpoint\label{sp:examp.altmon} 
(Examples and constructions of alternating monads.)
We see that we obtain a reasonable theory of exterior powers and determinants 
only over an alternating $\Fpm$-algebra~$K$. Therefore, we would like 
to know that there are sufficiently many of them.
\begin{itemize}
\item[a)] First of all, any classical commutative ring is an 
alternating $\Fpm$-algebra (cf.~\ptref{sp:alt.binop}). Thus 
$\bbZ$, $\bbQ$, $\bbR$, $\bbC$ are alternating.
\item[b)] Any $\Fpm$-subalgebra $K'$ of an alternating algebra~$K$ is 
also alternating, since the alternativity condition $\wedge t$ for an 
operation $t\in K'(n)$ in~$K'$ is equivalent to the corresponding condition 
in~$K$ (cf.~\ptref{sp:rem.altmon}). 
Hence $\Zinfty$, $\barZinfty$ and $\Zninfty$ are alternating, being 
$\Fpm$-subalgebras of~$\bbC$.
\item[c)] Recall that an $\Fpm$-algebra $K$ is alternating iff it is generated 
by a system of alternating operations (cf.~\ptref{sp:rem.altmon}), and that 
all operations of arity $\leq1$ are automatically alternating. Therefore, 
all $\Fpm$-algebras (absolutely) generated in arity $\leq1$ are alternating, 
e.g.\ $\Fpm$ itself, $\bbF_{1^n}$ for an even $n>0$, and $\bbF_{1^\infty}$.
\item[d)] Notice that any projective limit of alternating $\Fpm$-algebras 
is itself alternating. Indeed, if $K=\projlim K_\alpha$, then any 
$t\in K(n)=\projlim K_\alpha(n)$ is given by a compatible family of 
$t_\alpha\in K_\alpha(n)$, and the alternativity condition $\wedge t$ 
in $K(n(n-1)/2)^2$ is given by the compatible family of 
alternativity conditions $\wedge t_\alpha$, 
hence if all $\wedge t_\alpha$ hold, the same is true 
for~$\wedge t$. In particular, 
$\bbF_{\pm1^{-\infty}}=\projlim_{\text{$n$ odd}}\bbF_{1^{2n}}$ 
(cf.~\ptref{sp:examp.pres.comalg},c) is alternating.
\item[e)] Any strict quotient of an alternating algebra is clearly 
alternating. Hence $\Finfty$ is alternating, being a strict quotient 
of~$\Zinfty$. We might also check the alternativity condition 
$(0*\bu)*(-\bu*0)=0$ for the binary generator~$*$ of~$\Finfty$ directly, 
using $\Finfty=\Fpm[*^{[2]}\,|\,\bu*(-\bu)=0$, $\bu*\bu=\bu$, 
$x*y=y*x$, $(x*y)*z=x*(y*z)]$.
\item[f)] Since an alternating operation remains alternating after an 
application of a homomorphism of commutative $\Fpm$-algebras, we see that 
all sorts of inductive limits of alternating $\Fpm$-algebras are alternating.
For example, if $K$ is a commutative $\Fpm$-algebra, and 
$K_1$ and $K_2$ are alternating $K$-algebras, then $K_1\otimes_KK_2$ is 
also alternating, being generated by the operations of~$K_1$ and~$K_2$. 
\item[g)] Any (commutative) pre-unary algebra $K'$ over an alternating 
$\Fpm$-algebra~$K$ is alternating, since it is generated by the operations 
of~$K$ together with a family of unary, hence automatically alternating 
generators. For example, unary polynomial rings $\Fpm[T^{[1]},U^{[1]},\ldots]$
and $\Zinfty[T^{[1]},U^{[1]},\ldots]$ are alternating.
\item[h)] This is not true for binary polynomial rings: for example, 
$\Fpm[T^{[2]}]$ is not alternating, i.e.\ $T(T(0,\bu),T(-\bu,0))\neq0$ in this 
generalized ring (otherwise any binary operation would have been alternating!).
We can consider the largest alternating quotients 
$\Fpm[T_1^{[r_1]},\ldots]/_{alt}$ of these polynomial rings instead, 
thus obtaining free alternating $\Fpm$-algebras.
\end{itemize}

\nxsubpoint (Existence of $SL_{n,\Fpm}$.) 
Notice that we can define an affine group scheme~$SL_n$ over~$\Fpm$, such that
for any {\em alternating\/} $\Fpm$-algebra~$K$ we have 
\begin{equation}
SL_n(K)=\{g\in GL_n(K)\subset M(n,n;K)\,|\,\det(g)=\bu\}
\end{equation}
Of course, $SL_{n,\bbZ}$ coincides with the classical group scheme denoted 
in this way, since they both represent the same functor.

\nxpointtoc{Matrices with invertible determinant}\label{p:invdet.matr}
It would be nice to know that all matrices with invertible determinants 
are themselves invertible. However, we haven't managed to prove or disprove 
this statement so far. Therefore, we'll content ourselves with some partial 
results.

\nxsubpoint
We consider the following properties of an alternating $\Fpm$-algebra~$K$:

\noindent $(DET_r)$ Any $r\times r$-matrix $A\in M(r,r;K)$ with invertible 
determinant is invertible.

\noindent $(DET'_r)$ For any $r\times r$-matrix $A\in M(r,r;K)$ one can 
find an integer $N\geq1$ and a matrix $B\in M(r,r;K)$, such that 
$AB=BA=\det(A)^N$.

\noindent $(DET^*_r)$ Same as above, but with $N$ independent on the 
choice of~$A$.

\noindent $(DET^?_\infty)$ means that the properties $(DET^?_r)$ are fulfilled 
for all $r\geq0$; here the superscript $?$ can be replaced by $'$, $*$, 
or nothing at all.

Notice that these three properties trivially hold for $r\leq1$ (with $N=1$), 
so we might consider only $r\geq2$. Also note that if $(DET'_r)$ is satisfied 
for some $N\geq1$, we can replace $N$ by any $N'\geq N$: indeed, 
we just have to consider $B':=B\cdot\det(A)^{N'-N}$.

Clearly, $(DET^*_r)\Rightarrow(DET'_r)\Rightarrow(DET_r)$; hence 
$(DET^*_\infty)$ is the strongest among all conditions under consideration.

\nxsubpoint\label{sp:DET.decr.r}
An important point is that $(DET^?_{r+1})$ implies $(DET^?_r)$, hence 
$(DET^?_r)$ implies all $(DET^?_s)$ with $0\leq s\leq r$. Notice that this 
is true for $r=\infty$ as well.

So let's fix a matrix $A\in M(r,r;K)$ as required in~$(DET^?_r)$, and 
put $A':=A\oplus\bu\in M(r+1,r+1;K)$. Let's denote by $\lambda:K(r)\to K(r+1)$ 
and $\sigma:K(r+1)\to K(r)$ the natural embedding and projection ($\sigma$ 
maps $\{r+1\}$ into 0). Then $\sigma\lambda=I_r=\id_{K(r)}$, and obviously 
$\sigma A'\lambda=A$ and $\det(A')=\det(A)$. Now let's apply $(DET^?_{r+1})$ 
to~$A'$. We obtain a matrix $B'\in M(r+1,r+1;K)$, such that 
$A'B'=B'A'=\det(A')^N=\det(A)^N$ (with $N=0$ for the $(DET_r)$ property). 
Now put $B:=\sigma B'\lambda\in M(r,r;K)$. Clearly, $AB=A\sigma B'\lambda=
\sigma A'B'\lambda=\det(A)^N\cdot \sigma I_{r+1}\lambda=\det(A)^N\cdot I_r$, 
and similarly $BA=\sigma B'\lambda A=\sigma B'A'\lambda=\det(A)^N\cdot I_r$, 
i.e.\ $(DET^?_r)$ holds for~$A$.

\nxsubpoint
Notice that properties $(DET^*_r)$ and $(DET'_r)$ are stable under taking 
strict quotients, i.e.\ if $\rho:K\to\bar K$ is a strict epimorphism, and 
if one of these properties is valid for~$K$, it is valid for~$K'$ as well, 
with the same value of~$N=N(r)$ for the $(DET^*_r)$ property. 
Indeed, the induced maps $\rho_r^r:
M(r,r;K)=K(r)^r\to M(r,r;\bar K)$ are surjective, 
so we can start with an arbitrary matrix $\bar A\in M(r,r;\bar K)$, lift it 
arbitrarily to a matrix $A\in M(r,r;K)$, find a matrix $B\in M(r,r;K)$ 
that satisfies $(DET'_r)$, and consider $\bar B:=
\rho_r^r(B)\in M(r,r;\bar K)$. However, this reasoning is not applicable to 
the $(DET_r)$ properties, since in general we have no means to lift 
a matrix with invertible determinant $\bar A$ to a matrix~$A$ over~$K$ 
without violating the invertibility of the determinant.

An immediate application is that if we show that $\Zinfty$ satisfies 
$(DET^*_\infty)$ (the strongest of our properties), the same will be true 
for~$\Finfty$, and in particular any matrix over~$\Finfty$ with invertible 
determinant will be invertible. Actually, one can check $(DET_r)$ 
for~$\Finfty$ directly, by observing that the matrices with non-zero, i.e.\ 
invertible determinant in $M(r,r;\Finfty)$ are exactly the matrices from 
$GL_r(\Fpm)\subset M(r,r;\Finfty)$, i.e.\ the permutation matrices 
$(\pm\{\sigma_1\},\ldots,\pm\{\sigma_r\})$, with~$\sigma\in\gS_r$.

\nxsubpoint ($(DET^*_\infty)$ for $\Zinfty$ and~$\barZinfty$.)
One can try to show directly that $\Zinfty$ and $\barZinfty$  
satisfy~$(DET^*_\infty)$, as follows. A matrix $A\in M(r,r;\barZinfty)$ 
is actually a matrix $A=(a_{ij})\in M(r,r;\bbC)=\bbC^{\st r\times \st r}$ 
with the $L_1$-norm $\|A\|_1\leq 1$, i.e.\ with $\sum_i|a_{ij}|\leq 1$ for 
all~$j$. We have to check that the matrix $B:=A^*\cdot\Delta^{N-1}$ 
has $L_1$-norm $\leq1$ for some $N=N(r)>0$, where $A^*$ is the adjoint matrix 
of~$A$ and $\Delta:=\det A$, i.e.\ that $\|A^*\|_1\leq|\Delta|^{1-N}$ for 
all~$A\in M(r,r;\bbC)$ with $\|A\|_1\leq1$. We can remove the latter 
condition by rewriting our inequality:
\begin{equation}
\|A^*\|_1\cdot|\det A|^{N_r-1}\leq\|A\|_1^{rN_r-1}
\end{equation}
(The $\Zinfty$ case then follows, since if this inequality holds for 
matrices with complex coefficients, it holds {\em a fortiori\/} for 
matrices with real coefficients.) 
This inequality can be checked directly; for example, for $r=2$ we can 
take $N_r=2$, thus obtaining the following inequality:
\begin{equation}
\bigl(|a_{11}|+|a_{21}|\bigr)\cdot\bigl|a_{11}a_{22}-a_{12}a_{21}\bigr|\leq
\sup\bigl(|a_{11}|+|a_{12}|,|a_{21}|+|a_{22}|\bigr)^3
\end{equation}
We don't provide more details since we still hope to prove a more general 
result that would include unary polynomial rings over~$\Zinfty$ 
and~$\barZinfty$ as well. 

\nxsubpoint
Another simple observation is that all classical rings satisfy $(DET^*_r)$ 
with $N=1$. Indeed, we might take $B=A^*$, the adjoint matrix of~$A$; 
then $AA^*=A^*A=\det(A)$.

\nxsubpoint\label{sp:use.univ.matr} (Universal matrix.)
Fix an integer $r\geq0$ and consider the free alternating algebra
$\Lambda:=\Fpm[T_1^{[r]},\ldots,T_r^{[r]}]/_{alt}$ and the ``universal 
$r\times r$-matrix'' $A:=(T_1,\ldots,T_r)\in\Lambda(r)^r=M(r,r;\Lambda)$. 
Clearly, $(\Lambda,A)$ represents the functor $K\mapsto M(r,r;K)$ on the 
category of alternating $\Fpm$-algebras. We see that if 
$(DET'_r)$ is valid just for this matrix~$A$ over~$\Lambda$, with 
some value of $N=N_r$, then the $(DET'_r)$ and $(DET^*_r)$, hence also 
$(DET_r)$ properties are fulfilled for {\em all\/} alternating generalized 
rings, with the same value of~$N_r$. Therefore, the validity of $(DET'_r)$ 
for all alternating monads is equivalent to the validity of $(DET^*_r)$ for 
all alternating monads.

\nxsubpoint (Universal matrix with invertible determinant.)
In the above notations put $\Delta:=\det A\in|\Lambda|$ and consider the 
localization $\tilde\Lambda:=\Lambda[\Delta^{-1}]=\Lambda[\bar\Delta^{[1]}\,
|\,\Delta\bar\Delta=\bu]$; this is an alternating $\Fpm$-algebra, since 
it is unary over~$\Lambda$. Clearly, the image $\tilde A$ of $A$ in 
$M(r,r;\tilde\Lambda)$ is the ``universal $r\times r$-matrix with invertible 
determinant'', i.e.\ $(\tilde\Lambda,\tilde A)$ represents the functor 
$K\mapsto\{r\times r$-matrices over~$K$ with invertible determinant$\}$. 
Suppose that 
$(DET_r)$ is valid just for this matrix $\tilde A$ over~$\tilde\Lambda$, 
so that there is an inverse matrix $\tilde B$ to $\tilde A$. We'll see in the 
next chapter that the localization $\tilde\Lambda=\Lambda[\Delta^{-1}]$ admits 
the usual description, i.e.\ $\tilde\Lambda(n)$ consists of expressions 
$t/\Delta^k$, $t\in\Lambda(n)$, $k\geq0$, and $t/\Delta^k=t'/\Delta^{k'}$ iff 
$\Delta^{k'+m}\cdot t=\Delta^{k+m}\cdot t'$ for some $m\geq0$. Using this 
fact we can write $\tilde B\in M(r,r;\tilde\Lambda)=\tilde\Lambda(r)^r$ in form
$B/\Delta^m$ for some $B\in M(r,r;\Lambda)$. The equalities 
$(A/1)\cdot(B/\Delta^m)=I_r/1=(B/\Delta^m)\cdot(A/1)$ imply the existence of 
some $N\geq m$, such that $A\cdot(\Delta^{N-m}B)=\Delta^N\cdot I_r=
(\Delta^{N-m}B)\cdot A$ in $M(r,r;\Lambda)$. Therefore, in this case 
$(DET'_r)$ is fulfilled for the universal matrix~$A$ with this value of 
$N=N_r$, hence $(DET^*_r)$ is valid for all alternating monads with the 
same value of~$N_r$ by~\ptref{sp:use.univ.matr}. We see that 
{\em $(DET_r)$ holds over all alternating monads iff it holds 
for $\tilde A$ over~$\tilde\Lambda$ iff $(DET^*_r)$ holds over all 
alternating monads}, and in this case {\em an universal value of $N=N_r$ 
can be chosen for all alternating monads.}

\nxsubpoint
We see that proving or disproving the validity of $(DET^*_r)$, 
$(DET'_r)$ or $(DET_r)$ for all alternating monads is equivalent to 
proving or disproving $(DET'_r)$ for the universal matrix 
$A=(T_1,\ldots,T_r)$ over alternating generalized ring 
$\Lambda=\Fpm[T_1^{[r]},\ldots,T_r^{[r]}]/_{alt}$, 
so we'll concentrate our efforts on this case.

Let's introduce a technique that enables one to obtain lower bounds for 
the (universal) $N_r$, and to obtain good candidates for the matrix~$B$ 
of~$(DET'_r)$ for high values of~$N$. We consider for this the canonical 
map $\Lambda\to\Lambda_\bbZ=\bbZ\otimes_{\Fpm}\Lambda=\bbZ[T_1^{[r]},\ldots,
T_r^{[r]}]\cong\bbZ[T_{11},T_{12},\ldots,T_{1r},T_{21},\ldots,T_{rr}]$ 
(cf.~\ptref{sp:ex.aff.grp.sch}), where $(T_{i1},\ldots,T_{ir})$ 
is the image of~$T_i$ under the comparison map~$\pi_r:\Lambda(r)\to|\Lambda|^r$
(cf.~\ptref{sp:compar.maps}), e.g.\ $T_{i1}=T_i(\bu,0,\ldots,0)$. 
Clearly, the map $\Lambda\to\Lambda_\bbZ$ maps $T_i\in\Lambda(r)$ into 
$T_{i1}\{1\}+\cdots+T_{ir}\{r\}\in\Lambda_\bbZ(r)$; we'll write this 
in the following form:
\begin{equation}
T_i=T_{i1}\{1\}+T_{i2}\{2\}+\cdots+T_{ir}\{r\}\quad\text{(over~$\Lambda_\bbZ$)}
\end{equation}

Let us denote by $\Lambda'$ the image of $\Lambda\to\Lambda_\bbZ$, 
and by $A'$, $B'$ and $\Delta'$ the images under this map of 
the universal matrix~$A=(T_1,\ldots,T_r)$, the (hypothetical) matrix $B$
with property $(DET_r)$ with respect to~$A$ and some integer $N=N_r\geq1$, 
and the determinant $\Delta:=\det(A)$.

Clearly, $A'=(T_{ji})$, so that $\Delta'=\det A'=\det(T_{ij})$ is the 
determinant of the universal matrix over~$\bbZ$. If we have a matrix~$B$ 
as in~$(DET_r)$, we have $A'B'=(\Delta')^N=B'A'$ over 
$\Lambda'\subset\Lambda_\bbZ\subset\bbQ(T_{11},T_{12},\ldots,T_{rr})$. 
Since the determinant of~$A'$ is invertible in this field, $A'$ is itself 
invertible, hence $B'=(\Delta')^N\cdot(A')^{-1}=(\Delta')^{N-1}\cdot
(A')^*$ is completely determined by $N>0$.

On the other hand, the components $B'_j=(b'_{1j},b'_{2j},\ldots,b'_{rj})\in 
|\Lambda_\bbZ|^r=\Lambda_\bbZ(r)$ of the matrix~$B'$ have to lie in
$\Lambda'(r)$. Since $\Lambda'\subset\Lambda_\bbZ$ is generated by the 
images of the~$T_j$ in $\Lambda_\bbZ$, we see that $\Lambda'(r)\subset
\Lambda_\bbZ^r$ consists of expressions that can be obtained by applying 
the following rules finitely many times:
\begin{itemize}
\item $0$ and $\pm\{i\}$, for $1\leq i\leq r$, belong to $\Lambda'(r)$;
\item If $f_1$, \dots, $f_r$ lie in $\Lambda'(r)$, the same is true 
for $T_{i1}f_1+\cdots+T_{ir}f_r$, for any $1\leq i\leq r$.
\end{itemize}

If $(DET'_r)$ holds for $A$ with some $N\geq0$, then it holds for 
$A'$ and the same $N$ as well, hence the matrix $B'=(\Delta')^N\cdot(A')^{-1}$ 
lies in~$\Lambda'(r)$. Clearly, for $r\geq2$ this can happen only for~$N>0$, 
so that $B'=(\Delta')^{N-1}\cdot(A')^*$, i.e.\ $B=(b_{ij})$ with 
$b_{ij}=(-1)^{i+j}\Delta_{ji}\cdot{\Delta'}^{N-1}$, where $\Delta_{ji}$ is the 
corresponding principal minor of~$A'$. Therefore, we must have 
$(b_{11},b_{21},\ldots,b_{r1})\in\Lambda'(r)$ (and similarly for the 
other columns of~$B'$, but by symmetry it suffices to consider only the first 
column). 

Therefore, if for some value of~$N>0$ we can show that 
$(b_{11},\ldots,b_{r1})$ doesn't lie in~$\Lambda'(r)$, then 
the value of universal $N_r$ is at least $N+1$. On the other hand, if 
this vector lies in~$\Lambda'(r)$, then $(DET'_r)$ is fulfilled for
$A'$ and this value of~$N$. In this case we might consider different 
lifts~$B$ of $B'$ to $\Lambda$ and check whether $(DET'_r)$ is fulfilled 
for any of them. Notice that if we would show the injectivity of 
$\Lambda\to\Lambda_\bbZ$, the latter step would be unnecessary.

\nxsubpoint
Let's apply the above considerations to the first non-trivial case $r=2$; 
the case~$r\geq3$ seems to be more complicated only from the technical point 
of view, so all the interesting issues can be seen already for~$r=2$.

In this case we write $T$ and $U$ instead of $T_1$ and~$T_2$; so we have 
$\Lambda=\Fpm[T^{[2]},U^{[2]}]/_{alt}$, $A=(T,U)$, $\Lambda_\bbZ= 
\bbZ[T_1,T_2,U_1,U_2]$, where of course $T_1=T(\bu,0)$, $T_2=T(0,\bu)$ and 
so on. By~\eqref{eq:det.2x2} we obtain $\Delta=\det A=
T(U(0,\bu),U(-\bu,0))=U(T(0,-\bu),T(\bu,0))$; the image $\Delta'$ of 
$\Delta$ in $\Lambda_\bbZ$ is of course $T_1U_2-T_2U_1$, and 
$A'=\smallmatrix{T_1}{U_1}{T_2}{U_2}$, 
${A'}^*=\smallmatrix{U_2}{-U_1}{-T_2}{T_1}$.

We denote the components of $B$ by $(V,W)$; clearly, $V$ and $W\in\Lambda(2)$ 
have to satisfy
\begin{eqnarray}
T(V(\{1\},\{2\}),W(\{1\},\{2\}))&=&\Delta^N\cdot\{1\}\\
U(V(\{1\},\{2\}),W(\{1\},\{2\}))&=&\Delta^N\cdot\{2\}\\
V(T(\{1\},\{2\}),U(\{1\},\{2\}))&=&\Delta^N\cdot\{1\}\\
W(T(\{1\},\{2\}),U(\{1\},\{2\}))&=&\Delta^N\cdot\{2\}
\end{eqnarray}

We denote by $B'$ the image of $B$ in $\Lambda'\subset\Lambda_\bbZ$; 
its components $B'=(V',W')$ are given by $V'=U_2\Delta^{N-1}\{1\}-T_2
\Delta^{N-1}\{2\}$ and $W'=-U_1\Delta^{N-1}\{1\}+T_1\Delta^{N-1}\{2\}$. 
We want to find those values of $N\geq1$, for which $V'$ and~$W'$ lie 
indeed in~$\Lambda'$, and find possible values of~$V$ and~$W$ for such 
values of~$N$.

\nxsubpoint (Elements of~$\Lambda$.)
Notice that $\Lambda=\Fpm[T^{[2]},U^{[2]}]/_{alt}$ 
is generated by one constant~$0$, 
one unary operation~$-$, and two binary operations $T$ and~$U$. Therefore, 
elements of~$\Lambda(n)$ are formal expressions 
(terms, cf.~\ptref{sp:struct.ind}), 
usually written in prefix form, constructed with the aid of these operations 
from free variables $\{i\}_\stn$, $1\leq i\leq n$. Sometimes we omit 
braces around $\{i\}$, since positive integers have no other possible meaning 
in this context; thus we can write $T\,-1\,U\,2\,0$ instead of
$T(-\{1\},U(\{2\},0))$. Let us list all rules for the manipulation of 
these expressions, including implicit commutativity and alternativity 
relations; here $x$, $y$, $z$, \dots\ replace either free variables or 
arbitrary valid expressions:
\begin{enumerate}
\item $--x=x$, $-0=0$, $-Txy=T\,-x\,-y$, $-Uxy=U\,-x\,-y$;
\item $T00=0$, $U00=0$;
\item $T\,Txy\,Tzw=T\,Txz\,Tyw$;\quad
  $T\,Uxy\,Uzw=U\,Txz\,Tyw$;\\
  $U\,Uxy\,Uzw=U\,Uxz\,Uyw$;
\item $T\,T0x\,T-x0=0$;\quad $U\,U0x\,U-x0=0$.
\end{enumerate}
The relation $--x=x$ is nothing else than the relation $-^2=\bu$ from the 
definition of~$\Fpm$, the last group of relations are just the alternativity 
relations for~$T$ and~$U$, and the remaining relations are the commutativity 
relations between all pairs of generators.

Notice that the first set of relations allows us to rewrite any valid 
expression in such a way that unary $-$ is applied only to free variables.
This means that we can represent valid expressions by means of binary trees, 
with interior nodes marked by~$T$ and~$U$, and leaves marked by~$0$ 
or $\pm i$, $1\leq i\leq n$ (of course, $+i$ or~$i$ stands for $\{i\}_\stn$, 
and $-i$ for $-\{i\}_\stn$).

When we need to define a new operation in~$\Lambda(n)$, we usually write down 
either formulas like $Z\,x_1\ldots x_n=\langle$expression in 
variables $\pm x_i\rangle$, or formulas like $Z=\langle$expression in $\pm i$, 
$1\leq i\leq n\rangle$. For example, $\Delta$ can be defined either by
$\Delta\,z=T\;U\,0\,z\;U\,-z\,0$, or by $\Delta=T\;U\,0\,1\;U\,-1\,0$.

\nxsubpoint (Elements of $\Lambda'$.) 
Given two elements $f$ and $f'\in\Lambda(n)$, we write $f\sim f'$ and say that
{\em $f$ and $f'$ are weakly equivalent\/} if they have the same image~$\bar f$
in~$\Lambda'(n)$, or equivalenty, in~$\Lambda_\bbZ(n)$. 
Then $\Lambda'(n)$ 
is in one-to-one correspondence with weak equivalence classes of elements 
of~$\Lambda(n)$, and these elements can be represented by binary trees as 
described above. Of course, the image~$\bar f$ of~$f$ in~$\Lambda_\bbZ$ 
is computed by means of the following very simple rules:
\begin{itemize}
\item $0\in\Lambda(n)$ is mapped into $0=(0,\ldots,0)\in\Lambda_\bbZ(n)$.
\item $\pm i$, i.e.\ $\pm\{i\}_\stn$ is mapped into $\pm\{i\}=(0,\ldots,\pm1,
\ldots,0)$ in~$\Lambda_\bbZ(n)$.
\item If $f=T(f_1,f_2)$, i.e.\ if $f$ is described by a tree with root 
marked by~$T$, left subtree~$f_1$, and right subtree~$f_2$, then 
$\bar f=T_1\cdot\bar f_1+T_2\cdot \bar f_2$.
\item Similarly, $f=U(f_1,f_2)$ is mapped into 
$\bar f=U_1\cdot\bar f_1+U_2\cdot \bar f_2$
\end{itemize}
We see that any leaf contributes some monomial into~$\bar f$. More precisely, 
a leaf marked by zero contributes nothing, and a leaf marked by $\pm i$ 
contributes $\pm T_1^aT_2^bU_1^cU_2^d\{i\}$, where $a$ is the amount of 
left branches, taken after vertices marked by $T$ in the path 
from the root of the tree to the leaf under consideration, and 
integers $b$, $c$, $d\geq0$ are defined similarly.

An immediate consequence is that {\em replacing~$f$ by a weakly equivalent 
tree if necessary, we can assume that there is no cancellation between 
contributions of separate leaves.} Indeed, suppose that  
$+T_1^aT_2^bU_1^cU_2^d\{i\}$ cancels with $-T_1^aT_2^bU_1^cU_2^d\{i\}$. 
These two terms are contributed by two leaves, one of them marked by $i$, and 
the other by $-i$. Then we can simply replace both labels by~$0$; since 
this operation doesn't change~$\bar f$, but decreases the number of 
non-zero labels on leaves of~$f$, this process will stop after a finite 
number of steps.

Using relations $T\,0\,0=U\,0\,0=0$ if necessary, we see that {\em 
if $\bar f$ is homogeneous of degree~$d$, $f$ is weakly equivalent to 
a complete binary tree of depth~$d$, i.e.\ to a binary tree, in which all 
leaves are at distance~$d$ from the root. We can also assume that there 
are no cancellations between contributions of different leaves of~$f$.} 
Clearly, there are $2^d$ leaves 
in this case; each of these leaves can be labeled by~$0$ or $\pm i$. There 
are also $2^d-1$ intermediate nodes as well, labeled by letters~$T$ and~$U$.

\nxsubpoint (Lower bound for universal $N_2$.)
Let's apply these consideration to obtain a lower bound for the 
universal~$N_2$. So suppose that $(DET'_2)$ holds for the universal 
$2\times2$-matrix~$A$ for some $N>0$. Then $V'=U_2\Delta^{N-1}\{1\}-
T_2\Delta^{N-1}\{2\}$, $\Delta=T_1U_2-T_2U_1$, 
is of form $\bar V$ for some binary tree~$V$. Notice that $V'$ 
is homogeneous of degree~$2N-1$. Using the above results 
we can replace $V$ by a complete binary tree of depth~$2N-1$, such that 
there is no cancellation between contributions of individual leaves. 
Now let's consider individual monomials in $V'$. We see that all of them 
are of form $\pm T_1^aT_2^bU_1^cU_2^d\{i\}$ with $a+c=N-1$, $b+d=N$, 
and that $V'$ is the sum of~$2^N$ such monomials 
(one checks that there is no cancellation between individual monomials 
in the expression for~$V'$ by substituting $T_1:=U_1:=U_2:=1$, $T_2:=-1$). 
On the other hand, the only leaves in the complete binary tree of depth~$2N-1$ 
that might contribute such monomials are those that correspond to 
$(2N-1)$-digit binary numbers with $N-1$ zeroes and $N$ ones, if we agree to 
encode the path from the root to a leaf by means of a sequence of $2N-1$ 
binary digits, where $0$ corresponds to taking the left branch, and $1$ to 
the right branch. There are $\binom{2N-1}{N}$ such leaves, hence we can have 
at most $\binom{2N-1}{N}$ such monomials in $\bar V$. We have proved that the 
following inequality is necessary for the existance of~$V$ with $\bar V=V'$:
\begin{equation}\label{eq:lbound.N2}
\binom{2N-1}{N}\geq2^N
\end{equation}
One checks immediately that this condition is equivalent to $N\geq3$: indeed,
for $N=1$ we have $\binom11=1<2^1$, for $N=2$ we have $\binom32=3<2^2$, 
but for $N=3$ we get $\binom53=10>2^3$, and the inequality holds for 
$N>3$ as well. We conclude that {\em universal $N_2$ is at least~$3$.} 
Notice that for~$\Zinfty$ and~$\barZinfty$ we might take $N_2=2$, and for 
any classical ring we might even take $N_2=1$.

\nxsubpoint (Case $r=2$, $N=3$.)
Consider the case $r=2$, $N=3$. Then $V'=(T_1U_2-T_2U_1)^2(U_2\{1\}-T_2\{2\})=
T_1^2U_2^3\{1\}-2T_1T_2U_1U_2^2\{1\}+T_2^2U_1^2U_2\{1\}-T_1^2T_2U_2^2\{2\}+
2T_1T_2^2U_1U_2\{2\}-T_3^2U_1^2\{2\}$. We want to find a tree 
$V\in\Lambda(2)$, such that $\bar V=V'$. We know that if such a~$V$ exists, 
it can be replaced by a weakly equivalent complete binary tree of depth 5, 
without cancellation of contributions of different leaves. So in order to 
find such a~$V$ we can draw a complete binary tree of depth 5, mark by a 
``?'' those leaves that are encoded by a binary number with 3 ones and 2 zeroes
(there will be $\binom53=10$ such leaves), mark the remaining $32-10=22$ 
leaves by 0, and then try to replace eight question marks with labels 
$\pm i$, $i=1,2$ (each label must occur exactly twice) and two remaining 
question marks by zeroes, and fill in the intermediate nodes (there are 31 of 
them) by letters $T$ and $U$, so as to have the sum of contributions of 
leaves equal to~$V'$.

It is easy to find at least one such configuration; using rules $T00=U00=0$ 
we can replace subtrees with no non-zero leaves by a zero, thus obtaining 
a shorter expression for a~$V$ with~$\bar V=V'$. Here is one possible 
value for such~$V$:
\begin{equation}
V=TTU0U0U01U0UU0-1U-10TTU0U0-2UU020UUT0-2T20UU100
\end{equation}
Of course, we can obtain a $W$ with $\bar W=W'$ by replacing 
$Txy$ with $Uyx$ and $Uxy$ with $Tyx$ everywhere in the expression for~$V$, 
as well as $\pm\{i\}$ with~$\pm\{3-i\}$; 
in other words, we reflect the tree for~$V$ with respect to the vertical 
line and interchange $T\leftrightarrow U$ in the inner nodes and 
$\pm1\leftrightarrow\pm2$ in the leaves.

We can even find all possible values of~$V$, arising from a complete binary 
tree of depth~5; this is slightly more complicated than it seems since we 
might have some cancellation, i.e.\ we have to consider trees with 
more than eight non-zero leaf labels as well. However, this problem can 
be tackled with the aid of a computer.

Then we might take all values of $V$ thus obtained, try to check whether  
$VTxyUxy=\Delta^3x$ for these values of~$V$, then obtain the list of possible 
values of~$W$ by reflecting the list of possible values of~$V$ as described 
above. These $W$ would satisfy $WTxyUxy=\Delta^3y$. Then we might consider 
all pairs $(V,W)$ from the direct product of these two sets of possible values
and check whether $TVxyWxy=\Delta^3x$ and $UVxyWxy=\Delta^3y$. If this 
is true for some pair $(V,W)$, we can put $B:=(V,W)\in M(2,2;\Lambda)$ 
and conclude that $(DET'_2)$ holds for this matrix with $N=3$, and  
that universal~$N_2=3$. On the other hand, if these conditions would fail 
for all $(V,W)$, we might consider larger values of~$N$ in a similar manner.

Unfortunately, this procedure seems to be quite complicated to implement, 
on a computer or without it, since we don't have a reasonable way of checking 
whether two terms define the same element of $\Lambda(r)$ or not. This is 
the reason why we haven't managed to prove or disprove ``universal $N_2=3$'' 
so far.

Another interesting remark: one might try to use another procedure to find 
all possible values of~$V$. Namely, since we must have $VTxyUxy=\Delta^3x$, 
we might start by listing all complete binary tree of depth~$6$ corresponding 
to $\Delta^3\{1\}$, such that their lower level consists only of 
expressions like $T12$, $T-1-2$, $U12$, $U-1-2$, $T00$ and $U00$, and 
replace these expressions by $1$, $-1$, $2$, $-2$, $0$ and $0$, respectively; 
this would give us all candidates for~$V$, at least among complete 
binary trees of depth~$5$.

\nxsubpoint ($(DET^*_2)$ for torsion-free alternating algebras.)
Notice that we have just proved that $(DET'_2)$ is fulfilled for the 
matrix $A'$ over $\Lambda'$ for $N=3$, since we have shown that 
$V'$ and $W'$ lie in $\Lambda'(2)$, and not just in $\Lambda_\bbZ(2)$.
This result is sufficient by itself to make some interesting conclusions.

Namely, we say that an $\Fpm$-algebra~$K$ is {\em torsion-free\/} 
if the canonical homomorphism $K\mapsto K_\bbZ:=\bbZ\otimes_{\Fpm}K$ 
is injective, i.e.\ a monomorphism. Clearly, an $\Fpm$-algebra~$K$ is
torsion-free iff $K$ can be embedded into a (classical) $\bbZ$-algebra. 
In particular, $\Lambda'\subset\Lambda_\bbZ$ is torsion-free, and any 
$\Fpm$-algebra~$K$ admits a largest torsion-free quotient, namely, the 
image of $K\to K_\bbZ$. Another observation: any commutative torsion-free 
$\Fpm$-algebra is alternating, being a $\Fpm$-subalgebra of a classical 
commutative $\bbZ$-algebra~$K_\bbZ$.

Notice that any homomorphism $f:\Lambda=\Fpm[T^{[2]},U^{[2]}]/_{alt}\to K$ 
with a torsion-free $K$ factorizes through $\Lambda\to\Lambda'$; this 
follows from the commutativity of the following diagram together with the 
fact that $\Lambda'$ is the image of $\Lambda\to\Lambda_\bbZ$:
\begin{equation}
\xymatrix{
\Lambda\ar[r]\ar[d]^{f}&\Lambda_\bbZ\ar[d]^{f_\bbZ}\\
K\ar[r]&K_\bbZ}
\end{equation}

This implies that $\Lambda'$ and $A'\in M(2,2;\Lambda')$ represent the functor 
$K\mapsto M(2,2;K)$ on the category of {\em torsion-free\/} (alternating) 
$\Fpm$-algebras. Since we have shown that $(DET'_2)$ holds for $A'$ with 
$N=3$, we can conclude that {\em $(DET^*_2)$ is fulfilled for all 
(alternating) torsion-free $\Fpm$-algebras with $N=3$.}

\nxsubpoint (Torsion-free algebras over~$\Zinfty$ and~$\barZinfty$.)
Notice that $\Zinfty$ and $\barZinfty$ are torsion-free $\Fpm$-algebras, 
since they can be embedded into~$\bbC$. Another simple statement is 
that {\em any unary polynomial algebra~$A$ over a torsion-free 
$\Fpm$-algebra~$K$ is itself torsion-free.} Indeed, a filtered inductive limit 
argument reduces everything to the case of $A=K[T_1^{[1]},\ldots, T_n^{[1]}]$ 
(with finite~$n$), and then an induction argument shows that it is sufficient 
to treat the case $A=K[T^{[1]}]$. Then $K$ is torsion-free, hence 
$K\to K_\bbZ$ is a monomorphism, i.e.\ all $K(n)\to K_\bbZ(n)$ are injective; 
taking filtered inductive limits we see that $K(S)\to K_\bbZ(S)$ is injective 
for infinite sets~$S$ as well. Now consider the map $A(n)\to A_\bbZ(n)$. 
We know that $A(n)\cong K(\bbZ_{\geq0}\times\stn)$, and similarly 
$A_\bbZ(n)=(K_\bbZ[T])(n)\cong K_\bbZ(\bbZ_{\geq0}\times\stn)$. Then 
the map $A(n)\to A_\bbZ(n)$ is identified with the injective map 
$K(\bbZ_{\geq0}\times\stn)\to K_\bbZ(\bbZ_{\geq0}\times\stn)$, hence 
$A\to A_\bbZ$ is a monomorphism, i.e.\ $A$ is torsion-free.

We see that $(DET^*_2)$ holds for any unary polynomial algebra over~$\Zinfty$,
with $N=3$. Since this property is stable under strict quotients, we see that 
{\em $(DET^*_2)$ holds with $N=3$ for all pre-unary alternating $\Zinfty$- 
and $\barZinfty$-algebras, as well as for all strict quotients of torsion-free 
$\Zinfty$-algebras.} For example, $(DET^*_2)$ and $(DET_2)$ hold for 
$\Finfty[T^{[1]}]$.

\nxsubpoint (Open questions.)
Let us list some questions that have naturally arisen so far, and make some 
comments on each of them.
\begin{enumerate}
\item For which values of $r$ the property $(DET^*_r)$ holds over all 
alternating generalized rings? Or, equivalently, for which $r$'s 
$(DET'_r)$ is fulfilled for the universal matrix~$A$ over 
$\Lambda_r=\Fpm[T_1^{[r]},\ldots,T_r^{[r]}]/_{alt}$? What are the corresponding
(universal) values of~$N_r$? Is the property $(DET^*_\infty)$ universally 
true?
\item For which values of $r$ the property $(DET^*_r)$ holds over all 
torsion-free $\Fpm$-algebras? Or, equivalently, what are the values of~$r$, 
for which $(DET'_r)$ is fulfilled for the image of~$A$ in $\Lambda'_r$, 
the image of $\Lambda_r\to\Lambda_{r,\bbZ}$?
\item Is it true that {\em all\/} polynomial algebras over~$\Zinfty$ are 
torsion-free?
\end{enumerate}

\nxsubpoint (Some remarks.)
Notice that a positive answer for 1) 
would establish $(DET^*_\infty)$, the strongest of our 
properties, for all alternating monads, thus making all further questions 
irrelevant. That's why we actually don't believe that 1) is true, or even that 
$(DET^*_2)$ is universally true: this would be too nice to be true.

On the other hand, we think that 2) holds for all values of~$r$, and that 
our proof for $r=2$, $N=3$ given above can be extended to the general case, 
especially if we choose very large values of~$N$. Notice that 2) 
would imply $(DET^*_\infty)$ for all torsion-free $\Fpm$-, $\Zinfty$- and 
$\barZinfty$-algebras and their strict quotients, e.g.\ pre-unary alternating 
algebras over these generalized rings. This is more than enough for the 
applications in Arakelov geometry.

If both 2) and 3) are true, then $(DET^*_\infty)$ would hold for {\em all\/} 
alternating algebras over~$\Zinfty$. This statement would solve all our 
problems, but it seems again too strong to be true. So we don't believe in~3).

\nxsubpoint (Lower bounds for universal $N=N_r$.)
Notice that one can obtain lower bounds for universal $N=N_r$ in 2) and 1), 
essentially by the same reasoning as in the proof of~\eqref{eq:lbound.N2}. 
This yields the following inequality for $N=N_r$ of 2) and 1):
\begin{equation}
\frac{(Nr-1)!}{(N-1)!\cdot N!^{r-1}}\geq r!^N
\end{equation}
One can check that this inequality actually means 
$N\geq3$ for $2\leq r\leq 3$, and $N\geq2$ for $r\geq4$, so it doesn't 
really give us a lot of information. In fact even the inequality 
$N_r\leq N_{r+1}$ of~\ptref{sp:DET.decr.r} gives us a better lower bound
$N_r\geq N_2\geq 3$ for $r\geq2$.
However, we think that proving 2) for $r\geq3$ will be much simpler for 
values of $N$ much larger than any of these lower bounds. 
In fact, one might hope to find a regular way of constructing trees~$V$ 
with $\bar V=V'$ (this is sufficient to prove 2) for given $r$ and~$N$) 
when $N$ is very large.

\nxpointtoc{Complements}
We collect in this subsection some interesting observations that didn't 
fit naturally anywhere in this chapter. Some of them are quite important, 
but require more sophisticated techniques for their study than we have 
already developped (e.g.\ localization and homotopical algebra techniques).


\nxsubpoint\label{sp:tens.sqr.zinfty} 
(Tensor square of $\Zinfty$ and $\Finfty$.)
We have seen in~\ptref{sp:tens.sqr.z} that $\bbZ\otimes_{\Fone}\bbZ=
\bbZ\otimes_{\Fpm}\bbZ=\bbZ$, i.e.\ that $\Fone\to\bbZ$ is an epimorphism 
in the category of generalized rings. What can be said about tensor 
squares of $\Zinfty$, $\Zninfty$ and $\Finfty$? 

\nxsubpoint\label{sp:tens.sqr.finfty} (Tensor square of $\Finfty$.)
Recall that $\Finfty=\Fpm[*^{[2]}\,|\,\bu*\bu=\bu$, $\bu*(-\bu)=0$, 
$x*y=y*x$, $(x*y)*z=x*(y*z)]$ (cf.~\ptref{sp:examp.pres.comalg},e). 
We can use this to compute $\Sigma:=\Finfty\otimes_{\Fpm}\Finfty$. We see that 
$\Sigma$ is generated as a commutative $\Fpm$-algebra by two operations, 
say $*^{[2]}$ and $\vee^{[2]}$, each of them satisfies the conditions listed 
above for~$*$, i.e.\ is commutative, associative, idempotent, and 
$\bu*(-\bu)=\bu\vee(-\bu)=0$. We have to impose one additional condition as 
well, namely, the commutativity between $*$ and~$\vee$: 
\begin{equation}
(x\vee y)*(z\vee w)=(x*z)\vee(y*w)
\end{equation}
Putting here $z=y$, $w=x$ and applying relations $x*x=x=x\vee x$, 
$x*y=y*x$ and $x\vee y=y\vee x$, we obtain
\begin{equation}
x\vee y=(x\vee y)*(y\vee x)=(x*y)\vee(y*x)=x*y
\end{equation}
In other words, $*=\vee$, hence $\Finfty\otimes_{\Fpm}\Finfty=\Finfty$.
(This simple argument has been communicated to me by A.~Smirnov.)

\nxsubpoint\label{sp:tens.sqr.zninfty} (Tensor square of $\Zninfty$.)
The above argument can be easily generalized to prove $\Zninfty\otimes_{\Fpm}
\Zninfty=\Zninfty$. First of all, denote by $s_N\in\Zninfty(N)$, $N\geq1$,
the $N$-th averaging operation:
\begin{equation}
s_N=\frac1N\{1\}+\frac1N\{2\}+\cdots+\frac1N\{N\}
\end{equation}
One sees immediately that $\Zninfty$ is generated over $\Fpm$ by all such
operations $s_N$ (actually $s_{NN'}$ can be easily expressed in terms of
$s_N$ and $s_{N'}$, so $\{s_p\}_{p\in\bbP}$ would suffice).
Indeed, any rational octahedral combination
$t=\sum_{i=1}^n\lambda_i\{i\}=\sum_{i=1}^n (m_i/N)\cdot\{i\}\in\Zninfty(n)$,
where $N$, $m_i\in\bbZ$, $N>0$, and necessarily $\sum_i|m_i|\leq N$,
can be rewritten as $s_N$ applied to the list of arguments containing
$\pm\{i\}=\sgn m_i\cdot\{i\}$ exactly $|m_i|$ times, for all $1\leq i\leq n$,
and with the remaining $N-\sum_i|m_i|$ arguments put equal to zero:
\begin{equation}
t=\sum_{i=1}^n\frac{m_i}N\{i\}=s_N\bigl(
\overbrace{\pm\{1\},\ldots,\pm\{1\}}^{\text{$|m_1|$ times}},
\overbrace{\pm\{2\},\ldots}^{\text{$|m_2|$ times}},
\ldots,0,\ldots,0\bigr)
\end{equation}
Now $\Sigma:=\Zninfty\otimes_{\Fpm}\Zninfty$ is generated over~$\Fpm$ by two 
commuting sets of such averaging operations $s_N$, $s'_N$, $N\geq 1$,
and we are reduced to showing $s_N=s'_N$. Notice that
$s_N$ is invariant under any permutation of arguments, and that
$s_N(x,x,\ldots,x)=x$, and similarly for $s'_N$; applying
commutativity relation of~\ptref{def:comm} 
to $s_N$, $s'_N$, and $N\times N$-matrix
$x=(x_{ij})_{1\leq i,j\leq N}$, 
where $x_{ij}=\{((i+j-2)\bmod N)+1\}\in\Sigma(N)$
(i.e.\ $x$ is a Latin square: any row or column contains each $\{i\}$,
$1\leq i\leq N$, exactly once), we obtain $x_{.j}=
s_N(\{j\},\{j+1\}\ldots,\{j-1\})=s_N(\{1\},\ldots,\{N\})=s_N$, whence
$x'=s'_N(x_{.1},\ldots,x_{.N})=s'_N(s_N,\ldots,s_N)=s_N$;
a similar argument shows $x_{i.}=s'_N$ for any~$i$ and
$x''=s_N(x_{1.},\ldots,x_{N.})=s'_N$, so the commutativity relation
$x'=x''$ yields $s_N=s'_N$ as claimed.

\nxsubpoint\label{sp:tens.sqr.zninfty.fone} 
(Tensor square of $\Zinfty$, and tensor squares over $\Fone$.)
Notice that the above argument doesn't imply 
$\Zinfty\otimes_{\Fpm}\Zinfty=\Zinfty$,
since $\Zinfty$ is not generated over $\Fpm$ by the~$s_N$, and we actually
expect $\Zinfty\otimes_{\Fpm}\Zinfty\neq\Zinfty$.

It also leaves the question about tensor squares of $\Finfty$ 
and $\Zninfty$ over $\Fone$ open. Notice that 
$\Zninfty$ is generated over $\Fone$ by the symmetry $-^{[1]}$ subject to
$[-]^2=\bu$, and by the averaging operations $s_N$, and the previous argument
still shows $s_N=s'_N$ inside $\Sigma:=\Zninfty\otimes_{\Fone}\Zninfty$
since it never used the symmetry. Therefore, $\Sigma$ is generated over
$\Zninfty$ by an additional symmetry $\ominus^{[1]}$, commuting with
$[-]$ and all $s_N$, i.e.\ $\Zninfty\otimes_{\Fone}\Zninfty=
\Zninfty[\ominus^{[1]}]$. Notice, however, that apart from relation
$\ominus^2=\bu$, this new symmetry satisfies relations like
$s_2(x,\ominus x)=0$ since they are fulfilled by~$-$.
Similarly, $\Finfty\otimes_{\Fone}\Finfty=\Finfty[\ominus^{[1]}]$,
subject to relations $\ominus^2=\bu$ and $x*\ominus x=0$.
One can check that there are $\Finfty$-modules (hence also $\Zninfty$-modules)
$M$ with an involution $\ominus$ different from $-$: 
one can consider for example $M=\{0,a,b,-a,-b\}$ with $*$ given by
$x*x=x$, $x*y=0$ for $y\neq x$, and put $\ominus a:=b$, $\ominus b:=a$.
This implies $\Finfty\otimes_{\Fone}\Finfty\neq\Finfty$
and $\Zninfty\otimes_{\Fone}\Zninfty\neq\Zninfty$.

\nxsubpoint\label{sp:tens.z.zinfty}
Now let's show that $\bbZ\otimes_{\Fpm}\Zinfty=\bbR$ and 
$\bbZ\otimes_{\Fpm}\Zninfty=\bbQ$. To do this we use the localization theory 
developed in the next chapter, and in particular relations 
$\Zinfty[(1/2)^{-1}]\cong\bbR$ and $\Zninfty[(1/2)^{-1}]\cong\bbQ$; they 
will be also checked in the next chapter (cf.~\ptref{sp:examp.loc}). 
Now notice that the operation 
$1/2$ of~$\Zninfty$ becomes invertible in $\bbZ\otimes_{\Fpm}\Zninfty$: 
in fact, it admits an inverse $2:=\bu*\bu$, where we have temporarily denoted 
the addition of~$\bbZ$ by~$*$ to prevent confusion with octahedral linear 
combinations $\lambda_1\{1\}+\cdots+\lambda_n\{n\}$ of~$\Zninfty$.

Let us check that $2=\bu*\bu$ is indeed an inverse to $1/2$. We use for this 
the commutativity relation between $*$ and $s_2:=(1/2)\{1\}+(1/2)\{2\}$. 
We get $((1/2)x+(1/2)y)*((1/2)z+(1/2)w)=(1/2)(x*z)+(1/2)(y*w)$. Putting 
$y=z=0$, $x=w$, we obtain $2\cdot(1/2)x=((1/2)x)*((1/2)x)=(1/2)x+(1/2)x=x$, 
i.e.\ $2\cdot(1/2)=\bu$ in $\bbZ\otimes_{\Fpm}\Zninfty$ as claimed.

Once we know that $1/2$ becomes invertible in the tensor product, the rest 
follows immediately from universal properties of tensor products and 
localizations: $\bbZ\otimes_{\Fpm}\Zinfty=(\bbZ\otimes_{\Fpm}\Zinfty)
[(1/2)^{-1}]=\bbZ\otimes_{\Fpm}\Zinfty[(1/2)^{-1}]=\bbZ\otimes_{\Fpm}\bbR=
\bbZ\otimes_{\Fpm}\bbZ\otimes_\bbZ\bbR=\bbZ\otimes_\bbZ\bbR=\bbR$, since 
$\bbZ\otimes_{\Fpm}\bbZ=\bbZ$ by~\ptref{sp:tens.sqr.z}. Formula 
$\bbZ\otimes_{\Fpm}\Zninfty=\bbQ$ is deduced similarly from  
$\Zninfty[(1/2)^{-1}]=\bbQ$.

\nxsubpoint\label{sp:genfd} (Generalized fields.)
We would like to have a reasonable definition of generalized fields. 
Clearly, a generalized field is a generalized ring~$K$ with some additional 
properties. Our definition should coincide with the usual one for classical 
rings (i.e.\ a classical ring is a generalized field iff it is a classical 
field), and it shouldn't be too restrictive, otherwise we wouldn't be 
able to obtain a reasonable definition of ``residue fields'' of points of 
generalized schemes, studied in the next chapter.

In particular, we would like $\Finfty$ to be a generalized field, since it 
is the only reasonable candidate for the residue field of~$\Zinfty$. 
This rules out a potential definition of a generalized field as 
a generalized ring~$K$, such that all $K$-modules are free. Indeed, the 
subset $M:=\{0,\pm?\{1\}\pm?\{2\}\}\subset\Finfty(2)$ is a 5-element 
$\Finfty$-module, clearly not free, since $\card\Finfty(n)=3^n$ according 
to~\ptref{sp:examp.finfty}.

So we have to consider some other characteristic property of classical fields. 
For example, they have no non-trivial quotients. This leads to the following 
definition:

\begin{DefD}\label{def:gen.field}
We say that a generalized ring~$K$ is {\em subfinal} or 
{\em subtrivial} if it is isomorphic to a submonad of the final monad~$\st1$, 
i.e.\ to $\st1_+$ or $\st1$. According to~\ptref{sp:prop.triv.mon}, 
$K$ is subtrivial iff $\{1\}_\st2=\{2\}_\st2$.

We say that a generalized ring~$K$ is a {\em generalized field} if it is 
not subtrivial, but all its strict quotients different from~$K$ are 
subtrivial.
\end{DefD}

When~$K$ is a generalized ring with zero, $K$ is (sub)trivial iff 
$\bu=0$ in~$|K|$. Therefore, if $K$ is a generalized field with zero, we 
have $\bu\neq0$, i.e.\ $|K|=K(1)$ consists of at least two elements. This 
implies that the canonical homomorphism $K\to\END(|K|)$ is non-trivial, 
since~$0$ and $\bu$ act on~$|K|$ in different way, hence it has to be 
a monomorphism (otherwise the image of $K$ in $\END(|K|)$ would be a 
non-trivial strict quotient of~$K$), i.e.\ {\em if $K$ is a generalized field 
with zero, $|K|$ is an exact $K$-module.} More precisely, $|K|$ is an 
exact $K$-module whenever it consists of at least two elements; and $K(2)$ 
is an exact $K$-module in all cases. 

Another easy observation --- {\em $|K|$ doesn't admit any non-trivial strict 
quotient $A$ in $\catMod K$.} Indeed, if~$A$ is such a strict quotient, 
it is a quotient-monoid of~$|K|$ as well, hence we get an algebra structure 
on~$A$, and the corresponding unary algebra $K_A$ over~$K$ will be a 
non-trivial strict quotient of~$K$.

Conversely, {\em if $K$ is a generalized ring with zero, such that 
$|K|$ is non-trivial, exact, and doesn't admit any non-trivial strict quotients
in $\catMod K$, then $K$ is a generalized field.} Indeed, let 
$f:K\to K'$ be a strict epimorphism. Then $|K'|$ is a strict quotient of 
$|K|$ in $\catMod K$. If $|K'|=0$, $\bu=0$ in $|K'|$, and~$K'$ is itself 
trivial. On the other hand, if $|K'|\neq 0$, then necessarily $|K'|=|K|$, 
and the composite map $K\to K'\to\END(|K'|)=\END(|K|)$ is a monomorphism, 
$|K|$ being an exact $K$-module, hence $f:K\to K'$ is a monomorphism as well, 
hence $K'=K$. 

An immediate consequence is that {\em a classical ring is a generalized field 
iff it is a field in the usual sense}, and that {\em $\Finfty$ is a 
generalized field}.

\nxsubpoint (Examples of generalized fields.)
\begin{itemize}
\item We have seen that all classical fields, e.g.\ $\bbF_p$, $\bbQ$, 
$\bbR$ and $\bbC$, are generalized fields.
\item We have just seen that $\Finfty$ is a generalized field. This example 
shows that in general there can be modules over a generalized field that are 
not free.
\item It is easy to check with the aid of the criterion given above that 
$\Fempty$ and $\Fone$ are generalized fields, while $\Fpm$ and $\bbF_{1^n}$ 
for $n>1$ are not, since they admit $\Fone$ as a strict quotient. In fact, 
$\bbF_{1^n}$ is something like a local artinian ring with 
residue field~$\Fone$.
\end{itemize}

Notice that $|\bbF_{1^n}|$ is an exact $\bbF_{1^n}$-module, and all non-zero 
elements of $|\bbF_{1^n}|$ are invertible, but this is insufficient to 
conclude that $\bbF_{1^n}$ is a generalized field.

\nxsubpoint (Further properties of generalized fields.)
Applying Zorn lemma to the set of compatible equivalence relations $\equiv$ 
on a generalized ring~$K$, such that $\{1\}_2\not\equiv\{2\}_2$, we see that 
{\em any non-subtrivial generalized ring~$K$ admits a strict quotient that 
is a generalized field.} Another interesting property is obtained by 
inspecting the image of $|K|\stackrel s\to |K|$, $s\in|K|$, and 
the image of $K\to K[s^{-1}]$ as well: we see that 
{\em if~$K$ is a generalized field with zero, and $s\in|K|$, $s\neq0$, 
then the maps $[s]_{K(n)}:K(n)\to K(n)$ are injective for all $n\geq0$.} 
We don't know whether all non-zero elements of $|K|$ are invertible or not.

\nxsubpoint (Elementary $K$-theory.)
We can construct $K_0(\Sigma)$ for any generalized ring~$\Sigma$ by 
considering the quotient of the free abelian group generated by 
projective $\Sigma$-modules of finite type modulo relations 
$[M\oplus N]-[M]-[N]$. Clearly, the pullback functor defines a map 
$\rho^*:K_0(\Sigma)\to K_0(\Sigma')$ for any $\rho:\Sigma'\to\Sigma$.

We can try to construct $K_1(\Sigma)$ in the usual way, by considering the 
free abelian group generated by couples $[P,\phi]$, 
$P$ projective of finite type, 
$\phi\in\Aut_\Sigma(P)$, modulo certain relations. However, some of these 
relations have to involve something like ``short exact sequences''; 
a similar problem arises when we try to compute $K_0$ of larger categories 
of $\Sigma$-modules.

Probably the best way to deal with these problems is to construct a reasonable 
(triangulated) category of perfect complexes over~$\Sigma$. For now we'll 
content ourselves by indicating that a bicartesian square
\begin{equation}
\xymatrix{M\ar[r]\ar[d]&N\ar[d]\\M'\ar[r]&N'}
\end{equation}
is very much like a short exact sequence when we work over a classical ring, 
since in this case such a bicartesian square is essentially just a 
short exact sequence $0\to M\to M'\oplus N\to N'\to 0$, hence it yields 
a relation $[M]-[M']+[N]-[N']=0$ in~$K_0$. We can try to construct something 
similar in the general case, but it seems that we have to impose some 
additional conditions on such squares (e.g.\ require all horizontal
arrows to be constant perfect cofibrations in the sense of
Chapter~\ptref{chap:K0.Chow}) to obtain a reasonable theory.

We'll return to the algebraic $K$-theory for generalized rings
later in Chapter~\ptref{chap:K0.Chow}, where we'll apply a modification of
Waldhausen construction to define all algebraic $K$-groups of a generalized
ring~$\Sigma$.

\nxsubpoint (Theory of traces.)
We can construct an ``abstract theory of characteristic polynomials'' 
by considering the group $\tilde K_1(\Sigma)$, generated by pairs 
$[P,\phi]$ with $\phi\in\End_\Sigma(P)$ (now $\phi$ need not be invertible!) 
modulo relations similar to that of~$K_0$. We can also construct an 
``abstract theory of traces'' $K_1^+(\Sigma)$ by considering the 
$(\Sigma\otimes\bbZ)$-module generated by couples $[P,\phi]$ as above, 
subject to an additional family of relations, namely,
$$[P,t(\phi_1,\ldots,\phi_n)]=t\bigl([P,\phi_1],\ldots,[P,\phi_n])
\text{ for $t\in \Sigma(n)$.}$$

\nxsubpoint (A simpler construction of traces.)
However, we can construct a simpler theory of traces over a generalized 
ring~$\Sigma$. Namely, we can put $\tr\phi:=\tr\rho^*(\phi)$, where 
$\phi:\Sigma\to\bbZ\otimes\Sigma$ is the canonical map (the tensor product 
is taken here over~$\Fempty$, $\Fone$ or~$\Fpm$ in the most appropriate way). 
This gives us a theory of traces of matrices over~$\Sigma$ with values 
in $\Sigma\otimes\bbZ$, and this theory is quite sensible: for example, 
matrices over~$\Zinfty$ have traces in $\Zinfty\otimes_{\Fpm}\bbZ=\bbR$, 
and matrices over~$\bbF_{1^n}$ -- in $\bbZ[\sqrt[n]1]$. However, 
$\Finfty\otimes_{\Fpm}\bbZ=0$, so we don't obtain any theory of traces 
over~$\Finfty$ in this way.

\nxsubpoint\label{sp:val.rings} (Valuation rings.)
Recall that in~\ptref{sp:stab.of.subset} we have defined a submonad 
$(N:N)\subset\END(M)$ for any set $M$ and its subset $N$; 
$(N:N)(n)$ consists of all maps $f:M^n\to M$, such that $f(N^n)\subset N$. 
When $M$ is a $\Sigma$-module, we get a canonical homomorphism 
$\rho:\Sigma\to\END(M)$. Let us denote the preimage of $(N:N)$ under~$\rho$ 
by~$\Sigma_N$; clearly, $\Sigma_N$ is the largest algebraic submonad 
of~$\Sigma$, such that $N$ is a $\Sigma_N$-module.

Now suppose that $K$ is a classical field, and $\fnorm_v$ is a valuation 
on~$K$ (archimedian or not). Let's apply the above construction 
to $\Sigma=K$, $M=K$ (considered as a module over itself) and 
$N=N_v:=\{x\in K : |x|_v\leq 1\}$. We obtain an algebraic submonad 
$\cO_v:=K_{N_v}\subset K$; we say that {\em $\cO_v$ is the valuation ring 
of~$\fnorm_v$.} Clearly, $\cO_v$ is a classical ring iff 
$\fnorm_v$ is non-archimedian, and in this case it coincides with 
the classical valuation ring of $\fnorm_v$. In any case, 
the underlying set $|\cO_v|$ of~$\cO_v$ is easily seen to be equal to~$N_v$, 
and $\fnorm_v$ is completely determined by~$N_v$, hence also by~$\cO_v$,
up to equivalence. It is immediate that $\cO_v$ is a hypoadditive 
alternating $\Fpm$-subalgebra of~$K$.

Notice that the valuation rings of archimedian valuations on~$\bbR$, $\bbC$
and~$\bbQ$ are nothing else than $\Zinfty$, $\barZinfty$ and~$\Zninfty$. 
In general, it follows from the definitions that 
$\cO_v(n)$ consists of all $n$-tuples 
$\lambda=(\lambda_1,\ldots,\lambda_n)\in K^n$, such that 
$|\lambda_1x_1+\cdots+\lambda_nx_n|_v\leq1$ whenever all $|x_i|_v\leq1$. 
If $\fnorm_v$ is non-archimedian, this condition is equivalent to 
``all $|\lambda_i|_v\leq1$'', i.e.\ $\lambda\in N_v^n$. In the 
archimedian case one checks that this condition is equivalent to 
$|\lambda_1|_v+\cdots+|\lambda_n|_v\leq1$
(the verification is based on choosing rational numbers $y_i$ with 
$|\lambda_i|_v(1-\epsilon)\leq y_i\leq |\lambda_i|_v$ and considering 
$x_i:=y_i/\lambda_i$). Clearly, this is a generalization 
of the ``octahedral combinations'' used to construct $\Zinfty$, 
$\barZinfty$ and~$\Zninfty$.

It would be nice to have an intrinsic definition of valuation rings inside 
a classical field~$K$, and to transfer this definition to the case of 
generalized fields.


\cleardoublepage

\mysection{Localization, spectra and schemes}

Now we are going to develop a reasonable theory of spectra for generalized 
rings. More precisely, we develop {\em two\/} such theories, as two special 
cases of a general construction of theories of spectra. 
We construct the {\em unary} or {\em prime spectrum\/} 
$\Spec^u A$ or $\Spec^p A$ of a generalized ring~$A$. It admits a description 
in terms of prime ideals very similar to the classical case, and its basic 
properties are essentially the same.
However, this spectrum fails to have some more sophisticated properties,
e.g.\ finitely generated projective $A$-modules do not necessarily define 
locally free quasicoherent sheaves, and the image of a flat finitely presented 
morphism of such spectra is not necessarily open.

We deal with these problems by defining a more sophisticated theory of spectra,
namely, the {\em total spectra\/} $\Spec^t A$ or simply $\Spec A$ of 
generalized rings. We have a comparison map $\Spec^t A\to\Spec^p A$, that 
turns out to be an isomorphism of ringed spaces for a classical~$A$, but 
not in general.

Of course, once we construct the (total) spectrum of a generalized ring as 
a generalized (locally) ringed space or topos, we can immediately define  
the category of generalized schemes. We want to study their basic properties 
and discuss some examples, e.g.\ the projective spaces $\bbP^n_{\Fone}$
over~$\Fone$. 

\nxpointtoc{Unary localization}
The main definitions can be transferred verbatim from the classical case, 
since they depend (almost) only on the multiplicative structure of~$|A|$. 
However, we must be careful while dealing with the extra structures.

\begin{DefD} 
a) We say that a subset $S$ of a monoid~$M$ is a
{\em multiplicative system\/} if it is a submonoid of~$M$, i.e.\ 
if it is closed under multiplication and contains the identity of~$M$.

b) For any subset $S\subset M$ we denote by $\langle S\rangle$ 
the multiplicative system (i.e.\ the submonoid) generated by~$S$. 
If $S=\{f\}$, we write $S_f$ instead of $\langle\{f\}\rangle=
\{1,f,f^2,\ldots,f^n,\ldots\}$.

c) We say that a multiplicative system $S\subset M$ is 
{\em saturated\/} if $ax\in S$, $a\in M$ implies $x\in S$, i.e.\ if 
$S$ contains all divisors of any its element. The smallest saturated 
multiplicative system~$\tilde S$ containing a given multiplicative 
system~$S$ (i.e.\ the set of all divisors of all elements of~$S$) 
is called the {\em saturation} of~$S$.
\end{DefD}

\begin{DefD}\label{def:localiz}
a) Given a generalized ring $A$ and any subset~$S\subset|A|$, we denote by 
$(A[S^{-1}],i_A^S)$ the universal (i.e.\ initial) object in the category 
of couples $(B,\rho)$, where $B$ is a generalized ring, and $\rho:A\to B$ 
a homomorphism, such that all elements of $\rho_1(S)\subset|B|$ are 
invertible in $|B|$ with respect to the canonical commutative monoid 
structure on this set. We say that $A[S^{-1}]$ is the {\em localization 
of~$A$ with respect to~$S$}. When $S$ is a multiplicative system, we  
write $S^{-1}A$ instead of $A[S^{-1}]$ as well, and we write 
$A_f$ or $A[f^{-1}]$ instead of $S_f^{-1}A$ or~$A[\{f\}^{-1}]$.

b) Given a generalized ring~$A$, a subset~$S\subset|A|$, and an
$A$-module~$M$, we denote by $(M[S^{-1}],i_M^S)$ the universal (initial) 
object in the category of couples $(N,f)$, with $f:M\to N$ an $A$-module 
homomorphism, such that all elements $s\in S$ act bijectively on~$N$, i.e.\ 
such that all elements of~$S$ become invertible after the application 
of the canonical map $A\to\END(N)$.
We say that $M[S^{-1}]$ is the {\em localization of~$M$ with respect to~$S$}; 
when $S$ is a multiplicative system, we write $S^{-1}M$ instead of~$S$, 
and we write $M_f$ or $M[f^{-1}]$ instead of~$S_f^{-1}M$.
\end{DefD}

\nxsubpoint (Reduction to the multiplicative system case.)
Let us denote by $\Inv(B)$ the set of invertible elements of monoid~$|B|$, 
where $B$ is any algebraic monad. 
Notice that this is a multiplicative system in~$|B|$, hence its preimage 
$\rho^{-1}(\Inv(B))\subset|A|$ is a multiplicative system as well. 
We conclude that for any subset $S\subset|A|$ the condition 
$\rho(S)\subset\Inv(B)$, i.e.\ $S\subset\rho^{-1}(\Inv(B))$, is equivalent 
to $\rho(\langle S\rangle)\subset\Inv(B)$. This means that 
$A[S^{-1}]$ and $A[\langle S\rangle^{-1}]=\langle S\rangle^{-1}A$ solve 
the same universal problem, hence are isomorphic (when they are representable).
In other words, it suffices to show the existence of localizations with 
respect to multiplicative subsets $S\subset|A|$. 

Applying the above observation to $B:=\END_A(N)$ (this is a 
non-commu\-ta\-tive 
$A$-algebra, but we didn't use commutativity of~$B$ so far), we see that 
$S\subset|A|$ acts by bijections on an $A$-module~$N$ iff $\langle S\rangle$ 
acts by bijections on~$N$, hence $M[S^{-1}]=\langle S\rangle^{-1}N$, so 
we can restrict ourselves to the case when $S$ is a multiplicative system 
in the module case as well.

\nxsubpoint\label{sp:constr.localiz.mod} (Construction of $S^{-1}M$.)
Let $S\subset|A|$ be a multiplicative system and $M$ be an $A$-module. 
We construct the localization $S^{-1}M$ in the usual way. Namely, 
we consider the equivalence relation $\sim$ on the set $M\times S$, 
given by $(x,s)\sim(y,t)$ iff there is a $u\in S$, such that 
$ut\cdot x=us\cdot y$, and put $S^{-1}M:=(M\times S)/\sim$. One checks 
in the usual way that $\sim$ is indeed an equivalence relation on~$M\times S$; 
we denote by $x/s$ the class of $(x,s)$ in $S^{-1}M$ as usual.

We have constructed a set $S^{-1}M$ so far. Notice that this construction 
used only the monoid structure of~$|A|$ and the action of~$|A|$ on~$|M|$. 
Now consider small category $\cS$, defined as follows: $\Ob\cS=S$, 
i.e.\ the objects $[s]$ of $\cS$ are in one-to-one correspondence with 
elements $s\in S$, and $\Hom_\cS(s,t):=\{u\in S\,|\,su=t\}$; the composition 
of morphisms is given by the multiplication of~$S$.

We can use our $A$-module~$M$ to define a functor $\tilde M:\cS\to\catMod A$ 
by putting $\tilde M([s]):=M$ for all $s\in S$, and 
$\tilde M(u):=u_M:M\to M$ for any morphism $u:[s]\to[t]$ in~$\cS$ 
(notice that $u_M$ is an $A$-endomorphism of~$M$ because of the commutativity
of~$A$). Now it is immediate that $S^{-1}M=\injlim_\cS\tilde M$ in the 
category of sets. However, $\cS$ is obviously a filtered category, 
and filtered inductive limits in~$\catMod A$ and $\catSets$ coincide
(cf.~\ptref{sp:filtindlim.mod}); hence we obtain a natural 
$A$-module structure on the set $S^{-1}M$ by putting 
$S^{-1}M:=\injlim_\cS\tilde M$ in~$\catMod A$.

This $A$-module structure on~$S^{-1}M$ can be written down explicitly: 
for any operation $a\in A(n)$ and any elements $x_i/s\in S^{-1}M$ with 
common denominator~$s\in S$ we have
\begin{equation}
{[a]}_{S^{-1}M}(x_1/s,\ldots,x_n/s)=\bigl({[a]}_M(x_1,\ldots,x_n)\bigr)/s
\end{equation}

In particular, an element $t\in S$ acts on $S^{-1}M$ by mapping 
$x/s$ into $tx/s$; clearly, $x/s\mapsto x/st$ is the inverse map, i.e.\ 
$S$ acts bijectively on~$S^{-1}M$. We have a canonical $A$-linear map 
$i_M^S:M\to S^{-1}M$, $x\mapsto x/1$ as well; it is essentially given 
by the embedding of $M_{[1]}$ into $\injlim_\cS M_{[s]}$.

We have to show that for any $A$-module~$N$, such that $S$ acts bijectively 
on~$N$, we have $\Hom_A(S^{-1}M,N)\cong\Hom_A(M,N)$. This is quite clear: 
$\Hom_A(S^{-1}M,N)=\Hom_A(\injlim_\cS M_{[s]},N)\cong
\projlim_\cS\Hom(M_{[s]},N)$ consists of families $(f_s)_{s\in S}$, 
$f_s:M\to N$, such that $f_{st}=f_s\circ[t]_M$ for any $s$, $t\in S$; 
this is equivalent to $f_{st}=[t]_N\circ f_s$, and in particular 
$f_s=[s]_N^{-1}\circ f_1$, i.e.\ such families are in one-to-one correspondence
with $A$-homomorphisms $f_1:M\to N$, q.e.d.

\nxsubpoint\label{sp:ex.un.loc} (Existence and unarity of $A[S^{-1}]$.)
Notice that for any subset $S\subset|A|$ the localization 
$A[S^{-1}]$ can be constructed by considering the commutative $A$-algebra, 
generated by new unary operations $s^{-1}$, $s\in S$, subject to unary 
relations $s\cdot s^{-1}=\bu$:
\begin{equation}\label{eq:comm.expr.loc}
A[S^{-1}]=A[(s^{-1})_{s\in S}\,|\,s\cdot s^{-1}=\bu]
\end{equation}
This implies that {\em $A[S^{-1}]$ is a unary $A$-algebra}. Furthermore, 
we can define $A[S^{-1}]$ by a similar formula without imposing any 
implicit commutativity relations between the new generators and operations 
of~$A$:
\begin{equation}\label{eq:noncomm.expr.loc}
A[S^{-1}]=A\langle(s^{-1})_{s\in S}\,|\,s\cdot s^{-1}=\bu=s^{-1}\cdot s\rangle
\end{equation}
Indeed, we just have to check that the algebraic monad $B$ defined by the 
RHS is commutative, i.e.\ that the generators $s^{-1}$ commute 
between themselves and the operations of~$A$. The first statement is clear, 
since $s^{-1}t^{-1}=(ts)^{-1}=(st)^{-1}=t^{-1}s^{-1}$ in~$|B|$. 
Let's show that $s^{-1}$ commutes with all operations of~$A$. According 
to~\ptref{sp:descr.centr.el}, 
this is equivalent to ${[s^{-1}]}_{B(n)}:B(n)\to B(n)$ 
being an $A$-homomorphism for all $n\geq0$. But ${[s^{-1}]}_{B(n)}$ 
is the inverse of the automorphism ${[s]}_{B(n)}$ of~$B(n)$, and the 
latter is an $A$-homomorphism, since $s$ commutes with all elements of~$A$.

An immediate consequence of \eqref{eq:noncomm.expr.loc} is that 
{\em $i_A^S:A\to A[S^{-1}]$ is universal among all algebraic monad 
homomorphisms $\rho:A\to B$, such that $\rho(S)\subset\Inv(B)$,} 
i.e.\ $A[S^{-1}]$ is the ``non-commutative localization'' of $A$ 
with respect to~$S$ as well.

Another immediate consequence is that $i_A^S:A\to A[S^{-1}]$ is an 
NC-epimorphism (cf.~\ptref{sp:nc.epi}), i.e.\ any monad homomorphism 
$A\to B$ extends to at most one monad homomorphism $A[S^{-1}]\to B$, 
or equivalently, any $A$-module admits at most one $A[S^{-1}]$-module 
structure.

\nxsubpoint ($A[S^{-1}]$-modules.)
We see that algebraic monad homomorphisms $A[S^{-1}]\to B$ are in one-to-one 
correspondence with algebraic monad homomorphisms $\rho:A\to B$, such that 
all elements of $\rho(S)\subset|B|$ are invertible. Applying this to 
$B=\END(N)$ we see that {\em $\catMod{A[S^{-1}]}$ is 
isomorphic to the full subcategory of $\catMod A$, consisting of those 
$A$-modules~$N$, on which all elements $s\in S$ act by bijections.} 
Clearly, this isomorphism of categories is induced by the scalar 
restriction functor $(i_A^S)_*:\catMod{A[S^{-1}]}\to\catMod A$; 
in particular, {\em the scalar restriction functor $(i_A^S)^*$ is fully 
faithful.} Notice that it is also exact, and admits both a left and 
a right adjoint, since~$A[S^{-1}]$ is unary over~$A$
(cf.~\ptref{prop:cond.unarity}).

In this way we can identify $(i_A^S)_*$ with the natural embedding of 
the full subcategory of $\catMod A$, consisting of $A$-modules, on which 
$S$ acts bijectively, into $\catMod A$ itself. The left adjoint to 
this functor is clearly the functor $M\mapsto M[S^{-1}]$ 
of~\ptref{def:localiz}, i.e.\ 
{\em $(i_A^S)^*M=A[S^{-1}]\otimes_A M$ is canonically isomorphic 
to~$M[S^{-1}]$.} Another immediate consequence is that {\em any $M[S^{-1}]$ 
admits a canonical $A[S^{-1}]$-structure.}

\nxsubpoint\label{sp:class.descr.loc} (Description of $S^{-1}A$.)
Now we can combine our results to obtain an explicit description 
of $S^{-1}A$ when $S\subset|A|$ is a multiplicative system. Indeed, 
$(S^{-1}A)(n)=(i_A^S)^*(A(n))=S^{-1}A(n)$, hence according 
to~\ptref{sp:constr.localiz.mod} the set $(S^{-1}A)(n)$ of $n$-ary operations 
of~$S^{-1}A$ consists of fractions $a/s$, where $a\in A(n)$, $s\in S$, 
and $a/s=b/t$ iff $tu\cdot a=su\cdot b$ for some $u\in S$.

We have already described the $A$-module structure on each~$S^{-1}M$, 
hence also on~$S^{-1}A(m)$; 
using the fact that all $s\in S$ act on these modules by bijections, we 
obtain
\begin{multline}
{[a/s]}_{S^{-1}M}(x_1/t,\ldots,x_n/t)=\bigl([a]_M(x_1,\ldots,x_n)\bigr)\big/st
\\\text{ for any $a\in A(n)$, $x_i\in M$, $s$, $t\in S$}
\end{multline}
Putting here $M=A(m)$ we obtain the composition maps for the 
operations of~$S^{-1}A$. Of course, when we want to apply $a/s$ 
to some elements $x_i/t_i$ of~$S^{-1}M$ with different denominators, we have 
to reduce them to the common denominator, e.g.\ $t=t_1t_2\cdots t_n$. 
For example, for a binary operation $*^{[2]}$ we get
\begin{equation}
x_1/t_1*x_2/t_2=(t_2x_1*t_1x_2)/t_1t_2
\end{equation}
When $*=+$, we obtain the usual formula for the sum of two fractions.

\nxsubpoint (Flatness of localizations.)
Let's prove that {\em $A[S^{-1}]$ is a flat $A$-algebra for any $S\subset|A|$},
i.e.\ that $(i_A^S)^*:\catMod A\to\catMod{A[S^{-1}]}$ is (left) exact. Since 
$(i_A^S)_*$ is exact and conservative, it suffices to prove that 
$(i_A^S)_*(i_A^S)^*:\catMod A\to\catMod A$ is exact. But this functor 
is isomorphic to $M\mapsto M[S^{-1}]$; replacing $S$ by the multiplicative 
system generated by~$S$, we can assume that $S$ is a multiplicative system,
and then $M[S^{-1}]=S^{-1}M=\injlim_\cS M_{[s]}$, $M_{[s]}:=M$, where 
$\cS$ is the filtered category constructed in~\ptref{sp:constr.localiz.mod}.
Since filtered inductive limits commute with finite projective limits 
in $\catSets$, hence also in $\catMod A$, we see that $M\mapsto S^{-1}M$ 
commutes with finite projective and arbitrary inductive limits, and 
in particular is exact.

\nxsubpoint (Classical localizations.)
Suppose that $A$ is a classical commutative ring, i.e.\ a commutative 
$\bbZ$-algebra. Then any $A$-algebra is a $\bbZ$-algebra, hence a 
classical ring as well, so the universal property of $A[S^{-1}]$ 
implies that $A[S^{-1}]$ coincides with classical $A$-algebra 
denoted in the same way, 
i.e.\ our ``generalized'' unary localizations coincide with classical 
localizations when computed over classical rings. The same is true for 
the localizations of $A$-modules: if $M$ is a module over a classical ring~$A$,
then $S^{-1}M$ of~\ptref{def:localiz} coincides with the classical $A$-module
denoted in this way.

In this way our unary localization theory is completely compatible with 
the classical one.

\nxsubpoint (Invertible elements of localizations.)
Let $S\subset|A|$ be a multiplicative system, and $\tilde S$ its saturation. 
We claim that the invertible elements of $|S^{-1}A|$ are exactly the 
elements of the form $a/s$, with $a\in\tilde S$ and $s\in S$, and that 
$(i_A^S)^{-1}(\Inv(S^{-1}A))=\tilde S$. Indeed, the monoid $|S^{-1}A|$ is 
actually the localization of monoid $|A|$ with respect to~$S$, i.e.\ 
everything depends here only on the multiplicative structures, and 
we can repeat the classical proof verbatim.

\nxsubpoint\label{sp:submod.loc} (Submodules of localizations.)
Let us fix a generalized ring~$A$, an $A$-module~$M$, and a multiplicative 
system $S\subset|A|$. Any $A$-submodule $N\subset M$ defines an 
$S^{-1}A$-submodule $S^{-1}N\subset S^{-1}M$, since the localization functor 
is exact. We have an increasing map in the opposite direction as well: 
it transforms $N'\subset S^{-1}M$ into $(i_M^S)^{-1}(N')\subset M$. 
We say that an $A$-submodule $N$ of~$M$ is {\em saturated\/} 
(with respect to~$S$) if 
$sx\in N$ with $s\in S$, $x\in M$ implies~$x\in N$. For any 
$A$-submodule $N\subset M$ we define its {\em saturation~$\tilde N$} to be the 
minimal saturated submodule of~$M$ containing~$N$.

It is immediate that for any $N'\subset S^{-1}M$ the pullback 
$N:=(i_M^S)^{-1}(N')$ is a saturated submodule of~$M$, and 
in this case $S^{-1}N=N'$. Conversely, for a saturated $N\subset M$ 
we have $(i_M^S)^{-1}(S^{-1}N)=N$, and for a general~$N$ the LHS is 
exactly the saturation of~$N$; it consists of all elements $x\in M$, 
such that $sx\in N$ for some $s\in S$. In this way {\em there is an increasing 
bijection between saturated $A$-submodules of~$M$ and all $S^{-1}A$-submodules 
of $S^{-1}M$}. Of course, this bijection and its inverse are given by 
$N\mapsto S^{-1}N$ and $N'\mapsto(i_M^S)^{-1}(N')$. In this respect the 
situation is exactly the same as in the classical case.

In particular, the above applies to the {\em ideals} of~$A$, i.e.\ 
$A$-submodules of~$M:=|A|$. We see that the ideals of~$S^{-1}A$ are in
one-to-one correspondence with {\em saturated\/} ideals $\ga\subset |A|$.

\nxsubpoint\label{sp:strictquot.loc} (Strict quotients of~$S^{-1}M$.)
Now let's study the strict quotients of $S^{-1}M$. In the classical case we 
work in abelian categories, so the strict quotients of a module are in 
one-to-one correspondence with their kernels, and we don't obtain anything 
new compared to the situation just discussed. However, the general case 
is slightly more complicated. In any case any strict quotient
$M\twoheadrightarrow P$ has a kernel $K:=M\times_PM\subset M\times M$; 
it is an $A$-submodule of $M\times M$, and a compatible equivalence relation 
on~$M$ as well. Actually the compatibility condition means exactly 
that $K$ is a submodule of~$M\times M$, so we have just to require $K$ to 
be an equivalence relation. Conversely, any such $K$ defines a strict 
quotient $M/K$ of~$M$.

Notice that $S^{-1}K\subset S^{-1}M\times S^{-1}M$ is the kernel 
of~$S^{-1}M\to S^{-1}P$ as well because of the exactness of localizations. 
Conversely, if we start with a strict quotient~$P'$ of the $S^{-1}A$-module 
$S^{-1}M$, its kernel $K'\subset S^{-1}M\times S^{-1}M$ is both 
an $S^{-1}A$-submodule and an equivalence relation on~$S^{-1}M$, hence 
$K:=(i_{M\times M}^S)^{-1}(K')\subset M\times M$ is a (saturated) $A$-submodule
of~$M\times M$ and an equivalence relation on~$M$; this yields a 
strict quotient $M/K$ of~$M$. Another description: $K$ is the kernel and 
$M/K$ is the image of the composite map $M\to S^{-1}M\to P'$. 
The results of~\ptref{sp:submod.loc} 
imply $K'=S^{-1}K$, hence $P'=S^{-1}(M/K)=S^{-1}M/K'=P'$.

We have just shown that {\em there is a canonical increasing bijection 
between strict quotients of the $S^{-1}A$-module~$S^{-1}M$ and 
the strict quotients of $M$ with saturated kernels in $M\times M$.}

\nxsubpoint (Localization and base change.)
Let $A$ be a generalized ring, $B$ a commutative $A$-algebra, 
$\rho:A\to B$ the structural homomorphism, $S\subset|A|$ a multiplicative 
set, $M$ an $A$-module, and $N$ a $B$-module.

First of all, notice that $\rho(S)^{-1}N$ is computed by means of the 
same filtered inductive limit as $S^{-1}(\rho_*N)$, so we get a canonical 
bijection $\rho(S)^{-1}N\cong S^{-1}(\rho_*N)$, easily seen to be
an $S^{-1}A$-module isomorphism. That's why we usually denote both 
sides by~$S^{-1}N$.

Another easy observation: \eqref{eq:comm.expr.loc} implies 
$B[\rho(S)^{-1}]\cong B\otimes_AA[S^{-1}]$, i.e.\ $\rho(S)^{-1}B=
B\otimes_AS^{-1}A$ when $S$ is a multiplicative set. Since 
$S^{-1}M\cong S^{-1}A\otimes_A M$ for any $A$-module~$M$, we get 
\begin{equation}
\rho(S)^{-1}(B\otimes_A M)\cong \rho(S)^{-1}B\otimes_{S^{-1}A}S^{-1}M
\end{equation}
Usually we write $S^{-1}B$ instead of $\rho(S)^{-1}B$. Notice that 
we can deduce the above isomorphism $\rho(S)^{-1}N\cong S^{-1}(\rho_*N)$ 
from the ``base change theorem''~\ptref{th:aff.base.change} as well, 
using the flatness of $S^{-1}A$ over~$A$.

Final remark: if $A$ is a $K$-algebra and $K'$ is another $K$-algebra, then 
$S^{-1}(A\otimes_KK')\cong(S^{-1}A)\otimes_KK'$, so we can write simply 
$S^{-1}A\otimes_KK'$.

\nxsubpoint (Strict quotients of $S^{-1}A$.)
Let $A$ and $S$ be as above. Consider a strict quotient $\phi:A\to A'$ of~$A$.
We say that {\em $A'$ is saturated (with respect to~$S$)} if 
$A'(n)$ is a saturated strict quotient of~$A(n)$, i.e.\ the kernel 
$R(n):=A(n)\times_{A'(n)}A(n)$ is a saturated $A$-submodule of 
$A(n)\times A(n)$, for all $n\geq0$.

Given any strict quotient~$A'$ of~$A$ as above, we can construct a strict 
quotient $S^{-1}A'$ of $S^{-1}A$; since $(S^{-1}A')(n)=S^{-1}(A'(n))$ as 
an $S^{-1}A$-module, this construction coincides with that 
of~\ptref{sp:strictquot.loc} on the level of individual strict 
quotients $A(n)\to A'(n)$. Conversely, we can start with an arbitrary 
strict quotient $A''$ of $S^{-1}A$ and construct a {\em saturated\/} 
strict quotient $A'$ of $A$, namely, the image of homomorphism 
$A\to S^{-1}A\to A''$. On the level of individual components $A''(n)$ this 
construction coincides again with that of~\ptref{sp:strictquot.loc}. 
Therefore, {\em the two maps given above define a canonical increasing 
bijection between saturated strict quotients of~$A$ and all strict quotients 
of~$S^{-1}A$.}

\nxsubpoint (Application to generalized fields.)
The above statement immediately implies that {\em any localization 
$S^{-1}K$ of a generalized field~$K$ (cf.~\ptref{def:gen.field}) is 
either subtrivial or a generalized field; in the latter case 
$K\to S^{-1}K$ is a monomorphism.}

\begin{DefD}\label{def:reg.elem} (Regular elements.) 
We say that an element $f\in |A|$ is {\em regular with respect to an 
$A$-module~$M$}, or {\em $M$-regular}, if the map $f_M:M\to M$ is 
injective, or equivalently, if $M\to M_f$ is injective. We say that 
{\em $f$ is $A$-regular}, or simply {\em regular}, if it is 
$A(n)$-regular for all $n\geq0$, i.e.\ all maps $f_{A(n)}:A(n)\to A(n)$ 
are injective, or equivalently, if $A\to A_f$ is a monomorphism.
\end{DefD}
Clearly, $M$-regular elements constitute a saturated multiplicative 
system $S_M\subset|A|$, hence any element of a multiplicative system 
generated by $M$-regular elements is itself $M$-regular. It is easy to see 
that $i_M^S:M\to M[S^{-1}]$ is injective iff $S$ if $M$-regular, i.e.\ 
iff $S\subset S_M$.

Similarly, the set of all regular elements is a saturated multiplicative 
system $T\subset|A|$, and $A\to A[S^{-1}]$ is injective iff $S\subset T$. 
Therefore, $T^{-1}A$ is the largest localization of~$A$, such that 
$A\to T^{-1}A$ is injective. We say that $T^{-1}A$ is  
{\em the total ring of fractions\/} of~$A$.

When all non-zero elements of~$|A|$ are regular (i.e.\ all elements of~$|A|$, 
when $A$ is a monad without zero), we say that $A$ is a {\em domain}. 
For example, any generalized field~$K$ is a domain, its total 
fraction ring $T^{-1}K$ is a generalized field as well, $K\to T^{-1}K$ 
is a monomorphism, and any non-zero element of $|T^{-1}K|$ is invertible. 
In other words, {\em any generalized field can be embedded into a generalized 
field, in which all non-zero unary operations become invertible}. Clearly, 
$K\to T^{-1}K$ has a universal property among all such embeddings.

\nxsubpoint (Tensor products of localizations.)
Let $S$ and $T$ be two subsets of~$A$. Formulas~\eqref{eq:comm.expr.loc} 
and~\eqref{eq:noncomm.expr.loc} imply 
\begin{gather}
A[S^{-1}]\otimes_A A[T^{-1}]\cong A[(S\cup T)^{-1}]\\
A[S^{-1}]\boxtimes_A A[T^{-1}]\cong A[(S\cup T)^{-1}]
\end{gather}
In particular, the category of $A$-algebras of form $S^{-1}A$ 
is stable under coproducts and pushouts.

\nxsubpoint (Morphisms between localizations.)
Since all $A\to S^{-1}A$ are epimorphisms (and even NC-epimorphisms), 
there is at most one $A$-algebra homomorphism between two such localizations, 
i.e.\ the category of all localizations of~$A$ is actually a pre-ordered 
set. The universal property of $A[S^{-1}]$ shows that $A\to T^{-1}A$ 
factorizes through $A\to A[S^{-1}]$ iff all elements of $S$ become invertible 
in $T^{-1}A$, i.e.\ iff $S$ is contained in 
$\tilde T=(i_A^T)^{-1}(\Inv(T^{-1}A))$, the saturation on~$T$. When 
$S$ is a multiplicative system, this condition is also equivalent to 
$\tilde S\subset\tilde T$; since $\tilde S^{-1}A\cong S^{-1}A$, we see that 
{\em the category of localizations of~$A$ is equivalent to the category 
defined by the set of saturated multiplicative systems in~$|A|$, ordered 
by inclusion.}

In particular, $A\to A_f$ factorizes through $A\to A_g$ iff 
$\tilde S_g\subset\tilde S_f$ iff $g\in\tilde S_f$ iff $g$ divides $f^n$ 
for some $n\geq0$.

\nxsubpoint (Filtered inductive limits of localizations.)
If $S_\alpha$ is a collection of multiplicative subsets of $|A|$, 
filtered by inclusion, then $S:=\bigcup_\alpha S_\alpha$ is a multiplicative 
system as well, and $\injlim_\alpha S_\alpha^{-1}A\simto S^{-1}A$. This 
statement can be generalized to the case when we are given a filtered 
inductive system of generalized rings~$A_\alpha$ as well, $S_\alpha$ 
is a multiplicative system in~$|A_\alpha|$, and the transition 
morphisms $f_{\alpha\beta}:A_\alpha\to A_\beta$ satisfy 
$f_{\alpha\beta}(S_\alpha)\subset S_\beta$: we put $A:=\injlim A_\alpha$,
$S:=\injlim S_\alpha\subset |A|$, and then $\injlim S_\alpha^{-1}A_\alpha
\simto S^{-1}A$. We can also take an inductive system of $A_\alpha$-modules 
$M_\alpha$, put $M:=\injlim M_\alpha$, and obtain a canonical isomorphism
of $S^{-1}A$-modules $\injlim S_\alpha^{-1}M_\alpha\simto S^{-1}M$.

Notice that any multiplicative system $S\subset|A|$ is a filtered inductive 
limit of its finitely generated multiplicative subsystems 
$\langle f_1,\ldots,f_n\rangle$, $f_i\in S$. On the other hand, clearly 
$M[f_1^{-1},\ldots,f_n^{-1}]\cong M[f^{-1}]$, where $f:=f_1\cdots f_n\in S$, 
for any $A$-module~$M$. We conclude that
\begin{equation}
S^{-1}A\cong\injlim_{f\in S}A_f\quad\text{and}\quad 
S^{-1}M\cong\injlim_{f\in S}M_f
\end{equation}
The inductive limits are computed here with respect to the following preorder 
on~$S$: $f\prec g$ iff $\tilde S_f\subset\tilde S_g$ iff $f$ divides some 
power of~$g$. This condition guarantees the existence of maps 
$A_f\to A_g$ and $M_f\to M_g$ whenever $f\prec g$.

\nxsubpoint (Finitely presented localizations.)
Notice that $A[S^{-1}]$ is a finitely generated $A$-algebra iff we can find 
a finite subset $S_0\subset S$, such that $A[S_0^{-1}]\to A[S^{-1}]$ 
is surjective (cf.~\eqref{eq:comm.expr.loc}). Similarly, $A[S^{-1}]$ 
is a finitely presented $A$-algebra iff there is a finite subset 
$S_0\subset S$, such that $A[S_0^{-1}]\simto A[S^{-1}]$. In this case 
we can replace $S$ by $S_0$, and even by the product $f$ of all elements 
of~$S_0$, i.e.\ {\em finitely presented localizations of~$A$ are exactly 
the localizations of form~$A_f$, $f\in|A|$.}

\nxsubpoint (Localization, tensor products, and $\iHom$.)
If $S\subset|A|$ is a multiplicative system and $M$, $N$ two $A$-modules, 
we have canonical isomorphisms
\begin{equation}
S^{-1}(M\otimes_AN)\cong S^{-1}M\otimes_{S^{-1}A}S^{-1}N\cong
S^{-1}M\otimes_AS^{-1}N\cong S^{-1}M\otimes_AN
\end{equation}
In fact, the existence of the second and the third 
of these isomorphisms is due to the 
fact that $S$ acts by bijections on $S^{-1}M\otimes_A S^{-1}N$ 
and $S^{-1}M\otimes_AN$, hence these two modules admit a natural 
$S^{-1}A$-module structure. Now the first isomorphism is given 
by the isomorphism of~\ptref{sp:pullback.tensfunct}, valid for any 
$A$-algebra $A'$, when we put $A':=S^{-1}A$:
\begin{equation}
(M\otimes_AN)_{(A')}\cong M_{(A')}\otimes_{A'}N_{(A')}
\end{equation}

Similarly, for any two $A$-modules $M$ and $N$ we have a canonical 
$S^{-1}A$-module homomorphism
\begin{equation}
S^{-1}\Hom_A(M,N)\to\Hom_{S^{-1}A}(S^{-1}M,S^{-1}N)
\end{equation}
We claim that it is injective for a finitely generated~$M$, and 
an isomorphism for a finitely presented~$M$. In fact, we have a  
more general statement, valid for any {\em flat\/} $A$-algebra~$A'$
(hence also for $A'=S^{-1}A$), 
concerning the following canonical $A'$-module homomorphism:
\begin{equation}
\Hom_A(M,N)_{(A')}\stackrel\gamma\to\Hom_{A'}(M_{(A')},N_{(A')})
\end{equation}
To show the injectivity (resp.\ bijectivity) of~$\gamma$ 
for a finitely generated (resp.\ presented) $M$, we first show that 
it is an isomorphism for~$M=A(n)$: in this case the statement reduces to  
$(N^n)_{(A')}\cong(N_{(A')})^n$, an immediate consequence of the  
flatness of~$A'$. Then we choose a strict epimorphism $A(n)\to M$ 
(resp.\ a presentation $A(m)\rightrightarrows A(n)\to M$), and 
the general result follows from the special case just considered 
by a standard diagram chasing.

\nxsubpoint (Localizations and alternativity.)
Since $S^{-1}A$ is a unary $A$-algebra, $S^{-1}A$ is an alternating 
$\Fpm$-algebra whenever $A$ is one (cf.~\ptref{sp:examp.altmon},g). 
Moreover, any matrix $Z\in M(n,n;S^{-1}A)=(S^{-1}A)(n)^n\cong S^{-1}(A(n)^n)$ 
can be written in form $Z_0/s$, with $Z_0\in M(n,n;A)$ and $s\in S$, 
and $Z_0/s=Z'_0/s'$ iff $ts'\cdot Z_0=ts\cdot Z'_0$ for some $t\in S$. 
This immediately implies that {\em the properties $(DET'_r)$ and $(DET^*_r)$ 
of~\ptref{p:invdet.matr} hold for $S^{-1}A$ whenever they hold for~$A$.} 
We cannot say the same about properties~$(DET_r)$.

\nxsubpoint\label{sp:examp.loc} (Example.)
We know that localizations of classical rings coincide with classical 
localizations, so we don't obtain anything new in the classical case. 
Let's show that $\Zinfty[f^{-1}]\cong\bbR$, $\barZinfty[f^{-1}]
\cong\bbC$, and $\Zninfty[f^{-1}]\cong\bbQ$, for any 
$f\in|\Zinfty|$ (resp\dots), such that $0<|f|<1$. Let's 
treat the first case; the two remaining cases are proved similarly.

First of all, $\Zinfty\to\bbR$ is injective (i.e.\ a monomorphism), 
hence the same is true for $\Zinfty[f^{-1}]\to\bbR[f^{-1}]$ 
(this is a general fact: {\em if $A\to B$ is a monomorphism, then 
$S^{-1}A\to S^{-1}B$ is also a monomorphism}). Since $f\neq0$, 
$\bbR[f^{-1}]\cong\bbR$, and we can identify $B:=\Zinfty[f^{-1}]$ 
with a submonad of~$\bbR$: $\Zinfty\subset B=\Zinfty[f^{-1}]\subset\bbR$.
We want to show that $B=\bbR$, i.e.\ that any $\lambda=(\lambda_1,\ldots,
\lambda_n)\in\bbR(n)=\bbR^n$ lies in $B(n)$. Since $|f|<1$, we can find 
an integer $n\geq0$, such that $|\lambda_1|+\cdots+|\lambda_n|\leq f^{-n}$. 
Put $\mu:=f^n\lambda$; then $|\mu_1|+\cdots+|\mu_n|\leq 1$, hence 
$\mu$ lies in $\Zinfty(n)\subset B(n)$. On the other hand, $f^{-1}=1/f
\in B(1)$, hence $\lambda=\mu/f^n$ lies in $B(n)$ as well.

\nxsubpoint\label{sp:pseudoloc} (Axiomatic description of localizations. 
Pseudolocalizations.)
Consider the following set of properties of a commutative $A$-algebra~$B$, 
with the structural morphism $\rho:A\to B$:
\begin{itemize}
\item[1)] $\rho$ is an $NC$-epimorphism, i.e.\ any $A$-module admits 
at most one compatible $B$-module structure;
\item[2)] $\rho_*:\catMod B\to\catMod A$ is fully faithful;
\item[3)] $B$ is a flat $A$-algebra, i.e.\ $\rho^*$ is exact;
\item[4)] $B$ is a unary $A$-algebra;
\item[5)] $B$ is finitely presented over $A$.
\item[6)] $\rho_*$ is a $\otimes$-functor, i.e.\ 
$\rho_*M\otimes_A\rho_*N\to\rho_*(M\otimes_BN)$ is an isomorphism 
for any two $B$-modules $M$ and~$N$.
\end{itemize}
Notice that $2)\Rightarrow 1)$, that $1)$ implies that $\rho$ is an 
epimorphism, and that any (unary) localization $S^{-1}A$ satisfies 1)--4) 
and~6), and $A_f$ satisfies all of the above conditions. 
However, the unarity of $S^{-1}A$ hasn't 
been used in most of the previous considerations. Therefore, we might 
define the {\em pseudolocalizations of~$A$} to be the $A$-algebras~$B$ 
that satisfy 1)--3), and {\em open pseudolocalizations of~$A$} by requiring~5) 
as well. Most properties of localizations should be also true for 
pseudolocalizations. However, one can check that 6) is a formal consequence 
of 2) and~4), but not of~2) alone (consider $\rho:\Zinfty\to\Finfty$, 
$M=N=\Finfty(2)$ for this), 
so we cannot expect 6) to hold for all pseudolocalizations.

\nxsubpoint\label{sp:prop1.psloc} (Properties of pseudolocalizations.)
Clearly, the set of pseudolocalization morphisms is closed 
under composition and base change. For the latter property we use 
corollary~\ptref{sp:bc.stab.flat} of the ``base change theorem''
\ptref{th:aff.base.change} to show that flatness is preserved by any base 
change; to check stability of 2)
under some base change $f:A\to A'$ we observe that a $B':=A'\otimes_AB$-module
structure on an $A'$-module~$M$ is the same thing as a $B$-module structure 
on $f_*M$ (automatically commuting with the action of~$A'$ on~$M$: 
for any $t\in A'(n)$ the map ${[t]}_M:M^n\to M$ is an $A$-module 
homomorphism, hence also a $B$-module homomorphism, $\rho_*$ being 
left exact and fully faithful; notice that the flatness of~$\rho$ is 
not needed here), hence 
$\catMod{B'}$ can be identified with a full subcategory of $\catMod{A'}$, 
the preimage of $\catMod B\subset\catMod A$ under~$f_*$.

All pseudolocalizations of~$A$ constitute 
a full subcategory of the category of commutative $A$-algebras; 
property 1) implies that this subcategory is a preordered set, 
and it is stable under coproducts and pushouts. The same is true for 
the subcategory of open pseudolocalizations.

\nxsubpoint {\bf Questions.} 
Is it true that all unary pseudolocalizations are of form $S^{-1}A$? 
If not, is this true for a classical~$A$? Are there any non-unary 
pseudolocalizations of generalized rings? Of alternating monads with 
$(DET^*_\infty)$ property?

{\bf Answers.} No, no, yes, unknown.

\nxsubpoint (An example.)
Notice that generalized rings without zero often have non-unary 
pseudolocalizations. For example, if $A$ is a classical ring, and 
$S\subset|A|$ is a multiplicative system, the homomorphism of 
generalized rings $\Aff_A\to\Aff_{S^{-1}A}$ induced by $A\to S^{-1}A$ 
(cf.~\ptref{sp:examp.submon},g) is easily seen to be a pseudolocalization. 
On the other hand, $|\Aff_A|=|\Aff_{S^{-1}A}|=\{\bu\}$, so $\Aff_{S^{-1}A}$ 
is not unary over $\Aff_A$ when $S^{-1}A\neq A$. If $A\to S^{-1}A$ is 
not surjective, $\Aff_{S^{-1}A}$ cannot be even pre-unary over~$\Aff_A$.

Notice that there is only the trivial multiplicative system 
in $|\Aff_A|=\{\bu\}$, so we don't have any non-trivial localizations 
of~$\Aff_A$.

\nxsubpoint\label{sp:loc.matr} (Localization with respect to matrices.)
We can try to generalize our definition of unary localization by 
formally inverting families of square matrices 
$Z_i\in M(r_i,r_i;A)=A(r_i)^{r_i}$. Let's consider the case of one square 
matrix $Z\in M(r,r;A)$; the general case can be deduced from this one 
by (transfinite) induction.

So let's fix $Z=(Z_1,\ldots,Z_r)\in M(r,r;A)$ and define
\begin{gather}
A\langle Z^{-1}\rangle:=A\langle Z^-_1,\ldots, Z^-_r\,|\,
Z^-_i(Z_1,\ldots,Z_r)=\{i\}_r=Z_i(Z^-_1,\ldots,Z^-_r), \text{all $i$}\rangle\\
A[Z^{-1}]:=A[Z^-_1,\ldots, Z^-_r\,|\,
Z^-_i(Z_1,\ldots,Z_r)=\{i\}_r=Z_i(Z^-_1,\ldots,Z^-_r), \text{all $i$}]
\end{gather}
Clearly, $A\to A\langle Z^{-1}\rangle$ is universal among all monad 
homomorphisms $\rho:A\to B$, such that $\rho(Z)$ is invertible over~$B$, 
and $A\to A[Z^{-1}]$ has the same universal property among all 
commutative $A$-algebras.

We want to show that $B:=A\langle Z^{-1}\rangle$ is commutative; this 
would imply $A[Z^{-1}]\cong A\langle Z^{-1}\rangle$, 
similarly to what we had in~\ptref{sp:ex.un.loc}. So we have to 
check that the generators $Z^-_i$ commute among themselves and the operations 
of~$A$. For any $B$-module~$N$ we denote by ${[Z]}_N:N^r\to N^r$ the map 
with components ${[Z_j]}_N:N^r\to N$, and define ${[Z^-]}_N$ similarly. 
Clearly, all $Z^-_j$ commute with 
all operations from~$A$ iff all ${[Z^-_j]}_N$ are $A$-linear 
(for all $N$ or just for all $N=B(n)$) iff ${[Z^-]}_N:N^r\to N^r$ is 
$A$-linear for such $B$-modules~$N$. Now observe that ${[Z^-]}_N={[Z]}_N^{-1}$,
and ${[Z]}_N={[Z]}_{\rho_*N}$ is $A$-linear for any $N$ because of 
the commutativity of~$A$. 

It remains to check that $Z^-_j$ commute among themselves. This means 
equality of two maps $N^{\str\times\str}\to N$ for all $B$-modules~$N$ 
and all $1\leq i,j\leq r$. We can use the same trick as above and combine 
together these maps for all values of $i$ and $j$ into two maps 
$[Z^-\otimes Z^-]'_N, [Z^-\otimes Z^-]''_N:N^{\str\times\str}\to 
N^{\str\times\str}$; we have to check that these two maps coincide for 
all $B$-modules~$N$. Notice that these two maps can be described as
follows. The first of them takes a matrix $(z_{ij})\in N^{\str\times\str}$, 
applies ${[Z^-]}_N:N^r\to N^r$ to each of its rows, and then applies 
${[Z^-]}_N$ again to each of the columns, and the second map applies 
$[Z^-]_N$ first to columns, and then to rows. Hence 
$[Z^-\otimes Z^-]'_N=([Z\otimes Z]''_N)^{-1}$ and $[Z^-\otimes Z^-]''_N=
([Z\otimes Z]'_N)^{-1}$, so it remains to check 
$[Z\otimes Z]''_N=[Z\otimes Z]'_N$. We can replace here $N$ by $\rho_*N$, 
since $Z$ is defined over~$A$, and the resulting equality is true 
for any $A$-module $M$ (and in particular for $M=\rho_*N$) because of 
the commutativity of~$A$. 
One can also check that {\em $A[Z^{-1}]$ is 
alternating whenever~$A$ is alternating\/} essentially in the same way 
we have just verified the commutativity between the~$Z^-_j$.

An immediate consequence of $A\langle Z^{-1}\rangle=A[Z^{-1}]$ is that 
giving a monad homomorphism $A[Z^{-1}]\to B$ is equivalent to giving 
a monad homomorphism $\rho:A\to B$, such that $\rho(Z)$ becomes 
invertible over~$B$. Applying this to $B:=\END(M)$, we see that 
$A[Z^{-1}]$-modules are exactly the $A$-modules~$M$, such that 
${[Z]}_M:M^r\to M^r$ is bijective. In particular, the 
scalar restriction functor $(i_A^Z)_*$ is fully faithful, and 
$i_A^Z:A\to A[Z^{-1}]$ is an NC-epimorphism.

In this way we see that $A[Z^{-1}]$ has the properties 1), 2) and 5) 
of~\ptref{sp:pseudoloc}; if we localize with respect to an arbitrary 
set of matrices, we still obtain properties 1) and 2). However, 
we'll see in a moment that the property 3) (flatness) doesn't hold in general, 
so we don't obtain pseudolocalizations in this way.

\nxsubpoint (Example of a non-flat matrix localization.)
Put $A:=\bbZ_{\geq0}$, $Z:=\smallmatrix1011=(+,\{2\})\in 
M(2,2;A)=A(2)^2\cong\bbZ_{\geq0}^{2\times 2}$. Then $\catMod A$ 
is the category of commutative monoids, and ${[Z]}_M:M^2\to M^2$ 
is the map $(x,y)\mapsto(x+y,y)$. We see that ${[Z]}_M$ is bijective iff 
$M$ is an abelian group, hence $\catMod{A[Z^{-1}]}=\catMod\bbZ$, hence
$A[Z^{-1}]=\bbZ_{\geq0}[{\smallmatrix1011}^{-1}]=\bbZ$, i.e.\ 
$\bbZ$ is a matrix localization of~$\bbZ_{\geq0}$. 

We claim that $\bbZ$ is not flat over $\bbZ_{\geq0}$. Indeed, the base 
change functor $\bbZ\otimes_{\bbZ_{\geq0}}$ transforms any commutative 
monoid $M$ into its symmetrization $\tilde M$, and it is quite easy 
to construct a commutative monoid~$M$ and a submonoid $N\subset M$, 
such that $\tilde N\to\tilde M$ is not injective. For example, put 
$M:=\bbZ\times\bbZ_{\geq0}\times\{\pm1\}/\sim$, where the equivalence relation
$\sim$ is given by $(x,n,+1)\sim(x,n,-1)$ for any $x\in\bbZ$ and $n>0$, 
and $N:=\bbZ\times0\times\{\pm 1\}\subset M$. Then 
$\tilde N=\bbZ\times\{\pm1\}$ and $\tilde M=\bbZ\times\bbZ$, and the 
map $\tilde N\to\tilde M$ is given by $(x,\pm1)\to (x,0)$, hence it is 
not injective.

\nxsubpoint (Matrix localization and $(DET)$-properties.)
When $A$ is a classical ring, a square matrix $Z\in M(r,r;A)$ is 
invertible in an $A$-algebra~$B$ iff $\det(Z)$ is invertible.
More generally, suppose that $A$ is an alternating $\Fpm$-algebra, 
and $Z\in M(r,r;A)$ is a matrix, such that the property $(DET'_r)$ 
of~\ptref{p:invdet.matr} is fulfilled for~$Z$ (e.g.\ $A$ has the 
$(DET'_\infty)$-property). Suppose that $A[Z^{-1}]$ is also alternating 
(something that we haven't checked in detail in~\ptref{sp:loc.matr}). 
Then $\det(Z)$ is invertible in $A[Z^{-1}]$, hence $A\to A[Z^{-1}]$ 
factorizes through $A[\det(Z)^{-1}]$. On the other hand, $(DET'_r)$ 
means that we have a matrix $Z^*\in M(r,r;A)$ and an integer $N>0$, 
such that $ZZ^*=Z^*Z=\det(Z)^N$; this implies that $Z^*/\det(Z)^N$ is an 
inverse to $Z$ over $A[\det(Z)^{-1}]$, hence we have a homomorphism 
in the opposite direction $A[Z^{-1}]\to A[\det(Z)^{-1}]$. Since 
$A\to A[Z^{-1}]$ and $A\to A[\det(Z)^{-1}]$ are epimorphic, we conclude that 
under the above assumptions $A[Z^{-1}]\cong A[\det(Z)^{-1}]$, i.e.\ 
{\em matrix localization over alternating generalized rings with 
$(DET'_\infty)$-property shouldn't give anything new compared to  
unary localizations.}

\nxpointtoc{Prime spectrum of a generalized ring}
\nxsubpoint (Ideals in generalized rings.)
Recall that an {\em ideal\/~$\ga$} in a generalized ring~$A$ is by definition 
an $A$-submodule of~$|A|$. For any subset $S\subset|A|$ we denote by 
$(S)$ the ideal generated by~$S$, i.e.\ the smallest ideal in~$A$ 
containing~$S$; it is equal to the image of 
$A(S)\to|A|$ (cf.~\ptref{sp:subm.gen.subset}). When we have two subsets 
$S$, $T\subset|A|$, we usually write $(S,T)$ instead of~$(S\cup T)$, and 
similarly for larger finite families of subsets and/or elements of~$|A|$. 
In particular, for any two ideals $\ga$ and $\gb$ we can construct their 
upper bound $(\ga,\gb)$, denoted also by $\ga+\gb$, i.e.\ the smallest 
ideal containing both $\ga$ and~$\gb$; clearly, $\ga+\gb$ coincides with the 
image of $\ga\oplus\gb\to|A|$. On the other hand, the intersection of 
any family of ideals and the union of any filtering family of ideals is 
an ideal again; in particular, we have lower bounds $\ga\cap\gb$ in the 
lattice of ideals of~$|A|$ as well. Notice that the {\em unit ideal\/}
$(1)=|A|$, and the {\em initial ideal\/} $\emptyset_A$ (denoted also by 
$(0)$ or $0$ and called the {\em zero ideal\/} when $A$ is a generalized 
ring with zero) are the largest and the smallest elements of the 
ideal lattice of~$A$.

For any two ideals $\ga$ and $\gb$ we can construct their 
{\em product\/~$\ga\gb$}, i.e.\ the ideal generated by all products $xy$, 
with $x\in\ga$, $y\in\gb$. Another description: $\ga\gb$ is the image of 
$\ga\otimes_A\gb\to|A|$. Clearly, this multiplication of ideals is 
associative and commutative, with $(1)\cdot\ga=\ga$ and $\emptyset_A\cdot\ga=
\emptyset_A$. We can extend our definition and define 
$\ga_1\cdots\ga_nN\subset M$ (e.g.\ as the image of $\ga_1\otimes_A\cdots
\otimes_A\ga_n\otimes_A N\to M$), for any ideals $\ga_i\subset|A|$ and 
any $A$-submodule $N\subset M$. 

Finally, if $B$ is an $A$-algebra, and $\ga$ an ideal of~$A$, we denote 
by $B\ga$ the ideal of~$B$ generated by the image of~$\ga$, equal to 
the image of $B\otimes_A\ga\to B\otimes_A|A|=|B|$. 
Of course, all of the above constructions coincide 
with the usual ones over a classical~$A$.

\nxsubpoint (Prime ideals.)
We say that $\gp\subset|A|$ is a {\em prime ideal\/} if it an ideal, 
and its complement $S_\gp:=|A|-\gp$ is a multiplicative system, i.e.\ 
$\bu\not\in\gp$, and $xy\in\gp$ implies that either $x$ or $y$ is in~$\gp$.
The localization $S_\gp^{-1}A$ is usually denoted by $A_\gp$, and similarly 
for the localizations of an $A$-algebra~$B$ or a $B$-module~$N$.

\nxsubpoint (Maximal ideals.)
We say that $\gm\subset|A|$ is a {\em maximal ideal\/} if it is a maximal 
element in the ordered set of all ideals of $|A|$ distinct from~$(1)$, 
i.e.\ if $\gm\subset\ga$ implies $\ga=\gm$ or $\ga=(1)$. Notice that 
{\em any maximal ideal is prime}: indeed, by definition $1\not\in\gm$, 
and if $x$, $y\not\in\gm$, we have $(\gm,x)=(\gm,y)=(1)$, hence 
$(\gm,xy)\supset(\gm^2,x\gm,y\gm,xy)=(\gm,x)(\gm,y)=(1)$, and this is possible 
only if~$xy\not\in\gm$. Since any filtered union of ideals $\neq(1)$
is again an ideal $\neq(1)$, an application of Zorn lemma shows that 
{\em any ideal $\ga\neq(1)$ is contained in some maximal ideal\/}; 
more generally, any submodule $N$ of a finitely generated $A$-module~$M$ 
is either equal to~$M$ or contained in a maximal submodule of~$M$.
In particular, {\em any non-trivial generalized ring ($\st1_+$ included) 
contains at least one maximal (hence also prime) ideal.}

Applying this to principal ideals $(a)$ we see that {\em the union of 
all ideals $\neq(1)$ of~$A$ coincides with the union of all maximal ideals, 
equal to the complement of the set~$\Inv(A)$ of invertible elements of~$A$.}

\nxsubpoint\label{sp:un.local.rings} (Generalized local rings.)
It is natural to say that $A$ is {\em local\/} if it has exactly one
maximal ideal $\gm=\gm_A$. This is equivalent to saying that $\gm\subset|A|$
is an ideal, such that all elements of $|A|-\gm$ are invertible.
Furthermore, a homomorphism of generalized local rings $f:A\to B$
is {\em local\/} if $f^{-1}(\gm_B)=\gm_A$.

\nxsubpoint\label{sp:prime.spec} (Prime spectrum.)
The {\em prime spectrum\/} $\Spec A=\Spec^pA$ of a generalized ring~$A$ 
is simply the set of its prime ideals. Clearly, $\Spec A=\emptyset$ iff 
$A=\st1$. For any subset $M\subset|A|$ we put $V(M):=\{\gp\in\Spec A\,|\,
M\subset\gp\}$. We obtain immediately the usual formulas 
$V(M)=V(\ga)$ if $\ga=(M)$, $V(\emptyset)=\Spec A$, $V(1)=\emptyset$, 
$V(\bigcup_\alpha M_\alpha)=\bigcap_\alpha V(M_\alpha)$, and finally 
$V(\ga)\cup V(\gb)=V(\ga\cap\gb)=V(\ga\gb)$, that show that 
all these $V(M)$ are the closed subsets for a certain topology on 
$\Spec^pA$, called the {\em spectral\/} or {\em Zariski topology\/} 
on~$\Spec^pA$.

We write $V(f)$ instead of $V(\{f\})$ for any $f\in|A|$. Clearly, 
$V(M)=\bigcup_{f\in M}V(f)$ for any~$M$, hence the {\em principal 
open subsets} $D(f):=\{\gp\,|\,f\not\in\gp\}$ constitute a base of the 
Zariski topology on~$\Spec^pA$. We have the usual formula 
$D(fg)=D(f)\cap D(g)$.

For any homomorphism of generalized rings $\phi:A\to B$ and any 
prime ideal $\gq\subset|B|$ its pullback ${}^a\phi(\gq):=|\phi|^{-1}(\gq)$ 
is a prime ideal in~$A$, so we get a map ${}^a\phi:\Spec^pB\to\Spec^pA$. 
This map is continuous because of the usual formula 
$({}^a\phi)^{-1}(V(M))=V(\phi(M))$. In particular, 
$({}^a\phi)^{-1}(D(f))=D(\phi(f))$ for any $f\in|A|$.

\nxsubpoint\label{sp:specp.localiz} ($\Spec^pS^{-1}A$.)
Let $S$ be a multiplicative system in~$A$. We know that $\ga\mapsto S^{-1}\ga$ 
induces a bijection between $S$-saturated ideals $\ga\subset|A|$ and 
all ideals of~$S^{-1}A$ (cf.~\ptref{sp:submod.loc}), and a saturated~$\ga$
coincides with the pullback of~$S^{-1}\ga$. Clearly, the pullback of a 
prime ideal of $S^{-1}A$ is a saturated prime ideal $\gp\subset|A|$, 
and it is immediate that a prime ideal $\gp$ in~$A$ is $S$-saturated iff 
$\gp\cap S=\emptyset$. This implies that {\em ${}^ai_S^A:\Spec^pS^{-1}A\to
\Spec^pA$ induces a bijection between $\Spec^pS^{-1}A$ and the set of 
prime ideals $\gp\subset|A|$ that do not intersect~$S$.} Moreover, 
this bijection is actually a {\em homeomorphism\/} if we equip this set 
of prime ideals with the topology induced by that of~$\Spec^pA$. Indeed, 
$D(f/s)$, $f\in|A|$, $s\in S$, constitute a base of the topology of
$\Spec^pS^{-1}A$; but obviously $D(f/s)=D(f/1)=D(f)\cap\Spec^pS^{-1}A$. 
Henceforth we shall usually identify $\Spec^pS^{-1}A$ with this subspace 
of~$\Spec^pA$.

In particular, $\Spec^p A_f=\Spec^p A[f^{-1}]$ is homeomorphic to 
the principal open subset $D(f)\subset\Spec^p A$.

\nxsubpoint (An application.)
We can use the above constructions and definitions to show that 
{\em if $A$ is a generalized ring with zero, and $f\in|A|$ is not nilpotent, 
then $f\not\in\gp$ for some $\gp\in\Spec^pA$.} Therefore, the intersection 
of all prime ideals of a generalized ring with zero coincides with the 
set of nilpotent elements. On the other hand, if~$A$ is a generalized ring 
without zero, then the initial ideal $\emptyset_A$ is prime, hence 
the intersection of all prime ideals is empty.

\nxsubpoint\label{sp:opcov.specp} (Open covers and quasicompactness.)
Let $f_\alpha$ be a family of elements of~$|A|$, and let us denote by 
$\ga$ the ideal generated by this family. Then $\Spec^pA-\bigcup_\alpha 
D(f_\alpha)=\bigcap_\alpha V(f_\alpha)=V(\ga)$, hence 
$\Spec^p A=\bigcup_\alpha D(f_\alpha)$ iff $V(\ga)=\emptyset$ iff 
$\ga=(1)$, since any non-unit ideal is contained in some maximal ideal. 
In other words, {\em principal open subsets $D(f_\alpha)$ constitute 
a cover of~$\Spec^pA$ iff the ideal generated by $(f_\alpha)$ is equal 
to the unit ideal~$(1)$.} This means that 
$1=t(f_{\alpha_1},\ldots,f_{\alpha_n})$ for some $n\geq0$, $t\in A(n)$ 
and $\alpha_i$. In this case the finite set of elements $f_{\alpha_i}$ 
already generates the unit ideal, i.e.\ 
$\Spec^pA=\bigcup_{i=1}^nD(f_{\alpha_i})$. Using the fact that 
the principal open subsets are a base of topology of $\Spec^pA$, we see that 
{\em $\Spec^p A$ is quasicompact}. Of course, this is also true for any 
$D(g)\cong\Spec^p A_g$.

\nxsubpoint\label{sp:opcov.Dg} (Open covers of $D(g)$.)
Notice that $D(g)=\Spec^pA$ iff $V(g)=\emptyset$ iff $g$ is invertible in~$A$.
Therefore, $D(f)\subset D(g)$ iff $D(g)\cap\Spec^pA_f=\Spec^pA_f$ iff 
$D(g/1)=\Spec^pA_f$ iff $g$ is invertible in~$A_f$ iff $g$ divides some 
power of~$f$ iff $\tilde S_g\subset\tilde S_f$ iff $A\to A_f$ factorizes 
through~$A\to A_g$.

When $D(g)\subset\bigcup_\alpha D(f_\alpha)$? This condition is equivalent 
to $\bigcup_\alpha D(f_\alpha/1)=\Spec^pA_g$, and this is equivalent to
$\ga_g=(f_\alpha/1)\subset|A_g|$ being equal to~$|A_g|$. However, 
$\ga_g=A_g\cdot\ga$, where $\ga=(f_\alpha)\subset|A|$ is the ideal of~$A$ 
generated by~$f_\alpha$. Clearly, $\ga_g=(1)$ iff the saturation of 
$\ga$ with respect to $S_g=\{1,g,g^2,\ldots\}$ is equal to the unit ideal 
iff $S_g\cap\ga\neq\emptyset$ iff $g^n\in\ga=(f_\alpha)$ iff 
$g\in\rad(\ga):=\{g\,|\,g^n\in\ga$ for some $n>0\}$.

\nxsubpoint\label{sp:rad.ideal} (The radical of an ideal.)
For any ideal $\ga\subset|A|$ we denote by $\rad(\ga)$ the {\em radical\/}
of~$\ga$, i.e.\ the set of all elements $g\in|A|$, such that $g^n\in\ga$ 
for some $n>0$. We claim that {\em $\rad(\ga)$ is always an ideal, 
namely, the intersection of all prime ideals containing~$\ga$.} Clearly, 
$\rad(\ga)$ is contained in this intersection, since $g^n\in\ga\subset\gp$ 
implies $g\in\gp$. Conversely, suppose that $g\not\in\rad(\ga)$, i.e.\ 
$S_g\cap\ga=\emptyset$. Then the $S_g$-saturation of~$\ga$ doesn't 
contain~$\bu$, hence $\ga_g:=S_g^{-1}\ga\neq(1)$, hence this ideal is contained
in some maximal ideal $\gm$ of~$A_g$, and the preimage $\gp$ of~$\gm$ 
in~$A$ is a prime ideal of~$A$ that contains~$\ga$, but does not contain~$g$.

\nxsubpoint
Another simple statement is that for any $N>0$ a collection of 
elements $f_i\in|A|$ generates the unit ideal iff the same is true 
for $(f_i^N)$. 
Indeed, we know that $f_i$ generate the unit ideal iff $\bigcup_iD(f_i)=
\Spec^pA$, and similarly for $f_i^N$, and our statement follows from 
$D(f_i)=D(f_i^N)$.

\nxsubpoint
One can check that the irreducible closed subsets of $\Spec^pA$ are 
exactly the subsets of the form $V(\gp)$, with $\gp$ any prime ideal. 
Clearly, these closed subsets admit a unique generic point, namely, 
$\gp$ itself. This shows that $\Spec^pA$ is a sober topological space, 
and in particular a Kolmogorov space.

\nxsubpoint\label{sp:examp.prime.spec} 
(Examples.) Let's consider some examples of prime spectra. 
Examples b) and e) show that we shouldn't be too optimistic about the 
theory of prime spectra in the generalized context.
\begin{itemize}
\item[a)] Of course, for a classical commutative ring~$A$ we recover its 
prime spectrum $\Spec A$ in the usual sense.
\item[b)] On the other hand, $\Spec^p\Aff_A$ is a one-element set
for any classical ring~$A$, since $|\Aff_A|=\{\bu\}$, and the only prime 
ideal in~$\Aff_A$ is the initial ideal~$\emptyset$. This shows that 
the theory of prime spectra is not very nice for such generalized rings: 
one would rather expect $\Spec^p\Aff_A=\Spec^p A$ as a topological space 
(but with a different structure sheaf).
\item[c)] Similarly, $\Spec^p\Fempty$, $\Spec^p\Fone$, $\Spec^p\Fpm$
and $\Spec^p\bbF_{1^n}$ are one-element sets, since the initial ideal is 
the only prime ideal in any of these generalized rings.
\item[d)] We know already that $\Spec^p\Zinfty=\{0,\gm_\infty\}$ is a two-point
set with a generic and a closed point, hence it is similar 
to the spectrum of a discrete valuation ring.
\item[e)] The only prime ideals in $\Spec^p\Fone[T]$ are $(0)$ and $(T)$, and 
a similar result holds for all $\Spec^p\bbF_{1^n}[T]$, i.e.\ these spectra 
also look like spectra of DVRs. Notice that the canonical strict epimorphism 
$\phi:\Fone[T]\to\bbF_{1^n}=\Fone[\zeta\,|\,\zeta^n=\bu]$, $T\mapsto\zeta$, 
induces a map ${}^a\phi:\Spec^p\bbF_{1^n}\to\Spec^p\Fone[T]$ that is 
{\em not\/} closed: indeed, its image consists of the generic point 
of~$\Spec^p\Fone[T]$.
\end{itemize}

\nxpointtoc{Localization theories}
Now we want to define the {\em localization theories\/}, which will
be sometimes also called {\em theories of spectra\/}, giving
rise to (different) spectra of generalized rings, 
and study in some detail two of them.

\nxsubpoint\label{sp:cat.psloc} (Categories of pseudolocalizations.)
Recall that a {\em pseudolocalization\/} of a generalized ring~$A$ 
is a commutative $A$-algebra~$B$, or equivalently, a homomorphism of 
generalized rings $\rho:A\to B$, such that (cf.~\ptref{sp:pseudoloc}): 
1) the scalar restriction functor 
$\rho_*:\catMod B\to\catMod A$ is fully faithful; 
this condition implies that $\rho$ is an NC-epimorphism, 
hence an epimorphism as well;
2) $B$ is a flat $A$-algebra, i.e.\ the base change functor $\rho^*$ is exact. 
If in addition $B$ is finitely presented over~$A$, we say that it is 
an {\em open pseudolocalization\/} of~$A$.

Let's fix some $A$ and consider the category $\bar\cP_A$ of 
all pseudolocalizations of~$A$ 
(considered as a full subcategory of $A$-$\catAlg$), and its full subcategory 
$\cP_A$, consisting of all open pseudolocalizations of~$A$.

Recall that for any homomorphism of generalized rings $f:A\to A'$ the base 
change functor $f^*:=A'\otimes_A-$ transforms a pseudolocalization 
$\rho:A\to B$ of~$A$ into a pseudolocalization $\rho':A'\to B'$
(cf.~\ptref{sp:prop1.psloc}); we have actually shown that 
$\catMod{B'}\subset\catMod{A'}$ can be identified with the preimage of 
$\catMod B\subset\catMod A$ under the scalar restriction functor 
$f_*:\catMod{A'}\to\catMod A$, and that any compatible 
$B$-module structure on an $A'$-module $N$ automatically commutes with 
the $A'$-module structure, hence $B'=A'\otimes_AB\cong A'\boxtimes_AB$, i.e.\ 
$B'/A'$ is the non-commutative pullback of $B/A$ as well.

We can combine the above observations to define a fibered category 
$p:\bar\cP\to\catGenR$, with the fiber over~$A$ equal to~$\bar\cP_A$, 
and the pullback functors $f^*:\bar\cP_A\to\bar\cP_{A'}$ defined for 
any $f:A\to A'$ in the natural way. We have a fibered subcategory of 
open localizations $\cP\subset\bar\cP$ as well.

\nxsubpoint (Finite inductive limits of pseudolocalizations.)
Recall that both $\cP_A$ and $\bar\cP_A$ are essentially given by preordered 
sets, i.e.\ there is at most one morphism between two 
pseudolocalizations of~$A$. Moreover, all morphisms in these categories 
are epimorphisms and even NC-epimorphisms (i.e.\ they are epimorphic in the 
category of algebraic monads as well).
Another simple observation: both these categories 
are closed under coproducts, since pseudolocalization morphisms are stable 
under base change and composition. This implies immediately that 
both $\cP_A$ and $\bar\cP_A$ are closed under finite inductive limits.

\nxsubpoint (Filtered inductive limits.)
We claim that {\em $\bar\cP_A$ is stable under filtered inductive limits.} 
Together with our previous results on finite inductive limits this implies the 
stability of pseudolocalizations under arbitrary inductive limits.

Now suppose that $B=\injlim B_\alpha$, where all $A\to B_\alpha$ are 
pseudolocalizations. Put $\cA:=\catMod A$, $\cB_\alpha:=\catMod{B_\alpha}$, 
identified with a strictly full subcategory of~$\cA$. Notice that the 
image of the scalar restriction functor $\catMod B\to\cA$ is contained in 
each $\cB_\alpha$, hence it induces a functor 
$\catMod B\to\cB:=\bigcap_\alpha\cB_\alpha$. We claim 
that it is an equivalence; this would imply that $A\to B$ is a 
pseudolocalization, flatness of arbitrary inductive limits of flat algebras 
being evident. So let $M$ be an object of~$\cB$, i.e.\ an $A$-module that 
admits a (necessarily unique) compatible $B_\alpha$-module structure 
for each~$\alpha$. This means 
that $A\to\END(M)$ (uniquely) factorizes through each $B_\alpha$, hence 
also through $B=\injlim B_\alpha$, i.e.\ $M$ admits a unique $B$-module 
structure, compatible with its $A$-module structure. 
This proves the essential surjectivity of $\catMod B\to\cB$. 
Now it remains to check that $\Hom_B(N,N')=\Hom_A(N,N')$ for 
any two $B$-modules $N$ and $N'$. But this is clear: $\Hom_B(N,N')=
\bigcap_\alpha\Hom_{B_\alpha}(N,N')$, and $\Hom_{B_\alpha}(N,N')=\Hom_A(N,N')$ 
for all~$\alpha$.

{\bf Question.} We know that any localization $S^{-1}A$ is a filtered 
inductive limit of open (i.e.\ finitely presented) localizations $A_f$. 
Is it true that any pseudolocalization is an inductive (hence also filtered 
inductive) limit of open pseudolocalizations?

\nxsubpoint\label{sp:th.spec} (Localization theories.)
A {\em localization theory~$\cT$\/} (or a {\em theory of spectra}) 
of generalized rings (resp.\ of $K$-algebras, 
where $K$ is a fixed generalized ring) is a choice of strictly 
full subcategories $\cT_A\subset\cP_A$ for each generalized ring~$A$ 
(resp.\ each $K$-algebra~$A$), such that the following conditions are 
satisfied:
\begin{itemize}
\item[0)] All isomorphisms $(A\simto A')$ belong to~$\cT_A$.
\item[1)] $\cT_A$ is stable under base change, i.e.\ if $(A\to B)$ belongs 
to $\cT_A$, then $(A'\to A'\otimes_AB)$ belongs to $\cT_{A'}$ for any 
$f:A\to A'$. In other words, $\cT$ is a full fibered subcategory of~$\cP$.
\item[2)] $\cT$ is closed under composition, i.e.\ if $(A\to B)$ belongs to 
$\cT_A$, and $(B\to C)$ belongs to $\cT_B$, then $(A\to C)$ belongs to~$\cT_A$.
\item[3)] All open unary localizations $(A\to A_f)$, $f\in|A|$, belong 
to~$\cT_A$.
\end{itemize}

We denote by $\hat\cT_A$ the full subcategory of $\bar\cP_A$, consisting 
of filtering inductive limits of pseudolocalizations from~$\cT_A$. 
This gives rise to a full fibered subcategory $\hat\cT\subset\bar\cP$.

Clearly, $\cP$ is the maximal (``finest'') possible localization theory,
while $\cL$ (consisting of all open unary localizations $A\to A_f$) 
is the minimal (``coarsest'') one.

\nxsubpoint\label{sp:bprop.th.spec} (Basic properties.)
\begin{itemize}
\item[4)] $\cT_A$ is closed under coproducts, hence also under finite 
inductive limits of~$\cP_A$. (Notice that such inductive limits in~$\cP_A$
might be different from their values when computed in $A$-$\catAlg$.)
\item[5)] $\hat\cT_A$ is closed under arbitrary inductive limits. It contains 
$\hat\cL_A$, the category of all unary localizations $A\to S^{-1}A$.
\item[6)] If $(A\to A')$ and $(A\to A'')$ both belong to~$\cT_A$, and 
$f:A'\to A''$ is an $A$-algebra homomorphism, then $(A'\stackrel f\to A'')$ 
belongs to $\cT_{A'}$.
\item[7)] For any homomorphism of generalized rings $f:A\to B$ there is a 
pseudolocalization $A\to\bar B$ in $\hat\cT_A$, such that $f$ factorizes 
through a pseudolocalization $A\to A'$ in~$\hat\cT_A$ iff 
$A'\leq\bar B$ (i.e.\ iff there is an $A$-algebra homomorphism $A'\to\bar B$). 
In other words, the inclusion functor $\hat\cT_A\to\text{$A$-$\catAlg$}$ 
admits a right adjoint $B\mapsto\bar B$.
\item[8)] Arbitrary projective limits exist in $\hat\cT_A$ (not necessarily 
coinciding with those of~$A$-$\catAlg$).
\end{itemize}

\begin{Proof} Let's prove these properties. First of all, 1) and~2) imply~4):
indeed, if $A\to A'$ and $A\to A''$ both belong to $\cT_A$, then 
$A'\to A'\otimes_AA''\cong A'\boxtimes_AA''$ belongs to $\cT_A'$ by~1), 
hence $A\to A'\to A'\otimes_AA''$ belongs to $\cT_A$ by~2).
Now 5) is an immediate consequence of~4), and 6) follows from the fact that 
all pseudolocalizations are epimorphisms, hence $f:A'\to A''$ can be identified
with the pushout of $A\to A''$ with respect to $A\to A'$, and we 
just have to apply~1). Next, to show 7) we simply put $\bar B$ to be the 
inductive limit of all pseudolocalizations $A\to A'$, such that $A\to B$ 
factorizes through this pseudolocalization, and use 5) to show that 
$A\to\bar B$ lies indeed in $\hat\cT_A$. Finally, 8) is a formal consequence 
of 7) and of existence of arbitrary projective limits in $A$-$\catAlg$: 
we first compute the required projective limit~$B$ in $A$-$\catAlg$, and then 
apply $B\mapsto\bar B$ to it.
\end{Proof}

\nxsubpoint\label{sp:constr.spec} (Spectra.)
Let $\cT^?\subset\cP$ be any localization theory, $A$ a generalized ring 
(or a $K$-algebra, if $\cT^?$ is defined only for $K$-algebras). 
We construct the {\em (strong) $\cT^?$-spectrum of~$A$}, denoted by 
$\Spec^{\cT^?}_sA$, $\Spec^?_sA$, or even $\Spec^?A$, as follows.

a) We consider the opposite category $\cS^?_A:=(\cT^?_A)^{op}$ to $\cT^?_A$. 
Clearly, finite projective limits exist in~$\cS^?_A$, and all morphisms 
of $\cS^?_A$ are monomorphisms.

b) For any $A$-module~$M$ we define a presheaf of sets $\tilde M$ on~$\cS^?_A$,
i.e.\ a functor $\tilde M:(\cS^?_A)^{op}=\cT^?_A\to\catSets$, by putting 
$\tilde M(A'):=A'\otimes_AM$ (for now we consider only the underlying set 
of $A'\otimes_AM$). The functor $M\mapsto\tilde M$ is exact, all 
pseudolocalizations $A'$ being flat over~$A$.

c) We consider on $\cS^?_A$ the finest topology, for which all $\tilde M$ 
become sheaves. The results of SGA~4 tell us that $(U_\alpha\to U)$ is 
a covering in~$\cS^?_A$ iff
\begin{equation}\label{eq:cond.cover}
\xymatrix{\tilde M(U)\ar[r]&\prod_\alpha\tilde M(U_\alpha)
\ar@/_2pt/[r]\ar@/^2pt/[r]&\prod_{\alpha,\beta}
\tilde M(U_\alpha\cap U_\beta)}
\end{equation}
is left exact for all $A$-modules~$M$. 

d) We denote by $\Spec^?A$ the topos defined by the site~$\cS^?_A$ thus 
constructed, i.e.\ $\Spec^?A$ is the category of sheaves of sets on~$\cS^?_A$.

e) By construction all $\tilde M$ are sheaves on~$\Spec^?A$. In particular, 
this is true for all $\widetilde{A(n)}:A'\mapsto A'(n)$. 
Thus the collection of $(\widetilde{A(n)})$ defines a sheaf of generalized 
rings~$\tilde A$, characterized by $\tilde A(A')=A'$ 
(cf.~\ptref{p:algmon.over.topos}), hence $\Spec^?A:=
(\Spec^?A,\tilde A)$ is a generalized ringed topos 
(cf.~\ptref{sp:gen.ringsp.top}). We usually denote by $\sO_{\Spec^?A}$ 
the sheaf of generalized rings~$\tilde A$ just constructed, and say that 
it is {\em the structural sheaf of~$\Spec^?A$}. Since any $\tilde M(A')$ 
has a natural $A'$-module structure, we see that {\em all $\tilde M$ are 
$\sO_{\Spec^?A}$-modules.} 

f) Similarly, any (generalized) $A$-algebra~$B$ 
defines a $\sO_{\Spec^?_A}$-algebra~$\tilde B$, with individual components 
given by $\tilde B(n):=\widetilde{B(n)}$. Clearly, the value of~$\tilde B$ 
on a localization~$A'$ equals~$A'\otimes_AB$ 
(we use here the affine base change theorem~\ptref{th:aff.base.change} again).
Moreover, any $B$-module~$N$ can be treated as an~$A$-module, thus definiting 
a sheaf of~$\sO_{\Spec^?_A}$-modules~$\tilde N$, and it is easy to check that 
$\tilde N$ is in fact a sheaf of $\tilde B$-modules.

\nxsubpoint (Set-theoretical issues.)
Notice that if $A$ is $\univU$-small (this is what we've tacitly assumed
everywhere in this text), then the set of finitely presented $A$-algebras
is also $\univU$-small, since any such algebra admits a description
in terms of a finite list of generators (with appropriate arities)
and relations. In particular, the set of all open pseudolocalizations
of~$A$ is $\univU$-small, hence the same holds for 
$\cT^?_A$ and $\cS^?_A$, i.e.\ the site $\cS^?_A$ is $\univU$-small,
so we can safely consider corresponding topos $\Spec^?A$.

\nxsubpoint\label{sp:top.of.spec} (Topology of~$\cS^?_A$.)
The attentive reader may have noticed that SGA~4 actually gives 
a stronger requirement for a family $(U_\alpha\to U)$ 
to be a covering with respect to the finest topology on $~\cS^?_A$, 
for which all $\tilde M$ become sheaves. Namely, a family as above is 
a covering iff for any $V\in\Ob\cS^?_A$ and any $M\in\Ob\catMod A$ the 
following sequence of sets is left exact:
\begin{equation}\label{eq:cond.cover2}
\xymatrix{\tilde M(V\cap U)\ar[r]&\prod_\alpha\tilde M(V\cap U_\alpha)
\ar@/_2pt/[r]\ar@/^2pt/[r]&\prod_{\alpha,\beta}
\tilde M(V\cap U_\alpha\cap U_\beta)}
\end{equation}
However, in our case \eqref{eq:cond.cover2} is fulfilled for all $V$ and $M$ 
iff \eqref{eq:cond.cover} is fulfilled for all~$M$, so the description of 
covers in~$\cS^?_A$ given above in~\ptref{sp:constr.spec},c) turns out 
to be accurate.

Indeed, let's deduce \eqref{eq:cond.cover2} from~\eqref{eq:cond.cover}. 
An object $V\in\Ob\cS^?_A$ is given by some pseudolocalization 
$A\stackrel\rho\to A'$ of $\cT^?_A$. Put $N:=\rho_*\rho^*M$. 
According to~\ptref{th:aff.base.change},
$\rho_*$ commutes with any {\em flat\/} base change. This is also true 
for~$\rho^*$ for trivial reasons, and the definition of $\tilde M$ now implies 
$\tilde M(V\cap U)\cong\tilde N(U)$ for any $U\in\Ob\cS^?_A$. Therefore, 
\eqref{eq:cond.cover} for~$N$ implies \eqref{eq:cond.cover2} for~$M$ and~$V$.

\nxsubpoint 
We claim that {\em the topology of~$\cS^?_A$ is subcanonical\/}, i.e.\ 
all representable presheaves $\Hom(-,V)$ on $\cS^?_A$ are actually sheaves. 
This is easily reduced to proving the following: {\em if $(U_\alpha\to U)$ 
is a covering, and all $U_\alpha\subset V$, then $U\subset V$.}

{\bf Proof.} Let's denote $\sO_{\Spec^?A}$ simply by~$\sO$. Then for any 
$V\in\Ob\cS^?_A$ the generalized ring of sections $\sO(V)$ is nothing else 
than the pseudolocalization of~$A$ corresponding to~$V$, with the 
$A$-algebra structure given by the restriction map $A=\sO(e)\to\sO(V)$, 
where $e$ is the final object of~$\cS^?_A$ (corresponding to the trivial 
pseudolocalization $A\to A$). Recall that $\sO$ is a sheaf of generalized 
rings on~$\cS^?_A$; this implies the exactness of the following 
diagram of generalized rings (or $A$-algebras):
\begin{equation}
\xymatrix{\sO(U)\ar[r]&\prod_\alpha\sO(U_\alpha)\ar@/^2pt/[r]\ar@/_2pt/[r]&
\prod_{\alpha,\beta}\sO(U_\alpha\cap U_\beta)}
\end{equation}
Now $U_\alpha\subset V$ means the existence of a (necessarily unique) 
$A$-algebra homomorphism $f_\alpha:\sO(V)\to\sO(U_\alpha)$,
 and uniqueness shows that the two maps $\sO(V)\to\sO(U_\alpha\cap U_\beta)$ 
obtained from~$f_\alpha$ and~$f_\beta$ coincide. Left exactness of the 
above diagram now implies the existence of an $A$-algebra homomorphism 
$\sO(V)\to\sO(U)$, i.e.\ $U\subset V$, q.e.d.

\nxsubpoint\label{sp:funct.spec} (Functoriality.)
Let $f:A\to B$ be a homomorphism of generalized rings (or $K$-algebras, 
if $\cT^?$ is defined over $K$-$\catAlg$). We want to construct a 
morphism of ringed topoi $\af:\Spec^?B\to\Spec^?A$, enjoying all the 
usual properties. First of all, $f$ induces a base change functor 
$f^*:\cT^?_A\to\cT^?_B$, that can be considered as a functor 
$f^{-1}:\cS^?_A\to\cS^?_B$ between opposite categories. This functor 
$f^{-1}$ is left exact by construction; we have to check that it is a 
morphism of sites, i.e.\ that $(f^{-1}U_\alpha\to f^{-1}U)$ is a cover 
in $\cS^?_B$ whenever $(U_\alpha\to U)$ is a cover in~$\cS^?_A$. 

Fix a cover $(U_\alpha\to U)$ and let's prove that $(f^{-1}U_\alpha\to
f^{-1}U)$ is a cover as well. By definition, this means left exactness 
of the counterpart of~\eqref{eq:cond.cover} for $\tilde N$, 
$N$ any $B$-module, and the family $(f^{-1}U_\alpha\to f^{-1}U)$. 
Consider $M:=f_*N$, the scalar restriction of~$N$. We claim that 
$\tilde M(U)\cong\tilde N(f^{-1}U)$ for any $U\in\Ob\cS^?_A$. 
This would show that \eqref{eq:cond.cover} for $\tilde N$ and 
$(f^{-1}U_\alpha\to f^{-1}U)$ is equivalent to~\eqref{eq:cond.cover} 
for $\tilde M$ and $(U_\alpha\to U)$, valid by definition of the topology 
of~$\cS^?_A$.

So let's prove $\tilde M(U)\cong\tilde N(f^{-1}U)$. Notice that $U$ 
corresponds to some pseudolocalization $A\to A'$. Then $f^{-1}U$ 
corresponds to $B\to B':=A'\otimes_AB$, and 
$\tilde M(U)=A'\otimes_A M=A'\otimes_A f_*N$, $\tilde N(f^{-1}U)=
B'\otimes_BN$. Since $A'$ is flat over $A$, our statement now follows from 
the ``affine base change theorem''~\ptref{th:aff.base.change}.

Notice that we have a natural homomorphism $\theta:\sO_{\Spec^?A}\to
(\af)_*\sO_{\Spec^?B}$ of structural sheaves as well, defined by 
the canonical maps $\theta_U:\sO_{\Spec^?A}(U)=A'\to\sO_{\Spec^?B}(f^{-1}U)=
B\otimes_AA'$. Therefore, $\af=(\af,\theta):\Spec^?B\to\Spec^?A$ 
can be considered as a morphism of generalized ringed topoi.

\nxsubpoint\label{sp:spec.psloc.open} (Spectra of pseudolocalizations.)
Let $f:A\to A'$ be a pseudolocalization of~$\cT^?_A$. It corresponds to some 
object~$U$ in~$\cS^?_A$. Notice that $U$ is an {\em open object\/} 
of~$\cS^?_A$, i.e.\ a subobject of the final object $e_U$. Therefore, 
we can consider the open subtopos $\Spec^?A_{/U}$ of $\Spec^?A$, 
defined by site $\cS^?_A/U$. Notice that $\cS^?_A/U$ is 
canonically equivalent to~$\cS^?_{A'}$ (cf.\ property 5) 
of~\ptref{sp:bprop.th.spec}). We claim that this equivalence is compatible 
with the natural topologies on these sites. Indeed, the topology 
of $\cS^?_A/U$ is easily seen to be the finest topology, for which all 
$\tilde M|_U$, $M\in\Ob\catMod A$, become sheaves. On the other hand, 
the topology of $\cS^?_{A'}$ is the finest topology, for which all $\tilde N$, 
$N\in\Ob\catMod{A'}$, become sheaves. 
Now it is sufficient to check that these two families of sheaves of sets 
correspond to each other under equivalence $\cS^?_A/U\cong\cS^?_{A'}$. 
This is immediate: $\tilde M|_U$ is isomorphic to $\widetilde{f^*M}$; 
conversely, if we start with some $N\in\Ob\catMod{A'}$, then 
$\tilde N$ is isomorphic to $\widetilde{f_*N}|_U$.

Clearly, $\sO_{\Spec^?A}|_{U}$ is canonically isomorphic to $\sO_{\Spec^?A'}$, 
if we identify $\Spec^?A'$ with the open subtopos of~$\Spec^?A$ defined 
by~$U$. In other words, $\Spec^?A'$ is equivalent to the open ringed subtopos 
$\Spec^?A_{/U}$ of $\Spec^?A$. The morphism of ringed topoi 
$\Spec^?A'\simto\Spec^?A_{/U}\to\Spec^?A$ is easily seen to be isomorphic 
to~$\af$, i.e.\ {\em if $f:A\to A'$ is a pseudolocalization of~$\cT^?_A$, 
$\af:\Spec^?A'\to\Spec^?A$ is an open embedding of generalized ringed topoi.}

\nxsubpoint (Quasicoherent sheaves.)
Fix $X=(X,\sO_X):=\Spec^?A$ and consider the global sections functor 
$\Gamma=\Gamma_X:=\Gamma(X,-)$, given by $\Gamma(X,\sF):=\sF(e)$, where 
$e$ is the final object of $\Spec^?A$ or of~$\cS^?_A$. Since 
$\Gamma_X(\sO_X)=A$, $\Gamma_X$ induces a functor from the category of 
(sheaves of) $\sO_X$-modules into the category of $A$-modules: 
$\Gamma_X:\catMod{\sO_X}\to\catMod A$. We have a functor $\Delta=\Delta_X$ 
in the opposite direction as well, given by $\Delta_X(M):=\tilde M$.

It is easy to see that $\Delta_X$ is a left adjoint to~$\Gamma_X$, i.e.\ 
\begin{equation}
\Hom_{\sO_X}(\tilde M,\sF)\cong\Hom_A(M,\Gamma(X,\sF))
\end{equation}
Indeed, this follows immediately from the universal property of 
$\tilde M(U)=\sO_X(U)\otimes_AM$.

Next, it is immediate that $\Gamma(X,\tilde M)=M$, i.e.\ $\Gamma_X\circ\Delta_X
\cong\Id$; this means that $\Delta_X$ is fully faithful, hence {\em 
$\Delta_X$ induces an equivalence between $\catMod A$ and the essential image 
of~$\Delta_X$, and the quasiinverse equivalence is given by 
the restriction of $\Gamma_X$ to this essential image.}

We say that an $\sO_X$-module~$\sF$ is {\em quasicoherent\/} if it lies 
in the essential image of~$\Delta_X$, i.e.\ if it is isomorphic 
to some~$\tilde M$. We have just shown that {\em the category of quasicoherent
sheaves of $\Spec^?A$ is equivalent to $\catMod A$.} On the other hand, 
$\Delta_X$ is exact and commutes with arbitrary inductive limits, hence 
{\em the category of quasicoherent sheaves is closed under finite projective 
and arbitrary inductive limits of~$\catMod{\sO_X}$.}

\nxsubpoint\label{sp:qcoh.func.strtop} 
(Quasicoherent sheaves and functoriality.)
Let $f:A\to B$ be a homomorphism of generalized rings, 
$\af:Y:=\Spec^?B\to X:=\Spec^?A$ be the associated morphism of spectra. 
Notice that we have actually proved in~\ptref{sp:funct.spec} that 
$(\af)_*\circ\Delta_Y\cong\Delta_X\circ f_*$, i.e.\ 
$(\af)_*(\tilde N)\cong\widetilde{f_*N}$ for any $B$-module~$N$.
In particular, {\em the direct image functor $(\af)_*$ transforms 
quasicoherent sheaves on~$\Spec^?B$ into quasicoherent sheaves on~$\Spec^?A$.}

On the other hand, computing the left adjoint to $\Gamma_X\circ(\af)_*=
f_*\circ\Gamma_Y$ we obtain $(\af)^*\circ\Delta_X\cong\Delta_Y\circ f^*$, 
i.e.\ $(\af)^*\tilde M\cong\widetilde{f^*M}$ for any $A$-module~$M$. 
In particular, {\em the pullback functor $(\af)^*$ transforms 
quasicoherent $\sO_X$-modules into quasicoherent $\sO_Y$-modules.}

\nxsubpoint (Flat topology.)
We can generalize the above constructions as follows. We can consider 
category~$\cS_A^{fl}$, opposite to the category of flat $A$-algebras, 
and define there a topology by requiring all presheaves of sets 
$\tilde M:B\to B\otimes_AM$ to be sheaves. The description of covers given 
in~\ptref{sp:constr.spec},c) is still valid, since the reasoning 
of~\ptref{sp:top.of.spec} is still applicable. 
(Problem: not all fibered products exist in~$\cS_A^{fl}$.)

One can check that $\cS^A_{fl}$ corresponds to the small flat site 
of~$\Spec A$ (with respect to the (fpqc) topology) when $A$ is a classical 
commutative ring. That's why we say that {\em $\cS_A^{fl}$ is the small flat 
site of~$A$}, and the corresponding topos {\em $(\Spec A)_{fl}=\Spec^{fl}A:=
\tilde\cS_A^{fl}$ is the small flat topos of~$A$.}

If we choose some localization theory~$\cT^?$ over~$A$, we get a natural 
inclusion functor $\cS^?_A\to\cS_A^{fl}$, easily seen to be left exact 
and to transform covers into covers (more precisely, $(U_\alpha\to U)$ 
is a cover in~$\cS^?_A$ iff it is a cover in~$\cS^{fl}_A$), i.e.\ 
a morphism of sites. It induces a morphism of topoi $(\Spec A)_{fl}\to
\Spec^?A$. We have a natural sheaf of generalized rings $\sO$ or~$\sO_X$ on 
$(\Spec A)_{fl}$, that transforms any object of $\cS^{fl}_A$, given by 
some flat $A$-algebra~$B$, into~$B$ itself, and all $\tilde M$ are actually 
sheaves of $\sO$-modules on~$(\Spec A)_{fl}$.

\nxsubpoint\label{sp:comp.mor.strong.spec} (Comparison morphisms.)
Given two localization theories $\cT'$ and $\cT''$, we say that 
{\em $\cT''$ is a refinement of $\cT'$\/} if $\cT'\subset\cT''$, i.e.\ any 
open pseudolocalization of~$\cT'$ lies also in~$\cT''$. This defines a 
partial order between localization theories. Clearly, 
the theory $\cT^u$ consisting only of unary localizations of form $A\to A_f$ 
is the smallest element with respect to this order, while 
the theory $\cP$ consisting of all open pseudolocalizations is the largest 
element. The latter will be called the {\em total localization theory 
$\cT^t:=\cP$, and corresponding spectra will be denoted by~$\Spec^tA$.}

In any case, given two such theories $\cT'\subset\cT''$ as above, we obtain 
canonical inclusion functors $\cS_A'\to\cS_A''$; they are left exact and 
preserve covers, hence we have a morphism of sites and of topoi in the 
opposite direction $\Spec''A\to\Spec'A$, easily seen to be a morphism of 
generalized ringed topoi, depending functorially on~$A$. We say that 
{\em $\Spec''A\to\Spec'A$ is a comparison morphism with respect to these two 
theories.}

\nxsubpoint (Unary localization spectra $\Spec^uA$.)
Let us study spectra $\Spec^uA$ with respect to the smallest possible 
theory~$\cT^u$. Clearly, we have comparison morphisms $\Spec^?A\to\Spec^uA$ 
for any other theory~$\cT^?$; this explains the importance of~$\Spec^uA$.
First of all, the site $\cS^u_A$ is the opposite category to $\cT^u_A$, 
hence it consists (up to equivalence) of objects $D^u(f)$, $f\in|A|$, 
that correspond to unary localizations $A\to A_f$. We know that 
$\cS^u_A$ is a preordered set, i.e.\ between any two objects $D^u(f)$ and
$D^u(g)$ there is at most one morphism; when such a morphism exists, we write 
$D^u(f)\subset D^u(g)$. By definition $D^u(f)\subset D^u(g)$ iff 
$A\to A_f$ factorizes through $A\to A_g$ iff $f$ divides some power of~$g$ 
iff $D(f)\subset D(g)$ in the prime spectrum $\Spec^pA$ 
(cf.~\ptref{sp:opcov.Dg}). Therefore, $\cS^u_A$ can be identified with the 
ordered set of principal open subsets of the prime spectrum~$\Spec^pA$ 
of~\ptref{sp:prime.spec}.

\nxsubpoint\label{sp:covers.un.spec} (Covers in $\Spec^uA$.) 
We claim that {\em $(D^u(f_\alpha)\to
D^u(g))$ is a cover in $\cS^u_A$ iff $(D(f_\alpha)\to D(g))$ is a cover
in~$\Spec^pA$, i.e.\ iff $D(g)=\bigcup_\alpha D(f_\alpha)$.} Since the
principal open subsets form a base of topology of $\Spec^pA$, this would prove
that {\em $\Spec^uA$ coincides with the topos defined by topological
space~$\Spec^pA$}. Moreover, this would imply that $\Spec^uA$ has enough
points, and that the topological space associated to $\Spec^uA$ is canonically
homeomorphic to~$\Spec^pA$, the latter topological space being sober.

So let us show our statement about covers. Since replacing $A$ with~$A_g$
replaces both $\Spec^uA$ and $\Spec^pA$ with corresponding principal open
subsets or subtopoi (cf.~\ptref{sp:specp.localiz}
and~\ptref{sp:spec.psloc.open}), we can assume that $A=A_g$, hence
$D^u(g)=D^u(1)=e$ is the final object of~$\cS^u_A$, and $D(g)=D(1)=\Spec^pA$.

$\Rightarrow)$ Suppose that $(D^u(f_\alpha)\to D^u(1))$ is a cover 
in~$\cS^u_A$, but $\bigcup_\alpha D(f_\alpha)\neq\Spec^pA$. According 
to~\ptref{sp:opcov.specp}, this means that the ideal $\ga:=(f_\alpha)$ 
is $\neq(1)$. Then $\ga_{f_\alpha}=A_{f_\alpha}$ for any~$\alpha$, 
and the diagrams~\eqref{eq:cond.cover} are left exact for both $\tilde\ga$ 
and $\tilde A$. Consider the following diagram with left exact rows:
\begin{equation}
\xymatrix{\tilde\ga(U)\ar[r]\ar[d]&\prod_\alpha\tilde\ga(U_\alpha)
\ar[d]^{\sim}\ar@/_2pt/[r]\ar@/^2pt/[r]&\prod_{\alpha,\beta}
\tilde\ga(U_\alpha\cap U_\beta)\ar[d]^{\sim}\\
\tilde A(U)\ar[r]&\prod_\alpha\tilde A(U_\alpha)
\ar@/_2pt/[r]\ar@/^2pt/[r]&\prod_{\alpha,\beta}
\tilde A(U_\alpha\cap U_\beta)}
\end{equation}

The vertical arrows are induced by natural embeddings; we know that the
central and the right vertical arrows are isomorphisms, hence the same is true
for the remaining vertical arrow, i.e.\ $\ga=(1)$. This is absurd.

$\Leftarrow)$ Suppose that $\bigcup_\alpha D(f_\alpha)=\Spec^pA$. We have to
show that \eqref{eq:cond.cover} is exact for $(D(f_\alpha)\to D(1))$ and 
any $A$-module~$M$. We know that $\Spec^pA$ is quasicompact, so we can find 
a finite subcover, consisting of some $D(f_i):=D(f_{\alpha_i})$,
$1\leq i\leq n$. It is easy to see that it suffices to check left exactness 
of \eqref{eq:cond.cover} for this finite subcover. This means left exactness 
of the following diagram:
\begin{equation}\label{eq:cond.cov2}
\xymatrix{M\ar[r]&\prod_i M_{f_i}\ar@/_2pt/[r]\ar@/^2pt/[r]&
\prod_{i,j}M_{f_if_j}}
\end{equation}
provided the $f_i$ generate the unit ideal in~$A$
(cf.~\ptref{sp:opcov.specp}).

a) First of all, let us prove that if $x$, $y\in M$ are such that 
$x/1=y/1$ in any~$M_{f_i}$, then $x=y$. Indeed, $x/1=y/1$ in~$M_{f_i}$ means 
that $xf_i^N=yf_i^N$ for some $N>0$, and we can choose the same~$N$ for all 
values of~$i$, our cover being finite. We know that the elements 
$f_i^N$ still generate the unit ideal (since $\bigcup D(f_i^N)=\bigcup D(f_i)=
\Spec^pA$), so we can find an $H\in A(n)$, such that 
$H(f_1^N,\ldots,f_n^N)=\bu$. Then we obtain 
$x=\bu\cdot x=H(f_1^Nx,\ldots,f_n^Nx)=H(f_1^Ny,\ldots,f_n^Ny)=y$.

b) Now suppose we are given a family of elements $x_i/f_i^N$ in all $M_{f_i}$,
such that $x_i/f_i^N=x_j/f_j^N$ in $M_{f_if_j}$ for any $i$ and~$j$.  We use
again the finiteness of our cover to choose a common value of~$N$ for
all~$i$. Let us prove that $x_i/f_i^N=x/1$ for a suitably chosen $x\in M$.  We
know that $x_i/f_i^N=x_j/f_j^N$ means $x_if_i^Kf_j^{N+K}=x_jf_i^{N+K}f_j^K$
for some integer $K\geq 0$, chosen simultaneously for all values of~$i$
and~$j$. Replacing $x_i$ with $x_if_i^K$ and $N$ with $N+K$, we can assume 
$K=0$, i.e.\ $x_if_j^N=x_jf_i^N$.

By the same reasoning as before we can find some $H\in A(n)$, such
that $H(f_1^N,\ldots,f_n^N)=\bu$. Put $x:=H(x_1,\ldots,x_n)$.  For any~$i$ we
have $xf_i^N=H(x_1,\ldots,x_n)\cdot f_i^N=H(x_1f_i^N,\ldots,x_nf_i^N)=
H(f_1^Nx_i,\ldots,f_n^Nx_i)=H(f_1^N,\ldots,f_n^N)\cdot x_i=\bu x_i=x_i$, hence 
$x_i/f_i^N=x/1$, q.e.d.

\nxsubpoint (Open pseudolocalizations of classical rings.)  Let $f:A\to B$ be
a homomorphism of classical commutative rings. We claim that {\em $f$ is an
open pseudolocalization iff $\af:\Spec B\to\Spec A$ is an open embedding.} The
``if'' part is very easy to check, since $\af$ turns out to be flat and
finitely presented, and $\af_*$ induces a fully faithful functor from the
category of quasicoherent sheaves on~$\Spec B$ into the category of
quasicoherent sheaves on~$\Spec A$. Conversely, suppose that $f$ is an open 
pseudolocalization. This implies that $B$ is flat and finitely presented 
over~$A$, hence $\af$ is an open map, and in particular $U:=\af(\Spec B)$ is a
quasicompact open subset of~$\Spec A$. We have to check that the induced map 
$\Spec B\to U$ is an isomorphism; this is equivalent to checking that 
$\af_{D(s)}:\af^{-1}(D(s))\to D(s)$ is an isomorphism whenever $D(s)\subset
U$, $s\in A$, i.e.\ that $f$ becomes an isomorphism after any base change
$A\to A':=A_s$ for any such $s\in A$. Indeed, put $B':=A'\otimes_AB=B_s$. Then 
$f':A'\to B'$ is an open pseudolocalization again, but now $\af'$ is both flat 
and surjective, hence faithfully flat, i.e.\ $B'$ is a faithfully flat
$A'$-algebra. On the other hand, we know that $B'\otimes_{A'}B'\cong B'$,
since any pseudolocalization is an NC-epimorphism. In other words, $A'\to B'$
becomes an isomorphism after tensoring with $B'$ over~$A'$, hence it is itself
an isomorphism by faithful flatness, q.e.d.

\nxsubpoint\label{sp:str.spec.add.rings} 
(Comparison morphisms for spectra of classical rings.)  One can
easily check that a family of open pseudolocalizations $(B\to B_\alpha)$
defines a cover in the sense of~\ptref{sp:constr.spec},c) iff $\Spec B$ is the
union of its open subsets $\Spec B_\alpha$. An immediate consequence is that
{\em for any localization theory $\cT^?$ and any classical ring~$A$ the
comparison morphism $\Spec^?A\to\Spec^uA=\Spec A$ is an equivalence of ringed
topoi.} In other words, all localization theories yield the same result for
classical rings; we might have defined $\Spec A$ to be equal to $\Spec^tA$ as
well.

One can check essentially in the same way that whenever we have a
localization theory $\cT^?$ (say, for generalized $K$-algeras), such
that for any $K$-algebra $A$ and any finitely presented flat
$A$-algebra $B$ induced morphism $\Spec^?B\to\Spec^?A$ has open image
(i.e.\ decomposes into a faithfully flat morphism followed by an open
embedding), then the comparison morphisms $\Spec^tA\to\Spec^?A$ are
equivalences, i.e.\ any such theory can be actually replaced by the
total localization theory.

\nxsubpoint (Localizations vs.\ pseudolocalizations of classical rings.)
Notice that an open pseudolocalization $A\to B$ of classical rings need not be
an open localization $A\to A_s$, i.e.\ not all affine open subsets of $\Spec
A$ are principal open subsets. Let us illustrate the subtlety of the notion of
open pseudolocalizations by a seemingly simple question.

{\bf Question.} Let $A:=k[T_1,\ldots,T_n]$ be a polynomial ring over a
field~$k$. Is there a normal $k$-subalgebra $B\subset A$, distinct from~$A$,
such that $B\to A$ is a) an open localization; b) an open pseudolocalization?

While the answer to a) is obviously negative (all invertible elements of~$A$
lie in~$k$, so we cannot have $A=B_s$ for $s\not\in k$), the negative answer
to b) would in fact imply the jacobian conjecture, so the second part of the
question has to be very subtle.

\nxsubpoint (Points and quasicompactness of spectra.)
One can check directly or by invoking general results of SGA 4 that once we
know that all spectra for a certain theory~$\cT^?$ are quasicompact, they
admit enough points (being coherent; cf.\ SGA~4), and in any case they are
generated by open subobjects, hence they correspond to some sober topological
spaces, uniquely determined up to a canonical homeomorphism. These sober
topological spaces will be also denoted by $\Spec^?A$; they are generalized
ringed topological spaces, of course.

The only problem here is to prove the quasicompactness. I don't know whether
quasicompactness holds for all theories and all generalized rings,
so we have to circumvent this problem in another way.

\nxpointtoc{Weak topology and quasicoherent sheaves}
The ``strong'' topology just defined on $\cS^?_A=(\cT^?_A)^{op}$
and corresponding ``strong'' spectra $\Spec^?A=\Spec^?_sA$,
where $A$ is a generalized ring and $\cT^?$ is any localization theory,
has a serious drawback. Namely, when $X:=\Spec^?_A$ is not quasicompact
(we know already that this can happen only if $A$ is non-additive
and $\cT^?\neq\cT^u$), then the property of an $\sO_X$-module $\sF$
to be quasicoherent (i.e.\ of the form $\tilde M$ for some $A$-module~$M$)
is not necessarily local. Therefore, even if we construct ``strong'' schemes
by gluing together strong spectra (for some fixed theory $\cT^?$)
along open subsets, we won't obtain a reasonable notion of quasicoherent
sheaves over such schemes.

\zerosubpoint\nxsubpoint (Terminology and notations.)
In order to overcome this difficulty we consider two weaker topologies
on the same category $\cS^?_A=(\cT^?_A)^{op}$, namely, the
{\em weak\/} and the {\em finite\/} topology. Weak topology will
be weaker than the strong topology considered so far, and the finite topology
will be even weaker than the weak topology. Both these topologies have
correct properties with respect to quasicoherent sheaves, and the
finite topology has another advantage: resulting spectra turn out to
be quasicompact, hence also coherent in the sense of SGA~4, hence
they admit enough points, i.e.\ they are given by some sober topological
spaces, so if we restrict our attention to finite spectra, we can 
work with generalized ringed topological spaces and forget about
topoi.

{\bf Notation.} We denote by $\cS^?_{A,s}$, $\cS^?_{A,w}$ or $\cS^?_{A,f}$
the category $\cS^?_A$, considered as a site by means of its strong,
weak or finite topology, respectively. Corresponding topoi
(or topological spaces, when these topoi have enough points)
will be denoted by $\Spec^?_sA$, $\Spec^?_wA$ or $\Spec^?_fA$, 
respectively.

If $A$ is additive ({\em additive\/} actually means here ``additive with 
a compatible symmetry'', i.e.\ ``given by a classical ring''), 
or if $\cT^?$ is the unary localization theory $\cT^u$,
then all of the above topologies will coincide (cf.~\ptref{sp:comp.swf.spec}), 
so we'll still write $\cS^?_A$ and $\Spec^?A$, or even $\Spec A$.

\nxsubpoint (Fibered category of quasicoherent $\sO$-modules.)
Let us fix a generalized ring $A$ and a localization theory $\cT^?$. 
Put $\cS^?_A:=(\cT^?_A)^{op}$ as before. For any $U\in\Ob\cS^?_A$,
given by some open pseudolocalization $A\to A'$, denote by
$\stQCOH(U)$ the category of quasicoherent modules over $U$, or,
equivalently, the category of $A'$-modules: $\stQCOH(U)=\catMod{A'}$.
Whenever $V\subset U$, we have a base change (scalar extension) functor
$\stQCOH(U)\to\stQCOH(V)$, and these base change functors are transitive
up to a canonical isomorphism. This means that the collection of
categories $\{\stQCOH(U)\}$ together with these base change functors
defines a fibered category $\stQCOH=\stQCOH^?_A\to\cS^?_A$ over
category $\cS^?_A$ (cf.~\ptref{p:gen.on.stacks}, where definitions of
fibered categories and stacks are recalled, or
\cite[VI]{SGA1} and~\cite{Giraud}).

\nxsubpoint (Descent families.)
Let $\{U_\alpha\to U\}$ be any family of morphisms in category~$\cS^?_A$.
Denote by $\catDesc(\{U_\alpha\to U\},\stQCOH)$ the category of descent
data in $\stQCOH$ with respect to this family; such a descent datum
is a family of modules $M_\alpha\in\stQCOH(U_\alpha)$ together
with isomorphisms $\theta_{\alpha\beta}:M_\alpha|_{U_{\alpha\beta}}\simto
M_\beta|_{U_{\alpha\beta}}$ satisfying the cocycle relation, where
$U_{\alpha\beta}=U_\alpha\times_UU_\beta=U_\alpha\cap U_\beta$ in
our situation, and the ``restriction'' $M|_V$ denotes the image of $M$
under appropriate pullback, i.e.\ scalar extension functor.

We have a natural functor $\stQCOH(U)\to\catDesc(\{U_\alpha\to U\},
\stQCOH)$. Recall that $\{U_\alpha\to U\}$ is said to be 
a {\em family of descent\/} (resp.\ {\em of efficient descent})
for $\stQCOH$ if this functor is fully faithful (resp.\ an
equivalence of categories). If we fix a topology on $\cS^?_A$, then
$\stQCOH$ is said to be a {\em prestack\/} (resp.\ {\em stack})
if any cover for this topology is a family of descent (resp.\ of
efficient descent) for~$\stQCOH$. Finally, recall that all
{\em universal\/} families of descent (resp.\ of efficient descent)
for a fibered category, $\stQCOH$ in our case, constitute a topology
on $\cS^?_A$, called the descent (resp.\ efficient descent) topology
for $\stQCOH$ (cf.~\cite{Giraud}). This is the strongest topology
on $\cS^?_A$, for which $\stQCOH$ is a prestack (resp.\ a stack).

\nxsubpoint (Structural presheaf of $\cS^?_A$.)
Notice that some of the statements made before about sheaves on $\cS^?_A$
with respect to the strong topology are still valid as statements
about presheaves, without any choice of topology. For example,
we still have the {\em structural presheaf (of generalized rings)}
$\sO$ on $\cS^?_A$, which maps $U\in\Ob\cS^?_A$ given by some
open pseudolocalization $A\to A'$ into $A'$, and any $A$-module $M$
defines a {\em quasicoherent presheaf of $\sO$-modules $\tilde M$}
by $\tilde M(U):=M\otimes_AA'$. Furthermore, $\cS^?_A/U$ can be
identified with $\cS^?_{A'}$, and $\tilde M_{/U}$ is identified then
with $\widetilde{M|_U}$, where $M|_U=M\otimes_AA'$. 
Finally, $\Delta:M\mapsto\tilde M$
is still the left adjoint to global sections functor $\Gamma:\sF\mapsto
\sF(e)$, i.e.\
\begin{equation}
\Hom_\sO(\tilde M,\sF)\cong\Hom_A(M,\Gamma(\sF))
\quad\text{for any presheaf of $\sO$-modules~$\sF$.}
\end{equation}
Since $\Gamma\Delta=\Id_{\catMod A}$, we see that $\Delta$ is still fully
faithful, i.e.\ $\Hom_\sO(\tilde M,\tilde N)\cong\Hom_A(M,N)$ for
any two $A$-modules $M$ and~$N$.

\nxsubpoint (Descent topology on $\cS^?_A$.)
We claim that {\em the descent topology on $\cS^?_A$ with respect to
{\rm $\stQCOH$} coincides with the strong topology of~\ptref{sp:constr.spec}}.
More precisely, {\em a family $\{U_\alpha\to U\}$ in $\cS^?_A$ is
a family of descent for {\rm $\stQCOH$} iff 
it is a cover for the strong topology,
i.e.\ iff it satisfies the condition of \ptref{sp:constr.spec},c)}.
This means in particular that {\em any family of descent for {\rm $\stQCOH$}
is universal, i.e.\ remains such after any base change.}

(a) For any $U\in\Ob\cS^?_A$ and any $M$, $N\in\Ob\stQCOH(U)=\catMod{A'}$,
where $U$ is given by $A\to A'$, we define a presheaf of sets
$\iHom_U(M,N)$ on $\cS^?_A/U$ by mapping $V\to U$ into
$\Hom_{\stQCOH(V)}(M|_V,N|_V)$. It is a general fact that
$\{U_\alpha\to U\}$ is a family of descent for $\stQCOH$ iff
all $\iHom_U(M,N)$ satisfy the sheaf condition for this cover.
On the other hand, $(\iHom_U(M,N))(V)=\Hom_{\stQCOH(V)}(M|_V,N|_V)=
\Hom_{\sO(V)}(M|_V,N|_V)=\Hom_{\sO|_V}(\widetilde{M|_V},\widetilde{N|_V})=
\Hom_{\sO_V}(\tilde M|_V,\tilde N|_V)$, i.e.\ $\iHom_U(M,N)$ coincides with
the usual local Hom $\iHom_{\sO_U}(\tilde M,\tilde N)$.

(b) Now consider the strong topology on $\cS^?_A$. By definition any
$\tilde N$ is a sheaf over $\cS^?_A/U$, hence $\iHom_{\sO_U}(\tilde M,
\tilde N)$ is a sheaf as well, hence any cover $\{U_\alpha\to U\}$ 
for the strong topology is a descent family for~$\stQCOH$. 

(c) Conversely, let $\{U_\alpha\to U\}$ be a descent family for~$\stQCOH$.
Then by (a) the presheaf $\iHom_{\sO_U}(\tilde A|_U,\tilde M|_U)$
satisfies the sheaf condition for this family.
On the other hand, this presheaf is canonically isomorphic to $\tilde M|_U$
since $\tilde A=\sO$, i.e.\ $\tilde M$ satisfies the sheaf condition for
$\{U_\alpha\to U\}$. By \ptref{sp:constr.spec},c) this means that
this family is a cover for the strong topology, q.e.d.

\begin{DefD}\label{def:weak.spec}
Let $\cT^?$ be a localization theory and $A$ be a generalized ring. Define
$\cS^?_A:=(\cT^?_A)^{op}$ and {\rm $\stQCOH\to\cS^?_A$} as before. The
{\em weak topology\/} on $\cS^?_A$ is the efficient descent topology
with respect to~{\rm $\stQCOH$}. Category $\cS^?_A$, considered as a site with
respect to this topology, will be denoted by $\cS^?_{A,w}$.
The {\em weak $\cT^?$-spectrum of~$A$} is the corresponding topos
$\Spec^?_wA$, considered as a generalized ringed topos with structural
sheaf~$\sO$. A sheaf of $\sO$-modules $\sF$ is said to be {\em quasicoherent\/}
if it is isomorphic to $\tilde M$ for some $A$-module~$M$.
\end{DefD}
Since any universal family of efficient descent is a universal family of
descent, the weak topology is indeed weaker than the strong topology
of~\ptref{sp:constr.spec}. Therefore, any sheaf for the strong topology
on $\cS^?_A$ is a sheaf for the weak topology as well. This is
applicable in particular to $\sO$ and $\tilde M$, i.e.\ the structural
presheaf and the quasicoherent presheaves are indeed sheaves for the
weak topology.

\nxsubpoint\label{sp:weak.open.loc} (Restriction to open objects.)
Of course, for any $U\in\Ob\cS^?_A$ given by some open pseudolocalization
$A\to A'$ the weak spectrum $\Spec_w^?A'$ can be identified (as a 
generalized ringed topos) with the restriction $\Spec_w^?A_{/U}$ of
$\Spec_w^?A$ to its open object represented by~$U$. This observation
enables us to assume $U=e$, $A'=A$ while studying covers or families
of morphisms $\{U_\alpha\to U\}$.

\nxsubpoint\label{sp:qcoh.loc.for.weak} (Quasicoherence is a local property.)
Notice that the quasicoherence is a local property with respect to the weak
topology. In other words, if $\{U_\alpha\to e\}$ is a cover of the
final object of $\cS^?_A$ with respect to the weak topology, and
an $\sO$-module $\sF$ is such that each $\sF|_{U_\alpha}$ is quasicoherent,
i.e.\ $\sF|_{U_\alpha}\cong\tilde M_\alpha$ for some $\sO(U_\alpha)$-module
$M_\alpha$, then $\sF$ is itself quasicoherent. This statement 
actually follows immediately from the definition of weak topology:
$\{U_\alpha\to e\}$ is a family of efficient descent for $\stQCOH$,
and the family of $\sO(U_\alpha)$-modules $M_\alpha$ together with
isomorphisms $\theta_{\alpha\beta}:M_\alpha|_{U_{\alpha\beta}}\simto
M_\beta|_{U_{\alpha\beta}}$ induced by $\tilde M_\alpha|_{U_{\alpha\beta}}
\simto\sF|_{U_{\alpha\beta}}\simto\tilde M_\beta|_{U_{\alpha\beta}}$ 
obviously is a descent datum for $\stQCOH$ with respect to this cover.
We deduce existence of an $A$-module~$M$ such that 
$M|_{U_\alpha}\cong M_\alpha$ in a manner compatible with all 
$\theta_{\alpha\beta}$. In other words, we have a compatible family
of sheaf isomorphisms $\tilde M|_{U_\alpha}\cong\tilde M_\alpha\cong
\sF|_{U_\alpha}$, which can be glued into an isomorphism of $\sO$-modules 
$\tilde M\cong\sF$, i.e.\ $\sF$ is quasicoherent, q.e.d.

\nxsubpoint (Pullbacks of efficient descent families and 
functoriality of weak spectra.)
Let $f:A\to B$ be a homomorphism of generalized rings,
thus defining a left exact pullback functor $f^{-1}:\cS^?_A\to\cS^?_B$
as in~\ptref{sp:funct.spec}. We claim that {\em $f^{-1}$ 
preserves families of efficient descent for} $\stQCOH$. 
Since $f^{-1}$ is left exact, the same is true for universal 
efficient descent families (at least if $f$ is an open pseudolocalization
from $\cT^?_A$; however, this case already suffices 
to show~\ptref{sp:eff.desc.fam.are.univ},
which tells us that all efficient descent families are automatically
universal), 
i.e.\ $f^{-1}$ becomes a morphism of sites,
thus defining a morphism of generalized ringed topoi
${}^af:\Spec^?_wB\to\Spec^?_wA$ (the ``generalized ringed'' part
follows from corresponding strong topology result~\ptref{sp:funct.spec} 
once we show that $f^{-1}$ is a morphism of sites.)

(a) So let $\{U_\alpha\to U\}$ be an efficient descent family for $\stQCOH$
in~$\cS^?_A$. By~\ptref{sp:weak.open.loc} we may assume $U=e$.
Put $A_\alpha:=\sO(U_\alpha)$; then $f^{-1}U_\alpha$ is given by
$B\to B_\alpha:=B\otimes_AA_\alpha$, the pushout of $A\to A_\alpha$
with respect to $f:A\to B$. Let $(M_\alpha,\theta_{\alpha\beta})$
be a descent datum for $\stQCOH$ with respect to $\{f^{-1}U_\alpha\to e\}$.
Notice that $M_\alpha$ is a $B_\alpha$-module, and $\theta_{\alpha\beta}$
is an isomorphism of $B_{\alpha\beta}$-modules, where of course
$U_{\alpha\beta}:=U_\alpha\cap U_\beta$, $A_{\alpha\beta}:=
\sO(U_{\alpha\beta})=A_\alpha\otimes_AA_\beta$, $B_{\alpha\beta}=
B\otimes_AA_{\alpha\beta}$ and so on. By scalar restriction we
can treat $M_\alpha$ as $A_\alpha$-module and $\theta_{\alpha\beta}$ 
as an $A_{\alpha\beta}$-isomorphism, thus obtaining a descent datum for
$\{U_\alpha\to e\}$. By assumption any such descent datum is efficient,
so we obtain an $A$-module~$M$ together with $A_\alpha$-isomorphisms
$\gamma_\alpha:M_\alpha\simto M|_{U_\alpha}=M\otimes_AA_\alpha$. 
All we have to check
is that this $A$-module structure on~$M$ extends in a unique way to
a $B$-module structure, such that all $\gamma_\alpha$ become
$B_\alpha$-isomorphisms. (Notice that we are implicitly using here
the ``affine base change theorem''~\ptref{th:aff.base.change},
which enables us to identify $M\otimes_AA_\alpha$ with $M\otimes_BB_\alpha$
since $A\to A_\alpha$ is flat.)

(b) Let $t\in B(n)$ be any operation of~$B$. We want to show that
the $A$-homomorphism $[t]_M:M^n\to M$ is uniquely determined by our data.
Indeed, flatness of all $A\to A_\alpha$ implies $M^n|_{U_\alpha}=
(M|_{U_\alpha})^n\cong M_\alpha^n$, 
and the requirement of all $\gamma_\alpha$ to
be $B_\alpha$-homomorphisms means $[t]_M|_{U_\alpha}=[t]_{M_\alpha}$.
Therefore, the collection 
$\bigl\{[t]_{M_\alpha}:M_\alpha^n\to M_\alpha\bigr\}_\alpha$ defines
a morphism from the descent datum for $M^n$ into the descent datum for~$M$;
by descent we obtain a unique $A$-homomorphism $[t]_M:M^n\to M$ with
required properties, and all relations between operations $t$, $t'$, \dots\  
of~$B$ will be automatically fulfilled for $[t]_M$, $[t']_M$, \dots\ 
just because they hold for all $[t]_{M_\alpha}$, $[t']_{M_\alpha}$, \dots,
and because of the uniqueness of $[t]_M$. Thus $M$ admits a unique 
$B$-module structure with required properties, so our original descent
datum is efficient, q.e.d.

\nxsubpoint\label{sp:eff.desc.fam.are.univ} 
(All efficient descent families are universal.)
Applying the above result to an open pseudolocalization $A\to A'$ in~$\cT^?_A$,
corresponding to some $V\in\Ob\cS^?_A$, we see that {\em any
efficient descent family for {\rm $\stQCOH$} in $\cS^?_A$ is universal}.
Therefore, {\em covers for the weak $\cT^?$-topology on $\cS^?_A$ are
exactly the efficient descent families for} $\stQCOH$.

\nxsubpoint\label{sp:fin.desc.fam.eff} (Finite descent families are efficient.)
Clearly, any efficient descent family (for $\stQCOH$) is a descent family, 
i.e.\ any cover for the weak topology is a cover for the strong topology.
Let us prove a partial converse: {\em any {\bf finite}
descent family is efficient}, i.e.\ {\em any finite cover for the strong
topology is a cover for the weak topology.}

So let $\{U_i\to e\}_{1\leq i\leq n}$ be a finite cover of the final object
$e$ of $\cS^?_A$ with respect to the strong topology. Put
$U_{ij}:=U_i\cap U_j$, $A_i:=\sO(U_i)$ and $A_{ij}:=\sO(U_{ij})$ as before,
and consider any descent datum $(M_i,\theta_{ij})$ with respect to this
cover. For any $i$ and $j$ we have two maps $M_i\to M_{ij}=
M_i\otimes_{A_i}A_{ij}$: one is the natural scalar extension map, 
while the other is obtained by composing another scalar extension map
$M_i\to M_{ji}$ with $\theta_{ij}^{-1}:M_{ji}\simto M_{ij}$. Now we
can define an $A$-module~$M$ by means of the following left exact
diagram:
\begin{equation}
M\longto\prod_iM_i\rightrightarrows\prod_{i,j}M_{ij}
\end{equation}
Since this is a finite projective limit diagram, it remains such after
any flat base change, e.g.\ with respect to $A\to A_k$, i.e.\
$M\otimes_AA_k$ is canonically isomorphic to $\Ker(\prod_i M_{ik}
\rightrightarrows\prod_{i,j}M_{ijk})$. One easily sees that this
kernel is canonically isomorphic to $M_k=M_{kk}$, just 
because now we are solving a descent problem with respect to the
trivial cover $\{U_{ik}\to U_k\}_{1\leq i\leq n}$, hence
$M\otimes_AA_k\cong M_k$ for all $1\leq k\leq n$, i.e.\ $M$ is the
solution of our original descent problem, q.e.d.

This result leads us to the following definition:
\begin{DefD}\label{def:finite.spec} (Finite topology and spectra.)
The {\em finite topology\/} on $\cS^?_A$ is the topology generated by all
finite covers with respect to the strong (or equivalently, the weak)
topology on $\cS^?_A$. In other words, a family $\{U_\alpha\to U\}$
is a cover for the finite topology iff it admits a finite refinement
$\{V_i\to U\}$, which is an (efficient) descent family for {\rm $\stQCOH$}.
The corresponding site will be denoted by $\cS^?_{A,f}$, and the
corresponding generalized ringed topos by $\Spec^?_fA$; it will be called
the {\em finite $\cT^?$-spectrum\/} of~$A$.
\end{DefD}

\nxsubpoint (Properties of finite topology and spectra.)
All the nice properties of weak topology and weak spectra can be almost 
immediately transferred to the finite topology case since finite topology
is even weaker than the weak topology. In particular, all $\tilde M$
are sheaves of $\sO$-modules on $\Spec^?_fA$, quasicoherence is
a local property for the finite topology, and any $\phi:A\to B$
induces a morphism of generalized ringed topoi 
${}^a\phi:\Spec^?_fB\to\Spec^?_fA$.
Furthermore, all objects $U\in\Ob\cS^?_A$ are quasicompact for the finite
topology. Since $\cS^?_A$ is closed under finite projective limits,
this implies that $\Spec^?_fA$ is {\em coherent\/} in the sense
of SGA~4, hence it admits enough points by the general results of
{\em loc.\ cit.} On the other hand, $\Spec^?_fA$ is clearly generated by
its open objects (actually $U\in\Ob\cS^?_A$ would suffice), hence
the topos $\Spec^?_fA$ is defined by some (uniquely determined)
sober topological space, which will be usually denoted by $\Spec^?_fA$
as well. Therefore, the finite spectrum $\Spec^?_fA$ can be thought
of as a generalized ringed topological space, not a topos,
and this topological space will be automatically quasicompact.

\nxsubpoint\label{sp:comp.swf.spec} 
(Comparison of strong, weak and finite spectra.)
Clearly, we have morphisms of generalized ringed topoi
$\Spec^?_sA\to\Spec^?_wA\to\Spec^?_fA$, induced by the identity functor
on $\cS^?_A$. In some cases these arrows are equivalences:

(a) Suppose that any object $U\in\Ob\cS^?_A$ is quasicompact with respect
to the weak topology. Then by definition weak and finite topologies on
$\cS^?_A$ coincide, hence $\Spec^?_wA=\Spec^?_fA$. In particular,
these two topoi are given by some topological space.

(b) Similarly, if any object $U\in\Ob\cS^?_A$ is quasicompact with respect
to the strong topology, then strong, weak and finite topologies on
$\cS^?_A$ coincide, so we get $\Spec^?_sA=\Spec^?_wA=\Spec^?_fA$.

(c) This is applicable in particular for the unary localization theory
$\cT^u$ and any~$A$, since $\Spec^u_sA$ can be identified with the prime
spectrum $\Spec^pA$ (cf.~\ptref{sp:covers.un.spec}), 
which is quasicompact by~\ptref{sp:opcov.specp}, hence
$\Spec^u_fA=\Spec^u_wA=\Spec^u_sA=\Spec^pA$, i.e.\ we needn't
distinguish the three variants of spectra for the unary
localization theory.

(d) If $A$ is additive, then the comparison morphism $\Spec^?_sA\to\Spec^u_sA
=\Spec A$ is an isomorphism for any localization theory $\cT^?$ 
(cf.~\ptref{sp:str.spec.add.rings}), hence $\Spec^?A$
is quasicompact. If $A\to A'$ is any open pseudolocalization
from $\cT^?_A$, corresponding to some $U\in\Ob\cS^?_A$, 
then $A'$ is still additive, hence $U$ is quasicompact with respect
to the strong topology by the above reasoning. Applying (b)
we obtain $\Spec^?_fA=\Spec^?_wA=\Spec^?_sA=\Spec^uA=\Spec A$,
i.e.\ {\em any construction of spectra yields the prime spectrum
$\Spec A$ when applied to an additive (i.e.\ classical) commutative
ring~$A$.}

\nxsubpoint (Choice of localization theory and topology.)
Whenever we want to concider spectra of generalized rings, and
generalized schemes, obtained from such spectra by gluing, we have 
to fix a localization theory $\cT^?$ and a topology (weak or finite;
strong topology won't do for schemes since quasicoherence is not
local for strong topology). We have just seen that these choices do not
affect spectra of classical rings and classical schemes, so any
choice would be compatible with the classical theory of Grothendieck schemes
in this respect. Usually we do one of the extreme choices:

(a) Minimal choice: unary localization theory, and any topology.
We obtain $\Spec^uA=\Spec^pA$, the prime spectrum of~$A$
in the sense of~\ptref{sp:prime.spec}.

(b) Maximal choice: total localization theory $\cT^t$, consisting of
all open pseudolocalizations, and the finite topology. 
We obtain reasonable spectra $\Spec^t_fA$, which coincide with $\Spec A$
in the classical case, but do not coincide with $\Spec^uA$ for, say, 
$A=\Aff_\bbZ$ (cf.~\ptref{sp:examp.prime.spec}). 
In this way we still have quasicompactness and
nice properties of quasicoherent sheaves, these spectra are still
given by (generalized ringed) topological spaces, and we obtain
more natural results at least for generalized rings like $\Aff_R$.

The only obstruction to adopting (b) in all situations is the absence of
a direct description of all open pseudolocalizations of a generalized
ring~$A$, or of points of $\Spec^t_fA$. For such reasons we may want
to use unary spectra $\Spec^uA$ sometimes.

\nxsubpoint (Local rings.)
Once we fix a localization theory $\cT^?$ and a topology (usually we choose
the finite topology), we say that a generalized ring~$A$ is
{\em $\cT^?$-local\/} if the global sections functor~$\Gamma$ 
on $\Spec^?_fA$ is a point (or a fiber functor) on this topos. Another 
description: any cover of the final object $e$ in $\cS^?_f$ contains
the final object itself. For example, a generalized ring is $\cT^u$-local
iff it has a unique maximal ideal, i.e.\ is local in the sense 
of~\ptref{sp:un.local.rings}. Similarly, a homomorphism $\phi:A\to B$
of $\cT^?$-local rings is {\em local\/} if ${}^a\phi:\Spec^?_fB\to\Spec^?_fA$
maps the only closed point of $\Spec^?_f B$ (given by $\Gamma_B$) into
the only closed point of $\Spec^?_fA$. This is again consistent with
existing terminology.

\nxsubpoint\label{sp:gen.loc.rg.spc} 
(Generalized locally ringed spaces and topoi.)
Now we can define generalized locally ringed spaces (with respect to some
theory $\cT^?$) by requiring the stalks of the structural sheaf of
a generalized ringed space to be $\cT^?$-local, and local morphisms
of such spaces are then defined in the natural way. Furthermore,
we might extend these definitions to topoi, and show the usual
universal property of $\Spec^?_fA$ among all generalized locally ringed
spaces (or topoi) with the global sections of the structural sheaf isomorphic
to~$A$.

\nxpointtoc{Generalized schemes} 
Once we fix a localization theory~$\cT^?$ (say, of
commutative $K$-algebras), and choose a topology (weak or finite;
we'll usually choose finite) we can define a {\em scheme\/} for this theory, or
a {\em $\cT^?$-scheme}, to be a topological space or topos~$X$, together with
a structural sheaf~$\sO_X$ of (generalized) $K$-algebras, such that
$(X,\sO_X)$ is locally isomorphic to generalized ringed spaces or topoi of
form~$\Spec^?A$.  In the topos case this means that we have a cover
$(U_\alpha\to e)$ of the final object by open objects, such that each
$(X_{/U_\alpha},\sO_X|_{U_\alpha})$ is equivalent to some $\Spec^?A$. 
Of course,
once we know that all spectra for the chosen theory are quasicompact, we have
enough points, so we can work with sober topological spaces, and forget about
topoi. 

This is applicable in particular if we choose finite topology on all
our spectra. We are going to do this, so as to obtain a theory of schemes
with almost all basic classical properties fulfilled. Furthermore,
we fix some localization theory~$\cT^?$ and write $\Spec A$ instead of
$\Spec^?_fA$.

\nxsubpoint (Affine schemes.)
A {\em (generalized) affine scheme\/} $X$ is just a generalized (locally) 
ringed space isomorphic to $\Spec A$ for some generalized ring~$A$.
An {\em open subscheme\/} $U$ of a scheme~$X$ is an open subset $U\subset X$
with induced generalized ringed structure: $\sO_U=\sO_X|_U$. 
Notice that any generalized scheme can be covered by affine open subschemes,
essentially by definition of a scheme, and that affine open subsets constitute
a base of topology of any affine scheme $\Spec A$, essentially
by our construction of $\Spec A=\widetilde{\cS^?_{A,f}}$, hence of any scheme. 
This observation also implies that any open subscheme $U$ of a scheme~$X$
is indeed a scheme.

These properties would hold even if we had chosen the weak topology 
on $\cS^?_A$ instead of the finite one; however, the special property of
the finite topology is that {\em any affine scheme is quasicompact.}

\nxsubpoint\label{sp:mor.gen.sch} (Morphisms of generalized schemes.)
A {\em morphism\/} of generalized schemes $f:X\to Y$ is by definition
a {\em local\/} morphism of generalized locally ringed spaces or topoi,
cf.~\ptref{sp:gen.loc.rg.spc}. Another description: for any (generalized) 
affine open
subschemes $U\subset X$ and $V\subset Y$, such that $f(U)\subset V$, 
the restricted morphism $f_{U,V}:U\to V$ 
must be induced by a homomorphism of generalized rings
$\sO(V)\to\sO(U)$. One can either show the equivalence of these two
descriptions, or just adopt the second of them as a definition of
scheme morphism. One can also show that it would suffice to require
that the open subschemes $U$ with the above property (for variable~$V$)
cover~$X$ (cf.~\ptref{sp:loc.morph.gen.sch} below). An immediate consequence is
that {\em morphisms of affine schemes $f:\Spec B\to\Spec A$ are
in one-to-one correspondence with generalized ring homomorphisms
$A\to B$.} Covering an arbitrary scheme $X$ by affine open subschemes
$U_i=\Spec B_i$, and then each $U_i\cap U_j$ by affine open $U_{ij}^{(k)}$,
we deduce that {\em scheme morphisms $f:X\to\Spec A$ are in one-to-one
correspondence with homomorphisms $A\to\Gamma(X,\sO_X)$}, exactly as
in the classical situation.

\nxsubpoint\label{sp:loc.morph.gen.sch} (Localness of scheme morphisms.)
The statement implicitly used above is the following:
{\em if $f:\Spec B\to\Spec A$ is a generalized ringed space morphism,
such that all $f|_{U_i}:U_i=\Spec B_i\to\Spec A$  
are given by some homomorphisms $\phi_i:A\to B_i$, where
$\Spec B=\bigcup_iU_i$ is an affine open cover, then $f$ is itself given
by some homomorphism $\phi:A\to B$.}
Indeed, $U_{ij}=U_i\cap U_j\to U_i$ is obviously given
by a homomorphism, namely, $B_i\to B_{ij}:=B_i\otimes_BB_j$, 
essentially by definition of $\cS^?_{B,f}$ 
(cf.\ also~\ptref{sp:weak.open.loc}),
hence $f|_{U_{ij}}$ is given by $\phi_{ij}:A\stackrel{\phi_i}\to B_i\to 
B_{ij}$, equal to the map
induced by~$f|_{U_{ij}}$ on the global sections of structural sheaves.
In other words, we get maps $\phi_i:A\to B_i=\sO_{\Spec B}(U_i)$, and
the restrictions of $\phi_i$ and $\phi_j$ to $\sO_{\Spec B}(U_{ij})$ 
coincide. Since $\sO_{\Spec B}$ is a sheaf, the sheaf condition for 
cover $\{U_i\}$ implies existence and uniqueness of a generalized ring
homomorphism $\phi:A\to\sO_{\Spec B}(e)=B$, such that $A\stackrel\phi\to
B\to B_i$ coincide with $\phi_i$. This means ${}^a\phi|_{U_i}={}^a\phi_i=
f|_{U_i}$ for any~$i$, hence ${}^a\phi=f$, q.e.d.

\nxsubpoint\label{sp:fibprod.gen.sch} (Fibered products of 
generalized schemes.)
Previous statement immediately implies existence of fibered 
products of affine generalized schemes: $\Spec A\times_{\Spec C}\Spec B$
is given by usual formula $\Spec (A\otimes_CB)$, where this
tensor product of (generalized commutative) $C$-algebras is understood in the 
sense of~\ptref{sp:catgenr}. Once we know this, existence of fibered
products $X\times_SY$ of arbitrary generalized schemes is shown exactly in the same
way as in EGA~I, by gluing together fibered products of appropriate
affine open $U\subset X$ and $V\subset Y$ over $W\subset S$.
Most properties of fibered products of schemes can be generalized to
our case. Notice, however, that the map $|X\times_SY|\to|X|\times_{|S|}|Y|$,
where $|Z|$ denotes the underlying set of a generalized scheme~$Z$,
is not necessarily surjective, as illustrated by $X=Y=\Spec\bbZ$, 
$S=\Spec\Fempty$, $X\times_SX=\Spec\bbZ=X$ (cf.~\ptref{sp:tens.sqr.z}),
in contrast with the classical case.

\nxsubpoint\label{sp:final.gen.sch} (Final generalized scheme.)
Another consequence of~\ptref{sp:mor.gen.sch} is that
the category of generalized schemes admits a final object, namely,
$\Spec\Fempty$, simply because $\Fempty$ is initial in the category
of generalized rings (cf.~\ptref{sp:examp.submon},h). 
In this way any generalized scheme
admits a unique $\Fempty$-scheme structure, so we might say
{\em $\Fempty$-schemes\/} instead of {\em generalized schemes.}

\nxsubpoint\label{sp:gensch.with.0} (Generalized schemes with zero.)
Notice that ${}^a\phi:\Spec B\to\Spec A$ is a monomorphism 
in the category of generalized schemes whenever $\phi:A\to B$ is
an epimorphism of generalized rings. This is applicable in particular
to $\Fempty\to\Fone$ since $\Fone\otimes_{\Fempty}\Fone=\Fone$ 
(recall that $\Fone=\Fempty[0^{[0]}]$, hence
$\Fone\otimes\Fone=\Fempty[0^{[0]},\bar 0^{[0]}]$, and any two commuting
constants coincide). We see that $\Spec\Fone\to\Spec\Fempty$
is a monomorphism, i.e.\ {\em any generalized scheme~$X$ admits at most
one $\Fone$-structure.} If this be the case, we say that
$X$ is an {\em $\Fone$-scheme}, or a {\em (generalized) scheme with zero.}
Clearly, $X$ is a generalized scheme with zero iff $\Gamma(X,\sO_X)$ 
admits a constant, i.e.\ iff the set $\Gamma(X,\sO_X(0))$ is non-empty.

Another consequence of the fact that $\Spec\Fone\to\Spec\Fempty$ is
a monomorphism is that $X\times_{\Fone}Y=X\times_{\Fempty}Y$ for any
two $\Fone$-schemes $X$ and~$Y$, so we can write simply $X\times Y$.

\nxsubpoint\label{sp:add.gensch} (Additive generalized schemes.)
Similarly, $\Fone\to\bbZ$ and $\Fempty\to\bbZ$ are epimorphisms of
generalized rings since $\bbZ\otimes_{\Fempty}\bbZ=\bbZ$
(cf.~\ptref{sp:tens.sqr.z}). This means that $\Spec\bbZ\to\Spec\Fempty$
is a monomorphism of generalized schemes, i.e.\
{\em any generalized scheme~$X$ admits at most one $\bbZ$-structure.}
If this be the case, we say that $X$ is a {\em $\bbZ$-scheme},
or an {\em additive (generalized) scheme.} If $X$ is additive, we
have a homomorphism $\bbZ\to\sO_X(U)$ for any open $U\subset X$,
hence all $\sO_X(U)$ are additive generalized rings with symmetry $[-]$, 
i.e.\ are given by 
classical commutative rings, so $\sO_X$ can be thought of as a sheaf
of classical rings. Taking into account our previous results
\ptref{sp:alg.over.addring} and~\ptref{sp:comp.swf.spec}, 
we see that {\em the category of additive
generalized schemes is canonically equivalent to the category
of Grothendieck schemes.} Notice that $X\times_\bbZ Y=X\times_{\Fempty}Y$
for any additive schemes $X$ and~$Y$, so we can safely write $X\times Y$.

\nxsubpoint\label{sp:qcoh.sheaves} (Quasicoherent sheaves.)
Let $X=(X,\sO_X)$ be a generalized scheme. We say that an $\sO_X$-module
$\sF$ is {\em quasicoherent\/} if for any affine open $U=\Spec A\subset X$
the restriction $\sF|_U$ is isomorphic to $\sO_U$-module $\tilde M$ for
some $A$-module $M$, and then necessarily $M=\Gamma(U,\sF)$. Since
quasicoherence is a local property with respect to weak and finite
topologies (cf.~\ptref{sp:qcoh.loc.for.weak}), 
we see that it would suffice to require this condition 
for $U$ from some affine open cover of~$X$. Once we know this, all usual
properties of quasicoherent sheaves follow. In particular, quasicoherent
sheaves on~$X$ constitute a full subcategory $\catQCoh{\sO_X}$ in 
$\catMod{\sO_X}$, stable under arbitrary inductive and projective limits,
and for any scheme morphism $f:Y\to X$ the pullback functor $f^*:
\catMod{\sO_X}\to\catMod{\sO_Y}$ preserves quasicoherence. When both 
$X$ and~$Y$ are affine, the action of $f^*$ on quasicoherent sheaves
is given by scalar extension with respect to generalized ring
homomorphism~$\phi$ corresponding to~$f$. Furthermore, in this situation
(but not in general) $f_*$ also preserves quasicoherence, and its action
on quasicoherent sheaves is given by scalar restriction along~$\phi$
(we apply the ``affine base change theorem'' \ptref{th:aff.base.change} 
here; cf.~\ptref{sp:qcoh.func.strtop}). 

\nxsubpoint\label{sp:qcoh.ft} (Quasicoherent sheaves of finite type.)
Recall that a sheaf of $\sO_X$-modules $\sF$ on a generalized ringed space
$X=(X,\sO_X)$ is said to be {\em (locally) finitely generated,} or
{\em of finite type}, if one can find an open cover $\{U_\alpha\}$
and strict epimorphisms (i.e.\ surjections) $p_\alpha:
\sO_{U_\alpha}(n_\alpha)\twoheadrightarrow\sF|_{U_\alpha}$, where
$n_\alpha\geq0$ are integers, and $\sO_U(n)$ or $L_{\sO_U}(n)$ denotes 
the free $\sO_U$-module of rank~$n$.

Now suppose that $X$ is a generalized scheme, and that $\sF$ is quasicoherent;
since the above cover can be always chosen to consist of affine $U_\alpha$,
and functor $\Delta:M\mapsto\tilde M$, $\catMod{\sO(U)}\to\catMod{\sO_U}$
is exact for any open affine $U\subset X$, we see that 
{\em an $\sO_X$-module~$\sF$ is quasicoherent of finite type iff
one can find an affine open cover $\{U_\alpha\}$, such that each
$\sF|_{U_\alpha}$ is isomorphic to $\tilde M_\alpha$ for some
finitely generated $\sO(U_\alpha)$-module $M_\alpha$.} 
If we know already that $\sF$ is quasicoherent, this is equivalent to 
requiring all $\Gamma(U_\alpha,\sF)$ to be finitely generated 
as $\sO(U_\alpha)$-modules.

\nxsubpoint\label{sp:qcoh.ft.are.fg} (Global sections of quasicoherent
sheaves of finite type.)
We claim that {\em if $\sF$ is quasicoherent of finite type, then
$\Gamma(U,\sF)$ is a finitely generated $\sO(U)$-module for any affine
open $U\subset X$.} We can assume $X=U=\Spec A$, $\sF=\tilde M$ for
some $A$-module~$M$. Since $X$ is quasicompact, we can find a finite affine
open cover $X=\bigcup_{i=1}^nU_i$, $U_i=\Spec A_i$, such that
$\sF|_{U_i}=\tilde M_i$ for some finitely generated $A_i$-module $M_i$,
clearly equal to $M\otimes_AA_i$.
Now write $M$ as the filtered inductive limit of all its finitely generated
submodules: $M=\injlim_\alpha M_\alpha$, and put 
$M_{i\alpha}:=M_\alpha\otimes_AA_i$. 
Then $M_i=\injlim_\alpha M_{i\alpha}$, and $M_i$ is finitely generated,
whence $M_{i\alpha}=M_i$ for all $\alpha\geq\alpha(i)$. Taking
$\alpha\geq\alpha(i)$ simultaneously for $i=1,\ldots,n$, 
we get $M_i=M_{i\alpha}$
for all~$i$, hence $M=M_\alpha$ is finitely generated.

\nxsubpoint\label{sp:qcoh.fpres} (Finitely presented sheaves.)
Similarly, an $\sO_X$-module $\sF$ is said to be 
{\em (locally) finitely presented,} or {\em of finite presentation\/} 
if it becomes finitely presented 
on some open cover $\{U_\alpha\}$, i.e.\ if $\sF|_{U_\alpha}=
\Coker\bigl(p_\alpha,q_\alpha:
\sO_{U_\alpha}(m_\alpha)\rightrightarrows\sO_{U_\alpha}(n_\alpha)\bigr)$
for suitable $m_\alpha$, $n_\alpha\geq0$ and morphisms $p_\alpha$, $q_\alpha$.
Now if $X$ is a scheme, such a cover can be always assumed to be affine.
Clearly, each $\sF|_{U_\alpha}$ is quasicoherent, being an inductive limit
of quasicoherent sheaves, hence {\em any finitely presented $\sO_X$-module
$\sF$ over a scheme~$X$ is quasicoherent.} Furthermore, we see that
$M_\alpha=\Gamma(U_\alpha,\sF)$ is a finitely presented 
$\sO_X(U_\alpha)$-module, being a cokernel of two maps between free modules
of finite rank.

We claim that {\em if $\sF$ is finitely presented, then
$\Gamma(U,\sF)$ is a finitely presented $\sO(U)$-module for any
affine open $U\subset X$.} We can assume $X=U=\Spec A$ and $\sF=\tilde M$.
By~\ptref{sp:qcoh.ft.are.fg} $M$ is finitely generated, so we can find
a surjection $\pi:A(n)\to M$. Put $R:=A(n)\times_MA(n)$, and 
write $R=\injlim_\alpha R_\alpha$ where $R_\alpha$ are finitely generated
submodules of~$R$, and put $M_\alpha:=A(n)/R_\alpha=\Coker(R_\alpha
\rightrightarrows A(n))$. We get $M=\injlim_\alpha M_\alpha$, all
$M_\alpha$ are finitely presented, and all transition morphisms are
surjective. The proof is finished in the same way as 
in~\ptref{sp:qcoh.ft.are.fg}: we know that $M_i:=M\otimes_AA_i$ are
finitely presented $A_i$-modules, so, in the notations of {\em loc.cit.}\
$M_{i\alpha}\to M_i$ is an isomorphism for $\alpha\geq\alpha(i)$; choosing
$\alpha\geq\alpha(i)$ for all $i$, we see that $M=M_\alpha$ is finitely
presented.

\nxsubpoint\label{sp:qcoh.alg} (Sheaves of $\sO_X$-algebras.)
Notice that we have two different notions of an $\sO_X$-algebra $\sB$,
commutative or not, over a generalized scheme~$X$:

\begin{itemize}
\item[(a)]
We can define such an algebra as a sheaf of generalized rings~$\sB$
over~$X$ together with a central homomorphism $\sO_X\to\sB$.
\item[(b)]
We can define such an $\sB$ as an algebra in $\catMod{\sO_X}$,
considered as a $\otimes$-category with respect to~$\otimes_{\sO_X}$.
\end{itemize}
Similarly to what we had in~\ptref{sp:alg.over.genr} 
and~\ptref{sp:alg.in.Lmod}, any $\sO_X$-algebra in the sense of (b)
defines a {\em unary\/} $\sO_X$-algebra in the sense of (a), and the
functor thus defined is an equivalence of categories, essentially
for the same reasons as in~\ptref{sp:alg.in.Lmod}.

Furthermore, when $\sB$ is a quasicoherent $\sO_X$-algebra, then
for any affine open $U\subset X$ we have $\sB=\tilde B$ for some
$\sO_X(U)$-algebra~$B$ in the sense of~\ptref{sp:alg.over.genr},
where $\tilde B$ is understood as the sheaf of generalized rings over
$U$ defined by the family of sheaves of sets $\{\widetilde{B(n)}\}$,
with the structural maps $\tilde B(\phi)$ and $\mu_n^{(k)}$ induced
by corresponding structural maps for~$A$ (it is important here
that all structural maps involved are in fact $\sO(U)$-homomorphisms,
i.e.\ that $\sO(U)\to B$ is at least central, even if we don't
assume $B$ to be commutative). 

In this way we obtain a natural equivalence
between quasicoherent $\sO_X$-algebras~$\sB$ over an affine scheme
$X=\Spec A$ and $A$-algebras~$B$, given by $\sB\mapsto\Gamma(X,\sB)$,
$B\mapsto\tilde B$. This equivalence preserves unarity of algebras.

\nxsubpoint (Affine morphisms.)
Given a generalized scheme~$S$ and a quasicoherent $\sO_S$-algebra~$\sB$,
we can construct a generalized scheme $X=\Spec\sB$ over~$S$ by gluing
together $\Spec \sB(U)$ for all affine open $U\subset X$, in the
same way as in EGA~II 1.3.1. A scheme morphism $f:X\to S$ is said to be
{\em affine\/} if it is isomorphic to a morphism $\Spec\sB\to S$ of the above
form; in this case obviously $\sB=f_*\sO_X$, where $f_*$ is understood
as a direct image functor for sheaves of generalized rings over~$X$.

All general properties of affine morphisms from EGA~II~1 can be now
transferred to our situation. For example, $f:X\to S$ is affine iff
for any affine open $U\subset S$ its pullback $f^{-1}(U)$ is also affine
iff this condition holds for some open affine cover $\{U_\alpha\}$ of~$S$.

\nxsubpoint (Unary affine morphisms.)
We have a phenomenon specific to our situation: an affine morphism
$f:X\to S$ is said to be {\em unary\/} if corresponding quasicoherent
$\sO_S$-algebra $f_*\sO_X$ is unary. One can easily show that
unary affine morphisms are stable under composition, their diagonals
and graphs are also unary affine etc.

\nxsubpoint (Quasicompact morphisms.)
We say that a generalized scheme morphism $f:X\to S$ is
{\em quasicompact\/} if $f^{-1}(U)$ is quasicompact for any quasicompact
open subset $U\subset S$, i.e.\ if it can be covered by finitely many
affine open subsets of~$X$. It is sufficient to require this for all
open affine $U\subset S$, or just for open affine $U_\alpha$ from some
cover of~$S$.

\nxsubpoint (Quasiseparated morphisms.)
We say that a morphism $f:X\to Y$ is {\em quasiseparated\/} 
if the diagonal $\Delta_{X/Y}:X\to X\times_YX$ is quasicompact.
We say that a generalized scheme $X$ is {\em quasiseparated\/}
if $X\to\Spec\Fempty$ is quasiseparated. These notions have all
their usual properties known from EGA~IV~1. For example, $X$ is
quasiseparated iff it admits an affine open cover $\{U_\alpha\}$,
such that all $U_\alpha\cap U_\beta$ are quasicompact. All
affine schemes and morphisms are quasicompact and quasiseparated,
and $\Delta_{X/Y}$ is quasiseparated for any $X\stackrel f\to Y$.

\nxsubpoint (Open embeddings.)
We say that a generalized scheme morphism $f:Y\to X$ is
an {\em open embedding\/} or an {\em open immersion\/} 
if it is isomorphic to an embedding morphism
$U\to X$ for some open subset $U\subset X$ (clearly equal to
$f(Y)$). Open embeddings are always monomorphisms, hence
quasiseparated, but not necessarily quasicompact. An
{\em open subscheme\/ $U\subset X$} is a subobject of~$X$ 
in the category of schemes defined by an open embedding. Most of 
their elementary properties from EGA~I generalize immediately to
our situation.

\nxsubpoint (Closed embeddings?)
However, we don't seem to have a good notion of closed embeddings
for generalized schemes. This is partly due to the fact that
while the localization theory can be transferred to generalized rings
almost literally, this is not true for ideal/quotient ring theory.

Consider the following example. Let $A$ be a generalized ring with zero,
and put $X:=\Spec A$. One can consider some strict quotient $B$
of~$A$, e.g.\ by choosing two elements $f$, $g\in|A|$ and setting
$B:=A/\langle f=g\rangle$ (cf.~\ptref{sp:eqrel.gen.by.eqs}).
Since $A\to B$ is a (strict) epimorphism, $i:Y:=\Spec B\to X$
is a monomorphism of generalized schemes. In classical situation
$Y$ would have been a closed subscheme of~$X$, ``the locus where $f=g$'',
equal to $V(f-g)$.

However, in our situation $i$ will be a closed embedding in the topological
sense only if we divide by several equations of the form $f_k=0$
(then $Y$ will be homeomorphic to $V(\ga)$, where $\ga$ is the
ideal generated by $(f_k)$, if we work with prime spectra).
Since we don't have subtraction, we cannot replace equations $f_k=g_k$ 
with $f_k-g_k=0$ as in the classical case, so
in general $i:Y\to X$ can have a non-closed image,
as illustrated by the following example.

\nxsubpoint (Diagonal of $\bbA^1_{\Fone}$.)
Put $X:=\bbA^1_{\Fone}=\Spec\Fone[T]$ and consider
$\Delta_X:X\to X\times X$. Clearly, $\Delta_X$ is defined by
strict epimorphism $\Fone[T]\otimes\Fone[T]=\Fone[T,T']\to\Fone[T]$,
mapping both $T$ and~$T'$ into~$T$. Now let's compute the
corresponding prime spectra: $\Spec\Fone[T]=\{0,(T)\}$,
$\Spec\Fone[T,T']=\{0,(T),(T'),(T,T')\}$, and the image of $\Delta_X$
equals $\{0,(T,T')\}$, which is not closed in $\Spec\Fone[T,T']$,
i.e.\ the image of $\Delta_X$ is not closed, at least if we use~$\cT^u$.

\nxsubpoint (Separated morphisms?)
At this point we see that we don't have any reasonable notion of
separated morphisms $f:X\to Y$: if we require $\Delta_{X/Y}$ to
be a closed immersion (whatever that means), 
then even the affine line $\bbA^1_{\Fone}$ wouldn't be separated.

All we can do is to consider quasiseparated morphisms 
instead, whenever this is possible. It is not always enough:
for example, a section $\sigma:S\to X$ of an affine morphism $f:X\to S$ 
won't be a closed immersion, just because $\Delta_{X/S}$ is not.

\nxsubpoint\label{sp:imm.subsch} (Immersions and subschemes.)
However, some related notions can be still saved.
Namely, let us say that $f:Y\to X$ is an {\em immersion\/}
if it is a monomorphism, locally (in $Y$ and~$X$) given by
${}^a\phi:\Spec B\to\Spec A$, where $\phi:A\to B$ is a homomorphism
of a generalized ring~$A$ into a finitely generated $A$-algebra~$B$,
such that the scalar restriction functor $\phi_*:\catMod B\to\catMod A$
is fully faithful (hence $\phi$ is a NC-epimorphism, and in particular
an epimorphism, cf.~\ptref{sp:cat.psloc}), and retains this property after any
base change. Examples of such homomorphisms are given by
open pseudolocalizations and strict epimorphisms (i.e.\ surjective
homomorphisms).

A {\em subscheme $Y\subset X$}
is any subobject of~$X$ in the category of schemes, such that
$Y\to X$ is an immersion. Clearly, immersions are stable under composition
and base change, open embeddings are immersions, and any diagonal
morphism $\Delta_{X/S}$ is an immersion, hence the graph
$\Gamma_f:X\to X\times_SY$ of any morphism $f:X\to Y$ of
(generalized) $S$-schemes is an immersion as well.

\nxsubpoint (Monomorphisms are not injective!)
Notice that, in contrast to the classical case, a monomorphism
of generalized schemes can be not injective on points, as illustrated
by $\Spec\bbZ\to\Spec\Fempty$. It is an interesting question whether
an immersion can be non-injective.

\nxsubpoint\label{sp:closed.imm.sch} (``Closed'' immersions.)
We say that a morphism of generalized schemes $j:Y\to X$ is a 
``closed'' immersion if it is affine, and if the induced
homomorphism of sheaves of generalized rings $\sO_X\to f_*\sO_Y$
is surjective. Any ``closed'' immersion is an immersion in the
sense of~\ptref{sp:imm.subsch}; however, it is not necessarily closed in the
topological sense, whence the quotes. A ``closed'' immersion
is locally of the form $j:\Spec B\to\Spec A$, where $\phi:A\twoheadrightarrow
B$ is a surjective homomorphism of generalized rings. If we work in
$\cT^u$ (i.e.\ with prime spectra), then such a $j$ is injective on points
(since $\phi^{-1}(\gp)=\phi^{-1}(\gp')$ implies $\gp=\gp'$),
and the topology on $j(\Spec B)$ coincides with that induced from $\Spec A$,
since $j^{-1}(D(f))=D(\phi(f))$, open sets $D(\bar f)$ constitute
a base of topology of $\Spec B$, and any $\bar f\in|B|$ can be lifted
to some~$f\in|A|$, hence $D(\bar f)=j^{-1}(D(f))$.

These ``closed'' immersions are stable under composition and base change,
and the diagonal of an affine morphism is a ``closed'' immersion.
Therefore, it might make sense to define a separated (or rather ``separated'')
morphism $f:X\to Y$ by requiring its diagonal $\Delta_{X/Y}$ to be a ``closed''
immersion. Notice that ``closed'' immersions are always affine,
but not necessarily unary.

Of course, we define a ``closed'' subscheme $Y\subset X$ as a subobject
in the category of generalized schemes, such that $Y\to X$ is
a ``closed'' immersion.

\nxsubpoint (Disjoint union of affine schemes.)
Let $X_1=\Spec A_1$ and $X_2=\Spec A_2$ be two affine generalized schemes
with zero. We claim that $\Spec(A_1\times A_2)\cong X_1\sqcup X_2$;
in particular, the disjoint union $X_1\sqcup X_2$ is affine.

Indeed, according to~\ptref{th:mod.under.prod}, $\catMod{(A_1\times A_2)}$
is equivalent to $\catMod{A_1}\times\catMod{A_2}$; the inverse of
this equivalence transforms $(M_1,M_2)$ into $M_1\times M_2$,
with the action of $t=(t_1,t_2)\in(A_1\times A_2)(n)=A_1(n)\times A_2(n)$
defined componentwise. Now it is immediate that $A_1$ equals $A[\bu_1^{-1}]$,
where $A:=A_1\times A_2$ and $\bu_1=(\bu,0)\in|A|$, and similarly for
$A_2$, whence two open embeddings $\lambda_i:X_i=\Spec A_i\to X:=\Spec A$.
Next, $X_1\cap X_2$ equals $\Spec A[(\bu_1\bu_2)^{-1}]=\Spec A[0^{-1}]=
\emptyset$;
applying \ptref{th:mod.under.prod} once more, we see that
any quasicoherent sheaf $\tilde M$ on~$X$ satisfies the sheaf property
for family of open (pseudo)localizations $\{X_i\to X\}_{i=1,2}$, hence
$X=X_1\cup X_2$ by the definition of (strong) topology,
cf.~\ptref{sp:constr.spec} and~\ptref{sp:fin.desc.fam.eff}, q.e.d.

\nxsubpoint (Unarity is not local in the source.)
Notice that the unarity of affine schemes over a fixed base generalized
ring~$C$ is not a local property, as illustrated by
$C=\Zinfty$, $X=\Spec(\Zinfty\times\Zinfty)=\Spec\Zinfty\sqcup\Spec\Zinfty$:
we see that $X$ is a union of unary affine $C$-schemes, but it is
not unary itself.

\nxsubpoint (Morphisms of finite type and of finite presentation.)
We say that a generalized scheme morphism $f:X\to Y$
is {\em locally of finite type\/} if one can cover $Y$ by affine open
subschemes $V_i$, and each $f^{-1}(V_i)$ by affine open $U_{ij}$,
such that each $\sO(U_{ij})$ is a finitely generated $\sO(V_i)$-algebra
in the sense of~\ptref{sp:fingen.preun.alg}. Similarly,
$f$ is {\em locally of finite presentation\/} if one can find similar
covers with the property that each $\sO(U_{ij})$ be a finitely
presented $\sO(V_i)$-algebra in the sense of~\ptref{sp:finpres.un.alg}.

Furthermore, we say that $f$ is {\em of finite type\/}
if it is quasicompact and locally of finite type, and that
$f$ is {\em of finite presentation\/}
if it is quasicompact, quasiseparated and locally of finite presentation.

These notions have all formal properties of EGA~I~6 and EGA~IV~1,
modulo the following statement:

\nxsubpoint
{\em If $\Spec B\to\Spec A$ is locally of finite type (resp.\ presentation),
then $B$ is a finitely generated (resp.\ finitely presented) $A$-algebra.}
Unfortunately, we don't have a proof of this statement for $\cT^?\neq\cT^u$.
However, for the unary localization theory $\cT^u$ the proofs from EGA
work with minimal modifications.

Let us show for example the statement about finite generation. Suppose
that $s_1$, \dots, $s_n\in|B|$ are such that all $B_{s_i}$ are
finitely generated $A$-algebras, and choose some finite lists of
generators $f_{ij}\in B_{s_i}(r_{ij})$; we can write
$f_{ij}=g_{ij}/s_i^{k_{ij}}$ for some $k_{ij}\geq0$ and
$g_{ij}\in B(r_{ij})$, choose a common value $k$ for all $k_{ij}$
and replace $s_i$ by $s_i^k$, thus reducing to the case $k_{ij}=1$.
Next, since $\bigcup_{i=1}^nD(s_i)=\Spec B$, these elements
$s_i$ generate the unit ideal over~$B$ 
(cf.~\ptref{sp:covers.un.spec} and~\ptref{sp:opcov.specp}),
so we can find an operation $h\in B(n)$, such that $h(s_1,\ldots,s_n)=1$
(cf.~\ptref{sp:subm.gen.subset}). Now we put $B':=
A[s_i^{[1]},g_{ij}^{[r_{ij}]},h^{[n]}]\subset B$
and show that $B'=B$ in the same way as in EGA~I~6.5.

If all $B_{s_i}$ are finitely presented over~$A$, we can find a surjection
$\rho:B':=A[T_1^{[r_1]},\ldots,T_m^{[r_m]}]\to B$ from an algebra of
generalized polynomials, consider its kernel $R:=B'\times_BB'$,
choose finite lists of equations $(f_{ij}/\tilde s_i^N,g_{ij}/\tilde s_i^N)$
generating compatible equivalence relation $R_{\tilde s_i}$ 
on~$B_{\tilde s_i}$, the kernel of $B'_{\tilde s_i}\to B_{s_i}$, where
$\tilde s_i\in B'$ are some lifts of~$s_i$, and show that $R$ is generated 
by equations $\{f_{ij}=g_{ij}\}$ as 
a compatible equivalence relation on generalized ring~$B'$,
hence $B=B'/R$ is finitely presented over~$A$.

\nxsubpoint\label{sp:vb.proj} (Vector bundles and projective modules.)
A {\em vector bundle\/} $\sE$ over a generalized scheme~$S$ is 
defined as a {\em locally free $\sO_S$-module of finite type}, i.e.\
we require $\sE|_{U_i}\cong\sO_{U_i}(n_i)$
as an $\sO_{U_i}$-module, for some open cover $U_i$ of~$S$ and
some integers $n_i\geq0$. If all $n_i$ can be chosen to have the same
value~$n$, we say that {\em $\sE$ is a vector bundle of rank~$n$.}
A vector bundle is obviously a (locally) finitely presented
quasicoherent $\sO_S$-module.

\begin{Propz} 
If $\sE$ is a vector bundle over~$S$, then
$\Gamma(U,\sE)$ is a finitely generated projective module
over generalized ring~$\Gamma(U,\sO_S)$, for any affine open $U\subset S$.
\end{Propz}
\begin{Proof}
It suffices to consider the case $U=S=\Spec A$, $\sE=\tilde P$ for some
$A$-module~$P$, finitely presented by~\ptref{sp:qcoh.fpres}.
We have to check that~$P$ is projective, i.e.\ that the functor
$\Hom_A(P,-):\catMod A\to\catMod A$ preserves strict epimorphisms
(surjections of $A$-modules). Notice for this that
$\Hom_A(P,M)\widetilde{\phantom{M}}\cong\iHom_{\sO_S}(\tilde P,\tilde M)$
for any finitely presented $A$-module~$P$. The proof proceeds as in the
classical case of EGA~I: fix a finite presentation of~$P$, use
exactness properties of functors involved to reduce to the case of
a free~$P=A(n)$ of finite rank~$n$, and then we are reduced to proving
$\widetilde{M^n}\cong\tilde M^n$, true by exactness of $\Delta:M\to\tilde M$.
Applying this statement to our~$P$ and noticing that $f:M\to N$ is 
a strict epimorphism iff $\tilde f:\tilde M\to\tilde N$ is one,
we see that projectivity of~$P$ is equivalent to
$\iHom_{\sO_S}(\sE,\sF)\to\iHom_{\sO_S}(\sE,\sG)$ being a strict epimorphism
for any strict epimorphism $\phi:\sF\to\sG$ of quasicoherent $\sO_S$-modules.
This is a local property, obviously fulfilled for any free~$\sE$ of finite 
rank, hence for any locally free~$\sE$ of finite type, q.e.d.
\end{Proof}

\nxsubpoint (Converse is not true.)
Contrary to the classical case, $\tilde P$ is not necessarily a vector
bundle over~$\Spec A$ whenever $P$ is a finitely generated 
projective~$A$-module, at least if we use unary localization theory
$\cT^u$ to construct our spectra, as illustrated by the
non-free projective $\Finfty$-module~$P$ of~\ptref{sp:proj.finfty.not.free},
since $\Spec^u\Finfty$ is a one-point set.

\nxsubpoint\label{sp:lb.Pic} (Line bundles and Picard group.)
A {\em line bundle\/}~$\sL$ over~$S$ is defined as a vector bundle of rank~1,
i.e.\ an $\sO_S$-module, locally isomorphic to $|\sO_S|=\sO_S(1)$.
Tensor product of two line bundles is obviously again a line bundle,
and $\sL^{\otimes-1}:=\check\sL:=\iHom(\sL,|\sO_S|)$ is also a line bundle,
such that the canonical map $\sL\otimes_{\sO_S}\check\sL\to|\sO_S|$ is
an isomorphism.

This enables us to define the {\em Picard group\/} $\Pic(S)$ as the
set of isomorphism classes of line bundles over~$S$, with the group operation
given by the tensor product. If we denote by $|\sO_S|^\times$ the
sheaf of invertible elements in~$|\sO_S|$, then 
$\Pic(S)=H^1(S,|\sO_S|^\times)$, where the RHS is understood as the usual
cohomology of a sheaf of abelian groups.

\nxsubpoint (Flat morphisms.)
We say that a generalized scheme morphism $f:X\to S$ is
{\em flat\/} if for any affine open $U\subset X$ and $V\subset S$,
such that $U\subset f^{-1}(V)$, $\sO(U)$ is a flat $\sO(V)$-algebra.
This is equivalent to the existence of an affine open cover
$\{V_i\}$ of~$S$ and affine open covers $\{U_{ij}\}$ of $f^{-1}(V_i)$,
such that each $\sO(U_{ij})$ is a flat $\sO(V_i)$-algebra.

In order to check equivalence of these two descriptions we have to show
that an $A$-algebra $B$ is flat 
whenever there is an affine cover $\{U_i\to X\}$, such that
each $B_i$ is flat over~$A$, where we put $X:=\Spec B$, $S:=\Spec A$,
$U_i:=\Spec B_i$ and denote $\Spec B\to\Spec A$ by~$f$. 
This is immediate: indeed, functors
$f^*|_{U_i}:\catQCoh{\sO_S}\to\catQCoh{\sO_{U_i}}$ are exact,
hence $f^*:\catQCoh{\sO_S}\to\catQCoh{\sO_X}$ is exact, i.e.\
the scalar extension functor $\catMod A\to\catMod B$ is exact,
so $B$ is indeed a flat $A$-algebra.

\nxsubpoint (Faithfully flat morphisms.)
We say that $f:X\to S$ is {\em faithfully flat\/}
if it is flat, and exact functor $f^*:\catMod{\sO_S}\to\catMod{\sO_X}$
is fully faithful. This is equivalent to requiring all
$f_U^*:\catQCoh{\sO_U}\to\catQCoh{\sO_{f^{-1}(U)}}$ to be fully
faithful, where $U\subset S$ is any affine open subscheme, and
$f_U:f^{-1}(U)\to U$ is the restriction of~$f$.

\nxsubpoint (Flat topology (fpqc).)
Consider the category $\cS$ of affine generalized schemes, or, equivalently,
the opposite category of the category $\catGenR$ of generalized rings.
Introduce on~$\cS$ the topology generated by finite families of flat morphisms 
$\{\Spec A_i\to\Spec A\}$, which are (universal efficient) families
of descent for $\stQCOH$. The resulting topology is the {\em flat topology
(fpqc)}. Any scheme~$X$ represents a functor $\tilde X:\Spec A\mapsto
X(A)=\Hom(\Spec A,X)$, which is easily seen to be a sheaf for the
flat topology. Functor $h:X\mapsto\tilde X$ from generalized schemes
into the flat topos $\tilde\cS$ is fully faithful, so we can define
the flat topology on the category of all generalized schemes
simply by pulling back the canonical topology of $\tilde\cS$ with respect 
to~$h$. Of course, the topos defined by this larger site will be 
still~$\tilde\cS$.

One can describe this flat topology more explicitly. We don't want to
do it here; let us remark that a large portion of 
faithfully flat descent theory (cf.\ e.g.\ SGA~1) can be generalized 
to our case with the aid of flat topology just defined.

\nxsubpoint (\'Etale coverings, \'etale and smooth morphisms.)
Once we have a notion of flat topology, we can define for example
{\em \'etale coverings\/} $X\to S$ as morphisms, locally for the flat topology 
(over~$S$) isomorphic to ``constant \'etale coverings'' 
$S\times\stn\to S$, where $S\times\stn=S\sqcup\ldots\sqcup S$ denotes
the disjoint union of $n$ copies of~$S$. After that we can define
{\em \'etale morphisms\/} $f:X\to Y$ as morphisms that locally for the
flat topology on~$Y$, and after that --- locally for the Zariski topology
on~$X$ are of the form $X\stackrel j\to X'\stackrel\pi\to Y$,
where $j$ is an open embedding and $\pi$ an \'etale covering. Finally,
we can define {\em smooth morphisms\/} in a similar manner,
where this time $j$ is to be \'etale and $\pi:X'\to Y$ is to be
of the form $Y\times\bbA^n\to Y$.

However, we are not going to study or use these definitions in this work;
they've been mentioned just to demonstrate a possible way of extending
some classical notions to generalized schemes.

\nxsubpoint (Proper and projective morphisms.)
Unfortunately, the obvious definition of a proper morphism of
generalized schemes involving universal closedness 
doesn't seem to be useful, because
diagonals of affine morphisms are not proper in this sense. 
However, we have a reasonable notion of a projective morphism,
so we are going to study projective generalized schemes and morphisms
instead.

\nxpointtoc{Projective generalized schemes and morphisms}
Now we want to define projective generalized schemes and morphisms and
related notions. There are three major differences from the classical
case of EGA~II. Firstly, we don't have reasonable noetherian 
conditions, so we have to use finitely presented morphisms and
quasicoherent sheaves instead. Secondly, we have to deal with
{\em non-unary\/} generalized graded rings and their projective spectra,
together with various arising phenomena that have no counterparts
in the classical case. Thirdly, we have to consider ``closed''
immersions in the sense of~\ptref{sp:closed.imm.sch}, which are not
really closed as maps of topological spaces. This also  
complicates things considerably.

\nxsubpoint\label{sp:graded.alg.in.ACU} 
(Graded algebras and modules in a tensor category.)
Let $\Delta$ be a commutative monoid (e.g.\ $\bbN_0$ or $\bbZ$),
$\sA$ be an ACU $\otimes$-category, such that direct sums
(i.e.\ coproducts) indexed by subsets of $\Delta$ exist in $\sA$ and
commute with~$\otimes$. Suppose that $M\to M\oplus N$ is a monomorphism
for any $M$, $N\in\Ob\sA$ (condition automatically fulfilled
if $\sA$ has a zero object). Then we can define a
{\em $\Delta$-graded algebra~$S$ in~$\sA$} as an algebra
$S=(S,\mu,\epsilon)$ in~$\sA$ together with a direct sum decomposition
$S=\bigoplus_{\alpha\in\Delta}S^\alpha$, such that
$\epsilon:\Unit_\sA\to S$ factorizes through $S^0\subset S$,
and $S^\alpha\otimes S^\beta\to S\otimes S\stackrel\mu\to S$ factorizes through
$S^{\alpha+\beta}\subset S$, for any $\alpha$, $\beta\in\Delta$.
Next, we can define a {\em graded (left) $S$-module~$M$ in~$\sA$}
as an $S$-module $M=(M,\gamma)$ with a direct sum decomposition
$M=\bigoplus_{\alpha\in\Delta}M^\alpha$, such that
$S^\alpha\otimes M^\beta\to S\otimes M\stackrel\gamma\to M$ factorizes
through $M^{\alpha+\beta}\subset M$ for all $\alpha$ and~$\beta\in\Delta$.

\nxsubpoint\label{sp:alt.descr.gr.alg} (Alternative descriptions.)
The above situation admits an alternative description.
Namely, we can define a $\Delta$-graded algebra~$S$ in~$\sA$
as a family $\{S^\alpha\}_{\alpha\in\Delta}$ of objects of~$\sA$
together with identity map $\epsilon:\Unit_\sA\to S^0$ and
multiplication maps $\mu_{\alpha,\beta}:S^\alpha\otimes S^\beta\to 
S^{\alpha+\beta}$, satisfying natural associativity and unit relations.
A graded $S$-module $M$ can be described similarly, as a family
$\{M^\beta\}_{\beta\in\Delta}$ together with action maps
$\gamma_{\alpha,\beta}:S^\alpha\otimes M^\beta\to M^{\alpha+\beta}$,
satisfying obvious relations.

This approach has its advantages. It yields definitions equivalent to
those of~\ptref{sp:graded.alg.in.ACU} when $\Delta$-indexed
direct sums exist in~$\sA$ and commute with~$\otimes$, and
all $M\to M\oplus N$ are monic. However, these new definitions don't
need any of these assumptions and are valid in any ACU $\otimes$-category.

Yet another description: if $\Delta$-indexed direct sums exist in~$\sA$
and commute with~$\otimes$, we can consider the category
$\sA^\Delta$ of $\Delta$-indexed families of objects of~$\sA$,
and introduce a $\otimes$-structure on it by putting
$\{X^\alpha\}\otimes\{Y^\beta\}:=\{\bigoplus_{\alpha+\beta=\gamma}
X^\alpha\otimes Y^\beta\}_{\gamma\in\Delta}$. Then
$\Delta$-graded algebras and modules are nothing else than
algebras and modules in~$\sA^\Delta$.

\nxsubpoint (Application: $\Delta$-graded unary algebras and modules
over a generalized ring.)  
Let us fix a generalized base ring~$C$. Since $\catMod C$ is an ACU
$\otimes$-category with arbitrary direct sums commuting with the
tensor product, we can apply any of the above constructions to define
$\Delta$-graded $C$-algebras and modules over them.  Notice that the
condition ``$M\to M\oplus N$ injective'' is fulfilled for all
$C$-modules $M$ and~$N$ whenever~$C$ admits a zero,
so~\ptref{sp:graded.alg.in.ACU} is applicable at least over
such~$C$'s.  On the other hand, the alternative approach
of~\ptref{sp:alt.descr.gr.alg} is applicable without any restrictions
on~$C$, and has another convenient property. Namely, let $C_0\subset
C$ be a generalized subring; we would like to be able to treat any
graded $C$-algebra $S=\bigoplus S^\alpha$ as a graded $C_0$-algebra
via scalar restriction. However, if $C$ is not unary over $C_0$,
direct sums over $C$ and~$C_0$ differ, so the scalar restriction
of~$S$ wouldn't be equal to the direct sum of scalar restrictions of
its components~$S^\alpha$.

Nevertheless, if we adopt the ``family approach'' 
of~\ptref{sp:alt.descr.gr.alg}, this problem doesn't arise:
family $\{S^\alpha\}_{\alpha\in\Delta}$ together with $1\in S^0$
and $C$-bilinear maps $\mu_{\alpha,\beta}:S^\alpha\times S^\beta\to
S^{\alpha+\beta}$ retain all their properties after scalar restriction,
so we can work with $\Delta$-graded algebras and modules without
worrying too much about the choice of~$C$.

\nxsubpoint (Non-unary $\Delta$-graded $C$-algebras.)
The above description is somewhat incomplete. Namely, 
according to~\ptref{sp:alg.in.Lmod}, 
algebras in $\catMod C$ correspond to {\em unary\/} $C$-algebras,
so we've just described only unary graded $C$-algebras. If we
manage to describe non-unary graded $C$-algebras, then we can indeed
replace $C$ by $\Fempty$ and forget about it altogether, thus studying
$\Delta$-graded generalized rings and graded modules over them.

\nxsubpoint\label{sp:graded.alg.mon}
(Graded algebraic monads and generalized rings: non-unary case.)
A $\Delta$-graded algebraic monad or generalized ring~$S$ can
be defined as follows. We must have a collection of sets 
$\{S^\alpha(n)\}_{\alpha\in\Delta,n\geq0}$, together with
transition maps $S^\alpha(\phi):S^\alpha(n)\to S^\alpha(m)$ for
each map $\phi:\stn\to\stm$, such that 
$S^\alpha:\stn\mapsto S^\alpha(n)$ becomes a functor $\catN\to\catSets$,
which can be uniquely extended to an algebraic endofunctor $S^\alpha$ 
on~$\catSets$, if we wish to. Informally, $S^\alpha(n)$ is just the
degree~$\alpha$ component $S(n)^\alpha$ of free $S$-module~$S(n)$ of rank~$n$.

Next, we must have an identity $\bu\in S^0(1)$, or, equivalently,
a functorial morphism $\epsilon:\Id_{\catSets}\to S^0$,
and composition maps $\mu_{n,\beta}^{(k,\alpha)}:
S^\alpha(k)\times S^\beta(n)^k\to S^{\alpha+\beta}(n)$,
satisfying certain compatibility relations, which are just
``graded versions'' of those of~\ptref{sp:expl.algmon.cond}.
If we want to study $\Delta$-graded generalized rings, not just
algebraic monads, we must impose ``graded commutativity relations''
on $\mu$'s, involving commutative diagrams 
similar to~\eqref{eq:comm.comm.diagr}.

After this we can define a graded $S$-module~$M$ as a collection of
sets $\{M^\beta\}_{\beta\in\Delta}$ together with action maps
$\gamma_\beta^{(k,\alpha)}:S^\alpha(k)\times (M^\beta)^k\to M^{\alpha+\beta}$.
In this way we obtain a reasonable category of graded $S$-modules,
together with a monadic forgetful functor $\Gamma:\catGradMod S_\Delta\to
\catSets_{/\Delta}$ into the category of $\Delta$-graded sets.
This yields an algebraic monad over $\catSets_{/\Delta}$
in the sense of~\ptref{sp:catalgmon.mon}, 
i.e.\ $\Delta$-graded algebraic monads
just defined are special cases of those considered in {\em loc.cit.}

\nxsubpoint (Underlying set of a graded $R$-module.)
Notice that the underlying set $|M|=M$ of a graded $R$-module~$M$
is just the $\Delta$-graded set $|M|:=\bigsqcup_{\alpha\in\Delta}
M^\alpha$. In other words, any element $x\in M$ is automatically 
homogeneous. 

\nxsubpoint (Trivially graded generalized rings.)
Any generalized ring~$C$ can be treated as a trivially $\Delta$-graded
generalized ring $\bar C$ (usually denoted also by~$C$),
given by $\bar C^0(n):=C(n)$, $\bar C^\alpha(n):=\emptyset$ for $\alpha\neq0$,
i.e.\ $\bar C^0:=C$, $\bar C^\alpha:=\emptyset$ for $\alpha\neq0$
in terms of corresponding algebraic endofunctors.
Conversely, $R^0$ has a natural generalized ring structure for any
$\Delta$-graded generalized ring~$R$, and one easily checks that these
two functors $C\mapsto\bar C$ and $R\mapsto R^0$ are adjoint. 
In particular, a $C$-algebra structure
on a $\Delta$-graded generalized ring~$R$ is the same thing as 
a homomorphism $C\to R^0$, i.e.\ a $C$-algebra structure on~$R^0$.

\nxsubpoint\label{sp:deg.shift.mod} (Degree shift of graded modules.)
Let $M$ be a graded module over $\Delta$-graded generalized ring~$S$,
and $\lambda\in\Delta$ be an arbitrary element.
We denote by $M[\lambda]$, or sometimes by $M(\lambda)$, when no
confusion can arise, the corresponding {\em degree shift\/} of~$M$,
given by $M[\lambda]^\alpha:=M^{\lambda+\alpha}$,
$\gamma_{M[\lambda],\beta}^{(k,\alpha)}:=
\gamma_{M,\lambda+\beta}^{(k,\alpha)}$. These operations have all
their usual properties, e.g.\ $(M[\lambda])[\mu]=M[\lambda+\mu]$ and
$M[0]=M$.

\nxsubpoint (Free graded modules.)
We denote by $L_R(n)$ or simply by $R(n)$ the graded $R$-module
defined by collection $\{R^\alpha(n)\}_{\alpha\in\Delta}$,
with the $R$-action given by appropriate maps~$\mu^{(k,\alpha)}_{n,\beta}$. 
We can extend this definition to all sets $I$ by putting 
$L_R(I):=\injlim_{\stn\stackrel\phi\to I}L_R(n)$.
Another description: $L_R:\catSets\to\catGradMod R_\Delta$ is
left adjoint functor to $M\mapsto M^0$. Furthermore, we can always construct
a left adjoint to $\Gamma_R:\catGradMod R_\Delta\to\catSets_{/\Delta}$;
when $\Delta$ is a group, it can be defined
by transforming a family of sets $\{I_\alpha\}_{\alpha\in\Delta}$ into
$\bigoplus_\alpha L_R(I_\alpha)[-\alpha]$ (of course, arbitrary
inductive and projective limits of graded $R$-modules exist).
It makes sense to call such graded $R$-modules also {\em (graded) free}.

\nxsubpoint (Properties of graded rings and modules.)
One can transfer to graded case almost all our general considerations
done before in the non-graded case. For example, we can define
graded algebras given by generators and relations (with the
degree of each generator explicitly mentioned) by their universal property,
and show that they exist, and use this result to prove existence of
coproducts in the category of $\Delta$-graded generalized rings,
i.e.\ the ``graded tensor product'' $R'\otimes_RR''$ of two
graded algebras $R'$ and $R''$ over another graded generalized ring~$R$.

Furthermore, the categories of graded modules over graded
generalized rings also retain most of their nice properties,
e.g.\ existence of inductive and projective limits, 
scalar restriction and scalar extension functors and so on.
We can also define a {\em graded $R$-bilinear map\/}
$\phi:M\times N\to P$, where $M$, $N$, $P$ are graded $R$-modules, 
by requiring all $s_\phi(x):N\to P[\deg x]$, $y\mapsto\phi(x,y)$
to be graded $R$-homomorphisms for all $x\in M$,
and similarly for $d_\phi(y):M\to P[\deg y]$, for any $y\in N$
(here $M[\deg x]$ denotes the degree shift of~$M$ by~$\deg x$,
cf.~\ptref{sp:deg.shift.mod}). After this we can define and construct graded
tensor products $M\otimes_RN$ and inner Homs $\iHom_R(M,N)$;
one can check that $\iHom_R(M,N)^\alpha=\Hom_R(M,N[\alpha])$,
i.e.\ {\em graded component of $\iHom_R(M,N)$ of degree~$\alpha$
consists of homorphisms $M\to N$ of degree~$\alpha$.}

\nxsubpoint (Monoid of unary operations.)
Let $R$ be a $\Delta$-graded generalized ring. Maps
$\mu^{(1,\alpha)}_{1,\beta}:R^\alpha(1)\times R^\beta(1)\to 
R^{\alpha+\beta}(1)$ define on $|R|=|R(1)|:=\bigcup_\alpha R^\alpha(1)$
a (commutative) monoid structure with identity $\bu\in R^0(1)$,
and the degree map $\deg:|R|\to\Delta$, equal by definition
to~$\alpha$ on $R^\alpha(1)$, is obviously a monoid homomorphism.

\nxsubpoint (Graded multiplicative systems. Localization.)
A {\em (graded) multiplicative system} $S\subset|R|$ is simply a submonoid
of~$|R|$ with respect to monoid structure just discussed. Given
any subset $S\subset|R|$, such that $\deg S\subset\Delta$ consists
of invertible elements of~$\Delta$, we can consider the multiplicative system
(i.e.\ monoid) $\langle S\rangle\subset|R|$, generated by~$S$.
Next, we can consider graded $R$-algebras $R\stackrel\rho\to R'$
with the property that all elements of $\rho(S)\subset|R'|$ become
invertible in~$|R'|$, and define the {\em localization $R[S^{-1}]$}
as the universal (initial) $R$-algebra with this property.
Localizations always exist and can be defined by putting
$R[S^{-1}]:=R[s^{-1}|ss^{-1}=\bu]_{s\in S}$, where $\{s^{-1}\}_{s\in S}$
are new unary generators with $\deg s^{-1}:=-\deg s$.

Next, we can consider the full subcategory of $\catGradMod R_\Delta$,
consisting of those graded $R$-modules~$M$, on which all $s\in S$
act bijectively, i.e.\ the map $[s]_M:M\to M$ (or equivalently,
each of the maps $[s]^\alpha_M:M^\alpha\to M^{\alpha+\deg s}$)
is bijective. One checks in the usual manner that
this category is equivalent to $\catGradMod{R[S^{-1}]}$ by means of 
the scalar restriction functor, hence the scalar extension
$M[S^{-1}]$ of a graded $R$-module $M$ to $R[S^{-1}]$ can be
characterized by its universal property among all $R$-homomorphisms
from $M$ into $R$-modules $N$, on which $S$ acts bijectively.

\nxsubpoint (Classical description of localizations.)
One can always replace $S$ by $\langle S\rangle$ in the above considerations,
thus assuming $S$ to be a (graded) multiplicative system. In this case
we also write $S^{-1}R$ and $S^{-1}M$ instead of $R[S^{-1}]$ and $M[S^{-1}]$.
Then we can describe these objects in the classical way
(cf.~\ptref{sp:class.descr.loc}). For example,
$S^{-1}M=M\times S/\sim$, where $(x,s)\sim(y,t)$ iff there is some
$u\in S$, such that $utx=usy$. If we denote the class of $(x,s)$ by $x/s$,
then the grading on $S^{-1}M$ is given by $\deg(x/s):=\deg x-\deg s$
(recall that $\deg S\subset\Delta$ was supposed to consist of
invertible elements of~$\Delta$). Since $(S^{-1}R)(n)=S^{-1}(R(n))$,
$S^{-1}R(n)$ admits a similar description as $R(n)\times S/\sim$,
with equivalence relation $\sim$ and grading given by the same formulas.

\nxsubpoint (Flatness of localizations.)
Furthermore, graded localizations admit a filtered inductive limit
description similar to that of~\ptref{sp:constr.localiz.mod}.
Namely, let $S\subset|R|$ be a multiplicative system, such that
$\deg S\subset\Delta$ consists of invertible elements, and
let $M$ be a graded $R$-module. Consider small category $\cS$
of~\ptref{sp:constr.localiz.mod}, with $\Ob\cS=S$, and morphisms
given by $\Hom_\cS([s],[s']):=\{t\in S:ts=s'\}$. Define a functor
$\bar M:\cS\to\catGradMod R$ by putting $\bar M_{[s]}:=M[\deg s]$,
with the transition morphisms $t:M_{[s]}\to M_{[s']}$ given by
the action of $t\in S\subset|R|$ on~$M$. Then $S^{-1}M=\injlim_\cS M_{[s]}$,
and since filtered inductive limits are left exact,
$M\mapsto S^{-1}M$ is exact, i.e.\ $S^{-1}R$ is flat over~$R$. 

\nxsubpoint (Degree zero part of a localization.)
For any $f\in|R|$ with $\deg f$ invertible in $\Delta$
we denote by $R_{(f)}$ the degree zero component $(R_f)^0$
of localization $R_f=R[f^{-1}]$, and similarly $M_{(f)}:=(M_f)^0$
for any graded $R$-module~$M$. Clearly $M_{(f)}$ is a $R_{(f)}$-module;
it consists of fractions $x/f^n$, where $x\in M$, $n\geq0$ and 
$\deg x=n\deg f$. We have $x/f^n=y/f^m$ iff $f^{m+N}x=f^{n+N}y$ for
some $N\geq0$. Components $R_{(f)}(n)=R(n)_{(f)}$ of 
generalized ring $R_{(f)}$ admit a similar description. Notice that
$M_{(f)}$ depends only on components $M^\alpha$ with $\alpha=n\deg f$
for some $n\geq0$, and similarly for $R_{(f)}$.

\nxsubpoint\label{sp:loc.homog.loc} (Localization of $R_{(f)}$.)
Fix any element $\bar g=g/f^n\in|R_{(f)}|=R_{(f)}(1)$, where 
necessarily $g\in|R|$ is of degree $n\deg f$. One sees immediately
that $(R_{(f)})_{\bar g}=R_{(f)}[\bar g^{-1}]$ can be identified
with degree zero part of $(R_f)_{\bar g}=R[f^{-1},g^{-1}]=R_{fg}$,
i.e.\ with $R_{(fg)}$. More generally, if $f$, $g\in|R|$ are such that
$n\deg f=m\deg g$ for some integers $n$, $m\geq0$, then
$R_{(fg)}$ is canonically isomorphic both to $R_{(f)}[(g^m/f^n)^{-1}]$
and $R_{(g)}[(f^n/g^m)^{-1}]$.

\nxsubpoint (Cases $\Delta=\bbN_0$ and $\Delta=\bbZ$.)
Henceforth all graded rings and modules we consider will be
either $\bbZ$-graded or $\bbN_0$-graded. Furthermore,
any $\bbN_0$-graded ring~$R$ will be treated as a 
{\em positively\/} $\bbZ$-graded
algebra by putting $R_d(n):=\emptyset$ for $d<0$. Any $\bbN_0$-graded 
module~$M$ over such a ring~$R$ can be also treated as a positively
$\bbZ$-graded module by putting $M_d:=\emptyset$ for $d<0$.
Notice that now we use lower indices to denote graded components,
just to avoid confusion with direct power
$M^d$ (product of several copies of~$M$).

\nxsubpoint (Projective spectrum of a positively graded generalized ring.)
Let $R$ be a positively graded generalized ring, considered here
as a $\bbZ$-graded generalized ring. Consider monoid 
$|R|:=\bigsqcup_{n\in\bbZ}R_n(1)=\bigsqcup_{n\geq0}R_n(1)$, and 
its submonoid $|R|^+$ consisting of elements of strictly positive degree.
For any $f\in|R|^+$ we denote affine generalized scheme
$\Spec R_{(f)}$ by $D_+(f)$. If $g\in|R|^+$ is another such element,
then $D_+(f)$ contains principal open subset
$D_+(f)_g:=\Spec (R_{(f)})[(g^{\deg f}/f^{\deg g})^{-1}]$, 
canonically isomorphic
to $D_+(fg)=\Spec R_{(fg)}$ by~\ptref{sp:loc.homog.loc},
and similarly for $D_+(g)$. Gluing together generalized ringed spaces
$\{D_+(f)\}_{f\in|R|^+}$ along their isomorphic open subsets
$D_+(f)_g\cong D_+(fg)\cong D_+(g)_f$ in the usual fashion,
we obtain a generalized ringed space $\Proj R$, which
admits an open cover by subspaces isomorphic to $D_+(f)$ 
(so they will be identified with $D_+(f)=\Spec R_{(f)}$),
in such a way that $D_+(f)\cap D_+(g)=D_+(fg)$ inside $\Proj R$,
with induced morphisms $D_+(fg)\to D_+(f)$ equal to those considered above.

Since $\{D_+(f)\}_{f\in|R|^+}$ constitute an affine open cover,
$\Proj R$ is a generalized scheme (for the localization theory chosen).
It will be called the {\em projective spectrum\/} of~$R$.

\nxsubpoint (Quasicoherent sheaf defined by a graded module.)
Let $M$ be a graded module over a positively graded generalized ring~$R$.
We can construct a quasicoherent sheaf $\tilde M$ on $P:=\Proj R$ as follows.
Consider quasicoherent sheaves $\widetilde{M_{(f)}}$ on~$D_+(f)$;
canonical isomorphisms $(M_{(f)})_{g^{\deg f}/f^{\deg g}}\cong
M_{(fg)}$ induce canonical isomorphisms $\widetilde{M_{(f)}}
|_{D_+(f)\cap D_+(g)}\cong\widetilde{M_{(fg)}}\cong
\widetilde{M_{(g)}}|_{D_+(f)\cap D_+(g)}$, so these quasicoherent sheaves
can be patched together into a quasicoherent $\sO_P$-module
denoted by~$\tilde M$.

This construction is obviously functorial, so we get a functor
$\catGradMod R\to\catQCoh{\sO_{\Proj R}}$. This functor is left exact.
It even commutes with arbitrary inductive limits if $\Proj
R=\bigcup_{f\in R_1}D_+(f)$, individual functors $M\mapsto M_{(f)}$
having this property whenever $\deg f=1$.  Since the localization
functor $\catGradMod R\to\catGradMod{R_f}$, $M\mapsto M_f$, is exact
and commutes with arbitrary inductive limits, the only complicated
point is to show that left exact functor $\catGradMod{R_f}\to
\catMod{R_{(f)}}$, $N\mapsto N_0$, commutes with inductive limits.

One checks for this that the graded homomorphism $R_{(f)}[T^{\pm1}]\to
R_f$ of $R_{(f)}$-algebras, induced by $T\mapsto f$, where $T$ is a
unary generator of degree one, is in fact an isomorphism, i.e.\ $R_f$
is a {\em unary\/} graded $R_{(f)}$-algebra, hence the inductive
limits in $\catGradMod{R_f}$ can be computed degreewise in
$\catMod{R_{(f)}}$, hence $N\mapsto N_0$ commutes with arbitrary
inductive limits.

Another proof: if $\Proj R=\bigcup_{f\in R_1}D_+(f)$, we can prove
that $\Delta$ commutes with arbitrary inductive limits by constructing
explicitly a right adjoint functor $\Gamma_*$, cf.~\ptref{sp:adj.delta.gamma*}.

\nxsubpoint\label{sp:covers.of.projspec} (Covers of $\Proj R$.)
Let $(f_\alpha)$ be any family of elements of $|R|^+$.
Denote by $\ga\subset|R|$ the graded ideal
(i.e.\ $R$-submodule of $|R|=R(1)$) generated by $f_\alpha$.
Notice that $\ga$ is contained in $|R|^+$ since $|R|^+$ is itself
a graded ideal in $|R|$, $R$ being supposed to be positively graded.
We claim that {\em $D_+(g)\subset\bigcup_\alpha D_+(f_\alpha)$
iff $g^n\in\ga$ for some $n>0$, i.e.\ iff $g$ belongs to $\rad(\ga)$,
the radical of~$\ga$.} Since $\Proj R$ is by construction
the union of all such $D_+(g)$, we conclude that {\em
$\bigcup_\alpha D_+(f_\alpha)=\Proj R$ iff $\rad(\ga)=|R|^+$.}

Let us show our statement. Since $D_+(g)\cap D_+(f_\alpha)$ is
identified with $D(\bar f_\alpha)\subset D_+(g)=\Spec R_{(g)}$, where
$\bar f_\alpha=f_\alpha^{\deg g}/g^{\deg f_\alpha}$, we see that
$D_+(g)\subset\bigcup_\alpha D_+(f_\alpha)$ iff $D(\bar f_\alpha)$
constitute a cover of $D_+(g)=\Spec R_{(g)}$.  According
to~\ptref{sp:covers.un.spec} and~\ptref{sp:opcov.specp}, this is
equivalent to saying that the ideal in $R_{(g)}$ generated by the
$\bar f_\alpha$ contains the identity of $R_{(g)}$.  (Notice that such
a family constitutes a cover or not independently on the choice of
localization theory, since this issue can be reduced to whether this
is a family of efficient descent for $\stQCOH$ or not, a condition
independent on the choice of theory.)  One checks, more or less in the
classical fashion, that this ideal equals $\gb_{(g)}$, where $\gb$
denotes the graded ideal in $|R|$ generated by the $f_\alpha^{\deg g}$
(the key point here is that $\Delta:M\mapsto\tilde M$ always
preserves strict epimorphisms, even when it is not known to be right
exact); therefore, $1\in\gb_{(g)}$ iff $g^n\in\gb$ for some $n>0$ by
construction of $\gb_{(g)}$, i.e.\ iff $g\in\rad(\gb)=\rad(\ga)$,
q.e.d.

\nxsubpoint (Radicals of graded ideals. Graded prime ideals.)
The above reasoning used the following fact: {\em the radical
of a graded ideal in $|R|$ is itself a graded ideal.}
Indeed, we knew that $\gb\subset\ga$, so $\rad(\gb)\subset\rad(\ga)$;
on the other hand, all generators $f_\alpha$ of~$\ga$ belong to
$\rad(\gb)$, so $\ga\subset\rad(\gb)$ and $\rad(\ga)\subset\rad(\gb)$.
In order to conclude $\ga\subset\rad(\gb)$ we need to know that
$\rad(\gb)$ is a graded ideal. This can be shown in several different
ways. For example, we can introduce the notion of a 
{\em graded prime ideal\/} $\gp\subset|R|$, i.e.\ a graded ideal,
such that $|R|-\gp$ is a submonoid of $|R|$, and show that
$\rad(\gb)$ coincides with the intersection of all graded prime
ideals containing~$\gb$, more or less in the same way as 
in~\ptref{sp:rad.ideal}.

Furthermore, if we work in the unary localization theory $\cT^u$,
we can describe $\Proj R$ as the set of all graded prime ideals
$\gp\not\supset|R|^+$, and $D_+(f)$ is identified then with
the set of all $\gp\not\ni f$, exactly as in the classical 
theory of projective spectra of EGA~II~2.

\nxsubpoint\label{sp:un.proj.space} (Projective spaces.)
Now we are in position to define the {\em (unary) projective space
$\bbP^n_C$} over a generalized ring~$C$. By definition,
this is nothing else than 
$\Proj C[T_0^{[1]},T_1^{[1]},\ldots,T_n^{[1]}]$, the projective
spectrum of the unary polynomial algebra over~$C$ generated
by $n+1$ indeterminates of degree one. Our previous results show
that $D_+(T_i)$, $0\leq i\leq n$, constitute an affine cover
of $\bbP^n_C$, and that each of $D_+(T_i)$ is isomorphic to
$\bbA^n_C=\Spec[T_1^{[1]},\ldots,T_n^{[1]}]$, the unary
affine space over~$C$. For example, $\bbP^1_C$ is constructed
by gluing together two copies of the affine line in the classical fashion. 
In particular,
$\bbP^1_{\Fone}$ consists of three points (at least if 
we use unary localization theory): two closed points
$0$ and $\infty$, corresponding to graded prime ideals
$(T_0)$ and $(T_1)$ in $\Fone[T_0,T_1]$, and one generic point~$\xi$,
corresponding to the zero ideal $(0)$.

\nxsubpoint\label{sp:Serre.twists} (Serre twists.)
Suppose that $\Proj R$ is quasicompact, so that it can be 
covered by some $D_+(f_1)$, \dots, $D_+(f_n)$ with $\deg f_i>0$. 
We may even assume that all $\deg f_i$ are equal to some $N>0$,
by replacing $f_i$ with appropriate powers, and then that
all $\deg f_i=1$, replacing $R$ by $R^{[N]}$, given by
$R^{[N]}_d(n):=R_{Nd}(n)$ (this is possible at least if $R_0$ admits a zero), 
similarly to the constructions of EGA~II, 2.4.7.

One easily checks that $M_{(f)}\cong M[n]_{(f)}$ 
for any graded $R$-module $M$ and any $n\in\bbZ$ whenever $\deg f=1$.
Therefore, if we assume that $D_+(f)$ with $\deg f=1$ cover $\Proj R$,
we see that $\widetilde{M[n]}$ is locally isomorphic to~$\tilde M$.

Applying this to $M=|R|=R(1)$, we obtain {\em Serre line bundles\/}
$\sO[n]=\sO_{\Proj R}[n]:=\widetilde{R(1)[n]}$, for all $n\in\bbZ$.
They are indeed line bundles, since they are locally isomorphic to
$\sO$ by construction. They have their usual properties
$\sO[0]=\sO$, $\sO[n]\otimes_\sO\sO[m]\cong\sO[n+m]$ and
$\widetilde{M[n]}\cong\tilde M\otimes_\sO\sO[n]$. We define
{\em Serre twists\/} $\sF[n]:=\sF\otimes_\sO\sO[n]$
for any $\sO_{\Proj R}$-module~$\sF$.

When no confusion can arise, we use classical notation
$\sO(n)$ and $\sF(n)$ instead of $\sO[n]$ and $\sF[n]$.
Notice, however, that in our context $\sO(n)$ sometimes means
the sheaf of $n$-ary operations of~$\sO$, i.e.\ the free 
$\sO$-module $L_\sO(n)$ of rank~$n$.

\nxsubpoint\label{sp:constr.gamma*} (Functor $\Gamma_*$.)
Now suppose that we are given a generalized scheme~$X$
and a line bundle~$\sL$ on~$X$, i.e.\ a (necessarily quasicoherent)
$\sO_X$-module, locally isomorphic to $|\sO_X|$. Denote by
$\sL^{\otimes n}$, $n\in\bbZ$, the tensor powers of~$\sL$.
We can define
a positively graded generalized ring $R=\Gamma_*(\sO_X)$ by putting
\begin{equation}
R_d(n):=\Gamma(X,L_{\sO_X}(n)\otimes_{\sO_X}\sL^{\otimes d}),\quad n,d\geq0,
\end{equation}
where $L_{\sO_X}(n)$ denotes the free $\sO_X$-module of rank~$n$,
i.e.\ the trivial vector bundle of rank~$n$.
We define the ``multiplication maps'' $R_d(k)\times R_e(n)^k\to R_{d+e}(n)$
by ``twisting'' the ($\sO_X$-bilinear) multiplication maps 
$\mu_{\sO_X,n}^{(k)}:\sO_X(k)\times\sO_X(n)^k\to\sO_X(n)$
for~$\sO_X$ by appropriate powers of $\sL$ (notice that tensoring with
a line bundle is an exact functor) and taking global sections.
(Commutativity of $\sO_X$ is paramount here: we need to know that 
the multiplication maps we consider here are $\sO_X$-bilinear maps
of $\sO_X$-modules.) All relations implied in our definition of 
a graded generalized ring are easily seen to hold; actually
\ptref{sp:graded.alg.mon} was originally motivated exactly by this situation.

Next, for any $\sO_X$-module $\sF$ we can define a graded 
$\Gamma_*(\sO_X)$-module $M:=\Gamma_*(\sF)$ by putting
\begin{equation}
M_d:=\Gamma(X,\sF\otimes_{\sO_X}\sL^{\otimes d})
\end{equation}
Here we might use all $d\in\bbZ$, or restrict our attention only to $d\geq0$,
whatever we like most.

In this way we obtain a functor $\Gamma_*:\catQCoh{\sO_X}\to
\catGradMod R$, obviously left exact.

\nxsubpoint (Application to $X=\Proj R$, $\sL=\sO_X[1]$.)
Now let $R$ be a positively graded generalized ring, such that $X:=\Proj R$
is covered by $D_+(f)$ with $f\in R_1(1)$. Then $\sL:=\sO_X[1]$ is
a line bundle by~\ptref{sp:Serre.twists}, so we can apply the above
constructions to this case. Notice that for any graded $R$-module~$M$
we obtain a natural map of $\bbZ$-graded sets $\xi_M:M\to\Gamma_*(\tilde M)$.
Its component of degree~$d$, $\xi_{M,d}:M_d\to\Gamma_d(\tilde M)=
\Gamma(X,\tilde M\otimes\sL^{\otimes d})=\Gamma(X,\widetilde{M[d]})$
maps an element $x\in M_d$ into global section $\xi(x)\in\Gamma(X,
\widetilde{M[d]})$, such that $\xi(x)|_{D_+(f)}=x/1$ for any $f\in|R|^+$.

Applying this to free modules $M=R(n)$, we get maps of graded sets
$R(n)\to(\Gamma_*(\sO_X))(n)$, easily seen to combine together into
a graded generalized ring homomorphism $R\to\Gamma_*(\sO_X)$. 
Now if $\sF$ is any $\sO_X$-module, $\Gamma_*(\sF)$ can be
considered as an $R$-module by means of scalar restriction via
$R\to\Gamma_*(\sO_X)$. The map of graded sets $\xi_M:M\to\Gamma_*(\tilde M)$
is easily seen to be a homomorphism of graded $R$-modules, functorial
in~$M$, i.e.\ $\xi:\Id_{\catGradMod R}\to\Gamma_*\circ\Delta$
is a natural transformation, where $\Delta:\catGradMod R\to\catQCoh{\sO_X}$ 
denotes the functor $M\mapsto\tilde M$.

\nxsubpoint\label{sp:adj.delta.gamma*} 
(Adjointness of $\Delta$ and $\Gamma_*$.)
Let us show that in the above situation functors $\Delta$ and $\Gamma_*$,
considered as functors between $\catGradMod R$ and $\catQCoh{\sO_X}$
(or even $\catMod{\sO_X}$), are adjoint, and that 
$\xi_M:M\to\Gamma_*(\tilde M)$ is one of adjointness morphisms.
We show this by constructing the other adjointness morphism
$\eta_\sF:\widetilde{\Gamma_*(\sF)}\to\sF$, functorial in~$\sF$.
It would suffice to define $\eta_\sF|_{D_+(f)}:
\widetilde{\Gamma_*(\sF)}|_{D_+(f)}\to\sF|_{D_+(f)}$ for all $f\in|R|^+$.
Since the source of this map is the quasicoherent sheaf on affine scheme
$D_+(f)=\Spec R_{(f)}$ given by module $\Gamma_*(\sF)_{(f)}$,
we just have to define $R_{(f)}$-homomorphism
$\eta'_{\sF,f}:\Gamma_*(\sF)_{(f)}\to\Gamma(D_+(f),\sF)$.
We do this by mapping an element $x/f^n$, where $x\in\Gamma_{n\deg f}(\sF)=
\Gamma(X,\sF[n\deg f])$, into $x|_{D_+(f)}\otimes f^{-n}$, or equivalently,
into the preimage of $x|_{D_+(f)}$ under the isomorphism 
$\id_{\sF|_{D_+(f)}}\otimes f^n$, obtained by tensoring $\sF$ with 
$f^n:\sO_X|_{D_+(f)}\simto\sO_X[n\deg f]|_{D_+(f)}$.

Once $\eta$ is constructed, the verification of required compatibility
between $\xi$ and $\eta$ is done exactly in the same way as in EGA~II 2.6.5,
since related proofs of {\em loc.cit.} are ``multiplicative'',
i.e.\ don't use addition at all.

The following lemma is a counterpart of EGA~I~9.3.2:
\begin{LemmaD}\label{l:extend.sect} 
Let $X$ be quasicompact quasiseparated generalized
scheme, $\sL$ be a line bundle over~$X$, $f\in\Gamma(X,\sL)$ be
a global section, $X_f\subset X$ be the largest open subscheme of~$X$,
such that $f$ is invertible on~$X_f$. Let $\sF$ be a quasicoherent 
sheaf on~$X$. Then:
\begin{itemize}
\item[(a)] If $a$, $b\in\Gamma(X,\sF)$ are two global sections of
$\sF$, such that $a|_{X_f}=b|_{X_f}$, then $af^n=bf^n\in
\Gamma(X,\sF\otimes\sL^{\otimes n})$ for some integer $n\geq0$.
\item[(b)] For any section $a'\in\Gamma(X_f,\sF)$ one can
find an integer $n\geq0$ and a global section
$a\in\Gamma(X,\sF)$, such that $a'f^n=a|_{X_f}$.
\end{itemize}
\end{LemmaD}
\begin{Proof} Proof goes essentially in the same way as in {\em loc.cit.}
Since $X$ is quasicompact and $\sL$ locally trivial, we can find
a finite affine open cover $U_i=\Spec A_i$ of~$X$, $1\leq i\leq n$,
such that $\sL|_{U_i}\cong\sO_X|_{U_i}$. Let $f_i\in A_i=\Gamma(U_i,\sO_X)$
be the image of $f|_{U_i}$ under this isomorphism; then
$X_f\cap U_i=U_{i,f}=\Spec A_i[f_i^{-1}]=D(f_i)\subset\Spec A_i$.

(a) Since the cover $U_i$ is finite, it is enough to show 
existence of $n\geq0$, such that $af^n|_{U_i}=bf^n|_{U_i}$,
separately for each~$i$, i.e.\
we can replace $X$ by $U_i$ and assume $X=\Spec A$ affine, 
$\sL=\sO_X$ trivial, $f\in\Gamma(X,\sO_X)=A$, and $a$, $b\in M:=\Gamma(X,\sF)$
such that $a|_{D(f)}=b|_{D(f)}$, i.e. $a/1=b/1$ in $M_f$. This implies 
existence of $n\geq0$, such that $af^n=bf^n$, by definition of localizations.

(b) Reasoning as in (a), we can find an integer $n\geq0$, same for 
all~$i$ if we want, such that $a'f^n|_{U_i\cap X_f}$ extends to 
a section $a_i$ of $\sF$ over~$U_i$. However, $a_i|_{U_i\cap U_j}$
needn't coincide with $a_j|_{U_i\cap U_j}$; all we know is that
their restrictions to $U_i\cap U_j\cap X_f$ are equal, so applying (a) to
quasicompact generalized scheme $U_i\cap U_j$ we obtain
$a_if^m|_{U_i\cap U_j}=a_jf^m|_{U_i\cap U_j}$, and we can choose
an $m\geq0$ valid for all $1\leq i,j\leq n$. Then $a_if^m\in\Gamma(U_i,\sF)$ 
coincide on intersections of these open sets, so they can be glued to
a global section $a\in\Gamma(X,\sF)$, and obviously $a|_{X_f}=a'f^{n+m}$, 
q.e.d.
\end{Proof}

Once we have the above lemma, we can prove the following:
\begin{PropD}\label{prop:qcoh.on.proj}
Let $R$ be a positively graded generalized ring, such that
$P:=\Spec R$ is union of finitely many $D_+(f)$ with $\deg f=1$.
Then $\eta_{\sF}:\widetilde{\Gamma_*(\sF)}\to\sF$ 
is an isomorphism for any quasicoherent
$\sO_P$-module~$\sF$, functor $\Gamma_*:\catQCoh{\sO_P}\to\catGradMod R$
is fully faithful, and its left adjoint $\Delta:\catGradMod R\to
\catQCoh{\sO_P}$ is a localization functor
(in particular it is essentially surjective), i.e.\ $\catQCoh{\sO_P}$
is equivalent to the localization of $\catGradMod R$ with respect to
set of all homomorphisms in $\catGradMod R$, which become isomorphisms
after application of~$\Delta$, and this localization can be computed
by means of left fraction calculus.
\end{PropD}
\begin{Proof}
The statement about $\eta_{\sF}$ is deduced from~\ptref{l:extend.sect}
in the same way as in EGA~II 2.7.5, and the remaining statements
follow for general (category-theoretic) reasons, cf.~\cite{GZ}.
\end{Proof}

\nxsubpoint ($\Delta$ preserves finite type and presentation.)
Suppose that a positively graded ring~$R$ is such that
$P:=\Proj R=\bigcup_{f\in R_1}D_+(f)$, e.g.\ ideal $|R|^+$ is generated
by~$R_1$. Then $\widetilde{|R|[n]}=\sO_P[n]$ is a locally free
and in particular finitely presented $\sO_P$-module. Therefore,
if $M$ is a (graded) free $R$-module of finite rank~$k$, i.e.\
if $M=\bigoplus_{i=1}^k|R|[n_i]$, then $\tilde M$ is a vector bundle
of rank~$k$. Since $\Delta:M\mapsto\tilde M$ is right exact,
we see that $\tilde M$ is a quasicoherent $\sO_P$-module of
finite type (resp.\ finite presentation) whenever
$M$ is a graded $R$-module of finite type (resp.\ finite presentation).

\nxsubpoint\label{sp:fgp.qcoh.proj} (Partial converse.)
Let's suppose in addition that $P=\Proj R$ be quasicompact,
i.e.\ we are under the conditions of~\ptref{prop:qcoh.on.proj}.
Let $\sF$ be any quasicoherent $\sO_P$-module of finite type.
According to~\ptref{prop:qcoh.on.proj}, we can identify $\sF$ with
$\tilde M$ for a suitable graded $R$-module $M$, e.g.\ $M=\Gamma_*(\sF)$.
Next, we can write $M=\injlim_\alpha M_\alpha$, filtered inductive
limit of its graded submodules of finite type. Then $\sF=\tilde M=
\injlim_\alpha \tilde M_\alpha$, and using quasicompactness of $P$
we find an index $\alpha$, such that $\tilde M_\alpha=\tilde M=\sF$
(cf.~\ptref{sp:qcoh.ft.are.fg} for a similar reasoning),
i.e.\ {\em any quasicoherent $\sO_P$-module~$\sF$ of finite type
is isomorphic to $\tilde N$ for some graded $R$-module $N$ 
of finite type.} Reasoning further as in~\ptref{sp:qcoh.fpres},
we obtain a similar description of finitely presented $\sO_P$-modules:
{\em any finitely presented $\sO_P$-module~$\sF$ is isomorphic to
$\tilde N$ for some finitely presented graded $R$-module~$N$.}

Notice that we don't claim here $\Gamma_*(\sF)$ to be
a finitely generated or presented graded $R$-module.

\nxsubpoint (Finitely presented $\sO_P$-modules as a fraction category.)
Let's keep the previous assumptions. Let $\cC$ be the category
of finitely presented graded $R$-modules, $\cD$ be the category
of finitely presented $\sO_P$-modules, denote by $S\subset\Ar\cC$
the set of morphisms in $\cC$ that become isomorphisms after
an application of $\Delta:\cC\to\cD$. By universal property of
localizations of categories we obtain a functor
$\bar\Delta:\cC[S^{-1}]\to\cD$, essentially surjective
by~\ptref{sp:fgp.qcoh.proj}. We claim the following:

\begin{Propz} Subset $S\subset\Ar\cC$ admits a left fraction
calculus, so that morphisms $f':X\to Y$ in $\cC[S^{-1}]$ can be
described as $s^{-1}f$ for some $f:X\to Z$ and $s:Y\to Z$, $s\in S$
(cf.~{\rm \cite[I.2.3]{GZ}}), and $\bar\Delta:\cC[S^{-1}]\to\cD$
is an equivalence of categories.
\end{Propz}
\begin{Proof}
(a) Conditions a)--d) of \cite[I.2.2]{GZ} immediately follow from the
fact that $\cC$ is stable under finite inductive limits of $\catGradMod R$,
and that $\Delta:\cC\to\cD$ is right exact, i.e.\ commutes with such limits.
Therefore, $S$ admits left fraction calculus.

(b) Functor $\bar\Delta$ is faithful, i.e.\ injective on morphisms.
Indeed, suppose that $s^{-1}f$ and $t^{-1}g\in\Hom_{\cC[S^{-1}]}(X,Y)$
are such that $\tilde s^{-1}\tilde f=\tilde t^{-1}\tilde g$ in~$\cD$.
Constructing a pushout diagram on $s$ and~$t$ and using obvious
stability of $S$ under pushouts, we get a diagram
\begin{equation}
\xymatrix{
X\ar[rrd]_{f}\ar[rrrrd]^<>(.3){g}&&&Y\ar[ld]_{s}\ar[rd]^{t}\\
&&Y'\ar[rd]^{t'}&&Y''\ar[ld]_{s'}\\
&&&Z\ar@{-->}[r]^{u}&W}
\end{equation}
Put $h:=t'f$, $h':=s'g$, $W:=\Coker(h,h':X\to Z)$. Clearly
$\tilde{h}=\tilde{h'}$, hence $Z\stackrel u\to W$ belongs to~$S$,
$\Delta$ being exact, $W$ is finitely presented, being a finite inductive
limit of finitely presented graded $R$-modules, and obviously
$ut'f=uh=uh'=us'g$ and $ut's=us't$; we conclude that $s^{-1}f=t^{-1}g$ in
$\cC[S^{-1}]$.

(c) Functor $\bar\Delta$ is surjective on morphisms, hence 
fully faithful by~(b). Indeed, let $\phi:\tilde M\to\tilde N$ be
any $\sO_P$-morphism, for any $M$, $N\in\Ob\cC$. Put $N':=\Gamma_*(\tilde N)$.
By adjointness of $\Delta$ and $\Gamma_*$ we get graded homomorphisms
$f'=\phi^\flat:M\to N'$ and $s'=\id_{\tilde N}^\flat:N\to N'$, such that
$\tilde s'$ is an isomorphism by~\ptref{prop:qcoh.on.proj}, and
$\phi=(\tilde s')^{-1}\circ\tilde f'$. Now write $N'$ as a filtered
inductive limit of finitely presented modules $N_\alpha$ and
use that $\Hom(M,-)$ and $\Hom(N,-)$ commute with filtered inductive
limits, $M$ and $N$ being finitely presented. The conclusion is that
both $f'$ and $s'$ factorize through some finitely presented $N_\alpha$,
and by taking an even larger $\alpha$ we can assume $i_\alpha:N_\alpha\to N'$
to belong to~$S$. Then $s:N\to N_\alpha$ will also belong to~$S$,
since both $i_\alpha$ and $s'=i_\alpha s$ belong to~$S$, and
considering $f:M\to N_\alpha$, such that $f'=i_\alpha f$, we get
$\phi=\tilde s^{-1}\tilde f=\bar\Delta(s^{-1}f)$.

(d) Finally, $\bar\Delta$ is essentially surjective 
by~\ptref{sp:fgp.qcoh.proj}, hence an equivalence of categories, q.e.d.
\end{Proof}

\nxsubpoint (Almost isomorphisms of graded modules.)
We say that a homomorphism (of degree zero) $f:M\to N$ of
graded $R$-modules is {\em almost isomorphism\/} if
$f_n:M_n\to N_n$ is a bijection for $n\gg0$. Clearly,
any almost isomorphism $f$ induces an isomorphism 
$\tilde f:\tilde M\to\tilde N$ between corresponding quasicoherent sheaves.
The converse seems to hold only for finitely presented $M$ and~$N$,
and only under some additional restrictions, e.g.\ 
when $R$ is generated as a graded $R_0$-algebra
by a finite subset of~$R_1$, and in particular if $R$ is a finitely presented
graded unary $R_0$-algebra. This includes the projective spaces
$\bbP^n_{R_0}=\Proj R_0[T_0^{[1]},\ldots,T_n^{[1]}]$. However,
we don't want to provide more details here.

\nxsubpoint (Functoriality.)
Any homomorphism of graded generalized rings $\phi:R\to R'$ induces
an affine morphism of generalized schemes $\Proj(\phi):=\Phi:
G(\phi)\to\Proj R$, where
$G(\phi)\subset\Proj R'$ is the open subscheme of~$\Proj R'$ equal to
$\bigcup_{f\in|R|^+}D_+(\phi(f))$, essentially in the same way
as in EGA~II 2.8, i.e.\ $\Phi^{-1}(D_+(f))=D_+(\phi(f))$, and
the restriction of $\Phi$ to $D_+(\phi(f))$ is induced by canonical
homomorphism $R_{(f)}\to R'_{(\phi(f))}$. We don't need to
provide more details since the reasoning of {\em loc.cit.}\ is completely
``multiplicative'', i.e.\ doesn't use addition in any essential way.

\nxsubpoint (Base change.)
Let $C\to C'$ be an extension of generalized rings,
$R$ be a positively graded $C$-algebra. Consider $R':=R\otimes_CC'$
(understood for example as a pushout in the category of
graded generalized rings, where $C$ and $C'$ are considered with
their trivial grading). Graded homomorphism $\phi:R\to R'$
induces a morphism $\Phi:G(\phi)\to\Proj R$, where
$G(\phi)=\bigcup_{f\in|R|^+}D_+(\phi(f))\subset\Proj R'$. 
We obtain a commutative diagram
\begin{equation}
\xymatrix{
G(\phi)\ar[r]^{\Phi}\ar[d]&\Proj R\ar[d]\\
\Spec C'\ar[r]&\Spec C}
\end{equation}
One checks, essentially in the same way as in EGA~II, 2.8.9, that the
above diagram is cartesian (one just checks
$R'_{(f)}=R_{(f)}\otimes_CC'$ for this), and that in fact
$G(\phi)=\Proj R'$, at least if $R$ is a unary graded $C$-algebra
or if $C'$ is flat over~$C$ (this follows from the fact that 
in any of these cases the graded
$R'$-ideal $|R'|^+$ is generated by the image of $|R|^+$ in $|R'|$,
a consequence of the graded version of the ``affine base change''
theorem~\ptref{th:aff.base.change}).

Therefore, construction of $\Proj R$ always commutes with 
{\em flat\/} base change, and if $R$ is a unary graded $C$-algebra,
then it commutes with arbitrary base change.

\nxsubpoint (Quasicoherent sheaves of graded algebras.)
Now let $S$ be any generalized scheme. One can define
a {\em quasicoherent sheaf of graded $\sO_S$-algebras\/} $\sR$
by combining together~\ptref{sp:graded.alg.mon} and~\ptref{sp:qcoh.alg} 
in the natural way.
Any such $\sR$ has the property that $\sR|_U=\tilde R$ for any
open affine $U=\Spec C\subset S$, where $R$ is some graded $C$-algebra.
Reasoning further as in EGA II, 3.1, we construct a generalized 
$S$-scheme $\Proj\sR$ by gluing together $\Proj\Gamma(U,\sR)$ for
all open affine $U\subset S$. Graded quasicoherent sheaves of $\sR$-modules
on~$S$ define quasicoherent sheaves on $\Proj\sR$ by the same reasoning
as in {\em loc.cit.}, and further properties like the 
relative version of \ptref{prop:qcoh.on.proj} hold as well, at least if
we assume $S$ to be quasicompact and quasiseparated.

\nxsubpoint (Projective bundles.)
An important application is given by {\em projective bundles.}
Let $S$ be a generalized scheme, $\sE$ be a vector bundle over~$S$
(i.e.\ a locally free $\sO_S$-module of finite type), or even
any quasicoherent $\sO_S$-module. Consider
the symmetric algebra $\sR:=S_{\sO_S}(\sE)=\bigoplus_{n\geq0}S^n_{\sO_S}(\sE)$.
It is a $\bbZ$-graded algebra in the $\otimes$-category of 
quasicoherent $\sO_S$-modules, hence it can be considered as a
quasicoherent sheaf of unary graded $\sO_S$-algebras. Putting
$\bbP(\sE):=\Proj_{\sO_S}\sR$, we construct the
{\em projective bundle over~$S$ corresponding to vector bundle~$\sE$.}
If $\sE$ is a vector bundle of rank~$n+1$, 
then $\bbP(\sE)$ is locally (over~$S$) isomorphic to $\bbP^n_S/S$.

Some of basic properties of projective bundles known from EGA~II
generalize to our case. For example, canonical homomorphism
$\sE\otimes_{\sO_S}S_{\sO_S}(\sE)\to S_{\sO_S}(\sE)[1]$ induces
a surjection $\pi:p^*\sE\twoheadrightarrow\sO[1]$ on $\bbP(\sE)$,
where $p:\bbP(\sE)\to S$ denotes the structural morphism.

\nxsubpoint\label{sp:sect.proj.bdl} (Sections of projective bundles.)
Notice that sections $\Gamma_S(\bbP_S(\sE)/S)$ are in one-to-one
correspondence with strict quotients $\sE\twoheadrightarrow\sL$, which
are line bundles over~$S$. Namely, given an $S$-section
$\sigma:S\to\bbP_S(\sE)$, we obtain such a strict epimorphism
$\phi:=\sigma^*(\pi):\sE=\sigma^*p^*\sE
\twoheadrightarrow\sigma^*\sO[1]=:\sL$ by applying $\sigma^*$
to~$\pi:p^*\sE\twoheadrightarrow\sO[1]$.  Conversely, given
$\phi:\sE\twoheadrightarrow\sL$, we construct a section $\sigma:=\Proj
S_{\sO_S}(\phi): S\cong\bbP_S(\sL)\to\bbP_S(\sE)$. This morphism is
defined on the whole of $\bbP_S(\sL)=S$ since $S_{\sO_S}(\phi):
S_{\sO_S}(\sE)\to S_{\sO_S}(\sL)$ is surjective.  One checks that
these two maps $\sigma\leftrightarrow\phi$ are indeed inverse to each
other, essentially as in EGA~II 4.2.3. Noticing that the construction
of $\bbP_S(\sE)$ commutes with any base change $S'\stackrel f\to S$,
we see that the set of $S'$-valued points of $S$-scheme $\bbP_S(\sE)$
can be identified with the set of those strict quotients $\sL$ of
$\phi^*\sE$, which are line bundles on~$S'$, exactly as in the
classical case.

Another remark: in general $\sigma$ will be a ``closed'' immersion
in the sense of~\ptref{sp:closed.imm.sch}, 
for example because of surjectivity
of $S_{\sO_S}(\phi)$, or because of $\bbP_S(\sE)$ being ``separated'' 
over~$S$.

\nxsubpoint\label{sp:proj.spc.points} 
(Example: points of $\bbP^n$ over some generalized rings.)
This is applicable in particular to $\Fempty$-scheme
$\bbP^n=\bbP^n_{\Fempty}:=\Proj\Fempty[T_0^{[1]},\ldots,T_n^{[1]}]=
\bbP(\Fempty(n+1))$. We see that $R$-valued points of $\bbP^n$
correspond to strict quotients~$P$ of free $R$-module $L_R(n+1)=R(n+1)$,
which are projective ``of rank~$1$'' (i.e.\ $\tilde P$ has to be
a line bundle over $\Spec R$), for any generalized ring~$R$.
When $\Pic(R)=\Pic(\Spec R)$ is trivial (e.g.\ when $R$ is $\cT^?$-local), 
$P$ has to be isomorphic
to $|R|$, i.e.\ we have to consider the set of all surjective
$R$-homomorphisms $\phi:L_R(n+1)\twoheadrightarrow L_R(1)=|R|$,
modulo multiplication by invertible elements of~$\Aut_R(|R|)=|R|^{\times}$.
By the universal property of free modules $\phi$ is completely determined
by the set $(X_0,\ldots,X_n)\in|R|^n$ of its values $X_i:=\phi(\{i+1\})$
on the basis elements of~$L_R(n+1)$, and the image of $\phi$ is exactly
the ideal $\langle X_0,\ldots,X_n\rangle\subset|R|$. Therefore,
{\em if $\Pic(\Spec R)=0$, e.g.\ if $R$ is $\cT^?$-local, then
$\bbP^n(R)$ consists of $(n+1)$-tuples $(X_0,\ldots,X_n)\in|R|^{n+1}$ of
elements of $|R|$, generating the unit ideal in $|R|$, modulo
multiplication by invertible elements from~$|R|^\times$}
(cf.\ EGA~II, 4.2.6). Of course, the class of $(X_0,\ldots,X_n)$ modulo
this equivalence relation will be denoted by $(X_0:\ldots:X_n)$.

For example, $\bbP^n(\Fone)$ consists of $2^{n+1}-1$ points
$(X_0:X_1:\ldots:X_n)$ with not all $X_i\in|\Fone|=\{0,1\}$ equal to zero.
In particular, $\bbP^1_{\Fone}/\Spec\Fone$ has three sections
$\bbP^1(\Fone)=\{(0:1),(1:0),(1:1)\}$,
each supported at a different point of $\bbP^1_{\Fone}=
\Proj\Fone[X,Y]=\{0,\infty,\xi\}=\{(X),(Y),(0)\}$ 
(for the unary localization theory; cf.~\ptref{sp:un.proj.space}).
Notice that one of these sections is supported at the generic point~$\xi$!
Of course, this section is a ``closed'' immersion in the sense 
of~\ptref{sp:closed.imm.sch}, but not closed as a map of topological spaces.

\nxsubpoint\label{sp:def.ample.proj}
(Projective morphisms and ample line bundles.)
One can transfer more results from EGA~II to our case. For example,
one can start from a line bundle~$\sL$ over a generalized scheme~$X$
(better supposed to be quasicompact and quasicoherent),
construct graded generalized ring $\Gamma_*(\sO_X)$ by
the procedure of~\ptref{sp:constr.gamma*}, and a generalized
scheme morphism $r_\sL:G(\sL)\to\Proj\Gamma_*(\sO_X)$,
defined on an open subscheme $G(\sL)\subset X$ essentially in the
same way as in EGA~II, 3.7. It is natural to say that
$\sL$ is {\em ample\/} if $G(\sL)=X$ and $r_\sL$ is an
immersion in the sense of~\ptref{sp:imm.subsch}.

We say that $X$ is {\em projective\/} if it admits an ample $\sL$,
such that $r_\sL$ is an isomorphism of generalized schemes,
and that $X$ is {\em quasi-projective} if it admits an ample~$\sL$,
such that $r_\sL$ is an open embedding. We also need to assume~$X$
to be finitely presented over some affine base $S=\Spec C$ 
(e.g.\ $\Spec\Fempty$) in these definitions.

These definitions have relative versions, at least
over a quasicompact quasiseparated base~$S$. One can transfer
some of the classical results to our situation, e.g.\
$\sO[1]$ is $S$-ample on $\bbP_S(\sE)$.

\nxsubpoint\label{sp:proj.sch.vs.closed.subsch.projvb}
(Projective schemes vs.\ ``closed'' subschemes of $\bbP^N$.)
Notice that if $X=\Proj R$ is projective (say, over some $S=\Spec C$) 
in the above sense, we still cannot expect~$X$ to embed into~$\bbP^N_C$,
or even into an infinite-dimensional projective bundle $\bbP_C(\sE)$,
where $\sE:=\widetilde{|R|_1}$, simply because we cannot expect
the canonical graded $C$-algebra homomorphism $S_C(|R|_1)\to R$
to be surjective, even after replacing $R$ by some $R^{(d)}$ 
(i.e.\ ample line bundle $\sL$ by some power $\sL^{\otimes d}$),
since $R$ can happen not to be a {\em pre-unary\/} graded $C$-algebra. 

We'll distinguish these cases by saying that $X$ is {\em pre-unary\/}
or {\em unary projective\/} generalized scheme over~$S=\Spec C$, if 
$X=\Proj R$ for some pre-unary (resp.\ unary) graded~$C$-algebra~$R$.
For example, $\bbP^N_S$ is unary projective over~$S$, and any 
``closed'' subscheme of $\bbP^N_S$ is pre-unary projective over~$S$.

\nxsubpoint (Further properties.)
Some more classical properties of projective schemes and morphisms
can be transferred to our case. However, we don't seem to have
a reasonable notion of properness, and cohomological properties
can be partially transferred only with the aid of homotopic algebra of
Chapters \ptref{sect:homot.alg} and~\ptref{sect:homot.alg/topoi}, 
so there are remain not so many classical properties not discussed yet.

\cleardoublepage

\mysection{Arakelov geometry}
Now we want to discuss possible applications of generalized schemes
and rings to Arakelov geometry, thus fulfilling some of the
promises given in Chapters \ptref{sect:motivation} 
and~\ptref{sect:zinfty.lat.mod}. Our first aim is to construct
$\CompZ$, the ``compactification'' of $\Spec\bbZ$. After that 
we are going discuss models $\sX/\CompZ$ of algebraic varieties $X/\bbQ$.

We also present some simple applications, such as the classical
formula for the height of a rational point~$P$ on a projective variety~$X/\bbQ$
as the arithmetic degree $\deg\sigma_P^*\sO_\sX(1)$ of the pullback of
the Serre line bundle $\sO_\sX(1)$ on a projective model $\sX/\CompZ$
of~$X$ with respect to the section $\sigma_P:\CompZ\to\sX$ inducing
given point $P$ on the generic fiber.

\nxpointtoc{Construction of $\CompZ$}\label{p:constr.compz}
\nxsubpoint
(Notation: generalized rings $A_N$ and $B_N$.)
In order to construct $\CompZ$ we'll need two generalized subrings
$A_N\subset B_N\subset\bbQ$, defined for any integer $N>1$.
Namely, we put $B_N:=\bbZ[1/N]=\bbZ[N^{-1}]$ (so $B_N$ is actually 
a classical ring), and $A_N:=B_N\cap\Zninfty$. Therefore,
\begin{align}
A_N(n)=&\bigl\{(\lambda_1,\ldots,\lambda_n)\in\bbQ^n\,:\,
\lambda_i\in B_N, \sum_i|\lambda_i|\leq 1\bigr\},\quad\text{i.e.}\\
A_N(n)=&\Bigl\{\Bigl(\frac{u_1}{N^k},\ldots,\frac{u_n}{N^k}\Bigr)\,:\,
k\geq0, u_i\in\bbZ, \sum_i|u_i|\leq N^k\Bigr\}
\end{align}
In particular, $|A_N|=A_N(1)=\{\lambda\in B_N\,:\,|\lambda|\leq1\}$,
e.g.\ $1/N\in|A_N|$.

Another useful notation: we write $M\mid N^\infty$ if $M$ divides $N^k$ for
some $k\geq0$, i.e.\ the prime decomposition of~$M$ consists 
only of prime divisors of~$N$. Of course, we write $M\mid N$ when
$M$ is a divisor of~$N$, and $M\nmid N$ otherwise.

\nxsubpoint\label{sp:BN.as.localiz.AN}
($B_N$ is a localization of $A_N$.)
Image of $f:=1/N\in|A_N|$ under canonical embedding
$A_N\to B_N$ is obviously invertible in $|B_N|$,
hence we get a generalized ring homomorphism $A_N[f^{-1}]=A_N[(1/N)^{-1}]\to
B_N$.

\noindent
{\bf Claim.} {\em Map $A_N[(1/N)^{-1}]\to B_N$ is an isomorphism,
i.e.\ $A_N[(1/N)^{-1}]=B_N$ inside~$\bbQ$.}\\
{\bf Proof.} Since $A_N\to B_N$ is a monomorphism, i.e.\
all $A_N(n)\to B_N(n)$ are injective, the same is true for
induced map of localizations $A_N[f^{-1}]\to B_N[f^{-1}]=B_N$
(this is a general property of localizations). Since the composite
map $A_N\to A_N[f^{-1}]\to B_N$ is also a monomorphism, the same 
holds for the first arrow $A_N\to A_N[f^{-1}]$ as well, so we can
identify all generalized rings involved with their images in~$\bbQ$:
$A_N\subset A_N[f^{-1}]\subset B_N\subset\bbQ$. We have to show that
$A_N[f^{-1}]=B_N$, i.e.\ that any $\lambda=(\lambda_1,\ldots,\lambda_n)
\in B_N(n)=(B_N)^n$ belongs to $A_N[f^{-1}]$. This is clear:
choose large integer $k\geq0$, such that $\sum_i|\lambda_i|\leq N^k$
(we use $N>1$ here). Then $\mu:=f^k\lambda=N^{-k}\lambda$ belongs
to $A_N(n)$, hence $\lambda=\mu/f^k$ belongs to $A_N[f^{-1}](n)$, q.e.d.
(cf.~\ptref{sp:examp.loc} for a similar proof).\\
{\bf Remark.} Notice that the same reasoning shows that
$A_N[(1/d)^{-1}]=B_N$ for any divisor $d>1$ of~$N$.

\nxsubpoint
(Construction of $\CompZ\strut^{(N)}$.)
Since $B_N=A_N[(1/N)^{-1}]=\bbZ[N^{-1}]$, we see that both
$\Spec\bbZ$ and $\Spec A_N$ have principal open subsets isomorphic to
$\Spec B_N$, regardless of the localization theory $\cT^?$ 
chosen to construct spectra involved. Therefore, we can construct
a generalized scheme $\CompZ\strut^{(N)}$ by gluing together $\Spec\bbZ$
and $\Spec A_N$ along their open subschemes isomorphic to $\Spec B_N$.
Obviously $\CompZ\strut^{(N)}$ is quasicompact and quasiseparated. We are
going to describe its points, at least for the minimal localization
theory $\cT^u$ (i.e.\ if $\Spec=\Spec^u=\Spec^p$ are the prime spectra).

\nxsubpoint
(Prime ideals of $A_N$.)
Let us describe all prime ideals of $A_N$. First of all, we have ideal
$\gp_\infty$, preimage of $\gm_{(\infty)}\subset\Zninfty$ under the
canonical embedding $A_N\to\Zninfty$. Clearly, $\gp_\infty=
\{\lambda\in|B_N|\,:\,|\lambda|<1\}$. Its complement $|A_N|-\gp_\infty$
consists of $\{\pm1\}$, the invertible elements of $|A_N|$, hence
$A_N$ is local with maximal ideal $\gp_\infty$. Next, for any prime
$p\nmid N$ we get a prime ideal $\gp_p\subset|A_N|$, the preimage of
$pB_N\subset B_N$ under canonical map $A_N\to B_N$. Finally, we have
the zero ideal $0=\{0\}\subset|A_N|$.

We claim that this list consists of all prime ideals of $A_N$, i.e.\
$\Spec^p A_N=\{0,\gp_p,\ldots,\gp_\infty\}_{p\nmid N}$, where $p$ runs over all
primes $p\nmid N$. Indeed, since $B_N=A_N[(1/N)^{-1}]$,
the set $D(1/N)$ of prime ideals of $A_N$ not containing $1/N$
is in one-to-one correspondence to $\Spec B_N=\{0,(p),\ldots\}$,
for the same values of~$p$. Therefore, all we have to check is that
$\gp_\infty$ is the only prime ideal of~$A_N$, which contains~$1/N$.

So let $\gp\subset|A_N|$ be any prime ideal, such that $1/N\in\gp$.
Clearly, $\gp\subset\gp_\infty$. Conversely, suppose that 
$\lambda\in\gp_\infty$, i.e.\ $\lambda\in B_N$ and $|\lambda|<1$.
Then $|\lambda|^k<1/N$ for a sufficiently large integer $k>0$.
Now put $\mu:=N\lambda^k$. Then $\mu\in B_N$ and $|\mu|\leq 1$,
hence $\mu\in|A_N|$, hence $\lambda^k=\mu\cdot(1/N)\in|A_N|\cdot\gp=\gp$,
$\gp$ being an ideal in~$|A_N|$. Since $\gp$ is prime, $\lambda^k\in\gp$
implies $\lambda\in\gp$, so we have proved opposite inclusion 
$\gp_\infty\subset\gp$, q.e.d.

\nxsubpoint
(Prime spectrum of $A_N$.)
We have just shown that $\Spec^uA_N=\Spec^pA_N=\{\xi,\infty,p,\ldots\}$,
where $p$ runs over the set of primes $p\nmid N$. Of course, in this
notation $\xi$ corresponds to the zero ideal $0$, $p$ corresponds to
$\gp_p$, and $\infty$ to $\gp_\infty$. The topology on this set
is described as follows. Point $\xi$ is its generic point,
and point $\infty$ is its only closed point, $A_N$ being local.
The complement $\Spec^pA_N-\{\infty\}$ of this closed point is
homeomorphic to $\Spec B_N=\Spec\bbZ[1/N]$, hence each point $p$ is
closed in this complement, hence $\overline{\{p\}}=\{p,\infty\}$
in $\Spec^pA_N$. In other words, a non-empty subset $U\subset\Spec^pA_N$
is open iff its complement is either empty or 
consists of $\infty$ and of finitely many points $p$. 
In this respect $\Spec^pA_N$ looks like the spectrum of a
(classical) two-dimensional regular local ring. 

\nxsubpoint\label{sp:top.compuz.N}
(Topology of $\CompuZ^{(N)}$.)
Once we have a complete description of topological spaces
$\Spec^u\bbZ=\{\xi,p,\ldots\}_{p\in\bbP}$, $\Spec^u A_N=\{\xi,p,\ldots,\infty\}
_{p\nmid N}$ and $\Spec^u B_N=\{\xi,p,\ldots\}_{p\nmid N}$, we can
describe $\CompuZ^{(N)}$ as well. We see that $\CompuZ^{(N)}=\{\xi,p,\ldots,
\infty\}_{p\in\bbP}=\Spec\bbZ\cup\{\infty\}$, i.e.\ $\CompuZ^{(N)}$ has
an additional point $\infty$, corresponding to the archimedian valuation
of~$\bbQ$. The topology of $\CompuZ^{(N)}$ can be described as follows.
Point $\xi$ is generic, so it is contained in any non-empty open subset
of~$\CompuZ^{(N)}$. Points $\infty$ and $p$ for $p\mid N$ are closed,
and the remaining prime points $p\nmid N$ are not closed since
$\overline{\{p\}}=\{p,\infty\}$. A non-empty subset $U\subset\CompuZ^{(N)}$
is open iff it contains $\xi$, its complement is finite, and
if it either doesn't contain $\infty$ or contains all $p\nmid N$.

\nxsubpoint
(Local rings of $S_N:=\CompuZ^{(N)}$.)
Let's compute the stalks of $\sO=\sO_{S_N}$, where
$S_N:=\CompuZ^{(N)}$.  Since points $\xi$ and $p$ lie in open
subscheme $\Spec\bbZ\subset S_N$, we have $\sO_{S_N,\xi}=\bbQ$ and
$\sO_{S_N,p}=\bbZ_{(p)}$. On the other hand, $\Spec A_N$ is an open
neighborhood of $\infty$, hence
$\sO_{S_N,\infty}=A_{N,\gp_\infty}=A_N$, because $\gp_\infty$ is the
maximal ideal of local ring~$A_N$. In particular,
$\sO_{S_N,\infty}\neq\Zninfty$, contrary to what one might have
expected.

\nxsubpoint\label{sp:trans.maps.compz}
(Morphisms $f=f^{NM}_N:\CompZ\strut^{(NM)}\to\CompZ\strut^{(N)}$.)
Let $N$, $M>1$ be two integers. Denote by $U_1$ and $U_2$ the two open
subschemes of $\CompZ\strut^{(N)}$ isomorphic to $\Spec\bbZ$ and
$\Spec A_N$, respectively. Similarly, denote by $U'_1$ and $U'_2$
the open subschemes of $\CompZ\strut^{(NM)}$ isomorphic to
$\Spec\bbZ$ and $\Spec A_{NM}$, respectively. By construction
$\CompZ\strut^{(N)}=U_1\cup U_2$ and $\CompZ\strut^{(NM)}=U'_1\cup U'_2$.

Next, consider principal open subscheme $W=\Spec B_N=\Spec\bbZ[N^{-1}]$
of $U_1'=\Spec\bbZ$, and put $V_1:=U'_1$, $V_2:=U'_2\cup W$. We are going
to construct generalized scheme morphisms $f_i:V_i\to U_i$, $i:=1$, $2$,
such that $f_1|_{V_1\cap V_2}=f_2|_{V_1\cap V_2}$, and
$f_i^{-1}(U_1\cap U_2)=V_1\cap V_2$. Since
$\CompZ\strut^{(NM)}=V_1\cup V_2$, these two $f_i$ will define together
a generalized scheme morphism $f:\CompZ\strut^{(NM)}\to\CompZ\strut^{(N)}$,
such that $f^{-1}(U_i)=V_i$, $f_{U_i}=f_i:V_i\to U_i$, $i=1$, $2$.

(a) Since both $U_1=\Spec\bbZ$ and $V_1=U'_1=\Spec\bbZ$, we can take
$f_1:=\id_{\Spec\bbZ}$.

(b) Since $U_2=\Spec A_N$, generalized scheme morphisms
$f_2:V_2\to U_2=\Spec A_N$ are in one-to-one correspondence with
generalized ring homomorphisms $\phi:A_N\to\Gamma(V_2,\sO)$.
Since $V_2=U'_2\cup W$, and $U'_2\cap W=U'_2\cap U'_1\cap W=D(NM)\cap D(N)=
D(NM)$, we see that $\Gamma(V_2,\sO)$ equals the fibered product
of $\Gamma(U'_2,\sO)=A_{NM}$ and $\Gamma(W,\sO)=B_N$ over 
$\Gamma(D(NM),\sO_{\Spec\bbZ})=B_{NM}$, hence
$\Gamma(V_2,\sO)=A_{NM}\times_{B_{NM}}B_N=A_{NM}\times_\bbQ B_N=
A_{NM}\cap B_N=\Zninfty\cap B_{NM}\cap B_N=\Zninfty\cap B_N=A_N$,
so we can put $\phi:=\id_{A_N}:A_N\simto\Gamma(V_2,\sO)$.

(c) Notice that $V_1\cap V_2=U_1'\cap(U_2'\cup W)=(U_1'\cap U_2')\cup
(U_1'\cap W)=D_{U_1'}(NM)\cup W=W$, and $f_1^{-1}(U_1\cap U_2)=
f_1^{-1}(D_{U_1}(N))=D_{U_1'}(N)=W=V_1\cap V_2$. Furthermore, the map
$(f_1)_{U_1\cap U_2}:W\to U_1\cap U_2$ obtained from $f_1$ by base change
$U_1\cap U_2\to U_1$ is obviously the identity map of
$W=D_{U_1'}(N)=\Spec B_N=U_1\cap U_2$.

(d) Now compute $f_2^{-1}(U_1\cap U_2)=f_2^{-1}(D_{U_2}(1/N))=
(V_2)_{1/N}=D_{U_2'}(1/N)\cup D_W(1/N)=\Spec A_{NM}[(1/N)^{-1}]\cup W=
\Spec B_{NM}\cup\Spec B_N=\Spec B_N=W=V_1\cap V_2$. Here we have used
\ptref{sp:BN.as.localiz.AN}, which claims that $A_{NM}[(1/N)^{-1}]=B_{NM}$.

(e) It remains to check that the map $(f_2)_{U_1\cap U_2}:W\to U_1\cap U_2$
equals $(f_1)_{U_1\cap U_2}=\id_W$. Notice for this that
according to the construction of~$f_2$ given in~(b), the restriction
$f_2|_W:W=\Spec B_N\to U_2=\Spec A_N$ is induced by canonical embedding
$A_N\to B_N=A_N[(1/N)^{-1}]$, hence it is (up to a canonical isomorphism) 
the open immersion of $W=D_{\Spec\bbZ}(N)$ into $U_2=\Spec A_N$. This
implies $(f_2)_{U_1\cap U_2}=\id_W$, thus concluding the
proof of existence of~$f:\CompZ\strut^{(NM)}\to\CompZ\strut^{(N)}$ with
required properties.

(f) Notice that $f$ is quasicompact and quasiseparated, just because
both $f_1$ and $f_2$ are.

The reasoning used above to construct $f$ is valid regardless of the
localization theory $\cT^?$ chosen.

\nxsubpoint\label{sp:compz.trans.not.iso} (Map $f$ is not an isomorphism.)
Clearly, if $M\mid N^\infty$, numbers $NM$ and $N$ have same prime divisors,
hence $A_{NM}=A_N$, $B_{NM}=B_N$, $\CompZ\strut^{(NM)}=\CompZ\strut^{(N)}$,
and $f=f^{NM}_M:\CompZ\strut^{(NM)}\to\CompZ\strut^{(N)}$ is the identity map.

We claim that {\em if $M\nmid N^\infty$, then $f=f^{NM}_N$ is not an 
isomorphism.} If we use unary localization theory $\cT^u$, we can
just observe that $f$ is a continuous bijective map of underlying topological
spaces, but not a homeomorphism since any prime $p\mid M$, $p\nmid N$
defines a point, which is closed in $\CompuZ^{(NM)}$, but not in
$\CompuZ^{(N)}$. However, $f$ cannot be an isomorphism even if we 
choose another localization theory to construct our generalized schemes,
for example because $\Pic(\CompZ\strut^{(N)})$ is a free abelian group
of rank equal to the number of distinct prime divisors of~$N$,
regardless of the choice of~$\cT^?$ (cf.~\ptref{sp:pic.compz.N}),
and $NM$ has more prime divisors than $N$.

\nxsubpoint (Map $f$ is not affine.)
Furthermore, {\em if $M\nmid N^\infty$, then $f$ is not affine.}
Indeed, suppose that $f=f^{NM}_N$ is affine. Then $f_2=f_{U_2}:V_2\to 
U_2=\Spec A_N$ is affine, hence $V_2$ must be also affine. But
induced map $A_N\to \Gamma(V_2,\sO)=A_N$ is the identity map 
by~\ptref{sp:trans.maps.compz},(b), hence $f_2$ is an isomorphism.
Now $f_1=f_{U_1}:V_1\to U_1$ is always an isomorphism, so we conclude
that $f$ is an isomorphism, which is absurd by~\ptref{sp:compz.trans.not.iso}.

\nxsubpoint (Projective system $\CompZ\strut^{(\cdot)}$.)
The collection of all $\{\CompZ\strut^{(N)}\}_{N>1}$,
together with transition maps $f^N_M$, defined by~\ptref{sp:trans.maps.compz} 
whenever $M$ divides~$N$ (when $N=M$, we put $f_N^N:=\id$),
clearly constitutes a filtering projective system of generalized schemes over
the set of integers $N>1$ ordered by divisibility relation.

\nxsubpoint\label{sp:constr.compz} (Definition of $\CompZ$.)
Now we are tempted to define the ``true'' compactification of $\Spec\bbZ$
by putting $\CompZ:=\projlim_{N>1}\CompZ\strut^{(N)}$. 
However, we need to discuss the meaning of such a projective limit
before doing this. We suggest two possible approaches:

\nxsubpoint\label{sp:compZ.pro.sch} ($\CompZ$ as pro-generalized scheme.)
We can consider $\CompZ:=\quotprojlim_{N>1}\CompZ\strut^{(N)}$
as a pro-generalized scheme, i.e.\ a pro-object in the category of generalized
schemes (cf.\ SGA~4~I for a discussion of pro-objects). This means
that we put formally $\Hom(T,\CompZ):=\projlim_{N>1}\Hom(T,S_N)$ and
$\Hom(\CompZ,T):=\injlim_{N>1}\Hom(S_N,T)$, where $S_N:=\CompZ\strut^{(N)}$,
and $T$ is any generalized scheme. 

Furthermore, we can extend this to any stack~$\cC$ defined over the
category of generalized schemes by putting
$\cC(\CompZ):=\injLim_{N>1}\cC(S_N)$. Informally speaking, an object
of $\cC(\CompZ)$ is just a couple $(N,X)$, where $N>1$ and
$X\in\Ob\cC(S_N)$, and such a couple $(N,X)$ is identified with all
its pullbacks $(NM,(f^{NM}_N)^*X)$. In other words, we expect all
objects of $\cC(\CompZ)$ to come from a finite stage of the projective
limit.

However, this seems to be a reasonable approach only for stacks~$\cC$,
which already have similar property with respect to filtered projective
limits of affine generalized schemes, i.e.\ such that
$\cC(\Spec\injlim A_\alpha)\cong\injLim\cC(\Spec A_\alpha)$. For example,
$\cC=\stQCOH$ doesn't have this property, but stacks of finitely presented
objects (e.g.\ finitely presented sheaves of modules, finitely presented
schemes, or finitely presented sheaves over finitely presented schemes)
usually have this property, essentially by the classical reasoning
of EGA~IV~8. Therefore, we can say that a finitely presented
$\sO_{\CompZ}$-module $\sF$ is just a couple $(N,\sF_N)$, consisting
of an integer $N>1$ and a finitely presented $\sO_{S_N}$-module $\sF_N$.
In such a situation we'll write $\sF=f_N^*(\sF_N)$, where $f_N:\CompZ\to
S_N$ is the canonical projection (in the category of pro-generalized schemes).
Similarly, a finitely presented scheme $\sX$ over $\CompZ$,
and a finitely presented $\sO_\sX$-module $\sF$ can be described as 
``formal pullbacks'' of corresponding objects $X_N$, $\sF_N$,
defined over some finite stage $S_N=\CompZ\strut^{(N)}$ of 
the projective limit.

\nxsubpoint\label{sp:compZ.gen.rgsp} ($\CompZ$ as generalized ringed space.)
Another appoach is to compute $\CompZ:=\projlim_{N>1}S_N$, where
$S_N=\CompZ\strut^{(N)}$, in the category of generalized ringed
spaces.  This means that we first construct topological space $\CompZ$
as the projective limit of corresponding topological spaces $S_N$ (one
can also work on the level of corresponding sites or topoi, cf.\
SGA~4, but topological spaces are sufficient for our purpose), and
then define the structural sheaf by
$\sO_{\CompZ}:=\injlim_{N>1}f_N^{-1}\sO_{S_N}$, where $f_N:\CompZ\to
S_N$ are the natural projection maps, and $f_N^{-1}$ denotes the
``set-theoretic'' pullback of sheaves of generalized rings.

This approach has its advantages. We obtain a reasonable topological space,
and even a generalized ringed space, so we can study for example
line bundles, vector bundles or finitely presented sheaves of modules
over~$\CompZ$. However, this generalized ringed space $\CompZ$ is not a
generalized scheme, so in order to study schemes over $\CompZ$ in this
approach we would need to develop the notion of a ``relative generalized 
scheme''. This is possible, but complicates everything considerably,
so we prefer not to do this, and adopt the pro-generalized scheme approach
instead, whenever we need to study schemes over~$\CompZ$.

\nxsubpoint (Structure of generalized ringed space $\CompuZ$.)
Let us discuss the structure of~$\CompuZ$, considered as a 
generalized ringed space. First of all, each $S_N=\CompuZ^{(N)}=
\Spec\bbZ\cup\{\infty\}$ as a set, and all $f^{NM}_N:S_{NM}\to S_N$
induce identity maps on underlying sets, hence the underlying set of
$\CompuZ$ is also equal to $\Spec\bbZ\cup\{\infty\}$. Furthermore,
$\CompuZ$ contains $\Spec\bbZ$ as an open (generalized ringed) subspace,
simply because this is true for each $S_N=\CompuZ^{(N)}$. This determines
all neighborhoods in $\CompuZ$ of points $\neq\infty$; as to the neighborhoods
of $\infty$, by definition of projective limit of topological spaces,
a set $\infty\in U\subset\CompuZ$ is a neighborhood of~$\infty$ iff
it is one inside some $\CompuZ^{(N)}$. Applying our previous results
of~\ptref{sp:top.compuz.N}, we see that {\em $\xi$ is the generic point
of~$\CompuZ$, and all other points of $\CompuZ$ are closed. Non-empty
open subsets $U\subset\CompuZ$ are exactly the subsets containing $\xi$
and having a finite complement.} In this respect $\CompuZ$ is much more
similar to an algebraic curve, where all points except the generic one
are closed, than $\CompuZ^{(N)}$, where all primes $p\nmid N$
are ``entangled'' with~$\infty$. Furthermore, the stalk of $\sO_{\CompuZ}$
at $\infty$ equals $\injlim_{N>1}\sO_{\CompuZ\phantom{|}^{(N)},\infty}=
\injlim_{N>1}A_N=\Zninfty$, and the complement of $\infty$ is isomorphic to 
$\Spec\bbZ$ as a generalized ringed space.

In this way $\CompuZ$ looks very much like one's idea of
the ``correct'' compactification of $\Spec\bbZ$, apart from (not)
being a generalized scheme.

\nxsubpoint\label{sp:finpres.compuz.mod} 
(Finitely presented sheaves of $\sO_{\CompuZ}$-modules.)
Let $\sF$ be a finitely presented sheaf over~$\CompuZ$. Since
the complement of $\infty$ is isomorphic to $\Spec\bbZ$, and
$\sF|_{\Spec\bbZ}=\tilde M$ for some finitely generated $\bbZ$-module~$M$,
only the structure of~$\sF$ near~$\infty$ is interesting. 
By definition of finite presentation we can find an open neighborhood
$U$ of~$\infty$ and a right exact diagram
\begin{equation}
\xymatrix{
\sO_U(m)\ar@/_2pt/[r]\ar@/^2pt/[r]^{u,v}&\sO_U(n)\ar[r]^{p}&\sF|_U}
\end{equation}
Here $m$, $n\geq0$, and $u$, $v\in\Hom_{\sO_U}(\sO_U(m),\sO_U(n))=
\sO_U(n)^m$ are two $m\times n$-matrices over~$\sO_U$. By definition
of projective limit of topological spaces we may assume,
making~$U$ smaller if necessary, that $U=(f_N)^{-1}(U_N)$ for some
open neighborhood~$U_N$ of $\infty$ in $S_N:=\CompuZ^{(N)}$. Next,
$u$, $v\in\sO_U(n)^m=\injlim_M f_{MN}^{-1}\sO_{(f_N^{MN})^{-1}(U_N)}(n)^m$,
so replacing $N$ by a suitable multiple $NM$, we may assume that
$u$ and $v$ come from some $u_N$, $v_N\in\sO_{U_N}(n)^m$.
Consider finitely presented $\sO_{U_N}$-module $\sF'_N:=\Coker
(u_N,v_N:\sO_{U_N}(m)\rightrightarrows\sO_{U_N}(n))$. By construction
$f_N^*\sF'_N\cong\sF|_U$, and $\sF'_N$ agrees with $\sF|_{\Spec\bbZ}$
on $U_N\cap\Spec\bbZ$ since the restriction of $f_N$ to $U\cap\Spec\bbZ$
is an isomorphism from $U\cap\Spec\bbZ$ to $U_N\cap\Spec\bbZ$. Therefore,
$\sF'_N$ and $\sF|_{\Spec\bbZ}$ patch together into a finitely presented
$\sO_{S_N}$-module~$\sF_N$, such that $\sF=f_N^*(\sF_N)$.

We have just shown that {\em any finitely presented sheaf $\sF$ 
over~$\CompZ$ comes from a finitely presented sheaf $\sF_N$ defined
over a finite stage $S_N$ of the projective limit $\CompZ=\projlim_{N>1}S_N$.}
Therefore, the category of finitely presented sheaves over $\CompuZ$
coincides with that constructed in~\ptref{sp:compZ.pro.sch} 
via the pro-object approach.

\nxsubpoint\label{sp:equiv.2.appr} 
(Equivalence of the two approaches to $\CompZ$.)
We see that the pro-generalized scheme approach of~\ptref{sp:compZ.pro.sch}
and the generalized ringed space approach of~\ptref{sp:compZ.gen.rgsp}
yield equivalent results when we study {\em finitely presented\/}
sheaves of modules over $\CompZ$. This seems to be equally true
for all ``categories of finitely presented objects'' 
(sheaves of modules, relative generalized schemes,\dots),
so we can safely adopt any of these approaches (or use both at the same time)
whenever we consider only finitely presented objects.

\nxsubpoint (Quasicoherent sheaves over $\CompZ\strut^{(N)}$.)
Consider $\stQCOH(S_N)$, the category of quasicoherent $\sO_{S_N}$-modules,
where $S_N:=\CompZ\strut^{(N)}$ for some localization theory $\cT^?$.
Since $S_N=U_1\cup U_2=\Spec\bbZ\cup\Spec A_N$ is an affine open cover
of $S_N$ with $U_1\cap U_2=\Spec B_N$, we see that a quasicoherent
$\sO_{S_N}$-module~$\sF$ is completely determined by $\bbZ$-module
$M_\bbZ:=\Gamma(U_1,\sF)$, $A_N$-module $M_N:=\Gamma(U_2,\sF)$,
and $B_N$-module isomorphism $\theta_M:M_\bbZ\otimes_\bbZ B_N\simto
M_N\otimes_{A_N}B_N$, obtained by restricting to $U_1\cap U_2$
isomorphisms $\widetilde{M_\bbZ}\simto\sF|_{U_1}$ and $\sF|_{U_2}\simto
\widetilde{M_N}$, composing them, and taking global sections.

Furthermore, we obtain in this way an {\em equivalence\/} of $\stQCOH(S_N)$
with the category $\bar\cC_N$ of triples $(M_\bbZ,M_N,\theta_M)$,
consisting of a $\bbZ$-module~$M$, an $A_N$-module~$M_N$, and a
$B_N$-module isomorphism $\theta_N:(M_\bbZ)_{(B_N)}=M_\bbZ[N^{-1}]\simto
(M_N)_{(B_N)}=M_N[(1/N)^{-1}]$. In particular, 
$\stQCOH(\CompZ\strut^{(N)})$ {\em doesn't depend on the choice
of localization theory~$\cT^?$.} 

Pullback of quasicoherent sheaves with respect to $f_N^{NN'}:S_{NN'}\to S_N$
induces a functor $\bar\cC_N\to\bar\cC_{NN'}$, easily seen to transform
$(M_\bbZ,M_N,\theta)$ into $(M_\bbZ,(M_N)_{(A_{NN'})},
\theta_{(B_{NN'})})$.

\nxsubpoint\label{sp:finpres.compZ.N} 
(Finitely presented sheaves over $\CompZ\strut^{(N)}$.)
The essential image of the full subcategory of $\stQCOH(S_N)$ consisting
of finitely presented $\sO_{S_N}$-modules, where again 
$S_N=\CompZ\strut^{(N)}$, under the equivalence of categories
$\stQCOH(S_N)\to\bar\cC_N$, is obviously equal to the full subcategory
$\cC_N$ of~$\bar\cC_N$, consisting of triples $(M_\bbZ,M_N,\theta_N)$
as above, with $M_\bbZ$ a finitely presented $\bbZ$-module,
and $M_N$ a finitely presented $A_N$-module. In particular, 
this category doesn't depend on the choice of~$\cT^?$,
hence the same is true for the category of finitely presented sheaves
of modules over $\CompZ=\projlim_{N>1}S_N$ (for any understanding of 
this projective limit), this category being equivalent to 
$\injLim_{N>1}\cC_N$.

\nxsubpoint (``Point at $\infty$'' of $\CompZ$.)
Notice that for each $N>1$ we have a canonical embedding $A_N\to\Zninfty$,
which induces a morphism $\Spec\Zninfty\to \Spec A_N=U_2\subset
\CompZ\strut^{(N)}$. Let us denote by~$\eta_N$ the natural morphism
$\Spec\Zninfty\to S_N:=\CompZ\strut^{(N)}$. It is easy to see that
$\eta_N=f_N^{NM}\circ\eta_{NM}$ for any $M\geq1$, i.e.\
$\hat\eta=(\eta_N)_{N>1}$ is a morphism from $\Spec\Zninfty$ into the
projective system $(S_N)_{N>1}$. Therefore, it defines a morphism
$\hat\eta:\Spec\Zninfty\to\CompZ$, both 
in the category of pro-generalized schemes, and in the category
of generalized ringed spaces. If we adopt the generalized ringed space
approach, and use $\cT^?=\cT^u$, then $\hat\eta:\Spec\Zninfty\to\CompuZ$
is just the map with image $\{\xi,\infty\}$, inducing identity on
the local ring $\sO_{\CompuZ,\infty}=\Zninfty=\sO_{\Spec\Zninfty,\infty}$.

\nxsubpoint\label{sp:stalk.at.infty}
(Stalk at $\infty$ of a finitely presented $\sO_{\CompZ}$-module.)
Let $\sF$ be a finitely presented $\sO_{\CompZ}$-module. Then
$\hat\eta^*(\sF)$ is a finitely presented quasicoherent
$\sO_{\Spec\Zninfty}$-module, hence
$\sF_\infty:=\Gamma(\Spec\Zninfty,\hat\eta^*(\sF))$ is a finitely
presented $\Zninfty$-module. We will say that {\em $\sF_\infty$ is the
stalk of~$\sF$ at infinity}, regardless of the choice of $\cT^?$ and
of the understanding of $\CompZ$.

In the pro-generalized scheme approach, $\hat\eta^*(\sF)$
actually means $\eta_N^*(\sF_N)$, if $\sF$ is given by couple $(N,\sF_N)$.
On the other hand, in the generalized ringed space approach $\sF_\infty$
is indeed the stalk of~$\sF$ at~$\infty$, provided we use $\cT^?=\cT^u$.

\nxsubpoint\label{sp:cat.finpres.compzmod}
(Category of finitely presented $\sO_{\CompZ}$-modules.)
Let $\sF$ be a finitely presented $\sO_{\CompZ}$-module. It defines
a finitely generated $\bbZ$-module $M_\bbZ:=
\Gamma(\Spec\bbZ,\sF)$, a finitely presented 
$\Zninfty$-module $M_\infty:=\sF_\infty=\hat\eta^*\sF$,
together with canonical isomorphism
of $\bbQ$-vector spaces $\theta_M:M_\bbZ\otimes_\bbZ\bbQ\simto
M_\infty\otimes_{\Zninfty}\bbQ$, constructed for example by computing
the pullback of $\sF$ with respect to the ``generic point''
$\hat\xi:\Spec\bbQ\to\CompZ$ in two different ways, using factorizations
$\hat\xi:\Spec\bbQ\to\Spec\bbZ\to\CompZ$ and $\hat\xi:\Spec\bbQ\to
\Spec\Zninfty\stackrel{\hat\eta}\to\CompZ$. Of course,
$\bbQ$-vector space $\sF_\xi:=\hat\xi^*\sF\cong\sF_\infty\otimes_{\Zninfty}
\bbQ\cong\Gamma(\Spec\bbZ,\sF)\otimes_\bbZ\bbQ$ is called
the {\em generic fiber of~$\sF$.}

In this way we obtain a functor from the category of finitely presented
$\sO_{\CompZ}$-modules into the category $\cC$ of triples
$(M_\bbZ,M_\infty,\theta)$, consisting of a finitely generated
$\bbZ$-module $M_\bbZ$, finitely presented $\Zninfty$-module $M_\infty$,
and a $\bbQ$-vector space isomorphism $\theta:M_\bbZ\otimes_\bbZ\bbQ\simto
M_\infty\otimes_{\Zninfty}\bbQ$. It is easy to check that 
{\em this functor is in fact an equivalence of categories.} Let us show
for example the essential surjectivity of this functor. Since
the category of finitely presented $\sO_{\CompZ}$-modules is equivalent
to $\injLim_{N>1}\cC_N$ (cf.~\ptref{sp:finpres.compZ.N},
\ptref{sp:compZ.pro.sch} and~\ptref{sp:finpres.compuz.mod}), it suffices to
show that any triple $M:=(M_\bbZ,M_\infty,\theta)$ as above comes from
a triple $(M_\bbZ,M_N,\theta_N)\in\Ob\cC_N$ for some $N>1$
(i.e.\ $M_\infty=M_N\otimes_{A_N}\Zninfty$, 
and $\theta=\theta_N\otimes1_\bbQ$).
Indeed, $M_\infty$ is a finitely presented module over
filtering inductive limit $\Zninfty=\injlim_{N>1}A_N$, hence
one can find an $N>1$ and a finitely presented $A_N$-module
$M_N$, such that $M_\infty=M_N\otimes_{A_N}\Zninfty$
(the reasoning here is essentially the classical one of EGA~IV~8,
cf.\ also~\ptref{sp:finpres.compuz.mod}). Next, $\theta$ is an isomorphism
between $M_\bbZ\otimes_\bbZ\bbQ=:M_\bbQ=(M_\bbZ\otimes_\bbZ B_N)
\otimes_{B_N}\bbQ$ and
$M_\infty\otimes_{\Zninfty}\bbQ=(M_N\otimes_{A_N}B_N)\otimes_{B_N}\bbQ$.
Since $\bbQ=\injlim_{N'\geq1}B_{NN'}$, this isomorphism between
finitely presented $\bbQ$-modules comes from an isomorphism
$\theta_{NN'}:(M_\bbZ\otimes_\bbZ B_N)\otimes_{B_N}B_{NN'}\simto
(M_N\otimes_{A_N}B_N)\otimes_{B_N}B_{NN'}$ for some $N'\geq1$.
Then $(M_\bbZ,M_N\otimes_{A_N}A_{NN'},\theta_{NN'})\in\Ob\cC_{NN'}$
is a triple from $\Ob\cC_{NN'}$ inducing original triple $M\in\Ob\cC$ after
base change, q.e.d.

\nxsubpoint\label{sp:finpres.sch.compz}
(Finitely presented schemes over $\CompZ$.)
This argument can be extended to other categories of (relative)
finitely presented objects. For example, the category of
(relative) finitely presented schemes $\bar\sX/\CompZ$ is equivalent
to the category of triples $(\sX,\sX_\infty,\theta)$, where
$\sX$ is a finitely presented scheme over $\Spec\bbZ$,
$\sX_\infty$ is a finitely presented scheme over $\Spec\Zninfty$,
and $\theta:\sX_{(\bbQ)}\simto\sX_{\infty,(\bbQ)}$ is an isomorphism
of $\bbQ$-schemes. Indeed, any such triple comes from a triple
$(\sX,\sX_N,\theta_N)$, defined over a finite stage~$S_N$ of the projective
limit $\CompZ=\projlim S_N$, and then $\sX\to\Spec\bbZ=U_1\subset S_N$
and $\sX_N\to\Spec A_N=U_2\subset S_N$ can be patched together into a 
finitely presented generalized scheme $\bar\sX_N\to\CompZ\strut^{(N)}=S_N$.

\nxsubpoint\label{sp:finpres.models.compz}
(Implication for finitely presented models.)
An immediate consequence is that {\em constructing a 
finitely presented $\CompZ$-model $\bar\sX$
of a finitely presented $\bbQ$-scheme~$X$ is the same thing as
constructing a finitely presented $\bbZ$-model $\sX$ and
a finitely presented $\Zninfty$-model $\sX_\infty$ of given~$X$.}
Since the first half (existence and different properties of models 
over $\Spec\bbZ$) is quite well known, we'll concentrate 
our efforts on the second half (models over $\Zninfty$).
This agrees with our considerations of~\ptref{p:firstcompz},
arising from classical Arakelov geometry, provided we want to study
only {\em finitely presented\/} $\Zninfty$-models.

\nxsubpoint
($\CompZ\strut^{(N)}$ is finitely presented.)
Before going on, let us show that all $\CompZ\strut^{(N)}$ are
finitely presented (absolutely, i.e.\ over $\Fempty$,
or equivalently, over $\Fone=\Fempty[0^{[0]}]$). This will enable
us to speak about finitely presented schemes $\sX_N$ over
$\CompZ\strut^{(N)}$ without mentioning explicitly whether
$\sX_N$ is finitely presented absolutely, i.e.\ over $\Fempty$,
or relatively, i.e.\ over $\CompZ\strut^{(N)}$.

Since $S_N:=\CompZ\strut^{(N)}=U_1\cup U_2=\Spec\bbZ\cup\Spec A_N$,
with $U_1\cap U_2=\Spec B_N$, we see that $S_N$ is quasicompact and
quasiseparated. Next, $\bbZ=\Fpm[+^{[2]}\,|\,
\bu+0=\bu=0+\bu$, $\bu+(-\bu)=0]$ is finitely presented
over $\Fpm=\Fone[-^{[1]}\,|\,-^2=\bu]$, hence also over~$\Fone$,
so our statement is reduced to the following:

\begin{ThD}\label{th:AN.finpres}
Generalized ring $A_N=\Zninfty\cap\bbZ[1/N]$ is finitely presented
over $\Fpm$, hence also over $\Fone$ and~$\Fempty$. More precisely,
it is generated by the set of ``averaging operations''~$s_p\in A_N(p)$
given by
\begin{equation}
s_p(\{1\},\{2\},\ldots,\{p\})=\frac1p\{1\}+\frac1p\{2\}+\cdots+\frac1p\{p\}
\end{equation}
where $p$ runs over the finite set of all prime divisors of~$N>1$,
modulo finite list of relations \eqref{eq:avg.op.idemp}--\eqref{eq:avg.op.canc}
(for each of these operations~$s_p$) listed below.
\end{ThD}
\begin{Proof} (a)
Let us denote by $s_n\in\Zninfty(n)$ the {\em $n$-th averaging operation:}
\begin{equation}
s_n:=\frac1n\{1\}+\frac1n\{2\}+\cdots+\frac1n\{n\}\in\Zninfty(n)\subset
\Zinfty(n)\subset\bbR^n
\end{equation}
We have already seen in \ptref{sp:tens.sqr.zninfty} that these
operations $\{s_n\}_{n>1}$ generate $\Zninfty$, since any operation
$\lambda=(m_1/n)\{1\}+\cdots+(m_k/n)\{1\}\in\Zninfty(k)$, where
$n$, $m_i\in\bbZ$, $\sum_i|m_i|\leq n>0$, can be re-written as $s_n$,
applied to the list containing $\pm\{i\}$ exactly $|m_i|$ times,
and with the remaining $n-\sum|m_i|$ arguments set to zero. Furthermore,
if $\lambda\in A_N(k)$, then we can take $n=N^t$ for some integer $t>0$,
hence $A_N$ is generated by $\{s_{N^t}\}_{t>0}$. Next, the
operation $s_{nm}$ can be expressed in terms of $s_n$ and $s_m$:
\begin{equation}\label{eq:avg.mn}
s_{nm}=s_n(s_m(\{1\},\ldots,\{m\}),s_m(\{m+1\},\ldots,\{2m\}),\ldots)
\end{equation}
Therefore, $A_N$ is generated by one operation $s_N$, or by finite set
of operations $\{s_p\}_{p\mid N}$.

Notice that each $s_n$ satisfies following {\em idempotency,
symmetry\/} and {\em cancellation\/} relations:
\begin{align}
\label{eq:avg.op.idemp}
&s_n(\{1\},\{1\},\ldots,\{1\})=\{1\}\\
\label{eq:avg.op.symm}
&s_n(\{1\},\{2\},\ldots,\{n\})=s_n(\{\sigma(1)\},\ldots,\{\sigma(n)\}),
\quad\forall\sigma\in\gS_n\\
\label{eq:avg.op.canc}
&s_n(\{1\},\ldots,\{n-1\},-\{n-1\})=s_n(\{1\},\ldots,\{n-2\},0,0)
\end{align}
Furthermore, these relations for $s_{nm}$ follow from similar relations
for $s_n$ and $s_m$, once we take into account \eqref{eq:avg.mn} and
the implied commutativity relation between $s_n$ and~$s_m$. This
means that it suffices to impose the above relations only for 
$\{s_p\}_{p\mid N}$.

(b) Let us denote by $A'_N$ the (commutative) finitely presented
$\Fpm$-algebra generated by operations $S_p\in A'_N(p)$ for $p\mid N$,
subject to the above relations. It will be convenient to define
$S_n\in A'_N(n)$ for all $n\mid N^\infty$ by induction in~$n$, using
formula~\eqref{eq:avg.mn} for this, with $S_1:=\bu$ as the base of induction. 
One checks that the RHS
of~\eqref{eq:avg.mn} doesn't change when we interchange $m$ and~$n$,
using commutativity between $S_m$ and $S_n$ together
with~\eqref{eq:avg.op.symm}, so this inductive procedure does
define operations $\{S_n\}_{n\mid N^\infty}$,
satisfying~\eqref{eq:avg.mn} and
\eqref{eq:avg.op.idemp}--\eqref{eq:avg.op.canc}.

Denote by $\phi:A'_N\to A_N$ the generalized $\Fpm$-algebra 
homomorphism transforming
$S_p$ into $s_p$ for all primes $p\mid N$, hence also $S_n$ into $s_n$
for any $n\mid N^\infty$. It is a strict epimorphism (``surjective map'')
since the $\{s_p\}_{p\mid N}$ generate $A_N$ over $\Fpm$. Now we
are going to check that $\phi$ is a monomorphism (``injective''), since
this would imply that $A_N\cong A'_N$ is finitely presented with the
lists of generators and relations given above.
\end{Proof}

Now we need the following lemma:
\begin{LemmaD}
Any element $\lambda\in A'_N(k)$ can be written as
$S_n(u_1,\ldots,u_n)$, where each $u_i$ belongs to $\Fpm(k)\subset
A'_N(k)$, i.e.\ equals $0$ or $\pm\{x_i\}$ with $1\leq x_i\leq k$.
\end{LemmaD}
\begin{Proof}
Denote by $A''_N(k)\subset A'_N(k)$ the
set of all $k$-ary operations of $A'_N$, which can be represented in the 
above form. We have to check that $\|A''_N\|=\bigsqcup_k A''_N(k)$
contains all basis elements $\{x\}\in A'_N(k)$, $1\leq x\leq k$,
an obvious statement since $\{x\}=S_1(\{x\})$, and that $\|A''_N\|$
is stable under all maps $A'_N(\phi)$, $\phi:\stk\to\stk'$,
an obvious condition as well, and finally, that $\|A''_N\|$ is stable
under the action of any of generators $0^{[0]}$, $-^{[1]}$ and
$S_p^{[p]}$, $p\mid N$, of $A'_N$.

\begin{itemize}
\item
As to the zero $0^{[0]}$, it obviously lies in $A''_N(0)$,
as well as in all $A''_N(k)$.
\item
Stability under $-$ is also immediate once we use the (implied)
commutativity relation between $S_n$ and $-$:
\begin{equation}
-S_n(u_1,\ldots,u_n)=S_n(-u_1,\ldots,-u_n),
\end{equation}
combined with equalities $-(-u)=u$ and $-0=0$.
\item
Stability under $S_p$ is a bit more tricky. Let $z_1$, $\ldots$, $z_p
\in A''_N(k)$ be the list of arguments to~$S_p$. We have to show
that $z:=S_p(z_1,\ldots,z_p)$ lies in $A''_N(k)$ as well. By definition of
set $A''_N(k)$ each $z_i$ can be written as $S_{n_i}$ applied to some 
arguments $u_{ij}\in\Fpm(k)$, for some $n_i\mid N^\infty$. Notice that
any $n_i$ can be replaced by any its multiple $mn_i$ with $m\mid N^\infty$: 
all we need is to write the idempotency relations~\eqref{eq:avg.op.idemp}: 
$u_{ij}=S_m(u_{ij},\ldots,u_{ij})$, substitute these expressions into
$z_i=S_{n_i}(u_{i1},\ldots,u_{i,n_i})$, and apply~\eqref{eq:avg.mn}.
Therefore, we may replace $n_i$ by their product and assume all
$n_i$ to be equal to some~$n\mid N^\infty$. But then $z=S_p(z_1,\ldots,z_p)=
S_p(S_n(u_{11},\ldots,u_{1n}),\ldots,S_n(u_{p1},\ldots, u_{pn}))$ 
can be rewritten as $S_{pn}$ applied to the list of all $u_{ij}$, hence
$z$ belongs to $A''_N(k)$, q.e.d.
\end{itemize}
\end{Proof}

\nxsubpoint {\em (End of proof of~\ptref{th:AN.finpres}.)}
Let $z$, $z'\in A'_N(k)$ be two elements with $\phi(z)=\phi(z')$
in $A_N(k)$.
According to the above lemma, we can write $z=S_n(u_1,\ldots,u_n)$
and $z=S_m(u'_1,\ldots,u'_m)$ for some $m$, $n\mid N^\infty$,
$u_i$, $u'_j\in\Fpm(k)$. Furthermore, we can replace $m$ and~$n$
by their product and assume $m=n$, by the same argument as in the
proof of lemma. Next, using \eqref{eq:avg.op.canc} 
and~\eqref{eq:avg.op.symm}, we can assume that $\{x\}$ and $-\{x\}$
do not occur in the list $u_1$, $u_2$, $\ldots$, $u_n$ simultaneously,
simply because if $u_i=-u_j$, we can replace both $u_i$ and $u_j$ by
zero. Now we've obtained a {\em reduced\/} representation 
$z=S_n(u_1,\ldots,u_n)$. Notice that 
$\phi(z)=(m_1/n)\{1\}+\cdots+(m_k/n)\{k\}$, where $|m_x|$ is the multiplicity
of $\pm\{x\}$ in the list $u_1$, \dots, $u_n$, and the sign is chosen 
depending on whether $\{x\}$ or $-\{x\}$ is present in this list.
Since $\phi(z)=\phi(z')$, we see that $m_x=m'_x$ for all $1\leq x\leq k$,
hence the lists of arguments to $S_n$ in our expressions for $z$ and $z'$
coincide (up to a permutation), hence $z=z'$ by~\eqref{eq:avg.op.symm}.
This proves the injectivity of $\phi:A'_N\to A_N$, hence 
also~\ptref{th:AN.finpres}.

\begin{CorD} Morphisms $f^{NM}_N:\CompZ\strut^{(NM)}\to\CompZ\strut^{(N)}$
are finitely presented.
\end{CorD}
\begin{Proof}
This follows from the classical statement: ``$g\circ f$ finitely presented,
$g$ of finite type $\Rightarrow$ $f$ finitely presented'',
valid for generalized schemes as well.
\end{Proof}

\begin{RemD}\label{rem:gen.of.AN}
(a) Notice that the definition of $A_N=\Zninfty\cap\bbZ[1/N]$, used before
only for $N>1$, can be extended to $N=1$, yielding 
$A_1=\Zninfty\cap\bbZ=\Fpm$. Furthermore, \ptref{th:AN.finpres}
holds for $N=1$ as well: $A_1=\Fpm$ is a finitely presented $\Fpm$-algebra,
with an empty list of generators and relations.

(b) Theorem~\ptref{th:AN.finpres} holds for $N=\infty$, i.e.\ for
$\Zninfty=\injlim_{N>1}A_N$ as well. Indeed, this filtered inductive limit
representation immediately implies that $\Zninfty$ is generated over~$\Fpm$
by all averaging operations $\{s_n\}_{n>1}$, or just by all 
$\{s_p\}_{p\in\bbP}$, subject to relations 
\eqref{eq:avg.op.idemp}--\eqref{eq:avg.op.canc}.

(c) Essentially the same reasoning as in the proof of~\ptref{th:AN.finpres}
shows that $A_N$ is generated by {\em one\/} averaging operation~$s_N$,
subject to relations \eqref{eq:avg.op.idemp}--\eqref{eq:avg.op.canc}.

(d) Another interesting observation is that we've got a $p$-ary generator
$s_p$ of $A_p$ over~$\Fpm$. While $A_2$ and $A_3$ can be shown not to
be generated by operations of lower arities, $A_5$ can be
shown to be generated by quaternary operation 
$(2/5)\{1\}+(1/5)\{2\}+(1/5)\{3\}+(1/5)\{4\}$, 
so the beautiful system of generators given 
in~\ptref{th:AN.finpres} is not necessarily the one which minimizes the
arity of operations involved. 

(e) Notice that $\Zninfty$ is already generated over~$\Fpm$ 
by its {\em binary\/} operations $\lambda\{1\}+(1-\lambda)\{1\}$,
$\lambda\in[0,1]\cap\bbQ$, similarly to $\Zinfty$, so the phenomena
mentioned above in (d) are special for $N<\infty$.
\end{RemD}

\nxsubpoint (Special case: $N=2$.)
In particular, we obtain a presentation of~$A_2=\Zninfty\cap\bbZ[1/2]$.
It is generated by binary averaging operation $s_2$, which will be also denoted
by~$*$, when written in infix form. Thus $[*]=s_2=(1/2)\{1\}+(1/2)\{2\}$;
informally, $x*y=(x+y)/2$. We obtain $A_2=\Fpm[*^{[2]}\,|\,x*x=x$, $x*y=y*x$,
$x*(-x)=0]$. In this respect $A_2$ is somewhat similar to 
$\bbZ=\Fpm[+^{[2]}\,|\,x+0=x=0+x$, $x+(-x)=0]$: both are generated over~$\Fpm$
by one binary operation.

\nxsubpoint\label{sp:genrg.braid.repr} 
(Generalized rings and representations of braid groups.)
Let $A$ be a (commutative) generalized ring, $[*]\in A(2)$ be a 
binary operation, such that $\bu*\bu=\bu$. For example, we can take $A=A_2$,
or $A=\bbQ[q]$, $[*]=\{1\}*\{2\}:=q\{1\}+(1-q)\{2\}$, or the universal case
$A=\Fempty[*^{[2]}\,|\,\bu*\bu=\bu]$. Fix any integer $n\geq2$, and denote
by $t_k:A(n)\to A(n)$, where $1\leq k<n$, the endomorphism of 
free $A$-module~$A(n)$ defined by
\begin{equation}
t_k(\{i\})=\begin{cases}
\{i\},&\text{if $i\neq k,k+1$}\\
\{k+1\},&\text{if $i=k$}\\
\{k\}*\{k+1\},&\text{if $i=k+1$}
\end{cases}
\end{equation}
One checks directly that these matrices $t_k\in\End_A(A(n))\cong 
M(n\times n; A)=A(n)^n$ satisfy the braid relation:
\begin{equation}
t_k\circ t_{k+1}\circ t_k=t_{k+1}\circ t_k\circ t_{k+1},
\quad\text{for $1\leq k\leq n-2$.}
\end{equation}
Since $t_k$ and $t_{k+1}$ act non-identically only on basis elements 
$\{k\}$, $\{k+1\}$ and $\{k+2\}$, it suffices to check this relation for
$k=1$, $n=3$. We get $(t_1\circ t_2\circ t_1)(x,y,z)=
(t_1\circ t_2)(y,x*y,z)=t_1(y,z,(x*y)*z)=(z,y*z,(x*y)*z)$ and
$(t_2\circ t_1\circ t_2)(x,y,z)=(t_2\circ t_1)(x,z,y*z)=t_2(z,x*z,y*z)=
(x,y*z,(x*z)*(y*z))$, so we are reduced to showing $(x*y)*z=(x*z)*(y*z)$,
which is immediate from relation $z=z*z$ and the commutativity relation
$(x*y)*(z*w)=(x*z)*(y*w)$.

In this way we (almost) obtain a representation of the {\em braid group\/}
$B_n=\langle t_1,\ldots,t_{n-1}\,|\,t_kt_{k+1}t_k=t_{k+1}t_kt_{k+1}\rangle$
by $n\times n$-matrices over generalized ring~$A$. The only problem is
that our matrices $t_i$ may be not invertible; this can be circumvented
by considering the matrix localization $A[Z^{-1}]$ of $A$ with respect
to matrix~$Z:(x,y)\mapsto(y,x*y)$, cf.~\ptref{sp:loc.matr}.
One easily checks that $A_2[Z^{-1}]=B_2=\bbZ[1/2]$; as to the universal
case $A=\Fempty[*^{[2]}\,|\,\bu*\bu=\bu]$, I~don't know much about
$A[Z^{-1}]$.

Among other things, this example illustrates the complexity of non-unary
algebras: even $A=\Fempty[*^{[2]}\,|\,\bu*\bu=\bu]$ (or rather $A[Z^{-1}]$) 
is complicated enough to contain representations of the braid groups.
It is interesting whether one can generalize this construction to $p>2$,
using the special properties of operation $s_p$ instead of~$s_2$.

\begin{PropD}\label{prop:pic.AN.triv}
$\Pic(\Spec A_N)=0$, i.e.\ any line bundle over $\Spec A_N$ is trivial.
More precisely, if $P$ is a finitely generated projective $A_N$-module,
such that $\dim_\bbQ P_{(\bbQ)}=1$, then $P\cong|A_N|$.
\end{PropD}
\begin{Proof} (a)
Notice that if $\sL$ is a line bundle over $\Spec A_N$, then
$\sL=\tilde P$ for some finitely generated projective $A_N$-module~$P$
by~\ptref{sp:vb.proj}, and the generic fiber $\sL_\xi\cong P_{(\bbQ)}$
is obviously one-dimensional, so the last statement of the proposition
does imply the other ones. So let $P$ be as above. Choose
a surjection $\pi:A_N(n)\twoheadrightarrow P$, and
a section $j:P\to A_N(n)$ of~$\pi$, existing because of projectivity of~$P$.
Let $p:=j\circ\pi:A_N(n)\to A_N(n)$ be the corresponding projector.
Denote by $P'$, $\pi':B_N(n)\twoheadrightarrow P'$,
$j':P'\to B_N(n)$ and $p':B_N(n)\to B_N(n)$ 
the module and homomorphisms, obtained by scalar 
extension to~$B_N$. Now consider the following diagram:
\begin{equation}
\xymatrix@C+1pc{
A_N(n)\ar[r]^{\xi_{A_N(n)}}\ar@/_/[d]_{\pi}&B_N(n)\ar@/_/[d]_{\pi'}\\
P\ar[r]^{\xi_P}\ar@/_/[u]_{j}&P'\ar@/_/[u]_{j'}}
\end{equation}
Here $\xi_M:M\to M_{(B_N)}$ denotes the natural embedding of an
$A_N$-module~$M$ into its scalar extension to~$B_N$.

Notice that $\xi_{A_N(n)}$ is just the map $\phi_n:A_N(n)\to B_N(n)$
coming from the natural embedding $\phi:A_N\to B_N$, hence it is
injective, hence the same is true for its retract $\xi_P$. 

(b) Let us identify $A_N(n)$ with its image under~$\xi_{A_N(n)}$ in~$B_N(n)$,
which in turn can be identified with a subset of $\bbQ^n$. Identify
$P'$ with $j'(P')=p'(B_N(n))$, and $P$ with $j(P)=p(A_N(n))\subset
A_N(n)\subset B_N(n)$.  Then all modules involved are identified with
subsets of $B_N(n)\subset\bbQ^n$, and we have $P=P'\cap
A_N(n)=P'\cap\Zninfty(n)$.

Denote by $\vnorm$ the $L_1$-norm on~$\bbQ^n$, given by formula
$\|(\lambda_1,\ldots,\lambda_n)\|=\sum_i|\lambda_i|$. Then
$A_N(n)=\{\lambda\in B_N^n\,:\,\|\lambda\|\leq1\}$ and
$P=\{\lambda\in P'\,:\,\|\lambda\|\leq1\}$. Let $u_i:=p(\{i\})$, 
$1\leq i\leq n$, be the finite set of generators of $A_N$-module~$P$.
Put $C:=\max_i\|u_i\|$. Clearly $0<C\leq 1$, since $P\neq0$ and 
$P\subset A_N(n)$. Since $p(\lambda)=\sum_i\lambda_iu_i$ for any
$\lambda=(\lambda_1,\ldots,\lambda_n)\in A_N(n)$, we get
$\|p(\lambda)\|\leq C\cdot\sum_i|\lambda_i|=C\cdot\|\lambda\|\leq C$
for any $\lambda\in A_N(n)$, i.e.\ $C=\max_{x\in P}\|x\|$ since
$P=p(A_N(n))$. We may assume that $C=\|u_1\|$; otherwise we would
just renumber the basis elements of $A_N(n)$ and the generators $u_i$ of~$P$.

(c) We claim that $C=1$. Indeed, suppose that $C<1$. 
Then one can find a fraction
$v/N^k$, such that $1<v/N^k<C^{-1}$ (we use $N>1$ here).
Then $\mu:=v/N^k\in|B_N|$, hence $\mu u_1\in P'$, $P'$ being a $B_N$-module,
and by construction $\|\mu u_1\|\leq 1$, i.e.\ $\mu u_1\in A_N(n)$,
hence $\mu u_1\in P'\cap A_N(n)=P$ with $\|\mu u_1\|>\|u_1\|=C$, which
is absurd.

(d) Now we claim that $u_1$ freely generates $A_N$-module $P$, i.e.\ that
any $\lambda\in P$ can be uniquely written as $\lambda=c\cdot u_1$ 
with $c\in|A_N|$.
Uniqueness is clear, since we can always extend scalars to~$\bbQ$,
and there everything follows from $u_1\neq 0$. This argument also shows
existence of such $c\in\bbQ$, since $P_{(\bbQ)}\cong\bbQ$ is a line
in $A(n)_{(\bbQ)}=\bbQ^n$. We have to show that this $c$ actually lies in 
$|A_N|$. But this is clear: since $\|u_1\|=C=1$ and $\lambda=c\cdot u_1$,
we have $c=\pm\|\lambda\|=\pm\sum_i\lambda_i\in B_N$, because 
$\lambda\in P\subset B_N(n)$, i.e.\ all components $\lambda_i\in B_N$,
and $|c|=\|\lambda\|\leq 1$, just because $\lambda\in P\subset A_N(n)$.
We have just shown that $c\in|B_N|\cap|\Zninfty|=|A_N|$, q.e.d.
\end{Proof}

\begin{CorD}\label{cor:pic.zninfty.triv}
$\Pic(\Spec\Zninfty)=0$. More precisely, any 
finitely generated projective $\Zninfty$-module~$P$, such that
$\dim_\bbQ P_{(\bbQ)}=1$, is free of rank one.
\end{CorD}
\begin{Proof}
Since any finitely generated projective module is finitely
presented (being a direct factor $P=\Coker(p,\id)$, $p=p^2\in M(n,\Zninfty)$, 
of some $\Zninfty(n)$), the usual inductive limit argument shows that
$P$ comes from some finitely generated projective $A_N$-module~$P_N$,
for a suitable $N>1$. Then $P_{N,(\bbQ)}\cong P_{(\bbQ)}\cong\bbQ$,
hence $P_N$ is free of rank one by~\ptref{prop:pic.AN.triv}, 
hence this must be also true for $P=P_{N,(\Zninfty)}$.
\end{Proof}

\nxsubpoint\label{sp:pic.compz.N} 
(Line bundles over $\CompZ\strut^{(N)}$.)
Let $\sL$ be a line bundle over~$S_N:=\CompZ\strut^{(N)}=U_1\cup U_2=
\Spec\bbZ\cup\Spec A_N$. Since all line bundles over $\Spec\bbZ$
and $\Spec A_N$ are trivial (cf.~\ptref{prop:pic.AN.triv}),
we can choose trivializations $\phi_1:\sO_{\Spec\bbZ}\simto\sL|_{U_1}$
and $\phi_2:\sO_{\Spec A_N}\simto\sL|_{U_2}$. Comparing these 
two trivializations over $U_1\cap U_2=\Spec B_N$, we obtain an invertible
element $\lambda\in B_N^\times$, such that $\phi_2|_{U_1\cap U_2}=
\lambda\cdot\phi_1|_{U_1\cap U_2}$. This element~$\lambda$ determines
$\sL$ up to isomorphism, and we can always construct a line bundle
$\sL=\sO(\log\lambda)$ starting from any element $\lambda\in B_N^\times$,
simply by taking trivial line bundles over $U_1$ and $U_2$,
and gluing them over $U_1\cap U_2$ with the aid of~$\lambda$.

In this way we get a {\em surjection\/} of abelian groups
$B_N^\times\twoheadrightarrow\Pic(S_N)$. It is not injective only because
the trivializations $\phi_1$ and $\phi_2$ are not canonical: they are 
defined up to multiplication by some invertible elements 
$s_1\in\Gamma(U_1,\sO_{U_1}^\times)=\bbZ^\times=\{\pm1\}$ and
$s_2\in\Gamma(U_2,\sO_{U_2}^\times)=|A_N|^\times=\{\pm1\}$, i.e.\
up to sign, hence $\lambda$ is also defined up to sign. Identifying
$B_N^\times/\{\pm1\}$ with the group $B_{N,+}^\times$ of positive
invertible elements in $B_N=\bbZ[N^{-1}]$, we obtain the following statement:

\begin{Propz}
The Picard group $\Pic(\CompZ\strut^{(N)})$ is canonically isomorphic
to $\log B_{N,+}^\times$, the abelian group of positive invertible
elements of~$B_N=\bbZ[N^{-1}]$, written in additive form.
\end{Propz}

If $p_1$, \dots, $p_r$ is the set of all distinct prime divisors of~$N$,
then $B_{N,+}^\times$ consists of rational numbers $p_1^{k_1}\cdots p_r^{k_r}$,
where $k_i$ are arbitrary integers, hence $\log B_{N,+}^\times$
is a free abelian group with basis $\{\log p_i\}$, 
and $\Pic(\CompZ\strut^{(N)})$
is also isomorphic to $\bbZ^r$, with basis $\{\sO(\log p_i)\}_{1\leq i\leq r}$.
Notice that these line bundles $\sO(\log p_i)$ are in one-to-one
correspondence with $p\mid N$, the only primes ``not entangled''
with~$\infty$ in $\CompuZ^{(N)}$ (cf.~\ptref{sp:top.compuz.N}).

Among other things, this result completes the proof 
of~\ptref{sp:compz.trans.not.iso} for arbitrary localization
theory~$\cT^?$.

\nxsubpoint\label{sp:pic.compz0}
(Line bundles over $\CompZ$.)
Since $\CompZ=\projlim_{N>1}S_N$, where $S_N=\CompZ\strut^{(N)}$,
we see that the category of line bundles over~$\CompZ$ is $\injLim$
of categories of line bundles over~$S_N$, hence
$\Pic(\CompZ)=\injlim_{N>1}\Pic(S_N)=\injlim_{N>1}\log B_{N,+}^\times=
\log\bbQ^\times_+$, i.e.\ any line bundle over $\CompZ$ is isomorphic
to exactly one line bundle of form $\sO(\log\lambda)$, with $\lambda$
a positive rational number. In other words, $\Pic(\CompZ)$ is a free
abelian group generated by all $\sO(\log p)$, $p\in\bbP$.

\nxsubpoint\label{sp:deg.linbdl.compz}
(Degree of a line bundle over~$\CompZ$.)
Let $\sL$ be a line bundle over~$\CompZ$. According to~\ptref{sp:pic.compz0},
it is isomorphic to $\sO(\log\lambda)$ for exactly one 
$\lambda\in\bbQ_+^\times$. We will say that $\log\lambda\in\log\bbQ_+^\times$
is the {\em (arithmetic) degree of~$\sL$,} and write $\deg\sL=\log\lambda$.
Thus $\deg:\Pic(\CompZ)\to\log\bbQ_+^\times$ is the isomorphism between
$\Pic(\CompZ)$ and $\bbQ_+^\times$ constructed in~\ptref{sp:pic.compz0}.

Now we want to obtain an explicit description of $\deg\sL$.

Recall that, according to~\ptref{sp:cat.finpres.compzmod}, a line bundle
$\sL$ over $\CompZ$ is essentially the same thing as a triple
$L=(L_\bbZ,L_\infty,\theta_L)$, where $L_\bbZ$ is a free $\bbZ$-module
of rank one, $L_\infty$ is a free $\Zninfty$-module of rank one
(we apply here~\ptref{cor:pic.zninfty.triv}, $\tilde L_\infty=\hat\eta^*\sL$
being a line bundle over $\Spec\Zninfty$), and $\theta_L:L_\bbZ\otimes_\bbZ\bbQ
\simto L_\infty\otimes_{\Zninfty}\bbQ$ is the isomorphism of
one-dimensional $\bbQ$-vector spaces, arising from canonical isomorphisms
of both sides to the generic fiber $\sL_\xi=\hat\xi^*\sL$.

Now choose a free $\bbZ$-generator $e_1$ of $L_\bbZ$ and a free
$\Zninfty$-generator $e_2$ of $L_\infty$. Both are defined uniquely up to
sign, since $\bbZ^\times=\{\pm1\}=\Zninfty^\times$. Then
both $\theta_L(e_1\otimes 1)$ and $e_2\otimes 1$ are generators of 
one-dimensional $\bbQ$-vector space $L_\infty\otimes_{\Zninfty}\bbQ\cong
\sL_\xi$, hence $e_2\otimes1=\lambda\cdot\theta_L(e_1\otimes 1)$ for some
$\lambda\in\bbQ^\times$, defined uniquely up to sign. We can make
$\lambda$ unique by requiring $\lambda>0$. Then $\log\lambda\in
\log\bbQ^\times_+$ is exactly the degree of~$\sL$
(cf.\ \ptref{sp:pic.compz.N} and~\ptref{sp:cat.finpres.compzmod}).

\nxsubpoint
($\CompZ\strut^{(N)}$ is similar to $(\bbP^1)^r$.)
Let $N>1$, $S_N:=\CompZ\strut^{(N)}$ be as before. Denote by
$p_1$, \dots, $p_r$ all distinct prime divisors of~$N$. We have just seen
that $\Pic(S_N)\cong\bbZ^r$, similarly to the Picard group
$\Pic((\bbP^1)^r)$ of the product of~$r$ projective lines. In this respect
$S_N$ is quite similar to $P^r:=(\bbP^1)^r$. Notice, however, the following:
\begin{itemize}
\item[(a)] The automorphism group $\Aut(P^r)$ acts transitively on the
canonical basis elements of $\Pic(P^r)=\bbZ^r$, i.e.\ all
these line bundles $\sL_i:=\sO(0,\ldots,1,\ldots,0)$ over $P^r$ have the same
properties, while $\Aut(S_N)=1$ doesn't act transitively on the
generators $\sL'_i:=\sO(\log p_i)$ of $\Pic(S_N)$.
\item[(b)] This analogy between Picard groups shouldn't extend to higher
Chow groups: we expect that $c_1:\Pic(S_N)\to CH^1(S_N)$ is an 
isomorphism for~$S_N$, similarly to $c_1$ over $P^r$, but $CH^i(S_N)=0$
for $i>1$, in contrast with $CH^i(P^r)\cong\bigwedge^i\Pic(P^r)$.
\end{itemize}

\nxsubpoint
(Global sections of $\sO_{S_N}$.)
Let us compute the global sections of the structural sheaf of 
$S_N=\CompZ\strut^{(N)}$. Since $S_N=U_1\cup U_2$, with
$\Gamma(U_1,\sO)=\bbZ$, $\Gamma(U_2,\sO)=A_N$ and $\Gamma(U_1\cap U_2,\sO)=
B_N$, the sheaf condition for~$\sO$ yields $\Gamma(S_N,\sO)=\bbZ\times_{B_N}
A_N=\bbZ\cap A_N=\Fpm$, where the intersection is computed inside 
$B_N\subset\bbQ$. In order to check $\bbZ\cap A_N=\Fpm$ we write
$\bbZ\cap A_N=\bbZ\cap B_N\cap\Zninfty=\bbZ\cap\Zninfty$, and 
$(\bbZ\cap\Zninfty)(n)=\bbZ(n)\cap\Zninfty(n)=\{\lambda\in\bbZ^n\,:\,
\|\lambda\|\leq1\}=\Fpm(n)$, where $\vnorm$ is the $L_1$-norm on~$\bbQ^n$
as before. Actually this formula $\Gamma(\CompZ,\sO)=\bbZ\cap\Zinfty=\Fpm$
was exactly our original motivation for defining $\Fpm$ the way we've
defined it.

An immediate implication is that we can study all $S_N$ and their
projective limit $\CompZ$ as (pro-)generalized schemes over
$S_1:=\Spec\Fpm$, or over $\Spec\Fone$ or $\Spec\Fempty$, but these
cases almost exhaust all reasonable possibilities.

\nxsubpoint\label{sp:glob.sect.O.log.lambda}
(Global sections of $\sO(\log\lambda)$.)
Let $N>1$, $S_N$ be as above, $\lambda\in B_{N,+}^\times$. Let us
compute the set of global sections $\Gamma(S_N,\sO(\log\lambda))$.
Using trivializations $\phi_1$ and $\phi_2$ as in~\ptref{sp:pic.compz.N},
differing by multiplication by $\lambda$ over $U_1\cap U_2$,
together with the sheaf condition for $\sO(\log\lambda)$ with
respect to the cover $S_N=U_1\cap U_2$, we see that
$\Gamma(S_N,\sO(\log\lambda))\cong\{(x,y)\in\bbZ\times|A_N|\,:\,
x=\lambda y\}\cong\{x\in\bbZ\,:\,|\lambda^{-1}x|\leq1\}=\bbZ\cap
\lambda|A_N|=\bbZ\cap[-\lambda,\lambda]$.

For example, if $M\geq1$ is an integer, such that $M\mid N^\infty$,
then $\sO(\log M)$ admits $2M+1$ global sections, corresponding to integers
from $-M$ to $M$. The global sections corresponding to $\pm1$ are
invertible over $\Spec\bbZ$, while those corresponding to $\pm M$ are
invertible over $\Spec A_N$, or more precisely, outside the primes $p\mid M$.

Notice that $\Gamma(S_N,\sO(\log\lambda))$ doesn't change when we replace
$N$ by any its multiple~$N'$, hence $\Gamma(\CompZ,\sO(\log\lambda))=
\Hom(\sO(0),\sO(\log\lambda))$ equals 
$\injlim_{N'\geq1}\Hom_{\sO_{S_{NN'}}}(\sO(0),
\sO(\log\lambda))=\injlim_{N'\geq1}\Gamma(S_{NN'},\sO(\log\lambda))
\cong\bbZ\cap[-\lambda,\lambda]$ as well.

\nxsubpoint
(Global sections of twisted vector bundles.)
For any $\sO_{S_N}$-module~$\sF$ and any $\lambda\in B_{N,+}^\times$,
we denote by $\sF(\log\lambda)$ the corresponding {\em Serre twist\/}
$\sF\otimes_{\sO}\sO(\log\lambda)$. If $\sF$ was a vector bundle,
the same is true for $\sF(\log\lambda)$; in particular, we can start
from the trivial vector bundle $\sE_0:=\sO_{S_N}^{(r)}=L_\sO(r)$ of
rank~$r\geq0$, and consider $\sE:=\sE_0(\log\lambda)=L_\sO(r)\otimes_\sO
\sO(\log\lambda)$.

Let us compute the global sections of such twisted trivial vector bundles.
The same reasoning as in~\ptref{sp:glob.sect.O.log.lambda} yields
$\Gamma(\sO_{S_N},L_\sO(r)\otimes_\sO\sO(\log\lambda))=
\bbZ^r\cap\lambda A_N(r)=\bbZ^n\cap\lambda\cdot\Zninfty(r)=
\{\mu=(\mu_1,\ldots,\mu_r)\in\bbZ^n\,:\,\|\mu\|\leq\lambda\}$.
In other words, global sections of such vector bundles correspond
to integral points $(\mu_1,\ldots,\mu_r)\in\bbZ^n$, lying inside
octahedron $\sum_i|\mu_i|\leq\lambda$.

\nxsubpoint
(Ample line bundles over $\CompZ\strut^{(N)}$.)
Recall that the ample line bundles over $P^r=(\bbP^1)^r$ are 
$\sO(k_1,\ldots,k_r)$ with all $k_i>0$, so by analogy between
$\Pic(S_N)$ and $\Pic(P^r)$ we might expect that
the ample line bundles over $S_N=\CompZ\strut^{(N)}$ are
$\sO(\log M)$ with $M=p_1^{k_1}\cdots p_r^{k_r}$, where
all $k_i>0$, and $p_i$ are the distinct prime divisors of~$N$.

Since the above condition on~$M$ is equivalent to $M\mid N^\infty$
and $N\mid M^\infty$, and for such $M$ and~$N$ we have canonical
isomorphisms $S_N\cong S_{NM}\cong S_M$, we are reduced to checking
whether $\sO(\log N)$ is ample over~$S_N$.

\nxsubpoint
(Graded ring $\Gamma_*(\sO_{S_N})$ with respect to $\sL:=\sO_{S_N}(\log N)$.)
According to our general definitions~\ptref{sp:def.ample.proj}, 
we are to construct
the graded generalized ring $R:=\Gamma_*(\sO_{S_N})$, where
$\Gamma_d(\sO_{S_N})(n)=\Gamma(S_N,L_\sO(n)\otimes_\sO\sO(d\log N))$,
and check whether the natural map $S_N\to\Proj R$ is an isomorphism.

We know already that $R_d(n)=\{(\lambda_1,\ldots,\lambda_n)\in\bbZ^n\,|\,
\sum_i|\lambda_i|\leq N^d\}$. It will be convenient to write elements
of $R_d(n)$ in form $T^d(\lambda_1,\ldots,\lambda_n)$, or even
$\lambda_1 T^d\{1\}+\cdots+\lambda_n T^d\{n\}$, thus embedding $R$
into graded $\bbZ$-algebra $\bbZ[T]$; 
this embedding is actually compatible with the unit and multiplication
maps $\mu^{(k,d)}_{n,e}$, so it completely determines the structure of
generalized graded ring on $R$.

Sometimes it is convenient to think about formal variable~$T$
as if it were equal to $1/N$.

\nxsubpoint\label{sp:elems.f1f2}
(Two special elements of degree one.)
Notice that $|R|=R(1)=\bigsqcup_{d\geq0}R_d(1)$ consists of expressions
$uT^d$, with $d\geq0$, $u\in\bbZ$, $|u|\leq N^d$. In particular, we have
two special elements of degree one: $f_1:=T$ and $f_2:=NT$.

Denote by $\ga\subset|R|_+$ the graded ideal generated by $f_1$ and $f_2$.
It is clearly equal to the image of graded $R$-homomorphism $R(2)[-1]\to|R|^+$,
which maps $\{1\}\mapsto f_1$, $\{2\}\mapsto f_2$. Since $R_{d-1}(2)=
\{uT^{d-1}\{1\}+vT^{d-1}\{2\}\,:\,|u|+|v|\leq N^{d-1}\}$, we see that
$\ga_d\subset|R|_d$, where $d\geq1$, 
consists of elements $wT^d$, which can be written in form
$w=u+Nv$ for some $u$, $v\in\bbZ$ with $|u|+|v|\leq N^{d-1}$.
On the other hand, $|R|_{+,d}=|R|_d$ consists of all $wT^d$ with $|w|\leq N^d$,
so we cannot expect $\ga_d=|R|_{+,d}$ for $N>2$: indeed, $w=N^d-1$ cannot
be written as $w=u+Nv$ with integer $u$, $v$, such that $|u|+|v|\leq N^{d-1}$.

However, we claim that $\rad(\ga)=|R|_+$. Indeed, let $f:=wT^d$, $d>0$,
$|w|\leq N^d$ be any element of $|R|_+$. If $|w|=N^d$, we can take
$u=0$, $v=w/N=\pm N^{d-1}$ and write $f=wT^d=vT^{d-1}\cdot NT=vT^{d-1}
\cdot f_2\in\ga_d$, and we are done. If $|w|<N^d$, we can find an integer
$k\geq1$, such that $|N^{-d}w|^k\leq 1-1/N$ and $dk\geq3$. Then
$f^k=w^kT^{dk}$ with $|w^k|\leq(1-1/N)N^{dk}\leq N^{dk}-N^2$. Now write
$w^k=vN+u$ with $u$, $v\in\bbZ$, $|u|<N$, $uw^k\geq0$. Then
$|v|\leq|w^k|/N\leq N^{dk-1}-N$, hence $|u|+|v|\leq N^{dk-1}$, so
$f^k=w^k T^{dk}$ lies in $\ga$, and we are done.

\nxsubpoint\label{sp:loc.f1f2}
(Localizations of $R$ with respect to $f_1$, $f_2$ and $f_1f_2$.)
Now we want to compute generalized rings $R_{(f_1)}$, $R_{(f_2)}$
and $R_{(f_1f_2)}$.

(a) $R_{(f_1)}(n)=R_{(T)}(n)$ consists of fractions
$(u_1T^d\{1\}+\cdots+u_nT^d\{n\})/T^d$ with $u_i\in\bbZ$, 
$\sum_i|u_i|\leq N^d$, for all $d\geq0$. 
Therefore, $R_{(f_1)}(n)$ consists of all expressions $u_1\{1\}+\cdots+
u_n\{n\}$ with $u_i\in\bbZ$, i.e.\ $R_{(f_1)}(n)=\bbZ^n=\bbZ(n)$ 
and $R_{(f_1)}=\bbZ$.

(b) $R_{(f_2)}(n)=R_{(NT)}(n)$ consists of 
$(u_1T^d\{1\}+\cdots+u_nT^d\{n\})/(N^dT^d)$ with $u_i\in\bbZ$, 
$\sum|u_i|\leq N^d$, i.e.\ of expressions $v_1\{1\}+\cdots+v_n\{n\}$
with $v_i\in\bbZ[N^{-1}]=B_N$, $\sum_i|v_i|\leq1$, hence
$R_{(f_2)}=B_N\cap\Zninfty=A_N$.

(c) Finally, $R_{(f_1f_2)}=R_{(NT^2)}$ is isomorphic to $B_N=\bbZ[N^{-1}]$,
either by a similar direct computation, or because of formula
$R_{(f_1f_2)}\cong R_{(f_1)}[(f_2/f_1)^{-1}]=\bbZ[N^{-1}]$
shown in~\ptref{sp:loc.homog.loc}.

\nxsubpoint
(Projectivity of $\CompZ\strut^{(N)}=\Proj R$.)
We see that $\Proj R=D_+(f_1)\cup D_+(f_2)$ by
\ptref{sp:covers.of.projspec} and~\ptref{sp:elems.f1f2}. Furthermore,
$D_+(f_1)=\Spec R_{(f_1)}=\Spec\bbZ$, $D_+(f_2)=\Spec R_{(f_2)}=\Spec A_N$,
and $D_+(f_1)\cap D_+(f_2)=D_+(f_1f_2)=\Spec R_{(f_1f_2)}=\Spec B_N$,
i.e.\ $\Proj R$ can be constructed by gluing together $\Spec\bbZ$ and
$\Spec A_N$ along their open subsets isomorphic to $\Spec B_N$.

In other words, we have just shown that $\Proj R$ is isomorphic
to $\CompZ\strut^{(N)}$. This proves ampleness of $\sO(\log N)$
over $S_N:=\CompZ\strut^{(N)}$, as well as projectivity of $S_N$
(absolute or over~$\Fpm$).

Notice, however, that $S_N:=\CompZ\strut^{(N)}$ is a {\em non-unary\/}
projective generalized scheme,
cf.~\ptref{sp:proj.sch.vs.closed.subsch.projvb}.  In particular, even
if we know that $\Gamma(S_N,\sO(\log N))$ is a free $\Fpm$-module of
rank~$N$, we still shouldn't expect $S_N$ to be isomorphic to a
``closed'' subscheme of $\bbP^{N-1}_{\Fpm}$.

\nxsubpoint\label{sp:square.compz} (Square of $\CompZ$.)
We have seen in~\ptref{sp:tens.sqr.z} that $\bbZ\otimes\bbZ=\bbZ$,
when the tensor product (i.e.\ coproduct of generalized commutative
algebras) is taken over $\Fempty$, $\Fone$ or $\Fpm$, hence
$\Spec\bbZ\times_{\Spec\Fone}\Spec\bbZ=\Spec\bbZ$, i.e.\
$\Spec\bbZ\to\Spec\Fone$ is a monomorphism of generalized schemes.

Now consider $S_N:=\CompZ\strut^{(N)}=U_1\cup U_2$, where 
$U_1=\Spec\bbZ$ and $U_2=\Spec A_N$ as before. We claim that
{\em $S_N\times_{\Spec\Fpm}S_N=S_N$, i.e.\ $S_N\to\Spec\Fpm$ is 
a monomorphism of generalized schemes.} Since $\CompZ=\projlim_{N>1}S_N$,
we obtain $\CompZ\times_{\Spec\Fpm}\CompZ=\CompZ$ (e.g.\ in the category
of pro-generalized schemes) as well.

So let us show that $\Delta:=\Delta_{S_N}:S_N\to S_N\times S_N$ (all products
here are to be computed over $\Spec\Fpm$) is an isomorphism. Since
$\Delta^{-1}(U_i\times U_j)=U_i\cap U_j$, and open subschemes $U_i\times U_j$,
$1\leq i,j\leq 2$, cover $S_N\times S_N$, we are reduced to showing that
canonical morphisms $U_i\cap U_j\to U_i\times U_j$ are in fact isomorphisms.
If $i=j=1$, this is equivalent to $\bbZ\otimes_{\Fpm}\bbZ=\bbZ$,
valid by~\ptref{sp:tens.sqr.z}. If $i=j=2$, we have to show $A_N
\otimes_{\Fpm}A_N\cong A_N$. Since $A_N$ is generated over $\Fpm$ 
by one averaging operation $s_N$ (cf.~\ptref{rem:gen.of.AN},(c);
this is the only part of the proof valid over $\Fpm$, but not 
over $\Fone$ or $\Fempty$), 
all we have to check is that whenever
we have {\em two\/} commuting operations $s_N$, $s'_N$, both satisfying
\eqref{eq:avg.op.idemp}--\eqref{eq:avg.op.canc}, then necessarily $s'_N=s_N$.
This can be shown exactly in the same way as in~\ptref{sp:tens.sqr.zninfty},
applying commutativity relation between $s_N$ and $s'_N$ to a Latin square
filled by basis elements $\{1\}$, \dots, $\{N\}$.

Now only case $i\neq j$ remains. Since $U_1\cap U_2=\Spec B_N$, we are to 
check that $A_N\otimes\bbZ=B_N$. Put $\Sigma:=A_N\otimes\bbZ$, and denote
by $f$ the unary operation of $A_N$ and of~$\Sigma$, given by
$f:=s_N(\bu,0,\ldots,0)$, i.e.\ $f=1/N$. 
One easily checks that $f\cdot N=\bu$ in $|\Sigma|$,
by applying commutativity relation for $s_N$ and the $N$-tuple addition
operation $\{1\}+\dots+\{N\}$ to $N\times N$-matrix $(x_{ij})$ with
$x_{ii}=\bu$, $x_{ij}=0$ for $i\neq j$. This means that $f=1/N$ is invertible
in~$\Sigma$, hence $\Sigma=\Sigma[(1/N)^{-1}]=(A_N\otimes\bbZ)[(1/N)^{-1}]=
A_N[(1/N)^{-1}]\otimes\bbZ=B_N\otimes\bbZ=B_N\otimes_\bbZ\bbZ\otimes\bbZ=
B_N\otimes_\bbZ\bbZ=B_N$. We've used here $A_N[(1/N)^{-1}]=B_N$,
valid by~\ptref{sp:BN.as.localiz.AN}, as well as $\bbZ\otimes\bbZ=\bbZ$, shown
in~\ptref{sp:tens.sqr.z}.

This reasoning is a modification of that of~\ptref{sp:tens.z.zinfty},
used there to prove $\Zninfty\otimes\bbZ=\bbQ$ and $\Zinfty\otimes\bbZ=\bbR$.

\nxsubpoint\label{sp:smirnov.descr.compz}
(Projectivity of $f_N^{NM}$. $\CompZ$ as 
an infinite resolution of singularities.)
Since $f_N^{MN}:S_{NM}\to S_N$ is a morphism between two projective generalized
schemes, it is very reasonable to expect $f_N^{NM}$ to be 
a projective morphism. Let's accept this without proof for the time being.
Then we get the following picture. Each $S_N$ is something like a curve,
admitting $\Spec\bbZ$ as an open subset, and all projective maps $f_N^{NM}$ 
are isomorphisms over this open subset, i.e.\ they change something
only over~$\infty$, so they can be thought of as a sort of blow-up
(or a series of blow-ups) over~$\infty$. Since the complement of $\Spec\bbZ$
after this ``blow-up'' still consists of one point (at least for $\cT^u$),
this looks very much like blowing up a complicated 
cusp-like singularity over a curve.
This ``cusp'' appears to have infinite order, since it cannot be resolved
in any finite number of steps; we have an infinite sequence of
``blow-ups'', and the ``smooth model'' $\CompZ$ appears only
as a pro-object.

I would like to mention here that A.~Smirnov, when we discussed possible
meanings of $\CompZ$ several months ago, and this work was only in its
initial stage, mentioned to me that he expects something like an 
infinite resolution of singularities over a cusp-like point~$\infty$ to be
necessary to construct the ``correct'' (or ``smooth'') compactification
$\CompZ$. At that time I was somewhat skeptical about this.
However, this completely intuitive perception of $\CompZ$ turned out to
agree almost completely with the rigorous theory constructed here!
This remains a mystery for me.

\nxpointtoc{Models over $\CompZ$, $\Zinfty$ and $\Zninfty$}
Now we are going to discuss models of algebraic varieties $X/\bbQ$,
showing in particular that any such $X$ admits at least one finitely presented
model $\sX$ over $\CompZ$. We also discuss models over $\Zninfty$,
$\Zninfty$ and $\barZinfty$.

Then we are going to discuss the relation to Banach and even K\"ahler metrics,
e.g.\ by showing that the Fubiny--Study metric on $\bbP^N(\bbC)$ 
corresponds to a certain model $\bar\bbP^N/\barZinfty$ of $\bbP^N$.
Unfortunately,
models corresponding to ``nice'' metrics usually turn out to have
``bad'' algebraic properties (e.g.\ the model $\bar\bbP^N/\barZinfty$
just discussed is not finitely presented), and conversely. This probably
explains the complexity of (classical) Arakelov geometry, which might
be compared to the study of not necessarily finitely presented models
$\sX/\bbZ_p$ of an algebraic variety $X/\bbQ$. This also explains why
we can construct our variant of Arakelov geometry in a simpler fashion,
just by considering finitely presented models whenever possible.

\nxsubpoint 
(Reduction from $\CompZ$ to $\Spec\Zninfty$.)
Recall that constructing a finitely presented model $\bar\sX/\CompZ$ 
of an algebraic variety $X/\bbQ$ is essentially the same thing
as constructing finitely presented models $\sX/\bbZ$ and $\sX_\infty/\Zninfty$
of this~$X$ (cf.~\ptref{sp:finpres.models.compz}), 
so we need to concentrate our efforts only on finitely presented
models over~$\Zninfty$, assuming the first part (models over~$\bbZ$) 
to be well-known.

\nxsubpoint\label{sp:fix.arch.val.rg}
(Notations: $V\subset K$.)
Let us fix some field~$K$ and an ``archimedian valuation ring'' $V\subset K$
in the sense of~\ptref{sp:val.rings}. Usually we consider either $V=\Zninfty$,
$K=\bbQ$, or $V=\Zinfty$, $K=\bbR$, or $V=\barZinfty$, $K=\bbC$.
Denote by $\fnorm$ the archimedian valuation of~$K$ corresponding to~$V$,
and let $\tilde\bbQ\subset K$ be the closure of $\bbQ\subset K$. Obviously,
$\tilde\bbQ$ is canonically isomorphic to a subfield of~$\bbR$.
We'll consider only couples $(V,K)$, which satisfy the following extra
condition: {\em $|\lambda|$ belongs to $\tilde\bbQ\subset\bbR$ for any
$\lambda\in K$.} This enables us to consider $|\lambda|$ as an element of~$K$.

Notice that for any $f\in V$ with $0<|f|<1$ we have $V[f^{-1}]=K$,
hence $\Spec K$ is isomorphic to a principal open subset of~$\Spec V$.
The image of the only point of $\Spec K$ in $\Spec V$ is the
{\em generic\/} or {\em open point\/} $\xi$ of $\Spec V$. Any $\sX/V$
contains its generic fiber $\sX_\xi=\sX_{(K)}=\sX\otimes_VK$ as an
open generalized subscheme.

Now we are going to discuss models $\sX/V$ of a finitely presented
scheme $X/K$, i.e.\ an algebraic variety over~$K$.

\nxsubpoint
(Norm on $K[T_1,\ldots,T_k]$.)
We denote by $\vnorm$ the (coefficientwise) 
$L_1$-norm on the polynomial algebra
$A=K[T_1,\ldots,T_k]$, defined as follows:
\begin{equation}
\|P\|=\sum_{\alpha\in\bbN_0^k}|c_\alpha|,\quad\text{if }
P=\sum_{\alpha\in\bbN_0^k}c_\alpha T^\alpha
\end{equation}
Here we use the multi-index notation: $T^\alpha=
T_1^{\alpha_1}\cdots T_k^{\alpha_k}$. Furthermore, for any
element $\lambda=\sum_{i=1}^n\lambda_i\{i\}\in A(n)=A^n$ we put
$\|\lambda\|:=\sum_i\|\lambda_i\|$. Then generalized polynomial
ring $B:=V[T_1,\ldots,T_k]\subset K[T_1,\ldots,T_k]$ can be 
described as follows: $B(n)=\{\lambda\in A(n)\,:\,\|\lambda\|\leq 1\}$.

Notice that our assumptions on $(V,K)$ imply $\|\lambda\|\in K$ for any
$\lambda\in A(n)$, hence we can write any such $\lambda\neq 0$
as $c\cdot\mu$, where $c:=\|\lambda\|\in K^\times$, and $\mu:=c^{-1}\lambda
\in A(n)$ is of norm~one, and in particular belongs to~$B(n)$.

\nxsubpoint
(Finitely presented models of affine varieties.)
Let $X$ be an affine variety over~$K$, i.e.\ a finitely presented
affine scheme over~$\Spec K$. Then $X=\Spec A$, where
$A=K[T_1,\ldots,T_k]/(f_1,\ldots,f_m)$, for some finite number
of polynomials $f_j\in A_0:=K[T_1,\ldots,T_k]$, which can be
assumed to be $\neq0$. Furthermore, we are free to multiply 
$f_j$ by non-zero scalars from~$K$, so we can make $\|f_j\|\leq1$,
or even $\|f_j\|=1$, i.e.\ $f_j\in|B_0|$, where $B_0:=V[T_1,\ldots,T_k]\subset
A_0$ is the (generalized unary) polynomial algebra over~$V$.

Put $B:=B_0/\langle f_1=0,\ldots, f_m=0\rangle=V[T_1^{[1]},\ldots,T_k^{[1]}\,|
\,f_1=0,$ $\dots,$ $f_m=0]$. Then $\sX:=\Spec B$ is an affine
finitely presented generalized scheme over~$V$, such that
$\sX_{(K)}=\Spec(B\otimes_VK)=\Spec A=X$, i.e.\ {\em
any affine variety $X/K$ admits at least one finitely presented
affine model $\sX/V$.}

\nxsubpoint\label{sp:finpres.mod.projvar}
(Finitely presented models of projective varieties.)
Similarly, let $X=\Proj A$ be a projective variety over~$K$. Then
$A=K[T_0,\ldots,T_k]/(f_1,\ldots,f_m)$ for some homogeneous polynomials
$f_j\neq0$. Multiplying $f_j$ by suitable constants from $K^\times$,
we can assume $\|f_j\|=1$, i.e.\ $f_j$ lie in the generalized
graded polynomial algebra $R:=V[T_0,\ldots,T_k]$, with all free 
generators $T_i$ unary of degree one. Put $B:=R/
\langle f_1=0,\ldots,f_m=0\rangle$. This is a finitely presented
unary graded algebra over~$V$, such that $|B|_+$ is generated by
finitely many elements of degree one. We conclude that
$\sX:=\Proj B$ is a finitely presented projective model of $X$ over~$V$.
In particular, {\em any projective variety $X/K$ admits at least one
projective finitely presented model $\sX/V$,} and
{\em any projective variety $X/\bbQ$ admits a finitely presented model
$\bar\sX$ over~$\CompZ$.}

\nxsubpoint (Presence of $V$-torsion in finitely presented models.)
Notice that the finitely presented models $\sX$ constructed above usually
have some $V$-torsion, i.e.\ natural homomorphisms $\sO_{\sX}\to
\sO_{\sX}\otimes_VK$ are not necessarily monomorphisms of sheaves of
generalized rings. In fact, almost all finitely presented schemes
over $\Spec\Zninfty$, with the notable exception of $\bbA^n_{\Zninfty}$
and $\bbP^n_{\Zninfty}$, do have some ``torsion over $\infty$'',
hence the same is true for finitely presented models $\bar\sX/\CompZ$:
we can eliminate torsion over finite primes $p$ if we want,
but in general not over~$\infty$. 

One might think about this ``torsion over~$\infty$'' in models $\sX$
as a sort of ``built-in analytic torsion''. This explains why we
expect the Grothendieck--Riemann--Roch formula to hold for such (finitely
presented) models without any complicated correction terms: 
all analytic torsion is already accounted for by this torsion over~$\infty$,
built into the structure of the model chosen.

\nxsubpoint\label{sp:exist.tf.models} (Existence of torsion-free models.)
Suppose that our variety $X/K$ can be embedded as a closed subvariety
into some variety $P/K$, which admits a torsion-free model $\bar P/V$.
This condition is always fulfilled for any projective variety~$X$,
since we can take $P:=\bbP^N_K$, $\bar P:=\bbP^N_V$ for a suitable $N\geq0$.

Then we can construct a torsion-free model $\bar X/V$ of $X/K$,
simply by taking the ``scheme-theoretic closure'' of $X\subset P=\bar P_\xi
\subset\bar P$ inside $\bar P$. More formally, let $i:X\to P$ be the
closed immersion of $X$ into $P$, $j:P\to\bar P$ be the open immersion
of the generic fiber of~$\bar P$ into~$\bar P$. 
Consider following homomorphism of sheaves of generalized 
$\sO_{\bar P}$-algebras: 
$\phi:
\sO_{\bar P}\to j_*j^*\sO_{\bar P}=j_*\sO_P\twoheadrightarrow j_*i_*\sO_X$.
The first arrow here is the adjointness morphism for $j_*$ and $j^*$,
and the second one is obtained by applying $j_*$ to 
$\sO_P\twoheadrightarrow i_*\sO_X$. Since $j$ is quasicompact and 
quasiseparated (being affine), all $\sO_{\bar P}$-algebras involved
are quasicoherent, hence the same is true for $\sO_{\bar X}:=\phi(\sO_{\bar P})
\subset j_*i_*\sO_X$. Now let $\bar X$ be the ``closed'' subscheme
of~$\bar P$ defined by this quasicoherent strict quotient
$\sO_{\bar P}\twoheadrightarrow\sO_{\bar X}$. It is immediate 
that this $\bar X$ is indeed a torsion-free model of~$X$. 

Notice, however, that in general $\bar X$ will be of finite type over~$V$,
but not finitely presented. Therefore, such torsion-free models over $\Zninfty$
of varieties $X/\bbQ$ cannot be used (together with some flat model
$\sX/\bbZ$) to construct models $\bar\sX$ of $X$ over $\CompZ$.
All we can do is to write $\bar X$ as a projective limit of finitely
presented models $\bar X_\alpha$ of $X$ over~$V$, use these $\bar X_\alpha$
to obtain a projective system of finitely presented 
models $\bar\sX_\alpha/\CompZ$ of $X$, and put 
$\bar\sX:=\quotprojlim\bar\sX_\alpha$. In this way we obtain a
torsion-free (``flat'') finitely presented model $\bar\sX/\CompZ$
of~$X$ in the category of pro-generalized schemes. In the fancy language 
of~\ptref{sp:smirnov.descr.compz} one might say that in order to
eliminate ``embedded analytic torsion'' of a finitely presented
model over~$\CompZ$ we have to perform an infinite resolution of
singularities.

\nxsubpoint\label{sp:ext.ratpt.sect}
(Extending rational points to sections.)
One expects that any rational point $P\in X(\bbQ)$ of a {\em projective\/}
variety $X/\bbQ$ extends to a unique section $\sigma_P:\CompZ\to\bar\sX$,
for any projective model $\bar\sX/\CompZ$ of~$X$. According 
to~\ptref{sp:finpres.sch.compz}, finding such a $\sigma_P$ 
is equivalent to finding a section $\sigma'_P$ of a model
$\sX/\bbZ$, and a section $\sigma''_P$ of model $\sX_\infty/\Zninfty$,
both inducing~$P$ on the generic fiber. Existence and uniqueness
of $\sigma'_P$ is very well known, so we are reduced to
the case of $\Zninfty$-models.

In order to lift $K$-rational points of $\sX/V$ to $V$-rational points
one usually has to apply the valuative criterion of properness
in the easy direction. We would like, therefore, to discuss 
whether valuative criteria still work (at least in the easy direction) 
over archimedian valuation rings $V\subset K$ of~\ptref{sp:fix.arch.val.rg}.

\begin{DefD}
We say that a generalized scheme morphism $f:X\to S$ satisfies the
{\em valuative criterion of separability\/} (resp.\ {\em properness\/})
{\em (with respect to $(V,K)$)}, if, given any $u\in X(K)$ 
and $v\in S(V)$ fitting into commutative diagram
\begin{equation}\label{eq:diag.val.crit}
\xymatrix@C+1pc{
\Spec K\ar[r]^<>(.5){u}\ar[d]&X\ar[d]^{f}\\
\Spec V\ar[r]^<>(.5){v}\ar@{-->}[ru]|{w}&S}
\end{equation}
one can find at most one (resp.\ exactly one) $w\in X(V)$ completing
the above diagram.
\end{DefD}

It is clear that the class of morphisms satisfying the valuative criterion of
separability (resp.\ properness) with respect to $(V,K)$
is stable under composition and base change.

Notice that, given a diagram as above, we can always replace $X$
by $X_{(\Spec V)}=X\times_S\Spec V$, $f$ by $f_{(\Spec V)}$, and reduce the
above ``lifting problem'' to the case $S=\Spec V$, $v=\id_S$.

\begin{LemmaD}
If $f:X\to S$ is ``separated'' in the sense of~\ptref{sp:closed.imm.sch},
i.e.\ if $\Delta_{X/S}$ is a ``closed'' immersion 
(cf.\ {\em loc.~cit.}),
then $f$ satisfies the valuative criterion of separability with
respect to $(V,K)$.
\end{LemmaD}
\begin{Proof}
First of all, we can replace $f:X\to S$ by its base change with respect
to $v:\Spec V\to S$ from~\eqref{eq:diag.val.crit}, and assume
$S=\Spec V$, $v=\id_S$. Let $w$, $w':S\to X$ be two sections of $f:X\to S$,
coinciding on $\Spec K\subset S$, i.e.\ having the same generic fibers
$w_{(K)}=w'_{(K)}$. Let $T\subset S$ be the pullback of 
$\Delta_{X/S}:X\to X\times_S X$ with respect to $(w,w'):S\to X\times_SX$.
On one hand, $T\subset S$ is a ``closed'' subscheme of $S$,
since $f:X\to S$ has been supposed to be ``separated''. On the other hand,
$T=\Ker(w,w':S\rightrightarrows X)$, hence $T$ is the largest subobject
of~$S$, such that $w|_T=w'|_T$, and the construction
of~$T$ commutes with base change. Now $w_{(K)}=w'_{(K)}$ implies
$T_{(K)}=\Ker(w_{(K)},w'_{(K)})=S_{(K)}$. Applying 
lemma~\ptref{l:cl.subsch.cont.genpt} below, 
we see that $S=T$, hence $w=w'$, q.e.d.
\end{Proof}

\begin{LemmaD}\label{l:cl.subsch.cont.genpt}
If $T$ is a ``closed'' subscheme of $S=\Spec V$, such that $T_{(K)}=S_{(K)}$,
then $T=S$. Furthermore, this statement holds for any monomorphism of
generalized rings $V\to K$, not only for embeddings of 
an archimedian valuation ring~$V$ into its fraction field~$K$.
\end{LemmaD}
\begin{Proof}
By definition a ``closed'' subscheme $T\subset S=\Spec V$ equals
$\Spec W$ for some strict quotient $W$ of~$V$. On the other hand,
condition $T_{(K)}=S_{(K)}$ means that $\Spec K\to\Spec V$ factorizes 
through~$T=\Spec W$, i.e.\ that monomorphism $V\hookrightarrow K$ 
factorizes through strict epimorphism $\pi:V\twoheadrightarrow W$. 
This means that $\pi$ is both a strict epimorphism and a monomorphism,
hence an isomorphism, i.e.\ $W=V$ and $T=S$.
\end{Proof}

\begin{LemmaD}\label{l:cl.imm.are.proper}
Any ``closed'' immersion $f:X\to Y$ satisfies the valuative criterion of 
properness with respect to $(V,K)$.
\end{LemmaD}
\begin{Proof}
Applying base change with respect to $\Spec V\stackrel v\to Y$ as before,
we are reduced to case $Y=S=\Spec V$, $X=T\subset S$ a ``closed'' subscheme
of~$S$. We have to check that any point $u\in T_S(K)$ lifts to exactly
one point $w\in T_S(V)$. If $T_S(K)=\emptyset$, then $T_S(V)=\emptyset$ and
the statement is trivial. If $T_S(K)\neq\emptyset$, then
$T_S(K)\subset S_S(K)$ is a one-element set; this means that
$\Spec K\to S$ factorizes through $T$, i.e.\ $T_{(K)}=S_{(K)}$.
By~\ptref{l:cl.subsch.cont.genpt} we obtain $T=S$, hence $T_S(V)=S_S(V)=
\{\id_S\}$, and we are done.
\end{Proof}

\begin{ThD}\label{th:val.crit.prop.proj}
Let $\sF$ be a quasicoherent $\sO_S$-module of finite type over
some generalized scheme~$S$, and $X$ be a ``closed'' subscheme
of $\bbP_S(\sF)$. Then $X\to S$ satisfies the valuative criterion
of properness with respect to $(V,K)$.
\end{ThD}
\begin{Proof}
(a) By base change we can again assume $S=\Spec V$, $v=\id_S$.
Then $\sF=\tilde M$ for a finitely generated $V$-module~$M$, 
so we can find a strict epimorphism $\pi:V(n+1)\twoheadrightarrow M$,
which induces a strict epimorphism of unary graded $V$-algebras
$S_V(V(n+1))\twoheadrightarrow S_V(M)$, hence a ``closed'' immersion
$\bbP_S(\sF)\to\bbP_S(\widetilde{V(n+1)})=\bbP_S^n$. Now
$X\to\bbP_S^n$ is a ``closed'' immersion, hence it satisfies the valuative
criterion of properness for $(V,K)$ by~\ptref{l:cl.imm.are.proper},
so we are reduced to checking that $P:=\bbP_S^n\to S$ satisfies
the valuative criterion of properness with respect to $(V,K)$,
and we can still assume $v=\id_S$.

In other words, we are to check that the canonical map
$\bbP^n(V)\to\bbP^n(K)$ is bijective, for any $n\geq0$.

(b) Notice that $\Pic(V)=0$, for example because any non-trivial
finitely generated $V$-submodule $P\subset K$ is free of rank one:
indeed, if $u_1$, \dots, $u_n\in P$ are the $V$-generators of~$P$,
we just choose $u_k$ with maximal norm $|u_k|$, and notice that
$P$ is already generated by~$u_k$, since $u_i/u_k\in|V|$ for any other~$i$.

(c) Applying~\ptref{sp:proj.spc.points}, we see that elements of
$\bbP^n(V)$ are in one-to-one correspondence with strict quotients $L$
of $V(n+1)$, isomorphic to~$|V|$, and similarly for $\bbP^n(K)$.
So let $\sigma$ be an element of $\bbP^n(V)$, represented by a 
strict quotient $p:V(n+1)\twoheadrightarrow P$, where $P\cong|V|$.
Let $p_K:=p_{(K)}:K(n+1)=K^{n+1}\twoheadrightarrow P_K:=P_{(K)}\cong K$
be its scalar extension to~$K$; it defines the image of $\sigma$ 
in $\bbP^n(K)$. Now consider the following commutative diagram:
\begin{equation}\label{eq:diag.pts.proj.VK}
\xymatrix@C+1pc{
V(n+1)\ar[r]^<>(.5){\xi_{V(n+1)}}\ar[d]^{p}&K^{n+1}\ar[d]^{p_K}\\
P\ar[r]^{\xi_P}&P_K
}\end{equation}
Here $\xi_M:M\to M_{(K)}$ denotes the natural map from a $V$-module~$M$
into its scalar extension to~$K$. Notice that the horizontal arrows
are injective, since $P\cong|V|$ and $V\to K$ is a monomorphism of generalized
rings, and the vertical arrows are strict epimorphisms, i.e.\ surjective.

Therefore, the strict quotient $p:V(n+1)\twoheadrightarrow P$ is
completely determined by $p_K$, since $P$ can be identified with the image
of $p_K\circ\xi_{V(n+1)}$. Furthermore, if we start from any point 
of $\bbP^n(K)$, i.e.\ from a strict epimorphism 
$p_K:K^{n+1}\twoheadrightarrow K$, we can put $P:=\Im p_K\circ\xi_{V(n+1)}=
p_K(V(n+1))$. This is a non-trivial finitely generated $V$-submodule of~$K$,
hence it is free of rank one by~(b), so we get a strict quotient
$V(n+1)\twoheadrightarrow P\cong|V|$, easily seen to coincide with
original $p_K$ after base change to~$K$. In other words, we have just
shown that any point of $\bbP^n(K)$ uniquely lifts to a point of $\bbP^n(V)$,
q.e.d.
\end{Proof}

\begin{CorD}\label{cor:extend.ratpt.sect} 
(a) Let $\sX_\infty/V$ be a pre-unary projective model 
of a projective variety $X/K$ (this means exactly that $\sX_\infty$ 
can be embedded as a ``closed'' subscheme in $\bbP^n_V$,
cf.~\ptref{sp:proj.sch.vs.closed.subsch.projvb}). Then any rational point
$P\in X(K)$ extends to a unique section $\sigma_P:\Spec V\to\sX_\infty$.

(b) Let $\bar\sX/\CompZ$ be a finitely presented 
pre-unary projective model of a projective
variety $X/\bbQ$. Then any point $P\in X(\bbQ)$ extends to a unique
section $\sigma_P:\CompZ\to\bar\sX$.
\end{CorD}
\begin{Proof}
(a) is an immediate consequence of~\ptref{th:val.crit.prop.proj},
and (b) follows from (a), together with classical valuative criterion
of properness, which enables us to extend rational points to sections
of proper models over $\Spec\bbZ$, and~\ptref{sp:ext.ratpt.sect}, q.e.d.
\end{Proof}

Notice that this corollary is applicable in particular to the
finitely presented projective models $\bar\sX/\CompZ$ of
projective varieties $X/\bbQ$ constructed in~\ptref{sp:finpres.mod.projvar}.

\nxpointtoc{$\Zinfty$-models and metrics}
We want to sketch briefly the relationship between models over 
$\Zinfty$ or $\CompZ$ of an algebraic variety $X$ and a vector bundle
$\sE$ on~$X$, on one side,
and metrics on complex points $X(\bbC)$, $\sE(\bbC)$ of this algebraic
variety and vector bundle, on the other side. This will relate
our theory to ``classical'' Arakelov geometry,
based on consideration of metrics.
 
\nxsubpoint (Case of $\Zinfty$-lattices.)
Suppose that $E$ is a $\Zinfty$-lattice, i.e.\ a $\Zinfty$-module,
such that $E\to E_{(\bbR)}$ is injective and has compact image in
finite-dimensional real vector space $E_{(\bbR)}$.  Then the same
reasoning as in part (c) of the proof of~\ptref{th:val.crit.prop.proj}
shows that $\bbP_{\Zinfty}(E) \to\Spec\Zinfty$ satisfies the valuative
criterion of properness with respect to $(\Zinfty,\bbR)$. Indeed,
the only difference is that, in the notation of~\eqref{eq:diag.pts.proj.VK},
where we use $E$ and~$E_{(\bbR)}$ instead of $V(n+1)$ and $K^{n+1}$,
now $P:=p_K(E)=p_\bbR(E)\subset\bbR$ is not a finitely generated 
non-trivial $\Zinfty$-submodule of~$\bbR$; however, it is non-trivial
and compact, $E$ being a $\Zinfty$-lattice, and any non-trivial
compact $\Zinfty$-submodule of $\bbR$ is freely generated by any its
element with maximal absolute value.

Reasoning further as in the proof of~\ptref{th:val.crit.prop.proj}, we
see that, given any ``closed'' subscheme $\sX$ of a projective bundle
$\bbP_{\Zinfty}(E)$ defined by a $\Zinfty$-lattice~$E$, we obtain a
bijection between sections $\sigma:\Spec\Zinfty\to\sX$, and real
points of projective variety $X:=\sX_{(\bbR)}\subset\bbP_\bbR^{\dim
E_{(\bbR)}-1}$.

\nxsubpoint\label{sp:barzinfty.lat.proj.sp}
(Case of $\barZinfty$-lattices.)
Similarly, any non-trivial compact $\barZinfty$-sub\-mo\-dule of $\bbC$ is
freely generated by any its element with maximal absolute value,
so we can extend the above results to ``closed'' subschemes~$\sX$
of $\bbP_{\barZinfty}(E)$ for any $\barZinfty$-lattice~$E$, i.e.\
$\sX(\barZinfty)\to\sX(\bbC)$ is a bijection for any such~$\sX$.

Now consider projective complex algebraic variety $X:=\sX_{(\bbC)}$.
If $\sE$ is a vector bundle over $X$, which can be extended to a vector 
bundle (or maybe just a finitely presented sheaf) $\bar\sE$ on the
whole of $\sX\supset X$ (i.e.\ $\bar\sE$ is a ``model'' of $\sE$),
then we can start from any point $P\in X(\bbC)$, extend it to a section
$\sigma_P:\Spec\barZinfty\to X$, and compute $\sigma_P^*\sE$. This
will be a finitely presented $\sO_{\Spec\barZinfty}$-module,
with the generic fiber canonically isomorphic to $\bbC$-vector space $\sE(P)$.
In other words, we've got a $\barZinfty$-structure on each fiber $\sE(P)$,
i.e.\ essentially a {\em complex Banach norm\/} on $\sE(P)$
(cf.\ Chapter~\ptref{sect:zinfty.lat.mod}). 
This explains the relationship between 
$\barZinfty$-models in our sense, and metrics on vector bundles,
which are used instead of $\barZinfty$-structures in the classical
approach to Arakelov geometry.

\nxsubpoint (Fubini--Study metric.)
Let $Q(z_0,\ldots,z_n):=\sum_i|z_i|^2$ be the standard hermitian form
on $\bbC^{n+1}$, and consider corresponding $\barZinfty$-lattice
$E:=\{z\in\bbC^{n+1}\,:\,Q(z)\leq1\}$. Put $\tilde\bbP^n:=
\bbP_{\barZinfty}(E)$. Then $\tilde\bbP^n$ is a $\barZinfty$-model
of complex projective space $\bbP_\bbC^n$, and according 
to~\ptref{sp:barzinfty.lat.proj.sp}, any point $x\in\bbP^n(\bbC)$
extends to a unique section $\sigma_x\in\Gamma(\tilde\bbP^n/\Spec \barZinfty)$.

Suppose that we are able to construct a non-trivial (co)metric
on $\bbP^n(\bbC)$, by pulling back with respect to $\sigma_x$
some natural model of $\Omega^1_{\bbP^n/\bbZ}$, similarly to
what we did in~\ptref{sp:a1zinfty.infdata}
for $\bbA^1_{\Zinfty}=\Spec\Zinfty[T]$
(cf.\ also~\ptref{sp:proj.vb.zinflat}).
Then the only possibility for this metric is to be the Fubiny--Study
metric, up to multiplication by some constant. Indeed,
this metric must be equivariant under the action of
$\Aut(\tilde\bbP^n/\barZinfty)\supset\Aut_{\barZinfty}(E)=U(n+1)$,
and the only metric (up to a constant) on $\bbP^n(\bbC)$ 
equivariant under the (transitive) action of unitary group $U(n+1)$
is the Fubini--Study metric.

This argument extends to ``closed'' subvarieties $\sX\subset\tilde\bbP^n$:
if we manage to construct a metric on $\sX(\bbC)$, it must coincide
(at least on the smooth locus of $X:=\sX_{(\bbC)}$) with the
metric on $X$, induced from the Fubini--Study metric on~$\bbP^n_\bbC\supset X$.
(It makes sense to construct models of such $X\subset\bbP^n_\bbC$
by taking the ``scheme-theoretic closure'' of $X$ in $\tilde\bbP^n$,
cf.~\ptref{sp:exist.tf.models}.)

After discussing this relationship between our ``algebraic'' approach
to Arakelov geometry, based on generalized schemes and models over
$\CompZ$, and the classical approach, based on introducing suitable
metrics on the complex points of all varieties and vector bundles involved,
we would like to discuss a purely arithmetic application.

\nxpointtoc{Heights of rational points}
\nxsubpoint (Heights of rational points in projective space.)
Let $x\in\bbP^n(\bbQ)$ be a rational point of the $n$-dimensional projective
space~$\bbP^n$. Denote by $(x_0:x_1:\ldots:x_n)$ the homogeneous coordinates
of~$x$. These $x_i$ are rational numbers, not all equal to zero,
so we can define the {\em height\/} $H(x)=H_{\bbP^n}(x)$ of $x$ by
\begin{equation}\label{eq:class.height}
H(x):=\prod_{v\in\bbP\cup\{\infty\}}
\max\bigl({|x_0|}_v,{|x_1|}_v,\ldots,{|x_n|}_v\bigr)
\end{equation}
Here $v$ runs over the set of all valuations of~$\bbQ$, with the
$p$-adic valuation $\fnorm_p$ normalized by the usual requirement 
$|p|_p=p^{-1}$. Notice that this infinite product makes sense,
because $|x_i|_v=1$ for almost all (i.e.\ all but finitely many) 
valuations~$v$. Furthermore, this product doesn't depend on the choice
of homogeneous coordinates of~$x$: if we multiply $(x_0:x_1:\ldots:x_n)$
by any $\lambda\in\bbQ^\times$, the product~\eqref{eq:class.height}
is multiplied by $\prod_v|\lambda|_v$, equal to one by the product
formula~\eqref{eq:product.formula}.

One can describe $H(x)$ in another way. Namely, dividing the $x_i$ 
by their g.c.d., we can assume that the homogeneous coordinates
$x_i$ of point~$x$ are coprime integers. Then all terms 
in~\eqref{eq:class.height} for $v\neq\infty$ vanish, and we obtain
\begin{equation}\label{eq:class.height.coprime}
H(x)=\max\bigl(|x_0|,\ldots,|x_n|\bigr),
\quad\text{if $(x_i)$ are coprime integers.}
\end{equation}
In particular, $H(x)$ is a positive integer for any $x\in\bbP^n(\bbQ)$.

\nxsubpoint (Heights of rational points on a projective variety.)
Now suppose we have a projective variety $X/\bbQ$, embedded as a closed
subvariety into $\bbP^n_\bbQ$. Then we put $H_X(x):=H_{\bbP^n}(x)$
for any $x\in X(\bbQ)\subset\bbP^n(\bbQ)$. 

Notice that this height $H_X(x)$ actually depends not only on
$X/\bbQ$ and $x\in X(\bbQ)$, but also on the choice of closed embedding
$X\to\bbP^n$, or equivalently, on the choice of an ample line bundle
$\sL$ on~$X$ and a finite family of global sections of $\sL$.
However, it is a classical result that two such choices corresponding
to the same ample line bundle give rise to logarithmic height
functions $\log H_X$, $\log H'_X:X(\bbQ)\to\bbR_{\geq0}$ that
differ by a bounded amount.

\nxsubpoint (Logarithmic heights.)
One usually defines {\em logarithmic height\/} functions 
$h_X:X(\bbQ)\to\bbR_{\geq0}$ simply by putting $h_X:=\log H_X$.

\begin{ThD}\label{th:heights}
Let $\sX$ be a finitely presented ``closed'' subscheme of 
$\sP:=\bbP^n_{\CompZ}$.
Denote by $X\subset\bbP^n_\bbQ$ the generic fiber of~$\sX$, a closed
subvariety of projective space $\bbP^n_\bbQ$. Let $j:\sX\to\sP$
be the ``closed'' immersion,
and denote the ample (Serre) line bundle $j^*\sO_\sP(1)$ by $\sO_\sX(1)$.
Let $p:\sX\to\CompZ$ be the structural morphism. Then:
\begin{itemize}
\item[(i)] Any rational point $x\in X(\bbQ)\subset\bbP^n(\bbQ)$ extends
in a unique way to a section $\sigma_x$ of $\sX$ over~$\CompZ$.
\item[(ii)] Degree of line bundle $\sL:=\sigma_x^*\sO_\sX(1)$ over $\CompZ$
equals $h_X(x)=\log H_X(x)=\log H_{\bbP^n}(x)$, the logarithmic height
of~$x$:
\begin{equation}
h_X(x)=\deg\sigma_x^*\sO_\sX(1)
\end{equation}
\end{itemize}
\end{ThD}
\begin{Proof}
(i) has been already shown in~\ptref{cor:extend.ratpt.sect},
so we have to prove~(ii). Notice that $j\sigma_x$ is the 
section of $\sP/\CompZ$ extending $x\in X(\bbQ)\subset\bbP^n(\bbQ)$,
and $\sigma_x^*\sO_\sX(1)=\sigma_x^*j^*\sO_\sP(1)=(j\sigma_x)^*\sO_\sP(1)$.
Therefore, we may assume $\sX=\sP=\bbP^n_{\CompZ}$ and $X=\bbP^n_\bbQ$.

(a) Recall that sections of $\sP:=\bbP^n_{\CompZ}$ over $\CompZ$ are
in one-to-one correspondence with strict quotients~$\sL$
of trivial vector bundle $\sE:=\sO_{\CompZ}^{(n+1)}$, which are
line bundles. Furthermore, this correspondence transforms a 
section $\sigma:\CompZ\to\sP$ into $\sigma^*(\pi)$, where
$\pi:p^*\sE\twoheadrightarrow\sO_\sP(1)$ is the canonical surjection
over~$\sP$
(cf.~\ptref{sp:sect.proj.bdl}). In particular, section~$\sigma_x$ 
corresponds to strict epimorphism $\phi:=\sigma_x^*(\pi):\sE=\sigma_x^*p^*\sE
\twoheadrightarrow\sigma_x^*\sO_\sP(1)=\sL$, i.e.\ this strict quotient
of $\sE$ is actually the line bundle~$\sL$ from the statement of the
theorem.

(b) Let us study this strict epimorphism $\phi:\sE=\sO_{\CompZ}^{(n+1)}
\twoheadrightarrow\sL$ in more detail. Put $L_\bbZ:=\Gamma(\Spec\bbZ,\sL)$,
$L_\infty:=\sL_\infty$, $L_\bbQ:=\sL_\xi$, and define $E_\bbZ$, \dots,
$\phi_\bbQ:E_\bbQ\to L_\bbQ$ similarly. We obtain the following 
commutative diagram:
\begin{equation}
\xymatrix@C+10pt{
E_\bbZ\ar[d]^{\phi_\bbZ}\ar[r]^{u_E}&E_\bbQ\ar[d]^{\phi_\bbQ}&
E_\infty\ar[d]^{\phi_\infty}\ar[l]_{v_E}\\
L_\bbZ\ar[r]^{u_L}&L_\bbQ&L_\infty\ar[l]_{v_L}}
\end{equation}
All horizontal arrows are natural embeddings of $\bbZ$- or $\Zninfty$-modules
into their scalar extensions to~$\bbQ$. Since all modules involved are free,
and both $\bbZ\to\bbQ$ and $\Zninfty\to\bbQ$ are monomorphisms of 
generalized rings, all horizontal arrows are injective. On the other hand,
the vertical arrows are strict epimorphisms, i.e.\ surjective maps.

(c) Notice that $E_\bbQ=\bbQ^{n+1}$, and $L_\bbQ\cong\bbQ$. Fix any
such isomorphism $L_\bbQ\cong\bbQ$; then the images $x_i:=\phi_L(e_i)$
of the standard basis elements $e_i\in\bbQ^{n+1}$, $0\leq i\leq n$,
are exactly the homogeneous coordinates $(x_0:\ldots:x_n)$ of our
original point $x\in\bbP^n(\bbQ)$, just by definition of homogeneous
coordinates (cf.~\ptref{sp:proj.spc.points}). We can always fix
isomorphism $L_\bbQ\cong\bbQ$ so as to make all $x_i$ coprime integers.

(d) Now identify $E_\bbZ$ and $E_\infty$ with their images
$\bbZ^n$ and $\Zninfty(n)$ in $E_\bbQ=\bbQ^n$, and similarly
$L_\bbZ$ and $L_\infty$ with their images in~$L_\bbQ$. Furthermore,
identify $L_\bbQ=\bbQ$ as in~(c), so as to make the homogeneous coordinates
$x_i$ coprime integers. Then $\phi_\bbQ:\bbQ^n\to\bbQ$ is the map
$e_i\mapsto x_i$, $L_\bbZ\subset\bbQ$ is identified with
$\phi_\bbQ(\bbZ^n)$, and $L_\infty\subset\bbQ$ is identified with
$\phi_\bbQ(\Zninfty(n))$.

(e) Notice that $L_\bbZ\subset\bbQ$ is the $\bbZ$-submodule of~$\bbQ$, 
generated by coordinates $x_i$. Since these coordinates are coprime
integers, we get $L_\bbZ=\bbZ$. Similarly, $L_\infty\subset\bbQ$
is the $\Zninfty$-submodule of~$\bbQ$, generated by $x_0$, \dots, $x_n$.
Recall that any non-trivial finitely generated $\Zninfty$-submodule of $\bbQ$
is freely generated by its generator of largest absolute value
(cf.\ part~(c) of the proof of~\ptref{th:val.crit.prop.proj}),
hence $L_\infty=\lambda\cdot|\Zninfty|$ for 
$\lambda:=\max(|x_0|,\ldots,|x_n|)$.

(f) According to the description of $\deg\sL$ given 
in~\ptref{sp:deg.linbdl.compz},
we have to choose a free generator $f_1$ of $L_\bbZ$, a free generator
$f_2$ of $L_\infty$, compute $f_2/f_1$ in the generic fiber $L_\bbQ$,
and then $\deg\sL=\log(f_2/f_1)$. In our case $f_1=1$, $f_2=\lambda$,
hence $\deg\sL=\log\lambda$. On the other hand, $\lambda=H_{\bbP^n}(x)$
by~\eqref{eq:class.height.coprime}, hence $\deg\sigma_x^*\sO_\sX(1)=
\deg\sL=\log\lambda=\log H_X(x)=h_X(x)$, q.e.d.
\end{Proof}

\nxsubpoint (Comparison to classical Arakelov geometry.)
Classical Arakelov geometry has a similar formula relating height $\hat h(x)$
of a point $x\in\bbP^n(\bbQ)$ with arithmetic degree of the pullback
of corresponding ample (Serre) line bundle $\sO(1)$ on $\bbP^n(\bbQ)$,
equipped with its natural (``Fubini--Study'') metric. However,
the logarithmic height $\hat h(x)$ thus obtained is defined by formula
\begin{equation}
\hat h(x)=\log\bigl({|x_0|}^2+\cdots+{|x_n|}^2\bigr)^{1/2},
\quad\text{if $(x_i)$ are coprime integers.}
\end{equation}
On the other hand, we obtain in~\ptref{th:heights} 
the ``classical'' height $h(x)$ of~\eqref{eq:class.height.coprime},
used everywhere else in arithmetic geometry. Of course, this distinction
is not very significant, since $\hat h(x)$ and $h(x)$ differ
by a bounded amount, but it is quite amusing nonetheless.




\cleardoublepage

\mysection{Homological and homotopical algebra}
\label{sect:homot.alg}

The aim of this chapter is to develop a variant of homological algebra for 
the categories of modules over generalized rings. This achievement
will be extended in the next chapter to the categories 
of (sheaves of) modules over a generalized scheme or even a generalized 
ringed topos. Such homological (or rather homotopical) algebra
will be intensively used in the last chapter of this work to construct 
a reasonable arithmetic intersection theory in the framework of generalized 
schemes, suitable e.g.\ for a discussion of arithmetic Riemann--Roch theorem.

\zeropoint\nxpoint\label{p:great.plan} 
(Tale of the Great Plan.) Let us present a brief 
sketch of the layout of the remainder of this chapter. We also discuss
some of ideas that will be used in the following two chapters.
Thus one might consider this sketch as a common introduction to the
three ``homotopical'' chapters of this work. 

\nxsubpoint (Homotopical algebra and model categories.)
It is quite clear that classical homological algebra cannot work properly 
in a non-additive setup, and that one should use homotopical algebra 
instead, based on simplicial objects, model and homotopic categories, and 
weak equivalences,  
rather than on chain complexes, derived categories and quasi-isomorphisms. 
Therefore, we recall the basic definitions and constructions 
of homotopical algebra before using them. Our basic reference here is the 
foundational work of Quillen (\cite{Quillen}) that seems to be quite well 
suited for our purpose. For example, the results of 
{\em loc.cit.}\ enable us to construct a simplicial closed model category 
structure on $s(\catMod{\Sigma})$ almost immediately.
However, we sometimes insist on using 
a more modern terminology (e.g.\ we say ``acyclic fibrations'' instead of 
Quillen's ``trivial fibrations'').

\nxsubpoint 
(Comparison with the classical case: Dold--Kan and Eilenberg--Zilber theorems.)
Of course, we would like to compare our constructions and definitions with 
the classical ones when the generalized rings under consideration turn out 
to be additive. The key tool here is the {\em Dold--Kan correspondence\/} 
between the category~$s\cA$ of simplicial objects over 
an abelian category~$\cA$ and the category $\Ch(\cA)=\Ch_{\geq0}(\cA)$ 
of non-negative chain complexes over~$\cA$. 
More precisely, these two categories turn out 
to be equivalent, and chain homotopies in $\Ch(\cA)$ correspond exactly 
to the simplicial homotopies in $s\cA$.

We'd like to mention here the way the equivalence $K:\Ch(\cA)\to s\cA$ acts on 
objects:
\begin{equation}\label{eq:intr.simpl.from.cochain}
(KC_\cdot)_n=\bigoplus_{\eta:[n]\twoheadrightarrow [p]}C_p
\end{equation}
Here $[n]=\{0,1,\ldots,n\}$, and $\eta$ runs over all monotone surjective maps 
from $[n]$ onto~$[p]$. The quasi-inverse to $K$ is the {\em normalization 
functor~$N:s\cA\to\Ch(\cA)$}; it admits an explicit description as well.

Another important tool here, especially useful when dealing with 
bisimplical objects and chain bicomplexes over an abelian category, 
is the {\em Eilenberg--Zilber theorem\/} that asserts the following. 
Given a bisimplicial object $K_{..}$ over an 
abelian category~$\cA$, we can construct a chain bicomplex 
$N_IN_{II}K_{\cdot\cdot}$ and consider its totalization 
$\Tot(N_IN_{II}K_{\cdot\cdot})$; of course, 
this totalization corresponds to a certain simplicial object $\hat K_\cdot$. 
Now the Eilenberg--Zilber theorem asserts that $\hat K_\cdot$ is actually 
homotopy equivalent to the diagonal simplicial object $\tilde K_\cdot$ 
of~$K_{..}$, defined simply by $\tilde K_n:=K_{nn}$. Therefore, 
the diagonal simplicial object of a bisimplicial object is something like the 
totalization of a bicomplex, and we exploit this observation in the 
non-additive case as well.

\nxsubpoint (Model category structure on $s(\catMod\Sigma)$.)
For any algebraic monad~$\Sigma$ (i.e.\ a non-commutative generalized ring) 
we can obtain a simplicial closed model category structure on 
$s(\catMod\Sigma)$, by a direct application of a theorem of~\cite{Quillen}.
We check that when $\Sigma$ is additive (i.e.\ $\catMod\Sigma$ is abelian), 
$\Sigma=\Fempty$ or $\Fone$, we recover the well-known model category 
structures on $s\cA$ ($\cA$ an abelian category), simplicial sets, and 
simplicial sets with a marked point. In particular, for an additive $\Sigma$ 
the weak equivalences $f:X_\cdot\to Y_\cdot$ correspond via the Dold--Kan 
correspondence to the quasi-isomorphisms of non-negative chain complexes.
Therefore, the corresponding homotopic category $\Ho s(\catMod\Sigma)$, 
obtained by localizing $s(\catMod\Sigma)$ with respect to all 
weak equivalences, is a natural replacement for the subcategory 
$\cD^{\leq0}(\catMod\Sigma)$ of the derived category $\cD^-(\catMod\Sigma)$; 
if necessary, we can recover an analogue of $\cD^-(\catMod\Sigma)$ itself by 
formally inverting the translation functor in $\cD^{\leq0}$. We put 
$\cD^{\leq0}(\Sigma):=\Ho s(\catMod\Sigma)$.

Moreover, for any homomorphism of algebraic monads $\phi:\Sigma\to\Sigma'$ 
we show that the arising adjoint pair of functors $\phi^*$, $\phi_*$ 
between $s(\catMod\Sigma)$ and $s(\catMod{\Sigma'})$ constitute a Quillen 
functor. Therefore, we obtain the corresponding derived functors 
$\dL\phi^*:\cD^{\leq0}(\Sigma)\to\cD^{\leq0}(\Sigma')$ and 
$\dR\phi_*$ in the opposite direction, adjoint to each other. 
When both $\Sigma$ and $\Sigma'$ are additive, these derived functors 
reduce via the Dold--Kan correspondence to their classical counterparts.

\nxsubpoint (Derived tensor products, symmetric and exterior powers.)
Of course, we can derive the tensor product $\otimes_\Sigma$, thus 
obtaining a derived tensor product $\Lotimes:\cD^{\leq0}(\Sigma)\times
\cD^{\leq0}(\Sigma)\to\cD^{\leq0}(\Sigma)$. By Eilenberg--Zilber this 
derived tensor product turns out to coincide with the usual one for an 
additive~$\Sigma$. It has almost all properties one expects from it, 
e.g.\ commutativity, associativity, compatibility with $\dL\phi^*$, 
with the notable exception of the projection formula for a non-unary~$\phi$.

Moreover, when $\Sigma$ is alternating, we can try derive the exterior 
powers, thus obtaining some functors $\dL\bigwedge^r$, essentially in 
the same way Dold and Puppe constructed their homology of non-additive 
functors. However, it turns out to be technically more convenient to
derive the {\em symmetric powers\/} instead, similarly to what Dold
and Puppe actually did (notice, however, that they had a simpler situation
since they had to derive
{\em non-additive functors\/} between {\em additive categories}). 
Among other things, this
makes sense for any generalized rings, not just for alternating
$\Fpm$-algebras.

We show that these derived symmetric powers $\dL S^r$ satisfy at least 
some of their classical properties. These properties suffice to
introduce a pre-$\lambda$-ring structure on $K_0$ of
perfect complexes of modules later in Chapter~\ptref{chap:K0.Chow}.

\nxsubpoint (Global situation.)
Of course, we'd like to generalize the above constructions to the case 
of sheaves (quasicoherent or not) of $\sO_X$-modules on a generalized 
scheme~$X$, or even on a generalized ringed topos $(X,\sO_X)$. 
The obvious problem is that we don't have any reasonable model category 
structure on $s(\catMod{\sO_X})$. This is due to the fact we don't have 
enough projectives in these categories. However, we can still define the 
set of weak equivalences in such categories: for example, if $(X,\sO_X)$ 
has enough points (e.g.\ is a topological space), then $f:K_\cdot\to L_\cdot$ 
is a weak equivalence iff $f_p=p^*f$ is a weak equivalence 
in $s(\catMod{\sO_{X,p}})$ for any point $p$ of~$X$. 
For simplicial quasicoherent sheaves on a generalized scheme~$X$ 
it suffices to require $\Gamma(U,f)$ to be a weak equivalence of simplicial 
$\Gamma(U,\sO_X)$-modules for all $U$ from some affine cover of~$X$. 
In the general case we have to consider the pro-object of all hypercoverings 
of~$X$ instead; this presents some technical but not conceptual complications.

In any case, once we have a class of weak equivalences, we can define 
$\cD^{\leq0}(X,\sO_X)$ to be the localization of $s(\catMod{\sO_X})$ 
with respect to this class of morphisms; the only complication is to prove 
that we still get an $\univU$-category, i.e.\ that all sets of morphisms in 
$\cD^{\leq0}$ remain small. This might be achieved by considering only 
hypercoverings composed from objects of a small generating family 
for topos~$X$.

After that we can construct a derived functor $\dL f^*:\cD^{\leq0}(Y,\sO_Y)\to
\cD^{\leq0}(X,\sO_X)$ for any morphism of generalized ringed topoi 
$f:(X,\sO_X)\to(Y,\sO_Y)$, as well as the derived tensor product on these 
categories. For these constructions we need to consider the categories of 
presheaves as well: if $X=\tilde\cS$ for some small site~$\cS$, then 
the category of presheaves of $\sO_X$-modules on~$\cS$ does have a closed 
model category structure. Notice that we cannot obtain any $\dR f_*$ in this 
setup, so the projection formula cannot even be written down.

\nxsubpoint (Perfect complexes.)
Consider the case of a classical ring~$\Sigma$ first, so that we have the 
Dold--Kan correspondence. One checks directly that a chain complex $C_\cdot$ 
consists of finitely generated projective $\Sigma$-modules iff the same holds 
for the corresponding simplicial object $(KC_\cdot)_\cdot$, and that $C_\cdot$ 
has bounded homology iff the natural map $\sk_N KC_\cdot\to KC_\cdot$ is a 
weak equivalence for $N\gg0$. Therefore, both these conditions can be 
stated on the level of simplicial objects, and we can define the category of 
perfect complexes of $\Sigma$-modules (or rather perfect simplicial 
$\Sigma$-modules) for any algebraic monad~$\Sigma$. We can extend these 
definitions to the case of a generalized scheme~$X$ as well. 
The only potential problem here is whether we want to consider 
the locally projective or just the locally free $\sO_X$-modules, since 
in general we don't expect these two classes to coincide.

\nxsubpoint ($K$-theory.)
After that we are ready to construct $K_0(\Perf X)$. Of course, this will 
be the abelian group generated by weak equivalence classes 
of perfect simplicial $\sO_X$-modules modulo some relations. 
Clearly, we must require $[A_\cdot\oplus B_\cdot]=[A_\cdot]+[B_\cdot]$. 
A less trivial relation is this: 
if $A_\cdot$ and $B_\cdot$ are two perfect complexes, such that 
$A_n$ is isomorphic to $B_n$ for each $n\geq0$ (no compatibility with 
morphisms $A_n\to A_m$ required!), then $[A_\cdot]=[B_\cdot]$. 
These relations turn out to be sufficient in the additive case, 
so we use them in the general case as well. 

The derived tensor product defines a product on $K^\bullet(X):=K_0(\Perf X)$, 
and $K^\bullet(X)$ depends contravariantly on~$X$. 
After that we use derived symmetric powers $\dL S^r$ 
to define a pre-$\lambda$-ring structure on $K^\bullet(X)$.
We would like to check that this pre-$\lambda$-ring is actually 
a $\lambda$-ring. However, we didn't manage to prove all necessary relations
for our derived symmetric powers so far, so we replace $K^\bullet(X)$
with its largest $\lambda$-ring quotient $K^\bullet(X)_\lambda$ instead.

This means that we can define on $K^\bullet(X)_\lambda$ 
the $\gamma$-operations and the $\gamma$-filtration in the way known since 
Grothendieck, and define the Chow ring of~$X$ by $\Ch^\cdot(X):=\gr_\gamma
K^\bullet(X)_{\lambda,\bbQ}$. 
Grothendieck's argument gives us immediately a theory of 
Chern classes with values in $\Ch^\cdot(X)$ enjoying all the usual formal 
properties. An alternative approach, due to Soul\'e, 
is to define $\Ch^i(X)$ as the weight~$i$ part 
of $K^\bullet(X)_{\lambda,\bbQ}$ with respect to Adams operations.
We show that these two approaches yield the same result whenever the
$\gamma$-filtration is finite and separated.

\nxsubpoint (Geometric properties: projective bundles, regular embeddings 
and Gysin maps.) 
Of course, we'd like to prove the usual geometric properties of 
$\Ch^\cdot(X)$ and our Chern classes, e.g.\ the projective bundle theorem, 
existence of Gysin maps and so on. Such considerations are important for 
proving any sort of Riemann--Roch theorem; however, we postpone these 
``geometric'' considerations until a subsequent work.

\nxsubpoint (Covariant data and cosimplicial objects.)  Up to this
point we've been constructing only contravariant data, similar to the
cohomology ring, e.g.\ the perfect complexes, Chern classes, Chow
rings and so on. However, we'd like to have something covariant as
well, similar to the homology groups, e.g.\ the complexes with bounded
coherent cohomology, $K_\cdot(X):=K_0(\Coh X)$ etc. It is natural to
require these covariant objects to have an action of corresponding
contravariant objects, similar to the $\cap$-product action of
$H^*(X)$ on $H_*(X)$ or to the action of $K^\cdot(X)$ on $K_\cdot(X)$,
given by the derived tensor product in the classical case. Among other
things, this would enable us to write down a projection formula.

We see that we need something like cochain complexes. In the additive case 
the dual of Dold--Kan correspondence assures us that the category of 
non-negative cochain complexes over an abelian category~$\cA$ 
is equivalent to the category $c\cA$ of {\em cosimplicial\/} objects 
over~$\cA$. So we are tempted to use $c(\catMod\Sigma)$ and $c(\catMod{\sO_X})$
in the general case as well. Indeed, the action of $s(\catMod\Sigma)$ on 
$c(\catMod\Sigma)$ would be given by the derived inner Hom $\dR\iHom_\Sigma$, 
and similarly on any generalized ringed topos~$(X,\sO_X)$, and we can 
even hope to obtain a projection formula just from the adjointness of 
$\otimes$ and $\iHom$, and $f^*$ and $f_*$, {\em provided we are able to 
define a model category structure on~$c\cA$ and construct the derived 
direct image $\dR f_*$.}

Here is the problem. We cannot construct a reasonable model category structure 
on $c(\catMod\Sigma)$ or $c(\catMod{\sO_X})$ since these categories don't 
have enough injectives, and cannot derive $f_*$ for the same reason.
However, something can still be done. First of all, we can try to derive 
$\iHom:s(\catMod\Sigma)^{op}\times c(\catMod\Sigma)\to c(\catMod\Sigma)$ 
using only the model category structure on the first argument. Next, 
in the additive case
we can compute $\dR\Gamma$ and $\dR f_*$ by means of Verdier theorem 
(cf.~\cite{SGA4}, V) that asserts that all cohomology groups can be 
computed by something like \v Cech cohomology provided we consider all 
hypercoverings (which are simplicial objects themselves), not just the 
coverings, so we can mimic this construction in the non-additive case as well.

\nxsubpoint (Homotopy groups.)
Since any $A_\cdot$ in $s(\catMod\Sigma)$ can be considered as 
a simplicial set, we can compute its homotopy invariants $\pi_n(A_\cdot)=
\pi_n(\dR\Gamma_\Sigma A_\cdot)$, $n\geq0$, at least if we are given a 
base point in~$A_0$, e.g.\ if $\Sigma$ has a zero constant.
Of course, $\pi_0(A)$ is just a (pointed) set, 
$\pi_1(A)$ is a group, and $\pi_n(A)$ is an abelian group for $n\geq2$. 
Moreover, we get a $\Sigma$-structure on each of these sets, compatible with 
the group structure for $n\geq1$.
Recall that the category of abelian groups is the category of 
modules over~$\bbZ$, and the category of groups is defined by some 
non-commutative algebraic monad~$\bbG$ as well. All this means that 
$\pi_0(A)$ belongs to $\catMod{\Sigma}$, $\pi_1(A)$ to 
$\catMod{(\Sigma\otimes\bbG)}$, and 
$\pi_n(A)$ to $\catMod{(\Sigma\otimes\bbZ)}$. 
For example, if $\Sigma$ is additive, then $\Sigma\otimes\bbG\cong
\Sigma\otimes\bbZ\cong\Sigma$, and the $\Sigma$-modules $\pi_n(A)$ are 
canonically isomorphic to the homology groups $H_n(NA)$ of corresponding 
chain complex $NA$. Another example: if $\Sigma=\Zinfty$, then 
$\pi_0(A)$ is just a $\Zinfty$-module, $\pi_1(A)$ is a 
$\Zinfty\otimes\bbG$-module, and all higher $\pi_n(A)$, 
$n\geq2$, are $\bbR$-vector spaces.

Notice that the model category structure on $s(\catMod\Sigma)$ is defined in 
such a way that $f:A\to B$ is a weak equivalence iff all $\pi_n(f)$ are 
isomorphisms (more precisely, $\pi_n(f,x):\pi_n(A,x)\to\pi_n(B,f(x))$ for 
all $n\geq0$ and all choices of base point $x\in A_0$). However, two 
simplicial $\Sigma$-modules with isomorphic homotopy groups need not be 
weakly isomorphic or isomorphic in the homotopic category; such 
examples are very well known even in the additive case, e.g.\ $\Sigma=\bbZ$. 
In our case the homotopy groups carry even less information compared to the 
corresponding object of the homotopic category. For example, the 
``cones'' of $\Zinfty\stackrel{1/2}\to\Zinfty$ and 
$\Zinfty\stackrel{1/3}\to\Zinfty$ in $\cD^{\leq0}(\Zinfty)$ turn out to have 
the same homotopy groups ($\pi_n=0$ for $n>0$, $\pi_0=\Zinfty/\gm_\infty$), 
but they don't seem to be isomorphic in $\cD^{\leq0}(\Zinfty)$\dots\

All this means that we shouldn't expect to extract anything like 
the Euler characteristic of a simplicial $\Sigma$-module just from its 
homotopy groups.

\nxsubpoint\label{sp:Euler.char} (Euler characteristics.)
Of course, we'd like to be able to define and compute Euler characteristics 
of the (co)simplicial objects we construct, since this is  
the natural way to obtain arithmetic intersection numbers in our setup. 
Let's consider the case, say, of a perfect simplicial object~$A_\cdot$ 
over a generalized ring~$\Sigma$. 
Let $\ell$ be an additive function from the set of isomorphism 
classes of finitely generated projective modules over~$\Sigma$ into, say, 
a $\bbQ$-vector space~$V$, additivity meaning here just $\ell(P\oplus Q)=
\ell(P)+\ell(Q)$.

If $\Sigma$ is additive, we can consider the corresponding normalized chain 
complex $C_\cdot:=NA_\cdot$, and then $A.=KC_\cdot$ can be recovered according 
to formula~\eqref{eq:intr.simpl.from.cochain}; by additivity we get 
\begin{equation}
\ell(A_n)=\sum_{p=0}^n\binom np\cdot\ell(C_p)
\end{equation}
We can transform this identity formally to obtain
\begin{equation}
\ell(C_n)=\sum_{p=0}^n(-1)^{n-p}\binom np\cdot\ell(A_p)
\end{equation}
Now nothing prevents us from defining $\ell(C_n)$ by these formulas in the 
non-additive case and putting $\chi(A):=\sum_{n\geq0}(-1)^n\ell(C_n)$, 
provided we are able to prove that $\ell(C_n)=0$ for $n\gg0$. Of course, 
in order to obtain a reasonable Euler characteristic $\chi:K^\bullet(\Sigma)\to
V$ we'll have to check that $\chi(A)=\chi(A')$ for any weak equivalence 
$A\to A'$ as well.

We can consider the Euler--Poincar\'e series of~$A$ and~$C$, defined by the 
usual formula $P_t(A):=\sum_{n\geq0}\ell(A_n)t^n$, and similarly for~$C$.
We get $P_t(A)=P_{t/(1-t)}(C)$ and $P_u(C)=P_{u/(u+1)}(A)$, hence 
$\chi(A)=P_{-1}(C)=P_\infty(A)$. This seems absurd, but in our situation 
we expect $P_u(C)$ to be a polynomial, hence $P_t(A)$ to be a rational 
function, so we can indeed evaluate it at infinity. 
We see that some sort of regularization process is already being used here.

Another approach consists in considering the Hilbert polynomial $H_A(t)$ 
of~$A$, defined by the requirement $\ell(A_n)=H_A(n)$ for all integer 
$n\geq 0$. In our case we have $H_A(x)=\sum_{p=0}^N\ell(C_p)\binom xp$, 
so this is indeed a polynomial. On the other hand, $\binom{-1}p=(-1)^p$ for 
any $p\geq0$, so we get $\chi_A=H_A(-1)$, i.e.\ $\chi(A)$ is something like 
``$\ell(A_{-1})$'', and the Euler characteristic computation turns out to 
be a certain extrapolation problem.

The Hilbert polynomial point of view is also extremely useful. Suppose 
for example that we have a polynomial map $Q_r:V\to V$, such that 
$\ell(\wedge^rP)=Q_r(\ell(P))$ for any projective~$P$ of finite type. Then 
we get $H_{\lambda^rA}(n)=Q_r(H_A(n))$ for all $n\geq0$, hence 
$H_{\lambda^rA}(t)=Q_r(H_A(t))$ and $\chi(\lambda^rA)=Q_r(\chi(A))$.

\nxsubpoint 
(Intersection numbers and Euler characteristics.)
We'd like to be able to compute some intersection numbers, say, for 
any $a\in\Ch^d(X)$, $d=\dim X$. In order to do this we need a dimension 
theory on~$X$, so as to be able to define~$d$, as well as something like 
the ``fundamental class of variety~$X$'', $[X]\in H_{2d}(X)$; then the 
number we want to compute would be given by the $\cap$-product.

In our situation $a$ can be represented by a perfect simplicial object 
$A_\cdot$. 
And the most natural way to fix a dimension theory, at least for 
classical noetherian schemes, is to fix a dualizing (cochain) 
complex~$I^\cdot$, 
i.e.\ some cosimplicial object. Notice that in the classical 
case $I^\cdot$ provides both a dimension theory, since for any $x\in X$ 
we can define $d(x)$ to be the only index $i$, for which 
$\Ext^i(\kappa(x),I^\cdot)\neq0$, 
and something like a representative of the fundamental class, 
since $Z^\cdot:=\dR\Hom(A_\cdot,I^\cdot)$ is a complex (or a cosimplicial set),
trivial in degrees $\neq d=\dim X$, and having $\ell(H^d(Z))=\pm\chi(Z^\cdot)$ 
equal to the intersection number we are looking for.

It seems quite reasonable to mimic this construction in the non-additive 
case as well: fix some ``dualizing cosimplicial object'' $I^\cdot$ in 
$c(\catMod{\sO_X})$, and compute the cosimplicial set 
$Z^\cdot:=\dR\negthinspace\Hom(A_\cdot,I^\cdot)=\dR\Gamma\,\dR\negthinspace
\iHom(A_\cdot,I^\cdot)$ for any $A_\cdot$ from 
a suitable step of the $\gamma$-filtration of $K^\bullet(X)$. Then, if we 
are lucky, $\pm\chi(Z^\cdot)$ will give us the intersection number 
we are looking for. Of course, we have to choose $I^\cdot$ in a reasonable way;
if we think about $X$ as a smooth variety, then $\sO_X$ concentrated 
in degree~$0$, i.e.\ constant cosimplicial object with value~$\sO_X$, would do.

We refer to the above process as {\em the arithmetic integration}. 
Maybe we should even denote the resulting intersection number 
by~$\int_{I^\cdot}A_\cdot$.

\nxsubpoint (Euler characteristics of cosimplicial sets.)
First of all, we have to introduce a simplicial closed model category 
structure 
on the category $c\catSets$ of cosimplicial sets. We might do this with the 
aid of the dual of Quillen theorem already used to obtain a model category 
structure on $s(\catMod\Sigma)$, since $\catSets$ has sufficiently many 
injectives (something we don't have in an arbitrary $\catMod\Sigma$), 
and any set can be embedded into a group (e.g.\ the corresponding free group).
After that, given a cosimplicial set $Z_\cdot$, we can try to find a finite 
fibrant-cofibrant replacement $A_\cdot$ of~$Z_\cdot$, using the model category 
structure on~$c\catSets$. Since we work in the dual situation, in order to 
compute Euler characteristics we'll need a map $\ell$ from the category of 
finite sets into a vector space~$V$, such that 
$\ell(X\times Y)=\ell(X)+\ell(Y)$. The natural choice would be 
$V:=\bbR$, $\ell(X):=\log|X|$ (another option: 
$V=\bbQ\otimes_\bbZ\bbQ_{>0}^\times$). Then we can compute 
the Euler characteristic $\chi(Z_\cdot)=\chi(A_\cdot)$ by the procedure 
sketched above; it'll be a logarithm of a positive rational number. 
This definition makes sense at least for algebraic varieties $X$ over 
finite fields~$\bbF_q$: the arithmetic intersection numbers we obtain 
turn out to be equal to corresponding geometric intersection numbers 
multiplied by $\log q$.

\nxsubpoint (Mystery of good intersection numbers.)
It seems a bit strange that we obtain only logarithms of rational numbers 
in this way. However, this is not so strange once we recall that 
we cannot choose metrics on our arithmetic varieties and 
vector bundles arbitrarily. In fact, we have already seen 
in Chapter~\ptref{sect:zinfty.lat.mod} that 
arithmetic varieties of finite presentation over $\CompZ$ 
tend to have quite exotic 
singular (co)metrics, and, while more classical metrics, like the restriction 
of the Fubiny--Study metric on a closed subvariety of a projective space, 
can be described in our approach, 
they are not of finite presentation over~$\CompZ$, 
so we'll either need some regularization process while computing our 
Euler characteristics, or we'll have to represent the given metric as a 
limit of ``finitely presented metrics'' and compute the limit of intersection 
numbers. In any case we cannot expect the resulting number 
to be the logarithm of a rational number. 
This is something like inscribing convex polyhedra 
with rational vertices into the unit sphere and taking 
the limit of their volumes: their volumes will be rational, while the limit 
will be not.

\nxpointtoc{Model categories}
We would like to recall here the basic definitions and constructions of 
homotopic algebra that we are going to use later. Our basic reference here 
is \cite{Quillen}; however, our terminology is sometimes a bit different 
from that of Quillen. For example, all model categories we'll consider 
will be closed, so the axioms we give below for a model category 
(cf.\ e.g.\ \cite{Hovey}) are actually equivalent to 
Quillen's axioms for a closed model category.

\begin{DefD}\label{def:mod.cat}
A {\em (closed) model category~$\cC$} is a category $\cC$, 
together with three distinguinshed classes of morphisms, called 
{\bf fibrations}, {\bf cofibrations} and {\bf weak equivalences}, subject 
to conditions (CM1)--(CM5) listed below.
\end{DefD}

\nxsubpoint\label{sp:acyc.retr.liftp}
 (Acyclic (co)fibrations, retracts and lifting properties.)
Before listing the axioms of a model category, let's introduce some 
terminology. First of all, we say that a morphism $f:X\to Y$ in 
a category~$\cC$ with three distinguished classes of morphisms as above 
is an {\em acyclic fibration\/} (resp.\ {\em acyclic cofibration}) 
if it is both a fibration (resp.\ cofibration) and a weak equivalence.

Secondly, given two morphisms $i:A\to B$ and $f:X\to Y$ in an arbitrary 
category~$\cC$, we say that 
$i$ has the {\em left lifting property (LLP) with respect to~$f$}, or that 
$f$ has the {\em right lifting property (RLP) with respect to~$i$}, if 
for any two morphisms $u:A\to X$ and $v:B\to Y$, such that $v\circ i=f\circ u$,
one can find a morphism $h:B\to X$, such that $h\circ i=u$ and $f\circ h=v$:
\begin{equation}\label{eq:diag.lift}
\xymatrix{
A\ar[d]^{i}\ar@{-->}[r]^{u}&X\ar[d]^{f}\\
B\ar@{-->}[r]^{v}\ar@{.>}[ur]|{\exists h}&Y}
\end{equation}

Finally, we say that a morphism $g:Z\to T$ is a {\em retract\/} of 
another morphism $f:X\to Y$ iff there exist morphisms 
$(i,j):g\to f$ and $(p,q):f\to g$, such that $(p,q)\circ(i,j)=\id$, 
i.e.\ if there are morphisms $i:Z\to X$, $j:T\to Y$, $p:X\to Z$, $q:Y\to T$, 
such that $f\circ i=j\circ g$, $g\circ p=q\circ f$, $p\circ i=\id_Z$ and 
$q\circ j=\id_T$:
\begin{equation}\label{eq:diag.retr}
\xymatrix{
Z\ar[d]^{g}\ar[r]_{i}\ar@/^1pc/[rr]^{\id_Z}&X\ar[d]^{f}\ar[r]_{p}&Z\ar[d]^{g}\\
T\ar[r]^{j}\ar@/_1pc/[rr]_{\id_T}&Y\ar[r]^{q}&T}
\end{equation}
These notions have some obvious properties, e.g.\ if $i'$ is a retract of $i$, 
and $i$ has the LLP with respect to some~$f$, then $i'$ has the same LLP  
as well. Another simple property: 
the class of morphisms having the RLP with respect to some fixed morphism~$i$ 
is stable under base change and composition.

\nxsubpoint (Axioms of a closed model category.)

{\myindent{(CM1)} 
The category~$\cC$ is closed under inductive and 
projective limits. (Finite inductive and projective limits actually 
suffice for most applications.)
\smallskip

\myindent{(CM2)} 
Each of the three distinguished classes of morphisms in~$\cC$ 
is closed under retracts.
\smallskip

\myindent{(CM3)} 
(``2 out of 3'') Given two morphisms $f:X\to Y$ and $g:Y\to Z$, 
such that any two of $f$, $g$ and $g\circ f$ are weak equivalences, 
then so is the third.
\smallskip

\myindent{(CM4)} 
(Lifting axiom) Any fibration has the RLP with respect to all 
acyclic cofibrations, and any cofibration has the LLP with respect to all
acyclic fibrations, i.e.\ we can always complete a commutative diagram
\begin{equation}
\xymatrix{A\ar[r]\ar[d]^{i}&X\ar[d]^{f}\\B\ar@{..>}[ur]\ar[r]&Y}
\end{equation}

}{\setmyindent
provided $i$ is a cofibration, $f$ is a fibration, and either $i$ or $f$ 
is a weak equivalence.
\smallskip

\myindent{(CM5)}
(Factorization) Any morphism $f:X\to Y$ can be factored 
both into an acyclic cofibration followed by a fibration, or into 
a cofibration followed by an acyclic fibration: $f=p\circ j=q\circ i$, with 
$i$, $j$ cofibrations, $p$, $q$ fibrations, and $j$, $q$ weak equivalences.

}
\nxsubpoint (Immediate consequences.)
One observation is that the axioms are self-dual, if we interchange 
fibrations and cofibrations, so the opposite~$\cC^0$ of a 
model category~$\cC$ has a natural model category structure itself.
  
Another observation is that each of the lifting properties listed in (CM4) 
determines completely one distinguished class of morphisms 
in terms of another. 
For example, the acyclic cofibrations are {\em exactly\/} those morphisms 
of~$\cC$ that have the LLP with respect to all fibrations, and so on. 
This means that {\em any two of the three distinguished classes of morphisms 
in a model category completely determine the third}, i.e.\ the model 
category data is redundant. For example, once the 
fibrations and weak equivalences are known, the cofibrations are recovered 
as the morphisms with the LLP with respect to all acyclic fibrations. 
Similarly, weak equivalences can be described as those morphisms that 
can be decomposed into an acyclic cofibration followed by 
an acyclic cofibration, hence they are determined once one knows all 
fibrations and cofibrations.

Notice that any class of morphisms characterized by the RLP with 
respect to some other class of morphisms is automatically closed under 
retracts, base change and composition; this applies in particular to 
fibrations and acyclic fibrations.  Similarly, cofibrations and 
acyclic cofibrations are stable under retracts, composition and pushouts. 

\nxsubpoint (Fibrant and cofibrant objects.)
An object $X$ of a model category~$\cC$ is called {\em fibrant\/} if 
the morphism $X\to e_\cC$ from~$X$ into the final object is fibrant; 
dually, $X$ is called {\em cofibrant\/} if $\emptyset_\cC\to X$ is cofibrant.
Given any object~$X$, a {\em cofibrant replacement\/} for~$X$ is 
a weak equivalence $Z\to X$ with $Z$ cofibrant; dually, a {\em fibrant 
replacement\/} for~$X$ is a weak equivalence $X\to Z$ with a fibrant~$Z$. 
Fibrant and cofibrant replacements always exist because of (CM5) applied 
to $X\to e_\cC$ and $\emptyset_\cC\to X$, and we see that a cofibrant 
replacement $Z\to X$ for a fibrant~$X$ can be always chosen in such a way that
$Z$ is both fibrant and cofibrant, hence for an arbitrary~$X$ we can find 
weak equivalences $X\to Y\leftarrow Z$ with a fibrant-cofibrant~$Z$.
 
\nxsubpoint (Example.)
One of the most important examples for us is 
that of the category $\Ch(\cA)$ of {\em non-negative chain complexes\/} 
$K_\cdot=(\cdots\stackrel{\partial_1}\to K_1\stackrel{\partial_0}\to K_0)$ 
over an abelian category~$\cA$ with sufficiently 
many projective objects, e.g.\ $\cA=\catMod R$, $R$ a classical ring. 
In this case the weak equivalences are just the quasi-isomorphisms 
$f:K_\cdot\to L_\cdot$, i.e.\ chain maps, such that $H_n(f)$ is an 
isomorphism for all $n\geq0$. Moreover, $f$ is a fibration iff 
all $f_n$, $n>0$, are epimorphic, and $f$ is a cofibration iff 
all $f_n:K_n\to L_n$ are monomorphisms with $\Coker f_n$ a projective object, 
i.e.\ $f_n$ is actually the embedding $K_n\to L_n\cong K_n\oplus P_n$ with a 
projective~$P_n$. In this case all objects of $\Ch(\cA)$ turn out to be 
fibrant, the cofibrant objects are exactly the chain complexes with projective 
components, and a cofibrant replacement $P_\cdot\to K_\cdot$ of 
a complex~$K_\cdot$ is nothing else than a projective resolution of~$K_\cdot$.
This illustrates the importance of cofibrant replacements for the computation 
of left derived functors.

We'll construct later a counterpart of this model category over 
an arbitrary generalized ring~$\Sigma$, so as to recover $\Ch(\catMod\Sigma)$ 
(together with its model structure just described)
when $\Sigma$ is a classical ring.

\begin{DefD} (Homotopy category of a model category.)
Given any model category~$\cC$, its {\em model category} $\Ho\cC$ is 
the localization of $\cC$ with respect to the class of weak equivalences.
\end{DefD}
In other words, we have a {\em localization functor\/} $\gamma:\cC\to\Ho\cC$, 
such that $\gamma(f)$ is an isomorphism whenever $f$ is a weak equivalence, 
and this functor is universal among all functors $\cC\to\cD$ with 
this property. It is well-known that the localization of a category with 
respect to any set of morphisms always exists; however, in general one 
might need to enlarge the universe~$\univU$ since the set of morphisms between 
two objects of the localization is described as equivalence classes of 
paths in a certain graph, so this set needn't be $\univU$-small. 
This is not the case with the model categories: if $\cC$ is $\univU$-small, 
the same is true for $\Ho\cC$. This follows immediately from the 
``homotopic'' description of $\Ho\cC$ that we are going to recall later.

For example, $\Ho\Ch(\cA)$ for an abelian category~$\cA$ with sufficiently 
many projectives is equivalent to the full subcategory $\cD^{\leq0}(\cA)$ 
of the derived category $\cD^-(\cA)$.

Actually, the weak equivalences in a model category are crucial for 
defining the corresponding homotopy category and derived functors, 
while fibrations and cofibrations should be thought of as some means 
for auxiliary constructions in~$\cC$ and explicit computations of 
derived functors.

\nxsubpoint {\bf Notation.} $\Hom_{\Ho\cC}(X,Y)$ is sometimes denoted by 
$[X,Y]$; if $X$ and $Y$ are objects of~$\cC$, then 
$[X,Y]$ means $\Hom_{\Ho\cC}(\gamma X,\gamma Y)$.

\nxsubpoint 
Since our model categories are always closed, {\em a morphism 
$f:X\to Y$ in~$\cC$ is a weak equivalence if and only if $\gamma(f)$ is 
an isomorphism in~$\Ho\cC$} (cf.\ \cite[1.5]{Quillen}, prop.~1). 

\nxsubpoint\label{sp:cyl.path.obj} (Cylinder and path objects.)  
A {\em cylinder object\/}
$A\times I$ for an object~$A$ of a model category~$\cC$ is a diagram
$A\sqcup A\stackrel{\langle\partial_0,\partial_1\rangle}\longto 
A\times I\stackrel\sigma\to A$ with
$\sigma\circ\langle\partial_0,\partial_1\rangle=\nabla_A=
\langle\id_A,\id_A\rangle$,
such that $\langle\partial_0,\partial_1\rangle$ is a cofibration and
$\sigma$ is a weak equivalence. Dually, a {\em path object\/} $B^I$
for an object~$B$ is a diagram $B\stackrel s\to B^I\stackrel{(d_0,d_1)}
\longto B\times B$ with $(d_0,d_1)\circ s=\Delta_B=(\id_B,\id_B)$, 
where $s$ is a weak equivalence and $(d_0,d_1)$ is a fibration.

Notice that there can be non-isomorphic cylinder objects for the same 
object~$A$; in particular, $A\times I$ does {\em not\/} depend functorially 
on~$A$, and $A\times I$ is not the product of~$A$ and another object~$I$; 
it is just a notation. Same applies to path objects.

\nxsubpoint (Left and right homotopies.)
Given two maps $f,g:A\rightrightarrows B$, we say that 
{\em $f$ is left-homotopic to~$g$} and write $f\siml g$ if 
there is a morphism $h:\tilde A\to B$ from a cylinder object~$\tilde A$,  
$A\sqcup A\stackrel{\langle\partial_0,\partial_1\rangle}\longto\tilde A
\stackrel\sigma\to A$, such that $h\circ\langle\partial_0,\partial_1\rangle=
\langle f,g\rangle:A\sqcup A\to B$, i.e.\ $f=h\circ\partial_0$ and 
$g=h\circ\partial_1$. Dually, we say that {\em $f$ is right-homotopic to~$g$}
and write $f\simr g$ if there is a morphism $h:A\to\tilde B$ 
into a path object $\tilde B$ for~$B$, such that $d_0\circ h=f$ and 
$d_1\circ h=g$. In both cases we say that $h$ is the corresponding 
(left or right) homotopy between $f$ and~$g$.

For any $A$ and $B$ we denote by $\pi^\ell(A,B)$ (resp.\ by $\pi^r(A,B)$) 
the quotient of $\Hom(A,B)$ with respect to the equivalence relation 
generated by $\siml$ (resp.~$\simr$).

Let us list some basic properties of homotopies (cf.\ \cite[1.1]{Quillen}). 
If $A$ is cofibrant, then $\siml$ is an equivalence relation on 
$\Hom_\cC(A,B)$, and $f\siml g$ implies $f\simr g$. Moreover, in this case 
the composition induces a well-defined map $\pi^r(B,C)\times\pi^r(A,B)\to
\pi^r(A,C)$. If $A$ is cofibrant and $B$ is fibrant, $\simr$ and $\siml$ 
coincide on $\Hom(A,B)$; the resulting equivalence relation can be 
denoted simply by $\sim$, and $\pi(A,B):=\Hom(A,B)/\negthinspace\sim$. 
In this case given any $f\sim g$ and any cylinder object $A\times I$ for~$A$, 
we can find a left homotopy $h:A\times I\to B$ between $f$ and~$g$, 
and similarly for path objects~$B^I$ and right homotopies $h:A\to B^I$.

Another simple but important observation is that any of $f\siml g$ or 
$f\simr g$ implies $\gamma(f)=\gamma(g)$ in $\Ho\cC$. A slightly more 
complicated statement is that for cofibrant~$A$ and fibrant~$B$ the converse 
is true: $f\sim g$ iff $\gamma(f)=\gamma(g)$.

\nxsubpoint (Homotopy category and homotopies.)
Given a model category~$\cC$, let us denote by $\cC_c$, $\cC_f$ and $\cC_{cf}$ 
its full subcategories consisting of all cofibrant, fibrant and 
fibrant-cofibrant objects, respectively, and by $\Ho\cC_c$, $\Ho\cC_f$ and 
$\Ho\cC_{cf}$ the localizations of these categories with respect to those 
morphisms that are weak equivalences in~$\cC$. 

Denote by $\pi\cC_c$ the category with the same objects as $\cC_c$, but with 
morphisms given by $\Hom_{\pi\cC_c}(A,B)=\pi^r(A,B)$, and define similarly 
$\pi\cC_f$ and $\pi\cC_{cf}$ 
(in the latter case $\Hom_{\pi\cC_{cf}}=\pi(A,B)$). Since $f\siml g$ and 
$f\simr g$ both imply $\gamma(f)=\gamma(g)$, we get canonical functors 
$\bar\gamma_c:\pi\cC_c\to\Ho\cC_c$, $\bar\gamma_f:\pi\cC_f\to\Ho\cC_f$, 
and $\bar\gamma:\pi\cC_{cf}\to\Ho\cC$. 

{\bf Theorem.} (\cite[1.1]{Quillen}, th.~1)
{\em Functors $\Ho\cC_c\to\Ho\cC$, $\Ho\cC_f\to\Ho\cC$ and 
$\bar\gamma:\pi\cC_{cf}\to\Ho\cC$ are equivalences of categories.}

In other words, any object of $\cC$ becomes isomorphic to a fibrant-cofibrant 
object in the homotopy category (this is clear since any object can be 
connected by a chain of weak equivalences to the fibrant replacement of its 
cofibrant replacement), any morphism $\bar f\in[A,B]$ is representable by a 
morphism $f:A\to B$ if $A$ and $B$ are fibrant-cofibrant, and in this case 
$\bar f=\bar g$ iff $f\sim g$, i.e. $[A,B]\cong\pi(A,B)$. 
It is actually sufficient to require here $A$ to be cofibrant and 
$B$ to be fibrant; cf.~\cite[1.1]{Quillen}, cor.~1 of th.~1.

This theorem shows in particular that $\Ho\cC$ is a $\univU$-category whenever
$\cC$ is one.

\nxsubpoint (Example.)
In the situation $\cC=\Ch(\cA)$, $\cA$ an abelian category with sufficiently 
many projectives, the above theorem actually tells us that $\cD^{\leq0}(\cA)$ 
is equivalent to the category of non-negative chain complexes consisting of 
projective objects with morphisms considered up to (chain) homotopy.

\begin{DefD}\label{def:der.funct}
Given a functor $F:\cC\to\cD$ between two model categories, its 
{\em left derived} functor $\dL F:\Ho\cC\to\Ho\cD$ is the functor closest from 
the left to making the obvious diagram commutative. More formally, 
we must have a natural transformation $\epsilon:\dL F\circ\gamma_\cC\to
\gamma_\cD\circ F$, such that for any other functor $G:\Ho\cC\to\Ho\cD$ 
and natural transformation $\zeta:G\circ\gamma_\cC\to\gamma_\cD\circ F$ 
there is a unique natural transformation $\theta:G\to\dL F$, for which
$\zeta=\epsilon\circ(\theta\star\gamma_\cC)$. The {\em right derived functor}
$\dR F:\Ho\cC\to\Ho\cD$, $\eta:\gamma_\cD\circ F\to\dR F\circ\gamma_\cC$ 
is defined similarly.
\end{DefD}

Of course, the derived functors are defined uniquely by their universal 
properties, but the question of their existence is more subtle. Notice that 
we need only weak equivalences to define the derived functors.

\begin{DefD}\label{def:Q.funct} (Quillen functors.)
Given two adjoint functors $F:\cC\to\cD$, $G:\cD\to\cC$ between two model 
categories, we say that $(F,G)$ is a {\em Quillen pair}, or 
that $F$ is a {\em Quillen functor\/} if $F$ preserves 
weak equivalences between cofibrant objects of~$\cC$, 
and $G$ preserves weak equivalences between fibrant objects of~$\cD$.
\end{DefD}

The importance of Quillen functors is illustrated by the following theorem.

\begin{ThD}\label{th:Q.funct} {\rm (cf.~\cite[1.4]{Quillen}, prop.~2)}
Suppose that $F:\cC\to\cD$ and $G:\cC\to\cD$ constitute a Quillen pair. 
Then both derived functors $\dL F:\Ho\cC\to\Ho\cD$ and $\dR G:\Ho\cD\to\Ho\cC$ 
exist and are adjoint to each other. Moreover, $\dL F(\gamma_\cC X)$ 
is canonically isomorphic to $\gamma_\cD F(Z)$, where $X$ is any object 
of~$\cC$ and $Z\to X$ is any cofibrant replacement of~$X$. 
A similar description in terms of fibrant replacements can be given 
for~$\dR G$ as well.
\end{ThD}

For example, if $\rho:R\to S$ is any ring homomorphism, we have the 
base change and scalar restriction functors 
$\rho^*:\Ch(\catMod R)\to\Ch(\catMod S)$ and $\rho_*:\Ch(\catMod S)\to
\Ch(\catMod R)$. They constitute a Quillen pair, so the derived functors 
$\dL\rho^*$ and $\dR\rho_*$ exist and are adjoint to each other. 
We see that $\dR\rho_*$ can be computed 
componentwise since all objects of $\Ch(\catMod S)$ are already fibrant, 
and to compute $\dL\rho^*(X_\cdot)$ we need to apply $\rho^*$ componentwise to 
a cofibrant replacement, i.e.\ a projective resolution $P_\cdot\to X_\cdot$. 
Therefore, these functors coincide with those known from classical homological 
algebra. Of course, we want to generalize these results later to 
the case of an arbitrary algebraic monad homomorphism $\rho:\Sigma\to\Xi$.

Let us mention a criterion for existence of derived functors that doesn't 
require existence of adjoints:
\begin{PropD}\label{prop:ex.left.der} 
{(\rm \cite[1.4]{Quillen}, cor.\ of prop.~1)}
a) Let $F:\cC\to\cC'$ be a functor between model categories that transforms 
weak equivalences between cofibrant objects of~$\cC$ into weak equivalences 
in~$\cC'$. Then $\dL F:\Ho\cC\to\Ho\cC'$ exists, and for any cofibrant~$P$ 
the morphism $\epsilon_P:\dL F(\gamma_{\cC}P)\to \gamma_{\cC'}F(P)$ is an 
isomorphism, hence $\dL F(\gamma_\cC X)$ can be computed as 
$\gamma_{\cC'}F(P)$ for any cofibrant replacement $P\to X$.

b) If $F':\cC'\to\cC''$ is another functor preserving weak equivalences 
between cofibrant objects, and if~$F$ transforms cofibrant objects of~$\cC$ 
into cofibrant objects of~$\cC'$, then derived functors $\dL(F'\circ F)$, 
$\dL F'$ and $\dL F$ exist, and the functorial morphism 
$\dL F'\circ\dL F\to\dL(F'\circ F)$ arising from the universal property of
$\dL(F'\circ F)$ is an isomorphism.
\end{PropD}
\begin{Proof} 
The first part is proved in~\cite[1.4]{Quillen}, cor.\ of prop.~1, 
while the second one follows immediately from the description of derived 
functors given above.
\end{Proof}

\nxsubpoint (Suspension and loop functors; cf.\ \cite[1.2]{Quillen})
Given any two morphisms $f,g:A\rightrightarrows B$ with a cofibrant~$A$ and a 
fibrant~$B$, we can define a notion of a left homotopy~$H$ between two left 
homotopies $h:A\times I\to B$ and $h':A\times I'\to B$ from $f$ to~$g$. 
This defines a homotopy relation between homotopies from $f$ to~$g$, that 
turns out to be an equivalence relation, thus defining a set 
$\pi_1^\ell(A,B;f,g)$ of homotopy classes of homotopies $h:f\sim g$. 
We have a dual construction for right homotopies, but $\pi_1^r(A,B;f,g)$ 
turns out to be canonically isomorphic to $\pi_1^\ell(A,B;f,g)$, so we can 
write simply~$\pi_1(A,B;f,g)$. Furthermore, we have natural composition maps 
$\pi_1(A,B;f_2,f_3)\times\pi_1(A,B;f_1,f_2)\to\pi_1(A,B;f_1,f_3)$, coming 
from composition of homotopies; in particular, any $\pi_1(A,B;f,f)$ is a group,
and the collection of all $\pi_1(A,B;f,g)$ defines a groupoid.

Now suppose that $\cC$ is a pointed category, i.e.\ it has a zero object 
$0=0_\cC$. Then for any $A$ and $B$ we have a zero map $0=0_{AB}:A\to0\to B$, 
and we put $\pi_1(A,B):=\pi_1(A,B;0,0)$; this is a group for any
cofibrant~$A$ and fibrant~$B$. Next, we can extend $\pi_1(A,B)$ to a functor 
$A,B\mapsto[A,B]_1$, $(\Ho\cC)^0\times\Ho\cC\to\catGrps$, by requiring 
$[A,B]_1=\pi_1(A,B)$ whenever $A$ is cofibrant and $B$ is fibrant.

{\bf Theorem.} (\cite[1.2]{Quillen}, th.~2) 
There are two functors $\Ho\cC\to\Ho\cC$, called the {\em suspension functor 
$\Sigma$} and the {\em loop functor~$\Omega$}, such that
\begin{equation}
[\Sigma A,B]\cong[A,B]_1\cong[A,\Omega B]
\end{equation}

Recall that the suspension $\Sigma A$ of a cofibrant~$A$ can be computed as 
follows: choose any cylinder object~$A\times I$ for~$A$ and take the cofiber
(i.e.\ the pushout with respect to the map $A\sqcup A\to 0$) of the map 
$\langle\partial_0,\partial_1\rangle:A\sqcup A\to A\times I$. Similarly, 
$\Omega B$ is computed for a fibrant~$B$ as the fiber of 
$(d_0,d_1):B^I\to B\times B$.

Recall that $\Sigma^nX$ is a cogroup object and $\Omega^nX$ is a group 
object in~$\Ho\cC$ for $n\geq1$, which is commutative for $n\geq 2$. 
In particular, $\pi_n(X):=[\Sigma^n(*),X]\cong[*,\Omega^nX]$ is a group for 
$n\geq1$, commutative for $n\geq 2$, for any choice of~$*\in\Ob\Ho\cC$.

In the case $\cC=\Ch(\cA)$ the suspension functor $\Sigma$ is just the degree 
translation functor: $(\Sigma K_\cdot)_n=K_{n-1}$, $(\Sigma K_\cdot)_0=0$,
i.e.\ $\Sigma K_\cdot=K_\cdot[1]$,  
while its right adjoint $\Omega$ is given by the truncation of the opposite 
translation functor: $\Omega K_\cdot=\sigma_{\leq0}(K_\cdot[-1])$, i.e.\ 
$(\Omega K_\cdot)_n=K_{n+1}$ for $n>0$, $(\Omega K_\cdot)_0=\Ker
(\partial_0:K_1\to K_0)$.

\nxsubpoint (Cofibrantly generated model categories.)
Suppose that $I$ is a set of cofibrations and $J$ is a set of 
acyclic cofibrations of a model category~$\cC$, such that the fibrations 
are exactly the morphisms with the RLP with respect to all morphisms from~$J$, 
and the acyclic fibrations are exactly the morphisms with the RLP with respect
to~$I$. Then the model category structure of~$\cC$ is completely determined 
by these two sets~$I$ and~$J$, and we say that they {\em generate\/}
the model structure of~$\cC$.

We say that an object $X$ of~$\cC$ is {\em (sequentially) small\/} 
if $\Hom(X,-)$ commutes with sequential inductive limits, i.e.\ 
$$\Hom_\cC(X,\injlim_{n\geq0}Y_n)\cong\injlim_{n\geq0}\Hom_\cC(X,Y_n)\quad.$$ 
In particular, any finitely presented object~$X$ (characterized by the property
of~$\Hom(X,-)$ to commute with filtered inductive limits) is small. 

Finally, we say that a model category~$\cC$ is {\em cofibrantly generated\/} 
if it admits generating sets~$I$ and~$J$ consisting of morphisms with 
small sources. The importance of this notion is due to the fact that 
in a cofibrantly generated model category the factorizations of~(CM5) 
can be chosen to be functorial in~$f$.

Consider for example the case $\cC=\Ch(\catMod R)$, with~$R$ a classical ring. 
Denote by $D(n,R)$ the complex $K_\cdot$ with $K_n=K_{n+1}=R$, 
$\partial_n=\id_R$, and all other $K_i=0$, and by $S(n,R)$ the complex 
$R[n]$ consisting of~$R$ placed in (chain) degree~$n$. Then 
$I=\{S(n,R)\to D(n,R)\,|\,n\geq0\}$ and $J=\{0\to D(n,R)\,|\,n\geq0\}$ 
generate the model structure of~$\cC$, hence~$\Ch(\catMod R)$ is 
cofibrantly generated.

\nxpointtoc{Simplicial and cosimplicial objects}
Since the simplicial and cosimplical objects and categories will play a 
crucial role in the remaining part of this work, we'd like to list their 
basic properties and fix some notation. Our main reference here is~\cite{GZ}.

\begin{DefD} 
The {\em (non-empty finite) ordinal number category~$\catDelta$}
is by definition the full subcategory of the category of ordered sets, 
consisting of all standard finite non-empty ordered sets 
$[n]=\{0,1,\ldots,n\}$, $n\geq0$.
\end{DefD}
In other words, the objects $[n]$ of $\catDelta$ are in one-to-one 
correspondence with non-negative integers, and the morphisms 
$\phi:[n]\to[m]$ are simply the order-preserving maps: $\phi(x)\leq\phi(y)$ 
whenever $0\leq x\leq y\leq n$.

We consider two subcategories $\catDelta_+$ and $\catDelta_-$ of $\catDelta$ 
as well. They have the same objects as~$\catDelta$, but 
$\Hom_{\catDelta_+}([n],[m])$ consists of all injective order-preserving maps 
$\phi:[n]\to[m]$, while $\Hom_{\catDelta_-}([n],[m])$ consists of all 
surjective order-preserving maps.

\nxsubpoint\label{sp:face.deg} (Face and degeneracy maps.)
For any integers $0\leq i\leq n$ we define two morphisms in~$\catDelta$:
\begin{itemize}
\item The {\em face map\/} $\partial_n^i:[n-1]\to[n]$, the increasing 
injection that doesn't take value $i\in[n]$, $n>0$;
\item The {\em degeneracy map\/} $\sigma_n^i:[n+1]\to[n]$, the non-decreasing 
surjection that takes twice the value $i\in[n]$.
\end{itemize}
When no confusion can arise, we write simply $\partial^i$ and $\sigma^i$.

The compositions of these maps are subject to following relations:
\begin{align}
\label{eq:face.face}
\partial_{n+1}^j\partial_n^i&=\partial_{n+1}^i\partial_n^{j-1},
\quad{0\leq i<j\leq n+1}\\
\label{eq:deg.deg}
\sigma_n^j\sigma_{n+1}^i&=\sigma_n^i\sigma_{n+1}^{j+1},
\quad{0\leq i\leq j\leq n}\\
\label{eq:deg.face}
\sigma_{n-1}^j\partial_n^i&=
\begin{cases}
\partial_{n-1}^i\sigma_{n-2}^{j-1},&i<j\\
\id_{[n-1]},&j\leq i\leq j+1\\
\partial_{n-1}^{i-1}\sigma_{n-2}^j,&i>j+1
\end{cases}
\end{align}

It is well-known that any morphism $\phi:[n]\to[m]$ can be uniquely decomposed 
into a non-decreasing surjection $p:[n]\twoheadrightarrow[s]$ followed by an 
increasing injection $i:[s]\hookrightarrow[m]$. More precisely, any morphism 
$\phi:[n]\to[m]$ can be uniquely written in form
\begin{multline}
\phi=\partial_m^{i_{m-s}}\partial_{m-1}^{i_{m-s-1}}\cdots\partial_{s+1}^{i_1}
\sigma_s^{j_{n-s}}\sigma_{s+1}^{j_{n-s-1}}\cdots\sigma_{n-1}^{j_1}\quad,\\
\quad\text{where }m\geq i_{m-s}>\cdots>i_1\geq0,
\quad 0\leq j_{n-s}<\cdots<j_1<n
\end{multline}
An immediate consequence of this lemma is that $\catDelta$ can be identified 
with the category generated by the face and degeneracy maps subject to 
conditions listed above (cf.~\cite[2.2.2]{GZ}). Notice that the category 
generated by the face maps~$\partial_n^i$ subject to 
relations~\eqref{eq:face.face} coincides with~$\catDelta_+$, 
and the category generated by the degeneracy maps~$\sigma_n^i$ subject to 
relations~\eqref{eq:deg.deg} coincides with~$\catDelta_-$.

\begin{DefD}
Let~$\cC$ be an arbitrary category. A {\em simplicial object of~$\cC$} 
is a functor $X:\catDelta^0\to\cC$. The category of all simplicial objects 
of~$\cC$ will be denoted by $s\cC:=\catFunct(\catDelta^0,\cC)$. Dually, 
a {\em cosimplicial object of~$\cC$} is a functor $Y:\catDelta\to\cC$, 
and the category of all such cosimplicial objects will be denoted by 
$c\cC:=\catFunct(\catDelta,\cC)$.
\end{DefD}

Usually we write $X_n$ instead of $X([n])$, and for any $0\leq i\leq n$ put 
$d^{n,X}_i:=X(\partial_n^i):X_n\to X_{n-1}$ (if $n>0$) and
$s^{n,X}_i:=X(\sigma_n^i):X_n\to X_{n+1}$. 
When no confusion can arise, we write simply $d_i$ and $s_i$. In this way 
a simplicial object~$X$ of~$\cC$ is simply a collection of objects $X_n$, 
$n\geq0$, together with some morphisms $d^{n,X}_i:X_n\to X_{n-1}$, called 
{\em face operators}, and $s^{n,X}_i:X_n\to X_{n+1}$, called 
{\em degeneracy operators}, subject to the dual of relations 
of~\ptref{sp:face.deg}:
\begin{align}
\label{eq:coface.coface}
d^n_i d^{n+1}_j &= d^n_{j-1} d^{n+1}_i, \quad{0\leq i<j\leq n+1}\\
\label{eq:codeg.codeg}
s^{n+1}_i s^n_j &= s^{n+1}_{j+1} s^n_i,\quad{0\leq i\leq j\leq n}\\
\label{eq:codeg.coface}
d^n_i s^{n-1}_j &= \begin{cases}
  s^{n-2}_{j-1} d^{n-1}_i, &i<j\\
  \id_{X_{n-1}}, &j\leq i\leq j+1\\
  s^{n-2}_j d^{n-1}_{i-1}, &i>j+1
\end{cases}
\end{align}

Similarly, we write $Y^n$ instead of $Y([n])$ and put 
$\partial_{n,Y}^i:=Y(\partial_n^i):Y_{n-1}\to Y_n$, 
$\sigma_{n,Y}^i:=Y(\sigma_n^i):Y_{n+1}\to Y_n$; we sometimes write simply 
$\partial^i$ and $\sigma^i$. In this way a cosimplicial object~$Y$ is a 
collection of objects~$Y^n$, $n\geq0$, together with face and degeneracy 
operators $\partial_{n,Y}^i$, $\sigma_{n,Y}^i$, subject to 
relations~\ptref{sp:face.deg}.

Clearly, $c(\cC^0)=(s\cC)^0$ and $c\cC=(s(\cC^0))^0$, so statements about 
simplicial and cosimplicial objects can be obtained from each other 
by dualizing.

\nxsubpoint\label{sp:skel.coskel} (Skeleta and coskeleta.)
Consider the full subcategory $\catDelta_n\subset\catDelta$ consisting of 
objects $[p]$ with $p\leq n$. Functors $\catDelta_n^0\to\cC$ and 
$\catDelta_n\to\cC$ are called {\em truncated (co)simplicial objects of~$\cC$}.
Clearly, any (co)simplicial object can be truncated; if $I_n:\catDelta_n\to
\catDelta$ is the embedding functor, then the truncation functor is equal to 
$I_n^*:X\mapsto X\circ I_n$. Therefore, we can study its left and 
right adjoints, i.e.\ left and right Kan extensions $I_{n,!}$ and $I_{n,*}$ 
of $I_n$. Of course, in general they exist only if arbitrary inductive, resp.\ 
projective limits exist in~$\cC$. In our case it suffices to require existence 
of finite inductive, resp.\ projective limits, since $\catDelta_n/[m]$ and 
$[m]\backslash\catDelta_n$ are finite categories:
\begin{equation}
(I_{n,!}X)(m)=\injlim_{([m]\backslash\catDelta_n)^0}X
\end{equation}

Since $I_n$ is fully faithful, the same is true for its Kan extensions, hence 
$X\to I_n^*I_{n,!}X$ and $I_n^*I_{n,*}X\to X$ are isomorphisms for any 
truncated (co)sim\-pli\-cial object~$X$. We define the {\em $n$-th skeleton} 
(resp.\ {\em coskeleton}) of a (co)sim\-pli\-cial object~$X$ by 
$\sk_n X:=I_{n,!}I_n^*X$ (resp.\ $\cosk_n X:=I_{n,*}I_n^*X$). We have 
canonical morphisms $\sk_n X\to X$ and $X\to\cosk_n X$, that are isomorphisms 
in degrees $\leq n$; in fact, $\sk_n X\to X$ has a universal property among 
all morphisms $Y\to X$ that are isomorphisms in degrees $\leq n$, and 
similarly for~$\cosk_n X$.

\nxsubpoint (Simplicial sets.)
A very important example is given by the {\em category of simplicial sets}
$s\catSets=\catFunct(\catDelta^0,\catSets)$. Given a simplicial set~$X$, 
we say that elements $x\in X_n$ are the {\em $n$-simplices of~$X$}; 
$n$ is called the {\em degree\/} or {\em dimension\/} of~$x$.
If $x$ doesn't come from any $X_m$ with $m<n$, i.e.\ if $x$ doesn't lie in the 
union of the images of the degeneracy operators $s^{n,X}_i:X_{n-1}\to X_n$, 
we say that $x$ is a {\em non-degenerate $n$-simplex of~$X$.} 
Finally, we say that $X$ is a {\em finite simplicial set\/} if the set of
its non-degenerate simplices (in all degrees) is finite. It is easy to check 
that finite simplicial sets are exactly the finitely presented objects 
of~$s\catSets$. Notice that the $n$-th skeleton $\sk_n X$ of a simplicial 
set~$X$ is its smallest subobject containing all simplices of~$X$ of 
dimension $\leq n$. This explains why we say that {\em $X$ is of dimension 
$\leq n$} whenever $\sk_n X$ coincides with~$X$; this is equivalent to 
saying that all simplices of $X$ in dimensions $>n$ are degenerate. 

A simplicial set~$X$ is finite iff it has bounded dimension and all $X_n$ are 
finite. An important property is that {\em the product of two finite simplicial
sets is finite.}

\nxsubpoint (Examples of simplicial sets.)
Let us list some important examples of (finite) simplicial sets. 
First of all, for any $n\geq0$ we have the corresponding representable functor
$\Delta(n):=\Hom_{s\catSets}(-,[n])$, called the {\em standard $n$-simplex}.
By Yoneda $\Hom_{s\catSets}(\Delta(n),X)\cong X_n$. In this way we obtain 
a fully faithful functor $\Delta:\catDelta\to s\catSets$. It is easy to see 
that finite simplicial sets can be described as finite inductive limits of 
standard simplices.

Another important example is the {\em boundary\/} of the standard $n$-simplex: 
$\dot\Delta(n):=\sk_{n-1}\Delta(n)\subset\Delta(n)$. It is the largest 
simplicial subset of $\Delta(n)$ that doesn't contain 
$\id_{[n]}\in\Delta(n)_n$. Other description:
\begin{equation}
\dot\Delta(n)_m=\{\text{non-decreasing maps $\phi:[m]\to[n]$ with 
$\phi([m])\neq[n]$}\}
\end{equation}
In particular, $\dot\Delta(0)=\emptyset$.

We can define the {\em horns\/} $\Lambda_k(n)\subset\dot\Delta(n)\subset
\Delta(n)$, $0\leq k\leq n$, by removing one of the faces of $\dot\Delta(n)$. 
More precisely, $\Lambda_k(n)_m$ consists of all order-preserving maps 
$\phi:[m]\to[n]$, such that the complement of their image $[n]-\phi([m])$ 
is distinct from both $\emptyset$ and $\{k\}$. 

\nxsubpoint
(Geometric realization and singular simplices.)
For any set~$S$ we denote by $\Delta^{(S)}$ the convex hull of the standard 
base of $\bbR^{(S)}$, considered as a topological space. This is a 
covariant functor from $\catSets$ into $\catTop$. 
Restricting it to $\catDelta\subset\catSets$ we 
obtain a covariant functor $\Delta^{(-)}:\catDelta\to\catTop$, i.e.\ 
a cosimplicial topological space. Now we define the {\em singular functor\/}
$\Sing:\catTop\to s\catSets$ by putting $(\Sing X)([n]):=
\Hom_\catTop(\Delta^{([n])},X)$; clearly, $(\Sing X)_n$ is just the set 
of singular $n$-simplices of~$X$. We define the 
{\em geometric realization functor\/} $|-|:s\catSets\to\catTop$ as the 
left adjoint to $\Sing$; its uniqueness is evident, and its existence follows 
from the fact that $|\Delta(n)|$ must be equal to $\Delta^{([n])}$, that 
any simplicial set can be written down as an inductive limit of standard 
simplices, and that $|-|$ must commute with arbitrary inductive limits. 
This yields a formula for $|X|$:
\begin{equation}
|X|=\injlim\nolimits_{\catDelta/X}\Delta^{(-)}
\end{equation}
The geometric realization functor takes values in the full subcategory of 
$\catTop$ consisting of Kelley spaces (i.e.\ compactly generated Hausdorff 
spaces); if we consider it as a functor from $s\catSets$ into this 
subcategory, it preserves finite projective and arbitrary inductive 
limits (cf.~\cite[3.3]{GZ}).

\nxsubpoint 
(Closed model category structure on $\catTop$.)
There is a closed model category structure on the category of topological 
spaces, characterized by its classes of fibrations and weak equivalences: 
the fibrations are the Serre fibrations and weak equivalences are maps 
$f:X\to Y$ that preserve all homotopy groups, i.e.\ 
$\pi_n(f):\pi_n(X,x)\to\pi_n(Y,f(x))$ has to be an isomorphism for all 
$n\geq0$ and all base points $x\in X$ (cf.~\cite[2.3]{Quillen}).

\nxsubpoint
(Closed model category structure on $s\catSets$.)
The category of simplicial sets has a model category structure as well, 
that can be described as follows. The fibrations in~$s\catSets$ are the 
Kan fibrations, cofibrations are just the injective maps, and 
$f:X\to Y$ is a weak equivalence iff its geometric realization 
$|f|$ is one. In particular, all simplicial sets are cofibrant, but usually 
not fibrant; recall that $\Ch\cA$ had the opposite property.

Moreover, adjoint functors $|-|$ and $\Sing$ constitute a Quillen pair, 
and their derived functors $\dL|-|$ and $\dR\Sing$ are adjoint equivalences 
of categories; thus $\Ho s\catSets$ and $\Ho\catTop$ are naturally equivalent 
(cf.~\cite[2.3]{Quillen}).

\nxsubpoint\label{sp:cof.gen.ssets}
(Cofibrant generators for $s\catSets$.)
The model structure of $s\catSets$ happens to be cofibrantly generated. 
More precisely, $I:=\{\dot\Delta(n)\to\Delta(n)\,|\,n\geq0\}$ is a 
generating set of cofibrations with small sources, while 
$J:=\{\Lambda_k(n)\to\Delta(n)\,|\,0\leq k\leq n>0\}$ is a generating set of 
acyclic cofibrations.

Let's apply this to describe the acyclic fibrations $p:E\to B$, i.e.\ 
the maps with the RLP with respect to all inclusions 
$\dot\Delta(n)\to\Delta(n)$.  We know that 
$\Hom(\Delta(n),X)\cong X_n$, and using adjointness of 
$\sk_{n-1}$ and $\cosk_{n-1}$ we get 
$\Hom(\dot\Delta(n),X)=\Hom(\sk_{n-1}\Delta(n), X)\cong
\Hom(\Delta(n),\cosk_{n-1}X)\cong(\cosk_{n-1}X)_n$. We see that 
$p:E\to B$ is an acyclic fibration iff $E_n\to B_n\times_{(\cosk_{n-1}B)_n}
(\cosk_{n-1}E)_n$ is surjective for all $n>0$, and $p_0:E_0\to B_0$ is 
surjective as well. In particular, $X$ is fibrant iff 
$X_0\neq\emptyset$ and $X_n\to(\cosk_{n-1}X)_n$ is surjective for all~$n>0$.

\nxsubpoint
(Simplicial groups.)
The category $s\catGrps$ of simplical groups has a model category structure as 
well. Now we'd like to mention that $f:G\to H$ is a fibration in $s\catGrps$ 
iff it is a fibration in $s\catSets$, and that all simplicial groups are 
fibrant, both in $s\catGrps$ and in~$s\catSets$.

\nxpointtoc{Simplicial categories}
There are several ways of defining a simplicial structure on a category~$\cC$,
all of them more or less equivalent when arbitrary inductive and projective 
limits exist in~$\cC$. We'll use a definition based on the notion of 
an external $\otimes$-action.

\nxsubpoint (ACU $\otimes$-structure on $s\catSets$.)
First of all, notice that the direct product defines an ACU $\otimes$-structure
on $\catSets$, hence also on~$s\catSets$. We have 
$(K\otimes L)_n=(K\times L)_n=K_n\times L_n$ for this $\otimes$-structure.

\begin{DefD}
A {\em simplicial structure\/} on a category~$\cC$ is a right $\otimes$-action 
$\otimes:\cC\times s\catSets\to\cC$, $(X,K)\mapsto X\otimes K$
(cf.~\ptref{sp:ext.tens.act}) 
of~$s\catSets$ with the ACU $\otimes$-structure just considered on~$\cC$, 
required to commute with arbitrary inductive limits in~$K$.
A {\em simplicial category\/} is just a category with a simplicial structure.
A {\em simplicial functor} $F:\cC\to\cD$ between two simplicial categories 
is simply a {\em lax} external $\otimes$-functor in the sense 
of~\ptref{sp:tens.functor}. Recall that this involves giving a collection of 
morphisms $F(X)\otimes K\to F(X\otimes K)$, subject to certain relations.
\end{DefD}

Notice that the $\otimes$-structure on $s\catSets$ is commutative, so the 
distinction between left and right external $\otimes$-actions of $s\catSets$ 
is just a matter of convention. 

\nxsubpoint 
(A variant.) Since any simplicial set is an inductive (even filtered inductive)
limit of finite simplicial sets, and the category of finite simplicial sets is 
closed under products, we might have defined a simplicial structure as an 
external $\otimes$-action of the ACU $\otimes$-category of finite simplicial 
sets, required to commute with finite inductive limits. This definition is 
completely equivalent to the one given above when arbitrary inductive limits 
exist in~$\cC$.

Notice that any simplicial set is an inductive limit of standard 
simplices~$\Delta(n)$, so it might suffice to define all $X\otimes\Delta(n)$. 
However, the associativity conditions $(X\otimes\Delta(n))\otimes\Delta(m)\cong
X\otimes(\Delta(n)\times\Delta(m))$ cannot be naturally expressed in terms
of such tensor products since $\Delta(n)\times\Delta(m)$ is not a standard 
simplex itself.

\nxsubpoint\label{sp:simpl.hom.sets}
(Simplicial sets of homomorphisms.)
Let us fix two objects $X$ and $Y$ of a simplicial category~$\cC$. 
The functor $\Hom_\cC(X\otimes-,Y)$, $(s\catSets)^0\to\catSets$, transforms 
arbitrary inductive limits into corresponding projective limits; therefore, 
it is representable by an object $\iHom(X,Y)=\iHom_\cC(X,Y)$:
\begin{equation}
\Hom_{s\catSets}(K,\iHom_\cC(X,Y))\cong\Hom_\cC(X\otimes K,Y)
\end{equation}
In this way we get a functor $\iHom_\cC:\cC^0\times\cC\to s\catSets$.
These simplicial sets $\iHom(X,Y)$ will be called the 
{\em simplicial sets of homomorphisms}. Clearly, $\iHom_\cC(X,Y)_n=
\Hom_\cC(X\otimes\Delta(n),Y)$; in particular, the unit axiom for an external 
action implies that $\Hom_\cC(X,Y)\cong\iHom_\cC(X,Y)_0$ functorially in 
$X$ and~$Y$. A standard argument yields the existence of canonical
{\em composition maps}, satisfying natural associativity and unit conditions:
\begin{equation}
\circ:\iHom_\cC(Y,Z)\times\iHom_\cC(X,Y)\to\iHom_\cC(X,Z)\quad.
\end{equation}

Quillen has defined in~\cite{Quillen} a simplicial category as a category~$\cC$
together with a functor~$\iHom_\cC$, a collection of composition maps~$\circ$ 
and functorial isomorphisms $\iHom_\cC(X,Y)_0\cong\Hom_\cC(X,Y)$, subject to 
some relations (essentially the associativity and unit relations just 
mentioned). Of course, $\iHom_\cC$ determines the external action $\otimes$ 
uniquely up to a unique isomorphism, so Quillen's approach is essentially 
equivalent to one adopted here.

Yet another approach consists in describing the simplicial structure in terms 
of the functor $M:\cC^0\times(s\catSets)^0\times\cC\to\catSets$, 
$M(X,K;Y):=\Hom_\cC(X\otimes K,Y)$\dots

\nxsubpoint (Exponential objects.)
Given an object $Y\in\Ob\cC$ and a simplicial set~$K$, we define the 
corresponding {\em exponential object\/} $Y^K$ as the object representing 
the functor $\Hom_\cC(-\otimes K,Y)$, whenever it exists. Therefore, 
\begin{equation}\label{eq:adj.simpl.str}
\Hom_\cC(X\otimes K,Y)\cong\Hom_{s\catSets}(K,\iHom_\cC(X,Y))\cong
\Hom_\cC(X,Y^K)
\end{equation}

Clearly, $(Y^K)^L\cong Y^{K\times L}$, and $Y^{\Delta(0)}\cong Y$; 
of course, the functor $(Y,K)\mapsto Y^K$ determines completely 
the simplicial structure (actually this functor is an external $\otimes$-action
of $(s\catSets)^0$ on~$s\cC$), so 
we might have written the axioms in terms of this functor as well.

\nxsubpoint (Simplicial structure on categories of simplicial objects.)
Given any object~$X$ of a category~$\cC$ and any set~$E$, we denote by 
$X\times E$ the direct sum (i.e.\ coproduct) of $E$ copies of~$X$, i.e.\ 
$\Hom_\cC(X\times E,Y)\cong\Hom_\cC(X,Y)^E$.

Consider now the category $s\cC=\catFunct(\catDelta^0,\cC)$ of 
simplicial objects over a category~$\cC$, supposed to have arbitrary inductive 
and projective limits. 
Then $s\cC$ has a natural simplicial structure, given by
\begin{equation}
(X\otimes K)_n:=X_n\times K_n,\quad
\text{for any $X\in\Ob s\cC$, $K\in\Ob s\catSets$.}
\end{equation}
Moreover, one can prove (cf.~\cite[2.1]{Quillen}) that in this case
all exponential objects $Y^K$ are representable. Notice that the
simplicial extension $sF:s\cC\to s\cD$ of an arbitrary functor
$F:\cC\to\cD$ turns out to be automatically simplicial with respect to
this simplicial structure by means of canonical morphisms
$F(X_n)\times K_n\to F(X_n\times K_n)$.

If $\cC$ doesn't have any inductive or projective limits, this simplicial 
structure can be still described in terms of functors $\iHom_\cC$ or~$M$ 
as explained before. 

\nxsubpoint (Simplicial structure on categories of cosimplicial objects.)
Notice that the notion of a simplicial structure on a category~$\cC$ 
is essentially self-dual, if we don't forget to interchange the arguments of
$\iHom$, and $X\otimes K$ with $X^K$, i.e.\ the opposite of a simplicial 
category has a natural simplicial structure itself. 

Let $\cC$ be an arbitrary category with arbitrary inductive and projective 
limits. Since $c\cC=\catFunct(\catDelta,\cC)=
(s(\cC^0))^0$, we must have a simplicial structure on $c\cC$ as well. 
After dualizing and interchanging  
$X\otimes K$ and $X^K$, we see that this simplicial structure is 
easily characterized in terms of $X^K$. Namely, for any cosimplicial object~$X$
and any simplicial set~$K$ we have
\begin{equation}
(X^K)^n=\Hom(K_n,X^n)=(X^n)^{K_n}
\end{equation}
After this $\iHom_{c\cC}(X,Y)$ and $X\otimes K$ are completely determined 
by~\eqref{eq:adj.simpl.str}.

\nxsubpoint 
It seems a bit strange that the categories of simplicial and cosimplicial sets 
act on each other in an asymmetric way: $s\catSets$ acts both on itself 
and on $c\catSets$, while $c\catSets$ doesn't act on anything. 
This can be explained in the following fancy way. If we think of simplicial 
sets as chain complexes of $\Fempty$-modules, and of cosimplicial sets 
as cochain complexes, then the action of $s\catSets$ on itself is essentially 
the same thing as the tensor product over $\Fempty$ of chain complexes, while 
the action of $s\catSets$ on $c\catSets$ is something like $\Hom_{\Fempty}$ 
from a chain complex into a cochain complex, yielding a cochain complex.

\nxpointtoc{Simplicial model categories}
\begin{DefD} 
We say that a category~$\cC$ is a simplicial (closed) model category 
if it has both a simplicial and a closed model category structure, 
these two compatible in the following sense (cf.~\cite[2.2]{Quillen}):
\end{DefD}

{\myindent{(SM7)}
Whenever $i:A\to B$ is a cofibration and $p:X\to Y$ is a fibration in~$\cC$, 
the following map of simplicial sets is a fibration in~$s\catSets$:
\begin{equation}
\iHom(B,X)\stackrel{(i^*,p_*)}\longto
\iHom(i,p):=\iHom(A,X)\times_{\iHom(A,Y)}\iHom(B,X)
\end{equation}

\setmyindent
We require this fibration to be acyclic when either $i$ or $p$ is a weak 
equivalence.

}

This axiom is easily shown to be equivalent to any of the following:

\smallbreak
{\myindent{(SM7a)}
If $p:X\to Y$ is a fibration in~$\cC$ and $s:K\to L$ is a cofibration 
in~$s\catSets$, then
\begin{equation}
X^L\to X^K\times_{Y^K}Y^L
\end{equation}
\setmyindent

is a fibration in~$\cC$, acyclic if either $p$ or~$s$ is a weak equivalence.
}
\smallbreak
{\myindent{(SM7b)}
If $i:A\to B$ is a cofibration in~$\cC$ and $s:K\to L$ is a cofibration 
in~$s\catSets$, then
\begin{equation}
i\boxe s:A\otimes L\sqcup_{A\otimes K}B\otimes K\to B\otimes L
\end{equation}
\setmyindent

is a cofibration in~$\cC$, acyclic if either $i$ or~$s$ is a weak equivalence.
}

It is actually sufficient to require (SM7a) or (SM7b) when $s:K\to L$ runs 
through a system of generators for (acyclic or all) cofibrations 
in~$s\catSets$.

\nxsubpoint (Example.)
The category of simplicial sets $s\catSets$ happens to be a simplical 
model category.

\nxsubpoint (Simplicial homotopies.)
Given two objects $X$ and $Y$ of an arbitrary simplicial category~$\cC$, 
we say that $h\in\iHom_\cC(X,Y)_1$ is a {\em simplicial homotopy\/} 
between two morphisms $f,g\in\Hom_\cC(X,Y)=\iHom_\cC(X,Y)_0$ if 
$d_0(h)=f$ and $d_1(h)=g$. In this case we say that {\em $f$ and~$g$ 
are strictly simplicially homotopic} and write $f\simss g$. The equivalence 
relation on $\Hom_\cC(X,Y)$ generated by $\simss$ will be denoted by~$\sims$; 
if $f\sims g$, we say that {\em $f$ and $g$ are simplicially homotopic.} 
Notice that giving an element $h\in\iHom(X,Y)_1$ is equivalent to giving 
a morphism $h':X\otimes\Delta(1)\to Y$ or $h'':X\to Y^{\Delta(1)}$.

The quotient $\Hom(X,Y)/\negthinspace\sims$ will be denoted by 
$\pi_0\Hom(X,Y)$, or simply by $\pi_0(X,Y)$. For any full subcategory 
$\cC'\subset\cC$ we denote by $\pi_0\cC'$ the category with the same objects, 
but with morphisms given by $\pi_0(X,Y)$.

More generally, for any simplicial set $X$ we denote by $\pi_0 X$ the 
cokernel of the two face operators $d_0,d_1:X_1\rightrightarrows X_0$. 
We say that $\pi_0X$ is {\em the set of components\/} of~$X$.

\nxsubpoint (Homotopies in simplicial model categories.)
Recall the following proposition (cf.~\cite[2.2]{Quillen}, prop.~5):
\begin{Propz} 
If $f,g:X\rightrightarrows Y$ are two morphisms in a 
simplicial model category~$\cC$, then $f\sims g$ implies both 
$f\siml g$ and $f\simr g$. If $X$ is cofibrant and $Y$ is fibrant, 
then the homotopy relations $\simr$, $\siml$, $\sims$ and $\simss$ 
on~$\Hom_\cC(X,Y)$ coincide. In particular, $\Ho\cC\cong\pi\cC_{cf}\cong
\pi_0\cC_{cf}$, i.e.\ the homotopy category of~$\cC$ can be computed 
with the aid of simplicial homotopies.
\end{Propz}

The proof involves a very simple but useful observation. Namely, the
face and degeneracy maps induce two acyclic cofibrations
$\Delta(d_0),\Delta(d_1):\Delta(0)\rightrightarrows\Delta(1)$ and a
weak equivalence $\Delta(s_0):\Delta(1)\to\Delta(0)$. Using (SM7) in
its equivalent forms we see that induced diagrams $X\oplus X\to
X\otimes\Delta(1)\to X$ and $Y\to Y^{\Delta(1)}\to Y\times Y$ give us
a canonical way to construct a cofibrant cylinder for any
cofibrant~$X$, resp.\ a fibrant path object for any fibrant~$Y$.  If
we put $I:=\Delta(1)$ (``the simplicial segment''), then this cylinder
and path object will be indeed denoted by $X\times I=X\otimes I$,
resp.\ $Y^I$; this explains our previous notation for them.

\nxsubpoint\label{th:main.quillen} 
(Simplicial model category structure on~$s\cC$.)
Let~$\cC$ be a category closed under arbitrary inductive and projective limits.
We would like to discuss some conditions under which it is possible 
to endow the simplicial category~$s\cC$ with a compatible model category 
structure in a natural way (cf.~\cite[2.4]{Quillen}).

Recall that $f:X\to Y$ is said to be a {\em strict epimorphism\/} if 
it is the cokernel of its kernel pair $X\times_YX\rightrightarrows X$. 
An object $P\in\Ob\cC$ is {\em projective\/} if $\Hom_\cC(P,-)$ transforms 
strict epimorphisms into surjections, and {\em finitely presented\/} 
if $\Hom_\cC(P,-)$ preserves filtered inductive limits. We say that 
{\em $\cC$ has sufficiently many projectives\/} if for any object~$X$ 
there is a strict epimorphism $P\to X$ with a projective~$P$. 
A set $\cG\subset\Ob\cC$ is a set of {\em generators\/} of~$\cC$ if
for any object~$X$ of~$\cC$ one can find a strict epimorphism from a 
direct sum of objects of~$\cG$ into~$X$. 

Now we are able to state the following important theorem:
\begin{Thz} {\rm (\cite[2.4]{Quillen}, th.~4)}
Let $\cC$ be a category closed under arbitrary projective and inductive limits,
$s\cC$ be the corresponding category of simplicial objects. Define a morphism 
$f$ in~$s\cC$ to be a fibration (resp.\ weak equivalence) iff 
$\Hom(P,f)$ is a fibration (resp.\ weak equivalence) in~$s\catSets$ 
for each projective object~$P$ of~$\cC$, and a cofibration if it has the 
LLP with respect to all acyclic fibrations. Then $s\cC$ becomes a simplicial 
model category provided one of the following conditions holds:
\begin{itemize}
\item[a)] Every object of~$s\cC$ is fibrant.
\item[b)] Any object of~$\cC$ is a strict quotient of a cogroup object.
\item[c)] $\cC$ has a small set of finitely presented projective generators.
\end{itemize}
Moreover, these conditions are not independent: b) implies~a).
\end{Thz}

Notice that if $\cG$ is a set of projective generators of~$\cC$, then any 
projective object is a direct factor of a direct sum of objects from~$\cG$, 
hence $f$ is a fibration (resp.\ weak equivalence) in $s\cC$ iff 
$\Hom(P,f)$ is a fibration (resp.\ weak equivalence) for all $P\in\cG$.
(We use here that any product of weak equivalences in $s\catSets$ is again 
a weak equivalence!)

Another interesting observation is that in case c) the model structure 
of~$s\cC$ turns out to be cofibrantly generated. More precisely, if 
$\cG\subset\Ob\cC$ is a set of finitely presented projective generators 
as in~c), and $I$ and $J$ are the cofibrant generators for~$s\catSets$, 
then $\cG\otimes I=\{P\otimes f\,|\,P\in\cG$, $f\in I\}$ and 
$\cG\otimes J$ are cofibrant generators for $s\cC$.

\nxsubpoint\label{sp:appl.smc} (First applications.)
The above theorem is of course applicable to $\cC=\catSets$, $\st1$ being 
a finitely presented projective generator of~$\catSets$. The model structure 
we get coincides with the one we have already, since $\Hom(\st1,f)=f$. 

Another application: Let $R$ be a (classical) associative ring, 
$\cC:=\catMod R$. Then $R$ is a finitely presented projective generator 
of~$\cC$, so the condition c) of the theorem is fulfilled. Actually 
the codiagonal $\nabla_C:C\to C\oplus C$ defines a cogroup structure on 
any object of $\catMod R$, so conditions b) and a) are fulfilled 
as well.

\nxsubpoint (Application: category of cosimplicial sets.)
A less trivial example is given by the category~$\catSets^0$. In this 
case condition b), hence also a), is fulfilled since the strict monomorphisms 
of~$\catSets$ are just the injective maps, and any set can be embedded into 
a group, e.g.\ the corresponding free group. Therefore, 
$s(\catSets^0)=(c\catSets)^0$ has a simplicial model category structure. 
By dualizing we obtain such a structure on the category $c\catSets$ of 
cosimplicial sets. Since $\st2=\{1,2\}$ is an injective cogenerator 
of~$\catSets$, we see that {\em 
a map $f:X\to Y$ of cosimplicial sets is a cofibration 
(resp.\ a weak equivalence) iff the map of simplicial sets 
$f^*:\st2^Y\to\st2^X$ is a fibration (resp.\ a weak equivalence)}, 
and that {\em all cosimplicial sets are cofibrant.}

\nxsubpoint\label{sp:cat.simpl.mod} (Category of simplicial $\Sigma$-modules.)
Now let $\Sigma$ be an algebraic monad, $\cC:=\catMod\Sigma$. Clearly,
$|\Sigma|=\Sigma(1)$ is a finitely presented projective generator
for~$\cC$, so Theorem~\ptref{th:main.quillen},c) is applicable, hence
{\em the category $s(\catMod\Sigma)$ of simplicial $\Sigma$-modules
admits a natural simplicial model category structure.} Moreover, {\em this
model structure is cofibrantly generated by $L_\Sigma(I)$ and
$L_\Sigma(J)$, where $I$ and $J$ are any cofibrant generators
for~$s\catSets$.} Since $\Hom_\Sigma(|\Sigma|,-)$ is nothing else than
the forgetful functor $\Gamma_\Sigma:\catMod\Sigma\to \catSets$, we
see that {\em a morphism $f$ of simplicial $\Sigma$-modules is a
fibration (resp.\ weak equivalence) iff this is true for
$\Gamma_\Sigma(f)$, i.e.\ iff $f$ is a fibration (resp.\ weak
equivalence) when considered as a map of simplicial sets.}

When $\Sigma$ is a classical ring, we recover 
the simplical model category structure discussed in~\ptref{sp:appl.smc}. 
When $\Sigma=\Fempty$ we get the usual simplicial model category structure 
on $s(\catMod\Fempty)=s\catSets$, and when $\Sigma=\Fone$ we get the 
simplicial model category of pointed simplicial sets. Moreover, 
if we put $\Sigma=\bbG$, we recover the usual model category structure on 
the category of simplicial groups.

\nxsubpoint
(Lawvere theorem and Morita equivalence.)

\noindent
Quillen refers in~\cite[2.4]{Quillen}, rem.~1, to a theorem of Lawvere that 
asserts that a category closed under inductive limits and 
possessing a finitely presented projective generator is essentially a 
category of universal algebras, and conversely. Let us restate 
Lawvere theorem in our terms:

\begin{ThD} {\rm (Lawvere, cf.~\cite{Lawvere})}
Let $(\cC,P)$ be a couple consisting of a category~$\cC$ closed under
arbitrary inductive limits and a finitely presented projective
generator~$P$ of~$\cC$. Then the functor
$\Gamma:=\Hom_\cC(P,-):\cC\to\catSets$ admits a left adjoint $L$, thus
defining an {\em algebraic\/} monad~$\Sigma:=\Gamma L$ and a
comparison functor $I:\cC\to\catMod\Sigma$, transforming~$P$ into
$|\Sigma|$. This comparison functor is fully faithful; if all
equivalence relations in~$\cC$ are effective, i.e.\ if $R\cong
X\times_{X/R}X$ for any equivalence relation $R\rightrightarrows X$
in~$\cC$, then $I$ is an equivalence. Conversely, all these conditions
hold for $(\catMod\Sigma,|\Sigma|)$, for an arbitrary algebraic
monad~$\Sigma$.
\end{ThD}

In other words, Lawvere theorem essentially gives us an algebraic description 
of categories with one finitely presented projective generator. Quillen 
uses this theorem as a source of situations where 
theorem~\ptref{th:main.quillen} is applicable; we see that these situations 
correspond {\em exactly\/} to algebraic monads!

As to the proof of this theorem, it goes essentially in the same way
we have shown in~\ptref{sp:ff.zinfflat.mod} that $\ZinfFlat$ is
equivalent to a full subcategory of $\ZinfMod$; actually we might have
proved Lawvere theorem first and apply it afterwards to
deduce~\ptref{sp:ff.zinfflat.mod}.

Observe that if $(\cC,P)$ is a couple satisfying the conditions of Lawvere 
theorem, $(\cC,P^{(n)})$, $P^{(n)}=P^{\oplus n}$, is another. Therefore, 
we can start with an arbitrary algebraic monad~$\Sigma$ and put 
$(\cC,P):=(\catMod\Sigma,\Sigma(n))$; Lawvere theorem yields another 
algebraic monad $M_n\Sigma$, given by $(M_n\Sigma)(m)=
\Hom_\Sigma(\Sigma(n),\Sigma(nm))\cong\Sigma(nm)^n$, and asserts that 
$\catMod{M_n\Sigma}$ is equivalent to~$\catMod\Sigma$. We've got a generalized 
version of Morita equivalence here!

\nxsubpoint
(Free simplicial objects and morphisms.)
We can translate further Quillen results of~\cite[2.4]{Quillen} into 
the language of algebraic monads as well.

Namely, we say that a morphism $i:X\to Z$ in $s(\catMod\Sigma)$ is 
{\em free\/} if the degeneracy diagram of~$Z$, i.e.\ the restriction of $Z$ to 
$\catDelta_-\subset\catDelta$, is isomorphic to the direct sum of the 
degeneracy diagram of~$X$ and some free degeneracy diagram $L_\Sigma(C)$, 
$C:\catDelta_-^0\to\catSets$, and if the morphism of degeneracy diagrams 
$i|_{\catDelta_-}:X|_{\catDelta_-}\to Z|_{\catDelta_-}$ can be identified with 
the natural embedding $X_{\catDelta_-}\to X_{\catDelta_-}\oplus L_\Sigma(C)$.

Any free morphism $i:X\to Z$ is a cofibration in~$s(\catMod\Sigma)$; moreover, 
any morphism $f:X\to Y$ of this category can be factored into 
a free morphism followed by an acyclic fibration, and any cofibration 
is a retract of a free morphism. In particular, any cofibrant simplicial 
$\Sigma$-module is a rectract of a ``free'' simplicial $\Sigma$-module~$F$, 
if we agree to say that $F$ is free wheneven its degeneracy diagram is 
isomorphic to some $L_\Sigma(C)$, $C:\catDelta_-^0\to\catSets$.

\nxpointtoc{Chain complexes and simplicial objects over abelian categories}
Let $\cA$ be an abelian category, $\Ch(\cA)$ be the category of 
non-negative chain complexes over~$\cA$, and $s\cA=\catFunct(\catDelta^0,\cA)$ 
be the category of simplicial objects of~$\cA$. The main result here 
is the {\em Dold--Kan correspondence\/} that actually establishes 
an equivalence between these two categories (cf.~\cite{Dold}).

\nxsubpoint 
(Unnormalized chain complex defined by a simplicial object.)
The easiest way to obtain a chain complex $CX_\cdot$ from 
a simplicial object~$X_\cdot$ is the following. 
Put $(CX)_n:=X_n$ and define $\partial^{CX}_n:X_n\to X_{n-1}$ by 
\begin{equation}
\partial^{CX}_n:=\sum_{i=0}^n (-1)^i d^{n,X}_i\quad.
\end{equation}
Relation $\partial^{CX}_n\partial^{CX}_{n+1}=0$ follows immediately 
from~\eqref{eq:coface.coface}, and in this way we obtain the 
{\em unnormalized chain complex\/} of~$X$.

The whole construction is very classical. For example, 
$H_n(CL_\bbZ X)$ computes the simplicial homology of a simplicial set~$X$, 
considered here as a CW-complex, and $H_n(CL_\bbZ\Sing X)$ computes 
the singular homology~$H_n(X,\bbZ)$ of a topological space~$X$.

\nxsubpoint
(Normalized chain complex of a simplicial object.)
There are several other ways of obtaining a chain complex out of a simplicial 
object~$X$. Consider the {\em degenerate subcomplex\/} $DX\subset CX$, 
generated in each dimension by the images of degeneracy maps: 
\begin{equation}
(DX)_n:=\sum_{i=0}^{n-1} s^{n-1,X}_i(X_{n-1}),\quad\text{where $n>0$.}
\end{equation}
We put $(DX)_0=0$. One checks directly from~\eqref{eq:codeg.coface} that 
$DX$ is indeed a subcomplex of~$CX$. Moreover, it has an increasing filtration
by subcomplexes $F_pDX\subset DX$, defined by taking the sum of images 
of degeneracy maps $s_i$ with $i\leq p$. One can check that all 
$F_pDX/F_{p-1}DX$ are acyclic, hence $DX$ is acyclic, hence 
$H_n(CX)\cong H_n(CX/DX)$ for all $n\geq 0$.

We put $N'X:=CX/DX$; this is one of the equivalent descriptions of the 
{\em normalized chain complex\/} of~$X$. Another one is given by the formula
\begin{equation}
(NX)_n:=\bigcap_{i=1}^n\Ker d^{n,X}_i\subset X_n
\end{equation}
One can check that $NX$ is indeed a subcomplex of~$CX$, and that $DX$ is the 
complement of $NX$ in~$CX$, i.e.\ $X_n=(DX)_n\oplus(NX)_n$ for all $n\geq0$.
Therefore, $NX\cong N'X$, and it is this complex that is usually called 
the normalized chain complex of~$X$. Notice that $\partial^{NX}_n$ actually 
coincides with the restriction of $d^{X,n}_0:X_n\to X_{n-1}$ to 
$(NX)_n\subset X_n$, so we have a simpler formula in this case, and of course 
$H_n(NX)\simto H_n(CX)$ for all $n\geq0$.

\nxsubpoint\label{sp:simpl.obj.from.complex}
(Simplicial object defined by a non-negative chain complex.)
We have just constructed two functors $C,N:s\cA\to\Ch(\cA)$.
Now we'd like to define a functor $K:\Ch(\cA)\to s\cA$ in the opposite 
direction. Let's start with a non-negative chain complex~$A_\cdot$, with 
differentials $\partial^A_n:A_n\to A_{n-1}$, $n\geq1$. First of all, we 
construct a semisimplicial object $\tilde A:\catDelta_+^0\to\cA$ by putting
\begin{equation}
\tilde A_n:=A_n,\quad
d^{n,A}_k:=\begin{cases}\partial^A_n,&\text{if $k=0$}\\0,&\text{if $k>0$}
\end{cases}
\end{equation}
Since $\catDelta_+$ is generated by the face maps $\partial_n^k$, this 
completely determines $\tilde A$, provided relations~\eqref{eq:coface.coface} 
hold. This is indeed the case since the product of two face operators 
$d^{n,A}_kd^{n+1,A}_l$ is zero for all choices of $k$ and~$l$.

Now we define $KA$ to be $J_!\tilde A$, where $J_!$ is the left Kan extension 
of the embedding $J:\catDelta_+\to\catDelta$. This implies that 
$(KA)_n=(J_!\tilde A)([n])=\injlim_{([n]\backslash\catDelta_+)^0}\tilde A$. 
Using the fact that any morphism $\phi:[n]\to[m]$ in $\catDelta$ 
uniquely decomposes into a surjection $\eta:[n]\twoheadrightarrow[p]$ 
followed by an injection $i:[p]\to[m]$, i.e.\ a morphism of $\catDelta_+$, 
we see that the category $[n]\backslash\catDelta_+$ is a disjoint union of 
its subcategories, parametrized by the corresponding value of~$\eta$. 
Each of these subcategories has an initial object, namely, $\eta$ itself, 
and we arrive to the following classical formula:
\begin{equation}\label{eq:KA.n}
(KA)_n=\bigoplus_{\eta:[n]\twoheadrightarrow[p]}A_p
\end{equation}

The left Kan extension description given above defines the face and 
degeneracy operators on~$KA$ as well. More precisely, given a morphism 
$\phi:[m]\to[n]$ in~$\catDelta$ and a component $A_{p,\eta}=A_p$ of 
$(KA)_n$ indexed by some $\eta:[n]\twoheadrightarrow[p]$, the restriction of 
$(KA)(\phi):(KA)_n \to(KA)_m$ maps $A_{p,\eta}$ into $A_{q,\zeta}$, where 
$[m]\stackrel\zeta\twoheadrightarrow[q]\stackrel\psi\hookrightarrow[p]$ 
is the surjective-injective decomposition of $\eta\circ\phi$. The map 
$A_{p,\eta}\to A_{q,\zeta}$ itself is taken to be equal to $\tilde 
A(\psi)$, i.e.\ it is the identity if $p=q$, it is $\partial^A_p$ if 
$q=p-1$ and $\psi=\partial_p^0$, and it is zero otherwise.

\nxsubpoint (Lower degrees of $KA$.)
So if we start from $A=(\cdots\to A_2\stackrel{\partial_2}\to A_1
\stackrel{\partial_1}\to A_0)$, we get
\begin{equation}
\xymatrix{
\cdots A_0\oplus A_1\oplus A_1\oplus A_2
\ar[r]\ar@/_2pt/[r]\ar@/^2pt/[r]^<>(.5){d^2_i}
&A_0\oplus A_1\ar@/_1pt/[r]\ar@/^1pt/[r]^{d^1_i}\ar@/^1pc/@{.>}[l]|{s^1_i}
&A_0\ar@/^1pc/@{.>}[l]|{s^0_0}
}
\end{equation}
Here $s^0_0=\binom 10:A_0\to A_0\oplus A_1$, $d^1_0=(1,\partial^A_1)$, and
$d^1_1=(1,0):A_0\oplus A_1\to A_0$.

In general $A_p$ occurs in $(KA)_n$ exactly $\binom np$ times.

\begin{ThD}\label{th:DK} {\rm (Dold--Kan, cf.~\cite{Dold})}
For any abelian category~$\cA$ the functors $N:s\cA\rightleftarrows
\Ch_{\geq0}(\cA):K$ are adjoint equivalences of categories. In particular, 
$NKA_\cdot\cong A_\cdot$ for any $A_\cdot\in\Ch\cA$. 
Moreover, two chain maps $f,g:A\rightrightarrows B$ are chain homotopic iff 
$K(f),K(g):KA\rightrightarrows KB$ are simplicially homotopic.
\end{ThD}

\nxsubpoint
An immediate consequence is that $A\mapsto DKA$ is a universal 
construction that transforms any chain complex~$A$ over any abelian 
category~$\cA$ into an acyclic complex. One can show by applying this 
construction to an appropriate ``universal'' complex that this 
necessarily implies that $DKA$ is (chain) homotopic to zero, and, 
moreover, a chain homotopy $H:1_{DKA}\to 0$ can be defined by means of 
a universal formula. Therefore, $DX$ is (chain) homotopic to zero, 
and $CX$ is homotopic to $NX$ for any simplicial object $X$ of~$\cA$. 

\nxsubpoint (Dual Dold--Kan correspondence.)
Of course, we have a dual result: the category of cosimplicial objects 
$c\cA$ is equivalent to the category of non-negative cochain complexes 
over~$\cA$. The dual of $N'$ seems to be more convenient than that of~$N$: 
we get $(N'X^\cdot)^n=\bigcap_{k=0}^{n-1}\Ker\sigma_{n-1,X}^k\subset X^n$,
with $d_{N'X}^n$ given by the restriction of $\sum_{k=0}^{n+1}
(-1)^k\partial_{n+1,X}^k$.

\nxsubpoint (Model category structure on~$\Ch\cA$.)
Now suppose that $\cA$ is an abelian category satisfying the conditions 
of Quillen theorem~\ptref{th:main.quillen}, e.g.\ $\cA=\catMod R$ for 
a classical ring~$R$. Then we have a natural simplicial model category 
structure on $s\cA$, and we can transfer it to~$\Ch\cA$ by means of the 
Dold--Kan correspondence.

It turns out that a chain map $f:K\to L$ in $\Ch\cA$ is a weak equivalence 
iff all $H_n(f)$, $n\geq0$, are isomorphisms, a fibration iff 
all $f_n$, $n>0$, are epimorphisms, and a cofibration iff all $f_n$, 
$n\geq0$, are monomorphisms with a projective cokernel. This actually 
defines a model category structure on~$\Ch\cA$ whenever~$\cA$ is an abelian 
category with sufficiently many projective objects.

Moreover, we have a similar model category structure on the category 
$\Ch_{\gg-\infty}(\cA)$ of chain complexes bounded from below. In this case 
the weak equivalences are again the quasi-isomorphisms, the fibrations are 
the (componentwise) epimorphic maps, and the cofibrations are again 
the componentwise monic maps with projective cokernels. 

\nxsubpoint (Homology as homotopy groups.)
Given any simplicial group~$G$, one define its normalization $NG$ by 
$(NG)_n:=\bigcap_{i>0}\Ker d^n_i$, and define maps $d=d_n:(NG)_n\to(NG)_{n-1}$ 
by taking the restriction of~$d^n_0$ to $(NG)_n\subset G_n$. We obtain 
a complex of (non-abelian) groups $NG$, so we  can compute the 
{\em (Moore) homotopy groups\/} of~$G$ by
\begin{equation}
\pi_n(G):=\Ker\bigl(d:(NG)_n\to(NG)_{n-1}\bigr)\big/
\Im\bigl(d:(NG)_{n+1}\to(NG)_n\bigr)
\end{equation}
One can show that these groups are indeed canonically isomorphic to the 
homotopy groups $\pi_n(G,e)=\pi_n(|G|,e)$ of the underlying simplicial 
set of~$G$. Notice that in this case all $\pi_n(G)$ are groups, abelian 
if~$n>0$.

Now suppose $A$ is a simplicial abelian group, or more generally, a 
simplicial $R$-module, with $R$ a classical ring. We obtain immediately 
$\pi_n(A,0)\cong\pi_n(A)\cong H_n(NA)\cong H_n(CA)$, i.e.\ 
{\em the Dold--Kan correspondence induces a canonical isomorphism between 
homology groups of a chain complex and homotopy groups of the corresponding 
simplicial $R$-module.}

\nxsubpoint 
(Bisimplicial objects and bicomplexes.)
Now let's consider the categories of non-negative chain bicomplexes 
over an abelian category~$\cA$, and the category 
$ss\cA=\catFunct(\catDelta^0\times\catDelta^0,\cA)$ of bisimplicial objects 
of~$\cA$. Such bisimplicial objects $X_{\cdot\cdot}$ can be described as 
collections $\{X_{pq}\}_{p,q\geq0}$ of objects of~$\cA$, together with 
``vertical'' and ``horizontal'' face and degeneracy operators 
$d^{I,pq}_i:X_{pq}\to X_{p-1,q}$, $d^{II,pq}_i:X_{pq}\to X_{p,q-1}$ and 
so on, vertical and horizontal operators being required to commute between
themselves.

Clearly, the Dold--Kan correspondence extends to this situation 
(e.g.\ we can treat a bicomplex as a complex over the category of complexes, 
and similarly for the simplicial objects). In this case we have two 
commuting normalization functors $N_I$ and $N_{II}$, and 
the ``true'' normalization functor is their composite: $N=N_IN_{II}=N_{II}N_I$.

\nxsubpoint 
(Diagonal simplicial objects and totalizations of bicomplexes.)
Let's start with a bisimplicial object~$A=A_{\cdot\cdot}$ over $\cA$. 
We can compute the corresponding chain bicomplex $CA$ or $NA$, and 
take its totalization $\Tot(CA)$ or $\Tot(NA)$; these two chain complexes 
are homotopy equivalent. On the other hand, we can construct the 
{\em diagonal simplicial object\/} $\diag A$ of~$A$ by composing 
$A:\catDelta^0\times\catDelta^0\to\cA$ with the diagonal functor 
$\catDelta\to\catDelta\times\catDelta$. Clearly, $(\diag A)_n=A_{nn}$, 
and $d^{\diag A}_i=d^I_id^{II}_i$. After that we can consider the corresponding
unnormalized or normalized chain complex $C(\diag A)$ or $N(\diag A)$.

\begin{ThD}\label{th:EZ} {\rm (Eilenberg--Zilber; cf.~\cite[2]{DP})}
Let $A$ be a bisimplicial object over an abelian category~$\cA$. Then 
chain complexes $C(\diag A)$ and $\Tot CA$ are canonically 
quasi-isomorphic and even homotopy equivalent; in particular, their homology 
is canonically isomorphic: $\pi_n\diag A\cong H_n(C(\diag A))\cong 
H_n(\Tot CA)$ for all $n\geq0$.
\end{ThD}
Morally, this means that computing the diagonal of a bisimplicial object 
corresponds to computing the totalization of a bicomplex.

\nxsubpoint (Alexander--Whitney and shuffle maps.)
It's useful to have explicit constructions for the quasi-isomorphisms 
(even homotopy equivalences, because of the universality of all constructions)
$f:C(\diag A)\to\Tot(CA)$ and $\nabla:\Tot(CA)\to C(\diag A)$ 
implied in the above theorem.

The map $f:C(\diag A)\to\Tot(CA)$, $f_n:A_{nn}\to\bigoplus_{p+q=n}A_{pq}$ 
is called the {\em Alexander--Whitney map\/} (cf.~\cite[8.5.4]{Weibel} 
or~\cite[2]{DP}). 
Its component $f_{pq}:A_{nn}\to A_{pq}$ equals 
$d^I_{p+1}\cdots d^I_nd^{II}_0\cdots d^{II}_0$. The map 
$\nabla:\Tot(CA)\to C(\diag A)$ is called the {\em shuffle map}. 
Its components $\nabla_{pq}:A_{pq}\to A_{nn}$, $n=p+q$, are given by the sum
\begin{equation}
\nabla_{pq}=\sum_\mu\sgn(\mu) s^I_{\mu(n-1)}\cdots s^I_{\mu(p)}
s^{II}_{\mu(p-1)}\cdots s^{II}_{\mu(0)}
\end{equation}
over all $(p,q)$-shuffles $\mu$, i.e.\ all permutations $\mu$ of the set 
$\{0,1,\ldots,n-1\}$, such that $\mu(0)<\mu(1)<\cdots<\mu(p-1)$ and 
$\mu(p)<\mu(p+1)<\cdots<\mu(n-1)$.

It is possible to check that $f$ and $\nabla$ are indeed chain maps, 
inducing isomorphisms between homologies of corresponding complexes, inverse 
to each other. Once this is shown, the universality of all constructions 
involved implies that $f$ and $\nabla$ are homotopy equivalences between 
$C(\diag A)$ and $\Tot(CA)$, homotopy inverse to each other, so one might try 
to construct explicit homotopies $f\nabla\sim1$, $\nabla f\sim 1$ instead, 
and deduce the Eilenberg--Zilber theorem from this (cf.~\cite[2]{DP}). 

\nxsubpoint\label{sp:appl.class.der.tens} 
(Application: derived tensor product.)
Suppose we have two non-negative chain complexes $X$ and $Y$ over 
abelian category $\cA=\catMod R$, $R$ a commutative ring, and 
we want to compute the derived tensor product $X\Lotimes Y$. Classically 
we have to take projective resolutions, i.e.\ cofibrant replacements 
$P\to X$ and $Q\to Y$, compute the bicomplex $P\otimes Q$ given by 
$(P\otimes Q)_{ij}:=P_i\otimes Q_j$, and put $X\Lotimes Y:=\Tot(P\otimes Q)$.

However, Eilenberg--Zilber theorem implies that we can compute derived tensor 
products in another way. Namely, let's start from two simplicial 
$R$-modules $\tilde X=KX$ and~$\tilde Y=KY$, choose some cofibrant replacements
$\tilde P\to\tilde X$ and $\tilde Q\to\tilde Y$. 
Then $C\tilde P$ is chain homotopic to $N\tilde P$, and similarly for 
$C\tilde Q$, hence $\Tot(C\tilde P\otimes C\tilde Q)$ computes 
$X\Lotimes Y$ as well. Applying Eilenberg--Zilber we see that 
$CT$ with $T:=\diag(\tilde P\otimes\tilde Q)$ computes this derived tensor 
product as well, i.e.\ we just have to apply $\otimes_R$ componentwise: 
$T_n:=\tilde P_n\otimes_R\tilde Q_n$. This gives us a very convenient way of 
deriving bifunctors, additive or not, just by applying them 
componentwise to appropriate fibrant or cofibrant replacements (cf.~\cite{DP}).

\nxpointtoc{Simplicial $\Sigma$-modules}
Let $\Sigma$ be any algebraic monad, e.g.\ a generalized ring. 
Recall that the category $s(\catMod\Sigma)$ of simplicial $\Sigma$-modules 
has a cofibrantly generated simplicial (closed) model category structure, 
characterized as follows: a morphism $f:X\to Y$ in $s(\catMod\Sigma)$ 
is a fibration (resp.\ a weak equivalence) iff $\Gamma_\Sigma(f)$ is 
a fibration (resp.\ a weak equivalence) of underlying simplicial sets, 
and $f$ is a cofibration iff it has the LLP with respect to all acyclic 
fibrations (cf.~\ptref{sp:cat.simpl.mod}). When $\Sigma$ is a classical ring, 
$s(\catMod\Sigma)$ is equivalent to the category of non-negative 
chain complexes of $\Sigma$-modules (cf.~\ptref{th:DK}).

The corresponding homotopy category $\Ho s(\catMod\Sigma)$ will be also 
denoted by $\cD^{\leq0}(\catMod\Sigma)$ or $\cD^{\leq0}(\Sigma)$; when 
$\Sigma$ is a classical ring, this category is indeed equivalent to the 
subcategory of $\cD(\catMod\Sigma)$ consisting of (cochain) complexes 
concentrated in non-positive degrees.

\nxsubpoint (Algebraic simplicial monads.)
Notice that the simplicial extension $s\Sigma:s\catSets\to s\catSets$ of 
any monad $\Sigma$ on~$\catSets$ is a monad again, so we can consider the 
category $\catMod{(s\Sigma)}=(s\catSets)^{s\Sigma}$. However, since 
$((s\Sigma)(X))_n=\Sigma(X_n)$, a $s\Sigma$-structure on a simplicial 
set~$X$ is exactly a collection of $\Sigma$-structures on each~$X_n$, 
compatible with the face and degeneracy operators, i.e.\
$\catMod{(s\Sigma)}\cong s(\catMod\Sigma)$, so we can write simply 
$s\catMod\Sigma$.

One can define {\em algebraic simplicial monads\/} $\Xi$ on $s\catSets$ 
by requiring $\Xi$ to commute with filtered inductive limits of simplicial 
sets, and by imposing the following ``continuity condition'': 
the $n$-th truncation of $\Xi(X)$ should depend only on the $n$-th 
truncation of~$X$, i.e.\ $\sk_n\Xi(\sk_n X)\to\sk_n\Xi(X)$ has to 
be an isomorphism for all simplicial sets~$X$ and all $n\geq0$.

For example, $s\Sigma$ is an algebraic simplicial monad, 
for any algebraic monad $\Sigma$ on~$\catSets$. More generally, any simplicial 
algebraic monad defines an algebraic simplicial monad, 
but not the other way around.

We might extend to the case of algebraic simplicial monads $\Xi$ 
the theory developed in previous parts of this work for algebraic monads,
e.g.\ define ``simplicial sets of morphisms'' $\iHom_\Xi(X,Y)$, 
tensor products (for a commutative~$\Xi$), as well as construct a natural 
simplicial model category structure on $(s\catSets)^\Xi$ and prove some 
variants of most results that follow. We won't pursue this approach further, 
because it is too general for our purpose; 
however, some trace of it can be found in our 
notations, e.g.\ $\Hom_{s\Sigma}$ for the simplicial set of homomorphisms
of simplicial $\Sigma$-modules.

Notice that this approach might be useful for arithmetic geometry, even if 
it is beyond the scope of present work. For example, considering the 
cone of $\Zinfty\stackrel{f}\to\Zinfty$ as a new algebraic simplicial monad 
is a natural way to construct a quotient $\Zinfty/(f)$ different from 
$\Zinfty/\gm_\infty$. Therefore, if we want to consider divisors like 
$-\log|f|\cdot[\infty]$ as effective divisors given by closed subschemes of 
$\Spec\Zinfty$, we'll need such a way of constructing ``fine'' quotients. 

Now let's return to our much more specific situation.

\begin{PropD}\label{prop:ex.Lder.basech}
Let $\rho:\Sigma\to\Xi$ be an algebraic monad homomorphism. Consider the 
pair of adjoint functors $\rho^*=s\rho^*:s(\catMod\Sigma)\to s(\catMod\Xi)$ 
and $\rho_*:s(\catMod\Xi)\to s(\catMod\Sigma)$. We claim that this is 
a pair of Quillen functors (cf.~\ptref{def:Q.funct}), hence their derived 
functors $\dL\rho^*:\cD^{\leq0}(\Sigma)\to\cD^{\leq0}(\Xi)$ and
$\dR\rho_*:\cD^{\leq0}(\Xi)\to\cD^{\leq0}(\Sigma)$ exist and are adjoint to 
each other; they can be computed by applying $\rho^*$, resp.\ $\rho_*$ 
to arbitrary cofibrant, resp.\ fibrant replacements (cf.~\ptref{th:Q.funct}). 
Actually $\rho_*$ preserves all weak equivalences, hence $\dR\rho_*=\rho_*$.
\end{PropD}
{\bf Proof.} First of all, the definition of model category structure 
on $s(\catMod\Sigma)$ together with equality $\Gamma_\Xi=\Gamma_\Sigma\circ
\rho_*$ imply that $\rho_*$ preserves all fibrations and weak equivalences, 
hence also acyclic fibrations. Using adjointness of $\rho^*$ and $\rho_*$ 
we deduce that $\rho^*$ preserves arbitrary cofibrations and acyclic 
cofibrations, since $\rho^*(f)$ has the required LLP for any cofibration,
resp.\ acyclic cofibration $f:X\to Y$ in~$s(\catMod\Sigma)$.
Now it remains to check that $\rho^*$ preserves weak equivalences between 
cofibrant objects. We know this already for acyclic cofibrations, 
so we can finish the proof by applying the following lemma:

\begin{LemmaD}\label{l:acof.cof.suff}
a) Let $F:\cC\to\cD$ be a functor from a model category~$\cC$ that 
transforms acyclic cofibrations between cofibrant objects into isomorphisms.
Then the same is true for all weak equivalences between cofibrant objects
of~$\cC$.

b) Let $F:\cC\to\cD$ be a functor between two model categories that 
transforms acyclic cofibrations between cofibrant objects of~$\cC$ into 
weak equivalences in~$\cD$. Then $F$ preserves all weak equivalences between 
cofibrant objects of~$\cC$, hence $F$ has a left derived functor
$\dL F$ by~\ptref{prop:ex.left.der}.
\end{LemmaD}
{\bf Proof.}
b) follows immediately from a) applied to $\gamma_D\circ F:\cC\to\Ho\cD$ and 
from the fact that $\phi$ is a weak equivalence in~$\cD$ iff 
$\gamma_D(\phi)$ is an isomorphism, all our model categories being closed, 
so let's prove a).

So let $p:Q\to P$ be a weak equivalence between two cofibrant objects. By 
(CM5) and (CM3) we can decompose $p$ into an acyclic cofibration 
$i:Q\to P'$ followed by an acyclic fibration $\pi:P'\to P$. Clearly, $P'$ 
is itself cofibrant, and $F(i)$ is an isomorphism by assumption, so we are 
reduced to proving that $F(\pi)$ is an isomorphism for any acyclic fibration 
$\pi:P'\to P$ between cofibrant objects.

Since $P$ is cofibrant and $\pi$ is an acyclic fibration, any morphism 
$P\to P$ can be lifted to a morphism $P\to P'$; applying this to $\id_P$ 
we obtain 
a section $\sigma:P'\to P$. Clearly, $\sigma$ is a weak equivalence since both 
$\pi$ and $\pi\circ\sigma=\id_P$ are, so it can be decomposed into an 
acyclic cofibration $\tau:P\to P''$ followed by an acyclic fibration 
$\pi':P''\to P'$, and we have $\pi\circ\pi'\circ\tau=\id_P$. 
Since $P$ is cofibrant and $\tau:P\to P''$ is an acyclic cofibration, 
$P''$ is also cofibrant and $F(\tau)$ is an isomorphism by assumption, hence 
$F(\pi)\circ F(\pi')=F(\tau)^{-1}$ is an isomorphism as well.

Now $\pi':P''\to P'$ being an acyclic fibration between cofibrant objects, 
we can repeat the same construction again, obtaining a cofibrant object~$P'''$,
an acyclic fibration $\pi'':P'''\to P''$ and an acyclic cofibration 
$\tau':P'\to P'''$, such that $\pi'\circ\pi''\circ\tau'=\id_{P'}$, and 
$F(\tau')$ is an isomorphism, hence $F(\pi')\circ F(\pi'')$ as well.
We see that both $F(\pi)\circ F(\pi')$ and $F(\pi')\circ F(\pi'')$ are 
isomorphisms, hence $F(\pi)$ is an isomorphism by Lemma~\ptref{l:3iso} below, 
q.e.d.

\begin{LemmaD}\label{l:3iso}
If $\cdot\stackrel u\to\cdot\stackrel v\to\cdot\stackrel w\to\cdot$ 
are three composable morphisms of an arbitrary category, such that 
$vu$ and $wv$ are isomorphisms, then $u$, $v$ and $w$ are isomorphisms
themselves.
\end{LemmaD}
{\bf Proof.} Indeed, $v$ has both a left inverse $(wv)^{-1}w$ and a 
right inverse $u(vu)^{-1}$, hence it is invertible with $v^{-1}=(wv)^{-1}w=
u(vu)^{-1}$, hence $u=v^{-1}(vu)$ and $w=(wv)v^{-1}$ are invertible as well.

\nxsubpoint (Functoriality of derived base change.)
Now suppose that $\sigma:\Xi\to\Lambda$ is another algebraic 
monad homomorphism. Since $\rho^*$ transforms cofibrant objects into cofibrant 
objects, we see that $\dL(\sigma\circ\rho)^*=\dL(\sigma^*\circ\rho^*)\cong
\dL\sigma^*\circ\dL\rho^*$ (cf.~\ptref{prop:ex.left.der},b), and similarly 
$\dR(\sigma\circ\rho)_*=\dR(\rho_*\circ\sigma_*)\cong
\dR\rho_*\circ\dR\sigma_*$.

\nxsubpoint (Application: abelianization and cohomology.)
Let $\cC$ be a model category, $\cC_{ab}$ be the category of abelian group 
objects in~$\cC$, and suppose that the inclusion functor 
$i:\cC_{ab}\to\cC$ has a left adjoint~$ab$, and that $\cC_{ab}$ has a 
nice model category structure, satifying the conditions of~\cite[2.5]{Quillen}.
Recall that these conditions are fulfilled if $\cC_{ab}$ is equivalent to 
$s\cA$ for an abelian category~$\cA$, and if the derived functors 
$\dL ab$ and $\dR i$ exist. Then Quillen defines the 
{\em generalized cohomology\/}~$H^q(X,A)$ of an object~$X$ of~$\Ho\cC$ with 
coefficients in an object $A$ of $\Ho\cC_{ab}$ by
$$H^q(X,A):=[\dL ab(X),\Omega^{q+N}\Sigma^NA],\quad\text{for any $N\gg0$.}$$
Now if $\cC=s(\catMod\Sigma)$, then $\cC_{ab}=s\catMod{(\Sigma\otimes\bbZ)}$, 
and the existence of $ab$, $i$ and their derived follows from our general
results applied to algebraic monad homomorphism $\Sigma\to\Sigma\otimes\bbZ$.

\nxsubpoint (Base change and suspension.)
Recall that for any pointed model category we have a suspension functor 
$\Sigma:\Ho\cC\to\Ho\cC$; the suspension of a cofibrant~$X$ can be 
computed as the cofiber of $X\sqcup X\to X\times I$, where $X\times I$ 
is any cylinder object for~$X$. If $\cC$ is simplicial, $X\otimes\Delta(1)$ 
is a cylinder object for an acyclic~$X$, so we get a functorial way of 
computing $\Sigma X$ by means of the following cocartesian square:
\begin{equation}\label{eq:constr.susp}
\xymatrix@C+10pt{
X\sqcup X\ar[r]\ar[d]&X\otimes\Delta(1)\ar@{-->}[d]\\ 
0\ar@{-->}[r]&\Sigma X}
\end{equation}
Let's apply this to $s(\catMod\Sigma)$, for any algebraic monad~$\Sigma$ 
with zero. Then $(X\otimes\Delta(1))_n$ is the 
sum of $n+2$ copies of~$X_n$, parametrized by monotone maps 
$\eta:[n]\to[1]$, and the embeddings $X_n\rightrightarrows
(X\otimes\Delta(1))_n$ are just the embeddings into the components 
corresponding to constant~$\eta$'s. We deduce that
\begin{equation}
(\Sigma X)_n \cong \bigoplus_{\eta:[n]\twoheadrightarrow[1]} X_n\quad.
\end{equation}
Now if $\rho:\Sigma\to\Xi$ is a homomorphism of algebraic monads with zero, 
the base change functor $\rho^*$ preserves direct sums and cofibrant objects; 
hence $\Sigma\rho^*X\cong\rho^*\Sigma X$ for any cofibrant~$X$ in
$s(\catMod\Sigma)$, and $\Sigma\circ\dL\rho^*\cong\dL\rho^*\circ\Sigma$ 
on the level of derived functors (we might write $\Sigma$ for the functor 
on~$s(\catMod\Sigma)$ defined by~\eqref{eq:constr.susp}, and $\dL\Sigma$ 
for the suspension functor on the homotopic category). By adjointness we get 
$\dR\rho_*\circ\Omega\cong\Omega\circ\dR\rho_*$. This is actually a 
consequence of a general fact: 
Theorem~3 of \cite[1.4]{Quillen} implies that $\dL\rho^*$ preserves 
suspensions and cofibration sequences, while $\dR\rho_*$ preserves 
loop objects and fibration sequences.

\nxsubpoint\label{sp:mapping.cone.cyl} (Mapping cones and cylinders.)
Let $f:X\to Y$ be a morphism in any simplicial model category~$\cC$. 
We define the {\em (mapping) cylinder\/}~$Cyl(f)$ of~$f$ by 
the following cocartesian square:
\begin{equation}
\xymatrix{
X\otimes\{0\}\ar[r]^{f}\ar[d]&Y\otimes\{0\}\ar@{-->}[d]\\
X\otimes\Delta(1)\ar@{-->}[r]&Cyl(f)}
\end{equation}
Here $\{e\}$, $e=0,1$, denotes the corresponding vertex of~$\Delta(1)$, 
i.e.\ the subobject defined by $\Delta(\sigma_1^e):\Delta(0)\to\Delta(1)$. 
Clearly, in $s(\catMod\Sigma)$ we get $Cyl(f)_n\cong Y_n\oplus 
X_n^{\oplus(n+1)}$; we basically construct $(X\otimes\Delta(1))_n$ 
and replace the first copy of $X_n$ with~$Y_n$. 

When~$\cC$ is pointed, i.e.\ it has a zero object, 
we define the {\em (mapping) cone\/}~$C(f)$ of~$f$ as the cofiber of 
$X\cong X\otimes\{1\}\to X\otimes\Delta(1)\to Cyl(f)$. For $s(\catMod\Sigma)$ 
we get $C(f)_n\cong Y_n\oplus X_n^{\oplus(n)}$. One can check that these 
definitions of cone and cylinder correspond to those known from 
homological algebra via the Dold--Kan correspondence when $\Sigma$ is 
additive.

The above construction is especially useful when $X$ and $f$ are cofibrant; 
we can extend it for any morphism $f:X\to Y$ in $\Ho\cC$ by replacing $X$ 
and $Y$ by appropriate cofibrant representatives. Then $Y\to Cyl(f)$ is a 
weak equivalence for any~$f$, and $X\to Y\to C(f)$ is a cofibration sequence 
in~$\Ho\cC$ for any morphism~$f$, so for any $Z\in\Ob\Ho\cC$ we get a 
long exact sequence of homotopy:
\begin{multline}
\cdots\to[\Sigma^2X,Z]\to[\Sigma C(f),Z]\to[\Sigma Y,Z]\to\\
\to[\Sigma X, Z]\to[C(f),Z]\to[Y,Z]\to[X,Z]\to0
\end{multline}

Now an obvious but important statement is that {\em $\dL\rho^*$ preserves 
mapping cones and cylinders.} 

\nxsubpoint (Computation of $\pi_0$.)
Recall that for any simplicial set~$A$ we can compute its 
``set of connected components'' $\pi_0A=\pi_0(|A|)$ as the cokernel of 
$A_1\rightrightarrows A_0$ in the category of sets. We know that 
any weak equivalence induces an isomorphism of $\pi_0$'s, so we obtain a 
functor $\pi_0:\Ho s\catSets\to\catSets$.

Now suppose that $A$ is a simplicial $\Sigma$-module. Since weak equivalences 
of simplicial $\Sigma$-modules are defined in terms of underlying simplicial 
sets, we obtain a well-defined functor $\pi_0:\Ho s(\catMod\Sigma)\to\catSets$.
We claim that the surjection $A_0\twoheadrightarrow\pi_0A=A_0/\negthinspace
\sims$ is always compatible with the $\Sigma$-structure on~$A_0$, so that 
$\pi_0A$ inherits a canonical $\Sigma$-structure, and we get a functor 
$\pi_0:\Ho s(\catMod\Sigma)\to\catMod\Sigma$. Indeed, for any 
operation $t\in\Sigma(n)$, and any elements $x_1\simss y_1$, $x_2=y_2$, \dots, 
$x_n=y_n$ of~$A_0$ one can find by definition elements 
$h_i\in A_1$ with $d_0(h_i)=x_i$, $d_1(h_i)=y_i$ (when $x_i=y_i$ we take
$h_i:=s_0(x_i)$), put 
$h:=t(h_1,\ldots,h_n)$, and obtain $d_0(h)=x:=t(x_1,\ldots,x_n)$, 
$d_1(h)=y:=t(y_1,\ldots,y_n)$ since the $d_k$ are $\Sigma$-homomorphisms.
We conclude that $x\simss y$, so the equivalence relation $\sims$ generated 
by~$\simss$ is compatible with the $\Sigma$-structure as well.

On the other hand, any $\Sigma$-module~$A$ can be considered as a constant 
simplicial $\Sigma$-module that will be usually denoted by the same letter.
This yields a functor $\catMod\Sigma\to s(\catMod\Sigma)$ in the opposite 
direction, clearly right adjoint to~$\pi_0$, and $\pi_0A\cong A$ for any 
constant simplicial set. We deduce that the functor $\catMod\Sigma\to
s(\catMod\Sigma)$ is fully faithful. Notice that for any simplicial 
$\Sigma$-module we get a morphism $A\to\pi_0A$.

One can check that the higher homotopic invariants $\pi_nA$, $n\geq0$, 
also have a natural $\Sigma$-structure, commuting with their group structure,
at least if $\Sigma$ has a zero. One can use $\pi_n(A)\cong\pi_0(\Omega^nA)$ 
to show this.

\nxsubpoint
Since $\pi_0A=\Coker(A_1\rightrightarrows A_0)$ both in $\catMod\Sigma$ 
and $\catSets$, and $\rho^*$ is right exact, we see that 
$\pi_0\rho^*P\cong\rho^*\pi_0P$ 
for any cofibrant simplicial $\Sigma$-module~$P$
and any algebraic monad homomorphism $\rho:\Sigma\to\Xi$, hence 
$\pi_0\circ\dL\rho^*\cong\rho^*\circ\pi_0:\cD^{\leq0}(\Sigma)\to
\catMod\Xi$. In particular, for any $\Sigma$-module~$M$ we get
$\pi_0\dL\rho^*M\cong\rho^*M$, i.e.\ we have a canonical map $\dL\rho^*M\to
\rho^*M$ inducing an isomorphism of~$\pi_0$'s. This map is a weak equivalence 
iff all $\pi_n(\dL\rho^*M,x_0)=0$, $n\geq1$, for all choices
of base point $x_0\in\pi_0\dL\rho^*M=\rho^*\pi_0M$; then we can write 
$\dL\rho^*M\cong M$.

\nxsubpoint (Cohomological flatness.) 
We say that an algebraic monad homomorphism $\rho:\Sigma\to\Xi$ is 
{\em cohomologically flat\/} if $\rho^*$ preserves weak equivalences, i.e.\ 
if $\dL\rho^*=\rho^*$. This implies $\pi_n\dL\rho^*M=0$ for all $n\geq1$ 
and all $\Sigma$-modules~$M$. Notice that usual flatness 
(exactness of $\rho^*$) implies cohomological flatness. Indeed, since 
$\rho^*$ preserves acyclic cofibrations, all we have to check is that 
it preserves acyclic fibrations. But $p:X\to Y$ is an acyclic fibration iff 
all maps $X_n\to(\cosk_{n-1}X)_n\times_{(\cosk_{n-1}Y)_n}Y_n$ are surjective, 
and $(\cosk_{n-1}X)_n\cong\Hom(\dot\Delta(n),X)$ can be expressed in terms of 
finite projective limits of components of~$X$, $\dot\Delta(n)$ being a
finite inductive limit of standard simplices. Therefore, the property of~$p$ 
to be an acyclic fibration is preserved by any exact functor.

\nxsubpoint\label{sp:stab.ht.cat} (Stable homotopic category.)
There are several ways of transforming~$\Ho s\catSets=\cD^{\leq0}(\Fempty)$ 
(more precisely, the corresponding pointed category $\cD^{\leq0}(\Fone)$)
into a category on which
$\Omega$ and $\Sigma$ become autoequivalences, inverse to each other. 
The simplest way consists in considering the category of couples 
$(X,n)$, $X\in\Ho s\catSets$, $n\in\bbZ$, with $\Hom((X,n),(Y,m)):=
\injlim_{N\gg0}[\Sigma^{n+N}X,\Sigma^{m+N}Y]$. If we apply this construction 
to $\cD^{\leq0}(R)$ for a classical ring~$R$, we obtain $\cD^-(R)$, so 
we use the same notation $\cD^{-}(\Lambda)$ for the stable homotopy category
of $s(\catMod\Lambda)$ for any algebraic monad~$\Lambda$ with zero. 

Notice that 
$\cD^-(\Lambda)$ will be a {\em triangulated\/} category for any~$\Lambda$. 
The proof is exactly the same as for the stable homotopic category of 
simplicial sets: the abelian group structure on stable homomorphism 
sets comes from the abelian cogroup object structure on all $\Sigma^nX$, 
$n\geq2$, and the distinguished triangles arise from cofibration sequences. 

We get some {\em stable homotopy groups\/} 
$\pi_n(X):=[\Sigma^n\Delta(0),X]_{st}$; 
they constitute a homological functor on $\cD^-(\Lambda)$ with values in 
$\catMod\bbZ$ (actually in $\catMod{(\Lambda\otimes\bbZ)}$, if $\Lambda$ 
is commutative). However, these stable homotopy groups are not so easy to 
compute, and they do not reflect all cohomological structure of 
$\cD^-(\Lambda)$. For example, one can construct an object~$P$ in
$\cD^-(\Zinfty)$ which has $\pi_n(P)=0$ for all $n\in\bbZ$, but $P\neq0$.
(It is still unknown whether such examples exist over~$\Fone$.)

\nxpointtoc{Derived tensor product}\label{p:der.tensprod}
The aim of this subsection is to construct the {\em derived tensor product\/}
$\Lotimes=\dL\otimes_\Lambda:\cD^{\leq0}(\Lambda)\times
\cD^{\leq0}(\Lambda)\to\cD^{\leq0}(\Lambda)$
for any generalized ring, i.e.\ commutative algebraic monad~$\Lambda$. 
Before doing this we have to make some general remarks about 
deriving bifunctors.

\nxsubpoint (Simplicial extension of a bifunctor.)
Recall that any functor $F:\cC\to\cD$ has a natural simplicial extension
$sF:s\cC\to s\cD$, defined by composing functors $\catDelta^0\to\cC$ with~$F$. 
Now suppose that $\cC=\cC_1\times\cC_2$, i.e.\ we have a bifunctor 
$F:\cC_1\times\cC_2\to\cD$. Notice that $s(\cC_1\times\cC_2)=\catFunct
(\catDelta^0,\cC_1\times\cC_2)\cong\catFunct(\catDelta^0,\cC_1)\times
\catFunct(\catDelta^0,\cC_2)=s\cC_1\times s\cC_2$, so the simplicial extension
$sF$ of~$F$ is indeed a functor $s\cC_1\times s\cC_2\to s\cD$ as one would 
expect. Of course, we usually write just $F$ instead of~$sF$. For any two 
simplicial objects $X$ and $Y$ of $\cC_1$ and $\cC_2$, respectively, we get
$(sF(X,Y))_n=F(X_n,Y_n)$, with the face and degeneracy operators given by 
formulas like $d_i^{sF(X,Y)}=F(d_i^X, d_i^Y)$.

In particular, the tensor product $\otimes_\Lambda:\catMod\Lambda\times
\catMod\Lambda\to\catMod\Lambda$ extends to a bifunctor 
$s\catMod\Lambda\times s\catMod\Lambda\to s\catMod\Lambda$.

\nxsubpoint\label{sp:prod.mod.cat} (Product of model categories.)
The next thing we need to derive a bifunctor $F:\cC_1\times\cC_2\to\cD$, 
with $\cC_1$, $\cC_2$ and $\cD$ model categories, is a model category 
structure on the product $\cC_1\times\cC_2$. The natural choice is to 
declare $f=(f_1,f_2):(X_1,X_2)\to(Y_1,Y_2)$ a fibration (resp.\ cofibration, 
weak equivalence) in $\cC_1\times\cC_2$ iff both its components 
$f_i:X_i\to Y_i$, $i=1,2$, are fibrations (resp\dots) in~$\cC_i$. 
One checks immediately that this is indeed a model category structure, and 
that $\Ho(\cC_1\times\cC_2)\cong\Ho\cC_1\times\Ho\cC_2$, so we can safely 
consider the left and right derived functors of $F$ with respect to this 
model structure: they'll be some functors $\dL F,\dR F:\Ho\cC_1\times\Ho\cC_2
\to\Ho\cD$. Notice that if $\cC_1$ and $\cC_2$ are simplicial model categories,
we can define a compatible simplicial structure on $\cC_1\times\cC_2$ as well
by setting $(X_1,X_2)\otimes K:=(X_1\otimes K,X_2\otimes K)$ for any 
$X_i\in\Ob\cC_i$ and $K\in\Ob s\catSets$.

In particular, the above considerations enable us to define the 
left derived tensor product $\Lotimes:\cD^{\leq0}(\Lambda)\times
\cD^{\leq0}(\Lambda)\to\cD^{\leq0}(\Lambda)$. 
We just have to show its existence.

\nxsubpoint\label{sp:ex.lder.bif}
(Criterion for existence of left derived bifunctors.)
Now we'd like to combine~\ptref{prop:ex.left.der} and~\ptref{l:acof.cof.suff} 
with the above considerations. We obtain the following proposition:

\begin{Propz}
Let $\cC_1$, $\cC_2$ and $\cD$ be model categories, $F:\cC_1\times\cC_2\to\cD$ 
be a bifunctor. Then the following conditions are equivalent:
\begin{itemize}
\item[a)] $F$ transforms weak equivalences between cofibrant objects of
$\cC_1\times\cC_2$ into weak equivalences in~$\cD$.
\item[b)] $F$ transforms acyclic cofibrations between cofibrant objects of
$\cC_1\times\cC_2$ into weak equivalences in~$\cD$. In other words,
if $f_i:X_i\to Y_i$ are acyclic cofibrations between cofibrant objects 
of~$\cC_i$, then $F(f_1,f_2)$ is a weak equivalence.
\item[c)] For any cofibrant object~$P$ of~$\cC_1$ the functor 
$F(P,-):\cC_2\to\cD$ transforms acyclic cofibrations between cofibrant objects
into weak equivalences, and similarly for $F(-,Q):\cC_1\to\cD$, $Q$ any 
cofibrant object of~$\cC_2$.
\end{itemize}
If the above conditions hold, the functor $F$ admits a left derived 
$\dL F:\Ho\cC_1\times\Ho\cC_2\to\Ho\cD$, and $(\dL F)(\gamma X,\gamma Y)$ 
can be computed for any $X\in\Ob\cC_1$, $Y\in\Ob\cC_2$ by taking 
$\gamma F(P,Q)$, where $P\to X$ and $Q\to Y$ are arbitrary cofibrant
replacements.
\end{Propz}
{\bf Proof.} Immediate from~\ptref{prop:ex.left.der}, 
\ptref{l:acof.cof.suff} and the description of model category structure 
on $\cC_1\times\cC_2$ given in~\ptref{sp:prod.mod.cat}.

Now we'd like to verify the above conditions for the tensor product. 
However, it pays out to study a more general situation.

\begin{DefD}\label{def:comp.tens.mod.str}
Let $\otimes$ be a tensor structure on a model category~$\cC$. We say that 
these two structures are {\em compatible\/} if the following version 
of (SM7b) holds:

{\myindent{(TM)}
If $i:A\to B$ and $s:K\to L$ are cofibrations in~$\cC$, then
\begin{equation}
i\boxe s:A\otimes L\sqcup_{A\otimes K}B\otimes K\to B\otimes L
\end{equation}
\setmyindent

is a cofibration in~$\cC$ as well, 
acyclic if either $i$ or~$s$ is acyclic.
}

\noindent
In this case we say that $\cC$ is a {\em tensor model category.}

Similarly, if $\oslash:\cD\times\cC\to\cD$ is an external $\otimes$-action 
of a tensor model category~$\cC$ on another model category~$\cD$, we say that 
{\em $\oslash$ is compatible with the model structure of~$\cD$} if the 
following variant of (TM) holds:

{\myindent{(TMe)\ }
If $i:A\to B$ is a cofibration in~$\cD$ and $s:K\to L$ is a cofibration 
in~$\cC$, then
\begin{equation}
A\oslash L\sqcup_{A\oslash K}B\oslash K\to B\oslash L
\end{equation}
\setmyindent

is a cofibration in~$\cD$ as well, acyclic if either $i$ or~$s$ is acyclic.
}
\end{DefD}

\nxsubpoint (Examples.)
Clearly, the $\otimes$-structure on $s\catSets$ given by the direct product, 
as well as the $\otimes$-action of $s\catSets$ on any simplicial model 
category~$\cC$ are compatible with the model structures involved, 
essentially by definition. Another example will be given by $\otimes_\Lambda$
on $s(\catMod\Lambda)$ for any generalized ring~$\Lambda$ 
after we check (TM) in that case.

\nxsubpoint\label{sp:der.comp.tensprod} (Deriving $\otimes$.)
The importance of compatible tensor structures~$\otimes$ on 
a model category~$\cC$ is due to the fact they can be derived, thus yielding 
a tensor structure on $\Ho\cC$. More precisely:

\begin{Propz}
Let $\otimes$ be a compatible tensor structure on a model category~$\cC$. 
Suppose that $\emptyset\otimes X\cong\emptyset\cong X\otimes\emptyset$ 
for any object $X$ of~$\cC$, where $\emptyset$ denotes the initial object.
Then $P\otimes Q$ is a cofibrant object for any two cofibrant objects 
$P$ and $Q$ of~$\cC$, and $\otimes:\cC\times\cC\to\cC$ satisfies the 
conditions of~\ptref{sp:ex.lder.bif}; in particular, its left derived 
$\Lotimes=\dL\otimes$ exists and it can be computed by 
taking the tensor product of cofibrant replacements of both arguments. 
It defines a tensor structure on~$\Ho\cC$, 
with the same constraints as the original tensor structure on~$\cC$.
\end{Propz}

Similar statements can be made for a compatible $\otimes$-action 
$\oslash:\cD\times\cC\to\cD$ that preserves initial objects in each argument: 
it can be derived, thus yielding an external 
$\otimes$-action of $\Ho\cC$ on $\Ho\cD$.

\begin{Proof}
Applying (TM) to $\emptyset\to B$ and $\emptyset\to L$ we see immediately that 
$B\otimes L$ is cofibrant for any cofibrant $B$ and~$L$. Applying the same 
axiom to any acyclic cofibration $A\to B$ and any cofibration $\emptyset\to L$
we see that $A\otimes L\to B\otimes L$ is an acyclic cofibration as well. 
In particular, $-\otimes L$ transforms acyclic cofibrations between cofibrant 
objects into weak equivalences for any cofibrant~$L$, and we have a similar 
result for $B\otimes-$, for any cofibrant~$B$. So the condition
\ptref{sp:ex.lder.bif},c) holds and we can derive~$\otimes$. The case 
of an external $\otimes$-action is treated similarly with (TM) replaced 
by~(TMe).
\end{Proof}

\nxsubpoint\label{sp:comp.of.iHoms} (Compatibility and inner Homs.)
Suppose we are given some $\otimes$-structure on a model category~$\cC$, 
such that for any $A$ and $X$ in $\cC$ the functor $\Hom_\cC(A\otimes-,X)$ 
is representable by some inner Hom object $\iHom_\cC(A,X)$:
\begin{equation}
\Hom_\cC(K,\iHom_\cC(A,X))\cong\Hom_\cC(A\otimes K,X)
\end{equation}
Then we can state the compatibility axiom (TM) in the following equivalent way:

{\myindent{(TMh)}
Whenever $i:A\to B$ is a cofibration and $p:X\to Y$ is a fibration in~$\cC$, 
the following morphism is a fibration in~$\cC$:

}\vskip-1em
\begin{equation}
\iHom_\cC(B,X)\stackrel{(i^*,p_*)}\longto
\iHom_\cC(A,X)\times_{\iHom_\cC(A,Y)}\iHom_\cC(B,X)
\end{equation}

{\setmyindent
We require this fibration to be acyclic when either $i$ or $p$ is a weak 
equivalence.

}

The proof of equivalence of (TM) and (TMh) proceeds in the same way as 
the classical proof of (SM7)$\Leftrightarrow$(SM7b), using just the various 
lifting properties and adjointness isomorphisms involved, 
so we won't provide any further details.

Similarly, if we are given an external action $\oslash:\cD\times\cC\to\cD$, 
with both $\cC$ and $\cD$ model categories, and if the functor 
$\Hom_\cD(A\otimes-,X)$ is representable by some $\iHom_\cC(A,X)$, then 
the compatibility condition (TMe) has an equivalent form similar to (TMh) 
above, where $i:A\to B$ and $p:X\to Y$ are required now to be a cofibration 
and a fibration in~$\cD$.

\nxsubpoint (Deriving $\iHom$.)
Suppose the above conditions are fulfilled, so we have an inner Hom 
functor $\iHom_\cC:\cD^0\times\cD\to\cC$. Let's assume that 
$\iHom_\cC(\emptyset_\cD,X)\cong e_\cC$ and 
$\iHom_\cC(X,e_\cD)\cong e_\cC$ for 
any object~$X$ of~$\cD$; actually the second relation holds automatically, 
while the first is equivalent to $\emptyset_\cD\oslash K\cong\emptyset_\cD$
for all $K$ in~$\cC$. Then (TMh) implies the (variant of) condition
\ptref{sp:ex.lder.bif},c), sufficient for the existence of the 
right derived $\dR\iHom_\cC:\Ho\cD^0\times\Ho\cD\to\Ho\cC$, 
by the same reasoning as in~\ptref{sp:der.comp.tensprod}. In particular, 
$\dR\iHom_\cC(\gamma A,\gamma Y)$ can be computed as $\gamma\iHom_\cC(B,X)$, 
where $A\to B$ is any cofibrant replacement, and $X\to Y$ is any 
fibrant replacement; in this case $\iHom_\cC(B,X)$ will be automatically 
fibrant.

\nxsubpoint (Adjointness of $\Loslash=\dL\oslash$ and $\dR\iHom$.)
One can show that under the above conditions $\dR\iHom$ is an inner Hom
for $\Loslash:\Ho\cD\times\Ho\cC\to\Ho\cD$, so we get
\begin{equation}
\Hom_{\Ho\cD}(A\Loslash K,X)\cong\Hom_{\Ho\cC}(K,\dR\iHom_\cC(A,X))
\end{equation}
The proof goes as follows. We can assume $A$ and $K$ to be cofibrant and 
$X$ to be fibrant. Then $A\Loslash K=A\oslash K$ is again cofibrant and 
$\dR\iHom_\cC(A,X)$ is fibrant. Now let us fix a cofibrant $A$ and
consider the two adjoint functors $A\oslash-:\cC\to\cD$ and 
$\iHom_\cC(A,-):\cD\to\cC$. They are easily seen to be a Quillen pair, 
and for example $(\dL(A\oslash-))(K)$ clearly coincides with $A\Loslash K$, 
so the adjointness statement follows from 
Quillen's theorem~\ptref{th:Q.funct}. 

\nxsubpoint (Example.)
In particular, for any simplicial model category~$\cC$ we can derive 
$\otimes:\cC\times s\catSets\to\cC$ and $\iHom:\cC^0\times\cC\to s\catSets$.

\nxsubpoint (Inner Homs for $\otimes_\Lambda$.)
Let $\Lambda$ be a generalized ring, i.e.\ a commutative algebraic monad, 
and denote by $\otimes_\Lambda$ both the tensor product 
$\catMod\Lambda\times\catMod\Lambda\to\catMod\Lambda$ and its simplicial 
extension $s\otimes_\Lambda$. We'd like to show that 
{\em inner Homs $\iHom_{s\Lambda}$ exist for the simplicial extension 
of~$\otimes_\Lambda$}, and that {\em the underlying simplicial set of
$\iHom_{s\Lambda}(A,X)$ coincides with $\iHom(A,X)$ 
of~\ptref{sp:simpl.hom.sets}.}

Once these two statements are established, the rest is trivial: 
the compatibility condition (TMh) will follow from (SM7) and the fact that 
a morphism of simplicial $\Lambda$-modules is a fibration iff this is 
true for the underlying map of simplicial sets, and the additional conditions 
for $\otimes_\Lambda$ and $\iHom_\Lambda$ hold trivially, so we'll get 
the existence of $\Lotimes$ and $\dR\iHom_\Lambda$, their adjointness 
and so on.

\nxsubpoint ($\Lambda$-structure on sets of homomorphisms.)
First of all, for any two simplicial $\Lambda$-modules $X$ and $Y$ 
we get a $\Lambda$-structure on each set $\Hom_\Lambda(X_n,Y_n)$, 
$\Lambda$ being commutative (cf.~\ptref{sp:modstr.homsigma}),
thus obtaining an induced $\Lambda$-structure on
$\Hom_{s\Lambda}(X,Y)\subset\prod_{n\geq0}\Hom_\Lambda(X_n,Y_n)$.
If $Z$ is another simplicial $\Lambda$-module, the composition map 
$\Hom_{s\Lambda}(Y,Z)\times\Hom_{s\Lambda}(X,Y)\to\Hom_{s\Lambda}(X,Z)$ 
is clearly $\Lambda$-bilinear, this being true for 
individual composition maps 
$\Hom_\Lambda(Y_n,Z_n)\times\Hom_\Lambda(X_n,Y_n)\to\Hom_\Lambda(X_n,Z_n)$. 
In this sense $s\catMod\Lambda$ becomes a ``$\Lambda$-category''.

\nxsubpoint ($\Lambda$-structure on simplicial sets of homomorphisms.)
Let $X$ and $Y$ be two simplicial $\Lambda$-modules, and consider the 
following functor $\catDelta^0\to\catMod\Lambda$:
\begin{equation}
[n]\mapsto\Hom_{s\Lambda}(X\otimes\Delta(n),Y)\quad.
\end{equation}
This functor is a simplicial $\Lambda$-module that will be denoted by 
$\iHom_\Lambda(X,Y)$ or $\iHom_{s\Lambda}(X,Y)$. Notice that its underlying 
simplicial set coincides with the ``simplicial set of homomorphisms''
$\iHom(X,Y)$, characterized by the property $\Hom_{s\catSets}(K,
\iHom(X,Y))\cong\Hom_{s\Lambda}(X\otimes K, Y)$ 
(cf.~\ptref{sp:simpl.hom.sets}).

\nxsubpoint (Adjointness of $\iHom_{s\Lambda}$ and $s\otimes_\Lambda$.)
For any two simplicial $\Lambda$-modules $X$ and~$Y$ we have a canonical 
evaluation map of simplicial sets 
$\ev_{X,Y}:X\times\iHom_{s\Lambda}(X,Y)\to Y$, easily seen 
to be $\Lambda$-bilinear. Thus we get $\ev'_{X,Y}:X\otimes_\Lambda
\iHom_{s\Lambda}(X,Y)\to Y$. Now we use $\ev'_{X,Y}$ to construct 
a canonical map for any simplicial $\Lambda$-module~$Z$:
\begin{equation}
\Hom_{s\Lambda}(Z,\iHom_{s\Lambda}(X,Y))\to
\Hom_{s\Lambda}(X\otimes_\Lambda Z,Y)
\end{equation}
We claim that the above arrow is a functorial isomorphism, hence that 
$\iHom_{s\Lambda}$ is indeed an inner Hom for $s\otimes_\Lambda$.
First of all, this is true for $Z=L_\Lambda\Delta(n)$, since 
$\Hom_{s\Lambda}(L_\Lambda\Delta(n),X)\cong X_n$ and 
$X\otimes_\Lambda L_\Lambda K\cong X\otimes K$ for any simplicial set~$K$, 
in particular for $K=\Delta(n)$. Next, the simplicial $\Lambda$-modules
$L_\Lambda\Delta(n)$ constitute a set of generators for $s\catMod\Lambda$, 
so any $Z$ can be written as an inductive limit of objects of this form; 
our claim follows now from the fact that $\otimes_\Lambda$ commutes 
with arbitrary inductive limits.

\nxsubpoint (Existence of $\Lotimes=\dL\otimes_\Lambda$ 
and $\dR\iHom_{s\Lambda}$.)
Since the underlying simplicial set of $\iHom_{s\Lambda}(X,Y)$ coincides with 
the simplicial set of homomorphisms $\iHom(X,Y)$, and a morphism of 
simplicial $\Lambda$-modules is a fibration iff this is true for the 
underlying map of simplicial sets, we see that (SM7a) implies immediately 
(TMh), so {\em $s\otimes_\Lambda$ is a compatible ACU
$\otimes$-structure on $s\catMod\Lambda$ with inner Homs $\iHom_{s\Lambda}$.}
All the additional assumptions like $\emptyset\otimes X\cong\emptyset\cong
X\otimes\emptyset$ hold as well, so we can invoke our previous results and 
obtain the existence and adjointness of derived tensor products and inner
Homs:

\begin{Propz}
Let $\Lambda$ be a generalized ring. Then the simplicial extension of the 
tensor product $\otimes_\Lambda$ to the category of simplicial 
$\Lambda$-modules is compatible with the model structure on this category,
and $\iHom_{s\Lambda}$ is an inner Hom for this tensor structure.
Both these functors can be derived. The derived tensor product 
$\Lotimes=\dL\otimes_\Lambda$ induces an ACU $\otimes$-structure 
on the corresponding homotopy category $\cD^{\leq0}(\Lambda)$, and 
$\dR\iHom_\Lambda$ is an inner Hom for the derived tensor product.

Moreover, $\Lotimes$ can be computed by applying $\otimes_\Lambda$ to the 
cofibrant replacements of its arguments, yielding a cofibrant object as well, 
and $\dR\iHom_\Lambda$ can be computed by taking the cofibrant replacement 
of the first argument and fibrant replacement of the second argument 
and applying $\iHom_{s\Lambda}$.
\end{Propz}

\nxsubpoint (Further properties of derived tensor products.)
Now we can easily deduce some other natural properties of derived tensor 
products. For example, if $\Lambda$ is a generalized ring with zero, 
$\Lotimes$ commutes with suspension in each variable: 
$\Sigma X\Lotimes Y\cong\Sigma(X\Lotimes Y)\cong X\Lotimes\Sigma Y$. 
If $\rho:\Lambda\to\Xi$ is a homomorphism of generalized rings, we have
$\dL\rho^*(X\Lotimes\strut_\Lambda Y)\cong\dL\rho^* X\Lotimes\strut_\Xi
\dL\rho^* Y$. We get some formulas involving $\dR\iHom$ as well:
$\dR\iHom_\Lambda(X,\dR\rho_*Y)\cong\dR\iHom_\Xi(\dL\rho^*X,Y)$.

\nxsubpoint (Stable version.)
Since both $\Lotimes$ and $\dL\rho^*$ commute with suspension, 
we can construct their stable versions.
Thus $\Lotimes$ defines an ACU $\otimes$-structure on $\cD^-(\Lambda)$ for any 
generalized ring~$\Lambda$ with zero, and we get a $\otimes$-functor 
$\dL\rho^*:\cD^-(\Lambda)\to\cD^-(\Xi)$ for any homomorphism 
$\rho:\Lambda\to\Xi$.

\nxsubpoint (Derived tensor powers.)
Fix any integer $k\geq0$ and consider the $k$-fold tensor power functor 
$T^k:X\mapsto X^{\otimes k}$ on $\catMod\Lambda$ and $s\catMod\Lambda$. 
We see that this functor preserves cofibrant objects and weak equivalences 
between them, so it has a left derived $\dL T^k$, that can be computed 
with the aid of cofibrant replacements. Clearly, $\dL T^k$ coincides with 
the $k$-fold derived tensor product on $\cD^{\leq0}(\Lambda)$. However, 
$T^k$ is not additive even in the classical case, so we've got here a 
natural example of a derived non-additive functor (cf.~\cite{DP}).
Notice that $T^k$ is not ``linear'' in several other respects, for example 
$T^k\circ\Sigma\cong\Sigma^k\circ T^k$ and $T^k(X\otimes K)\cong
T^k(X)\otimes (K^{\times n})$, i.e.\ $T^k$ is ``homogeneous of degree~$k$''.

\nxsubpoint (Classical case.)
Notice that when the generalized ring $\Lambda$ is additive, i.e.\ 
just a classical commutative ring, then the derived tensor product and 
derived pullbacks we constructed here coincide with their classical 
counterparts via Dold--Kan and Eilenberg--Zilber
(cf.~\ptref{sp:appl.class.der.tens}).

\cleardoublepage
\mysection{Homotopic algebra over topoi}\label{sect:homot.alg/topoi}

Now we'd like to extend our constructions to the case of a category of 
modules over a generalized ringed site or topos. The main idea here is to 
develop ``local'' variants of all homotopic notions discussed so far, 
and try to convert existing proofs in a more or less systematic fashion; 
this would enable us to re-use most statements and proofs without much work.

The dictionary to translate mathematical statements from the ``global'' 
or ``point'' context (based on constructions in $\catSets$) to the 
``local'' or ``topos'' context (based on similar constructions in an 
arbitrary topos~$\cE$) is given by the so-called {\em Kripke--Joyal 
semantics}, that essentially asserts that topoi are models of intuitionistic 
set theory.

We are going to present the main ideas of this semantics in 
an informal fashion after we give all necessary preliminary definitions; 
let us content ourselves for now with the following remark. While 
transferring from the ``point'' to the ``topos'' situation, we have to replace 
sets by sheaves of sets (i.e.\ objects of our topos $\cE$), elements of 
sets by sections of sheaves, sets of maps by sheaves of maps (i.e.\ 
local Homs $\iHom_\cE$), maps of sets by morphisms of sheaves, 
small categories by inner categories, and arbitrary categories by stacks
and/or prestacks. 
Moreover, all constructions and situations we consider have to be 
compatible with arbitrary pullbacks $i_X^*:\cE\to\cE_{/X}$, for any object
$X$ of~$\cE$, and the properties of maps, sets etc.\ we consider have to 
be {\em local}. Finally, the last but not least important ingredient of 
Kripke--Joyal semantics is that we have to replace the classical logic 
in all proofs by the {\em intuitionistic logic.} This means that we cannot
use the law of excluded middle ($A\vee\neg A$), the proofs based on 
{\em reductio ad absurdum\/} ($A\wedge\neg B\Rightarrow0$ doesn't imply 
$A\Rightarrow B$, and $\neg\forall x.\neg A[x]$ doesn't imply 
$\exists x.A[x]$), as well as the axiom of choice. Our hope here is that 
most algebra is essentially intuitionistic (that's why the modules 
over a ringed topos have almost all properties of modules over a ring, 
for example), even if most calculus is not (as established by Brouwer and 
Weyl in their attempts to build all mathematics intuitionistically). 

Since we always try to keep things as algebraic as possible, we might 
hope to be able to transfer the definitions and constructions of 
model categories to the topos case, thus obtaining the notion of a
{\em model stack.}

\nxpointtoc{Generalities on stacks}\label{p:gen.on.stacks}
Now we are going to recall the basic definitions of fibered categories, 
cartesian functors, prestacks and stacks 
(cf.\ \cite[VI]{SGA1} and~\cite{Giraud}).

\nxsubpoint ($\cS$-categories; cf.\ SGA~1.)
When we are given a functor $p:\cC\to\cS$, we can say that {\em $\cC$ is an
$\cS$-category,} or {\em a category over~$\cS$.} 
Then for any $S\in\Ob\cS$ we denote by $\cC_S$ or $\cC(S)$ the 
fiber of $p$ over $S$, considered here as a point subcategory of~$\cS$. 
For any $\phi:T\to S$, $X\in\Ob\cC_T$, $Y\in\Ob\cC_S$ we denote by 
$\Hom_\phi(X,Y)$ the set of all {\em $\phi$-morphisms $f:X\to Y$}, i.e.\ 
all morphisms $f:X\to Y$, such that $p(f)=\phi$. If $\phi=\id_S$, we write 
$\Hom_S$ instead of~$\Hom_{\id_S}$; clearly, $\Hom_S$ is just the 
set of morphisms inside $\cC_S$.

An {\em $\cS$-functor\/} $F:\cC\to\cC'$ between two $\cS$-categories is 
simply any functor, such that $p'\circ F=p$; an 
{\em $\cS$-natural transformation $\eta:F\to G$} between two $\cS$-functors 
$F$, $G:\cC\to\cC'$ is a natural transformation $\eta:F\to G$, such that 
$p'\star\eta=\Id_p$. In this way we obtain a strictly associative 
2-category $\catCat/\cS$, so we can define ``$\cS$-adjoint functors'', 
``$\cS$-equivalences'' and so on with the usual properties (cf.~\cite{SGA1}). 

Informally, one should think of an $\cS$-category~$\cC$ as ``collection of 
categories $\cC_S$ parametrized by $\cS$''; $\cS$-functors $F:\cC\to\cC'$, 
$\cS$-equivalences etc.\ should be thought of as ``compatible'' collections of 
fiberwise functors $F_S:\cC_S\to\cC'_S$ etc.

Given any functor $H:\cS'\to\cS$, we can consider the pullback $p':\cC'\to\cS'$
of any $\cS$-category $p:\cC\to\cS$; clearly, $\cC'_S=\cC_{H(S)}$.

\nxsubpoint (Fibered categories and cartesian functors; cf.\ SGA~1.)
A $\phi$-morphism $h:Y\to X$ in an $\cS$-category $\cC$, where
$X\in\Ob\cC_S$, $Y\in\Ob\cC_T$, $\phi:T\to S$, is said to be 
{\em cartesian\/} if it induces a bijection $h_*:\Hom_T(Z,Y)\to
\Hom_f(Z,X)$, for all $Z\in\Ob\cC_T$. We say that {\em $\cC$ is a fibered 
category over~$\cS$} if for any $\phi:T\to S$ and any $X\in\Ob\cC_S$ one 
can find a cartesian $\phi$-morphism $Y\to X$ with target~$X$, and if 
the composite of any two cartesian morphisms is again cartesian.

When $\cC$ is fibered over~$\cS$, we choose for each $\phi:T\to S$ and 
$X\in\Ob\cC_S$ some cartesian $\phi$-morphism $\phi^*X\to X$, thus obtaining a 
{\em pullback functor\/} $\phi^*:\cC_S\to\cC_T$, defined uniquely up to 
a unique isomorphism. We obtain natural isomorphisms 
$(\phi\psi)^*\cong\psi^*\circ\phi^*:\cC_S\to\cC_U$ for any
$U\stackrel\psi\to T\stackrel\phi\to S$, satisfying the pentagon axiom. 
In this way {\em a fibered category over $\cS$ is essentially 
a contravariant 2-functor $\cS^0\to\catCat$.} The difference is that 
we don't fix any particular choice of pullback functors $\phi^*$.

An $\cS$-functor $F:\cC\to\cC'$ between two fibered categories over~$\cS$ 
is called {\em cartesian\/} if it transforms cartesian morphisms of $\cC$ into 
cartesian morphisms of~$\cC'$. This essentially means that we have a collection
of functors $F_S:\cC_S\to\cC'_S$, commuting with all pullback functors 
$\phi^*$ in $\cC$ and $\cC'$ up to an isomorphism. 

We denote by $\catFunct_{cart,\cS}(\cC,\cC')$ or $\catCart_\cS(\cC,\cC')$  
the category of all cartesian functors $\cC\to\cC'$. Since $\cS$ is fibered 
over itself, we can consider ``the category of cartesian sections'' 
$\Gamma_{cart,\cS}(\cC):=\catCart_\cS(\cS,\cC)$; we denote it by 
$\projLim_\cS\cC$ as well. Clearly, when $\cS$ has a final object~$e$, 
we have an equivalence $\cC_e\simto\Gamma_{cart}(\cC/\cS)$.

Notice that for any object $F$ of~$\cS$, or, more generally, for any presheaf 
$F$ over~$\cS$, the natural functor $\cS_{/F}\to\cS$ defines a fibered 
category over~$\cS$, so we can 
consider $\catCart_\cS(\cS_{/F},\cC)$; when $F$ lies in~$\cS$, this
category is equivalent to $\cC_F$. If we replace each $\cC_S$ with 
equivalent category $\catCart_\cS(\cS_{/F},\cS)$, we obtain an 
$\cS$-equivalent $\cS$-fibered category with a natural choice of pullback 
functors, such that $(\psi\phi)^*=\phi^*\psi^*$. Actually this construction, 
when considered for all presheaves $F\in\Ob\hat\cS$,
yields a canonical extension $\cC^+/\hat\cS$ of $\cC/\cS$. Notice that 
$S\mapsto\Ob\cC^+_S$ and $S\mapsto\Ar\cC^+_S$ are presheaves over~$\cS$, 
i.e.\ we obtain a presheaf of categories over~$\cS$, or an inner category
in $\hat\cS_\univV$, for a suitable universe $\univV\ni\univU$.

Finally, if $H:\cS'\to\cS$ is any functor, then the pullback with respect 
to~$H$ transforms cartesian functors between $\cS$-fibered categories 
into cartesian functors between $\cS'$-fibered categories. 
For example, $\cC^+\times_{\hat\cS}\cS$ is $\cS$-equivalent to~$\cC$.
When $\cS'=\cS_{/F}$, we put $\cC_{/F}:=\cC\times_\cS\cS_{/F}$.

\begin{NotatD} 
Whenever a pullback functor $\phi^*:\cC_S\to\cC_T$ in an 
$\cS$-fibered category $\cC$ has a left adjoint, it will be denoted by 
$\phi_!:\cC_T\to\cC_S$. Similarly, the right adjoint to $\phi^*$ will be 
denoted by~$\phi_*$.
\end{NotatD}

\nxsubpoint (Hom-presheaves.)
Given a fibered $\cS$-category $\cC$ and two objects $X$, $Y\in\Ob\cC_S$, 
$S\in\Ob\cC$, we define the {\em local Hom\/} 
$\iHom_{\cC_{/S}}(X,Y):(\cS_{/S})^0\to\catSets$ by 
$(T\stackrel\phi\to S)\mapsto\Hom_{\cC_T}(\phi^*X,\phi^*Y)$. 
Clearly, $\iHom_{\cC_{/S}}(X,Y)$ 
is a presheaf over 
$\cS_{/S}$, with global sections equal to $\Hom_{\cC_S}(X,Y)$, and 
$\phi^*\negthinspace\iHom_{\cC_{/S}}(X,Y)\cong
\iHom_{\cC_{/T}}(\phi^*X,\phi^*Y)$, for any $\phi:T\to S$ in~$\cS$.

\nxsubpoint\label{def:prestack} (Prestacks.)
We say that a fibered category $\cC$ over a {\em site\/} $\cS$ is 
a {\em prestack\/} if for any $S\in\Ob\cS$, $X$, $Y\in\Ob\cC_S$ the 
Hom-presheaf $\iHom_{\cC_{/S}}(X,Y)$ is actually a sheaf over~$\cS_{/S}$. 
{\em Morphisms\/} of prestacks are just cartesian $\cS$-functors.

Notice that the above definition applies when $\cS$ is actually a topos~$\cE$,
considered as a site with respect to its canonical topology. In this case 
$\widetilde{\cS_{/S}}$ is equivalent to $\cE_{/S}$, so we get a local 
Hom-object $\iHom_{\cC_{/S}}(X,Y)\in\Ob\cE_{/S}$.

\nxsubpoint (Descent data.)
Given any object~$S\in\Ob\cS$ and any sieve $R\subset S$, i.e.\ a subpresheaf 
of~$S$ in~$\hat\cS$, we obtain a natural functor for any fibered 
$\cS$-category~$\cC$:
\begin{equation}\label{eq:desc.comp.funct}
\cC(S)=\cC_S\cong\catCart(\cS_{/S},\cC)\to\cC^+(R)=\catCart(\cS_{/R},\cC)
\end{equation}
We say that $\cC^+(R)$ is {\em the category of descent data with 
respect to~$R$}, and denote it by $\catDesc(R;\cC)$. When $R$ is generated 
by some cover $(S_\alpha\stackrel{f_\alpha}\to S)$, and double and triple 
fibered products of $S_\alpha$ over $S$ do exist, we denote by 
$\catDesc((S_\alpha\to S);\cC)$ the category of descent data with respect 
to this cover; by definition, an object of this category is a descent datum
with respect to $(S_\alpha\to S)$ with values in~$\cC$, i.e.\ a collection 
of objects $X_\alpha\in\Ob\cC(S_\alpha)$ and isomorphisms 
$\theta_{\alpha\beta}:\pr_1^*X_\alpha\simto\pr_2^*X_\beta$ in 
$\cC(S_\alpha\times_SS_\beta)$, subject to the usual cocycle relation in 
$\cC(S_\alpha\times_SS_\beta\times_SS_\gamma)$. One checks immediately that 
the naturally arising functor $\catDesc(R;\cC)\to\catDesc((S_\alpha\to S);\cC)$
is indeed an equivalence of categories. This justifies the above terminology.

\begin{DefD}\label{def:stack} {\rm (cf.~\cite[II,1.2]{Giraud})}
A fibered category $\cC$ over a site $\cS$ is called a {\em stack}
(resp. {\em prestack}) {\em over~$\cS$} if for any object $S\in\Ob\cS$ and 
any covering sieve $R$ of~$S$ the functor $\cC(S)\to\cC^+(R)$ 
of~\eqref{eq:desc.comp.funct} is an equivalence (resp.\ fully faithful),
i.e.\ if all covers in~$\cS$ have the effective descent (resp.\ descent)
property with respect to~$\cC$.
\end{DefD}

One checks immediately that $\cC(S)\to\cC^+(R)$ is fully faithful iff
all local Hom-presheaves $\iHom_{\cC_{/S}}(X,Y)$ satisfy the sheaf condition 
for $R\subset S$; therefore, {\em a fibered category $\cC/\cS$ is a prestack 
iff all $\iHom_{\cC_{/S}}(X,Y)$ are sheaves}, so \ptref{def:stack} 
and~\ptref{def:prestack} are compatible in this respect.

\nxsubpoint (Examples.)
(a) For any site $\cS$ we have a ``stack of sets'' (or sheaves)
$\stSETS_{\cS}$, given by $\stSETS_\cS(X):=\widetilde{\cS_{/X}}$, with 
the pullback functors given by restriction of sheaves. We might have denoted 
this stack by $\stSHEAVES_{\cS}$ as well, but we wanted to illustrate that 
this stack will indeed play the role of $\catSets$ in our further 
local considerations.

(b) For any category $\cE$ with fibered products the category of arrows 
$\Ar\cE$ together with the target functor $t:\Ar\cE\to\cE$ defines a 
fibered category structure, such that $(\Ar\cE)(X)=\cE_{/X}$. When~$\cE$ 
is a topos, $\widetilde{\cE_{/X}}$ is naturally equivalent to $\cE_{/X}$, 
hence $\Ar\cE$ is $\cE$-equivalent to $\stSETS_{\cE}$; in particular, 
it is an $\cE$-stack. Actually, we'll usually put $\stSETS_{\cE}:=\Ar\cE$ 
when working over a topos.

(c) When $\Sigma$ is an algebraic monad over a topos~$\cE$, we can define 
the {\em stack of $\Sigma$-modules\/} $\stMOD\Sigma=\stMODo_{\cE,\Sigma}$
by $\stMOD\Sigma(X):=\catMod{\Sigma_{|X}}$. Of course, we can do the same 
construction over a site~$\cS$, if we define $\stMOD\Sigma(X)$ to be 
the category of sheaves of $\Sigma_{|X}$-modules over~$\cS_{/X}$. When 
$\Sigma$ is the constant algebraic monad $\Fempty$, we recover $\stSETS_\cE$.

\nxsubpoint\label{sp:stk.over.site.topos} 
(Stacks over a site and over the corresponding topos.)
If we have a stack $\cC$ over a site $\cS$, we can extend to a stack 
$\tilde\cC$ over topos $\tilde\cS$, simply by restricting $\cC^+$ to 
$\tilde\cS$:
$\tilde\cC:=\cC^+\times_{\hat\cS}\tilde\cS$, i.e.\ $\tilde\cC(F)=\cC^+(F)=
\catCart(\cS_{/F},\cC)$ for any sheaf~$F$. Conversely, any stack 
$\cC'$ over $\tilde\cS$ defines a stack $\cC'\times_{\tilde\cS}\cS$ 
over~$\cS$ by restriction. Now we claim that {\em these two functors 
establish a 2-equivalence between the 2-categories of stacks over~$\cS$ 
and over~$\tilde\cS$} (cf.~\cite[II,3.3]{Giraud}). We'll use this observation 
to work with stacks over topoi whenever possible.

The main idea here is that any sheaf~$F$ can be covered in~$\tilde\cS$ by
a small family of objects $S_\alpha$, and each $S_\alpha\times_FS_\beta$ 
has a similar cover by some $S_{\alpha\beta}^{(\nu)}$ as well, 
hence for any stack $\cC'/\tilde\cS$ the category $\cC'(F)$ can be described 
in terms of some sort of descent data, consisting of collections 
of objects $X_\alpha\in\cC'(S_\alpha)$ and isomorphisms 
$\theta_{\alpha\beta}^{(\nu)}$ in $\cC'(S_{\alpha\beta}^{(\nu)})$, hence 
$\cC'$ is completely determined by its restriction to~$\cS$.

\nxsubpoint (Fibers over direct sums.)
If $\cC$ is a stack over a topos~$\cE$, then $\cC(\bigsqcup_\alpha S_\alpha)$ 
is naturally equivalent to the product $\prod_\alpha\cC(S_\alpha)$,
just because this product is the category of descent data for 
the cover $(S_\alpha\to S)_\alpha$, $S:=\bigsqcup S_\alpha$, direct sums 
in~$\cE$ being disjoint.

\nxsubpoint (Stacks over the point site and topos.)
The above considerations apply in particular to the point site $\st1$
(the final category) and the point topos $\catSets$. We see that 
{\em stacks over~$\catSets$ are essentially just categories.} More precisely, 
a stack $\cC$ over $\catSets$ defines a category $\cC(\st1)$, and conversely,
any category $\cD$ defines a stack~$\cD^+$ over~$\catSets$ by 
$\cD^+(I):=\cD^I$ for any set~$I$.

\nxsubpoint\label{sp:prestack.from.fib}
(Prestack associated to a fibered category.)
Given a fibered category $\cC$ over a site~$\cS$, the {\em stack\/} 
(resp.\ {\em prestack}) 
{\em associated to~$\cC$} is defined by the following requirement: 
it is a stack (resp.\ prestack) $\bar\cC/\cS$ with a cartesian $\cS$-functor 
$I:\cC\to\bar\cC$, such that for any stack (resp.\ prestack) $\cD/\cS$ the 
functor $I^*:\catCart_\cS(\bar\cC,\cD)\to\catCart_\cS(\cC,\cD)$ is an 
equivalence.

The prestack $\bar\cC$ associated to $\cC$ can be constructed as follows. 
Consider the sheafifications $a\iHom_{\cC_{/S}}(X,Y)$ of local Hom-presheaves,
where $S\in\Ob\cS$, $X,Y\in\cC_S$.
Since $a:\hat\cS\to\tilde\cS$ commutes with finite products, we get 
canonical maps $a\iHom_{\cC_{/S}}(Y,Z)\times a\iHom_{\cC_{/S}}(X,Y)\to
a\iHom_{\cC_{/S}}(X,Z)$ induced by composition. Now we put 
$\Ob\bar\cC_S:=\Ob\cC_S$, 
$\Hom_{\bar\cC_S}(X,Y):=\Gamma(S,a\iHom_{\cC_{/S}}(X,Y))$, and define 
for any $\phi:T\to S$ the pullback functors $\phi^*:\bar\cC_S\to\bar\cC_T$ 
on morphisms $f\in\Hom_{\bar\cC_S}(X,Y)=\sF(S)$, $\sF:=a\iHom_{\cC_{/S}}(X,Y)$,
by means of $\sF(\phi):\sF(S)\to\sF(T)$. One checks immediately that 
this defines indeed a fibered category $\bar\cC/\cS$, and by 
construction all $\iHom_{\bar\cC_{/S}}(X,Y)=
a\negthinspace\iHom_{\cC_{/S}}(X,Y)$ are indeed sheaves,
i.e.\ $\bar\cC$ is a prestack. The universal property of the ``identity'' 
functor $I:\cC\to\bar\cC$ is also immediate, once we observe that 
any cartesian functor $F:\cC\to\cD$ induces maps of presheaves
$F_{X,Y}:\iHom_{\cC_{/S}}(X,Y)\to\iHom_{\cD_{/S}}(F(X),F(Y))$.

\nxsubpoint\label{sp:stack.from.fib} (Stack associated to a fibered category.)
Now we'd like to present a construction of the stack $\hat\cC$ associated
to a fibered category~$\cC$ over a site~$\cS$, different from 
that given in~\cite{Giraud}. First of all, 
we construct the associated prestack $\bar\cC$ by sheafifying all 
Hom-presheaves as above, and observe that a stack associated to $\bar\cC$ 
will be also associated to~$\cC$, i.e.\ we can assume $\cC$ to be 
a prestack.

Next, we define $\hat\cC(S)$ to be the category of descent data in~$\cC$
($=\bar\cC$) with respect to all covers of~$S$. More formally, 
we consider the ordered set $J(S)$ of covering sieves of~$S$
as a category, and put
\begin{equation}
\hat\cC(S):=\injLim_{R\in J(S)}\cC^+(R)\quad.
\end{equation}
Recall that for any fibered category $\cF$ over another category $\cI$ 
the pseudolimit $\injLim_\cI\cF$ can be constructed simply as the 
localization of $\cF$ with respect to all cartesian morphisms; in particular,
we get canonical functors $\cF_i\to\injLim_\cI\cF$, $i\in\Ob\cI$, 
commuting with the pullback functors up to some natural isomorphisms.

In our situation all pullback (=restriction) functors 
$\cC^+(R)\to\cC^+(R')$ will be fully faithful for any covering sieves
$R'\subset R\subset S$, $\cC$ being a prestack, 
hence the functors $\cC^+(R)\to\hat\cC(S)$ and $J_S:\cC(S)\to\hat\cC(S)$ 
will be fully faithful as well.

This construction defines a fibered category~$\hat\cC/S$ and a 
fully faithful cartesian functor $J:\cC\to\hat\cC$; we leave to the reader 
the verification of the universal property of $J$ and 
of $\hat\cC$ being a stack.

\nxsubpoint (Direct image of a stack.)
Let $\cC$ be a stack over a site $\cS$, and $u:\cS\to\cS'$ be a morphism of 
sites, given by some pullback functor $u^*:\cS'\to\cS$ (if $\cS$ and 
$\cS'$ are closed under finite projective limits, this means simply that 
$u^*$ is left exact and preserves covers). Then we can define a new fibered 
category~$\cC'$ over $\cS'$, namely, the pullback $\cC':=\cC\times_\cS\cS'$ 
of $\cC$ with respect to $u^*$. One checks immediately (at least when 
fibered products exist in both sites) that $\cC'$ is a {\em stack\/} 
over~$\cS'$. It is usually called {\em the direct image of $\cC$ with 
respect to~$u$} and is denoted by $u^{st}_*\cC$.

\nxsubpoint (Inverse image of a stack.)
Given a morphism of sites $u:\cS\to\cS'$, and two stacks $\cC\stackrel p\to\cS$
and $\cC'\stackrel{p'}\to\cS'$, a {\em $u$-functor\/} $F:\cC'\to\cC$ 
is simply any functor, such that $p\circ F=u^*\circ p'$, and the induced
functor $\bar F:\cC'\to\cC\times_\cS\cS'=u^{st}_*\cC$ is {\em cartesian}. 
Informally, $F$ is a collection of functors $F_{S'}:\cC'(S')\to\cC(u^*S')$, 
parametrized by $S'$ in~$\cS'$.

If a $u$-functor $I:\cC'\to u_{st}^*\cC'$ has a universal property among 
all $u$-functors with source $\cC'$, i.e.\ if for any stack $\cC$ over~$\cS$ 
the functor $I^*:\catCart(u_{st}^*\cC',\cC)\to\catFunct_u(\cC',\cC)$ is 
an equivalence of categories, then we say that the stack $u_{st}^*\cC'$ 
is the {\em inverse image\/} or {\em pullback\/} of $\cC'$ with respect to~$u$.

One can construct the inverse image $u_{st}^*\cC'$ as follows. First of all, 
define a fibered category $u^\bullet\cC'$ over $\cS$ with the required 
universal property among all fibered categories $\cC/\cS$, by the
``left Kan extension formula'' 
$(u^\bullet\cC')(T):=\injLim_{\cS'_{/T}}\cC'$, for all $T\in\Ob\cS$. Then 
the stack associated to $u^\bullet\cC'$ has the property required 
from~$u_{st}^*\cC'$.

{\bf Example.} The pullback of $\catSets$, considered here as a stack over 
the point topos/site, to an arbitrary topos~$\cE$, is the substack 
$\stLCSETS_\cE\subset\stSETS_\cE$ of locally constant objects of~$\cE$.

\nxpointtoc{Kripke--Joyal semantics}
Now we describe a version of Kripke--Joyal semantics, suitable for transferring
intuitionistic statements about categories into statements about stacks over 
a site or topos.

We fix a site~$\cS$; in some cases we assume $\cS$ to have finite projective 
limits. When $\cS$ is a topos, we denote it by~$\cE$ as well.

\nxsubpoint (Contexts.)
We define a {\em context\/} $\sitV$ (over a ``base object'' $S$ of $\cS$) as
a finite collection of {\em variables\/} (usually denoted by lowercase 
Latin letters), together with their {\em types\/} (usually denoted by 
uppercase letters) and {\em values\/} of appropriate type. Strictly 
speaking, the assignment of values is needed only for semantics
(``evaluation''), but not for checking the syntax of an expression, so 
we should distinguish {\em signatures\/} $\sitV_0$ (lists of 
variables together with their types, but without any values), and 
{\em contexts}, but usually we'll mix them for shortness.

Sometimes we have two contexts $\sitV$ and $\sitW$ over the same base object
$S$, such that all variables of~$\sitV$ appear in~$\sitW$ with the same 
types and values; then we write $\sitW\geq\sitV$ or $\sitV\leq\sitW$ and say 
that {\em $\sitW$ extends $\sitV$}.

Notice that all variables, bound or free, must have some specified type; 
when $x$ is of type $X$, we usually write $x:X$ or even $x\in X$.

Once we have a context~$\sitV$, we can construct different {\em terms\/} and 
{\em propositions\/}, involving the variables of $\sitV$ and some basic 
logical operations, and evaluate them. Notice that all terms we construct 
have a well-defined type, and their values, when defined at all, 
are of corresponding type. As to the propositions, all we can say is whether 
they are (universally) valid in a context~$\sitV$ or not. When
a proposition~$A$ is valid in context~$\sitV$, we write $\sitV\models A$. 
The value of a term $t$ in context $\sitV$ will be denoted by~$t(\sitV)$
or $t_\sitV$; this applies to variables as well.

\nxsubpoint (Simple interpretation.)
The main idea here is that types correspond to some sheaves over~$\cS_{/S}$, 
or to sets in classical case, while variables and terms correspond to 
sections of these sheaves, or to elements of sets in the classical case.
As to the propositions, they correspond to {\em local\/} properties
of these sheaves and sections.

\nxsubpoint (Pullback of contexts.)
The most important property of contexts is that any context 
$\sitV/S$ can be pulled back with respect to any morphism $\phi:T\to S$ 
in~$\cS$, yielding a context $\phi^*\sitV/T$. The way this is done 
is usually immediate from the context. For example, if a type $X$ corresponds 
to a sheaf $X_\sitV$ over $\cS_{/S}$, then the value $X_{\phi^*\sitV}$ is 
simply the restriction of $X_\sitV$ to $\cS_{/T}$, and if a variable $x:X$
has value $x_\sitV\in\Gamma_S(X_\sitV)=X_\sitV(S)$, then the value 
$x_{\phi^*\sitV}$ is simply the restriction $\phi^*x_\sitV=
(X_\sitV(\phi))(x_\sitV)$ of~$x_\sitV$ to~$T$.

\nxsubpoint (Localness of all properties.)
Any proposition $A$, depending on a context $\sitV$ (of some fixed 
signature~$\sitV_0$, which lists all variables of~$A$ with appropriate 
type), defines a {\em local property\/} of the objects involved.
This means the following:
\begin{itemize}
\item If $A$ holds in~$\sitV$, it holds in all pullbacks $\phi^*\sitV$.
\item If $\{S_\alpha\stackrel{\phi_\alpha}\to S\}$ is a cover in~$\cS$, and
$\sitV$ is a context over~$S$, such that $A$ holds in all 
$\phi_\alpha^*\sitV$, then $A$ holds in $\sitV$ itself.
\end{itemize}

\nxsubpoint (Constants.)
Apart from variables, we can also have some {\em constants\/}. These are 
essentially some variables with predefined values, always the same in 
all contexts we consider; these values come from ``outside'', and the 
constants are usually denoted by the same letter they've been denoted in the 
``outside'' reasoning. For example, if we have somehow obtained a section
$x\in X(S)$ of some sheaf~$X$, we can use $x$ as a constant of type~$X$.
An important thing about constants is that they are {\em never\/} considered
free variables.

Notice that we have {\em constant types\/} as well, as illustrated by the
constant type~$X$ given by sheaf~$X$ in the above example.

As to the {\em constant properties\/} $P(\sitW)$ of a context of 
signature $\sitW_0$, the main requirement for them to be allowed to 
appear in our expressions is to be {\em local}. Then $P$ can be used 
as a proposition of any signature $\geq\sitW_0$.

\nxsubpoint (Types.)
These are our constructions of {\em types}, in a context $\sitW/S$. 
We list the ``sheaf types'' first, the values of which are sheaves over 
$\cS_{/S}$:
\begin{itemize}
\item Any sheaf $X$ over $\cS$ or $\cS_{/S}$ defines a {\em constant type},
denoted also by~$X$.
\item If $X$ and $Y$ are sheaf types, then $Y^X$ or $X\to Y$ denotes the 
(sheaf) type $\iHom_{\cS_{/S}}(X,Y)$.
\item If $A$ is a proposition depending on the variables of~$\sitW$ as well 
as on a variable $x:X$, then $\{x|A\}$ denotes the largest subsheaf 
$X'\subset X$, such that $A$ holds for all pullbacks $\phi^*\sitW$,
$\phi:T\to S$, whenever we choose the value of $x$ belonging to 
$X'(T)\subset X(T)$.
\item If $X$ and $Y$ are sheaf types, then $X\times Y$ is another.
\item One can extend this list by means of any sheaf construction compatible 
with pullbacks, using the same notation.
\end{itemize}
We have some additional constructions involving ``large'' types:
\begin{itemize}
\item $\Ob\cC$ is a type for any stack~$\cC$ over~$\cS$. Values of 
terms of type $\Ob\cC$ are objects of $\cC(S)$.
\item If $x$ and $y$ are terms of type~$\Ob\cC$, then 
$x\to y$ or $\iHom_\cC(x,y)$ is the sheaf type given by 
$\iHom_{\cC_{/S}}(x_\sitV,y_\sitV)$.
\end{itemize}

\nxsubpoint (Terms.)
We write here $x:X$ to denote that $x$ is a term (i.e.\ expression) 
of type~$X$.
\begin{itemize}
\item Any variable $x:X$ is a valid term of type~$X$.
\item Any section $x\in X(S)$ defines a constant $x$ of (sheaf) type~$X$.
\item Any object $x\in\Ob\cC(S)$ defines a constant $x\typ \Ob\cC$.
\item Tuples: $(x_1,x_2,\ldots,x_n)$ is a term of type 
$X_1\times X_2\times\cdots\times X_n$, whenever $x_i\typ X_i$.
\item Application of functions:
$f\,x$ or $f(x)$ is a term of type $Y$ for any $f:X\to Y$ and $x\typ X$.
\item $\lambda$-abstraction:
$\lambda x.t$ or $x\mapsto t$ is a term of type~$X\to Y$ whenever 
$t$ is a valid term of type~$Y$ in any context obtained from $\sitV$ 
by pulling back and assigning to $x$ any value of type~$X$ 
\item Hilbert's iota: $\iota_xA$ or $\iota x.A$ is a term of type~$X$
whenever $A$ is a valid proposition in any context obtained from $\sitV$
by pulling back and assigning to $x$ any value of type~$X$, provided 
$\sitV\models\exists!x:X.A$.
\item $\id_x$ is a term of type $x\to x$ for any $x:\Ob\cC$.
\item $g\circ f$ is a term of type $x\to z$ for any $f:x\to y$, $g:y\to z$,
$x,y,z:\Ob\cC$.
\item We use the standard conventions for binary operations. For example, 
if $+:X\times X\to X$, $x,y\typ X$, we write $x+y$ instead of $+(x,y)$.
\end{itemize}

\nxsubpoint (Propositions.)
We explain here the {\em syntax\/} of propositions in a context $\sitV$ 
(or a signature $\sitV_0$); the semantics will be explained later.
\begin{itemize}
\item Conjuction $A\& B=A\wedge B$, disjunction $A\vee B$, 
implication $A\Rightarrow B$, equivalence $A\Leftrightarrow B$ and 
negation $\neg A$ are (valid) propositions whenever $A$ and $B$ are.
\item Logical constants $\st1$ (true) and $\st0$ (false) are propositions.
\item $x=y$ is a proposition whenever $x$ and $y$ are variables of the same 
sheaf type~$X$.
\item $(\forall x)A$, $(\forall x\typ X)A$, $\forall x.A$ and 
$\forall x\typ X.A$ 
are (equivalent) propositions for any proposition~$A$, valid in the signature
obtained from $\sitV_0$ by adding a variable~$x$ of type~$X$ 
(if a variable named $x$ is already present in~$\sitV_0$, it has to be removed
first).
\item $(\exists x)A$, $\exists x\typ X.A$ and so on are valid propositions for 
any~$A$ as above.
\item Same applies to $(\exists!x\typ X)A$, $\exists!x.A$ etc.
\end{itemize}

\nxsubpoint (Semantics of propositions.)
Now we describe when a proposition~$A$ holds in a context~$\sitV/S$
(notation: $\sitV\models A$ or $A(\sitV)$), 
which assignes appropriate values to all variables from~$A$.

Before listing the rules recall that all our propositions must be local, i.e.:
\begin{itemize}
\item If $A$ holds in $\sitV$, it holds in any pullback $\phi^*\sitV$.
\item If $\{\phi_\alpha:S_\alpha\to S\}$ is a cover in~$\cS$, and 
$\phi_\alpha^*\sitV\models A$ for all~$\alpha$, then $\sitV\models A$.
\end{itemize}

Below we denote by $\sitV\times_ST$ the pullback $\phi^*\sitV$ 
of $\sitV$ with respect to a morphism $\phi:T\to S$. If we have a cover 
$\{S_\alpha\to S\}$, we put $\sitV_\alpha:=\sitV\times_SS_\alpha$.
Now the rules:
\begin{itemize}
\item Conjunction: $A\& B=A\wedge B$ holds in $\sitV$ iff 
both $A$ and $B$ hold in~$\sitV$.
\item Disjunction: $A\vee B$ holds in $\sitV/S$ iff there is a cover
$\{S_\alpha\to S\}$, such that for each $\alpha$ at least one of $A$ and 
$B$ holds in~$\sitV_\alpha/S_\alpha$.
\item Negation: $\neg A$ holds in $\sitV/S$ iff $A$ does not hold in any 
pullback $\phi^*\sitV$ of~$\sitV$.
\item Implication: $A\Rightarrow B$ holds in~$\sitV/S$ iff in any pullback
$\phi^*\sitV$, $\phi^*\sitV\models B$ whenever $\phi^*\sitV\models A$.
\item Equivalence: $A\Leftrightarrow B$ holds in~$\sitV/S$ iff in any 
pullback of $\sitV$ each of~$A$ and~$B$ holds whenever the other holds. 
In other words, $\sitV\models A\Leftrightarrow B$ is equivalent to
$\sitV\models (A\Rightarrow B)\&(B\Rightarrow A)$. 
\item Constants: $\st1$ is always true, and $\st0$ is always false.
\item $\sitV\models x=y$, where $x$ and $y$ are terms of sheaf type~$X$, 
iff $x_\sitV=y_\sitV$ in $X(S)$.
\item Universality: $\sitV\models(\forall x\typ X)A$ means that $A$ holds 
in any context $\sitW/T$, obtained by pulling $\sitV$ back with respect to
any morphism $\phi:T\to S$, and assigning to a new variable~$x\typ X$ an
arbitrary value~$x_\sitW\in X(T)$.
\item Existence: $\sitV\models(\exists x\typ X)A$ means that there is a 
cover $\{S_\alpha\to S\}$ and some elements $x_\alpha\in X(S_\alpha)$, 
such that $A$ holds in contexts $\sitW_\alpha$, obtained by pulling 
$\sitV$ back to $S_\alpha$ and assigning to $x\typ X$ the value~$x_\alpha$.
\item Uniqueness: $(\exists!x\typ X)A[x]$ is equivalent to 
$(\exists x:X.A[x])\&(\forall x\typ X.\forall y\typ X.A[x]\& A[y]\Rightarrow 
x=y)$.
The brackets $[]$ are used here to point out some free variables entering 
in~$A$, as well as to denote the result of substituting these variables.
\item Closure: If $A$ involves some free variables $x_1\typ X_1$, \dots, 
$x_n\typ X_n$, then $\sitV\models A$ actually means 
$\sitV\models\forall x_1\typ X_1\ldots\forall x_n\typ X_n.A$.
\item When no context is given, $\models A$ (i.e. ``$A$ holds'') means 
that $A$ holds in the empty context over any object $S$ of~$\cS$ (or just 
over the final object).
\end{itemize}

\nxsubpoint (a) Notice that existence in Kripke--Joyal semantics 
always actually means {\em local existence}, and that universality means 
universality after arbitrary pullbacks.

(b) Example: if $X$ and $Y$ are sheaves, and $f:X\to Y$ is a morphism of 
sheaves, thus defining a constant $f$ of type $X\to Y$, then the 
``surjectivity condition'' $\forall y\typ Y.\exists x\typ X.f(x)=y$ means 
(in Kripke--Joyal semantics) that for any object~$S$ of~$\cS$
and any section $y\in Y(S)$ one can find a cover $\{S_\alpha\to S\}$
and sections $x_\alpha\in X(S_\alpha)$, such that $f_{S_\alpha}(x_\alpha)=
y|_{S_\alpha}$, i.e.\ the surjectivity of~$f$ as a map of sheaves.

(c) Suppose that $\sitV\models\exists!x\typ X.A[x]$, 
i.e.\ ``such an $x\in X$ exists and is unique''. First of all, existence
$\exists x.A[x]$ means local existence, i.e.\ we have a cover 
$\{S_\alpha\to S\}$ and elements $x_\alpha\in X(S_\alpha)$, having the 
property expressed by~$A$. Now the uniqueness implies that the restrictions 
of $x_\alpha$ and $x_\beta$ to $S_\alpha\times_SS_\beta$ coincide, so we 
can glue them to section $x\in X(S)$ with property $A$, $X$ being a sheaf 
and $A$ being a local property. In other words, {\em we don't need to 
pass to a cover to find an object that exists and is unique.} 
This section $x\in X(S)$ is actually the value of term $\iota x\typ X.A[x]$ 
in context~$\sitV$.

(d) Notice that the proposition $A\vee\neg A$ doesn't usually hold,
i.e.\ we have indeed to forget the law of excluded middle and its 
consequences. One can check that the logical laws still applicable 
in these situations are exactly those of intuitionistic logic. For example, 
we have the {\em modus ponens:} $(A\Rightarrow B)\& A\Rightarrow B$. 
We have also to forget about the axiom of choice and Hilbert's tau 
($\tau_x A$ chooses any $x$, for which $A[x]$ holds, whenever this is possible,
thus implying the axiom of choice). All we have is Hilbert's iota that
singles out objects that exist and are unique.

(e) We abbreviate $\forall x_1\typ X_1.\forall x_2\typ X_2.\forall\ldots$ into 
$\forall x_1\typ X_1,x_2\typ X_2,\ldots$; and we write $x,y,z\typ X$ instead 
of $x\typ X,y\typ X,z\typ X$.

(f) Notice that our system is actually redundant. For example, 
$\lambda x\typ X.A$ 
can be replaced with $\iota f\typ X\to Y.\forall x\typ X.f(x)=A$. However,
none of basic logical operations $\vee$, $\wedge$ and $\Rightarrow$ can 
be expressed in terms of the others (but $\neg A$ is actually equivalent to
$A\Rightarrow\st0$).

\nxsubpoint (Global and local notions.)
When we need to distinguish the intuitionistic notions arising 
from Kripke--Joyal semantics
from those used in classical sense,
we call them {\em local}, e.g.\ local existence, local disjunction etc.,
as opposed to {\em global\/} notions: global existence, global disjuction
etc.

Notice that $\sitV\models A$ is a classical statement, obeying the rules of 
classical logic, while $A$ itself is intuitionistic.

\nxsubpoint (Topos case.)
Of course, Kripke--Joyal semantics can be applied directly to any topos~$\cE$,
considered as a site with respect to its canonical topology. In this case 
$\widetilde{\cE_{/S}}\cong\cE_{/S}$, i.e.\ the ``sheaf types'' of~$\cE$ 
are actually objects of $\cE$ or $\cE_{/S}$. Some rules can be simplified 
by using the existence of (small) sums $\bigsqcup_\alpha S_\alpha$, 
together with $X(\bigsqcup_\alpha S_\alpha)\cong\prod_\alpha X(S_\alpha)$, 
for any sheaf~$X$. For example, ``local existence'' $\exists x.A[x]$ in
some context $\sitV/S$ can be understood as follows: ``there is an 
{\em epimorphism\/} $S'\twoheadrightarrow S$ and an element $x\in X(S')$
having the required property in the pullback $\sitV\times_SS'$''. Another 
example: $\sitV\models A_1\vee A_2$ iff there is an epimorphism 
$S_1\sqcup S_2\twoheadrightarrow S$, such that $\sitV\times_SS_i\models A_i$,
$i=1,2$.

\nxsubpoint (Localness of all propositions.)
Notice that the localness of all properties expressed by propositions~$A[x]$ 
of Kripke--Joyal semantics is actually a consequence of logical 
quantifier elimination and introduction rules. For example,
\begin{equation}
(\forall_R)\;\frac{\displaystyle B\Rightarrow A[a]}
{\displaystyle B\Rightarrow \forall x\typ X.A[x]}
\quad\text{if $a$ is free in~$A$ and doesn't occur in~$B$}
\end{equation}
when applied to $B=A[x]$ and a new variable $z:\st1$, 
yields $A[x]\Rightarrow\forall z.A[x]$, hence 
that $A[x]$ holds in all pullbacks of the original context whenever 
it holds in this context itself. The second 
localness condition follows similarly from $\exists z.A[x]\Rightarrow A[x]$,
a consequence of rule
\begin{equation}
(\exists_L)\;\frac{\displaystyle A[a]\Rightarrow C}
{\displaystyle (\exists x\typ X.A[x])\Rightarrow C}
\quad\text{if $a$ is free in~$A$ and doesn't occur in $C$}
\end{equation}

\nxsubpoint\label{sp:presh.vs.sh} (Presheaves vs.\ sheaves.)
Notice that almost all constructions of Kripke--Joyal semantics appear to
be applicable to the case when the ``sheaf types'' are allowed to be 
arbitrary presheaves~$X$. However, the localness of proposition $x=y$,
where $x,y\typ X$, implies separability of presheaf~$X$; and the
axiom for the iota-symbol
\begin{equation}
(\iota)\;\frac{\displaystyle B\Rightarrow \exists!x\typ X.A[x]}
{\displaystyle B\Rightarrow A[\iota x\typ X.A[x]]}
\end{equation}
actually expresses the sheaf condition for~$X$.

In other words, the sheaf condition is necessary to pick up elements 
of ``sets'' characterized by some property $A[x]$, once it is shown that 
such an element exists and is unique.

\nxsubpoint\label{sp:prestk.vs.stk} (Prestacks vs.\ stacks.)
Similarly, one might think that we might use prestacks $\cC$ instead of 
stacks in these considerations, since all we formally need is that 
all $\iHom_{\cC_{/S}}(x,y)$, $x,y\in\Ob\cC_S$, be sheaves. However,
the stack condition for $\cC$ is actually equivalent to the ability to pick 
up objects of a ``category'' $\cC$, characterized by some property $A[x]$ 
uniquely up to a {\em unique\/} isomorphism, once the (local) existence
of such an $x$ is shown.

Indeed, once we know that such an $x$ exists (locally) in some context
$\sitV/S$, we can find a cover $\{S_\alpha\to S\}$ and objects 
$x_\alpha\in\cC(S_\alpha)$ with property~$A$. Using uniqueness, we obtain
isomorphisms $\theta_{\alpha\beta}$ between pullbacks of $x_\alpha$ and
$x_\beta$ to $\cC(S_\alpha\times_SS_\beta)$. Finally, pulling everything back
to $S_\alpha\times_SS_\beta\times_SS_\gamma$, and using uniqueness of these
isomorphisms, we obtain the cocycle relation for the $\theta$s, i.e.
we've got a descent datum. Now if $\cC$ is a stack, this descent datum is 
effective, and we get an object $x\in\cC(S)$ with required property.

In this way the stack condition is required in our ``intuitionistic 
category theory'' to be able to pick up objects, defined uniquely up to a 
unique isomorphism. This is important since we don't have the axiom of choice
or Hilbert's $\tau$ to do this in the usual way. For example, this ability
is important to construct initial or final objects, adjoint functors, 
or to show that any equivalence (defined here as a fully faithful 
essentially surjective functor) admits a quasi-inverse (adjoint) equivalence.

\nxsubpoint (Transitivity of pullbacks and stacks.)
The reader may have noticed a certain problem with our constructions.
Namely, if a context~$\sitV$ involves some variable $x$ of ``large'' 
or ``stack'' type $\Ob\cC$, then the pullbacks are not transitive:
$\psi^*\phi^*\sitV\neq(\phi\psi)^*\sitV$ since $(\phi\psi)^*(x)$ is known
just to be isomorphic to $\psi^*\phi^*(x)$. Actually, a problem appears 
even if we construct one pullback $\phi^*\sitV$, since the pullback 
functors $\phi^*:\cC_S\to\cC_T$, $\phi:T\to S$, are defined only up to 
an isomorphism, and the stack structure of~$\cC$ doesn't provide a 
canonical choice of~$\phi^*$.

Of course, one might tackle with these problems by replacing $\cC$ by the
$\cS$-equivalent stack~$\cC^+$, so as to have $(\phi\psi)^*=\psi^*\phi^*$. 
However, we don't like this solution for several reasons, one of them being 
that the embedding $\cC\to\cC^+$ is itself defined by a large-scale 
application of the axiom of choice.

Instead, we prefer to think of objects of $\cC$ as ``defined up to 
an isomorphism'', using the fact that we cannot distinguish two isomorphic 
objects of a category by internal means: they must have exactly the 
same properties. We have already seen in~\ptref{sp:prestk.vs.stk} that
this is a viable point of view. Therefore, all our considerations should 
not change if we replace some objects (i.e.\ values of some variables
of type~$\Ob\cC$) by other objects isomorphic to them, for example if 
we choose $\phi^*(x)$ in another way.

This is actually the reason why $x=y$ is {\em not\/} a proposition 
in our system when $x$ and $y$ are of type $\Ob\cC$, while 
$x\simeq y$ (``$x$ is isomorphic to $y$'') is one:
$x\simeq y:=(\exists f:x\to y.\exists g:y\to x.
(f\circ g=\id_y\& g\circ f=\id_x))$.

\nxsubpoint\label{sp:loc.prop.obj.mor} 
(Local properties of objects and morphisms of stacks.)
Another consequence of this philosophy is that {\em local properties 
of objects of a stack should be stable under isomorphisms.} In other 
words, a subset $P\subset\Ob\cC$ is a {\em local set\/} or {\em class
of objects of~$\cC$} if: 
(a) $\phi^*(X)\in P$ whenever $X\in P$;
(b) $X\in P$ whenever all $\phi_\alpha^*X\in P$ for some cover 
$\{\phi_\alpha:S_\alpha\to S\}$;
and (c) $X\simeq Y$ in $\cC_S$ and $X\in P$ implies $Y\in P$.

Indeed, without property (c) properties (a) and (b) do not make any sense, 
just because of different possible choices of $\phi^*X$. On the other 
hand, if $P$ satisfies these properties, we can easily transfer~$P$ 
to any stack $\cS$-equivalent to~$\cC$.

Similar remarks apply to {\em local\/} sets of morphisms~$Q$ 
in fibers of~$\cC$: they have to be closed under isomorphisms as well, i.e.\ 
if $f$ in $\Ar\cC_S$ belongs to $Q$, and $u$, $v$ are isomorphisms in~$\cC_S$,
such that $ufv$ is defined, then $ufv$ has to belong to $Q$ as well.

\nxsubpoint
Of course, our description of Kripke--Joyal semantics is by no means 
rigorous and complete from a logician's point of view. 
For example, we didn't pay enough attention
to distinguish between syntax and semantics, or between axioms, axiom 
schemata, and deduction rules, and we haven't in fact said anything about 
proofs at all (a suitable extension of Gentzen's LJ system combined with 
$\lambda$-calculus would do). 
Another obvious gap is that we haven't explained how one 
can incorporate into the system described above cartesian functors 
$F:\cC\to\cD$ (this is quite clear anyway: if $x\typ\Ob\cC$, then 
$F(x)\typ\Ob\cD$, and if $f\typ x\to y$, then $F(f)\typ F(x)\to F(y)$), and, 
more importantly, how one can construct new stacks from existing ones
(e.g.\ $\cC\times_\cS\cD$, corresponding to 
``product of intuitionistic categories'', or $\catCart_\cS(\cC,\cD)$, 
or $\cC_{/x}$, $x\in\Ob\cC_e$\dots). We think these things are already 
clear enough, and that we'll be able to perform them by ``external'' means
if necessary.

\nxpointtoc{Model stacks}
Now we are going to present the definition of a model stack over a site or 
topos, and the construction of its homotopic category. The main idea 
here is to deal with a stack as a ``intuitionistic category'', transfer 
Quillen's definitions and proofs to the intuitionistic case, and interpret 
the result in Kripke--Joyal semantics, thus regaining descriptions in terms
of stacks over sites.

\nxsubpoint\label{sp:loc.cls.mor} (Local classes of morphisms.)
We fix a stack $\cC$ over a site $\cS$ for the most part of this subsection.
Sometimes we'll assume for simplicity that $\cS$ is closed under finite 
projective limits and its topology is subcanonical, even if this is 
inessential for most statements. 
Usually we'll have three {\em local\/} classes of morphisms in fibers 
of~$\cC$ (cf.~\ptref{sp:loc.prop.obj.mor}), called {\em fibrations},
{\em cofibrations\/} and {\em weak equivalences\/} (cf.~\ptref{def:mod.cat}).
Of course, a morphism $f:X\to Y$ in $\cC_S$ is said to be an {\em acyclic 
fibration\/} (resp.\ {\em acyclic cofibration})
if it is both a fibration (resp.\ cofibration) and a weak equivalence, 
and we say that an object $X\in\Ob\cC_S$ is {\em cofibrant\/} 
(resp.\ {\em fibrant}) if $\emptyset_{\cC_S}\to X$ is a cofibration,
resp.\ if $X\to e_{\cC_S}$ is a fibration. Clearly, these are local 
classes of objects and morphisms in the sense of~\ptref{sp:loc.prop.obj.mor}.

\nxsubpoint (Retracts and local retracts.)
Recall that a morphism $g:Z\to T$ is a retract (or a {\em global retract\/})
of another morphism $f:X\to Y$ in $\cC_S$ iff  
there are morphisms $i:Z\to X$, $j:T\to Y$, $p:X\to Z$, $q:Y\to T$, 
such that $f\circ i=j\circ g$, $g\circ p=q\circ f$, $p\circ i=\id_Z$ and 
$q\circ j=\id_T$ (cf.~\eqref{eq:diag.retr}).

We can express this condition by formula $\exists i:Z\to X.\exists j:T\to Y.
\exists p:X\to Z.\exists q:Y\to T.(f\circ i=j\circ g\&g\circ p=q\circ f\&
p\circ i=\id_Z\&q\circ j=\id_T)$ of signature $X$, $Y$, $Z$, $T:\Ob\cC$, 
$f:X\to Y$, $g:Z\to T$. Interpreting this formula in Kripke--Joyal 
semantics, we obtain the notion of a {\em local retract}:
$g$ is a local retract of $f$ iff such morphisms $i$, $j$, $p$, $q$ as 
above exist locally, i.e.\ iff $g|_{S_\alpha}$ is a retract of $f|_{S_\alpha}$
in $\cC_{S_\alpha}$ for some cover $\{S_\alpha\to S\}$ in~$\cS$.

Usually we introduce the ``local'' or ``intuitionistic'' counterparts 
of usual (i.e.\ ``global'' or ``classical'') notions by a similar procedure,
that needn't be explicit each time.

For example, if $(P)$ is a local property of morphisms in fibers of~$\cC$, 
we can say that $f:X\to Y$ in $\cC_S$ is a local retract of 
an (unspecified) morphism with property $(P)$ iff 
there is some cover $\{S_\alpha\to S\}$,
such that each $f|_{S_\alpha}$ is a retract in $\cC_{S_\alpha}$ of 
a morphism $u_\alpha:Z_\alpha\to T_\alpha$ having property $(P)$.

\nxsubpoint (Local lifting properties.)
Similarly, given two morphisms $i:A\to B$ and $p:X\to Y$ in~$\cC_S$, 
we say that
{\em $i$ has the local left lifting property (local LLP or lLLP) with 
respect to $p$}, or that {\em $p$ has the lRLP with respect to~$i$},
if for any $T\to S$ in $\cS$ and any morphisms $u:A|_T\to X|_T$ and 
$v:B|_T\to Y|_T$, such that $v\circ i|_T=p|_T\circ v$, there is a cover
$\{T_\alpha\to T\}$ and morphisms $h_\alpha:B|_{T_\alpha}\to X|_{T_\alpha}$,
such that $h_\alpha\circ i|_{T_\alpha}=u|_{T_\alpha}$ and
$p|_{T_\alpha}\circ h_\alpha=v|_{T_\alpha}$. Again, this condition is 
obtained simply by interpreting the usual definition of lifting properties
in Kripke--Joyal semantics. We can express it by means of the 
following diagram:
\begin{equation}\label{eq:diag.loclift}
\xymatrix{
A\ar[d]^{i}\ar@{-->}[r]^{\forall u}&X\ar[d]^{p}\\
B\ar@{-->}[r]^{\forall v}\ar@{.>}[ur]|{\exists h}&Y}
\end{equation}
This diagram is understood as follows. First, we are given or can construct 
the solid arrows ($i$ and~$p$ in this case) in $\cC_S$, without pulling 
anything back. Then we pull back with respect to any $T\to S$ 
and choose arbitrarily the dashed arrows $u$ and $v$, 
so as to make the diagram commutative. 
After this the dotted arrow $h$ exists only after some other pullbacks
$\{T_\alpha\to T\}$; they must constitute a cover of~$T$ because of the 
$\exists$ sign.

\begin{DefD}
A stack~$\cC$ over a site~$\cS$ with three distinguished {\em local} 
classes of morphisms in fibers as in~\ptref{sp:loc.cls.mor} is said to be 
a {\bf model stack} if the conditions (MS1)--(MS5) hold:
\end{DefD}

{\myindent{(MS1)} All fiber categories $\cC_S$ are closed under 
arbitrary (small) projective and inductive limits, and all pullback 
functors $\phi^*:\cC_S\to\cC_T$ admit both left and right adjoints
$\phi_!$, $\phi_*:\cC_T\to\cC_S$.

\myindent{(MS2)} Each distinguished class of morphisms is (local and) 
stable under global retracts (in each fiber~$\cC_S$).

\myindent{(MS3)} (``2-out-of-3'') Given $X\stackrel f\to Y\stackrel g\to Z$ 
in $\cC_S$, such that two of $f$, $g$ and $gf$ are weak equivalence, 
then so is the third.

\myindent{(MS4)} (Lifting.) Any cofibration $i:A\to B$ in~$\cC_S$ has 
the local LLP with respect to all acyclic fibrations $p:X\to Y$, and 
any acyclic cofibration $i:A\to B$ has the local LLP with respect to all
fibrations $p:X\to Y$. In other words, for any cofibration $i:A\to B$
and fibration $p:X\to Y$ in~$\cC_S$, one of them being acyclic, and 
any morphisms $u:A\to X$, $v:B\to Y$ in~$\cC_S$, such that
$vi=pu$, one can {\em locally\/} find $h_\alpha:B|_{S_\alpha}\to 
X|_{S_\alpha}$, making~\eqref{eq:diag.loclift} commutative.

\myindent{(MS5)} (Factorization.) Any morphism $f:X\to Y$ in $\cC_S$ 
can be {\em globally\/} factorized into 
$X\stackrel\alpha\to Z\stackrel\beta\to Y$ and 
$X\stackrel\gamma\to W\stackrel\delta\to Y$ (in~$\cC_S$), where
$\alpha$ and $\gamma$ are cofibrations, $\beta$ and $\delta$ are fibrations,
and $\alpha$ and $\delta$ are weak equivalences.

}
\nxsubpoint {\bf Remarks.}
(a) One can show that (MS1) is actually equivalent to the more natural 
requirement of~$\cC$ to be closed under arbitrary {\em local\/}
inductive and projective limits (taken over abritrary inner categories $\cI$
in~$\tilde\cS$, i.e.\ essentially sheaves of small categories). 
However, the axiom (MS1) is technically 
simpler to verify, and it doesn't require any knowledge of local limits.

(b) Notice that (MS2) actually implies that all distinguished classes are 
closed under {\em local\/} retracts as well, and that for example even 
if $f:X\to Y$ is local retract of an unspecified fibration, then it is 
a fibration itself. Therefore, we might replace (MS2) by the more naturally 
looking condition of each distinguished class to be stable under local 
retracts.

(c) Similarly, the explicit description given in (MS4) seems too weak,
because the statement ``$i:A\to B$ has the local LLP with respect to all 
acyclic cofibrations'' actually means that we are free to choose 
an acyclic cofibration $p: X\to Y$ after making an arbitrary pullback, 
and then $u$ and $v$ can be chosen after another arbitrary pullback. 
However, all these additional pullbacks 
are not necessary, once we know that all distinguished classes are local
and in particular stable under pullbacks.

\nxsubpoint
(d) The most interesting is our form of (MS5). Indeed, the natural local form
would be a (strictly) weaker statement:

{\myindent{(MS5l)} (Local factorization.) 
Any morphism $f:X\to Y$ in $\cC_S$ 
can be {\em locally\/} factorized into $X|_{S_i}\stackrel{\alpha_i}\to 
Z_i\stackrel{\beta_i}\to Y|_{S_i}$
and $X|_{S_i}\stackrel{\gamma_i}\to W_i\stackrel{\delta_i}\to Y|_{S_i}$ 
(in~$\cC_{S_i}$ for some cover $\{S_i\to S\}$), where
$\alpha_i$ and $\gamma_i$ are cofibrations, $\beta_i$ and $\delta_i$ 
are fibrations, and $\alpha_i$ and $\delta_i$ are weak equivalences.

}
The reason why we require the stronger axiom (MS5) is that in fact all
model categories people really use (e.g.\ all cofibrantly generated 
model categories) admit functorial factorizations (e.g.\ constructed by 
means of Quillen's small object argument) in (M5), and actually some 
authors {\em require\/} the existence of such functorial factorizations 
in the definition of a model category (cf.\ e.g.\ \cite{Hovey}). 
This means that all model stacks we consider actually satisfy a
{\em stronger\/} version (MS5f) of (MS5), which in particular implies (MS5):

{\myindent{(MS5f)} (Functorial factorization.)
One can choose factorizations of (MS5) functorially in $f\in\Ar\cC_S$,
i.e.\ we actually have functors $W$, $Z:\Ar\cC_S\to\cC_S$, 
$\alpha$, $\beta$, $\gamma$, $\delta:\Ar\cC_S\to\Ar\cC_S$, such that 
for any $f:X\to Y$ in $\Ar\cC_S$ we get 
$X\stackrel{\alpha(f)}\longto Z(f)\stackrel{\beta(f)}\longto Y$, 
$X\stackrel{\gamma(f)}\longto W(f)\stackrel{\delta(f)}\longto Y$ 
with the properties listed in (MS5). Moreover, these functors are
compatible with base change functors $\phi^*:\cC_S\to\cC_T$, i.e.\ they
extend to cartesian functors between appropriate stacks.

}
However, (MS5) seems to be sufficient for almost all our constructions,
so we don't insist on requiring (MS5f). On the other hand, (MS5) 
allows us to construct (globally) fibrant and cofibrant replacements, 
something we wouldn't be able to do having only~(MS5l). This somewhat 
simplifies the exposition and allows us to construct globally more than 
it would be possible by only local means.

\nxsubpoint 
(Basic properties of distinguished classes.)
One checks, essentially in the classical way, that for example the 
acyclic cofibrations are {\em exactly\/} the morphisms in fibers of~$\cC$
that have the local LLP with respect to all fibrations (of course, 
if we fix some $i:A\to B$ in~$\cC(S)$, we are free to choose a fibration
$p:X\to Y$ and morphisms $u:A\to X$ and $v:B\to Y$ after any pullback). 
Furthermore, it is easy to see that the class of morphisms $P'$ in fibers of a 
stack~$\cC$, characterized by the local LLP with respect to some other
{\em local\/} class of morphisms $P$, is itself local, and stable under 
composition, pushouts, {\em finite\/} direct sums and 
local retracts (but not sequential inductive limits!) in the fibers of~$\cC$,
as well as under all pullback functors $\phi^*$. 
This applies in particular to cofibrations and acyclic cofibrations
in a model stack,
and the classes of fibrations and acyclic fibrations have dual properties,
e.g.\ stability under pullbacks, finite products, and~$\phi^*$.

Similarly, one sees immediately from (MS5) and (MS3) that the weak equivalences
are exactly those morphisms that can be factorized into an acyclic cofibration
followed by an acyclic fibration. Therefore, any two of distinguished local
classes of a model stack~$\cC$ completely determine the remaining class.

\nxsubpoint
(Absence of stability under sequential inductive limits.)
Notice that the class of morphisms $P'$ characterized by the local LLP 
with respect to some other local class $P$ needn't be closed under 
sequential limits. In other words, if $A_0\stackrel{i_0}\to A_1
\stackrel{i_1}\to A_2\to\cdots$ is an inductive system in $\cC(S)$,
such that each $i_n:A_n\to A_{n+1}$ has the local LLP with respect to
some $p:X\to Y$ in $\cC(S)$ (and all its pullbacks), we cannot conclude
that $i:A:=A_0\to B:=\injlim_nA_n$ has the same property. Indeed, let us 
try to repeat the usual proof. Fix some $v:B\to Y$ and denote by
$Z_n:=\iHom_Y(A_n,X)$ the subsheaf of $\iHom(A_n,X)$ consisting of all
local liftings of~$v$, i.e.\ $Z_n=\{f:A_n\to X\,|\,p\circ f=v\circ j_n\}$
in Kripke--Joyal semantics, where $j_n:A_n\to B$ is the natural embedding.
Put $Z:=\iHom_Y(B,X)=\iHom_Y(\injlim_nA_n,X)=\projlim_n Z_n$. Then the 
local LLP of $i_n$ means that each $i_n^*:Z_{n+1}\to Z_n$ is an epimorphism;
and hypothetical local LLP of $i$ would mean that $i^*:Z\to Z_0$ is an
epimorphism as well.

In other words, we have a projective system $\cdots\to Z_2\to Z_1\to Z_0$ 
with epimorphic transition morphisms in a topos $\cE=\tilde\cS$, and 
we would like to conclude that $Z=\projlim_nZ_n\to Z_0$ is also 
epimorphic. However, the usual proof of this statement for $\cE=\catSets$ 
invokes the (dependent countable) axiom of choice, so it cannot be transferred
to the topos case. In fact, one can construct examples where $Z\to Z_0$ is
not epimorphic while all $Z_{n+1}\to Z_n$ are. Put $\cE:=\cB_{\hat\bbZ}$; 
this is the topos of discrete $\hat\bbZ$-sets, i.e.\ sets~$X$ with an
action of $\hat\bbZ$, such that the stabilizer of any point $x\in X$ is 
open in~$\hat\bbZ$. Epimorphisms in this topos are just the surjective
maps of discrete $\hat\bbZ$-sets, and projective limits can be computed 
by taking the subset of all points of the usual projective limit having 
an open stabilizer.
Next, put $Z_n:=\bbZ/p^n\bbZ$, and let the generator
$\sigma$ of $\hat\bbZ$ act on $Z_n$ by adding one: $\sigma(x)=x+1$. 
Clearly, these $\hat\bbZ$-sets $Z_n$ together with canonical projections
$Z_{n+1}\to Z_n$ define a projective system in $B_{\hat\bbZ}$ with 
epimorphic transition morphisms. However, the projective limit $Z$ of 
this system in $\cB_{\hat\bbZ}$ is an empty set, as well as the product
$\prod_nZ_n\supset Z$, the stabilizer of any point of corresponding 
set-theoretical product being of infinite index in~$\hat\bbZ$.

\nxsubpoint (Fiberwise dual of a model stack.)
Given a fibered category $\cC\to\cS$, we define its {\em fiberwise dual\/} 
or {\em opposite\/} $\cC^{fop}\to\cS$ by replacing each 
fiber category $\cC(S)$ by its opposite $\cC(S)^0$ 
while preserving the pullback functors $\phi^*:\cC(S)\to\cC(T)$. Of course, 
this definition doesn't actually depend on the choice of pullback functors
$\phi^*$. Now we see immediately that {\em the axioms of a model stack are 
self-dual}, i.e.\ $\cC^{fop}$ is a model stack over~$\cS$ whenever
$\cC$ is one (of course, one has to interchange fibrations and cofibrations 
while dualizing). This remark enables us to deduce new statements about
model stacks by duality.

\nxsubpoint (Flat stacks.)
We say that a fibered category $\cC\to\cS$ satisfying (MS1) 
over a category $\cS$ with fibered products is {\em flat\/} if the
following condition is fulfilled:

{\myindent{(MS1+)} For any cartesian square in~$\cS$
\begin{equation}
\xymatrix{
T'\ar[r]^{g}\ar[d]^{q}&T\ar[d]^{p}\\
S'\ar[r]^{f}&S}
\end{equation}
\setmyindent the canonical morphism $q_!g^*\to f^*p_!$ is an isomorphism 
of functors $\cC(T)\to\cC(S')$, hence the same is true for its adjoint
$p^*f_*\to g_*q^*$ as well.

}
Notice that we have actually {\em two\/} canonical morphisms $q_!g^*\to
f^*p_!$. The first is deduced by adjointness from $g^*\to q^*f^*p_!\cong
 g^*p^*p_!$, itself obtained by applying $g^*\star-$ to the unit 
$\Id_{C(T)}\to p^*p_!$. The second one is deduced by adjointness from 
$f_!q_!g^*\cong p_!g_!g^*\to p_!$, obtained by applying $p_!\star-$ to
the counit $g_!g^*\to\Id_{C(S')}$. In most cases it is clear that these two 
coincide; let us include for simplicity the equality of these two canonical 
morphisms $q_!g^*\to f^*p_!$ as an additional requirement in (MS1+).

\nxsubpoint\label{sp:equiv.flat.stk} (Equivalent formulations of flatness.)
Let $\cC$ be a stack over~$\cS$, and 
$T\stackrel\phi\to S$ be a morphism in~$\cS$.
We denote by $\cC^{T/S}$ or $\cC_{/S}^T$ (or simply $\cC^T$ when $S=e_\cS$) 
the direct image $\phi_*^{st}\cC_{/T}$ of stack $\cC_{/T}:=\cC\times_\cS
\cS_{/T}$ over $\cS_{/T}$. Clearly, $\cC_{/S}^T(S')=
\cC(T\times_SS')$ for any $S'\in\Ob\cS_{/S}$, so the pullback functors 
$\phi_{S'}^*:\cC(S')\to\cC(T\times_SS')$ for $\phi_{S'}:T\times_SS'\to S'$ 
combine together to a {\em cartesian\/} functor $\bphi^*:\cC_{/S}\to\cC_{/S}^T$
over~$\cS_{/S}$. Now the flatness condition (MS1+) can be interpreted as 
a requirement for all $\bphi^*$ to admit {\em cartesian\/} $\cS_{/S}$-adjoint
functors $\bphi_!$, $\bphi_*:\cC_{/S}^T\to\cC_{/S}$. This adjointness
can be interpreted in yet another way as a functorial $\cS_{/S}$-isomorphism
of local Hom-sheaves, i.e.\ an $\cS_{/S}$-isomorphism of cartesian functors
$(\cC_{/S}^T)^{fop}\times\cC_{/S}\to\stSETS_{\cS_{/S}}$:
\begin{multline}\label{eq:locadj.phi!.phi*}
\iHom_{\cC_{/S}^T|S'}(X,\bphi^*Y)\cong\iHom_{\cC_{/S}|S'}(\bphi_!X,Y)
\quad\text{for all $S'\in\Ob\cS_{/S}$,}\\
\text{$X$ in $\cC_{/S}^T(S')=\cC(T\times_SS')$,
$Y$ in $\cC_{/S}(S')=\cC(S')$.}
\end{multline}
We obtain a similar ``local adjointess'' interpretation for $\bphi^*$ and
$\bphi_*$ of course.

Notice that the Kripke--Joyal philosophy insists that we should always use
Hom-sheaves, never Hom-sets; therefore, \eqref{eq:locadj.phi!.phi*} is 
the only correct way to discuss adjoint functors to $\bphi^*$ from this point
of view. Therefore, it might be very natural to include (MS1+) into the list of
axioms for model stacks and consider only flat model stacks.

From the intuitionistic point of view $\cC_{/S}$ is an 
``intuitionistic category'' (over base site $\cS_{/S}$), $T\in\Ob\cC_{/S}$
is an ``intuitionistic set'', and $\cC_{/S}^T$ is ``the category of 
families of objects of $\cC_{/S}$ indexed by~$T$'', i.e.\ some sort of
(local) product of categories. Furthermore, $\bphi^*:\cC_{/S}\to\cC_{/S}^T$ 
is the ``constant family functor'', and its left and right adjoints
$\bphi_!=:\coprod_{T/S}$ and $\bphi_*=:\prod_{T/S}:\cC_{/S}^T\to\cC_{/S}$ 
should be thought of as ``(local) coproducts and products of families of 
objects of category $\cC_{/S}$ indexed by~$T$''. 
Then (MS1+) assures us that these local coproducts and products are indeed
``local'' or ``universal'', i.e.\ compatible with pullbacks.

\nxsubpoint\label{sp:loc.sum.prod} (Local sums and products in~$\cC$.)
When we have $T\stackrel\phi\to S$ in~$\cS$, and an object $X\in\Ob\cC(T)$ 
``of $\cC$ over $T$'', we denote $\phi_!X$ by $\coprod_{T/S}X$ or
$\bigoplus_{T/S}X$, 
and $\phi_*X$ by $\prod_{T/S}X$. When $S$ is the final object $e_\cS$,
we write simply $\coprod_TX$ or $\bigoplus_TX$, and $\prod_TX$, 
respectively. These notations are motivated by usual notations for sheaves 
and presheaves (cf. SGA~3), as well as the case $\cS=\catSets$: then
$\cC(I)=\cC(\st1)^I$, and $\prod_I$, $\coprod_I:\cC(\st1)^I\to\cC(\st1)$ 
are just the usual product and coproduct functors.

Furthermore, we denote $\phi_!\phi^*X$ by $X\times_ST$ or $X\otimes_ST$ 
(omitting $S$ if $S=e_\cS$), and $\phi_*\phi^*X$ by $\iHom_S(T,X)$
($\iHom(T,X)$ or even $X^T$ for $S=e_\cS$). Again, the motivation comes from
the case $\cS=\catSets$. In any case $X\mapsto X\otimes_ST$ and
$Y\mapsto\iHom_S(T,Y)$ are adjoint: $\Hom_S(X\otimes_ST,Y)\cong
\Hom_T(X|_T,Y|_T)\cong\Hom_S(X,\iHom_S(T,Y))$.

\nxsubpoint (Local limits.)
Suppose we are given inner category $C$ in a topos~$\cE$, given by its 
objects of morphisms $C_1$ and of objects $C_0\in\Ob\cE$, together with the 
source and target morphisms $s$, $t:C_1\to C_0$, the identity morphism
$i:C_0\to C_1$, and the composition morphism 
$\mu:C_2:=C_1\times_{s,t}C_1\to C_1$,
subject to usual identity and associativity conditions. Suppose that 
we have an inner functor $F$ from $C$ to a model stack $\cD\stackrel p\to\cE$, 
i.e.\ a inner category $F$ in $\cD$, such that $p(F)=C$, and
$s_F:F_1\to F_0$ is cartesian. Then we can compute the corresponding 
inner or local limits $\injlim_CF$ and $\projlim_CF$ in the classical way:
$\injlim_CF:=\Coker(\coprod_{C_1}F_1\rightrightarrows\coprod_{C_0}F_0)$,
and $\projlim_CF:=\Ker(\prod_{C_0}F_0\rightrightarrows\prod_{C_1}F_1)$.

When we have the flatness condition (MS1+), all these local sums, products 
and limits commute with arbitrary base change, so these notions fit 
nicely into the Kripke--Joyal philosophy: they might be thought of 
as ``limits in an intuitionistic category~$\cD$ along a small 
intuitionistic index category~$C$''.

\nxsubpoint\label{sp:modstk.sites.topoi} 
(Model stacks over sites and topoi.)
Recall that we have a correspondence (a 2-equivalence of 2-categories, 
actually) between stacks $\cC$ over a site $\cS$, and stacks $\tilde\cC$
over the corresponding topos~$\tilde\cS$. Under this correspondence
$\cC$ is simply the ``restriction'' of~$\tilde\cC$ to~$\cS$, i.e.\ it is
$\tilde\cC\times_{\tilde\cS}\cS$, and conversely, 
$\tilde\cC$ is recovered from $\cC$ by restricting $\cC^+$ from $\hat\cS$
to $\tilde\cS$: $\tilde\cC:=\cC^+\times_{\hat\cS}\tilde\cS$
(cf.~\ptref{sp:stk.over.site.topos}). Now, if we have for example a 
local class of morphisms $\tilde P$ in the fibers of~$\tilde\cC\to\tilde\cS$,
we obtain a local class of morphisms $P$ in the fibers of $\cC\to\cS$
simply by ``restriction'' since $\cC(S)=\tilde\cC(\tilde S)$, where 
$\tilde S$ is the sheafification of the presheaf represented by~$S$, and
conversely, starting from $P\subset\Ar\cC$, we can extend it uniquely 
to a local class $\tilde P\subset\Ar\tilde\cC$ as follows:
a morphism $f:X\to Y$ in $\tilde\cC(X)$, $X\in\Ob\tilde\cS$, belongs to
$\tilde P$ iff its pullbacks 
$\phi^*(f)\in\Ar\tilde\cC(\tilde S)\cong\Ar\cC(S)$ 
with respect to all $\phi:\tilde S\to X$ belong to~$P$. It is immediate that 
these constructions are inverse to each other, once we take into account
that local properties can be also transferred along any $\cS$-equivalence 
of stacks (cf.~\ptref{sp:loc.prop.obj.mor}).

In particular, we can transfer the three distinguished classes entering into 
the structure of a model stack from stacks $\cC\to\cS$ to stacks 
$\tilde\cC\to\tilde\cS$, and conversely. Now it is easy to see that 
the axioms of a model stack hold for $\cC$ iff they hold for $\tilde\cC$, 
with the following two exceptions. Firstly, (MS5) cannot be transfered 
from $\cC$ to $\tilde\cC$, but both (MS5l) and (MS5f) can, and since 
all model stacks we consider satisfy the stronger condition (MS5f), 
this is not really a big problem for our considerations. Secondly, 
it is not immediate that (MS1) extends from $\cC$ to $\tilde\cC$
(the other direction is clear). We'll check this in a moment; 
for now, let us state

\begin{PropD}
The 2-category of model stacks over a site~$\cS$, satisfying (MS5f),
is 2-equivalent to the 2-category of model stacks over the corresponding 
topos~$\tilde\cS$, satisfying (MS5f). The same applies if we 
replace in the definition of model stacks (MS5) with a weaker local 
condition (MS5l). Moreover, any model stack over topos~$\tilde\cS$ 
defines by restriction a model stack over $\cS$ without any additional 
requirements.
\end{PropD}

In this way it is essentially the same thing to consider model stacks over 
a site or over the corresponding topos.

\begin{Proof}
The only problem is to transfer (MS1) from $\cC$ to~$\tilde\cC$.
(a) Let's prove that all $\tilde\cC(X)$, $X\in\Ob\tilde\cS$, are closed 
under arbitrary projective limits. By definition $\tilde\cC(X)=\cC^+(X)=
\catCart_\cS(\cS_{/X},\cC)$, so suppose we want to compute a projective
limit $F=\projlim F_i$ of cartesian functors $F_i:\cS_{/X}\to\cC$. To this 
end we define $F$ fiberwise, i.e.\ for any $S\stackrel\xi\to X$ in $\cS_{/X}$
we compute $F(\xi):=\projlim F_i(\xi)$ in the fiber $\cC(S)$. Now we observe 
that for all $\phi:T\to S$ the functor $\phi^*:\cC(S)\to\cC(T)$ commutes 
with arbitrary projective limits, having a left adjoint $\phi_!$ by~(MS1),
hence $\phi^*F(\xi)\cong F(\xi\phi)$, i.e.\ we've indeed constructed a 
cartesian functor $F:\cS_{/X}\to\cC$, clearly having the universal properly 
required from $\projlim F_i$. The case of inductive limits is dealt with 
similarly, using that $\phi^*$ has a right adjoint $\phi_*$, hence commutes 
with arbitrary inductive limits as well.

(b) Now we have to show that for any morphism of sheaves $f:X\to Y$
the pullback functor $f^*:\tilde\cC(Y)\to\tilde\cC(X)$ has both left and 
right adjoints $f_!$ and $f_*$. We can replace $\cS$ by any its small 
generating category $\cS'$ since in this case $\tilde\cS'\cong\tilde\cS$,
and $\tilde\cC$ is equivalent to the extension to $\tilde\cS'$ of 
$\cC':=\cC\times_\cS\cS'$. Therefore, we can assume $\cS$ to be small. 

In this case one can actually show the statement for the pullback
functor $f^*:\cC^+(Y)\to\cC^+(X)$ for any morphism of presheaves $f:X\to Y$.
The main idea here is that $\cC^+(X)=\catCart_\cS(\cS_{/X},\cC)$, and
similarly for $\cC^+(Y)$, and under this indentification $f^*$ is the 
precomposition functor $\bar f^*$ with respect to $\bar f:\cS_{/X}\to\cS_{/Y}$.
We want to show that $\bar f^*$ admits, say, a left adjoint $\bar f_!$
(existence of $\bar f_*$ can be obtained then by fiberwise duality). 
We can construct first a ``left Kan extension functor''
$\bar f^0_!:\catFunct_\cS(\cS_{/X},\cC)\to\catFunct_\cS(\cS_{/Y},\cC)$ 
by putting $(\bar f^0_!F)(T\stackrel\eta\to Y):=
\injlim_{\cS_{/X\times_YT}}\phi_!F(\xi)$, 
where the limit is computed in~$\cC(T)$ along all diagrams
\begin{equation}
\xymatrix{
S\ar@{-->}[r]^{\xi}\ar@{-->}[d]^{\phi}&X\ar[d]^{f}\\
T\ar[r]^{\eta}&Y}
\end{equation}
The problem with $\bar f^0_!$ is that $\bar f^0_!F$ needn't be cartesian for 
a cartesian~$F$, at least if we don't have additional flatness conditions
like (MS1+). However, if we show that the inclusion 
$J:\catCart_\cS(\cS_{/Y},\cC)\to\catFunct_\cS(\cS_{/Y},\cC)$ admits a 
left adjoint $J_!$, we can put $f_!=\bar f_!:=J_!\circ\bar f^0_!$. 
To show existence of $J_!$ we define $(J_!F)(S)$ for any $\cS$-functor 
$F:\cS_{/Y}\to\cC$ by putting it equal to the inductive limit in $\cC(S)$ 
of $u_{1,!}v^*_1u_{2,!}v^*_2\cdots v^*_nF(S_n)$, taken along the category 
of all diagrams $S=S_0\stackrel{u_1}\leftarrow T_0\stackrel{v_1}\rightarrow
S_1\stackrel{u_2}\leftarrow T_2\stackrel{v_2}\to\cdots\stackrel{v_n}\to S_n$,
all integer $n\geq0$, where we allow to split any arrow in two of the same
direction, and insert couples of identity morphisms at any point. When
we have flatness (MS1+), the construction of $J_!F$ is considerably simplified:
it suffices to take the inductive limit along all diagrams with $n=1$.
\end{Proof}

We don't provide more details because in all our applications the model stack
can be constructed directly over the topos~$\tilde\cS$ and then restricted
to $\cS$ if necessary, 
so we'll need to apply the above proposition only in the obvious direction.

\nxsubpoint (Arbitrary limits in a model stack over a topos.)
One can easily show that whenever $p:\cC\to\cS$ is a fibered category 
satisfying (MS1), and certain types of (say, projective) limits exist both 
in $\cS$ and all fibers $\cC(S)$, then these types of limits exist 
in~$\cC$ as well, and $p$ commutes with these limits. In particular, 
{\em arbitrary (small) projective and inductive limits exist in a model
stack $\cC$ over a topos~$\cE$, and the projection $\cC\stackrel p\to\cE$ 
commutes with these limits.} The proof goes as follows. Suppose we want to
compute $X=\projlim_\cI X_i$ in $\cC$. We put $S_i:=p(X_i)$ and compute
$S:=\projlim_\cI S_i$ in $\cS$. Then we put $X:=\projlim_\cI\pi_i^*X_i$ 
in~$\cC(S)$, where $\pi_i:S\to S_i$ are the natural morphisms, and 
easily check (using the existence of left adjoints $\phi_!$) that
$X$ is indeed $\projlim_\cI X_i$ in~$\cC$. To compute inductive limits 
we put similarly $S:=\injlim_\cI S_i$ and $X:=\injlim_\cI\lambda_{i,!}X_i$ 
in~$\cC(S)$, where $\lambda_i:S_i\to S$ are the natural morphisms.

\nxsubpoint (Final section of $\cC\to\cS$.)
Let~$\cS$ be a category with finite projective limits, 
$\cC\stackrel p\to\cS$ be a fibered category with finite projective limits 
in each fiber~$\cC(S)$. Suppose all pullback functors $\phi^*:\cC(S)\to\cC(T)$
to be left exact. Then the above reasoning shows that finite projective limits
exist in~$\cC$ itself, and $p$ is left exact. Now let us choose a final
object $[S]$ in each fiber $\cC(S)$. Clearly, $\phi^*[S]\cong[T]$ for 
any $\phi:T\to S$ in~$\cS$; moreover, $\Hom_\cC(X,[S])\cong\Hom_\cS(p(X),S)$
for any $X\in\Ob\cC$, and in particular $\Hom_\cC([T],[S])\cong
\Hom_\cS(T,S)$. This means that the restriction of $p$ to the fibered 
subcategory $[\cS]$ of $\cC$ consisting of objects $[S]$ is an 
{\em isomorphism\/} of categories; taking its inverse we obtain a cartesian 
section $\sigma:\cS\to\cC$, clearly the final object of $\projLim_\cS\cC$.
Since $[S]$ is a final object in~$\cC(S)$, we obtain a left exact 
faithful functor $\cC(S)\to\cC_{/[S]}$, so we can write $\Hom_{[S]}(X,Y)$
or $\Hom_S(X,Y)$ instead of $\Hom_{\cC(S)}(X,Y)$, for any $X$, $Y\in\Ob\cC(S)$.
Moreover, given any $\phi:T\to S$ in $\cS$ and
any $X\to[S]$ in $\cC(S)$, the pullback $\phi^*X$ is easily seen to
be the fibered product $X\times_{[S]}[T]$ in~$\cC$. This explains our 
alternative notations $X\times_{[S]}[T]$ and $X\times_ST$ for the 
pullback $\phi^*X$, and allows us to mix fibered products in fibers of $\cC$
with fibered products corresponding to pullbacks $\phi^*X$ without risk
of confusion. Notice, however, that the product $X\times Y$ of two objects
in $\cC(S)$ corresponds to $X\times_{[S]}Y$ in~$\cC$, so it might be 
convenient to denote this product by $X\times_SY$.

\begin{DefD} (Homotopic stack of a model stack.)
Given a homotopic stack $\cC$ over a site~$\cS$, or, more generally, 
a stack $\cC$ with a local class of weak equivalences, and a cartesian functor 
$\gamma:\cC\to\HO\cC$ into another stack $\HO\cC$, we say that 
{\em $\HO\cC$ is the homotopic stack of~$\cC$} if for any stack $\cD$ over
$\cS$ the induced functor 
$\gamma^*:\catCart_\cS(\HO\cC,\cD)\to\catCart_\cS(\cC,\cD)$ is fully faithful
and its essential image consists of all cartesian functors $F:\cC\to\cD$
that transform weak equivalences into isomorphisms. We define the 
{\em homotopic prestack $\HO_p\cC$} and {\em homotopic fibered category
$\HO_0\cC$} by similar requirements, where $\cD$ runs through all 
prestacks (resp.\ fibered categories) over~$\cS$.
\end{DefD}

Notice that the above requirements determine $\HO\cC$, $\HO_p\cC$ and
$\HO_0\cC$ up to an $\cS$-equivalence. Moreover, in the definition of
$\HO_0\cC$ and $\HO_p\cC$ we can require $\gamma^*$ to induce an {\em
isomorphism} between $\catCart_\cS(\HO_?\cC,\cD)$ and the full subcategory 
of $\catCart_\cS(\cC,\cD)$ described above. Then $\HO_0\cC$ and $\HO_p\cC$
become determined up to an $\cS$-isomorphism; in this case we speak about the
{\em strict\/} homotopic fibered category or prestack of~$\cC$.

\nxsubpoint (Homotopic category of a model stack.)
We denote by $\Ho\cC$ the category $\projLim_\cS\HO\cC$ 
of global (cartesian) sections of $\HO\cC$. If $\cS$ has a final object~$e$,
we usually replace $\projLim_\cS\HO\cC$ by equivalent category $(\HO\cC)(e)$.
Categories $\Ho_p\cC$ and $\Ho_0\cC$ are defined similarly. Notice, however,
that $\Ho\cC$ and $\Ho_p\cC$ depend on the whole stack~$\cC$, not just on 
the category $\projLim_\cS\cC$ or $\cC(e)$ of its global sections.

\nxsubpoint (Existence of homotopic stacks.)
We want to show that $\HO_0\cC$, $\HO_p\cC$ and $\HO\cC$ always exist for any 
fibered category $\cC\to\cS$ with a local class of weak equivalences,
at least if we don't mind enlarging the universe. Indeed, $\HO_0\cC$ can 
be constructed fiberwise by putting $(\HO_0\cC)(S)$ equal to the localization
of $\cC(S)$ with respect to weak equivalences lying in this category; 
since all pullback functors $\phi^*:\cC(S)\to\cC(T)$ transform 
weak equivalences into weak equivalences, they induce functors 
$\phi^*:(\HO_0\cC)(S)\to(\HO_0\cC)(T)$, thus defining a fibered category 
$\HO_0\cC\to\cS$ and a cartesian functor $\gamma:\cC\to\HO_0\cC$, having 
required universal property in its strict form.

Next, we can construct $\HO_p\cC$ as the prestack associated to~$\HO_0\cC$;
recall that this can be done by leaving all objects intact but replacing 
the $\Hom$-sets in fibers by the global sections of sheafifications 
of corresponding $\iHom$-presheaves (cf.~\ptref{sp:prestack.from.fib}). 
If we combine this construction of associated prestack together 
with the previous construction of $\HO_0\cC$, 
we see that the prestack thus obtained satisfies the 
universal property of $\HO_p\cC$ in its {\em strict\/} form.

Finally, we can construct the homotopic stack $\HO\cC$ as the stack
associated to $\HO_p\cC$; according to the construction given 
in~\ptref{sp:stack.from.fib}, the natural embedding $\cS$-functor 
$\HO_p\cC\to\HO\cC$ is fully faithful. Therefore, we can recover
$\HO_p\cC$ from $\HO\cC$ by taking the essential image of $J:\cC\to\HO\cC$.
We can even recover the strict form of $\HO_p\cC$ by putting 
$\Ob\HO_p\cC:=\Ob\cC$, $\Hom_{\HO_p\cC}(X,Y):=\Hom_{\HO\cC}(J(X),J(Y))$.

\nxsubpoint\label{sp:homstk.sites.topoi} 
(Homotopic stacks over sites and topoi.)
Let~$\cC$ be a model stack over a site $\cS$. Let us denote by $\tilde\cC$
its extension to the corresponding topos $\tilde\cS$; then $\cC$ is
equivalent to the ``restriction'' $\tilde\cC|_\cS=\tilde\cC\times_{\tilde\cS}
\cS$ of $\tilde\cC$ to~$\cS$ (cf.~\ptref{sp:stk.over.site.topos}). 
Now we can construct $\HO_0\tilde\cC$, $\HO_p\tilde\cC$ and $\HO\tilde\cC$
over~$\tilde\cS$, and compare their ``restrictions'' to $\cS$ with
$\HO_0\cC$, $\HO_p\cC$ and $\HO\cC$. We claim that {\em in each of these 
cases we obtain $\cS$-equivalent categories.} Indeed, this is clear 
for $\HO_0$, because it can be computed fiberwise, hence $\cC\mapsto\HO_0\cC$
commutes with base change $\cS\to\tilde\cS$. As for $(\HO\tilde\cC)|_\cS$ and
$\HO\cC$, their equivalence follows from the 2-equivalence of 
categories of $\cS$-stacks and $\tilde\cS$-stacks 
(cf.~\ptref{sp:stk.over.site.topos}), once we take into account that
for any stack $\cD/\cS$
the natural ``restriction'' functor $\catCart_{\tilde\cS}(\tilde\cC,\tilde\cD)
\to\catCart_\cS(\cC,\cD)$ induces an equivalence between the subcategories
of cartesian functors that transform weak equivalences into isomorphisms,
the class of weak equivalences being {\em local.}

Now only the case of $\HO_p\cC$ remains. However, $\HO_p\tilde\cC$ can be 
computed fiberwise as the essential image of $\tilde\cC\to\HO\tilde\cC$, 
so this case follows from that of~$\HO\cC$ already considered.

Therefore, we can always work over a topos if we want to, restricting 
to the original site at the end if necessary. Notice that {\em a priori\/}
one would rather expect only the stack~$\HO\cC$ to depend essentially only
on the topos~$\tilde\cS$ but not the particular site $\cS$ chosen to 
represent this topos, but, surprisingly, this independence extends to
$\HO_p\cC$ and $\HO_0\cC$. However, $\HO_p\tilde\cC/\tilde\cS$ 
cannot be recovered from $\HO_p\cC/\cS$ alone: 
we need the cartesian functor $\cC\to\HO_p\cC$ as well, and the same applies 
to~$\HO_0\tilde\cC$.

\nxsubpoint (Derived cartesian functors between model stacks.)
Not surprisingly, the definition of derived cartesian functors between 
model stacks is completely similar to the classical one 
(cf.~\ptref{def:der.funct}). Namely, given a cartesian functor 
$F:\cC\to\cD$ of model categories over the same site~$\cS$, its
{\em left derived\/} $\dL F:\HO\cC\to\HO\cD$ is a cartesian functor
$\dL F$ between corresponding homotopy categories, together with a 
natural transformation $\epsilon:\dL F\circ\gamma_\cC\to
\gamma_\cD\circ F$ (over~$\cS$), such that for any other cartesian functor 
$G:\HO\cC\to\HO\cD$ 
and natural transformation $\zeta:G\circ\gamma_\cC\to\gamma_\cD\circ F$ 
there is a unique natural transformation $\theta:G\to\dL F$, for which
$\zeta=\epsilon\circ(\theta\star\gamma_\cC)$. The {\em right derived functor}
$\dR F:\HO\cC\to\HO\cD$, $\eta:\gamma_\cD\circ F\to\dR F\circ\gamma_\cC$ 
is defined similarly.

When $\dL F$ maps the full subcategory $\HO_p\cC$ of~$\HO\cC$ into the
full subcategory $\HO_p\cD$ of~$\HO\cD$, the induced functor
$\HO_p\cC\to\HO_p\cD$ will be denoted by~$\dL_p F$ or even $\dL F$. 
Furthermore, 
the induced functors between the corresponding categories of global sections
will be also denoted by $\dL F$, e.g.\ $\dL F:\Ho\cC\to\Ho\cD$.

Of course, the above definitions and constructions are essentially invariant 
under extending everything from $\cS$ to the corresponding topos~$\tilde\cS$,
so we always can work over topoi if we need to.

\nxpointtoc{Homotopies in a model stack}
Now we want to construct the homotopic stack $\HO\cC$ in a manner 
as close as possible to classical Quillen's approach, based on equivalence 
$\Ho\cC\cong\pi\cC_{cf}$, and prove or at least state corresponding 
criteria for the existence and adjointness of derived cartesian functors. 
Among other things, this will demonstrate that $\HO\cC$ is actually a
$\univU$-category, so we don't need to enlarge the universe.

Let us fix a model stack $\cC$ over a site~$\cS$.

\nxsubpoint (Fibrant and cofibrant objects.)
Recall that an object $X\in\Ob\cC(S)$ is cofibrant (resp.\ fibrant)
iff $\emptyset_{\cC(S)}\to X$ is a cofibrant morphism 
(resp.\ iff $X\to e_{\cC(S)}$ is fibrant). Since the property of a morphism
to be a (co)fibration is local, and all pullback functors $\phi^*$ preserve
the initial and final objects by (MS1), the full subcategories $\cC_c$, $\cC_f$
and $\cC_{cf}$ of $\cC$ consisting of cofibrant, fibrant, and 
fibrant-cofibrant objects are strictly full cartesian subcategories 
and substacks of~$\cC/\cS$.

\nxsubpoint (Fibrant and cofibrant replacements.)
Given an object $X\in\Ob\cC(S)$, its {\em cofibrant replacement\/} is any
acyclic fibration (or in some cases just any weak equivalence) $Q\to X$ with
a cofibrant~$Q$. Applying (MS5) to $\emptyset\to X$, we see that 
cofibrants replacements always exist ``globally'', i.e.\ in $\cC(S)$ itself. 
When we have functorial factorizations (MS5f), we can choose the cofibrant
replacements functorially, i.e.\ construct a cartesian functor $Q:\cC\to\cC_c$
and a natural transformation $\xi:Q\to\Id_\cC$, such that $\xi_X:Q(X)\to X$
is a cofibrant replacement for all $X\in\Ob\cC$.

Similar remarks apply of course to {\em fibrant replacements\/} 
$X\to R$, defined as acyclic cofibrations (or just weak equivalences) with 
a fibrant target.

\nxsubpoint (Cylinder and path objects.)
Cylinder and path objects are defined fiberwise in the classical way
(cf.~\ptref{sp:cyl.path.obj}). For example, a cylinder object 
$A\times I$ of an object $A\in\Ob\cC(S)$ is a diagram
$A\sqcup A\stackrel{\langle\partial_0,\partial_1\rangle}\longto 
A\times I\stackrel\sigma\to A$ with
$\sigma\circ\langle\partial_0,\partial_1\rangle=\nabla_A=
\langle\id_A,\id_A\rangle$,
such that $\langle\partial_0,\partial_1\rangle$ is a cofibration and
$\sigma$ is a weak equivalence. Notice that the global existence of cylinder 
and path objects is a consequence of (MS5).

\nxsubpoint (Local, semilocal and global homotopies.)
Let us fix two parallel morphisms $f,g:A\rightrightarrows X$ in a 
fiber~$\cC(S)$. A {\em left homotopy $h$ from $f$ to $g$} is defined in 
the classical way: it is a morphism $h:A\times I\to X$ from any cylinder 
object $A\times I\in\Ob\cC(S)$ of~$A$ as above, such that 
$f=h\partial_0$ and $g=h\partial_1$. When such a homotopy (globally) 
exists, we say that {\em $f$ is globally left homotopic to~$g$} and 
write $f\Hl\approx g$. We denote by $\pi^\ell_\approx(A,X)$ the 
quotient of $\Hom_{\cC(S)}(A,X)$ with respect to (the equivalence relation
generated by) $\Hl\approx$. Next, we say that {\em $f$ is semilocally left 
homotopic to~$g$}, if the cylinder object $A\times I$ for~$A$ can be 
still chosen globally, i.e.\ in~$\cC(S)$, but the left homotopy 
$h:A\times I\to X$ from~$f$ to~$g$ exists only locally, 
i.e.\ there is a cover $\{S_\alpha\to S\}$
and some left homotopies $h_\alpha:A\times I|_{S_\alpha}\to X|_{S_\alpha}$
between $f|_{S_\alpha}$ and $g|_{S_\alpha}$. In this case we write 
$f\Hl\simeq g$, and denote by $\pi^\ell_\simeq(A,X)$ the corresponding 
quotient of $\Hom_{\cC(S)}(A,X)$.

Finally, we say that {\em $f$ is locally left homotopic to~$g$} and 
write $f\siml g$ if both the cylinder object and the homotopy exist only 
locally, i.e.\ if there is a cover $\{S_\alpha\to S\}$, such that 
$f|_{S_\alpha}\Hl\approx g|_{S_\alpha}$. Notice that this is exactly 
the interpretation of the classical definition in the Kripke--Joyal 
semantics. The corresponding quotient of $\Hom_S(A,X)$ will be denoted 
by $\pi^\ell(A,X)$.

Now the formula $(T\stackrel\phi\to S)\mapsto\pi^\ell(\phi^*A,\phi^*X)$ 
defines a presheaf on $\cS_{/S}$, denoted by $\bpi^\ell(A,X)$. 
We denote its sheafification by $\tilde\bpi^\ell(A,X)$, and the sections 
of this sheafification over~$S$ by $\tilde\pi^\ell(A,X)$. From the
Kripke--Joyal point of view this sheaf $\tilde\bpi^\ell(A,X)$ is the correct
generalization of classical $\pi^\ell(A,X)$.

Of course, all the notions introduced above have their right counterparts.
Clearly, $\Hl\approx\Rightarrow\Hl\simeq\Rightarrow\siml$ and
$\Hr\approx\Rightarrow\Hr\simeq\Rightarrow\simr$.

In the following several lemmas $A$, $B$, $C$, $X$, $Y$ are objects 
of~$\cC(S)$. These lemmas are natural counterparts of Quillen's lemmas 4--8 of
\cite[1.1]{Quillen}, and are shown essentially in the same way 
(in the case of $\siml$ and $\simr$ {\em exactly\/} in the same way,
if we interpret the proofs in Kripke--Joyal semantics).
\begin{LemmaD}
If $A$ is cofibrant, then $\siml$ and $\Hl\simeq$ are equivalence relations 
on $\Hom_{\cC(S)}(A,B)$, $\bpi^\ell(A,B)$ is a separated presheaf, and 
$\bpi^\ell(A,B)\to\tilde\bpi^\ell(A,B)$ is injective.
\end{LemmaD}
\begin{Proof} Same as in \cite[1.1]{Quillen}, Lemma~4. Once we know that
$\siml$ is a {\em local\/} equivalence relation, we see that 
it is representable by a subsheaf $R\subset\iHom(A,B)\times\iHom(A,B)$,
hence the presheaf quotient $\bpi^\ell(A,B)=\iHom(A,B)/R$ is separated. 
The last statement follows from $\tilde\bpi^\ell(A,B)=a\bpi^\ell(A,B)$.
\end{Proof}

\begin{LemmaD}\label{l:left.homot.to.right}
Let $A$ be cofibrant, and $f$, $g\in\Hom_{\cC(S)}(A,B)$. Then:
\begin{itemize}
\item[(i)] $f\Hl\simeq g\Rightarrow f\siml g\Rightarrow f\Hr\simeq g\Rightarrow
f\simr g$.
\item[(ii)] If $f\Hr\simeq g$, there (locally) exists 
a right homotopy $k:A\to B^I$ from $f$ to~$g$ with $s:B\to B^I$ an acyclic 
cofibration in~$\cC(S)$. 
If $f\simr g$, same conclusion with $B^I$ existing only locally.
\item[(iii)] If $u:B\to C$, then $f\simr g\Rightarrow uf\simr ug$, and
$f\Hr\simeq g\Rightarrow uf\Hr\simeq ug$.
\end{itemize}
\end{LemmaD}
\begin{Proof} Identical to \cite[1.1]{Quillen}, Lemma~5. Notice that to prove
the only non-trivial implication $\siml\Rightarrow\Hr\simeq$ in (i) we use
(MS5) to construct a global path object $B^I$ for~$B$.
\end{Proof}

\begin{LemmaD}\label{l:homot.comp}
If $A$ is cofibrant, then composition in~$\cC(S)$ induces maps of
sets $\pi^r(B,C)\times\pi^r(A,B)\to\pi^r(A,C)$, and similarly for 
$\pi^r_\simeq$ and $\tilde\pi^r$, as well as maps of presheaves
$\bpi^r(B,C)\times\bpi^r(A,B)\to\bpi^r(A,C)$ and corresponding sheaves
$\tilde\bpi^r(B,C)\times\tilde\bpi^r(A,B)\to\tilde\bpi^r(A,C)$. All these maps 
are compatible with pullbacks $\phi^*$, for all $\phi:T\to S$.
\end{LemmaD}
\begin{Proof} Same as in \cite[1.1]{Quillen}, Lemma~6, for the first two 
statements, using (iii) of the previous lemma~\ptref{l:left.homot.to.right}. 
The remaining statements follow immediately.
\end{Proof}

\begin{LemmaD}\label{l:afib.to.monom}
Let $A$ be cofibrant and $p:X\to Y$ be an acyclic fibration. Then the maps
of sets $p_*:\pi^\ell(A,X)\to\pi^\ell(A,Y)$, $\pi^\ell_\simeq(A,X)\to
\pi^\ell_\simeq(A,X)$ and of presheaves $\bpi^\ell(A,X)\to\bpi^\ell(A,Y)$
are injective, while the maps of sheaves 
$\tilde\bpi^\ell(A,X)\to\tilde\bpi^\ell(A,Y)$ and of sets 
$\tilde\pi^\ell(A,X)\to\tilde\pi^\ell(A,Y)$ are bijective.
\end{LemmaD}
\begin{Proof}
The proof of injectivity is that of \cite[1.1]{Quillen}, Lemma~7,
while the proof of surjectivity given in {\em loc.cit.}\ actually shows only 
{\em local\/} surjectivity, i.e.\ that $\bpi^\ell(A,p)$ becomes surjective 
after sheafification. The remaining statements are now immediate.
\end{Proof}

\nxsubpoint
Notice that in the above notations (after appropriate identifications) 
we obtain canonical inclusions 
$\pi^\ell(A,X)\subset\pi^\ell(A,Y)\subset\tilde\pi^\ell(A,X)=\tilde\pi^\ell
(A,Y)$. Therefore, if we fix $A$ and $X$, but let $p:X\to Y$ run over 
all acyclic fibrations, we get an embedding of the union of all $\pi^\ell(A,Y)$
into $\tilde\pi^\ell(A,X)$. It is an interesting question whether this union 
equals $\tilde\pi^\ell(A,X)$. We'll see later that when this 
is true, we can compute $\HO_p\cC$ by a right fraction calculus in
$\pi\cC_c$, similar to the constriction of the derived category of 
an abelian category.

\begin{LemmaD}\label{l:homot.into.iso}
(i) Let $F:\cC\to\cB$ be a cartesian functor from a model stack~$\cC$ 
into a prestack (or, more generally, a fibered category with separated
$\iHom$-presheaves)~$\cB$, transforming weak equivalences into isomorphisms.
If $f\siml g$ or $f\simr g$, then $F(f)=F(g)$.

(ii) Let a cartesian $F:\cC_c\to\cB$ carry weak equivalences into 
isomorphisms, with $\cC$ and $\cB$ as above. Then $f\simr g$ implies
$F(f)=F(g)$.

(iii) Let $F:\cC_{cf}\to\cB$ be a cartesian functor as above, then
$f\sim g$ implies $F(f)=F(g)$.
\end{LemmaD}
\begin{Proof} Same as in \cite[1.1]{Quillen}, Lemma~8. Indeed, to show (i)
assume first that $f\Hr\approx g$, i.e.\ we have a right homotopy 
$k:A\to B^I$ and a path object $B\stackrel s\to B^I\stackrel{(d_0,d_1)}\longto
B\times B$ for~$B$, such that $d_0k=f$, $d_1k=g$. Then $F(d_0)F(s)=\id=
F(d_1)F(s)$, hence $F(d_0)=F(d_1)$, $F(s)$ being an isomorphism, hence
$F(f)=F(d_0)F(k)=F(d_1)F(k)=F(g)$. Now if we know just $f\simr g$, then
$f|_{S_\alpha}\Hr\approx g|_{S_\alpha}$ on some cover $\{S_\alpha\to S\}$,
so we obtain $F(f)|_{S_\alpha}=F(g)|_{S_\alpha}$, hence $F(f)=F(g)$, 
the presheaf $\iHom_{\cB_{/S}}(F(A),F(B))$ being separated. 
The statement for $f\siml g$ is shown dually. 

Now to show (ii) and (iii) for some $f,g:A\rightrightarrows B$ 
we reason similarly; the only problem is that we need $B^I$ to lie in
$\cC_c$ (resp.\ $\cC_{cf}$). This is achieved by invoking
\ptref{l:left.homot.to.right}, (ii): we can (locally) choose the path object
in such a way that $s:B\to B^I$ be an acyclic cofibration; then $B^I$ will lie
in~$\cC_c$ whenever $B$ is itself cofibrant. On the other hand, $B^I\to
B\times B$ is always a fibration, hence when $B$ is fibrant, so is $B^I$.
\end{Proof}

\nxsubpoint
We denote by $\pi\cC_c$ the fibered $\cS$-category with the same objects 
as $\cC_c$, but with morphisms given by $\pi^r$, and define
$\pi\cC_f$ dually. According to~\ptref{l:homot.comp}, the composition
in $\pi\cC_c$ is well-defined and compatible with pullbacks, so we've
indeed described a fibered category. We define $\tilde\pi\cC_c$ and
$\tilde\pi\cC_f$ in a similar manner, using $\tilde\pi$ instead of~$\pi$.
Clearly, $\tilde\pi\cC_c$ and $\tilde\pi\cC_f$ are the prestacks 
associated to $\pi\cC_c$ and $\pi\cC_f$ (cf.~\ptref{sp:prestack.from.fib}). 

We construct $\pi\cC_{cf}$ and $\tilde\pi\cC_{cf}$ similarly, taking into 
account that for a cofibrant~$A$ and a fibrant~$B$ the equivalence relations
$\siml$, $\Hl\simeq$, $\simr$, $\Hr\simeq$ coincide on $\Hom(A,B)$
by~\ptref{l:left.homot.to.right}, so we can write 
simply $\pi(A,B)$ and~$\tilde\pi(A,B)$. 
Notice that by \ptref{l:afib.to.monom} and its dual both acyclic cofibrations
and acyclic fibrations become isomorphisms in $\tilde\pi\cC_{cf}$; since 
any weak equivalence in $\cC_{cf}$ can be decomposed  by (MS5) and~(MS3) into 
an acyclic cofibration followed 
by an acyclic fibration, with the intermediate object lying automatically 
in $\cC_{cf}$, we conclude that {\em $\cC_{cf}\to\tilde\pi\cC_{cf}$ 
transforms weak equivalences into isomorphisms.} Therefore, we obtain
a cartesian functor $\bar\gamma_{cf}^+:\HO_p\cC_{cf}\to\tilde\pi\cC_{cf}$,
the target being a prestack.

Finally, \ptref{l:homot.into.iso} implies existence of cartesian functors 
$\pi\cC_c\to\tilde\pi\cC_c\stackrel{\bar\gamma_c}\to\HO_p\cC_c\to\HO\cC_c$,
similar functors $\pi\cC_f\to\cdots\to\HO\cC_f$, as well as cartesian 
functors $\pi\cC_{cf}\to\tilde\pi\cC_{cf}\stackrel{\bar\gamma}\to\HO_p\cC
\to\HO\cC$ and $\bar\gamma_{cf}:\tilde\pi\cC_{cf}\to\HO_p\cC_{cf}$. 
Our next goal is to show that $\bar\gamma:\tilde\pi\cC_{cf}\to
\HO_p\cC$ is an $\cS$-equivalence.

\begin{LemmaD}\label{l:afib.to.isom}
Cartesian functor $\cC_c\to\tilde\pi\cC_c$ transforms acyclic fibrations
into isomorphisms. Dually, $\cC_f\to\tilde\pi\cC_f$ transforms acyclic
cofibrations into isomorphisms.
\end{LemmaD}
\begin{Proof}
Let $p:X\to Y$ be an acyclic fibration of cofibrant objects. Using (MS4)
we obtain {\em local\/} existence of a section $\sigma:Y\to X$, 
$p\sigma=\id_Y$. Now if $\sigma'$ is another (local) section of~$p$
(over~$S$ or over some $S_\alpha\to S$), then $p\sigma=p\sigma'$ 
implies $\sigma\siml\sigma'$ by \ptref{l:afib.to.monom}, hence
$\sigma\simr\sigma'$ by \ptref{l:left.homot.to.right},(i), so
the image $\bar\sigma$ of $\sigma$ in $\tilde\pi^r(Y,X)$ exists locally
and is unique, so it exists globally as well. In other words,
$\bar\sigma:Y\to X$ is a section of $\bar p$ in $\pi\cC_c$. On the other hand,
consider (locally) $\sigma p$ and $\id_X$; they have equal images under
$p_*$, hence by the same lemmas we obtain $\sigma p\siml\id_X$ and
$\sigma p\simr\id_X$, i.e.\ $\bar\sigma\bar p=\id_X$ in~$\tilde\pi\cC_c$
(locally, hence also globally, $\tilde\pi\cC_c$ being a prestack),
so $\bar\sigma$ is a two-sided inverse to~$\bar p$.
\end{Proof}

\begin{ThD}\label{th:homot.mod.stk} {\rm (cf.~\cite[1.1]{Quillen}, Theorem~1)}
The cartesian functor $\bar\gamma:\tilde\pi\cC_{cf}\to\HO_p\cC$ is 
an $\cS$-equivalence, hence $\HO_p\cC$ (resp.\ $\HO\cC$) is 
$\cS$-equivalent to the prestack (resp.\ stack) associated to $\cS$-fibered 
category $\pi\cC_{cf}$.

In particular, $\HO_p\cC$ and $\HO\cC$ are $\univU$-categories.
\end{ThD}
\begin{Proof}
(a) First of all, notice that $\bar\gamma_{cf}$ and $\bar\gamma_{cf}^+$ are
adjoint $\cS$-equivalences between prestacks $\tilde\pi\cC_{cf}$ and
$\HO_p\cC_{cf}$, and even $\cS$-isomorphisms, inverse to each other, if 
we choose the {\em strict\/} version of $\HO_p\cC_{cf}$. Indeed, all we have
to check is that $\cC_{cf}\to\tilde\pi\cC_{cf}$ transforms weak equivalences 
into isomorphisms, something that we know already, and that this cartesian
functor is universal among all cartesian functors $F:\cC_{cf}\to\cB$ into
$\cS$-prestacks $\cB$ that transform weak equivalences into isomorphisms,
something we know from \ptref{l:homot.into.iso},(iii), $\tilde\pi\cC_{cf}$ 
being the prestack associated to $\pi\cC_{cf}$.

(b) Now let us construct a cartesian functor $\bar Q:\cC\to\tilde\pi\cC_c$.
Let's choose for this arbitrary {\em strict\/} cofibrant replacements 
$Q(X)\to X$ for all $X\in\Ob\cC$, strictness understood here as 
the requirement of $Q(X)\to X$ to be an acyclic fibration, 
not just a weak equivalence,
and put $\bar Q(X)$ equal to the image of $Q(X)$ in $\tilde\pi\cC_c$. Notice 
that {\em $\bar Q(X)$ is independent (up to a canonical isomorphism)
on the choice of $Q(X)\to X$.} Indeed, whenever we have two such strict
cofibrant replacements $p:Q\to X$, $p':Q'\to X$, 
we can find another strict replacement $p'':Q''\to X$ and 
acyclic fibrations $\sigma:Q''\to Q$ and $\sigma':Q''\to Q'$, such that 
$p\sigma=p''=p'\sigma'$. To achieve this we simply take a strict
cofibrant replacement $Q''\to Q\times_XQ'$, using the fact that acyclic 
fibrations are stable under base change and composition. 
Now acyclic fibrations $\sigma$ and $\sigma'$ become isomorphisms in
$\tilde\pi\cC_c$ by \ptref{l:afib.to.isom}, hence an isomorphism
$\bar\sigma'\circ\bar\sigma^{-1}:\bar Q\stackrel\sim\leftarrow\bar 
Q''\simto\bar Q'$. 
The independence of this isomorphism on the choice of $Q''\to Q\times_XQ'$ 
is checked similarly.

(c) So far we have constructed $\bar Q$ only on objects. Notice, however,
that the independence on the choice of $Q(X)\to X$ implies the compatibility
with pullback functors $\phi^*$, since $Q(\phi^*X)\to \phi^*X$ and 
$\phi^*Q(X)\to\phi^*X$ are two strict cofibrant replacements of~$\phi^*X$. 
Thus it remains to define $\bar Q$ on morphisms in fibers of~$\cC$.
If $f:X\to Y$ is a morphism in $\cC(S)$, and $p_X:Q(X)\to X$, 
$p_Y:Q(Y)\to Y$ are two strict cofibrant replacements, we can apply 
(MS4) to lift {\em locally\/} $fp_X:Q(X)\to Y$ to a morphism 
$Q(f):Q(X)\to Q(Y)$, which is unique up to a left homotopy by 
\ptref{l:afib.to.monom}, hence {\em a fortiori\/} unique up to a right
homotopy by \ptref{l:left.homot.to.right},(i), 
hence we obtain locally a well-defined morphism
$\bar Q(f):\bar Q(X)\to \bar Q(Y)$ in $\tilde\pi\cC_c$, hence
$\bar Q(f)$ is defined globally as well, $\tilde\pi\cC_c$ being a prestack. 
For any $g:Y\to Z$ we get $Q(g)Q(f)\siml Q(gf)$ and $Q(\id_X)\siml\id_{Q(X)}$
by \ptref{l:afib.to.monom}, hence $\bar Q(g)\bar Q(f)=\bar Q(gf)$ and
$\bar Q(\id_X)=\id_{\bar Q(X)}$ by~\ptref{l:left.homot.to.right},(i) again.
This finishes the construction of cartesian functor 
$\bar Q:\cC\to\tilde\pi\cC_c$.

(d) We construct a cartesian functor $\bar R:\cC\to\tilde\pi\cC_f$ similarly,
by putting $\bar R(X)$ to be the image of any strict fibrant replacement 
$i_X:X\to R(X)$. If $X$ is cofibrant, so will be $R(X)$, $i_X$ being an 
acyclic cofibration. Hence $\bar R$ induces a functor $\bar R_c:
\cC_c\to\tilde\pi\cC_{cf}$. Notice that $f\simr g:X\to Y$ in $\cC_c$ implies
$i_Yf\simr i_Yg$ by~\ptref{l:left.homot.to.right},(iii), 
hence (locally) $R(f)i_X\simr R(g)i_X$, i.e.\ the images of these two 
elements coincide in $\tilde\pi^r(X,R(Y))$; the dual of
\ptref{l:afib.to.monom} now implies $\bar R(f)=\bar R(g)$ (locally,
hence also globally) in $\tilde\pi^r(R(X),R(Y))=\tilde\pi(R(X),R(Y))$, i.e.\ 
$\bar R_c(f)=\bar R_c(g)$. This means that $\bar R_c$ factorizes through
$\pi\cC_c$, hence also through $\tilde\pi\cC_c$, the target category 
being a prestack. The induced functor will be denoted 
$\tilde R_c:\tilde\pi\cC_c\to\tilde\pi\cC_{cf}$.

(e) Now consider the composite cartesian functor 
$\gamma':\cC\stackrel{\bar Q}\to
\tilde\pi\cC_c\stackrel{\tilde R_c}\to\tilde\pi\cC_{cf}$. Its target category
is a prestack, and it is immediate from (MS3) that the image under $\gamma'$ 
of a weak equivalence $f:X\to Y$ is locally representable by a
weak equivalence $RQ(f):RQ(X)\to RQ(Y)$, hence $\gamma'(f)$ is 
locally an isomorphism, hence this is true globally. Therefore,
$\gamma'$ induces a cartesian functor 
$\tilde\gamma:\HO_p\cC\to\tilde\pi\cC_{cf}$.

(f) Consider cartesian functors between prestacks $\HO_p\cC_{cf}
\stackrel I\to\HO_p\cC\stackrel{\tilde\gamma}\to\tilde\pi\cC_{cf}
\stackrel{\bar\gamma_{cf}}\to\HO_p\cC_{cf}$. Clearly, the composite 
functor is $\cS$-isomorphic to the identity, since we can choose $Q(X)=X=R(X)$ 
for any $X\in\Ob\cC_{cf}$ in the construction of $\gamma'=\tilde R_c\bar Q$.
The last arrow $\bar\gamma_{cf}$ is an $\cS$-equivalence by (a), hence 
$\HO_p\cC_{cf}\stackrel I\to\HO_p\cC\stackrel{\tilde\gamma}\to
\tilde\pi\cC_{cf}$ is
an $\cS$-equivalence as well; in particular, $\tilde\gamma$ is 
surjective on morphisms. On the other hand, it is faithful, 
the composite functor $\HO_p\cC\stackrel{\tilde\gamma}\to\tilde\pi\cC_{cf}
\cong\HO_p\cC_{cf}\stackrel I\to\HO_p\cC$ 
being isomorphic to the identity functor as well: indeed, it transforms 
$\gamma X$ into $\gamma RQ(X)$, which is canonically isomorphic to $\gamma X$
via $\gamma(i_{Q(X)})\gamma(p_X)^{-1}$. This proves that $I$ is 
essentially surjective, and that $\tilde\gamma$ 
is fully faithful; since the same is true for $\bar\gamma_{cf}$ and
$\bar\gamma_{cf}\tilde\gamma I\cong\Id$, we see that $I$ is 
fully faithful as well, hence an $\cS$-equivalence. From this and (a)
we deduce immediately that $\tilde\gamma:\HO_p\cC\to\tilde\pi\cC_{cf}$
is an $\cS$-equivalence as well, with $\bar\gamma=I\bar\gamma_{cf}$
its quasi-inverse, q.e.d.
\end{Proof}

\begin{CorD} {\rm (cf.~cor.~1 of th.~1, \cite[1.1]{Quillen}) }
If $A$ is cofibrant and $B$ fibrant in $\cC(S)$, then
$[A,B]_S:=\Hom_{(\HO\cC)(S)}(\gamma A,\gamma B)$ is canonically
isomorphic to $\tilde\pi(A,B)$.
\end{CorD}

\nxsubpoint
Notice that if we would have (SM5l) instead of (SM5), the whole construction
of $\bar R$ and $\bar Q$ would work only on the level of associated stacks,
fibrant and cofibrant replacements being defined only locally.

It is also interesting to note that all the above constructions and 
statements are still valid when $\cC$ is a model {\em pre}stack over~$\cS$.

\begin{PropD}
(a) A fibration $p:X\to Y$ in $\cC_{cf}(S)$ is acyclic iff 
$\gamma(p)$ is an isomorphism in $\HO_p\cC$ iff locally $p$ is a dual of 
a strong deformation retract, i.e.\ iff after pulling back to some cover
$\{S_\alpha\to S\}$ one can find a section $t:Y\to X$, $pt=\id_X$, 
and a left homotopy $h:X\times I\to X$ from $tp$ to $\id_Y$ with $ph=p\sigma$.

(b) A morphism $f:X\to Y$ in $\cC(S)$ is a weak equivalence iff
$\gamma(f)$ is an isomorphism in $\HO\cC$ (or in its full subcategory
$\HO_p\cC$).
\end{PropD}
{\bf Proof.} 
Same as that of Lemma~1 and Proposition~1 of \cite[1.5]{Quillen}.

\nxsubpoint (Loop and suspension functors.)
When $\cC$ is a {\em pointed\/} model stack over~$\cS$, i.e.\ each fiber
$\cC(S)$ admits a zero object, we can define loop and suspension objects 
in the usual way. For example, the suspension $\Sigma A$ of a cofibrant 
object $A\in\Ob\cC_c(S)$ can be defined as the cofiber of $A\sqcup A\to
A\times I$, for any cylinder object $A\times I$ for~$A$. One checks 
in the usual manner that $\gamma\Sigma A$ depends (up to a canonical 
isomorphism in $(\HO_p\cC)(S)$) only on the weak equivalence class of~$A$,
so we get well-defined cartesian functors $\Sigma=\dL\Sigma$
and $\Omega=\dR\Omega:\HO_p\cC\to\HO_p\cC$, and corresponding functors
$\HO\cC\to\HO\cC$, denoted in the same way. When we have (SM5f), we can
construct cartesian functors $\Omega:\cC_f\to\cC_f$ and $\Sigma:\cC_c\to\cC_c$
as well.

\nxsubpoint (Fibration and cofibration sequences.)
One can define fibration and cofibration sequences in $\HO_p\cC$ for any
pointed model category~$\cC$
essentially in the same way as in \cite[1.3]{Quillen}, and extend these 
definitions to $\HO\cC$ ``by descent'' (cf.~\ptref{sp:stack.from.fib}). 
Then one can construct associated long exact sequences of 
sheaves of abelian groups/groups/sets on $\cS_{/S}$, 
again by repeating the reasoning of {\em loc.cit.}\ for
the homotopic prestack $\HO_p\cC$, and extending to the associated stack
$\HO\cC$ ``by descent''.

\nxsubpoint (Criteria for existence of derived functors.)
Of course, we have the usual criteria for existence of derived functors.
Suppose for example that $F:\cC\to\cD$ is a cartesian functor between
model stacks over~$\cS$, transforming acyclic cofibrations between fibrant
objects (in fibers of~$\cC$) into weak equivalences in~$\cD$. Then its 
left derived functors $\dL F=\dL_pF:\HO_p\cC\to\HO_p\cD$ as well as
$\dL F:\HO\cC\to\HO\cD$ exist, and the natural morphism
$\eta_X:\dL F(\gamma_\cC X)\to\gamma_\cD F(X)$ is an isomorphism for any
cofibrant~$X$, i.e.\ $\dL F(X)$ can be computed as $F(P)$ for any
cofibrant replacement $P\to X$. The proof is essentially that of
\ptref{prop:ex.left.der} and \ptref{l:acof.cof.suff} (cf.\ also 
prop.~1 of \cite[1.4]{Quillen}), at least for
$\dL_p F:\HO_p\cC\to\HO_p\cD$. Then we construct $\dL F:\HO\cC\to\HO\cD$
by extending $\dL_p F$ to associated stacks. In the sequel we'll
refer to these criteria as ``model stack'' or ``local'' versions
of~\ptref{prop:ex.left.der} and \ptref{l:acof.cof.suff}. Notice that
the possibility to construct $\dL_p F$ on the level of homotopic prestacks 
is due to the existence of global factorizations (MS5): its weaker version
(MS5l) would enable us to construct $\dL F$, but not $\dL_p F$.

\nxsubpoint (Adjoint cartesian Quillen functors.)
We have a notion of {\em cartesian Quillen functors\/} or {\em pairs\/}:
these are $\cS$-adjoint cartesian functors $F:\cC\to\cD$, $G:\cD\to\cC$ 
between model stacks over~$\cS$, such that $F$ transforms weak equivalences
(or just acyclic cofibrations) in fibers of~$\cC_c$ into weak equivalences
in~$\cD$, and dually $G$ preserves weak equivalences in $\cD_f$.
In this case derived functors $\dL_p F:\HO_p\cC\to\HO_p\cD$ and
$\dR_p G:\HO_p\cD\to\HO_p\cC$ exist, can be computed by means of 
cofibrant (resp.\ fibrant) replacements, and are $\cS$-adjoint. Same 
applies to their stack extensions $\dL F:\HO\cC\to\HO\cD$ and
$\dR G:\HO\cD\to\HO\cC$.

\nxsubpoint (Fibered products of model stacks.)
Whenever $\cC_1$ and $\cC_2$ are fibered categories over~$\cS$, 
their fibered product $\cC_1\times_\cS\cC_2$ is another one; 
this is clearly the {\em strict\/} product of $\cC_1$ and $\cC_2$ in
the 2-category of fibered $\cS$-categories, i.e.\ for any fibered
$\cS$-category $\cD$ we have an {\em isomorphism\/} of categories
$\catCart_\cS(\cD,\cC_1\times_\cS\cC_2)\cong\catCart_\cS(\cD,\cC_1)\times
\catCart_\cS(\cD,\cC_2)$. All this means that this fibered product is 
a natural replacement of product of categories; in fact, it is the 
``fiberwise product of two families of categories parametrized by~$\cS$''
since $(\cC_1\times_\cS\cC_2)(S)=\cC_1(S)\times\cC_2(S)$, and 
the corresponding pullback functors $\phi^*$ can be constructed 
as $\phi_{\cC_1}^*\times\phi_{\cC_2}^*$. Therefore, a natural replacement 
for bifunctors is given by cartesian functors like 
$F:\cC_1\times_\cS\cC_2\to\cD$; if we need, say, a bifunctor contravariant 
in the first argument and covariant in the second, we use the 
fiberwise opposite: $G:\cC_1^{fop}\times_\cS\cC_2\to\cD$.

Now similarly to what we had in~\ptref{sp:prod.mod.cat},
{\em $\cC_1\times_\cS\cC_2$ admits a natural model stack structure 
for any two model stacks $\cC_1$ and $\cC_2$ over a site~$\cS$.} 
Of course, we declare a morphism $f=(f_1,f_2)\in\Ar(\cC_1\times_\cS\cC_2)(S)=
\Ar(\cC_1(S)\times\cC_2(S))$ to be a (co)fibration or a weak equivalence 
iff both $f_1$ and $f_2$ are. This enables us to derive functors 
$F:\cC_1\times_\cS\cC_2\to\cD$ and $G:\cC_1^{fop}\times_\cS\cC_2\to\cD$ in 
the obvious manner. We obtain criteria of existence of derived ``bifunctors''
similar to~\ptref{sp:ex.lder.bif} as well, where the corresponding requirements
(that of $F(P,-)$ to transform acyclic cofibrations between cofibrant objects
of~$\cC_2$ into weak equivalences in~$\cD$ for any cofibrant object $P$ 
of~$\cC_1$, and the symmetric one) have to be understood fiberwise
(i.e.\ $F(P,f)$ has to be a weak equivalence in $\cD(S)$ 
for any $P\in\Ob\cC_{1,c}(S)$ and any acyclic cofibration $f:X\to Y$ 
in $\cC_{2,c}$, together with the symmetric condition). In this case 
$\dL_p F(\gamma X_1,\gamma X_2)\cong \gamma F(P_1,P_2)$, where 
$P_i\to X_i$ are arbitrary cofibrant replacements in $\cC_i(S)$.

\nxsubpoint\label{sp:cart.tens.struct} (Cartesian $\otimes$-structures.)
One can define the notion of a $\otimes$-structure and of constraints for 
such a structure for objects of an arbitrary strictly associative
2-category with strict products, so as to recover the usual definitions
in the case of the 2-category of all categories. On the other hand,
applying these definitions to the 2-category of fibered categories and 
cartesian functors over a fixed category~$\cS$, we obtain the notion 
of a {\em cartesian $\otimes$-structure\/} on a fibered category $\cC/\cS$.
By definition, this is just a cartesian functor 
$\otimes:\cC\times_\cS\cC\to\cC$; and, for example, an associativity 
constraint is an $\cS$-isomorphism of appropriate cartesian functors 
$\cC\times_\cS\cC\times_\cS\cC\to\cC$. A {\em fibered $\otimes$-category
over~$\cS$} is by definition a fibered category $\cC\to\cS$ with a
cartesian $\otimes$-structure; we speak about AU, ACU etc.\ fibered
$\otimes$-categories depending on the presence and compatibility of 
appropriate constraints.

If we think of $\cC\to\cS$ as ``a family of categories parametrized by~$\cS$'',
then, say, a cartesian ACU $\otimes$-structure on $\cC$ is essentially the
same thing as an ACU $\otimes$-structure on each fiber $\cC(S)$, together 
with a $\otimes$-functor structure on each pullback functor
$\phi^*:\cC(S)\to\cC(T)$ (i.e.\ some functorial isomorphisms
$\phi^*(X\otimes Y)\cong\phi^*X\otimes\phi^*Y$, compatible with constraints),
the canonical isomorphisms $(\phi\psi)^*\cong\psi^*\phi^*$ being required
to be isomorphisms of $\otimes$-functors.

Sometimes we denote by $X\otimes_S Y$ the image of $(X,Y)$ under $\otimes$,
when both $X$ and $Y$ lie over $S\in\Ob\cS$.

\nxsubpoint (Compatible cartesian $\otimes$-structures on a model stack.)
Now let $\otimes:\cC\times_\cS\cC\to\cC$ be a cartesian $\otimes$-structure
on a model stack $\cC$ over a site~$\cS$. We say that this cartesian 
$\otimes$-structure is {\em compatible\/} with the model structure of~$\cC$
if the condition (TM) of \ptref{def:comp.tens.mod.str} holds in each fiber
$\cC(S)$ of~$\cC$.

If a compatible cartesian $\otimes$-structure (fiberwise) preserves 
initial objects in each argument, it admits a left derived 
$\Lotimes:\HO\cC\times_\cS\HO\cC\to\HO\cC$ or 
$\Lotimes:\HO_p\cC\times_\cS\HO_p\cC\to\HO_p\cC$, that can be computed 
by means of cofibrant replacements in each argument, by the same 
reasoning as in~\ptref{sp:der.comp.tensprod}. We denote also by $\Lotimes$
the induced functors between corresponding categories of global (cartesian) 
sections $\Ho\cC$ and $\Ho_p\cC$, thus obtaining $\otimes$-structures 
on each of these categories, with the same constraints as the original one.

\nxsubpoint (Compatible external cartesian $\otimes$-actions.)
Of course, we can define a {\em right external cartesian $\otimes$-action\/} 
of a fibered AU $\otimes$-category $\cC$ on another fibered category 
$\cD$ over the same base category~$\cS$ as a cartesian functor 
$\oslash:\cD\times_\cS\cC\to\cD$, together with external associativity and 
unity constraints compatible with those of~$\cC$. A left external cartesian
$\otimes$-action $\obslash:\cC\times_\cS\cD\to\cD$ is defined similarly.
Again, this is essentially equivalent to giving a family of external
$\otimes$-actions $\oslash_S:\cD(S)\times\cC(S)\to\cC(S)$, compatible with 
the pullback functors.

When in addition $\cC$ and $\cD$ are model stacks over a site~$\cS$, we
can define a {\em compatible\/} (say, right) cartesian $\otimes$-action
$\oslash:\cD\times_\cS\cC\to\cD$, by requiring 
(TMe) of~\ptref{def:comp.tens.mod.str} to hold fiberwise. 
Of course, if $\oslash$ (fiberwise) preserves initial
objects in each argument, it can be derived with the aid of cofibrant 
replacements in each position, thus defining an external cartesian 
$\otimes$-action of $\HO_p\cC$ on $\HO_p\cD$ and of $\HO\cC$ on $\HO\cD$,
as well as external $\otimes$-actions between corresponding global section
categories.

\nxsubpoint (Equivalent conditions in terms of inner Homs.)
Compatibility conditions (TM) and (TMe) have equivalent formulations 
similar to (TMh) and (SM7a) in terms of appropriate inner Homs, 
whenever these exist and are cartesian. More precisely, let us
consider the case of a cartesian external $\otimes$-action 
$\oslash:\cD\times_\cS\cC\to\cD$, with $\cC$ and $\cD$ model stacks over
a site~$\cS$. We say that a {\em cartesian\/} functor 
$\iHom_\cC:\cD^{fop}\times_\cS\cD\to\cC$ is an inner Hom for $\oslash$
if we have functorial isomorphisms
\begin{multline}
\Hom_{\cD(S)}(X\oslash_S K,Y)\cong\Hom_{\cC(S)}(K,\iHom_\cC(X,Y))\\
\quad\text{for all $S\in\Ob\cS$, $X$, $Y\in\Ob\cD(S)$, $K\in\Ob\cC(S)$}
\end{multline}
compatible with all pullback functors, i.e.\ an $\cS$-isomorphism of
cartesian functors $\cD^{fop}\times_\cS\cC^{fop}\times_\cS\cD\to
\catSets\times\cS$. These functorial isomorphisms can be rewritten in terms
of corresponding local Hom-sheaves, thus yielding an $\cS$-isomorphism
of cartesian functors $\cD^{fop}\times_\cS\cC^{fop}\times_\cS\cD\to
\stSETS_\cS$:
\begin{multline}
\iHom_{\cD_{/S}}(X\oslash_S K,Y)\cong\iHom_{\cC_{/S}}(K,\iHom_\cC(X,Y))\\
\quad\text{for all $S\in\Ob\cS$, $X$, $Y\in\Ob\cD(S)$, $K\in\Ob\cC(S)$}
\end{multline}

Now $\oslash$ is compatible with model structures involved, i.e.\ 
$\oslash$ satisfies (TMe) fiberwise, iff $\iHom_\cC$ satisfies (TMh) fiberwise.
The proof is still the same classical interplay of adjointness and lifting 
properties as before in~\ptref{sp:comp.of.iHoms} or in the classical proof
of (SM7)$\Leftrightarrow$(SM7b). Of course, now we have only local lifting 
properties, so we have to reason locally; Kripke--Joyal semantics provides 
a formal way to transfer this proof to the local case as usual.

On the other hand, we might also have another flavor of inner Hom, namely, 
a cartesian functor $\iHom_\cD:\cC^{fop}\times_\cS\cD\to\cD$, characterized by
\begin{multline}
\iHom_{\cD_{/S}}(X\oslash_S K,Y)\cong\iHom_{\cD_{/S}}(X,\iHom_\cD(K,Y))\\
\quad\text{for all $S\in\Ob\cS$, $X$, $Y\in\Ob\cD(S)$, $K\in\Ob\cC(S)$}
\end{multline}
We usually denote $\iHom_\cD(K,X)$ by $X^K$; when such inner Homs exist,
(fiberwise) condition (TMe) for $\oslash$ is equivalent to (fiberwise) 
condition (SM7a) for~$\iHom_\cD$.

\nxsubpoint (Deriving inner Homs.)
Suppose that $\oslash$ is a compatible external cartesian $\otimes$-action
as above, (fiberwise) preserving initial objects in each variable, and that
one or both of the above inner Homs exist.
Then it is easy to see that inner Homs $\iHom_\cC$ and
$\iHom_\cD$ (fiberwise) transform initial objects in first variable 
and final objects in second variable into final objects, hence the 
same criterion used to derive $\oslash$ is applicable to these inner Homs,
once we replace appropriate stacks by their fiberwise duals and use
(TMh) or (SM7a) instead of (TMe). In this way we obtain right derived functors
$\dR\iHom_\cC:\HO\cD^{fop}\times_\cS\HO\cD\to\HO\cC$ and
$\dR\iHom_\cD:\HO\cC^{fop}\times_\cS\HO\cD\to\HO\cD$, and their 
homotopic prestack counterparts, that can be computed 
by taking cofibrant replacements of the first argument and 
fibrant replacements of the second one. Now observe that for any cofibrant
$Q\in\Ob\cD(S)$ cartesian functors $K\mapsto \phi^*Q\oslash K$ and
$X\mapsto\iHom_\cC(\phi^*Q,X)$, where $T\stackrel\phi\to S$, $K\in\Ob\cC(T)$,
$X\in\Ob\cD(T)$, constitute a Quillen $\cS_{/S}$-adjoint pair of 
cartesian functors $\cC_{/S}\leftrightarrows\cD_{/S}$, 
hence their derived are $\cS_{/S}$-adjoint as well. 
This implies that {\em $\dR\iHom_\cC$ and $\dR\iHom_\cD$ are
inner Homs for $\Loslash:\HO\cD\times_\cS\HO\cC\to\HO\cD$.} Of course, 
we obtain a similar statement on the level of homotopic prestacks,
and for corresponding categories of global (cartesian) sections $\Ho_p\cD$,
$\Ho\cD$ etc.\ as well.

\nxsubpoint (Possible applications.)
The above results would immediately enable us to derive tensor products and
inner Homs on a generalized (commutatively) ringed topos $(\cE,\sO)$, 
and base change with respect to any morphism of generalized ringed topoi, 
provided we manage to construct appropriate model structures on 
stacks $\stks\stSETS_\cE$ and $\stks\stMOD\sO$ over~$\cE$.

\nxpointtoc{Pseudomodel stack structure on simplicial sheaves}
Throughout this subsection we fix a topos~$\cE$ and consider the stack
$\stSETS=\stSETS_\cE\to\cE$ of ``sets'' or ``sheaves'' over~$\cE$,
characterized by $\stSETS_\cE=\Ar\cE$, $\stSETS_\cE(S)=\cE_{/S}$,
and the stack of simplicial sets or sheaves 
$\stks\stSETS_\cE\cong\stCART_\cE(\catDelta^0\times\cE,\stSETS_\cE)$, 
characterized by $(\stks\stSETS_\cE)(S)=s(\cE_{/S})$. We would like to
define a reasonable model stack structure on $\stks\stSETS_\cE$. 
In particular, this would enable us to speak about (co)fibrations and
weak equivalences in its fiber $(\stks\stSETS_\cE)(e)=s\cE$. When
$\cE=\catSets$, we must recover the classical model category structure on
$s\catSets$; but in general we won't obtain a model category structure
on $s\cE$, and not even a (local) model stack structure, 
lifting properties being fulfilled only locally and under some 
additional restrictions. 

Unfortunately, classical construction of a model category structure on
$s\catSets$ is not intuitionistic, so we cannot transfer it immediately
to the topos case, i.e.\ to stack $\stks\stSETS_\cE$ over~$\cE$. 
It turns out that we have to weaken either the (local) lifting axiom or
the factorization axiom of a model stack. We chose to sacrifice part
of the (local) lifting axiom, thus obtaining the notion of a 
{\em pseudomodel stack}. Then we can obtain a {\em pseudo}model structure on
$\stks\stSETS_\cE$.

Some of our statements will be proved only for a topos~$\cE$ with sufficiently
many points; while this additional condition shouldn't affect the validity 
of our statements, it considerably shortens the proofs. The general case 
might be treated either by transferring results from $\hat\cS$ to $\tilde\cS$
in the usual manner, or by using suitable local limits instead of points.
This will be done elsewhere.

Notice that we usually don't consider the site case here; 
of course, we don't lose much, 
since we can always recover corresponding constructions 
over a site~$\cS$ by doing everything over the topos $\tilde\cS$ first 
and then pulling back resulting (pre)stacks with respect to $\cS\to\tilde\cS$
(cf.\ \ptref{sp:stk.over.site.topos}, \ptref{sp:modstk.sites.topoi} 
and~\ptref{sp:homstk.sites.topoi}).

\nxsubpoint (Simplicial objects in a fibered category.)
Given a fibered category~$p_\cC:\cC\to\cS$, we define
the corresponding {\em fibered category of simplicial objects\/}
$\stks\cC$ by means of the following cartesian square of categories
\begin{equation}
\xymatrix{
\stks\cC\ar[r]\ar[d]^{p_{\stks\cC}}&s\cC\ar[d]^{s(p_\cC)}\\
\cS\ar[r]^{I}&s\cS}
\end{equation}
Here $s\cC=\catFunct(\catDelta^0,\cC)$ and $s\cS$ are the usual categories 
of simplicial objects, and $I:\cS\to s\cS$ is the constant simplicial object
embedding. By definition $(\stks\cC)(S)=s(\cC(S))$ for any $S\in\Ob\cS$,
i.e.\ the fibers of $\stks\cC$ consist of simplicial objects over the 
corresponding fibers of~$\cC$, and the pullback functors $\phi^*$ of
$\stks\cC$ are just the simplicial extensions $s(\phi^*)$.

Now if $\cC\to\cS$ is a (pre)stack over a site~$\cS$, then $\stks\cC$ is
obviously another one. In particular, for any topos~$\cE$ we obtain 
a well-defined stack $\stks\stSETS_\cE$ over~$\cE$, characterized by
$(\stks\stSETS_\cE)(X)=s(\cE_{/X})$.

One can extend cartesian functors $F:\cC\to\cD$ over~$\cS$ to cartesian 
functors $\stks F:\stks\cC\to\stks\cD$ in the usual manner; when no 
confusion can arise, we denote this simplicial extension by the same 
letter~$F$. Moreover, this construction is compatible with fibered 
products: $\stks(\cC_1\times_\cS\cC_2)=\stks\cC_1\times_\cS\stks\cC_2$,
so for example a ``cartesian bifunctor'' $F:\cC_1\times_\cS\cC_2\to\cD$ 
admits a simplicial extension $F=\stks F:\stks\cC_1\times_\cS\stks\cC_2\to
\stks\cD$.

Notice that in these constructions $\catDelta^0$ can be replaced by 
an arbitrary small index category~$\cI$. For example, $\cI=\catDelta$ 
yields the stack of cosimplicial objects~$\stkc\cC$.

\nxsubpoint
(Constant simplicial objects.)
Let us denote by $q:\cE\to\catSets$ the canonical morphism of $\cE$ into 
the point topos $\catSets$. Then $q_*X=\Gamma(X)=\Hom_\cE(e_\cE,X)$, and
$q^*K$ is just the constant sheaf or object of~$\cE$, usually denoted
by $\cst{K}_\cE$. Moreover, $q^*$ extends to a functor 
$q^*=s(q^*):s\catSets\to s\cE$; the image under this functor of a simplicial
set $K$ will be also denoted by $\cst{K}_\cE$, or even by~$K$ when 
no confusion can arise.

In particular, we obtain the standard simplices $\cst{\Delta(n)}_\cE$, for 
all $n\geq0$; they still have their characteristic properties 
$\Hom_{s\cE}(\cst{\Delta(n)}_\cE,X)=\Gamma(X_n)$ for any $X\in\Ob s\cE$, and
$\iHom_{\stks\cE|S}(\cst{\Delta(n)}_\cE,X)=X_n$ for any 
$X\in\Ob(\stks\stSETS_\cE)(S)=\Ob s(\cE_{/S})$.

\nxsubpoint 
(Cartesian ACU $\otimes$-structure on $\stSETS$.)
Clearly, the fibered products $\times_S:\cE_{/S}\times\cE_{/S}\to\cE_{/S}$ 
are compatible with arbitrary pullbacks, thus defining a cartesian
``bifunctor'' $\otimes=\times_\cE:\stSETS_\cE\times_\cE\stSETS_\cE\to
\stSETS_\cE$, easily seen to be an ACU cartesian $\otimes$-structure 
on $\stSETS/\cE$ (cf.~\ptref{sp:cart.tens.struct}). Of course, this
cartesian $\otimes$-structure admits (cartesian) inner Homs
$\iHom=\iHom_\cE:\stSETS_\cE^{fop}\times_\cE\stSETS_\cE\to\stSETS_\cE$, 
given by usual local Hom-objects: $\iHom_{\cE|S}(X,Y)=\iHom_S(X,Y)=
\iHom_{\cE_{/S}}(X,Y)$ is just the object of $\cE_{/S}$ representing 
the Hom-sheaf $\iHom_{\cE_{/S}}(X,Y)$ over~$\cE_{/S}$, for any
$X$, $Y\in\Ob\cE_{/S}$.

\nxsubpoint\label{sp:act.SETS.flat.stk}
(Cartesian $\otimes$-action of $\stSETS$ on any flat stack~$\cC$.)
Let $\cC$ be any {\em flat\/} stack over~$\cE$, i.e.\ we require (MS1) and
(MS1+). Then $(X,T)\mapsto X\times_ST=\phi_!\phi_*X$, where
$(T\stackrel\phi\to S)\in\Ob\cE_{/S}=\Ob\stSETS_\cE(S)$, $X\in\Ob\cC(S)$, 
defines a cartesian ``bifunctor'' $\oslash:\cC\times_\cE\stSETS_\cE\to
\stSETS_\cE$, easily seen to be a cartesian $\otimes$-action of $\stSETS$
on~$\cC$ (cf.\ \ptref{sp:equiv.flat.stk} and~\ptref{sp:loc.sum.prod}).
Moreover, both flavors of inner Homs exist in this situation:
$\iHom_\stSETS:\cC^{fop}\times_\cE\cC\to\stSETS_\cE$ is given simply 
by the local $\iHom$-objects, i.e.\ $\iHom_{\stSETS|S}(X,Y)$ represents 
the Hom-sheaf $\iHom_{\cC_{/S}}(X,Y)$, while $\iHom_\cC:\stSETS_\cE^{fop}
\times_\cE\cC\to\cC$ is given by $(T,X)\mapsto X^{T/S}=\phi_*\phi^*X$,
for any $(T\stackrel\phi\to S)\in\Ob\cE_{/S}$, $X\in\Ob\cC(S)$
(cf.~\ptref{sp:loc.sum.prod}).

In the sequel the $\stSETS$-valued inner Homs will be called {\em local},
just to distinguish them from all other sorts of inner Homs that will 
appear. We'll usually denote local Homs simply by $\iHom$ or $\iHom_S$,
or by $\iHom_{\cC|S}$ when we want to indicate the stack $\cC$ as well.

\nxsubpoint
(Cartesian ACU $\otimes$-structure on $\stks\stSETS$.)
Of course, the above cartesian ACU $\otimes$-structure
$\otimes=\times:\stSETS\times_\cE\stSETS\to\stSETS$ extends to a cartesian ACU
$\otimes$-structure $\otimes:=s(\otimes)$ on $\stks\stSETS$. Whenever 
we have two simplicial objects $X$, $Y\in\Ob s(\cE_{/S})$, their 
``tensor product'' $X\otimes_S Y$ is simply the componentwise fibered product
$X\times_SY$, given by $(X\times_SY)_n=X_n\times_SY_n$.

We claim that {\em inner Homs {\rm $\iHom_{\stks\stSETS}:
\stks\stSETS^{fop}\times_\cE\stks\stSETS\to\stks\stSETS$} do exist for this
cartesian $\otimes$-structure.} Indeed, the requirement for this inner
Hom $\iHom_{\stks\stSETS|S}(X,Y)$, where $X$, $Y$ are two simplicial objects
over~$\cE_{/S}$, is
\begin{equation}\label{eq:univ.simpl.ihom}
\iHom_S(K\otimes_SX,Y)\cong\iHom_S(K,\iHom_{\stks\stSETS|S}(X,Y))
\end{equation}
Of course, we are free to choose a simplicial object $K$ after any pullback
$S'\to S$ as well. Putting here $K=\cst{\Delta(n)}$, we see 
that $\iHom_{\stks\stSETS|S}(X,Y)$ has to be computed by the classical formula
\begin{equation}\label{eq:constr.simpl.ihom}
\iHom_{\stks\stSETS|S}(X,Y)_n:=\iHom_S(\cst{\Delta(n)}\otimes_SX,Y)
\end{equation}
Now the verification of \eqref{eq:univ.simpl.ihom} for all 
$K\in\Ob s(\cE_{/T})$ can be done directly in the classical manner, 
the classical proof being intuitionistic, 
or by using an obvious ``local Yoneda lemma''.

\nxsubpoint
($\otimes$-action of $\stks\stSETS$ on $\stks\cC$.)
Now let $\cC$ be a {\em flat\/} stack over~$\cE$. The external cartesian
$\otimes$-action $\oslash:\cC\times_\cE\stSETS\to\cC$ 
of~\ptref{sp:act.SETS.flat.stk} extends to an external cartesian
$\otimes$-action of $\stks\stSETS$ on $\stks\cC$, usually denoted by
$\oslash$ or~$\otimes$. Existence of both flavors of inner Homs with 
respect to this action can be shown again in more or less classical fashion
(cf.~\cite[2.1]{Quillen}); 
in particular, $\iHom_{\stks\stSETS}:(\stks\cC)^{fop}\times_\cE\stks\cC\to
\stks\stSETS$ is still given by~\eqref{eq:constr.simpl.ihom}.

\nxsubpoint\label{sp:cof.gen.ssheaves}
(Standard cofibrant generators of $s\catSets$ and $\stks\stSETS_\cE$.)
Recall that $s\catSets$ admits two standard sets of small cofibrant generators
(cf.~\ptref{sp:cof.gen.ssets}):
\begin{align}
I:=&\{\dot\Delta(n)\to\Delta(n)\,|\,n\geq0\}\\
J:=&\{\Lambda_k(n)\to\Delta(n)\,|\,0\leq k\leq n>0\}
\end{align}
Here, of course, $I$ generates the cofibrations, and $J$ the acyclic 
cofibrations.

Now it is reasonable to expect that for any morphism of topoi $f:\cE'\to\cE$
the pullback functor $f^*:s\cE\to s\cE'$ will preserve cofibrations and 
acyclic cofibrations. Let's apply this to $q:\cE\to\catSets$; 
we obtain two following sets of maps of constant simplicial sheaves:
\begin{align}
\cst{I}=\cst{I}_\cE:=&\{\cst{\dot\Delta(n)}\to\cst{\Delta(n)}\,|\,n\geq0\}\\
\cst{J}=\cst{J}_\cE:=&\{\cst{\Lambda_k(n)}\to\cst{\Delta(n)}\,
|\,0\leq k\leq n>0\}
\end{align}

Clearly, for any reasonable model stack structure on $\stSETS_\cE$ these
sets of morphisms $\cst{I}$, $\cst{J}\subset\Ar (\stks\stSETS_\cE)(e)$,
as well as all their pullbacks,
have to be cofibrations, resp.\ acyclic cofibrations. Moreover, we are tempted
to choose these sets $\cst{I}$ and $\cst{J}$ as cofibrant generators for
this model stack structure. Unfortunately, in general this doesn't
define a model stack structure on $\stks\stSETS_\cE$ unless topos~$\cE$
is {\em sequential\/}, i.e.\ $\projlim_nZ_n\to Z_0$ is an epimorphism
for any projective system $\cdots\to Z_2\to Z_1\to Z_0$ in~$\cE$ with 
epimorphic transition morphisms. This condition is almost never fulfilled:
for example, a metrizable topological space~$X$ defines a sequential topos
iff $X$ is discrete.

\nxsubpoint
(Cofibrant generators of a model stack.)
Given a model stack $\cC\to\cE$ and two sets of morphisms  
$I$, $J\subset\Ar\cC(e)$ in its final fiber, we say that {\em $I$ and $J$ 
are cofibrant generators for~$\cC$} if the acyclic fibrations 
(resp.\ fibrations) of $\cC$ are exactly the morphisms with the local 
RLP with respect to (all pullbacks of) all morphisms from $I$ (resp.~$J$).
Notice that this condition automatically implies that $I$ (resp.~$J$)
consists of cofibrations (resp.\ acyclic cofibrations).

We say that $\cC$ is {\em cofibrantly generated\/} if it admits sets of
cofibrant generators $I$, $J$ as above. In this case the model stack structure
of $\cC$ is completely determined by these two sets. Indeed, we can
define acyclic fibrations (resp.\ fibrations) as morphisms with 
the local RLP with respect to all pullbacks
of morphisms from $\cst{I}$ (resp.~$\cst{J}$), then define acyclic cofibrations
(resp.\ cofibrations) as the morphisms with local LLP with respect to all
fibrations (resp.\ acyclic fibrations), and finally define weak equivalences 
as morphisms which can be decomposed into an acyclic cofibration followed 
by an acyclic fibration.

If we start from two arbitrary sets $I$, $J\subset\Ar\cC(e)$ and define 
cofibrations, fibrations and weak equivalences as above, we'll still have
to check some of the axioms (MS1)--(MS5) to show that we've indeed obtained
a model stack structure; in particular, we need to prove (MS5) 
(factorization), (MS2) (localness and stability of weak equivalences under 
retracts) and (MS3) (2-out-of-3 for weak equivalences). 
In addition, we have to check that e.g.\ acyclic fibrations are indeed 
the morphisms that are both fibrations and acyclic.

\begin{DefD}\label{def:dist.morph.sSETS}
{\rm (Distinguished classes of morphisms in $\stks\stSETS_\cE$.)}
We define subsets of fibrations, acyclic fibrations, {\bf strong} 
cofibrations and {\bf strong} acyclic cofibrations 
in sets of morphisms in fibers of $\stks\stSETS_\cE$ 
by means of the above construction, applied 
to standard cofibrant generators $\cst{I}$, $\cst{J}$ 
of~\ptref{sp:cof.gen.ssheaves}.
\end{DefD}

This doesn't define a model stack structure on $\stks\stSETS_\cE$
unless $\cE$ is sequential: in general 
we have to sacrifice either the local lifting axiom (MS4)
or the factorization axiom (MS5). This is due to the fact that the class of
{\em strong\/} cofibrations defined above by the local LLP with respect to
all acyclic fibrations is too small, and in particular it is not closed 
under infinite sums and sequential inductive limits.

We deal with this problem by embedding the strong cofibrations into a larger
class of cofibrations, stable under small sums, composition, pushouts,
sequential inductive limits, and pullbacks with respect to arbitrary morphisms
of topoi. In other words, we drop the local lifting axiom (MS4) to be able
to obtain factorization (MS5) and even~(MS5f).

\begin{DefD}\label{def:closed.class}
We say that a class of morphisms $\propP$ in fibers of a 
{\em flat\/} stack~$\cC$ over~$\cS$ is {\em closed\/} if:
\begin{enumerate}
\item $\propP$ is local, and in particular stable under all pullbacks~$\phi^*$;
\item $\propP$ contains all isomorphisms in fibers of~$\cC$;
\item $\propP$ is stable under pushouts;
\item $\propP$ is stable under composition;
\item $\propP$ is stable under ``sequential composition'', i.e.\ 
whenever $A_0\stackrel{i_0}\to A_1\stackrel{i_1}\to A_2\to\cdots$ is a
composable infinite sequence of morphisms from $\propP$, its ``composition''
$A_0\stackrel i\to A_\infty=\injlim_nA_n$ is also in $\propP$.
\item $\propP$ is stable under retracts (hence also under local retracts);
\item $\propP$ is stable under arbitrary (small) direct sums;
\item $\propP$ is stable under all $\phi_!$, hence under all 
$-\otimes_ST=\phi_!\phi^*$; notice that this condition makes 
sense only for a flat $\cC$. (In fact, stability under $-\otimes_ST$ would 
suffice for our applications.)
\end{enumerate}

Similarly, a class of morphisms $\propP$ in a category $\cC$ is 
{\em (globally) closed\/} if it satisfies those of the above conditions that 
are applicable, i.e.\ all with exception of the first and the last. 
If we omit condition 6) from the above list, we obtain the notions of a
{\em semiclosed\/} and {\em (globally) semiclosed\/} class of morphisms.
\end{DefD}

For example, the classes of cofibrations and acyclic cofibrations in a model
category~$\cC$ are (globally) closed. Unfortunately, this is not true 
for model stacks $\cC/\cS$: in general we don't have 5), 7) and~8).

\begin{DefD}\label{def:closure.class}
Given a flat stack $\cC/\cS$ and 
any set of morphisms $I\subset\Ar\cC(e)$ in the 
final fiber of~$\cC$, we denote by $\Cl I$ the {\em closure\/} of $I$
in~$\cC$, i.e.\ the smallest closed class of morphisms in fibers of~$\cC$
that contains~$I$. We define the {\em semiclosure\/} $\SCl I$ similarly;
and if $\cC$ is just a category, we obtain corresponding global notions
$\GCl I$ and $\GSCl I$, which will be usually denoted simply by $\Cl I$
and~$\SCl I$.
\end{DefD}

Clearly, any intersection of closed sets of morphisms is still closed, so
the closure $\Cl I$ always exists and can be constructed as the intersection
of all closed sets of morphisms containing~$I$.

If $\cC$ is a model category cofibrantly generated by 
sets $I$ and $J$ of morphisms with small sources, then Quillen's small object 
argument shows that $\Cl I$ is the set of cofibrations, and 
$\Cl J$ the set of acyclic cofibrations of~$\cC$. Unfortunately, this is 
not true for model stacks.

\begin{DefD}\label{def:acof.simpl.sh} 
(Acyclic cofibrations of simplicial sheaves.) 
The set of cofibrations (resp.\ acyclic cofibrations) in $\stks\stSETS_\cE/\cE$
is by definition the closure in $\stks\stSETS_\cE$ 
of $\cst{I}_\cE$ (resp.\ $\cst{J}_\cE$) of~\ptref{sp:cof.gen.ssheaves}.
\end{DefD}

Now we are going to check the factorization axiom (MS5f) for these sets of
(acyclic) cofibrations and (acyclic) fibrations, and to prove that all 
strong (acyclic) cofibrations in the sense of~\ptref{def:dist.morph.sSETS} are
(acyclic) cofibrations in the sense of above definition. For this 
we'll need a version of Quillen's small object argument.

\begin{DefD} 
(Sequentially small and finitely presented objects.)
An object $X$ of a category~$\cC$ is {\em (sequentially) small\/} 
(resp.\ {\em finitely presented\/}) 
if $\Hom_\cC(X,-):\cC\to\catSets$ commutes with sequential inductive limits
(resp.\ filtered inductive limits). Similarly, an object $X\in\cC(S)$ of 
a fiber of a stack $\cC/\cS$ is {\em locally (sequentially) small\/}
(resp.\ {\em locally finitely presented\/}) if $\iHom_{\cC_{/T}}(X|_T,-):
\cC(T)\to\stSETS_\cS(T)=\widetilde{\cS_{/T}}$ commutes with 
sequential (resp.\ filtered) inductive limits in $\cC(T)$, for any 
$T\to S$ in~$\cS$.
\end{DefD}

Evidently, these properties are indeed local in the sense 
of~\ptref{sp:loc.prop.obj.mor}. When $\cS=\st1$ or $\catSets$, we recover 
the corresponding ``global'' notions. Notice that the substacks of 
locally small and locally finitely presented objects are stable under 
finite inductive limits in fibers of~$\cC$. 

\nxsubpoint (Small objects in $\stks\stSETS_\cE$.)
First of all, notice that the {\em standard simplices\/}
$\cst{\Delta(n)}_\cE$ are locally finitely presented, hence also locally small,
in $s\cE$ (or, more precisely, in $(\stks\stSETS_\cE)(e)$). Indeed, 
this follows from the formula $\iHom(\cst{\Delta(n)},X)=X_n$, valid for any
$X\in\Ob s\cE$, and the fact that arbitrary inductive limits in $s\cE$
are computed componentwise, hence $\iHom(\cst{\Delta(n)},-):X\mapsto X_n$ 
commutes with arbitrary inductive limits.

Since any finite inductive limit of locally finitely presented objects 
is again finitely presented, and the constant simplicial sheaf functor 
$sq^*:s\catSets\to s\cE$ commutes with arbitrary inductive limits,
we see that $\cst A_\cE$ is locally finitely presented for any {\em finite\/}
simplicial set~$A$, such simplicial sets being exactly the finite inductive 
limits of standard simplices.

Now all simplicial sets involved in standard cofibrant generators 
$I$ and $J$ of $s\catSets$ are finite, hence the sources and targets of 
all morphisms of $\cst I_\cE$ and $\cst J_\cE$ are locally finitely presented 
and locally small.

\nxsubpoint ($\iHom(\cst{A},-)$ for a finite simplicial set~$A$.)
Before going on, let's consider the functors $\iHom(\cst{A},-):s\cE\to\cE$,
where $A$ is any finite simplicial set. We can represent $A$ as a cokernel 
of a couple of morphisms between finite sums of standard simplices:
$A=\Coker\bigl(\coprod_{j=1}^t\Delta(m_j)\rightrightarrows
\coprod_{i=1}^s\Delta(n_i)\bigr)$ in $s\catSets$. Applying exact 
functor $sq^*:A\mapsto\cst A$ we obtain a similar relation in~$s\cE$.
Now $\iHom$ is clearly left exact in both arguments, hence 
\begin{equation}
\iHom(\cst{A},X)\cong \Ker\bigl(\prod_{i=1}^s X_{n_i}\rightrightarrows
\prod_{j=1}^t X_{m_j}\bigr)
\end{equation}
In other words, $\iHom(\cst{A},X)$ is a finite projective limit of the 
components of~$X$. 

\nxsubpoint (Fibrations and acyclic fibrations.)
This is applicable in particular to $A=\dot\Delta(n)=
\sk_{n-1}\Delta(n)$, and we obtain $\iHom(\cst{\dot\Delta(n)},X)\cong
(\cosk_{n-1}X)_n$ for any $X\in\Ob s\cE$, similarly to the classical 
case~\ptref{sp:cof.gen.ssets}. On the other hand, $p:X\to Y$ has the 
local RLP with respect to some $i:A\to B$ in~$\cC(S)$, $\cC/\cE$ any stack,
iff $\iHom(B,X)\to\iHom(A,X)\times_{\iHom(A,Y)}\iHom(B,Y)$ is an epimorphism
in~$\cE_{/S}$. Applying this to all $i$ from $\cst I_\cE$, we see that 
{\em $p:X\to Y$ is an acyclic fibration iff the natural morphisms
$X_n\to(\cosk_{n-1}X)_n\times_{(\cosk_{n-1}Y)_n}Y_n$ are epimorphic 
in~$\cE_{/S}$ for all~$n\geq0$,} i.e.\ a condition completely similar to
its classical counterpart~\ptref{sp:cof.gen.ssets}.

Similarly, applying this to all $i:A\to B$ from $\cst J_\cE$, we obtain 
a description of fibrations $p:X\to Y$ in~$s\cE$ in terms 
of a countable set of conditions each requiring a certain morphism 
between finite projective limits of components of~$X$ and~$Y$ to be 
an epimorphism in $\cE_{/S}$.

\nxsubpoint\label{sp:stab.fibr.pullbacks}
(Stability of fibrations under pullbacks with respect to any morphism 
of topoi.)
An immediate consequence is that {\em for any finite simplicial set~$A$ 
functors $\iHom(\cst A,-)$ commute with pullbacks $f^*$ with respect to 
any morphism of topoi $f:\cE'\to\cE$}, i.e.\ 
$f^*\iHom(\cst A_\cE, X)\cong\iHom(\cst A_{\cE'},f^*X)$ for any 
$X\in\Ob s\cE$. Indeed, this follows from the fact that $\iHom(\cst A,X)$ 
is just some finite projective limit of components of~$X$, and
the exactness of~$f^*$. Therefore, the morphisms between finite projective 
limits of components of~$X$ and~$Y$, the epimorphicity of which is 
equivalent to $p:X\to Y$ being a fibration (resp.\ an acyclic fibration), 
are preserved by~$f^*$. Now $f^*$ is exact, and in particular preserves 
epimorphisms. We conclude that {\em fibrations and acyclic fibrations in
$s\cE$ are preserved by pullbacks with respect to arbitrary morphisms of 
topoi $f:\cE'\to\cE$.}

\nxsubpoint
(One step of Quillen's construction.)
Let $I=\{i_\alpha:A_\alpha\to B_\alpha\}_{\alpha\in I}$ 
be a (small) set of morphisms in a fiber $\cC(S)$ of a 
flat stack~$\cC/\cE$. Then for any morphism $X\stackrel f\to Y$ 
in~$\cC(S)$ we construct a functorial decomposition $X\stackrel{\xi_f}\longto
R_{I,Y}(X)\stackrel{R(f)}\longto Y$ as follows. 
For any $\alpha\in I$ denote by 
$H_\alpha=H_\alpha(f)\in\Ob\cE_{/S}$ the {\em sheaf\/} of diagrams
\begin{equation}\label{eq:cdiag.of.Halpha}
\xymatrix{
A_\alpha\ar[r]^{i_\alpha}\ar@{-->}[d]^{u_\alpha}&
B_\alpha\ar@{-->}[d]^{v_\alpha}\\
X\ar[r]^{f}&Y}
\end{equation}
In other words, $H_\alpha$ is just the fibered product
$\iHom(A_\alpha,X)\times_{\iHom(A_\alpha,Y)}\iHom(B_\alpha,X)$ in~$\cE_{/S}$.
Now $R_{I,Y}(X)$ is constructed by means of the following cocartesian square:
\begin{equation}\label{eq:def.Qstep}
\xymatrix{
\coprod_\alpha A_\alpha\otimes H_\alpha\ar[r]^{(i_\alpha)}\ar[d]_{(u_\alpha)}&
\coprod_\alpha B_\alpha\otimes H_\alpha\ar@{-->}[d]\ar[ddr]^{(v_\alpha)}\\
X\ar@{-->}[r]^{\xi_f}\ar[rrd]_{f}&R_{I,Y}(X)\ar@{-->}[rd]|{R(f)}\\
&&Y}
\end{equation}
The top horizontal arrow $(i_\alpha)$ is simply $\coprod_\alpha 
i_\alpha\otimes\id_{H_\alpha}$. As to the left vertical arrow 
$(u_\alpha):\coprod_\alpha A_\alpha\otimes H_\alpha\to X$, its components
$u'_\alpha:A_\alpha\otimes H_\alpha\to X$ are constructed by
applying $\id_{A_\alpha}\otimes-$ to the natural projection $H_\alpha\to
\iHom(A_\alpha,X)$ and composing the result with the canonical ``evaluation''
morphism $\ev_{A_\alpha,X}:A_\alpha\otimes\iHom(A_\alpha,X)\to X$. 
The arrow $(v_\alpha)$ is constructed similarly with the aid of
the second projection $H_\alpha\to\iHom(B_\alpha,Y)$ and $\ev_{B_\alpha,Y}$, 
and the commutativity of the outer circuit of the above diagram follows 
from the definition of~$H_\alpha$ and the explicit description of 
evaluation morphisms involved.

So let us recall the construction of $\ev_{A,X}:A\otimes_S\iHom_S(A,X)\to X$, 
for any $A$, $X\in\Ob\cC(S)$. Put $H:=\iHom_S(A,X)$, and let $\phi:H\to S$
be the structural morphism of $H\in\Ob\cE_{/S}$. Then by definition 
$A\otimes_S H=\phi_!\phi_*A$, hence giving $\ev_{A,X}:\phi_!\phi^*A\to X$ 
is equivalent to giving some $\ev_{A,X}^\flat:A|_H=\phi^*A\to X|_H$, i.e.\ 
an element $\ev_{A,X}^\flat\in\Hom_{\cC(H)}(A|_H,X|_H)=(\iHom_S(A,X))(H)$,
i.e.\ an $S$-morphism from $H$ into $\iHom_S(A,X)=H$. Of course, we put
$\ev_{A,X}^\flat:=\id_H$.

\nxsubpoint (Sequential closure of Quillen's construction.)
Notice that $R_{I,Y}$ can be treated as a functor $\cC(S)_{/Y}\to\cC(S)_{/Y}$,
and $\xi_f:X\to R_{I,Y}(X)$ is a functorial morphism $\xi:\Id_{\cC(S)_{/Y}}
\to R_{I,Y}$. Therefore, we can iterate the construction:
$X\to R_{I,Y}(X)\to R_{I,Y}^2(X)\to\cdots$. Put $R_{I,Y}^\infty(X):=
\injlim_n R_{I,Y}^n(X)$; when no confusion can arise, we denote 
$R_{I,Y}^\infty(X)$ simply by $R_I(X)$ or $R_I(f)$. By construction,
we obtain a functorial decomposition $X\to R_I(f)\to Y$ of any $X\stackrel f
\to Y$. Furthermore, the construction of $R_{I,Y}(X)$ commutes with arbitrary 
pullbacks $\phi^*:\cC(S)\to\cC(T)$, hence the same is true for $R_I(X)$, 
the stack~$\cC$ being flat. We conclude that {\em functors 
$R_I:\Ar\cC(S)\to\cC(S)$ define together a cartesian functor $R_I:
\stkAr\cC\to\cC$ of stacks over~$\cE$.}

\begin{PropD}\label{prop:Quillen.fact} (Quillen's factorization.)
If $I\subset\Ar\cC(e)$ is a small set of morphisms in the final fiber
of a flat stack $\cC$ over a topos $\cE$, and all morphisms from~$I$ have 
locally small sources, then $X\to R_I(f)$ belongs to the closure $\Cl I$
of~$I$ in~$\cC$, and $R_I(f)\to Y$ has the local RLP with respect to all
morphisms from~$I$, for any $f:X\to Y$ in $\cC(S)$. In other words,
$R_I:\stkAr\cC\to\cC$ together with appropriate natural transformations 
provides a functorial factorization for all $f:X\to Y$,
similar to that of~(MS5f).
\end{PropD}
\begin{Proof}
The fact that $X\to R_{I,Y}(X)$ belongs to $\Cl I$ is immediate from
\eqref{eq:def.Qstep}, \ptref{def:closed.class} and~\ptref{def:closure.class}.
Since $X\to R_I(f)=R_{I,Y}^\infty(X)$ is a sequential composition of morphisms
of the above form, it lies in $\Cl I$ as well. Now let us prove that 
$R_I(f)\to Y$ has the local RLP with respect to any 
$i_\alpha:A_\alpha\to B_\alpha$ from~$I$. Since the construction of~$R_I(f)$ 
is compatible with all pullbacks $\phi^*$, we may assume that our lifting 
problem $u:A_\alpha\to R_I(f)$, $v:B_\alpha\to Y$ is given inside $\cC(S)$ 
itself, and, moreover, may even assume $S=e$, replacing $\cC$ and $\cE$ with 
$\cC_{/S}$ and $\cE_{/S}$ if necessary. Now $A_\alpha$ is locally small
by assumption, hence $\iHom(A_\alpha,R_I(f))=\iHom(A_\alpha,\injlim_n
R_{I,Y}^n(X))\cong\injlim_n\iHom(A_\alpha, R_{I,Y}^n(X))$, so we can
find (locally, i.e.\ on some cover) some $n\geq0$ and a morphism
$u':A_\alpha\to R_{I,Y}^n(X)$ that represents~$u$. Notice that 
$u'$ and $v$ define a commutative diagram~\eqref{eq:cdiag.of.Halpha} 
with $f:X\to Y$ replaced by $f_n:R_{I,Y}^n(X)\to Y$, 
i.e.\ $w:=(u',v)$ is a section of $H_\alpha=H_\alpha(f_n)$. 
Now $\id_{B_\alpha}\otimes w$ defines a morphism 
$B_\alpha\to B_\alpha\otimes H_\alpha$; composing it with the natural 
embedding of $B_\alpha\otimes H_\alpha$ into the coproduct
$\coprod_{\beta\in I}B_\beta\otimes H_\beta$ and the right vertical arrow
of~\eqref{eq:def.Qstep}, we obtain a morphism 
$h':B_\alpha\to R_{I,Y}^{n+1}(X)$, making the following diagram commutative:
\begin{equation}
\xymatrix{
A_\alpha\ar[r]^{i_\alpha}\ar[d]^{u'}&B_\alpha\ar[d]^{h'}\\
R_{I,Y}^n(X)\ar[r]^{\xi_{f_n}}&R_{I,Y}^{n+1}(X)}
\end{equation}
Composing $h'$ with the natural embedding of $R_{I,Y}^{n+1}(X)$ into 
the inductive limit $R_I(f)$, we obtain a (local) solution $h:B_\alpha\to
R_I(f)$ of the lifting problem $(u,v)$, q.e.d.
\end{Proof}

\begin{CorD}\label{cor:funct.decomp.sSETS}
Any morphism $f:X\to Y$ in $s\cE$ can be functorially 
decomposed both into a cofibration $X\to R_{\cst I}(X)$ followed by an acyclic 
fibration $R_{\cst I}(X)\to Y$, and into an acyclic cofibration 
$X\to R_{\cst J}(X)$ followed by a fibration $R_{\cst J}(X)\to Y$. In other 
words, (MS5f) holds in $\stks\stSETS_\cE/\cE$ for this choice of 
(acyclic) fibrations and cofibrations.
\end{CorD}
\begin{Proof} Immediate from~\ptref{prop:Quillen.fact} and
definitions \ptref{def:acof.simpl.sh} and~\ptref{def:dist.morph.sSETS}.
\end{Proof}

\nxsubpoint
Unfortunately, in our case the cofibrations (resp.\ acyclic cofibrations) 
are not necessarily strong, i.e.\ they needn't have the local LLP with 
respect to all acyclic fibrations (resp.\ fibrations). The reason for this
is that, say, the class of strong cofibrations still contains $\cst I$,
but is not closed (e.g.\ it is not closed under infinite direct sums).
On the other hand, the opposite inclusion is still true: 
{\em any strong cofibration is a cofibration, and similarly for 
strong acyclic cofibrations.} Indeed, we can apply the classical reasoning 
of \cite{Quillen}: any strong cofibration $f:X\to Y$ can be factorized into
a cofibration $i:X\to Z$ and an acyclic fibration $p:Z\to Y$, and
$f$ has the local LLP with respect to~$p$, hence we can locally find a section
$\sigma:Y\to Z$ of~$p$, hence $f$ is a local retract of~$i$, hence a 
cofibration, the class of cofibrations being closed.

\begin{DefD}
We say that a morphism $f:X\to Y$ in $s\cE$, or in $s\cE_{/S}=
(\stks\stSETS_\cE)(S)$ is a {\bf weak equivalence} if it can be factorized 
into an acyclic cofibration followed by an acyclic fibration.
\end{DefD}

It is not evident from this definition that being a weak equivalence is 
a local property, or that weak equivalences satisfy the 2-out-of-3 axiom
(MS3). We are going to check this in a moment, at least for topoi with enough
points.

\begin{PropD}
Any acyclic cofibration in $s\cE$ is both a cofibration and a weak 
equivalence, and any acyclic fibration is both a fibration and a weak 
equivalence.
\end{PropD}
\begin{Proof} (a) Acyclic cofibrations and acyclic fibrations are 
weak equivalences by definition. (b) Let us show that any acyclic fibration
in $s\cE$ is a fibration. Indeed, any morphism $\Lambda_k(n)\to\Delta(n)$ 
of~$J$ can be decomposed into $\Lambda_k(n)\to\dot\Delta(n)\to\Delta(n)$,
and both morphisms of this decomposition are pushouts of morphisms from~$I$. 
Clearly, this is still true for $\cst I_\cE$ and $\cst J_\cE$. Since strong
cofibrations in~$s\cE$ contain $\cst I_\cE$ and are stable under pushouts 
and composition, we see that $\cst J_\cE$ consists of strong cofibrations,
i.e.\ any acyclic fibration $p:X\to Y$ in $s\cE$ has the local RLP with 
respect to all morphisms from $\cst J_\cE$, i.e.\ is a fibration.

(c) Now let's prove that acyclic cofibrations are cofibrations, i.e.\ 
$\Cl\cst J_\cE\subset\Cl\cst I_\cE$. It is enough to show for this that 
$\cst J_\cE\subset\Cl\cst I_\cE$, i.e.\ that $\cst J_\cE$ consists of 
cofibrations. But this is shown by the same reasoning as above: any morphism
of $\cst J_\cE$ is a composition of pushouts of morphisms from $\cst I_\cE$,
hence lies in $\Cl\cst I_\cE$.
\end{Proof}

\begin{PropD}\label{prop:stab.distcl.topos.pb}
Let $f:\cE'\to\cE$ be any morphism of topoi. Then $f^*=sf^*:s\cE\to s\cE'$ 
preserves cofibrations, acyclic cofibrations, fibrations, acyclic fibrations
and weak equivalences.
\end{PropD}
\begin{Proof}
(a) The stability of fibrations and acyclic fibrations under $f^*$ 
has been already shown in~\ptref{sp:stab.fibr.pullbacks}, 
since the property of $p:X\to Y$ to be a fibration or an
acyclic fibration is equivalent to epimorphicity of certain morphisms between
finite projective limits of components $X_n$ and~$Y_n$.
(b) The stability of weak equivalences will immediately follow from that 
of acyclic cofibrations and acyclic fibrations.
(c) Let's deal with the case of, say, cofibrations. We want to show that
$f^*$ maps $\Cl\cst I_\cE$ into $\Cl\cst I_{\cE'}$. Consider for this 
the class of morphisms $\propP$ in fibers of $\stks\stSETS_\cE/\cE$ defined as 
follows: a morphism $i:A\to B$ in $(\stks\stSETS_\cE)(S)=s(\cE_{/S})$ 
belongs to $\propP$ iff $f_S^*(i)$ is a cofibration in $s(\cE'_{/f^*S})$, where
$f_S:\cE'_{/f^*S}\to\cE_{/S}$ denotes the induced morphism of topoi. 
One checks immediately that $\propP$ is closed and contains $I_\cE$, hence it 
has to contain $\Cl\cst I_\cE$, i.e.\ all cofibrations of $\stks\stSETS_\cE$.
This proves the required stability for cofibrations; the case of acyclic 
cofibrations is dealt with similarly, with $\cst I$ replaced by~$\cst J$.
\end{Proof}

We have already remarked that when $\cE=\catSets$, we recover the usual 
model category structure on~$\catSets$, i.e.\ the (acyclic) cofibrations, 
fibrations and weak equivalences in $s\catSets$ obtained from our definitions 
coincide with their classical counterparts.
\begin{PropD}\label{prop:f.af.w.pointwise}
Let $f:X\to Y$ be a morphism in $s\cE$. If $f$ is an acyclic cofibration,
cofibration, acyclic fibration, fibration or weak equivalence in $s\cE$,
then $f_p=p^*(f):X_p\to Y_p$ has the same property in $s\catSets$ for 
any point $p:\catSets\to\cE$ of topos~$\cE$. Conversely, if $\cE$
has enough points, and all $f_p:X_p\to Y_p$ are fibrations, acyclic fibrations
or weak equivalences in $s\catSets$, then the same is true for $f$, i.e.\ 
these properties can be checked pointwise.
\end{PropD}
\begin{Proof}
(a) The first statement is a special case of \ptref{prop:stab.distcl.topos.pb},
once we take into account the remark made just before this proposition. 

(b) Let's prove the second statement. 
The case of fibrations and acyclic fibrations
is evident: indeed, $f:X\to Y$ is a fibration (resp.\ an acyclic fibration)
iff certain morphisms between finite projective limits of components of
$X$ and~$Y$ are epimorphic (cf.~\ptref{sp:stab.fibr.pullbacks}). Now if
$\cE$ has enough points, then $u:Z\to W$ is an epimorphism in $\cE$ iff
all $u_p:Z_p\to W_p$ are surjective; since $p^*$ is exact, we obtain 
our statement for fibrations and acyclic fibrations. 

(c) Suppose that all $f_p:X_p\to Y_p$ are weak equivalences in $s\catSets$.
Consider the decomposition $X\stackrel i\to Z\stackrel s\to Y$ of~$f$
into an acyclic cofibration followed by a fibration, existing 
by~\ptref{cor:funct.decomp.sSETS}. For any point $p$ the composition
$s_p\circ i_p=f_p$ is a weak equivalence by assumption, as well as $i_p$,
hence $s_p$ is a weak equivalence as well by the 2-out-of-3 axiom in
$s\catSets$, hence $s_p$ is an acyclic fibration for all points $p$ of~$\cE$,
and we can conclude that $s$ is an acyclic fibration in~$\cE$ by (b).
This proves that $f=s\circ i$ is a weak equivalence, q.e.d.
\end{Proof}

\begin{CorD}
If topos~$\cE$ has enough points, then weak equivalences in 
$\stks\stSETS_\cE$ constitute a local class, stable under retracts and 
composition, and satisfying the 2-out-of-3 axiom (MS3). Moreover,
in this case $f:X\to Y$ is an acyclic fibration in $s\cE$ iff it is 
both a fibration and a weak equivalence.
\end{CorD}

Of course, these statements must be true without any assumption on~$\cE$. 
However, when~$\cE$ has enough points, we are able to deduce the above 
statements easily from the classical (highly non-trivial) fact 
that $s\catSets$ is indeed a model category. 

\nxsubpoint\label{q:c.ac.ptwise?} {\bf Question.} 
Can cofibrations and acyclic cofibrations 
be checked pointwise? We believe that the answer is negative, and even 
that not all objects of $s\cE$ are cofibrant. The reason for our disbelief
is the following. The classical proof of this fact (in $s\catSets$)
goes as follows: we write $A=\injlim_n\sk_nA$, and then write each
$\sk_{n-1}A\to\sk_nA$ as a pushout of cofibration 
$\dot\Delta(n)\otimes\Sigma_n\to\Delta(n)\otimes\Sigma_n$, where
$\Sigma_n\subset A_n$ is the subset of {\em non-degenerate simplices},
i.e.\ the complement of the union of images of all degeneracy maps 
$s_n^i:A_{n-1}\to A_n$. Now we observe that this proof is not intuitionistic
since it uses the fact that any simplex is either degenerate or non-degenerate,
so we see that this proof doesn't work in $s\cE$, i.e.\ ``there is no reason''
for $A$ to be cofibrant. Of course, this doesn't prove anything from the 
formal point of view. We have good candidates for such a non-cofibrant
$A$ over topos $\cE=\hat\cS$, $\cS=[1]=\{0\to 1\}$: consider for 
example the pushout of $\Delta(1)\otimes U\to\Delta(1)$ and
$\Delta(1)\otimes U\to\Delta(0)\otimes U$, where $U$ is the only non-trivial
open object of~$\cE$, but the problem here is 
that it is quite complicated to show that something is {\em not\/}
a cofibration\dots

\nxsubpoint {\bf Question.} Is it true that if $f:X\to Y$ in $s\cE$ is 
both a weak equivalence and a cofibration, then it is an acyclic cofibration?
We expect the answer to this question to be positive, due to the fact that
we have postulated stability of cofibrations and acyclic cofibrations under
all ``local sums'' $\phi_!$, not just under all $-\otimes_ST=\phi_!\phi^*$.

\nxpointtoc{Pseudomodel stacks}
Let us summarize the properties of distinguished classes of morphisms 
in $s\cE$ and $\stks\stSETS_\cE$. We see that these classes do not 
satisfy the axioms for a model category or a model stack.
The following definition describes a self-dual subset of properties 
that we have in $s\cE$ and $\stks\stSETS_\cE$:

\begin{DefD}\label{def:psmod.cat} (Pseudomodel categories.)
We say that a category~$\cC$ with five distinguished local classes of morphisms
(localness understood here as closedness under isomorphisms), called
{\em cofibrations, acyclic cofibrations, fibrations, acyclic fibrations} 
and {\em weak equivalences} is a {\bf pseudomodel category} if
the following conditions hold:
\end{DefD}
{

\myindent{(PM1)} 
The category~$\cC$ is closed under finite inductive and 
projective limits.

\myindent{(PM2)} 
Each of the five distinguished classes of morphisms in~$\cC$ 
is closed under composition and retracts. Acyclic cofibrations and 
cofibrations are stable under pushouts and finite direct sums.
Acyclic fibrations and fibrations are stable under pullbacks and
finite products.
\smallskip

\myindent{(PM3)} 
(``2 out of 3'') Given two morphisms $f:X\to Y$ and $g:Y\to Z$, 
such that any two of $f$, $g$ and $g\circ f$ are weak equivalences, 
then so is the third.
\smallskip

\myindent{(PM5)}
(Factorization) Any morphism $f:X\to Y$ can be factored 
both into an acyclic cofibration followed by a fibration, or into 
a cofibration followed by an acyclic fibration: $f=p\circ j=q\circ i$, with 
$i$ a cofibration, $j$ an acyclic cofibration, $p$ a fibration, 
and $q$ an acyclic fibration.
\smallskip

\myindent{(PM8)}
Any acyclic cofibration is both a weak equivalence and a cofibration.
Any acyclic fibration is both a weak equivalence and a fibration.
Any weak equivalence can be factorized into an acyclic cofibration followed
by an acyclic fibration.

}
\nxsubpoint
Notice the absence of the lifting axiom (MS4). All other axioms but (PM8) are 
more or less similar to those of a closed model category; however, 
we have to include some stabilities in (PM2), which are consequences of 
other axioms for model categories; and we mention {\em five\/} distinguished 
classes in this axiom. As to (PM8), notice that we don't require any of 
the opposite implications, summarized in the following stronger condition:

\smallskip
{\myindent{(PM8+)} Acyclic cofibrations are exactly the morphisms which 
are both cofibrations and weak equivalences. Acyclic fibrations are exactly
the morphisms which are both fibrations and weak equivalences.
\smallbreak

}
The reason for omitting this stronger condition is that we have shown only 
its second half for $s\cE$, and we want the axioms of a pseudomodel category
to be self-dual.

\begin{DefD}
A {\bf pseudomodel stack} $\cC$ over a topos~$\cE$ is a flat stack
$\cC/\cE$ together with five local classes of morphisms in its fibers,
called as in~\ptref{def:psmod.cat}, satisfying the above axioms
(PM1)--(PM3), (PM5), (PM8) in each fiber $\cC(S)$.
\end{DefD}
Notice that each fiber $\cC(S)$ of a pseudomodel stack~$\cC$ is a pseudomodel
category.

\nxsubpoint
We know that $s\cE$ has a natural pseudomodel category structure,
and $\stks\stSETS_\cE$ a natural pseudomodel stack structure, at least
if $\cE$ has enough points. These pseudomodel structures enjoy some 
additional properties, e.g.\ closedness of sets of cofibrations and
acyclic cofibrations, cartesian functorial factorization (MS5f), 
and the second half of (PM8+).

\nxsubpoint (Homotopic category of a pseudomodel category.)
Once we have a class of weak equivalences in a pseudomodel category~$\cC$,
we can define the corresponding {\em homotopic category\/}
$\Ho\cC=\Ho_0\cC$ as the localization of~$\cC$ with respect to the 
weak equivalences, at least if we don't mind enlarging the universe~$\univU$.
In most cases one can show {\em a posteriori\/} that the homotopic category
will be still a $\univU$-category, so this enlargement of universe turns out
to be unnecessary. Therefore, one can define left and right derived functors
$\dL F$, $\dR F:\Ho\cC\to\Ho\cD$ of any functor $F:\cC\to\cD$ between 
pseudomodel categories (cf.~\ptref{def:der.funct}). When $\cD$ is an arbitrary
category, without any pseudomodel structure, we apply this definition 
considering isomorphisms in~$\cD$ as weak equivalences. For example,
$\dL F:\Ho\cC\to\Ho\cD$ is also the left derived functor of
$\gamma_\cD\circ F:\cC\to\Ho\cD$.

When it comes to considering homotopic fibered categories, (pre)stacks etc.,
we'll consider only the fibered category $\HO_0\cC$, given by
$(\HO_0\cC)(S):=\Ho\cC(S)$. 

\nxsubpoint (Cofibrant objects, replacements,\dots)
We define cofibrant and fibrant objects, (strict) (co)fibrant replacements
inside a pseudomodel category in the same way we did it before for 
model categories. Notations like $\cC_c$, $\cC_f$ and $\cC_{cf}$ will be 
used to denote the full subcategories consisting of cofibrant objects etc.\ 
as before. Notice that $\cC_c$, $\cC_f$ and $\cC_{cf}$ satisfy all the
axioms of pseudomodel categories except (PM1), if we agree 
to consider in (PM2) only those direct sums, pushouts etc.\ which lie inside 
these subcategories.

\begin{ThD}\label{th:ex.left.der.psmod}
Let $\cC$ be a pseudomodel category, $F:\cC\to\cB$ any functor transforming
weak equivalences between cofibrant objects into isomorphisms. Then 
$F$ admits a left derived $\dL F:\Ho\cC\to\cB$, which can be computed 
with the aid of cofibrant replacements, i.e.\ 
$\eta_P:\dL F(\gamma P)\to F(P)$ is an isomorphism for all 
$P\in\Ob\cC_c$, and $\dL F(\gamma X)\cong F(P)$ for any cofibrant 
replacement $P\to X$ of~$X$.
\end{ThD}
Applying this theorem to $\cB=\Ho\cD$, we obtain
\begin{CorD}
If $F:\cC\to\cD$ is a functor between pseudomodel categories, transforming 
weak equivalences between cofibrant objects of $\cC$ into weak equivalences
in~$\cD$, then $F$ admits a left derived $\dL F$, which can be computed 
by means of cofibrant replacements.
\end{CorD}
\begin{Proof} (of \ptref{th:ex.left.der.psmod})
The proof is somewhat similar to that of~\ptref{th:homot.mod.stk}.

(a) First of all, we construct a functor $\tilde F:\cC\to\cB$ as follows.
For any $X\in\Ob\cC$ we choose any its strict cofibrant replacement,
i.e.\ any acyclic fibration $Q\stackrel p\to X$ with a cofibrant~$Q$; 
such cofibrant replacements exist by (PM5). Now we put $\tilde F(X):=F(Q)$. 
We have to check that $\tilde F(X)$ doesn't depend on the choice of
$Q\stackrel p\to X$ up to a canonical isomorphism. We reason as in
the proof of~\ptref{th:homot.mod.stk}: if $Q'\stackrel{p'}\to X$ is another
strict cofibrant replacement, we consider any strict cofibrant replacement
$Q''\stackrel r\to Q\times_XQ'$; composing $r$ with the projections of
$Q\times_XQ'$ and taking (PM2) into account, we obtain two acyclic fibrations
$\sigma:Q''\to Q$ and $\sigma':Q''\to Q'$ between cofibrant objects,
such that $p\sigma=p'\sigma'$. Hence $F(\sigma)$ and $F(\sigma')$ are 
isomorphisms in~$\cB$, so we obtain an isomorphism 
$\phi_{Q''}:=F(\sigma')\circ F(\sigma)^{-1}:F(Q)\simto F(Q')$, compatible with 
$F(p)$ and~$F(p')$. Next, we need to check independence of this isomorphism
of the choice of $Q''\stackrel r\to Q\times_XQ'$; this is done similarly:
if $Q''_1\stackrel{r_1}\to Q\times_XQ'$ is another such choice, 
we choose any strict cofibrant replacement
$Q''_2$ of the fibered product of $Q''$ and $Q''_1$ over $Q\times_XQ'$, thus
obtaining acyclic fibrations $Q''_2\stackrel u\to Q''$ and 
$Q''_2\stackrel v\to Q''_1$. Now $\phi_{Q''_2}=F(\sigma'\circ u)\circ
F(\sigma\circ u)^{-1}=F(\sigma')\circ F(u)\circ F(u)^{-1}\circ F(\sigma)^{-1}=
\phi_{Q''}$, and $\phi_{Q''_2}=\phi_{Q''_1}$ for similar reasons.

(b) Once we have constructed $\tilde F$ on objects and have shown  
independence of $\tilde F(X)$ on the choice of strict cofibrant replacements,
we can easily define $\tilde F$ on morphisms. Indeed, let $f:X\to Y$ be
any morphism in $\cC$. Choose any strict cofibrant replacement 
$Q\stackrel p\to Y$; since acyclic fibrations are stable under base change,
$Q\times_YX\to X$ is still an acyclic fibration, so if we choose any
strict cofibrant replacement $Q'\to Q\times_YX$, the composite map
$p':Q'\to Q\times_YX\to X$ will be still a cofibrant replacement,
and the other composition $\bar f:Q'\to Q\times_YX\to Q$ has the property
$p\circ \bar f = f\circ p'$. Now $F(\bar f):F(Q')\cong\tilde F(X)\to
F(Q)\cong\tilde F(Y)$ is a natural candidate for $\tilde F(f)$. 
Its independence on the choice of $Q'\to Q\times_YX$ and $Q\to Y$ 
is shown as in~(a).
After this $\tilde F(\id_X)=\id$ is checked immediately (put $Q'=Q$),
and $\tilde F(g\circ f)=\tilde F(g)\circ\tilde F(f)$ for some
$X\stackrel f\to Y\stackrel g\to Z$ is checked as follows: we first choose 
some strict cofibrant replacement $Q\to Z$, then a strict cofibrant replacement
$Q'\to Q\times_ZY$, and finally a strict cofibrant replacement 
$Q''\to Q'\times_YX$. Then we get morphisms 
$\bar g:Q'\to Q$, $\bar f:Q''\to Q'$, such that $\tilde F(g\circ f)$,
$\tilde F(g)$ and $\tilde F(f)$ are canonically identified with
$F(\bar g\circ\bar f)$, $F(\bar g)$ and $F(\bar f)$:
\begin{equation}
\xymatrix{
Q''\ar[r]^{\bar f}\ar[d]&Q'\ar[r]^{\bar g}\ar[d]&Q\ar[d]\\
X\ar[r]^{f}&Y\ar[r]^{g}&Z}
\end{equation}

(c) Notice that the 2-out-of-3 axiom (PM3) together with (PM8) imply that
$\bar f:Q'\to Q$ is a weak equivalence whenever $f:X\to Y$ is one.
Since $Q'$ and $Q$ are cofibrant, we conclude that $\tilde F(f)=F(\bar f)$ 
is an isomorphism in~$\cB$, i.e.\ {\em $\tilde F$ transforms weak equivalences
into isomorphisms.} By definition, this means that $\tilde F$ can be uniquely
factorized through $\cC\stackrel\gamma\to\Ho\cC$. We denote the arising functor
$\Ho\cC\to\cB$ by $\dL F$. Thus $(\dL F)(\gamma X)=\tilde F(X)\cong F(Q)$
for any strict cofibrant replacement $Q\to X$.

(d) We need a natural transformation $\eta:\dL F\circ\gamma=\tilde F\to F$.
It can be defined as follows: if $Q\stackrel p\to X$ is any strict 
cofibrant replacement of~$X$, so that $\tilde F(X)\cong F(Q)$,
we put $\eta_X:=F(p)$. The independence of $\eta_X$ on the choice of~$Q\to X$
and its compatibility with morphisms follows immediately from (a) and~(b).
Notice that $\eta_X$ is an isomorphism for any cofibrant~$X$ since we can
take $Q=X$. We conclude that $(\dL F)(\gamma X)\cong(\dL F)(\gamma Q)\cong
F(Q)$ for any cofibrant replacement $Q\to X$ (not necessarily strict).

(e) Now we want to show that $(\dL F,\eta)$ is the left derived functor of~$F$.
Let $(G,\xi)$, $G:\Ho\cC\to\cB$, $\xi:G\circ\gamma\to F$ be another such pair.
We have to show existence and uniqueness of $\zeta:G\to\dL F$, such that 
$\xi=\eta\circ(\zeta\star\gamma)$. By the universal property of 
localizations 
$\Hom(G,\dL F)\cong\Hom(G\circ\gamma, \dL F\circ\gamma)$, so we have to
show existence and uniqueness of $\tilde\zeta:=\zeta\star\gamma
:G\circ\gamma\to\tilde F=\dL F\circ\gamma$, such that 
$\xi=\eta\circ\tilde\zeta$. 

(f) Uniqueness of $\tilde\zeta_X:G(\gamma X)\to\tilde F(X)$ 
is shown as follows. Choose any strict cofibrant replacement
$Q\stackrel p\to X$ and consider the following commutative diagram:

\begin{equation}\label{eq:constr.of.zeta}
\xymatrix{
G(\gamma Q)\ar[r]_{\tilde\zeta_Q}\ar@/^1pc/[rr]^{\xi_Q}
\ar[d]_{G(\gamma(p))}^{\sim}&
\tilde F(Q)\ar[r]_{\eta_Q}^{\sim}\ar[d]_{\tilde F(p)}^{\sim}&
F(Q)\ar[d]^{F(p)}\\
G(\gamma X)\ar[r]^{\tilde\zeta_X}\ar@/_1pc/[rr]_{\xi_X}&
\tilde F(X)\ar[r]^{\eta_X}&F(X)}
\end{equation}
We see that $\tilde\zeta_X=\tilde F(p)\circ\eta_Q^{-1}\circ\xi_Q
\circ G(\gamma(p))^{-1}$, i.e.\ it is completely determined by~$\xi$.

(g) Now it remains to show existence of $\tilde\zeta$, i.e.\ we need to show 
that the morphisms $\tilde\zeta_X$ defined as above do not depend on the 
choice of $Q\stackrel p\to X$ and are compatible with morphisms in~$\cC$.
Independence on choice of~$Q$ is shown as in (a): we reduce to the case of
$Q''\stackrel{\sigma}\to Q\stackrel p\to X$, with $\sigma$ an acyclic 
fibration between cofibrant objects, and a $3\times 3$-diagram similar to
\eqref{eq:constr.of.zeta}, with one extra horizontal line corresponding 
to~$Q''$ shows the required independence. As to the compatibility of
$\tilde\zeta$ with morphisms $f:X\to Y$ in~$\cC$, it is checked by
choosing $Q\to Y$, $Q'\to X$ and $\bar f:Q'\to Q$ as in~(b), q.e.d.
\end{Proof}

\nxsubpoint (Construction of $\bar Q$.)
Notice that the first part of the above proof, namely, construction of
$\tilde F:\cC\to\cB$ and $\dL F:\Ho\cC\to\cB$ in (a)--(c), 
works for any functor 
$F:\cC_c\to\cB$ transforming weak equivalences into isomorphisms.
In particular, we can apply it to $\gamma_c:\cC_c\to\Ho\cC_c$, thus
obtaining a functor $\bar Q:\Ho\cC\to\Ho\cC_c$, such that 
$\bar Q(\gamma X)$ is canonically isomorphic to $\gamma_c Q$ for any
strict cofibrant replacement $Q\stackrel p\to X$. On the other hand, we have a 
natural functor $I:\Ho\cC_c\to\Ho\cC$ induced by the embedding
$\cC_c\to\cC$. Clearly, $\bar Q\circ I$ is canonically isomorphic to
$\Id_{\Ho\cC_c}$, since one can always take $Q=X$ while computing
$\bar Q I(\gamma_c X)=\bar Q(\gamma X)$ for a cofibrant~$X$. On the other 
hand, we have a canonical natural transformation $\eta:I\circ\bar Q\to
\Id_{\Ho\cC}$, given by $\eta_{\gamma X}:=\gamma(p):
\gamma Q=I\bar Q(\gamma X)\to\gamma X$ as in part (d) of the above proof.
Now it is immediate that this $\gamma(p)$ is an isomorphism in $\Ho\cC$,
$p$ being an acyclic fibration, hence $I\circ\bar Q$ is also canonically 
isomorphic to the identity functor, i.e.\ {\em 
$I$ and $\bar Q$ are adjoint equivalences between $\Ho\cC$ and $\Ho\cC_c$.}
By duality $\Ho\cC_f$ is also equivalent to $\Ho\cC$.

\nxsubpoint (Construction of $\bar R_c$.)
Notice that the above construction of equivalence $\bar Q:\Ho\cC_c\to\Ho\cC$
still works if we weaken the pseudomodel category axioms for~$\cC$ 
by requiring just existence of pullbacks of acyclic fibrations instead 
of (PM1). Since $\cC_c$ satisfies the dual of this condition
(it is stable under pushouts of acyclic cofibrations), we can apply the 
above reasoning to $\cC_c^0$, thus obtaining an equivalence of categories
$\bar R_c:\Ho\cC_c\to\Ho\cC_{cf}$, adjoint to the natural
``embedding'' $\Ho\cC_{cf}\to\Ho\cC_c$. We conclude that 
{\em similarly to the classical case of model categories,
$\Ho\cC\cong\Ho\cC_c\cong\Ho\cC_{cf}\cong\Ho\cC_f$ for any pseudomodel
category~$\cC$.}
\begin{LemmaD}\label{l:acof.cof.suff2}
Let $\cC$ be pseudomodel category satisfying the first half of (PM8+),
i.e.\ such that the acyclic cofibrations are exactly the cofibrations 
which are also weak equivalences. Suppose that $F:\cC_c\to\cB$ transforms 
acyclic cofibrations between objects of~$\cC_c$ into isomorphisms in~$\cB$.
Then $F$ transforms weak equivalences in $\cC_c$ into isomorphisms in~$\cB$.
In particular, if $F:\cC\to\cB$ transforms acyclic cofibrations between
cofibrant objects into isomorphisms, then~\ptref{th:ex.left.der.psmod}
is applicable.
\end{LemmaD}
\begin{Proof}
The statement is essentially that of~\ptref{l:acof.cof.suff}, but we 
need to replace our original argument by another one, 
attributed in~\cite{DS} to K.~Brown. This argument goes as follows.
Let $f:X\to Y$ be a weak equivalence in~$\cC_c$.
Consider the factorization of the ``cograph'' 
$\langle f,\id_Y\rangle:X\sqcup Y\to Y$ into a cofibration
$q:X\sqcup Y\to Z$ followed by an acyclic fibration
$p:Z\to Y$. The natural embeddings $i:X\to X\sqcup Y$, 
$j:Y\to X\sqcup Y$ are also cofibrant, $X$ and $Y$ being cofibrant,
and $pqi=f$, $pqj=\id_Y$ and $p$ are weak equivalences, hence 
$qi$ and $qj$ are both weak equivalences and cofibrations, hence
acyclic cofibrations. By assumption $F(qi)$ and $F(qj)$ are isomorphisms
in~$\cB$, as well as $F(pqj)=F(\id_Y)$, hence $F(p)$ and
$F(f)=F(pqi)$ are also isomorphisms, q.e.d.
\end{Proof}

\nxsubpoint (Computation of $\Ho\cC$ and $\Ho\cC_{cf}$ by means of 
mixed fraction calculus.)
Quillen remarks in \cite[1.1]{Quillen} that $\Ho\cC$ and $\Ho\cC_{cf}$ 
do not admit left or right fraction calculus, but can be computed 
by a mixture of both. We claim that a similar statement is true 
in any of $\Ho\cC$, $\Ho\cC_c$, $\Ho\cC_f$ and $\Ho\cC_{cf}$. Namely,
any morphism $\phi\in\Hom_{\Ho\cC}(\gamma X,\gamma Y)$ can be written 
in form $\gamma(i)^{-1}\circ\gamma(f)\circ\gamma(p)^{-1}$, where
$X'\stackrel p\to X$ is a strict cofibrant replacement, 
$Y\stackrel i\to Y'$ is a strict fibrant replacement,
and $f:X'\to Y'$ a morphism in $\cC$, and similarly for any of the smaller
categories. In order to show this one simply checks that the set of all
such $\phi$'s is closed under composition and contains all morphisms
from $\cC$ and the inverses of all weak equivalences of~$\cC$.
Next, two such expressions $\gamma(i_1)^{-1}\gamma(f_1)\gamma(p_1)^{-1}$ 
and $\gamma(i_2)^{-1}\gamma(f_2)\gamma(p_2)^{-1}$ are equal, where
$X\stackrel{p_k}\leftarrow X'_k\stackrel{f_k}\to Y'_k\stackrel{i_k}
\leftarrow Y$ are as above, iff one can find a strict cofibrant replacement
$X'\to X'_1\times_XX'_2$, a strict fibrant replacement
$Y'_1\coprod_YY'_2\to Y'$ and a morphism $f:X'\to Y'$ making the 
obvious diagram commutative.

This is almost all we can obtain without the lifting properties,
which in the case of model categories enable us to show that any morphism
in $\Ho\cC_{cf}\cong\Ho\cC$ comes from some morphism in $\cC_{cf}$.

\nxsubpoint\label{sp:ptwise.psms.sE}
(Pointwise pseudomodel structure on $s\cE$ and $\stks\stSETS_\cE$.)
Notice that, while the five distinguished sets of morphisms in $s\cE$
and $\stks\stSETS_\cE$ have been defined without any reference to points, 
we have checked that we do obtain a pseudomodel structure in this way
only when $\cE$ has enough points. But in this case we can replace 
acyclic cofibrations and cofibrations by (potentially larger) classes 
of {\em pointwise acyclic cofibrations\/} and {\em pointwise cofibrations\/}
(of course, a morphism $f:X\to Y$ in $s\cE$ is a pointwise cofibration
iff $f_p:X_p\to Y_p$ is a cofibration in $s\catSets$ for all points 
$p$ of~$\cE$). Enlarging the three other classes in this manner wouldn't
produce anything new because of~\ptref{prop:f.af.w.pointwise}.
In this way we obtain another pseudomodel structure on $s\cE$ and
$\stks\stSETS_\cE$, which will be called the {\em pointwise pseudomodel 
structure}. In fact, all axioms of pseudomodel structure but (PM5) 
can be checked now pointwise, and factorization (PM5),
and even the functorial factorization (MS5f), follows from our
previous results~\ptref{cor:funct.decomp.sSETS}. Furthermore, 
the classes of pointwise (acyclic) cofibrations are closed,
and we even have (PM8+). 

We see that in some respects this pointwise
pseudomodel structure is even better than the one we had before; 
it would be nice to obtain a description of pointwise (acyclic) cofibrations
that would be valid in topoi without enough points. For example, a positive 
answer to~\ptref{q:c.ac.ptwise?} would yield such a description, 
thus actually showing that our original pseudomodel structure on $s\cE$
is even better than we thought before.

In any case, the weak equivalences for these two structures are the same,
hence the corresponding homotopic categories $\Ho s\cE$ and derived functors 
also coincide, i.e.\ {\em we can freely use pointwise pseudomodel structures 
while computing derived functors.}

\nxpointtoc{Pseudomodel structure on $\stks\stMOD\sO_\cE$ and
$s\catMod\sO$}\label{p:psmod.str.sO-mod}
Let $(\cE,\sO)$ be a generalized ringed topos, i.e.\ 
$\sO$ be an arbitrary (inner) algebraic monad over a topos~$\cE$
(cf.~\ptref{sp:gen.ringsp.top} and~\ptref{p:algmon.over.topos}; 
notice that we don't need commutativity of $\sO$ here).
We want to construct 
a pseudomodel category structure on the category $s\catMod\sO$ 
of simplicial $\sO$-modules in~$\cE$, as well as a pseudomodel stack
structure on the flat stack $\stks\stMOD\sO_\cE$.

\nxsubpoint\label{sp:distcl.sO.mod} (Basic definitions.)
Since $L_\Sigma(I)$ and $L_\Sigma(J)$ are cofibrant generators of 
our standard model category structure on $s\catMod\Sigma$, for any 
algebraic monad~$\Sigma$, we are tempted to take the sets
$L_\sO(\cst I_\cE)$ and $L_\sO(\cst J_\cE)$ of morphisms inside
$(\stks\stMOD\sO_\cE)(e)=s\catMod\sO$ as the ``cofibrant generators''
for a pseudomodel structure on $\stks\stMOD\sO_\cE$. Here
$L_\sO$ denotes both the free $\sO$-module functor $\cE\to\catMod\sO$
and its simplicial extension $s\cE\to s\catMod\sO$.

More precisely:
\begin{itemize}
\item A morphism $f:X\to Y$ in $(\stks\stMOD\sO_\cE)(S)=s\catMod{\sO|_S}$
is a {\em fibration} (resp.\ {\em acyclic fibration}) iff it has the local RLP
with respect to all morphisms from $L_\sO(\cst J)$ 
(resp.\ $L_\sO(\cst I)$).
\item A morphism $f:X\to Y$ in $s\catMod{\sO|_S}$ is a {\em cofibration\/}
(resp.\ {\em acyclic cofibration}) iff it belongs to the closure 
$\Cl L_\sO(\cst I)$ (resp.\ $\Cl L_\sO(\cst J)$) in the sense 
of~\ptref{def:closure.class}.
\item Finally, a morphism $f:X\to Y$ in $s\catMod{\sO|_S}$ is a
{\em weak equivalence\/} iff it can be factorized into an acyclic cofibration
followed by an acyclic fibration.
\end{itemize} 

\nxsubpoint
Notice that adjoint functors $L_\sO$ and $\Gamma_\sO$ between
$\catMod\sO$ and $\cE$ extend to adjoint cartesian functors 
${\bm L}_\sO$ and ${\bm\Gamma}_\sO$ between corresponding stacks
$\stMOD\sO_\cE$ and $\stSETS_\cE$, and the same applies to their simplicial
extensions $\stks\stMOD\sO_\cE\leftrightarrows\stks\stSETS_\cE$; resulting
adjoint cartesian functors will be also denoted by ${\bm L}_\sO$ and
${\bm\Gamma}_\sO$, or even by $L_\sO$ and $\Gamma_\sO$, when no confusion
can arise.

In particular, we obtain the following formulas:
\begin{align}
\Hom_{s\catMod\sO}(L_\sO(A),X)\cong&\Hom_{s\cE}(A,\Gamma_\sO(X))\\
\label{eq:innadj.LG}
\iHom_{\stks\stMOD\sO_\cE}(L_\sO(A),X)\cong&\iHom_{\stks\stSETS_\cE}
(A,\Gamma_\sO(X))
\end{align}

\nxsubpoint\label{sp:f.af.smod.O}
(Fibrations and acyclic fibrations in $\catMod\sO$.)
The second of the above formulas immediately implies that 
some $p:X\to Y$ in $s\catMod\sO$ has the local RLP with respect to
$L_\sO(i)$ for some $i:A\to B$ in $s\cE$ iff $\Gamma_\sO(p)$ has the local
RLP with respect to $i$. Considering here all $i$ from 
$\cst I_\cE$ (resp.\ $\cst J_\cE$), we conclude that
{\em $p:X\to Y$ is an acyclic fibration (resp.\ fibration) in
$s\catMod\sO$ iff $\Gamma_\sO(p)$ is an acyclic fibration 
(resp.\ fibration) in $s\cE$.} Furthermore, the latter condition can be 
expressed in terms of epimorphicity of certain morphisms between 
finite projective limits of components of $\Gamma_\sO(X)$ and 
$\Gamma_\sO(Y)$. Now notice that $\Gamma_\sO:\catMod\sO\to\cE$ is 
left exact, and $\Gamma_\sO(f)$ is an epimorphism in $\cE$ iff
$f$ is a strict epimorphism in $\catMod\sO$. We conclude that
{\em $p:X\to Y$ is an acyclic fibration (resp.\ fibration) in
$s\catMod\sO$ iff certain morphisms between finite projective limits
of components of $X$ and~$Y$ are strictly epimorphic}, and that
{\em $p:X\to Y$ is an acyclic fibration in $s\catMod\sO$ iff
all $X_n\to Y_n\times_{(\cosk_{n-1}Y)_n}(\cosk_{n-1}X)_n$ are strict 
epimorphisms in $\catMod\sO$.}

\begin{PropD}\label{prop:stab.distcl.topos.Opb}
Let $f:\cE'\to\cE$ be a morphism of topoi, $\sO$ a generalized ring in~$\cE$,
$\sO':=f^*\sO$ its pullback to $\cE'$. Then 
all five distinguished classes in $s\catMod\sO$ are stable under
the induced pullback functor $f^*:s\catMod\sO\to s\catMod{\sO'}$.
\end{PropD}
\begin{Proof}
The proof is completely similar to that of~\ptref{prop:stab.distcl.topos.pb}.
(a) Stability of acyclic cofibrations and cofibrations is shown as in 
{\em loc.cit.}: we consider the class $\propP$ of morphisms in fibers
of $\stks\stMOD\sO_\cE$, consisting of those $\xi:X\to Y$ in
$s\catMod{\sO|_S}$ which become cofibrations (resp.\ acyclic cofibrations)
in $s\catMod{\sO'|_{f^*S}}$ after application of $f_S^*$, easily check 
that $\propP$ is closed and that it contains $L_\sO(\cst I_\cE)$
(resp.\ $L_\sO(\cst J_\cE)$) since $f^*\circ L_\sO=L_{\sO'}\circ f^*$,
hence $\propP$ has to contain $\Cl L_\sO(\cst I_\cE)$ (resp.\dots),
i.e.\ all cofibrations (resp.\ acyclic cofibrations) of $\stks\stMOD\sO_\cE$.
(b) Stability of acyclic fibrations and fibrations immediately follows from
{\em loc.cit.}\ and \ptref{sp:f.af.smod.O}, since 
$\Gamma_{\sO'}\circ f^*=f^*\circ\Gamma_{\sO}$.
(c) Stability of weak equivalences follows from that of acyclic cofibrations
and acyclic fibrations.
\end{Proof}

\nxsubpoint (Simple cases.)
(a) When $\cE=\catSets$, hence $\sO$ is just a (constant) algebraic monad,
then the five distinguished classes in $s\catMod\sO$ defined 
in~\ptref{sp:distcl.sO.mod} coincide with those given by the
model category structure on $s\catMod\sO$ defined in~\ptref{sp:cat.simpl.mod}.
Indeed, this is evident for fibrations and acyclic fibrations; 
as to cofibrations and acyclic cofibrations with respect to 
this model structure, they clearly constitute closed sets of morphisms 
containing $L_\sO(I)$ (resp.\ $L_\sO(J)$), 
hence also the closures of these sets.
The opposite inclusion is shown by Quillen's small object argument, 
which establishes that any cofibration (resp.\ acyclic cofibration) is
a retract of a morphism of $L_\sO(I)$ (resp.\ $L_\sO(J)$), (CM4) being 
fulfilled in $s\catMod\sO$, hence lies itself in $L_\sO(I)$ (resp.\dots).
Finally, the case of weak equivalences follows from those of 
acyclic cofibrations and acyclic fibrations.

(b) On the other hand, we can take $\sO=\Fempty$, considered as a constant
algebraic monad over any topos~$\cE$. Then $s\catMod\sO=s\cE$,
$L_\sO(\cst I_\cE)=\cst I_\cE$, $L_\sO(\cst J_\cE)=\cst J_\cE$, and
we recover the pseudomodel structure on $s\cE$ considered
before in~\ptref{p:psmod.str.sO-mod}.

\begin{PropD}\label{prop:f.af.weq.sO-mod}
If $f:X\to Y$ belongs to any of five distinguished classes of morphisms
in $s\catMod\sO$, then $f_p:X_p\to Y_p$ belongs to the same class in
$s\catMod{\sO_p}$ for any point $p$ of~$\cE$. Conversely, if $\cE$
has enough points, and if all $f_p:X_p\to Y_p$ are fibrations 
(resp.\ acyclic fibrations, weak equivalences) in $s\catMod{\sO_p}$,
then $f:X\to Y$ is itself a fibration (resp.\dots) in $s\catMod\sO$.
\end{PropD}
\begin{Proof} Completely similar to that of~\ptref{prop:f.af.w.pointwise}.
\end{Proof}

\begin{CorD}\label{cor:indeed.psmod.sO-mod}
If $\cE$ has enough points, then \ptref{sp:distcl.sO.mod} defines a 
pseudomodel structure on $s\catMod\sO$ and $\stks\stMOD\sO_\cE$,
satisfying the second half of (PM8+), the functorial factorization (MS5f),
and having closed classes of cofibrations and acyclic cofibrations.
\end{CorD}
\begin{Proof}
All statements are deduced pointwise from corresponding statements
for $s\catMod\sO_p$, valid since $s\catMod\sO_p$ is known to be 
a model category by~\ptref{sp:cat.simpl.mod}. The only exception is
the (functorial) factorization (MS5f), which can be again shown by 
``local'' Quillen's small object argument \ptref{prop:Quillen.fact}, 
applied to $L_\sO(\cst I)$ and $L_\sO(\cst J)$, provided we show that 
these morphisms have locally small sources in $s\catMod\sO$. But this is
immediate from the fact that $\Gamma_\sO$ commutes with filtered 
inductive limits, combined with~\eqref{eq:innadj.LG}, 
which actually implies that $L_\sO$ transforms
locally small objects of $s\cE$ into locally small objects of~$s\catMod\sO$.
\end{Proof}

\nxsubpoint (Pointwise pseudomodel structure on $s\catMod\sO$.)
When $\cE$ has enough points, we can replace cofibrations and 
acyclic cofibrations with (potentially larger) classes of pointwise 
(acyclic) cofibrations, similarly to what we did in~\ptref{sp:ptwise.psms.sE}.
In this way we obtain a ``better'' pseudomodel structure on $s\catMod\sO$,
which has closed classes of cofibrations and acyclic cofibrations,
and satisfies (MS5f) and (PM8+).  We call it the {\em pointwise pseudomodel
structure on~$s\catMod\sO$}.
Of course, this pseudomodel structure has the same weak equivalences 
as before, hence it gives rise to the same homotopic category
$\cD^{\leq0}(\cE,\sO)=\Ho s\catMod\sO$ and same derived functors.

\nxsubpoint (Derived category of $(\cE,\sO)$.)
We denote by $\cD^{\leq0}(\cE,\sO)$ or by $\cD^{\leq0}(\sO)$ the
homotopic category $\Ho s\catMod\sO$, i.e.\ the localization of
$s\catMod\sO$ with respect to weak equivalences. If $\sO$ admits a zero
(i.e.\ a central constant), we have well-defined suspension functors
$\Sigma$ on $(s\catMod\sO)_c$ and $\dL\Sigma$ on $\cD^{\leq0}(\cE,\sO)$;
formally inverting this functor in the usual manner
(cf.~\ptref{sp:stab.ht.cat}), we obtain the
{\em stable\/} homotopic category $\cD^-(\cE,\sO)$ and a natural functor
$\cD^{\leq0}(\cE,\sO)\to\cD^-(\cE,\sO)$.

\nxsubpoint\label{sp:der.sO-mod.addc} (Additive case.)
Now suppose that $\sO$ is additive, i.e.\ it admits a central zero and
addition. Then $\sO$ is nothing else than a classical ring in $\cE$,
i.e.\ $(\cE,\sO)$ is a ringed topos in the classical sense.
In this case $\catMod\sO$ is additive, hence the Dold--Kan correspondence
(cf.~\ptref{th:DK}) establishes an equivalence between the category 
$s\catMod\sO$ of simplicial $\sO$-modules and the category
$\Ch^{\geq0}(\catMod\sO)$ of non-negative chain complexes of $\sO$-modules.
Furthermore, weak equivalences in $s\catMod\sO$ correspond exactly
to quasi-isomorphisms in $\Ch^{\geq0}(\catMod\sO)$, at least if
$\cE$ has enough points: indeed, in this case both weak equivalences and 
quasi-isomorphisms can be checked pointwise, and we are reduced to
check this statement over $\catSets$, where it is classical. We conclude that
{\em $\Ho s\catMod\sO$ is indeed equivalent to the full subcategory
$\cD^{\leq0}(\cE,\sO)$ of $\cD(\cE,\sO):=\cD(\catMod\sO)$,} hence
{\em the corresponding stable homotopic category $\cD^-(\cE,\sO)$ is 
equivalent to $\cD^-(\catMod\sO)$.}

\nxsubpoint\label{sp:cofrepl.add.sO-mod} 
(Cofibrant replacements in the additive case.)
Furthermore, one checks immediately that the class of injective chain maps
with flat cokernels is closed in $\Ch^{\geq0}(\catMod\sO)$ 
(more precisely, in corresponding stack), and it contains 
complexes corresponding to simplicial objects $L_\sO(\cst I_\cE)$. This
implies that {\em cofibrations in $s\catMod\sO$ correspond to (some set of)
injective chain maps in $\Ch^{\geq0}(\catMod\sO)$ with flat cokernels}.
In particular, cofibrant replacements $Q\to X$ in $s\catMod\sO$ yield
flat resolutions $N(Q)\to N(X)$ in $\Ch(\catMod\sO)$. This means that
deriving functors by means of cofibrant replacements corresponds to deriving
functors by means of flat resolutions. 

We'll use this remark later to
compare our derived tensor products and pullbacks with their classical
additive counterparts.

\begin{ThD}\label{th:ex.lder.pb} (Derived pullbacks.)
Let $f=(\phi,\theta):(\cE',\sO')\to(\cE,\sO)$ be a morphism of generalized
ringed topoi, i.e.\ $\phi:\cE'\to\cE$ is a morphism of topoi and
$\theta:\sO\to\phi_*\sO'$ is a morphism of algebraic inner monads on~$\cE$.
Suppose that $\cE'$ and $\cE$ have enough points. Then the pullback 
functor $f^*:s\catMod\sO\to s\catMod{\sO'}$ admits a left derived
$\dL f^*:\cD^{\leq0}(\cE,\sO)\to\cD^{\leq0}(\cE',\sO')$, which can be computed
by means of cofibrant replacements. When both $\sO$ and $\sO'$ are additive,
$\dL f^*$ is identified with its classical counterpart via Dold--Kan
correspondence.
\end{ThD}
\begin{Proof}
Let's endow $s\catMod\sO$ and $s\catMod{\sO'}$ with their 
pointwise pseudomodel structures. According to~\ptref{th:ex.left.der.psmod},
all we have to check is that $f^*$ preserves weak equivalences between
pointwise cofibrant objects. But for any such $u:A\to B$ in $s\catMod\sO$
and any point $p$ of $\cE'$ we have $f^*(u)_p=(\theta^\sharp_p)^*(u_{f(p)})$,
where $(\theta^\sharp_p)^*$ denotes the scalar extension with respect to
algebraic monad homomorphism $\theta^\sharp_p:\sO_{f(p)}\to\sO'_p$.
Now everything follows from the fact that weak equivalences can be checked 
pointwise and that $(\theta^\sharp_p)^*$ preserves weak equivalences between
cofibrant objects by~\ptref{prop:ex.Lder.basech}.
\end{Proof}

\nxsubpoint (Functoriality of derived pullbacks.)
One can check directly that $f^*$ preserves cofibrations and acyclic 
cofibrations, using the same trick as in the proof 
of~\ptref{prop:stab.distcl.topos.Opb}. In particular, it preserves cofibrant 
objects and pointwise cofibrant objects. This immediately implies that
whenever $g:(\cE'',\sO'')\to(\cE',\sO')$ is another morphism of 
generalized ringed topoi with enough points, we have
$\dL(fg)^*=\dL g^*\circ\dL f^*$.

\nxsubpoint (Pullbacks with respect to flat morphisms.)
Suppose that $f=(\phi,\theta):(\cE',\sO')\to(\cE,\sO)$ is {\em flat},
i.e.\ $\theta^\sharp:\phi^*\sO\to\sO'$ is a flat extension of algebraic monads
in~$\cE'$, i.e.\ $(\theta^\sharp)^*:\catMod{\phi^*\sO}\to\catMod{\sO'}$ is 
exact, and in particular preserves finite projective limits and strict 
epimorphisms. Since the property of some $p:X\to Y$ to be a fibration
or an acyclic fibration in $s\catMod\sO$ can be expressed in terms of 
such limits and strict epimorphisms, we see that {\em $f^*$ preserves 
fibrations and acyclic fibrations for any flat~$f$}. We have already the 
same property for (acyclic) cofibrations, hence for weak equivalences as well.
Now a functor $f^*$ that preserves weak equivalences can be derived in the 
trivial manner, usually expressed by the formula $\dL f^*=f^*$.

\nxpointtoc{Derived local tensor products}
Henceforth we assume $(\cE,\sO)$ to be a generalized 
{\em commutatively\/} ringed topos.
Then we have a tensor product functor $\otimes=\otimes_\sO:
\catMod\sO\times\catMod\sO\to\catMod\sO$, which canonically extends 
to a cartesian functor $\otimes:\stMOD\sO_\cE\times_\cE\stMOD\sO_\cE
\to\stMOD\sO_\cE$. Moreover, the simplicial extensions of these functors 
define an ACU $\otimes$-structure on category $s\catMod\sO$ and on
stack $\stks\stMOD\sO_\cE$. We want to prove the following statement:

\begin{ThD}\label{th:ex.loc.der.tensprod}
Suppose that $\cE$ has enough points. 
Then functor $\otimes:s\catMod\sO\times s\catMod\sO\to s\catMod\sO$ admits a 
left derived functor $\Lotimes=\dL\otimes$, which can be 
computed with the aid of cofibrant replacements in each variable,
i.e.\ $X\Lotimes X'\cong Q\otimes Q'$, for any cofibrant replacements
$Q\to X$ and $Q'\to X'$ in $s\catMod\sO$. Moreover, $\Lotimes$ defines an
ACU $\otimes$-structure on $\Ho s\catMod\sO=\cD^{\leq0}(\cE,\sO)$.
\end{ThD}
\begin{Proof}
Of course, the ``product'' pseudomodel structure on 
a product of pseudomodel categories $\cC_1\times\cC_2$
or on a fibered product of pseudomodel stacks $\cC_1\times_\cE\cC_2$ 
is defined as in~\ptref{sp:prod.mod.cat}, i.e.\ $f=(f_1,f_2):(X_1,X_2)\to
(Y_1,Y_2)$ belongs to one of five distinguished classes of morphisms
in $\cC_1\times\cC_2$ iff both $f_i:X_i\to Y_i$, $i=1,2$, 
belong to the same distinguished class in~$\cC_i$. Therefore, everything
will follow from~\ptref{th:ex.left.der.psmod}, provided we manage to
show that $X\otimes Y$ is cofibrant when both $X$ and $Y$ are cofibrant
(in $s\catMod\sO$), and that a tensor product of weak equivalences between
cofibrant objects is again a weak equivalence. Now this is very easily done
for the pointwise pseudomodel structure on $s\catMod\sO$, for which both 
statements follow immediately from our former results 
of~\ptref{p:der.tensprod}. We know that both pseudomodel structures on
$s\catMod\sO$ have same weak equivalences, hence same homotopic categories
and derived functors, so we obtain all statements of the theorem for the
``smaller'' pseudomodel structure as well, especially if we take into account
that cofibrant replacements for the ``smaller''
pseudomodel structure on $s\catMod\sO$ are automatically 
pointwise cofibrant replacements, so they can be used to compute~$\Lotimes$.
\end{Proof}

\nxsubpoint
Notice that we didn't actually show that the tensor product of cofibrant
simplicial $\sO$-modules is cofibrant, since the corresponding statement
for pointwise cofibrant objects suffices for 
proving~\ptref{th:ex.loc.der.tensprod}. However, we'll show in a moment
that this statement is true, and, moreover, that a tensor product
of acyclic cofibrations between cofibrant objects is again an
acyclic cofibration.

\nxsubpoint (Compatible $\otimes$-structures and $\otimes$-actions.)
Suppose that $\cC_1$, $\cC_2$ and $\cD$ are three pseudomodel categories, and
$\otimes:\cC_1\times\cC_2\to\cD$ is a bifunctor, commuting with inductive
limits in each argument. We say that $\otimes$ is {\em compatible
with given pseudomodel structures\/} if the axiom (TM)
of~\ptref{def:comp.tens.mod.str} holds. This definition is in particular 
applicable to $\otimes$-structures on pseudomodel categories and
external $\otimes$-actions of one pseudomodel category on another.
Moreover, it can be easily extended to the case of pseudomodel stacks
$\cC_1$, $\cC_2$ and $\cD$ over a topos~$\cE$ and a cartesian bifunctor
$\otimes:\cC_1\times_\cE\cC_2\to\cD$: in this situation we just require
(TM) to hold fiberwise.

Notice that whenever $\otimes$ is a compatible ACU $\otimes$-structure on 
a pseudomodel category, we show as in~\ptref{sp:der.comp.tensprod} 
that $A\otimes X$ is cofibrant whenever
$A$ and $X$ are cofibrant, and that a tensor product of acyclic cofibrations
between cofibrant objects is an acyclic cofibration.

\begin{PropD}\label{prop:comp.tens.sO-mod}
Tensor product on $s\catMod\sO$ and on $\stks\stMOD\sO_\cE$ 
is compatible with their pseudomodel
structure. This applies in particular to $\otimes_{\Fempty}=\times$ on
$s\cE$.
\end{PropD}
{\bf Proof.} 
When $\cE$ has enough points, this compatibility for the pointwise pseudomodel 
structure is immediate. Let us show the statement for the ``smaller''
classes of cofibrations and acyclic cofibrations, given by 
$\Cl L_\sO(\cst I_\cE)$, $\Cl L_\sO(\cst J_\cE)$, without assuming 
$\cE$ to have enough points (strictly speaking, we haven't shown that
we obtain a pseudomodel structure in this case, but this isn't necessary
to discuss the validity of (TM)). 

(a) Let us denote by $Z(i,s)$ the source of
$i\boxe s$, i.e.\ the pushout of $i\otimes\id_K:A\otimes K\to B\otimes K$
and $\id_A\otimes s:A\otimes K\to A\otimes L$. Then 
$Z$ is a bifunctor $\Ar\cC\times\Ar\cC\to\cC$, where $\cC$ is either
$s\catMod\sO$ or a fiber of $\stks\stMOD\sO_\cE$, 
$B\otimes L$ can be identified with $Z(\id_B,s)$ and $Z(i,\id_L)$, and
$i\boxe s:Z(i,s)\to Z(\id_B,s)$ equals 
$Z((i,\id_B),(\id_K,\id_L))=:Z(i,\id_B;s)$.

(b) Fix any cofibration $s$ and consider the class $\propP$, 
consisting of morphisms $i$, such that $i\boxe s$ is a cofibration
(resp.\ an acyclic cofibration). According to Lemma~\ptref{l:cl.box.div.class}
below, this class $\propP$ is closed, so it contains all
cofibrations (resp.\ acyclic cofibrations) iff it
contains $L_\sO(\cst I)$ (resp.\ $L_\sO(\cst J$)). 
In other words, it suffices to
check (TM) for $i$ from $L_\sO(\cst I)$ (resp.\ $L_\sO(\cst J)$).

(c) Interchanging the arguments in the reasoning of (b), we see that 
we can assume that $s$ lies in $L_\sO(\cst I)$ as well, i.e.\ we have to check
that $i\boxe s$ is a cofibration (resp.\ acyclic cofibration) whenever
$i$ lies in $L_\sO(\cst I)$ (resp.\ $L_\sO(\cst J)$) and $s$ lies in 
$L_\sO(\cst I)$. 

(d) Now notice that $L_\sO:s\cE\to s\catMod\sO$ commutes with $\boxe$
and preserves cofibrations and acyclic cofibrations, hence we are reduced to
proving the statement in $s\cE$ for $s\in\cst I$, $i\in\cst I$ or $\cst J$.
Furthermore, $q^*:s\catSets\to s\cE$ also commutes with $\boxe$ and preserves 
cofibrations and acyclic cofibrations, so we are reduced to proving a 
statement about $I$ and $J$ in $s\catSets$, which follows immediately 
from the compatibility of $s\catSets$ with 
its tensor product $\otimes=\times$, q.e.d.

\begin{LemmaD}\label{l:cl.box.div.class}
Let $\propP$ be a closed class in category $\cC=s\catMod\sO$ or in stack
$\cC=\stks\stMOD\sO_\cE$, and $s:K\to L$ be any morphism in this category or 
in the final fiber of this stack. 
Then the class $\propQ:=(\propP:s)$ of morphisms
$i$, such that $i\boxe s$ lies in $\propP$, is itself closed.
\end{LemmaD}
\begin{Proof}
We have to check conditions of~\ptref{def:closed.class} for 
$\propQ=(\propP:s)$ one by one.

(a) Condition 1) (localness of $\propQ$) is evident, as well as condition 
2) (all isomorphisms do lie in $\propQ$), 3) (if $i'$ is a pushout of
$i$, then $i'\boxe s$ is easily seen to be a pushout of $i\boxe s$),
6) (if $i'$ is a retract of $i$, then $i'\boxe s$ is a retract of $i\boxe s$)
and 7) (since $-\boxe s$ commutes with arbitrary direct sums). 
So only conditions 4), 5) and 8) remain.

(b) Let us check condition 4), i.e.\ stability under composition. Let
$A\stackrel i\to B\stackrel j\to C$ be a composable sequence of morphisms 
in~$\propQ$.
Then $ji\boxe s=Z(ji,\id_C;s):Z(ji,s)\to Z(\id_C,s)$ can be decomposed into
$Z(i,\id_C;s):Z(ji,s)\to Z(j,s)$, followed by $Z(j,\id_C;s)=j\boxe s:
Z(j,s)\to Z(\id_C,s)$. Now notice that $Z(i,\id_C;s)$ is a pushout
of $Z(i,\id_B;s)=i\boxe s:Z(i,s)\to Z(\id_B,s)$, just because $Z(-,s)$ is 
right exact and the following square is cocartesian in $\Ar\cC$:
\begin{equation}
\xymatrix@C+10pt{
i\ar[r]^{(i,\id_B)}\ar[d]_{(\id_A,j)}&\id_B\ar[d]^{(\id_B,j)}\\
ji\ar[r]^{(i,\id_C)}&j}
\end{equation}
We conclude that $ji\boxe s$ is a composite of $j\boxe s$ and a pushout of
$i\boxe s$, hence it lies in $\propP$ whenever 
both $j\boxe s$ and $i\boxe s$ do.

(c) Now let us check 5), i.e.\ stability under sequential composition.
Let $A_0\stackrel{i_0}\to A_1\stackrel{i_1}\to A_2\to\ldots$ be a composable 
sequence of morphisms from~$\propQ$. Put $B:=A_\infty:=\injlim_n A_n$, and
denote by $j_n:A_n\to B$ the natural embeddings, and by 
$k_n:A_0\to A_n$ the compositions $i_{n-1}\cdots i_1i_0$, so as to have
$j_nk_n=j_0$. We must check that $j_0:A_0\to B$ lies in $\propQ$, i.e.\
that $j_0\boxe s=Z(j_0,\id_B;s)$ lies in $\propP$. Notice that
$(j_0,\id_B)$ is the inductive limit of $(k_n,\id_B)$, i.e.\ 
the sequential composition of
$j_0\stackrel{(i_0,\id_B)}\longto j_1\stackrel{(i_1,\id_B)}\longto
j_2\to\cdots$. Since $Z(-;s)$ commutes with arbitrary inductive limits, we see
that $j_0\boxe s$ is the sequential composition of $Z(i_n,\id_B;s):
Z(j_n,s)\to Z(j_{n+1},s)=Z(i_nj_n,s)$. But we have seen in (b) that each
$Z(i_n,\id_B;s)$ is a pushout of $Z(i_n,\id_{A_{n+1}};s)=i_n\boxe s$,
hence lies in $\propP$, hence the same is true for their sequential composition
$j_0\boxe s$, $\propP$ being closed.

(d) It remains to check 8), i.e.\ the stability under $\phi_!$ in 
the stack case. This stability follows from the canonical isomorphism
$(\phi_!A)\otimes K\cong\phi_!(A\otimes\phi^*K)$, which enables us to
identify $\phi_!(i)\boxe s$ with $\phi_!(i\boxe\phi^*(s))$, q.e.d.
\end{Proof}

\nxsubpoint
One checks exactly in the same way as above that the external $\otimes$-action
of $s\cE$ on $s\catMod\sO$ is compatible with pseudomodel structures
and pointwise pseudomodel structures, and that this $\otimes$-action can 
be derived with the aid of cofibrant replacements (or even pointwise 
cofibrant replacements).

\nxsubpoint\label{sp:fancy.lang}
(Fancy description of this method of proof as ``devissage''.)
We have already used several times the following method to prove 
stability of cofibrations or acyclic cofibrations under some functor $F$:
we show that the class of morphisms which become (acyclic) cofibrations
after applying~$F$ is closed, checking conditions of \ptref{def:closed.class}
one by one, and then check that the generators $L_\sO(\cst I)$ 
(resp.\ $L_\sO(\cst J)$) belong to this class; this is usually done by
observing that functor $F$ essentially ``commutes'' with $L_\sO\circ q^*:
s\catSets\to s\catMod\sO_\cE$, so we are finally reduced to checking the 
statetement for morphisms from $I$ (resp.\ $J$) in $s\catSets$.

Up to now we have applied this method to prove stability of cofibrations
and acyclic cofibrations under topos pullbacks and base change,
and we have also applied it to bifunctor $\otimes$ just now 
in~\ptref{prop:comp.tens.sO-mod}, where we have been able to use separately 
the above reasoning first in one variable (fixing the other), and then 
in the second, thus reducing (TM) to checking that
$i\boxe s$ is a cofibration for any $i$, $s\in L_\sO(\cst I)$, and 
an acyclic cofibration for $i\in L_\sO(\cst I)$, $s\in L_\sO(\cst J)$,
something that we can check in $s\catSets$ and apply $L_\sO\circ q^*$ 
afterwards. However, we applied this reasoning to 
an ``additive'' or right exact functor~$F$ 
(right exact in each variable in the case of $\otimes$) so far.
Since we are going to apply it shortly to the much more complicated case
of (definitely non-additive) symmetric powers, we'd like to give a fancy
description of this method of proof.

Namely, one can imagine that cofibrations are just monomorphisms
$i:A\to B$ in some abelian category (e.g.\ modules over some classical
commutative ring~$\Lambda$, or complexes of such modules), 
having a projective cokernel~$P=B/A$
(one can check that such monomorphisms are just the closure of
one-element set $\{0\to\Lambda\}$ in $\catMod\Lambda$),
and acyclic cofibrations are cofibrations with ``negligible'' cokernel
$P\approx 0$ (something like a complex of projectives with trivial cohomology).
We might think of $A\stackrel i\to B$ as a ``presentation'' of $P$ in form
$B/A$. If $i':A'\to B'$ is a pushout of $i:A\to B$, then $B'/A'\cong B/A$,
i.e.\ a pushout just provides another presentation of the same object. 
Next, when we decompose a cofibration $A\to B$ into a composition of
``simpler'' cofibrations $A=A_0\to A_1\to\cdots\to A_n=B$, this essentially 
means that we define on $P=B/A$ an increasing filtration $F_iP=A_i/A_0$ with
``simpler'' quotients $P_i=A_i/A_{i-1}$, and try to deduce some property
of $P$ from similar properties of $P_i$, i.e.\ do some sort of 
{\em devissage.} When we consider direct sums of morphisms, they correspond
to direct sums of cokernels, and retracts of $i:A\to B$ correspond to 
direct factors of $P=B/A$. Therefore, the above method of proving 
statements about cofibrations or acyclic cofibrations can be thought of
as some sort of devissage, where the role of simple objects is played by the
``cokernels'' of standard generators from $L_\sO(\cst I)$ or $L_\sO(\cst J)$.

For example, if $i:A\to B$ and $s:K\to L$ are two such ``cofibrations'', with
projective cokernels $A/B=P$ and $K/L=Q$, 
then $i\boxe s:\ldots\to B\otimes L$ of (TM) is injective with cokernel 
$P\otimes Q$, i.e.\ $\boxe$ is just ``the tensor product of projective modules
in terms of their presentations'', and (TM) can be interpreted as saying
that $P\otimes Q$ is projective whenever $P$ and $Q$ are projective, and
that $P\otimes Q\approx 0$ if in addition either $P\approx0$ or $Q\approx0$,
a very natural statement indeed. Moreover, in classical situation we would
probably show such statements by writing both $P$ and $Q$ as direct factors
of free modules, and decomposing these free modules further into direct sums
of free modules of rank 1; in our situation we reduce to the case where 
$P$ and $Q$ are ``cokernels'' of morphisms from $L_\sO(\cst I)$ or
$L_\sO(\cst J)$ instead.

\nxsubpoint
Once we know that derived local tensor product $\Lotimes=\Lotimes_\sO$ exists
and can be computed by means of cofibrant replacements, we deduce the
usual properties (associativity, commutativity, compatibility with pullbacks
and base change etc.)\ in the same way as before. In particular, we have
$X\Lotimes K\cong X\Lotimes_\sO\dL L_\sO(K)$, for any 
$X\in\Ob\cD^{\leq0}(\cE,\sO)$ and $K\in\Ob\Ho s\cE=
\Ob\cD^{\leq0}(\cE,\Fempty)$. Notice, however, that we cannot derive inner
Homs as we did in~\ptref{sp:comp.of.iHoms}, because inner Homs don't
commute with fiber functors $p^*$, so we cannot check required properties 
pointwise. Therefore, we don't obtain any inner Homs for the ACU
$\otimes$-structure on $\cD^{\leq0}(\cE,\sO)$ given by~$\Lotimes$.

\nxsubpoint
If $\sO$ admits a zero, then $\Lotimes$ commutes with the suspension functor 
in each variable, hence induces an ACU $\otimes$-structure $\Lotimes$
on the corresponding stable category $\cD^-(\cE,\sO)$. All formulas and 
properties mentioned above extend to the stable case as well.

\nxsubpoint
When $\sO$ is additive, cofibrant replacements in $s\catMod\sO$ 
correspond via Dold--Kan to some flat resolutions in the
category of chain complexes 
$\Ch^{\geq0}(\catMod\sO)$ (cf.~\ptref{sp:cofrepl.add.sO-mod}), 
hence our derived local 
tensor product corresponds to its classical counterpart, once we
take into account Eilenberg--Zilber theorem (cf.~\ptref{th:EZ}).

\nxpointtoc{Derived symmetric powers}
This subsection is dedicated to the proof of existence of
derived symmetric powers $\dL S^n=\dL S^n_\sO:\cD^{\leq0}(\cE,\sO)\to
\cD^{\leq0}(\cE,\sO)$, for any generalized commutatively ringed topos
$(\cE,\sO)$ with sufficiently many points (cf.~\ptref{th:ex.der.symmprod}).
Notice that, unlike our previous considerations, 
the results of this subsection seem to be new even in the
additive case.

\zerosubpoint\nxsubpoint 
(Importance of derived symmetric powers.)
Existence of derived symmetric powers is crucial for the goal of this work,
since we want to define Chow rings by means of the $\gamma$-filtration,
exactly the way Grothendieck did it in his proof of Riemann--Roch theorem, 
and a $\lambda$-ring structure on $K_0$ of perfect complexes (or rather perfect
simplicial modules) is necessary for this. A $\lambda$-structure can be 
described either in terms of exterior power operations $\{\lambda^n\}_{n\geq0}$
or in terms of symmetric powers $\{s^n\}_{n\geq0}$. We choose the second
approach for the following reasons:
\begin{itemize}
\item Symmetric powers have good properties in all situations, while
exterior powers are nice only for an {\em alternating\/}~$\sO$
(cf.~\ptref{def:altop}), e.g.\ $S_\sO^n(X\oplus Y)=\bigoplus_{p+q=n}
S_\sO^p(X)\otimes_\sO S_\sO^q(Y)$ is true unconditionally, 
while a similar formula for $\bigwedge^n$ requires alternativity of~$\sO$.
\item We don't want to take care of the $\pm$ signs in all our proofs.
\item The proof for symmetric powers is ultimately reduced to the 
basic case of morphisms from $I$ or $J$ in $s\catSets$, which can be dealt
with the aid of topological realization, while the basic case for exterior 
powers would be $L_{\Fpm}(I)$ or $L_{\Fpm}(J)$ in $s\catMod\Fpm$, 
more difficult to tackle with.
\item Symmetric powers have a long and venerable history of being derived,
and they actually were the first example considered by Dold and Puppe
in their work \cite{DP} 
on derivation of non-additive functors between abelian categories.
\end{itemize}

\nxsubpoint
(Notations.)
Let us fix a generalized commutatively ringed topos $(\cE,\sO)$. 
We denote by $\cC$ the category $s\catMod\sO$, endowed with its 
pseudomodel structure discussed in~\ptref{sp:distcl.sO.mod}
and~\ptref{cor:indeed.psmod.sO-mod} 
(when $\cE$ doesn't have enough
points, we cannot speak about a pseudomodel structure, but we still have
well-defined cofibrations and acyclic cofibrations in~$\cC$)
and its compatible ACU $\otimes$-structure $\otimes=\otimes_\cE$
(cf.~\ptref{prop:comp.tens.sO-mod}). On some occasions $\cC$ will denote
the corresponding stack $\stks\stMOD\sO_\cE$, and we'll deliberately confuse
$\cC$ with its final fiber $\cC(e)=s\catMod\sO$.

Notice that almost all our considerations will be actually valid for
any pseudomodel category or stack with a compatible ACU $\otimes$-structure,
if we assume in addition that $\otimes:\cC\times\cC\to\cC$ commutes with
arbitrary inductive limits in each variable (condition automatically fulfilled
whenever $\otimes$ admits inner Homs).

When a finite group~$G$ acts (say, from the right) 
on an object~$X$ of~$\cC$, we denote by $X/G$ or $X_G$ its {\em object of 
coinvariants}, i.e.\ the largest strict quotient of~$X$, on which $G$
acts trivially. We can also describe $X_G$ as the cokernel of 
$\pr_1,\alpha:X\times G\rightrightarrows X$, where 
$\alpha:X\times G=\bigsqcup_{g\in G}X\to X$ is the group action. 
Right exactness of $\otimes$ implies $(X\otimes Y)_G\cong X_G\otimes Y$ for
any $Y$ with a trivial action of~$G$, hence
\begin{equation}\label{eq:coinv.tens.prod}
(X\otimes Y)_{G\times H}\cong ((X\otimes Y)_G)_H\cong
(X_G\otimes Y)_H\cong X_G\otimes Y_H
\end{equation}
for any $G$-object $X$ and $H$-object $Y$.

We denote tensor powers by $T^n(A)$ or $A^{\otimes n}$, and the symmetric 
powers $T^n(A)/\gS_n=T^n(A)_{\gS_n}$ by $S^n_\sO(A)$ or simply $S^n(A)$.
We use these notations both for these functors in $\catMod\sO$ and for
their simplicial extensions to $s\catMod\sO$.

\begin{ThD}\label{th:ex.der.symmprod}
Let $(\cE,\sO)$ be a generalized commutatively ringed topos,
$n\geq0$ be an integer. Then:
\begin{itemize}
\item[(a)] Symmetric power functor $S^n:s\catMod\sO\to s\catMod\sO$
preserves (pointwise) cofibrations and acyclic cofibrations between 
(pointwise) cofibrant objects.%
\footnote{\rm We give a complete proof of (a) only in the pointwise case,
since the general case would require a development of certain techniques
we want to avoid in this work, and there are some indications one would
need to extend generating sets $L_\sO(\cst I_\cE)$ and $L_\sO(\cst J_\cE)$
in order to prove it, thus replacing cofibrations and acyclic cofibrations
with something more similar to corresponding pointwise notions anyway.
In any case, this doesn't affect validity of (b) and (c), 
any cofibration being pointwise cofibrant as well.}
\item[(b)]
If $\cE$ has enough points, $S^n$ preserves weak equivalences between cofibrant
objects, and the same is true even for weak equivalences between pointwise
cofibrant objects.
\item[(c)] If $\cE$ has enough points, the symmetric power functor 
admits a left derived functor $\dL S^n:\cD^{\leq0}(\cE,\sO)\to
\cD^{\leq0}(\cE,\sO)$, which can be computed with the aid of cofibrant 
replacements or pointwise cofibrant replacements.
\end{itemize}
\end{ThD}
\begin{Proof}
We want to show that everything follows from (a), which will be proved 
in the remaining part of this subsection by induction in~$n$.
(a)$\Rightarrow$(b): Follows immediately from the fact that weak equivalences
can be checked pointwise (cf.~\ptref{prop:f.af.weq.sO-mod}),
Lemma~\ptref{l:acof.cof.suff}, 
and commutativity of symmetric powers with topos pullbacks
and fiber functors $p^*$ in particular. (b)$\Rightarrow$(c): 
Apply~\ptref{th:ex.left.der.psmod}.
\end{Proof}

Before proving (a), or actually a more precise statement
\ptref{prop:symmprod.steps.cof},
we want to introduce some auxilliary constructions.

\nxsubpoint 
(Hypercube construction.)
Let us fix an ordered set~$I$; usually it will be a subset of $\bbZ$ with 
induced order, e.g.\ $[1]=\{0,1\}$, or some other $[m]=\{0,1,\ldots,m\}$,
or $\omega=\{m:m\geq0\}$ -- in short, any ordinal $\leq\omega$.
We consider ordered sets as categories with at most one morphism
between any two objects in the usual manner. Let us also fix some $n\geq1$,
and put $\stn:=\{1,2,\ldots,n\}$ as an (unordered) set. Then 
the symmetric group $\gS_n=\Aut(\stn)$ acts from the right on
$I^n=\catFunct(\stn,I)$, hence it also acts from the left on functors
$I^n\to\cC$.

More precisely, if $\epsilon=(\epsilon_1,\ldots,\epsilon_n)\in I^n$ and
$\sigma\in\gS_n$ is a permutation, then $\epsilon\sigma=
(\epsilon_{\sigma(1)},\ldots,\epsilon_{\sigma(n)})$, and if 
$P:I^n\to\cC$ is a functor, then $\sigma P:I^n\to\cC$ is given by
$(\sigma P)(\epsilon):=P(\epsilon\sigma)$, i.e.\ 
$(\sigma P)(\epsilon_1,\ldots,\epsilon_n)=P(\epsilon_{\sigma(1)},\ldots,
\epsilon_{\sigma(n)})$.

Now let $\{X_i:I\to\cC\}_{1\leq i\leq n}$ be any family of functors. 
For example, if $I=[1]=\{0\to 1\}$, then each $X_i$ is just a morphism
$X_i(0)\to X_i(1)$ in~$\cC$. Define $\tilde X:I^n\to\cC$ by taking the
tensor product of values of $X_i$:
\begin{equation}
\tilde X(\epsilon_1,\ldots,\epsilon_n)=X_1(\epsilon_1)\otimes\cdots
\otimes X_n(\epsilon_n)
\end{equation}
When $I=[1]$, we obtain an $n$-dimensional hypercube of arrows in~$\cC$,
hence the name for our construction.

\nxsubpoint (Extension of $\tilde X$ to subsets of~$I^n$.)
Let $R\subset I^n$ be any subset of $I^n$, considered here as an ordered set
with respect to induced order, or as a full subcategory of $I^n$.
We extend $\tilde X:I^n\to\cC$ from $I^n$ to the set $\gP(I^n)$ of subsets
of~$I^n$, ordered by inclusion, as follows:
\begin{equation}\label{eq:ext.Xtil.to.subsets}
\tilde X_!(R):=\injlim_{\epsilon\in R}\tilde X(\epsilon)
\end{equation}
Clearly, $\tilde X_!(\{\epsilon\})=\tilde X(\epsilon)$. 
More generally, if $R$ has a largest
element $\epsilon$, then $X_!(R)=X(\epsilon)$. For example, if $I=[m]$, then
$X_!(I^n)=\tilde X(m,\ldots, m)=X_1(m)\otimes\cdots\otimes X_n(m)$.

Notice that the right action of $\gS_n$ on $I^n$ naturally extends to a 
right action on $\gP(I^n)$, and we have
\begin{equation}
\tilde X_!(R\sigma)=\injlim_{\epsilon\in R}\tilde X(\epsilon\sigma)=
(\sigma\tilde X)_!(R)
\end{equation}
In other words, $\sigma(\tilde X_!)=(\sigma\tilde X)_!$, so we can write
simply $\sigma\tilde X_!$.

\nxsubpoint\label{sp:mult.box.prod}
(Application: multiple box products.)
Let $I=[1]$, and $\{X_i\}_{1\leq i\leq n}$ be a family of functors
$I\to\cC$, i.e.\ a collection of morphisms $\{X_i(0)\stackrel{u_i}\to
X_i(1)\}$ in~$\cC$. Consider ordered subset $R:=I^n-\{(1,1,\ldots,1)\}\subset
I^n$, obtained by removing the ``final'' vertex of our hypercube. 
Then we get a morphism 
\begin{equation}
u_1\boxe u_2\boxe\cdots\boxe u_n:
\tilde X_!(R)\to\tilde X_!(I^n)=X_1(1)\otimes X_2(1)\otimes\cdots\otimes
X_n(1)
\end{equation}
We say that $u_1\boxe\cdots\boxe u_n$ is the {\em multiple box product of
$u_1$, \dots, $u_n$.} When $n=2$, we recover the usual $\boxe$-product 
$u\boxe v$ of \ptref{def:comp.tens.mod.str}, used in particular to state~(TM).

Now it is an easy exercise to check that $(u_1\boxe u_2)\boxe u_3=
u_1\boxe u_2\boxe u_3=u_1\boxe(u_2\boxe u_3)$; the verification is quite 
formal and uses only the right exactness of~$\otimes$. In other words,
{\em the box-product $\boxe:\Ar\cC\times\Ar\cC\to\Ar\cC$ is associative.}
It is of course also commutative, $\otimes$ being commutative.
(All this would seem not so unexpected if we use the ``fancy language''
of~\ptref{sp:fancy.lang} and recall that $\boxe$ ``corresponds to tensor
product of cokernels''.) 

We conclude by an easy induction using (TM), valid in $\cC=s\catMod\sO$
by~\ptref{prop:comp.tens.sO-mod}, that {\em $u_1\boxe u_2\boxe\cdots\boxe u_n$
is a cofibration whenever all $u_i$ are cofibrations, and it is an acyclic
cofibration if in addition at least one of $u_i$ is an acyclic cofibration.}

\nxsubpoint\label{sp:symm.hyperc}
(Symmetric hypercube construction.)
Let us fix an ordered set~$I$ and an integer $n>0$ as above, and
any functor $X:I\to\cC$. Put all $X_i:=X$, $1\leq i\leq n$, and
construct $\tilde X:I^n\to\cC$ with the aid of tensor product in~$\cC$ 
as before:
\begin{equation}
\tilde X(\epsilon_1,\ldots,\epsilon_n)=X(\epsilon_1)\otimes\cdots\otimes
X(\epsilon_n)
\end{equation}
The associativity and commutativity constraints for $\otimes$ yield
canonical isomorphisms $\eta_{\sigma,\epsilon}:\tilde X(\epsilon)\to
\tilde X(\epsilon\sigma)=(\sigma\tilde X)(\epsilon)$, which can be 
described on the level of elements as the maps 
$X(\epsilon_1)\otimes\cdots\otimes X(\epsilon_n)\to 
X(\epsilon_{\sigma(1)})\otimes\cdots\otimes X(\epsilon_{\sigma(n)})$,
given by
$x_1\otimes\cdots\otimes x_n\mapsto x_{\sigma(1)}\otimes\cdots
\otimes x_{\sigma(n)}$.

Clearly, $\eta_{\sigma,\epsilon}:\tilde X(\epsilon)\to
(\sigma\tilde X)(\epsilon)$ are functorial in $\epsilon$, i.e.\ they 
define a natural transformation $\eta_\sigma:\tilde X\to\sigma\tilde X$,
for any $\sigma\in\gS_n$. Moreover, these $\eta_\sigma$
{\em satisfy a cocycle relation $\eta_{\sigma\tau}=\sigma(\eta_\tau)\circ
\eta_\sigma:\tilde X\to\sigma\tilde X\to\sigma\tau\tilde X$}.

By taking inductive limits $\injlim_R$ along subsets $R\subset I^n$,
we get induced morphisms $\eta_{\sigma,R}:\tilde X_!(R)\simto
\tilde X_!(R\sigma)=(\sigma\tilde X_!)(R)$, i.e.\ $\eta_\sigma$
extends to a functorial isomorphism $\tilde X_!\to\sigma\tilde X_!$,
where $\tilde X_!:\gP(I^n)\to\cC$ is the extension of $\tilde X$ defined
by~\eqref{eq:ext.Xtil.to.subsets}. It is immediate that these extended
$\eta_\sigma$ still satisfy the cocycle relation 
$\eta_{\sigma\tau}=\sigma(\eta_\tau)\circ\eta_\sigma$, i.e.\ 
$\eta_{\sigma\tau,R}=\eta_{\tau,R\sigma}\circ\eta_{\sigma,R}$.

In particular, for a {\em symmetric\/} $R\subset I^n$ (i.e.\ stable under
$\gS_n$: $R\sigma=R$ for all $\sigma\in\gS_n$) we get isomorphisms
$\eta_{\sigma,R}:\tilde X_!(R)\simto\tilde X_!(R\sigma)=\tilde X_!(R)$,
and the cocycle relation implies that these 
$\{\eta_{\sigma,R}\}_{\sigma\in\gS_n}$ define a {\em right\/} action
of $\gS_n$ on $\tilde X_!(R)$, hence we can consider the 
coinvariants of $\tilde X_!(R)$, which will be denoted by $\tilde X_!^s(R)$:
\begin{equation}
\tilde X_!^s(R):=\tilde X_!(R)_{\gS_n}\quad\text{for any symmetric 
$R\subset I^n$}
\end{equation}
We say that $X_!^s:\gP(I^n)^{\gS_n}\to\cC$ is the {\em symmetric extension\/}
of $X$ or~$\tilde X$.

Notice that the above construction applies to any $R\subset I^n$, fixed
by some subgroup $G\subset\gS_n$: we still obtain a right action of
$G$ on $\tilde X_!(R)$, hence can consider the coinvariants $\tilde X_!(R)_G$.

\nxsubpoint\label{sp:can.dec.symm.pow}
(Application: canonical decomposition of symmetric powers.)
Let us apply the above construction to any morphism $u:X(0)\to X(1)$,
considered here as a functor $X:[1]\to\cC$. 
We obtain some symmetric extension $\tilde X_!^s:\gP([1]^n)^{\gS_n}\to\cC$.
Put $R_k:=[1]^n_{\leq k}=\{\epsilon\in [1]^n : |\epsilon|\leq k\}$,
where $|\epsilon|$ denotes the sum $\sum_{i=1}^n\epsilon_i$ of components 
of~$\epsilon\in [1]^n\subset\bbZ^n$. Clearly, 
$\{0\}=R_0\subset R_1\subset\cdots\subset R_n=[1]^n$ 
is an increasing filtration of $[1]^n$ by symmetric subsets. 
Put $F_kS^n(u):=\tilde X_!^s(R_k)=\tilde X_!(R_k)_{\gS_n}$; taking into account
obvious identifications of $F_0S^n(u)=\tilde X_!^s(\{0\})$ with 
$(X(0)^{\otimes n})_{\gS_n}=S^nX(0)$ and of $F_nS^n(u)=\tilde X_!^s([1]^n)$
with $S^nX(1)$, we obtain a canonical decomposition of $S^n(u)$:
\begin{multline}
S^n(u):S^nX(0)=F_0S^n(u)\stackrel{\rho_1^{(n)}(u)}\longto
F_1S^n(u)\stackrel{\rho_2^{(n)}(u)}\longto\cdots\\
\cdots\longto F_{n-1}S^n(u)\stackrel{\rho_n^{(n)}(u)}\longto F_nS^n(u)=S^nX(1)
\end{multline}
Sometimes we consider also $R_{-1}=[1]^n_{\leq-1}=\emptyset$, together with
$\rho_0^{(n)}(u):F_{-1}S^n(u)=\emptyset\to F_0S^n(u)=S^nX(0)$. 

We also introduce a shorter notation $\rho_n(u)$ for the
``$n$-fold symmetric box product''
$\rho_n^{(n)}(u):F_{n-1}S^n(u)\to S^nX(1)$; we'll see in a moment that 
$\rho_n(u)$ ``computes the $n$-th symmetric power of $\Coker u$'', hence
its importance for our considerations.

\nxsubpoint (Fancy explanation.)
The meaning of the above decomposition of $S^n(u)$ is better explained in
the ``fancy language'' of~\ptref{sp:fancy.lang}. One can check that the 
above construction, when applied to an injective homomorphism
$X(0)\stackrel u\to X(1)$ with projective cokernel~$P$ (in $\catMod\Lambda$,
$\Lambda$ a classical commutative ring), yields a decomposition of
$S^n(u)$ into $n$ morphisms $\rho_k^{(n)}$, which correspond to the canonical
increasing filtration of $\Coker S^n(u)=S^nX(1)/S^nX(0)$ with 
quotients isomorphic to $S^kP\otimes S^{n-k}X(0)$. In particular,
$\rho_n(u)=\rho_n^{(n)}(u)$ corresponds to the last step of this filtration,
with quotient $S^nP$, i.e.\ $\Coker\rho_n(u)\cong S^n(\Coker u)$ :
{\em $\rho_n$ computes symmetric powers of the cokernel}, similarly 
to the relationship between the box-product and the tensor product.

Now the statement (a) of~\ptref{th:ex.der.symmprod} will be a consequence
of the following more precise statement:
\begin{PropD}\label{prop:symmprod.steps.cof}
Whenever $u:A\to B$ is a cofibration (resp.\ acyclic cofibration)
between cofibrant objects of $\cC=s\catMod\sO$, the same is true for
all $\rho_k^{(n)}(u):F_{k-1}S^n(u)\to F_kS^n(u)$, $1\leq k\leq n$, 
hence also for their composition 
$S^n(u)=\rho_n^{(n)}\rho_{n-1}^{(n)}\cdots\rho_1^{(n)}:S^n(A)\to S^n(B)$.
\end{PropD}
{\bf Proof.} Induction in~$n\geq0$, cases $n=0$ and $n=1$ being trivial.
The induction step will occupy the rest of this subsection; its main
idea is of course our usual sort of ``devissage'' (cf.~\ptref{sp:fancy.lang}).
However, we want to work with acyclic cofibrations and cofibrations between
cofibrant objects only. The possibility of doing this is due to the fact
that {\em acyclic cofibrations (resp.\ cofibrations) 
between cofibrant objects coincide with the closure of $L_\sO(\cst I)$ 
(resp.\ $L_\sO(\cst J)$) inside the full subcategory $\cC_c$ of 
$\cC=s\catMod\sO$, or, more precisely, inside the full substack $\cC_c$
of $\cC=\stks\stMOD\sO_\cE$,} shown in the following lemma:

\begin{LemmaD}
Let $\cC$ be a flat stack over a topos~$\cE$, $I$, $J$ be (small) 
sets of morphisms in the final fiber of $\cC$, such that $J\subset\Cl I$.
Let us call morphisms from $\Cl I$ (resp.\ $\Cl J$) cofibrations 
(resp.\ acyclic cofibrations), and denote by $\cC_c$ the full substack
of~$\cC$ consisting of cofibrant objects. Denote by $\Cl_c I$ 
(resp.\ $\Cl_c J$) the closure of $I$ (resp.\ $J$) in $\cC_c$. Then
$\Cl_c I=\Cl I\cap\Ar\cC_c$ and $\Cl_c J=\Cl J\cap\Ar\cC_c$, i.e.\ 
$\Cl_c I$ (resp.\ $\Cl_c J$) consists exactly of cofibrations 
(resp.\ acyclic cofibrations) between cofibrant objects.
\end{LemmaD}
{\bf Proof.}
It suffices to show $\Cl_c J=\Cl J\cap\Ar\cC_c$, since the other formula
can be then shown by invoking this lemma with $J:=I$. 
Clearly, $\Cl J\cap\Ar\cC_c$ is a closed class of morphisms in $\cC_c$,
containing~$J$, hence also $\Cl_c J$, so one inclusion is trivial.

To prove the opposite inclusion we consider the class $\propP$ 
of morphisms $A\stackrel u\to B$ in (fibers of) $\cC$ having the following
property: all pushouts $A'\stackrel{u'}\to B'$ 
of $\phi^*(u)$ with respect to morphisms $\phi^*A\to A'$ in $\cC(T)$ 
with cofibrant target lie in~$\Cl_c J$, where $\phi:T\to S$ runs over all 
morphisms in~$\cE$ with target~$S$, and $u$ lies in $\Ar\cC(S)$.

Clearly, $J\subset\propP$; if we show that $\propP$ is closed in $\cC$,
we would conclude $\Cl J\subset\propP$, hence also the required inclusion
$\Cl J\cap\Ar\cC_c\subset\propP\cap\Ar\cC_c\subset\Cl_c J$, applying the
characteristic property of morphisms $A\stackrel u\to B$ from $\propP$
to the pushout with respect to $\id_A$, whenever~$A$ is cofibrant.

\nxsubpoint (Closedness of $\propP$.)
So we are reduced to proving closedness of~$\propP$, i.e.\ verifying
conditions 1)--8) of~\ptref{def:closed.class}. 

(a) Condition 1) 
(localness of $\propP$) is obvious from the definition of~$\propP$;
since all our constructions are also compatible with pullbacks, we
need to check the basic property of morphisms from $\propP$ only for pushouts
computed in the same fiber, i.e.\ assume $\phi=\id_S$ while checking 2)--7),
but not 8).

(b) Conditions 2) (all isomorphisms do lie in $\propP$) and 3)
(stability of $\propP$ under pushouts) are evident. Conditions 4) and 5)
(stability under composition and sequential composition) are also evident,
once we take into account that $\Cl J\subset\Cl I$ implies that any pushout
$u':A'\to B'$ with cofibrant source~$A'$ of a morphism $u:A\to B$ 
from~$\propP$ has a cofibrant target~$B'$ as well.

(c) Condition 6) (stability under retracts) is somewhat more complicated.
Let $v:X\to Y$ be a retract of $u:A\to B$ from $\propP$, 
i.e.\ we have a commutative diagram
\begin{equation}\label{eq:retr.diag2}
\xymatrix{
X\ar[d]^{v}\ar[r]_{i}\ar@/^1pc/[rr]^{\id_X}&A\ar[d]^{u}\ar[r]_{p}&X\ar[d]^{v}\\
Y\ar[r]^{j}\ar@/_1pc/[rr]_{\id_Y}&B\ar[r]^{q}&Y}
\end{equation}
According to Lemma~\ptref{l:repl.retr} below, we may assume $A=X$, 
$i=p=\id_A$, at the cost of replacing $u$ by its pushout~$\tilde u$ 
with respect to~$p$. Since $\tilde u$ still lies in $\propP$, we can safely
do this. Now the pushout $v':A'\to Y'$ of $v:A\to Y$ with respect to
any morphism $A\to A'$ with cofibrant target~$A'$ is a retract of 
the corresponding
pushout $u':A'\to B'$ of $u$, lying in $\Cl_c J$ by assumption. In particular
$u'\in\Cl_c J\subset\Cl J\subset\Cl I$ is a cofibration in~$\cC$, hence the
same is true for its retract~$v'$, hence the target of~$v'$ lies in~$\cC_c$,
i.e.\ $v'$ is a retract of~$u'$ inside $\cC_c$, hence it lies in $\Cl_c J$ 
as well.

(d) Let us prove condition 7) (stability under direct sums).
Let $\{u_\alpha:A_\alpha\to B_\alpha\}_{\alpha\in K}$ 
be a family of morphisms from~$\propP$, and put $A:=\bigoplus A_\alpha$,
$B:=\bigoplus B_\alpha$, $u:=\bigoplus u_\alpha:A\to B$. We have to check that
$u$ lies in $\propP$. Let $\sigma:A=\bigoplus A_\alpha\to A'$ be any morphism
with cofibrant target. Consider the components $\sigma_\alpha:A_\alpha\to A'$
of~$\sigma$, together with their direct sum 
$\tilde\sigma:=\bigoplus\sigma_\alpha:A\to A'\otimes K=\bigoplus_\alpha A'$.
Denote by $\nabla$ the codiagonal map $A'\otimes K=\bigoplus_\alpha A'\to A'$,
and consider the following commutative diagram with cocartesian squares:
\begin{equation}
\xymatrix@C+10pt{
A\ar[d]^{u}\ar[r]^<>(.5){\tilde\sigma}&
A'\otimes K\ar[d]^{v}\ar[r]^<>(.5){\nabla}&
A'\ar[d]^{u'}\\
B\ar[r]&\bigoplus B'_\alpha\ar[r]&B'
}
\end{equation}
Here $v$ is the direct sum of pushouts $v_\alpha:A'\to B'_\alpha$ of morphisms
$u_\alpha$. Each such pushout lies in $\Cl_c J$ by assumption, hence the
same is true for their direct sum $v$, and in particular the source and the
target of~$v$ lie in~$\cC_c$. On the other hand, the right square of the 
above diagram is cocartesian in $\cC_c$, hence $u'$ lies in $\Cl_c J$ as well.
Now we just observe that $u'$ is exactly the pushout of $u$ with respect to
$\sigma=\nabla\tilde\sigma$.

(e) It remains to check 8), i.e.\ stability under $\phi_!$ or 
``local direct sums''. So let us fix some $\phi:T\to S$ in $\cE$; we can assume
$S=e$, replacing $\cE$ and $\cC$ by $\cE_{/S}$ and $\cC_{/S}$ if necessary.
Let $u:A\to B$ be a morphism of $\cC(T)$ lying in $\propP$. We want to show 
that any pushout of $\phi_!(u):\phi_!A\to\phi_!B$ with respect to a morphism
$\sigma:\phi_!A\to A'$ in~$\cC(e)$ with cofibrant target lies in $\Cl_c J$.
Denote by $\sigma^\flat:A\to\phi^*A'$ the morphism in~$\cC(T)$ defined by
adjointness of $\phi_!$ and $\phi^*$. Let $v:\phi^*A'\to Z$ be the pushout
in $\cC(T)$ of $u$ with respect to $\sigma^\flat$; since $\phi^*A'$ is 
cofibrant and $u$ lies in $\propP$, $v$ lies in $\Cl_cJ$, hence the same
is true for $\phi_!(v)$. Now consider the following commutative diagram
in $\cC(e)$ with cocartesian squares:
\begin{equation}
\xymatrix@C+10pt{
\phi_!A\ar[d]^{\phi_!(u)}\ar[r]^<>(.5){\phi_!(\sigma^\flat)}&
\phi_!\phi^*A'\ar[d]^{\phi_!(v)}\ar[r]^<>(.5){\nabla}&A'\ar[d]^{u'}\\
\phi_!B\ar[r]&\phi_!Z\ar[r]&B'
}
\end{equation}
Clearly, $u':A'\to B'$ is a pushout of $\phi_!(v)$ in $\cC_c(e)$, hence it
lies in $\Cl_c J$ as well. On the other hand, $u'$ is exactly the pushout
in~$\cC(e)$ of~$\phi_!(u)$ 
with respect to $\nabla\circ\phi_!(\sigma^\flat)=\sigma$, q.e.d.

\begin{LemmaD}\label{l:repl.retr}
Suppose that $v:X\to Y$ is a retract of $u:A\to B$ in some category~$\cC$,
by means of a retract diagram~\eqref{eq:retr.diag2}. Suppose that the pushout 
$\tilde u:X\to Z$ of $u$ with respect to $p:A\to X$ exists in~$\cC$. Then
$v$ is a retract of $\tilde u$, and one can find a retract diagram
$(\tilde\imath,\tilde\jmath):v\to\tilde u$, $(\tilde p,\tilde q):\tilde u\to v$
with property $\tilde\imath=\tilde p=\id_X$.
\end{LemmaD}
\begin{Proof}
Consider the following diagram induced by~\eqref{eq:retr.diag2}, where the
right square is cocartesian:
\begin{equation}
\xymatrix{
X\ar[d]^{v}\ar[r]^{i}&A\ar[d]^{u}\ar[r]^{p}&
X\ar@{-->}[d]_{\tilde u}\ar[ddr]^{v}\\
Y\ar[r]^{j}&B\ar@{-->}[r]^{r}\ar[rrd]_{q}&
Z\ar@{-->}[rd]|{\tilde q}\\&&&Y}
\end{equation}
Put $\tilde\imath:=\tilde p:=\id_X$, $\tilde\jmath:=rj$. Then these morphisms
$(\tilde\imath,\tilde\jmath):v\to\tilde u$ and 
$(\tilde p,\tilde q):\tilde u\to v$ do
satisfy the retract conditions, thus proving that $v$ is a retract of
$\tilde u$ with required property. Indeed,
$\tilde u\tilde\imath=\tilde u\id_X=\tilde upi=rjv=\tilde\jmath v$,
$v\tilde p=v=\tilde q\tilde u$, $\tilde p\tilde\imath=\id_X$, and
$\tilde q\tilde\jmath=\tilde qrj=qj=\id_Y$, q.e.d.
\end{Proof}

Before proving the induction step of~\ptref{prop:symmprod.steps.cof},
we need to construct some decompositions of $\rho_k^{(n)}(u)$, where
$u=\cdots u_2u_1$ is a finite or sequential composition of morphisms,
into some ``simple'' morphisms.

\nxsubpoint (Sieves, symmetric sieves and their extensions.)
Let $I$ be an ordered set as before. We say that a subset $R\subset I^n$
is a {\em sieve\/} if $x\leq y$ and $y\in R$ implies $x\in R$. We say that
$R$ is a {\em symmetric sieve\/} if it is fixed by $\gS_n$: $R=R\sigma$
for all $\sigma\in\gS_n$. An {\em extension of sieves\/} is simply a
couple $(R,R')$ of sieves in $I^n$, such that $R\subset R'$.
Finally, an extension $R\subset R'$ of symmetric sieves will be said to
be {\em simple\/} if $R\neq R'$, but there are no symmetric sieves $R''$,
such that $R\subsetneq R''\subsetneq R'$.

\nxsubpoint\label{sp:struct.simple.syms}
(Structure of simple extensions of symmetric sieves.)
Let $I$ be an ordinal $\leq\omega$, i.e.\ either a finite set $[m]$,
or the infinite set $\omega$ of non-negative integers. Let $R\subsetneq R'$
be any non-trivial extension of symmetric sieves. Choose any 
minimal element $\gamma=(\gamma_1,\ldots\gamma_n)$ 
of the complement $R'-R$ of $R$ in~$R'$
(e.g.\ choose $\gamma$ with minimal value of $|\gamma|$).
Clearly, $\gamma\not\in R$, but $\gamma-e_k=(\gamma_1,\ldots,\gamma_k-1,
\ldots,\gamma_n)\in R$, for any $k$, such that $\gamma_k>0$.
One easily checks that whenever we are given a symmetric sieve $R\subset I^n$ 
and an element $\gamma\in I^n$ with the above property, 
then $R'':=R\cup\gamma\gS_n$ is a symmetric sieve in~$I^n$, and, 
moreover, the extension $R\subset R''$ is simple.

We conclude that any non-trivial extension of symmetric sieves $R\subset R'$
contains a simple extension $R\subset R''=R\cup\gamma\gS_n\subset R'$ of the
above form; if $R\subset R'$ is itself simple, then $R'=R''$, i.e.\ 
{\em all simple extensions of symmetric sieves are of form 
$R\subset R\cup\gamma\gS_n$.} Since we are free to replace $\gamma$ by any
$\gamma\sigma$, {\em we can assume that $\gamma_1\leq\gamma_2\leq\cdots\leq
\gamma_n$.}

\nxsubpoint (Relation to the symmetric hypercube construction.)
Now let us fix any functor $X:I\to\cC$, i.e.\ a finite or infinite sequence
of composable morphisms, and apply the symmetric hypercube construction.
Let $R\subset R'=R\cup\gamma\gS_n$ be a simple extension of symmetric 
sieves as above, $\gamma=(\gamma_1,\ldots,\gamma_n)\in I^n$, and
denote by~$G$ the stabilizer $\Stab_{\gS_n}(\gamma)$ of~$\gamma$.
Put $W':=\{\gamma_1-1,\gamma_1\}\times\cdots\times
\{\gamma_n-1,\gamma_n\}\cap I^n$ -- the ``small hypercube at~$\gamma$'',
and $W:=W'\cap R=W'-\{\gamma\}$, the latter equality being true because of 
the characteristic property of $\gamma$. Clearly, $W\subset R$ and
$W'\subset R'$, so we get a commutative square
\begin{equation}
\xymatrix{
X_!(W)\ar[r]\ar[d]&X_!(W')\ar[d]\\X_!(R)\ar[r]&X_!(R')
}
\end{equation}
Notice that $W$ and $W'$ are stabilized by $G$, hence $G$ acts from the right
on $X_!(W)$ and $X_!(W')$, as well as on $X_!(R)$ and $X_!(R')$,
on which in fact a larger group $\gS_n\supset G$ acts
(cf.~\ptref{sp:symm.hyperc}), 
and all arrows are compatible with these group actions,
whence a commutative diagram of coinvariants:
\begin{equation}\label{eq:simple.step.square}
\xymatrix{
X_!(W)_G\ar[r]^{v_\gamma}\ar[d]&X_!(W')_G\ar[d]\\
X^s_!(R)\ar[r]^{u_{R,\gamma}}&X^s_!(R')
}
\end{equation}
We want to prove that this square is cocartesian, and find an expression
for $v_\gamma$ in terms of box-products of morphisms $u_k:X(k-1)\to X(k)$
and their ``symmetric box-powers'' $\rho_d(u_k)$.

\begin{PropD}\label{prop:cocart.step.sq}
Square \eqref{eq:simple.step.square} is cocartesian.
\end{PropD}
\begin{Proof}
(a) First of all, the above commutative square can be constructed componentwise
in terms of components $X(i)_k$ of simplicial $\sO$-modules $X(i)$,
so we can check its couniversality componentwise, i.e.\ replace all
$X(i)$ with $X(i)_k$ and consider the corresponding symmetric hypercube
construction in $\cC:=\catMod\sO$, not in $s\catMod\sO$. This reduction step
is in fact not indispensable (all we actually need is a functor 
$\tilde X:I^n\to\cC$ together with a ``cocycle'' $\{\eta_\sigma:\tilde X\simto
\sigma\tilde X\}$, so we can work in any~$\cC$), 
but it can simplify the understanding of the proof.

(b) We have to check that \eqref{eq:simple.step.square} becomes a cartesian
square of sets after an application of contravariant functor 
$h_V:=\Hom_\cC(-,V)$, for any $V\in\Ob\cC$. 
Let us fix such a~$V$ and denote by $P$, $\tilde P$ and $\tilde P_*$ the
compositions of $X$, $\tilde X$ and $\tilde X_!$ with $h_V$. Denote by
$\sigma_*$ the maps $h_V(\eta_\sigma):\tilde P(\epsilon\sigma)\to
\tilde P(\epsilon)$ and $\tilde P_*(S\sigma)\to\tilde P_*(S)$, for any
$S\subset I^n$.
Taking into account that $h_V$ transforms inductive limits into projective 
limits and coinvariants into invariants, we see that we have to
prove that the following commutative square of sets is cartesian:
\begin{equation}
\xymatrix{
\tilde P_*(R')^{\gS_n}\ar[r]^{u^*}\ar[d]^{q^*}&
\tilde P_*(R)^{\gS_n}\ar[d]^{p^*}\\
\tilde P_*(W')^G\ar[r]^{v^*}&
\tilde P_*(W)^G}
\end{equation}
where $\tilde P_*(S)=\projlim_{\epsilon\in S}\tilde P(\epsilon)$ is 
a projective limit of sets, i.e.\ it is just the set of collections
$x=\{x_\epsilon\}_{\epsilon\in S}$, $x_\epsilon\in\tilde P(\epsilon)$,
compatible in the following sense: $x_\epsilon|_\delta=x_\delta$ for 
any $\delta\leq\epsilon$ in $S$, where of course the ``restriction''
$z|_\delta$ denotes the image of $z\in\tilde P(\epsilon)$ under the
canonical ``restriction'' map $\tilde P(\epsilon)\to\tilde P(\delta)$.

(c) Let us check that an element $x=\{x_\epsilon\}_{\epsilon\in R'}$
of $\tilde P_*(R')^{\gS_n}$ is completely determined by its restrictions
$u^*(x)$ and $q^*(x)$ to $R$ and $W'$, respectively. Indeed,
$u^*(x)$ determines all $x_\epsilon$ for $\epsilon\in R$, and
$q^*(x)$ determines $x_\gamma$, hence also all
$x_{\gamma\sigma}=\sigma_*^{-1}(x_\gamma)$, $x$ being $\gS_n$-invariant.
We see that all components of $x$ are thus determined 
since $R'=R\cup\gamma\gS_n$.

(d) Now it remains to check that given any two such collections
$y=\{y_\epsilon\}_{\epsilon\in R}\in\tilde P_*(R)^{\gS_n}$ and
$z=\{z_\epsilon\}_{\epsilon\in W'}\in\tilde P_*(W')^G$, such that 
$p^*(y)=v^*(z)$, i.e.\ $y_\epsilon=z_\epsilon$ for all 
$\epsilon\in W=W'\cap R$, one can find a compatible 
collection $x=\{x_\epsilon\}_{\epsilon\in R'}$ in $\tilde P_*(R')^{\gS_n}$,
which restricts to $y$ on $R$ and to $z$ on $W'$. Of course, we put
$x_\epsilon:=y_\epsilon$ for $\epsilon\in R$ and $x_\epsilon:=z_\epsilon$
for $\epsilon\in W'$, thus defining $x$ on $R\cup W'=R\cup\{\gamma\}$,
these two definitions being compatible on $R\cap W'=W$.
Furthermore, the element $x$ being constructed has to be $\gS_n$-invariant,
i.e.\ we must have $x_{\gamma\sigma}=\sigma_*^{-1}(x_\gamma)=
\sigma_*^{-1}(z_\gamma)$ for all
$\sigma\in\gS_n$. This indeed defines $x$ on $\gamma\gS_n$ since if
$\gamma\sigma=\gamma\tau$, then $\tau=g\sigma$ for some $g\in G$,
hence $\tau_*^{-1}(z_\gamma)=\sigma_*^{-1}(g_*^{-1}(z_\gamma))=
\sigma_*^{-1}(z_\gamma)$ since $z$ was supposed to be $G$-invariant.
We have thus determined $x=(x_\epsilon)$ for all $\epsilon\in R'=R\cup
\gamma\gS_n$.

(e) Notice that this collection $x=(x_\epsilon)$ is $\gS_n$-invariant 
as required. Indeed, for $\epsilon\in R$ this is evident since $y$ was
supposed to be $\gS_n$-invariant, and for $\epsilon=\gamma\sigma$ we get
$x_{\epsilon\tau}=x_{\gamma\sigma\tau}=(\sigma\tau)_*^{-1}(z_\gamma)=
\tau_*^{-1}(\sigma_*^{-1}(z_\gamma))=\tau_*^{-1}(x_\epsilon)$ as required.

(f) Now let us check that this collection $x\in(\prod_{\epsilon\in R'}
\tilde P(\epsilon))^{\gS_n}$ is compatible, i.e.\ that it lies in
$\projlim_{R'}\tilde P(\epsilon)=\tilde P_*(R')$. This would finish the proof.
So let's prove $x_{\epsilon'}|_\epsilon=x_\epsilon$ for any
$\epsilon\leq\epsilon'\in R'=R\cup\gamma\gS_n$. If $\epsilon'$ lies in $R$,
the same is true for $\epsilon\leq\epsilon'$, $R$ being a sieve, and
the required compatibility is evident since $x_\epsilon=y_\epsilon$,
$x_{\epsilon'}=y_{\epsilon'}$, and $y$ was supposed to be a compatible family.
So assume $\epsilon'\in\gamma\gS_n$. If $\epsilon$ also lies in $\gamma\gS_n$,
then $|\epsilon'|=|\epsilon|$, hence $\epsilon'=\epsilon$ and the compatibility
is trivial, so we assume $\epsilon\in R$. Now the $\gS_n$-invariance of~$x$
shows that $x_{\epsilon'}|_\epsilon=x_\epsilon$ iff
the same is true for $\epsilon\sigma^{-1}\leq\epsilon'\sigma^{-1}$,
i.e.\ we can assume $\epsilon'=\gamma$.
Since $\epsilon\leq\gamma$ and $\epsilon\neq\gamma$, $\epsilon_i<\gamma_i$
at least for one index~$i$; put $\delta:=\gamma-e_i=(\gamma_1,\ldots,
\gamma_i-1,\ldots,\gamma_n)\in W$. Now $\epsilon\leq\delta\leq\gamma$,
$x_\gamma|_\delta=x_\delta$ because both $\delta$ and $\gamma$ lie in $W'$
and $z$ was supposed to be compatible, and $x_\delta|_\epsilon=x_\epsilon$
because both $\delta$ and $\epsilon$ lie in~$R$, q.e.d.
\end{Proof}

\begin{PropD}\label{prop:comput.symm.step}
Let $I\subset\bbZ$, $n\geq0$, $X:I\to\cC$ and 
$\gamma=(\gamma_1,\ldots,\gamma_n)$ be as above, and suppose that
$\gamma_1\leq\cdots\leq\gamma_n$. Put $G:=\Stab_{\gS_n}(\gamma)\subset\gS_n$,
$\tilde W':=\{\gamma_1-1,\gamma_1\}\times
\cdots\times\{\gamma_n-1,\gamma_n\}\subset\bbZ^n$, $\tilde W:=\tilde W'-
\{\gamma\}$, $W':=\tilde W'\cap I^n$, $W:=\tilde W\cap I^n=W'-\{\gamma\}$,
and denote by $v_\gamma$ the morphism $\tilde X_!(W)_G\to\tilde X_!(W')_G$
induced by inclusion $W\subset W'$. Then:
\begin{itemize}
\item[(a)] We can extend $X$ to a functor $\bbZ\to\cC$, without changing
$v_\gamma$. Moreover, we can do it in such a way that $\tilde X_!(W)\cong
\tilde X_!(\tilde W)$ and $\tilde X_!(W')\cong\tilde X_!(\tilde W')$.
\item[(b)] If $\gamma=0=(0,0,\ldots,0)$, then $v_\gamma$ is the morphism
$\emptyset\to S^nX(0)$.
\item[(c)] If $\gamma_1=\gamma_2=\cdots=\gamma_n=k$, then
$v_\gamma$ coincides with $\rho_n(u_k)$, where $u_k:X(k-1)\to X(k)$ is the
transition morphism, provided we put $X(-1)=\emptyset$.
\item[(d)] Suppose that $\gamma_d<\gamma_{d+1}$ for some $0<d<n$, and put
$e:=n-d$. Write $\gamma=\gamma'\star\gamma''$, where $\star$
denotes concatenation of sequences $\star:\bbZ^d\times\bbZ^e\to\bbZ^n$, i.e.\
$\gamma'=(\gamma_1,\ldots,\gamma_d)$ and $\gamma''=(\gamma_{d+1},\ldots,
\gamma_n)$. Then $v_\gamma$ can be identified with $v_{\gamma'}\boxe
v_{\gamma''}$.
\item[(e)] Let $\{k_1,\ldots,k_s\}$ be the set of all distinct components of
$\gamma$, listed in increasing order, and $d_i>0$ be the multiplicity
of $k_i$ in~$\gamma$, so that $\sum_id_i=n$. 
Then $v_\gamma$ can be identified with
$\rho_{d_1}(u_{k_1})\boxe\cdots\boxe\rho_{d_s}(u_{k_s})$, where
$u_{k_i}:X(k_i-1)\to X(k_i)$ are the transition morphisms of~$X$, and
$X(-1):=\emptyset$.
\end{itemize}
\end{PropD}
\begin{Proof}
Notice that (e) follows from (d) and (c) by induction in~$s$, and
(b) is a special case of (c), or just evident by itself. 
In order to show (a) we first extend
$I=[m]$ to $I=\omega$ if $I$ was finite by putting $X(k):=X(m)$ for $k>m$,
and $u_k:X(k-1)\to X(k)$ is of course taken equal to $\id_{X(m)}$.
Then we extend $X$ to $\bbZ$ by putting $X(k):=\emptyset$ for $k<0$.
Now (a) is trivial except for statement $\tilde X_!(W)=\tilde X_!(\tilde W)$; 
but if we compare inductive limits involved we see that the diagram of
the second inductive limit is obtained from that of the first 
by adding some initial
objects and some morphisms, all of them coming from these new initial objects.
Of course, this cannot change the inductive limit. Therefore, we can suppose
$I=\bbZ$ while proving (c) and (d).

Now (c) is immediate since it involves only the restriction of $X:\bbZ\to\cC$
to the subset $\{k-1,k\}\cong[1]$, given by the morphism $u_k:X(k-1)\to X(k)$,
and $W\subset W'$ corresponds exactly to $[1]^n_{\leq n-1}\subset[1]^n$,
hence $v_\gamma=\rho_n(u_k)$ by definition \ptref{sp:can.dec.symm.pow}
of~$\rho_n$.

Let us show the remaining statement (d). Put $G':=\Stab_{\gS_d}(\gamma')$,
$G'':=\Stab_{\gS_e}(\gamma'')$; then clearly $G=G'\times G''$ if
we identify $\gS_d\times\gS_e$ with a subgroup of $\gS_n$ in the usual way.
Next, consider $W_1\subset W'_1\subset I^d$ and
$W_2\subset W'_2\subset I^e$, constructed from $\gamma'$ and~$\gamma''$
in the same way $W\subset W'$ have been constructed from~$\gamma$.
Then $W'=W'_1\times W'_2$ and $W=W'_1\times W_2\cup W_1\times W'_2$.
Now $W'_1\times W_2$ and $W_1\times W'_2$ are two subsieves of 
$W'_1\times W'_2$ with intersection equal to $W_1\times W_2$; 
this implies formally (using only properties of inductive limits) that
the following square is cocartesian:
\begin{equation}
\xymatrix{
\tilde X_!(W_1\times W_2)\ar[r]\ar[d]&\tilde X_!(W'_1\times W_2)\ar[d]\\
\tilde X_!(W_2\times W'_2)\ar[r]&\tilde X_!(W)}
\end{equation}
Clearly, $\tilde X_!(R_1\times R_2)=\tilde X_!(R_1)\otimes\tilde
X_!(R_2)$ for any subsets $R_1\subset I^d$, $R_2\subset I^e$. Applying
this observation to the above diagram we see that $\tilde
X_!(W)\to\tilde X_!(W')= \tilde X_!(W'_1)\times\tilde X_!(W'_2)=\tilde
X(\gamma')\otimes\tilde X(\gamma'')$ is exactly the box-product of
$\tilde X_!(W_1)\to\tilde X_1(W'_1)$ and $\tilde X_!(W_2)\to\tilde
X_2(W'_2)$.  (One might have also seen from the very beginning that
$\tilde X_!(W)\to\tilde X_!(W')$ equals the multiple box product
$u_{\gamma_1}\boxe\cdots\boxe u_{\gamma_n}$ and use the associativity
of box products mentioned in~\ptref{sp:mult.box.prod}; in fact, we
have just proved this associativity.)  Taking coinvariants under
$G=G'\times G''$ with the aid of~\eqref{eq:coinv.tens.prod}, we obtain
$v_\gamma=v_{\gamma'}\boxe v_{\gamma''}$ as required, q.e.d.
\end{Proof}

\begin{CorD}\label{cor:prev.steps.symmpr}
Let $u:X(0)\to X(1)$ be any morphism in~$\cC$, $0<k\leq n$.
Then $\rho_k^{(n)}(u):F_{k-1}S^n(u)\to F_kS^n(u)$ is a pushout of
$\id_{S^{n-k}X(0)}\otimes\rho_k(u)=i\boxe\rho_k(u)$, where
$i$ is the only morphism $\emptyset\to S^{n-k}X(0)$.
\end{CorD}
\begin{Proof}
Put $I:=[1]^n$, $R:=[1]^n_{\leq k-1}$, $R':=[1]^n_{\leq k}$,
$\gamma:=(0,\ldots,0,1,1,\ldots 1)$, $|\gamma|=k$. Then $R\subset
R'=R\cup\gamma\gS_n$ is a simple extension of symmetric sieves,
so we are in position to apply~\ptref{prop:cocart.step.sq}, which tells us that
$u_{R,\gamma}=\rho_k^{(n)}(u):\tilde X_!^s(R)\to\tilde X_!^s(R')$
is a pushout of $v_\gamma:\tilde X_!(W)_G\to\tilde X_!(W')_G$ for
$G=\Stab_{\gS_n}(\gamma)$. Now \ptref{prop:comput.symm.step} yields
$v_\gamma=\id_{S^{n-k}X(0)}\otimes\rho_k(u)$ as required.
\end{Proof}

Now we are ready to formulate the statement we are really going to prove:
\begin{ThD}\label{th:main.prop.rhon}
For any $n\geq1$ functors $\rho_n$ preserve (pointwise) cofibrations and
acyclic cofibrations in $\cC=s\catMod\sO$ and $\stks\stSETS_\cE$.%
\footnote{\rm We provide a complete proof only in the pointwise case}
\end{ThD}
\nxsubpoint (This statement implies the main theorem.)
Before proving the above statement, let us show how it will imply
\ptref{prop:symmprod.steps.cof},
hence also~\ptref{th:ex.der.symmprod}.
Indeed, applying \ptref{th:main.prop.rhon}
to any cofibration $u:\emptyset\to A$ and taking into account that in this
case $F_{n-1}S^n(u)=\emptyset$, and $\rho_n(u):F_{n-1}S^n(u)\to S^n(A)$ can
be identified with $S^n(u):\emptyset\to S^n(A)$, we see that
{\em $S^n(A)$ is cofibrant for any cofibrant~$A$.} Next, let us take
any cofibration (resp.\ acyclic cofibration) $u:A\to B$ between 
cofibrant objects $A$ and $B$. Then all $S^{n-k}A$ are cofibrant,
and all $\rho_k(u)$ are cofibrations (resp.\ acyclic cofibrations)
by~\ptref{th:main.prop.rhon}, hence the same is true for all
$S^{n-k}A\otimes\rho_k(u)$ by (TM), hence also for their pushouts
$\rho_k^{(n)}(u)$, $0<k\leq n$. Since $F_0S^n(u)=S^nA$ is already known to
be cofibrant, we prove by induction in~$k$ that all $F_kS^n(u)$ are
also cofibrant, i.e.\ all $\rho_k^{(n)}$ are cofibrations (resp.\
acyclic cofibrations) between cofibrant objects, exactly as
claimed in~\ptref{prop:symmprod.steps.cof}.

\nxsubpoint\label{sp:proof.main.rhon} {\bf Proof} of \ptref{th:main.prop.rhon}.
Of course, we use our usual sort of devissage, 
explained in~\ptref{sp:fancy.lang}.
Let us denote by $\propP$ the set of morphisms $u:A\to B$ in 
$\cC=s\catMod\sO$ (or rather in the fibers of $\stks\stMOD\sO_\cE$, but 
we can work in a fixed fiber throughout all steps of the proof but one),
such that $\rho_n(u)$ is a cofibration (resp.\ acyclic cofibration)
for all $n\geq 1$ (in particular, $u=\rho_1(u)$ is itself a cofibration,
resp.\ acyclic cofibration). We are going to show that $\propP$ is closed and
that it contains the standard generators 
$L_\sO(\cst I_\cE)$ (resp.\ $L_\sO(\cst J_\cE)$); this would imply the
theorem. So let's check conditions 1)--8) of~\ptref{def:closed.class}
one by one.

(a) Condition 1) (localness of $\propP$) is evident, as well as condition 2),
since $\rho_n(u)$ is an isomorphism whenever $u$ is one. Condition 6)
is also evident: if $u'$ is a retract of $u$, then $\rho_n(u')$ is a retract
of $\rho_n(u)$, retracts being preserved by arbitrary functors 
$\Ar\cC\to\Ar\cC$. Now let's check 3). Suppose that $u'$ is a pushout of~$u$;
then $u'\boxe v$ is easily seen to be a pushout of $u\boxe v$ for any $v$,
hence by induction the $n$-fold box power 
${u'}^{\boxe n}=u'\boxe\cdots\boxe u'$ is a pushout of $u^{\boxe n}$, i.e.\
we have a cocartesian square relating these two morphisms. Applying the
right exact functor of taking $\gS_n$-coinvariants to this square, we
see that $\rho_n(u')$ is a pushout of $\rho_n(u)$, thus proving 3).

(b) Let's check condition 5) (stability under sequential composition).
Let $X(0)\stackrel{u_1}\to X(1)\stackrel{u_2}\to X(2)\to\cdots$ be a composable
sequence of morphisms from $\propP$, considered here as a functor 
$X:I\to\cC$, where $I=\omega$. Denote the sequential composition of 
this sequence by $w:X(0)\to X(\omega):=\injlim_k X(k)$.
We apply the symmetric hypercube construction
to functor~$X$, introducing notations $\tilde X$, $\tilde X_!$ and 
$\tilde X_!^s$ as before. 
Next, we construct an infinite sequence of simple extensions of
symmetric sieves $R_0\subset R_1\subset R_2\subset\cdots$ in~$I^n$ as follows.
Put $R_0:=\{\epsilon\in I^n\,|\,\exists i:\epsilon_i=0\}$, and define
$R_N$ for $N>0$ by induction. Namely, we choose any element $\gamma^{(N)}\in 
I^n-R_{N-1}$ with minimal possible $|\gamma^{(N)}|$; 
since $R_{N-1}$ is symmetric, we can assume that $\gamma^{(N)}_1\leq\cdots\leq
\gamma^{(N)}_n$. Now we put $R_N:=R_{N-1}\cup\gamma^{(N)}\gS_n$; according
to \ptref{sp:struct.simple.syms}, $R_{N-1}\subset R_N$ is indeed a simple
extension of symmetric sieves. Notice that $X_!^s(R_0\cap[m]^n)\to
X_!^s([m]^n)$ can be identified with $\rho_n(u_mu_{m-1}\cdots u_1)$;
taking $\injlim_m$, we see that $X_!^s(R_0)\to X_!^s(I^n)$ can be identified
with $\rho_n(w)$. On the other hand, by construction $\bigcup_NR_N=I^n$, hence
$\injlim_N X_!^s(R_N)=X_!^s(I^n)=S^nX(\omega)$, i.e.\ 
$\rho_n(w):X_!^s(R_0)\to X_!^s(I^n)$ is the sequential composition of morphisms
$w_N:X_!^s(R_{N-1})\to X_!^s(R_N)$. 
According to \ptref{prop:cocart.step.sq} and \ptref{prop:comput.symm.step},
$w_N$ is a pushout of $v_{\gamma^{(N)}}=\rho_{d_1}(u_{k_1})\boxe\cdots
\boxe\rho_{d_s}(u_{k_n})$, where $0<k_1<\cdots<k_s$ are distinct components of
$\gamma^{(N)}$ (all necessarily $>0$ since $\gamma^{(N)}\not\in R_0$)
and $d_i>0$ are their multiplicities, $\sum d_i=n$. (One can even check that
each such box product appears exactly for one value of~$N$.)
Since all $u_k$ were supposed to be in $\propP$, all $\rho_d(u_k)$ are
cofibrations (resp.\ acyclic cofibrations), hence the same is true for
their box products $v_{\gamma^{(N)}}$ by (TM), hence for their pushouts
$w_N$ and the sequential composion $\rho_n(w)$ of $w_N$ as well. Therefore,
$w\in\propP$.

(c) Condition 4) (stabilitity under compositions) can be either reduced to
5) by extending $X(0)\stackrel{u}\to X(1)\stackrel{v}\to X(2)$ to the right
by infinitely many copies of $\id_{X(2)}$, or proved directly by the same
reasoning applied to $I=[2]$. In the latter case the only difference is that
the sequence of symmetric sieves $R_N$ stops after a finite number of 
extensions, i.e.\ we get $R_N=[2]^n$ for some $N$ (actually for $N=n+1$),
thus obtaining a decomposition of $\rho_n(vu)$ into pushouts of morphisms
$\rho_n(u)$, $\rho_{n-k}(u)\boxe\rho_k(v)$, where $0<k<n$, and $\rho_n(v)$.

(d) Let $\{u_i:A_i\to B_i\}_{i\in T}$ be any family of morphisms from~$\propP$.
We want to check that $u:=\bigoplus u_i$ also lies in $\propP$.
Consider the following formula for $S^n(B)$:
\begin{equation}\label{eq:symm.pow.dirsum}
S^n(B)=S^n\bigl(\bigoplus_{i\in T} B_i\bigr)=
\bigoplus\limits_{\alpha\in\bbN_0(T)\atop|\alpha|=n}
\bigotimes\limits_{\alpha_i>0}S^{\alpha_i}B_i
\end{equation}
This formula can be deduced for example from a similar formula for tensor
powers $B^{\otimes n}$:
\begin{equation}\label{eq:tens.pow.dirsum}
B^{\otimes n}=\bigl(\bigoplus_{i\in T}B_i\bigr)^{\otimes n}=
\bigoplus_{\phi\in T^n}\bigotimes_{1\leq i\leq n}B_i^{\otimes\phi_i}
\end{equation}
The latter formula can be shown first for a finite~$T$ by induction, and then
the case of an arbitrary~$T$ is shown by taking a filtered inductive limit 
along finite subsets $T_0\subset T$. Then the formula for $S^n(B)$ is
obtained by re-grouping the direct summands according to $\gS_n$-orbit
decomposition of $T^n$, and computing the $\gS_n$-coinvariants, 
taking into account that $T^n/\gS_n=\{\alpha\in\bbN_0(T):|\alpha|=n\}$
(another natural notation for $T^n/\gS_n$ is $S^nT$ --- symmetric power
in~$\catSets$).

Reasoning in essentially the same way with box products, box powers and
``symmetric box powers'' $\rho_n$, we obtain a similar formula for $\rho_n(u)$:
\begin{equation}
\rho_n(u)=\rho_n\bigl(\bigoplus_{i\in T}u_i\bigr)=
\bigoplus\limits_{\alpha\in\bbN_0(T)\atop|\alpha|=n}
\Box_{\alpha_i>0}\rho_{\alpha_i}(u_i)
\end{equation}
This formula shows immediately that $\rho_n(u)$ is an (acyclic) cofibration
whenever all $\rho_d(u_i)$ are, i.e.\ $\propP$ is closed under direct sums.

(e) Now it remains to check 8), i.e.\ stability of $\propP$ under
``local direct sums'' $\phi_!$, for any $\phi:T\to S$ in~$\cE$. We can 
assume $S=e$ as usual. If we would be able to check cofibrations and 
acyclic cofibrations pointwise, we would simply apply
\begin{equation}
(\phi_!X)_p=\bigoplus_{q\in T_p}X_q
\end{equation}
and reduce everything to (d). Notice, however, that this remark actually 
proves our theorem for the pointwise cofibrant structure on $s\catMod\sO$,
thus allowing us to derive symmetric powers with the aid of pointwise cofibrant
replacements, hence they can be derived with the aid of cofibrant replacements
as well. In general case we need to apply some ``local variant'' of
formula~\eqref{eq:symm.pow.dirsum}; this is not so easy to do because the
proof of \eqref{eq:tens.pow.dirsum}$\Rightarrow$\eqref{eq:symm.pow.dirsum}
is not intuitionistic, so we postpone this verification
for the time being, since we have just explained how to 
prove~\ptref{th:ex.der.symmprod}(b) and (c) without this verification,
using pointwise pseudomodel structures instead.

(f) Once we have shown that $\propP$ is closed, all we need to check is that
it contains the standard generators $L_\sO(\cst I_\cE)$ 
(resp.\ $L_\sO(\cst J_\cE))$. Since symmetric and tensor powers, 
cofibrations and acyclic cofibrations, and
all our other constructions as well, are preserved by base change functors
$L_\sO:s\cE\to s\catMod\sO$ and $q^*:s\catSets\to s\cE$, we are 
reduced to check the statement for $I$ and $J$ in $s\catSets$.
This is done in Lemma~\ptref{l:symmpow.gen.ssets} below, so we
can finish the proof of~\ptref{th:main.prop.rhon}, which is already known
to imply the main theorem~\ptref{th:ex.der.symmprod}.

\begin{LemmaD}\label{l:symmpow.gen.ssets}
(a) Let $u:A\to B$ be any cofibration in $s\catSets$
(e.g.\ a generator $\dot\Delta(n)\to\Delta(n)$ from~$I$). Then
$\rho_m(u)$ is a cofibration in $s\catSets$, for any $m>0$.

(b) Let $u:\Lambda_k(n)\to\Delta(n)$ be a standard acyclic cofibration 
from~$J$, $0\leq k\leq n>0$. Then $\rho_m(u)$ is an acyclic cofibration
in $s\catSets$, for any $m>0$.
\end{LemmaD}
\begin{Proof}
(a) First of all, $u^{\boxe n}$ is a cofibration in $s\catSets$,
and $\rho_n(u)$ is obtained from $u^{\boxe n}$ by taking $\gS_n$-coinvariants.
Recall that cofibrations in $s\catSets$ are exactly the (componentwise)
injective maps of simplicial sets. Now the statement follows from the following
fact: if $i:X\to Y$ is an injective map of $\gS_n$-sets, then 
the induced map of coinvariants (i.e.\ orbit sets) 
$i_{\gS_n}:X/\gS_n\to Y/\gS_n$ is also injective.

(b) We prove the statement by induction in $m>0$, case $m=1$ being trivial.
We know by (a) that $\rho_m(u)$ is a cofibration, and all
$\rho_k(u)$, $0<k<m$, are acyclic cofibrations by induction hypothesis, 
hence the same is true
for all $S^{m-k}A\otimes\rho_k(u)$ and for their pushouts 
$\rho_k^{(m)}(u):F_{k-1}S^m(u)\to F_kS^m(u)$ 
(cf.~\ptref{cor:prev.steps.symmpr}), and we conclude that
$w:=\rho_{m-1}^{(m)}\cdots\rho_2^{(m)}\rho_1^{(m)}:S^mA\to F_{m-1}S^m(u)$
is an acyclic cofibration, and in particular a weak equivalence.
Now $S^m(u)=\rho_m(u)w:S^nA\to S^nB$ by~\ptref{sp:can.dec.symm.pow},
so by the 2-out-of-3 axiom it would suffice to show that $S^m(u)$ is
a weak equivalence.

(c) Let us show that $S^m(u)$ is a weak equivalence, for any 
$u:\Lambda_k(n)\to\Delta(n)$ from~$J$. Consider for this the embedding
$h:\Delta(0)\to\Lambda_k(n)$, defined by the $k$-th vertex of $\Lambda_k(n)$,
and put $h':=hu$:
\begin{equation}
\xymatrix{
\Delta(0)\ar[d]^{h}\ar[rd]^{h'}\\
\Lambda_k(n)\ar[r]^{u}&\Delta(n)}
\end{equation}
By the 2-out-of-3 axiom it would suffice to show that both $S^n(h)$ and
$S^n(h')$ are weak equivalences. Now a morphism in $s\catSets$ is a 
weak equivalence iff this is true for its geometric realization;
using the fact that $|X\times Y|=|X|\times|Y|$ for any two finite simplicial
sets $X$ and~$Y$, and that $|\cdot|$ is right exact, we see that
$|S^nX|=S^n|X|=|X|^n/\gS_n$ for any finite simplicial set~$X$,
so, say, $|S^n(h)|$ can be identified with $S^n(|h|)$.

(d) Notice that both $|h|$ and $|h'|$ are homotopy equivalences.
More precisely, if we denote by $\pi:Z:=|\Lambda_k(n)|\to\pt=|\Delta(0)|$
the only map from $Z$ into a point, then $\pi\circ|h|=\id_{\pt}$, and
$|h|\circ\pi$ is homotopic to $\id_Z$. In fact $Z=|\Lambda_k(n)|$
equals $\{(\lambda_0,\ldots,\lambda_n):$ all $\lambda_i\geq0$,
$\sum_i\lambda_i=1$, $\exists i\neq k:\lambda_i=0\}$, and $|\Delta(n)|$
is described similarly without the last condition. The $k$-th vertex of
$\Delta(n)$, i.e.\ the image of $|h|$ or $|h'|$, is just the
$k$-th basis vector $(0,\ldots,1,\ldots,0)$. Now we define maps
$H_t:|Z|\to|Z|$ for all real $0\leq t\leq 1$ by
\begin{equation}
H_t(\lambda_0,\lambda_1,\ldots,\lambda_k,\ldots,\lambda_n):=
(t\lambda_0,t\lambda_1,\ldots,1-t(1-\lambda_k),\ldots,t\lambda_n)
\end{equation}
Clearly, this formula defines a continuous family of maps $H_t:|\Delta(n)|\to
|\Delta(n)|$, respecting $|\Lambda_k(n)|\subset|\Delta(n)|$,
such that $H_0$ is just the constant map onto the $k$-th vertex, and
$H_1=\id$, i.e.\ $H=\{H_t\}_{0\leq t\leq 1}$ yields a homotopy from
$\id$ to the constant map both for $|\Lambda_k(n)|$ and for $|\Delta(n)|$.

(e) Now we put $\tilde H_t:=S^n(H_t):S^nZ\to S^nZ$ for all $t\in[0,1]$, 
where $Z$ is either $|\Lambda_k(n)|$ or $|\Delta(n)|$. Clearly,
$\tilde H_t$ depends continuously on~$t$, hence it defines a homotopy
$\tilde H$ between $\id_{S^nZ}$ and a constant map, hence 
$S^nZ\to\pt$ and $\pt\to S^nZ$ are weak equivalences, both for
$Z=|\Lambda_k(n)|$ and $Z=|\Delta(n)|$, q.e.d.
\end{Proof}

\cleardoublepage

\mysection{Perfect cofibrations and Chow rings}
\label{chap:K0.Chow}

The aim of this chapter is to define perfect cofibrations and perfect 
simplicial objects in $s\catMod\sO$, where $\cE=(\cE,\sO)$ is an
arbitrary generalized (commutatively) ringed topos, usually assumed to
have enough points. Then we define $K_0$ of the set of isomorphism classes
of perfect objects in $\cD^{\leq0}(\cE,\sO)=\Ho s\catMod\sO$ and
introduce on this $K^0(\cE):=K_0(\Perf(\cE,\sO))$ 
a pre-$\lambda$-ring structure.
We postpone the verification of the fact that this is a $\lambda$-ring 
structure; instead, we provisionally replace $K^0(\cE)$ by its largest
quotient which is a $\lambda$-ring, and use it to construct $\Ch\cE$ as
the associated graded with respect to the $\gamma$-filtration. 
In this way we construct a reasonable intersection and Chern class theory
on any generalized ringed topos, and in particular any generalized scheme.

We finish this chapter by a computation of $K^0(\CompZ)$ and $\Ch(\CompZ)$,
obtaining ``correct'' answers. The case of a projective bundle will be
treated elsewhere.

\zeropoint\nxpoint (Notations.) 
We fix a generalized (commutatively) ringed
topos $\cE=(\cE,\sO)$, usually supposed to have enough points. We denote
by $\cC$ either the pseudomodel category $s\catMod\sO$ of 
simplicial $\sO$-modules, or the corresponding stack $\stks\stMOD\sO_\cE$.
We denote standard finite ordered sets by $[n]=\{0,1,\ldots,n\}$,
and the category of all such finite ordered sets by $\catDelta$; thus
a simplicial object over~$\cD$ is a functor $X:\catDelta^0\to\cD$.
The full subcategory of $\catDelta$ consisting of $[n]$ with $n\leq m$ 
will be denoted by $\catDelta_{\leq m}$ (cf.~\ptref{sp:skel.coskel},
where it has been denoted $\catDelta_m$), and $X_{\leq m}:\catDelta_{\leq m}^0
\to\cD$ will denote the $m$-th truncation of~$X$.

The standard generating sets for cofibrations and acyclic cofibrations
in $s\catSets$ will be still denoted by $I$ and~$J$. Thus
$I=\{\dot\Delta(n)\to\Delta(n):n\geq0\}$ and $J=\{\Lambda_k(n)\to\Delta(n):
0\leq k\leq n>0\}$, and $L_\sO(\cst I_\cE)$ and $L_\sO(\cst J_\cE)$ are
our generators for the pseudomodel structure of $s\catMod\sO$.

\nxpointtoc{Simplicial dimension theory}\label{p:dim.theory}
Up to now we didn't use much the fact that $\cC=s\catMod\sO$, usually only
using the fact that we can reason componentwise, and sometimes reducing
the statements being proved to basic case $s\catSets$, where the
needed statements are usually quite well known. Now we'd like to 
construct a {\em simplicial\/} or {\em homotopic dimension theory\/} for
objects and morphisms in~$\cC$, together with canonical dimension 
decompositions, similar to the theory of brute truncations of a complex
of modules over an additive sheaf of rings. 

Several first definitions will be in fact valid in any category $\cC=s\cD$
of simplicial objects over a category $\cD$ with finite inductive limits:

\nxsubpoint 
(Morphisms concentrated in dimension $>n$ and $\leq n$.)
Let $f:X\to Y$ be any morphism in $\cC=s\cD$. We say that
{\em $f$ is concentrated in dimension $>n$}, if $f$ induces an isomorphism
$f_{\leq n}:X_{\leq n}\to Y_{\leq n}$ of $n$-truncated objects.
This is equivalent to saying that all $f_i:X_i\to Y_i$ are isomorphisms for
$0\leq i\leq n$, or that $\sk_n(f):\sk_nX\to\sk_nY$ is an isomorphism in~$\cC$.

We say that {\em $f$ is concentrated in dimension $\leq n$}, or that
{\em $f$ has (relative) dimension $\leq n$}, if the following square is 
cocartesian, where $\sk_nX\to X$ and $\sk_nY\to Y$ are the canonical 
morphisms (cf.~\ptref{sp:skel.coskel}):
\begin{equation}
\xymatrix{
\sk_nX\ar[d]\ar[r]^{\sk_nf}&\sk_nY\ar[d]\\
X\ar[r]^{f}&Y}
\end{equation}
Since $\sk_n(f)$ is completely determined by the truncation $f_{\leq n}$,
we see that {\em any $f:X\to Y$ of relative dimension $\leq n$ is completely
determined (up to a unique isomorphism) by $X$ and~$f_{\leq n}$, i.e.\ by the
collection $f_i:X_i\to Y_i$, $i\leq n$, if $X$ and~$Y$ are already fixed.}

We extend the above notions to any $n\in\bbZ$, possibly negative, by putting
$\catDelta_n:=\emptyset$, $\sk_nX:=\emptyset_\cC$ for all $n<0$. 
When $f$ is concentrated in dimension $>n$, we also say that it is 
concentrated in dimension $\geq n+1$, and so on. If $f$ is concentrated 
in dimension $\leq b$ and $\geq a$, we say that {\em $f$ is concentrated in
dimensions $[a,b]$}. When $a=b$, we say that {\em $f$ is pure(ly) 
of dimension $a$} instead. 
We allow here $a=-\infty$ and $b=+\infty$ as well.

When $n<0$, any $f$ is concentrated in dimension
$>n$, but only isomorphisms are concentrated in dimension $\leq n$.

One checks very easily using formula $\sk_n\circ\sk_m=\sk_{\min(m,n)}$ for
any $m$, $n\in\bbZ$, itself an immediate consequence of the universal
property of $\sk_nX$,
that whenever $m\leq n$, condition
``$f$ is concentrated in dimension $\leq m$'' implies that
``$f$ is concentrated in dimension $\leq n$'', and 
``$f$ is concentrated in dimension $>n$'' implies ``$f$ is concentrated in
dimension $>m$''. In particular, we can define the 
{\em (relative) dimension $\dim f$ of~$f$} as the infimum of all
$n\in\bbZ$, such that $f$ is concentrated in dimension $\leq n$. Then
{$f$ is concentrated in dimension $\leq n$ iff $\dim f\leq n$.}
Another immediate consequence is that {\em if $f$ is concentrated in
dimensions $[a,b]$ and $[a,b]\subset[a',b']$, then $f$ is concentrated in
dimensions $[a',b']$}.

Finally, we extend all the above notions to objects $X\in\Ob\cC$ by 
considering $\nu_X:\emptyset\to X$, e.g.\ $\dim X:=\dim\nu_X$ and so on.
The only object of dimension $<0$ is the initial object $\emptyset$,
which has $\dim\emptyset=-\infty$, and the simplicial objects of
dimension $\leq 0$ (hence pure of dimension~$0$) are the constant simplicial
objects, since $\sk_0X$ is the constant simplicial object $X_0$ for any~$X$.

\nxsubpoint 
(Compatibility with composition and inductive limits.)
One checks immediately from definitions that the above properties are stable
under composition, and in one case even satisfy the 2-out-of-3 property:
{\em if two of $f$, $g$, $g\circ f$ are concentrated in dimension
$>n$, then so is the third. If $g\circ f$ and $f$ are of relative dimension
$\leq n$ (resp.\ are concentrated in $[a,b]$), then the same holds for~$g$.}
In particular, if $X$ and $Y$ are of dimension $\leq n$, the same is true
for any $f:X\to Y$.

Since $\sk_n$ commutes with arbitrary inductive limits
(and functors $\phi^*$ and $\phi_!$ in stack $\stks\stMOD\sO_\cE$ for
any $\phi:T\to S$, if we want to consider stacks here), we see that the above
properties of simplicial objects are stable under all sorts 
of inductive limits, and the property of a morphism $f$ to be 
concentrated in $[a,b]$ is stable under pushouts, composition, sequential
composition, retracts, direct sums (and local direct sums $\phi_!$ in the stack
situation), hence {\em morphisms concentrated in dimension $[a,b]$
constitute a closed class in $\cC=s\catMod\sO$ or $\cC=\stks\stMOD\sO_\cE$}
in the sense of~\ptref{def:closed.class}.

\nxsubpoint
(Example: standard generators of $s\catSets$.)
It is easy to see that the standard simplex $\Delta(n)\in\Ob s\catSets$ 
is of dimension~$n$ (i.e.\ $\dim X=n$),
hence any finite simplicial set has finite dimension, equal in fact to
the highest dimension in which it has non-degenerate simplices.
 
Consider a standard cofibration $i_n:\dot\Delta(n)\to\Delta(n)$ from~$I$.
It is of dimension $\leq n$ because both sides are. On the other hand,
$\sk_{n-1}\Delta(n)=\dot\Delta(n)$, i.e.\ $i_n$ is of dimension $>n-1$,
hence {\em $i_n$ is pure of dimension~$n$.} One checks similarly that
the standard acyclic cofibrations $\Lambda_k(n)\to\Delta(n)$,
$0\leq k\leq n>0$, are concentrated in dimension $[n-1,n]$.

{\bf Notation.} We denote by $I_{\leq n}$, $I_{\geq n}$, $I_{[a,b]}$, $I_S$
(where $S\subset\bbZ$) and $I_{=n}=I_{[n,n]}$ 
the subsets of~$I$ consisting of generators of specified dimensions.

\nxsubpoint
(Functorial dimensional decompositions.)
We claim that {\em any morphism $f:X\to Y$ can be uniquely 
(up to a unique isomorphism) factorized into $X\stackrel g\to Z\stackrel h\to 
Y$, where $g$ is of dimension $\leq n$ and $h$ is concentrated in
dimension $>n$.} Indeed, $\sk_n h$ has to be an isomorphism, i.e.\ 
all $h_i:Z_i\to Y_i$, $i\leq n$, are isomorphisms; we can assume, replacing
$Z_i$ by isomorphic objects $Y_i$ if necessary, that $h_i=\id_{Y_i}$,
for all $i\leq n$. Then $h_{\leq n}=\id$, hence $\sk_n h=\id$ and
$\sk_n f=\sk_n h\sk_n g=\sk_ng$, i.e.\ $\sk_n g$ is known; since
$g$ is of dimension $\leq n$, this completely determines $g$, which
can be identified with the pushout of $\sk_ng=\sk_nf$ with respect to
$\sk_nX\to X$. In this way we obtain uniqueness; conversely, we can
always define $g=\sk_{\leq n}(f)$ and $h=\sk_{>n}(f)$ with required 
properties by means of the following cocartesian square:
\begin{equation}\label{eq:constr.dim.decomp}
\xymatrix{
\sk_nX\ar[r]^{\sk_nf}\ar@{-->}[d]&\sk_nY\ar[d]\ar@/^/[ddr]\\
X\ar@{-->}[r]^{\sk_{\leq n}f}\ar[rrd]_{f}&Z\ar@{-->}[rd]|{\sk_{>n}f}\\
&&Y}
\end{equation}
Therefore, the above factorization $f=\sk_{>n}(f)\circ\sk_{\leq n}(f)$
is functorial in~$f$, and commutes with arbitrary inductive limits.

Now suppose that $n<m$. We can construct a decomposition of $f:X\to Y$
into $X\stackrel u\to W\stackrel v\to Z\stackrel w\to Y$, where
$u$ is concentrated in $\leq n$, $v$ in $[n+1,m]$ and $w$ in $>m$, for example
by putting $u:=\sk_{\leq n}(f)$, $v:=\sk_{\leq m}\sk_{>n}(f)$, 
$w:=\sk_{>m}\sk_{>n}(f)=\sk_{>m}(f)$. Such decompositions are also unique
(apply previous result first to $f=(wv)u$, then to $wv$), and since
$u=\sk_{\leq n}\sk_{\leq m}(f)=\sk_{\leq n}(f)$, $v=\sk_{>n}\sk_{\leq m}(f)$,
$w=\sk_{>m}(f)$ constitute another such decomposition, we have
$\sk_{\leq m}\sk_{>n}=\sk_{>n}\sk_{\leq m}$ whenever $n<m$. We put
\begin{align}
\sk_{[n,m]}(f):=&\sk_{\leq m}\sk_{\geq n}(f)=\sk_{\geq n}\sk_{\leq m}(f),
\quad\text{for any $n\leq m$,}\\
\sk_{=n}(f):=&\sk_{[n,n]}(f)\quad\text{for any $n$.}
\end{align}
Since $\sk_{\leq n}(f)=\sk_{=n}(f)\circ\sk_{\leq n-1}(f)$,
$\sk_{<0}(f)=\id$, and $\injlim_n\sk_{\leq n}(f)=f$, we obtain a 
canonical decomposition of $f$ into a sequential composition of 
morphisms $\sk_{=n}(f)$, each of pure dimension~$n$:
\begin{equation}\label{eq:full.dim.decomp}
X=F_0X\stackrel{\sk_{=0}f}\longto F_1X\stackrel{\sk_{=1}f}\longto F_2X
\to\cdots F_\infty X=Y
\end{equation}
The above decomposition is functorial in~$f$, unique up to an isomorphism,
and we have e.g.\ $\sk_{[a,b]}(f)=\sk_{=b}(f)\circ\sk_{=b-1}(f)\circ\cdots
\circ\sk_{=a}(f)$. In the additive case this decomposition corresponds
to an increasing filtration $F_kP$ on a chain complex~$P$, ``cokernel'' of~$f$,
namely, the brute filtration: $(F_kP)_n=0$ for $n\geq k$, $(F_kP)_n=P_n$ for 
$n<k$.

\nxsubpoint 
(Extension of terminology.)
Given any subset $S\subset\bbZ$, we say that {\em $f$ is concentrated in
dimension(s)~$S$} if $\sk_{=n}(f)$ are isomorphisms for all $n\not\in S$,
i.e.\ if the corresponding steps of~\eqref{eq:full.dim.decomp} are trivial.
In this case $f$ is a (finite or sequential) composition of morphisms
isomorphic to $\sk_{=n}(f)$, for elements $n\in S$ listed in increasing order.

Clearly, this definition is compatible with our previous terminology
when $S=[a,b]$. 

\begin{PropD}\label{prop:cof.bounded.dim}
(a) If a morphism $f:X\to Y$ belongs to $\Cl L_\sO(\cst I_S)$, where 
$S\subset\bbZ$ is any subset, then $\sk_n(f)$ and its pushout
$\sk_{\leq n}(f)$ belong to $\Cl L_\sO(\cst I_{S\cap[0,n]})$, and
$\sk_{>n}(f)$ belongs to $\Cl L_\sO(\cst I_{S\cap[n+1,\infty)})$,
hence $\sk_{[a,b]}(f)$ belongs to $\Cl L_\sO(\cst I_{S\cap[a,b]})$.

(b) If $f:X\to Y$ is a cofibration, then $\sk_n f$, $\sk_{\leq n}f$, 
$\sk_{>n}f$, $\sk_{[a,b]}f$ and $\sk_{=n}f$ are also cofibrations.
Moreover, $\sk_{[a,b]}f$ lies in the closure of $L_\sO(\cst I_{[a,b]})$.

(c) The set of cofibrations concentrated in dimension $[a,b]$ coincides
with the closure (in~$\cC=\stks\stMOD\sO_\cE$) of $L_\sO(\cst I_{[a,b]})$.

(d) More generally, the set of cofibrations concentrated in dimensions 
$S\subset\bbZ$ coincides with the closure of $\cst I_S$.
\end{PropD}
\begin{Proof} (a) is reduced by our usual ``devissage'' 
(cf.~\ptref{sp:fancy.lang}) to the case of a standard cofibration 
$i_n:\dot\Delta(n)\to\Delta(n)$ in $s\catSets$ with $n\in S$, where it is
immediate, $i_n$ being purely of dimension~$n$. The ability to apply 
devissage to $\sk_n$ is due to the fact that 
$\sk_n$ commutes with any inductive limits,
composition, sequential composition, local direct sums ($\phi_!$) and so on.
Case of $\sk_{>n}$ is not much more complicated: the only difference is that
now $\sk_{>n}(gf)=\sk_{>n}g\circ($some pushout of $\sk_{>n}f)$, similarly to
what we had for $u\mapsto u\boxe i$ in~\ptref{l:cl.box.div.class}.
We can emphasize this similarity by means of a statement that contains 
both these situations, cf.~\ptref{l:cocart.defect} below.

(b) is now immediate from (a), and (c) follows from (b): one inclusion is due
to the fact that morphisms concentrated in dimensions $[a,b]$ constitute
a closed class, and the other inclusion follows from (b) and the fact that
$f$ is concentrated in dimensions $[a,b]$ iff $f$ is isomorphic to
$\sk_{[a,b]}f$.

Let's show (d). If $f:X\to Y$ lies in~$\Cl L_\sO(I_S)$, 
then by (a) $\sk_{=n}(f)$ lies in $\Cl L_\sO(I_{S\cap\{n\}})$; 
for $n\not\in S$ this means that $\sk_{=n}(f)$ is
an isomorphism, i.e.\ the $n$-th step of dimensional decomposition
is trivial, hence $f$ is concentrated in dimensions from $S$ by definition.
Conversely, if $f$ is a cofibration concentrated in dimensions from~$S$,
then it is a (finite or sequential) composition of morphisms $\sk_{=n}(f)$
for $n\in S$; each of these morphisms lies in $\Cl L_\sO(I_{=n})\subset
\Cl L_\sO(I_S)$ by (b), hence the same is true for their composition~$f$.
\end{Proof}

\begin{LemmaD}\label{l:cocart.defect}
Let $\cC$ and $\cD$ be any two categories closed under (small) 
inductive limits, $F$, $G:\cC\to\cD$ be two functors commuting with 
arbitrary inductive limits, and $\eta:F\to G$ be a natural transformation.
For any $u:X\to Y$ in $\cC$ denote by $Z(u)$ the fourth vertex of the
cocartesian square built on $F(u)$ and $\eta_X$, and by $\eta[u]:Z(u)\to G(B)$
the natural morphism, thus defining functors $Z:\Ar\cC\to\cD$ and
$\eta[\cdot]:\Ar\cC\to\Ar\cD$. Then:
\begin{itemize}
\item[(a)] $Z(\id_Y)=F(Y)$, and $Z((u,\id_Y)):Z(u)\to Z(\id_Y)$ can be 
identified with $\eta[u]$.
\item[(b)] For any composable $X\stackrel u\to Y\stackrel v\to W$ 
the morphism $Z(u,\id_W):Z(vu)\to Z(v)$ is a pushout of 
$Z(u,\id_Y)=\eta[u]$, and $\eta[vu]=\eta[v]\circ Z(u,\id_W)$.
\item[(c)] If $u':X'\to Y'$ is a pushout of $u:X\to Y$, then
$\eta[u']$ is a pushout of $\eta[u]$.
\item[(d)] $\eta\bigl[\bigoplus_\alpha u_\alpha\bigr]=
\bigoplus_\alpha\eta[u_\alpha]$.
\item[(e)] If $u$ is a sequential composition $\cdots u_2u_1u_0$, then
$\eta[u]$ is sequential composition of pushouts of $\eta[u_n]$.
\end{itemize}
\end{LemmaD}
\begin{Proof}
Statement (a) is immediate, and (b) is shown by a simple diagram chase, 
not involving any properties of~$F$ or $G$. Statement (c) is shown 
by another diagram chase, using that $F$ and~$G$ preserve cocartesian 
squares. Remaining statements are then shown in the same way
as in~\ptref{l:cl.box.div.class}.
\end{Proof}

\begin{PropD}\label{prop:dim.tens.prod}
(a) If $i:A\to B$ and $s:X\to Y$ are cofibrations of dimension $\leq n$ and
$\leq m$, respectively, then $i\boxe s$ is a cofibration of dimension 
$\leq n+m$. In particular, $\dim(i\boxe s)\leq\dim i+\dim s$ for any two 
cofibrations $i$ and~$s$.

(b) If $X$ and $Y$ are cofibrant, then $\dim(X\otimes Y)\leq\dim X+\dim Y$.
\end{PropD}
\begin{Proof} (b) is immediate from (a) since $\nu_{X\otimes Y}=\nu_X
\boxe\nu_Y$. Let's prove (a). 
By \ptref{prop:cof.bounded.dim},(c) we see that
$i$ lies in the closure of $L_\sO(\cst I_{\leq n})$, and $s$ in that
of $L_\sO(\cst I_{\leq m})$. Applying devissage for the box product,
first in one variable and then in the other (application of devissage
is made possible by lemma~\ptref{l:cl.box.div.class}), we are finally reduced
to checking that $w:=i_n\boxe i_m$ is a cofibration of dimension $\leq n+m$ in
$s\catSets$, for any $n$, $m\geq0$. First of all, we see that the target 
$\Delta(n)\otimes\Delta(m)=\Delta(n)\times\Delta(m)$ of~$w$ is of dimension
$n+m$: its simplices are order-preserving maps $[k]\to[n]\times[m]$,
and non-degenerate simplices correspond to injective such maps, which 
can exist only for $k\leq m+n$. Next, $w$ is a cofibration, hence componentwise
injective; this immediately implies that it preserves non-degenerate simplices,
hence its source can have non-degenerate simplices only in dimension 
$\leq n+m$, i.e.\ both the source and the target of~$w$ are of dimension 
$\leq n+m$, hence the same holds for~$w$ itself.
\end{Proof}

\nxsubpoint\label{sp:dim.pb.cof}
One shows even simpler that {\em for any morphism of generalized ringed topoi
$f:(\cE',\sO')\to(\cE,\sO)$, functor $f^*$ preserves cofibrations 
and cofibrant objects concentrated in dimension $[a,b]$.} This is actually
true for arbitrary morphisms or objects of $s\catMod\cE$ 
(not just cofibrant), just because the components of $\sk_n$, 
$\sk_{\leq n}$ and $\sk_{>n}$ are computed with the aid of finite inductive
limits only.

\nxsubpoint\label{sp:const.cof}
(Constant cofibrations. Components of cofibrations.)
We say that a morphism $f_0:X_0\to Y_0$ in $\catMod\sO$ is a 
{\em (constant) cofibration\/}
if it becomes a cofibration in $s\catMod\sO$ when considered as a morphism
of constant simplicial objects. Since constant objects and morphisms
are exactly those of dimension $\leq0$, we see that $f_0$ is a cofibration
iff it lies in $\Cl L_\sO(\cst I_{=0})$. Since $Z\mapsto Z_0$ commutes
with any limits, our usual devissage argument shows that (constant) 
cofibrations in $\catMod\sO$ are {\em exactly\/} morphisms from the
closure in~$\catMod\sO$ of one-element set $\emptyset_\sO\to L_\sO(1)$. 
For example, cofibrations in $\catSets$ are just the injective maps.

Now functor $P_n:s\catMod\sO\to\catMod\sO$, $Z\mapsto Z_n$, 
also preserves any limits, so we can apply  
devissage again and prove that {\em whenever $f:X\to Y$ is a cofibration in
$s\catMod\sO$, all its components $f_n:X_n\to Y_n$ are cofibrations
in $\catMod\sO$, i.e.\ they belong to $\Cl\{\emptyset_\sO\to L_\sO(1)\}$.}

\begin{PropD}
Let $(\cE,\sO)$ be a generalized ringed topos, not necessarily commutative.
Then:
\begin{itemize}
\item[(a)] All components $f_n:X_n\to Y_n$ of a cofibration $f:X\to Y$
in $s\catMod\sO$ are monomorphisms.
\item[(b)] Any constant cofibration $f_0:X_0\to Y_0$ in $\catMod\sO$
is a monomorphism. 
\end{itemize}
\end{PropD}
\begin{Proof} Statement (a) follows immediately from (b), since we have
just shown that all components of a cofibration are constant cofibrations.
Let's prove (b). We assume for simplicity that $\cE$ has enough points, 
so it is enough to prove (b) over $\cE=\catSets$. Then we can assume
that $\sO$ has at least one constant, i.e.\ $\sO(0)=\sO^{(0)}=\emptyset_\sO$ 
is non-empty: otherwise
we can do scalar extension with respect to 
$r:\sO\mapsto\sO\boxtimes\bbF_1=\sO\langle c^{[0]}\rangle$ 
(here we use 
the non-commutative tensor product $\boxtimes$ of~\ptref{sp:nc.tensprod},
i.e.\ $c$ is not required to commute with any operations of~$\sO$).
Then $M\to r_*r^*M$ is injective for any $\sO$-module~$M$, since it can
be identified with embedding $M\to M\oplus L_\sO(1)$, which admits 
a section when $M$ is non-empty, and is injective for trivial reasons if~$M$
is empty (this reasoning wouldn't be valid over an arbitrary topos~$\cE$, 
where we would
have to decompose $\cE$ into the open subtopos over which $M$ locally admits
a section, and its closed complement). Therefore, injectivity of
$r^*(i)$ would imply injectivity of $i$, and $r^*$ preserves cofibrations,
so the reduction step is done.

So let us assume that $\sO$ admits a constant $c$. We claim that then
{\em any constant cofibration $f_0:X_0\to Y_0$ admits a left inverse},
hence is injective. This statement is shown again by ``devissage'',
using that constant cofibrations coincide with the closure of 
$i_0:L_\sO(0)\to L_\sO(1)$ in $\catMod\sO$, and that $i_0$ admits a left
inverse~$\pi_0$, defined e.g.\ by the constant $c\in L_\sO(0)$ (elements of $M$
are in one-to-one correspondence with $\sO$-homomorphisms $L_\sO(1)\to M$,
and condition $\pi_0 i_0=\id$ is automatic, 
$L_\sO(0)=\emptyset_\sO$ being initial).
\end{Proof}

\begin{PropD}\label{prop:compw.crit.cof}
Morphism $f:X\to Y$ is a cofibration iff all $(\sk_{\geq n}f)_n$ are constant
cofibrations iff all $(\sk_{=n}f)_n$ are constant cofibrations, where
$n\geq0$ runs through all non-negative integers.
\end{PropD}
\begin{Proof}
We know that $\sk_{\geq n}f$ and $\sk_{=n}f$ are cofibrations whenever
$f$ is one (cf.~\ptref{prop:cof.bounded.dim}), and that all components of
a cofibration are constant cofibrations (cf.~\ptref{sp:const.cof}),
hence the first condition implies the other two. On the other hand,
$\sk_{=n}f=\sk_{\leq n}\sk_{\geq n}f$ is a pushout of
$\sk_n\sk_{\geq n}f$, hence $(\sk_{=n}f)_n$ is a pushout 
of~$(\sk_{\geq n}f)_n$, i.e.\ the second condition implies the third.

Now suppose that $(\sk_{=n}f)_n$ is a constant cofibration. Since
$f$ is a sequential composition of the $\sk_{=n}f$, it would suffice to show
that $\sk_{=n}f$ is a cofibration, for all $n\geq0$. Putting $g:=\sk_{=n}f$,
we are reduced to checking the following statement: {\em if
$g:X'\to Y'$ is purely of dimension $n$, and $g_n:X'_n\to Y'_n$ is a 
constant cofibration, then $g$ is itself a cofibration.}
According to Lemma~\ptref{l:pshout.dec.pure.mor} below, any 
$g:X'\to Y'$ purely of dimension~$n$ is a pushout of $g_n\boxe i_n=
g_n\boxe L_\sO(i_n)$, where $i_n:\dot\Delta(n)\to\Delta(n)$ is the
standard cofibrant generator. We know already that $\boxe$ respects
cofibrations (cf.~\ptref{prop:comp.tens.sO-mod}), whence the statement.
\end{Proof}

During the proof we have obtained the following interesting statement:
\begin{CorD}
A morphism $f:X\to Y$ is a cofibration purely of dimension~$n$ iff
it is a pushout of $u\boxe i_n$, for some constant cofibration~$u$.
\end{CorD}

\begin{LemmaD}\label{l:pshout.dec.pure.mor}
If $f:X\to Y$ is purely of dimension~$n\geq0$, then it can be canonically
identified with a pushout of $f_n\boxe i_n=f_n\boxe L_\sO(i_n)$,
where $f_n:X_n\to Y_n$ is the $n$-th component of~$f$, considered here as
a morphism of constant simplicial objects, and $i_n:\dot\Delta(n)\to
\Delta(n)$ is the standard cofibrant generator of $s\catSets$.
\end{LemmaD}
\begin{Proof}
First of all, $f$ is a pushout of $\sk_n f$, being of dimension $\leq n$,
$f':=\sk_n f$ is still purely of dimension~$n$, and $f_n=(\sk_n f)_n$, so
we can replace $f$ with $\sk_n f$ and assume that $X$ and $Y$ are of
dimension $\leq n$. Now consider the following commutative diagram:
\begin{equation}
\xymatrix{
X_n\otimes\dot\Delta(n)\ar[rr]^{f_n\otimes\id_{\dot\Delta(n)}}
\ar[d]^{\id_{X_n}\otimes i_n}
&&Y_n\otimes\dot\Delta(n)\ar[dddl]_(.7){w}\ar@{-->}[d]^{\tilde\imath}
\ar[rdd]^{\id_{Y_n}\otimes i_n}\\
X_n\otimes\Delta(n)\ar@{-->}[rr]^{\tilde f}
\ar[rrrd]_(.7){f_n\otimes\id_{\Delta(n)}}\ar[ddr]_{\ev_{\Delta(n),X}}
&&Z\ar@{-->}[rd]|<>(.5){f_n\boxe i_n}\ar@{-->}[ddl]^{h}\\
&&&Y_n\otimes\Delta(n)\ar[d]^{\ev_{\Delta(n),Y}}\\
&X\ar[rr]^{f}&&Y}
\end{equation}
Here $\ev_{\Delta(n),X}:X_n\otimes\Delta(n)\to X$ is the ``evaluation map'',
defined by adjointness from the canonical isomorphism 
$\iHom(\Delta(n),X)\cong X_n$, and $\ev_{\Delta(n),Y}$ is defined similarly,
so the only solid arrow left unexplained is $w:Y_n\otimes\dot\Delta(n)\to X$.
By adjointness it suffices to construct $w^\flat:Y_n\to\iHom(\dot\Delta(n),X)$.
First of all, we have a canonical morphism $f_*:\iHom(\dot\Delta(n),X)\to
\iHom(\dot\Delta(n),Y)$; it is even an {\em iso}morphism since 
$\iHom(\dot\Delta(n),X)$ is a certain finite projective limit of components
$X_k$ of~$X$ with $k\leq n-1$, and all $f_k:X_k\to Y_k$, $0\leq k\leq n-1$,
are isomorphisms, $f$ being purely of dimension~$n$ by assumption.
Composing the inverse of this isomorphism $f_*$ with 
$i_n^*:Y_n=\iHom(\Delta(n),Y)\to\iHom(\dot\Delta(n),Y)$, we obtain 
our~$w^\flat$.

Now it is immediate that the subdiagram of the above diagram consisting of
solid arrows is commutative, so we can construct $Z$ by means of a cocartesian
square, and define remaining arrows $\tilde\imath$, $\tilde f$,
$f_n\boxe i_n$ and $h$, thus completing the construction of the diagram.

We want to show that the square of this diagram with vertices in $Z$,
$Y_n\otimes\Delta(n)$, $X$ and~$Y$ is cocartesian. Since all simplicial 
objects involved in the diagram are of dimension $\leq n$, it would
suffice to show that this square becomes cocartesian after truncation
to dimension $\leq n$, i.e.\ that its components in dimensions $k\leq n$
constitute a cocartesian square. When we truncate to dimensions $\leq n-1$,
$i_n$ becomes an isomorphism, hence the same is true for $\id_{X_n}\otimes
i_n$, $\id_{Y_n}\otimes i_n$, $\tilde\imath$, hence $(f_n\boxe i_n)_{\leq n-1}$
is also an isomorphism, and $f_{\leq n-1}$ is an isomorphism as well,
$f$ being concentrated in dimension $>n-1$, so our square is cocartesian
in dimensions $<n$ for trivial reasons.

Now it remains to compute the components in dimension~$n$. First of all,
$(\Delta(n))_n$ consists of all order-preserving maps $\phi:[n]\to[n]$,
and $(\dot\Delta(n))_n$ consists of all such non-surjective maps,
i.e.\ $(\Delta(n))_n=(\dot\Delta(n))_n\sqcup\{\sigma\}$, where
$\sigma=\id_\stn$ is the only non-degenerate $n$-dimensional simplex 
of~$\Delta(n)$. This immediately implies that $(\id_{X_n}\otimes i_n)_n$
can be identified with the canonical embedding $X_n\otimes\dot\Delta(n)_n\to
\bigl(X_n\otimes\dot\Delta(n)_n\bigr)\oplus X_n$, hence $\tilde\imath$ can be
identified with $Q:=Y_n\otimes\dot\Delta(n)_n\to Q\oplus X_n\cong Z_n$, 
and a similar description is valid for $(\id_{Y_n}\otimes i_n)_n$,
with the second summand equal to~$Y_n$ instead of~$X_n$. This means that
we can identify $(f_n\boxe i_n)_n$ with $\id_Q\oplus f_n$, i.e.\
we get the following commutative diagram:
\begin{equation}
\xymatrix@C+13pt{
X_n\ar[d]^{f_n}\ar[r]^<>(.5){\otimes\sigma}&
Q\oplus X_n\ar[d]^{(f_n\boxe i_n)_n}\ar[r]^<>(.5){h_n}&X_n\ar[d]^{f_n}\\
Y_n\ar[r]^<>(.5){\otimes\sigma}&
Q\oplus Y_n\ar[r]^<>(.5){(\ev_{\Delta(n),Y})_n}&Y_n}
\end{equation}
The composites of horizontal arrows are equal to identity 
(because $\sigma=\id_\stn$), hence the outer circuit is cocartesian,
and the left square is cocartesian because $(f_n\boxe i_n)_n=\id_Q\oplus f_n$,
hence the right square is also cocartesian, q.e.d.
\end{Proof}

\nxpointtoc{Finitary closures and perfect cofibrations}
Now we are going to define and study the basic properties of 
{\em perfect cofibrations\/} and {\em perfect (simplicial) objects},
over any generalized commutatively ringed topos $\cX=(\cX,\sO)$, 
sometimes supposed to have enough points. The basic idea behind these 
definitions is the following. We define {\em finitarily closed\/}
sets of morphisms in a category or in fibers of a stack by imposing
``finitary'' versions of conditions 1)--8) of~\ptref{def:closed.class},
e.g.\ we allow only finite direct sums and finite compositions. Then
we can define {\em finitary closure\/} of any set of morphisms; 
and the {\em perfect cofibrations\/} will be exactly morphisms from
the finitary closure of the set $L_\sO(\cst I_\cX)$ of standard cofibrant
generators for~$\stks\stMOD\sO_\cX$, and {\em perfect\/} or
{\em perfectly cofibrant objects\/~$X$} will be the objects for which
$\nu_X:\emptyset\to X$ is a perfect cofibration.

Then we are going to prove some basic properties of perfect cofibrations
and perfect objects, essentially applying the same (or even simpler)
``devissage'' as before (cf.~\ptref{sp:fancy.lang}). Since most of the proofs
would just repeat some of the arguments applied before to prove 
properties of cofibrations, we are going to give references to our 
former proofs instead of re-proving the statements for the perfect case.

\zerosubpoint\nxsubpoint\label{sp:retract.ver}
(Retract and retract-free versions of the theory.)
The theory developped below can be actually constructed in two different 
fashions, depending on whether we allow retracts in the definition of
finitary closed sets ({\em ``retract version''} of the theory)
or not ({\em ``retract-free version''}). We'll usually develop only
the retract version, leaving the retract-free version to the reader.
Usually the retract condition in the definition of finitarily closed sets
will be used only to propagate itself, so the proofs for the retract-free
theory are just sub-proofs of those considered by us. 

However, when there will be some subtle points where the retract and
the retract-free variants differ, we'll point this out explicitly.
One of such examples will be encountered during the consideration of
affine generalized schemes, since a locally free $\sO$-module needn't
be globally free, while a locally projective $\sO$-module is 
globally projective. There is another important point affected 
by the choice of theory:

{\bf Terminology.} A {\em vector bundle~$\sE$} over $(\cX,\sO)$ will be
just any (locally) finitely presented locally projective $\sO$-module
(resp.\ locally free in the retract-free version of the theory).

\begin{DefD}\label{def:finclosed.class}
Let $\cC$ be a stack over a site~$\cS$, such that finite inductive limits
exist in each fiber of $\cC/\cS$ and all pullback functors 
$\phi^*:\cC(S)\to\cC(T)$ are right exact, and $\propP\subset\Ar\cC$ be 
any set of morphisms lying in the fibers of~$\cC$. We say that
$\propP$ is {\em finitarily closed\/} if the following conditions hold
(cf.~\ptref{def:closed.class}):
\begin{itemize}
\item[1)] $\propP$ is local, and in particular stable under base change
$\phi^*:\cC(S)\to\cC(T)$ with respect to any $\phi:T\to S$ in~$\cS$.
\item[2)] $\propP$ contains all isomorphisms in fibers of $\cC/\cS$.
\item[3)] $\propP$ is stable under pushouts.
\item[4)] $\propP$ is stable under composition (of two morphisms in a fiber
of~$\cC$).
\item[6)] $\propP$ is stable under retracts (in the retract-free version of the theory this condition is omitted).
\end{itemize}
Similarly, a set $\propP\subset\cC$ of morphisms in an arbitrary category~$\cC$
with finite inductive limits 
is said to be {\em (globally) finitarily closed\/} if it fulfills conditions
2), 3), 4) and 6) (in the retract-free version the latter condition 
is omitted).

When we want to emphasize that we work in the retract-free version, we
speak about {\em finitarily semiclosed\/} sets.
\end{DefD}

\nxsubpoint\label{sp:stab.fin.dir.sum} (Stability under finite direct sums.)
Notice that conditions 2), 3) and 4) actually imply a weakened version of~7):
\begin{itemize}
\item[$7^w)$] $\propP$ is stable under finite direct sums.
\end{itemize}
Indeed, an empty direct sum is an isomorphism $\emptyset\to\emptyset$,
hence lies in $\propP$ by 2), and $u\oplus v=(u\oplus\id)\circ(\id\oplus v)$
is a composition of a pushout of~$u$ and a pushout of~$v$, hence belongs
to~$\propP$ by 3) and 4).

\nxsubpoint (Examples.)
(a) Notice that any closed class of morphisms is also finitarily closed,
e.g.\ the sets of cofibrations or acyclic cofibrations in pseudomodel
stack $\stks\stMOD\sO_\cX$ are finitarily closed.

(b) Another example is given by {\em strong cofibrations\/} and
{\em strong acyclic cofibrations\/}, i.e.\ morphisms with the local LLP
with respect to all acyclic fibrations (resp.\ all fibrations)
in a pseudomodel stack. More generally, any class of morphisms characterized
by local LLP with respect to another local class of morphisms is
finitarily closed.

\begin{DefD}
Let $I$ be any set of morphisms in a category~$\cC$ with finite inductive 
limits. Its {\em (global) finitary closure\/}, 
denoted by $\GFinCl_\cC I$, $\GFinCl I$ or just
$\FinCl I$, is the smallest finitarily closed subset of $\Ar\cC$ 
containing~$I$. Similarly, suppose $\cC/\cS$ be a stack, 
such that finite inductive limits exist in each fiber of $\cC/\cS$,
all $\phi^*:\cC(S)\to\cC(T)$ are right exact,
and that $\cS$ admits a final object~$e$. Then for any set $I\subset\Ar\cC(e)$
we define its {\em finitary closure\/} $\FinCl I$ or $\FinCl_\cC I$ as
the smallest finitarily closed subset $\propP\subset\Ar\cC$ containing~$I$.

When we want to emphasize that we work in the retract-free theory,
we speak about {\em finitary semiclosures\/}, and denote them by~$\FinSCl I$.
\end{DefD}

Any intersection of finitarily closed sets is again finitarily closed,
hence $\FinCl_\cC I$ exist for any~$I$: we just have to take the
intersection of all finitarily closed sets containing~$I$.

\begin{LemmaD}\label{l:glob.fin.cl}
Let $\cC$ be a category with finite inductive limits, $I\subset\Ar\cC$. Then:
\begin{itemize}
\item[(a)] The global finitary semiclosure $\GFinSCl I$ of~$I$ consists
of the morphisms $f\in\Ar\cC$ that are (isomorphic to) finite
(possibly empty) compositions of pushouts of morphisms from~$I$.
\item[(b)] The global finitary closure $\GFinCl I$ consists of retracts 
of morphisms from the global finitary semiclosure~$\GFinSCl I$. Moreover,
by~\ptref{l:repl.retr} it is enough to consider here 
only retracts with fixed source.
\end{itemize}
\end{LemmaD}
\begin{Proof} (a) Clearly, all such morphisms have to belong to $\GFinSCl I$;
conversely, the set of all such morphisms is easily seen to be stable under
pushouts and composition, and to contain all isomorphisms and~$I$, whence
the opposite inclusion.

(b) We may replace $I$ with $\GFinSCl I$ and assume $I$ to be already 
finitarily semiclosed. Denote by $\tilde I$ the set of all fixed source
retracts of morphisms from~$I$; since $I$ is closed under pushouts,
lemma~\ptref{l:repl.retr} implies that $\tilde I$ contains {\em all\/}
retracts of morphisms from~$I$. All we need to check is that 
$\tilde I\supset I$ is finitarily closed. It is obviously stable under 
pushouts and retracts, so all we need to check is that $\tilde I$ is 
closed under composition.

(c) Recall that $f'$ is a fixed source retract of $f$ iff there exist 
morphisms $i$ and $\sigma$, such that $if'=f$, $\sigma f=f'$, and 
$\sigma i=\id$. Now suppose that composable morphisms $f'$ and $g'$
are fixed source retracts of morphisms $f$ and $g$ from~$I$. First of all,
if $(j,\tau)$ are the retract morphisms for~$g'$, then the same morphisms
satisfy the retract property for $gf'$ and $g'f'$, i.e.\ $g'f'$ is a retract
of~$gf'$, so it would suffice to prove $gf'\in\tilde I$, i.e.\ we can
assume $g'=g\in I$. Now let $(i,\sigma)$ be the retract morphisms for~$f'$,
and denote by $\tilde g$ the pushout of $g$ with respect to~$i$;
since $\sigma i=\id$, the pushout of $\tilde g$ with respect to~$\sigma$ 
can be identified with $g$,
i.e.\ we get the following commutative diagram with {\em two\/} cocartesian
squares, one built on $i$ and $\tilde\imath$, and the other built on
$\sigma$ and $\tilde\sigma$:
\begin{equation}
\xymatrix@C+3pc{
&B\ar@{-->}[r]^{\tilde g}\ar@/_/[d]_{\sigma}&
\tilde C\ar@/_/@{-->}[d]_{\tilde\sigma}\\
A\ar[ru]^{f}\ar[r]^{f'}&B'\ar[r]^{g}\ar@/_/[u]_{i}&
C\ar@{-->}@/_/[u]_{\tilde\imath}}
\end{equation}
Now $(\tilde\imath,\tilde\sigma)$ satisfy the retract relations for
$gf'$ and $\tilde g f$, i.e.\ $gf'$ is a retract of $\tilde gf$, and the
latter morphism belongs to~$I$ since $\tilde g$ is a pushout of~$g$.
\end{Proof}

\begin{LemmaD}\label{l:loc.fin.cl}
Let $\cS$ be a site with final object~$e$, $\cC/\cS$ be a stack with 
right exact pullback functors, and $I\subset\Ar\cC(e)$ 
be any set of morphisms. Then a morphism $f:X\to Y$ in $\cC(S)$, 
$S\in\Ob\cS$, belongs to the finitary closure $\Cl_\cC I$ iff it locally
belongs to the global finitary closure of $\phi^*I$, i.e.\ iff
one can find a cover $\{S_\alpha\to S\}$, such that each $f|_{S_\alpha}$
belongs to the global finitary closure in $\cC(S_\alpha)$ 
of the pullback $I|_{S_\alpha}$ of~$I$ to $S_\alpha$.
\end{LemmaD}
\begin{Proof}
Clearly, any morphism from $\GFinCl_{\cC(S_\alpha)}I|_{S_\alpha}$
lies in $\FinCl_\cC I$, and the latter class of morphisms is local, hence
the second condition (``$f$ locally belongs to the global finitary closure'') 
implies 
the first. To show the opposite inclusion we have to check that the set
of morphisms $f$ having this property is finitarily closed. This fact is 
immediate, once we take into account that all $\phi^*$ are right exact:
the only non-trivial case is that of compositions, but if $f$ belongs
to the global finitary closure on one cover $\{S_\alpha\to S\}$, and $g$
on another $\{S'_\beta\to S\}$, then all we have to do is to choose a
common refinement of these two covers. (Notice that such an argument
wouldn't work for sequential compositions and infinite direct sums,
i.e.\ finitarity is essential here.)
\end{Proof}

\begin{DefD}\label{def:perf.cof.obj}
Let $\cX=(\cX,\sO)$ be a generalized ringed topos. We say that a morphism
$f:X\to Y$ in $\stks\stMOD\sO_\cX(S)=s\catMod{\sO|_S}$ is {\em perfectly 
cofibrant}, or a {\em perfect cofibration}, or just {\em perfect\/} 
if it belongs to the finitary
closure $\FinCl L_\sO(\cst I_\cX)$ of the set of standard cofibrant 
generators of $\stks\stMOD\sO_\cX$. We say that $f$ is a {\em globally
perfect cofibration\/} or just {\em globally perfect\/} 
if it belongs to the global finitary closure
$\GFinCl L_{\sO|_S}(\cst I_{\cX_{/S}})$ of the set of standard cofibrant 
generators in $s\catMod{\sO|_S}$.

An object $X\in\Ob\stks\stMOD\sO_\cX(S)=\Ob s\catMod{\sO|_S}$ is said to
be {\em perfect\/} or {\em perfectly cofibrant\/} 
(resp.\ {\em globally perfect}) if $\emptyset\to X$ 
is a perfect cofibration (resp.\ globally perfect cofibration).
\end{DefD}

\begin{DefD}\label{def:perf.obj.der.cat}
An object $\bar X$ of the homotopic category $\cD^{\leq0}(\cX,\sO)=\Ho
s\catMod\sO$ is {\em perfect\/} if it is isomorphic to an object of the
form $\gamma X$, for some perfect $X$ in $s\catMod\sO$. Similarly,
a morphism in $\cD^{\leq0}(\cX,\sO)$ is {\em perfect\/} if it is isomorphic
to $\gamma(f)$ for some perfect cofibration in~$s\catMod\sO$.

If $\sO$ admits a zero, we extend these definitions to the stable homotopic
category $\cD^-(\cX,\sO)$ in the natural way, e.g.\ an object $(\bar X,n)$ 
of~$\cD^-(\cX,\sO)$ is perfect iff for some $N\geq0$ the object
$\Sigma^N\bar X$ is perfect in $\cD^{\leq0}(\cX,\sO)$.
\end{DefD}
Notice that these homotopic category notions are not local, i.e.\ 
if $\bar X\in\Ob\cD^{\leq0}$ becomes perfect on some cover, this doesn't
necessarily imply perfectness of $\bar X$ itself. This is one of the
reasons we prefer to work with perfect objects and morphisms inside
$s\catMod\sO$ and not $\cD^{\leq0}(\cX,\sO)$ in almost all situations.

\nxsubpoint 
(Direct description of perfect objects and cofibrations.)
According to~\ptref{l:loc.fin.cl}, an object or a morphism in $s\catMod\sO$ 
is perfect iff locally it is globally perfect, i.e.\ iff it becomes 
globally perfect on some cover. According to~\ptref{l:glob.fin.cl},
a morphism is globally perfect iff it is a (fixed source) retract of a 
finite composition of pushouts of morphisms from $L_\sO(\cst I_\cX)$;
in the retract-free theory, of course, we don't even have to consider retracts.

This explicit description, together with the fact that strong cofibrations 
constitute a finitarily closed set of morphisms containing $L_\sO(\cst I)$,
immediately implies the properties of perfect objects and cofibrations 
summarized in the following proposition:

\begin{PropD}\label{prop:basic.prop.perf}
(a) All perfect morphisms are strong cofibrations (i.e.\ have the
local LLP with respect to all acyclic fibrations) and  
cofibrations (i.e.\ lie in $\Cl L_\sO(\cst I)$). All perfect objects are
cofibrant and strongly cofibrant.

(b) Any globally perfect cofibration has bounded dimension. Any perfect 
cofibration or perfect object has locally bounded dimension; if $\cX$
is quasicompact, any perfect cofibration or object has bounded dimension.
\end{PropD}

\nxsubpoint\label{sp:dim.dec.perf}
(Dimensional decomposition of perfect cofibrations.)
Reasoning in the same way as in~\ptref{prop:cof.bounded.dim}, 
with ``finitary devissage'' used
instead of ``devissage'', we obtain the following statement:
\begin{Propz}
(a) Functors $\sk_n$, $\sk_{\leq n}$, $\sk_{>n}$ and $\sk_{[a,b]}$
of~\ptref{p:dim.theory} transform perfect cofibration into perfect 
cofibrations.

(b) A perfect cofibration is concentrated in dimensions $S\subset\bbZ$ 
iff it belongs to $\FinCl L_\sO(\cst I_S)$. In particular, 
perfect cofibrations purely of dimension~$n$ coincide with 
$\FinCl L_\sO(\cst I_{=n})$.

(c) A perfect cofibration $f$ of bounded dimension $N<+\infty$ 
(condition automatically verified for a quasicompact~$\cX$) can be canonically
decomposed into a finite composition of perfect cofibrations $\sk_{=k}(f)$,
$0\leq k\leq N$, each of them purely of dimension~$k$.
\end{Propz}

\nxsubpoint\label{sp:const.perf}
(Constant perfect objects and cofibrations.)
Similarly to~\ptref{sp:const.cof}, 
we say that a morphism $f_0:X_0\to Y_0$ in $\catMod\sO$ is a 
{\em (constant) perfect cofibration\/}
if it becomes a perfect cofibration in $s\catMod\sO$ when considered as 
a morphism of constant simplicial objects. An object $X_0$ of $\catMod\sO$
is said to be {\em perfect\/} if it is perfect as a constant object of 
$s\catMod\sO$.

One checks, reasoning as in~\ptref{sp:const.cof}, that 
{\em the constant perfect cofibrations coincide with the finitary closure 
$\FinCl\{\emptyset_\sO\to L_\sO(1)\}$}. Applying \ptref{l:loc.fin.cl}
and~\ptref{l:glob.fin.cl}, we see that a constant morphism is a 
(constant) perfect cofibration iff it can be locally represented as a retract
of a standard embedding $X_0\to X_0\oplus L_\sO(n)$.
Therefore, constant perfect objects are just local retracts of free 
$\sO$-modules of finite rank, i.e.\ the {\em vector bundles}
(cf.~\ptref{sp:retract.ver}).

Clearly, all components $f_n:X_n\to Y_n$ of a perfect cofibration $f:X\to Y$
in $s\catMod\sO$ are constant perfect cofibrations (same reasoning as 
in~\ptref{sp:const.cof}). In particular, {\em all components $X_n$ of 
a perfect simplicial object~$X$ are vector bundles.}

\nxsubpoint\label{sp:perf.cof.ss}
(Perfect cofibrations of simplicial sets.)
Let's consider the situation $\cX=\catSets$, $\sO=\Fempty$, i.e.\ 
$\catMod\sO=\catSets$, $s\catMod\sO=s\catSets$.

(a) First of all, (constant) perfect cofibrations in $\catSets$ are just 
injective maps of sets $f_0:X_0\to Y_0$, such that $Y_0-f_0(X_0)$ is finite.
In particular, constant perfect sets are just the finite sets 
(``vector bundles over~$\Fempty$'').

(b) Therefore, a perfect cofibration $f:X\to Y$ in $s\catSets$ has the 
following properties: it is a cofibration (i.e.\ a componentwise injective 
map) of bounded dimension $N<+\infty$, and all components 
$f_n:X_n\to Y_n$ are constant perfect cofibrations, i.e.\ 
all sets $Y_n-f_n(X_n)$ are finite. Conversely, one easily checks that 
whenever $f:X\to Y$ is a cofibration in $s\catSets$ of bounded dimension 
$N<+\infty$, such that $Y_n-f_n(X_n)$ is finite (for all $n$ or just for
$n\leq N$), then one can represent $f:X\to Y$ as a finite composition 
of pushouts of standard generators~$i_n$. In order to do this just 
choose a simplex $\sigma\in Y_n-f_n(X_n)$ of minimal dimension~$n$; 
then all its faces already 
lie in $X$, so we can ``glue in'' this simplex to~$X$ along its boundary, 
i.e.\ decompose $f$ into
$X\to X'\to Y$, where $X\to X'$ is the pushout of 
$i_n:\dot\Delta(n)\to\Delta(n)$ with respect to the map $\dot\Delta(n)\to X$
defined by the faces of~$\sigma$; one can check then that $X'\to Y$ is
still a cofibration (i.e.\ componentwise injective map) satisfying the 
above conditions for the same $N$ but with smaller value of 
$\sum_{n\leq N}|Y_n-f_n(X_n)|$. The proof is concluded by induction in this
number.

(c) In particular, a simplicial set~$X$ is perfect iff it has bounded dimension
and all its components~$X_n$ are finite. This is equivalent to saying that
the set of non-degenerate simplices of~$X$ is finite, i.e.\ that 
$X$ is a finite simplicial set. Hence {\em perfect simplicial sets are
exactly the finite simplicial sets.}

(d) Now let $f:X\to Y$ be any cofibration of simplicial sets, and suppose
$Y$ to be perfect, i.e.\ finite. Then $f$ maps non-degenerate simplices of
$X$ into non-degenerate simplices of~$Y$, all components of~$f$ being 
injective, hence $X$ also has finitely many non-degenerate simplices,
i.e.\ it is finite, or perfect. Now it is immediate that $f$ has bounded 
dimension, and all $Y_n-f_n(X_n)\subset Y_n$ are finite, hence $f$ is 
a perfect cofibration. We have just shown that {\em any cofibration
of simplicial sets with perfect target is itself perfect, and its source is 
also perfect.}

(e) This is applicable in particular to cofibrations $i_n\boxe i_m$ 
with target $\Delta(n)\otimes\Delta(m)=\Delta(n)\times\Delta(m)$, i.e.\
{\em all $i_n\boxe i_m$ are perfect cofibrations in $s\catSets$.}

\begin{PropD}\label{sp:pullb.perf}
Let $f:\cX'=(\cX',\sO')\to\cX=(\cX,\sO)$ be any morphism of generalized 
commutatively ringed topoi. Then $f^*:s\catMod\sO\to s\catMod{\sO'}$ 
preserves perfect cofibrations and objects. Moreover, 
$\dim f^*(u)\leq \dim u$ for any (perfect) cofibration $u$, and
$\dim f^*X\leq\dim X$ for any (perfectly) cofibrant object~$X$.
\end{PropD}
\begin{Proof}
The statement about dimensions doesn't need perfectness, and has been 
already discussed in~\ptref{sp:dim.pb.cof}. The fact about perfect objects 
follows immediately from that about perfect cofibrations, and the latter 
is shown by an obvious finitary devissage.
\end{Proof}

\begin{PropD}\label{prop:tens.prod.perf}
Let $\cX=(\cX,\sO)$ be a generalized (commutatively) ringed topos. Then:
\begin{itemize}
\item[(a)] $i\boxe s$ is a perfect cofibration in $s\catMod\sO$ 
whenever $i$ and $s$ are. Furthermore, $\dim(i\boxe s)\leq\dim i+\dim s$.
\item[(b)] $X\otimes Y=X\otimes_\sO Y$ is perfect for any two perfect 
simplicial objects $X$ and~$Y$ of~$s\catMod\sO$. Furthermore,
$\dim(X\otimes Y)\leq\dim X+\dim Y$.
\end{itemize}
\end{PropD}
\begin{Proof}
(b) is an obvious consequence of (a) since $\nu_{X\otimes Y}=\nu_X\boxe\nu_Y$.
The second statement of (a) has been shown in~\ptref{prop:dim.tens.prod}.
The first statement is reduced by finitary version of devissage used 
in~\ptref{prop:comp.tens.sO-mod} to the case of $i_n\boxe i_m$ in $s\catSets$,
already dealt with in \ptref{sp:perf.cof.ss},(e), q.e.d.
\end{Proof}

\begin{CorD}\label{cor:compw.crit.perf}
(Componentwise criterion of perfectness.)
Let $f:X\to Y$ be a morphism in $s\catMod\sO$ of bounded dimension
$N<+\infty$. Then $f$ is a perfect cofibration iff all
$(\sk_{\geq n}f)_n$, $0\leq n\leq N$, are constant perfect cofibrations iff all
$(\sk_{=n}f)_n$, $0\leq n\leq N$, are constant perfect cofibrations.
\end{CorD}
\begin{Proof}
These conditions are necessary by  
\ptref{sp:dim.dec.perf} and~\ptref{sp:const.perf}.
The proof of the opposite implications goes exactly in the same way 
as in~\ptref{prop:compw.crit.cof},
taking into account that in our case the dimensional decomposition of~$f$
is finite: $f=\sk_{=N}(f)\circ\cdots\circ\sk_{=1}(f)\circ\sk_{=0}(f)$,
and that $u\boxe i_n=u\boxe L_\sO(i_n)$ is a perfect cofibration whenever
$u$ is a constant perfect cofibration by~\ptref{prop:tens.prod.perf}.
\end{Proof}

\begin{PropD}\label{prop:symm.pow.perf}
Let $(\cX,\sO)$ be as above. Denote by $S^n=S_\sO^n$ the symmetric power 
functors, and by $\rho_n$ the ``$n$-th symmetric box power'' functor
of~\ptref{sp:can.dec.symm.pow}. Then:
\begin{itemize}
\item[(a)] $\rho_n(u)$ is a perfect cofibration in $s\catMod\sO$ whenever
$u$ is one, for any~$n\geq0$. Furthermore, $\dim\rho_n(u)\leq n\dim(u)$.
\item[(b)] $S^n(X)$ is a perfect simplicial object whenever $X$ is one.
Furthermore, $\dim S^n(X)\leq n\dim X$.
\end{itemize}
\end{PropD}
\begin{Proof}
Statement (b) is a consequence of (a) since $\nu_{S^nX}=\rho_n(\nu_X)$.
Statement (a) is shown simultaneously for all values of~$n\geq1$ 
(case $n=0$ is trivial since $\rho_0(u)=\nu_{L_\sO(1)}$) by 
our usual finitary devissage. The only non-trivial case is that of
$\rho_n(vu)$, but according to~\ptref{sp:proof.main.rhon}(c),
this $\rho_n(vu)$ is a finite composition of pushouts of 
$\rho_k(v)\boxe\rho_{n-k}(u)$, $0\leq k\leq n$.
This accomplishes the devissage step once we take \ptref{prop:tens.prod.perf}
into account. The statement about dimensions is shown similarly,
considering only perfect cofibrations of relative dimension $\leq d$, i.e.\ 
elements of the finitary closure $\FinCl L_\sO(\underline I_{\leq d})$:
then from $\dim\rho_k(v)\leq k\dim v\leq kd$, and similarly for
$\rho_{n-k}(u)$, we deduce $\dim(\rho_k(v)\boxe\rho_{n-k}(u))\leq nd$,
hence the same is true for the composition $\rho_n(vu)$ of pushouts of 
these morphisms. Therefore, our finitary devissage reduces everything to the
case of $\rho_n(i_m)$ in $s\catSets$; the statement is immediate there
by~\ptref{sp:perf.cof.ss},(d) and~(e). Since the target of $\rho_n(i_m)$ 
is $S^n\Delta(m)$, itself a quotient of~$\Delta(m)^n$, hence $nm$-dimensional,
we see that the same holds for $\rho_n(i_m)$, i.e.\ $\dim\rho_n(i_m)\leq nm$.
This proves the dimension statement.
\end{Proof}

\nxsubpoint (Derived category case.)
Since derived functors $\Loplus$, $\Lotimes$ and $\dL S^n$ exist and 
can be computed by means of cofibrant replacements (cf.\
\ptref{th:ex.loc.der.tensprod} and~\ptref{th:ex.der.symmprod}), 
we can extend our previous results to $\cD^{\leq0}(\cX,\sO)=\Ho s\catMod\sO$:
whenever $\bar X$ and $\bar Y$ are perfect in $\cD^{\leq 0}$, the same
holds for $\bar X\Loplus\bar Y$, $\bar X\Lotimes\bar Y$ and $\dL S^n\bar X$.
Furthermore, left derived pullbacks $\dL f^*$ with respect to 
morphisms of generalized ringed topoi preserve perfectness as well
(cf.~\ptref{th:ex.lder.pb}).

One might try to extend the above results to morphisms in $\cD^{\leq0}$.
The idea is to observe that any morphism in $\cD^{\leq0}$ can be
represented or replaced by a cofibration between cofibrant objects, 
and to ``derive'' $\boxe$ and $\rho_n$ with the aid of such 
``doubly cofibrant replacements'', using a more sophisticated version
of~\ptref{th:ex.left.der.psmod}. We don't want to provide more details 
for the time being.

\nxsubpoint (Perfect cofibrations between perfect objects.)
Denote by $\cC_p$ the full substack of $\cC=\stks\stMOD\sO_\cX$, consisting 
of perfect objects. Since $L_\sO(\cst I)\subset\Ar\cC_p$, we might ask whether
the finitary closure of $L_\sO(\cst I)$ inside $\cC_p$ coincides with the 
set of perfect cofibrations between perfect objects, i.e.\ whether it is
sufficient to consider pushouts with perfect target while constructing 
a perfect cofibration between perfect objects. In fact, 
\ptref{l:glob.fin.cl} and~\ptref{l:loc.fin.cl} immediately imply that
the answer to this question is positive.

\nxsubpoint\label{sp:perf.cones.cyl} (Perfectness of cones and cylinders.)
Suppose $\sO$ to have a zero. Then for any map $f:X\to Y$ of simplicial
$\sO$-modules we can construct its {\em cylinder\/} $Cyl(f)$ and 
{\em cone\/} $C(f)$ by means of the
following cocartesian squares (cf.~\ptref{sp:mapping.cone.cyl}):
\begin{equation}\label{eq:diag.perf.cone.cyl}
\xymatrix{
X\otimes\{0\}\ar[r]^{f}\ar[d]&Y\otimes\{0\}\ar@{-->}[d]\\
X\otimes\Delta(1)\ar@{-->}[r]&Cyl(f)}
\quad
\xymatrix{
X\otimes\{1\}\ar[r]\ar[d]&Cyl(f)\ar@{-->}[d]\\
0\ar@{-->}[r]&C(f)}
\end{equation}
If $X$ is perfect, then $\id_X\otimes\Delta(\partial^0):
X\otimes\{0\}\to X\otimes\Delta(1)$ is a perfect cofibration, since
$\Delta(\partial^0)$ is a perfect cofibration of simplicial sets. Therefore,
$Y\cong Y\otimes\{0\}\to Cyl(f)$ is a perfect cofibration whenever 
$X$ is perfect. If $f:X\to Y$ is a perfect cofibration between 
perfect objects, then $X\otimes\Delta(1)\to Cyl(f)$ is a perfect cofibration, 
being a pushout of $f$, hence $X\otimes\{1\}\to X\otimes\Delta(1)\to Cyl(f)$
is also perfect, hence the same holds for its pushout $0\to C(f)$, i.e.\ 
{\em the cone of a perfect cofibration between cofibrant objects is 
perfect.}

\nxsubpoint\label{sp:perf.susp} (Perfectness of suspensions.)
Suppose $\sO$ still has a zero, and let $X$ be any perfect simplicial
$\sO$-module. Then $X\oplus X\to X\otimes\Delta(1)$ is a perfect cofibration,
being equal to $\id_X\otimes j=\nu_X\boxe j$, since $j:\Delta(0)\sqcup
\Delta(0)\to\Delta(1)$ is a perfect cofibration of simplicial sets.
Hence the cofiber $0\to\Sigma X$ 
of $\id_X\otimes j$, i.e.\ its pushout with respect to
$X\oplus X\to 0$, is a perfect cofibration, i.e.\ {\em the 
suspension $\Sigma X$ of any perfect object is perfect.}

In particular, this justifies the definition of perfectness in the stable
case, given in~\ptref{def:perf.obj.der.cat}.

\nxsubpoint\label{sp:perf.cof.addc} (Additive case.)
Now suppose that $\sO$ is additive, i.e.\ is a classical sheaf of 
commutative rings. Then:

(a) Constant perfect cofibrations are just injective maps
$f:\sF'\to\sF$ of $\sO$-modules with the cokernel $\sE:=\Coker f$ a 
vector bundle, i.e.\ a local retract of a free $\sO$-module (resp.\ 
a locally free $\sO$-module in the retract-free theory). 
Indeed, it is evident that
any constant perfect cofibration has this property; conversely, if this is 
the case, $\sF\to\sE=\Coker f$ locally splits, $\sE$ being locally projective,
hence $\sF'\to\sF$ is locally isomorphic to $\sF'\to\sF'\oplus\sE$, i.e.\ 
to a retract of a map $\sF'\to\sF'\oplus\sO(n)$.

(b) Simplicial $\sO$-modules of dimension $\leq n$ correspond via
Dold--Kan to chain complexes concentrated in (chain) degrees $[0,n]$. 
In order to see this we just recall that by definition of dimension and 
of $\sk_n$ the simplicial $\sO$-modules of dimension $\leq n$ coincide with 
the essential image of the left Kan extension $I_{n,!}$ of the functor
$I_n:\catDelta_{\leq n}\to\catDelta$ (cf.~\ptref{sp:skel.coskel}), 
and the associated simplicial object 
functor $K:\Ch_{\geq0}(\catMod\sO)\to s\catMod\sO$ is also defined by means of 
a certain left Kan extension $J_!$ for $J:\catDelta_+\to\catDelta$, 
cf.~\ptref{th:DK} and~\ptref{sp:simpl.obj.from.complex}. Then our 
statement follows almost immediately from transitivity of left Kan 
extensions and definitions, once we consider the following commutative square
of categories
\begin{equation}
\xymatrix{
\catDelta_{+,\leq n}\ar[r]^{I_{+,n}}\ar[d]^{J_n}&\catDelta_+\ar[d]^{J}\\
\catDelta_{\leq n}\ar[r]^{I_n}&\catDelta}
\end{equation}
and observe that $((I_{+,n})_!X)_m=\injlim_{[m]\rightarrowtail[p],p\leq n}X_p$
equals $0$ for $m>n$ and $X_m$ for $m\leq n$.

(c) Since $\sk_nX\to X$ can be described as the universal (final) 
object in the category of morphisms from simplicial objects of 
dimension~$\leq n$ into~$X$, we see that $\sk_nX$ corresponds under
Dold--Kan equivalence to the ``brute truncation'' $\sigma_{\geq-n}P$ 
of corresponding chain complex~$P$. This enables us to compute the 
counterparts of $\sk_{\leq n}$, $\sk_{>n}$, $\sk_{=n}$ for chain complexes,
since adjoint equivalences 
$K:\Ch_{\geq0}(\catMod\sO)\leftrightarrows s\catMod\sO:N$
have to preserve arbitrary inductive limits. For example, a chain map
$f':X\to Y$ corresponds to a morphism of simplicial objects of dimension
$\leq n$ iff $f'_k:X_k\to Y_k$ are isomorphisms for all $k>n$.

(d) In particular, if a map of simplicial $\sO$-modules~$f$ corresponds
to a chain map $f'=Nf:X\to Y$, then $\sk_{=n}(f)$ corresponds to 
$f'_{=n}:(\cdots\to X_{n+1}\to X_n\to Y_{n-1}\to\cdots\to Y_0)\to
(\cdots\to X_{n+1}\to Y_n\to Y_{n-1}\to\cdots\to Y_0)$. Since $\sk_{=n}(f)$
can be identified with $K(f'_{=n})$, formula \eqref{eq:KA.n} yields 
an identification of $(\sk_{=n}(f))_n$ with $\id_Q\oplus f_n:Q\oplus X_n\to 
Q\oplus Y_n$, where $Q$ is a certain direct sum of $Y_k$ with $k<n$.
In particular, description given in (a) implies that 
$(\sk_{=n}(f))_n=\id_Q\oplus f_n$ is a constant 
perfect cofibration iff $f_n:X_n\to Y_n$ is one.

(e) Combining (b) and (d) with the dimensionwise criterion of 
perfectness~\ptref{cor:compw.crit.perf}, we see that {\em perfect cofibrations
of bounded dimension correspond via Dold--Kan to chain maps
$f:X\to Y$ with all components $f_n:X_n\to Y_n$ as in (a), such that
all $f_n$ are isomorphisms for $n\gg0$.} Taking (a) into account we see that
{\em perfect cofibrations of bounded dimension correspond via Dold--Kan
to injective chain maps $f:\sF_\cdot\to\sF'_\cdot$ with the cokernel 
a finite complex of vector bundles, i.e.\ a perfect complex 
in classical sense.}

(f) Applying the above result to morphisms $0\to X$, we see that
{\em perfect simplicial objects of bounded dimension correspond via
Dold--Kan to bounded chain complexes of vector bundles, i.e.\ to
perfect complexes in the classical sense}. If $\cX$ is quasicompact,
any perfect simplicial object has bounded dimension, hence perfect simplicial
objects correspond exactly to perfect complexes in the classical sense.

\nxpointtoc{$K_0$ of perfect morphisms and objects}
Our nearest goal is to construct $K_0$ of perfect morphisms and/or objects
over a generalized ringed topos $\cX=(\cX,\sO)$, sometimes assumed to have
enough points. We do this in quite a natural way, suggested by our definition
of perfect objects, which may be thought of as a certain modification of 
Waldhausen's construction to the case when the basic category is not required 
to have a zero object (cf.~\ptref{sp:comp.Waldhausen} for a more detailed 
comparison). In fact, we define {\em three\/} variants of such $K_0$:
\begin{itemize}
\item $K_0$ of all perfect cofibrations, denoted by $K^0_{big}(\cX,\sO)$.
\item $K_0$ of perfect cofibrations between perfect objects 
$K^0_{perf}(\cX,\sO)$ or simply $K^0(\cX,\sO)$. It will be
also called ``$K_0$ of perfect simplicial $\sO$-modules'' or
``$K_0$ of perfect objects''.
\item $K_0$ of constant perfect objects, i.e.\ vector bundles over~$\cX$, 
denoted by $K^0_{vect}(\cX,\sO)$ or $\hat K^0(\cX,\sO)$.
\end{itemize}
Besides, we can construct the above notions in the retract and retract-free
settings, cf.~\ptref{sp:retract.ver}, 
and the six arising abelian groups (four of them  
actually are commutative pre-$\lambda$-rings) 
are related to each other by some canonical homomorphisms.

\begin{DefD}\label{def:K0.perf.cof}
The {\em $K_0$ of perfect cofibrations}, or {\em ``the big $K_0$''} of
$\cX=(\cX,\sO)$, denoted by $K^0_{big}(\cX,\sO)$, is the free abelian group
generated by all perfect cofibrations (or just by their isomorphism classes)
$[u]$ in $s\catMod\sO$, modulo following relations:
\begin{itemize}
\item[0)] If $u$ is an isomorphism, then $[u]=0$.
\item[1)] If $X\stackrel u\to Y\stackrel v\to Z$ are composable perfect 
cofibrations, then $[vu]=[v]+[u]$.
\item[2)] If $X'\stackrel{u'}\to Y'$ is a pushout of a perfect cofibration
$X\stackrel u\to Y$, then $[u']=[u]$.
\item[3)] If $X\stackrel u\to Y$ is isomorphic to $Z\stackrel v\to W$ in the
derived category $\cD^{\leq0}(\cX)=\Ho s\catMod\sO$, then $[v]=[u]$.
\end{itemize}
\end{DefD}

Notice that condition 0) is actually superfluous, being a special case of 3),
since 1) already implies $[\id_X]+[\id_X]=[\id_X]$, hence $[\id_X]=0$. 
On the other hand, 0) and 1) already imply that $[u]$ actually depends on the 
isomorphism class of $u$ in $s\catMod\sO$, something that is also immediate 
from 2). Finally, the last condition 3) implies that $[u]$ depends only
on the isomorphism class of $u$ in $\cD^{\leq0}(\cX,\sO)$, i.e.\ we might 
say that $K^0_{big}(\cX,\sO)$ is generated by isomorphism classes
of perfect morphisms in $\cD^{\leq0}(\cX,\sO)$ 
(cf.~\ptref{def:perf.obj.der.cat}), and 
use only morphisms in derived category in the above definition
(actually condition 2) is not so easy to state in the derived category). 
However, we prefer to work in $s\catMod\sO$ for technical reasons.

Notice that it is not clear whether $K^0_{big}(\cX,\sO)$ admits a small set 
of generators, i.e.\ whether $K^0_{big}(\cX,\sO)$ is a ($\univU$-)small 
abelian group. That's why we introduce a ``smaller'' version of $K^0$,
which can be shown to coincide with the previous one in ``good situations''.

\begin{DefD}
The {\em $K_0$ of perfect cofibrations between perfect objects},
or simply {\em $K_0$ of perfect objects}, denoted by $K^0_{perf}(\cX,\sO)$ 
or $K^0(\cX,\sO)$, is the free abelian group generated by isomorphism
classes $[u]$ of perfect cofibrations $u:X\to Y$ between perfect objects 
in $s\catMod\sO$, subject to the same relations 0)--3) as above, where 
of course we consider only perfect cofibrations between perfect objects
(e.g.\ $X'$ in 2) has to be perfect).
\end{DefD}

Clearly, $[u]\mapsto[u]_{big}$ defines a homomorphism
$K^0(\cX,\sO)\to K^0_{big}(\cX,\sO)$, which can be shown to be 
an isomorphism in all ``good'' situations (e.g.\ if $\cX$ is quasicompact
and $\sO$ additive, but not only in this case). On the other hand,
it is easy to deduce from \ptref{l:glob.fin.cl} and \ptref{l:loc.fin.cl}
that the set of isomorphism classes of perfect objects of $s\catMod\sO$ 
is small, hence $K^0(\cX,\sO)$ admits a small set of generators and 
in particular is itself small. Therefore, we usually consider this 
$K^0(\cX,\sO)$, for example to avoid all possible 
set-theoretical complications.

\nxsubpoint\label{sp:perf.obj.in.K0} (Image of a perfect object in $K^0$.)
Given any perfect object $X$ of $s\catMod\sO$, we denote by $[X]$ the
element $[\nu_X]$ of $K^0(\cX,\sO)$ or $K^0_{big}(\cX,\sO)$ (if we 
want to distinguish these two situations, we write $[X]$ and $[X]_{big}$,
respectively), where $\nu_X:\emptyset\to X$ is the only morphism from the
initial object as before.

Now if $X\stackrel u\to Y$ is any perfect cofibration between perfect objects,
then $\nu_Y=u\circ\nu_X$, hence by 1) we get $[Y]=[\nu_Y]=[u]+[\nu_X]=[u]+[X]$,
i.e.\ $[u]=[Y]-[X]$. In particular, {\em $K^0(\cX,\sO)$ is generated by 
elements $[X]$ corresponding to isomorphism classes of perfect objects.}
(Of course, this is not necessarily true for $K^0_{big}(\cX,\sO)$.)
That's why we call $K^0(\cX,\sO)$ the $K_0$ of perfect objects.

Since $[u]=[Y]-[X]$, we might define $K^0(\cX,\sO)$ as the free abelian group
generated by isomorphism classes $[X]$ (either in $s\catMod\sO$ or 
$\cD^{\leq0}=\Ho s\catMod\sO$) of perfect objects, modulo following relations:
\begin{itemize}
\item[$0')$] $[\emptyset]=0$, where $\emptyset$ denotes the initial object
of $s\catMod\sO$.
\item[$2')$] Whenever we have a cocartesian square \eqref{eq:perf.cocart.sq} 
of perfect objects
in $s\catMod\sO$ with perfect cofibrations for horizontal arrows, 
$[Y]-[X]=[Y']-[X']$.
\begin{equation}\label{eq:perf.cocart.sq}
\xymatrix{
X\ar[r]^{u}\ar[d]&Y\ar[d]\\X'\ar[r]^{u'}&Y'}
\end{equation}
\item[$3')$] If two perfect objects $X$ and $Y$ become isomorphic in 
$\cD^{\leq0}(\cX,\sO)$, then $[X]=[Y]$.
\end{itemize}
Notice that we don't need to write down relations corresponding to 1),
since $[vu]=[Z]-[X]=([Z]-[Y])+([Y]-[X])=[v]+[u]$ is automatic for 
any composable couple of perfect cofibrations $X\stackrel u\to Y\stackrel v\to
Z$ between perfect objects. On the other hand, condition $3')$ 
implies following weaker conditions:
\begin{itemize}
\item[$3^w)$] If $f:X\to Y$ is a weak equivalence between perfect objects,
then $[X]=[Y]$.
\item[$3^a)$] If $f:X\to Y$ is an acyclic cofibration between perfect objects,
then $[X]=[Y]$.
\end{itemize}
Notice that $3^w)$ doesn't imply $3')$, since under the condition of $3')$ 
we know only that $X$ and $Y$ can be connected by a path of weak equivalences,
with intermediate nodes not necessarily perfect.

\nxsubpoint\label{sp:dirsum.in.K0} (Image of direct sums in $K_0$.)
Recall that any finitarily closed set of morphisms, e.g.\ the set of 
perfect cofibrations, is closed under finite direct sums since
$u\oplus v=(u\oplus\id)\circ(\id\oplus v)$ (cf.~\ptref{sp:stab.fin.dir.sum}).
Relations 1) and 2) immediately imply
\begin{itemize}
\item[4)] $[u\oplus v]=[u]+[v]$
\end{itemize}
Similarly, relations $0')$ and $2')$, where we put $X:=\emptyset$,
$Y'=X'\oplus Y$, imply $[X\oplus Y']=[X]+[Y']$. We can also deduce 
this formula from 4) applied to $\nu_X$ and $\nu_{Y'}$.

\nxsubpoint\label{sp:K0.fincl.set} 
($K_0$ of any finitarily closed set of morphisms.)
More generally, let $\cC$ be any category and $\propP\subset\Ar\cC$ be any 
(globally) finitarily (semi)closed set of morphisms in~$\cC$. Notice that
in particular this includes existence in $\cC$ of all pushouts of morphisms
from~$\propP$. Then we can define $K_0(\propP)$ as the free abelian group
generated by (isomorphism classes of) morphisms from~$\propP$, modulo 
relations 1) and 2), hence also automatically 0). Moreover, this $K_0$ is 
functorial in the following sense: if $\propP'\subset\Ar\cC'$ is another
finitarily semiclosed set as above, and $F:\cC\to\cC'$ is any right exact 
functor such that $F(\propP)\subset\propP'$, then $F$ induces a homomorphism
$F_*:K_0(\propP)\to K_0(\propP')$.

However, this construction is insufficient to describe our $K^0_{perf}$ and
$K^0_{big}$ because of relations 3), which involve isomorphisms in another 
(namely, derived) category. In order to cover these cases we consider
{\em triples} $(\cC,\propP,\gamma)$, where $\propP\subset\Ar\cC$ is as above,
and $\gamma:\cC\to\tilde\cC$ is an arbitrary functor. Then we can define 
$K_0(\propP,\gamma)$ by taking the quotient of $K_0(\propP)$ modulo 
relations 3), i.e.\ $[u]=[u']$ whenever $\gamma(u)$ is isomorphic to 
$\gamma(u')$ in $\tilde\cC$. This construction is also functorial,
if we define a morphism $(F,\tilde F,\eta):
(\cC,\propP,\gamma)\to(\cC',\propP',\gamma')$ as a triple consisting of 
a right exact functor $F:\cC\to\cC'$, such that $F(\propP)\subset\propP'$,
its ``derived functor'' $\tilde F:\tilde\cC\to\tilde\cC'$, and a functorial
isomorphism $\eta:\tilde F\circ\gamma\simto\gamma'\circ F$ (usually we 
just write $\tilde F\circ\gamma=\gamma'\circ F$).

For example, taking $\cC:=s\catMod\sO$, $\propP:=\{$perfect cofibrations in
$\cC\}$, $\tilde\cC:=\Ho\cC=\cD^{\leq0}(\cX,\sO)$, we recover 
$K^0_{big}(\cX,\sO)$. Abelian group $K^0_{perf}(\cX,\sO)$ can be also 
constructed in this manner if we take the full subcategory 
$\cC_{perf}\subset\cC$ instead of~$\cC$, but keep same $\tilde\cC=
\cD^{\leq0}(\cX,\sO)$. Finally, group homomorphism $K^0_{perf}(\cX,\sO)\to
K^0_{big}(\cX,\sO)$ is a special case of functoriality just discussed.

\nxsubpoint\label{sp:K0.vect.bdl} 
($K_0$ of constant perfect objects, i.e.\ vector bundles.)
In particular, we can consider $K_0$ of (constant) perfect cofibrations
between constant perfect objects, i.e.\ vector bundles. The resulting 
abelian group will be denoted by $K^0_{vect}(\cX,\sO)$ or $\hat K^0(\cX,\sO)$.

Reasoning as in~\ptref{sp:perf.obj.in.K0}, we see that 
$\hat K^0(\cX,\sO)$ is generated by isomorphism classes $[X]$ of vector 
bundles over $\cX$ modulo relations $0')$ and $2')$, 
where of course we consider 
only {\em constant\/} cocartesian squares~\eqref{eq:perf.cocart.sq},
i.e.\ all four objects involved must be vector bundles, and $X\stackrel u\to Y$
must be a constant perfect cofibration between vector bundles, i.e.\ 
it must lie in the finitary closure of $\{\emptyset_\sO\to L_\sO(1)\}$
(cf.~\ptref{sp:const.perf}).

Notice the absence of relations $3')$. Such relations are actually 
unnecessary due to the following fact: {\em if two constant simplicial
objects become isomorphic in $\Ho s\catMod\sO$, then they are already 
isomorphic in $\catMod\sO$.} We are going to check this in a moment.

\nxsubpoint ($\bpi_0$ of a simplicial $\sO$-module.)
(a) Given any $X\in\Ob s\catMod\sO$, we define $\bpi_0X=
\bpi_0(X)\in\Ob\catMod\sO$
as the cokernel of $d^X_0,d^X_1:X_1\rightrightarrows X_0$ in $\catMod\sO$.
Clearly, $\bpi_0$ is a functor $s\catMod\sO\to\catMod\sO$, and
$\bpi_0X_0=X_0$ for any constant simplicial object~$X_0$.

(b) It is immediate that the set $\propP$ of morphisms $f$ in fibers 
of $\stks\stMOD\sO_\cX$, such that $\bpi_0(f)$ is an isomorphism, is closed
(cf.~\ptref{def:closed.class}), 
$\bpi_0X$ being defined with the aid of a finite inductive limit
of components of~$X$. It is equally obvious that $\bpi_0$ commutes with
generalized ringed topos pullbacks. Now the standard acyclic cofibrations 
$\Lambda_k(n)\to\Delta(n)$, $0\leq k\leq n>0$ from~$J$ clearly induce
isomorphisms on $\pi_0$, both $\Lambda_k(n)$ and $\Delta(n)$ being 
connected and non-empty, hence $L_\sO(\cst J)$ lies in $\propP$, hence
$\propP$ contains all acyclic cofibrations, i.e.\ {\em $\bpi_0(f)$
is an isomorphism for any acyclic cofibration~$f$.}

(c) Now let $f:X\to Y$ be an acyclic fibration in $s\catMod\sO$. 
According to \ptref{sp:f.af.smod.O}, this means that 
$X_n\to Y_n\times_{(\cosk_{n-1}Y)_n}(\cosk_{n-1}X)_n$ are strict epimorphisms
for all $n\geq0$; in particular, $f_0:X_0\to Y_0$ and
$X_1\to Y_1\times_{(Y_0\times Y_0)}(X_0\times X_0)$ have to be 
strict epimorphisms, being just special cases of the above morphisms for
$n=0$ and~$n=1$. Notice that $\bpi_0X=\Coker(X_1\rightrightarrows X_0)$
actually depends only on $X_R\subset X_0\times X_0$, the image of
$X_1\to X_0\times X_0$, and in our case 
$X_R$ is the preimage of $Y_R$ under strict 
epimorphism $X_0\times X_0\to Y_0\times Y_0$. This immediately implies that
$X_0\twoheadrightarrow Y_0\twoheadrightarrow Y_0/Y_R=\bpi_0Y$ is the cokernel
of $X_R\rightrightarrows X_0$ as well (indeed, $X_R$ contains $X_0\times_{Y_0}
X_0$, hence $X_0/X_R$ has to factorize through $Y_0$, $X_0\to Y_0$ being 
a strict epimorphism, and the kernel of $Y_0\to X_0/X_R$ obviously contains
the image of $X_R$ in $Y_0\times Y_0$, i.e.\ $Y_R$; the rest is trivial).
In other words, $\bpi_0X\cong\bpi_0Y$, i.e.\ {\em
$\bpi_0(f)$ is an isomorphism for any acyclic fibration~$f$.}

(d) Combining (b) and (c) together, we see that $\bpi_0$ transforms 
weak equivalences into isomorphisms, hence it can be derived in a trivial 
manner, yielding a functor $\bpi_0:\cD^{\leq0}(\cX,\sO)\to\catMod\sO$.

(e) An interesting consequence is that {\em $\bpi_0$ commutes with the 
forgetful functor $\Gamma_\sO:\catMod\sO\to\cX$}, i.e.\ $\bpi_0X=\Coker
(X_1\rightrightarrows X_0)$ in $\catMod\sO$ coincides with 
the same cokernel computed in topos~$\cX$. Indeed, one can find 
a weak equivalence (even an acyclic cofibration) $f:X\to X'$,
and {\em $\Gamma_\sO(f)$ is a weak equivalence iff $f$ is one} by 
\ptref{prop:f.af.weq.sO-mod}, \ptref{prop:f.af.w.pointwise} and 
\ptref{sp:cat.simpl.mod}, hence both $\bpi_0(f)$ and $\bpi_0\Gamma_\sO(f)$
are isomorphisms, and we are reduced to the case of a fibrant~$X$
(fibrant in $s\catMod\sO$ or $s\cX$ -- this is the same thing). 
Consider $\bpi_0\Gamma_\sO X:=\Coker(X_1\rightrightarrows X_0)$ in $\cX$;
all we have to check is that $X_0\to\bpi_0\Gamma_\sO X$ is compatible 
with the $\sO$-module structure on $X_0$, i.e.\ that it is compatible with the
structure maps
$\sO(n)\times X_0^n\to X_0$. This is immediate, once we take into 
account that $\bpi_0(X\times Y)\cong\bpi_0X\times\bpi_0Y$ for any
fibrant sheaves of sets $X$ and~$Y$.

(f) In any case, we have constructed a functor $\bpi_0:\cD^{\leq0}(\cX,\sO)
\to\catMod\sO$, having property $\bpi_0\gamma A_0\cong A_0$ for any 
constant simplicial $\sO$-module $A_0$, hence $\gamma A_0\cong\gamma B_0$
implies $A_0\cong B_0$ in $\catMod\sO$ for any two $\sO$-modules $A_0$ and
$B_0$, exactly as claimed in~\ptref{sp:K0.vect.bdl}.

\nxsubpoint (Morphisms between different versions of $K_0$.)
Of course, we get canonical homomorphisms of abelian groups
\begin{equation}
K^0_{vect}(\cX,\sO)=\hat K^0(\cX)\longrightarrow
K^0_{perf}(\cX,\sO)=K^0(\cX)\longrightarrow
K^0_{big}(\cX,\sO)
\end{equation}
which can be deduced from general functoriality of~\ptref{sp:K0.fincl.set}.
Furthermore, we obtain similar homomorphisms for the retract-free versions
of the above groups, as well as canonical homomorphisms from the 
retract-free into the retract versions, fitting together into a commutative
$3\times 2$-diagram.

\nxsubpoint\label{sp:K0.vect.addc} ($K^0_{vect}$ in the additive case.)
Now suppose $\sO$ to be additive, so that $\cX=(\cX,\sO)$ is just a 
classical ringed topos. In this case constant perfect cofibrations 
are just injective maps of $\sO$-modules $\sF'\to\sF$ with the cokernel 
a vector bundle (cf.~\ptref{sp:perf.cof.addc},(a)). 
Consider any short exact sequence of vector bundles on~$\cX$:
\begin{equation}
0\longto \sF'\stackrel u\longto\sF\stackrel v\longto\sE\longto0
\end{equation}
According to the description just recalled, $u$ is a perfect cofibration
between constant perfect objects, and we get a cocartesian square 
as in~\eqref{eq:perf.cocart.sq}:
\begin{equation}
\xymatrix{
\sF'\ar[r]^{u}\ar[d]&\sF\ar[d]^v\\
0\ar[r]^{u'}&\sE}
\end{equation}
Relations $0')$ and $2')$ yield $[\sF]-[\sF']=[\sE]-[0]=[\sE]$, i.e.\ 
we have following relations in $K^0_{vect}(\cX,\sO)$:
\begin{itemize}
\item[$5^+)$] $[\sF]=[\sF']+[\sE]$ for any short exact sequence 
of vector bundles\\$0\to\sF'\to\sF\to\sE\to 0$.
\end{itemize}
Conversely, relations $5^+)$ imply $0')$ and $2')$. As to $0')$,
putting $\sF'=\sF=\sE=0$ immediately yields $[0]=0$. In order to show
$2')$ consider any square \eqref{eq:perf.cocart.sq} and put $E:=\Coker u$.
Since $u'$ is a pushout of $u$, $\Coker u'\cong\Coker u=E$; on the other 
hand, $E$ has to be a vector bundle, $u$ being a perfect cofibration. 
Relations $5^+)$ for $0\to X\stackrel u\to Y\to E\to 0$ and 
$0\to X'\stackrel{u'}\to Y'\to E\to 0$ yield $[Y']-[X']=[E]=[Y]-[X]$ 
as required in $2')$.

We conclude that {\em in the additive case $K^0_{vect}(\cX,\sO)$ coincides 
with the free abelian group generated by isomorphism classes of vector 
bundles modulo relations $5^+)$, i.e.\ with the classical (Grothendieck's)
$K_0$ of vector bundles over $\cX$.}

\nxsubpoint\label{sp:K0.perf.addc} ($K^0_{perf}$ in the additive case.)
Now let $\sO$ still be additive, and consider $K^0=K^0_{perf}(\cX)$.
Let's suppose in addition all perfect objects over~$\cX$ to have bounded
dimension (quasicompactness of $\cX$ would suffice for this).
In this case perfect cofibrations between perfect objects correspond
via Dold--Kan equivalence to monomorphic chain maps $f:X\to Y$ of
perfect non-negatively graded chain complexes (i.e.\ finite chain complexes
consisting of vector bundle) with $E:=\Coker f$ also a perfect complex,
cf.~\ptref{sp:perf.cof.addc},(e). Reasoning as above we see that in this
case $K^0_{perf}(\cX)$ is the free abelian group generated by 
isomorphism classes $[E]$ of perfect non-negatively graded chain complexes,
modulo relations $5^+)$, where we consider all short exact sequences of 
perfect complexes, and~$3')$. This is {\em almost\/} the classical
$K_0$ of the triangulated category of perfect complexes 
(since all distinguished triangles are isomorphic to triangles defined 
by short exact sequences of complexes) up to some minor points.

\nxsubpoint (Contravariance of $K^0$.)
Let $f:(\cX',\sO')\to(\cX,\sO)$ be any morphism of generalized ringed topoi.
Recall that according to~\ptref{sp:pullb.perf}), 
the pullback functor $f^*:s\catMod\sO\to s\catMod{\sO'}$ 
preserves perfect cofibrations and perfect objects
hence also perfect cofibrations between perfect objects. Furthermore,
$f^*$ admits a left derived $\dL f^*:\cD^{\leq0}(\cX,\sO)\to
\cD^{\leq0}(\cX',\sO')$, which can be computed with the aid of cofibrant 
replacements (cf.~\ptref{th:ex.lder.pb}), i.e.\ 
$\dL f^*\gamma X\cong\gamma f^*X$ for any cofibrant and in particular for 
any perfect~$X$. This means that we are in position to apply general 
functoriality statements of~\ptref{sp:K0.fincl.set}, provided we consider only 
perfect cofibrations between perfect objects, thus obtaining the following
statement:
\begin{Propz}
For any morphism of generalized ringed topoi $f:\cX'\to\cX$ we have 
well-defined canonical homomorphisms $f^*:K^0_{perf}(\cX)\to K^0_{perf}(\cX')$
and $K^0_{vect}(\cX)\to K^0_{vect}(\cX')$, uniquely determined by rule
\begin{equation}
f^*[X]=[f^*X],\quad f^*[u]=[f^*(u)]
\end{equation}
Furthermore, these homomorphisms $f^*$ are functorial in $f$, i.e.\ 
$(f\circ g)^*=g^*\circ f^*$.
\end{Propz}

Notice that we don't obtain a similar statement for $K^0_{big}$; this is
probably due to the fact that a more ``correct'' definition of $K^0_{big}$
should involve only perfect cofibrations between cofibrant objects.

\begin{ThD} {\rm (Multiplication on $K^0$.)}
Let $\cX=(\cX,\sO)$ be a generalized commutatively ringed topos. 
Then the tensor product $\otimes=\otimes_\sO$ induces a commutative
multiplication law (i.e.\ symmetric biadditive map) 
on $K^0_{perf}(\cX)$ and $K^0_{vect}(\cX)$ by the rules
\begin{equation}\label{eq:mult.on.K0}
[X]\cdot[Y]=[X\otimes Y],\quad [u]\cdot[v]=[u\boxe v]
\end{equation}
where $X$ and $Y$ are perfect objects, and $u$ and $v$ are perfect cofibrations
between perfect objects (resp.\ constant perfect objects and perfect 
cofibrations between such for $K^0_{vect}(\cX)$), hence the same 
is true for $X\otimes Y$ and $u\boxe v$ by~\ptref{prop:tens.prod.perf}. 
Furthermore,
$K^0_{perf}(\cX)$ and $K^0_{vect}(\cX)$ become commutative rings under 
this multiplication, and canonical maps $K^0_{vect}(\cX)\to K^0_{perf}(\cX)$
are ring homomorphisms,
as well as the maps $f^*$ induced by generalized ringed topos pullbacks.
\end{ThD}
\begin{Proof} (a) All we have to check is that \eqref{eq:mult.on.K0}
really defines a product on $K^0_{perf}$ or $K^0_{vect}$; all remaining 
statements will follow immediately from associativity and commutativity of
tensor products and their compatibility with topos pullbacks/scalar
extension. We'll discuss only the case of $K^0_{perf}$; $K^0_{vect}$ is
treated similarly if we consider only constant perfect objects. Notice
that the formula $[X]\cdot[Y]=[X\otimes Y]$ is a special case of
$[u]\cdot[v]=[u\boxe v]$ for $u=\nu_X$, $v=\nu_Y$, so we'll concern ourselves 
only with the latter formula.

(b) We have to check that the map $(u,v)\mapsto[u\boxe v]$ is compatible
with the relations 0)--3) in each argument. By symmetry it suffices to 
check this compatibility only in first argument, while assuming $v$ to 
be fixed. Now relation 0) is trivial: if $u$ is an isomorphism, so is
$u\boxe v$, hence $[u\boxe v]=0$. Relation 2) is also simple since
$u'\boxe v$ is a pushout of $u\boxe v$ whenever $u'$ is a pushout of~$u$,
hence $[u'\boxe v]=[u\boxe v]$ as expected. Relation 1) follows from the fact
that $u'u\boxe v$ can be decomposed into $u'\boxe v$ and a pushout of 
$u\boxe v$ (cf.~\ptref{l:cocart.defect} or~\ptref{l:cl.box.div.class}), 
hence $[u'u\boxe v]$ is indeed equal to $[u'\boxe v]+[u\boxe v]$
by 1) and 2). So only relations 3) remain.

(c) Before checking compatibility with relations 3) let's show the following
formula: {\em if $X\stackrel u\to Y$ and $Z\stackrel v\to W$ are perfect 
cofibrations between perfect objects, then $[u\boxe v]=
[Y\otimes W]+[X\otimes Z]-[X\otimes W]-[Y\otimes Z]$.} Indeed, consider 
the diagram used to define $u\boxe v$:
\begin{equation}
\xymatrix{
X\otimes Z\ar[r]^{u\otimes\id_Z}\ar[d]&Y\otimes Z\ar@{-->}[d]\ar[ddr]\\
X\otimes W\ar@{-->}[r]^{u'}\ar[rrd]&T\ar@{-->}[rd]|{u\boxe v}\\
&&Y\otimes W}
\end{equation}
All objects involved in this diagram are perfect, and all morphisms are
perfect cofibrations by~\ptref{prop:tens.prod.perf}; 
furthermore, 2) or $2')$ is applicable to this 
cocartesian square, yielding $[u']=[u\otimes\id_Z]=[Y\otimes Z]-[X\otimes Z]$.
On the other hand, $[X\otimes W]+[u']+[u\boxe v]=[Y\otimes W]$ by 1); 
combining these formulas together, we obtain the announced formula for 
$[u\boxe v]$.

(d) Now suppose that perfect cofibrations between perfect objects 
$X\stackrel u\to Y$ and $X'\stackrel{u'}\to Y'$ become isomorphic in
$\cD^{\leq0}(\cX)$, i.e.\ $\gamma(u')\cong\gamma(u)$, and in particular
$\gamma X\cong\gamma X'$ and $\gamma Y\cong\gamma Y'$. Applying the formula 
of (c) to $[u\boxe v]$ and $[u'\boxe v]$, we see that it would suffice 
to show the following statement: {\em if $\gamma X\cong\gamma X'$ and
$\gamma Y\cong\gamma Y'$ in $\cD^{\leq0}(\cX)$, then $[X\otimes Y]=
[X'\otimes Y']$.}

(e) This statement follows from existence of derived tensor products
$\Lotimes$ and the fact that they can be computed with the aid of
cofibrant replacements (cf.~\ptref{th:ex.loc.der.tensprod}): indeed, we have
$\gamma(X\otimes Y)\cong\gamma X\Lotimes\gamma Y\cong\gamma X'\Lotimes
\gamma Y'\cong\gamma(X'\otimes Y')$, hence $[X\otimes Y]=[X'\otimes Y']$
in $K^0(\cX)$ by~3).

(f) The above reasoning shows that for any fixed $v$ the rule 
$[u]\mapsto[u\boxe v]$ determines a well-defined map 
$h_v:K^0(\cX)\to K^0(\cX)$;
interchanging arguments, we obtain compatibility with relations 0)--3) in
$v$ as well, i.e.\ show that $[v]\mapsto h_v$ is a well-defined homomorphism 
$K^0(\cX)\to\End(K^0(\cX))$, which corresponds by adjointness to a 
$\bbZ$-linear map $K^0(\cX)\otimes_\bbZ K^0(\cX)\to K^0(\cX)$, 
i.e.\ to a $\bbZ$-bilinear map $K^0(\cX)\times K^0(\cX)\to K^0(\cX)$
with required property $([u],[v])\mapsto[u\boxe v]$, q.e.d.
\end{Proof}

\nxsubpoint (Identity of $K^0(\cX)$.)
It is worth mentioning that the identity of $K^0(\cX)$ under the multiplication
just defined is $1:=[L_\sO(1)]=[\nu_{L_\sO(1)}]$. This is due to the fact that
$L_\sO(1)\otimes X\cong X$ (i.e.\ $L_\sO(1)$ is the unit for $\otimes$) and
$K^0(\cX)$ is generated by $[X]$, or can be equally easy deduced directly from
$\nu_{L_\sO(1)}\boxe u\cong u$ (i.e.\ $\nu_{L_\sO(1)}$ is a unit for $\boxe$).

\nxsubpoint (General remarks on pre-$\lambda$-rings.)
Given any commutative ring~$K$, we denote by $\hat G^t(K)$ or
by $1+K[[t]]^+$ the set of formal power series over~$K$ in one 
indeterminate~$t$ with constant term equal to one, 
considered as an abelian group under multiplication. We write the 
group law of $\hat G^t(K)$ additively; when we consider a formal series $f(t)$
of $1+K[[t]]^+$ as an element of this group, we denote it by $\{f(t)\}$;
thus $\{f(t)\}+\{g(t)\}=\{f(t)g(t)\}$. In this way we obtain a functor
$\hat G^t$ from commutative rings into abelian groups.

Of course, $t$ can be replaced here by
any other letter, e.g.\ $\hat G^u(K)=1+K[[u]]^+$. 

Now recall the following definition (cf.\ SGA~6, V.2.1):
\begin{DefD}\label{def:pre.lambda.ring}
A pre-$\lambda$-ring is a commutative ring~$K$, endowed with an abelian
group homomorphism $\lambda_t:K\to\hat G^t(K)=1+K[[t]]^+$, such that
$\lambda_t(x)=1+xt+\cdots$ for any $x\in K$.
\end{DefD}
Recall that the coefficient at $t^n$ of $\lambda_t(x)$ is usually denoted 
by~$\lambda^n(x)$ or $\lambda^nx$; thus a pre-$\lambda$-ring is a 
commutative ring~$K$ together with a family of unary operations 
$\lambda^n:K\to K$, $n\geq0$, called {\em exterior power operations,} 
such that $\lambda^0(x)=1$ and $\lambda^1(x)=x$ for any $x\in K$, and
\begin{equation}
\lambda^n(x+y)=\sum_{p+q=n}\lambda^p(x)\lambda^q(y)
\end{equation}

\nxsubpoint (Symmetric operations on a pre-$\lambda$-ring.)
Given any pre-$\lambda$-ring $K$, one defines {\em symmetric (power)
operations} $s^n$ on $K$ by means of the following generating series:
\begin{equation}\label{eq:symm.pow.ops}
s_t(x)=\sum_{n\geq0}s^n(x)t^n:=\lambda_{-t}(x)^{-1}
\end{equation}
Of course, one can write arising relations explicitly:
\begin{equation}
\sum_{p+q=n}(-1)^p\lambda^p(x)s^q(x)=0
\end{equation}

Since $f(t)\mapsto f(-t)^{-1}$ is an automorphism of $1+K[[t]]^+$,
$s_t:K\to 1+K[[t]]^+$ is an abelian group homomorphism as well,
and $s_t(x)=\lambda_{-t}(x)^{-1}=(1-xt+\cdots)^{-1}=1+xt+\cdots$,
i.e.\ we still have $s^0(x)=1$, $s^1(x)=x$, and
\begin{equation}
s^n(x+y)=\sum_{p+q=n}s^p(x)s^q(y)
\end{equation}
In other words, symmetric power operations define another 
pre-$\lambda$-structure on $K$.

Conversely, given any abelian group homomorphism $s_t:K\to 1+K[[t]]^+$,
such that $s_t(x)=1+xt+\cdots$, or equivalently, a collection of
symmetric power operations $\{s^n:K\to K\}_{n\geq0}$ satisfying the 
above relations, we can recover $\lambda_t$ from \eqref{eq:symm.pow.ops},
thus obtaining a pre-$\lambda$-ring structure on~$K$. In other words,
{\em pre-$\lambda$-rings admit an equivalent description in terms of
symmetric operations~$s^n$.}

\begin{ThD}
Let $\cX=(\cX,\sO)$ be as above. There are unique $\bbZ$-linear maps 
$s^n:K^0(\cX)\to K^0(\cX)$, $n\geq0$, where $K^0(\cX)$ denotes either
$K^0_{perf}(\cX)$ or $K^0_{vect}(\cX)$, such that
\begin{equation}\label{eq:symmpow.on.K0}
s^n[X]=[S^n(X)],\quad s^n[u]=[\rho_n(u)]
\end{equation}
where $S^n=S_\sO^n$ is the $n$-th symmetric power functor, and
$\rho_n(u)$ is defined in~\ptref{sp:can.dec.symm.pow}.
Furthermore, these maps $s^n$ are compatible with ring homomorphisms 
$f^*$ induced by generalized ringed topos morphisms as well as 
with ring homomorphism $K^0_{vect}(\cX)\to K^0_{perf}(\cX)$, 
and they satisfy following relations:
\begin{align}
s^0(\xi)=&1,\\
s^1(\xi)=&\xi,\\
s^n(\xi+\eta)=&\sum_{p+q=n}s^p(\xi)\cdot s^q(\eta)
\end{align}
In other words, $s^n$ are the symmetric power operations for a unique
pre-$\lambda$-ring structure on $K^0(\cX)$.
\end{ThD}
\begin{Proof}
We deal again only with the case $K^0(\cX)=K^0_{perf}(\cX)$; the case of
$K^0_{vect}$ is treated similarly, restricting all considerations to
constant perfect objects and cofibrations. Also notice that
$\rho_n(\nu_X)=\nu_{S^nX}$, hence formula $s^n[u]=[\rho_n(u)]$ implies
$s^n[X]=[S^n(X)]$, i.e.\ it is enough to consider only the first formula.

(a) First of all, notice that $[S^0(X)]=[L_\sO(1)]=1$ for any perfect~$X$,
and $[\rho_0(u)]=[\nu_{L_\sO(1)}]=1$ for any $u$ by definition, hence
$s^0$ is well-defined and $s^0(\xi)=1$ for all $\xi\in K^0(\cX)$. 
Similarly, $[S^1(X)]=[X]$ and $[\rho_1(u)]=[u]$ just because $\rho_1(u)=u$,
hence $s_1$ is also well-defined and $s^1(\xi)=\xi$ for all~$\xi$.

(b) For any perfect cofibration $u$ between perfect objects 
we put
\begin{equation}
s_t(u):=\sum_{n\geq0}[\rho_n(u)]\in 1+K^0(\cX)[[t]]^+
\end{equation}
In order for this expression to make sense we need to know that all
$\rho_n(u)$ are perfect cofibrations between perfect objects;
this is true by~\ptref{prop:symm.pow.perf}, \ptref{sp:can.dec.symm.pow}
and \ptref{cor:prev.steps.symmpr} (in fact, we obtain by induction in~$k$ that
all $\rho_k^{(n)}(u):F_{k-1}S^n(u)\to F_kS^n(u)$ are perfect cofibrations 
between perfect objects). Moreover, the free term $s^0(u)$ of $s_t(u)$ 
equals one by (a), hence $s_t(u)$ has indeed free term equal to one.
We write $s_t(X)$ instead of $s_t(\nu_X)$ for any perfect~$X$ as usual.

(c) Now all we have to check is that $[u]\to s_t(u)$ is a well-defined map
$K^0(\cX)\to 1+K^0(\cX)[[t]]^+$, i.e.\ that the $s_t(u)$ satisfy relations
0)--3) of \ptref{def:K0.perf.cof}. This is evident for 0); for 2) we just
use the fact that $\rho_n(u')$ is a pushout of $\rho_n(u)$ whenever
$u'$ is a pushout of $u$ (cf.~\ptref{sp:proof.main.rhon},(a)). As to 1),
we use that $\rho_n(vu)$ can be decomposed into a composition of pushouts
of $\rho_k(v)\boxe\rho_{n-k}(u)$, $0\leq k\leq n$ 
(cf.~\ptref{sp:proof.main.rhon},(c)), hence 
$[\rho_n(vu)]=\sum_{p+q=n}[\rho_p(v)][\rho_q(u)]$, i.e.\ $s_t(vu)=s_t(v)s_t(u)$
as required by 1).

(d) So only condition~3) remains. Since for any perfect cofibration between
perfect objects $u:X\to Y$ we have $s_t(u)=s_t(Y)\cdot s_t(X)^{-1}$, 
it would suffice to show that $\gamma X\cong\gamma X'$ in $\cD^{\leq0}(\cX)$ 
implies $s_t(X)=s_t(X')$. This is immediate from existence of derived 
symmetric powers and the fact that they can be computed by means of 
cofibrant replacements (cf.~\ptref{th:ex.der.symmprod}): we obtain
$\gamma(S^nX)\cong\dL S^n(\gamma X)\cong\dL S^n(\gamma X')\cong\gamma(S^nX')$,
hence $[S^nX]=[S^nX']$ by~3). This completes the proof of existence of
symmetric power operations $s^n$ with required property, and the 
construction of the pre-$\lambda$-structure on~$K^0(\cX)$.

(e) We still have to check that the symmetric power operations $s^n$ 
are compatible with maps $f^*$ and $K^0_{vect}(\cX)\to K^0_{perf}(\cX)$;
the latter statement is evident, and the former is an immediate consequence 
of commutativity of generalized ringed topos pullbacks with symmetric 
powers of $\sO$-modules, q.e.d.
\end{Proof}

\nxsubpoint ($K_0$ of an algebraic monad $\Lambda$.)
All of the above constructions are applicable to any algebraic monad $\Lambda$,
since it can be considered as an algebraic monad over the point topos 
$\catSets$. We'll usually write $K^0(\Lambda)$ and $K^0_{vect}(\Lambda)=
\hat K^0(\Lambda)$ instead of 
$K^0(\catSets,\Lambda)$ and $K^0_{vect}(\catSets,\Lambda)$. 
These are abelian groups for an arbitrary $\Lambda$, 
and commutative pre-$\lambda$-rings for a commutative~$\Lambda$. 
Furthermore, any algebraic monad homomorphism
$\rho:\Lambda\to\Lambda'$ defines a generalized ringed topos morphism
$(\catSets,\Lambda')\to(\catSets,\Lambda)$, whence canonical homomorphisms
(even pre-$\lambda$-rings homomorphisms for commutative $\Lambda$ and
$\Lambda'$) $\rho_*:K^0(\Lambda)\to K^0(\Lambda')$. Of course,
$\rho_*$ is induced by scalar extension, e.g.\ $\rho_*[P]=[P_{(\Lambda')}]$
for any projective $\Lambda$-module $P$ of finite type. When $\Lambda$ 
is additive, we recover classical $K^0$ of a ring according 
to~\ptref{sp:K0.vect.addc}.

\nxsubpoint (Computation of $\hat K^0(\Fempty)$.)
In particular, the pre-$\lambda$-rings $\hat K^0(\Fempty)=K^0_{vect}
(\catSets,\Fempty)$ and $K^0(\Fempty)$ are quite important since 
$(\catSets,\Fempty)$ is the 2-final object in the 2-category of generalized
ringed topoi, and therefore we obtain pre-$\lambda$-ring homomorphisms
$\hat K^0(\Fempty)\to\hat K^0(\cX)$ and $\hat K^0(\Fempty)\to K^0(\Fempty)
\to K^0(\cX)$ for any generalized ringed topos $\cX=(\cX,\sO)$.

We'll show in several steps that $\hat K^0(\Fempty)\cong K^0(\Fempty)\cong 
\bbZ$ considered as a pre-$\lambda$-ring with respect to its only 
$\lambda$-ring structure given by $\lambda^k(n)=\binom nk$,
$\lambda_t(n)=(1+t)^n$, $s_t(n)=(1-t)^{-n}$, $s^k(n)=\binom{n+k-1}k$
(cf.\ SGA~6 V~2.5).

Let's show the statement about $\hat K^0(\Fempty)=K^0_{vect}(\Fempty)$.
A perfect (constant) $\Fempty$-module is just a finite set~$X$, hence
$\hat K^0(\Fempty)$ is generated by $[X]$, for all finite sets~$X$.
Since $[X\oplus Y]=[X]+[Y]$ by~\ptref{sp:dirsum.in.K0}, 
and any finite set $X$ is 
a direct sum (i.e.\ disjoint union) of $|X|$ copies of $\st1$, we get
$[X]=|X|\cdot[\st1]=|X|\cdot1$ since $[\st1]$ is the identity of 
$\hat K^0(\Fempty)$. In particular, the canonical ring homomorphism
$\bbZ\to\hat K^0(\Fempty)$ is surjective. On the other hand, the
map $X\mapsto|X|$ obviously satisfies relations $0')$ and $2')$ since
a constant perfect cofibration is just an injective map $u:X\to Y$ with
finite $Y-u(X)$, and the number of elements in the complement of
the image of~$u$ is obviously preserved under any pushouts. 
This yields a map $\hat K^0(\Fempty)\to
\bbZ$, $[X]\mapsto|X|$ in the opposite direction, clearly inverse to the
previous one, hence $\hat K^0(\Fempty)\cong\bbZ$. This isomorphism is 
obviously compatible with multiplication; as to the $\lambda$-structures,
$S^nX$ is just the $n$-th symmetric power $X^{(n)}:=X^n/\gS_n$ 
of a set~$X$, hence $|S^nX|=\binom{|X|+n-1}{n}$ as claimed.

\nxsubpoint (Computation of $[i_n]$ in $K^0(\Fempty)$ and $K^0(\cX)$.)
Now let's compute $[i_n]$ in $K^0(\Fempty)$, where $i_n:\dot\Delta(n)\to
\Delta(n)$ denotes a standard cofibrant generator of $s\catSets$. 
Consider for this the standard acyclic cofibration $j:\Lambda_n(n)\to
\Delta(n)$. On one hand, it can be decomposed into $\Lambda_n(n)\stackrel{u}\to
\dot\Delta(n)\stackrel{i_n}\to\Delta(n)$, where $u$ is a pushout of
$i_{n-1}$, hence $[j]=[i_n]+[u]=[i_n]+[i_{n-1}]$. On the other hand,
$[j]=0$ by 3) and 0) since $\gamma(j)$ is an isomorphism in the derived 
category, $j$ being an acyclic cofibration. We conclude that 
$[i_n]+[i_{n-1}]=0$ for any $n>0$, and clearly $[i_0]=1$, hence
\begin{equation}\label{eq:class.i_n}
[i_n]=(-1)^n\quad\text{for any $n\geq0$.}
\end{equation}
Notice that the above formula holds in any $K^0(\cX)$ as well,
$K^0(\Fempty)\to K^0(\cX)$ being a ring homomorphism, where of course
$[i_n]$ is understood as $L_\sO(\cst{i_n})$.

Now we are going to deduce from \eqref{eq:class.i_n} some sort of 
Euler characteristic formula:
\begin{PropD}
Let $\cX=(\cX,\sO)$ be a generalized (commutatively) ringed topos,
$u:X\to Y$ be a perfect cofibration of bounded dimension 
$\dim u\leq N<+\infty$ between perfect objects. Then all $\sk_{=n}u$ are 
also perfect cofibrations between perfect objects, all $(\sk_{=n}u)_n$ 
are constant perfect cofibrations between vector bundles, and
\begin{equation}\label{eq:dim.dec.K0}
[u]=\sum_{n=0}^N[\sk_{=n}u]=\sum_{n=0}^N(-1)^n[(\sk_{=n}u)_n]
\end{equation}
In particular, $[u]$ lies in the image of $\hat K^0(\cX)\to K^0(\cX)$.
\end{PropD}
\begin{Proof}
Indeed, the dimensional decomposition~\eqref{eq:full.dim.decomp} of $u$
into $\sk_{=n}u$ consists of finitely many steps, all $\sk_{=n}u$ being
isomorphisms for $n>N$, and all $\sk_{=n}u:F_n(u)\to F_{n+1}(u)$ 
are perfect cofibrations by~\ptref{sp:dim.dec.perf}, hence by induction 
all $F_n(u)$ are also perfect, $F_0(u)=X$ being perfect by assumption,
hence $[\sk_{=n}u]$ is an element of $K^0(\cX)$, and the first equality 
of~\eqref{eq:dim.dec.K0} follows from relations~1). Furthermore,
all components of $\sk_{=n}u$, and in particular $w_n:=(\sk_{=n}u)_n$,
are constant perfect cofibrations between vector bundles 
by~\ptref{sp:const.perf}, hence $[w_n]$ makes sense both in
$\hat K^0(\cX)$ and $K^0(\cX)$. Finally, $\sk_{=n}u$ is a pushout of
$w_n\boxe i_n$ by~\ptref{l:pshout.dec.pure.mor}, hence $[\sk_{=n}u]=
[w_n][i_n]=(-1)^n[w_n]$ by~\eqref{eq:class.i_n}. Computing the sum over all
$n\leq N$ we obtain the second half of~\eqref{eq:dim.dec.K0}, q.e.d.
\end{Proof}

In particular, we obtain the following interesting statement:
\begin{CorD}\label{cor:surj.K0.to.K0}
If all perfect simplicial objects over $\cX=(\cX,\sO)$ have bounded dimension,
e.g.\ if $\cX$ is quasicompact, then the canonical pre-$\lambda$-ring
homomorphism $\hat K^0(\cX)\to K^0(\cX)$ is surjective, i.e.\ 
$K^0(\cX)$ is generated by elements $[P]$, where $P$ runs over isomorphism
classes of vector bundles over~$\cX$.
\end{CorD}

For example, $\hat K^0(\Lambda)\to K^0(\Lambda)$ is surjective for 
any generalized ring~$\Lambda$.

\nxsubpoint ($\hat K^0$ and $K^0$ of a classical field $k$.)
Now let $k$ be a classical (i.e.\ additive) field. Then
$\hat K^0(k)\to K^0(k)$ is surjective, and $\bbZ\to\hat K^0(k)$ is
an isomorphism, the inverse being given by $[V]\mapsto\dim_kV$ 
(cf.~\ptref{sp:K0.vect.addc}). We claim that $\hat K^0(k)\to K^0(k)$ is 
also an isomorphism, i.e.\ $K^0(k)=K^0_{perf}(K)$ is also isomorphic to
$\bbZ$ with its standard $\lambda$-structure. Indeed, \ptref{sp:K0.perf.addc}
immediately implies that $\chi:[X]\mapsto\sum_{n\geq0}(-1)^n\dim_k H^n(NX)$
is a well-defined homomorphism $K^0(k)\to\bbZ$, such that $\chi(1)=1$,
where $NX$ denotes the normalized chain complex associated to a simplical
$k$-vector space. This $\chi$ is clearly an inverse to the surjective map
$\bbZ=\hat K^0(k)\to K^0(k)$.

\nxsubpoint ($K^0(\Fempty)=\bbZ$.)
We already know that $\hat K^0(\Fempty)\cong\bbZ$ and that the homomorphism
$\hat K^0(\Fempty)\to K^0(\Fempty)$ is surjective. We claim that it is 
in fact an isomorphism, i.e.\ $K^0(\Fempty)\cong\bbZ$ with the standard
$\lambda$-structure. Indeed, since $\bbQ$ is an extension of $\Fempty$,
we get a map $K^0(\Fempty)\to K^0(\bbQ)$, $[X]\mapsto[L_\bbQ(X)]$,
and $K^0(\bbQ)$ is isomorphic to $\bbZ$ via the map~$\chi$, hence
the composite map $\chi\circ L_\bbQ$ is an additive map $K^0(\Fempty)\to\bbZ$
mapping $1$ into $1$, hence the surjection $\bbZ\to K^0(\Fempty)$ must 
be an isomorphism.

\nxsubpoint ($K^0(\bbF_{1^n})=\hat K^0(\bbF_{1^n})=\bbZ$.)
Now let us compute $\hat K^0(\bbF_{1^n})$ and $K^0(\bbF_{1^n})$,
where $\bbF_{1^n}$ is the generalized ring defined 
in~\ptref{sp:examp.pres.comalg}. Recall that $\bbF_{1^n}=\Fone[\zeta^{[1]}\,|
\,\zeta^n=\bu\}$ and $\Fone=\Fempty[0^{[0]}]$, and 
$\catMod{\bbF_{1^n}}$ consists of sets $X$ together with a marked point
$0=0_X\in X$ and a bijection $\zeta=\zeta_X:X\to X$ respecting~$0_X$, 
such that $\zeta_X^n=\id_X$, 
i.e.\ an action of the cyclic group $C_n:=\bbZ/n\bbZ$. 
We claim that {\em the canonical 
$\lambda$-homomorphisms $\bbZ=\hat K^0(\Fempty)\to\hat K^0(\bbF_{1^n})\to
K^0(\bbF_{1^n})$ are actually isomorphisms}, i.e.\ {\em
both $\hat K^0(\bbF_{1^n})$ and $K^0(\bbF_{1^n})$ are canonically isomorphic
to~$\bbZ$.}

(a) First of all, a free $\bbF_{1^n}$-module of rank~$r$ is just the
$nr+1$-element set $0\sqcup \str\times C_n$, 
where $C_n=\bbZ/n\bbZ$ is the cyclic group
of order $n$ with the obvious action of~$\zeta$. Conversely, if $X$ is 
a finite set with a marked point~$0$, such that the stabilizer in $C_n$ 
of any non-zero $x\in X$ is trivial, then the $C_n$-orbit decomposition of~$X$
is of the above form, hence $X$ is a free $\bbF_{1^n}$-module of finite rank.

(b) An immediate consequence is that {\em any projective $\bbF_{1^n}$-module
of finite type is free}. Indeed, any projective module is a retract,
i.e.\ both a quotient and a submodule of a free module, and (a) shows that
any submodule of a free $\bbF_{1^n}$-module is free. Furthermore,
all constant perfect cofibrations are easily seen to be of the form
$X\to X\oplus L_{\bbF_{1^n}}(r)$, i.e.\ $X\to X\sqcup \str\times C_n$.
This means that $\hat K^0(\bbF_{1^n})$ is generated by $[L_{\bbF_{1^n}}(1)]=1$,
i.e.\ $\bbZ=\hat K^0(\Fempty)\to\hat K^0(\bbF_{1^n})$ is surjective,
and $\hat K^0(\bbF_{1^n})\to K^0(\bbF_{1^n})$ is surjective 
by~\ptref{cor:surj.K0.to.K0}. 

(c) It remains to check that these two surjections have trivial kernel.
It would suffice for this to construct a $\bbZ$-linear map 
$\chi:K^0(\bbF_{1^n})\to\bbZ$, such that $\chi(1)=1$. We define such a map
as follows: put $A:=\bbZ\otimes_{\Fone}\bbF_{1^n}=\bbZ[\zeta]/(\zeta^n-1)$,
choose any maximal ideal $\gm\subset A$ (such an ideal exists since
$A$ is a free $\bbZ$-module of rank $n>0$, hence a non-trivial commutative
ring), and put $k:=A/\gm$. Then $k$ is a classical field, and we have
a canonical generalized ring homomorphism $\rho:\bbF_{1^n}\to A\to k$,
which induces a map $\rho_*:K^0(\bbF_{1^n})\to K^0(k)\cong\bbZ$ with 
the required property.

\nxsubpoint\label{sp:poly.func} (Polynomial functions.)
Let $A$ be any abelian group. We denote by $A^{\bbN_0}$ the set of all
$A$-valued functions $f:\bbN_0\to A$ defined on the set of 
non-negative integers. Of course, $A^{\bbN_0}$ has a canonical abelian
group structure given by pointwise addition: $(f+g)(n)=f(n)+g(n)$.

For any $f\in A^{\bbN_0}$ we denote by $\Delta f\in A^{\bbN_0}$ the 
$A$-valued function given by
\begin{equation}
(\Delta f)(n):=f(n+1)-f(n)\quad\text{for all $n\geq0$.}
\end{equation}
We say that $f$ is a {\em $A$-valued polynomial function\/} if
$(\Delta^nf)(0)=0$ for all $n>N$, where $N\geq0$ is some integer.
(The smallest~$N$ with this property is called the {\em degree\/} of~$f$.)
An easy induction shows that this condition is equivalent to $\Delta^nf=0$
for all $n>N$, or to $\Delta^{N+1}f=0$. Furthermore, another induction
in~$x\geq0$ yields
\begin{equation}\label{eq:newton.formula}
f(x)=\sum_{n\geq0}\binom xn\cdot(\Delta^n f)(0)
\end{equation}
for any $x\geq0$ and any $f\in A^{\bbN_0}$. Notice that this sum is
finite for any $x\geq0$ since $\binom xn=0$ for $n>x$.

Now if $f\in A^{\bbN_0}$ is an $A$-valued polynomial function, all terms
corresponding to $n>N$ in the above sum are automatically zero, hence
\begin{equation}\label{eq:gen.poly.func}
f(x)=\sum_{n=0}^N\binom xn\cdot c_n\quad\text{for some $c_n\in A$.}
\end{equation}
Conversely, if $f(x)$ is given by the above formula for some $c_n\in A$,
$0\leq n\leq N$, then $(\Delta f)(x)$ is easily seen to be equal to
$\sum_{n=0}^{N-1}\binom xn\cdot c_{n+1}$, and an obvious induction yields
$(\Delta^nf)(0)=c_n$, where of course $c_n:=0$ for $n>N$, 
thus proving that $f$ is an $A$-valued polynomial function, as well as
showing the uniqueness of coefficients~$c_n$.

In other words, we see that the abelian group $A[\binom xn]\subset A^{\bbN_0}$ 
of $A$-valued polynomials is canonically isomorphic to 
$A\otimes\bbZ[\binom xn]$, where $\bbZ[\binom xn]\subset\bbQ[x]$ is the set
of all integer-valued polynomials, a free $\bbZ$-module with basis $\binom xn$,
$n\geq0$.

Another immediate consequence of the above constructions is that any
polynomial function $f:\bbN_0\to A$ can be uniquely extended 
by \eqref{eq:gen.poly.func} to a polynomial function $f:\bbZ\to A$,
where the coefficients $c_n$ are necessarily equal to $(\Delta^nf)(0)$,
i.e.\ the extension of~$f$ is given by the Newton extrapolation 
formula~\eqref{eq:newton.formula}. In particular, we can consider
\begin{equation}\label{eq:eval.at.-1}
f(-1)=\sum_{n=0}^N(-1)^nc_n=\sum_{n\geq0}(-1)^n(\Delta^nf)(0)
\end{equation}

\nxsubpoint (Hilbert function of a perfect cofibration.)
Let $u:X\to Y$ be a perfect cofibration in $s\catMod\sO$. 
We denote by $P_u:\bbN_0\to\hat K^0(\cX,\sO)$ the
{\em Hilbert function\/} of~$u$, given by
\begin{equation}
P_u(n):=[u_n],\quad\text{for any $n\geq0$.}
\end{equation}
Here $u_n:X_n\to Y_n$ denotes the corresponding component of~$u$; 
according to \ptref{sp:const.perf}, all $u_n$ are indeed constant cofibrations
between vector bundles, hence $[u_n]$ indeed makes sense in $\hat K^0(\cX)$.
Furthermore, we have $P_{u\boxe v}=P_uP_v$ (i.e.\ $P_{u\boxe v}(n)=
P_u(n)\cdot P_v(n)$ for all $n\geq0$) since $(u\boxe v)_n$ obviously equals
$u_n\boxe v_n$. Similarly, $P_{\rho_n(u)}=s^nP_u$.

If $X$ is a perfect simplicial object, we define its 
{\em Hilbert function\/} $P_X$ by $P_X:=P_{\nu_X}$, i.e.\ 
$P_X(n)=[X_n]$. Then $X\mapsto P_X$ satisfies $0')$ and~$2')$,
$P_{X\otimes Y}=P_XP_Y$, and $P_{S^nX}=s^nP_X$, i.e.\ all natural relations
of $K^0(\cX)$ except $3')$ are fulfilled.

Notice that $u\mapsto P_u$ satisfies relations 0)--2), since for example 
if $u'$ is a pushout of~$u$, then all $u'_n$ are pushouts of $u_n$, hence
$P_{u'}(n)=[u'_n]=[u_n]=P_u(n)$ by 2) in $\hat K^0(\cX)$. However, relations
3) are not satisfied, since components of a weak equivalence needn't be
isomorphisms.

\begin{PropD}\label{prop:HP.perf.cof} 
{\rm (Hilbert polynomials and Euler characteristic.)}
Let $\cX=(\cX,\sO)$ be a generalized ringed topos as usual,
$u:X\to Y$ be a perfect cofibration of bounded dimension $\dim u\leq N<+\infty$
between perfect objects (e.g.\ if $\cX$ is quasicompact, any such~$u$ has
bounded dimension). Then $P_u\in\hat K^0(\cX)^{\bbN_0}$ is a polynomial 
function of degree $\leq N$; more precisely,
\begin{equation}\label{eq:coeff.of.HP}
P_u(n)=\sum_{k=0}^N c_k\binom nk,\quad\text{where $c_k:=[(\sk_{=k}u)_k]$}
\end{equation}
Furthermore, the image of $P_u(-1)=\sum_{k=0}^N(-1)^kc_k$ under the canonical
$\lambda$-homomorphism $\hat K^0(\cX)\to K^0(\cX)$ equals $[u]$:
\begin{equation}\label{eq:HP.at.-1}
[u]=P_u(-1)=\sum_{k=0}^N(-1)^k(\Delta^kP_u)(0)=
\sum_{k=0}^N(-1)^kc_k=\sum_{k=0}^N(-1)^k[(\sk_{=k}u)_k]
\end{equation}
Similar formulas are valid for a perfect object~$X$ of dimension $\leq N$,
e.g.\ $P_X$ is a polynomial function of degree $\leq N$, and $P_X(-1)=[X]$.
\end{PropD}
\begin{Proof}
Of course, the statements about $P_X$ follow from those about $P_u$
applied to $u=\nu_X$, so we'll treat only the case of $P_u$. Notice that
it would suffice to show \eqref{eq:coeff.of.HP}: the remaining statements
would then follow from \eqref{eq:newton.formula}, \eqref{eq:eval.at.-1}
and \eqref{eq:dim.dec.K0}. Since $u\mapsto P_u$ satisfies relations 0)--2) 
and the dimensional decomposition \eqref{eq:full.dim.decomp} of~$u$
is finite, it is enough to show \eqref{eq:coeff.of.HP} for $u$ purely of
dimension~$k$. In this case $u$ is a pushout of $w\boxe i_k$, where
$w:=u_k=(\sk_{=k}u)_k$ (cf.~\ptref{l:pshout.dec.pure.mor}), hence
$P_u=P_{w\boxe i_k}=P_w\cdot P_{i_k}=wP_{i_k}$, $w$ being constant, 
i.e.\ we are reduced to proving
\begin{equation}
P_{i_k}(n)=\binom nk \quad\text{in $\hat K^0(\Fempty)=\bbZ$.}
\end{equation}
This equality is shown by direct computation: by definition,
$P_{i_k}(n)=[i_{k,n}]=|\Delta(k)_n|-|\dot\Delta(k)_n|$, $\Delta(k)_n$ consists
of all non-decreasing maps $\phi:[n]\to [k]$, and $\dot\Delta(k)_n$ of 
all non-surjective such maps, hence $P_{i_k}(n)$ equals to the number
of surjective non-decreasing maps $\phi:[n]\twoheadrightarrow[k]$, easily seen
to be equal to $\binom nk$, such maps being in one-to-one correspondence
to $k$-element subsets $I=\{n_1,\ldots,n_k\}\subset\{1,2,\ldots,n\}$
by the following rule: $I=\{j\in[n]\,:\,j>0$, $\phi(j-1)<\phi(j)\}$.
This finishes the proof of~\ptref{prop:HP.perf.cof}.
\end{Proof}

\nxsubpoint (Relation to Euler characteristic.)
Recall that in the additive case $P_{KA}(-1)=\sum_{k\geq0}(-1)^k[A_k]$ 
for any perfect complex $A$, according to the computation already 
mentioned in~\ptref{sp:Euler.char}. Therefore, the above formula
$[u]=P_u(-1)$ should be thought of as a counterpart of the Euler 
characteristic formula.

\nxsubpoint (Suspensions, cones and cylinders.)
Suppose that $\sO$ admits a zero, so that we can construct suspensions,
cones and cylinders.

(a) Notice that $\Delta(0)\to\Delta(1)$ is an acyclic cofibration of finite
simplicial sets, hence $[\Delta(1)]=[\Delta(0)]=1$ in $K^0(\Fempty)$.
We can use this to conclude $[i_1]=-1$ for $i_1:\Delta(0)\sqcup\Delta(0)=
\dot\Delta(1)\to\Delta(1)$, or just apply~\eqref{eq:class.i_n}.

(b) Let $X$ be a perfect simplicial object, hence its suspension
$\Sigma X$ is also perfect by~\ptref{sp:perf.susp}. Recall that
$\Sigma X$ is the cofiber of $\id_X\otimes i_1=\nu_X\boxe i_1:X\oplus X\to
X\otimes\Delta(1)$, i.e.\ $0\to\Sigma X$ is a pushout of this morphism,
hence $[\Sigma X]=[\nu_X]\cdot[i_1]=-[X]$. If $u:X\to Y$ is a perfect 
cofibration between perfect objects, the same holds for $\Sigma(u)$,
this morphism being a pushout of $u\boxe i_1$, and $[\Sigma(u)]=[u]\cdot[i_1]=
-[u]$.

(c) Now let $f:X\to Y$ be a perfect cofibration between perfect objects.
Applying 2) to cocommutative squares of~\eqref{eq:diag.perf.cone.cyl},
we obtain $[Y\otimes\{0\}\to Cyl(f)]=[f]\cdot[\Delta(0)\to\Delta(1)]=
[f]\cdot 0=0$, hence $[Cyl(f)]=[Y]$. The second square yields
$[C(f)]=[Cyl(f)]-[X]=[Y]-[X]$, exactly what one usually has in the classical
(additive) situation.

\nxsubpoint\label{sp:comp.Waldhausen} 
(Relation to Waldhausen's construction.)
Our definition of $K_0$ might be thought of as a modification of 
Waldhausen's construction of $K_0$ presented in \cite{Waldhausen}. 
Since Waldhausen defines higher algebraic $K$-theory as well, this 
observation might be used to define higher $K$-functors in our situation
by modifying Waldhausen's construction in a suitable way.

Recall that Waldhausen considers {\em categories~$\cC$ with cofibrations 
and weak equivalences}, sometimes called {\em Waldhausen categories}.
A Waldhausen category is a pointed category $\cC$ 
(i.e.\ it has a zero object $0$) with two classes of morphisms
$\propP$ and $\propW\subset\Ar\cC$, called {\em cofibrations\/} and
{\em weak equivalences}, satisfying certain conditions. Namely, both classes
are supposed to contain all isomorphisms and be closed under composition.
Secondly, the following conditions have to be fulfilled:
\begin{itemize}
\item[Co2)] Morphisms $0\to X$, for all $X\in\Ob\cC$, are cofibrations.
\item[Co3)] Cofibrations are stable under pushouts (in particular,
pushouts of cofibrations are required to exist).
\item[We2)] (Gluing lemma) If the left horizontal arrows $u$ and $v$
in the diagram
\begin{equation}
\xymatrix{
B\ar[d]^{\sim}&A\ar[l]_{u}\ar[d]^{\sim}\ar[r]&C\ar[d]^{\sim}\\
B'&A'\ar[r]\ar[l]_{v}&C'}
\end{equation}
are cofibrations, and the vertical arrows are weak equivalences, then
the induced map of pushouts $B\vee_AC\to B'\vee_{A'}C'$ is also
a weak equivalence.
\end{itemize}
For any cofibration $X\stackrel u\to Y$ we define its cofiber $E$ as 
$Y\vee_X0$. Clearly, if $u'$ is any pushout of $u$, then the cofiber of $u'$
is also equal to~$E$. Now Waldhausen's $K_0(\cC)=K_0(\cC,\propP,\propW)$
can be defined as the free abelian group generated by objects $[X]$ of~$\cC$
modulo relations $(W_1)$ $[Y]=[X]+[E]$ for any cofibration $X\stackrel u\to Y$ with
cofiber~$E$, and $(W_2)$ $[X]=[X']$ for any weak equivalence $X\to X'$.

Now we can see a certain similarity to our construction, where we consider
the category $\cC$ of all perfect objects,
$\propP$ is the set of perfect cofibrations, and $\propW$ is the
set of weak equivalences in~$\cC$,
at least in the case when $\sO$ admits a zero. 
Reasoning as in~\ptref{sp:K0.vect.addc} (with cofibers instead of cokernels), 
we can show that
relations $(W_1)$ are in fact equivalent to relations $0')$ and $2')$,
and $(W_2)$ correspond to our $3^w)$, but not $3')$ (i.e.\ we impose
more relations than Waldhausen).

It is an interesting question whether the axiom We2), which
makes sense for $\sO$ without a zero as well, holds for our choice
of $\cC$, $\propP$ and $\propW$. The answer is extremely likely to be 
positive (for example, it is positive for the classical case of 
pointed simplicial sets, i.e.\ simplicial $\bbF_1$-modules).
If the answer {\em is\/} positive, we might try to apply Waldhausen's
construction directly to our situation, thus immediately defining 
higher algebraic $K$-groups as well.

\nxpointtoc{Projective modules over $\Zinfty$}
Our next aim is to compute $\hat K^0(\CompZ)$. However, we need to know more
about projective bundles $\sE$ over $\CompZ$ to do this. Since 
$\sE|_{\Spec\bbZ}$ is given by a finitely generated projective (i.e.\ free)
$\bbZ$-module, we just have to study finitely generated projective modules
over $\Zninfty$ or $\Zinfty$. In fact, we are going to show that such 
projective modules are free as well, and that any constant perfect cofibration
between vector bundles over $\Zinfty$ is of the form $\Zinfty^{(n)}\to
\Zinfty^{(n)}\oplus\Zinfty^{(m)}\cong\Zinfty^{(n+m)}$, and similarly 
for~$\Zninfty$ and $\barZinfty$. This would imply 
$\hat K^0(\Zinfty)= K^0(\Zinfty)=\bbZ$, and similarly for $\Zninfty$ 
and~$\barZinfty$.

\nxsubpoint (Strictly convex archimedian valuation rings.)
We are going to treat our three cases $\Zinfty$, $\Zninfty$ and 
$\barZinfty$ simultaneously. In order to do this we fix following data:
\begin{itemize}
\item A (classical) field~$K$ together with a non-trivial archimedian valuation
$\fnorm$. Notice that necessarily $\charact K=0$, i.e.\ $K\supset\bbQ$.
We consider $K$ as a topological field with respect to the topology given
by its norm~$\fnorm$.
\item $V\subset K$ denotes the (generalized) valuation ring of $\fnorm$ 
in~$K$, defined as in~\ptref{sp:val.rings}. In other words, 
$V(n)\subset K^{(n)}=K^n$ consists of $\lambda=(\lambda_1,\ldots,\lambda_n)\in
K^n$ with $\sum_i|\lambda_i|\leq 1$.
\item $\tilde\bbQ$ denotes the closure of $\bbQ\subset K$ in $K$, and
$\tilde\bbQ_+\subset\tilde\bbQ$ denotes the closure of $\bbQ_+$. Since 
any element of $\tilde\bbQ$ is a limit of a Cauchy sequence in~$\bbQ$, 
we obtain a natural embedding $\tilde\bbQ\subset\bbR$, i.e.\ we can treat
elements of $\tilde\bbQ\subset K$ as real numbers.
\item We require $(K,\fnorm)$ to be {\em strictly convex\/}, meaning the
following:
\begin{equation}
|x+y|=|x|+|y|\quad\text{for $x$, $y\in K$ }\Rightarrow
\text{ $y=0$ or $xy^{-1}\in\tilde\bbQ_+$.}
\end{equation}
\item Finally, we require $|x|\in\tilde\bbQ_+$ for any $x\in K$.
\end{itemize}

The above conditions are easily seen to be satisfied by $K=\bbQ$, $\bbR$ and
$\bbC$ with their standard archimedian norms, given by absolute value.
In these cases we have $V=\Zninfty$, $\Zinfty$ and $\barZinfty$, 
respectively.

\begin{ThD}\label{th:proj.over.archv} 
{\rm (Projective modules over strictly convex archimedian valuation rings.)}
Let $(K,V)$ be as above, i.e.\ $V$ be a strictly convex archimedian valuation
ring, e.g.\ $V=\Zinfty$, $\Zninfty$ or $\barZinfty$. Denote by 
$V^{(n)}$ or $L_V(n)$ the free $V$-module of rank~$n$. Then:
\begin{itemize}
\item[(a)] Any projective $V$-module $P$ of finite type is free of finite
rank, i.e.\ isomorphic to some $L_V(n)$, $n\geq0$.
\item[(b)] If $L_V(n)$ is isomorphic to $L_V(m)$, then $n=m$.
\item[(c)] Any constant perfect cofibration of $V$-modules is of the form
$M\to M\oplus L_V(n)$. In particular, any constant perfect cofibration
between vector bundles over~$V$ is of the form $L_V(m)\to L_V(m)\oplus L_V(n)
\cong L_V(m+n)$.
\end{itemize}
\end{ThD}
\begin{Proof}
Notice that (b) follows immediately from existence of embedding
$\rho:V\to K$: applying scalar extension functor $\rho^*=K\otimes_V-$ to
$L_V(n)$ we obtain $L_K(n)=K^n$, a $K$-vector space of dimension~$n$,
and $K^n\simeq K^m$ is possible only for $n=m$.

As to the proof of (a) and (c), it will occupy the rest of this subsection.
We are going to obtain matrix descriptions of objects involved in (a) and (c),
reduce (c) to~(a) and prove~(a). Therefore, we fix $K$ and $V$ as above
until the end of this subsection.
\end{Proof}

\nxsubpoint (Finitely generated projective modules over any~$\Sigma$.)
Let $P$ be a finitely generated projective module over an arbitrary 
algebraic monad~$\Sigma$. Fix any finite system $f_1$, \dots, $f_n$
of generators of~$P$, or equivalently, a strict epimorphism 
$\phi:L_\Sigma(n)\twoheadrightarrow P$, 
such that $\phi(e_i)=f_i$, where $e_i=\{i\}$ denotes
the $i$-th basis element of $L_\Sigma(n)$. Since $P$ is projective,
$\phi$ admits a section $\sigma:P\to L_\Sigma(n)$, and $p:=\sigma\phi\in
\End_\Sigma(L_\Sigma(n))=M(n;\Sigma)=\Sigma(n)^n$ is a projector
($p^2=p$) with $P=\Coker(p,\id)=\Ker(p,\id)=p(L_\Sigma(n))$. Conversely,
if $p\in\End_\Sigma(L_\Sigma(n))$ is a projector, then the corresponding 
direct factor $P$ of $L_\Sigma(n)$ is projective and finitely generated,
being a strict quotient of $L_\Sigma(n)$. Therefore, any finitely generated
projective $\Sigma$-module $P$ with $n$ generators can be described by
means of a projector $p\in M(n,\Sigma)=\Sigma(n)^n$, and any such projector
determines such a projective $\Sigma$-module with $n$ fixed generators.

\nxsubpoint (Constant perfect cofibrations over any~$\Sigma$.)
Similarly, let $f:M\to N$ be a constant perfect cofibration 
of $\Sigma$-modules. According to~\ptref{sp:const.perf},
such an $f$ is a fixed-source retract of a morphism of the form 
$i_M:M\to M\oplus L_\Sigma(n)$, i.e.\ there are $\Sigma$-homomorphisms
$\sigma:M\oplus L_\Sigma(n)\twoheadrightarrow N$ and $j:N\to
M\oplus L_\Sigma(n)$, such that $\sigma j=\id_N$, $jf=i_M$ and $\sigma i_M=j$.
Putting $p:=j\sigma$, we obtain a projector 
$p=p^2\in\End_\Sigma(M\oplus L_\Sigma(n))$, such that $p i_M=i_M$; 
conversely, any such projector $p$ defines a direct factor $N$ of
$M\oplus L_\Sigma(n)$, such that $M\to N$ is a retract of $i_M$, i.e.\ 
a constant perfect cofibration with source~$M$. In this way we obtain a
description of constant perfect cofibrations in terms of projectors as well.

Notice that $\sigma:M\oplus L_\Sigma(n)\twoheadrightarrow N$ is necessarily
a strict epimorphism, i.e.\ $N$ is generated by $M$ and elements
$u_i:=\sigma(e_i)$, $1\leq i\leq n$. Conversely, 
let $f:M\to N$ be a constant perfect cofibration, and 
$u_1$, \dots, $u_m\in N$ be any system of elements which generates $N$ together
with $f(M)$, thus defining a strict epimorphism $\sigma':M\oplus L_\Sigma(m)
\twoheadrightarrow N$.
Recall that any constant perfect cofibration is strongly cofibrant and in
particular has the LLP with respect to strict epimorphisms (i.e.\ 
surjective maps) of $\Sigma$-modules; applying this to $f$ and $\sigma'$,
we see that $f$ is a fixed-source retract of $M\to M\oplus L_\Sigma(m)$.

\nxsubpoint\label{sp:proj.over.v} (Projective modules over $V$.)
Now let $P$ be a finitely generated projective module over~$V$.
Choose a system of generators $u_1$, \dots, $u_n$ of~$P$ with {\em minimal\/}
$n\geq0$. According to the above considerations, this yields a description
of $P$ as a direct factor (retract) of $L_V(n)$, given by a certain projector
$a=a^2\in\End_V(L_V(n))=M(n;V)$. We will ultimately show that the
minimality of $n$ implies $a=\id$, i.e.\ $P\cong L_V(n)$, thus proving
\ptref{th:proj.over.archv},(a).

Since $V\subset K$, $M(n;V)\subset M(n;K)=K^{n\times n}$, i.e.\ 
$a$ is given by a certain matrix $(a_{ij})_{1\leq i,j\leq n}$ 
with coefficients in~$K$, such that $a(e_j)=\sum_i a_{ij}e_i$.
Denote by $\vnorm$ the $L_1$-norm in $K^n$, i.e.\ 
$\|\lambda_1e_1+\cdots+\lambda_ne_n\|=|\lambda_1|+\cdots+|\lambda_n|$.
Then $V(n)$ can be identified with $\{\lambda\in K^n:\|\lambda\|\leq 1\}$,
and a matrix $a\in M(n;K)$ lies in $M(n;V)=V(n)^n$ iff all $a(e_j)$ lie
in $V(n)$, i.e. iff
\begin{equation}
\sum_{i=1}^n|a_{ij}|\leq1,\quad\text{for all $1\leq j\leq n$.}
\end{equation}

We fix $P$ and a matrix $a=a^2$ as above, and put $u_j:=a(e_j)$; these elements
generate a $V$-submodule of $L_V(n)$ isomorphic to~$P$, which will be
usually identified with~$P$. Notice that elements $u_j$ are identified 
then with our original system of generators of~$P$.

\nxsubpoint\label{sp:const.perf.cof.over.v} 
(Constant perfect cofibrations over $V$.)
Similarly, let $f:M\to N$ be a constant perfect cofibration over~$V$.
Choose a system of generators $u'_1$, \dots, $u'_n\in N$ of $N$ ``over $M$''
(i.e.\ $u'_i$ together with $f(M)$ generate $N$) with a {\em minimal\/} 
$n\geq0$. Then $f$ is a retract of $i_M:M\to M\oplus L_V(n)$, given by a
certain projector $p=p^2\in\End_V(M\oplus L_V(n))$, such that $p i_M=i_M$. 
Denote by $a$ the composite map $L_V(n)\to M\oplus L_V(n)\stackrel p\to
M\oplus L_V(n)\to 0\oplus L_V(n)=L_V(n)$, or equivalently, the map induced
by~$p$ on strict quotient $(M\oplus L_V(n))/M\cong L_V(n)$; 
clearly, $a=a^2$ is a projector
in $\End_V(L_V(n))$ defining projective module~$P$, the cofiber of~$f$.

We will ultimately show that $p=\id_M\oplus a$, hence $f$ will be 
identified with $M\to M\oplus P$; this will enable us to deduce
\ptref{th:proj.over.archv},(c) from \ptref{th:proj.over.archv},(a).

Let us denote $u'_j:=p(e_j)\in M\oplus L_V(n)$; these elements can be still
identified with our original generators of $N$, considered here as a 
submodule of $M\oplus L_V(n)$. The projections $u_j:=a(e_j)$ of
$u'_j$ to $L_V(n)$ will then generate~$P$.

\nxsubpoint (Elements of direct sums over~$V$.)
Let $M$ and $M'$ be two arbitrary $V$-modules. Since $M\oplus M'$ is
generated by $M\cup M'$, any element $z\in M\oplus M'$ can be written in
form $z=t(x_1,\ldots,x_m,y_1,\ldots,y_n)$ for some $t\in V(m+n)$, 
$x_i\in M$, $y_j\in M'$, $m$, $n\geq0$. Using definition of $V$ we obtain
\begin{equation}
z=\sum_i\lambda_ix_i+\sum_j\mu_jy_j,\quad
\text{where $\sum_i|\lambda_i|+\sum_j|\mu_j|\leq 1$, $x_i\in M$, $y_j\in M'$.}
\end{equation}
Put $\lambda:=\sum_i|\lambda_i|$, $\mu:=\sum_j|\mu_j|$, 
$x:=\sum_i(\lambda^{-1}\lambda_i)x_i\in M$, $y:=\sum_j(\mu^{-1}\mu_j)y_j
\in M'$; if $\lambda=0$, then all $\lambda_i=0$, so we put $x:=0$, and 
similarly $y=0$ if $\mu=0$. Then $|\lambda|+|\mu|=\lambda+\mu\leq 1$, and
$z=\lambda x+\mu y$, i.e.\ {\em any element $z$ of $M\oplus M'$ can be written
in form $\lambda x+\mu y$ for $x\in M$, $y\in M'$, $|\lambda|+|\mu|\leq1$.}
This decomposition is not unique; however, elements $\lambda x\in M$
and $\mu y\in M'$ are completely determined by $z$, being its images
under projections $M\oplus M'\to M$ and $M\oplus M'\to M'$.

This is applicable in particular to $z\in M\oplus L_V(n)$. In this case we
put $\|z\|:=\|\mu y\|\leq 1$. This number is well-defined, $\|z\|=0$ iff
$z\in M$, and $\|z\|=1$ implies $|\mu|\geq 1$ since $\|y\|\leq 1$,
$y$ being an element of $L_V(n)=V(n)\subset K^n$, hence $|\mu|=1$,
$\lambda=0$ in view of $|\lambda|+|\mu|\leq 1$, i.e.\ {\em
$\|z\|=1$ implies $\lambda=0$ and $z\in L_V(n)$.}

\nxsubpoint\label{sp:mingen.proj.v} (First consequence of minimality of~$n$.)
Now let us return back to the situation of~\ptref{sp:proj.over.v}, i.e.\
we still have a projective $V$-module $P$, determined by a projector
$a=a^2\in M(n;V)\subset M(n;K)$ with minimal possible~$n$. Put
$I:=\{i\in\stn:\|u_i\|=1\}$, where $u_i=a(e_i)$ as before.

Let $x\in K^n$ be any vector. Then $\|a(x)\|=\sum_i|\sum_j a_{ij}x_j|\leq
\sum_i\sum_j|a_{ij}|\cdot|x_j|=\sum_j|x_j|\cdot(\sum_i|a_{ij}|)\leq
\sum_j|x_j|=\|x\|$ since $\sum_i|a_{ij}|=\|a(e_j)\|\leq 1$ for any~$j$.
Now suppose that we have an equality: $\|a(x)\|=\|x\|$. This implies 
in particular an equality in $\sum_j|x_j|\cdot(\sum_i|a_{ij}|)\leq
\sum_j|x_j|$, which is possible only if for any~$j$ we have either
$x_j=0$ or $\sum_i|a_{ij}|=1$, i.e.\ $j\in I$. Therefore, {\em
if $\|a(x)\|=\|x\|$, then $x_i=0$ for $i\not\in I$.}

Since $a^2=a$, we have $a(u_j)=a^2(e_j)=a(e_j)=u_j$ for any~$j$, 
and in particular
$\|a(u_j)\|=\|u_j\|$, so by the above observation $a_{ij}=0$ for
$i\not\in I$ and any~$j$. This means that $u_j$ lies in the $V$-submodule
generated by $\{e_i\}_{i\in I}$; applying $a$ we see that any
$u_j=a(u_j)$ lies in the $V$-span of $\{u_i=a(e_i)\}_{i\in I}$, i.e.\
$P$ is generated already by $\{u_i\}_{i\in I}$. Since $n$ was chosen
to be minimal possible, we must have $I=\stn$, i.e.
\begin{equation}
\|u_i\|=\|a(e_i)\|=\sum_{j=1}^n|a_{ji}|=1\quad
\text{for all $1\leq i\leq n$.}
\end{equation}

\nxsubpoint\label{sp:mingen.perfcof.v} 
(Consequence for constant perfect cofibrations.)
Now let $f:M\to N$ be a constant perfect cofibration, $u'_1$, \dots,
$u'_n$ be a minimal system of relative generators of $N$ over $M$ as
before, and define $p=p^2\in\End_V(M\oplus L_V(n))$, $a=a^2\in\End_V(L_V(n))$,
$P$ as in~\ptref{sp:const.perf.cof.over.v}. Then $u'_i=p(e_i)\in N\subset
M\oplus L_V(n)$, hence it can be written in form
\begin{equation}
u'_i=\lambda_iv_i+\mu_iw_i,\quad\text{where $|\lambda_i|+|\mu_i|\leq1$,
$v_i\in M$, $w_i\in L_V(n)$.}
\end{equation}
Furthermore, $\mu_i\neq0$, otherwise $u'_i$ would lie in $M$ and we could
remove it from our system of relative generators, thus obtaining a 
contradiction with minimality of~$n$, and $\mu_iw_i$ equals 
$u_i=a(e_i)$, the projection of $u'_j\in M\oplus L_V(n)$ to~$L_V(n)$.

Now put $I:=\{i:\|u_i\|=1\}=\{i:\|u'_i\|=1\}$ as before. For any~$i$
we have $a(u_i)=u_i$ and $u_i=\mu_iw_i$ for $\mu_i\neq0$, hence
$a(w_i)=w_i$ and $\|a(w_i)\|=\|w_i\|$ as well, and 
by~\ptref{sp:mingen.proj.v} we obtain $(w_i)_j=0$ for $j\not\in I$,
i.e.\ $w_i=\sum_{j\in I}(w_i)_je_j$, and $\sum_j|(w_i)_j|\leq 1$ since
$w_i\in L_V(n)$. Now $p(w_i)=\sum_{j\in I}(w_i)_ju'_j$ belongs to the
$V$-submodule $N'$ of $M\oplus L_V(n)$ generated by $M$ and 
$\{u'_j\}_{j\in I}$, hence the same is true for
$u'_i=p(u'_i)=p(\lambda_i v_i+\mu_iw_i)=\lambda_i v_i+\mu_i p(w_i)$,
taking into account that $p|_M=\id_M$ and $|\lambda_i|+|\mu_i|\leq 1$.
In other words, $\{u'_i\}_{i\in I}$ is a smaller system of relative generators
of $N$ over $M$ unless $I=\stn$. Since $n$ was supposed to be minimal,
this means $\|u_i\|=\|u'_i\|=|\mu_i|\cdot\|w_i\|=1$ for all $i$, hence
$|\mu_i|=1$ and $\lambda_i=0$, i.e.\ {\em all $u'_i$ actually lie in
$L_V(n)\subset M\oplus L_V(n)$.} This means exactly that
$p=\id_M\oplus a$, i.e.\ {\em any constant perfect cofibration over~$V$
is of the form $M\to M\oplus P$ for a finitely generated projective~$P$.}

Therefore, we have reduced \ptref{th:proj.over.archv},(c) to
\ptref{th:proj.over.archv},(a) as promised before, so it remains to show
the latter statement.

\nxsubpoint (Finite Markov chains defined by projective $V$-modules.)
Now let us forget about perfect cofibrations and deal with a projective
$V$-module~$P$ with a minimal system of generators $u_1$, \dots, $u_n$,
given by a matrix $a=a^2=(a_{ij})\in M(n;V)\subset M(n;K)$ 
as in~\ptref{sp:proj.over.v}. We have already seen that minimality of~$n$
implies $\|u_j\|=\sum_i|a_{ij}|=1$ for all~$j$. Put $p_{ij}:=|a_{ij}|$,
$0\leq p_{ij}\leq 1$, thus obtaining a matrix $p=(p_{ij})\in 
M(n;\tilde\bbQ_+)\subset M(n;\bbR)$. We claim that $p$ is still a projector,
i.e.\ $p^2=p$. Indeed, $a^2=a$ means $a_{ik}=\sum_j a_{ij}a_{jk}$, hence
$p_{ik}=|a_{ik}|\leq \sum_j|a_{ij}|\cdot|a_{jk}|=\sum_j p_{ij}p_{jk}$.
On the other hand, $1=\sum_i p_{ik}\leq\sum_{i,j}p_{ij}p_{jk}=\sum_j p_{jk}=1$,
hence all inequalities under consideration have to be equalities, meaning that
$p_{ik}=\sum_j p_{ij}p_{jk}$, i.e.\ $p^2=p$ as claimed.

Notice that one can associate to $p=(p_{ij})$ a finite Markov chain with
$n$ states by putting the probability of going from state $j$ into state $i$
in one step equal to $p_{ij}$. Condition $\sum_ip_{ij}=1$ now means that
the sum of probabilities to end up into some state starting from state~$j$
equals one, and $p^2=p$ means that the probability of going from $i$ into $j$
in two steps equals the probability to go from $i$ into $j$ in one step,
so we've got some very special Markov chains here.

We consider also the directed graph $\Gamma$ with $n$ vertices,
such that the edge $i\to j$ belongs to $\Gamma$ iff $p_{ij}>0$, i.e.\
iff $a_{ij}\neq 0$.

\nxsubpoint (Transitivity of $\Gamma$.)
We claim that {\em if $p_{ij}>0$ and $p_{jk}>0$, then $p_{ik}>0$.}
Indeed, we know that all $p_{ij}\geq0$, and if $p_{ij}>0$ and $p_{jk}>0$,
then at least one summand in $p_{ik}=\sum_j p_{ij}p_{jk}$ is strictly
positive, and the others are non-negative, hence their sum $p_{ik}>0$.

\nxsubpoint\label{sp:refl.Gamma} (Reflexivity of $\Gamma$.)
Another immediate consequence of minimality of~$n$ is that {\em for 
all $i$ the edge $i\to i$ belongs to $\Gamma$, i.e.\ $a_{ii}\neq 0$ 
and $p_{ii}>0$.} Indeed, if $a_{ii}=0$, then $u_i$ is a $V$-linear
combination of elements $\{e_j\}_{j\neq i}$, hence $u_i=a(u_i)$ is
a $V$-linear combination of $\{u_j\}_{j\neq i}$, hence $\{u_j\}_{j\neq i}$
already generate $P$ as a $V$-module. This contradicts the minimality of~$n$.

\nxsubpoint (Connectedness of $\Gamma$.)
We can suppose that $\Gamma$ is connected. Indeed, otherwise we can split
our index set $\stn$ into $S\sqcup T$, such that $a_{ij}=0$ unless 
$(i,j)\in S\times S\cup T\times T$. Renumbering our generators in such
a way that all elements of $S$ are listed before all elements of~$T$,
we obtain a decomposition $a=a'\oplus a''$, where $a'$ and $a''$ are two
projectors of smaller size, defining a decomposition $P=P'\oplus P''$;
now it would be sufficient to show that both $P'$ and $P''$ are free to 
conclude that $P$ is free, so we proceed by an induction in~$n$, and 
we have to consider only the case of a connected~$\Gamma$.

\nxsubpoint (Strong connectedness of $\Gamma$.)
Now we assume $\Gamma$ to be connected and show that it can be even
supposed to be {\em strongly connected}, i.e.\ $p_{ij}>0$ for any $i$ and~$j$.
(This is equivalent to the ergodicity of the corresponding Markov chain.) 
Indeed, for any $j\in\stn$ consider $S_j:=\{i:p_{ij}>0\}$. Reflexivity
of $\Gamma$ means $j\in S_j$, and transitivity means $S_k\subset S_j$
for any $k\in S_j$. Now choose $i_0$, such that $|S_{i_0}|$ is minimal
possible, and put $S:=S_{i_0}$. Clearly, $S\neq\emptyset$ since $i_0\in S$,
and $S_i\subset S$ for any $i\in S$ implies $S_i=S$ for any $i\in S$
in view of the minimality of $|S|$. In other words, for any $i$ and 
$j\in S$ we have $p_{ij}>0$, and for $i\not\in S$, $j\in S$ we have $p_{ij}=0$.
If $S=\stn$, $\Gamma$ is strongly connected and we are done. Otherwise
let's reorder our indices in $\stn$ so that $S=\stk=\{1,2,\ldots,k\}$,
$0<k<n$. Since $a_{ij}=0$ for $i>k$, $j\leq k$, we see that the matrix $a$
is block-diagonal: $a=\smallmatrix{a'}{c}{0}{a''}$. Clearly $a'=(a')^2$,
and $a'\in M(k;V)$ since the sums of valuations of elements taken along columns
of $a'$ coincide with such sums computed along corresponding columns 
of~$a$. Similarly, $a''=(a'')^2$ and $a''\in M(n-k;V)$. Denote by
$P'$ and $P''$ projective $V$-modules defined by $a'$ and $a''$.
Then $(a',a)$ describe $f:P'\to P$ as a retract of $L_V(k)\to L_V(n)$,
hence $f:P'\to P$ is a perfect cofibration between projective $V$-modules
with cofiber equal to~$P''$. Clearly the $(u_i)_{1\leq i\leq k}$ generate
$P'$ as a $V$-module, and the $(u_i)_{k<i\leq n}$ constitute a system of
relative generators of $P$ over $P'$. Since the union of any system of
generators of $P'$ and any system of relative generators of $P$ over $P'$
consititutes a system of generators of~$P$, minimality of~$n$ for~$P$ implies 
minimality of $k$ for~$P'$ and of $n-k$ for~$P$ relatively to $P'$.
We have already shown in~\ptref{sp:mingen.perfcof.v} 
that in this situation we have $P=P'\oplus P''$; since $P'$ is generated by
$k<n$ elements and $P''$ by $n-k<n$ elements, they are free by induction
hypothesis, hence $P$ is also free, and we are done.

\nxsubpoint (All columns of $p$ are equal.)
Now we assume $\Gamma$ to be strongly connected and 
prove that all columns of $p$ are
equal, i.e.\ $p_{ij}=p_{ii}=:\delta_i>0$ for all $i$ and~$j$. First of all,
denote by $v\in\bbR^n$ the vector with all components equal to~$1$.
Then $p^tv=v$: indeed, $(p^tv)_j=\sum_ip_{ij}v_i=\sum_ip_{ij}=1$. Let
$x\in\bbR^n$ be another vector with property $p^tx=x$. We claim that 
{\em $x$ is necessarily proportional to $v$, i.e.\ all components $x_i$ 
of~$x$ are equal.} Indeed, put $m:=\min_i x_i$; replacing $x$ with
$x-mv$, we may assume $m=0$, i.e.\ all $x_i\geq0$, and $x_k=0$ for at least
one~$k$. We are going to show that $x=0$. Indeed, we know that
all $x_i\geq0$, $x_k=0$, all $p_{ik}>0$ by strong connectedness of~$\Gamma$,
and $x_k=\sum_ip_{ik}x_i$. This is possible only if all $x_i=0$, i.e.\ 
if $x=0$.

Now let's fix any index $i$ and put $x$ equal to the transpose of
$i$-th row of~$p$: $x_j:=p_{ij}$. Since $p^2=p$, we have $p^tx=x$, hence
all components of $x$ have to be equal, i.e.\ $p_{ij}=p_{ii}$ for all
$j$ and $i$ as claimed.

\nxsubpoint (Renormalization of $a$.)
Up to now we haven't used strict convexity, but only relation 
``$|x|\in\tilde\bbQ_+$ for all $x\in K$''. Now we are going to use
strict convexity. For any $x$, $y\in K^*$ let's write $x\sim y$ 
and say that {\em $x$ and $y$ are proportional\/} if $yx^{-1}\in\tilde\bbQ_+$.
Strict convexity now asserts that $|x+y|=|x|+|y|$ implies $x\sim y\sim x+y$.
An easy induction shows that $|x_1+\cdots+x_k|=|x_1|+\cdots+|x_k|$ 
for $x_i\in K^*$ implies $x_i\sim x_1+\cdots+x_k$ as well.

Consider now equality $a_{ik}=\sum_j a_{ij}a_{jk}$. Since $p_{ik}=\sum
p_{ij}p_{jk}$ as well, and $|a_{ij}|=p_{ij}>0$ by strong connectedness 
of~$\Gamma$, we can apply strict convexity and obtain 
$a_{ik}\sim a_{ij}a_{jk}$ for any $i$, $j$, $k$. In particular,
$a_{ii}\sim a_{ii}^2$, hence $1\sim a_{ii}$, i.e.\ $a_{ii}\in\tilde\bbQ_+$;
we'll write this as ``$a_{ii}>0$'' and ``$a_{ii}=|a_{ii}|=p_{ii}$''. 

Put $\epsilon_i:=a_{1i}/|a_{1i}|$. Clearly, $|\epsilon_i|=1$, hence
$\epsilon_i$ is invertible in~$V$, so we can replace our original
system of generators $u_i$ by $\epsilon_iu_i$. This amounts to
replacing $a$ with its conjugate $a'$ by the diagonal matrix with 
entries $\epsilon_i$, i.e.\ $a'_{ij}=\epsilon_ia_{ij}\epsilon_j^{-1}$.
This matrix $a'$ describes the same projective module~$P$, and we have
$|a'_{ij}|=|a_{ij}|=p_{ij}$. On the other hand, obviously $\epsilon_1=1$
and $a'_{1j}=|a_{1j}|=p_{1j}\in\tilde\bbQ_+$ by construction. Therefore,
replacing $a$ with $a'$, we may assume all $a_{1i}\sim 1$.
Since $a_{1i}a_{ij}\sim a_{1j}$ for any $i$ and $j$, we obtain $a_{ij}\sim 1$,
i.e.\ $a_{ij}\in\tilde\bbQ_+$, or, equivalently, $a_{ij}=|a_{ij}|=p_{ij}$,
i.e.\ we can identify $a$ with~$p$. In particular, all columns of~$a$
are now equal.

\nxsubpoint (End of proof.)
We have just shown that we can assume all columns of~$a$ to be equal,
i.e.\ $u_1=u_2=\cdots=u_n$. Since $\{u_i\}_{1\leq i\leq n}$ was supposed
to be a minimal system of generators of~$P$, this is possible only if $n=1$.
Then $a_{11}=1$ and $P=L_V(1)$ is indeed free as claimed
in~\ptref{th:proj.over.archv}.

\nxsubpoint (Minimal number of generators.)
Notice that the number $n$ of elements in a minimal system of generators of
a projective $V$-module $P$ must be actually equal to the rank of~$P$,
i.e.\ it is the only integer $m\geq0$, for which $P\simeq L_V(m)$.
Indeed, uniqueness of~$m$ follows from \ptref{th:proj.over.archv},(b),
$L_V(m)$ obviously admits a system of $m$ generators, and doesn't
admit a smaller system of generators just because the $K$-vector space
$L_K(m)=K^m=L_V(m)\otimes_VK$ doesn't.

\nxsubpoint (Consequence for $K^0(V)$.)
An immediate consequence of~\ptref{th:proj.over.archv} 
is that $\hat K^0(V)=K^0(V)=\bbZ$ for a strictly convex $V$ as above, 
e.g.\ $V=\Zninfty$, $\Zinfty$ or $\barZinfty$. 
Indeed, since any finitely generated projective
$V$-module is free, $\bbZ=\hat K^0(\Fempty)\to\hat K^0(V)$ is surjective,
and $\hat K^0(V)\to K^0(V)$ is surjective for any generalized ring~$V$.
On the other hand, we have a map $K^0(V)\to K^0(K)=\bbZ$ induced 
by the embedding of~$V$ into its ``fraction field'' $K$. Since the composite
of all these ring homomorphisms equals $\id_\bbZ$, we see that all of them
have to be isomorphisms as claimed.

\nxsubpoint\label{sp:proj.finfty.not.free} (Example over~$\Finfty$.)
We have just shown that any finitely generated projective $\Zinfty$-module
is free, i.e.\ $\Zinfty$ behaves exactly like a local ring in this respect. 
However, its ``residue field'' $\Finfty$ is not so nice, and it actually
admits not only non-free modules (cf.~\ptref{sp:genfd}),
but also non-free projective modules of finite type!
Indeed, consider the endomorphism $p\in\End_{\Finfty}(L_{\Finfty}(2))$,
given on the basis elements by $p:\{1\}\mapsto\{1\}$, $\{2\}\mapsto
?\{1\}+?\{2\}$ in the notation of~\ptref{sp:examp.finfty}. Then it
is easy to see that $p^2=p$, and the corresponding projective module
$P=p(L_{\Finfty}(2))$ is the $\Finfty$-submodule of $L_{\Finfty}(2)$
generated by $\{1\}$ and $?\{1\}+?\{2\}$. It consists of 5 elements,
namely, 0, $\{1\}$, $?\{1\}+?\{2\}$, $-\{1\}$ and $-?\{1\}-?\{2\}$. Since
$|\Finfty(n)|=3^n$ for any $n\geq0$ (cf.~\ptref{sp:examp.finfty}),
$P$ cannot be a free $\Finfty$-module.
(This example has been communicated to the author by A.~Smirnov.)

\nxpointtoc{Vector bundles over $\CompZ$}
Now we want to obtain a reasonable description of vector bundles~$\sE$
over the compactified $\CompZ$. It will be used in~\ptref{th:K0.compZ} to
compute $\hat K^0(\CompZ)$.

\nxsubpoint (Compactified $\CompZ$.)
Recall that $\CompZ$, the ``compactification of $\Spec\bbZ$'',
can be constructed as follows (cf.~\ptref{p:constr.compz}).
Fix any integer $N>1$, put
$B_N:=\bbZ[1/N]$ and $A_N:=B_N\cap\Zninfty$; then
$B_N$ is a unary localization both of $\bbZ$ and $A_N$ since
$B_N=\bbZ[N^{-1}]=A_N[(1/N)^{-1}]$, i.e.\ $\Spec B_N$ is a principal 
open subscheme both in $\Spec\bbZ$ and $\Spec A_N$. 
Then generalized scheme $\CompZ\strut^{(N)}$ is obtained by gluing
$\Spec\bbZ$ and $\Spec A_N$ along their open subschemes isomorphic
to $\Spec B_N$. 

This generalized scheme $\CompZ\strut^{(N)}$ does depend on the
choice of $N>1$; however, we get natural ``transition morphisms''
$f^M_N:\CompZ\strut^{(M)}\to\CompZ\strut^{(N)}$ whenever $N\mid M$
(cf.~\ptref{sp:trans.maps.compz}),
so these $\CompZ\strut^{(N)}$ constitute a projective system of
generalized schemes, indexed by the filtered set of integers $N>1$,
ordered by divisibility.

Then we define the ``true'' compactification $\CompZ$
as $\projlim_{N>1}\CompZ\strut^{(N)}$. There are two understandings of
this projective limit: we can compute it either in the
category of pro-generalized schemes (cf.~\ptref{sp:compZ.pro.sch}),
or in the category of generalized ringed spaces 
(cf.~\ptref{sp:compZ.gen.rgsp}).
Both approaches are equivalent when we consider only finitely presented
objects over~$\CompZ$ (cf.~\ptref{sp:equiv.2.appr}); 
since we are now concerned only with
perfect simplicial $\sO_{\CompZ}$-modules, all components of which
are vector bundles, hence finitely presented, we may freely
choose between these two approaches.

\nxsubpoint\label{sp:finpres.compz.mod} 
(Finitely presented $\sO_{\CompZ}$-modules.)
Recall that we have an open embedding $\Spec\bbZ\to\CompZ$,
as well as ``point at infinity'' $\hat\eta:\Spec\Zninfty\to\CompZ$
and ``generic point'' $\hat\xi:\Spec\bbQ\to\CompZ$
(cf.~\ptref{p:constr.compz}).

Now let $\sF$ be any finitely presented $\sO_{\CompZ}$-module. 
Then we put $M_\bbZ:=\Gamma(\Spec\bbZ,\sF|_{\Spec\bbZ})$,
$M_\infty:=\sF_\infty=\Gamma(\Spec\Zninfty,\hat\eta^*\sF)$,
and finally $M_\bbQ:=\sF_\xi:=\Gamma(\Spec\bbQ,\hat\xi^*\sF)$. Clearly,
$M_\bbZ$ is a finitely presented $\bbZ$-module, $M_\infty$ is
a finitely presented $\Zninfty$-module, and $M_\bbZ\otimes_\bbZ\bbQ
\cong M_\bbQ\cong M_\infty\otimes_{\Zninfty}\bbQ$, for example because
$\hat\xi:\Spec\bbQ\to\CompZ$ can be factorized both through $\Spec\bbZ\subset
\CompZ$ and $\hat\eta:\Spec\Zninfty\to\CompZ$. Therefore,
we obtain a canonical isomorphism of $\bbQ$-vector spaces
$\theta_M:M_\bbZ\otimes_\bbZ\bbQ\simto M_\infty\otimes_{\Zninfty}\bbQ$.

In this way we obtain a functor $F$ from the category of finitely presented
$\sO_{\CompZ}$-modules into the category $\cC$ of triples $(M_\bbZ,M_\infty,
\theta_M)$, consisting of a finitely generated $\bbZ$-module~$M_\bbZ$,
a finitely presented $\Zninfty$-module $M_\infty$, and an isomorphism
of $\bbQ$-vector spaces $\theta_M:M_\bbZ\otimes_\bbZ\bbQ\simto
M_\infty\otimes_{\Zninfty}\bbQ$. Furthermore, we have seen
in~\ptref{sp:cat.finpres.compzmod} that {\em this functor~$F$ is
actually an equivalence of categories.} In particular, {\em the category 
of finitely presented $\sO_{\CompZ}$-modules doesn't depend on the
localization theory $\cT^?$ used to construct spectra involved.}

\nxsubpoint (Vector bundles over $\CompZ$.)
Since any vector bundle $\sE$ over $\CompZ$ is in particular a finitely
presented $\sO_{\CompZ}$-module, we can apply the above considerations,
thus obtaining a triple $E=(E_\bbZ,E_\infty,\phi)$, where $E_\bbZ$
is a finitely generated projective (hence free) $\bbZ$-module, 
$E_\infty$ is a finitely generated projective (hence free 
by~\ptref{th:proj.over.archv}) $\Zninfty$-module, and $\phi:E_\bbZ
\otimes_\bbZ\bbQ\to E_\infty\otimes_{\Zninfty}\bbQ$ is an isomorphism of their
scalar extensions to~$\bbQ$. Conversely, any such triple 
$(E_\bbZ,E_\infty,\phi)$ defines a finitely presented $\sO_{\CompZ}$-module,
which is locally free over $\Spec\bbZ$, and free in some neighbourhood of
$\infty$ as well by the inductive limit argument, since $\sE_\infty=E_\infty
\simto L_{\Zninfty}(r)$ is free. Therefore, $\sE$ is a vector bundle.
Among other things, we have just shown that {\em any vector bundle
$\sE$ over~$\CompZ$ is locally free}, i.e.\ {\em vector bundles over~$\CompZ$
in the retract and retract-free versions of the theory coincide.}
Therefore, we may limit ourselves to the consideration of the retract-free
version of vector bundles. Using \ptref{th:proj.over.archv},(c), we see
that this is equally true for constant perfect cofibrations between 
vector bundles.

\nxsubpoint (Rank of a vector bundle over $\CompZ$.)
Another immediate consequence is that {\em any vector bundle~$\sE$ over
$\CompZ$ has a well-defined rank $r\geq0$, equal to $\dim_\bbQ\sE_\xi$.}
In other words, $\sE$ is locally isomorphic to $\sO_{\CompZ}^{(r)}=
L_{\sO_{\CompZ}}(r)$. Indeed, this is obvious over~$\Spec\bbZ$, and true
over $\sO_{\CompZ,\infty}=\Zninfty$ by~\ptref{th:proj.over.archv}, hence
true over some neighbourhood of infinity by the inductive limit argument.

\nxsubpoint\label{sp:matr.descr.vb.compz} 
(Matrix description of a vector bundle over~$\CompZ$.)
Let $\sE$ be a vector bundle over $\CompZ$ of rank~$r$. We know that
$\sE$ is completely determined by $(E_\bbZ,E_\infty,\phi)$, where
$E_\bbZ\simeq\bbZ^r$, $E_\infty\simeq\Zninfty^{(r)}$, and
$\phi$ is a $\bbQ$-isomorphism between $E_\bbZ\otimes_\bbZ\bbQ$ and
$E_\infty\otimes_{\Zninfty}\bbQ$. Another equivalent description:
$E=\sE_\xi$ is a $\bbQ$-vector space of dimension~$r$, $E_\bbZ\subset E$
is a $\bbZ$-structure on~$E$ (i.e.\ a lattice in~$E$), and $E_\infty\subset E$
is a free $\Zninfty$-structure on~$E$ (i.e.\ an octahedron 
inside~$E$ centered at the origin).

Such a situation admits a description in terms of bases and matrices.
Namely, let us choose any base $(e_i)_{1\leq i\leq r}$ of $\Zninfty$-module 
$E_\infty$ (i.e.\ $\{\pm e_i\}$ are the vertices of octahedron $E_\infty$),
and any base $(f_i)_{1\leq i\leq r}$ of $\bbZ$-lattice $E_\bbZ\subset E$.
Since both $(e_i)$ and $(f_i)$ are bases of $\bbQ$-vector space $E$,
they are related to each other by means of a matrix $A=(a_{ij})\in GL_r(\bbQ)$:
\begin{equation}
e_i=\sum_{j=1}^r a_{ij}f_j,\quad a_{ij}\in\bbQ
\end{equation}

\nxsubpoint\label{sp:vb.as.cosets} (Dependence on the choice of bases.)
If we choose another base $(f'_i)$ for $E_\bbZ$, it is related to
$(f_i)$ by means of a matrix $B=(b_{ij})\in GL_r(\bbZ)$, where
$f_i=\sum_j b_{ij}f'_j$. Then $e_i=\sum_{j,k}a_{ij}b_{jk}f'_k$, i.e.\ 
$A$ is replaced by $A'=AB$. Similarly, if we replace $(e_i)$ by another base
$(e'_i)$ of $E_\infty$, these two bases are related to each other by means
of $C=(c_{ij})\in GL_r(\Zninfty)=\Oct_r$, where $e'_i=\sum_j c_{ij}e_j$,
and we obtain $A'=CA$. Therefore, {\em 
multiplying $A$ from the left by matrices from $GL_r(\Zninfty)=\Oct_r$
and from the right by matrices from $GL_r(\bbZ)$ doesn't change corresponding
vector bundle~$\sE$.} Conversely, if two matrices $A$ and $A'$ define 
isomorphic vector bundles, we can assume (by transporting all structure)
these vector bundles to coincide, i.e.\ $A$ and $A'$ would arise from different
choices of bases in same modules $E_\infty$ and $E_\bbZ$, hence
$A'=CAB$ for some $C\in GL_r(\Zninfty)$ and $B\in GL_r(\bbZ)$. We have just
shown the following statement:

{\em Isomorphism classes of vector bundles $\sE$ of rank~$r$ over~$\CompZ$
are in one-to-one correspondence with double cosets $GL_r(\Zninfty)\backslash
GL_r(\bbQ)/ GL_r(\bbZ)$.}

This description of isomorphism classes of vector bundles in terms of
double cosets of matrices is completely similar to the classical description
of vector bundles over a projective curve over a field.

\nxsubpoint (Operations on vector bundles in terms of matrices.)
All reasonable operations with vector bundles, such as direct sums,
tensor products, symmetric, tensor and exterior powers, 
can be computed by means
of corresponding operations with matrices. For example, let $\sE$ and
$\sE'$ be two vector bundles of rank $r$ and $r'$, respectively,
and $(e_i)_{1\leq i\leq r}$, $(f_i)_{1\leq i\leq r}$, 
$(e'_j)_{1\leq j\leq r'}$, $(f'_j)_{1\leq j\leq r'}$ be the bases used 
to construct matrices $A\in GL_r(\bbQ)$ and $A'\in GL_{r'}(\bbQ)$ 
representing $\sE$ and $\sE'$. Then $\sE\otimes_\sO\sE'$ is
given by $E_\bbZ\otimes_\bbZ E'_\bbZ$ and $E_\infty\otimes_{\Zninfty}E'_\infty$
in $E\otimes_\bbQ E'$, and $(f_i\otimes f'_j)_{i,j}$, $(e_i\otimes e'_j)_{i,j}$
constitute natural bases for these tensor products of free modules. The
matrix relating these two bases will be exactly the Kronecker product
$A\otimes A'$ of matrices $A$ and $A'$ as claimed.

As to the other operations, they can be dealt with in a completely similar
fashion, by considering naturally arising bases. Notice that for 
the exterior powers we use alternativity of $\Zninfty$ (valid since
$\Zninfty\subset\bbQ$) to assure that an exterior power of a free module
is still free, with the base given in the classical way.

\nxsubpoint\label{sp:pic.compz} (Picard group of $\CompZ$.)
Let's compute Picard group $\Pic(\CompZ)$ using
the description of vector bundles given in~\ptref{sp:vb.as.cosets}.
We get $\Pic(\CompZ)=GL_1(\Zninfty)\backslash GL_1(\bbQ)/GL_1(\bbZ)=
\{\pm1\}\backslash\bbQ^*/\{\pm1\}=\bbQ^*_+$, 
(cf.\ also~\ptref{sp:pic.compz0}), 
i.e.\ {\em isomorphism classes
of line bundles over $\CompZ$ are parametrized by positive rational
numbers.} We have just seen that tensor product of vector bundles corresponds
to Kronecker product of matrices; for line bundles this means product of
rational numbers, i.e.\ {\em $\Pic(\CompZ)\cong\bbQ^*_+$ as an abelian group.}
Notice that $\bbQ^*_+=\bbZ^{(\bbP)}$, where $\bbP$ is the set of prime
numbers, i.e.\ {\em $\Pic(\CompZ)$ is a free abelian group with countably
many generators.}

Since $\Pic(\CompZ)$ is often written additively, and $\bbQ^*_+$ is written
multiplicatively, we denote by $\log\bbQ^*_+$ the group $\bbQ^*_+$ written
in additive form; its elements will be often denoted by $\log\lambda$,
$\lambda\in\bbQ^*_+$ (this is just a formal notation, but we can identify
$\log\bbQ^*_+$ with a subgroup of $\bbR^+$ by means of the ``true''
logarithm if we like). Thus $\log\lambda+\log\mu=\log (\lambda\mu)$.

We denote by $\sO(\log\lambda)=\sO_{\CompZ}(\log\lambda)$ the line bundle
on $\CompZ$ given by $1\times 1$-matrix $(\lambda)$. Thus
$x\mapsto[\sO(x)]$ is an isomorphism $\log\bbQ^*_+\simto\Pic(\CompZ)$.

\nxsubpoint (``Serre twists'' on $\SpecZ$.)
Given any $\sO_{\CompZ}$-module $\sF$ and any $x\in\log\bbQ^*_+$, 
we denote by $\sF(x)$ the corresponding ``Serre twist'' 
$\sF\otimes_\sO\sO(x)$. Since $\sO(x+y)=\sO(x)\otimes_\sO\sO(y)$,
we have $\sF(x+y)=(\sF(x))(y)$. This operation has nearly all the properties
one has for Serre twists on a projective line over a field
(with the obvious difference that the Picard group of a projective line
is $\bbZ$, not $\log\bbQ^*_+$).

For example, if $\sE$ is a vector bundle of rank~$r$, so is $\sE(\log\lambda)$;
furthermore, if $\sE$ is given by matrix $A$, then $\sE(\log\lambda)$ is
given by~$\lambda A$.

\nxsubpoint (Determinant line bundles.)
Let $\sE$ be a vector bundle of rank~$r$ on $\CompZ$, given by a matrix~$A$.
In the classical case $\bbP^1_k$ of a projective line over a field we have
$c_1(\sE)=c_1(\det\sE)$, and the first Chern class $c_1:\Pic(\bbP^1)\to
CH^1(\bbP^1)$ induces an isomorphism, both sides being isomorphic to~$\bbZ$.
(In this case $c_1(\sE)$ is actually the {\em degree\/} of~$\sE$.)
Therefore, over $\bbP^1$ we might identify $c_1(\sE)\in CH^1(\bbP^1)$
with $[\det\sE]\in\Pic(\bbP^1)$. We can try to do the same over $\CompZ$.
In our case $\det\sE$ will be given by $\det A$, and its class in 
$\Pic(\CompZ)$ is given by $\log|\det A|\in\log\bbQ^*_+$. We'll see later
that this is indeed a valid computation of $\deg\sE=c_1(\sE)$ over~$\CompZ$;
we'll just remark now that this computation already shows us how
the logarithms of volumes appear as arithmetic intersection
numbers or arithmetic degrees of vector bundles in Arakelov geometry.

\nxsubpoint (Homomorphisms of vector bundles.)
Since the functor $F:\sE\mapsto(\Gamma(\Spec\bbZ,\sE),\sE_\infty,\phi_\sE)$
of~\ptref{sp:finpres.compz.mod} is an equivalence of categories, we can compute
$\Hom_{\CompZ}(\sE,\sE')$ for any two vector bundles $\sE$ and $\sE'$
in terms of their matrices. In fact, if $\sE$ is given by $E_\bbZ$ and
$E_\infty$ inside $\bbQ$-vector space~$E$, and $\sE'$ is given similarly
by $E'_\bbZ$ and $E'_\infty\subset E'$, then $\Hom_{\CompZ}(\sE,\sE')$
can be identified with $\bbQ$-linear maps $u:E\to E'$, such that
$u(E_\bbZ)\subset E'_\bbZ$ and $u(E_\infty)\subset E'_\infty$. If we 
choose bases as in~\ptref{sp:matr.descr.vb.compz}, thus obtaining matrices
$A\in GL_r(\bbQ)$, $A'\in GL_{r'}(\bbQ)$, then we can consider the 
(transposed) matrix
$U\in M(r,r';\bbQ)$ of $u$ with respect to bases $(e_i)$ and $(e'_j)$,
defined by $u(e_i)=\sum_j u_{ij}e'_j$.
Then the conditions for $U$ to define a homomorphism $\sE\to\sE'$ are
the following:
\begin{itemize}
\item $\sum_j|u_{ij}|\leq1$ for any $1\leq i\leq r$ (this is equivalent to
$u(E_\infty)\subset E'_\infty$).
\item $A^{-1}UA'\in M(r,r';\bbZ)$ 
(this expresses $u(E_\bbZ)\subset E'_\bbZ$).
\end{itemize}

\nxsubpoint (Global sections of~$\sE$ and Mumford-regular vector bundles.)
The above considerations can be applied to
$\Gamma(\CompZ,\sE)=\Hom_{\sO_{\CompZ}}(\sO,\sE)$. We obtain immediately
that the $\Fpm$-module $\Gamma(\CompZ,\sE)$ can be identified with
$E_\bbZ\cap E_\infty\subset E$. Notice that this intersection is
always a finite set, more or less for compactness (or boundedness) reasons.

Similarly, $\Gamma_{\CompZ}(\sE(\log\lambda))\cong\Hom_\sO(\sO(-\log\lambda),
\sE)$ can be identified with $E_\bbZ\cap\lambda E_\infty$, i.e.\ 
$f_\sE:\log\lambda\mapsto\card(E_\bbZ\cap\lambda E_\infty)$ is
the ``Hilbert function'' of~$\sE$. We see that Hilbert functions
of vector bundles over~$\CompZ$ are closely related to the problem of
counting lattice points inside polyhedra (octahedra in this case) with
rational vertices. Another interesting thing is that asymptotically
$f_\sE(\log\lambda)=r\log\lambda+\log|\det A|+\log(2^r/r!)+o(1)$
when $\lambda\to+\infty$,
i.e.\ the rank and the ``arithmetic degree'' $c_1(\sE)=\log|\det A|$ 
determine the asymptotic behavior of~$f_\sE$.

We say that a vector bundle $\sE$ is {\em Mumford-regular\/} if 
all the vertices of octahedron $E_\infty\subset E$ lie in lattice $E_\bbZ$,
i.e.\ if $\sE$ is given by a matrix~$A$ with integer coefficients.
(Usually a vector bundle $\sE$ over $\bbP^N$ is said to be Mumford-regular
if higher cohomology groups $H^q(\bbP^N,\sE(-q))$ vanish for all $q>0$; 
any such $\sE$ is known to be generated
by its global sections.) However, in our case a Mumford-regular vector
bundle needn't be generated by its global sections (i.e.\ 
$E_\infty\cap E_\bbZ$ doesn't necessarily generate $E_\bbZ$), as
illustrated by the vector bundle of rank 3 defined by
\begin{equation}
A=\left(\begin{matrix}1\\&1\\-1&-1&2\end{matrix}\right)
\end{equation}

Clearly, for any vector bundle $\sE$ we can find a Serre twist
$\sE(\log\lambda)$, which is Mumford-regular. More precisely, there 
is a minimal $\lambda_0$ with this property (namely, the inverse of the 
g.c.d.\ 
of all entries $a_{ij}$ of~$A$, i.e.\ the generator of the fractional
ideal generated by the $a_{ij}$), and all other $\lambda$s are precisely
the multiples of~$\lambda_0$.

\nxsubpoint\label{sp:mumford.antireg} 
(Mumford-antiregularity and homomorphisms into line bundles.)
Let's consider the dual problem. Given any vector bundle $\sE$
of rank~$r$,
we can consider $\Hom_\sO(\sE,\sO(\log\lambda))$. It can be identified with
the set of linear forms $u\in E_\bbZ^*\subset E^*$, such that
$u(E_\infty)\subset[-\lambda,\lambda]$, 
i.e.\ $u\in E_\bbZ^*\cap \lambda E_\infty^*$, where
$E_\infty^*$ is the dual of $E_\infty$, a cube with rational vertices 
in $E^*$ centered at the origin. In other words, the ``dual Hilbert function''
\begin{equation}
p_\sE(\lambda):=\card\Hom_\sO(\sE,\sO(\log\lambda))=|E_\bbZ^*\cap\lambda 
E_\infty^*|
\end{equation}
counts the number of lattice points of $E_\bbZ^*$ inside cubes $\lambda 
E_\infty^*$. This function has somewhat better properties than
the ``ordinary Hilbert functions'' considered before. For example,
\begin{equation}
p_{\sE\oplus\sE'}(\lambda)=p_\sE(\lambda)\cdot p_{\sE'}(\lambda)
\end{equation}
an equality that has no counterpart for ``ordinary Hilbert functions''.
Furthermore, we have $\log p_\sE(\lambda)=r\log(2\lambda)-\log|\det A|
+o(1)$ asymptotically in $\lambda\to+\infty$, without any additional constants,
where $A$ is any matrix definining~$\sE$.

We say that $\sE$ is {\em Mumford-antiregular\/} if the dual base $(e_i^*)$
of the base $(e_i)$ consisting of vertices of octahedron $E_\infty$,
i.e.\ the centers of faces of cube $E_\infty^*$, lie in dual lattice
$E_\bbZ^*$. Clearly, this is equivalent to requiring all coefficients 
of~$A^{-1}$ to be integer.

Of course, there are $\lambda\in\bbQ^*_+$, such that $\sE(-\log\lambda)$
is Mumford-antiregular, and such $\lambda$'s are multiples of minimal
such $\lambda_0$, equal to the g.c.d.\ of coefficients of~$A^{-1}$.
In particular, if $\sE$ is Mumford-antiregular, so is $\sE(-\log\lambda)$ for
all ``positive'' $\lambda$, i.e.\ for all $\lambda\in\bbN$.

It turns out that for any $\sE$ one can find a polynomial $\tilde p_\sE\in
\bbQ[T]$, such that $p_\sE(\lambda)=\tilde p_\sE(\lambda)$ for all
$\lambda\in\bbQ^*_+$, for which $\sE(-\log\lambda)$ is Mumford-antiregular.
We'll return to this question later, indicating a method for
computation of polynomial~$\tilde p_\sE$.

\nxsubpoint\label{sp:matr.perfcof.vb} 
(Matrix description of perfect cofibrations of vector bundles.)
Let $q:\sE'\to\sE$ be a perfect cofibration with cofiber $\sE''$, i.e.\
$0\to\sE'\stackrel q\to\sE\stackrel\pi\to\sE''\to0$ 
is a cofibration sequence, and
$[\sE]=[\sE']+[\sE'']$ in $\hat K^0(\CompZ)$. We want to obtain a matrix
description of such a situation. Let $r$ be the rank of $\sE'$,
$s$ be the rank of $\sE''$; then the rank of $\sE$ equals $r+s$ since
we have a short exact sequence of generic fibers. Let us denote
$E:=\sE_\xi$, $E_\bbZ:=\Gamma(\Spec\bbZ,\sE)\subset E$ and
$E_\infty:=\sE_\infty$ as before, and similarly for $\sE'$ and $\sE''$.
We identify $E'$ with a $\bbQ$-vector subspace of $E$ by means of $q_\xi$;
then $E''$ can be identified with $E/E'$.
Since $\Spec\bbZ$ is additive, a cofibration sequence of free $\bbZ$-modules
$0\to E'_\bbZ\to E_\bbZ\to E''_\bbZ\to 0$ is nothing else than a short 
exact sequence. Therefore, if $(f'_i)_{1\leq i\leq r}$ is any base of
$E'_\bbZ$, $(f''_j)_{1\leq j\leq s}$ any base of $E''_\bbZ$,
$f_i:=f'_i$ for $1\leq i\leq r$, and $f_{r+i}\in E_\bbZ$ are any lifts of
$f''_i$ to $E_\bbZ$, $1\leq i\leq s$, then $(f_i)_{1\leq i\leq r+s}$ is
a base of $E_\bbZ$. Similarly, $E'_\infty\to E_\infty$ is a perfect 
cofibration of vector bundles over $\Zninfty$, hence by
\ptref{th:proj.over.archv},(c) it is isomorphic to $L_{\Zninfty}(r)\to
L_{\Zninfty}(r)\oplus L_{\Zninfty}(s)=L_{\Zninfty}(r+s)$. In other words,
we can choose a $\Zninfty$-base $(e_i)_{1\leq i\leq r+s}$ of $E_\infty$,
such that its first $r$ elements $(e'_i:=e_i)_{1\leq i\leq r}$ constitute
a base of $E'_\infty$; then the images $e''_i$ of $e_{r+i}$, $1\leq i\leq s$,
constitute a base of $E''_\infty$.

Now let $A\in GL_{r+s}(\bbQ)$ be the matrix relating $(e_i)$ to $(f_i)$,
and define $A'\in GL_r(\bbQ)$ and $A''\in GL_s(\bbQ)$ similarly. Then
$A$, $A'$, $A''$ describe vector bundles $\sE$, $\sE'$ and $\sE''$,
respectively, and by construction $A=\smallmatrix{A'}{0}{*}{A''}$ is
block-triangular with diagonal blocks $A'$ and $A''$. Conversely,
if a vector bundle $\sE$ admits a description in terms of a block-triangular
matrix $A$ of the above form, we obtain morphisms of corresponding 
vector bundles $\sE'\to\sE\to\sE''$; this sequence is a perfect cofibration
sequence over $\Spec\bbZ$, being just a short exact sequence of free
$\bbZ$-modules there, and a perfect cofibration sequence in some neighbourhood
of $\infty$, since $\sE'_\infty\to\sE_\infty\to\sE''_\infty$ is a perfect
cofibration sequence of free $\Zninfty$-modules by construction.

Therefore, {\em perfect cofibrations $\sE'\to\sE$ between vector bundles
and perfect cofibration sequences $\sE'\to\sE\to\sE''$ over $\CompZ$ 
correspond to block-triangular decompositions $A=\smallmatrix{A'}{0}{*}{A''}$
of corresponding matrices.} Of course, we obtain such block-triangular
decompositions only for special choices of bases for $E_\infty$ and $E_\bbZ$,
i.e.\ they hold only for some special representatives $A\in GL_{r+s}(\bbQ)$
of double coset in $GL_{r+s}(\Zninfty)\backslash GL_{r+s}(\bbQ)
/GL_{r+s}(\bbZ)$. However, $A'$ and $A''$ might have been chosen arbitrarily.

\begin{PropD}\label{prop:can.repr.GlQ/GlZ}
{\rm (Canonical representatives for $GL_n(\bbQ)/GL_n(\bbZ)$.)}\\
Any right coset $GL_n(\bbQ)/GL_n(\bbZ)$ contains exactly one matrix
$A=(a_{ij})\in GL_n(\bbQ)$ satisfying following conditions:
\begin{itemize}
\item $A$ is lower-triangular, i.e.\ $a_{ij}=0$ for $i<j$.
\item $A$ has positive diagonal elements, i.e.\ $a_{ii}>0$.
\item $0\leq a_{ij}<a_{ii}$ for all $1\leq j<i\leq n$.
\end{itemize}
\end{PropD}
\begin{Proof}
(a) First of all, notice that the above conditions are invariant under
multiplication of~$A$ by any positive rational number~$\gamma$. Let
us choose any matrix $A$ in given right coset; replacing $A$ by $\gamma A$
if necessary, we can assume $A\in M(n;\bbZ)$, $\det A\neq 0$.
Then the right coset of~$A$ consists only of matrices with integer 
coefficients.

(b) We are free to replace $A$ with $A'=AB$ for any $B\in GL_n(\bbZ)$.
This amounts to applying the same matrix~$B$ to all rows of~$A$,
or to performing some linear operations with columns of~$A$. In particular,
we are allowed to permute two or more columns, change the sign of a column, 
and add to a column any integer multiple of another one.

(c) Consider the matrices~$A$ from given right coset with $a_{11}>0$,
and choose one with minimal $a_{11}$ (this is possible since $a_{11}\in\bbZ$).
Then all elements $a_{1i}$ of the first row of~$A$ must be divisible 
by~$a_{11}$: otherwise we would subtract from the $i$-th column a multiple
of the first column and make $0<a_{1i}<a_{11}$; then we would permute the
first and the $i$-th column and obtain a matrix~$A'$ in the same right coset
with $0<a'_{11}<a_{11}$. This contradicts the minimality of $a_{11}$.
Therefore, any $a_{1i}$ is divisible by $a_{11}$; subtracting from the
$i$-th column $a_{1i}/a_{11}$ times the first column, we obtain a matrix~$A$
with $a_{11}>0$, $a_{1i}=0$ for $1<i\leq n$.

(d) Now consider the submatrix of $A$ formed by its rows and columns
with indices from $2$ to~$n$. Using only operations with these columns,
i.e.\ considering the coset of $A$ in $GL_n(\bbQ)/GL_{n-1}(\bbZ)$,
we can make $a_{22}>0$ and $a_{2i}=0$ for $i>2$ by the same reasoning as
in~(c). Furthermore, subtracting from the first column $[a_{21}/a_{22}]$
times the second column, we can make $0\leq a_{21}<a_{22}$. Notice that
all these manipulations do not affect the first row.

(e) Next, we can make $a_{33}>0$ and $a_{3i}=0$ for $i>3$ by manipulations
with columns $\geq3$, and then subtract suitable multiples of the third 
column from the first and the second column so as to have
$0\leq a_{3i}<a_{33}$ for $i<3$. These operations don't affect the 
two first rows. We proceed further by induction, making the fourth,
\dots, $n$-th row of~$A$ of the form indicated in the proposition.
This proves the existence of a representative with required properties.

(f) Now let us show the uniqueness. Let $A=(a_{ij})$ and $A'=(a'_{ij})$ 
be two matrices with the above properties, 
and suppose that they lie in the same
right coset with respect to $GL_n(\bbZ)$, i.e.\ $A'=AB$ for some
$B\in GL_n(\bbZ)$. First of all, $B$ is lower-triangular with positive
diagonal entries, just because $A$ and $A'$ are such. Next, all coefficients
of $B$ are integers, hence $b_{ii}\geq1$, and the product of $b_{ii}$ equals
$\det B=\pm1$, hence all $b_{ii}=1$, i.e.\ $B$ is lower-unitriangular,
and $a'_{ii}=a_{ii}>0$ for all~$i$. Now suppose that $A\neq A'$, i.e.\ 
$B\neq E$. Then $b_{ij}\neq 0$ for some $i>j$. Choose indices $(i,j)$,
such that $b_{ij}\neq 0$, $i>j$, with minimal $i$.
Then $a'_{ij}=\sum_k a_{ik}b_{kj}=a_{ij}+a_{ii}b_{ij}$
since $b_{kj}=0$ for $k<j$, $b_{jj}=1$, $b_{kj}=0$ for $j<k<i$, and
$a_{ik}=0$ for $k>i$. On the other hand, $0\leq a_{ij}<a_{ii}$ and
$0\leq a'_{ij}<a'_{ii}=a_{ii}$, hence $|a_{ii}b_{ij}|=|a'_{ij}-a_{ij}|<
a_{ii}$. This contradicts $b_{ij}\in\bbZ$, $b_{ij}\neq0$, q.e.d.
\end{Proof}

\nxsubpoint\label{sp:can.descr.vb} (Vector bundles of rank~$n$.)
The above result enables us to choose some canonical representatives
$A\in GL_n(\bbQ)$ for the isomorphism class of a vector bundle $\sE$
of rank~$n$. Namely, we know (cf.~\ptref{sp:vb.as.cosets}) that
the set of isomorphism classes of vector bundles of rank~$n$ over~$\CompZ$
can be identified with $\Oct_n\backslash GL_n(\bbQ)/GL_n(\bbZ)$, where
$\Oct_n=GL_n(\Zninfty)$ is the octahedral group, and in particular is finite.
Since $\Oct_n$ is finite, any such double coset is a finite union of
right cosets $GL_n(\bbQ)/GL_n(\bbZ)$, and these right cosets can be
identified with the set of their canonical representatives in~$GL_n(\bbQ)$
by~\ptref{prop:can.repr.GlQ/GlZ}. Therefore, any isomorphism class of
vector bundles canonically corresponds to a finite set of ``canonical''
matrices as in~\ptref{prop:can.repr.GlQ/GlZ}.

We would like to choose one canonical matrix from this finite set. 
In order to do this notice that when we study representatives~$A$ of such
double cosets as above, we are allowed not only to multiply $A$ from the
right by matrices $B\in GL_n(\bbZ)$, i.e.\ to permute columns of~$A$
and add multiples of one column of~$A$ to any other, but to multiply
$A$ from the left from matrices $C\in\Oct_n$, i.e.\ by permutation
matrices with non-zero entries equal to $\pm1$. In other words, now
we are allowed to permute rows and change their sign as well.

If we would find a canonical way of ordering the rows of~$A$, preserved under
multiplication by matrices from $GL_n(\bbZ)$ from the right, we would
re-order the rows in this way by means of a matrix from $\Oct_n$, 
and then all the freedom left would be to change the sign of a row,
something that must be immediately compensated by changing the sign in the
corresponding column if we insist to mantain all matrices in the canonical
form of~\ptref{prop:can.repr.GlQ/GlZ}.

Denote by $d_i>0$ the g.c.d.\ of all elements $(a_{ij})_{1\leq j\leq n}$ of
the $i$-th row of~$A$. Clearly, the sequence $(d_1,\ldots,d_n)$ doesn't
change when we replace $A$ by $AB$ with $B\in GL_n(\bbZ)$. Therefore,
``generically'' (if all $d_i$ are distinct), we obtain a natural ordering
of the rows of~$A$, and we can require $d_1<d_2<\cdots<d_n$ in order to
make the choice of~$A$ more canonical. If not all $d_i$ are distinct,
we can still require $d_1\leq d_2\leq\cdots\leq d_n$, but now there might
be some degree of freedom left: we may still permute rows with equal values
of~$d_i$.

\nxsubpoint\label{sp:dec.vb.lb} 
(Decomposition of a vector bundle into line bundles.)
Let $\sE$ be a vector bundle of rank~$n$. According to
\ptref{sp:vb.as.cosets} and \ptref{prop:can.repr.GlQ/GlZ},
$\sE$ can be described by a matrix $A\in GL_n(\bbQ)$ of canonical form.
In particular, we can suppose $A$ to be lower-triangular with $a_{ii}>0$.
On the other hand, any block-triangular decomposition of~$A$ corresponds
by~\ptref{sp:matr.perfcof.vb} to a perfect cofibration of vector bundles.
Since $A$ is completely triangular, we obtain a finite filtration
$0=\sE_0\subset\sE_1\subset\cdots\subset\sE_n=\sE$ of $\sE$ by vector bundles,
such that each $\sE_{k-1}\to\sE_k$ is a perfect cofibration of vector
bundles with cofiber equal to line bundle $\sO(\log a_{kk})$.

In particular, $[\sE]=\sum_{k=1}^n[\sO(\log a_{kk})]$, i.e.\ 
{\em $\hat K^0(\CompZ)$ is generated by line bundles.} Furthermore,
$\hat K^0(\CompZ)\to K^0(\CompZ)$ is surjective by~\ptref{cor:surj.K0.to.K0},
$\CompZ$ being quasicompact, hence {\em the same is true for $K^0(\CompZ)$.}

\nxsubpoint (Semistable filtration.)
If all $d_i=\gcd(a_{i1},\ldots,a_{in})$ of \ptref{sp:can.descr.vb} are
distinct, we can make the above filtration on~$\sE$ canonical by choosing
a canonical representative $A$ in such a way that $d_1<d_2<\cdots<d_n$.
In general if $d_1=\cdots=d_{i_1-1}<d_{i_1}=\cdots=d_{i_2-1}<d_{i_2}=\cdots$,
then only the sub-filtration consisting of the $\sE_{i_k}$ is canonical.
Let's say that {\em $\sE$ is semistable of slope $\log d$} if all $d_i=d$.
Then the sub-filtration we've just discussed has semistable cofibers,
i.e.\ is some sort of Harder--Narasimhan filtration. 
Notice, however, that in general the slope $\log d$ 
of a semistable $\sE$ needn't equal $\deg\sE/\rank\sE$:
all we can say is that $d^n$ is a divisor of $\det A$, i.e.\ 
$\deg\sE-\rank\sE\cdot\log d$ is ``positive'' ($\log$ of a natural number).

\begin{LemmaD}\label{l:rel.lb.K0.compz}
{\rm (Relations between line bundles in $\hat K^0(\CompZ)$.}
\begin{itemize}
\item[(a)] Let $a$, $b$, $a'$, $b'\in\bbQ^*_+$ be such that $ab=a'b'$. Then
\begin{equation}\label{eq:rel.lb.K0.compz}
[\sO(\log a)]+[\sO(\log b)]=[\sO(\log a')]+[\sO(\log b')]
\quad\text{in $\hat K^0(\CompZ)$.}
\end{equation}
\item[(b)] The map $\log\lambda\mapsto[\sO(\log\lambda)]-1$ is an
abelian group homomorphism $\log\bbQ^*_+=\Pic(\CompZ)\to\hat K^0(\CompZ)$.
\end{itemize}
\end{LemmaD}
\begin{Proof}
(a) First of all, notice that it suffices to prove \eqref{eq:rel.lb.K0.compz}
under the additional assumption $b/b'\in\bbZ$. Indeed, in the general case
we can always find a common divisor $b''>0$ of both $b$ and $b'$,
and put $a'':=ab/b''$; then it would suffice to show \eqref{eq:rel.lb.K0.compz}
for $a$, $b$, $a''$, $b''$ and for $a'$, $b'$, $a''$, $b''$. 

So let us assume $b/b'\in\bbN$. Consider the line bundle $\sE$ of rank 2
over $\CompZ$ defined by matrix $A=\smallmatrix{a}{}{b'}{b}$. Since this
matrix is triangular, we have $[\sE]=[\sO(\log a)]+[\sO(\log b)]$ 
by~\ptref{sp:dec.vb.lb}. On the other hand, we are free to multiply $A$ 
by matrices $C\in\Oct_2$ from the left; taking $C=\smallmatrix0110$, we get
$A'=CA=\smallmatrix{b'}{b}{a}{}$. Let $A''=A'B=\smallmatrix{\tilde b'}{}
{c}{\tilde a'}$ be the canonical form of $A'$ in the sense 
of~\ptref{prop:can.repr.GlQ/GlZ}. Since the row g.c.d.s of $A'$ and $A''$
must be equal, we get $\tilde b'=\gcd(b',b)=b'$; since $b'\tilde a'=\det A''=\pm\det A=\pm ab=\pm a'b'$ and all numbers involved are positive, we have
$\tilde a'=a'$ as well, i.e.\ $A''=\smallmatrix{b'}{}{c}{a'}$. Since
$A''$ still defines the same vector bundle~$\sE$, we get 
$[\sE]=[\sO(\log b')]+[\sO(\log a')]$ by~\ptref{sp:dec.vb.lb}, whence
\eqref{eq:rel.lb.K0.compz}.

(b) This statement is immediate from (a): we have to check
$[\sO(\log\lambda+\log\mu)]-1=[\sO(\log\lambda)]-1+[\sO(\log\mu)]-1$,
and this is a special case of (a) for $a=\lambda$, $b=\mu$, $a'=\lambda\mu$,
$b'=1$.
\end{Proof}

\nxsubpoint\label{sp:l.ring.from.abgr} 
($\lambda$-ring $\tilde A$ defined by an abelian group~$A$.)
Let $A$ be an arbitrary abelian group. Denote by $\tilde A$ the ring
$\bbZ\times A$, where the multiplication is determined by the requirement
$A^2=0$ in $\tilde A$, i.e.\ $(m,x)\cdot (n,y)=(mn,nx+my)$,
and in particular $(1,x)\cdot(1,y)=(1,x+y)$. The addition
is of course defined componentwise: $(m,x)+(n,y)=(m+n,x+y)$.
We identify $\bbZ$ with subring $\bbZ\times 0$ in $\tilde A$, and $A$ 
with $0\times A\subset\tilde A$, so we can write $n+x$ instead of $(n,x)$
if we like. Clearly, $\tilde A$ is a commutative ring with unity $(1,0)$.
Furthermore, it has a unique pre-$\lambda$-structure, such that
$\lambda^k(1,x)=0$ for $k\geq2$: indeed, uniqueness is evident since
elements $(1,x)$ generate $\tilde A$ as an abelian group, and
writing $(n,x)=(n-1)+(1,x)$ we get $\lambda_t(n,x)=(1+t)^{n-1}\lambda_t(1,x)=
(1+t)^{n-1}(1+t(1+x))$, whence $\lambda^k(n+x)=\binom nk+\binom{n-1}{k-1}x$
for any $k\geq1$, and of course $\lambda^0(n+x)=1$. It is immediate that
these formulas do define a pre-$\lambda$-ring structure on $\tilde A$.
Furthermore, this is actually a $\lambda$-ring structure
(cf.\ SGA~6~V): indeed, we have to check that $\lambda_t:\tilde A\to
1+\tilde A[[t]]^+$ is a $\lambda$-homomorphism; it is already an abelian group 
homomorphism, and elements $1+x$ generate $\tilde A$, so by linearity 
it is enough to check compatibility of $\lambda_t$ with 
$\lambda$-operations and multiplication 
on these generators, where it is evident since $\lambda^k(1+x)=0$ 
for $k>1$.

\begin{LemmaD}\label{l:pic.to.K0}
Let $X=(X,\sO_X)$ be a generalized ringed space. Put $A:=\Pic(X)$,
and construct $\lambda$-ring $\tilde A=\bbZ\times A$ as above.
Suppose that the map $c_1:A=\Pic(X)\to K^0:=\hat K^0(X)$ given by
$\sL\mapsto[\sL]-1$ is a homomorphism of abelian groups. Then:
\begin{itemize}
\item[(a)] The map $\phi:\tilde A\to K^0$ given by $(n,\sL)\mapsto n-1+[\sL]=
n+c_1(\sL)$ is a ring homomorphism and even a $\lambda$-homomorphism.
\item[(b)] If $K^0$ is generated by line bundles, then $\phi$ is surjective,
and $K^0$ is a $\lambda$-ring.
\item[(c)] If both $\bbZ\to K^0$ and $c_1:A\to K^0$ are injective, then
$\phi$ is injective. If in addition $K^0$ is generated by line bundles,
then $\phi:\tilde A\to K^0$ is an isomorphism.
\item[(d)] All of the above applies if we put $K^0:=K^0(X)$ instead
of $\hat K^0(X)$.
\end{itemize}
\end{LemmaD}
\begin{Proof}
(a) Map $\phi$ is an abelian group homomorphism by construction; since the 
multiplication and $\lambda_t$ are $\bbZ$-(bi)linear, it is enough to check
compatibility of $\phi$ with these maps on any system of generators 
of $\tilde A$ as an abelian group. Elements $1+x\in\tilde A$
where $x=\cl\sL$ constitute such a system of generators, we have
$\phi(1+x)=[\sL]$, $(1+x)(1+x')=1+(x+x')$ and $\phi(1+x+x')=[\sL\otimes\sL']=
[\sL][\sL']=\phi(1+x)\phi(1+x')$ whence the compatibility with multiplication.
As to the $\lambda$-operations, we notice that $s^n\sL=\sL^{\otimes n}$ for
any $n\geq0$, whence $s_t([\sL])=\sum_{n\geq0}[\sL]^nt^n$ and
$\lambda_t([\sL])=s_{-t}([\sL])^{-1}=1+[\sL]t$; since $\lambda_t(1+x)=
1+(1+x)t$ in $1+\tilde A[[t]]^+$, we obtain compatibility with 
$\lambda$-operations as well.

(b) Since $\phi(\tilde A)$ contains all $[\sL]$, the first statement is
obvious. Now $\tilde A$ is a $\lambda$-ring, $K^0$ is a pre-$\lambda$-ring
and $\phi:\tilde A\to K^0$ is a surjective $\lambda$-homomorphism,
hence all necessary relations between $\lambda$-operations are fulfilled
in $K^0$ as well, i.e.\ $K^0$ is a $\lambda$-ring.

(c) Let $\ga\subset\tilde A$ be the kernel of $\phi$. Suppose that
$\ga\neq 0$, and let $0\neq(n,x)\in\ga$. Clearly $n\neq0$ since
$\phi(0,x)=c_1(x)\neq0$ for $x\neq0$. If $n<0$, then $-(n,x)=(-n,-x)$
also lies in~$\ga$, so we can assume $n>0$. Since $\phi$ is a
$\lambda$-homomorphism, $\ga$ is a $\lambda$-ideal, hence
$\lambda^n(n,x)\in\ga$. Clearly $\lambda^n(n,x)=(1,y)$ for some $y\in A$.
Now $(1,0)=(1,-y)(1,y)$ also lies in ideal $\ga$, hence $\phi(1)=\phi(1,0)=0$.
This contradicts injectivity of $\bbZ\to K^0$.

(d) Obvious since we used only that $K^0$ is pre-$\lambda$-ring and the
formulas for multiplication and symmetric powers of images of line
bundles in $K^0$, valid for $K^0(X)$ as well.
\end{Proof}

\begin{ThD}\label{th:K0.compZ}
Let $A=\log\bbQ_+^*$, denote by $\tilde A=\bbZ\times\log\bbQ_+^*$ the
$\lambda$-ring constructed in~\ptref{sp:l.ring.from.abgr}, and define
$\phi:\tilde A\to \hat K^0(\CompZ)$ by $n+\log\lambda\mapsto
n-1+[\sO(\log\lambda)]$. Then $\phi$ is an isomorphism of
$\lambda$-rings, and in particular $\hat K^0(\CompZ)$ and
$K^0(\CompZ)$ are $\lambda$-rings.
\end{ThD}
\begin{Proof}
(a) Recall that $\Pic(\CompZ)=\log\bbQ_+^*=A$, and the isomorphism
$A\to\Pic(\CompZ)$ is given by $\log\lambda\to[\sO(\log\lambda)]$,
cf.~\ptref{sp:pic.compz}. Furthermore, the map $c_1:A\to K^0:=
\hat K^0(\CompZ)$ given by $\log\lambda\mapsto[\sO(\log\lambda)]-1$
is a homomorphism by~\ptref{l:rel.lb.K0.compz}, so we are in position
to apply~\ptref{l:pic.to.K0}. Thus $\phi:\tilde A\to K^0$ is a
$\lambda$-homomorphism. It is surjective by \ptref{l:pic.to.K0},(b)
and \ptref{sp:dec.vb.lb}. Since $\CompZ$ is quasicompact,
$\hat K^0(\CompZ)\to K^0(\CompZ)$ is also surjective; this already implies
that pre-$\lambda$-rings $\hat K^0(\CompZ)$ and $K^0(\CompZ)$ are 
actually $\lambda$-rings.

(b) Now it remains to show injectivity of~$\phi$. According to
\ptref{l:pic.to.K0},(c), it would suffice to show injectivity of
$\bbZ\to K^0$ and $c_1:A\to K^0$. The injectivity of the first map
follows from existence of map $\Spec\bbQ\to\CompZ$, which induces a
map $K^0\to\hat K^0(\bbQ)=\bbZ$. In order to show injectivity of $c_1$
we construct a degree map $\deg:\hat K^0(\CompZ)\to\log\bbQ^*_+$, such that
$\deg\circ c_1=\id$. Namely, for any vector bundle $\sE$ 
(of some rank~$n$) over~$\CompZ$,
described by a matrix~$A$, we put $\deg\sE:=\log|\det A|$. This number
is well-defined since matrices from $GL_n(\bbZ)$ and $\Oct_n=GL_n(\Zninfty)$
have determinants $\pm1$. Next, for any cofibration sequence $\sE'\to\sE\to
\sE''$ we have $A=\smallmatrix{A'}{}{*}{A''}$ for a suitable choice of
bases by~\ptref{sp:matr.perfcof.vb}, hence $\log|\det A|=\log|\det A'|+
\log|\det A''|$, i.e.\ $\deg\sE=\deg\sE'+\deg\sE''$. In other words,
$\sE\mapsto\deg\sE$ is an additive function on vector bundles, hence
it induces an abelian group homomorphism $\deg:\hat K^0(\CompZ)\to
\log\bbQ^*_+$, $\hat K^0(\CompZ)$ being the abelian group generated by
elements $[\sE]$ and relations $[\sE]=[\sE']+[\sE'']$, and
obviously $\deg\sO(\log\lambda)=\log\lambda$, i.e.\ $\deg\circ c_1=\id$
as claimed, q.e.d.
\end{Proof}

\nxsubpoint ($K^0(\CompZ)$.)
We have just shown that $\hat K^0(\CompZ)=\bbZ\oplus\log\bbQ^*_+$,
and that $K^0(\CompZ)$ is a quotient of this $\lambda$-ring. We would
like to prove that $K^0(\CompZ)$ is also isomorphic to $\bbZ\oplus
\log\bbQ^*_+$.
Since $\bbZ\to K^0(\CompZ)$ is injective for the same reason as before,
all we have to do for this is to construct a degree map
$\deg:K^0(\CompZ)\to\log\bbQ^*_+$. It is quite easy to define
$\deg X$ for any perfect simplicial $\sO_{\CompZ}$-module~$X$:
we might for example consider the function $u_X:\bbN_0\to\log\bbQ^*_+$
given by $u_X(n)=\deg X_n$, all components $X_n$ being vector bundles
by~\ptref{sp:const.perf}. Then we can show that $u_X$ is a polynomial
function of degree $\leq\dim X$ by the same reasoning as 
in~\ptref{prop:HP.perf.cof}, and put $\deg X:=u_X(-1)$.
The map $X\mapsto\deg X$ thus defined is easily seen to satisfy relations
0)--2) of~\ptref{def:K0.perf.cof}, but relation 3), which asserts that
$\deg X=\deg Y$ whenever $\gamma X$ is isomorphic to $\gamma Y$ in the
derived category, seems much harder to prove. If we would have used
Waldhausen's definition of $K_0$ instead, we would have to show that
$\deg X=\deg Y$ whenever there is a weak equivalence $X\to Y$,
a weaker but still complicated statement. Therefore, we leave this
topic for the time being, even if we think equality $K^0(\CompZ)=
\bbZ\oplus\log\bbQ^*_+$ to be highly plausible.

\nxpointtoc{Chow rings, Chern classes and intersection theory}
Our next step is to recover the Chow ring $CH(X)$ 
of a generalized ringed space or topos~$X$ (e.g.\ a generalized scheme) 
from $K^0(X)$ or $\hat K^0(X)$, considered as an augmented $\lambda$-ring,
and obtain a theory of Chern classes $c_i:K^0(X)\to CH^i(X)$ satisfying 
all classical relations. We do this by an application of the classical
procedure due to Grothendieck: define $\gamma$-operations, use them to
define the $\gamma$-filtration on~$K^0$, and finally define $CH(X)$ as the
associated graded $\gr_\gamma K^0(X)$. Chern classes can be defined then
by $c_i(\xi):=\cl\gamma^i(\xi)\in CH^i(X)$. The intersection theory 
thus constructed can be shown to be universal among all intersection theories
admitting Chern classes (i.e.\ ``orientable''), at least after tensoring 
with~$\bbQ$, and, when applied to smooth algebraic varieties, 
coincides (again after tensoring with~$\bbQ$) with any other intersection
theory with rational coefficients, admitting Chern classes and such that
the Riemann--Roch formula is valid. These results of Grothendieck show that
$CH(X,\bbQ)=\gr_\gamma K^0(X)_\bbQ$ is in all respects a reasonable definition
of the Chow ring of~$X$ with rational coefficients.

\zerosubpoint\nxsubpoint
We are going to recall some of the constructions involved here; we refer
to SGA~6~V for the proofs. However, there are some complications
in our situation that we need to deal with:
\begin{itemize}
\item We have shown that $K^0(X)$ and $\hat K^0(X)$ are pre-$\lambda$-rings,
but in order to apply Grothendieck's construction we need a $\lambda$-ring.
We deal with this problem by replacing $K^0(X)$ with its largest quotient
$K^0(X)_\lambda$ which is a $\lambda$-ring. In all situations where we
are able to prove that $K^0(X)$ is already a $\lambda$-ring, $K^0(X)_\lambda$
coincides with $K^0(X)$, hence our result will coincide with the
classical one.
\item We need in fact an augmented $\lambda$-ring, i.e.\ we need an 
augmentation $\epsilon:K^0(X)\to H^0(X,\bbZ)$ given by the rank of 
vector bundles. For this we need to know that the rank of a vector bundle 
over~$X$ is a well-defined locally constant integer-valued function.
Sometimes we can achieve this by considering the retract-free version of $K^0$
(then all vector bundles are automatically locally free), but this is
not sufficient (we might have $L_\sO(n)\cong L_\sO(m)$ for $n\neq m$).
We'll discuss what can done about it.
\end{itemize}

\nxsubpoint\label{sp:lambda.str.formser} 
($\lambda$-structure on $\hat G^t(K)$.)
Let $K$ be any commutative ring with unity.
Recall that one defines a pre-$\lambda$-ring structure on
the set $\hat G^t(K)=1+K[[t]]^+$ of formal series with free term equal to one,
considered as an abelian group under multiplication, by imposing the
following requirements:
\begin{itemize}
\item $\hat G^t(K)$ depends functorially on~$K$.
\item $\{1+Xt\}\circ\{1+Yt\}=\{1+XYt\}$ in $\hat G^t(\bbZ[X,Y])$, where
$\circ$ denotes the new multiplication on $\hat G^t(K)$, as opposed to the
usual multiplication of formal series, which corresponds now to the
addition of ring $\hat G^t(K)$.
\item $\lambda^k\{1+Xt\}=\{1\}=0$ in $\hat G^t(\bbZ[X])$ for any $k\geq2$.
\item Operations $\circ$ and $\lambda^k$ are ``continuous'' with respect to
the natural filtration on $\hat G^t(K)$. In other words, 
the first $n$ coefficients of $f\circ g$ and $\lambda^k(f)$ are 
completely determined by the first $N$ coefficients of $f$ and~$g$
(resp.~$f$) for some $N=N(n)>0$.
\end{itemize}
One can show that these conditions uniquely determine operations $\circ$ and
$\lambda^k$ on $\hat G^t$, and that $\hat G^t(K)$ is a pre-$\lambda$-ring
with respect to these operations (cf.\ SGA~6~V). Furthermore,
these operations can be expressed by means of certain universal polynomials
with integer coefficients, which can be found by means of 
elementary symmetric polynomial calculus:
\begin{multline}
\{1+x_1t+x_2t^2+\cdots\}\circ\{1+y_1t+y_2t^2+\cdots\}=
1+P_1(x_1;y_1)t+\\
+P_2(x_1,x_2;y_1,y_2)t^2+\cdots+P_n(x_1,\ldots,x_n;y_1,\ldots,
y_n)t^n+\cdots
\end{multline}
Each $P_n=P_n(x_1,\ldots,x_n;y_1,\ldots,y_n)$ is isobaric of weight~$n$
both in $(x_i)$ and $(y_i)$, where the weight of $x_i$ and $y_i$ is set
equal to~$i$. For example, $P_1(x_1;y_1)=x_1y_1$, and
$P_2(x_1,x_2;y_1,y_2)=x_2y_1^2+y_2x_1^2-2x_2y_2$ (cf.\ {\em loc.cit.}).
The $\lambda$-operations are also given by some universal polynomials:
\begin{equation}
\lambda^k\{1+x_1t+x_2t^2+\cdots\}=1+Q_{k,1}(x_1,\ldots,x_k)t+\cdots+
Q_{k,n}(x_1,\ldots,x_{kn})t^n+\cdots
\end{equation}
These universal polynomials $Q_{i,j}(x_1,x_2,\ldots,x_{ij})$ are 
isobaric of weight~$ij$. For example, $Q_{k,1}=x_k$, and
$Q_{2,2}(x_1,x_2,x_3,x_4)=x_1x_3-x_4$ (cf.\ SGA~6~V 2.3).

\begin{DefD} {\rm ($\lambda$-rings.)}
A pre-$\lambda$-ring $K$ is a {\em $\lambda$-ring\/} if
the abelian group homomorphism $\lambda_t:K\to\hat G^t(K)$ is
a $\lambda$-homomorphism with respect to the $\lambda$-structure on 
$\hat G^t(K)$ just discussed.
\end{DefD}
In other words, $\lambda^n(xy)$ and $\lambda^j(\lambda^i(x))$ in~$K$
must be given by the same universal polynomials:
\begin{align}\label{eq:lambda.prod.rels}
\lambda^n(xy)&=P_n(x,\lambda^2x,\ldots,\lambda^nx;y,\lambda^2y,\ldots,
\lambda^ny)\\\label{eq:lambda.lambda.rels}
\lambda^j(\lambda^ix)&=Q_{i,j}(x,\lambda^2x,\ldots,\lambda^{ij}x)
\end{align}

Almost all pre-$\lambda$-rings considered in practice are $\lambda$-rings.
For example, $\hat G^t(K)$ is actually a $\lambda$-ring for any commutative
ring~$K$. Another example: $K^0(X)$ is a $\lambda$-ring for any classical
ringed space or topos~$X$. Unfortunately, we cannot generalize the proof
of this classical statement to our case, at least until we compute
$K^0$ of projective bundles.

\nxsubpoint\label{sp:assoc.lambda.ring} 
($\lambda$-ring associated to a pre-$\lambda$-ring.)
Given any pre-$\lambda$-ring~$K$, we denote by $K_\lambda$ the quotient
of $K$ modulo the $\lambda$-ideal generated by 
relations \eqref{eq:lambda.prod.rels} and~\eqref{eq:lambda.lambda.rels}. 
Clearly, $K_\lambda$ is the largest $\lambda$-ring 
quotient of $K$; if $K$ is already a $\lambda$-ring,
then $K_\lambda=K$.

\nxsubpoint (Binomial rings and their $\lambda$-structure.)
We say that a commutative ring~$K$ is {\em binomial\/}
(cf.\ SGA~6 V~2.7) if it is torsion-free as an abelian group, and
if for any $x\in K$ and any $n\geq 1$ the element $x(x-1)\cdots(x-n+1)$
is divisible by $n!$ in~$K$. The quotient is then denoted by $\binom xn$.

We endow a binomial ring by its canonical $\lambda$-ring structure, given 
by $\lambda^nx:=\binom xn$. In other words,
$\lambda_t(x)=(1+t)^x$, hence $s_t(x)=(1-t)^{-x}$ and $s^n(x)=(-1)^n\binom{-x}n
=\binom{x+n-1}n$. For example, $\bbZ$ and 
$H^0(X,\bbZ)$ (for any topological space or topos~$X$) are binomial rings.

\nxsubpoint (Adams operations on a pre-$\lambda$-ring.)
Given any pre-$\lambda$-ring~$A$, we define the {\em Adams operations\/}
$\Psi^n:A\to A$, $n\geq1$, by means of the following identity:
\begin{equation}\label{eq:def.adams.op}
\sum_{n=1}^\infty(-1)^{n-1}\Psi^n(x)t^{n-1}=\frac{d}{dt}\lambda_t(x)
\big/\lambda_t(x)=\frac{d}{dt}\log\lambda_t(x)
\end{equation}
We can rewrite the above relation coefficientwise:
\begin{equation}
n\lambda^n(x)=\sum_{k=1}^n(-1)^{k-1}\Psi^k(x)\lambda^{n-k}(x)
\end{equation}
These relations allow one to express easily the Adams operations in terms
of exterior power operations, or conversely. In any case,
$\Psi^n(x)$ is given by a universal polynomial with integer coefficients
in $x$, $\lambda^2(x)$, \dots, $\lambda^n(x)$, isobaric of weight~$n$.

Symmetric and exterior power operations are related by
$s_t(x)\lambda_{-t}(x)=1$, so we have
\begin{equation}
\sum_{n=1}^\infty\Psi^n(x)t^{n-1}=-\frac{d}{dt}s_t(x)\big/s_t(x)=
-\frac{d}{dt}\log s_t(x)
\end{equation}
This equality enables one to express symmetric operations in terms of
Adams operations and conversely.

The usual properties of logarithmic derivatives together with
$\lambda_t(x+y)=\lambda_t(x)\lambda_t(y)$ imply
\begin{equation}
\Psi^n(x+y)=\Psi^n(x)+\Psi^n(y),\quad\Psi^1(x)=x
\end{equation}
Conversely, if $A\supset\bbQ$, then any family of abelian group
homomorphisms $\{\Psi^n:A\to A\}_{n\geq1}$, such that $\Psi^1=\Id_A$, 
defines a pre-$\lambda$-ring
structure on~$A$ by means of~\eqref{eq:def.adams.op}. A $\lambda$-homomorphism
$f:B\to A$ can be described then as a ring homomorphism commuting with all
Adams operations~$\Psi^n$.

\nxsubpoint (Adams operations on a $\lambda$-ring.)
One easily shows, starting from
\begin{equation}
\Psi^n\{1+Xt\}=1+X^nt\quad\text{in $\hat G^t(\bbZ[X])$}
\end{equation}
that whenever~$A$ is a $\lambda$-ring, its Adams operations satisfy
\begin{equation}
\Psi^n(1)=1,\quad\Psi^n(xy)=\Psi^n(x)\Psi^n(y),\quad
\Psi^i(\Psi^j(x))=\Psi^{ij}(x)
\end{equation}
For example, $\Psi^n=\id_K$ on any binomial ring~$K$.

Conversely, if a pre-$\lambda$-ring $A$ is $\bbZ$-torsion-free and 
if its Adams operations satisfy the above relations, $A$ can be shown
to be a $\lambda$-ring. Therefore, a $\lambda$-ring structure on
a commutative ring $A\supset\bbQ$ can be described as a family of
ring endomorphisms $\{\Psi^n:A\to A\}_{n\geq1}$, such that
$\Psi^n\circ\Psi^m=\Psi^{nm}$ and $\Psi^1=\id$. This is probably the
simplest way to define a $\bbQ$-$\lambda$-algebra.

\nxsubpoint ($\gamma$-operations.)
Given a pre-$\lambda$-ring~$A$, we define $\gamma$-operations 
$\gamma^n:A\to A$, $n\geq0$, by means of the following generating series:
\begin{equation}
\gamma_t(x)=\sum_{n=0}^\infty\gamma^n(x)t^n:=\lambda_{t/(1-t)}(x)
\end{equation}
Since $\lambda_t(x)=\gamma_{t/(t+1)}(x)$, the $\lambda$-operations completely
determine the $\gamma$-opera\-tions, and conversely. Formal manipulations
with power series, using $\lambda_t(x+n-1)=\lambda_t(x+n)/(1+t)$, yield
\begin{equation}
\gamma^n(x)=\lambda^n(x+n-1)
\end{equation}

\nxsubpoint (Augmented $K$-$\lambda$-algebras.)
Let $K$ be any binomial ring. An {\em augmented $K$-$\lambda$-algebra\/}
is by definition a $K$-$\lambda$-algebra $A$ together with an
{\em augmentation\/} $\epsilon:A\to K$, supposed to be a 
$K$-$\lambda$-homomorphism. If $A$ is $\bbZ$-torsion-free, the latter
condition can be expressed with the aid of Adams operations;
taking into account that all $\Psi_K^n=\id_K$, we get
$\epsilon\circ\Psi^n=\epsilon$ for any $n\geq1$. Let's identify $K$ 
with a subring of~$A$. One can put $\Psi^0:=\epsilon$, considered
as a ring homomorphism $A\to A$. Then $\Psi^n\circ\Psi^0=\Psi^0=\Psi^0\circ
\Psi^n$ for all $n\geq0$, i.e.\ {\em an augmentation can be considered 
as an extension $\Psi^0$ of Adams operations $\Psi^n$ to $n=0$.}
Conversely, $K$ can be recovered from~$A$ and $\Psi^0$ as the image of
$\Psi^0$, or as the set of all $x$, such that $x=\Psi^0(x)$,
and relations $\Psi^n\circ\Psi^0=\Psi^n$ imply that $K$
is a binomial ring, at least if $A$ is torsion-free.

\nxsubpoint (Augmented $K$-pre-$\lambda$-algebras.)
If $A$ is just an augmented $K$-pre-$\lambda$-algebra
over a binomial ring~$K$, $\lambda$-homomorphism
$\epsilon:A\to K$ factorizes through the largest $\lambda$-ring quotient
$A_\lambda$ of~$A$, discussed in~\ptref{sp:assoc.lambda.ring}, since
$K$ is a $\lambda$-ring. Therefore, $A_\lambda$ becomes an
augmented $K$-$\lambda$-algebra, and we can apply to it all the
constructions that follow.

\nxsubpoint\label{sp:cond.augm.K0} (Important example: augmented~$K^0$.)
Let $X=(X,\sO_X)$ be a generalized ringed space or topos.
Put $A:=\hat K^0(X)$ or $K^0(X)$ and consider binomial ring
$K:=H^0(X,\bbZ)$. We have natural $\lambda$-homomorphisms 
$K\to\hat K^0(X)\to K^0(X)$, given by $n\mapsto[L_\sO(n)]$, i.e.\ 
$A$ is a pre-$\lambda$-algebra.
Let's assume the following:
\begin{itemize}
\item Any vector bundle $\sE$ over~$X$ is locally free.
Notice that this holds over any classical locally ringed topos.
In general we might either prove it or use the retract-free version of the
theory, where this statement is automatic.
\item The rank of a vector bundle $\sE$ over~$X$ is locally well-defined,
i.e.\ for any non-empty open subset $U\subset X$ (object of~$X$ in the
topos case) and any $n\neq m$ the free $\sO_X|U$-modules 
$L_{\sO_X|U}(n)$ and $L_{\sO_X|U}(m)$ are not isomorphic. Notice that this
condition is automatic if $\sO_X$ is alternating (cf.~\ptref{def:altop})
and if $\sO_X(U)$ is non-subtrivial for all $U\neq\emptyset$, since then
the rank of $\sE=L_{\sO_X|U}(n)$ can be recovered as the largest integer~$k$,
such that $\bigwedge^k\sE\neq0$.
\end{itemize}
Under these conditions the rank of a vector bundle is well-defined and
additive, hence it defines an additive map 
$\epsilon:\hat K^0(X)\to H^0(X,\bbZ)=K$.
In some cases it can be extended to a map $K^0(X)\to K$, for example if
$X$ is quasi-compact and has a dense open subset $U\subset X$ isomorphic 
to a classical scheme, such that the restriction map $H^0(X,\bbZ)\to
H^0(U,\bbZ)$ is an isomorphism, 
a condition usually fulfilled by all our ``arithmetic models''
$\sX/\CompZ$ of algebraic varieties~$X/\bbQ$.

Therefore, we obtain an augmented $K$-pre-$\lambda$-algebra~$A$. If it is
not a $\lambda$-ring, we can consider augmented $K$-$\lambda$-algebra 
$A_\lambda$ instead.

\nxsubpoint ($\gamma$-filtration.)
Let $A$ be any augmented $K$-$\lambda$-algebra over a binomial ring~$K$.
We define a decreasing filtration $F_\gamma^i$ on $A$, called the
{\em $\gamma$-filtration}, as follows: $F_\gamma^nA$ is the ideal in~$A$
generated by elements
\begin{equation}
\gamma^{k_1}(x_1)\gamma^{k_2}(x_2)\cdots\gamma^{k_s}(x_s),\quad
\text{where $\sum_i k_i\geq n$, $x_i\in A$, $\epsilon(x_i)=0$}
\end{equation}
In particular, $F_\gamma^0A=A$ and $F_\gamma^1A=\Ker\epsilon$.

Since $F_\gamma$ depends functorially on~$A$, we've obtained a functorial
construction $A\mapsto\gr_\gamma A$ of a graded commutative ring from
an augmented $K$-$\lambda$-algebra~$A$. We'll see in a moment that
under some additional restrictions $K$ and $A$ can be recovered from
graded commutative ring $\gr_\gamma A$, thus establishing an equivalence
of categories.

\nxsubpoint (Chow ring.)
In particular, if $X=(X,\sO_X)$ satisfies the conditions 
of~\ptref{sp:cond.augm.K0}
we can apply the above construction to $K=H^0(\bbZ,X)$ and 
$A=\hat K^0(X)_\lambda$ or $A=K^0(X)_\lambda$, thus obtaining
the {\em Chow ring of~$X$:}
\begin{equation}
CH(X):=\gr_\gamma\hat K^0(X)_\lambda,\quad
CH(X,\bbQ):=CH(X)_\bbQ=CH(X)\otimes_\bbZ\bbQ
\end{equation}
Usually we cannot expect to obtain a reasonable intersection theory
with integer coefficients unless $X$ is an algebraic variety over a field,
so we prefer to consider only $CH(X)_\bbQ$.

By construction $CH(X)$ and $CH(X)_\bbQ$ are contravariant in~$X$,
i.e.\ if $f:Y\to X$ is a morphism of generalized ringed spaces or topoi
satisfying the conditions of~\ptref{sp:cond.augm.K0}, we obtain
pullback maps $f^*:CH(X)\to CH(Y)$. Under some very special conditions 
for~$f$ (e.g.\ being a ``regular immersion'' or ``projective locally 
complete intersection'') 
one can hope to construct ``Gysin maps'' $f_!:CH(Y)\to CH(X)$
in the opposite direction.

\nxsubpoint (Example: Chow ring of $\CompZ$.)
Let $A$ be any abelian group (say, $A=\log\bbQ^*_+$). 
Consider the $\lambda$-ring $\tilde A=\bbZ\times A$ 
of~\ptref{sp:l.ring.from.abgr}, given by $(n,x)\cdot(m,y)=
(nm,ny+mx)$, $\lambda^k(n,x)=\bigl(\binom nk,\binom{n-1}{k-1}x\bigr)$.
Clearly $\tilde A$ is an augmented $\bbZ$-$\lambda$-algebra 
with augmentation $\epsilon:(n,x)\mapsto n$, so we can apply the above
construction to it. 

First of all, we know that $F_\gamma^0\tilde A=\tilde A$ 
and $F_\gamma^1\tilde A=\Ker\epsilon=0\times A\cong A$.
Since $\gamma^k(0,x)=\lambda^k(k-1,x)=
(0,\binom{k-2}{k-1}x)=0$ for $k\geq 2$ and any $x\in A$, we see that
$F_\gamma^2\tilde A$ is the ideal generated by products $\gamma^1(x)
\gamma^1(y)$ for $x$, $y\in\Ker\epsilon=A$; since $\gamma^1=\id$ and
$A^2=0$ in $\tilde A$, we obtain $F_\gamma^2\tilde A=0$. This implies
$\gr_\gamma\tilde A=\bbZ\oplus A$, which is isomorphic as a ring
to~$\tilde A$ itself.

Now let's apply this to $\hat K^0(\CompZ)$, equal to $\tilde A$ for
$A=\log\bbQ^*_+$ by~\ptref{th:K0.compZ}. We obtain
\begin{equation}
CH(\CompZ)=CH^0(\CompZ)\oplus CH^1(\CompZ)=\bbZ\oplus\log\bbQ_*^+=
\bbZ\oplus\Pic(\CompZ)
\end{equation}
We can compute $CH^0(\CompZ)_\bbQ=\bbQ$ and $CH^1(\CompZ)_\bbQ=
\bbQ\otimes_\bbZ\log\bbQ^*_+$ as well, but in this special situation we
obtain the ``correct'' result even without tensoring with~$\bbQ$, just
because we have $F_\gamma^2A=0$ in this situation, i.e.\ 
``the intersection theory of $\CompZ$ is one-dimensional''.

\nxsubpoint (Chern ring and formal computations with Chern classes.)
Let $C=\oplus_{n\geq0}C_n$ be any graded commutative $K$-algebra
(e.g.\ a Chow ring $CH(X)$) over a binomial ring~$K$ (in most cases
we'll have $K=C_0$). 
Consider the completion
$\hat C:=\prod_{n\geq0}C_n=\{c_0+c_1+\cdots+c_n+\cdots\,|\,c_i\in C_i\}$
and multiplicative subgroup $1+\hat C^+=1\times\prod_{n\geq1}C_n=
\{1+c_1+c_2+\cdots+c_n+\cdots\}\subset\hat C$. Define the
{\em Chern ring\/} $\Ch(C)$ by
\begin{equation}
\Ch(C):=K\times(1+\hat C^+)
\end{equation}
The operations on $\Ch(C)$ are defined as follows:
\begin{itemize}
\item Addition is defined componentwise: $(a,f)+(b,g)=(a+b,fg)$.
\item Action of $K$ on $\Ch(C)$ (i.e.\ restriction of multiplication
to $K\times\Ch(C)$) is given by
\begin{equation}
(a,1)*(b,1+f)=(ab,(1+f)^a),\quad\text{$a$, $b\in K$, $1+f\in 1+\hat C^+$,}
\end{equation}
where the power $(1+f)^a$ is defined by the Newton binomial formula
\begin{equation}
(1+f)^a=\sum_{n=0}^\infty\binom anf^n
\end{equation}
\item Multiplication satisfies
\begin{equation}
(1,1+x)*(1,1+y)=(1,1+x+y)\quad\text{for any $x$, $y\in C_1$.}
\end{equation}
\item $\lambda$-operations satisfy
\begin{align}
\lambda^n(a,1)&=\bigl(\binom an,1\bigr)\quad\text{for any $a\in K$, $n\geq0$,}\\
\lambda^n(1,1+x)&=0=(0,1)\quad\text{for any $x\in C_1$, $n\geq2$.}
\end{align}
\item $\Ch(C)$ is a commutative pre-$\lambda$-ring depending functorially
on graded $K$-algebra~$C$. Its operations $*$ and $\lambda^n$ are
continuous in the same sense as in~\ptref{sp:lambda.str.formser}.
\end{itemize}
One can check that the above conditions determine a unique $\lambda$-ring
structure on $\Ch(C)$, and that the coefficients of products and
exterior powers can be expressed with the aid of some universal polynomials
with integer coefficients (cf.\ SGA~6 V.6).

If we think of an element
$\tilde c(\sE)=(n,1+c_1+c_2+\cdots)$ of~$\Ch(C)$ as 
the rank~$n$ and the collection
$c_i$ of Chern classes of a vector bundle~$\sE$, then these universal
polynomials are nothing else than the classical rules for the
Chern classes of tensor product of two vector bundles or the exterior power
of a vector bundle, given in~\cite{Gr1}. 
In other words, if $C=CH(X)$, we can expect
the completed Chern class $\tilde c:\hat K^0(X)\to\Ch(C)$ to be a 
$\lambda$-homomorphism.

\nxsubpoint (Formal definition of Chern classes.)
Let $A$ be again an augmented $K$-$\lambda$-algebra over a binomial ring~$K$.
(Our principal example is still $K=H^0(X,\bbZ)$, $A=\hat K^0(X)_\lambda$
for a suitable generalized ringed space~$X$.) Put $C:=\gr_\gamma A$
and define the {\em Chern classes $c_i:A\to C_i$, $i\geq1$,} by
\begin{equation}
c_i(x):=\cl_{\gr_\gamma^iA}\gamma^i(x-\epsilon(x))\in C_i,\quad
\text{for any $x\in A$.}
\end{equation}
We define the {\em completed total Chern class $\tilde c:A\to\Ch(C)$} by
\begin{equation}
\tilde c(x):=\bigl(\epsilon(x),1+c_1(x)+c_2(x)+\cdots+c_n(x)+\cdots\bigr)
\end{equation}
Then {\em $\tilde c:A\to\Ch(C)$ is a $\lambda$-homomorphism of augmented 
$K$-$\lambda$-algebras} (cf.\ SGA~6 V.6.8), i.e.\ the formal Chern classes
thus defined satisfy all classical relations of~\cite{Gr1}.

In particular, applying this to $K=H^0(X,\bbZ)$, $A=\hat K^0(X)_\lambda$
or $K^0(X)_\lambda$
as above, we obtain reasonable Chern classes $c_i:A=\hat K^0(X)\to CH^i(X)=
\gr_\gamma^i\hat K^0(X)_\lambda$ and the completed total Chern class
$\tilde c:\hat K^0(X)\to\Ch(CH(X))$. When $\sE$ is a vector bundle over~$X$,
we write $c_i(\sE)$ and $\tilde c(\sE)$ instead of $c_i([\sE])$ and
$\tilde c([\sE])$.

\nxsubpoint (Example: Chern classes over $\CompZ$.)
Applying the above construction to $X=\CompZ$, we obtain exactly the
answer we expect: $c_i(\sE)=0$ for $i\geq2$, and $c_1(\sE)=\deg\sE$
for any vector bundle $\sE$ over~$\CompZ$, where the degree of a 
vector bundle is defined as in the proof of~\ptref{th:K0.compZ}.

\nxsubpoint (Additivity of $c_1$ on line bundles.)
One has classical formula for $c_1$ of a product:
\begin{equation}
c_1(xy)=\epsilon(x)c_1(y)+\epsilon(y)c_1(x)
\end{equation}
This implies that the restriction of $c_1$ onto the multiplicative subgroup
$\{x:\epsilon(x)=1\}$ is a homomorphism of abelian groups. Applying this
to our favourite example $A=\hat K^0(X)_\lambda$, 
we see that {\em the first Chern class induces a homomorphism $c_1:
\Pic(X)\to CH^1(X)$.}

\nxsubpoint (Chern character.)
Let $C$ be a graded commutative algebra over a binomial ring~$K$ as before.
We define a functorial homomorphism $\ch:\Ch(C)\to\hat C_\bbQ$
(or to $K\oplus\hat C^+_\bbQ$), where
$\hat C_\bbQ:=\widehat{C\otimes_\bbZ\bbQ}$,
by requiring $\ch$ to be ``continuous'' and imposing following relations:
\begin{align}
\ch(a,1)&=a,\quad a\in K,\\
\ch(1,1+x)&=\exp(x)=\sum_{n\geq0}\frac{x^n}{n!},\quad x\in C_1\\
\ch(u+v)&=\ch(u)+\ch(v),\quad u,v\in\Ch(C)
\end{align}
We say that $\ch$ is the {\em Chern character;} it satisfies
$\ch(u*v)=\ch(u)\ch(v)$ for any $u$, $v\in\Ch(C)$, 
and it can be also defined by
\begin{equation}
\ch(a,\phi)=a+\eta(\log\phi)
\end{equation}
where $\log:1+\hat C^+\to\hat C^+_\bbQ$ is defined by the usual series,
and $\eta$ is the additive endomorphism of $\hat C^+_\bbQ$ multiplying
the degree~$k$ component~$a_k$ of an element $a=a_1+a_2+\cdots$ of 
$\hat C^+_\bbQ$ by $(-1)^{k-1}/(k-1)!$ (cf.\ SGA~6 V.6.3).

One can check that $\ch$ is compatible with the multiplication on $\Ch(C)$,
i.e.\ $\ch$ is a ring homomorphism (cf.~{\em loc.cit.}). Furthermore,
if $K\supset\bbQ$ (e.g.\ if we tensorise everything with~$\bbQ$ from the 
very beginning), and if $C_0=K$, then $\ch$ is a $K$-algebra {\em iso}morphism
(since both $\eta$ and $\log$ are isomorphisms in this case).

Applying this construction to $C=\gr_\gamma A$ obtained from an
augmented $K$-$\lambda$-algebra~$A$, we obtain a $K$-algebra homomorphism
$\ch:A\stackrel{\tilde c}\to\Ch(C)\stackrel{\ch}\to\hat C_\bbQ$, e.g.\ 
$\ch:\hat K^0(X)\to\widehat{CH(X)}_\bbQ$, called {\em Chern character}. Thus
$\ch(x+y)=\ch(x)+\ch(y)$ and $\ch(xy)=\ch(x)\ch(y)$ as in the classical case.
Of course, if $A=\hat K^0(X)$, we write $\ch(\sE)$ or $\ch\sE$ instead of
$\ch([\sE])$. In this way we obtain a reasonable Chern character on
$K^0(X)$ or $\hat K^0(X)$ in a completely formal fashion.

\nxsubpoint (Adams operations on $\Ch(C)$.)
One can compute the Adams operations $\Psi^n$ on $\Ch(C)$.
Since $K$ is binomial, we have $\Psi^n(a,1)=(a,1)$ for any $n\geq1$ and
any $a\in K$. Next, let's compute over $C=\bbZ[T_1,\ldots,T_n]$, $K=\bbZ$, 
graded by putting
$\deg T_i=1$. Since $\lambda_t(1,1+T_i)=1+(1,1+T_i)t$, we get
$\Psi^k(1,1+T_i)=(1,1+T_i)^{*k}=(1,1+kT_i)$ by definition of multiplication
on $\Ch(C)$ and by general rule $\lambda^t(x)=1+xt\Rightarrow\Psi^k(x)=x^k$. 
Since $\Psi^k$ are additive, we can conclude $\Psi^k(0,1+T_i)=(0,1+kT_i)$,
hence $\Psi^k(0,\prod_i(1+T_i))=(0,\prod_i(1+kT_i))$. Taking $\gS_n$-invariants
we obtain $\Psi^k(0,1+X_1+X_2+\cdots+X_n)=(0,1+kX_1+k^2X_2+\cdots+k^nX_n)$
over $C=\bbZ[T_1,\ldots,T_n]^{\gS_n}=\bbZ[X_1,\ldots,X_n]$, graded by
$\deg X_i=i$. Using universality of such algebras $\bbZ[X_1,\ldots,X_n]$
together with ``continuity'' of $\Psi^k$ we get
\begin{equation}
\Psi^k(a,1+c_1+c_2+\cdots+c_n+\cdots)=(a,1+kc_1+k^2c_2+\cdots+k^nc_n+\cdots)
\end{equation}
inside any $\Ch(C)$. Notice that this equality extends to $k=0$ if we define
$\Psi^0:=\epsilon$. 

\nxsubpoint (Adams operations on $\hat C$ via $\ch$.)
Let's still suppose $K\supset\bbQ$, $C$ be a graded $K$-algebra with $C_0=K$,
Then $\ch:\Ch(C)\to\hat C$ is a ring isomorphism, so we can transfer the
$\lambda$-ring structure from $\Ch(C)$ to $\hat C\supset C$. Since $\hat C
\supset K\supset\bbQ$, this $\lambda$-structure is completely determined by
its Adams operations.
\begin{Propz}
Adams operations $\Psi^k$ on $\hat C$ with respect to
the $\lambda$-structure just discussed are given by
\begin{equation}\label{eq:Adams.on.Chat}
\Psi^k(c_0+c_1+c_2+\cdots+c_n+\cdots)=c_0+kc_1+k^2c_2+\cdots+k^nc_n+\cdots
\end{equation}
In other words, $\Psi^k$ are continuous on $\hat C$ and restrict to $k^n$
on $C_n$.
\end{Propz}
\begin{Proof} Both sides are ``continuous'' and additive, and $\Psi^k$ 
restricted to $C_0=K$ is trivial just because $K$ is binomial, so 
we can suppose $c_0=0$. Next, $\ch$ is an isomorphism, hence it is enough
to check the statement on ``universal elements'' $\ch(0,1+X_1+\cdots+X_n)$
of $C=\bbQ[X_1,\ldots,X_n]$, graded by $\deg X_i=i$. Let's embed
$\bbQ[X_1,\ldots,X_n]$ into $C'=\bbQ[T_1,\ldots,T_n]$, $\deg T_i=1$,
by means of elementary symmetric polynomials. Then $(0,1+X_1+\cdots+X_n)=
(0,1+T_1)+\cdots+(0,1+T_n)$, so by additivity we are reduced to check
\eqref{eq:Adams.on.Chat} on elements $\ch(0,1+t)$, $t\in C_1$,
or on elements $\ch(1,1+t)=1+\ch(0,1+t)$.
By definition $\ch(1,1+t)=\exp(t)=\sum_{n\geq0}t^n/n!$, and
$\lambda^k(1,1+t)=0$ in $\Ch(C)$ for $k\geq2$, again by definition of
$\lambda$-operations in $\Ch(C)$, whence $\Psi^k(1,1+t)=(1,1+t)^{*k}=
(1,1+kt)$ as before, and $\ch\circ\Psi^k(1,1+t)=\ch(1,1+kt)=
\exp(1+kt)=\sum_{n\geq0}k^nt^n/n!=\Psi^k\circ\ch(1,1+t)$, if
$\Psi^k$ on $\hat C$ are defined by \eqref{eq:Adams.on.Chat}. This proves
our statement.
\end{Proof}

Notice that these Adams operations $\Psi^k:\hat C\to\hat C$ 
respect $C\subset\hat C$, i.e.\ $\Psi^k(C)\subset C$, hence 
$C$ is a $K$-$\lambda$-subalgebra of $\hat C$ in a natural fashion.

\nxsubpoint (Weight decomposition with respect to Adams operations.)
Notice that the homogeneous component $C_n$ can be easily recovered from
$\lambda$-ring $C$ or $\hat C$ as the subgroup $\hat C_{(n)}$ 
of elements of weight $n$ with respect to Adams operations:
\begin{equation}
\hat C_{(n)}=\{x\in\hat C\,:\,\Psi^k(x)=k^nx\text{ for any $k\geq1$}\}.
\end{equation}
Furthermore, the above formula holds for $k\geq0$ if we put $\Psi^0:=\epsilon$.
On the other hand, if we fix any $k\geq2$, $\hat C^{(n)}$ 
is already completely determined by requirement $\Psi^k(x)=k^nx$ for this
fixed value of~$k$ since all $k^n$ are distinct. In particular,
binomial ring $K=C_0=\hat C_{(0)}$ is completely determined by the 
$\lambda$-structure of~$C$ or $\hat C$, and is the largest binomial ring
contained in this $\lambda$-ring.

The above construction is valid for any $\bbQ$-$\lambda$-algebra~$A$:
we obtain weight subgroups $A_{(n)}\subset A$, and the sum
$A_{(\cdot)}:=\sum_{n\in\bbZ}A_{(n)}\subset A$ 
will be always direct (e.g.\ because of
the non-degeneracy of Vandermond matrices), but in general it needn't
coincide with~$A$. In any case, $A_{(n)}\cdot A_{(m)}\subset A_{(n+m)}$
for any $n$, $m\in\bbZ$ just because all $\Psi^k$ are ring homomorphisms,
i.e.\ $A_{(\cdot)}$ is an augmented graded $K$-$\lambda$-subalgebra of~$A$ over
binomial ring $K=A_{(0)}$.

Since $\ch:\Ch(C)\to\hat C$ is a $\lambda$-ring isomorphism 
(we still assume $C_0=K\supset\bbQ$), 
we see that $C_n$ might be equally easily recovered from
$\Ch(C)$ by $C_n=\Ch(C)_{(n)}$.

\nxsubpoint (Classification of augmented $K$-$\lambda$-algebras with
discrete $\gamma$-filtration over $K\supset\bbQ$.)
Let fix a binomial ring~$K\supset\bbQ$. If $C$ is a graded $K$-algebra with
$C_0=K$, then the $\gamma$-filtration on $\Ch(C)$ is easily seen to
coincide with the natural one (e.g.\ using the Chern character isomorphism):
\begin{equation}
F_\gamma^n\Ch(C)=\{(0,1+c_n+c_{n+1}+\cdots),\,c_i\in C_i\}
\end{equation}
Hence $\gr_\gamma\Ch(C)\cong C$, at least as a graded $K$-module. In fact,
this is an isomorphism of $K$-algebras if we define $\gr_\gamma\Ch(C)\simto C$
with the aid of $\ch$, i.e.
\begin{equation}
\cl_{\gr_\gamma^n\Ch(C)}(0,1+c_n+c_{n+1}+\cdots)\mapsto
\frac{(-1)^{n-1}}{(n-1)!}\cdot c_n
\end{equation}
In other words, we have a functorial isomorphism $\gr_\gamma\Ch(C)\simto C$.

Conversely, if $A$ is a graded $K$-$\lambda$-algebra, the completed total
Chern class $\tilde c:A\to\Ch(C)=\Ch(\gr_\gamma A)$ is a homomorphism
of augmented $K$-$\lambda$-algebras, which can be shown to induce an
isomorphism on associated graded with respect to the $\gamma$-filtration. 
Modulo some additional verifications we obtain the following statement:

\begin{Thz} For any integer $N\geq0$ and any $K\supset\bbQ$ 
functors $A\mapsto\gr_\gamma A$ and $C\mapsto\Ch(C)$ are quasi-inverse
equivalences between the category of augmented $K$-$\lambda$-algebras~$A$, 
such that $F_\gamma^{N+1}A=0$, and the category of graded $K$-algebras~$C$,
such that $C_0=K$ and $C_n=0$ for $n>N$.
\end{Thz}
{\bf Proof.} More details can be found in SGA~6 V~6.11.

Therefore, the category of augmented $K$-$\lambda$-algebras~$A$ with discrete
$\gamma$-filtration is equivalent to the category of graded $K$-algebras~$C$,
such that $C_0=K$, $C_n=0$ for $n\gg0$, i.e.\ any such~$A$ is isomorphic to
the Chern ring $\Ch(C)$ of a suitable~$C$. Furthermore, under these conditions
the Chern character $\ch:\Ch(C)\to\hat C=C$ is an isomorphism,
i.e.\ {\em the Chern character induces an isomorphism between 
$A$ and a graded $K$-algebra~$C$, considered as a $\lambda$-ring via
Adams operations $\Psi^k$ given by $\Psi^k|_{C_n}=k^n$.}
Since the graded components $C_n$ of~$C$ come from the weight decomposition
with respect to Adams operations, and $C$ is isomorphic to~$A$ as a
$\lambda$-ring, we see that {\em $\gr_\gamma^nA=C_n=C_{(n)}$ is 
canonically isomorphic to abelian group $A_{(n)}$.} This is essentially
the construction of Soul\'e (cf.~\cite{Soule1}):
$CH^i(X)_\bbQ=K^0(X)_{\bbQ,(n)}$. Actually Soul\'e considers weight 
decomposition of higher algebraic $K$-groups as well, something we
don't discuss in this work.

If $A$ is complete with respect to its $\gamma$-filtration, i.e.\ 
$A=\projlim_n A/F_\gamma^nA$, then $A\cong\Ch(C)$ by taking projective limits,
and $\ch:\Ch(C)\to\hat C$ is still an isomorphism, i.e.\
{\em augmented $K$-$\lambda$-algebras, complete with respect 
to $\gamma$-filtration, are isomorphic to $\Ch(C)$ or $\hat C$ for
graded $K$-algebras~$C$ with $C_0=K$.}
In this case we still have $C_n=\gr_\gamma^nA\cong A_{(n)}$.

If the $\gamma$-filtration on $A$ is just separated, we can still
embed $A$ into its completion $\hat A_\gamma$ and obtain $\hat A_\gamma
\cong\Ch(C)\cong\hat C$ for $C=\gr_\gamma A_\gamma=\gr_\gamma A$.
Then $C_n$ is still equal to $\hat A_{\gamma,(n)}$, but $A_{(n)}$ might 
be smaller, i.e.\ Soul\'e's construction might give different (smaller)
result.

Finally, if we don't know anything about the $\gamma$-filtration on~$A$,
we still have a $K$-$\lambda$-algebra homomorphism $A\to\hat A_\gamma$,
which induces maps $A_{(n)}\to C_n=\gr_\gamma^nA$, as well as a
homomorphism of graded $K$-algebras $A_{(\cdot)}\to C=\gr_\gamma A$.
In general we cannot expect this homomorphism to be injective or
surjective, so we get two different ``Chow rings'' $A_{(\cdot)}$ and
$\gr_\gamma A$.

\nxsubpoint (Consequences for~$K^0$.)
If $X$ is a ``nice'' generalized ringed space, so that 
$A:=K^0(X)$ is an augmented $K$-$\lambda$-algebra for $K=H^0(X,\bbZ)$
with discrete $\gamma$-filtration, then the Chow ring
$C=\gr_\gamma A_\bbQ=CH(X,\bbQ)$ is concentrated in bounded degrees
($C_n=0$ for $n\gg0$), and $C$ can be recovered from~$A$ via
Soul\'e's construction: $C_n\cong A_{\bbQ,(n)}\subset A_\bbQ$, 
$A_\bbQ=\bigoplus_{n=0}^N A_{\bbQ,(n)}$.

We know that under the above conditions the completed Chern class
$\tilde c:K^0(X)_\bbQ\to\Ch(C)$ is an isomorphism, i.e.\
{\em an element $\xi\in K^0(X)$ is completely determined
(up to $\bbZ$-torsion) by its rank $\epsilon(\xi)$ and its Chern classes
$c_i(\xi)\in C_i=CH^i(X,\bbQ)$, $0<i\leq N$, and all combinations of
rank/Chern classes are possible.}

\nxsubpoint (Equivalent description of augmented $\lambda$-algebras
with discrete $\gamma$-filtration.)
We have just seen that whenever $A$ is an augmented $K$-$\lambda$-algebra
over a binomial ring $K\supset\bbQ$, and the $\gamma$-filtration of~$A$
is discrete (i.e.\ finite), then $A=\bigoplus_{n=0}^NA_{(n)}$ for some
$N>0$, and $C_n=\gr_\gamma^nA\cong A_{(n)}$. In particular, $K=A_{(0)}$
is also determined by the $\lambda$-structure of~$A$.

Conversely, if $A$ is a $\bbQ$-$\lambda$-algebra, such that
$A=\bigoplus_{n=0}^NA_{(n)}$ for some $N>0$, then $A$ is isomorphic
to graded $K$-algebra~$C$, where we put $K:=A_{(0)}$, $C_n:=A_{(n)}$
and define the Adams operations on~$C$ by $\Psi^k|_{C_n}=k^n$ as before.
This implies $A\cong C\cong\hat C\cong\Ch(C)$, i.e.\ $A$ is an augmented
$K$-$\lambda$-algebra with discrete $\gamma$-filtration.

\nxsubpoint (Relation to Chern character.)
Let $A$ be as above, i.e.\ a $\bbQ$-$\lambda$-algebra with weight decomposition
$A=\bigoplus_{n=0}^NA_{(n)}$, or, equivalently, an augmented 
$K$-$\lambda$-algebra with $F^{N+1}_\gamma A=0$ over some $K\supset\bbQ$.
One might ask what the components $x_{(n)}$ 
in $A_{(n)}\cong C_n=\gr^\gamma_n A$ of
an element~$x\in A$ are. The answer is very simple: the isomorphism
between $A$ and $C$ (or $A_{(n)}$ and $C_n$) was given by the Chern character,
hence the $x_{(n)}$ are identified with the components of the Chern character
$\ch(x)$, and $x$ is identified with $\ch(x)$ itself.

\nxsubpoint\label{sp:form.duality.K0} (Duality and Adams operations.)
Let $A$ be an augmented $K$-$\lambda$-algebra. A {\em duality\/} on~$A$
is just an involution $x\mapsto x^*$ on~$A$ compatible with all structure.
If we have such a duality, we can put $\Psi^0:=\epsilon$, $\Psi^{-1}(x):=x^*$,
$\Psi^{-n}(x)=\Psi^n(x^*)$, thus defining ring endomorphisms $\Psi^n:A\to A$
for all $n\in\bbZ$, such that $\Psi^1=\id_A$ and $\Psi^{nm}=\Psi^n\circ\Psi^m$
for all $m$, $n\in\bbZ$. 

For example, if we start from a graded $K$-algebra
$C$, we can define a duality on $\Ch(C)$ by requiring it to be ``continuous'',
functorial, and such that $(1,1+x)^*=(1,1-x)$ for any $x\in C_1$,
$(a,0)^*=(a,0)$ for any $a\in K$. If $K\supset\bbQ$, we can transfer this
duality to $\hat C$ via $\ch$; we obtain $x^*=(-1)^nx$ for $x\in C_n$,
i.e.\ the formula $\Psi^k(x)=k^nx$ for $x\in C_n$ is valid for all $k\in\bbZ$,
not just for $k\geq1$. Then $C=\bigoplus C_n$ is a weight decomposition
with respect to the family of all Adams operations $\{\Psi^k\}_{k\in\bbZ}$.

In our favourite situation $A=\hat K^0(X)$, $K=H^0(X,\bbZ)$ we have a
natural duality, at least if $X$ is additive, induced by duality of vector
bundles: $[\sE]^*:=[\check\sE]$. One checks immediately (e.g.\
reducing to the case of line bundles by Grothendieck projective bundle
argument) that $c_n(x^*)=(-1)^nc_n(x)$, i.e.\ this duality on $A$
is compatible via the completed Chern class $\tilde c$ with the duality 
on~$\Ch(C)$ (and on $\hat C_\bbQ$ via Chern character) just discussed.

If $X$ is non-additive, $\check\sE=\iHom(\sE,\sO)$ 
needn't be a vector bundle. However,
if the $\gamma$-filtration on $A_\bbQ$ (or $A_{\bbQ,\lambda}$, if 
$A$ is not a $\lambda$-ring) is discrete, we get $\ch:A_\bbQ\cong C_\bbQ$,
so we can still transfer the duality just discussed from $C_\bbQ$ to
$A_\bbQ$ (or $A_{\bbQ,\lambda}$) in a completely formal fashion.

\nxpointtoc{Vector bundles over $\CompZ$: further properties}
Vector bundles over $\CompZ$ and related generalized schemes appear to
possess some very interesting number-theoretic properties by themselves.
We cannot discuss them here in much detail, but we would like to mention
at least some of their properties.

\nxsubpoint (Formal duality and dual parametrization.)
For any $A\in GL_n(\bbQ)$ put $A^*:=(A^t)^{-1}$. Then $A\mapsto A^*$ is 
an involution on group $GL_n(\bbQ)$, compatible with multiplication and
preserving subgroups $GL_n(\bbZ)$ and $GL_n(\Zninfty)=\Oct_n$. Therefore,
$A\mapsto A^*$ induces an involution on
$\Oct_n\backslash GL_n(\bbQ)/GL_n(\bbZ)$, i.e.\ the moduli space 
of vector bundles of rank~$n$ over~$\CompZ$. 
This involution can be used in two different ways:
\begin{itemize}
\item We can use it to define a ``formal duality operation'' $\sE\mapsto\sE^*$
on vector bundles: if $\sE$ is given by some matrix~$A$, then $\sE^*$
will be the vector bundle of the same rank given by~$A^*$. In general
$\sE\mapsto\sE^*$ is not a contravariant functor, and $\sE^*\not\cong
\check\sE=\iHom_\sO(\sE,\sO)$ as one would expect in the additive case
(indeed, $\check\sE_\infty$ is the dual of octahedron $\sE_\infty$,
i.e.\ a cube, hence $\check\sE$ is not even a vector bundle if $\rank\sE>2$.) 
However, $\sE\mapsto\sE^*$ is compatible with all isomorphisms of 
vector bundles, and, furthermore, if $\sE'\to\sE\to\sE''$ is a cofibration
sequence, then $A=\smallmatrix{A'}{}{*}{A''}$, whence
$A^*=\smallmatrix{{A'}^*}{*}{}{{A''}^*}$, so we get a dual cofibration
sequence ${\sE''}^*\to\sE^*\to{\sE'}^*$. This means that $[\sE]\mapsto[\sE^*]$
induces an involution on $\hat K^0(\CompZ)=\bbZ\oplus\log\bbQ^*_+$,
easily seen to coincide with the ``formal duality'' $n+\log\lambda\mapsto
n-\log\lambda$ discussed in~\ptref{sp:form.duality.K0}.
\item On the other hand, we can use $\sE\mapsto A^*$ as another parametrization
of vector bundles of rank~$n$ over $\CompZ$ by double coset space
$\Oct_n\backslash GL_n(\bbQ)/GL_n(\bbZ)$. For example, we can
reduce $A^*$ to its canonical form of~\ptref{prop:can.repr.GlQ/GlZ},
consider the set of row g.c.d.s of~$A^*$, construct the ``dual
Harder--Narasimhan filtration'' using the canonical form of~$A^*$,
and define the ``cosemistable'' vector bundles over~$\CompZ$ 
of slope $\log\lambda$
by requiring all row g.c.d.s of~$A^*$ to be equal to~$\lambda$.
This dual parame\-tri\-za\-tion and dual notions seem to be more convenient
in some cases.
\end{itemize}

\nxsubpoint (Application: $\Hom_\sO(\sE,\sO)$ of a Mumford-antiregular~$\sE$.)
Let $\sE$ be a Mumford-antiregular vector bundle of rank~$n$ over~$\CompZ$,
i.e.\ let's suppose all elements of $A^{-1}$ (or its transposed matrix
$A^*$) to lie in~$\bbZ$ (cf.~\ptref{sp:mumford.antireg}). Let's denote by $(e_i)$ and $(f_i)$ the
bases of $E_\infty=\sE_\infty\subset E=\sE_\xi$ and 
$E_\bbZ=\Gamma(\Spec\bbZ,\sE)\subset E$ used to construct $A$ from~$\sE$,
cf.~\ptref{sp:matr.descr.vb.compz}. Let $(e^*_i)$ and $(f^*_i)$ be the
dual bases in~$E^*$; since $e_i=\sum_ja_{ij}f_j$, the dual bases are related
by $e^*_i=\sum_ja^*_{ij}f^*_j$, where $A^*=(a^*_{ij})$. Clearly,
$(f^*_i)$ is a base of dual lattice $E_\bbZ^*\subset E^*$. On the other
hand, since $E_\infty$ was an octahedron with vertices $\pm e_i$,
its dual $E^*_\infty$ consists of all linear forms $u=\sum u_je^*_j$,
such that $|\langle\pm e_i,u\rangle|\leq 1$ for all~$i$. This condition
is obviously equivalent to $|u_j|\leq 1$ for all~$j$, i.e.\ $E^*_\infty$
is the cube with rational vertices $\pm e^*_1\pm e^*_2\pm\cdots\pm e^*_n$.
Since all $a^*_{ij}\in\bbZ$, all $e^*_i$ lie in $E^*_\bbZ$, hence
the same is true for all vertices of cube $E^*_\infty$.

Now let's compute the number $p_\sE(0)=\card\Hom_\sO(\sE,\sO)=|E^*_\bbZ\cap
E^*_\infty|$, discussed in~\ptref{sp:mumford.antireg}. Since all $e^*_i$
lie in lattice~$E^*_\bbZ$, we can replace $E^*_\infty$ by cube
$E^*_\infty+e^*_1+\cdots+e^*_n=[0,2e^*_1]\times\cdots\times[0,2e^*_n]$.
This cube is a fundamental domain for lattice $2\Lambda$ 
generated by $\{2e^*_i\}_i$,
and this lattice contains $E^*_\bbZ$, hence its fundamental domain contains
exactly $(E^*_\bbZ:2\Lambda)=|\det 2A^*|=2^n|\det A^*|$ points of $E^*_\bbZ$.

However, this reasoning is slightly imprecise since we consider the closed
cube $[0,2e^*_1]\times\cdots\times[0,2e^*_n]$, while a fundamental domain
for $\Lambda$ would be given e.g.\ by product of semi-open intervals 
$(0,2e^*_1]\times\cdots\times(0,2e^*_n]$. In order to compensate for this
we introduce for any point $u=\sum_i u_ie^*_i$ its {\em support\/}
$I=\supp u:=\{i\in\stn:u_i\neq0\}$ and write
\begin{multline}
p_\sE(0)=|E^*_\bbZ\cap E^*_\infty|=\sum_{I\subset\stn}
\card\{u\in E^*_\bbZ: 0\leq\langle e_i,u\rangle\leq 2, \supp u=I\}=\\
\sum_{I\subset\stn}2^{|I|}(E^*_I\cap E^*_\bbZ:\Lambda_I)
=\sum_I 2^{|I|}c_I(A)
\end{multline}
Here $E^*_I\subset E^*$ denotes the span of $\{e^*_i\}_{i\in I}$,
and $\Lambda_I=E^*_I\cap\Lambda$ is the $\bbZ$-span of the same set.
We have put
\begin{equation}
c_I(A)=c_I(\sE):=(E^*_I\cap E^*_\bbZ:\Lambda_I)
\end{equation}
This expression makes sense even if $\sE$ is not Mumford-antiregular,
i.e.\ if $E^*_I\cap E^*_\bbZ$ doesn't contain $\Lambda_I$: indeed,
we can always take the absolute value of the determinant of 
any matrix relating any two bases of these two lattices in $E^*_I$.

\nxsubpoint (Computation of $c_I(A)$.)
We want to explain how $c_I(A)=(E^*_I\cap E^*_\bbZ:\Lambda_I)$ can
be computed in terms of matrix~$A$.
Consider the embedding $\kappa_I:E^*_I\to E^*$. Its dual $\kappa_I^*$
is the canonical projection $E\to E_I$, $\sum x_ie_i\mapsto
\sum_{i\in I^c}x_ie_i$, where 
$E_I\subset E$ denotes the $\bbQ$-span of $\{e_i\}_{i\in I}$.
Then $E^*_I\cap E^*_\bbZ=\kappa_I^{-1}(E^*_\bbZ)=\{u\in E^*_I\,|\,
\langle x,\kappa_I(u)\rangle\in\bbZ$ for all $x\in E_\bbZ\}=\{u\in E^*_I\,|\,
\langle \kappa_I^*(x),u\rangle\in\bbZ$ for all $x\in E_\bbZ\}=
(\kappa_I^*(E_\bbZ))^*$. In other words, for any $I\subset\stn$,
$|I|=r$, the positive rational number $c_I(A)$ can be computed as follows:
\begin{itemize}
\item Consider the $n\times r$-submatrix $A^*_I$ of $A^*$, consisting
of rows of~$A^*$ with indices in~$I$. Since $A^*\in GL_n(\bbQ)$,
the $r$ rows of $A^*_I$ are linearly independent, i.e.\
$\rank A^*_I=r$.
\item Consider the $\bbZ$-sublattice $\kappa_I^*(E_\bbZ)$
in $\bbQ^r$ generated by the $n$ columns of $A^*_I$. Since these rows
span $\bbQ^r$ as a $\bbQ$-vector space, this is indeed a sublattice,
so we can find a matrix $B\in GL_r(\bbQ)$, the rows of which
constitute a base of this lattice. In fact, the usual ``integer
Gauss elimination process'' (similar to that used 
in~\ptref{prop:can.repr.GlQ/GlZ}) yields a {\em triangular\/} matrix~$B$ with
this property.
\item Now $c_I(A)=|\det B^*|^{-1}=|\det B|$, and we are done.
\end{itemize}

Notice that the above algorithm yields a positive rational number $c_I(A)$
for any matrix $A\in GL_n(\bbQ)$, and $c_I(\lambda A)=
\lambda^{-|A|}c_I(A)$ for any $\lambda\in\bbQ^*_+$.

\nxsubpoint (Dual Hilbert polynomial of a vector bundle.)
Let $\sE$ be a vector bundle over~$\CompZ$, given by a matrix 
$A\in GL_n(\bbQ)$. Recall that we have defined in~\ptref{sp:mumford.antireg}
the {\em dual Hilbert function} $p_\sE:\log\bbQ^*_+\to\bbN$ 
of~$\sE$ (or~$A$) by
\begin{equation}
p_\sE(\log\lambda)=\card\Hom_\sO(\sE,\sO(\log\lambda))=
\card|E^*_\bbZ\cap\lambda E^*_\infty|
\end{equation}
If $\lambda$ is such that $\sE(-\log\lambda)$ is Mumford-antiregular
(we know that such $\lambda$s are exactly the integer multiples of 
$\lambda_0$, the g.c.d.\ of all $a^*_{ij}$), then
\begin{multline}
p_\sE(\log\lambda)=\card\Hom_\sO(\sE(-\log\lambda),\sO)=\\
=\sum_{I\subset\stn}2^{|I|}c_I(\lambda^{-1}A)=
\sum_{I\subset\stn}2^{|I|}\lambda^{|I|}c_I(A)=\tilde p_\sE(2\lambda)
\end{multline}
where $\tilde p_\sE=\tilde p_A\in\bbQ[T]$ is the polynomial given by
\begin{equation}
\tilde p_\sE(T)=\sum_{I\subset\stn}c_I(A)T^{\card I}
\end{equation}
It is natural to call this $\tilde p_\sE=\tilde p_A$ the {\em
dual Hilbert polynomial\/} of $\sE$ or~$A$.

Clearly, $\tilde p_\sE(T)$ is a polynomial of degree $n=\rank\sE$,
all its coefficients are positive rational numbers, its free term
$\tilde p_\sE(0)=c_{\emptyset}(A)=1$, and its leading coefficient is
$c_\stn(A)=|\det A^*|=|\det A|^{-1}$.

\nxsubpoint (Examples of dual Hilbert polynomials.)
The dual Hilbert polynomial $\tilde p_\sE=\tilde p_A$ can be easily computed
for any matrix $A\in GL_n(\bbQ)$, at least with the aid of a computer. If $A$
is in canonical (or just lower-triangular) form, the computation
is even more simple. For example, $\tilde p_{\sO(\log\lambda)}(T)=
\tilde p_{(\lambda)}(T)=1+\lambda^{-1}T$. Another example:
\begin{equation}
\tilde p_A(T)=1+(a+\gcd(c,d))T+adT^2\quad
\text{for }A^*=\left(\begin{matrix}a\\c&d\end{matrix}\right)
\end{equation}
Consider the vector bundle $\sE$ given by $A=\smallmatrix{1}{1/2}{0}{1}$,
$A^*=\smallmatrix{1}{0}{-1/2}{1}$.
Then $\tilde p_\sE(T)=1+\frac32T+T^2\neq(1+T)^2=\tilde p_\sO(T)^2$, regardless
of the existence of cofibration sequence $\sO\to\sE\to\sO$. Therefore,
{\em $\sE\mapsto\tilde p_\sE$ is not additive.} Actually we could expect this
from our knowledge of $\hat K^0(\CompZ)$: all additive functions of
$\sE$ can be expressed with the aid of the rank $\epsilon(\sE)$ and
its degree $\deg\sE=c_1(\sE)$, but clearly $\tilde p_\sE(T)$ contains 
much more information than that.

Since $p_\sE(\log\lambda)$ coincides with $p_\sE(2\lambda)$ for
infinitely many values of~$\lambda$, we see that {\em $\sE\mapsto p_\sE$
is also non-additive.}

\nxsubpoint (Properties of dual Hilbert functions and polynomials.)
By definition of coproduct $\oplus$ for any two vector bundles 
$\sE$ and $\sE'$ and any $\lambda\in\bbQ^*_+$ we have
\begin{equation}
\Hom_\sO(\sE\oplus\sE',\sO(\lambda))=\Hom_\sO(\sE,\sO(\lambda))\times
\Hom_\sO(\sE,\sO(\lambda))
\end{equation}
Counting elements of these sets we obtain
\begin{equation}
p_{\sE\oplus\sE'}(\log\lambda)=p_\sE(\log\lambda)\cdot p_{\sE'}(\log\lambda)
\end{equation}
Since we can find a rational $\lambda_0>0$, such that 
$p_\sE(\log\lambda)=\tilde p_\sE(2\lambda)$
for all integer multiples $\lambda$ of $\lambda_0$, and similarly
for $\sE'$ and $\sE\oplus\sE'$, we obtain equality of polynomials
\begin{equation}
\tilde p_{\sE\oplus\sE'}=\tilde p_\sE\cdot\tilde p_{\sE'}
\end{equation}
One might also prove this equality directly, by showing
\begin{equation}
c_I(A)=c_{I\cap\stk}(A')\cdot c_{I\cap\stk^c}(A'')
\quad\text{for }A=A'\oplus A''=\smallmatrix{A'}{0}{0}{A''}.
\end{equation}
Here $A'\in GL_k(\bbQ)$, $A''\in GL_{n-k}(\bbQ)$, and $\stk^c:=\stn-\stk$.

When $A=\smallmatrix{A'}{0}{*}{A''}$, essentially the same computation shows
that $c_I(A)$ divides $c_{I\cap\stk}(A')\cdot c_{I\cap\stk^c}(A'')$,
hence {\em we have a coefficientwise inequality of polynomials 
$\tilde p_\sE(T)\leq\tilde p_{\sE'}(T)\tilde p_{\sE''}(T)$ for any
cofibration sequence $\sE'\to\sE\to\sE''$.}

\nxsubpoint ($\tilde p_{\sE\otimes\sE'}$.)
One might ask whether it is possible to compute $\tilde p_{\sE\otimes\sE'}(T)$
knowing only $\tilde p_\sE(T)$ and $\tilde p_{\sE'}(T)$, and similarly
for exterior and symmetric powers. If $\sE'=\sO(\log\mu)$, then
$\sE\otimes_\sO\sE'=\sE(\log\mu)$ is given by matrix $\mu A$, and we 
get $\tilde p_{\sE(\log\mu)}(T)=\tilde p_\sE(T/\mu)$. However,
in general $\tilde p_\sE(T)$ and $\tilde p_{\sE'}(T)$ are insufficient
to compute $\tilde p_{\sE\otimes\sE'}(T)$, as illustrated by the following
example. Let $\sE$ and $\sE'$ be vector bundles of rank~$2$ given by
matrices $A^*=\smallmatrix{2}{0}{2}{6}$ and ${A'}^*=\smallmatrix{1}{0}{3}{12}$.
Then $\tilde p_\sE(T)=1+4T+12T^2=\tilde p_{\sE'}(T)$; however, 
direct computation shows that polynomials
$\tilde p_{\sE\otimes\sE}(T)$ and $\tilde p_{\sE\otimes\sE'}(T)$ differ.

\nxsubpoint (Difference between $p_\sE$ and $\tilde p_\sE$.)
We know that the deviation $\delta_\sE(\lambda):=p_\sE(\log\lambda)-
\tilde p_\sE(2\lambda)$ equals zero for all $\lambda$ divisible by some
$\lambda_0$. However, values of $\delta_\sE(\lambda)$ for $\lambda\not\in
\lambda_0\bbZ$ might have also some interesting number-theoretical properties,
e.g.\ if we study $\delta_\sE(\lambda)$ for $\lambda\in\bbZ$ when 
$\lambda_0=N>1$. It might be even possible to obtain $\delta_\sE(a)=
\legendre ap$ for all $a>0$, $a\in\bbZ$, as illustrated
by one of elementary proofs of the quadratic reciprocity law given by
Gauss, based on a problem of counting points in a triangle:
\begin{align}
\legendre ap&=(-1)^{d(a,p)},\quad\text{where}\\
d(a,p)&=\sum_{0<k<p/2}\left[\frac{2ak}{p}\right]=
\card\bigl\{(x,y)\in\bbN^2\,:\,x\equiv 1\bmod2,\,\frac xp+\frac ya\leq1\bigr\}
\end{align}
In general it might be interesting 
to study the algebra of functions $\bbQ^*_+\to\bbQ$
generated by all dual Hilbert functions~$p_\sE$.

\nxsubpoint (Torsion-free finitely presented $\sO_{\CompZ}$-modules.)
Let $\sF$ be a torsion-free finitely presented $\sO_{\CompZ}$-module.
By~\ptref{sp:finpres.compz.mod}, the category of such~$\sF$'s is equivalent
to the category of triples $(F,F_\bbZ,F_\infty)$, where $F$ is a 
finite-dimensional $\bbQ$-vector space, $F_\bbZ\subset F$ is a lattice in~$F$,
and $F_\infty\subset F$ is a finitely presented torsion-free 
$\Zninfty$-submodule generating~$F$, i.e.\ a symmetric convex polyhedron
inside~$F$ (with rational vertices), not contained in any hyperplane of~$F$.
We see that in this respect $\CompZ$ differs from the projective line $\bbP^1$:
{\em a torsion-free finitely presented $\sO_{\CompZ}$-module $\sF$ 
needn't be a vector bundle.} Notice that the problem of computing
$\card\Gamma(\CompZ,\sF(\log\lambda))$ is nothing else than the classical
problem of counting points of lattice~$F_\bbZ$ inside a convex 
polyhedron~$\lambda F_\infty$.

If we cease to require finite presentation at $\infty$, and use $\Zinfty$
instead of~$\Zninfty$, we obtain the problem of counting lattice points
inside any symmetric convex set $F_\infty\subset F$, e.g.\ a ball.
Even if we don't obtain an equivalence of categories, we still get
a sheaf of $\sO_{\CompZ}$-modules $\sF$ given by
$\sF(\Spec\bbZ[1/N])=F_\bbZ[1/N]$, $\sF(\Spec(\Zninfty\cap\bbZ[1/N]))=
F_\bbZ[1/N]\cap F_\infty$, and $\sF|_{\Spec\bbZ}$ will be a vector bundle.
However, in general $\sF$ won't be even quasi-coherent in the neighbourhood
of $\infty$.

\nxsubpoint (Relation to euclidean metrics.)
In particular, euclidean lattices, i.e.\ lattices $E_\bbZ$ embedded into
a finite-dimensional real space~$E_\bbR$, equipped with a positive definite
quadratic form~$Q$, can be expressed in this form, by using
$E_\infty:=\{x\in E_\bbR\,:\,Q(x)\leq1\}$. This corresponds to the
classical understanding of Arakelov geometry and archimedian structure
as explained in~\ptref{ss:bunoncompz}.

In general, however, there are much more $\Zinfty$-structures 
(i.e.\ norms) on $E_\bbR$ in our sense than quadratic forms, and the
euclidean structures just discussed are not given by finitely presented
$\Zinfty$-modules. We can compensate this as follows. Given
a $\Zinfty$-structure $E_\infty$, i.e.\ a compact convex body
inside a finite-dimensional real space~$E_\bbR$ (cf.~\ptref{def:zinflat}),
we can construct an euclidean structure, i.e.\ a positive definite
quadratic form~$Q$ on~$E_\bbR$ by one of the following methods:
\begin{itemize}
\item If $E_\infty$ is free, i.e.\ if $E_\infty=\Zinfty^{(n)}$ is a 
octahedron centered at the origin, we can choose any its basis $(e_i)$
and declare it an orthonormal basis for~$Q$. Since $\Oct_n\subset O_n(\bbR)$,
$Q$ does not depend on the choice of basis for $E_\infty$.
\item We can define $Q$ by averaging over the dual convex set $E^*_\infty
=\{u\in E^*:u(E_\infty)\subset[-1,1]\}$:
\begin{equation}
Q(x)=\frac{c_n}{\mu(E^*_\infty)}\int_{E^*_\infty}\langle x,u\rangle^2\,du
\end{equation}
for some positive constant $c_n>0$ depending on $n=\dim E_\bbR$.
\item Dually, we can define $Q^*:E^*\to\bbR$ by averaging over~$E_\infty$,
and then consider the dual quadratic form of~$Q^*$.
\end{itemize}
In general these methods give different results. The second and 
the third method have the property to reproduce the original quadratic
form~$Q'$ (up to a constant) if $E_\infty$ was already quadratic 
(i.e.\ if $E_\infty=\{x:Q'(x)\leq1\}$), so the constants may be chosen so as
to have $Q=Q'$ in this case.

On the other hand, if $E_\infty$ is a free $\Zinfty$-module, and 
we fix any basis $(e_i)$ for $E_\infty$, then the resulting quadratic form
will be necessarily invariant under $\Aut_{\Zinfty}(E_\infty)\cong
\Aut_{\Zinfty}(\Zinfty^{(n)})=\Oct_n$, and any quadratic form $Q$ invariant
under $\Oct_n$ is easily seen to be given by a diagonal matrix. Therefore,
we can always choose the constants in the second and the third method
so as to make them coincide with the first method for free $\Zinfty$-modules.
We'll adopt this approach for simplicity.

If $E_{(\infty)}$ is a $\Zninfty$-structure inside a $\bbQ$-vector space~$E$,
we can extend it to a $\Zinfty$-structure $E_\infty$ inside $E_\bbR=
E\otimes_\bbQ\bbR$ 
(e.g.\ by taking the closure of $E_{(\infty)}$ in $E_\bbR$) and apply
any of the above constructions. Notice that if $E_{(\infty)}$ was 
finitely presented, then $E_\infty$ is a symmetric convex polyhedron
with rational vertices, and any of the above methods will produce 
a quadratic form with rational coefficients.

\nxsubpoint (Relation to Shimura varieties.)
The above construction induces a map on moduli spaces
\begin{equation}
\nu:\Oct_n\backslash GL_n(\bbQ)/GL_n(\bbZ)\to
O_n(\bbR)\backslash GL_n(\bbR)/GL_n(\bbZ)
\end{equation}
with dense image, so, for example, a continuous function $f$ on the RHS is
completely determined by the function $g:=f\circ\nu$, continuous and constant
on the fibers of~$\nu$. Conversely, any continuous function $g$ on the LHS
determines a function on the RHS: we first extend it by continuity to
$\Oct_n\backslash GL_n(\bbR)/GL_n(\bbZ)$, and then integrate it over the 
fibers. If we consider vector bundles $\sE$ up to a Serre twist, and
euclidean lattices up to similarity, we obtain another map with dense image
\begin{equation}
\nu':\Oct_n\cdot\bbQ^*_+\backslash GL_n(\bbQ)/GL_n(\bbZ)\to
O_n(\bbR)\cdot\bbR^*_+\backslash GL_n(\bbR)/GL_n(\bbZ)
\end{equation}
The target is the Shimura variety of $GL_n$ without level structure,
so we obtain some relation between e.g.\ automorphic forms with respect 
to~$GL_n$ and sections of certain bundles over the moduli space of
vector bundles of rank~$n$ over~$\CompZ$.

We can introduce a level~$N$ structure on this moduli space, so as to make 
it rigid. This means that we fix a $\bbZ/N\bbZ$-base of $E_\bbZ/NE_\bbZ$,
i.e.\ an isomorphism $E_\bbZ/NE_\bbZ\cong(\bbZ/N\bbZ)^n$. Equivalently,
we consider vector bundles $\sE$ over $\CompZ$ of rank~$n$ together
with a trivialization of $\sE|_{\Spec\bbZ/N\bbZ}$. Then the target of
corresponding map $\nu''$ will be the Shimura variety of $GL_n$ with
respect to the (full) level~$N$ structure. We can consider a 
``level~$\infty$ structure'' on vector bundles $\sE/\CompZ$ if we like:
this corresponds to choosing a $\Zninfty$-base of $E_\infty$. The
corresponding moduli space is $GL_n(\bbQ)/GL_n(\bbZ)$; it is rigid.

\nxsubpoint (Operations with euclidean lattices.)
Choosing a lifting~$\sE$ of an euclidean lattice $\Lambda\subset E$
with respect to
the map $\nu:\Oct_n\backslash GL_n(\bbR)/GL_n(\bbZ)\to O_n(\bbR)\backslash
GL_n(\bbR)/GL_n(\bbZ)$ roughly corresponds to fixing an orthonormal
base in the euclidean space~$E$. Then we might do some operations with~$\sE$
(which is either a vector bundle over~$\CompZ$, or can be approximated
by such), and in some occasions the image of the result under~$\nu$
would not depend on the choice of lifting~$\sE$, thus defining an operation
with euclidean lattices. For example, $\sE\mapsto\sE^*$ corresponds to
the operation of taking dual euclidean lattice. Operations 
$S^k\sE$, $\bigwedge^k\sE$, $\sE\oplus\sE'$ and $\sE\otimes\sE'$ 
have this compatibility property as well,
giving rise to corresponding operations with euclidean lattices.
This probably explains why it was possible to establish a working
Arakelov geometry using vector bundles with suitable hermitian metrics,
and also establishes some relationship between classical Arakelov geometry
and its version discussed in this work.

\nxsubpoint (Automorphisms of vector bundles.)
We have already remarked that the moduli space $\Oct_n\backslash
GL_n(\bbQ)/GL_n(\bbZ)$ is not rigid, i.e.\ a vector bundle $\sE/\CompZ$
may have a non-trivial automorphism group. However, this group is
always finite since $\Aut_\sO(\sE)\subset\Aut_{\Zninfty}(E_\infty)=\Oct_n$.
Fixing a ``level~$\infty$ structure'' makes the moduli space 
$GL_n(\bbQ)/GL_n(\bbZ)$ rigid.

\nxsubpoint\label{sp:comp.OK} (Generalization to other number fields.)
Let $K$ be any number field. Then we might attempt to 
construct $\widehat{\Spec\cO_K}$ in different ways. 
The most canonical of them is to consider the 
``semilocal ring'' $\cO_{K,\infty}$, the intersection of all archimedian
valuation rings $\cO_v$ inside~$K$. Then $\Spec\cO_{K,\infty}$ is
one-dimensional, and its closed points correspond to some families of 
$v|\infty$ (some or all archimedian valuations can be glued together;
this happens already for real quadratic fields $K/\bbQ$; but
if $K$ is a CM-field, this never happens, thus suggesting another 
connection to Shimura varieties).
In any case, we can construct
$\widehat{\Spec\cO_K}\strut^{(N)}$ by choosing any integer $N>1$ and gluing
$\Spec(\cO_K)$ with $\Spec(\cO_{K,\infty}\cap\cO_K[1/N])$ along their 
common open subset $\Spec(\cO_K[1/N])$. The generalized schemes thus 
constructed depend on the choice of $N>1$, but they constitute a
projective system over the set of integers $N>1$ ordered by divisibility,
so we can still define $\widehat{\Spec\cO_K}$ as the projective limit
of this system, either in the category of pro-generalized schemes,
or in the category of generalized ringed spaces. 

This $\widehat{\Spec\cO_K}$ has one generic
point, and its closed points correspond to all valuations of~$K$, but
with some archimedian valuations glued together. This suggests that
the ``true'' $\widehat{\Spec\cO_K}$ should be constructed in a
more sophisticated way. In the fancy language of ``infinite resolution
of singularities'' introduced in~\ptref{sp:smirnov.descr.compz}
one might say that this time the singularity over $\infty$ is even more
complicated than a cusp of infinite order, so we need more complicated
resolution of singularities, which will ``disentagle'' different
points lying over~$\infty$, thus yielding the ``smooth'' model
of~$\widehat{\Spec\cO_K}$.

It would be interesting to classify vector bundles over $\widehat{\Spec\cO_K}$
and to compute~$K^0$. Unfortunately, we cannot even 
apply~\ptref{th:proj.over.archv}, since the condition ``$|x|\in\tilde\bbQ_+$
for all $x\in K$'' is not fulfilled. This problem might become simpler
after the ``correct'' $\widehat{\Spec\cO_K}$ is constructed.

\nxsubpoint (Cyclotomic extensions.)
Another possibility is to study the ``cyclotomic extensions''
$Z_n:=\CompZ\otimes_{\Fone}\bbF_{1^n}$. Such generalized schemes are not
``integral'' (their generic fiber is $\Spec\bbQ[T]/(T^n-1)$), but
their structure (especially over~$\infty$) 
seems to be simpler than that discussed in~\ptref{sp:comp.OK}.
They have the following nice property: a vector bundle of rank~$r$ over~$Z_n$
corresponds to a vector bundle $\sE$ of rank~$rn$ over~$\CompZ$ together
with an action $\sigma:\sE\to\sE$ of the cyclic group 
$C_n=\langle\sigma\rangle$, i.e.\ an element of order~$n$ inside
$\Aut(\sE)\subset\Oct_{nr}$. However, not all such couples $(\sE,\sigma)$
correspond to a vector bundle over~$Z_n$, i.e.\ we have just a fully
faithful functor, not an equivalence of categories. In any case,
the computation of $\hat K^0(Z_n)$ shouldn't be too difficult.


\nxsubpoint (Intersection theory of $\CompZ\strut^{(N)}$.)
Another interesting possibility is to study vector bundles and perfect
cofibrations over $S_N:=\CompZ\strut^{(N)}$ for some $N>1$,
and to compute the Chow ring $CH(S_N)$ afterwards.
If any vector bundle over $\Spec A_N$, where $A_N=\Zninfty\cap\bbZ[N^{-1}]$,
had been trivial, then we would obtain results similar to those obtained
above for $\CompZ$: vector bundles over $S_N$ would be
parametrized by $GL_n(A_N)\backslash GL_n(B_N)/GL_n(\bbZ)$,
where $B_N:=\bbZ[N^{-1}]$, and 
$CH^0(S_N)=\bbZ$, $CH^1(S_N)\cong\Pic(S_N)\cong\log B_{N,+}^\times\cong\bbZ^r$,
where $p_1$, \dots, $p_r$ are distinct prime divisors of~$N$
(cf.~\ptref{sp:pic.compz.N}), and $CH^i(S_N)=0$ for $i\geq2$. However,
even proving $\Pic(A_N)=0$ has been much more complicated than proving
$\Pic(\Zninfty)=0$ (cf.~\ptref{prop:pic.AN.triv}), and we see no reason
for all finitely generated projective modules over $A_N$ to be free,
since $S_N$ is sort of a ``non-smooth version'' of $\CompZ$
(cf.~\ptref{sp:smirnov.descr.compz}). 

\cleardoublepage


\pagestyle{plain}

\bigbreak

\end{document}